\documentclass[11pt]{article}
\usepackage{latexsym,amsmath,amssymb,mathrsfs,graphicx,hyperref,mathpazo,soul,color,amsthm,mathtools,xcolor,multirow}
\hypersetup{pdftex,colorlinks=true,linkcolor=blue,citecolor=blue,filecolor=blue,urlcolor=blue}

\usepackage{url}
\usepackage{newtxtext} 
\usepackage[varvw]{newtxmath} 
\usepackage{caption}
\def\tr{{\raise0pt\hbox{$\scriptscriptstyle\top$}}}

\DeclarePairedDelimiter\floor{\lfloor}{\rfloor}

\usepackage{longtable}

\usepackage{rotating}
\newtheorem{theorem}{Theorem}[section]

\newtheorem{lemma}[theorem]{Lemma}

\newtheorem{proposition}[theorem]{Proposition}

\newtheorem{problem}[theorem]{Problem}
\newtheorem{definition}[theorem]{Definition}

\newtheorem{remark}[theorem]{Remark}
\numberwithin{equation}{section}
\numberwithin{table}{section}
\title{A systematic approach to Diophantine equations: two thousand solved examples}
\author{Ashleigh Ratcliffe}
\begin{document}
\maketitle

\begin{abstract}
	Monograph ``B. Grechuk, Polynomial Diophantine equations. A systematic approach'' suggests solving Diophantine equations systematically in certain order. Many hundreds of the equations are left to the reader. Here, we provide complete solutions to all these equations. The difficulties of solved equations range from elementary to research level. In the last section, we present a summary table of all solved equations ordered by their size, which makes them easy to find. As a result, this document may be used as a database of solved Diophantine equations of small size, to which more difficult equations can be reduced.
\end{abstract}

\tableofcontents

\section*{Introduction}

A polynomial Diophantine equation is an equation of the form 
$$
P(x_1,x_2,...,x_n)=0
$$ 
where $P$ is a polynomial with integer coefficients. 
If the equation (in reduced form) has $k$ monomials of degrees $d_1, \dots, d_k$ with coefficients $a_1, \dots, a_k$, then its size $H$ suggested in \cite{Z2018} is defined as
$$
H(P)= \sum_{j=1}^k |a_j| 2^{d_j}.
$$
As an example, let us look at the equation $3x^2 -y+7xy=1$.
$$
H(3x^2 -y+7xy-1) = 3(2^2)+2+7(2)(2)+1 = 43.
$$

Book \cite{mainbook} suggests a systematic approach to solving Diophantine equations, this means working through equations in order of size.

This document presents solutions to exercises from book \cite{mainbook}. The exercises in Chapter 3 which ask the reader to investigate which equations they can solve are not included in this document and are left open for the reader. In the last section, we present a summary table of all equations, ordered by $H$, which are solved in \cite{mainbook} and in this document.

\section{Chapter 1}
\subsection{Exercise 1.2}

Equations are considered \emph{equivalent} if one equation can be transformed to the other by a permutation of variables, multiplication by a non-zero constant, or by the change of variables $x_i \to -x_i$. 

\textbf{\emph{List all polynomial Diophantine equations of size $H \leq 6$, up to equivalence.}}

See Table \ref{tab:H6all}.

\vspace{10pt}

\begin{minipage}{\textwidth}
\begin{center}
\begin{tabular}{ |c|c|c|c|c|c|c|c|c|c| } 
 \hline
 $H$ & Equation & $H$ & Equation & $H$ & Equation& $H$ & Equation& $H$ & Equation \\ 
 \hline\hline
 $2$ & $x=0$ & $4$ & $xy=0$& $5$ & $x+y+1=0$   & $6$ & $x^2+x=0$& $6$ & $x^2+y=0$\\ 
 \hline
 $3$ & $x+1=0$ & $5$ & $x+3=0$ & $5$ & $xy+1=0$ & $6$ &  $x+y+2=0$& $6$ & $x+y+z=0$  \\ 
 \hline
 $4$ & $x+2=0$ & $5$ & $2x+1=0$ &  $6$ & $x+4=0$ & $6$ &  $2x+y=0$& $6$ & $xy+z=0$  \\ 
 \hline
 $4$ &$x^2=0$ & $5$ & $x^2+1=0$  & $6$& $x^2-2=0$& $6$ &$xy+2=0$  &&\\ 
 \hline
$4$ &   $x+y=0$&$5$ & $x^2-1=0$ & $6$ & $x^2+2=0$ &$6$& $x+xy=0$ &&\\\hline
\end{tabular}
\captionof{table}{All equations of size $H\leq 6$ up to equivalence.\label{tab:H6all}}
\end{center} 
\end{minipage}

\subsection{Exercise 1.3}\label{ex:1var}
\textbf{\emph{For each equation listed in Table \ref{tab:H12onevarlist}, find all its integer solutions. }}

\begin{center}
\begin{tabular}{ |c|c|c|c|c|c| } 
 \hline
 $H$ & Equation & $H$ & Equation & $H$ & Equation \\ 
 \hline\hline
 $8$ & $x^3=0$ & $11$ & $x^3+3=0$ & $12$ & $x^3+x+2=0$ \\ 
 \hline
 $9$ & $x^3+1=0$ & $11$ & $x^3-x+1=0$ & $12$ & $x^3-2x=0$ \\ 
 \hline
 $10$ & $x^3+2=0$ & $11$ & $x^3+x+1=0$ & $12$ & $x^3+2x=0$ \\ 
 \hline
 $10$ & $x^3-x=0$ & $12$ & $x^3+4=0$ & $12$ & $x^3+x^2=0$ \\ 
 \hline
 $10$ & $x^3+x=0$ & $12$ & $x^3-x+2=0$ & & \\ 
 \hline
\end{tabular}
\captionof{table}{\label{tab:H12onevarlist} Equations in one variable of degree $d>2$ and size $H \leq 12$.}
\end{center} 

Any one-variable equation with integer coefficients can be written in the form 
$$
	a_d x^d + \dots + a_{m+1} x^{m+1} + a_m x^m = 0,
$$
for some integers $0\leq m \leq d$ and integer coefficients $a_m, a_{m+1}, \dots, a_d$ such that $a_m \neq 0$ and $a_d\neq 0$. Then any non-zero integer solution must be a divisor of $a_m$. By substituting $x=0$ and the (finitely many) divisors of $a_m$ instead of $x$ into the equation, we can check which of them (if any) are the integer solutions.

For example, for equation $x^3-x+2=0$ we have $a_m=2$, hence we need to check values $x=0$ and $x=-2,-1,1,2$. Direct substitution shows that none of the listed values satisfy the equation, hence it has no integer solutions. 

The other equations can be solved similarly. Table \ref{tab:H12onevarlistsol} lists, for each equation, the value of $a_m$, its divisors to check, and the final answer.

\begin{center}
	\begin{tabular}{ |c|c|c|c|c|c| } 
		\hline
		Equation & $a_m$ & Values to check ($0$ and divisors of $a_m$) & Integer Solutions \\ 
		\hline\hline
		 $x^3=0$ & $1$ & $0,-1,1$ & $x=0$ \\\hline
		$x^3+1=0$ & $1$ & $0,-1,1$  & $x=-1$ \\\hline
		 $x^3+2=0$ & $2$ & $0,-1,1,-2,2$  &  No integer solutions \\\hline 
		 $x^3-x=0$ & $-1$ & $0,-1,1$  &  $x=0$ or $x=\pm 1$ \\\hline
		 $x^3+x=0$ & $1$ & $0,-1,1$  &  $x=0$ \\\hline
		 $x^3+3=0$ & $3$ & $0,-1,1,-2,2,-3,3$ &  No integer solutions \\\hline
		 $x^3-x+1=0$ & $1$ & $0,-1,1$  &  No integer solutions \\\hline
		 $x^3+x+1=0$ & $1$ & $0,-1,1$  &  No integer solutions \\\hline
		 $x^3+4=0$ & $4$  &$0,-1,1,-2,2,-3,3,-4,4$  &  No integer solutions \\\hline
		$x^3-x+2=0$ & $2$ & $0,-1,1,-2,2$  &  No integer solutions \\\hline
		 $x^3+x+2=0$ & $2$ & $0,-1,1,-2,2$  &  $x=-1$ \\  \hline
		 $x^3-2x=0$ &$-2$ & $0,-1,1,-2,2$ &  $x=0$ \\  \hline
		 $x^3+2x=0$ & $2$ & $0,-1,1,-2,2$  &  $x=0$  \\ \hline
		 $x^3+x^2=0$ & $1$ & $0,-1,1$  &  $x=-1$ or $x=0$ \\  \hline
	\end{tabular}
	\captionof{table}{\label{tab:H12onevarlistsol} Integer solutions to the equations listed in Table \ref{tab:H12onevarlist}.}
\end{center} 

\subsection{Exercise 1.5}\label{ex:sepvar}
\textbf{\emph{Table \ref{tab:H7sepvarlist} lists all equations of the form
\begin{equation}\label{eq:sepvar}
x_i + P(x_1,\dots,x_{i-1},x_{i+1},\dots,x_n) = 0,
\end{equation}
where $P$ is a polynomial with integer coefficients, of size $H\leq 7$. For each equation, express a variable as a polynomial of other variables, and then represent the solution set as a polynomial family. }}
\begin{center}
\begin{tabular}{ |c|c|c|c|c|c| } 
 \hline
 $H$ & Equation & H & Equation & H & Equation \\ 
 \hline\hline
  $4$ & $x+y=0$ & $6$ & $x+y+z=0$ & $7$ & $x^2+y-1=0$ \\ 
 \hline
 $5$ & $x+y+1=0$ & $6$ & $xy+z=0$ & $7$ & $x+y+z+1=0$\\ 
 \hline
 $6$ & $x+y+2=0$ & $7$ & $x+y+3=0$ & $7$ & $xy+z+1=0$ \\ 
 \hline
 $6$ & $2 x+y=0$ & $7$ & $2 x+y+1=0$ & & \\ 
 \hline
 $6$ & $x^2+y=0$ & $7$ & $x^2+y+1=0$ & & \\ 
 \hline
\end{tabular}
\captionof{table}{\label{tab:H7sepvarlist} Equations of the form \eqref{eq:sepvar} of size $H\leq 7$.}
\end{center} 

For all equations of the form \eqref{eq:sepvar}, 
we can describe all integer solutions as 
$$
(x_i,x_1,\dots,x_{i-1},x_{i+1},\dots,x_n)=( -P(u_1,\dots,u_{i-1},u_{i+1},\dots,u_n),u_1,\dots,u_{i-1},u_{i+1},\dots,u_n)
$$
for arbitrary integers $u_1,\dots,u_{i-1},u_{i+1},\dots,u_n$.

For example, let us consider equation $xy+z=0$. We can rearrange this equation to $z=-xy$, then taking $x=u$ and $y=v$ where $u$ and $v$ are arbitrary integers, we can describe all integer solutions as $(x,y,z)=(u,v,-uv)$ for integers $u,v$.

\begin{center}
\begin{tabular}{ |c|c||c|c|} 
 \hline
 Equation & $(x,y)$ or $(x,y,z)$ & Equation & $(x,y)$ or $(x,y,z)$ \\ 
 \hline\hline
$x+y=0$&$(u,-u)$ &$x+y+3=0$&$(u,-u-3)$ \\\hline
$x+y+1=0$&$(u,-u-1)$ & $2x+y+1=0$ &$(u,-2u-1)$ \\\hline
$x+y+2=0$&$(u,-u-2)$ & $x^2+y+1=0$&$(u,-u^2-1)$ \\\hline
$2x+y=0$&$(u,-2u)$ &$x^2+y-1=0$&$(u,1-u^2)$ \\\hline
$x^2+y =0$ &$(u,-u^2)$ & $x+y+z+1=0$&$(u,v,-u-v-1)$  \\\hline
$x+y+z=0$ & $(u,v,-u-v)$ & $xy+z+1=0$&$(u,v,-uv-1)$ \\\hline
$xy+z=0$&$(u,v,-uv)$  & &\\\hline

\end{tabular}
\captionof{table}{\label{tab:Ex1.5} Integer solutions to the equations listed in Table \ref{tab:H7sepvarlist}. Assume that $u,v$ are arbitrary integers.}
\end{center}

\subsection{Exercise 1.7}\label{ex:comvar}
\textbf{\emph{Table \ref{tab:H8comvar} lists the equations of size $H\leq 8$ in which all monomials share a variable. For each equation, write it in the form 
\begin{equation}\label{eq:xiR}
P=x_i \cdot R(x_1,\dots, x_n) = 0
\end{equation}
and use this representation to describe all integer solutions. }}
\begin{center}
\begin{tabular}{ |c|c|c|c|c|c| } 
 \hline
 $H$ & Equation &  $H$ & Equation &  $H$ & Equation \\ 
 \hline\hline
 $4$ & $xy=0$ &  $8$ & $x^2+xy=0$ &  $8$ & $xyz=0$ \\ 
 \hline
 $6$ & $xy+x=0$ &  $8$ & $x^2 y=0$ &  &  \\ 
 \hline
 $8$ & $xy+2x=0$ &  $8$ & $xy+xz=0$ &   &  \\ 
 \hline
\end{tabular}
\captionof{table}{\label{tab:H8comvar} Equations of size $H\leq 8$ in which all monomials share a variable.}
\end{center} 

For all equations \eqref{eq:xiR} either $x_i=0$ or $R(x_1,\dots, x_n) = 0$. 

As an example, let us consider the equation $xy+xz=0$. We can rearrange this equation to $x(y+z)=0$. Hence, either $x=0$ (and then $y,z$ are arbitrary), or $y+z=0$. In the latter case, $(y,z)=(v,-v)$ for some integer $v$, and we may take $x$ arbitrary. Finally, all integer solutions to $xy+xz=0$ are $(x,y,z)=(0,u,v),(u,v,-v)$ where $u,v$ are arbitrary integer parameters. 

All other equations in Table \ref{tab:H8comvar} can be solved similarly. Table \ref{tab:H8comvar_sol} provides for each equation, the representation in the form \eqref{eq:xiR}, and the set of integer solutions.

\begin{center}
\begin{tabular}{ |c|c|c|c|c|c| } 
 \hline
 Equation & Form \eqref{eq:xiR} & Solution $(x,y)$ or $(x,y,z)$ \\ 
 \hline\hline
$xy=0$ & $(x)(y)=0$ & $(0,u),(u,0)$  \\\hline
$xy+x=0$ & $x(y+1)=0$ & $(0,u),(u,-1)$ \\\hline
 $xy+2x=0$ & $x(y+2)=0$ & $(0,u),(u,-2)$ \\\hline
 $x^2+xy=0$ & $x(x+y)=0$ & $(0,u),(u,-u)$  \\\hline
 $x^2 y=0$ & $(x^2)(y)=0$ & $(0,u),(u,0)$  \\\hline
 $xy+xz=0$ & $x(y+z)=0$ & $(0,u,v),(u,v,-v)$ \\\hline
$xyz=0$ & $(x)(y)(z)=0$ & $(0,u,v),(u,0,v),(0,u,v)$ \\ \hline
\end{tabular}
\captionof{table}{\label{tab:H8comvar_sol} Integer solutions to the equations listed in Table \ref{tab:H8comvar}. Assume that $u,v$ are arbitrary integers.}
\end{center}

\subsection{Exercise 1.8}\label{ex:reducible}
\textbf{\emph{For each equation in Table \ref{tab:H15redQ}, write it in the form 
\begin{equation}\label{eq:QR}
Q \cdot R = 0,
\end{equation}
where $Q$ and $R$ are non-constant polynomials with integer coefficients, and use this representation to describe all its integer solutions. }}
\begin{center}
\begin{tabular}{ |c|c|c|c|c|c| } 
 \hline
 $H$ & Equation &  $H$ & Equation &  $H$ & Equation \\ 
 \hline\hline
 $8$ & $x^2-y^2=0$ & $13$ & $x^2-y^2+2x+1=0$ &  $15$ & $x^2y+x^2-y-1=0$ \\ 
 \hline
 $9$ & $xy+x+y+1=0$ & $13$ & $x^2y+x-y+1$ &  $15$ & $x^2y+x^2+y+1=0$ \\ 
 \hline
 $11$ & $x^2+xy+y-1=0$ & $14$ & $x^2+xy+x-y-2=0$ & $15$ & $x^3+xy+y+1=0$ \\ 
 \hline
 $12$ & $xy+2x+y+2=0$ & $15$ & $xy+3x+y+3=0$ & $15$ & $xy+xz+x+y+z+1=0$  \\ 
 \hline
  $12$ & $xy+xz+y+z=0$ & $15$ & $2xy+2x+y+1=0$ & $15$ & $xyz+xy+z+1=0$ \\ 
 \hline
  $12$ & $x^2-y^2+x+y=0$ & $15$ & $x^2+xy+2x+y+1=0$ & &   \\ 
 \hline
  $12$ & $x^2+xy+x+y=0$ & $15$ & $x^2y+xy+x+1=0$ & &   \\ 
 \hline
\end{tabular}
\captionof{table}{\label{tab:H15redQ} Equations of size $H\leq 15$ that can be factorised.}
\end{center} 

For equations of the form \eqref{eq:QR} we can solve factors $Q=0$ and $R=0$ separately. When solving $Q=0$ (or $R=0$, respectively), any variables that do not participate in $Q$ (or $R$, respectively) can be set as arbitrary integer parameters. 

For example let us consider the equation 
\begin{equation}\label{eq:x2pxypym1}
x^2+xy+y-1=0,
\end{equation}
which can be rearranged to $(x+1)(y+x-1)=0$. Then either $x+1=0$ (and then $y$ is arbitrary) or $y+x-1=0$. Hence all integer solutions to \eqref{eq:x2pxypym1} are $(x,y)=(-1,u),(u,1-u)$ for arbitrary integer $u$. 

All other equations in Table \ref{tab:H15redQ} can be solved similarly. Table \ref{exercise1.8} provides for each equation, the representation in the form \eqref{eq:QR}, and the set of integer solutions.

\begin{center}
\begin{tabular}{ |c|c|c|c|c|c| } 
 \hline
Equation & Form \eqref{eq:QR} & Integer solutions $(x,y)$ or $(x,y,z)$ \\
\hline \hline
$x^2-y^2=0$& $(x+y)(x-y)=0$&$(u,\pm u)$ \\\hline
$xy+x+y+1=0$&$(x+1)(y+1)=0$  &$(-1,u),(u,-1)$ \\\hline
$x^2+xy+y-1=0$&$(x+1)(y+x-1)=0$&$(-1,u),(u,1-u)$ \\\hline
$xy+2x+y+2=0$ & $(x+1)(y+2)=0$ & $(-1,u),(u,-2)$ \\\hline
$xy+xz+y+z=0$&$(x+1)(y+z)=0$&$(-1,u,v),(u,v,-v)$ \\\hline
$x^2-y^2+x+y=0$&$(-y+x+1)(y+x)=0$&$(u,u+1),(u,-u)$ \\\hline
$x^2+xy+x+y=0$ &$(x+1)(x+y)=0$ & $(-1,u),(u,-u)$ \\\hline
$x^2-y^2+2x+1=0$ & $(1+x-y)(1+x+y)=0$ & $(u,-1-u),(u,1+u)$ \\\hline
$x^2y+x-y+1=0$ & $(1+x)(1-y+xy)=0$&$(-1,u),(0,1),(2,-1)$ \\\hline
$x^2+xy+x-y-2=0$ &$(-1+x)(2+x+y)=0$&$(u,-2-u),(1,u)$\\\hline
$xy+3x+y+3=0$&$(1+x)(3+y)=0$&$(u,-3),(-1,u)$\\\hline
$2xy+2x+y+1=0$&$(2x+1)(1+y)=0$&$(u,-1)$\\\hline
$x^2+xy+2x+y+1=0$&$(1+x)(1+x+y)=0$&$(-1,u),(u,-1-u)$\\\hline
$x^2y+xy+x+1=0$&$(1+x)(1+xy)=0$&$(1,-1),(-1,u)$\\\hline
$x^2y+x^2-y-1=0$&$(x-1)(x+1)(1+y)=0$&$(\pm 1,u),(u,-1)$\\\hline
$x^2y+x^2+y+1=0$&$(1+x^2)(1+y)=0$&$(u,-1)$\\\hline
$x^3+xy+y+1=0$&$(1+x)(1-x+x^2+y)=0$&$(-1,u),(u,-u^2+u-1)$\\\hline
$xy+xz+x+y+z+1=0$&$(1+x)(1+y+z)=0$&$(-1,u,v),(u,v,-1-v)$\\\hline
$xyz+xy+z+1=0$&$(1+xy)(1+z)=0$&$(-1,1,u),(1,-1,u),(u,v,-1)$\\\hline
\end{tabular}
\captionof{table}{Integer solutions to the equations listed in Table \ref{tab:H15redQ}. Assume that $u,v$ are arbitrary integers. \label{exercise1.8}}
\end{center} 

\subsection{Exercise 1.10}\label{ex:H9comconst}
\textbf{\emph{For each equation in Table \ref{tab:H9comconst}, write it in the form 
\begin{equation}\label{eq:xiRf}
x_i \cdot R(x_1,\dots,x_n) = -f, \quad f\neq 0,
\end{equation}
and use this representation to find all its integer solutions. }}

\begin{center}
\begin{tabular}{ |c|c|c|c|c|c| } 
 \hline
 $H$ & Equation &  $H$ & Equation & $H$ & Equation \\ 
 \hline\hline
 $5$ & $xy+1=0$ &  $8$ & $xy+x+2=0$ & $9$ & $x^2+xy-1=0$ \\ 
 \hline
 $6$ & $xy+2=0$ &  $9$ & $xy+5=0$ & $9$ & $x^2+xy+1=0$ \\ 
 \hline
 $7$ & $xy+3=0$ &  $9$ & $xy+x+3=0$ & $9$ & $x^2y+1=0$ \\ 
 \hline
 $7$ & $xy+x+1=0$ &  $9$ & $xy+2x+1=0$ & $9$ & $xy+xz+1=0$ \\ 
 \hline
 $8$ & $xy+4=0$ &  $9$ & $2xy+1=0$ & $9$ & $xyz+1=0$ \\ 
 \hline
\end{tabular}
\captionof{table}{\label{tab:H9comconst} Equations of the form \eqref{eq:xiRf} of size $H\leq 9$.}
\end{center} 

In any integer solution to \eqref{eq:xiRf}, $x_i$ must be a divisor of $f$. We can then substitute all divisors and solve the resulting equation in $n-1$ variables.

For example, let us consider the equation $xy+x+1=0$ which can be rewritten in form $x(y+1)=-1$. 
Then $x$ is a divisor of $-1$, hence either $x=-1$ (and then $y+1=1$) or $x=1$ (and then $y+1=-1$). 
We then obtain that all integer solutions to the original equation are $(x,y)=(-1,0),(1,-2)$.

All other equations in Table \ref{tab:H9comconst} can be solved similarly. Table \ref{exercise1.10} provides for each equation, the representation in the form \eqref{eq:xiRf}, and the set of integer solutions.

\begin{center}
\begin{tabular}[c]{| c | c | c |}
\hline
Equation & Form \eqref{eq:xiRf} & Integer solutions $(x,y)$ or $(x,y,z)$ \\
\hline \hline
$xy+1=0$&$xy=-1$&$\pm(1,-1)$ \\\hline
$xy+2=0$&$xy=-2$&$\pm(1,-2),\pm(2,-1)$ \\\hline
$xy+3=0$ & $xy=-3$&$\pm(1,-3),\pm(3,-1)$\\\hline
$xy+x+1=0$&$x(y+1)=-1$&$(-1,0),(1,-2)$ \\\hline
$xy+4=0$ & $xy=-4$&$\pm(1,-4),\pm(2,-2),\pm(4,-1)$\\\hline
$xy+x+2=0$&$x(1+y)=-2$&$(-2,0),(-1,1),(1,-3),(2,-2)$\\\hline
$xy+5=0$&$xy=-5$&$\pm(1,-5),\pm(5,-1)$\\\hline
$xy+x+3=0$&$x(y+1)=-3$&$(-3,0),(-1,2),(1,-4),(3,-2)$\\\hline
$xy+2x+1=0$&$x(y+2)=-1$&$(-1,-1),(1,-3)$\\\hline
$2xy+1=0$&$x(2y)=-1$&No integer solutions. \\\hline
$x^2+xy-1=0$&$x(x+y)=1$&$(\pm1,0)$\\\hline
$x^2+xy+1=0$&$x(x+y)=-1$&$\pm(1,-2)$\\\hline
$x^2y+1=0$&$x^2y=-1$&$(\pm 1,-1)$\\\hline
$xy+xz+1=0$&$x(y+z)=-1$&$(-1,u,1-u),(1,u,-1-u)$\\\hline
$xyz+1=0$&$xyz=-1$&$(-1,-1,-1),(-1,1,1),(1,-1,1),(1,1,-1)$\\\hline

\end{tabular}
\captionof{table}{Integer solutions to the equations listed in Table \ref{tab:H9comconst}. Assume that $u$ is an arbitrary integer.\label{exercise1.10}}
\end{center}

\subsection{Exercise 1.11}\label{ex:bilin}

\textbf{\emph{Table \ref{tab:H12bilinlist} lists all equations of the form 
\begin{equation}\label{eq:bilinear2var}
axy+bx+cy+d=0,
\end{equation}
where $a,b,c,d$ are integer coefficients with $a\neq 0$, of size $H\leq 12$. For each equation, represent it in the form \begin{equation}\label{eq:bilinearfac}
(ax+c)(ay+b)=bc-ad.
\end{equation}  
and use this representation to find all its integer solutions.}}
\begin{center}
\begin{tabular}{ |c|c|c|c|c|c| } 
 \hline
 $H$ & Equation &  $H$ & Equation & $H$ & Equation \\ 
 \hline\hline
 $8$ & $xy+x+y=0$ &  $11$ & $xy+x+y-3=0$ & $12$ & $xy+x+y+4=0$ \\ 
 \hline
 $9$ & $xy+x+y-1=0$ &  $11$ & $xy+x+y+3=0$ & $12$ & $xy+2x+y-2=0$ \\ 
 \hline
 $10$ & $xy+x+y-2=0$ &  $11$ & $xy+2x+y-1=0$ & $12$ & $xy+2x+2y=0$ \\ 
 \hline
 $10$ & $xy+x+y+2=0$ &  $11$ & $xy+2x+y+1=0$ & $12$ & $xy+3x+y=0$ \\ 
 \hline
 $10$ & $xy+2x+y=0$ &  $12$ & $xy+x+y-4=0$ & $12$ & $2xy+x+y=0$ \\ 
 \hline
\end{tabular}
\captionof{table}{\label{tab:H12bilinlist} Equations of the form \eqref{eq:bilinear2var} of size $H\leq 12$.}
\end{center} 

For each equation of the form \eqref{eq:bilinear2var}, we can rewrite the equation as \eqref{eq:bilinearfac}, and solve it by solving the systems $ax+c=d_i$ and $ay+b=d_j$ where $d_i,d_j$ are integers satisfying $d_i \cdot d_j=bc-ad$. 

For example, let us consider the equation $xy+x+y-1=0$, which is equation \eqref{eq:bilinear2var} with $a=b=c=1$ and $d=-1$. Hence in form \eqref{eq:bilinearfac} we have $(x+1)(y+1)=2$. Hence we need to solve the systems
\begin{itemize}
\item $x+1=-1$ and $y+1=2$,
\item $x+1=1$ and $y+1=-2$,
\item $x+1=-2$ and $y+1=1$,
\item $x+1=2$ and $y+1=-1$.
\end{itemize}
We then obtain that all integer solutions to the original equation are $(x,y)=(-3,-2),(-2,-3),(0,1),(1,0)$.

All other equations in Table \ref{tab:H12bilinlist} can be solved similarly. Table \ref{exercise1.11} provides for each equation, the representation in the form \eqref{eq:bilinearfac} and the set of integer solutions.

\begin{center}
\begin{tabular}{ |c|c|c|c|c|c| } 
 \hline
 Equation & Form \eqref{eq:bilinearfac} & Integer solutions $(x,y)$ \\
\hline \hline
$xy+x+y=0$&$(x+1)(y+1)=1$&$(-2,-2),(0,0)$ \\\hline
$xy+x+y-1=0$&$(x+1)(y+1)=2$&$(-3,-2),(-2,-3),(0,1),(1,0)$\\\hline
$xy+x+y-2=0$&$(x+1)(y+1)=3$&$(-4,-2),(-2,-4),(0,2),(2,0)$\\\hline
$xy+x+y+2=0$&$(x+1)(y+1)=-1$&$(-2,0),(0,-2)$\\\hline
$xy+2x+y=0$&$(x+1)(y+2)=2$&$(-3,-3),(-2,-4),(0,0),(1,-1)$\\\hline
$xy+x+y-3=0$&$(x+1)(y+1)=4$&$(-5,-2),(-3,-3),(-2,-5),(0,3),(1,1),(3,0)$\\\hline
$xy+x+y+3=0$&$(x+1)(y+1)=-2$&$(-3,0),(-2,1),(0,-3),(1,-2)$\\\hline
$xy+2x+y-1=0$&$(x+1)(y+2)=3$&$(-4,-3),(-2,-5),(0,1),(2,-1)$\\\hline
$xy+2x+y+1=0$&$(x+1)(y+2)=1$&$(-2,-3),(0,-1)$\\\hline
$xy+x+y-4=0$&$(x+1)(y+1)=5$&$(-6,-2),(-2,-6),(0,4),(4,0)$\\\hline
$xy+x+y+4=0$&$(x+1)(y+1)=-3$&$(-4,0),(0,-4),\pm (2,-2)$\\\hline
$xy+2x+y-2=0$&$(x+1)(y+2)=4$&$(-5,-3),(-3,-4),(-2,-6),(0,2),(1,0),(3,-1)$\\\hline
$xy+2x+2y=0$&$(y+2)(x+2)=4$&$(-6,-3),(-4,-4),(-3,-6),(-1,2),(0,0),(2,-1)$\\\hline
$xy+3x+y=0$&$(x+1)(y+3)=3$&$(-4,-4),(-2,-6),(0,0),(2,-2)$\\\hline
$2xy+x+y=0$&$(2x+1)(2y+1)=1$&$(-1,-1),(0,0)$\\\hline
\end{tabular}
\captionof{table}{Integer solutions to the equations listed in Table \ref{tab:H12bilinlist}.\label{exercise1.11}}
\end{center} 

\subsection{Exercise 1.12}\label{ex:PQc}
\textbf{\emph{Table \ref{tab:H12PQclist} lists equations representable in the form 
\begin{equation}\label{eq:PQc}
Q \cdot R = c
\end{equation}
of size $H\leq 12$. For each equation, write down its representation \eqref{eq:PQc} and use it to find all its integer solutions.  }}

\begin{center}
\begin{tabular}{ |c|c|c|c|c|c| } 
 \hline
 $H$ & Equation & $H$ & Equation & $H$ & Equation \\ 
 \hline\hline
  $9$ & $x^2-y^2-1=0$ & $11$ & $x^2+x-y^2-1=0$ & $12$ & $x^2+xy+y-2=0$ \\ 
 \hline
 $10$ & $x^2-y^2-2=0$ & $11$ & $x^2+x-y^2+1=0$ & $12$ & $x^2+xy+y+2=0$ \\ 
 \hline
 $10$ & $x^2+x-y^2=0$ & $12$ & $x^2-y^2-4=0$ &  $12$ & $x^2+xy+2y=0$ \\ 
 \hline
 $10$ & $x^2+xy+y=0$ & $12$ & $x^2+x-y^2-2=0$ & $12$ & $x^2+xy+x-y=0$ \\ 
 \hline
 $11$ & $x^2-y^2-3=0$ & $12$ & $x^2+x-y^2+2=0$ & $12$ & $x^2y+x-y=0$ \\ 
 \hline
 $11$ & $x^2+xy+y+1=0$ & $12$ & $x^2+2x-y^2=0$ & & \\ 
 \hline
\end{tabular}
\captionof{table}{\label{tab:H12PQclist} Equations representable in the form \eqref{eq:PQc} of size $H\leq 12$.}
\end{center} 

Equations of the form \eqref{eq:PQc} can be solved using a similar method to that one used in Section \ref{ex:bilin}, by solving the systems of equations $Q=d_i$ and $R=d_j$ where $d_i,d_j$ are integers satisfying $d_i \cdot d_j = c$. 

However, for some equations $P(x,y)=0$, it is not clear how to find a constant $c$ such that $P+c$ can be factorised as $Q\cdot R$. To do this, we can consider the system $\frac{\partial P}{\partial x}=\frac{\partial P}{\partial y}=0$, and if $(x_0,y_0)$ is a solution to this system, then we can take $c=-P(x_0,y_0)$.

For example, let us consider the equation
\begin{equation}\label{eq:x2pxmy2m1}
	x^2+x-y^2-1=0.
\end{equation}
Let $P(x,y)=x^2+x-y^2-1$ be the left-hand side of \eqref{eq:x2pxmy2m1}. Then $\frac{\partial P}{\partial x}=2x+1$, $\frac{\partial P}{\partial y}=-2y$, and $\frac{\partial P}{\partial x}=\frac{\partial P}{\partial y}=0$ for $x=-\frac{1}{2}$ and $y=0$. Let $c=-P\left(-\frac{1}{2},0\right)=\frac{5}{4}$. Then 
$$
P(x,y)+c=x^2+x-y^2+\frac{1}{4}=\left(x+\frac{1}{2}+y\right)\left(x+\frac{1}{2}-y\right).
$$ 
Hence \eqref{eq:x2pxmy2m1} is equivalent to $\left(x+\frac{1}{2}+y\right)\left(x+\frac{1}{2}-y\right)=\frac{5}{4}$, or, $(2x+1+2y)(2x+1-2y)=5$. 
We can then find solutions to this equation by using the factors of 5, giving four systems of equations to solve:
\begin{itemize}
	\item $2x+1+2y=5$ and $2x+1-2y=1$
	\item $2x+1+2y=1$ and $2x+1-2y=5$
	\item $2x+1+2y=-1$ and $2x+1-2y=-5$
	\item $2x+1+2y=-5$ and $2x+1-2y=-1$
\end{itemize}
By solving these systems, we conclude that all integer solutions to equation \eqref{eq:x2pxmy2m1} are $(x,y)=(-2,\pm 1)$, $(1,\pm 1)$.

\vspace{10pt}

Another example we will look at in detail is
\begin{equation}\label{eq:x2ypxmy}
	x^2y+x-y=0.
\end{equation}
Let $P(x,y)=x^2y+x-y$, then we have $\frac{\partial P}{\partial x}=2xy+1$, $\frac{\partial P}{\partial y}=x^2-1$, and $\frac{\partial P}{\partial x}=\frac{\partial P}{\partial y}=0$ for $x=1$ and $y=- \frac{1}{2}$. Let $c=-P\left(1, - \frac{1}{2}\right)=-1$. So, $P(x,y)+c=x^2y+x-y-1=(xy+y+1)(x-1)$. Then \eqref{eq:x2ypxmy} is equivalent to $(xy+y+1)(x-1)=-1$. We can then find solutions to this equation by using the factors of $-1$, giving two systems of equations to solve: 
\begin{itemize}
	\item $x-1=-1$ and $xy+y+1=1$
	\item $x-1=1$ and $xy+y+1=-1$
\end{itemize}
The solution to the first system is $(x,y)=(0,0)$, while the second system only has a non-integer solution $(x,y)=\left(2,- \frac{2}{3}\right)$. Hence, the unique integer solution to equation \eqref{eq:x2ypxmy} is $(x,y)=(0,0)$.

All other equations in Table \ref{tab:H12PQclist} can be solved similarly. Table \ref{exercise1.12} provides for each equation, the representation in the form \eqref{eq:PQc}, and the set of integer solutions.

\begin{center}
\begin{tabular}{ |c|c|c|c|c|c| } 
\hline
Equation & Form \eqref{eq:PQc}  & Integer solutions $(x,y)$ \\
\hline\hline
$x^2-y^2-1=0$&$(x+y)(x-y)=1$&$(\pm 1,0)$\\\hline
$x^2-y^2-2=0$&$(x+y)(x-y)=2$&No integer solutions. \\\hline
$x^2+x-y^2=0$&$(2x+1-2y)(2x+1+2y)=1$&$(-1,0),(0,0)$ \\\hline
$x^2+xy+y=0$&$(x+1)(y+x-1)=-1$&$(-2,4),(0,0)$ \\\hline
$x^2-y^2-3=0$&$(x+y)(x-y)=3$&$(\pm 2,\pm 1)$\\\hline
$x^2+xy+y+1=0$&$(x+1)(x+y-1)=-2$&$(-3,5),(-2,5),(0,-1),(1,-1)$\\\hline
$x^2+x-y^2-1=0$&$(2x+1+2y)(2x+1-2y)=5$&$(-2,\pm 1),(1,\pm 1)$\\\hline
$x^2+x-y^2+1=0$&$(2x+1+2y)(2x+1-2y)=-3$&$(-1,\pm 1),(0,\pm 1)$\\\hline
$x^2-y^2-4=0$&$(x+y)(x-y)=4$&$(\pm2,0)$\\\hline
$x^2+x-y^2-2=0$&$(2x+1+2y)(2x+1-2y)=9$&$(-3,\pm 2),(-2,0),(1,0),(2,\pm 2)$ \\ \hline
$x^2+x-y^2+2=0$&$(2x+1-2y)(2x+1+2y)=-7$&$(-2,\pm 2),(1,\pm 2)$\\\hline
$x^2+2x-y^2=0$&$(x+1-y)(x+1+y)=1$&$(-2,0),(0,0)$\\\hline
$x^2+xy+y-2=0$&$(x+1)(x-1+y)=1$&$(-2,2),(0,2)$\\\hline
$x^2+xy+y+2=0$&$(x+1)(x-1+y)=-3$&$(-4,6),(-2,6),(0,-2),(2,-2)$\\\hline
$x^2+xy+2y=0$&$(x+2)(x-2+y)=-4$&$(-6,9),(-4,8),(-3,9),$ \\ && $(-1,-1),(0,0),(2,-1)$\\\hline
$x^2+xy+x-y=0$&$(x-1)(x+2+y)=-2$&$(-1,0),(0,0),(2,-6),(3,-6)$\\\hline
$x^2y+x-y=0$&$(x+1)(xy-y+1)=1$&$(0,0)$\\\hline
\end{tabular}
\captionof{table}{Integer solutions to the equations listed in Table \ref{tab:H12PQclist}. \label{exercise1.12}}
\end{center} 

\subsection{Exercise 1.13}\label{ex:noreal}
\textbf{\emph{For each equation in Table \ref{tab:H16noreallist}, observe that the polynomial in the left-hand side can be represented as a sum of squares of polynomials. Conclude that the equation has no real solutions and therefore no integer solutions.   }}

\begin{center}
\begin{tabular}{ |c|c|c|c| } 
 \hline
 $H$ & Equation & $H$ & Equation \\ 
 \hline\hline
  $\geq 9$ & $x^2+y^2=c$, $\,\,c<0$ & $\geq 14$ &  $x^2+2x+y^2+1=c$, $\,\,c<0$ \\ 
 \hline
  $\geq 11$ & $x^2+x+y^2=c$, $\,\,c<0$ & $\geq 15$ & $x^2+xy+y^2+x=c$, $\,\,c<0$ \\ 
 \hline
 $\geq 13$ & $2x^2+y^2=c$, $\,\,c<0$ & $\geq 15$ & $x^2+y^2+z^2+x=c$, $\,\,c<0$ \\ 
 \hline
 $\geq 13$ & $x^2+x+y^2+y=c$, $\,\,c<0$ & $\geq 15$ & $2x^2+y^2+y=c$, $\,\,c<0$ \\ 
 \hline
 $\geq 13$ & $x^2+xy+y^2=c$, $\,\,c<0$ & $\geq 15$ & $2x^2+x+y^2=c$, $\,\,c<0$ \\ 
 \hline
 $\geq 13$ & $x^2+y^2+z^2=c$, $\,\,c<0$ & $\geq 16$ & $x^2+2x+y^2+y+1=c$, $\,\,c<0$ \\ 
 \hline
\end{tabular}
\captionof{table}{\label{tab:H16noreallist} Equations with no real solutions of size $H\leq 16$.}
\end{center} 

To solve equations in this section, we will express one side of the equation as a sum of squares, which may require multiplying the equation by a constant and completing the square. As we have $c<0$, any sum of squares equal to $c$ is impossible. 

Let us now consider one equation in detail, for example,
\begin{equation}\label{eq:x2pxypy2pxmc}
x^2+xy+y^2+x=c.
\end{equation}
Multiplying by $4$, we obtain $4x^2+4xy+4y^2+4x=4c$, which can be rearranged to $(x+2)^2+(x+2y)^2+2x^2=4c+4$. Because $c<0$ is an integer, we have $c\leq -1$, hence $4c+4\leq 0$, thus this equality may hold only if $(x+2)^2+(x+2y)^2+2x^2=4c+4=0$. But then we have $(x+2)^2=0$, $(x+2y)^2=0$ and $2x^2=0$. For $(x+2)^2=0$ to be true, we must have $x=-2$, however, for $2x^2=0$, we must have $x=0$. This is a contradiction.  Therefore, equation \eqref{eq:x2pxypy2pxmc} has no integer solutions. 

All other equations in Table \ref{tab:H16noreallist} can be solved similarly. Table \ref{tab:H16noreallistsol} represents these equations in the form $P(x,y)=f(c)$ where $P(x,y) \geq 0$ as a sum of squares, and $f$ is a function such that $f(c)\leq 0$ for $c\leq -1$. The case $P(x,y)=f(c)=0$ is impossible by direct analysis, therefore the original equation has no integer solutions. 

\begin{center}
\begin{tabular}{ |c|c|c|c| } 
\hline
Equation ($c<0$ is an integer parameter) & $P(x,y)=f(c)$ \\
\hline \hline
$x^2+y^2=c$ & $x^2+y^2=c$ \\\hline
$x^2+x+y^2=c$ & $(2x+1)^2 +4y^2=4c+1$ \\\hline
$2x^2+y^2=c$ & $2x^2+y^2=c$ \\\hline
$x^2+x+y^2+y=c$ & $(2x+1)^2+(2y+1)^2=4c+2$ \\\hline
$x^2+xy+y^2=c$ & $(2x+y)^2+3y^2=4c$ \\\hline
$x^2+y^2+z^2=c$ & $x^2+y^2+z^2=c$ \\\hline
$x^2+2x+y^2+1=c$ & $(x+1)^2+y^2=c$ \\\hline
$x^2+xy+y^2+x=c$ & $(x+2)^2+(x+2y)^2+2x^2=4c+4$ \\\hline
$x^2+y^2+z^2+x=c$ & $(2x+1)^2+4y^2+4z^2=4c+1$ \\\hline
$2x^2+y^2+y=c$ & $8x^2+(2y+1)^2=4c+1$ \\\hline
$2x^2+y^2+x=c$ & $(4x+1)^2+8y^2=8c+1$ \\\hline
$x^2+2x+y^2+y+1=c$  & $(2x+2)^2+(2y+1)^2=4c+1$ \\\hline
\end{tabular}
\captionof{table}{Necessary information to prove that the equations listed in Table \ref{tab:H16noreallist} have no real solutions, and therefore no integer solutions. In the second column, $P(x,y) \geq 0$ and $f(c)\leq 0$, with at least one inequality being strict.\label{tab:H16noreallistsol}}
\end{center}

\subsection{Exercise 1.14}\label{ex:boundedreal}
\textbf{\emph{For each equation $P=0$ in Table \ref{tab:H12boundedreal}, prove that the solution set to the inequality $P\leq 0$ is a bounded region with finitely many integer points, and use this information to find all integer solutions to $P=0.$}}

\begin{minipage}{\textwidth}
\begin{center}
\begin{tabular}{ |c|c|c|c|c|c| } 
 \hline
 $H$ & Equation & $H$ & Equation & $H$ & Equation \\ 
 \hline\hline
  $8$ & $x^2+y^2=0$ & $11$ & $x^2+x+y^2-1=0$ & $12$ & $x^2+xy+y^2=0$ \\ 
 \hline
  $9$ & $x^2+y^2-1=0$ & $12$ & $x^2+y^2-4=0$ & $12$ & $2x^2+y^2=0$ \\ 
 \hline
 $10$ & $x^2+y^2-2=0$ & $12$ & $x^2+x+y^2-2=0$ &  $12$ & $x^2+y^2+z^2=0$ \\ 
 \hline
 $10$ & $x^2+x+y^2=0$ & $12$ & $x^2+x+y^2+y=0$ &  &  \\ 
 \hline
 $11$ & $x^2+y^2-3=0$ & $12$ & $x^2+2x+y^2=0$ &  &  \\ 
 \hline
\end{tabular}
\captionof{table}{\label{tab:H12boundedreal} Equations with bounded set of real solutions of size $H\leq 12$.}
\end{center} 
\end{minipage}

\vspace{10pt}

The equations in Table \ref{tab:H12boundedreal} have a bounded set of real solutions. Due to this, we may check all integer values in this region and see which are integer solutions to the original equation. 

For example, let us consider equation
\begin{equation}\label{eq:x2pxpy2py}
x^2+x+y^2+y=0.
\end{equation}
Multiplying this equation by $4$ and adding $2$ to both sides, we obtain $(2x+1)^2+(2y+1)^2=2$. Hence any real solution to this equation must satisfy $(2x+1)^2\leq 2$ and $(2y+1)^2\leq 2$. These inequalities define a bounded region, in which the only integer points are $(x,y)=(-1,-1),(-1,0),(0,-1),(0,0)$. A direct substitution shows that all these points are the integer solutions to equation \eqref{eq:x2pxpy2py}. 

Table \ref{tab:H12boundedrealsol} gives, for each equation in Table \ref{tab:H12boundedreal}, the bounded region containing all its real solutions, and the set of all its integer solutions.

\begin{center}
\begin{tabular}{ |c|c|c|c|c|c| } 
 \hline
Equation & Region implied by equation & Integer solutions $(x,y)$ or $(x,y,z)$ \\ 
 \hline\hline
$x^2+y^2=0$ & $x^2 \leq 0$ and $y^2 \leq 0$ &  $(0,0)$ \\  \hline
$x^2+y^2-1=0$ & $x^2 \leq 1$ and $y^2 \leq 1$ &  $(\pm 1,0),(0,\pm 1)$ \\  \hline
$x^2+y^2-2=0$ & $x^2 \leq 2$ and $y^2 \leq 2$ & $(\pm 1, \pm 1)$ \\  \hline
$x^2+x+y^2=0$ & $(2x+1)^2 \leq 1$ and $4y^2 \leq 1$ & $(-1,0),(0,0)$ \\  \hline
$x^2+y^2-3=0$ & $x^2 \leq 3$ and $y^2 \leq 3$  & No integer solutions  \\  \hline
$x^2+x+y^2-1=0$ & $(2x+1)^2 \leq 5$ and $4y^2 \leq 5$ & $(-1,\pm 1),(0,\pm 1)$  \\  \hline
$x^2+y^2-4=0$ & $x^2\leq 4$ and $y^2 \leq 4$ & $(\pm 2,0),(0,\pm 2)$ \\  \hline
$x^2+x+y^2-2=0$ &   $(2x+1)^2 \leq 9$ and $4y^2 \leq 9$ & $(-2,0),(1,0)$ \\  \hline
$x^2+x+y^2+y=0$ &  $(2x+1)^2 \leq 2$ and $(2y+1)^2 \leq 2$ & $(-1,-1),(-1,0),(0,-1),(0,0)$ \\  \hline
$x^2+2x+y^2=0$ &  $(2x+2)^2 \leq 4$ and $4y^2 \leq 4$ & $(-2,0),(-1,\pm 1),(0,0)$ \\  \hline
$x^2+xy+y^2=0$ & $(2x+y)^2 \leq 0$ and $3y^2 \leq 0$ & $(0,0)$ \\ \hline
$2x^2+y^2=0$ & $2x^2 \leq 0$ and $y^2 \leq 0$ & $(0,0)$ \\  \hline
$x^2+y^2+z^2=0$ & $x^2 \leq 0$ and $y^2 \leq 0$ and $z^2 \leq 0$ & $(0,0,0)$  \\  \hline
\end{tabular}
\captionof{table}{\label{tab:H12boundedrealsol} Integer solutions to the equations listed in Table \ref{tab:H12boundedreal}.}
\end{center}

\subsection{Exercise 1.15}\label{ex:boundedreal2}
\textbf{\emph{Prove that each equation in Table \ref{tab:H14boundedreal2} does not have a real solution with all variables simultaneously large in the absolute value. Then use this information to find all its integer solutions.}}
\begin{center}
\begin{tabular}{ |c|c|c|c|c|c| } 
 \hline
 $H$ & Equation & $H$ & Equation & $H$ & Equation \\ 
 \hline\hline
  $12$ & $x^2y+x+y=0$ & $14$ & $x^2 y+x+y+2=0$ & $14$ & $xyz+x+y+2=0$ \\ 
 \hline
  $12$ & $xyz+x+y=0$ & $14$ & $x^2 y+2 x+y=0$ & $14$ & $xyz+x+y+z=0$ \\ 
 \hline
 $13$ & $x^2y+x+y+1=0$ & $14$ & $x^2 y+x+2y=0$ &  $14$ & $xyz-x-y-z=0$ \\ 
 \hline
 $13$ & $xyz+x+y+1=0$ & $14$ & $x^2 y+x-2y=0$ & $14$ & $x y z+2 x+y=0$ \\ 
 \hline
\end{tabular}
\captionof{table}{\label{tab:H14boundedreal2} Equations of size $H\leq 14$ whose real solutions cannot have all variables large.}
\end{center} 

Equations in Table \ref{tab:H14boundedreal2} cannot have all variables being arbitrarily large in the absolute value. Some equations may have a variable which cannot be large in the absolute value, while other equations may be unbounded for each variable, requiring us to restrict variables simultaneously. Let us demonstrate this with examples.

We will first consider the equation
\begin{equation}\label{eq:x2ypxpy}
x^2y+x+y=0.
\end{equation}
This equation has integer solutions where $y=-\frac{x}{x^2+1}$ is integer. However, $(|x|-1)^2 \geq 0$ implies that $x^2+1 \geq 2|x|$, hence $|y|=\frac{|x|}{x^2+1}\leq \frac{1}{2}$, so the solution is bounded in $y$. The only integer satisfying this inequality is $y=0$. Substituting $y=0$ into the original equation, we obtain that the only integer solution to equation \eqref{eq:x2ypxpy} is $(x,y)=(0,0)$.  

\vspace{10pt}

We will next consider the equation
\begin{equation}\label{eq:xyzpxpy}
xyz+x+y=0.
\end{equation}
The set of real solutions for this equation is unbounded for each variable, however, we can restrict variables simultaneously. For example, if $|x| \geq 1$ and $|y| \geq 1$, then 
$$
|z|=\frac{|x+y|}{|xy|} \leq \frac{1}{|x|}+\frac{1}{|y|} \leq 2.
$$
Hence, for any real solution, we must have either $|x| <1$ or $|y|<1$, or $|z| \leq 2$. After checking these values, we can then conclude that all integer solutions to equation \eqref{eq:xyzpxpy} are $(x,y,z)=(u,-u,0)$, $(0,0,u)$, $\pm(1,1,-2)$, $\pm(2,2,-1)$, where $u$ is an arbitrary integer.

\vspace{10pt}

We will next consider the equation
\begin{equation}\label{eq:x2ypxpyp1}
x^2y+x+y+1=0.
\end{equation}
This equation has integer solutions where $y=-\frac{x+1}{x^2+1}$ is integer. However, if $|x|>1$, then 
$$
x^2+1>|x|+1\geq |x+1|>0,
$$
but then $0<|y|<1$, a contradiction. Hence, $|x|\leq 1$, and we may then conclude that all integer solutions to \eqref{eq:x2ypxpyp1} are $(x,y)=(-1,0),(0,-1),(1,-1)$.

\vspace{10pt}

We will next consider the equation
\begin{equation}\label{xyzpxpyp1}
xyz+x+y+1=0.
\end{equation}
Assume that $(x,y,z)$ is an integer solution with $xy\neq 0$. Then $z=\frac{-(x+y+1)}{xy}$ is an integer. So, 
$$
|z|=\frac{|x+y+1|}{|xy|} \leq \frac{1}{|y|} + \frac{1}{|x|}+\frac{1}{|xy|} \leq 3.
$$
Hence, either $|z|\leq 3$ or $x=0$ or $y=0$. By checking these cases, we find that all integer solutions to equation \eqref{xyzpxpyp1} are 
\begin{equation}\label{xyzpxpyp1sol}
	\begin{aligned}
	(x,y,z)=(u,-1,1),(u, -1 - u, 0),(-1, 0, u),(-1, u, 1),(0, -1, u),(1, 1, -3), (1, 2,-2), \\ (2, 1, -2),(2, 3, -1),(3, 2, -1), \quad \text{where $u$ is an arbitrary integer.}
	\end{aligned}
	\end{equation}

\vspace{10pt}

We will next consider the equation
\begin{equation}\label{eq:x2ypxpyp2}
x^2y+x+y+2=0.
\end{equation}
Let us first assume that $|x|>2$ and $|y|>2$. Then 
$$
|x^2y|>2|xy|=|xy|+|xy|>2|x|+2|y|=|x|+|x|+|y|+|y|>|x|+|y|+4>|x+y+2|,
$$
which is a contradiction. Hence, we must have either $|x| \leq 2$ or $|y| \leq 2$. By substituting these values into the original equation, we may conclude that all integer solutions to equation \eqref{eq:x2ypxpyp2} are $(x,y)=(-2,0),(0,-2)$.

\vspace{10pt}

We will next consider the equation
\begin{equation}\label{eq:x2yp2xpy}
	x^2y+2x+y=0.
\end{equation}
Let us first assume that $|x|>2$ and $|y|>2$. Then 
$$
|x^2y|>2|xy|=|xy|+|xy|>|2x|+2|y| > |2x|+|y| \geq |2x+y|,
$$
which is a contradiction. Hence, we must have either $|x| \leq 2$ or $|y| \leq 2$. By substituting these values into the original equation, we may conclude that the integer solutions to equation \eqref{eq:x2yp2xpy} are $(x,y)=(0,0),\pm(1,-1)$.

\vspace{10pt}

We will next consider the equation
\begin{equation}\label{eq:x2ypxp2y}
	x^2y+x+2y=0.
\end{equation}
Let us first assume that $|x|>2$ and $|y|>2$. Then 
$$
|x^2y|>2|xy|=|xy|+|xy|>2|x|+|2y| > |2y|+|x| \geq |2y+x|,
$$
which is a contradiction. Hence, we must have either $|x| \leq 2$ or $|y| \leq 2$. By substituting these values into the original equation, we may conclude that the unique integer solution to equation \eqref{eq:x2ypxp2y} is $(x,y)=(0,0)$.

\vspace{10pt}

We will next consider the equation
\begin{equation}\label{eq:x2ypxm2y}
x^2y+x-2y=0.
\end{equation}
Let us first assume that $|x|>2$ and $|y|>2$. Then 
$$
|x^2y|>2|xy|=|xy|+|xy|>|2y|+2|x| > |2y|+|x| \geq |x-2y|,
$$
which is a contradiction. Hence, we must have either $|x| \leq 2$ or $|y| \leq 2$. By substituting these values into the original equation, we may conclude that the integer solutions to equation \eqref{eq:x2ypxm2y} are $(x,y)=$, $(0,0)$, $\pm(1,1)$, $\pm(2,-1)$.

\vspace{10pt}

We will next consider the equation
\begin{equation}\label{eq:xyzpxpyp2}
xyz+x+y+2=0.
\end{equation}
Let us first assume that $|x|>2$, $|y|>2$ and $|z|>2$. Then 
$$
|xyz| > |2xy| = |xy|+|xy| > 2|x|+2|y| =|x|+|x|+|y|+|y|>|x|+|y|+2 \geq |x+y+2|,
$$ 
which is a contradiction. Hence, we must have either $|x| \leq 2$ or $|y| \leq 2$ or $|z| \leq 2$. By substituting these values into the original equation, we may conclude that the integer solutions to equation \eqref{eq:xyzpxpyp2} are $(x,y,z)=$ ${(u,-2-u,0)}$, $(-2,0,u)$, $(0,-2,u)$,  $(-1,1,2)$, $(1,-1,2)$, $(1,1,-4)$, $(1,3,-2)$, $(2,4,-1)$, $(3,1,-2)$, $(4,2,-1)$.

\vspace{10pt}

We will next consider the equation
\begin{equation}\label{xyzpxpypz}
xyz+x+y+z=0.
\end{equation}
This equation is equivalent to $2xyz+2x+2y+2z=0$. Let us first assume that $|x|>2$, $|y|>2$ and $|z|>2$. Then, 
$$
\begin{aligned}
|2xyz| =&  |xyz|+|xyz| > 2|xy| + 4|z| \\ =& |xy|+|xy|+4|z|>  2|x|+2|y|+2|z| \geq |2x+2y+2z|,
\end{aligned}
$$
which is a contradiction. By checking values of $|x|\leq 2$, $|y|\leq 2$ and $|z|\leq 2$, we obtain that integer solutions to equation \eqref{xyzpxpypz} are 
\begin{equation}\label{eq:xyzmxmymzsol}
	(x,y,z)=(-1,1,u),(0,u,-u) \quad \text{and their permutations, where $u$ is an arbitrary integer.}
\end{equation}

\vspace{10pt}

We will next consider the equation
\begin{equation}\label{eq:xyzmxmymz}
xyz-x-y-z=0.
\end{equation}
This equation is equivalent to $2xyz-2x-2y-2z=0$. Let us first assume that $|x|>2$, $|y|>2$ and $|z|>2$. Then, 
$$
\begin{aligned}
	|2xyz| = & |xyz|+|xyz| > 2|xy| + 4|z|  =|xy|+|xy|+4|z| \\ & > 2|x|+2|y|+2|z| \geq |2x+2y+2z|,
\end{aligned}
$$ 
which is a contradiction. By checking values of $|x|\leq 2$, $|y|\leq 2$ and $|z|\leq 2$, we obtain that all integer solutions to equation \eqref{eq:xyzmxmymz} are $(x,y,z)=(0,u,-u),(-1,-2,-3),(1,2,3)$ and their permutations, where $u$ is an arbitrary integer.

\vspace{10pt}

The final equation we will consider is
\begin{equation}\label{eq:xyzp2xpy}
xyz+2x+y=0.
\end{equation}
Let us first assume that $|x|>2$, $|y|>2$ and $|z|>2$. Then, 
$$
|xyz| > 2|xy| = |xy|+|xy| > 2|x|+2|y| =2|x|+|y|+|y|>2|x|+|y| \geq |2x+y|,
$$
which is a contradiction. By checking values of $|x|\leq 2$, $|y|\leq 2$ and $|z|\leq 2$, we obtain that all integer solutions to equation \eqref{eq:xyzp2xpy} are $(x,y,z)= (u,-2u,0)$, $(0,0,u)$, $\pm(1,-1,1)$, $\pm(1,1,-3)$, $\pm(1,2,-2)$, $\pm(2,4,-1)$, $\pm(3,3,-1)$, where $u$ is an arbitrary integer.

Table \ref{tab:H14boundedreal2sol} provides, for each equation in Table \ref{tab:H14boundedreal2}, the restriction on the variables derived in the proofs above, and the final set of all integer solutions.

\begin{center}
\begin{tabular}{ |c|c|c|c|c|c| } 
 \hline
Equation & Restriction & Integer solutions $(x,y)$ or $(x,y,z)$ \\ 
 \hline\hline
 $x^2y+x+y=0$ & $y=0$ & $(0,0)$  \\\hline
  $xyz+x+y=0$ & $|z| \leq 2$ or $x=0$ or $y=0$ & $(u,-u,0),(0,0,u),\pm(1,1,-2),\pm(2,2,-1)$ \\\hline
$x^2y+x+y+1=0$ & $|x| \leq 1$ & $(-1,0),(0,-1),(1,-1)$ \\\hline
  $xyz+x+y+1=0$ & $|z| \leq 3$ or $x=0$ or $y=0$ & $(u,-1,1),(u,-1-u,0),(-1,0,u),$ \\ && $(-1,u,1),(0,-1,u),(1,1,-3),(1,2,-2),$ \\ && $(2,1,-2),(2,3,-1),(3,2,-1)$ \\\hline
 $x^2 y+x+y+2=0$ & $|x| \leq 2$ or $|y| \leq 2$ & $(-2,0),(0,-2)$ \\\hline
$x^2 y+2 x+y=0$ & $|x| \leq 2$ or $|y| \leq 2$ & $(0,0),\pm(1,-1)$ \\\hline
 $x^2 y+x+2y=0$ & $|x| \leq 2$ or $|y| \leq 2$ &  $(0,0)$ \\\hline
 $x^2 y+x-2y=0$ & $|x| \leq 2$ or $|y| \leq 2$ & $(0,0),\pm(1,1),\pm(2,-1)$ \\\hline
 $xyz+x+y+2=0$ & $|x| \leq 2$ or $|y| \leq 2$ or $|z| \leq 2$ & $(u,-2-u,0),(-2,0,u), (0,-2,u),$ \\ && $(-1,1,2), (1,-1,2), (1,1,-4),(1,3,-2),$ \\ && $ (2,4,-1),(3,1,-2),(4,2,-1)$ \\  \hline
 $xyz+x+y+z=0$ & $|x| \leq 2$ or $|y| \leq 2$ or $|z| \leq 2$ & $(-1,1,u),(0,u,-u)$ \\ && and their permutations \\ \hline
  $xyz-x-y-z=0$ & $|x| \leq 2$ or $|y| \leq 2$ or $|z| \leq 2$ & $(0,u,-u),(-1,-2,-3),(1,2,3)$ \\ && and their permutations \\  \hline
  $x y z+2 x+y=0$ & $|x| \leq 2$ or $|y| \leq 2$ or $|z| \leq 2$ & $(u,-2u,0), (0,0,u),\pm(1,-1,1), \pm(1,1,-3), $ \\ && $\pm(1,2,-2),\pm(2,4,-1), \pm(3,3,-1)$ \\ 
 \hline
\end{tabular}
\captionof{table}{\label{tab:H14boundedreal2sol} Integer solutions to the equations listed in Table \ref{tab:H14boundedreal2}. Assume that $u$ is an arbitrary integer.}
\end{center} 

\subsection{Exercise 1.19}\label{ex:sepvara}
\textbf{\emph{Table \ref{tab:H11axPlist} lists non-linear equations of the form 
\begin{equation}\label{eq:sepvara}
a x_i + P(x_1,\dots,x_{i-1},x_{i+1},\dots,x_n) = 0,
\end{equation}
where $a\neq 0$ is an integer and $P$ is a polynomial with integer coefficients, of size $H\leq 10$. For each equation, find all its integer solutions.  }}
\begin{center}
\begin{tabular}{ |c|c|c|c|c|c| } 
 \hline
 $H$ & Equation &  $H$ & Equation & $H$ & Equation \\ 
 \hline\hline
 $8$ & $x^2+2y=0$ & $9$ & $xy+2z+1=0$ &  $10$ & $x^2+3 y = 0$  \\ \hline
  $8$ & $xy+2z=0$ & $10$ & $x^2+2y-2=0$ & $10$ & $x y+2 z+2 = 0$ \\ \hline
 $9$ & $x^2+2y-1=0$ &  $10$ & $x^2+2y+2=0$ & $10$ & $x y+x+2 z = 0$ \\ \hline
  $9$ & $x^2+2y+1=0$ & $10$ & $x^2+x+2 y=0$ & $10$ & $x y+3 z = 0$ \\ \hline
\end{tabular}
\captionof{table}{\label{tab:H11axPlist} Non-linear equations of the form \eqref{eq:sepvara} of size $H\leq 10$.}
\end{center} 

It is clear that equation \eqref{eq:sepvara} may have integer solutions only if $P(x_1,\dots,x_{i-1},x_{i+1},\dots,x_n)$ is divisible by $a$. We must then consider the cases for which values of the variables $x_1,\dots,x_{i-1},x_{i+1},\dots,x_n$ this is true for. 

As an example, let us consider the equation
$$
x^2+2y-1=0.
$$
We can see from this equation that $x^2$ is odd, so $x$ must be odd. Let $x=2u+1$, then integer solutions to this equation are $(x,y,z)=(2u+1,-2u^2-2u)$ where $u$ is an arbitrary integer.

\vspace{10pt}

As another example, let us consider the equation
\begin{equation}\label{eq:xypxp2z}
xy+x+2z=0.
\end{equation}
We can see from this equation that we must have $xy+x$ even, which implies the following cases, (i) $x$ is even, or (ii) $y$ is odd. 
\begin{itemize}
\item[(i)] Let $x=2u$ and $y$ arbitrary, then solutions to this equation are $(x,y,z)=(2u,v,-uv-u)$, where $u,v$ are arbitrary integers. 
\item[(ii)] Let $y=2v-1$ and $x$ arbitrary, then solutions to this equation are $(x,y,z)=(u,2v-1,-uv)$, where $u,v$ are arbitrary integers. 
\end{itemize}
Therefore, the integer solutions to this equation \eqref{eq:xypxp2z} can be written as $(x,y,z)=$ $(2u,v,-uv-u)$, $(u,2v-1,-uv)$, where $u,v$ are arbitrary integers.

The other equations in Table \ref{tab:H11axPlist} can be solved similarly, and their solutions are presented in Table \ref{tab:H11axPlistsol}.

\begin{center}
\begin{tabular}{ |c|c|c|c|c|c| } 
 \hline
Equation & Cases & Solution $(x,y,z)$ \\ 
 \hline\hline
$x^2+2y=0$ & $x$ is even & $(2u,-2u^2)$  \\ \hline
 $xy+2z=0$ & $x$ or $y$ is even & $(2u,v,-uv),(u,2v,-uv)$ \\ \hline
  $x^2+2y-1=0$ & $x$ is odd & $(2u+1,-2u^2-2u)$  \\ \hline
 $x^2+2y+1=0$ & $x$ is odd&  $(2u+1,-2u^2-2u-1)$ \\ \hline
$xy+2z+1=0$ & $x$ and $y$ are odd & $(2u+1,2v+1,-2uv-v-u-1)$ \\ \hline
$x^2+2y-2=0$ & $x$ is even & $(2u,1-2u^2)$ \\ \hline
  $x^2+2y+2=0$ & $x$ is even & $(2u,-1-2u^2)$ \\ \hline
  $x^2+x+2 y=0$ & $x$ is either even or odd & $(2u,-2u^2-u),(2u+1,-2u^2-3u-1)$ \\ \hline
   $x^2+3 y = 0$ & $x$ is divisible by $3$ & $(3u,-3u^2)$ \\ \hline
$x y+2 z+2 = 0$ & $x$ or $y$ is even & $(2u,v,-1-uv),(u,2v,-1-uv)$ \\ \hline
   $x y+x+2 z = 0$ & $x$ is even, or $y$ is odd & $(2u,v,-uv-u),(u,2v-1,-uv)$ \\ \hline
$x y+3 z = 0$ & $x$ or $y$ is divisible by $3$ & $(u,3v,-uv),(3u,v,-uv)$ \\ \hline
\end{tabular}
\captionof{table}{\label{tab:H11axPlistsol} Integer solutions to the equations listed in Table \ref{tab:H11axPlist}. Assume that $u,v$ are arbitrary integers.}
\end{center} 

\subsection{Exercise 1.23}\label{ex:lin2var}
\textbf{\emph{Table \ref{tab:H14lin2varlist} lists 2-variable linear equations 
\begin{equation}\label{eq:lin2var}
ax+by+c=0,
\end{equation}
where $a,b$ are non-zero integers, of size $H\leq 14$. For each equation, use Theorem \ref{th:lin2var} below 
to describe all its integer solutions.  }}

\begin{theorem}\label{th:lin2var}[Theorem 1.22 in the book]
	Equation \eqref{eq:lin2var} has an integer solution if and only if $d=\gcd(a,b)$ is a divisor of $c$. In this case, an integer solution $(x_0,y_0)$ to \eqref{eq:lin2var} can be efficiently computed by the extended Euclidean algorithm. Then all the solutions to \eqref{eq:lin2var} are given by
	$$
		x = x_0 + \frac{b}{d}u, \quad y = y_0 - \frac{a}{d} u, \quad\quad u \in {\mathbb Z}.
	$$
\end{theorem} 

\begin{center}
\begin{tabular}{ |c|c|c|c|c|c| } 
 \hline
 $H$ & Equation &  $H$ & Equation & $H$ & Equation \\ 
 \hline\hline
 $9$ & $2x+2y+1=0$ & 13 & $2x+2y+5=0$ & 14 & $3x+3y+2 = 0$ \\ 
 \hline
 $10$ & $3x+2y=0$ & 13 & $3x+2y+3=0$ & 14 & $4x+3y = 0$ \\ 
 \hline
 $11$ & $2x+2y+3=0$ & 13 & $3x+3y+1=0$ & 14 & $5x+2y = 0$ \\ 
 \hline
 $11$ & $3x+2y+1=0$ & 13 & $4x+2y+1=0$ & &  \\ 
 \hline
 $12$ & $3x+2y+2=0$ & 14 & $3x+2y+4=0$ & &  \\ 
 \hline
\end{tabular}
\captionof{table}{\label{tab:H14lin2varlist} Linear equations in $2$ variables of size $H\leq 14$.}
\end{center} 

Theorem \ref{th:lin2var} is easy to apply. For each equation in Table \ref{tab:H14lin2varlist}, Table \ref{exercise1.22} provides the necessary steps to apply the Theorem and summarises the equation's integer solutions.

\begin{center}
\begin{tabular}{ |c|c|c|c|c|c| } 
\hline
Equation & $d=\gcd (a,b)$ & Does $d$ divide $c$? & $(x_0,y_0)$ & Integer solutions $(x,y)$ \\
\hline \hline
$2x+2y+1=0$ &2&No&&No integer solutions \\\hline
$3x+2y=0$&1&Yes&$(0,0)$&$(2u,-3u)$ \\\hline
$2x+2y+3=0$ &2&No&&No integer solutions \\\hline
$3x+2y+1=0$&1&Yes&$(1,-2)$&$(1+2u,-2-3u)$ \\\hline
$3x+2y+2=0$ & 1 & Yes & $(0,-1)$ & $(2u,-1-3u)$ \\\hline
$2x+2y+5=0$ & 2 & No &&No integer solutions \\\hline
$3x+2y+3=0$&1&Yes&$(-1,0)$&$(-2u-1,3u)$\\\hline
$3x+3y+1=0$&3&No&&No integer solutions\\\hline
$4x+2y+1=0$&2&No&&No integer solutions \\\hline
$3x+2y+4=0$&1&Yes&$(0,-2)$&$(2u,-2-3u)$\\\hline
$3x+3y+2=0$&3&No&&No integer solutions \\\hline
$4x+3y=0$&1&Yes&$(0,0)$&$(3u,-4u)$\\\hline
$5x+2y=0$&1&Yes&$(0,0)$&$(2u,-5u)$\\\hline
\end{tabular}
\captionof{table}{Integer solutions to the equations listed in Table \ref{tab:H14lin2varlist}. Assume that $u$ is an arbitrary integer.\label{exercise1.22}}
\end{center} 

\subsection{Exercise 1.26}\label{ex:lingen}
\textbf{\emph{Table \ref{tab:H17lingen} lists linear equations in $n>2$ variables of size $H\leq 17$. For each equation, use Theorem \ref{th:linnvar} below to describe all its integer solutions.   }}

\begin{theorem}\label{th:linnvar}[Theorem 1.25 in the book] 
	Equation 
	\begin{equation}\label{eq:linear}
		a_1 x_1 + \dots + a_n x_n + b = 0,
	\end{equation}
	where $a_1,a_2\dots, a_n$ are non-zero integers, and $b$ is an arbitrary integer,
	 has an integer solution if and only if $d=\text{gcd}(a_1,\dots,a_n)$ is a divisor of $b$. Further, there is a polynomial-time algorithm which, given $a_1,\dots, a_n$ and $b$ with $\text{gcd}(a_1,\dots,a_n)$ dividing $b$, returns an integer solution $x^0=(x_1^0, \dots, x_n^0)$ to \eqref{eq:linear}. Then each of the (infinitely many) solutions to \eqref{eq:linear} can be written as 
	$$
	x = x^0 + x_h
	$$ 
	where $x_h$ is a solution to the homogeneous equation \eqref{eq:linhom} described in Proposition \ref{prop:linhom} below.
\end{theorem}

\begin{proposition}\label{prop:linhom}[Proposition 1.24 in the book]  
	Equation 
	\begin{equation}\label{eq:linhom}
		a_1 x_1 + \dots + a_n x_n = 0
	\end{equation}
 has infinitely many integer solutions. Further, one can choose $n-1$ solutions\footnote{Here, $j$ is just an upper index, not a power.} 
	$$
	x^j=(x_1^j, \dots, x_n^j), \quad j=1,\dots, n-1
	$$
	such that each solution $x=(x_1, \dots, x_n)$ to \eqref{eq:linhom} is an integer linear combination of those $n-1$ solutions:
	$$
	x = \sum_{j=1}^{n-1} u_j x^j \quad\quad u_j \in {\mathbb Z}, \,\, j=1,\dots, n-1.
	$$
	Moreover, there is a polynomial-time algorithm which, given coefficients $a_1,\dots, a_n$ of \eqref{eq:linhom}, returns $x^1, \dots, x^{n-1}$.
\end{proposition}

\vspace{10pt}

\begin{minipage}{\textwidth}
\begin{center}
\begin{tabular}{ |c|c|c|c|c|c| } 
 \hline
 $H$ & Equation &  $H$ & Equation & $H$ & Equation \\ 
 \hline\hline
 $13$ & $2x+2y+2z+1=0$ & 16 & $3x+2y+2z+2=0$ & 17 & $3x+3y+2z+1 = 0$ \\ 
 \hline
 $14$ & $3x+2y+2z=0$ & 16 & $3x+3y+2z=0$ & 17 & $4x+2y+2z+1 = 0$ \\ 
 \hline
 $15$ & $2x+2y+2z+3=0$ & 17 & $2x+2y+2z+5=0$ & 17 & $2x+2y+2z+2t+1 = 0$ \\ 
 \hline
 $15$ & $3x+2y+2z+1=0$ & 17 & $3x+2y+2z+3=0$ & &  \\ 
 \hline
\end{tabular}
\captionof{table}{\label{tab:H17lingen} Linear equations in $n>2$ variables of size $H\leq 17$.}
\end{center} 
\end{minipage}

\vspace{10pt}

All equations in this exercise are of the form \eqref{eq:linear}. These equations have integer solutions if and only if $d=\text{gcd}(x_1,\dots,x_n)$ divides $b$. 
If the solution set is non-empty, it can be written as
\begin{equation}\label{eq:linearsol}
x = x_0 + x_h,
\end{equation}
where $x=(x_1,\dots, x_n)$, $x_0$ is any solution to \eqref{eq:linear}, and $x_h$ is a solution to the homogeneous equation 
\eqref{eq:linhom}.
In turn, the solution set to \eqref{eq:linhom} can be written as a linear combination with integer coefficients of $n-1$ solutions to \eqref{eq:linhom}.

As a first simple example, consider the equation
$$
2x+2y+2z+3=0.
$$
Because the left-hand side is always odd, this equation has no integer solutions. 

\vspace{10pt}

As another example, we will consider the equation
\begin{equation}\label{eq:3xp3yp2z}
3x+3y+2z=0.
\end{equation}
From this equation it is easy to see that $2z$ must be divisible by $3$, so let $z=3z'$, then we obtain $x+y+2z'=0$ which is an equation of the form \eqref{eq:sepvar} whose integer solution is $(x,y,z')=(-u-2v,u,v)$ with $u,v$ arbitrary integers. Therefore, we obtain that integer solutions to \eqref{eq:3xp3yp2z} are
$$
(x,y,z)=(-u-2v,u,3v) \quad \text{ with } u,v \quad \text{ arbitrary integers.}
$$

\vspace{10pt}

As another example, we will consider the equation
\begin{equation}\label{eq:3xp3yp2zp1}
3x+3y+2z+1=0.
\end{equation}
By trial and error, it is easy to find a solution $(x^0,y^0,z^0)=(-1,0,1)$. The integer solutions to the homogenous equation \eqref{eq:3xp3yp2z} are $(x,y,z)=(-u-2v,u,3v)$ with $u,v$ arbitrary integers. Therefore, by \eqref{eq:linearsol} we have that the integer solutions to equation \eqref{eq:3xp3yp2zp1} are 
$$
(x,y,z)=(-u-2v-1,u,3v+1) \quad \text{ with } u,v \quad \text{ arbitrary integers.}
$$

All other equations in Table \ref{tab:H17lingen} can be solved similarly. Table \ref{exercise1.25} provides the necessary details to use formula \eqref{eq:linearsol} and summarises all integer solutions to the equations.

\begin{center}
\begin{tabular}{ |c|c|c|c|c|c| } 
\hline
Equation & $d$ & Does $d$ divide $b$? & $x^0$ & Integer solutions $(x,y,z)$ or $(x,y,z,t)$ \\ \hline
\hline
$2x+2y+2z+1=0$&2&No&&No integer solutions. \\\hline
$3x+2y+2z=0$&1&Yes&$(0,0,0)$&$(2u,v,-3u-v)$ \\\hline
$2x+2y+2z+3=0$&2&No& &No integer solutions. \\\hline
$3x+2y+2z+1=0$&1&Yes&$(1,-2,0)$&$(2u+1,v-2,-3u-v)$ \\\hline
$3x+2y+2z+2=0$&1&Yes&$(0,0,-1)$&$(2u,v,-1-v-3u)$\\\hline
$3x+3y+2z=0$&1&Yes&$(0,0,0)$&$(v,-2u-v,3u)$\\\hline
$2x+2y+2z+5=0$&2&No&&No integer solutions.\\\hline
$3x+2y+2z+3=0$&1&Yes&$(-1,0,0)$&$(2u-1,v,-v-3u)$\\\hline
$3x+3y+2z+1=0$&1&Yes&$(-1,0,1)$&$(-1+v,-2u-v,1+3u) $ \\\hline
$4x+2y+2z+1=0$&2&No&&No integer solutions. \\\hline
$2x+2y+2z+2t+1=0$&2&No&&No integer solutions. \\\hline
\end{tabular}
\captionof{table}{Integer solutions to the equations listed in Table \ref{tab:H17lingen}. Assume that $u,v$ are arbitrary integers.\label{exercise1.25}}
\end{center}

\subsection{Exercise 1.29}\label{ex:nomod}
\textbf{\emph{For each of the equations listed in Table \ref{tab:H16nomod}, find an integer $m\geq 2$ such that the equation has no solutions modulo $m$. Conclude that the equation has no integer solutions.  }}
\begin{center}
\begin{tabular}{ |c|c|c|c|c|c| } 
 \hline
 $H$ & Equation &  $H$ & Equation & $H$ & Equation \\ 
 \hline\hline
 $14$ & $x^2+xy-y^2+2=0$ & 15 & $x^2+x+1+2 y z=0$ & 15 & $2 x^2-y^2+y+1 = 0$ \\ 
 \hline
 $15$ & $x^2+xy-y^2+3=0$ & 15 & $2 x^2-y^2-3=0$ & 15 & $x^3-x+y^2+1 = 0$ \\ 
 \hline
 $15$ & $x^2+x y-y^2+x-1=0$ & 15 & $2 x^2-y^2+3=0$ & 16 & $x^3-x+y^2-2 = 0$ \\ 
 \hline
 $15$ & $x^2+x-1+2 y z=0$ & 15 & $2 x^2-y^2+y-1 = 0$ & 16 & $2 x^2-y^2 + x + 2 = 0$ \\ 
 \hline
\end{tabular}
\captionof{table}{\label{tab:H16nomod} Equations of size $H \leq 16$ with no solutions modulo some integer.}
\end{center} 

If an equation has no integer solutions modulo some integer $m$, then the equation has no integer solutions. 

Table \ref{exercise1.26} proves that all equations in Table \ref{tab:H16nomod} have no integer solutions, as they have no integer solutions modulo some integer $m$. 

\begin{longtable}[c]{| c | c | p{6.5 cm} |}
 \hline
Equation & $m$ & Proof \\
 \hline
 \endfirsthead
 \multicolumn{3}{|c|}{$(H=14)$}\\ \hline
$2+x^2+xy-y^2=0$ & 4 & \begin{tabular}{c| c c c c}
$x \, \backslash \, y$ & 0 & 1&2&3 \tabularnewline \hline
0 & 2&1&2&1 \tabularnewline
1&3&3&1&1 \tabularnewline
2&2&3&2&3 \tabularnewline
3&3&1&1&3 \tabularnewline
\end{tabular} \\\hline
\multicolumn{3}{|c|}{$(H=15)$}\\ \hline
$3+x^2+xy-y^2=0$ & 5 & \begin{tabular}{c| c c c c c}
$x \, \backslash \, y$ & 0 & 1&2&3 &4\tabularnewline \hline
0 &3&2&4&4&2 \tabularnewline
1&4&4&2&3&2 \tabularnewline
2&2&3&2&4&4 \tabularnewline
3&2&4&4&2&3 \tabularnewline
4&4&2&3&2&4 \tabularnewline
\end{tabular}\\\hline
$-1+x+x^2+xy-y^2=0$ & 4 &\begin{tabular}{c| c c c c}
$x \, \backslash \, y$ & 0 & 1&2&3 \tabularnewline \hline
0 &3&2&3&2 \tabularnewline
1&1&1&3&3 \tabularnewline
2&1&2&1&2 \tabularnewline
3&3&1&1&3 \tabularnewline
\end{tabular} \\\hline 

$3+2x^2-y^2=0$ & 8 &\begin{tabular}{c| c c c c c c c c}
$x \, \backslash \, y$ & 0 & 1&2&3 &4 &5&6&7\tabularnewline \hline
0&3&2&7&2&3&2&7&2 \tabularnewline
1&5&4&1&4&5&4&1&4 \tabularnewline
2&3&2&7&2&3&2&7&2 \tabularnewline
3&5&4&1&4&5&4&1&4 \tabularnewline
4&3&2&7&2&3&2&7&2 \tabularnewline
5&5&4&1&4&5&4&1&4 \tabularnewline
6&3&2&7&2&3&2&7&2 \tabularnewline
7&5&4&1&4&5&4&1&4 \tabularnewline
\end{tabular} \\\hline

$-3+2x^2-y^2=0$ & 8 &\begin{tabular}{c| c c c c c c c c}
$x \, \backslash \, y$ & 0 & 1&2&3 &4 &5&6&7\tabularnewline \hline
0& 5 & 4&1&4&5&4&1&4 \tabularnewline
1& 7&6&3&6&7&6&3&6\tabularnewline
2&5&4&1&4&5&4&1&4 \tabularnewline
3&7&6&3&6&7&6&3&6 \tabularnewline
4&5&4&1&4&5&4&1&4 \tabularnewline
5&7&6&3&6&7&6&3&6 \tabularnewline
6&5&4&1&4&5&4&1&4 \tabularnewline
7&7&6&3&6&7&6&3&6 \tabularnewline
\end{tabular} \\\hline

$1+2x^2+y-y^2=0$ & 2 &\begin{tabular}{c| c c}
$x \, \backslash \, y$ & 0 & 1\tabularnewline \hline
0&1&1 \tabularnewline
1&1&1\tabularnewline
\end{tabular} \\
\end{longtable}

In Table \ref{tab:H16nomod}, the non-existence of integer solution modulo $m$ has been proved by considering $m^2$ cases. In many examples, the number of cases can be greatly reduced by smarter analysis. For example, for equations 
$$
\pm 3+2x^2-y^2=0
$$
it is clear that $y$ is odd, and therefore $y^2=1$ modulo $8$. But then $2x^2=1\pm 3$ modulo $8$, hence $x^2$ is either $(1-3)/2=-1$ or $(1+3)/2=2$ modulo $4$, which is impossible.

\subsection{Exercise 1.40}\label{ex:2mona}
\textbf{\emph{Table \ref{tab:H16mon2a} lists two-monomial equations of the form 
		\begin{equation}\label{eq:2mon1signed}
			x_1^{r_1} x_2^{r_2} \dots x_m^{r_m} = \pm x_{m+1}^{r_{m+1}} \dots x_n^{r_n}.
		\end{equation}
	 of size $H\leq 16$. For each equation, describe all its integer solutions. }}

	\begin{center}
	\begin{tabular}{ |c|c|c|c|c|c| } 
		\hline
		$H$ & Equation &  $H$ & Equation & $H$ & Equation \\ 
		\hline\hline
		$8$ & $x^2-yz=0$ & $12$ & $xyz-t^2=0$  & $16$ & $x^2y-z^2t=0$ \\ 
		\hline
		$8$ & $xy-zt=0$ & $12$ & $x^2y-zt=0$  & $16$ & $x_1x_2x_3-x_4^2x_5=0$ \\
		\hline 
		$12$ & $x^3-y^2=0$ & $12$ & $x_1x_2x_3-x_4x_5=0$   & $16$ & $x_1x_2x_3-x_4x_5x_6=0$  \\ 
		\hline
		$12$ & $x^2y-z^2=0$ & $16$ & $x^3-yzt=0$ &  & \\
		\hline
		$12$ & $x^3-yz=0$ & $16$ & $x^3-y^2z=0$ &  &  \\
		\hline
	\end{tabular}
	\captionof{table}{\label{tab:H16mon2a} Equations of the form \eqref{eq:2mon1signed} of size $H\leq 16$.}
\end{center} 

Two-monomial equations are of the form 
\begin{equation}\label{eq:2mon}
a x_1^{r_1} x_2^{r_2} \dots x_m^{r_m} = b x_{m+1}^{r_{m+1}} \dots x_n^{r_n},
\end{equation}
where $r_1,\dots,r_n$ are positive integers, and $a,b$ are non-zero coefficients.

A trivial integer solution is one with $x_1 x_2 \cdots x_n=0$ and these solutions are easy to find.
Theorem \ref{th:2mon1}, Proposition \ref{prop:linmin} and Remark \ref{remark:dxkmyk} below reduce the task of describing all non-trivial integer solutions to \eqref{eq:2mon} to the task of finding all minimal solutions to equation
\begin{equation}\label{eq:linnonneghom}
	\sum_{i=1}^m r_i z_i - \sum_{i=m+1}^n r_i z_i = 0, \quad \quad z_i \geq 0, \,\, i=1,\dots, n.
\end{equation}
A solution $z=(z_1,\dots,z_n)\neq (0,\dots,0)$ to this equation is called \emph{minimal} if there is no other solution $z'=(z'_1,\dots,z'_n)\neq (0,\dots,0)$ to \eqref{eq:linnonneghom} such that $z' < z$ (which means that $z'_i\leq z_i$ for $i=1,\dots,n$ with at least one inequality being strict).

\begin{theorem}\label{th:2mon1}[Theorem 1.39 in the book] 
Let $0<m<n$ be integers. Let $r_1,\dots,r_n$ be positive integers at least one of which is odd. Let 
$$
c^k = (c_{k1}, \dots, c_{kn}), \quad k=1,2,\dots, N,
$$
be the complete list of minimal solutions to the equation \eqref{eq:linnonneghom}.
Then for any integers $u_1, \dots, u_N$
\begin{equation}\label{eq:2mon1sol}
x_i = \prod_{k=1}^N (u_k)^{c_{ki}}, \quad i=1,2,\dots,n
\end{equation}
is a solution to \begin{equation}\label{eq:2mon1}
x_1^{r_1} x_2^{r_2} \dots x_m^{r_m} = x_{m+1}^{r_{m+1}} \dots x_n^{r_n}.
\end{equation}
Conversely, any non-trivial integer solution to \eqref{eq:2mon1} must be of the form \eqref{eq:2mon1sol}.
\end{theorem}

\begin{proposition}\label{prop:linmin}[Proposition 1.34 in the book] 
For any minimal solution $(z_1,\dots,z_n)$ to 
$$
\sum_{i=1}^m r_i z_i - \sum_{i=m+1}^n r_i z_i = b, \quad \quad z_i \geq 0, \,\, i=1,\dots, n,
$$
where $r_1, \dots, r_n,b$ are positive integers, we have 
$$
\max\{z_1,\dots,z_m\} \leq \max\{r_{m+1},\dots,r_n,b\} \quad \text{and} \quad \max\{z_{m+1},\dots,z_n\} \leq \max\{r_1,\dots,r_m\}.
$$
\end{proposition} 

\begin{remark}\label{remark:dxkmyk}[Remark follows Proposition 1.31 in the book]
For any integers $d$ and $k\geq 2$, if $\sqrt[k]{d}$ is not an integer then it cannot be rational. This is equivalent to the statement that $(x,y)=(0,0)$ is the only solution to the equation
\begin{equation}\label{eq:dxkmyk}
dx^k-y^k=0, \quad k\geq 2 \quad \text{and} \quad d \,\, \text{is not a perfect k-th power}.
\end{equation}
\end{remark}

Now let us apply this theory to solving equations. Equation
\begin{equation}\label{eq:x2myz}
x^2-yz=0
\end{equation}
is solved in Section 1.6.1 in the book, and its integer solutions are
\begin{equation}\label{eq:x2myzsol}
(x,y,z)=(uvw,uv^2,uw^2) \quad \text{for some integers} \quad u,v,w.
\end{equation}

\vspace{10pt}

Equation
\begin{equation}\label{eq:xymzt}
xy-zt=0
\end{equation}
is solved in Section 1.6.1 in the book, and its integer solutions are
\begin{equation}\label{eq:xymztsol}
(x,y,z,t)=(uv,wr,uw,vr) \quad \text{for some integers} \quad u,v,w,r.
\end{equation}

\vspace{10pt}

Equation
\begin{equation}\label{eq:x3my2}
x^3-y^2=0
\end{equation}
is solved in Section 1.6.1 in the book, and its integer solutions are
\begin{equation}\label{eq:x3my2sol}
(x,y)=(u^2,u^3) \quad \text{for some integer} \quad u.
\end{equation}

\vspace{10pt}

Equation
\begin{equation}\label{eq:x1x2x3mx4x5}
x_1x_2x_3-x_4x_5=0
\end{equation}
is solved in Section 1.6.3 in the book, and its integer solutions are
\begin{equation}\label{eq:x1x2x3mx4x5sol}
(x,y)=(u_1u_2,u_3u_4,u_5u_6,u_1u_3u_5,u_2u_4u_6) \quad \text{for some integers} \quad u_1,u_2,u_3,u_4,u_5,u_6.
\end{equation}

\vspace{10pt}

Let us now consider the equation
\begin{equation}\label{eq:x2ymz2}
	x^2y-z^2=0.
\end{equation}
If $x=0$ then $z=0$, and $y$ can be arbitrary, so that $(x,y,z)=(0,u,0)$ for $u \in {\mathbb Z}$. If $x\neq 0$, then $\sqrt{y}=|z/x|$ is rational, hence it is integer by Remark \ref{remark:dxkmyk}.
Let $z/x=v$, then $y=v^2$. Let us take $x=u$ arbitrary, and then $z=xv=uv$. In conclusion, all integer solutions to \eqref{eq:x2ymz2} are
\begin{equation}\label{eq:x2ymz2sol}
	(x,y,z)=(0,u,0) \quad \text{or} \quad (u,v^2,uv) \quad \text{for some integers} \quad u,v. 
\end{equation}

\vspace{10pt}

The next equation to consider is
\begin{equation} \label{eq:x3myz}
x^3-yz=0.
\end{equation}
This is an equation of the form \eqref{eq:2mon1}, and we will apply Theorem \ref{th:2mon1}. For this equation, the corresponding equation \eqref{eq:linnonneghom} is
$$
3x_1-y_1-z_1=0.
$$
By Proposition \ref{prop:linmin}, the minimal solutions to this equation satisfy $x_1\leq 1$ and $\max\{y_1,z_1\}\leq 3$, and easy trial and error returns a list of four minimal solutions $(x_1,y_1,z_1)=(1,0,3), (1,3,0), (1,2,1), $ and $(1,1,2)$. Let us organize these minimal solutions as the rows of the following matrix:
$$
c = 
\begin{bmatrix}
	c_{11}&c_{12}&c_{13}\\
	c_{21}&c_{22}&c_{23}\\
	c_{31}&c_{32}&c_{33}\\
	c_{41}&c_{42}&c_{43}\\
\end{bmatrix}
=
\begin{bmatrix}
	1&0&3\\
	1&3&0\\
	1&2&1\\
	1&1&2\\
\end{bmatrix}.
$$
Using these minimal solutions and Theorem \ref{th:2mon1}, we can find the non-trivial integer solutions to the initial equation by \eqref{eq:2mon1sol}:
 \[ x=   \prod_{k=1}^{N} (u_k)^{c_{k1}} = u_1^{c_{11}}u_2^{c_{21}}u_3^{c_{31}}u_4^{c_{41}} = u_1u_2u_3u_4,\] 
\[ y=   \prod_{k=1}^{N} (u_k)^{c_{k2}} = u_1^{c_{12}}u_2^{c_{22}}u_3^{c_{32}}u_4^{c_{42}}  = u_2^3\,u_3^2\,u_4,\] 
\[ z=   \prod_{k=1}^{N} (u_k)^{c_{k3}} = u_1^{c_{13}}u_2^{c_{23}}u_3^{c_{33}}u_4^{c_{43}} = u_1^3\,u_3 u_4^2.\]

Therefore, non-trivial integer solutions to equation \eqref{eq:x3myz} are of the form 
\begin{equation} \label{eq:x3myzsol}
(x,y,z)=\left(u_1u_2u_3u_4,u_2^3 \, u_3^2\, u_4,u_1^3\,u_3 u_4 ^2\right) \quad \text{for some integers} \quad u_1,u_2,u_3,u_4.
\end{equation}

As Theorem \ref{th:2mon1} only finds non-trivial solutions, we must now check that the trivial solutions are included. The trivial solutions of this equation are $(x,y,z)=(0,u,0)$ and $(0,0,u)$, where $u$ is an arbitrary integer. These solutions are of the form \eqref{eq:x3myzsol} with $(u_1,u_2,u_3,u_4)=(0,1,1,u)$ and $(1,0,u,1)$, respectively.

Hence, all integer solutions to equation \eqref{eq:x3myz} are of the form \eqref{eq:x3myzsol}.
 
 \vspace{10pt}

The next equation to consider is
\begin{equation} \label{eq:xyzmt2}
xyz-t^2=0.
\end{equation} 
This is an equation of the form \eqref{eq:2mon1}, and we will apply Theorem \ref{th:2mon1}. For this equation, the corresponding equation \eqref{eq:linnonneghom} is
$$
2t_1-x_1- y_1-z_1=0.
$$
By Proposition \ref{prop:linmin}, the minimal solutions to this equation satisfy $t_1\leq 1$ and $\max\{x_1,y_1,z_1\}\leq 2$, and easy trial and error returns a list of six minimal solutions $(x_1,y_1,z_1,t_1)$, which are the rows of the following matrix
$$c = \begin{bmatrix}
	1&0&1&1\\
	1&1&0&1\\
	0&1&1&1\\
	2&0&0&1\\
	0&2&0&1\\
	0&0&2&1\\
\end{bmatrix}.$$

Using these minimal solutions and Theorem \ref{th:2mon1}, we can find non-trivial integer solutions to the initial equation by \eqref{eq:2mon1sol}:
 \[ x=   \prod_{k=1}^{N} (u_k)^{c_{k1}} = u_1^{c_{11}}u_2^{c_{21}}u_3^{c_{31}}u_4^{c_{41}}u_5^{c_{51}}u_6^{c_{61}} = u_1 u_2 u_4 ^2,\] 
\[ y=   \prod_{k=1}^{N} (u_k)^{c_{k2}} = u_1^{c_{12}}u_2^{c_{22}}u_3^{c_{32}}u_4^{c_{42}} u_5^{c_{52}}u_6^{c_{62}} = u_2 u_3  u_5^2,\] 
\[ z=   \prod_{k=1}^{N} (u_k)^{c_{k3}} = u_1^{c_{13}}u_2^{c_{23}}u_3^{c_{33}}u_4^{c_{43}}u_5^{c_{53}}u_6^{c_{63}} = u_1 u_3 u_6^2,\]
\[ t=   \prod_{k=1}^{N} (u_k)^{c_{k4}} = u_1^{c_{14}}u_2^{c_{24}}u_3^{c_{34}}u_4^{c_{44}}u_5^{c_{54}}u_6^{c_{64}} = u_1 u_2 u_3 u_4 u_5 u_6.\]

Therefore any non-trivial integer solutions to equation \eqref{eq:xyzmt2} must be of the form 
\begin{equation} \label{eq:xyzmt2sol}
(x,y,z,t)=(u_1 u_2 u_4 ^2 \, ,u_2 u_3  u_5^2 \, ,u_1 u_3 u_6^2 \, ,u_1 u_2 u_3 u_4 u_5 u_6) \quad \text{for some integers} \quad u_1, \dots, u_6.
\end{equation} 
As Theorem \ref{th:2mon1} only finds non-trivial solutions, we must now check that the trivial solutions are included. The trivial solutions of \eqref{eq:xyzmt2} are $(x,y,z,t)=$ $(0,u,v,0)$, $(u,0,v,0)$ and $(u,v,0,0)$, where $u,v$ are arbitrary integers. These solutions are of the form \eqref{eq:xyzmt2sol} with $(u_1,u_2,u_3,u_4,u_5,u_6)=$ $(v,u,1,0,1,1)$, $(1,u,v,1,0,1)$ and $(u,1,v,1,1,0)$, respectively.

Hence all integer solutions to equation \eqref{eq:xyzmt2} are of the form \eqref{eq:xyzmt2sol}.

\vspace{10pt}

The next equation to consider is 
\begin{equation}\label{eq:x2ymzt}
	x^2y-zt=0.
\end{equation}
This is an equation of the form \eqref{eq:2mon1}, and we will apply Theorem \ref{th:2mon1}. For this equation, the corresponding equation \eqref{eq:linnonneghom} is
$$
2x_1+y_1=z_1+t_1.
$$
By Proposition \ref{prop:linmin}, the minimal solutions to this equation satisfy $\max\{x_1,y_1\} \leq 1$ and $\max\{z_1,t_1\}\leq 2$, and easy trial and error returns a list of five minimal solutions $(x_1,y_1,z_1,t_1)$, which are the rows of the following matrix
$$
c = \begin{bmatrix}
	1&0&1&1\\
	0&1&1&0\\
	0&1&0&1\\
	1&0&2&0\\
	1&0&0&2\\
\end{bmatrix}.
$$

Using these minimal solutions and Theorem \ref{th:2mon1}, we can find non-trivial integer solutions to the initial equation by \eqref{eq:2mon1sol}:
 \[ x=   \prod_{k=1}^{N} (u_k)^{c_{k1}} = u_1^{c_{11}}u_2^{c_{21}}u_3^{c_{31}}u_4^{c_{41}}u_5^{c_{51}} = u_1u_4u_5, \] 
\[ y=   \prod_{k=1}^{N} (u_k)^{c_{k2}} = u_1^{c_{12}}u_2^{c_{22}}u_3^{c_{32}}u_4^{c_{42}} u_5^{c_{52}} = u_2u_3, \] 
\[ z=   \prod_{k=1}^{N} (u_k)^{c_{k3}} = u_1^{c_{13}}u_2^{c_{23}}u_3^{c_{33}}u_4^{c_{43}}u_5^{c_{53}} = u_1u_2u_4^2, \]
\[ t=   \prod_{k=1}^{N} (u_k)^{c_{k4}} = u_1^{c_{14}}u_2^{c_{24}}u_3^{c_{34}}u_4^{c_{44}}u_5^{c_{54}} = u_1u_3u_5^2. \]

Therefore any non-trivial integer solutions to equation \eqref{eq:xyzmt2} must be of the form 
\begin{equation} \label{eq:x2ymztsol} 
	(x,y,z,t)=(u_1u_4u_5, u_2u_3, u_1u_2u_4^2 \, , u_1u_3u_5^2 \, ) \quad \text{for some integers} \quad u_1, \dots, u_{5}.
 \end{equation} 
The trivial integer solutions to \eqref{eq:x2ymzt} are also included in this description. Hence, all integer solutions to equation \eqref{eq:x2ymzt} are of the form \eqref{eq:x2ymztsol}.

\vspace{10pt}

The next equation to consider is 
\begin{equation} \label{trivial2}
x^3-yzt=0.
\end{equation}
This is an equation of the form \eqref{eq:2mon1}, and we will apply Theorem \ref{th:2mon1}. For this equation, the corresponding equation \eqref{eq:linnonneghom} is
$$
3x_1=y_1+z_1+t_1.
$$
By Proposition \ref{prop:linmin}, the minimal solutions to this equation satisfy $t_1\leq 1$ and $\max\{x_1,y_1,z_1\}\leq 3$, and easy trial and error returns a list of ten minimal solutions $(x_1,y_1,z_1,t_1)$, which are the rows of the following matrix

$$
c = \begin{bmatrix}
	1&1&1&1\\
	1&1&2&0\\
	1&1&0&2\\
	1&2&1&0\\
	1&2&0&1\\
	1&0&1&2\\
	1&0&2&1\\
	1&0&0&3\\
	1&3&0&0\\
	1&0&3&0\\
\end{bmatrix}.
$$
 Using these minimal solutions and Theorem \ref{th:2mon1}, we can find non-trivial integer solutions to the initial equation by \eqref{eq:2mon1sol}:
  \[ x=   \prod_{k=1}^{N} (u_k)^{c_{k1}} = u_1^{c_{11}}u_2^{c_{21}}u_3^{c_{31}}u_4^{c_{41}}u_5^{c_{51}}u_6^{c_{61}}u_7^{c_{71}}u_8^{c_{81}}u_9^{c_{91}}u_{10}^{c_{10\,1}} = u_1u_2u_3u_4u_5u_6u_7u_8u_9u_{10},\] 
\[ y=   \prod_{k=1}^{N} (u_k)^{c_{k2}} = u_1^{c_{12}}u_2^{c_{22}}u_3^{c_{32}}u_4^{c_{42}} u_5^{c_{52}}u_6^{c_{62}}u_7^{c_{72}}u_8^{c_{82}}u_9^{c_{92}}u_{10}^{c_{10\,2}}  = u_1u_2u_3 \, u_4^2 \, u_5^2 \, u_9^3, \] 
\[ z=   \prod_{k=1}^{N} (u_k)^{c_{k3}} = u_1^{c_{13}}u_2^{c_{23}}u_3^{c_{33}}u_4^{c_{43}}u_5^{c_{53}}u_6^{c_{63}}u_7^{c_{73}}u_8^{c_{83}}u_9^{c_{93}}u_{10}^{c_{10\,3}} = u_1u_2^2 \, u_4u_6 \, u_7^2 \, u_{10}^3, \]
\[ t=   \prod_{k=1}^{N} (u_k)^{c_{k4}} = u_1^{c_{14}}u_2^{c_{24}}u_3^{c_{34}}u_4^{c_{44}}u_5^{c_{54}}u_6^{c_{64}}u_7^{c_{74}}u_8^{c_{84}}u_9^{c_{94}}u_{10}^{c_{10\,4}} = u_1u_3^2 \, u_5u_6^2 \, u_7u_8^3. \]

 Therefore, non-trivial integer solutions to equation \eqref{trivial2} must be of the form 
 \begin{equation}\label{triv2} 
 	(x,y,z,t)=(u_1u_2u_3u_4u_5u_6u_7u_8u_9u_{10},u_1u_2u_3 \, u_4^2 \, u_5^2 \, u_9^3,u_1u_2^2 \, u_4u_6 \, u_7^2 \, u_{10}^3,u_1u_3^2 \, u_5u_6^2 \, u_7u_8^3 \,) 
 \end{equation} 
for some integers $u_1, \dots, u_{10}$.
 The trivial integer solutions to \eqref{trivial2} are also included in this description. Hence, all integer solutions to equation \eqref{trivial2} are of the form \eqref{triv2}.

\vspace{10pt}

The next equation to consider is 
\begin{equation}\label{trivial6}
	x^3-y^2z=0.
\end{equation}
This is an equation of the form \eqref{eq:2mon1}, and we will apply Theorem \ref{th:2mon1}. For this equation, the corresponding equation \eqref{eq:linnonneghom} is
$$
3x_1=2y_1+z_1.
$$
By Proposition \ref{prop:linmin}, the minimal solutions to this equation satisfy $x_1 \leq 2$ and $\max\{y_1,z_1\}\leq 3$, and easy trial and error returns a list of three minimal solutions $(x_1,y_1,z_1)$, which are the rows of the following matrix
$$
c = \begin{bmatrix}
	1&1&1\\
	1&0&3\\
	2&3&0\\
\end{bmatrix}.
$$
Using these minimal solutions and Theorem \ref{th:2mon1}, we can find non-trivial integer solutions to the initial equation by \eqref{eq:2mon1sol}:
 \[ x=   \prod_{k=1}^{N} (u_k)^{c_{k1}} = u_1^{c_{11}}u_2^{c_{21}}u_3^{c_{31}} = u_1 u_2 u_3 ^2, \] 
\[ y=   \prod_{k=1}^{N} (u_k)^{c_{k2}} = u_1^{c_{12}}u_2^{c_{22}}u_3^{c_{32}}= u_1 u_3^3, \] 
\[ z=   \prod_{k=1}^{N} (u_k)^{c_{k3}} = u_1^{c_{13}}u_2^{c_{23}}u_3^{c_{33}} = u_1 u_2^3. \]
Therefore any non-trivial integer solution must be of the form 
\begin{equation}\label{trivial6sol}
(x,y,z)=(u_1u_2u_3^2 \, ,u_1u_3^3 \, ,u_1u_2^3 \, ) \quad \text{for some integers} \quad u_1, u_2, u_3.
\end{equation} 
The trivial integer solutions to \eqref{trivial6} are also included in this description. Hence, all integer solutions to equation \eqref{trivial6} are of the form \eqref{trivial6sol}. 

\vspace{10pt}

The next equation to consider is 
\begin{equation}\label{trivial7}
	x^2y-z^2t=0.
\end{equation}
This is an equation of the form \eqref{eq:2mon1}, and we will apply Theorem \ref{th:2mon1}. For this equation, the corresponding equation \eqref{eq:linnonneghom} is
$$
2x_1+y_1-2z_1-t_1=0.
$$
By Proposition \ref{prop:linmin}, the minimal solutions to this equation satisfy $\max\{x_1,y_1\} \leq 2$ and $\max\{z_1,t_1\}\leq 2$, and easy trial and error returns a list of four minimal solutions $(x_1,y_1,z_1,t_1)$, which form the rows of the following matrix
$$
c = \begin{bmatrix}
	0&1&0&1\\
	1&0&1&0\\
	0&2&1&0\\
	1&0&0&2\\
\end{bmatrix}.
$$
Using these minimal solutions and Theorem \ref{th:2mon1}, we can find non-trivial integer solutions to the initial equation by \eqref{eq:2mon1sol}:
 \[ x=   \prod_{k=1}^{N} (u_k)^{c_{k1}} = u_1^{c_{11}}u_2^{c_{21}}u_3^{c_{31}}u_4^{c_{41}} = u_2u_4, \] 
\[ y=   \prod_{k=1}^{N} (u_k)^{c_{k2}} = u_1^{c_{12}}u_2^{c_{22}}u_3^{c_{32}}u_4^{c_{42}}= u_1u_3^2, \] 
\[ z=   \prod_{k=1}^{N} (u_k)^{c_{k3}} = u_1^{c_{13}}u_2^{c_{23}}u_3^{c_{33}}u_4^{c_{43}} =u_2u_3, \]
\[ t=   \prod_{k=1}^{N} (u_k)^{c_{k4}} = u_1^{c_{14}}u_2^{c_{24}}u_3^{c_{34}}u_4^{c_{44}} = u_1u_4^2. \]
Therefore any non-trivial integer solution must be of the form 
\begin{equation}\label{trivial7sol}
(x,y,z,t)=(u_2u_4 ,u_1u_3^2 \, ,u_2u_3, u_1u_4^2 \, ) \quad \text{for some integers} \quad u_1, \dots, u_4.
\end{equation}
The trivial integer solutions to \eqref{trivial7} are also included in this description. Hence, all integer solutions to equation \eqref{trivial7} are of the form \eqref{trivial7sol}.

\vspace{10pt}

The next equation to consider is 
\begin{equation} \label{trivial4}
x_1x_2x_3- x_4^2 \, x_5=0.
\end{equation}
This is an equation of the form \eqref{eq:2mon1}, and we will apply Theorem \ref{th:2mon1}. The corresponding equation \eqref{eq:linnonneghom} is $$
z_1+z_2+z_3-2z_4-z_5=0.
$$
 By Proposition \ref{prop:linmin}, the minimal solutions to this equation satisfy $\max\{z_4,z_5\}\leq 1$ and $\max\{z_1,z_2,z_3\}\leq 2$, and easy trial and error returns a list of nine minimal solutions $(z_1,z_2,z_3,z_4,z_5)$ which form the rows of the following matrix
 $$
 c = \begin{bmatrix}
 	1&1&0&1&0\\
 	1&0&1&1&0\\
 	0&1&1&1&0\\
 	1&0&0&0&1\\
 	0&1&0&0&1\\
 	0&0&1&0&1\\
 	2&0&0&1&0\\
 	0&2&0&1&0\\
 	0&0&2&1&0\\
 \end{bmatrix}.
 $$
Using these minimal solutions and Theorem \ref{th:2mon1}, we can find non-trivial integer solutions to the initial equation  by \eqref{eq:2mon1sol}:
 \[ x_1=   \prod_{k=1}^{N} (u_k)^{c_{k1}} = u_1^{c_{11}}u_2^{c_{21}}u_3^{c_{31}}u_4^{c_{41}}u_5^{c_{51}}u_6^{c_{61}}u_7^{c_{71}}u_8^{c_{81}}u_9^{c_{91}} = u_1u_2u_4 u_7^2, \] 
\[ x_2=   \prod_{k=1}^{N} (u_k)^{c_{k2}} = u_1^{c_{12}}u_2^{c_{22}}u_3^{c_{32}}u_4^{c_{42}} u_5^{c_{52}}u_6^{c_{62}}u_7^{c_{72}}u_8^{c_{82}}u_9^{c_{92}} = u_1u_3u_5 u_8^2, \] 
\[ x_3=   \prod_{k=1}^{N} (u_k)^{c_{k3}} = u_1^{c_{13}}u_2^{c_{23}}u_3^{c_{33}}u_4^{c_{43}}u_5^{c_{53}}u_6^{c_{63}}u_7^{c_{73}}u_8^{c_{83}}u_9^{c_{93}} =u_2u_3u_6u_9^2, \]
\[ x_4=   \prod_{k=1}^{N} (u_k)^{c_{k4}} = u_1^{c_{14}}u_2^{c_{24}}u_3^{c_{34}}u_4^{c_{44}}u_5^{c_{54}}u_6^{c_{64}}u_7^{c_{74}}u_8^{c_{84}}u_9^{c_{94}}= u_1u_2u_3u_7u_8u_9, \]
\[ x_5=   \prod_{k=1}^{N} (u_k)^{c_{k5}} = u_1^{c_{15}}u_2^{c_{25}}u_3^{c_{35}}u_4^{c_{45}}u_5^{c_{55}}u_6^{c_{65}}u_7^{c_{75}}u_8^{c_{85}}u_9^{c_{95}}= u_4u_5u_6. \]
Therefore, any non-trivial integer solution must be of the form 
\begin{equation}\label{triv4} 
(x_1,x_2,x_3,x_4,x_5)=(u_1u_2u_4 u_7^2 \, ,u_1u_3u_5 u_8^2 \, ,u_2u_3u_6u_9^2\, ,u_1u_2u_3u_7u_8u_9,u_4u_5u_6),
\end{equation} 
for some integers $u_1, \dots, u_9$. 
The trivial integer solutions to \eqref{trivial4} are also included in this description. Hence, all integer solutions to equation \eqref{trivial4} are of the form \eqref{triv4}.

\vspace{10pt}

The final equation to consider is 
\begin{equation}\label{trivial3}
	x_1x_2x_3-x_4x_5x_6=0.
\end{equation}
This is an equation of the form \eqref{eq:2mon1}, and we will apply Theorem \ref{th:2mon1}. The corresponding equation \eqref{eq:linnonneghom} is 
$$
z_1+z_2+z_3-z_4-z_5-z_6=0.
$$
 By Proposition \ref{prop:linmin}, the minimal solutions to this equation satisfy $\max\{z_4,z_5,z_6\}\leq 1$ and $\max\{z_1,z_2,z_3\}\leq 1$, and easy trial and error returns a list of nine minimal solutions $(z_1,z_2,z_3,z_4,z_5,z_6)$, which form the rows of the following matrix
 $$
 c = 
 \begin{bmatrix}
 	1&0&0&1&0&0\\
 	1&0&0&0&1&0\\
 	1&0&0&0&0&1\\
 	0&1&0&1&0&0\\
 	0&1&0&0&1&0\\
 	0&1&0&0&0&1\\
 	0&0&1&1&0&0\\
 	0&0&1&0&1&0\\
 	0&0&1&0&0&1\\
 \end{bmatrix}.
 $$
Using these minimal solutions and Theorem \ref{th:2mon1}, we can find non-trivial integer solutions to the initial equation by \eqref{eq:2mon1sol}:
 \[ x_1=   \prod_{k=1}^{N} (u_k)^{c_{k1}} = u_1^{c_{11}}u_2^{c_{21}}u_3^{c_{31}}u_4^{c_{41}}u_5^{c_{51}}u_6^{c_{61}}u_7^{c_{71}}u_8^{c_{81}}u_9^{c_{91}} = u_1u_2u_3,\] 
\[ x_2=   \prod_{k=1}^{N} (u_k)^{c_{k2}} = u_1^{c_{12}}u_2^{c_{22}}u_3^{c_{32}}u_4^{c_{42}} u_5^{c_{52}}u_6^{c_{62}}u_7^{c_{72}}u_8^{c_{82}}u_9^{c_{92}} = u_4u_5u_6, \] 
\[ x_3=   \prod_{k=1}^{N} (u_k)^{c_{k3}} = u_1^{c_{13}}u_2^{c_{23}}u_3^{c_{33}}u_4^{c_{43}}u_5^{c_{53}}u_6^{c_{63}}u_7^{c_{73}}u_8^{c_{83}}u_9^{c_{93}} =u_7u_8u_9,  \]
\[ x_4=   \prod_{k=1}^{N} (u_k)^{c_{k4}} = u_1^{c_{14}}u_2^{c_{24}}u_3^{c_{34}}u_4^{c_{44}}u_5^{c_{54}}u_6^{c_{64}}u_7^{c_{74}}u_8^{c_{84}}u_9^{c_{94}}= u_1u_4u_7, \]
\[ x_5=   \prod_{k=1}^{N} (u_k)^{c_{k5}} = u_1^{c_{15}}u_2^{c_{25}}u_3^{c_{35}}u_4^{c_{45}}u_5^{c_{55}}u_6^{c_{65}}u_7^{c_{75}}u_8^{c_{85}}u_9^{c_{95}}= u_2u_5u_8, \]
\[ x_6=   \prod_{k=1}^{N} (u_k)^{c_{k6}} = u_1^{c_{16}}u_2^{c_{26}}u_3^{c_{36}}u_4^{c_{46}}u_5^{c_{56}}u_6^{c_{66}}u_7^{c_{76}}u_8^{c_{86}}u_9^{c_{96}}= u_3u_6u_9. \]
Therefore, any non-trivial integer solution must be of the form 
\begin{equation}\label{triv3}
		(x_1,x_2,x_3,x_4,x_5,x_6)=(u_1u_2u_3,u_4u_5u_6,u_7u_8u_9,u_1u_4u_7,u_2u_5u_8,u_3u_6u_9), 
\end{equation} 
for some integers $u_1, \dots, u_9$. 
The trivial integer solutions to \eqref{trivial3} are also included in this description. Hence, all integer solutions to equation \eqref{trivial3} are of the form \eqref{triv3}.

Table \ref{table1.38} presents all integer solutions to the two-monomial equations listed in Table \ref{tab:H16mon2a}.

\begin{center}
\begin{tabular}{ |c|c|c|} 
		\hline
		Equation & Integer solutions $(x,y)$ or $(x,y,z)$ or $(x,y,z,t)$ or $(x_1,\dots, x_n)$ \\
		\hline\hline
		$x^2-yz=0$&$(u_1 u_2 u_3,u_1 u_2^2,u_1 u_3^2)$ \\\hline
		$xy-zt=0$&$(u_1 u_2, u_3 u_4, u_1 u_3, u_2 u_4)$ \\\hline
		$x^3-y^2=0$&$(u_1 ^2,u_1 ^3)$ \\\hline
		$x^2y-z^2=0$&$(0,u_1,0),(u_1,u_2^2, u_1 u_2)$\\\hline
		$x^3-yz=0$&$(u_1u_2u_3u_4,u_2^3 \, u_3^2\, u_4,u_1^3\,u_3 u_4 ^2)$\\\hline
		$xyz-t^2=0$&$(u_1 u_2 u_4 ^2 \, ,u_2 u_3  u_5^2 \, ,u_1 u_3 u_6^2 \, ,u_1 u_2 u_3 u_4 u_5 u_6)$\\\hline
		$x^2y-zt=0$&$(u_1u_4u_5, u_2u_3, u_1u_2u_4^2 \, , u_1u_3u_5^2 \, )$\\\hline
		$x_1x_2x_3-x_4x_5=0$&$(u_1u_2,u_3u_4,u_5u_6,u_1u_3u_5,u_2u_4u_6)$ \\\hline
		$x^3-yzt=0$ & $(u_1u_2u_3u_4u_5u_6u_7u_8u_9u_{10},u_1u_2u_3 \, u_4^2 \, u_5^2 \, u_9^3,$ \\ & $u_1u_2^2 \, u_4u_6 \, u_7^2 \, u_{10}^3,u_1u_3^2 \, u_5u_6^2 \, u_7u_8^3 \,)$\\\hline
		$x^3-y^2z=0$&$(u_1u_2u_3^2 \, ,u_1u_3^3 \, ,u_1u_2^3 \, )$\\\hline
		$x^2y-z^2t=0$&$(u_2u_4 ,u_1u_3^2 \, ,u_2u_3, u_1u_4^2 \, )$\\\hline
		$x_1x_2x_3- x_4^2 \, x_5=0$&$(u_1u_2u_4 u_7^2 \, ,u_1u_3u_5 u_8^2 \, ,u_2u_3u_6u_9^2\, ,u_1u_2u_3u_7u_8u_9,u_4u_5u_6) $\\\hline
		$x_1x_2x_3-x_4x_5x_6=0$&$(u_1u_2u_3,u_4u_5u_6,u_7u_8u_9,u_1u_4u_7,u_2u_5u_8,u_3u_6u_9)$\\\hline
	\end{tabular}
	\captionof{table}{Integer solutions to the equations listed in Table \ref{tab:H16mon2a}. Assume that $u_1,\dots,u_{10}$ are integers. \label{table1.38}}
\end{center}

\subsection{Exercise 1.42}\label{ex:2monb}
\textbf{\emph{Table \ref{tab:H16mon2b} lists two-monomial equations of size $H\leq 16$. For each equation, describe all its integer solutions.  }}

\begin{center}
	\begin{tabular}{ |c|c|c|c|c|c| } 
		\hline
		$H$ & Equation &  $H$ & Equation & $H$ & Equation \\ 
		\hline\hline
		$12$ & $2x^2-y^2=0$ & $16$ & $x^3-2y^2=0$ & $16$ & $xyz-2t^2=0$ \\ 
		\hline
		$12$ & $2x^2-yz=0$ & $16$ & $3x^2-yz=0$ & $16$ & $x^2y-2zt=0$ \\
		\hline
		$12$ & $x^2-2yz=0$  & $16$ & $x^2-3yz=0$ & $16$ & $3xy-zt=0$ \\
		\hline
		$12$ & $2xy-zt=0$ & $16$ & $x^2y-2z^2=0$ & $16$ &  $x_1x_2x_3-2x_4x_5=0$  \\
		\hline
		$16$ & $3x^2-y^2=0$ & $16$ & $x^3-2yz=0$ &  & \\
		\hline
	\end{tabular}
	\captionof{table}{\label{tab:H16mon2b} The remaining two-monomial equations of size $H\leq 16$.}
\end{center}

Equation
$$
	2x^2-y^2=0
$$
is solved in Section 1.6.1 of the the book, and its unique integer solution is
$$
	(x,y)=(0,0).
$$

\vspace{10pt}

Equation
$$
	x^2-2yz=0
$$
is solved in Section 1.6.3 of the book, and its integer solutions are
$$
	(x,y,z)=(2uvw,2uv^2,uw^2) \quad \text{or} \quad (2uvw,uv^2,2uw^2) \quad \text{for some integers} \quad u,v,w.
$$

\vspace{10pt}

Let us now consider the equation
\begin{equation}\label{eq:2x2myz}
	2x^2-yz=0.
\end{equation}
We can see from this equation that $yz$ is even, so either $y$ is even, or $z$ is even. 
In the first case, after substituting $y=2y'$, for some integer $y'$, into \eqref{eq:2x2myz} and cancelling $2$, we obtain, up to the names of variables, equation \eqref{eq:x2myz}. By \eqref{eq:x2myzsol}, its integer solutions are 
$
(x,y',z)=(uvw,uv^2,uw^2) 
$
for some integers $u,v,w$, or in the original variables, $(x,y,z)=(uvw,2uv^2,uw^2)$. The second case is similar. Finally, all integer solutions to \eqref{eq:2x2myz} are
$$
	(x,y,z)=(uvw,2uv^2,uw^2) \quad \text{or} \quad (uvw,uv^2,2uw^2) \quad \text{for some integers} \quad u,v,w. 
$$

\vspace{10pt}

The next equation to consider is
\begin{equation}\label{eq:2xymzt}
	2xy-zt=0.
\end{equation}
We can see from this equation that $zt$ is even, so either $z$ is even, or $t$ is even. In the first case, after substituting $z=2z'$, for some integer $z'$, into \eqref{eq:2xymzt} and cancelling $2$, we obtain, up to the names of variables, equation \eqref{eq:xymzt}. By \eqref{eq:xymztsol}, its integer solutions are
$
(x,y,z',t)=(uv,wr,uw,vr) 
$
for some integers $u,v,w,r$, or in the original variables, $(x,y,z,t)=(uv,wr,2uw,vr)$. The second case is similar. Hence, all integer solutions to \eqref{eq:2xymzt} are 
$$
	(x,y,z,t)=(uv,wr,2uw,vr) \quad \text{or} \quad (uv,wr,uw,2vr) \quad \text{for some integers} \quad u,v,w,r. 
$$

\vspace{10pt}

The next equation to consider is
$$
	3x^2-y^2=0.
$$
This equation is of the form $dx^k-y^k=0$. Because $\sqrt{3}$ is not integer, Remark \ref{remark:dxkmyk} implies that this equation has only the trivial solution 
$$
(x,y)=(0,0). 
$$ 

\vspace{10pt}

The next equation to consider is
\begin{equation}\label{eq:x3m2y2}
	x^3-2y^2=0.
\end{equation}
From the equation it is clear that $x$ is even. After substituting $x=2x'$ for some integer $x'$ into \eqref{eq:x3m2y2} and cancelling $2$, we obtain $4(x')^3-y^2=0$. From this, it is clear that $y$ is even, so let $y=2y'$ for some integer $y'$. After substituting this into the previous equation and cancelling $4$, we obtain, up to the names of variables, equation \eqref{eq:x3my2}. By \eqref{eq:x3my2sol}, its integer solutions are $(x',y')=(u^2,u^3)$ for some integer $u$. Then, in the original variables we have that the integer solutions to \eqref{eq:x3m2y2} are 
$$
(x,y)=(2u^2,2u^3) \quad \text{for some integer} \quad u. 
$$

\vspace{10pt}

The next equation to consider is
\begin{equation}\label{eq:3x2myz}
	3x^2-yz=0.
\end{equation}
We can see from this equation that $yz$ is a multiple of 3, so either $y$ is a multiple of 3, or $z$ is a multiple of 3. 
In the first case, after substituting $y=3y'$ and cancelling $3$, up to the names of variables, we obtain equation \eqref{eq:x2myz}. By \eqref{eq:x2myzsol}, its integer solutions are $(x,y',z)=(uvw,uv^2,uw^2)$ for some integers $u,v,w$, or in the original variables, $(x,y,z)=(uvw,3uv^2,uw^2)$. The second case is similar. Hence, all integer solutions to \eqref{eq:3x2myz} are
\begin{equation}\label{eq:3x2myzsol}
(x,y,z)=(uvw,3uv^2,uw^2) \quad \text{or} \quad (uvw,uv^2,3uw^2) \quad \text{for some integers} \quad u,v,w. 
\end{equation}

\vspace{10pt}

The next equation to consider is
\begin{equation}\label{eq:x2m3yz}
	x^2-3yz=0.
\end{equation}
From the equation it is clear that $x$ is a multiple of 3, so let $x=3x'$ for some integer $x'$. After substituting this into \eqref{eq:x2m3yz} and cancelling $3$, we obtain $3(x')^2-yz=0$. Up to the names of variables, this is equation \eqref{eq:3x2myz}, and by \eqref{eq:3x2myzsol}, its integer solutions are $(x',y,z)=(uvw,3uv^2,uw^2)$ or $(uvw,uv^2,3uw^2)$ for some integers $u,v,w$, or in the original variables, we have that the integer solutions to equation \eqref{eq:x2m3yz} are
$$
(x,y,z)=(3uvw,uv^2,3uw^2) \quad \text{or} \quad (3uvw,3uv^2,uw^2) \quad \text{for some integers} \quad u,v,w. 
$$

\vspace{10pt}

The next equation to consider is 
\begin{equation}\label{eq:x2ym2z2}
	x^2y-2z^2=0.
\end{equation}
If $x=0$ then $z=0$, and we can take $y=u$ arbitrary, so that $(x,y,z)=(0,u,0)$ for integer $u$. Now assume that $x\neq 0$, and let us prove that $y$ is even. If $z=0$ then $y=0$ is even. If $z\neq 0$, then $y\neq 0$. Let $x_2,y_2,z_2$ be the exponents with which prime $p=2$ enters the prime factorizations of $x,y,z$, respectively. Then $x^2y=2z^2$ implies that $2x_2+y_2=1+2z_2$, from which it is clear that $y_2$ is odd, and in particular $y_2\neq 0$. Hence, $y$ is even.
So, let $y=2y'$ for some integer $y'$, then up to the names of variables, we obtain equation \eqref{eq:x2ymz2}. By \eqref{eq:x2ymz2sol}, its integer solutions are $(x,y',z)=(0,u,0)$ or $(u,v^2,uv)$ for some integers $u,v$. Hence, in the original variables, we obtain that all integer solutions to \eqref{eq:x2ym2z2} are
$$
(x,y,z)=(0,u \, ,0) \quad \text{or} \quad (u,2v^2 \, ,uv) \quad \text{for some integers} \quad u,v . 
$$ 

\vspace{10pt}

The next equation to consider is
\begin{equation}\label{eq:x3m2yz}
	x^3-2yz=0.
\end{equation}
From this equation, we can see that $x^3$ is even, therefore $x$ is even. Let $x=2x'$ for some integer $x'$. After substituting this into \eqref{eq:x3m2yz} and cancelling $2$, we obtain, 
\begin{equation}\label{eq:x3m2yzred}
4(x')^3-yz=0,
\end{equation}
 this gives three cases, (i) $y$ is even and $z$ is even, (ii) $y$ is a multiple of 4, or, (iii) $z$ is a multiple of 4. Let us look at the first case. Let $y=2y'$ and $z=2z'$ for some integers $y'$ and $z'$, substituting these into \eqref{eq:x3m2yzred} and cancelling $4$, we obtain, up to the names of variables, equation \eqref{eq:x3myz}. By \eqref{eq:x3myzsol}, its integer solutions are $(x',y',z')=(u_1u_2u_3u_4,u_2^3 \, u_3^2 \, u_4,u_1^3 \, u_3 u_4^2)$ for integers $u_1,u_2,u_3,u_4$, or in the original variables, $(x,y,z)=(2u_1u_2u_3u_4,2u_2^3 \, u_3^2 \, u_4,2u_1^3 \, u_3 u_4^2)$. The second and third cases are similar. Hence, all integer solutions to equation \eqref{eq:x3m2yz} are
$$
\begin{aligned}
	(x,y,z)=(2u_1u_2u_3u_4,2u_2^3 \, u_3^2 \, u_4,2u_1^3 \, u_3 u_4^2) \quad \text{or} \quad (2u_1u_2u_3u_4,4u_2^3 \, u_3^2 \, u_4,u_1^3 \, u_3 u_4^2) \\ \text{or} \quad (2u_1u_2u_3u_4,u_2^3 \, u_3 ^2 \, u_4,4u_1^3 \, u_3 u_4 ^2) \quad \text{for some integers} \quad u_1,u_2,u_3 ,u_4. 
\end{aligned}
$$

\vspace{10pt}

The next equation to consider is
\begin{equation}\label{eq:xyzm2t2}
	xyz-2t^2=0.
\end{equation}
From this equation, we can see that $xyz$ is even, so either $x$ is even, $y$ is even, or $z$ is even. For the first case, let $x=2x'$. After substituting this into \eqref{eq:xyzm2t2} and cancelling $2$, we obtain, up to the names of variables, equation \eqref{eq:xyzmt2}. By \eqref{eq:xyzmt2sol}, its integer solutions are $(x',y,z,t)=(u_1 u_2 u_4 ^2 \, ,u_2 u_3  u_5^2 \, ,u_1 u_3 u_6^2 \, ,u_1 u_2 u_3 u_4 u_5 u_6)$ for integers $u_1,\dots,u_6$, or in the original variables 
$$
(x,y,z,t)=(2u_1 u_2 u_4 ^2 \, ,u_2 u_3  u_5^2 \, ,u_1 u_3 u_6^2 \, ,u_1 u_2 u_3 u_4 u_5 u_6).
$$ 
The second and third cases are similar. Hence, all integer solutions to \eqref{eq:xyzm2t2} are 
$$
\begin{aligned}
	(x,y,z,t)=& \, (2u_1 u_2 u_4 ^2 \, ,u_2 u_3  u_5^2 \, ,u_1 u_3 u_6^2 \, ,u_1 u_2 u_3 u_4 u_5 u_6), \quad \text{or,} \\ & (u_1 u_2 u_4 ^2 \, ,2u_2 u_3  u_5^2 \, ,u_1 u_3 u_6^2 \, ,u_1 u_2 u_3 u_4 u_5 u_6), \quad \text{or,} \\ & (u_1 u_2 u_4 ^2 \, ,u_2 u_3  u_5^2 \, ,2u_1 u_3 u_6^2 \, ,u_1 u_2 u_3 u_4 u_5 u_6) \quad \text{for some integers} \quad u_1,\ldots, u_6. 
\end{aligned}
$$

\vspace{10pt}

The next equation to consider is 
\begin{equation} \label{eq:x2ym2zt}
	x^2y-2zt=0.
\end{equation}
We can see from this equation that $x^2y$ is even. So either $y$ is even, or $x^2$ is even. Let us look at when $y$ is even, so let $y=2y'$ for some integer $y'$. Substituting this into \eqref{eq:x2ym2zt} and cancelling $2$, we obtain, up to the names of variables, equation \eqref{eq:x2ymzt}. By \eqref{eq:x2ymztsol}, its integer solutions are $(x,y',z,t)=(u_1u_4u_5, u_2u_3, u_1u_2u_4^2 \, ,$ $ u_1u_3u_5^2  \, )$ for some integers $u_1, \ldots, u_5,$ or in the original variables, $(x,y,z,t)=(u_1u_4u_5, 2u_2u_3, u_1u_2u_4^2 \, , u_1u_3u_5^2  \, )$.

Next let us look at when $x^2$ is even. If $x^2$ is even, then $x$ is even, so let $x=2x'$ for some integer $x'$, then we have $2(x') ^2 y-zt=0$. Now we can see from this equation that $zt$ is even, so either $z$ is even, or $t$ is even. Let $z$ be even, so $z=2z'$ for some integer $z'$, after substituting this into the previous equation and cancelling $2$, we again obtain, up to the names of variables, equation \eqref{eq:x2ymzt}. Then in the original variables we obtain $(x,y,z,t)=(2u_1u_4u_5, u_2u_3, 2u_1u_2u_4^2 \, , u_1u_3u_5^2 \, )$. 
The case for when $t$ is even is similar. Finally, we may conclude that all integer solutions to equation \eqref{eq:x2ym2zt} are
$$
	\begin{aligned}
(x,y,z,t)=(u_1u_4u_5, 2u_2u_3, u_1u_2u_4^2 \, , u_1u_3u_5^2  \, ) \quad \text{or} \quad (2u_1u_4u_5, u_2u_3, u_1u_2u_4^2 \, , 2u_1u_3u_5^2  \, ) \quad \text{or} \\ (2u_1u_4u_5, u_2u_3, u_1u_2u_4^2 \, , 2u_1u_3u_5^2  \, ) \quad \text{for some integers} \quad u_1, \ldots, u_5.
\end{aligned}
$$

\vspace{10pt}

The next equation to consider is
\begin{equation}\label{eq:3xymzt}
	3xy-zt=0.
\end{equation}
We can see from this equation that $zt$ is a multiple of 3, so either $z$ is a multiple of 3, or $t$ is a multiple of 3. In the first case, after substituting $z=3z'$ for some integer $z'$ into \eqref{eq:3xymzt} and cancelling $3$, we obtain, up to the names of variables, equation \eqref{eq:xymzt}. By \eqref{eq:xymztsol}, its solutions are $(x,y,z',t)=(uv,wr,uw,vr)$ for some integers $u,v,w,r$, or in the original variables, $(x,y,z,t)=(uv,wr,3uw,vr)$. The second case is similar. Hence, all integer solutions to equation \eqref{eq:3xymzt} are
$$
(x,y,z,t)=(uv,wr,3uw,vr) \quad \text{or} \quad (uv,wr,uw,3vr) \quad \text{for some integers} \quad u,v,w,r. 
$$

\vspace{10pt}

The final equation to consider is
\begin{equation}\label{eq:x1x2x3m2x4x5}
	x_1x_2x_3-2x_4x_5=0.
\end{equation}
We can see from this equation that either (i) $x_1$ is even, (ii) $x_2$ is even, or (iii) $x_3$ is even. In the first case, substituting $x_1=2x_1'$ for some integer $x_1'$ into \eqref{eq:x1x2x3m2x4x5} and cancelling $2$, we obtain, up to the names of variables, equation \eqref{eq:x1x2x3mx4x5}. By \eqref{eq:x1x2x3mx4x5sol}, its integer solutions are $(x_1',x_2,x_3,x_4,x_5) = (u_1u_2,u_3u_4,u_5u_6,u_1u_3u_5,u_2u_4u_6)$, for some integers $u_1,\ldots, u_6$, or in the original variables, 
$$
(x_1,x_2,x_3,x_4,x_5) = (2u_1u_2,u_3u_4,u_5u_6,u_1u_3u_5,u_2u_4u_6) \quad \text{for some integers} \quad u_1,\ldots, u_6.
$$
The second and third cases are similar. Hence, all integer solutions to \eqref{eq:x1x2x3m2x4x5} are 
$$
	\begin{aligned}
		(x_1,x_2,x_3,x_4,x_5)=(2u_1u_2,u_3u_4,u_5u_6,u_1u_3u_5,u_2u_4u_6) \quad \text{or} \quad (u_1u_2,2u_3u_4,u_5u_6,u_1u_3u_5,u_2u_4u_6) \\ \text{or} \quad (u_1u_2,u_3u_4,2u_5u_6,u_1u_3u_5,u_2u_4u_6) \quad \text{for some integers} \,\, u_1,\ldots, u_6. 
	\end{aligned}
$$

Table \ref{tab:H16mon2bsol} presents all integer solutions to the two-monomial equations listed in Table \ref{tab:H16mon2b}.

	\begin{center}
	\begin{tabular}{ |c|c|c|} 
		\hline
		Equation & Integer solutions $(x,y)$ or $(x,y,z)$ or $(x,y,z,t)$ or $(x_1,\dots, x_5)$ \\
		\hline\hline
		$2x^2-y^2=0$&$(0,0)$ \\\hline
		$2x^2-yz=0$& $(u_1 u_2 u_3,2u_1 u_2^2,u_1 u_3^2),(u_1 u_2 u_3,u_1u_2^2,2u_1 u_3^2) $\\\hline
		$x^2-2yz=0$&$(2u_1 u_2 u_3,2u_1u_2^2,u_1u_3^2),(2u_1 u_2 u_3,u_1 u_2^2,2u_1 u_3^2)$ \\\hline
		$2xy-zt=0$&$(u_1 u_2,u_3 u_4,2u_1 u_3,u_2 u_4),(u_1 u_2,u_3 u_4,u_1 u_3,2u_2 u_4)$\\\hline
		$3x^2-y^2=0$&$(0,0)$\\\hline
		$x^3-2y^2=0$&$(2u_1^2,2u_1^3)$\\\hline
		$3x^2-yz=0$&$(u_1 u_2 u_3,3u_1 u_2^2,u_1 u_3^2), (u_1 u_2 u_3,u_1 u_2^2,3u_1 u_3^2)$\\\hline
		$x^2-3yz=0$&$(3u_1 u_2 u_3,u_1 u_2^2,3u_1 u_3^2), (3u_1 u_2 u_3,3u_1 u_2^2,u_1u_3^2)$\\\hline
		$x^2y-2z^2=0$ & $(0,u_1,0),(u,2u_2^2 ,u_1 u_2)$\\\hline
		$x^3-2yz=0$&$(2u_1u_2u_3u_4,2u_2^3 \, u_3^2 \, u_4,2u_1^3 \, u_3 u_4^2),$\\
		&$  (2u_1u_2u_3u_4,4u_2^3 \, u_3^2 \, u_4,u_1^3 \, u_3 u_4^2),$\\
		&$(2u_1u_2u_3u_4,u_2^3 \, u_3 ^2 \, u_4,4u_1^3 \, u_3 u_4 ^2) $\\\hline
		$xyz-2t^2=0$&$(2u_1 u_2 u_4 ^2 \, ,u_2 u_3  u_5^2 \, ,u_1 u_3 u_6^2 \, ,u_1 u_2 u_3 u_4 u_5 u_6), $\\
		&$(u_1 u_2 u_4 ^2 \, ,2u_2 u_3  u_5^2 \, ,u_1 u_3 u_6^2 \, ,u_1 u_2 u_3 u_4 u_5 u_6),$\\
		&$(u_1 u_2 u_4 ^2 \, ,u_2 u_3  u_5^2 \, ,2u_1 u_3 u_6^2 \, ,u_1 u_2 u_3 u_4 u_5 u_6)$\\\hline
		$x^2y-2zt=0$&$(u_1u_4u_5, 2u_2u_3, u_1u_2u_4^2 \, , u_1u_3u_5^2  \, ),$\\
		&$(2u_1u_4u_5, u_2u_3, 2u_1u_2u_4^2 \, , u_1u_3u_5^2 \, ),$\\
		&$(2u_1u_4u_5, u_2u_3, u_1u_2u_4^2 \, , 2u_1u_3u_5^2  \, )$\\\hline
		$3xy-zt=0$&$(u_1 u_2,u_3 u_4,3u_1 u_3,u_2 u_4),(u_1 u_2 ,u_3 u_4 ,u_1 u_3,3u_2 u_4)$\\\hline
		$x_1x_2x_3-2x_4x_5=0$&$(2u_1u_2,u_3u_4,u_5u_6,u_1u_3u_5,u_2u_4u_6),$\\
		&$(u_1u_2,2u_3u_4,u_5u_6,u_1u_3u_5,u_2u_4u_6),$\\
		&$ (u_1u_2,u_3u_4,2u_5u_6,u_1u_3u_5,u_2u_4u_6) $\\\hline
	\end{tabular}
	\captionof{table}{Integer solutions to the equations listed in Table \ref{tab:H16mon2b}. Assume that $u_1,\dots,u_{6}$ are integers. \label{tab:H16mon2bsol}}
\end{center}

\section{Chapter 2}

The following theorem proved in \cite{vaserstein2010polynomial} is a useful result for this chapter.
\begin{theorem}\label{th:vaserstein2010}
[Theorem 2.1 in the book] The set of all integer solutions to
\begin{equation}\label{eq:xypztp1}
xy-zt=1
\end{equation} 
is a polynomial family with $46$ parameters. In other words, there exist polynomials $U(u)$, $V(u)$, $W(u)$ and $R(u)$ in $46$ variables $u=(u_1,\dots,u_{46})$ with integer coefficients such that 
$$
U(u)V(u) - W(u)R(u) = 1 \quad\quad \forall u \in {\mathbb Z}^{46},
$$
and, conversely, for every integer solution $(x,y,z,t)$ of \eqref{eq:xypztp1}, there exists $u\in {\mathbb Z}^{46}$ such that $x=U(u)$, $y=V(u)$, $z=W(u)$ and $t=R(u)$. 
\end{theorem}

An example of an equation which can be reduced to equation \eqref{eq:xypztp1} is 
\begin{equation}\label{eq:x2mxmyz}
x^2-x-yz=0
\end{equation}
which is solved in the book, and its integer solutions are
\begin{equation}\label{eq:x2mxmyzsol}
(x,y,z)=(uv,uw,vr) \quad \text{for some integers} \quad u,v,w,r \quad \text{such that} \quad uv-wr=1.
\end{equation}

\subsection{Exercise 2.6}\label{ex:x2pyzpc}
\textbf{\emph{Solve equation 
\begin{equation}\label{eq:x2pyzpc}
yz-x^2=c
\end{equation} 
for $c=\pm 2, \pm 3, \pm 4$ and $\pm 5$.}}

To solve these equations, we must first look at the following definition.
\begin{definition}\label{def:symreduced}[Definition 2.4 in the book]
A $2\times 2$ symmetric matrix 
\begin{equation}\label{eq:symreduced}
A=\begin{bmatrix}
a & b \\
b & d
\end{bmatrix}
\end{equation}
with integer entries 
is called \emph{reduced} if either $a=d=0$ and $b\geq 0$, or the following conditions are satisfied:
\begin{itemize}
\item[(i)] $(1-|a|)/2\leq b \leq |a|/2$, and
\item[(ii)] $|d|\geq |a|>0$, and
\item[(iii)] if $|d|=|a|$ then $b\geq 0$, and
\item[(iv)] if $|d|=|a|$ and $b=0$ then $a\geq d$.
\end{itemize}
\end{definition}

Let $c$ be the determinant of the matrix, then conditions (i) and (ii) mean that
\begin{equation}\label{eq:symreducedbounds}
|c|=|ad-b^2| \geq |a||d|-|b|^2 \geq a^2 - \frac{a^2}{4} = \frac{3a^2}{4} \geq 3b^2,
\end{equation}
hence, for any $c$, we have only finitely many options for $a$ and $b$, and then $d=\frac{c+b^2}{a}$ is determined uniquely.

\begin{proposition}\label{prop:symactions}[Proposition 2.5 in the book]
Every symmetric $2\times 2$ matrix with integer entries and non-zero determinant is equivalent to some reduced matrix. 
\end{proposition}

Let $M_j, j=1,\dots,k$ be the full list of all reduced matrices with determinant $c$, and let $(a_j,b_j,d_j)$, $j=1,\dots,k$ be the corresponding entries. 
Then Proposition \ref{prop:symactions} implies that for any solution $(x,y,z)$ to \eqref{eq:x2pyzpc} we have in component form
\begin{equation}\label{eq:x2pyzpcsol}
(x,y,z) = (a_j uw+b_j(uv+wr)+d_j vr, \,\, a_j u^2 + 2 b_j r u + d_j r^2, \,\, a_j w^2 + 2 b_j v w + d_j v^2) 
\end{equation}
for some $j \in \{1,\dots,k\}$ and some integers $u,v,w,r$ such that $uv-wr=1$.

Equation
\begin{equation}\label{eq:x2pyzp1}
	yz-x^2=1
\end{equation}
is solved in Section 2.1.3 of the book and its integer solutions are
\begin{equation}\label{eq:x2pyzp1sol}
	(x,y,z) = \pm(uw+vr, \,\, u^2+r^2, \,\, v^2 + w^2), 
	\quad \text{where} \,\, uv-wr=1.
\end{equation}

\vspace{10pt}

Equation
\begin{equation}\label{eq:x2pyzm1}
	yz-x^2=-1
\end{equation}
is solved in Section 2.1.3 of the book and its integer solutions are
\begin{equation}\label{eq:x2pyzm1sol}
	(x,y,z) = (uw-vr, \,\, u^2-r^2, \,\, w^2-v^2) \quad \text{or} \quad (uv+wr, \,\,2 r u, \,\, 2 v w), \quad \text{where} \,\, uv-wr=1.
\end{equation}

\vspace{10pt}

Let us first consider equation \eqref{eq:x2pyzpc} with $c=2$, or
\begin{equation}\label{eq:yzmx2m2}
	yz-x^2=2.
\end{equation}
Let $A$ given by \eqref{eq:symreduced} be a reduced matrix with determinant $ad-b^2=2$. If $a=0$, then $b^2=-2$, a contradiction, hence $a \neq 0$, and conditions (i)-(iv) in Definition \ref{def:symreduced} must be satisfied. Then by \eqref{eq:symreducedbounds} we have $b^2\leq 2/3$, and, because $b$ is an integer, this is possible only if $b=0$. But then $2=ad-b^2=ad$ implies that $(a,d)=\pm(1,2)$ or $\pm(2,1)$. The only solutions satisfying conditions (i)-(iv)  are $(a,d)=\pm(1,2)$. Substituting these into \eqref{eq:x2pyzpcsol}, we obtain that the integer solutions to equation \eqref{eq:yzmx2m2} are
\begin{equation}\label{eq:x2pyzp2sol}
	(x,y,z) = \pm(uw+2vr, \, u^2+2r^2, \, 2v^2 + w^2), 
	\quad \text{where} \,\, uv-wr=1.
\end{equation}

\vspace{10pt}

The next equation we will consider is \eqref{eq:x2pyzpc} with $c=-2$, or
\begin{equation} \label{eq:yzmx2p2}
	yz-x^2=-2.
\end{equation}
Let $A$ given by \eqref{eq:symreduced} be a reduced matrix with determinant $ad-b^2=-2$. If $a=0$, then $b^2=2$, a contradiction, hence $a \neq 0$, and conditions (i)-(iv) in Definition \ref{def:symreduced} must be satisfied. Then by \eqref{eq:symreducedbounds} we have $b^2\leq 2/3$, and, because $b$ is an integer, $b=0$. But then $-2=ad-b^2=ad$ implies that $(a,d)=\pm(-1,2)$ or $\pm(2,-1)$. The only solutions satisfying conditions (i)-(iv)  are $(a,d)=\pm(-1,2)$. Substituting these into \eqref{eq:x2pyzpcsol}, we obtain that the integer solutions to equation \eqref{eq:yzmx2p2} are
\begin{equation}\label{eq:x2pyzm2sol}
	(x,y,z) = \pm(-uw+2vr, \, -u^2+2r^2, \, 2v^2 - w^2), 
	\quad \text{where} \,\, uv-wr=1.
\end{equation}

\vspace{10pt}

The next equation we will consider is \eqref{eq:x2pyzpc} with $c=3$, or
\begin{equation} \label{eq:yzmx2m3}
	yz-x^2=3.
\end{equation}
Let $A$ given by \eqref{eq:symreduced} be a reduced matrix with determinant $ad-b^2=3$. If $a=0$ then $b^2=-3$, a contradiction, hence $a \neq 0$, and conditions (i)-(iv) in Definition \ref{def:symreduced} must be satisfied. Then by \eqref{eq:symreducedbounds} we have $b^2\leq 1$, and, because $b$ is an integer, this is possible only if $b=\pm 1$ or $b=0$. If $b=\pm 1$, then $3=ad-1$, so $ad=4$, which implies that $(a,d)=\pm(1,4), \pm(4,1)$ or $\pm(2,2)$. If $b=-1$, these solutions do not satisfy conditions (i)-(iv). If $b=1$ then the only solutions satisfying conditions (i)-(iv) are $(a,d)=\pm(2,2)$. Substituting these into \eqref{eq:x2pyzpcsol}, we obtain
\begin{equation}\label{eq:yzmx2m3solb}
	\begin{aligned}
		(x,y,z) =&  (-2uw+uv+wr-2vr, \, -2u^2+2ru-2r^2, \, -2w^2 +2vw - 2v^2),  \quad \text{or} \\ & (2uw+uv+wr+2vr, \, 2u^2+2ru+2r^2, \, 2w^2 +2vw + 2v^2), \quad \text{where} \,\, uv-wr=1.
	\end{aligned}
\end{equation}
If $b=0$, then $3=ad-b^2=ad$, which implies that $(a,d)=\pm(1,3)$ or $\pm(3,1)$. The only solutions satisfying conditions (i)-(iv)  are $(a,d)=\pm(1,3)$. Substituting these into \eqref{eq:x2pyzpcsol}, we obtain
\begin{equation}\label{eq:yzmx2m3sola}
	(x,y,z) = \pm(uw+3vr, \, u^2+3r^2, \, 3v^2 + w^2), \,\, \text{where} \,\, uv-wr=1.
\end{equation}
Therefore, we have that the integer solutions to equation \eqref{eq:yzmx2m3} are either of the form \eqref{eq:yzmx2m3solb} or \eqref{eq:yzmx2m3sola}. 

\vspace{10pt}

The next equation we will consider is \eqref{eq:x2pyzpc} with $c=-3$, or
\begin{equation} \label{eq:yzmx2p3}
	yz-x^2=-3.
\end{equation}
Let $A$ given by \eqref{eq:symreduced} be a reduced matrix with determinant $ad-b^2=-3$. If $a=0$, then $b^2=3$, a contradiction, hence $a \neq 0$, and conditions (i)-(iv) in Definition \ref{def:symreduced} must be satisfied. Then by \eqref{eq:symreducedbounds} we have $b^2\leq 1$, and, because $b$ is an integer, this is possible only if $b=\pm 1$ or $b=0$. If $b=\pm 1$, then $-3=ad-b^2$ gives $ad=-2$ implies that $(a,d)=\pm(-1,2)$ or $\pm(2,-1)$, however, these solutions do not satisfy conditions (i)-(iv). If $b=0$, then $-3=ad-b^2=ad$ implies that $(a,d)=\pm(-1,3)$ or $\pm(3,-1)$. The only solutions satisfying conditions (i)-(iv)  are $(a,d)=\pm(-1,3)$. Substituting these into \eqref{eq:x2pyzpcsol}, we obtain that the integer solutions to equation \eqref{eq:yzmx2p3} are
\begin{equation}\label{eq:yzmx2p3sol}
	(x,y,z) = \pm(-uw+3vr, \, -u^2+3r^2, \, 3v^2 - w^2), 
	\quad \text{where} \,\, uv-wr=1.
\end{equation}

\vspace{10pt}

The next equation we will consider is \eqref{eq:x2pyzpc} with $c=4$, or
\begin{equation} \label{eq:yzmx2m4}
	yz-x^2=4.
\end{equation}
Let $A$ given by \eqref{eq:symreduced} be a reduced matrix with determinant $ad-b^2=4$. If $a=0$, then $b^2=-4$, a contradiction, hence $a \neq 0$, and, conditions (i)-(iv) in Definition \ref{def:symreduced} must be satisfied. Then by \eqref{eq:symreducedbounds} we have $b^2\leq \frac{4}{3}$, and, because $b$ is an integer, this is possible only if $b=\pm 1$, or $b=0$. If $b= \pm 1$, then $ad=5$, which implies that $(a,d)=\pm(1,5)$ or $\pm(5,1)$. However, these solutions do not satisfy conditions (i)-(iv). Hence, $b=0$, so $ad=4$, which implies that $(a,d)=\pm(1,4), \pm(4,1)$ or $\pm(2,2)$. The only solutions satisfying conditions (i)-(iv)  are $(a,d)=\pm(1,4)$ or $\pm(2,2)$. Substituting these into \eqref{eq:x2pyzpcsol}, we obtain that the integer solutions to \eqref{eq:yzmx2m4} are
\begin{equation}\label{eq:yzmx2m4sol} 
	\begin{aligned}
		(x,y,z)= & \pm(uw+4vr, \, u^2+4r^2, \, 4v^2 + w^2), \quad \text{or} \\ & \pm(2uw+2vr, \, 2u^2+2r^2, \, 2v^2 + 2w^2),  \text{where} \,\, uv-wr=1. 
	\end{aligned} 
\end{equation}

\vspace{10pt}

The next equation we will consider is \eqref{eq:x2pyzpc} with $c=-4$, or
\begin{equation} \label{eq:yzmx2p4}
	yz-x^2=-4. 
\end{equation}
Let $A$ given by \eqref{eq:symreduced} be a reduced matrix with determinant $ad-b^2=-4$. If $a=0$, then Definition \ref{def:symreduced} implies that $d=0$ and $b \geq 0$, hence $b=2$. Substituting this into \eqref{eq:x2pyzpcsol}, we obtain
\begin{equation} \label{eq:yzmx2p4sola}
	(x,y,z) = (2uv+2wr, \, 4ru, \, 4vw), 
	\quad \text{where} \,\, uv-wr=1.
\end{equation}
If $a \neq 0$, conditions (i)-(iv) in Definition \ref{def:symreduced} must be satisfied. Then by \eqref{eq:symreducedbounds} we have $b^2\leq \frac{4}{3}$, and, because $b$ is an integer, this is possible only if $b=\pm 1$ or $b=0$. If $b= \pm 1$, then $ad=-3$, which implies that $(a,d)=\pm(-1,3)$ or $\pm(3,-1)$. However, these solutions do not satisfy conditions (i)-(iv). Hence, $b=0$. Then $ad=-4$, which implies that $(a,d)=\pm(-1,4),\pm(4,-1)$ or $\pm(2,-2)$. The only solutions satisfying conditions (i)-(iv)  are $(a,d)=\pm(-1,4)$ and $(2,-2)$. Substituting these into \eqref{eq:x2pyzpcsol}, we obtain
\begin{equation}\label{eq:yzmx2p4solb}
	(x,y,z) = \pm(-uw+4vr, \, -u^2+4r^2, \, 4v^2 - w^2), 
	\quad \text{where} \,\, uv-wr=1,
\end{equation} or 
\begin{equation}\label{eq:yzmx2p4solc}
	(x,y,z) = (2uw-2vr, \, 2u^2-2r^2, \, -2v^2 +2w^2), 
	\quad \text{where} \,\, uv-wr=1.
\end{equation}
Therefore, the integer solutions to equation \eqref{eq:yzmx2p4} are of the form \eqref{eq:yzmx2p4sola}, \eqref{eq:yzmx2p4solb} or \eqref{eq:yzmx2p4solc}.

\vspace{10pt}

The next equation we will consider is \eqref{eq:x2pyzpc} with $c=5$, or
\begin{equation} \label{eq:yzmx2m5}
	yz-x^2=5.
\end{equation}
Let $A$ given by \eqref{eq:symreduced} be a reduced matrix with determinant $ad-b^2=5$. If $a=0$, then $b^2=-5$, a contradiction, hence $a \neq 0$, and, conditions (i)-(iv) in Definition \ref{def:symreduced} must be satisfied. Then by \eqref{eq:symreducedbounds} we have $b^2\leq \frac{5}{3}$, and, because $b$ is an integer, this is possible only if $b=\pm 1$, or $b=0$. If $b=\pm 1$, then $ad=6$, which implies that $(a,d)=\pm(1,6), \pm(6,1), \pm(2,3)$ or $\pm(3,2)$. If $b=-1$, then these solutions do not satisfy conditions (i)-(iv). If $b=1$, the only solutions which satisfy conditions (i)-(iv) are $(a,d)=\pm(2,3)$. Substituting these into \eqref{eq:x2pyzpcsol}, we obtain 
\begin{equation} \label{eq:yzmx2m5solb}
	(x,y,z)=(2uw+uv+wr+3vr, \, 2u^2+2ru+3r^2, \, 2w^2 +2vw + 3v^2) \quad \text{where} \,\, uv-wr=1,
\end{equation}
or 
\begin{equation}\label{eq:yzmx2m5solc}
	(x,y,z)=(-2uw+uv+wr-3vr, \, -2u^2+2ru-3r^2, \, -2w^2 +2vw -3v^2) \quad \text{where} \,\, uv-wr=1. 
\end{equation}
If $b=0$ then $ad=5$, which implies that $(a,d)=\pm(1,5)$ or $\pm(5,1)$. The only solutions satisfying conditions (i)-(iv)  are $(a,d)=\pm(1,5)$. Substituting these into \eqref{eq:x2pyzpcsol}, we obtain 
\begin{equation}\label{eq:yzmx2m5sola}
	(x,y,z) = \pm(uw+5vr, \, u^2+5r^2, \, 5v^2 + w^2), \quad \text{where} \,\, uv-wr=1.
\end{equation}
Therefore, we have that the integer solutions to equation \eqref{eq:yzmx2m5} are of the form \eqref{eq:yzmx2m5solb}, \eqref{eq:yzmx2m5solc} or \eqref{eq:yzmx2m5sola}.

\vspace{10pt}

The next equation we will consider is \eqref{eq:x2pyzpc} with $c=-5$, or
\begin{equation} \label{eq:yzmx2p5}
	yz-x^2=-5.
\end{equation}
Let $A$ given by \eqref{eq:symreduced} be a reduced matrix with determinant $ad-b^2=-5$. If $a=0$, then $b^2=5$, a contradiction, hence $a \neq 0$, and, conditions (i)-(iv) in Definition \ref{def:symreduced} must be satisfied. Then by \eqref{eq:symreducedbounds} we have $b^2\leq \frac{5}{3}$, and, because $b$ is an integer, this is possible only if $b=\pm 1$ or $b=0$. If $b=\pm 1$, then $ad=-4$, which implies that $(a,d)=\pm(-1,4), \pm(4,-1)$ or $\pm(2,-2)$. If $b=-1$, then these solutions do not satisfy conditions (i)-(iv). If $b=1$ then the only solutions satisfying conditions (i)-(iv)  are $(a,d)=\pm(2,-2)$. Substituting these into \eqref{eq:x2pyzpcsol}, we obtain
\begin{equation}\label{eq:yzmx2p5solb}
	(x,y,z) = (2uw+uv+wr-2vr, \, 2u^2+2ru-2r^2, \, -2v^2+2vw +2w^2), 
	\quad \text{where} \,\, uv-wr=1,
\end{equation}
or
\begin{equation}\label{eq:yzmx2p5solc}
	(x,y,z) = (-2uw+uv+wr+2vr, \, -2u^2+2ru+2r^2, \, 2v^2+2vw -2w^2), 
	\quad \text{where} \,\, uv-wr=1.
\end{equation}
If $b=0$ then $ad=-5$ which implies that $(a,d)=\pm(-1,5)$ or $\pm(5,-1)$. The only solutions satisfying conditions (i)-(iv)  are $(a,d)=\pm(-1,5)$. Substituting these into \eqref{eq:x2pyzpcsol}, we obtain
\begin{equation}\label{eq:yzmx2p5sola}
	(x,y,z) = \pm(-uw+5vr, \, -u^2+5r^2, \, 5v^2 - w^2), 
	\quad \text{where} \,\, uv-wr=1.
\end{equation}
Therefore, the integer solutions to equation \eqref{eq:yzmx2p5} are of the form \eqref{eq:yzmx2p5solb},  \eqref{eq:yzmx2p5solc} or \eqref{eq:yzmx2p5sola}.

Table \ref{table2.16} summarises the integer solutions to equation \eqref{eq:x2pyzpc} for $1\leq |c| \leq 5$.

	\begin{center}
\begin{tabular}{| c  | c |c|}
			\hline
			Equation  & Integer solutions $(x,y,z)$ \\
			\hline \hline
			$yz-x^2=1$ & $\pm(uw+vr, \,\, u^2+r^2, \,\, v^2 + w^2)$ \\\hline
			$yz-x^2=-1$ & $(uw-vr, \,\, u^2-r^2, \,\, w^2 - v^2)$, or \\
			& $(uv+wr,2ru,2vw)$ \\\hline
			$yz-x^2=2$ & $\pm(uw+2vr, \,\, u^2+2r^2, \,\, 2v^2 + w^2)$ \\\hline
			$yz-x^2=-2$ & $\pm(-uw+2vr, \, -u^2+2r^2, \, 2v^2 - w^2)$ \\\hline
			$yz-x^2=3$ & $\pm(uw+3vr, \,\, u^2+3r^2, \,\, 3v^2 + w^2)$, or \\
			& $(2uw+uv+wr+2vr, \,\, 2u^2+2ru+2r^2, \,\, 2w^2 +2vw + 2v^2)$, or \\
			& $ (-2uw+uv+wr-2vr, \,\, -2u^2+2ru-2r^2, \,\, -2w^2 +2vw - 2v^2)$ \\\hline
			$yz-x^2=-3$ & $\pm(-uw+3vr, \, -u^2+3r^2, \, 3v^2 - w^2)$ \\\hline
			$yz-x^2=4$ & $\pm(uw+4vr, \,\, u^2+4r^2, \,\, 4v^2 + w^2)$, or \\
			& $\pm(2uw+2vr, \,\, 2u^2+2r^2, \,\, 2v^2 + 2w^2)$ \\\hline
			$yz-x^2=-4$& $\pm(-uw+4vr, \, -u^2+4r^2, \, 4v^2 - w^2), \text{or,}$ \\
			&$(2uw-2vr, \, 2u^2-2r^2, \, -2v^2 +2w^2), \text{or,}$\\
			&$ (2uv+2wr, \, 4ru, \, 4vw) $ \\\hline
			$yz-x^2=5$ &$ \pm(uw+5vr, \,\, u^2+5r^2, \,\, 5v^2 + w^2)$, or \\
			& $(2uw+uv+wr+3vr, \,\, 2u^2+2ru+3r^2, \,\, 2w^2 +2vw + 3v^2)$, or \\
			& $ (-2uw+uv+wr-3vr, \,\, -2u^2+2ru-3r^2, \,\, -2w^2 +2vw -3v^2)$ \\\hline 
			$yz-x^2=-5$ &$ \pm(-uw+5vr, \, -u^2+5r^2, \, 5v^2 - w^2), \, \text{or}$ \\
			& $(2uw+uv+wr-2vr, \, 2u^2+2ru-2r^2, \, -2v^2+2vw +2w^2), \text{or,}$ \\
			& $ (-2uw+uv+wr+2vr, \, -2u^2+2ru+2r^2, \, 2v^2+2vw -2w^2)$ \\\hline
		\end{tabular}
		\captionof{table}{Integer solutions to equations of the form \eqref{eq:x2pyzpc} for $1\leq |c| \leq 5$. Assume that $u,v,w,r$ are integers satisfying $uv-wr=1$.\label{table2.16}}
	\end{center}

\subsection{Exercise 2.8}
\textbf{\emph{Use the solutions to the equations in Section \ref{ex:x2pyzpc} to prove the following proposition.
\begin{proposition}\label{prop:divsofx2pc}[Proposition 2.7 in the book]
	For every integer $x$, the following statements are true: 
	\begin{center}
		\begin{tabular}{ |c|c|c|c|c|c|c| } 
			\hline
			 &  & are of the form \\ 
			\hline
			All positive divisors of & $x^2+1$ & $u^2+r^2$ \\ 
			\hline
			All positive divisors of & $x^2+2$ & $u^2+2r^2$ \\ 
			\hline
			All divisors of & $x^2-2$ & $u^2-2r^2$ and $2u^2-r^2$ \\
			\hline
			All positive odd divisors of & $x^2+3$ & $u^2+3r^2$ \\  
			\hline
			All divisors of & $x^2-3$ & $u^2-3r^2$ or $3u^2-r^2$ \\ 
			\hline
			All positive divisors of & $x^2+4$ & $u^2+4r^2$ or $2u^2+2r^2$ \\ 
			\hline
			All positive divisors of & $x^2+5$ & $u^2+5r^2$ or $2u^2+2ur+3r^2$ \\ 
			\hline
			All odd divisors of & $x^2-5$ & $u^2-5r^2$ and $5u^2-r^2$ \\ 
			\hline
		\end{tabular}
	\end{center}
	for some integers $u,r$ such that $\text{gcd}(u,r)=1$.
\end{proposition}}}

	All integer solutions in Exercise \ref{ex:x2pyzpc} have the condition $uv-wr=1$. This implies that $\text{gcd}(u,r)=1$. 
	
	The description \eqref{eq:x2pyzp1sol} of all integer solutions to equation $yz=x^2+1$ implies that every divisor $y$ of $x^2+1$ is of the form $\pm(u^2+r^2)$ for some coprime integers $u$ and $r$. 
	In particular, all \emph{positive} divisors of $x^2+1$ can be represented as a sum of two squares of coprime integers. 
	
	The description \eqref{eq:x2pyzp2sol} of all integer solutions to equation $yz=x^2+2$ implies that every divisor $y$ of $x^2+2$ is of the form $\pm(u^2+2r^2)$ for some coprime integers $u$ and $r$.
	In particular, all \emph{positive} divisors of $x^2+2$ can be represented as the sum $u^2+2r^2$ for coprime $u,r$.
	
	The description \eqref{eq:x2pyzm2sol} of all integer solutions to equation $yz=x^2-2$ implies that every divisor $y$ of $x^2-2$ is of the form $\pm(u^2-2r^2)$ for some coprime integers $u$ and $r$. In other words, all divisors of $x^2-2$ can be represented as the difference $u^2-2r^2$ or $2r^2-u^2$ for coprime $u,r$. 
	However, the equalities
	$$
	u^2-2r^2=2(r+u)^2-(u+2r)^2 \quad \text{and} \quad 2r^2-u^2=(2r-u)^2-2(u-r)^2
	$$
	show that these representations are equivalent, so, all divisors of $x^2-2$ can be represented as the difference $u^2-2r^2$  for coprime $u,r$, and also as the difference $2r^2-u^2$ (or, equivalently $2u^2-r^2$) for coprime $u,r$.
	
	All integer solutions to equation $yz=x^2+3$ are given by \eqref{eq:yzmx2m3solb} or \eqref{eq:yzmx2m3sola}. In \eqref{eq:yzmx2m3solb}, $y$ and $z$ are both even, hence every \emph{odd} divisor of $x^2+3$ is of the form $\pm(u^2+3r^2)$ for some coprime integers $u$ and $r$. 
	In particular, all \emph{positive odd} divisors of $x^2+3$ can be represented as the sum $u^2+3r^2$ for coprime $u,r$.
	
	The description \eqref{eq:yzmx2p3sol} of all integer solutions to equation $yz=x^2-3$ implies that every divisor $y$ of $x^2-3$ is of the form $\pm(u^2-3r^2)$ for some coprime integers $u$ and $r$. 

	The description \eqref{eq:yzmx2m4sol} of all integer solutions to equation $yz=x^2+4$ implies that every divisor $y$ of $x^2+4$ is of the form $\pm(2u^2+2r^2)$ or $\pm(u^2+4r^2)$ for some coprime integers $u$ and $r$. 
	In particular, all \emph{positive} divisors of $x^2+4$ can be represented as either $u^2+4r^2$ or $2u^2+2r^2$ for coprime $u,r$.
	
	The description \eqref{eq:yzmx2m5solb}-\eqref{eq:yzmx2m5sola} of all integer solutions to the equation $yz=x^2+5$ implies that every \emph{positive} divisor $y$ of $x^2+5$ is of the form $u^2+5r^2$ or $2u^2+2ru+3r^2$ for some coprime integers $u$ and $r$. 
	
	Finally, all integer solutions to equation $yz=x^2-5$ are given by \eqref{eq:yzmx2p5solb}-\eqref{eq:yzmx2p5sola}.
	In \eqref{eq:yzmx2p5solb} and \eqref{eq:yzmx2p5solc}, $y$ and $z$ are both even, hence  
	every \emph{odd} divisor $y$ of $x^2-5$ is of the form $\pm(u^2-5r^2)$ for some coprime integers $u$ and $r$. 
	In other words, all odd divisors of $x^2-5$ can be represented as the difference $u^2-5r^2$ or $5r^2-u^2$ for coprime $u,r$. However, the equalities
	$$
	u^2-5r^2=5(u+2r)^2-(5r+2u)^2 \quad \text{and} \quad 5r^2-u^2=(5r-2u)^2-5(u-2r)^2
	$$
	show that these representations are equivalent, so, all odd divisors of $x^2-5$ can be represented as the difference $u^2-5r^2$ for coprime $u,r$, and also as the difference $5r^2-u^2$ (or, equivalently, $5u^2-r^2$) for coprime $u,r$.

\subsection{Exercise 2.10}\label{ex:xypztpc}
\textbf{\emph{Solve equation 
\begin{equation}\label{eq:xypztpc}
xy-zt=c
\end{equation}
where $x,y,z,t$ are variables and $c$ is a constant for $c=3,4$ and $5$.}}

To solve these equations, we will use the following proposition.

\begin{proposition}\label{prop:sl2actions}[Proposition 2.9 in the book]
Let $c\neq 0$ and let $(x,y,z,t)$ be any solution to \eqref{eq:xypztpc}. Then there exists a matrix $A \in SL_2(\mathbb Z)$ such that 
$$
\begin{bmatrix}
x & z \\
t & y
\end{bmatrix}
=
A \cdot
\begin{bmatrix}
a & b \\
0 & d
\end{bmatrix},
$$
for some integers $a,b,d$ such that $0\leq b < d$ and $ad=c$.
\end{proposition}

It is easy to see that for every $c\neq 0$ there are only finitely many triples $a,b,d$ satisfying the conditions of Proposition \ref{prop:sl2actions}. Denote $(a_j,b_j,d_j)$, $j=1,\dots,k$, as the full list of these triples. Then Proposition \ref{prop:sl2actions} implies that for any solution $(x,y,z,t)$ to \eqref{eq:xypztpc} we have
$$
\begin{bmatrix}
x & z \\
t & y
\end{bmatrix}
=
\begin{bmatrix}
u & w \\
r & v
\end{bmatrix}
\cdot 
\begin{bmatrix}
a_j & b_j \\
0 & d_j
\end{bmatrix}
=
\begin{bmatrix}
a_j u & b_j u + d_j w \\
a_j r & b_j r + d_j v
\end{bmatrix},
$$
or in component form
 \begin{equation}\label{eq:xypztpcsol}
(x,y,z,t) = (a_j u, \,\, b_j r + d_j v, \,\, b_j u + d_j w, \,\, a_j r) \quad \text{for some} \quad j \in \{1,\dots,k\}, 
\end{equation}
for some integers $u,v,w,r$ such that $uv-wr=1$. 

A similar equation, \eqref{eq:xypztpc} with $c=2$, or
\begin{equation}\label{eq:xypztp2}
xy-zt=2,
\end{equation}
is solved in the book, and its integer solutions are
\begin{equation}\label{eq:xypztp2sol}
(x,y,z,t)=(2u,v,w,2r) \quad \text{or} \quad (u,2v,2w,r) \quad \text{or} \quad (u,2v+r,2w+u,r) \quad \text{where} \,\, uv-wr=1.
\end{equation} 

\vspace{10pt}

The first equation we will solve is
\begin{equation}\label{ex2.8_3}
xy-zt=3.
\end{equation}
This is equation \eqref{eq:xypztpc} with $c=3$. From Proposition \ref{prop:sl2actions} we have that $c=ad=3$, so we have the following values for $(a,d)=(1,3),$ and $(3,1)$. Then, possible triples $(a,b,d)$ in Proposition \ref{prop:sl2actions} are $(1,0,3), (1,1,3),$ $(1,2,3),$ and $(3,0,1)$. By \eqref{eq:xypztpcsol}, we obtain that all integer solutions to equation \eqref{ex2.8_3} are
$$
(x,y,z,t) = (u,3v,3w,r), (u, r+3v,u+3w,r), (u,2r+3v, 2u+3w,r), (3u,v,w,3r), \text{where } uv-wr=1.
$$

\vspace{10pt}

The next equation we will solve is
\begin{equation}\label{ex2.8_4}
xy-zt=4.
\end{equation}
This is equation \eqref{eq:xypztpc} with $c=4$. From Proposition \ref{prop:sl2actions} we have that $c=ad=4$, so we have the following values for $(a,d)=(1,4), (4,1) $ and $(2,2)$. Then, possible triples $(a,b,d)$ in Proposition \ref{prop:sl2actions} are $(1,0,4), (1,1,4),$ $(1,2,4),(1,3,4), (2,0,2), (2,1,2)$ and $(4,0,1)$. By \eqref{eq:xypztpcsol}, we obtain that all integer solutions to equation \eqref{ex2.8_4} are
$$
\begin{aligned}
(x,y,z,t) = (u,4v,4w,r), (u, r+4v,u+4w,r), (u,2r+4v, 2u+4w,r), (u, 3r+4v, 3u+4w,r), \\(2u,2v,2w,2r), (2u, r+2v, u+2w,2r),(4u,v,w,4r), \text{where } uv-wr=1.
\end{aligned}
$$

\vspace{10pt}

The next equation we will solve is
\begin{equation}\label{ex2.8_5}
xy-zt=5.
\end{equation}
This is equation \eqref{eq:xypztpc} with $c=5$. From Proposition \ref{prop:sl2actions} we have that $c=ad=5$, so we have the following values for $(a,d)=(1,5),$ and $(5,1)$. Then, possible triples $(a,b,d)$ in Proposition \ref{prop:sl2actions} are $(1,0,5), (1,1,5),$ $(1,2,5),(1,3,5), (1,4,5), $ and $(5,0,1)$. By \eqref{eq:xypztpcsol}, we obtain that all integer solutions to equation \eqref{ex2.8_5} are
$$
\begin{aligned}
(x,y,z,t) = (u,5v,5w,r), (u, r+5v,u+5w,r), (u,2r+5v, 2u+5w,r), (u, 3r+5v, 3u+5w,r),\\ (u, 4r+5v,4u+5w,r), (5u,v,w,5r), \text{where } uv-wr=1.
\end{aligned}
$$

Table \ref{table2.8} summarises the integer solutions to the equations solved in this exercise.
\begin{center}
\begin{tabular}{| c  | c |c|}
\hline
Equation & Integer solutions $(x,y,z,t)$\\
\hline \hline
$xy-zt=2$ & $(2u,v,w,2r), (u,2v,2w,r) \quad \text{or} \quad (u,2v+r,2w+u,r)$ \\\hline
$xy-zt=3$ & $(u,3v,3w,r), (u, r+3v,u+3w,r), (u,2r+3v, 2u+3w,r), \text{or, } (3u,v,w,3r)$ \\\hline
$xy-zt=4$ & $(u,4v,4w,r), (u, r+4v,u+4w,r), (u,2r+4v, 2u+4w,r), (u, 3r+4v, 3u+4w,r),$ \\
& $ (2u,2v,2w,2r), (2u, r+2v, u+2w,2r), \text{or}, \, (4u,v,w,4r)$ \\\hline
$xy-zt=5$ & $ (u,5v,5w,r), (u, r+5v,u+5w,r), (u,2r+5v, 2u+5w,r), (u, 3r+5v, 3u+5w,r),$ \\
&$(u, 4r+5v,4u+5w,r),  \text{or}, \, (5u,v,w,5r)$ \\\hline
\end{tabular}
\captionof{table}{Integer solutions to equations of the form \eqref{eq:xypztpc} with $2 \leq c \leq 5$. Assume that $u,v,w,r$ are integers satisfying $uv-wr=1$.\label{table2.8}}
\end{center}

\subsection{Exercise 2.13}\label{ex:linsub}
\textbf{\emph{For each of the equations listed in Table \ref{tab:H12linsub}, find linear substitutions in the forms 
$$
x'_i = x_i + b
$$
or
$$
x'_i = x_i + b x_j, \quad x'_j = x_j
$$
for some $i,j \in \{1,\dots,n\}$ and some integer $b$, that reduce it to a simpler equation with smaller $H$. Use this to solve the equation.}}
\begin{center}
\begin{tabular}{ |c|c|c|c|c|c| } 
 \hline
 $H$ & Equation &  $H$ & Equation & $H$ & Equation \\ 
 \hline\hline
 $10$ & $x^2+yz+y=0$ & $12$ & $xy+zt+2x=0$ & $12$ & $x^2+yz+y+z=0$ \\ 
 \hline
 $10$ & $xy+zt+x=0$ & $12$ & $xy+xz+y-z=0$ & $12$ & $x^2+yz+2y=0$ \\ 
 \hline
 $10$ & $xy+xz+y=0$ & $12$ & $xy+xz+y+2=0$ & $12$ & $x^2+yz+x+y=0$ \\ 
 \hline
 $11$ & $xy+zt+x+1=0$ & $12$ & $xy+xz+2y=0$ & $12$ & $x^2+yz+2x=0$ \\ 
 \hline
 $11$ & $xy+xz+y+1=0$ & $12$ & $xy+xz+x+y=0$ & $12$ & $x^2+yz+yt=0$ \\ 
 \hline
 $11$ & $x^2+yz+y-1=0$ & $12$ & $xy+xz+yt=0$ & $12$ & $x^2+xy+yz=0$ \\ 
 \hline
 $11$ & $x^2+yz+y+1=0$ & $12$ & $x^2+xy-z^2=0$ & $12$ & $x^2+xy+zt=0$ \\ 
 \hline
 $12$ & $xy+zt+x+2=0$ & $12$ & $x^2+yz+y-z=0$ & $12$ & $x^2+xy+z^2=0$ \\ 
 \hline
 $12$ & $xy+zt+x+y=0$ & $12$ & $x^2+yz+y-2=0$ & $12$ & $xy+xz+yz=0$ \\ 
 \hline
 $12$ & $xy+zt+x+z=0$ & $12$ & $x^2+yz+y+2=0$ & $12$ & $x_1x_2+x_1x_3+x_4x_5=0$ \\ 
 \hline
\end{tabular}
\captionof{table}{\label{tab:H12linsub} Equations of size $H \leq 12$ solvable by invertible linear substitutions.}
\end{center} 

In this section, we will solve equations by finding linear substitutions in order to reduce the equation to a simpler one (one with a smaller value of H).
Linear substitutions that we may use will be of the following forms:
\begin{itemize}
\item $x'_i = x_i +b$
\item $x'_i = x_i+bx_j+c$
\end{itemize}
where $b,c$ are integers, and $x_i,x_j, x'_i$ are variables.

Equation
$$
x^2+yz+y=0
$$
is equivalent to $x^2-yz-y=0$ which is solved in Section 2.2.1 of the book, and its integer solutions are
$$
(x,y,z)=(u v w,-uv^2,uw^2-1) \quad \text{for some integers} \quad u,v,w.
$$

Equation
$$
xy+zt+x=0
$$
is equivalent to $xy-zt+x=0$ which is solved in Section 2.2.1 of the book, and its integer solutions are
$$
(x,y,z,t)=(uv,wr-1,-uw,vr) \quad \text{for some integers} \quad u,v,w,r.
$$

Equation
$$
xy+xz+y=0
$$
is solved in Section 2.2.1 of the book, and its integer solutions are
$$
(x,y,z)=(u,-uv,uv+v) \quad \text{for some integers} \quad u,v.
$$

Equation
$$
x^2+yz+y-1=0
$$
is equivalent to $x^2-yz-y-1=0$ which is solved in Section 2.2.1 of the book, and its integer solutions are
$$
(x,y,z)=(uw - vr, r^2 - u^2, w^2 - v^2 - 1) \quad \text{or} \quad (uv + wr, -2ru, 2vw - 1), \quad \text{where} \quad uv - wr = 1.
$$

Equation
$$
xy+yz+xz=0
$$
is solved in Section 2.2.1 of the book, and its integer solutions are
$$
(x,y,z,t)=(uvw,uv^2-uvw,uw^2-uvw) \quad \text{for some integers} \quad u,v,w.
$$

The first equation we will consider is
$$
xy+zt+x+1=0.
$$
By making the substitutions $y'=y+1$ and $z'=-z$, this equation is reduced to $xy'-z't=-1$, which is equivalent to $z't-xy'=1$, which up to the names of variables is \eqref{eq:xypztp1}, whose solution is a polynomial family with 46 parameters. Let us write $(x,y',z',t) = (w,r,u,v)$ with $uv-wr=1$, or in the original variables, 
$$
(x,y,z,t)=(w,r-1,-u,v) \quad \text{with} \quad uv-wr=1.
$$

\vspace{10pt}

The next equation we will consider is
\begin{equation}\label{2.11_b}
xy+xz+y+1=0.
\end{equation}
By making the substitution $z'=y+z$, the equation is reduced to $xz'+y+1=0$ which is of the form \eqref{eq:sepvar} and its integer solution is $(x,y,z')=(u,-1-uv,v)$, where $u,v$ are integers. Hence the integer solution to equation \eqref{2.11_b} is 
$$
(x,y,z)=(u,-1-uv,1+v+uv) \quad \text{for some integers} \quad u,v. 
$$

\vspace{10pt}

The next equation we will consider is
\begin{equation} \label{2.11_c}
x^2+yz+y+1=0.
\end{equation}
By making the substitutions $z'=z+1$ and $y'=-y$, the equation is reduced to $y'z'-x^2=1$, which up to the names of variables is \eqref{eq:x2pyzp1}, and whose solution is $(x,y',z')=\pm(uw+vr,u^2+r^2,v^2+w^2)$ with $uv-wr=1$. Hence the integer solutions to equation \eqref{2.11_c} are
$$
(x,y,z)=(uw+vr,-u^2-r^2,v^2+w^2-1) \,\, \text{or} \,\, (-uw-vr,u^2+r^2,-v^2-w^2-1) \quad \text{with} \quad uv-wr=1.
$$

\vspace{10pt}

The next equation we will consider is
\begin{equation} \label{2.11d}
xy+zt+x+2=0.
\end{equation}
By making the substitutions $y'=y+1$ and $x'=-x$, the equation is reduced to $x'y'-zt=2$, which up to the names of variables is equation \eqref{eq:xypztp2}, whose solution is $(x',y',z,t)=(2u,v,w,2r)$, or $(u,2v,2w,r)$, or $(u,2v+r,2w+u,r)$, with $uv-wr=1$. Hence, the integer solutions to equation \eqref{2.11d} are
$$ 
\begin{aligned}
(x,y,z,t)=(-2u,v-1,w,2r), \,\, \text{or} \,\, (-u,2v-1,2w,r), \,\, \text{or} \,\, (-u,2v+r-1,2w+u,r), \\ \text{with} \quad uv-wr=1.
\end{aligned}
$$

\vspace{10pt}

The next equation we will consider is
\begin{equation}\label{eq:xypztpxpy}
xy+zt+x+y=0.
\end{equation}
By making the substitutions $y'=y+1$, $x'=x+1$ and $z'=-z$ this equation is reduced to $x'y'-z't=1$, which up to the names of variables is \eqref{eq:xypztp1}, whose solution is $(x',y',z',t) = (u,v,w,r)$ with $uv-wr=1$, or in the original variables we obtain that all integer solutions to \eqref{eq:xypztpxpy} are 
$$
(x,y,z,t)=(u-1,v-1,-w,r) \quad \text{with} \quad uv-wr=1.
$$

\vspace{10pt}

The next equation we will consider is
\begin{equation}\label{eq:xypztpxpz}
xy+zt+x+z=0.
\end{equation}
By making the substitutions $y'=y+1$, $t'=t+1$ and $z'=-z$, the equation is reduced to $xy'-z't'=0$, which up to the names of variables is \eqref{eq:xymzt}, whose solution is $(x,y',z',t')=(uv,wr,uw,vr)$ or in the original variables we obtain that all integer solutions to equation \eqref{eq:xypztpxpz} are 
$$
(x,y,z,t)=(uv,wr-1,-uw,vr-1) \quad \text{for some integers} \quad u,v,r. 
$$

\vspace{10pt}

The next equation we will consider is
\begin{equation}\label{eq:xypztp2x}
xy+zt+2x=0.
\end{equation}
By making the substitutions $y'=y+2$ and $z'=-z$, the equation is reduced to $xy'-z't=0$, which up to the names of variables is \eqref{eq:xymzt}, whose solution is $(x,y',z',t)=(uv,wr,uw,vr)$ for some integers $u,v,w,r$, or in the original variables we obtain that the integer solutions to equation \eqref{eq:xypztp2x} are 
$$
(x,y,z,t)=(uv,wr-2,-uw,vr) \quad \text{for some integers} \quad u,v,w,r. 
$$

\vspace{10pt}

The next equation we will consider is 
\begin{equation}\label{2.11_a}
xy+xz+y-z=0.
\end{equation}
By making the substitutions $y'=y+z$ and $x'=x+1$, the equation is reduced to $x'y'=2z$ which is equivalent to an equation in Table \ref{tab:H11axPlistsol} whose solution is $(x',y',z)=(2u,v,uv)$ or $(u,2v,uv)$ for some integers $u,v$. Hence, in the original variables, we have that the integer solutions to equation \eqref{2.11_a} are
$$
(x,y,z)=(u-1,2v-uv,uv), \,\, \text{or} \,\, (2u-1,v-uv,uv) \quad \text{for some integers} \quad u,v. 
$$

\vspace{10pt}

The next equation we will consider is
\begin{equation}\label{2.11_ai}
xy+xz+y+2=0.
\end{equation}
By making the substitutions $y'=y+z$ and $x'=x+1$, the equation is reduced to $y'x'+2=z$ which is of the form \eqref{eq:sepvar} and its integer solutions are $(x',y',z)=(u,v,uv+2)$ for some integers $u,v$. Hence, in the original variables we obtain that the integer solutions to equation \eqref{2.11_ai} are
$$
(x,y,z)=(u-1,v-uv-2,uv+2) \quad \text{for some integers} \quad u,v.
$$

\vspace{10pt}

The next equation we will consider is
\begin{equation}\label{2.11_aii}
xy+xz+2y=0.
\end{equation}
By making the substitutions $y'=y+z$ and $x'=x+2$, the equation is reduced to $x'y'=2z$ which is equivalent to an equation in Table \ref{tab:H11axPlistsol} whose solution is $(x'y',z)=(2u,v,uv)$ or $(u,2v,uv)$ for some integers $u,v$. Hence, in the original variables, we obtain that all integer solutions to equation \eqref{2.11_aii} are
$$
(x,y,z)=(2u-2,v-uv,uv), \,\, \text{or} \,\, (u-2,2v-uv,uv) \quad \text{for some integers} \quad u,v. 
$$ 
 
 \vspace{10pt}
 
 The next equation we will consider is
\begin{equation}\label{2.11_aiii}
xy+xz+x+y=0.
\end{equation}
By making the substitution $z'=y+z+1$, the equation is reduced to $xz'+y=0$, which up to the names of variables is an equation in Table \ref{tab:Ex1.5} whose solution is $(x,y,z')=(u,-uv,v)$ for some integers $u,v$. Hence, in the original variables, we obtain that all integer solutions to equation \eqref{2.11_aiii} are
$$
(x,y,z)=(u,-uv,v+uv-1) \quad \text{for some integers} \quad u,v. 
$$

\vspace{10pt}
 
 The next equation we will consider is
\begin{equation}\label{2.11_aiv}
xy+xz+yt=0.
\end{equation}
By making the substitutions $z'=y+z$ and $y'=-y$, the equation is reduced to $xz'-y't=0$, which up to the names of variables is equation \eqref{eq:xymzt}, whose solution is $(x,y',z',t)=(uv,uw,wr,vr)$ for some integers $u,v,w,r$. Hence, in the original variables, we obtain that all integer solutions to equation \eqref{2.11_aiv} are
$$
(x,y,z,t)=(uv,-uw,wr+uw,vr) \quad \text{for some integers} \quad u,v,w,r. 
$$ 
 
 \vspace{10pt}
 
  The next equation we will consider is
\begin{equation}\label{2.11_av}
x^2+xy-z^2=0.
\end{equation}
By making the substitution $y'=x+y$, the equation is reduced to $z^2-xy'=0$, which up to the names of variables is equation \eqref{eq:x2myz}, whose solution is $(x,y',z)=(uw^2,uv^2,uvw)$ for some integers $u,v,w$. Hence, in the original variables, we obtain that all integer solutions to equation \eqref{2.11_av} are
$$
(x,y,z)=(uw^2,uv^2-uw^2,uvw) \quad \text{for some integers} \quad u,v,w. 
$$ 

\vspace{10pt}

The next equation we will consider is
\begin{equation}\label{2.11_avi}
x^2+yz+y-z=0,
\end{equation}
By making the substitutions $t=-z-1$ and $y'=y-1$, the equation is reduced to $y't-x^2=1$, which up to the names of variables is \eqref{eq:x2pyzp1}, whose solution is $(x,y',t)=\pm(uw+vr,u^2+r^2,v^2+w^2) \text{ with } uv-wr=1$. Hence, in the original variables, we obtain that all integer solutions to equation \eqref{2.11_avi} are
$$
\begin{aligned} 
(x,y,z)=(uw+vr,u^2+r^2+1,-v^2-w^2-1) \,\, \text{or} \,\, (-uw-vr, -u^2-r^2+1,v^2+w^2-1)\\ 
\text{with} \quad uv-wr=1. 
\end{aligned}
$$ 

\vspace{10pt}

The next equation we will consider is
\begin{equation}\label{2.11_avii}
x^2+yz+y-2=0.
\end{equation}
By making the substitutions $z'=z+1$ and $y'=-y$, the equation is reduced to $y'z'-x^2=-2$ which up to the names of variables is equation \eqref{eq:yzmx2p2}, whose solution is $(x,y',z)=\pm(-uw+2vr,-u^2+2r^2,2v^2-w^2) \text{ with } uv-wr=1$. Hence, in the original variables, we obtain that all integer solutions to equation \eqref{2.11_avii} are 
$$
\begin{aligned} (x,y,z)=(-uw+2vr,u^2-2r^2,2v^2-w^2-1) \,\, \text{or} \,\, (uw-2vr,-u^2+2r^2,-2v^2+w^2-1) \\
 \text{with} \quad uv-wr=1.  
 \end{aligned}
 $$ 

\vspace{10pt}

The next equation we will consider is
\begin{equation}\label{2.11_aviii}
x^2+yz+y+2=0.
\end{equation}
By making the substitutions $z'=z+1$ and $y'=-y$, the equation is reduced to $y'z'-x^2=2$ which up to the names of variables is equation \eqref{eq:yzmx2m2}, whose solution is $(x,y',t)=\pm(uw+2vr,u^2+2r^2,2v^2+w^2) \text{ with } uv-wr=1$. Hence, in the original variables, we obtain that all integer solutions to equation \eqref{2.11_aviii} are
$$
\begin{aligned} 
(x,y,z)= (uw+2vr,-u^2-2r^2,2v^2+w^2-1) \,\, \text{or} \,\, (-uw-2vr,u^2+2r^2,-2v^2-w^2-1)\\ 
\text{with} \quad uv-wr=1. 
\end{aligned}
$$ 

\vspace{10pt}

The next equation we will consider is
\begin{equation}\label{2.11_aviv}
x^2+yz+y+z=0.
\end{equation}
By making the substitutions $z'=z+1$, $y'=y+1$ and $t=-z'$, the equation is reduced to $ty'-x^2=-1$, which up to the names of variables is equation \eqref{eq:x2pyzm1}, whose solution is $(x,y',t)=(uw-vr,u^2-r^2,w^2-v^2)$ or $(uv+wr,2ru,2vw)$ with  $uv-wr=1$. Hence, in the original variables, we obtain that the integer solutions to equation \eqref{2.11_aviv} are
$$
\begin{aligned} 
(x,y,z)= (uw-vr,u^2-r^2-1,v^2-w^2-1) \,\, \text{or} \,\, (uw+vr,2ru-1,-2vw-1) \,\, \text{with} \,\, uv-wr=1. 
\end{aligned}
$$

\vspace{10pt}

The next equation we will consider is
\begin{equation}\label{2.11_aX}
x^2+yz+2y=0.
\end{equation}
By making the substitutions $z'=z+2$ and $y'=-y$, the equation is reduced to $x^2-y'z'=0$, which up to the names of variables is equation \eqref{eq:x2myz}, whose solution is $(x,y',z)=(uvw, uv^2,uw^2)$ where $u,v,w$ are integers. Hence, in the original variables, we obtain that the integer solutions to equation \eqref{2.11_aX} are
$$
(x,y,z)=(uvw,-uv^2,uw^2-2) \quad \text{for some integers} \quad u,v,w. 
$$   

\vspace{10pt}

The next equation we will consider is
\begin{equation}\label{2.11_aXi}
x^2+yz+x+y=0.
\end{equation}
By making the substitutions $z'=z+1$ and $y'=-y$, the equation is reduced to $x^2+x-y'z'=0$, which is equivalent to equation \eqref{eq:x2mxmyz}, up to the names of variables. By \eqref{eq:x2mxmyzsol} its integer solutions are $(x,y',z')=(uv,uw,vr)$ for some integers $u,v,w,r$ such that $wr-uv=1$. Therefore, the integer solutions to equation \eqref{2.11_aXi} are
\begin{equation}\label{2.11_aXisol}
(x,y,z)=(uv,-uw,vr-1) \quad \text{with} \quad wr-uv=1.
\end{equation}
For completeness, we will now check that this solution is correct. If we substitute the solution into equation \eqref{2.11_aXi}, we have 
$$
x^2+yz+x+y=(uv)^2+(-uw)(vr-1)+uv+(-uw)=u^2v^2+uv-uvrw,
$$
from this it is not immediately clear that this is equal to 0. However, due to the constraint on the variables $wr-uv=1$, we may substitute $wr=1+uv$, and we obtain
$$
u^2v^2+uv-uvrw=u^2v^2+uv-uv(1+uv)=u^2v^2+uv-uv-u^2v^2=0.
$$ 
Hence, \eqref{2.11_aXisol} is a solution to equation \eqref{2.11_aXi}. 

\vspace{10pt}

The next equation we will consider is
\begin{equation}\label{2.11_aXii}
x^2+yz+2x=0.
\end{equation}
By making the substitutions $x'=x+1$ and $y'=-y$, the equation is reduced to $y'z-(x')^2=-1$, which up to the names of variables is \eqref{eq:x2pyzm1}, whose solution is $(x',y',z)=(uw-vr, u^2-r^2,w^2-v^2)$ or $(uv+wr,2ru,2vw)$ with $uv-wr=1$, see \eqref{eq:x2pyzm1sol}. Hence, in the original variables, we obtain that the integer solutions to equation \eqref{2.11_aXii} are
$$
(x,y,z)=(uw-vr-1, r^2-u^2,w^2-v^2) \,\, \text{or} \,\, (uv+wr-1,-2ru,2vw) \quad \text{with} \quad uv-wr=1.
$$  

\vspace{10pt}

The next equation we will consider is
\begin{equation}\label{2.11_aXiii}
x^2+yz+yt=0.
\end{equation}
By making the substitutions $z'=z+t$ and $y'=-y$, the equation is reduced to $x^2-y'z'=0$, which up to the names of variables is \eqref{eq:x2myz}, whose solution is $(x,y',z)=(uvw, uv^2,uw^2)$ for some integers $u,v,w$. Hence, in the original variables, we obtain that the integer solutions to equation \eqref{2.11_aXiii} are
$$
(x,y,z,t)=(uvw,-uv^2,uw^2-r,r) \quad \text{for some integers} \quad u,v,w,r. 
$$   
 
 \vspace{10pt}
 
The next equation we will consider is
\begin{equation}\label{2.11_aXiv}
x^2+xy+yz=0.
\end{equation}
By making the substitutions $z'=x+z$ and $y'=-y$, the equation is reduced to $x^2-y'z'=0$, which up to the names of variables is \eqref{eq:x2myz}, whose solution is $(x,y',z)=(uvw, uv^2,uw^2)$ for some integers $u,v,w$. Hence, in the original variables, we obtain that the integer solutions to equation \eqref{2.11_aXiv} are
$$
(x,y,z)=(uvw,-uv^2,uw^2-uvw) \quad \text{for some integers} \quad u,v,w. 
$$  

\vspace{10pt}

The next equation we will consider is
\begin{equation}\label{2.11_aXv}
x^2+xy+zt=0.
\end{equation}
By making the substitutions $y'=x+y$ and $z'=-z$, the equation is reduced to $xy'-z't=0$, which up to the names of variables is \eqref{eq:xymzt}, whose solution is $(x,y',z',t)=(uv,wr,uw,vr)$ for some integers $u,v,w,r$. Hence, in the original variables, we obtain that the integer solutions to equation \eqref{2.11_aXv} are
$$
(x,y,z,t)=(uv,wr-uv,-uw,vr) \quad \text{for some integers} \quad u,v,w,r. 
$$ 

\vspace{10pt}

The next equation we will consider is
\begin{equation}\label{2.11_aXvi}
x^2+xy+z^2=0.
\end{equation}
By making the substitutions $y'=x+y$ and $x'=-x$, the equation is reduced to $z^2-y'x'=0$, which up to the names of variables is \eqref{eq:x2myz}, whose solution is $(x',y',z)=(uw^2, uv^2,uvw)$ for some integers $u,v,w$. Hence, in the original variables, we obtain that the integer solutions to equation \eqref{2.11_aXvi} are
$$
(x,y,z)=(-uw^2,uv^2+uw^2,uvw) \quad \text{for some integers} \quad u,v,w. 
$$  

\vspace{10pt}

The next equation we will consider is
\begin{equation}\label{2.11_aXvii}
x_1x_2+x_1x_3+x_4x_5=0.
\end{equation}
By making the substitutions $x'=x_2+x_3$ and $x_4'=-x_4$, the equation is reduced to $x'x_1-x_4'x_5=0$, which up to the names of variables is \eqref{eq:xymzt}, whose solution is $(x_1,x',x_4',x_5)=(uv,wr,uw,vr)$ for some integers $u,v,w,r$. Hence, in the original variables, we obtain that the integer solutions to equation \eqref{2.11_aXvii} are
$$
(x_1,x_2,x_3,x_4,x_5)=(uv,wr-s,s,-uw,vr) \quad \text{for some integers} \quad u,v,w,r,s. 
$$ 

Table \ref{table2.11} summarises the integer solutions to the equations solved in this exercise.

\begin{center}
\begin{tabular}{| c  | c |c|}
\hline
Equation & Integer solutions $(x,y,z)$ or $(x,y,z,t)$ or $(x_1,\dots,x_n)$ \\
\hline \hline
$x^2+yz+y=0$&$(uvw,-uv^2,uw^2-1)$ \\\hline
$xy+zt+x=0$&$(uv,wr-1,-uw,vr)$ \\\hline
$xy+xz+y=0$&$(u,-uv,uv+v)$ \\\hline
$xy+zt+x+1=0$&$(w,r-1,-u,v) \, \text{with} \, uv-wr=1.$ \\\hline
$xy+xz+y+1=0$&$(u,-1-uv,1+v+uv)$ \\\hline
$x^2+yz+y-1=0$&$(uw-vr,r^2-u^2,w^2-v^2-1) \text{ or }$\\
&$(uv+wr,-2ru,2vw-1)\text{ with } uv-wr=1$ \\
\hline
\end{tabular}
\captionof{table}{Integer solutions to the equations listed in Table \ref{tab:H12linsub}. Assume that $u,v,w,r,s$ are integers. \label{table2.11}}
\end{center}

\subsection{Exercise 2.28}
\textbf{\emph{Extend the proofs of Propositions \ref{prop:GausEuc} and \ref{prop:QuadIntEuc357} to prove Theorem \ref{th:Euclist} for $d=-11,-2,2,3$, and $13$. Conclude that the statement of Proposition \ref{prop:x2pyzpxp1} remains valid for $k=\pm 3$.}}

We will start bt reminding the reader the statements of Theorem \ref{th:Euclist} and Propositions \ref{prop:GausEuc} and \ref{prop:QuadIntEuc357}.
 
\begin{theorem}\label{th:Euclist}[Theorem 2.26 in the book]
${\mathbb Z}[w_d]$ is an Euclidean domain if and only if
$$
d=-11,-7,-3,-2,-1, 2, 3, 5, 6, 7, 11, 13, 17, 19, 21, 29, 33, 37, 41, 57, \,\text{or}\,\, 73.
$$
\end{theorem}

\begin{proposition}\label{prop:GausEuc}[Proposition 2.15 in the book]
The Gaussian integers form a Euclidean domain with $f=N$, where $N$ is given by $N(x)=N(a+bi)=a^2+b^2$. 
\end{proposition}

\begin{proposition}\label{prop:QuadIntEuc357}[Proposition 2.25 in the book]
For $d=5,-3,-7$, the set ${\mathbb Z}[w_d]$ with the usual addition and multiplication forms a Euclidean domain with Euclidean valuation $f(\alpha)=|N(\alpha)|$. 
\end{proposition}

\begin{proposition}\label{prop:x2pyzpxp1}[Proposition 2.27 in the book]
	Let $k=-1,1$ or $2$. Then integers $x,y,z$ satisfy equation $x^2+x+k=yz$
	if and only if 
	$$
		(x,y,z)=(kur+vw+wr,\, \pm (ku^2 + uw + w^2),\, \pm (v^2 + vr + kr^2)) \quad \text{where} \quad uv-wr=1,
	$$
	where the choice of signs for $y$ and $z$ is either ``$+$'' for both variables or ``$-$'' for both variables.
\end{proposition}  

We next present the proofs of Propositions \ref{prop:GausEuc} and \ref{prop:QuadIntEuc357} given in the book.

\begin{proof}[Proof of Proposition \ref{prop:GausEuc}]
	Let $a$ and $b\neq 0$ be Gaussian integers. Then $N(b)$ is a positive integer, hence $N(b)\geq 1$. Thus, $N(ab)=N(a)N(b)\geq N(a)$, and property (i) follows. To prove (ii), write the complex number $a/b$ in the form $u+iv$ for some rational numbers $u,v$. Let us choose integers $x$ and $y$ that are the nearest to $u$ and $v$, respectively, that is, $|u-x|\leq 1/2$ and $|v-y|\leq 1/2$. Then $q=x+iy$ and $r=a-bq$ are Gaussian integers. Also,
	$$
	N(r)=N(a-bq)=N(b)N(a/b-q)=N(b)N((u-x)+(v-y)i)
	$$
	$$
	=N(b)((u-x)^2+(v-y)^2) \leq N(b)((1/2)^2+(1/2)^2) = N(b)/2 < N(b).
	$$
\end{proof}
Proposition \ref{prop:GausEuc} proves Theorem \ref{th:Euclist} for $d=-1$. 

\begin{proof}[Proof of Proposition \ref{prop:QuadIntEuc357}]
	Clearly, ${\mathbb Z}[w_d] \subset {\mathbb Q}(\sqrt{d})$. For any $\alpha=a_1+b_1w_d \in {\mathbb Z}[w_d]$ and $\beta=a_2+b_2w_d \in {\mathbb Z}[w_d]$ we obviously have $\alpha+\beta \in {\mathbb Z}[w_d]$ and $-\alpha \in {\mathbb Z}[w_d]$. Also, each of the listed values of $d$ is of the form $d=-4k+1$ for integer $k$, hence
	$$
		w_d^2 = \left(\frac{1+\sqrt{d}}{2}\right)^2=\frac{1+2\sqrt{d}+d}{4}=\frac{d-1}{4}+w_d = -k + w_d,
	$$
	and
	$$
	\alpha \beta = (a_1+b_1 w_d)\cdot(a_2+b_2 w_d) = a_1a_2+(a_1b_2+a_2b_1)w_d+b_1b_2(-k + w_d) \in {\mathbb Z}[w_d].
	$$
	Thus, the set ${\mathbb Z}[w_d]$ is closed under addition and multiplication and therefore forms a commutative ring. Next, if $\alpha\beta=0$, then $0=N(\alpha\beta)=N(\alpha)N(\beta)$, hence either $N(\alpha)=0$ or $N(\beta)=0$, which implies that $\alpha=0$ or $\beta=0$. Hence, ${\mathbb Z}[w_d]$ is an integral domain. 
	
	We next prove that ${\mathbb Z}[w_d]$ is a Euclidean domain. We follow the proof of Example 12.2 in \cite{dujella2021number}. First observe that for every $\alpha=a+bw_d \in {\mathbb Z}[w_d]$, we have $\alpha=(a+b/2)+(b/2)\sqrt{d}$, hence
	$$
		N(\alpha) = (a+b/2)^2-d(b/2)^2 = (a+b/2)^2-(-4k+1)(b/2)^2 = a^2 + ab + k b^2,
	$$
	hence $N(\alpha)$ is an integer. Now, fix $\alpha,\beta \in {\mathbb Z}[w_d]$ with $\beta\neq 0$. Then $f(b)=|N(b)|$ is a positive integer, hence
	$$
	f(ab) = |N(ab)| = |N(a)||N(b)| \geq |N(a)| = f(a),
	$$
	which is the property (i) in the definition of Euclidean domain. To prove (ii), note that the ratio $a/b$ belongs to ${\mathbb Q}(\sqrt{d})$, and can therefore be written in the form $u+v\sqrt{d}$ for some rational numbers $u,v$. Let $y$ be the nearest integer to $2v$, that is, an integer such that $|2v-y|\leq 1/2$, or, equivalently, $|s|\leq 1/4$, where $s=v-y/2$. Also, let $x$ be the nearest integer to $u-y/2$, that is, an integer such that $|t|\leq 1/2$, where $t=u-y/2-x$. Then $q=x+yw_d$ and $r=a-bq$ are elements of ${\mathbb Z}[w_d]$. Further, 
	$$
	a/b - q = (u+v\sqrt{d})-\left(x+y\frac{1+\sqrt{d}}{2}\right) = t + s\sqrt{d},
	$$
	hence
	$$
	f(r)=|N(r)|=|N(a-bq)|=|N(b)||N(a/b-q)|=f(b)|t^2-ds^2|.
	$$
	Because $|d|<12$, we have $|t^2-ds^2|\leq t^2+|d|s^2 < (1/2)^2+12(1/4)^2 = 1$, hence $f(r)<f(b)$, and property (ii) in the definition of Euclidean domain follows.
\end{proof} 

The proof of Proposition \ref{prop:QuadIntEuc357} is directly applicable to $d=-11$. 

We will now prove Theorem \ref{th:Euclist} for $d=13$. As $|d| \not < 12$, the proof of Proposition  \ref{prop:QuadIntEuc357} is not directly applicable as we would have $|t^2-13 s^2| \leq t^2 + 13 s^2 < (1/2)^2+13(1/4)^2 \not < 1$. However, we can follow the proof up to and including the line
$$
f(r)=|N(r)|=|N(a-bq)|=|N(b)||N(a/b-q)|=f(b)|t^2-ds^2|
$$
and then we must extend the proof, in order to improve the inequality. The inequality $|t^2-13 s^2| \leq t^2 + 13 s^2$ can be improved to $|t^2-13 s^2| \leq \max(t^2, 13 s^2)$, because if $t^2-13 s^2 \geq 0$ then $|t^2-13 s^2| = t^2-13 s^2 \leq t^2$, and if $t^2-13 s^2 \leq 0$ then $|t^2-13 s^2| = 13 s^2-t^2 \leq 13 s^2$.
Hence, we have the following cases to consider, (i) $\max(t^2, 13 s^2) = t^2$ and (ii) $\max(t^2, 13 s^2)=13s^2$.
Let us first consider case (i). Then we have $|t^2-13 s^2| \leq  t^2 < \left( \frac{1}{2} \right) ^2 = \frac{1}{4}$, hence $f(r) < \frac{1}{4}f(b)$ so $f(r) < f(b)$. Let us now look at case (ii). We have $|t^2-13 s^2| \leq  13s^2 < 13 \left( \frac{1}{4} \right) ^2 = \frac{13}{16}$ hence $f(r)< \frac{13}{16}f(b)$ hence $f(r) <f(b)$.
Therefore, property (ii) in the definition of Euclidean domain holds. 

In order to prove Theorem \ref{th:Euclist} for $d=-2,2$ and $3$, we will use the absolute value of the field norm as a Euclidean valuation. The field norm is given by $N(x)=N(a+b\sqrt{d})=(a+b\sqrt{d})(a-b\sqrt{d})=a^2-db^2$. 
Let $a$ and $b \neq 0$ be of the form $m+n\sqrt{d}$ with $m,n$ integers. Then $|N(b)|$ is a positive integer, hence $|N(b)| \geq 1$. Thus $|N(ab)|=|N(a)||N(b)| \geq |N(a)|$, and property (i) in the definition of Euclidean domain follows. To prove (ii), write the number $a/b$ in the form $u+v\sqrt{d}$ for some rational numbers $u,v$. Let us choose integers $x$ and $y$ that are the nearest to $u$ and $v$, respectively, that is, $|u-x|\leq 1/2$ and $|v-y|\leq 1/2$. Then $q=x+y\sqrt{d}$ and $r=a-bq$ are of the form $m+n\sqrt{d}$ with $m,n$ integers.  Also,
$$
|N(r)|=|N(a-bq)|=|N(b)||N(a/b-q)|=|N(b)||N((u-x)+(v-y)\sqrt{d})|
$$
$$
=|N(b)|\cdot |(u-x)^2-d(v-y)^2| \leq |N(b)|\max(|u-x|^2,|d||v-y|^2).
$$
This then gives two cases, (a) $\max(|u-x|^2,|d||v-y|^2) = |u-x|^2$ and (b) $\max(|u-x|^2,|d||v-y|^2) = |d||v-y|^2$.
Let us first look at case (a). We then have 
$$
|N(r)| \leq |N(b)|\max(|u-x|^2,|d||v-y|^2) = |N(b)||u-x|^2 \leq |N(b)|(1/2)^2 = |N(b)|/4 < |N(b)|.
$$ 
Let us now look at case (b). We then have 
$$
|N(r)| \leq |N(b)|\max(|u-x|^2,|d||v-y|^2) = |N(b)|(|d||v-y|^2) 
$$
$$
\leq |N(b)|(|d|(1/2)^2) = |d||N(b)|/4 < |N(b)|,
$$ 
where the last inequality follows from $|d|<4$. In both cases, we have $|N(r)| < |N(b)|$.

Hence, we have shown that Theorem \ref{th:Euclist} holds for $d=-2,2,3$.

Finally, as we have proved Theorem \ref{th:Euclist} for $d=-11$ and $d=13$, we can then conclude that the proof of Proposition \ref{prop:x2pyzpxp1} is directly applicable for $k = \pm 3$.

\subsection{Exercise 2.31}\label{ex:divfunct}
\textbf{\emph{Use the divisor function $D$ to describe the sets of integer solutions to the equations listed in Table \ref{tab:H13open}.}}

\begin{center}
\begin{tabular}{ |c|c|c|c|c|c| } 
 \hline
 $H$ & Equation & $H$ & Equation & $H$ & Equation \\ 
 \hline\hline
 $13$ & $x^2+y^2+zt+1=0$ & $13$ & $xyz+t^2+1=0$ & $13$ & $x^2y+z^2+1=0$ \\ 
 \hline
 $13$ & $x^2+y^2+zt-1=0$ & $13$ & $xyz+t^2-1=0$ & $13$ & $x^2y+z^2-1=0$ \\ 
 \hline
 $13$ & $x^3+yz+1=0$ & $13$ & $x^2y+zt+1=0$ & $13$ & $x_1x_2x_3+x_4x_5+1=0$ \\ 
 \hline
\end{tabular}
\captionof{table}{\label{tab:H13open} Equations of size $H=13$ for which existence of polynomial parametrization is unknown.}
\end{center} 

In this exercise, we will describe all integer solutions to equations in $n \geq 3$ variables (of size $H=13$), 
using the divisor function. This is, for any integer $k \geq 1$ and $m$, we will use the notation 
\begin{equation}\label{Divisor Function}
D_k(m) =\{x \in \mathbb{Z}, \quad \text{such that } x^k \text{ is a divisor of } m \}.
\end{equation}  
If $k=1$, we will simplify the notation, instead of writing $D_1(m)$, we will write $D(m)$. 

Let us consider one equation in detail, for example, the equation
\begin{equation}\label{eq:x3pyzp1}
	x^3+yz+1=0.
\end{equation}
If $y=0$, then $x^3=-1$ and $z$ can be arbitrary, therefore integer solutions to \eqref{eq:x3pyzp1} with $y=0$ are described by $(x,y,z)=(-1,0,u)$ where $u$ is an arbitrary integer. Now assume that $y \neq 0$. Then \eqref{eq:x3pyzp1} implies that $z=-\frac{x^3+1}{y}$. Hence, we may let $x=u$ to be an arbitrary integer, then let $y$ to be any divisor $v$ of $u^3+1$, and then $z=-\frac{1+u^3}{v}$.
Hence, we may summarise all integer solutions to equation \eqref{eq:x3pyzp1} as
 \begin{equation*}\label{eq:x3pyzp1sol}
 	(x,y,z)=(-1,0,u), \,\, \text{or} \,\, \left(u,v,-\frac{1+u^3}{v}\right), \quad u \in \mathbb{Z}, \quad v \in D(1+u^3).
 \end{equation*}

All other equations in Table \ref{tab:H13open} can be solved similarly and their solutions are presented in Table \ref{tab:H13opensol}.

\begin{center}
\begin{tabular}[c]{| c  | c |c|c|}
\hline
Equation & Solution $(x,y,z)$ or $(x,y,z,t)$, or $(x_1,x_2,x_3,x_4,x_5)$ \\\hline \hline
 $x^2+y^2+zt+1=0$ & $\left(u,v,w,-\frac{1+u^2+v^2}{w}\right), \quad u,v \in \mathbb{Z}, \quad w \in D(1+u^2+v^2)$ \\ \hline
$x^2+y^2+zt-1=0$ &$(\pm 1,0,0,u), (0,\pm 1,0,u), \left(u,v,w,\frac{1-u^2-v^2}{w}\right), $ \\
& $ u,v \in \mathbb{Z}, \quad w \in D(1-u^2-v^2),$ \\\hline
$x^3+yz+1=0$ & $(-1,0,u), \left(u,v,-\frac{1+u^3}{v}\right), \quad u \in \mathbb{Z}, \quad v \in D(1+u^3)$ \\\hline
$xyz+t^2+1=0$ & $\left(u,v, - \frac{w^2+1}{uv},w \right), \quad w \in \mathbb{Z}, \quad u \in D(w^2+1), \quad v \in D\left(\frac{w^2+1}{u}\right)$ \\\hline
$xyz+t^2-1=0$ & $\left(u,v, \frac{1-w^2}{uv},w \right), \quad w \in \mathbb{Z}, \quad u \in D(1-w^2), \quad v \in D\left(\frac{1-w^2}{u}\right)$ \\
& $(0,u,v,\pm 1), (u,0,v,\pm1), \quad u,v \in \mathbb{Z}$ \\
\hline
$x^2y+zt+1=0$& $(\pm 1,-1,0,u), \left(u,v, w, - \frac{1+u^2v}{w} \right), \quad u,v \in \mathbb{Z}, \quad w \in D(1+u^2v)$ \\\hline
 $x^2y+z^2+1=0$ & $\left(v,-\frac{u^2+1}{v^2},u \right), \quad u \in \mathbb{Z}, \quad v \in D_2(u^2+1)$ \\\hline
$x^2y+z^2-1=0$ & $(0,u,\pm 1), \left(v,\frac{1-u^2}{v^2},u \right), \quad u \in \mathbb{Z}, \quad v \in D_2(1-u^2)$ \\\hline
$x_1x_2x_3+x_4x_5+1=0$ &$(-1,-1,-1,0,u),(-1,1,1,0,u),(1,1,-1,0,u),(1,-1,1,0,u), $ \\
& $\left(u,v,w,t,- \frac{1+uvw}{t}\right), \quad u,v,w \in \mathbb{Z}, \quad t \in D(1+uvw)$  \\\hline
\end{tabular}
\captionof{table}{\label{tab:H13opensol} Integer solutions to the equations listed in Table \ref{tab:H13open}.}
\end{center}

\section{Chapter 3}
\subsection{Exercise 3.3}\label{ex:Vieta2var}
\textbf{\emph{Solve each of the equations listed in Table \ref{tab:H18Vietavar2}, preferably with human proofs. 
	}}

\begin{center}
\begin{tabular}{ |c|c|c|c|c|c| } 
 \hline
 $H$ & Equation & $H$ & Equation & $H$ & Equation \\ 
 \hline\hline
 $13$ & $x^2+xy-y^2-1=0$  & $16$ & $x^2+xy-y^2+4=0$ &  $17$ & $x^2+x+x y-y^2-3=0$  \\ 
 \hline
 $13$ & $x^2+xy-y^2+1=0$ & $16$ & $x^2+x+x y-y^2-2=0$  & $18$ & $x^2+x+x y-y^2+4=0$ \\ 
 \hline
 $14$ & $x^2+x+x y-y^2 = 0$ & $16$ & $x^2+x+x y-y^2+2=0$ & $18$ & $x^2+3x+x y-y^2 = 0$ \\ 
 \hline
 $15$ & $x^2+x+x y-y^2+1 = 0$ & $17$ & $x^2+xy-y^2-5=0$ & &  \\ 
 \hline
$16$ & $x^2+xy-y^2-4=0$  & $17$ & $x^2+xy-y^2+5=0$ &  & \\ 
 \hline
\end{tabular}
\captionof{table}{\label{tab:H18Vietavar2} Two-variable quadratic equations of size $H\leq 18$ solvable by Vieta jumping.}
\end{center}

This section solves equations of the form 
\begin{equation}\label{eq:quad2vargen2}
ax^2+bxy+cy^2+dx+ey+f=0 
\end{equation} 
with $a,b,c,d,e,f$ integers using Vieta jumping. To solve these equations, we will make a jump in $x$ and a jump in $y$, and use these to make a contradiction for integers $x$ and $y$. We have that for a solution $(x,y)$ to equation \eqref{eq:quad2vargen2}, $\left(- \frac{by+d}{a}-x,y\right)$ and $\left(x,-\frac{bx+e}{c}-y\right)$ are also solutions. To solve equations, we will use $x_1=- \frac{by+d}{a}-x$ and $y_1=-\frac{bx+e}{c}-y$. Let (a) be the transformation $\left(- \frac{by+d}{a}-x,y\right)$, and let (b) be the transformation $\left(x,-\frac{bx+e}{c}-y\right)$. We will also use transformations (ab), first apply (a) to $(x,y)$ and then apply (b), and the transformation (ba), first apply (b) to $(x,y)$ and then apply (a). Applying these transformations to $(x,y)$ will give a chain of solutions, hence we can represent solutions to equation \eqref{eq:quad2vargen2} recursively. We will call an integer solution \textbf{minimal} if neither of the transformations (a) and (b) strictly decrease the norm $N$, where $N=|x|+|y|$.

Equation 
$$
x^2+xy-y^2-1=0
$$
is solved in Section 3.1.1 of the book, and its integer solutions are 
$$
(x,y)=\pm (F_{2n-1},F_{2n}) \quad \text{or} \quad \pm (F_{2n+1},-F_{2n}), \quad n=0,1,2,\dots,
$$ 
where $F_n$ is the $n^{th}$ term of the Fibonacci sequence, defined by
\begin{equation}\label{eq:Fibonacci}
	F_0 = 0, \quad F_1 = 1, \quad F_{n+2} = F_{n+1} + F_n, \quad n=0,1,2,\dots,
\end{equation} 
and we denote $F_{-1}=1$ by convention.

Equation 
$$
x^2+xy-y^2+1=0
$$
is solved in Section 3.1.1 of the book, and its integer solutions are
$$
(x,y)=\pm (F_{2n-1},F_{2n}) \quad \text{or} \quad \pm (F_{2n+1},-F_{2n}) \quad \text{where} \quad F_n \quad \text{is given by \eqref{eq:Fibonacci}}.
$$

Equation 
$$
x^2+x+xy-y^2=0
$$
is solved in Section 3.1.2 of the book, and its integer solutions are given by $(x,y)=(x_n,y_n)$, where
$$
	(x_0,y_0)=(0,0) \quad \text{and} \quad (x_{n+1},y_{n+1})=(-x_n-y_n-1, -x_n - 2y_n-1), \quad n=0,1,2\dots,
$$ 
and $(x,y)=(x'_n,y'_n)$, where
$$
	(x'_0,y'_0)=(0,0) \quad \text{and} \quad (x'_{n+1},y'_{n+1})=(-2x'_n+y'_n-1, x'_n - y'_n), \quad n=0,1,2\dots.
$$

The first equation we will consider is
\begin{equation}\label{vieta_1}
x^2+x+xy-y^2+1=0
\end{equation}
For any solution $(x,y)$ to \eqref{vieta_1}, the transformations (a) $(x,y)\to(-x-y-1,y)$ and (b) $(x,y)\to(x,x-y)$ return integer solutions to the same equation. 
Let $(x,y)$ be a minimal solution of equation \eqref{vieta_1}. Assume that $x$ and $y$ are non-zero. Then, by Vieta jumping $(x_1,y)$ is also an integer solution where $x_1=-1-y-x=-\frac{y^2-1}{x}$. Because transformation (a) does not decrease $|x|+|y|$, we must have $|x_1| \geq |x|$ or $\left| \frac{1-y^2}{x}\right| \geq |x|$ or $|y^2-1|\geq |x|^2$. Because $y \neq 0$, $|y^2-1|=|y|^2-1$ hence, $|y|^2-1 \geq |x|^2$ hence $|y| > |x|$ so $|y| \geq |x|+1$. Now do Vieta jumping in $y$ and obtain the integer solution $(x,y_1)$ where $y_1=x-y=-\frac{x^2+x+1}{y}$ because transformation (b) does not decrease $|x|+|y|$, $|y_1| \geq |y|$ or $\left|-\frac{x^2+x+1}{y}\right|\geq |y|$ or $x^2+|x|+1\geq |x^2+x+1| \geq |y|^2 \geq (|x|+1)^2=|x|^2+2|x|+1$ which is a contradiction. Hence, either $x=0$ or $y=0$. Therefore, the minimal solutions to equation \eqref{vieta_1} are $(x,y)=(0,-1)$ and $(0,1)$.
Using transformations (a) and (b) starting from $(0,-1)$, we obtain the chain of solutions
$$
\dots \xleftrightarrow{\text{(b)}} (-2,1) \xleftrightarrow{\text{(a)}}  (0,1) \xleftrightarrow{\text{(b)}} (0,-1) \xleftrightarrow{\text{(a)}} (0,-1) \xleftrightarrow{\text{(b)}} (0,1) \xleftrightarrow{\text{(a)}} (-2,1) \xleftrightarrow{\text{(b)}} (-2,-3) \xleftrightarrow{\text{(a)}} \dots 
$$
Because solution $(x,y)=(0,1)$ is also in this chain, it contains all integer solutions to \eqref{vieta_1}.
This chain is symmetric, so we may only look at solutions on the right from $(0,-1)$. 
We can represent any solution $(x,y)$ to equation \eqref{vieta_1} using the transformations (ab) and (ba) $n$ times for some $n\geq0$. In this case, (ab) is $(x,y) \to (-x-y-1,-x-2y-1)$ and (ba) is $(x,y) \to (-2x+y-1,x-y)$. Therefore, the integer solutions to \eqref{vieta_1} can be written as $(x,y)=(x_n,y_n)$ where
$$
(x_0,y_0)=(0,-1) \text{ and } (x_{n+1},y_{n+1})=(-x_n-y_n-1,-x_n-2y_n-1), \quad n=0,1,2,\ldots,
$$
or, $(x,y)=(x'_n,y'_n)$, where
$$
(x'_0,y'_0)=(0,-1) \text{ and } (x'_{n+1},y'_{n+1})=(-2x'_n+y'_n-1,x'_n-y'_n), \quad n=0,1,2,\ldots.
$$

\vspace{10pt}

The next equation we will consider is
\begin{equation}\label{vieta_2}
x^2+xy-y^2-4=0.
\end{equation}
For this equation, the transformations (a) and (b) are $(x,y) \to (-x-y,y)$ and $(x,y) \to (x,x-y)$. 
Let $(x,y)$ be a minimal solution of equation \eqref{vieta_2}. Assume that $x$ and $y$ are non-zero. Then, by Vieta jumping $(x_1,y)$ is also an integer solution where $x_1=-y-x=\frac{-y^2-4}{x}$. Because transformation (a) does not decrease $|x|+|y|$, we must have $|x_1| \geq |x|$ or $\left| \frac{4+y^2}{x}\right| \geq |x|$ or $|y^2+4|\geq |x|^2$, or $y^2+4 \geq x^2$ hence, $y^2 \geq x^2-4$. Now do Vieta jumping in $y$ and obtain the integer solution $(x,y_1)$ where $y_1=x-y=\frac{-x^2+4}{y}$. Because transformation (b) does not decrease $|x|+|y|$, $|y_1| \geq |y|$ or $\left|\frac{x^2-4}{y}\right|\geq |y|$ or $|x^2-4| \geq |y|^2$. $|x^2-4| = x^2-4$ provided that $|x| \geq 2$. So, we have that $x^2-4 \geq y^2 \geq x^2-4$, hence $y^2=x^2-4$, which implies that $(x,y)=(-2,0),(2,0)$. We must also check $|x| \leq 1$, which does not give integer solutions to equation \eqref{vieta_2}. Therefore, the minimal solutions to equation \eqref{vieta_2} are $(x,y)=(-2,0),(2,0)$.
Using transformations (a) and (b), starting from $(-2,0)$ we obtain the chain of solutions
$$
\dots \xleftrightarrow{\text{(b)}} (4,-2) \xleftrightarrow{\text{(a)}}  (-2,-2) \xleftrightarrow{\text{(b)}} (-2,0) \xleftrightarrow{\text{(a)}} (2,0) \xleftrightarrow{\text{(b)}} (2,2) \xleftrightarrow{\text{(a)}} (-4,2) \xleftrightarrow{\text{(b)}} (-4,-6) \xleftrightarrow{\text{(a)}} \dots 
$$
Because solution $(2,0)$ is also in this chain, it contains all integer solutions to \eqref{vieta_2}. The left-hand side of the chain, from $(-2,0)$ is a negation of the right-hand side. So, we can represent any solution $\pm (x,y)$ to equation \eqref{vieta_2} using the transformations (ab) and (ba) $n$ times for some $n\geq0$. In this case, (ab) is $(x,y) \to (-x-y,-x-2y)$ and (ba) is $(x,y) \to (-2x+y,x-y)$. Therefore, the integer solutions to \eqref{vieta_2} can be written as  $(x,y)=\pm (x_n,y_n)$ where 
$$
(x_0,y_0)=(-2,0) \text{ and } (x_{n+1},y_{n+1})=(-x_n-y_n,-x_n-2y_n), \quad n=0,1,2,\ldots,
$$
or, $(x,y)=\pm (x'_n,y'_n)$, where
$$
(x'_0,y'_0)=(-2,0) \text{ and } (x'_{n+1},y'_{n+1})=(-2x'_n+y'_n,x'_n-y'_n), \quad n=0,1,2,\ldots.
$$

\vspace{10pt}

The next equation we will consider is
\begin{equation}\label{vieta_3}
x^2+xy-y^2+4=0.
\end{equation}
For this equation, the transformations (a) and (b) are $(x,y) \to (-x-y,y)$ and $(x,y) \to (x,x-y)$. Let $(x,y)$ be a minimal solution of equation \eqref{vieta_3}. Assume that $x$ and $y$ are non-zero. Then, by Vieta jumping $(x_1,y)$ is also an integer solution where $x_1=-y-x=-\frac{y^2-4}{x}$. Because transformation (a) does not decrease $|x|+|y|$, we must have $|x_1| \geq |x|$ or $\left| \frac{y^2-4}{x}\right| \geq |x|$ or $|y^2-4|\geq |x|^2$, or $y^2-4 \geq x^2$ provided that $|y| \geq 2$. Hence, $y^2-4 \geq x^2$. Now do Vieta jumping in $y$ and obtain the integer solution $(x,y_1)$ where $y_1=x-y=\frac{x^2+4}{y}$. Because transformation (b) does not decrease $|x|+|y|$, $|y_1| \geq |y|$ or $\left|\frac{x^2+4}{y}\right|\geq |y|$ or $|x^2+4| \geq |y|^2$. $|x^+4| = x^2+4$, so, we have that $x^2 \geq y^2-4 \geq x^2$ which only has solutions for $y^2-4=x^2$, so we have integer solutions $(x,y)=(0,2)$ and $(0,-2)$. We must also check $|y| \leq 1$, which does not give integer solutions to equation \eqref{vieta_3}. Therefore, the minimal solutions to equation \eqref{vieta_3} are $(x,y)=(0,-2),(0,2)$.
Using transformations (a) and (b), starting from $(0,-2)$, we obtain a chain of solutions
$$
\dots \xleftrightarrow{\text{(b)}} (6,-4) \xleftrightarrow{\text{(a)}}  (-2,-4) \xleftrightarrow{\text{(b)}} (-2,2) \xleftrightarrow{\text{(a)}} (0,2) \xleftrightarrow{\text{(b)}} (0,-2) \xleftrightarrow{\text{(a)}} (2,-2) \xleftrightarrow{\text{(b)}} (2,4) \xleftrightarrow{\text{(a)}} \dots 
$$
Because solution $(0,2)$ is also in this chain, it contains all integer solutions to \eqref{vieta_3}. The left-hand side of the chain from $(0,2)$ is a negation of the right-hand side. So, we can represent any solution $\pm (x,y)$ to equation \eqref{vieta_3} using the transformations (ab) and (ba) $n$ times for some $n\geq0$. In this case, (ab) is $(x,y) \to (-x-y,-x-2y)$ and (ba) is $(x,y) \to (-2x+y,x-y)$. Therefore, the integer solutions to \eqref{vieta_3} can be written as $(x,y)= \pm (x_n,y_n)$ where 
$$
(x_0,y_0)=(0,-2) \text{ and } (x_{n+1},y_{n+1})=(-x_n-y_n,-x_n-2y_n), \quad n=0,1,2,\ldots,
$$
or, $(x,y)= \pm (x'_n,y'_n)$, where
$$
(x'_0,y'_0)=(0,-2) \text{ and } (x'_{n+1},y'_{n+1})=(-2x'_n+y'_n,x'_n-y'_n), \quad n=0,1,2,\ldots.
$$

\vspace{10pt}

The next equation we will consider is
\begin{equation}\label{vieta_4}
x^2+x+xy-y^2-2=0.
\end{equation}
For this equation, the transformations (a) and (b) are $(x,y) \to (-x-y-1,y)$ and $(x,y) \to (x,x-y)$.
Let $(x,y)$ be a minimal solution of equation \eqref{vieta_4}. Assume that $x$ and $y$ are non-zero. Then, by Vieta jumping $(x_1,y)$ is also an integer solution where $x_1=-1-y-x=-\frac{y^2+2}{x}$. Because transformation (a) does not decrease $|x|+|y|$, we must have $|x_1| \geq |x|$ or $\left| \frac{y^2+2}{x}\right| \geq |x|$ or $|y^2+2|\geq |x|^2$, or $y^2+2 \geq x^2$. Now do Vieta jumping in $y$ and obtain the integer solution $(x,y_1)$ where $y_1=x-y=-\frac{x^2+x-2}{y}$. Because transformation (b) does not decrease $|x|+|y|$, $|y_1| \geq |y|$ or $\left|-\frac{x^2+x-2}{y}\right|\geq |y|$ or $|x^2+x-2| \geq y^2$. Also, $|x^2+x-2| = x^2+x-2$ provided that $x\geq 1$ or $x \leq -2$. So we have that $x^2+x-2 \geq y^2 \geq x^2-2$. Let $|y| \geq |x|+1$ so $x^2+x-2 \geq y^2 \geq x^2+2|x|+1$ which is impossible. Let $|y| \leq |x|-1$, then we have $x^2-2 \leq y^2 \leq x^2-2|x|+1$ which is only true for $|x| \leq \frac{3}{2}$. So, we have the following cases to check, $|y|=|x|$ or $x=-1, 0,$ or $1$.  These cases give the following integer solutions $(x,y)=(-2,-2),(1,0),(1,1)$. 
Using transformations (a) and (b), starting from $(1,0)$ we obtain a chain of solutions
\begin{equation}\label{eq:vieta_4_chain}
\dots \xleftrightarrow{\text{(a)}}  (-3,-4) \xleftrightarrow{\text{(b)}} (-3,1) \xleftrightarrow{\text{(a)}} (1,1) \xleftrightarrow{\text{(b)}} (1,0) \xleftrightarrow{\text{(a)}} (-2,0) \xleftrightarrow{\text{(b)}} (-2,-2) \xleftrightarrow{\text{(a)}} (3,-2) \xleftrightarrow{\text{(b)}} \dots 
\end{equation}
From this chain, we can see that the solutions $(-2,-2)$ and $(1,1)$ can be transformed into smaller ones, hence they are not minimal, thus the only minimal solution is $(1,0)$. Hence, \eqref{eq:vieta_4_chain} 
contains all integer solutions to \eqref{vieta_4}. Let us represent any solution $(x,y)$ to equation \eqref{vieta_4} using the transformations (ab) and (ba) $n$ times for some $n\geq0$. In this case, (ab) is $(x,y) \to (-x-y-1,-x-2y-1)$ and (ba) is $(x,y) \to (-2x+y-1,x-y)$. We can represent the left-hand side of the chain by applying (ba) to $(1,0)$ and applying (ab) to $(1,1)$. We can represent the right-hand side of the chain by applying (ab) to $(1,0)$ and applying (ba) to $(1,1)$. 

Therefore, all integer solutions can be described as $(x,y)=(x_n,y_n)$ where 
$$
(x_0,y_0)=(1,0) \text{ and } (x_{n+1},y_{n+1})=(-x_n-y_n-1,-x_n-2y_n-1), \quad n=0,1,2,\ldots,
$$
or, $(x,y)=(x'_n,y'_n)$, where
$$
(x'_0,y'_0)=(1,0) \text{ and } (x'_{n+1},y'_{n+1})=(-2x'_n+y'_n-1,x'_n-y'_n), \quad n=0,1,2,\ldots,
$$
or, $(x,y)=(x''_n,y''_n)$ where 
$$
(x''_0,y''_0)=(1,1) \text{ and } (x''_{n+1},y''_{n+1})=(-x''_n-y''_n-1,-x''_n-2y''_n-1), \quad n=0,1,2,\ldots,
$$
or, $(x,y)=(x'''_n,y'''_n)$, where
$$
(x'''_0,y'''_0)=(1,1) \text{ and } (x'''_{n+1},y'''_{n+1})=(-2x'''_n+y'''_n-1,x'''_n-y'''_n), \quad n=0,1,2,\ldots.
$$

\vspace{10pt}

The next equation we will consider is
\begin{equation}\label{vieta_5}
x^2+x+xy-y^2+2=0.
\end{equation}
For this equation, the transformations (a) and (b) are $(x,y) \to (-x-y-1,y)$ and $(x,y) \to (x,x-y)$.
Let $(x,y)$ be a minimal solution of equation \eqref{vieta_5}. Assume that $x$ and $y$ are non-zero. Then, by Vieta jumping $(x_1,y)$ is also an integer solution where $x_1=-1-y-x=-\frac{y^2-2}{x}$. Because transformation (a) does not decrease $|x|+|y|$, we must have $|x_1| \geq |x|$ or $\left| \frac{2-y^2}{x}\right| \geq |x|$ or $|y^2-2|\geq |x|^2$. If $|y| \geq 2$, $|y^2-2|=|y|^2-2$ hence, $|y|^2-2 \geq |x|^2$ hence $|y| > |x|$ so $|y| \geq |x|+1$. Now do Vieta jumping in $y$ and obtain the integer solution $(x,y_1)$ where $y_1=x-y=-\frac{x^2+x+2}{y}$. Because transformation (b) does not decrease $|x|+|y|$, $|y_1| \geq |y|$ or $\left|-\frac{x^2+x+2}{y}\right|\geq |y|$ or $x^2+4|x|+4 > x^2+|x|+2 \geq |x^2+x+2| \geq |y|^2 \geq (|x|+1)^2=|x|^2+2|x|+1$ which is only true for $1 \geq |x|$. Hence, either $|x| \leq 1$ or $|y|\leq1$. We then find the integer solutions to equation \eqref{vieta_5} with $|x| \leq 1$ or $|y|\leq1$ are $(x,y)=(-1,-2), (-1,1)$.
Using transformations (a) and (b), starting with $(-1,1)$, we obtain the chain of solutions
$$
\dots \xleftrightarrow{\text{(b)}} (2,-2) \xleftrightarrow{\text{(a)}}  (-1,-2) \xleftrightarrow{\text{(b)}} (-1,1) \xleftrightarrow{\text{(a)}} (-1,1) \xleftrightarrow{\text{(b)}} (-1,-2) \xleftrightarrow{\text{(a)}} (2,-2) \xleftrightarrow{\text{(b)}} (2,4) \xleftrightarrow{\text{(a)}} \dots 
$$
Because $(-1,-2)$ is also in this chain, it contains all integer solutions to \eqref{vieta_5}. This chain is symmetric, so we may only look at solutions on the right from $(-1,1)$. 
We can represent any solution $(x,y)$ to equation \eqref{vieta_5} using the transformations (ab) and (ba) $n$ times for some $n\geq0$. In this case, (ab) is $(x,y) \to (-x-y-1,-x-2y-1)$ and (ba) is $(x,y) \to (-2x+y-1,x-y)$. Therefore, the integer solutions to \eqref{vieta_5} can be written as $(x,y)=(x_n,y_n)$ where
$$
(x_0,y_0)=(-1,1) \text{ and } (x_{n+1},y_{n+1})=(-x_n-y_n-1,-x_n-2y_n-1), \quad n=0,1,2,\ldots,
$$
or, $(x,y)=(x'_n,y'_n)$, where
$$
(x'_0,y'_0)=(-1,1) \text{ and } (x'_{n+1},y'_{n+1})=(-2x'_n+y'_n-1,x'_n-y'_n), \quad n=0,1,2,\ldots.
$$

\vspace{10pt}

The next equation we will consider is
\begin{equation}\label{vieta_6}
x^2+xy-y^2-5=0.
\end{equation}
For this equation, the transformations (a) and (b) are $(x,y) \to (-x-y,y)$ and $(x,y) \to (x,x-y)$. 
Let $(x,y)$ be a minimal solution of equation \eqref{vieta_6}. Assume that $x$ and $y$ are non-zero. Then, by Vieta jumping $(x_1,y)$ is also an integer solution where $x_1=-y-x=-\frac{y^2+5}{x}$. Because transformation (a) does not decrease $|x|+|y|$, we must have $|x_1| \geq |x|$ or $\left| \frac{y^2+5}{x}\right| \geq |x|$ or $|y^2+5|\geq |x|^2$, or $y^2+5 \geq x^2$. Now do Vieta jumping in $y$ and obtain the integer solution $(x,y_1)$ where $y_1=x-y=\frac{-x^2+5}{y}$. Because transformation (b) does not decrease $|x|+|y|$, $|y_1| \geq |y|$ or $\left|\frac{x^2-5}{y}\right|\geq |y|$ or $|x^2-5| \geq |y|^2$. Also, $|x^2-5| = x^2-5$, provided that $|x| \geq 3$, so, we have that $x^2 \geq y^2+5 \geq x^2$ which only has solutions if $x^2=y^2+5$, whose only integer solutions are $(x,y)=(\pm 3,\pm 2)$. Substituting these values of $x$ and $y$ into \eqref{vieta_6} we do not obtain integer solutions to equation \eqref{vieta_6}. We must now check $|x| \leq 2$, which gives minimal solutions to equation \eqref{vieta_6} $(x,y)=(-2,-1),(2,1)$. 
Using transformations (a) and (b), starting with $(-2,-1)$ we obtain the chain of solutions
\begin{equation}\label{eq:vieta_6_chain_i}
\dots \xleftrightarrow{\text{(b)}} (-7,4) \xleftrightarrow{\text{(a)}}  (3,4) \xleftrightarrow{\text{(b)}} (3,-1) \xleftrightarrow{\text{(a)}} (-2,-1) \xleftrightarrow{\text{(b)}} (-2,-1) \xleftrightarrow{\text{(a)}} (3,-1) \xleftrightarrow{\text{(b)}} (3,4) \xleftrightarrow{\text{(a)}}  \dots 
\end{equation}
We can see that this chain is symmetric, but $(2,1)$ does not appear in this chain. 
So we can create another chain using transformations (a) and (b) starting with $(2,1)$ and we obtain the chain of solutions
$$
\dots \xleftrightarrow{\text{(b)}}  (7,-4) \xleftrightarrow{\text{(a)}}  (-3,-4) \xleftrightarrow{\text{(b)}} (-3,1) \xleftrightarrow{\text{(a)}} (2,1) \xleftrightarrow{\text{(b)}} (2,1) \xleftrightarrow{\text{(a)}} (-3,1) \xleftrightarrow{\text{(b)}} (-3,-4) \xleftrightarrow{\text{(a)}} \dots 
$$
We can see that this chain, is a negation of chain \eqref{eq:vieta_6_chain_i}. Therefore, we can represent any solution $\pm (x,y)$ to equation \eqref{vieta_6} using the transformations (ab) and (ba) $n$ times for some $n\geq0$. In this case, (ab) is $(x,y) \to (-x-y,-x-2y)$ and (ba) is $(x,y) \to (-2x+y,x-y)$. Therefore, the integer solutions to \eqref{vieta_6} can be written as $(x,y)=\pm(x_n,y_n)$ where
$$
(x_0,y_0)=(2,1) \text{ and } (x_{n+1},y_{n+1})=(-x_n-y_n,-x_n-2y_n), \quad n=0,1,2,\ldots,
$$
or, $(x,y)=\pm(x'_n,y'_n)$, where
$$
(x'_0,y'_0)=(2,1) \text{ and } (x'_{n+1},y'_{n+1})=(-2x'_n+y'_n,x'_n-y'_n), \quad n=0,1,2,\ldots.
$$

\vspace{10pt}

The next equation we will consider is
\begin{equation}\label{vieta_7}
x^2+xy-y^2+5=0.
\end{equation}
For this equation, the transformations (a) and (b) are $(x,y) \to (-x-y,y)$ and $(x,y) \to (x,x-y)$.  
Let $(x,y)$ be a minimal solution of equation \eqref{vieta_7}. Assume that $x$ and $y$ are non-zero. Then, by Vieta jumping $(x_1,y)$ is also an integer solution where $x_1=-y-x=\frac{-y^2+5}{x}$. Because transformation (a) does not decrease $|x|+|y|$, we must have $|x_1| \geq |x|$ or $\left| \frac{y^2-5}{x}\right| \geq |x|$ or $|y^2-5|\geq |x|^2$, or $y^2-5 \geq x^2$ provided that $|y| \geq 3$. Hence, $y^2-5 \geq x^2$. Now do Vieta jumping in $y$ and obtain the integer solution $(x,y_1)$ where $y_1=x-y=\frac{x^2+5}{y}$. Because transformation (b) does not decrease $|x|+|y|$, $|y_1| \geq |y|$ or $\left|\frac{x^2+5}{y}\right|\geq |y|$ or $|x^2+5| \geq |y|^2$. Also, $|x^2+5| = x^2+5$, so, we have that $y^2 \geq x^2+5 \geq y^2$ which only has solutions if $y^2=x^2+5$, and this equations only integer solutions are $(x,y)=(\pm2,\pm3)$. After substituting these values into \eqref{vieta_7}, we see that these are not integer solutions. We must now check $|y| \leq 2$, and we obtain a list of integer solutions, including the minimal solutions, to equation \eqref{vieta_7} is $(x,y)=(-1,2),(1,-2)$. 
Using transformations (a) and (b), starting with $(-1,2)$ we obtain the chain of solutions
\begin{equation}\label{eq:vieta_7_chain_i}
\dots \xleftrightarrow{\text{(b)}} (4,-3) \xleftrightarrow{\text{(a)}}  (-1,-3) \xleftrightarrow{\text{(b)}} (-1,2) \xleftrightarrow{\text{(a)}} (-1,2) \xleftrightarrow{\text{(b)}} (-1,-3) \xleftrightarrow{\text{(a)}} (4,-3) \xleftrightarrow{\text{(b)}} (4,7) \xleftrightarrow{\text{(a)}} \dots 
\end{equation}
We can see that this chain is symmetric, but the solution $(1,-2)$ does not appear in this chain. 
So we can create another chain using transformations (a) and (b) starting with $(1,-2)$ and we obtain the chain of solutions
$$
\dots \xleftrightarrow{\text{(b)}}  (-4,3) \xleftrightarrow{\text{(a)}}  (1,3) \xleftrightarrow{\text{(b)}} (1,-2) \xleftrightarrow{\text{(a)}} (1,-2) \xleftrightarrow{\text{(b)}} (1,3) \xleftrightarrow{\text{(a)}} (-4,3) \xleftrightarrow{\text{(b)}} (-4,-7) \xleftrightarrow{\text{(a)}} \dots 
$$
We can see that this chain, is a negation of chain \eqref{eq:vieta_7_chain_i}. Therefore, we can represent any solution $\pm (x,y)$ to equation \eqref{vieta_7} using the transformations (ab) and (ba) $n$ times for some $n\geq0$. In this case, (ab) is $(x,y) \to (-x-y,-x-2y)$ and (ba) is $(x,y) \to (-2x+y,x-y)$. Therefore, the integer solutions to \eqref{vieta_7} can be written as $(x,y)=\pm(x_n,y_n)$ where
$$
(x_0,y_0)=(1,-2) \text{ and } (x_{n+1},y_{n+1})=(-x_n-y_n,-x_n-2y_n), \quad n=0,1,2,\ldots,
$$
or, $(x,y)=\pm(x'_n,y'_n)$, where
$$
(x'_0,y'_0)=(1,-2) \text{ and } (x'_{n+1},y'_{n+1})=(-2x'_n+y'_n,x'_n-y'_n), \quad n=0,1,2,\ldots.
$$

\vspace{10pt}

The next equation we will consider is
\begin{equation}\label{vieta_8}
x^2+x+xy-y^2-3=0.
\end{equation}
For this equation, the transformations (a) and (b) are $(x,y) \to (-x-y-1,y)$ and $(x,y) \to (x,x-y)$.
Let $(x,y)$ be a minimal solution of equation \eqref{vieta_8}. Assume that $x$ and $y$ are non-zero. Then, by Vieta jumping $(x_1,y)$ is also an integer solution where $x_1=-y-x-1=\frac{-y^2-3}{x}$. Because transformation (a) does not decrease $|x|+|y|$, we must have $|x_1| \geq |x|$ or $\left| \frac{y^2+5}{x}\right| \geq |x|$ or $|y^2+3|\geq |x|^2$, or $y^2+3 \geq x^2$. Now do Vieta jumping in $y$ and obtain the integer solution $(x,y_1)$ where $y_1=x-y=-\frac{x^2+x-3}{y}$. Because transformation (b) does not decrease $|x|+|y|$, $|y_1| \geq |y|$ or $\left|\frac{x^2+x-3}{y}\right|\geq |y|$ or $|x^2+x-3| \geq |y|^2$. So we have that $|x^2+x-3| \geq y^2  \geq x^2-3$. If $|y| \leq |x|-1$ then we have $x^2-3 \leq y^2 \leq  x^2-2|x|+1$ which only holds for $|x| \leq 2$. If $|y| \geq |x|+1$ then we have $|x^2+x-3| \geq y^2 \geq x^2+2|x|+1$ which is only possible for $|x|\leq 2$. So we have the following cases to check, $|x|=|y|$, $y=0$ or $|x| \leq 2$. These cases give the integer solutions $(x,y)=(-2,-1),(2,-1),(2,3)$. 
Using transformations (a) and (b), starting with $(-2,-1)$, we obtain a chain of solutions
$$
 \dots \xleftrightarrow{\text{(a)}}  (2,3) \xleftrightarrow{\text{(b)}} (2,-1) \xleftrightarrow{\text{(a)}} (-2,-1) \xleftrightarrow{\text{(b)}} (-2,-1) \xleftrightarrow{\text{(a)}} (2,-1) \xleftrightarrow{\text{(b)}} (2,3) \xleftrightarrow{\text{(a)}} (-6,3) \xleftrightarrow{\text{(b)}} \dots 
$$
We can see that the solutions found previously, $(2,-1)$ and $(2,3)$ are in the same chain, therefore it contains all solutions. This chain is symmetric and so we may only look at the solutions on the right from $(-2,-1)$.  We can represent any solution $(x,y)$ to equation \eqref{vieta_8} using the transformations (ab) and (ba) $n$ times for some $n\geq0$. In this case, (ab) is $(x,y) \to (-x-y-1,-x-2y-1)$ and (ba) is $(x,y) \to (-2x+y-1,x-y)$. Therefore, the integer solutions to \eqref{vieta_8} can be written as $(x,y)=(x_n,y_n)$ where
$$
(x_0,y_0)=(-2,-1) \text{ and } (x_{n+1},y_{n+1})=(-x_n-y_n-1,-x_n-2y_n-1), \quad n=0,1,2,\ldots,
$$
or, $(x,y)=(x'_n,y'_n)$, where
$$
(x'_0,y'_0)=(-2,-1) \text{ and } (x'_{n+1},y'_{n+1})=(-2x'_n+y'_n-1,x'_n-y'_n), \quad n=0,1,2,\ldots.
$$

\vspace{10pt}

The next equation we will consider is
\begin{equation}\label{vieta_9}
x^2+x+xy-y^2+4=0.
\end{equation}
For this equation, the transformations (a) and (b) are $(x,y) \to (-x-y-1,y)$ and $(x,y) \to (x,x-y)$.
Let $(x,y)$ be a minimal solution of equation \eqref{vieta_9}. Assume that $x$ and $y$ are non-zero. Then, by Vieta jumping $(x_1,y)$ is also an integer solution where $x_1=-y-x-1=-\frac{y^2-4}{x}$. Because transformation (a) does not decrease $|x|+|y|$, we must have $|x_1| \geq |x|$ or $\left| \frac{y^2-4}{x}\right| \geq |x|$ or $|y^2-4|\geq |x|^2$. $|y^2-4| = y^2-4$ provided that $|y| \geq 2$, so $y^2-4 \geq x^2$ for $|y|\geq 2$. Now do Vieta jumping in $y$ and obtain the integer solution $(x,y_1)$ where $y_1=x-y=-\frac{x^2+x+4}{y}$. Because transformation (b) does not decrease $|x|+|y|$, $|y_1| \geq |y|$ or $\left|\frac{x^2+x+4}{y}\right|\geq |y|$ or $|x^2+x+4| \geq |y|^2$. So we have that $x^2+x+4 \geq y^2 \geq x^2+4$. 
If $|y| \leq |x|-2$ then we have $x^2+4 \leq y^2 \leq x^2-4|x|+4$, which is only true for $x=0$. If $|y| \geq |x|+2$ then $x^2+x+4 \geq y^2 \geq x^2+4|x|+4$ which is only true for $x=0$. Therefore, we have the following cases to check, $|y|=|x|+1$,  $|y|=|x|$, $|y|=|x|-1$, $x=0$, or $|y|\leq 1$. These cases give the integer solutions $(x,y)=(0,\pm 2),(1,-2),(-3,2)$.
Using transformations (a) and (b), starting with $(0,-2)$, we obtain a chain of solutions
$$
 \dots \xleftrightarrow{\text{(b)}} (7,-5) \xleftrightarrow{\text{(a)}} (-3,-5) \xleftrightarrow{\text{(b)}} (-3,2) \xleftrightarrow{\text{(a)}} (0,2) \xleftrightarrow{\text{(b)}} (0,-2) \xleftrightarrow{\text{(a)}} (1,-2) \xleftrightarrow{\text{(b)}} (1,3) \xleftrightarrow{\text{(a)}}  \dots 
$$

Because $(0,2),(1,-2)$ and $(-3,2)$ are also in this chain, it contains all integer solutions to \eqref{vieta_9}. We can represent any solution $(x,y)$ to equation \eqref{vieta_9} using the transformations (ab) and (ba) $n$ times for some $n\geq0$. In this case, (ab) is $(x,y) \to (-x-y-1,-x-2y-1)$ and (ba) is $(x,y) \to (-2x+y-1,x-y)$. We can create the left-hand side of the chain by applying (ba) to $(0, -2)$ and applying (ab) to $(0, 2)$. We can create the right-hand side of the chain by applying (ab) to $(0, -2)$ and applying (ba) to $(0, 2)$.

Therefore, the integer solutions to \eqref{vieta_9} can be written as $(x,y)=(x_n,y_n)$ where
$$
(x_0,y_0)=(0,2) \text{ and } (x_{n+1},y_{n+1})=(-x_n-y_n-1,-x_n-2y_n-1), \quad n=0,1,2,\ldots,
$$
or, $(x,y)=(x'_n,y'_n)$, where
$$
(x'_0,y'_0)=(0,2) \text{ and } (x'_{n+1},y'_{n+1})=(-2x'_n+y'_n-1,x'_n-y'_n), \quad n=0,1,2,\ldots,
$$
or,
$$
 (x''_0,y''_0)=(0,-2) \text{ and } (x''_{n+1},y''_{n+1})=(-x''_n-y''_n-1,-x''_n-2y''_n-1), \quad n=0,1,2,\ldots,
$$
or, $(x,y)=(x'''_n,y'''_n)$, where
$$
(x'''_0,y'''_0)=(0,-2) \text{ and } (x'''_{n+1},y'''_{n+1})=(-2x'''_n+y'''_n-1,x'''_n-y'''_n), \quad n=0,1,2,\ldots.
$$

\vspace{10pt}

The final equation we will consider is
\begin{equation}\label{vieta_10}
x^2+3x+xy-y^2=0.
\end{equation}
For this equation, the transformations (a) and (b) are $(x,y) \to (-x-y-3,y)$ and $(x,y) \to (x,x-y)$. 
Let $(x,y)$ be a minimal solution of equation \eqref{vieta_10}. Assume that $x$ and $y$ are non-zero. Then, by Vieta jumping $(x_1,y)$ is also an integer solution where $x_1=-y-x-3=-\frac{y^2}{x}$. Because transformation (a) does not decrease $|x|+|y|$, we must have $|x_1| \geq |x|$ or $\left| -\frac{y^2}{x}\right| \geq |x|$ or $y^2\geq x^2$. Now do Vieta jumping in $y$ and obtain the integer solution $(x,y_1)$ where $y_1=x-y=-\frac{x^2+3x}{y}$. Because transformation (b) does not decrease $|x|+|y|$, $|y_1| \geq |y|$ or $\left|\frac{x^2+3x}{y}\right|\geq |y|$ or $|x^2+3x| \geq |y|^2$, $|x^2+3x|=x^2+3x$ provided that $x\geq 0$ or $x \leq -3$. So we have that $x^2+3x \geq y^2 \geq x^2$. 
If $|y| \leq |x|-1$ then we have $x^2 \leq y^2 \leq x^2-2|x|+1$, which is only true for $x=0$. If $|y| \geq |x|+2$ then $x^2+3x \geq y^2 \geq x^2+4|x|+4$ which is impossible. Therefore, we have the following cases to check, $|y|=|x|+1$, $|y|=|x|$, $x=0,-1,-2$, or $y=0$. These cases give the integer solutions $(x,y)=(0,0),(-3,0),(-3,-3),(3,-3)$. Using transformations (a) and (b), starting with $(0,0)$, we obtain the chain of solutions
$$
 \dots \xleftrightarrow{\text{(a)}} (-3,-3) \xleftrightarrow{\text{(b)}} (-3,0) \xleftrightarrow{\text{(a)}} (0,0) \xleftrightarrow{\text{(b)}} (0,0) \xleftrightarrow{\text{(a)}} (-3,0) \xleftrightarrow{\text{(b)}} (-3,-3) \xleftrightarrow{\text{(a)}} (3,-3) \xleftrightarrow{\text{(b)}}  \dots 
$$
All solutions found previously are included in this chain, so it contains all solutions. This chain is symmetric, so we may only look at the solutions on the right from $(0,0)$. We can represent any solution $(x,y)$ to equation \eqref{vieta_10} using the transformations (ab) and (ba) $n$ times for some $n\geq0$. In this case, (ab) is $(x,y) \to (-x-y-3,-x-2y-3)$ and (ba) is $(x,y) \to (-2x+y-3,x-y)$. 
Therefore, the integer solutions to \eqref{vieta_10} can be written as $(x,y)=(x_n,y_n)$ where
$$
(x_0,y_0)=(0,0) \text{ and } (x_{n+1},y_{n+1})=(-x_n-y_n-3,-x_n-2y_n-3), \quad n=0,1,2,\ldots,
$$
or, $(x,y)=(x'_n,y'_n)$, where
$$
(x'_0,y'_0)=(0,0) \text{ and } (x'_{n+1},y'_{n+1})=(-2x'_n+y'_n-3,x'_n-y'_n), \quad n=0,1,2,\ldots
$$

Table \ref{table3.3} summarises the minimal integer solutions and transformations to use to find all solutions to the equations solved in this section.

\begin{center}
\begin{tabular}[c]{| c  | c |c|c|}

\hline
Equation & Minimal solutions $(x,y)$ &  Transformation (a) & Transformation (b) \\\hline \hline
$x^2+xy-y^2=0$&$(0,0)$&$(-y-x,y)$&$(x,x-y)$ \\\hline
$x^2+xy-y^2-1=0$&$(\pm1,0)$&$(-y-x,y)$&$(x,x-y)$ \\\hline
$x^2+xy-y^2+1=0$&$(0,\pm1)$&$(-y-x,y)$&$(x,x-y)$ \\\hline
$x^2+x+xy-y^2=0$&$(0,0),(1,-1)$&$(-x-y-1,y)$&$(x,x-y)$\\\hline
$x^2+x+xy-y^2+1=0$ & $(0,-1)$ &$(-x-y-1,y)$&$(x,x-y)$\\\hline
$x^2+xy-y^2-4=0$ & $(-2,0),(2,0)$&$(-x-y,y)$&$(x,x-y)$ \\\hline
$x^2+xy-y^2+4=0$ & $(-2,0),(2,0)$ &$(-x-y,y)$&$(x,x-y)$\\\hline
$x^2+x+xy-y^2-2=0$ &$(1,0),(1,1)$&$(-x-y-1,y)$&$(x,x-y)$ \\\hline
$x^2+x+xy-y^2+2=0$ &$(-1,1)$ &$(-x-y-1,y)$&$(x,x-y)$\\\hline
$x^2+xy-y^2-5=0$ & $(-2,-1),(2,1)$ &$(-x-y,y)$&$(x,x-y)$\\\hline
$x^2+xy-y^2+5=0$ &$(-1,2),(1,-2)$&$(-x-y,y)$&$(x,x-y)$ \\\hline
$x^2+x+xy-y^2-3=0$ &$(-2,-1)$&$(-x-y-1,y)$&$(x,x-y)$ \\\hline
$x^2+x+xy-y^2+4=0$ & $(0,-2),(0,2)$&$(-x-y-1,y)$&$(x,x-y)$\\\hline
$x^2+3x+xy-y^2=0$ &$(0,0)$&$(-x-y-3,y)$&$(x,x-y)$ \\\hline
\end{tabular}
\captionof{table}{\label{table3.3} Minimal solutions and the transformations to generate all integer solutions to the equations in Table \ref{tab:H18Vietavar2}.}
\end{center}

\subsection{Exercise 3.5} \label{ex:Pell}
\textbf{\emph{For integers $6 \leq d \leq 12$, $d\neq 9$, find the fundamental solution to equation 
\begin{equation}\label{eq:x2mdy2m1}
x^2 - dy^2 = 1.
\end{equation} 
Use it to describe the set of all its integer solutions.
}}

 Recall that the fundamental solution $(x,y)=(a,b)$ to \eqref{eq:x2mdy2m1} is its solution in positive integers, such that \eqref{eq:x2mdy2m1} has no integer solutions $(x,y)$ with $0 <y<b$.  To find it, we may try positive integer values of $y$ in increasing order, until we find a solution where $x$ is also integer. We will then write the complete solution set $(x,y)=(\pm x_n, \pm y_n)$ using the recursive formula 
 \begin{equation}\label{eq:x2mdy2m1sol}
 	(x_0,y_0)=(1,0), \quad \text{and} \quad (x_{n+1}, y_{n+1}) = (ax_n + bdy_n, bx_n + ay_n), \quad n=0,1,2,\ldots.
 \end{equation}

We will first look at equation \eqref{eq:x2mdy2m1} with $d=6$, which is,
\begin{equation} \label{pell6}
x^2-6y^2=1
\end{equation}
Substituting $y=1$ into \eqref{pell6} gives $x^2=7$, hence $x$ is not integer. Next, substituting $y=2$ into \eqref{pell6} gives $x^2=25$, therefore $x$ is integer and we have found the fundamental solution to equation \eqref{pell6} as $(a,b)=(5,2)$. Note that $(a,b)=(-5,2)$ is not the fundamental solution because $a,b$ must be positive. Then by \eqref{eq:x2mdy2m1sol} we have that the complete set of integer solutions to equation \eqref{pell6} is $(x,y)=(\pm x_n, \pm y_n)$ 
$$
 (x_0,y_0)=(1,0) \quad \text{and} \quad (x_{n+1},y_{n+1})= (5x_n+12y_n,2x_n+5y_n), \quad n=0,1,2,\dots.
$$

The remaining equations can be solved similarly, and Table \ref{table3.5} presents the equation, its fundamental solution and the recursive formula obtained by \eqref{eq:x2mdy2m1sol} that describe all the integer solutions.

\begin{center}
\begin{tabular}[c]{| c  | c |c|c|}
\hline
Equation & Fundamental solution $(a,b)$ &  Recursive formula $(x_{n+1},y_{n+1})$ \\
\hline \hline
$x^2-2y^2=1$ & $(3,2)$ & $(3x_n+4y_n,2x_n+3y_n)$ \\\hline
$x^2-3y^2=1$ & $(2,1)$ & $(2x_n+3y_n,x_n+2y_n)$ \\\hline
$x^2-5y^2=1$ & $(9,4)$ & $(9x_n+20y_n,4x_n+9y_n)$ \\\hline
$x^2-6y^2=1$ & $(5,2)$ & $(5x_n+12y_n,2x_n+5y_n)$ \\\hline
$x^2-7y^2=1$ & $(8,3)$ & $(8x_n+21y_n,3x_n+8y_n)$ \\\hline
$x^2-8y^2=1$  & $(3,1)$ & $(3x_n+8y_n,x_n+3y_n)$ \\\hline
$x^2-10y^2=1$ & $(19,6)$ & $(19x_n+60y_n,6x_n+19y_n)$ \\\hline
$x^2-11y^2=1$ & $(10,3)$ & $(10x_n+33y_n,3x_n+10y_n)$ \\\hline
$x^2-12y^2=1$ & $(7,2)$ & $(7x_n+24y_n,2x_n+7y_n)$ \\\hline
\end{tabular}
\captionof{table}{Fundamental solutions to equations of the form \eqref{eq:x2mdy2m1} for $2 \leq d \leq 12$, $d \neq 4,9$. The complete set of integer solutions to the equations are described by $(x,y)=(\pm x_n,\pm y_n)$, $n=0,1,2, \dots,$ where $(x_0,y_0)=(1,0)$ and the sequence $(x_n,y_n)$ is described by the final column. \label{table3.5}}
\end{center}

\subsection{Exercise 3.7}\label{ex:Pell2}
\textbf{\emph{Use the continued fraction expansion of $\sqrt{d}$ to compute the fundamental solution to equation \eqref{eq:x2mdy2m1} for all $14\leq d \leq 30$, $d\neq 16,25$. }}

This exercise is finding the fundamental solution to the equation \eqref{eq:x2mdy2m1} for $14 \leq d \leq 30$, $d \neq 16, 25$. Unlike in Section \ref{ex:Pell} we will be using continued fractions to find the fundamental solution, rather than trying values of $y$ in order.  We will use the iterative formula $x_{n+1}=\frac{1}{x_n-\floor*{x_n}}$, $n\geq 1$, with $x_0=\sqrt{d}$, so then, $\sqrt{d} = \floor*{x_0} + \frac{1}{\floor*{x_1}+\frac{1}{\floor*{x_2}+\frac{1}{\floor*{x_3}+\dots}}}$. 

Let us consider, for example, equation \eqref{eq:x2mdy2m1} with $d=14$, or specifically
\begin{equation} \label{pell14}
	x^2-14y^2=1.
\end{equation}
To find the fundamental solution, we will use rational approximations of $\sqrt{14}$. 
Let $x_0 = \floor*{\sqrt{14}} = 3$. This then gives the (trivial) approximation $\sqrt{14} \approx \frac{3}{1}$. We next check if $(x,y)=(3,1)$ is a solution to \eqref{pell14}, which it is not, so, we need to find a better approximation. We can do this by writing $\sqrt{14}=3+\frac{1}{x_1}$. We then have $\floor*{x_1}=1$, which gives $\sqrt{14} \approx 3+\frac{1}{1} = \frac{4}{1}$. We next check if $(x,y)=(4,1)$ is a solution to \eqref{pell14}, which it is not. Repeating this method, we obtain a sequence of rational approximations to $\sqrt{14}$ as $\frac{3}{1}, \frac{4}{1}, \frac{3}{2}, \frac{11}{3}, \frac{15}{4}, \ldots$. Direct substitution shows that $(x,y)=(15,4)$ is the first term of this sequence which is a solution to equation \eqref{pell14}. Hence, this must be the fundamental solution. 

The fundamental solutions to the equations with $14 < d \leq 30$, $d\neq 16,25$, can be found similarly. Table \ref{table3.7} presents the equation, the sequence of rational approximations found by the iterative process as above and the fundamental solution.

\begin{center}
	\begin{tabular}[c]{| c  | c |c|c|}
		\hline
		Equation & Rational Approximation of $\sqrt{d}$ & Fundamental solution $(x,y)$  \\
		\hline \hline
		$x^2-13y^2=1$ &  $\frac{3}{1}, \frac{4}{1}, \frac{7}{2}, \frac{11}{3}, \frac{18}{5},\frac{119}{33},\frac{137}{38}, \frac{256}{71}, \frac{393}{109},\frac{649}{180}, \ldots$ & $(649,180)$     \\\hline
		$x^2-14y^2=1$ &  $\frac{3}{1}, \frac{4}{1}, \frac{3}{2}, \frac{11}{3}, \frac{15}{4}, \ldots$ & $(15,4)$     \\\hline
		$x^2-15y^2=1$ &  $\frac{3}{1},\frac{4}{1}, \ldots$ & $(4,1)$    \\\hline
		$x^2-17y^2=1$ & $\frac{4}{1},\frac{33}{8}, \ldots$ & $(33,8)$    \\\hline
		$x^2-18y^2=1$ & $\frac{4}{1},\frac{17}{4}, \ldots$ & $(17,4)$   \\\hline
		$x^2-19y^2=1$ & $\frac{4}{1},\frac{9}{2},\frac{13}{3}, \frac{48}{11},\frac{61}{14},\frac{170}{39}, \ldots$ & $(170,39)$  \\\hline
		$x^2-20y^2=1$ & $\frac{4}{1},\frac{9}{2}, \ldots$ & $(9,2)$  \\\hline
		$x^2-21y^2=1$ & $\frac{4}{1},\frac{5}{1},\frac{9}{2}, \frac{23}{5},\frac{32}{7},\frac{55}{12}, \ldots$ & $(55,12)$   \\\hline
		$x^2-22y^2=1$ & $\frac{4}{1},\frac{5}{1},\frac{14}{3}, \frac{61}{13},\frac{136}{29},\frac{197}{42}, \ldots$ & $(197,42)$ \\\hline
		$x^2-23y^2=1$ & $\frac{4}{1},\frac{5}{1},\frac{19}{4}, \frac{24}{5}, \ldots$ & $(24,5)$ \\\hline
		$x^2-24y^2=1$ & $\frac{4}{1},\frac{5}{1}, \ldots$ & $(5,1)$ \\\hline
		$x^2-26y^2=1$ & $\frac{5}{1},\frac{51}{10}, \ldots$ & $(51,10)$ \\\hline
		$x^2-27y^2=1$ & $\frac{5}{1},\frac{26}{5}, \ldots$ & $(26,5)$ \\\hline
		$x^2-28y^2=1$ & $\frac{5}{1},\frac{16}{3},\frac{37}{7}, \frac{127}{24}, \ldots$ & $(127,24)$ \\\hline
		$x^2-29y^2=1$ & $\frac{5}{1},\frac{11}{2},\frac{16}{3}, \frac{27}{5},\frac{70}{13}, \frac{727}{135},$ $ \frac{1524}{283}, \frac{2251}{418}, \frac{3775}{701}, \frac{9801}{1820}, \ldots$ & $(9801,1820)$ \\\hline
		$x^2-30y^2=1$ & $\frac{5}{1},\frac{11}{2}, \ldots$ & $(11,2)$ \\\hline
	\end{tabular}
	\captionof{table}{Rational approximations for $\sqrt{d}$, used to find fundamental solutions for equations of the form \eqref{eq:x2mdy2m1}, for $13 \leq d \leq 30$, $d \neq 16,25$. \label{table3.7}}
\end{center}

\subsection{Exercise 3.11} \label{ex:genPell}
\textbf{\emph{For each of the equations listed in Table \ref{tab:H20Pell}, use Theorem \ref{th:x2mdy2mm} below (or reduction to a smaller equation of the same type) to describe all its integer solutions.
 \begin{theorem}\label{th:x2mdy2mm}[Theorem 3.9 in the book]
Let $m\neq 0,1$ be an integer, and let $d$ be a positive integer that is not a perfect square. Let $(a,b)$ be the fundamental solution of the equation $x^2-dy^2=1$. Then for each fundamental solution $(x^*,y^*)$ of the equation $x^2-dy^2=m$ where $d$ and $m$ are integer parameters, we have
$$
0 \leq y^* \leq \frac{b}{\sqrt{2(a+\epsilon)}}\sqrt{|m|}
$$
where $\epsilon=1$ if $m>0$ and $\epsilon = -1$ if $m<0$. 
\end{theorem}}}

\begin{center}
\begin{tabular}{ |c|c|c|c|c|c| } 
 \hline
 $H$ & Equation & $H$ & Equation & $H$ & Equation \\ 
 \hline\hline
 $13$ & $x^2-2y^2=1$ & $16$ & $x^2-2y^2=-4$ & $19$ & $x^2-3y^2=-3$ \\ 
 \hline
 $13$ & $x^2-2y^2=-1$ & $17$ & $x^2-3y^2=1$ & $20$ & $x^2-2y^2=8$ \\ 
 \hline
 $14$ & $x^2-2y^2=2$ & $18$ & $x^2-3y^2=-2$ & $20$ & $x^2-2y^2=-8$ \\ 
 \hline
 $14$ & $x^2-2y^2=-2$ & $19$ & $x^2-2y^2=7$ & $20$ & $x^2-3y^2=4$ \\ 
 \hline
 $16$ & $x^2-2y^2=4$ & $19$ & $x^2-2y^2=-7$ &  & \\ 
 \hline
\end{tabular}
\captionof{table}{\label{tab:H20Pell} General Pell's equations of size $H\leq 20$.}
\end{center} 

Recall that a solution $(x^*,y^*)$ to 
\begin{equation}\label{eq:x2mdy2mm} 
x^2-dy^2=m
\end{equation} 
 is called \textbf{fundamental} if (i) $y^*\geq 0$, (ii) there is no solution $(x,y)$ associated to $(x^*,y^*)$ with $y^*>y\geq 0$, and (iii) if $(-x^*,y^*)$ is associated to $(x^*,y^*)$, then $x^*\geq 0$.

The following proposition defines what it means for two solutions to be associated. 
\begin{proposition}\label{prop:Pellassociated}[Proposition 3.8 in the book]
Assume that $m\neq 0$, and $d>0$ is not a perfect square. Then two solutions $(x_1,y_1)$ and $(x_2,y_2)$ to equation $x^2-dy^2=m$, where $d$ and $m$ are integer parameters, are associated if and only if
$$
x=\frac{x_1x_2-dy_1y_2}{m} \quad \text{and} \quad y=\frac{x_1y_2-x_2y_1}{m}
$$
are both integers.
\end{proposition}

Let $F$ be the set of fundamental solutions to  \eqref{eq:x2mdy2mm}, then the set of all integer solutions $(x,y)$ is given by
\begin{equation}\label{eq:x2mdy2mmsol}
x + y \sqrt{d} = \pm (x^* + y^* \sqrt{d})(a + b\sqrt{d})^n \quad \text{for some} \quad n \in {\mathbb Z} \quad \text{and some} \quad (x^*,y^*) \in F.
\end{equation}
We will now prove that \eqref{eq:x2mdy2mmsol} can be written as $(x,y)=(\pm x_n,\pm y_n)$ where $(x_n,y_n)$ is given by, 
\begin{equation}\label{recursion_3.10}
 (x_0,y_0)=(x^*,y^*) \text{ and } (x_{n+1},y_{n+1})=(ax_n+bdy_n,bx_n+ay_n) \text{ with } n=0,1,2, \ldots. 
 \end{equation}
\begin{proof}
Starting with \eqref{eq:x2mdy2mmsol}, we will first look at the situation that $n \geq 0$.
For $n=0$, we have $(x,y)=\pm (x^*,y^*)$. Now, assume that $n > 0$. 
Then $x+y\sqrt{d} = \pm (ax^*+bdy^*  +\sqrt{d}(ay^*+bx^*))(a+b\sqrt{d})^{n-1}$. This can then be represented using the recursive formula, $(x,y)=\pm (x_n,y_n)$ where $(x_n,y_n)$, $n \geq 0$ is given by  
$$
(x_0,y_0)=(x^*,y^*) \text{ and } (x_{n+1},y_{n+1})=(ax_n+bdy_n,ay_n+bx_n).
$$

Next, we will look at \eqref{eq:x2mdy2mmsol} where $n$ is negative, i.e. $n < 0$. Let $k= |n|$. Then, 
$$
\begin{aligned}
x+\sqrt{d} = & \pm \frac{(x^*+y^*\sqrt{d})(a+b\sqrt{d})^k}{a+b\sqrt{d}}  = \pm \frac{(x^*+y^*\sqrt{d})(a+b\sqrt{d})^k(a-b\sqrt{d})}{(a+b\sqrt{d})(a-b\sqrt{d})} \\ = & \pm \frac{(ax^*-bdy^*  +\sqrt{d}(ay^*-bx^*))(a+b\sqrt{d})^k}{a^2-db^2}   =  \pm (ax^*-bdy^*  +\sqrt{d}(ay^*-bx^*))(a+b\sqrt{d})^k
\end{aligned}
$$
as $a^2-db^2=1$. This can be written recursively as $(x,y)=\pm (x'_n,y'_n)$ where $(x'_n,y'_n)$, $n \geq 0$ is given by  
$$
(x'_0,y'_0)=(x^*,y^*) \text{ and } (x'_{n+1},y'_{n+1})=(ax'_n-bdy'_n,ay_n-bx'_n).
$$
We must next check that $(x,y)=(-x^*,y^*)$ and $(x,y)=(x^*,-y^*)$ are solutions of \eqref{eq:x2mdy2mm}. By direct substitution, we can see that they are solutions. To summarise, we can therefore represent all integer solutions $(x,y)$ to equation \eqref{eq:x2mdy2mm} recursively using the following formula
$(x,y)=(\pm x_n, \pm y_n)$ where $(x_n,y_n)$, $n \geq 0$ is given by \eqref{recursion_3.10}.
\end{proof}

Equation
\begin{equation}\label{eq:x2m2y2m1}
x^2-2y^2=1
\end{equation}
is solved in Section 3.1.3 of the book, and its integer solutions are $(x,y)=(\pm x_n, \pm y_n)$, where
\begin{equation}\label{eq:x2m2y2m1sol}
(x_0,y_0)=(1,0), \quad \text{and} \quad (x_{n+1}, y_{n+1}) = (3x_n + 4y_n, 2x_n + 3y_n), \quad n=0,1,2,\ldots.
\end{equation}

\vspace{10pt}

Equation
\begin{equation}\label{eq:x2m2y2p1}
x^2-2y^2=-1
\end{equation}
is solved in Section 3.1.5 of the book, and its integer solutions are $(x,y)=(\pm x_n, \pm y_n)$ where
\begin{equation}\label{eq:x2m2y2p1sol}
 (x_0,y_0)=(1,1) \text{ and } (x_{n+1},y_{n+1})=(3x_n+4y_n,2x_n+3y_n), \quad n=0,1,2,\ldots.
\end{equation}

\vspace{10pt}

Equation
$$
x^2-2y^2=-2
$$
is solved in Section 3.1.5 of the book, and its integer solutions are $(x,y)=(\pm x_n, \pm y_n)$ where
$$
 (x_0,y_0)=(0,1) \text{ and } (x_{n+1},y_{n+1})=(3x_n+4y_n,2x_n+3y_n), \quad n=0,1,2,\ldots.
$$

\vspace{10pt}

Equation
\begin{equation}\label{eq:x2m3y2m1}
x^2-3y^2=1
\end{equation}
is solved in Section 3.1.3 of the book, and its integer solutions are $(x,y)=(\pm x_n, \pm y_n)$, where
\begin{equation}\label{eq:x2m3y2m1sol}
 (x_0,y_0)=(1,0) \text{ and } (x_{n+1},y_{n+1})=(2x_n+3y_n,x_n+2y_n), \quad n=0,1,2,\ldots.
\end{equation}

\vspace{10pt}

The first equation we will consider is
\begin{equation}\label{eq:x2m2y2m2}
	x^2-2y^2=2.
\end{equation}
We can see from this equation that $x^2$ is even, so, let $x=2x'$, substituting this into \eqref{eq:x2m2y2m2} and cancelling $2$, we obtain
$$
y^2-2(x')^2=-1,
$$
which up to the names of variables is equation \eqref{eq:x2m2y2p1}. By \eqref{eq:x2m2y2p1sol}, we know that the integer solutions to the equation $X^2-2Y^2=-1$ can be described as $(X,Y)=(\pm X_n, \pm Y_n)$ where
$$
	(X_0,Y_0)=(1,1) \text{ and } (X_{n+1},Y_{n+1})=(3X_n+4Y_n,2X_n+3Y_n), \quad n=0,1,2,\ldots.
$$
Hence, to obtain solutions to equation \eqref{eq:x2m2y2m2}, we must make the substitutions $X=y$ and $Y=x'=\frac{x}{2}$ or equivalently, $X_n=y_n$ and $2Y_n=x_n$. Therefore, we obtain the initial point $(x_0,y_0)=(2Y_0,X_0)=(2,1)$, with the recursive formula
$$
(x_{n+1},y_{n+1})=(2Y_{n+1},X_{n+1})=(2(2X_n+3Y_n),3X_n+4Y_n)=(4y_n+3x_n,3y_n+2x_n), \quad n=0,1,2,\ldots.
$$

To conclude, we may describe all integer solutions to equation \eqref{eq:x2m2y2m2} as $(x,y)=(\pm x_n, \pm y_n)$ where
\begin{equation}\label{eq:x2m2y2m2sol}
	(x_0,y_0)=(2,1) \text{ and } (x_{n+1},y_{n+1})=(3x_n+4y_n,2x_n+3y_n), \quad n=0,1,2,\ldots.
\end{equation}

\vspace{10pt}

The next equation we will consider is
\begin{equation}\label{eq:3.10_ni}
x^2-2y^2=4.
\end{equation}
We can see from this equation that $x^2$ is even, so, let $x=2x'$ for some integer $x'$. After substituting this into \eqref{eq:3.10_ni} and cancelling $2$, we obtain
 $$
 2(x')^2-y^2=2.
 $$ 
Therefore, $y^2$ must also be even, so, let $y=2y'$ for some integer $y'$. After substituting this into the previous equation and cancelling $2$, we obtain \eqref{eq:x2m2y2m1}, up to the names of variables, whose solution is given by \eqref{eq:x2m2y2m1sol}. Hence, in the original variables, we obtain that all integer solutions to \eqref{eq:3.10_ni} are 
$(x,y)=(\pm x_n, \pm y_n)$ where 
\begin{equation}\label{eq:3.10_nisol}
(x_0,y_0)=(2,0) \text{ and } (x_{n+1},y_{n+1})=(3x_n+4y_n,2x_n+3y_n), \quad n=0,1,2,\ldots.
\end{equation}

\vspace{10pt}

The next equation we will consider is
\begin{equation}\label{eq:3.10_nii}
x^2-2y^2=-4.
\end{equation}
Similar to the previous equation, we must have that both $x$ and $y$ are even. So let $x=2x'$ and $y=2y'$ for integers $x'$ and $y'$. After substituting these into \eqref{eq:3.10_nii} and cancelling $4$, we obtain \eqref{eq:x2m2y2p1}, up to the names of variables, whose solution is given by \eqref{eq:x2m2y2p1sol}. Hence, in the original variables, we obtain that the integer solutions to \eqref{eq:3.10_nii} are
$(x,y)=(\pm x_n, \pm y_n)$ where
\begin{equation}\label{eq:3.10_niisol}
(x_0,y_0)=(2,2) \text{ and } (x_{n+1},y_{n+1})=(3x_n+4y_n,2x_n+3y_n), \quad n=0,1,2,\ldots.
\end{equation}

\vspace{10pt}

The next equation we will consider is
\begin{equation}\label{eq:3.10_ii}
x^2-3y^2=-2.
\end{equation}
To solve this equation we will use Theorem \ref{th:x2mdy2mm}. The fundamental solution to the corresponding Pell's equation $x^2-3y^2=1$ is $(a,b)=(2,1)$. By Theorem \ref{th:x2mdy2mm}, we have $0 \leq y^* \leq \frac{1}{\sqrt{2(2-1)}}\sqrt{|2|}=1$. Hence, either $y^*=0$ or $y^*=1$. Substituting these values into \eqref{eq:3.10_ii} we find the integer solutions $(x,y)=(\pm1,1)$. By Proposition \ref{prop:Pellassociated}, we have that $(-1,1)$ and $(1,1)$ are associated. Therefore the fundamental solution to \eqref{eq:3.10_ii} is $(x_0,y_0)=(1,1)$. To describe all integer solutions to equation \eqref{eq:3.10_ii}, we can substitute these values into \eqref{recursion_3.10} and express this as the recursive formula $(x,y)=(\pm x_n, \pm y_n)$ where 
\begin{equation}\label{eq:x2m3y2p2sol}
(x_0,y_0)=(1,1) \text{ and } (x_{n+1},y_{n+1})=(2x_n+3y_n,x_n+2y_n), \quad n=0,1,2,\ldots.
\end{equation}

\vspace{10pt}

The next equation we will consider is
\begin{equation}\label{eq:x2m2y2m7}
x^2-2y^2=7.
\end{equation}
To solve this equation we will use Theorem \ref{th:x2mdy2mm}. The fundamental solution to the corresponding Pell's equation $x^2-2y^2=1$ is $(a,b)=(3,2)$. By Theorem \ref{th:x2mdy2mm}, we have $0 \leq y^* \leq \frac{2}{\sqrt{2(3+1)}}\sqrt{|7|}$. Hence, either $y^*=0$ or $y^*=1$. Substituting these values into \eqref{eq:x2m2y2m7} we find the integer solutions $(x,y)=(\pm 3,1)$. By Proposition \ref{prop:Pellassociated}, we have that $(-3,1)$ and $(3,1)$ are not associated. Therefore, we will have two initial points for the recursive solution. To describe all integer solutions to equation \eqref{eq:x2m2y2m7}, we can substitute these values into \eqref{recursion_3.10} and express this as the recursive formula $(x,y)=(\pm x_n, \pm y_n)$ where
\begin{equation}\label{eq:x2m2y2m7sol1}
(x_0,y_0)=(3,1) \text{ and } (x_{n+1},y_{n+1})=(3x_n+4y_n,2x_n+3y_n), \quad n=0,1,2,\ldots,
\end{equation}
and $(x,y)=(\pm x'_n, \pm y'_n)$ where, 
\begin{equation}\label{eq:x2m2y2m7sol2}
(x'_0,y'_0)=(-3,1) \text{ and } (x'_{n+1},y'_{n+1})=(3x'_n+4y'_n,2x'_n+3y'_n), \quad n=0,1,2,\ldots.
\end{equation}

\vspace{10pt}

The next equation we will consider is
\begin{equation}\label{eq:x2m2y2p7}
x^2-2y^2=-7.
\end{equation}
To solve this equation we will use Theorem \ref{th:x2mdy2mm}. The fundamental solution to the corresponding Pell's equation $x^2-2y^2=1$ is $(a,b)=(3,2)$. By Theorem \ref{th:x2mdy2mm}, we have $0 \leq y^* \leq \frac{2}{\sqrt{2(3-1)}}\sqrt{|-7|}$. Hence, either $y^*=0$, $y^*=1$ or $y^*=2$. Substituting these values into \eqref{eq:x2m2y2p7} we find the integer solutions $(x,y)=(\pm 1,2)$. By Proposition \ref{prop:Pellassociated}, we have that $(1,2)$ and $(-1,2)$ are not associated. Therefore, we will have two initial points for the recursive solution. To describe all integer solutions to equation \eqref{eq:x2m2y2p7}, we can substitute these values into \eqref{recursion_3.10} and express this as the recursive formula $(x,y)=(\pm x_n, \pm y_n)$ where
\begin{equation}\label{eq:x2m2y2p7sol1}
(x_0,y_0)=(1,2) \text{ and } (x_{n+1},y_{n+1})=(3x_n+4y_n,2x_n+3y_n), \quad n=0,1,2,\ldots,
\end{equation}
and $(x,y)=(\pm x'_n, \pm y'_n)$ where, 
\begin{equation}\label{eq:x2m2y2p7sol2}
(x'_0,y'_0)=(-1,2) \text{ and } (x'_{n+1},y'_{n+1})=(3x'_n+4y'_n,2x'_n+3y'_n), \quad n=0,1,2,\ldots.
\end{equation}

\vspace{10pt}

The next equation we will consider is
\begin{equation}\label{eq:3.10_iii}
x^2-3y^2=-3.
\end{equation}
We can see from this equation that $x^2$ is a multiple of 3, so, let $x=3x'$ for some integer $x'$. After substituting this into \eqref{eq:3.10_iii} and cancelling $3$, we obtain \eqref{eq:x2m3y2m1} up to the names of variables, whose solution is given by \eqref{eq:x2m3y2m1sol}. Hence, in the original variables, we can describe all integer solutions to \eqref{eq:3.10_iii} as
$(x,y)=(\pm x_n,\pm  y_n)$ where 
$$
(x_0,y_0)=(0,1) \text{ and } (x_{n+1},y_{n+1})=(2x_n+3y_n,x_n+2y_n), \quad n=0,1,2,\ldots.
$$

\vspace{10pt}

The next equation we will consider is
\begin{equation}\label{eq:3.10_niii}
x^2-2y^2=8.
\end{equation}
We can see from this equation that $x^2$ is even, so, let $x=2x'$ for some integer $x'$. After substituting this into \eqref{eq:3.10_niii} and cancelling $2$, we obtain \eqref{eq:3.10_nii}, up to the names of variables, whose solution is given by \eqref{eq:3.10_niisol}. Hence, in the original variables, we obtain that all integer solutions to \eqref{eq:3.10_niii} are
 $(x,y)=(\pm x_n, \pm y_n)$ where
$$
(x_0,y_0)=(4,2) \text{ and } (x_{n+1},y_{n+1})=(3x_n+4y_n,2x_n+3y_n), \quad n=0,1,2,\ldots.
$$

\vspace{10pt}

The next equation we will consider is
\begin{equation}\label{eq:3.10_niv}
x^2-2y^2=-8.
\end{equation}
We can see from this equation that $x^2$ is even, so, let $x=2x'$ for some integer $x'$. After substituting into the original equation and cancelling $2$, we obtain \eqref{eq:3.10_ni}, up to the names of variables, whose solution is given by \eqref{eq:3.10_nisol}. Hence, in the original variables, we obtain that all integer solutions to \eqref{eq:3.10_niv} are
$(x,y)=(\pm x_n, \pm y_n)$ where
$$
(x_0,y_0)=(0,2) \text{ and } (x_{n+1},y_{n+1})=(3x_n+4y_n,2x_n+3y_n), \quad n=0,1,2,\ldots.
$$

\vspace{10pt}

The final equation we will consider is
\begin{equation}\label{eq:3.10_vi}
x^2-3y^2=4.
\end{equation}
To solve this equation we will use Theorem \ref{th:x2mdy2mm}. The fundamental solution to the corresponding Pell's equation $x^2-3y^2=1$ is $(a,b)=(2,1)$. By Theorem \ref{th:x2mdy2mm}, we have $0 \leq y^* \leq \frac{1}{\sqrt{2(2+1)}}\sqrt{|4|}$. Hence $y^*=0$. Substituting this into \eqref{eq:3.10_vi} we find the integer solutions $(x,y)=(\pm 2,0)$. By Proposition \ref{prop:Pellassociated}, we have that $(2,0)$ and $(-2,0)$ are associated. Hence the fundamental solution is $(x^*,y^*)=(2,0)$. 
To describe all integer solutions to equation \eqref{eq:3.10_vi}, we can substitute these values into \eqref{recursion_3.10} and express this as the recursive formula $(x,y)=(\pm x_n, \pm y_n)$ where
$$
(x_0,y_0)=(2,0) \text{ and } (x_{n+1},y_{n+1})=(2x_n+3y_n,x_n+2y_n), \quad n=0,1,2,\ldots.
$$

 Table \ref{table3.10} summarises all integer solutions to the equations presented in Table \ref{tab:H20Pell} by listing their initial solutions and recursive formulas.

	\begin{center}
	\begin{tabular}[c]{|c| c  | c |}
		\hline
		Equation & $(x_0,y_0)$ & $(x_{n+1},y_{n+1})$   \\
		\hline \hline
		$x^2-2y^2=1$&$(1,0)$&$(3x_n +4y_n,2x_n +3y_n)$ \\\hline
		$x^2-2y^2=-1$&$(1,1)$&$(3x_n +4y_n,2x_n +3y_n)$\\\hline
		$x^2-2y^2=-2$&$(0,1)$ & $(3x_n+4y_n,2x_n+3y_n)$  \\\hline
		$x^2-2y^2=2$ & $(2,1)$ & $(3x_n+4y_n,2x_n+3y_n)$ \\\hline
		$x^2-2y^2=4$ & $(2,0)$ & $(3x_n+4y_n,2x_n+3y_n)$ \\\hline
		$x^2-2y^2=-4$ & $(2,2)$ & $(3x_n+4y_n,2x_n+3y_n)$ \\\hline
		$x^2-3y^2=1$&$(1,0)$&$(2x_n +3y_n,x_n +2y_n)$ \\\hline
		$x^2-3y^2=-2$ & $(1,1)$ & $(2x_n+3y_n,x_n+2y_n)$ \\\hline
		$x^2-2y^2=7$&$(-3,1)$&$(3x_n+4y_n,2x_n+3y_n)$\\
		& $(3,1)$&$(3x_n+4y_n,2x_n+3y_n)$\\\hline
		$x^2-2y^2=-7$&$(1,2)$& $(3x_n+4y_n,2x_n+3y_n)$\\
		&$(-1,2)$& $(3x_n+4y_n,2x_n+3y_n)$\\\hline
		$x^2-3y^2=-3$ & $(0,1)$ & $(2x_n+3y_n,x_n+2y_n)$ \\\hline
		$x^2-2y^2=8$&$(4,2)$& $(3x_n+4y_n,2x_n+3y_n)$\\\hline
		$x^2-2y^2=-8$&$(0,2)$& $(3x_n+4y_n,2x_n+3y_n)$\\\hline
		$x^2-3y^2=4$ & $(2,0)$ & $(2x_n+3y_n,x_n+2y_n)$ \\\hline
	\end{tabular}
	\captionof{table}{Integer solutions to the equations listed in Table \ref{tab:H20Pell}. All integer solutions can be described as $(x,y)=(\pm x_n,\pm y_n)$, with $n=0,1,2, \dots$ \label{table3.10}}
\end{center}

\subsection{Exercise 3.15}\label{ex:genquad2var}
\textbf{\emph{Solve all equations listed in Table \ref{tab:H17quad2var}.}}

\begin{center}
\begin{tabular}{ |c|c|c|c|c|c| } 
 \hline
 $H$ & Equation & $H$ & Equation & $H$ & Equation \\ 
 \hline\hline
 $14$ & $x^2+x-2y^2=0$ & $16$ & $2x^2+x-y^2+y=0$ & $17$ & $2x^2+x-y^2-3=0$ \\ 
 \hline
 $14$ & $2x^2+x-y^2=0$ & $16$ & $2x^2+2x-y^2=0$ & $17$ & $2x^2+x-y^2+3=0$ \\ 
 \hline
 $15$ & $2x^2+x-y^2-1=0$ & $16$ & $2x^2-y^2+y-2=0$ & $17$ & $2x^2+x-y^2+y-1=0$ \\ 
 \hline
 $15$ & $2x^2+x-y^2+1=0$ & $16$ & $2x^2-y^2+y+2=0$ &  $17$ & $2x^2+x-y^2+y+1=0$ \\ 
 \hline
 $16$ & $2x^2+x-y^2-2=0$ & $17$ & $2x^2+2x-y^2+1=0$ & & \\ 
 \hline
\end{tabular}
\captionof{table}{\label{tab:H17quad2var} Two-variable quadratic equations of size $H\leq 17$.}
\end{center}

Equation
\begin{equation}\label{eq:x2pxm2y2}
x^2+x-2y^2=0
\end{equation}
is solved in Section 3.1.6 of the book, and its integer solutions are $(x,y)=(x_n,\pm y_n)$ or $(x,y)=(x'_n,\pm y'_n)$,  where $(x_n,y_n)$ is given by
\begin{equation}\label{eq:x2pxm2y2sol1}
(x_0,y_0)=(0,0), \quad \text{and} \quad (x_{n+1}, y_{n+1}) = (3x_n + 4y_n+1, 2x_n + 3y_n+1), \quad n=0,1,2,\ldots,
\end{equation}
and $(x'_n,y'_n)$  is given by
\begin{equation}\label{eq:x2pxm2y2sol2}
(x'_0,y'_0)=(-1,0), \quad \text{and} \quad (x'_{n+1}, y'_{n+1}) = (3x'_n + 4y'_n+1, 2x'_n + 3y'_n+1), \quad n=0,1,2,\ldots
\end{equation}

\vspace{10pt}

Equation
$$
2x^2+x-y^2=0
$$
is solved in Section 3.1.6 of the book, and its integer solutions are $(x, y) = \left(\frac{x_{2n}}{2}, \pm y_{2n}\right)$, where $(x_n, y_n)$ is given by \eqref{eq:x2pxm2y2sol1}, or $(x, y) = (\frac{x_{2n+1}}{2}, \pm y_{2n+1})$, where $(x_n, y_n)$ is given by \eqref{eq:x2pxm2y2sol2}.
 
 \vspace{10pt}

Equation 
$$
 2x^2 + x - y^2 - 1 = 0
$$
is solved in Section 3.1.7 of the book, and its integer solutions are $(x,y)=\left(\frac{x_{2k+1}-1}{4} ,\pm \frac{y_{2k+1}}{2}\right)$ or $ \left(\frac{-x_{2k}-1}{4} ,\pm \frac{y_{2k}}{2}\right)$, where the sequence $(x_n, y_n)$, $n \geq 0$, is given by
\begin{equation}\label{2x2pxmy2m1redsol}
	(x_0,y_0) = (3,0), \text{ and } (x_{n+1},y_{n+1}) = (3x_n +4y_n,2x_n +3y_n), \quad n = 0,1,2,\ldots.
\end{equation}

\vspace{10pt}

Let us consider the equation 
\begin{equation}\label{eq:2x2pxmy2p1}
	2x^2+x-y^2+1=0.
\end{equation}
After multiplication by $8$ and rearranging, we obtain the general Pell's equation
\begin{equation}\label{eq:2x2pxmy2p1red}
	X^2 - 2Y^2 = -7,
\end{equation}  
with the restrictions that (i) $X=4x+1$ is $1$ modulo $4$ and (ii) $Y=2y$ is even. We next observe that in any integer solution to \eqref{eq:2x2pxmy2p1red} $X^2=2Y^2-7$ is odd, hence $X$ must be odd. This implies that $2Y^2=X^2+7$ is divisible by $4$, so $Y$ must be even, hence condition (ii) is not a restriction. 
Equation \eqref{eq:2x2pxmy2p1red} is \eqref{eq:x2m2y2p7}, up to the names of variables, and its integer solutions are 
$(X,Y)=(\pm x_n, \pm y_n)$ and $(X,Y)=(\pm x'_n, \pm y'_n)$, where $(x_n,y_n)$, $n\geq0$, is given by \eqref{eq:x2m2y2p7sol1}, and $(x'_n,y'_n)$, $n \geq 0$, is given by \eqref{eq:x2m2y2p7sol2}.
We must next investigate for which of the solutions we have that $X$ is $1$ modulo $4$. Recurrence relation $x_{n+1}=3x_n+4y_n$ implies that if $x_n \equiv 1 (\text{mod } 4)$ then $x_{n+1} \equiv 3 (\text{mod } 4)$ and vice versa. This together with the initial condition $x_0=1$ implies that $x_n \equiv 1 (\text{mod } 4)$ if and only if $n$ is even, that is, $n=2k$ for some integer $k$. Similarly, $-x_n \equiv 1 (\text{mod } 4)$ if and only if $n$ is odd, that is, $n=2k+1$ for some integer $k$. In the same way, we obtain that $x'_n \equiv 1 (\text{mod } 4)$ if and only if $n$ is odd, that is, $n=2k+1$ for some integer $k$ and $-x'_n \equiv 1 (\text{mod } 4)$ if and only if $n$ is even, that is, $n=2k$ for some integer $k$.

Finally, we conclude that the integer solutions to \eqref{eq:2x2pxmy2p1} are 
\begin{equation}\label{eq:2x2pxmy2p1finalsol1}
	(x,y)=\left(\frac{x_{2k}-1}{4}, \pm \frac{y_{2k}}{2}\right) \quad \text{or} \quad (x,y)=\left(\frac{-x_{2k+1}-1}{4}, \pm \frac{y_{2k+1}}{2}\right), \quad k=0,1,2,\ldots,
\end{equation}
where $(x_n,y_n)$ is given by \eqref{eq:x2m2y2p7sol1}, and
\begin{equation}\label{eq:2x2pxmy2p1finalsol}
	(x,y)=\left(\frac{x'_{2k+1}-1}{4}, \pm \frac{y'_{2k+1}}{2}\right) \quad \text{or} \quad (x,y)=\left(\frac{-x'_{2k}-1}{4}, \pm \frac{y'_{2k}}{2}\right), \quad k=0,1,2,\ldots,
\end{equation}
where $(x'_n,y'_n)$ is given by \eqref{eq:x2m2y2p7sol2}.

\vspace{10pt}

The next equation we will consider is 
\begin{equation}\label{eq:2x2pxmy2m2}
2x^2+x-y^2-2=0.
\end{equation}
After multiplication by $8$ and rearranging, we obtain the general Pell's equation
\begin{equation}\label{eq:2x2pxmy2m2red}
X^2 - 2Y^2 = 17,
\end{equation}  
with the restrictions that (i) $Y=2y$ is even and (ii) $X=4x+1$ is $1$ modulo $4$. We next observe that in any integer solution to \eqref{eq:2x2pxmy2m2red}, $X^2=2Y^2+17$ is odd, so $X$ must be odd. This implies that $2Y^2=X^2-17$ is divisible by $4$, hence $Y$ is even, and condition (i) is not a restriction. By the same method used in Section \ref{ex:genPell}
we obtain that the integer solutions to equation \eqref{eq:2x2pxmy2m2red} are $(X,Y)=(\pm x_n, \pm y_n)$ and $(X,Y)=(\pm x'_n, \pm y'_n)$, where 
\begin{equation}\label{eq:x2m2y2m17sol}
(x_0,y_0)=(-5,2), \quad \text{and} \quad (x_{n+1}, y_{n+1}) = (3x_n + 4y_n, 2x_n + 3y_n), \quad n=0,1,2,\ldots,
\end{equation}
and
\begin{equation}\label{eq:x2m2y2m17sol_i}
(x'_0,y'_0)=(5,2), \quad \text{and} \quad (x'_{n+1}, y'_{n+1}) = (3x'_n + 4y'_n, 2x'_n + 3y'_n), \quad n=0,1,2,\ldots
\end{equation}

We next investigate for which of the solutions we have that $X$ is $1$ modulo $4$. Recurrence relation $x_{n+1}=3x_n+4y_n$ implies that if $x_n \equiv 1 (\text{mod } 4)$ then $x_{n+1} \equiv 3 (\text{mod } 4)$ and vice versa. This together with the initial condition $x_0=-5$ implies that $x_n \equiv 1 (\text{mod } 4)$ if and only if $n$ is odd, that is, $n=2k+1$ for some integer $k$. Similarly, $-x_n \equiv 1 (\text{mod } 4)$ if and only if $n$ is even, that is, $n=2k$ for some integer $k$. So, the integer solutions to \eqref{eq:2x2pxmy2m2} are
\begin{equation}\label{eq:x2m2y2m17sol1}
(x,y)=\left(\frac{x_{2k+1}-1}{4}, \pm \frac{y_{2k+1}}{2}\right) \quad \text{or} \quad (x,y)=\left(\frac{-x_{2k}-1}{4}, \pm \frac{y_{2k}}{2}\right), \quad k=0,1,2,\ldots,
\end{equation}
where $(x_n,y_n)$ is given by \eqref{eq:x2m2y2m17sol}.

Recurrence relation $x'_{n+1}=3x'_n+4y'_n$ implies that if $x'_n \equiv 1 (\text{mod } 4)$ then $x'_{n+1} \equiv 3 (\text{mod } 4)$ and vice versa. This together with the initial condition $x'_0=5$ implies that $x'_n \equiv 1 (\text{mod } 4)$ if and only if $n$ is even, that is, $n=2k$ for some integer $k$. Similarly, $-x'_n \equiv 1 (\text{mod } 4)$ if and only if $n$ is odd, that is, $n=2k+1$ for some integer $k$. So, the integer solutions to \eqref{eq:2x2pxmy2m2} are 
\begin{equation}\label{eq:x2m2y2m17sol2}
(x,y)=\left(\frac{x'_{2k}-1}{4}, \pm \frac{y'_{2k}}{2}\right) \quad \text{or} \quad (x,y)=\left(\frac{-x'_{2k+1}-1}{4}, \pm \frac{y'_{2k+1}}{2}\right), \quad k=0,1,2,\ldots,
\end{equation}
where $(x'_n,y'_n)$ is given by \eqref{eq:x2m2y2m17sol_i}.

To summarise, all integer solutions to \eqref{eq:2x2pxmy2m2} are given by \eqref{eq:x2m2y2m17sol1} and \eqref{eq:x2m2y2m17sol2}, where the sequences $(x_n,y_n)$, $n\geq 0$, and $(x'_n,y'_n)$, $n\geq 0$, are given by \eqref{eq:x2m2y2m17sol} and \eqref{eq:x2m2y2m17sol_i}, respectively.

\vspace{10pt}

The next equation we will consider is
\begin{equation}\label{eq:2x2pxmy2py}
	2x^2+x-y^2+y=0.
\end{equation}
After multiplication by $8$ and rearranging we obtain the general Pell's equation
\begin{equation}\label{eq:2x2pxmy2pyred}
	X^2 - 2Y^2 = -1,
\end{equation}  
with the restrictions that (i) $Y=2y-1$ is odd and (ii) $X=4x+1$ is $1$ modulo $4$. 
We next observe that in any integer solution to \eqref{eq:2x2pxmy2pyred}, $X^2=2Y^2-1$ is odd hence $X$ must be odd. Also, $2Y^2=X^2+1$ is even but not divisible by $4$, hence $Y$ is odd, so condition (i) is not a restriction.
Equation \eqref{eq:2x2pxmy2pyred} is \eqref{eq:x2m2y2p1}, up to the names of variables, and its integer solutions are $(X,Y)=(\pm x_n, \pm y_n)$, where $(x_n,y_n)$ is given by \eqref{eq:x2m2y2p1sol}.
We next investigate for which of the solutions we have that $X$ is $1$ modulo $4$.
Recurrence relation $x_{n+1}=3x_n+4y_n$ implies that if $x_n \equiv 1 (\text{mod } 4)$ then $x_{n+1} \equiv 3 (\text{mod } 4)$ and vice versa. This together with the initial condition $x_0=1$ implies that $x_n \equiv 1 (\text{mod } 4)$ if and only if $n$ is even, that is, $n=2k$ for some integer $k$. Similarly, $-x_n \equiv 1 (\text{mod } 4)$ if and only if $n$ is odd, that is, $n=2k+1$ for some integer $k$. Hence, all integer solutions to \eqref{eq:2x2pxmy2py} are described by
\begin{equation}\label{eq:2x2pxmy2pysol_}
	(x,y)=\left(\frac{x_{2k}-1}{4},  \frac{ 1 \pm y_{2k}}{2}\right) \quad \text{or} \quad (x,y)=\left(\frac{-x_{2k+1}-1}{4},  \frac{ 1 \pm y_{2k+1}}{2}\right), \quad k=0,1,2,\ldots,
\end{equation}
where the sequence $(x_n,y_n)$, $n\geq 0$, is given by \eqref{eq:x2m2y2p1sol}.

\vspace{10pt}

The next equation we will consider is
\begin{equation}\label{eq:2x2p2xmy2}
2x^2+2x-y^2=0.
\end{equation}
We can easily see from this equation that $y^2$ is even, and so $y$ is even. Let $y=2Y$, then the equation is reduced to $x^2+x-2Y^2=0$, which is equation \eqref{eq:x2pxm2y2}, up to the names of variables. By \eqref{eq:x2pxm2y2sol1} and \eqref{eq:x2pxm2y2sol2}, its integer solutions are $(x,Y)=(x_n,\pm Y_n)$ or $(x,Y)=(x'_n,\pm Y'_n)$, where $(x_n,Y_n)$ is given by
$$
(x_0,Y_0)=(0,0), \quad \text{and} \quad (x_{n+1}, Y_{n+1}) = (3x_n + 4Y_n+1, 2x_n + 3Y_n+1), \quad n=0,1,2,\ldots,
$$
and $(x'_n,Y'_n)$ is given by
$$
(x'_0,Y'_0)=(-1,0), \quad \text{and} \quad (x'_{n+1}, Y'_{n+1}) = (3x'_n + 4Y'_n+1, 2x'_n + 3Y'_n+1), \quad n=0,1,2,\ldots.
$$
Therefore, we have that the integer solutions to equation \eqref{eq:2x2p2xmy2} are $(x,y)=(x_n, \pm y_n)$ or $(x,y)=(x'_n, \pm y'_n)$, where $(x_n,y_n), n \geq 0$ is given by
\begin{equation}\label{eq:2x2p2xmy2sol1}
(x_0,y_0)=(0,0), \quad \text{and} \quad (x_{n+1}, y_{n+1}) = (3x_n + 2y_n+1, 4x_n + 3y_n+2), \quad n=0,1,2,\ldots,
\end{equation}
and $(x'_n,y'_n), n \geq 0$, is given by
\begin{equation}\label{eq:2x2p2xmy2sol2}
(x'_0,y'_0)=(-1,0), \quad \text{and} \quad (x'_{n+1}, y'_{n+1}) = (3x'_n + 2y'_n+1, 4x'_n + 3y'_n+2), \quad n=0,1,2,\ldots
\end{equation}
To summarise, all integer solutions to \eqref{eq:2x2p2xmy2} are given by $(x,y)=(x_n,\pm y_n)$ and $(x,y)=(x'_n, \pm y'_n)$, where the sequences $(x_n,y_n)$ and $(x'_n,y'_n)$, $n\geq 0$, are given by \eqref{eq:2x2p2xmy2sol1} and \eqref{eq:2x2p2xmy2sol2}, respectively.

\vspace{10pt}

The next equation we will consider is
\begin{equation}\label{eq:2x2pymy2m2}
2x^2-y^2+y-2=0.
\end{equation}
After multiplication by $4$ and rearranging we obtain the general Pell's equation
\begin{equation}\label{eq:2x2pymy2m2red}
Y^2 - 2X^2 = -7,
\end{equation}  
with the restrictions that (i) $Y=2y-1$ is odd and (ii) $X=2x$ is even.
We next observe that in any integer solution to \eqref{eq:2x2pymy2m2red}, $Y^2=2X^2-7$ is odd hence $Y$ must be odd, so condition (i) is not a restriction. Up to the names of variables, \eqref{eq:2x2pymy2m2red} is \eqref{eq:x2m2y2p7} whose solution is $(Y,X)=(\pm x_n, \pm y_n)$, where $(x_n,y_n)$ is given by \eqref{eq:x2m2y2p7sol1}, and $(Y,X)=(\pm x'_n, \pm y'_n)$, where $(x'_n,y'_n)$ is given by \eqref{eq:x2m2y2p7sol2}.
We next investigate for which of the solutions we have that $X$ is $0$ modulo $2$.
Recurrence relation $y_{n+1}=2x_n+3y_n$ implies that if $y_n \equiv 0 (\text{mod } 2)$ then $y_{n+1} \equiv 0 (\text{mod } 2)$ and vice versa. This together with the initial condition $y_0=2$ implies that $y_n \equiv 0 (\text{mod } 2)$ for all $n$. 
Hence, all integer solutions to \eqref{eq:2x2pymy2m2} are given by
$$
(x,y)=\left(\pm \frac{y_{n}}{2}, \frac{1 \pm x_{n}}{2}\right) \text{ or } (x,y)=\left(\pm \frac{y'_{n}}{2}, \frac{1 \pm x'_{n}}{2}\right),
$$
where the sequence $(x_n,y_n)$ and $(x'_n,y'_n)$, $n\geq 0$, are given by \eqref{eq:x2m2y2p7sol1} and \eqref{eq:x2m2y2p7sol2}, respectively.
We may write this equivalently as
$(x,y)=(\pm x_n, y_n)$ or $(\pm x_n, 1-y_n)$, where $(x_n,y_n)$, $n\geq 0$, is given by
\begin{equation}\label{eq:2x2pymy2m2finalsol3}
	(x_0,y_0)=(1,0) \text{ and } (x_{n+1},y_{n+1})=(3x_n+2y_n-1,4x_n+3y_n-1),
\end{equation}
or $(x,y)=(\pm x'_n, y'_n)$, or $(\pm x'_n, 1-y'_n)$, where $(x'_n,y'_n)$, $n\geq 0$, is given by
\begin{equation}\label{eq:2x2pymy2m2finalsol}
	(x'_0,y'_0)=(1,1) \text{ and } (x'_{n+1},y'_{n+1})=(3x'_n+2y'_n-1,4x'_n+3y'_n-1).
\end{equation}

\vspace{10pt}

The next equation we will consider is 
\begin{equation}\label{eq:2x2my2pyp2}
2x^2-y^2+y+2=0. 
\end{equation} 
After multiplication by $8$ and rearranging
we obtain the general Pell's equation
\begin{equation}\label{eq:2x2my2pyp2red}
X^2 - 2Y^2 = 9,
\end{equation}  
with the restrictions that (i) $X=2y-1$ is odd and (ii) $Y=2x$ is even. 
We next observe that in any integer solution to \eqref{eq:2x2my2pyp2red}, $X^2=2Y^2+9$ is odd hence $X$ must be odd. Also, $2Y^2=X^2-9$ is divisible by 4, hence $Y$ is even, so condition (i) is not a restriction and condition (ii) is not a restriction.
We can find solutions to \eqref{eq:2x2my2pyp2red} by the same method in Section \ref{ex:genPell}, which gives solutions $(X,Y)=(\pm x_n, \pm y_n)$, where
\begin{equation}\label{eq:2x2my2pyp2sol}
(x_0,y_0)=(3,0) \text{ and } (x_{n+1},y_{n+1})=(3x_n+4y_n,2x_n+3y_n) \quad n=0,1,2,\ldots.
\end{equation}
Hence, all integer solutions to \eqref{eq:2x2my2pyp2} are described by
$$
(x,y)=\left(\pm \frac{y_{n}}{2},  \frac{ 1 \pm x_{n}}{2}\right)  \quad n=0,1,2,\ldots,
$$
where the sequence $(x_n,y_n)$, $n\geq 0$, is given by \eqref{eq:2x2my2pyp2sol}. This can be written equivalently as  $(x,y)=(\pm x_n,  y_n)$ or $(\pm x_n,1-y_n)$, where $(x_n,y_n)$ is given by
\begin{equation}\label{eq:2x2my2pyp2finalsol}
	(x_0,y_0)=(0,2) \text{ and } (x_{n+1},y_{n+1})=(3x_n+2y_n-1,4x_n+3y_n-1).
\end{equation}

\vspace{10pt}

The next equation we will consider is
\begin{equation}\label{eq:2x2p2xmy2p1}
	2x^2+2x-y^2+1=0. 
\end{equation}
After multiplication by $8$ and rearranging, we obtain the general Pell's equation
\begin{equation}\label{eq:2x2p2xmy2p1red}
	X^2 - 2Y^2 = 2,
\end{equation}  
with the restrictions that (i) $Y=2x+1$ is odd and (ii) $X=2y$ is even. 
We next observe that in any integer solution to \eqref{eq:2x2p2xmy2p1red}, $X^2=2Y^2+2$ is even hence $X$ must be even so condition (ii) is not a restriction.
Equation \eqref{eq:2x2p2xmy2p1red} is \eqref{eq:x2m2y2m2}, up to the names of variables, and the integer solutions to \eqref{eq:2x2p2xmy2p1red} are $(X,Y)=(\pm x_n, \pm y_n)$, where $(x_n,y_n)$ is given by \eqref{eq:x2m2y2m2sol}.
We next investigate for which of the solutions we have that $Y$ is $1$ modulo $2$.
Recurrence relation $y_{n+1}=2x_n+3y_n$ implies that if $y_n \equiv 1 (\text{mod } 2)$ then $y_{n+1} \equiv 1 (\text{mod } 2)$ and vice versa. This together with the initial condition $y_0=1$ implies that $y_n \equiv 1 (\text{mod } 2)$ for all $n$. So, the integer solutions to \eqref{eq:2x2p2xmy2p1} are
$$
(x,y)=\left(\frac{-1 \pm  y_{n}}{2},  \frac{  \pm x_{n}}{2}\right)  \quad n=0,1,2,\ldots,
$$
where the sequence $(x_n,y_n)$, $n\geq 0$, is given by \eqref{eq:x2m2y2m2sol}. This can be written equivalently as $(x,y)=(X_n, \pm Y_n)$ or $(x,y)=(-X_n-1,\pm Y_n)$, where $(X_n,Y_n)$ is given by
\begin{equation}\label{eq:2x2p2xmy2p1finalsol}
	(X_0,Y_0)=(0,1) \text{ and } (X_{n+1},Y_{n+1})=(3X_n+2Y_n+1,4X_n+3Y_n+2).
\end{equation}

\vspace{10pt}

The next equation we will consider is
\begin{equation}\label{eq:2x2pxmy2m3}
	2x^2+x-y^2-3=0. 
\end{equation}
After multiplication by $8$ and rearranging, we obtain the general Pell's equation
\begin{equation}\label{eq:2x2pxmy2m3red}
	X^2 - 2Y^2 = 25,
\end{equation}  
with the restrictions that (i) $Y=2y$ is even and (ii) $X=4x+1$ is 1 modulo 4. 
We next observe that in any integer solution to \eqref{eq:2x2pxmy2m3red}, $X^2=2Y^2+25$ is odd hence $X$ must be odd. Also, $2Y^2=X^2-25$ is divisible by 4, hence $Y$ is even, so condition (i) is not a restriction.
We can find solutions to \eqref{eq:2x2pxmy2m3red} by the same method from Section \ref{ex:genPell} which gives the integer solutions $(X,Y)=(\pm x_n, \pm y_n)$, where
\begin{equation}\label{eq:2x2pxmy2m3sol}
	(x_0,y_0)=(5,0) \text{ and } (x_{n+1},y_{n+1})=(3x_n+4y_n,2x_n+3y_n) \quad n=0,1,2,\ldots.
\end{equation}
We next investigate for which of the solutions we have that $X$ is $1$ modulo $4$.
Recurrence relation $x_{n+1}=3x_n+4y_n$ implies that if $x_n \equiv 1 (\text{mod } 4)$ then $x_{n+1} \equiv 3 (\text{mod } 4)$ and vice versa. This together with the initial condition $x_0=5$ implies that $x_n \equiv 1 (\text{mod } 4)$ if and only if $n$ is even, that is, $n=2k$ for some integer $k$. Similarly, $-x_n \equiv 1 (\text{mod } 4)$ if and only if $n$ is odd, that is, $n=2k+1$ for some integer $k$. Hence, all integer solutions to \eqref{eq:2x2pxmy2m3} are described by
\begin{equation}\label{eq:2x2pxmy2m3sol_}
	(x,y)=\left(\frac{x_{2k}-1}{4},  \frac{ \pm y_{2k}}{2}\right) \quad \text{or} \quad (x,y)=\left(\frac{-x_{2k+1}-1}{4},  \frac{ \pm y_{2k+1}}{2}\right), \quad k=0,1,2,\ldots,
\end{equation}
where the sequence $(x_n,y_n)$, $n\geq 0$, is given by \eqref{eq:2x2pxmy2m3sol}.

\vspace{10pt}

The next equation we will consider is
\begin{equation}\label{eq:2x2pxmy2p3}
	2x^2+x-y^2+3=0 .
\end{equation}
After multiplication by $8$ and rearranging, we obtain the general Pell's equation
\begin{equation}\label{eq:2x2pxmy2p3red}
	X^2 - 2Y^2 = -23,
\end{equation}  
with the restrictions that (i) $Y=2y$ is even and (ii) $X=4x+1$ is 1 modulo 4. 
We next observe that in any integer solution to \eqref{eq:2x2pxmy2p3red}, $X^2=2Y^2-23$ is odd hence, $X$ must be odd. Also, $2Y^2=X^2+23$ is divisible by $4$, hence $Y$ is even, so condition (i) is not a restriction.
We can find solutions to \eqref{eq:2x2pxmy2p3red} by the same method from Section \ref{ex:genPell} which gives the integer solutions $(X,Y)=(\pm x_n, \pm y_n)$ or $(X,Y)=(\pm x'_n, \pm y'_n)$, where $(x_n,y_n)$, $n \geq 0$, is given by
\begin{equation}\label{eq:2x2pxmy2p3sol}
	(x_0,y_0)=(-3,4) \text{ and } (x_{n+1},y_{n+1})=(3x_n+4y_n,2x_n+3y_n),
\end{equation}
and $(x'_n,y'_n)$, $n \geq 0$, is given by
\begin{equation}\label{eq:2x2pxmy2p3soli}
	(x'_0,y'_0)=(3,4) \text{ and } (x'_{n+1},y'_{n+1})=(3x'_n+4y'_n,2x'_n+3y'_n). 
\end{equation}
We next investigate for which of the solutions we have that $X$ is $1$ modulo $4$.
Recurrence relation $x_{n+1}=3x_n+4y_n$ implies that if $x_n \equiv 1 (\text{mod } 4)$ then $x_{n+1} \equiv 3 (\text{mod } 4)$ and vice versa. This together with the initial condition $x_0=-3$ implies that $x_n \equiv 1 (\text{mod } 4)$ if and only if $n$ is even, that is, $n=2k$ for some integer $k$. Similarly, $-x_n \equiv 1 (\text{mod } 4)$ if and only if $n$ is odd, that is, $n=2k+1$ for some integer $k$. So, the integer solutions to \eqref{eq:2x2pxmy2p3} are 
\begin{equation}\label{eq:2x2pxmy2p3finalsol}
	(x,y)=\left(\frac{x_{2k}-1}{4},  \frac{ \pm y_{2k}}{2}\right) \quad \text{or} \quad (x,y)=\left(\frac{-x_{2k+1}-1}{4},  \frac{ \pm y_{2k+1}}{2}\right), \quad k=0,1,2,\ldots,
\end{equation}
where $(x_n,y_n)$ is given by \eqref{eq:2x2pxmy2p3sol}.
Recurrence relation $x'_{n+1}=3x'_n+4y'_n$ implies that if $x'_n \equiv 1 (\text{mod } 4)$ then $x'_{n+1} \equiv 3 (\text{mod } 4)$ and vice versa. This together with the initial condition $x'_0=3$ implies that $x'_n \equiv 1 (\text{mod } 4)$ if and only if $n$ is odd, that is, $n=2k+1$ for some integer $k$. Similarly, $-x'_n \equiv 1 (\text{mod } 4)$ if and only if $n$ is even, that is, $n=2k$ for some integer $k$. So, the integer solutions to  \eqref{eq:2x2pxmy2p3} are
\begin{equation}\label{eq:2x2pxmy2p3finalsol2}
	(x,y)=\left(\frac{x'_{2k+1}-1}{4},  \frac{ \pm y'_{2k+1}}{2}\right) \quad \text{or} \quad (x,y)=\left(\frac{-x'_{2k}-1}{4},  \frac{\pm y'_{2k}}{2}\right), \quad k=0,1,2,\ldots,
\end{equation}
where $(x'_n,y'_n)$ is given by \eqref{eq:2x2pxmy2p3soli}.

To summarise, all integer solutions to  \eqref{eq:2x2pxmy2p3} are described by \eqref{eq:2x2pxmy2p3finalsol} and \eqref{eq:2x2pxmy2p3finalsol2} where $(x_n,y_n)$ and $(x'_n,y'_n)$, $n\geq 0$, are given by \eqref{eq:2x2pxmy2p3sol} and \eqref{eq:2x2pxmy2p3soli}, respectively.                                                           

\vspace{10pt}

The next equation we will consider is
\begin{equation}\label{eq:2x2pxpymy2m1}
2x^2+x-y^2+y-1=0.
\end{equation}
After multiplication by $8$ and rearranging, we obtain the general Pell's equation
\begin{equation}\label{eq:2x2pxpymy2m1red}
X^2 - 2Y^2 = 7,
\end{equation}  
with the restrictions that (i) $Y=2y-1$ is odd and (ii) $X=4x+1$ is 1 modulo 4. 
We next observe that in any integer solution to \eqref{eq:2x2pxpymy2m1red}, $X^2=2Y^2+7$ is odd hence $X$ must be odd. Also, $2Y^2=X^2-7$ is 2 modulo 4, hence $Y$ is odd, so the condition (i) is not a restriction.
Up to the names of variables, equation \eqref{eq:2x2pxpymy2m1red} is \eqref{eq:x2m2y2m7} and its integer solutions are $(X,Y)=(\pm x_n, \pm y_n)$ or $(X,Y)=(\pm x'_n, \pm y'_n)$, where $(x_n,y_n)$ is given by \eqref{eq:x2m2y2m7sol1} and $(x'_n,y'_n)$ is given by \eqref{eq:x2m2y2m7sol2}.
We next investigate for which of the solutions we have that $X$ is $1$ modulo $4$.
Recurrence relation $x_{n+1}=3x_n+4y_n$ implies that if $x_n \equiv 1 (\text{mod } 4)$ then $x_{n+1} \equiv 3 (\text{mod } 4)$ and vice versa. This together with the initial condition $x_0=3$ implies that $x_n \equiv 1 (\text{mod } 4)$ if and only if $n$ is odd, that is, $n=2k+1$ for some integer $k$. Similarly, $-x_n \equiv 1 (\text{mod } 4)$ if and only if $n$ is even, that is, $n=2k$ for some integer $k$. So, the integer solutions to \eqref{eq:2x2pxpymy2m1} are 
\begin{equation}\label{eq:2x2pxpymy2m1finalsol}
(x,y)=\left(\frac{x_{2k+1}-1}{4},  \frac{1 \pm y_{2k+1}}{2}\right) \quad \text{or} \quad (x,y)=\left(\frac{-x_{2k}-1}{4},  \frac{1 \pm y_{2k}}{2}\right), \quad k=0,1,2,\ldots,
\end{equation}
where $(x_n,y_n)$ is given by \eqref{eq:x2m2y2m7sol1}.
Recurrence relation $x'_{n+1}=3x'_n+4y'_n$ implies that if $x'_n \equiv 1 (\text{mod } 4)$ then $x'_{n+1} \equiv 3 (\text{mod } 4)$ and vice versa. This together with the initial condition $x'_0=-3$ implies that $x'_n \equiv 1 (\text{mod } 4)$ if and only if $n$ is even, that is, $n=2k$ for some integer $k$. Similarly, $-x'_n \equiv 1 (\text{mod } 4)$ if and only if $n$ is odd, that is, $n=2k+1$ for some integer $k$. So, the integer solutions to \eqref{eq:2x2pxpymy2m1} are
\begin{equation}\label{eq:2x2pxpymy2m1finalsol2}
	(x,y)=\left(\frac{x'_{2k}-1}{4},  \frac{1 \pm y'_{2k}}{2}\right) \quad \text{or} \quad (x,y)=\left(\frac{-x'_{2k+1}-1}{4},  \frac{ 1 \pm y'_{2k+1}}{2}\right), \quad k=0,1,2,\ldots,
\end{equation}
where $(x'_n,y'_n)$ is given by \eqref{eq:x2m2y2m7sol2}.

To summarise, all integer solutions to \eqref{eq:2x2pxpymy2m1} are given by  \eqref{eq:2x2pxpymy2m1finalsol} and \eqref{eq:2x2pxpymy2m1finalsol2} where the sequences $(x_n,y_n)$, $n \geq 0$, and $(x'_n,y'_n)$, $n \geq 0$, are given by \eqref{eq:x2m2y2m7sol1} and \eqref{eq:x2m2y2m7sol2}, respectively.

\vspace{10pt}

The next equation we will consider is
\begin{equation}\label{eq:2x2pxmy2pyp1}
2x^2+x-y^2+y+1=0. 
\end{equation}
After multiplication by $8$ and rearranging, we obtain the general Pell's equation
\begin{equation}\label{eq:2x2pxmy2pyp1red}
X^2 - 2Y^2 = -9,
\end{equation}  
with the restrictions that (i) $Y=2y-1$ is odd and (ii) $X=4x+1$ is $1$ modulo $4$. 
We next observe that in any integer solution to \eqref{eq:2x2pxmy2pyp1red}, $X^2=2Y^2-9$ is odd hence $X$ must be odd. Also, $2Y^2=X^2+9$ is even but not divisible by $4$, hence $Y$ is odd, so condition (i) is not a restriction.
We can find solutions to \eqref{eq:2x2pxmy2pyp1red} by the same method from Section \ref{ex:genPell} which gives the integer solutions $(X,Y)=(\pm x_n, \pm y_n),$ where
\begin{equation}\label{eq:2x2pxmy2pyp1sol}
(x_0,y_0)=(3,3) \text{ and } (x_{n+1},y_{n+1})=(3x_n+4y_n,2x_n+3y_n) \quad n=0,1,2,\dots.
\end{equation}
We next investigate for which of the solutions we have that $X$ is $1$ modulo $4$.
Recurrence relation $x_{n+1}=3x_n+4y_n$ implies that if $x_n \equiv 1 (\text{mod } 4)$ then $x_{n+1} \equiv 3 (\text{mod } 4)$ and vice versa. This together with the initial condition $x_0=3$ implies that $x_n \equiv 1 (\text{mod } 4)$ if and only if $n$ is odd, that is, $n=2k+1$ for some integer $k$. Similarly, $-x_n \equiv 1 (\text{mod } 4)$ if and only if $n$ is even, that is, $n=2k$ for some integer $k$. Hence, all integer solutions to \eqref{eq:2x2pxmy2pyp1} are described by
\begin{equation}\label{eq:2x2pxmy2pyp1sol_}
(x,y)=\left(\frac{-x_{2k}-1}{4},  \frac{ 1 \pm y_{2k}}{2}\right) \quad \text{or} \quad (x,y)=\left(\frac{x_{2k+1}-1}{4},  \frac{ 1 \pm y_{2k+1}}{2}\right), \quad k=0,1,2,\dots,
\end{equation}
where the sequence $(x_n,y_n)$, $n\geq 0$, is given by \eqref{eq:2x2pxmy2pyp1sol}.

 The integer solutions to the equations from Table \ref{tab:H17quad2var} are summarised in Table \ref{table3.14}. 
 
\begin{center}
\begin{tabular}[c]{| c  | c |c|c|}
		\hline
		Equation & Solution $(x,y)$, $n \geq 0$ \\\hline \hline
		$x^2+x-2y^2=0$& $(x_n, \pm y_n)$, where $(x_n,y_n)$ is given by either \eqref{eq:x2pxm2y2sol1} or \eqref{eq:x2pxm2y2sol2}. \\\hline
		$2x^2+x-y^2=0$&  $ (\frac{x_{2k}}{2}, \pm y_{2k})$, where $(x_n, y_n)$ is given by \eqref{eq:x2pxm2y2sol1} or \\ & $ (\frac{x_{2k+1}}{2}, \pm y_{2k+1})$, where $(x_n, y_n)$ is given by \eqref{eq:x2pxm2y2sol2}.  \\\hline
		$2x^2+x-y^2-1=0 $& $\left(\frac{x_{2k+1}-1}{4}, \pm \frac{y_{2k+1}}{2}\right)$ or $\left(\frac{-x_{2k}-1}{4}, \pm \frac{y_{2k}}{2}\right)$ \\ &  where $(x_n,y_n)$ is given by \eqref{2x2pxmy2m1redsol}. \\\hline
		$2x^2+x-y^2+1=0 $& \eqref{eq:2x2pxmy2p1finalsol1} and \eqref{eq:2x2pxmy2p1finalsol} where $(x_n,y_n)$ and $(x'_n,y'_n)$ are given by \\ & \eqref{eq:x2m2y2p7sol1}  and  \eqref{eq:x2m2y2p7sol2}, respectively.\\
\hline
		
	\end{tabular}
	\captionof{table}{Integer solutions to the equations listed in Table \ref{tab:H17quad2var}.}\label{table3.14}
\end{center}

\subsection{Exercise 3.26}\label{ex:Thue}
\textbf{\emph{(i) Use the {\tt Reduce} command in Mathematica to find all integer solutions to the equations 
\begin{equation}\label{2x3py3mm}
2x^3 + y^3 = m,
\end{equation}
\begin{equation}\label{x3px2ymy3mm}
x^3+x^2y-y^3 = m,
\end{equation}
and
\begin{equation}\label{x3px2ypy3mm}
x^3+x^2y+y^3 = m,
\end{equation}
with $m=1,2,3$.
\\
(ii) Investigate for which of these equations you can produce a direct proof (without computer assistance) that they have no other integer solutions. }}

\subsubsection{Exercise 3.26 (i)}

First, let us find all integer solutions to equations \eqref{2x3py3mm}, \eqref{x3px2ymy3mm} and \eqref{x3px2ypy3mm} where $m=1,2,3$. These equations are \emph{Thue equations}, as they are equations of the form 
$$
P(x,y)=m,
$$
where $P(x,y)$ is a homogeneous polynomial of degree $\geq 3$ and $m$ is integer. We can solve Thue equations in Mathematica using the code
$$
{\tt Reduce [P(x,y) == m, \{x,y\}, Integers] }
$$
The output lists all integer solutions to the equation. These lists are then presented in Table \ref{tab:3.25}. 

\begin{center}
\begin{tabular}{ |c|c| } 
 \hline
Equation & Integer solutions $(x,y)$ \\ 
 \hline\hline
$2x^3+y^3=1$ & $(0,1),(1,-1)$\\\hline
$2x^3+y^3=2$ & $(1,0)$\\\hline
$2x^3+y^3=3$ & $(1,1),(4,-5)$\\\hline
$x^3+x^2y-y^3=1$ & $(-3,-4),(0,-1),(1,-1),(1,0),(1,1) $ \\\hline
$x^3+x^2y-y^3=2$ & No integer solutions \\\hline
$x^3+x^2y-y^3=3$ & No integer solutions \\\hline
$x^3+x^2y+y^3=1$ & $(-1,1),(0,1),(1,0),(3,-2)$ \\\hline
$x^3+x^2y+y^3=2$ & No integer solutions \\\hline
$x^3+x^2y+y^3=3$ & $(1,1),(2,-1)$ \\\hline
\end{tabular}
\captionof{table}{\label{tab:3.25} Integer solutions to equations of the form \eqref{2x3py3mm}, \eqref{x3px2ymy3mm} and \eqref{x3px2ypy3mm} with $m=1,2,$ and $3$.}
\end{center} 

\subsubsection{Exercise 3.26 (ii)}

We will next provide a computer-free solutions of the first, fifth, sixth, and eighth equation from Table \ref{tab:3.25}. 

The first equation we will consider is
\begin{equation}\label{3.25i_m1}
2x^3+y^3=1.
\end{equation}
After substituting $x=-v$ and $y=u$ into \eqref{3.25i_m1} we obtain $u^3-2v^3=1$, which has been solved in Section 3.2.1 of the book, and its only integer solutions are $(u,v)=(-1,-1)$ and $(u,v)=(1,0)$. Hence, in the original variables, we obtain that the only integer solutions to \eqref{3.25i_m1} are $(x,y)=(1,-1)$ and $(x,y)=(0,1)$.

Next let us look at the equations
\begin{equation}\label{3.25ii_m2}
x^3+x^2y \pm y^3=2.
\end{equation}
We can see that the right-hand side is even. Let us look at the left-hand side of these equations modulo 2.

\begin{center}
	\begin{tabular}{ c| c c }
	$y / x$ & $0$ & $1$ \\\hline
		$0$ & 0 &1 \\
$1$ &1&1 \\   
	\end{tabular}
\captionof{table}{$x^3+x^2y \pm y^3$ modulo 2.\label{tab:3.25ii_m2}}	
\end{center}

We can see from Table \ref{tab:3.25ii_m2} that integer solutions only exist if $x$ and $y$ are both even. 
So, let $x=2x'$ and $y=2y'$ where $x'$ and $y'$ are some integers. 
After substituting into equation \eqref{3.25ii_m2} and cancelling $2$, we obtain $4((x')^3+(x')^2y' \pm (y')^3)=1$. This equation has no solutions, as the left-hand side is even, while the right-hand side is odd. Therefore, equations \eqref{3.25ii_m2} have no integer solutions.

\vspace{10pt}

Next, let us look at the equation
\begin{equation}\label{3.25ii_m3}
x^3+x^2y-y^3=3.
\end{equation}
We can see from this equation that the right-hand side is 0 modulo 3. Let us look at the left-hand side of the equation modulo 3.

\begin{center}
	\begin{tabular}{ c| c c c }
	$y / x$ & $0$ & $1$  & $2$ \\\hline
$0$ & 0 &1&2 \\
$1$ &2&1&2 \\
$2$ &1&1&2\\
	\end{tabular}
\captionof{table}{$x^3+x^2y-y^3$ modulo 3.\label{tab:3.25ii_m3}}	
\end{center}

We can see from Table \ref{tab:3.25ii_m3} that integer solutions only exist if $x \equiv 0 \text{ (mod 3)}$ and $y \equiv 0 \text{ (mod 3)}$. So, let $x=3x'$ and $y=3y'$ where $x'$ and $y'$ are some  integers. 
After substituting and cancelling $3$, we obtain

$9((x')^3+(x')^2y'-(y')^3)=1$. We can see that this equation has no integer solutions, as the left-hand side is divisible by $3$, while the right-hand is not. 
Therefore, equation \eqref{3.25ii_m3} has no integer solutions.

\subsection{Exercise 3.32}\label{ex:H16elliptic}
\textbf{\emph{Table \ref{tab:H16elliptic} lists equations defining elliptic curves in Weierstrass form, that is, equations of the form 
\begin{equation}\label{eq:Weiform}
y^2 + a x y + c y = x^3 + b x^2 + d x + e,
\end{equation} 
with $a,b,c,d,e$ integer coefficients, whose coefficients satisfy 
\begin{equation}\label{eq:Weiformcond}
-f^2 g^2 + 32 g^3 + f^3 h - 36 f g h + 108 h^2 \neq 0,
\end{equation} 
with $f=a^2+4b$, $g=ac+2d$ and $h=c^2+4e$.
For each equation (i) use the {\tt EllipticCurve} command in SageMath to find all its integer solutions, (ii) check that LMFDB lists the same integer solutions for these equations, and (iii) investigate for which equations you can produce a direct proof (without computer assistance and without referring to LMFDB) that they have no other integer solutions. }}

\begin{center}
\begin{tabular}{ |c|c|c|c|c|c| } 
 \hline
 $H$ & Equation & $H$ & Equation & $H$ & Equation \\ 
 \hline\hline
 $13$ & $y^2=x^3-1$ & $15$ & $y^2=x^3-x+1$ & $16$ & $y^2=x^3+x+2$ \\ 
 \hline
 $13$ & $y^2=x^3+1$ & $15$ & $y^2=x^3+x-1$ & $16$ & $y^2=x^3-2x$ \\ 
 \hline
 $14$ & $y^2=x^3-2$ & $15$ & $y^2=x^3+x+1$ & $16$ & $y^2=x^3+2x$ \\ 
 \hline
 $14$ & $y^2=x^3+2$ & $15$ & $y^2+y=x^3-1$ & $16$ & $y^2+y=x^3-2$ \\ 
 \hline
 $14$ & $y^2=x^3-x$ & $15$ & $y^2+y=x^3+1$ & $16$ & $y^2+y=x^3+2$ \\ 
 \hline
 $14$ & $y^2=x^3+x$ & $16$ & $y^2=x^3-4$ & $16$ & $y^2+y=x^3-x$ \\ 
 \hline
 $14$ & $y^2+y=x^3$ & $16$ & $y^2=x^3+4$ & $16$ & $y^2+y=x^3+x$ \\ 
 \hline
 $15$ & $y^2=x^3-3$ & $16$ & $y^2=x^3-x-2$ & & \\ 
 \hline
 $15$ & $y^2=x^3+3$ & $16$ & $y^2=x^3+x-2$ & & \\ 
 \hline
\end{tabular}
\captionof{table}{\label{tab:H16elliptic} Equations of size $H\leq 16$ defining elliptic curves in Weierstrass form.}
\end{center} 

\subsubsection{Exercise 3.32 (i)-(ii)}
Elliptic curves in Weierstrass form are equations of the form \eqref{eq:Weiform} satisfying \eqref{eq:Weiformcond}.
We can find all integer solutions to these equations using the SageMath command 
$$
{\tt EllipticCurve([a,b,c,d,e]).integral\_points()}
$$

However, it is important to note that this command does not output all integer solutions, it only outputs all values of $x$ for which the equation has integer solutions. If there are two integer solutions to the equation $(x,y_1)$ and $(x,y_2)$ with the same $x$, then only one of them is displayed, in the format $(x : y_1 : 1)$. Hence, for each output $(x : y : 1)$, we need to solve the equation in $x$, and find $y$. Let us consider, for example, the equation
\begin{equation}\label{eq:y2pymx3mx}
y^2+y=x^3+x.
\end{equation}
The SageMath command 
$$
	{\tt EllipticCurve([0,0,1,1,0]).integral\_points()}
$$
outputs 
$$
[(0 : 0 : 1), (1 : 1 : 1), (3 : 5 : 1)],
$$
which means that we should check values $x=0$, $x=1$ and $x=3$. For $x=0$, the equation reduces to $y^2+y=0$ and its integer solutions are $y=0$ and $y=-1$. For $x=1$, we obtain $y^2+y=2$, hence $y=1$ or $y=-2$.  Finally, for $x=3$, we obtain $y^2+y=30$, hence $y=5$ or $y=-6$. In conclusion, all integer solutions to \eqref{eq:y2pymx3mx} are
$$
(x,y)=(0, -1), (0,0), (1,-2), (1,1), (3,-6) \quad \text{and} \quad (3,5).
$$

The other equations in Table \ref{tab:H16elliptic} can be solved similarly. Table \ref{tab:3.31} presents all integer solutions to the equations listed in Table \ref{tab:H16elliptic}.

\begin{center}
\begin{tabular}{ |c|c||c|c| } 
		\hline
		Equation & Integer solutions $(x,y)$ & Equation & Integer solutions $(x,y)$ \\ 
		\hline\hline
		$y^2=x^3-1$& $(1,0)$ &$y^2=x^3-4$& $(2,\pm2),(5,\pm11)$ \\\hline
		$y^2=x^3+1$& $(-1,0),(0,\pm1),(2,\pm3)$ & $y^2=x^3+4$ &$(0,\pm2)$ \\\hline
		$y^2=x^3-2$& $(3,\pm 5)$ & $y^2=x^3-x-2$& $(2,\pm2)$ \\\hline
		$y^2=x^3+2$& $(-1,\pm1)$ & $y^2=x^3+x-2$& $(1,0)$ \\\hline
		$y^2=x^3-x$& $(\pm1,0),(0,0)$ & $y^2=x^3+x+2$& $(-1,0),(1,\pm2)$ \\\hline
		$y^2=x^3+x$& $(0,0)$ & $y^2=x^3-2x$& $(-1,\pm1),(0,0),(2,\pm2),$ \\ \cline{1-2} 
		$y^2+y=x^3$& $(0,-1),(0,0)$ && $(338,\pm6214) $\\\hline
		$y^2=x^3-3$& No integer solutions & $y^2=x^3+2x$& $(0,0)$\\\hline
		$y^2=x^3+3$& $(1,\pm2)$ & $y^2+y=x^3-2$& $(2,-3),(2,2)$\\\hline
		$y^2=x^3-x+1$ &$(0,\pm1),(\pm 1,\pm1),$ & $y^2+y=x^3+2$& $(0,-2),(0,1)$ \\\cline{3-4}
		& $(3,\pm5),(5,\pm11),(56,\pm419)$ & $y^2+y=x^3-x$& $(\pm 1,-1),(\pm 1,0),$ \\\cline{1-2}
		$y^2=x^3+x-1$& $(1,\pm1),(2,\pm3),(13,\pm47)$ & &$(0,-1),(0,0),(2,-3),$ \\\cline{1-2}
		$y^2=x^3+x+1$ &$(0,\pm1),(72,\pm611)$ & &$(2,2),(6,-15),(6,14)$ \\\hline 
		$y^2+y=x^3-1$& $(1,-1),(1,0),(7,-19),(7,18) $ & $y^2+y=x^3+x$& $(0,-1),(0,0),(1,-2),$\\\cline{1-2}
		$y^2+y=x^3+1$& $\pm (1,1),(-1,0),(1,-2),$& & $(1,1),(3,-6),(3,5) $\\
		&$(11,-37),(11,36)$&&\\\hline
	\end{tabular}
	\captionof{table}{\label{tab:3.31} Integer solutions to the equations in Table \ref{tab:H16elliptic}.}
\end{center} 

\subsubsection{Exercise 3.32 (iii)}
The proof that equation 
$$
y^2=x^3-3
$$ 
has no integer solutions can be found in Section 7.1.2 of the book.

\subsection{Exercise 3.33}\label{ex:H24Weideg}
\textbf{\emph{Solve all equations from Table \ref{tab:H24Weideg}. }}

\begin{center}
\begin{tabular}{ |c|c|c|c|c|c| } 
 \hline
 $H$ & Equation & $H$ & Equation & $H$ & Equation \\ 
 \hline\hline
 $16$ & $y^2=x^3-x^2$ & $20$ & $y^2=x^3-3x-2$ & $24$ & $y^2+3xy=x^3$ \\ 
 \hline
 $16$ & $y^2=x^3+x^2$ & $20$ & $y^2=x^3-3x+2$ & $24$ & $y^2+xy=x^3+2x^2$ \\ 
 \hline
 $16$ & $y^2+xy=x^3$ & $20$ & $y^2+xy=x^3-x^2$ & $24$ & $y^2+xy-y=x^3-2x-2$ \\ 
 \hline
 $19$ & $y^2=x^3-x^2-x+1$ & $20$ & $y^2+xy=x^3+x^2$ & & \\ 
 \hline
 $19$ & $y^2=x^3+x^2-x-1$ & $23$ & $y^2+xy-y=x^3+x^2-1$ & & \\ 
 \hline
\end{tabular}
\captionof{table}{\label{tab:H24Weideg} Degenerate Weiestrass equations of size $H\leq 24$.}
\end{center} 

Let us call an equation a Degenerate Weierstrass equation if it is of the form \eqref{eq:Weiform} and the condition \eqref{eq:Weiformcond} fails.

Equation
$$
y^2=x^3-x^2
$$
is solved in Section 3.2.5 of the book, and its integer solutions are
$$
(x,y)=(0,0), \quad \text{or} \quad (u^2+1,u^3+u) \quad \text{for some integer} \quad u.
$$

Equation
$$
y^2=x^3+x^2
$$
is solved in Section 3.2.5 of the book, and its integer solutions are
$$
(x,y)=(u^2-1,u^3-u) \quad \text{for some integer} \quad u.
$$

Equation
$$
y^2+xy=x^3
$$
is solved in Section 3.2.5 of the book, and its integer solutions are
$$
(x,y)=(u^2+u,u^3+u^2) \quad \text{for some integer} \quad u.
$$

The first equation we will consider is
\begin{equation}\label{y2ex3mx2mxp1}
y^2=x^3-x^2-x+1,
\end{equation}
which can be factorised as $y^2=(x-1)^2(x+1)$. If $y=0$ then $x=1$ or $x=-1$. Now assume that $y \neq 0$. This suggests that $x+1$ is a perfect square. So, let $x+1=u^2$ for some integer $u$, and we may assume that $u(u^2-2)$ has the same sign as $y$. By substituting $x=u^2-1$, we obtain $y^2=u^2(u^2-2)^2$. Hence
$$
y=u(u^2-2)=u^3-2u.
$$
Note that solution $(x,y)=(-1,0)$ is described by $(x,y)=(u^2-1,u^3-2u)$ with $u=0$. 
To conclude, the integer solutions to equation \eqref{y2ex3mx2mxp1} are
$$
(x,y)=(1,0) \quad \text{and} \quad (u^2-1,u^3-2u) \quad \text{for some integer} \quad u.
$$

\vspace{10pt}

The next equation we will consider is
\begin{equation}\label{y2ex3px2mxm1}
y^2=x^3+x^2-x-1,
\end{equation}
which can be factorised as $y^2=(x+1)^2(x-1)$. If $y=0$ then $x=1$ or $x=-1$. Now assume that $y \neq 0$. This suggests that $x-1$ is a perfect square, so let $x-1=u^2$ for some integer $u$, and we may assume that $u$ has the same sign as $y$.  By substituting $x=u^2+1$, we obtain $y^2=u^2(u^2+2)^2$. Hence
$$
y=u(u^2+2)=u^3+2u.
$$
Note that solution $(x,y)=(1,0)$ is described by $(x,y)=(u^2+1,u^3+2u)$ with $u=0$.
To conclude, the integer solutions to equation \eqref{y2ex3px2mxm1} are 
$$
(x,y)=(-1,0) \quad \text{and} \quad (u^2+1,u^3+2u)  \quad \text{for some integer} \quad u.
$$

\vspace{10pt}

The next equation we will consider is
\begin{equation}\label{y2ex3m3xm2}
y^2=x^3-3x-2,
\end{equation}
which can be factorised as $y^2=(x+1)^2(x-2)$. If $y=0$ then $x=-1$ or $x=2$. Now assume that $y \neq 0$. This suggests that $x-2$ is a perfect square, so let $x-2=u^2$ for some integer $u$, where $u$ has the same sign as $y$. By substituting $x=u^2+2$, we obtain $y^2=u^2(u^2+3)^2$. Hence,
$$
y=u(u^2+3)=u^3+3u.
$$
Note that solution $(x,y)=(2,0)$ is described by $(x,y)=(u^2+2,u^3+3u)$ with $u=0$.
To conclude, the integer solutions to equation \eqref{y2ex3m3xm2} are
$$
(x,y)=(-1,0) \quad \text{and} \quad (u^2+2,u^3+3u) \quad \text{for some integer} \quad u. 
$$

\vspace{10pt}

The next equation we will consider is
\begin{equation}\label{y2ex3m3xp2}
y^2=x^3-3x+2,
\end{equation}
which be factorised as $y^2=(x-1)^2(x+2)$. If $y=0$ then $x=-2$ or $x=1$. Now assume that $y \neq 0$. This suggests that $x+2$ is a perfect square, so let $x+2=u^2$ for some integer $u$, where $u(u^2-3)$ has the same sign as $y$. After substituting $x=u^2-2$, we obtain $y^2=u^2(u^2-3)^2$. Hence,
$$
y=u(u^2-3)=u^3-3u.
$$
Note that solution $(x,y)=(-2,0)$ is described by $(x,y)=(u^2-2,u^3-3u)$ with $u=0$.
To conclude, the integer solutions to equation \eqref{y2ex3m3xp2} are 
$$
(x,y)=(1,0) \quad \text{and} \quad (u^2-2,u^3-3u) \quad \text{for some integer} \quad u. 
$$

\vspace{10pt}

The next equation we will consider is
\begin{equation}\label{y2pxyex3mx2}
y^2+xy=x^3-x^2,
\end{equation}
which is equivalent to $(2y+x)^2=x^2(4x-3)$. If $2x+y=0$ then we must have $(x,y)=(0,0)$. Now assume that $2y+x \neq 0$. This suggests that $4x-3$ is the square of an odd integer. Let the integer be $2u+1$, where $2u+1$ has the same sign as $2y+x$. Then $4x-3=(2u+1)^2=4u^2+4u+1$, so $x=u^2+u+1$. Also, $(2y+x)^2=x^2(2u+1)^2$, hence, 
$$
2y+x=x(2u+1)=(u^2+u+1)(2u+1) = 2u^3+3u^2+3u+1.
$$ 
So, $y= u^3+u^2+u$. 
To conclude, the integer solutions to equation \eqref{y2pxyex3mx2} are
$$
(x,y)=(0,0) \quad \text{and} \quad (u^2+u+1,u^3+u^2+u) \quad \text{for some integer} \quad u.
$$

\vspace{10pt}

The next equation we will consider is
\begin{equation}\label{y2pxyex3px2}
y^2+xy=x^3+x^2,
\end{equation}
which is equivalent to $(2y+x)^2=x^2(4x+5)$. If $2y+x=0$ then we must have $(x,y)=(0,0)$. Now assume that $2y+x \neq 0$. This suggests that $4x+5$ is the square of an odd integer, let this integer be $2u+1$ where $x(2u+1)$ has the same sign as $2y+x$. Then, $4x+5=(2u+1)^2=4u^2+4u+1$, and so $x=u^2+u-1$. Also, $(2y+x)^2=x^2(2u+1)^2$, hence 
$$
2y+x=x(2u+1)=(u^2+u-1)(2u+1) = 2u^3+3u^2-u-1.
$$ 
So, $y= u^3+u^2-u$. 
To conclude, the integer solutions to equation \eqref{y2pxyex3px2} are
$$
(x,y)=(0,0) \quad \text{and} \quad (u^2+u-1,u^3+u^2-u) \quad \text{for some integer} \quad u. 
$$

\vspace{10pt}

The next equation we will consider is
\begin{equation}\label{y2pxymyex3px2m1}
y^2+xy-y=x^3+x^2-1,
\end{equation}
which is equivalent to $(2y+x-1)^2=(x+1)^2(4x-3)$. If $2y+x-1=0$ then $(x,y)=(-1,1)$. Now assume that $2y+x-1 \neq 0$. This suggests that $4x-3$ is the square of an odd integer, let this integer be $2u+1$ where $2u+1$ has the same sign as $2y+x-1$. Then, $4x-3=(2u+1)^2=4u^2+4u+1$, and so $x=u^2+u+1$. Also, $(2y+x-1)^2=(x+1)^2(2u+1)^2$, hence,
$$
2y+x-1=(x+1)(2u+1)=(u^2+u+2)(2u+1) = 2u^3+3u^2+5u+2.
$$ 
So, $y=u^3+u^2+2u+1$. 
To conclude, the integer solutions to equation \eqref{y2pxymyex3px2m1} are 
$$
(x,y)=(-1,1) \quad \text{and} \quad (u^2+u+1,u^3+u^2+2u+1) \quad \text{for some integer} \quad u. 
$$

\vspace{10pt}

The next equation we will consider is
\begin{equation}\label{y2p3xyex3}
y^2+3xy=x^3,
\end{equation}
which is equivalent to $(2y+3x)^2=x^2(4x+9)$. If $2y+3x=0$ then $(x,y)=(0,0)$. Now assume that $2y+3x \neq 0$. This suggests that $4x+9$ is the square of an odd integer, let this integer be $2u+1$, where $x(2u+1)$ has the same sign as $2y+3x$. Then, $4x+9=(2u+1)^2=4u^2+4u+1$, and so $x=u^2+u-2$. Also, $(2y+3x)^2=x^2(2u+1)^2$, hence, 
$$
2y+3x=x(2u+1) = (u^2+u-2)(2u+1) =2u^3+3u^2-3u-2.
$$ 
So, $y= u^3-3u+2$. 
To conclude, the integer solutions to equation \eqref{y2p3xyex3} are
$$
(x,y)=(u^2+u-2, u^3-3u+2) \quad \text{for some integer} \quad u. 
$$

\vspace{10pt}

The next equation we will consider is
\begin{equation}\label{y2pxyex3p2x2}
y^2+xy=x^3+2x^2,
\end{equation}
which is equivalent to $(2y+x)^2=x^2(4x+9)$. If $2y+x=0$ then $(x,y)=(0,0)$. Now assume that $2y+x \neq 0$. This suggests that $4x+9$ is the square of an odd integer, let this integer be $2u+1$ where $x(2u+1)$ has the same sign as $2y+x$. Then, $4x+9=(2u+1)^2=4u^2+4u+1$, and so $x=u^2+u-2$. Also, $(2y+x)^2=x^2(2u+1)^2$, hence,
$$
2y+x=x(2u+1) = (u^2+u-2)(2u+1) =2u^3+3u^2-3u-2.
$$
So, $y=u^3+u^2-2u$.
To conclude, the integer solutions to equation \eqref{y2pxyex3p2x2} are
$$
(x,y)=(u^2+u-2,u^3+u^2-2u) \quad \text{for some integer} \quad u.
$$ 

\vspace{10pt}

The final equation we will consider is
\begin{equation}\label{y2pxymyex3m2xm2}
y^2+xy-y=x^3-2x-2,
\end{equation}
which is equivalent to $(2y+x-1)^2=(x+1)^2(4x-7)$. If $2y+x-1=0$ then $(x,y)=(-1,1)$. Now assume that $2y+x-1\neq 0$. This suggests that $4x-7$ is the square of an odd integer, let this integer be $2u+1$ where $2u+1$ has the same sign as $2y+x-1$. Then, $4x-7=(2u+1)^2=4u^2+4u+1$, and so $x=u^2+u+2$. Also, $(2y+x-1)^2=(x+1)^2(2u+1)^2$, hence, 
$$
2y+x-1=(x+1)(2u+1)=(u^2+u+3)(2u+1) =2u^3+3u^2+7u+3.
$$ 
So, $y=u^3+u^2+3u+1$. 
To conclude, the integer solutions to equation \eqref{y2pxymyex3m2xm2} are 
$$
(x,y)=(-1,1) \quad \text{and} \quad (u^2+u+2,u^3+u^2+3u+1) \quad \text{for some integer} \quad u.
$$

Table \ref{Table:ex3.32} presents the integer solutions to the equations listed in Table \ref{tab:H24Weideg}.

\begin{center}
	\begin{tabular}{|c|c|c|}
		\hline
		Equation & Solution $(x,y)$ \\\hline \hline
		$y^2=x^3-x^2$&$(0,0), (u^2 +1,u^3 +u)$ \\\hline
		$y^2=x^3+x^2$&$ (u^2 -1,u^3 -u)$ \\\hline
		$y^2+xy=x^3$&$(u^2 + u, u^3 + u^2)$ \\\hline
		$y^2=x^3-x^2-x+1$&$(1,0),(u^2-1,u^3-2u)$ \\\hline
		$y^2=x^3+x^2-x-1$&$(-1,0),(u^2+1,u^3+2u)$ \\\hline
		$y^2=x^3-3x-2$&$(-1,0),(u^2+2,u^3+3u)$ \\\hline
		$y^2=x^3-3x+2$&$(1,0),(u^2-2,u^3-3u)$ \\\hline
		$y^2+xy=x^3-x^2$&$(0,0),(u^2+u+1,u^3+u^2+u)$ \\\hline
		$y^2+xy=x^3+x^2$&$(0,0),(u^2+u-1,u^3+u^2-u)$ \\\hline
		$y^2+xy-y=x^3+x^2-1$&$(-1,1),(u^2+u+1,u^3+u^2+2u+1)$ \\\hline
		$y^2+3xy=x^3$&$(u^2+u-2,u^3-3u+2)$ \\\hline
		$y^2+xy=x^3+2x^2$&$(u^2+u-2,u^3+u^2-2u)$ \\\hline
		$y^2+xy-y=x^3-2x-2$& $(-1,1),(u^2+u+2,u^3+u^2+3u+1)$\\\hline
	\end{tabular}
	\captionof{table}{Integer solutions to the equations in Table \ref{tab:H24Weideg}. Assume that $u$ is an arbitrary integer.}\label{Table:ex3.32}
\end{center}

\subsection{Exercise 3.35}\label{ex:H20almWei}
\textbf{\emph{Solve all equations from Table \ref{tab:H20almWei}. }}

\begin{center}
\begin{tabular}{ |c|c|c|c|c|c| } 
 \hline
 $H$ & Equation & $H$ & Equation & $H$ & Equation \\ 
 \hline\hline
 $17$ & $2y^2=x^3-1$ & $19$ & $2y^2+y=x^3-1$ & $20$ & $2y^2+y=x^3+2$ \\ 
 \hline
 $17$ & $2y^2=x^3+1$ & $19$ & $2y^2+y=x^3+1$ & $20$ & $2y^2+y=x^3-x$ \\ 
 \hline
 $18$ & $2y^2=x^3+2$ & $20$ & $2y^2=x^3-2x$ & $20$ & $2y^2+y=x^3+x$ \\ 
 \hline
 $18$ & $2y^2=x^3-x$ & $20$ & $2y^2=x^3+2x$ & $20$ & $2y^2+2y=x^3$ \\ 
 \hline
 $18$ & $2y^2=x^3+x$ & $20$ & $2y^2=x^3-x+2$ & $20$ & $2y^2+xy=x^3$ \\ 
 \hline
 $18$ & $2y^2+y=x^3$ & $20$ & $2y^2=x^3+x-2$ & $20$ & $2y^2=x^3-x^2$ \\ 
 \hline
 $19$ & $2y^2=x^3-3$ & $20$ & $2y^2=x^3+x+2$ & $20$ & $2y^2=x^3+x^2$ \\ 
 \hline
 $19$ & $2y^2=x^3+3$ & $20$ & $2y^2+y=x^3-2$ & &  \\ 
 \hline
\end{tabular}
\captionof{table}{\label{tab:H20almWei} Almost Weierstrass equations \eqref{eq:Weiformgen} with non-zero $f,g$ of size $H\leq 20$.}
\end{center} 

Let us call equations \emph{Almost Weierstrass equations} if they are of the form 
\begin{equation}\label{eq:Weiformgen}
fy^2 +axy+cy = gx^3 +bx^2 +dx+e \quad \text{where} \quad a,b,c,d,e,f,g \in \mathbb{Z}.
\end{equation} 
To solve equations of this type, we will be using a linear change of variables to put equations into Weierstrass form \eqref{eq:Weiform}. If the reduced equation is non-degenerate (that is, of the form \eqref{eq:Weiform} and the condition \eqref{eq:Weiformcond} holds), we can solve it using the method in Section \ref{ex:H16elliptic}, and if the reduced equation is degenerate (that is, of form \eqref{eq:Weiform} and the condition \eqref{eq:Weiformcond} fails), we can use the method in Section \ref{ex:H24Weideg}. 

Equation 
$$
2y^2=x^3-1
$$
is solved in Section 3.2.6 of the book and its integer solutions are
$$
(x,y)=(1,0).
$$

Equation 
$$
2y^2=x^3+1
$$
is solved in Section 3.2.6 of the book and its integer solutions are
$$
(x,y)=(-1,0),(1,\pm 1), (23,\pm 78).
$$

The first equation we will consider is
\begin{equation}\label{2y2ex3p2}
2y^2=x^3+2.
\end{equation}
After multiplication by $8$, we obtain
$$
Y^2=X^3+16,
$$ 
where $X=2x$ and $Y=4y$.
The SageMath command ${\tt EllipticCurve([0,0,0,0,16]).integral\_points();}$ returns that all integer solutions to the reduced equation are $(X,Y)=(0,\pm 4)$. Hence, in the original variables we obtain that the integer solutions to equation \eqref{2y2ex3p2} are
$$
(x,y)=(0,\pm 1).
$$

\vspace{10pt}

The next equation we will consider is
\begin{equation}\label{2y2ex3mx}
2y^2=x^3-x.
\end{equation}
After multiplication by $8$, we obtain
$$
Y^2=X^3-4X,
$$ 
where $X=2x$ and $Y=4y$.
The SageMath command ${\tt EllipticCurve([0,0,0,-4,0]).integral\_points();}$ returns that all integer solutions to the reduced equation are $(X,Y)=(0,0),(\pm 2,0)$. Hence, in the original variables we obtain that the integer solutions to
equation \eqref{2y2ex3mx} are
$$
(x,y)=(0,0),(\pm 1,0).
$$

\vspace{10pt}

The next equation we will consider is
\begin{equation}\label{2y2ex3px}
2y^2=x^3+x.
\end{equation}
After multiplication by $8$, we obtain
$$
Y^2=X^3+4X,
$$
where $X=2x$ and $Y=4y$. 
The SageMath command ${\tt EllipticCurve([0,0,0,4,0]).integral\_points();}$ returns that all integer solutions to the reduced equation are $(X,Y)=(0,0),(2,\pm 4)$. Hence, in the original variables we obtain that the integer solutions to equation
\eqref{2y2ex3px} are
$$
(x,y)=(0,0),(1, \pm 1).
$$

\vspace{10pt}

The next equation we will consider is
\begin{equation}\label{2y2pyex3}
2y^2+y=x^3.
\end{equation}
After multiplication by $8$, we obtain
$$
Y^2+2Y=X^3,
$$ 
where $X=2x$ and $Y=4y$.
The SageMath command ${\tt EllipticCurve([0,0,2,0,0]).integral\_points();}$ returns that all integer solutions to the reduced equation are $(X,Y)=(-1,-1),(0,-2),(0,0),(2,-4),(2,2)$. We now need to select solutions such that $X=2x$ is even and $Y=4y$ is divisible by $4$. Such solutions are $(X,Y)=(0,0)$ and $(2,-4)$. Hence, in the original variables we obtain that the integer solutions to equation
\eqref{2y2pyex3} are 
$$
(x,y)=(0,0),(1, -1).
$$

\vspace{10pt}

The next equation we will consider is
\begin{equation}\label{2y2ex3m3}
2y^2=x^3-3.
\end{equation}
After multiplication by 8, we obtain 
$$
Y^2=X^3-24,
$$ 
where $X=2x$ and $Y=4y$. The SageMath command ${\tt EllipticCurve([0,0,0,0,-24]).integral\_points();}$ returns that the reduced equation has no integer solutions and therefore the original equation \eqref{2y2ex3m3} has no integer solutions.

\vspace{10pt}

The next equation we will consider is
\begin{equation}\label{2y2ex3p3}
2y^2=x^3+3.
\end{equation}
After multiplication by $8$, we obtain
$$
Y^2=X^3+24,
$$ 
where $X=2x$ and $Y=4y$. The SageMath command
${\tt EllipticCurve([0,0,0,0,24]).integral\_points();}$ returns that
the integer solutions to the reduced equation are $(X,Y)=$ $(-2,\pm 4),$ $(1,\pm 5),$ $(10,\pm32),$ \\$(8158,\pm 736844).$ Hence, in the original variables we obtain that the integer solutions to 
equation \eqref{2y2ex3p3} are
$$
(x,y)=(-1, \pm 1),(5,\pm 8),(4079,\pm 184211).
$$

\vspace{10pt}

The next equation we will consider is
\begin{equation}\label{2y2pyex3m1}
2y^2+y=x^3-1.
\end{equation}
After multiplication by $8$, we obtain
$$
Y^2+2Y=X^3-8,
$$ 
where $X=2x$ and $Y=4y$. The SageMath command
${\tt EllipticCurve([0,0,2,0,-8]).integral\_points();}$ returns that
the integer solutions to the reduced equation are $(X,Y)=$ $(2,-2),$ $(2,0),$ $(32,-182),$ $(32,180)$. Hence, in the original variables we obtain that the integer solutions to
equation \eqref{2y2pyex3m1} are 
$$
(x,y)=(1,0),(16,45).
$$

\vspace{10pt}

The next equation we will consider is
\begin{equation}\label{2y2pyex3p1}
2y^2+y=x^3+1.
\end{equation}
After multiplication by $8$, we obtain
$$
Y^2+2Y=X^3+8,
$$ 
where $X=2x$ and $Y=4y$. The SageMath command
${\tt EllipticCurve([0,0,2,0,8]).integral\_points();}$ returns that the integer solutions to the reduced equation are
$(X,Y)=$ $(-2,-2),$ $(-2,0),$ $(0,-4),$ $(0,2),$ $(3,-7),$ $(3,5),$ $(6,-16),$ $(6,14),$ $(40,-254),$  $(40,252).$ Hence, in the original variables we obtain that the integer solutions to
equation \eqref{2y2pyex3p1} are
$$
(x,y)=(-1,0),(0,-1),(3,-4),(20,63).
$$

\vspace{10pt}

The next equation we will consider is
\begin{equation}\label{2y2ex3m2x}
2y^2=x^3-2x.
\end{equation}
After multiplication by $8$, we obtain
$$
Y^2=X^3-8X,
$$ 
where $X=2x$ and $Y=4y$. The SageMath command
${\tt EllipticCurve([0,0,0,-8,0]).integral\_points();}$ returns that the only integer solution to the reduced equation is
$(X,Y)=(0,0)$. Hence, in the original variables we obtain that the unique integer solution to
equation \eqref{2y2ex3m2x} is
$$
(x,y)=(0,0).
$$

\vspace{10pt}

The next equation we will consider is
\begin{equation}\label{2y2ex3p2x}
2y^2=x^3+2x.
\end{equation}
After multiplication by $8$, we obtain
$$
Y^2=X^3+8X,
$$ 
where $X=2x$ and $Y=4y$. The SageMath command
${\tt EllipticCurve([0,0,0,8,0]).integral\_points();}$ returns that the integer solutions to the reduced equation are
$(X,Y)=$ $(0,0),$ $(1,\pm 3),$ $(8,\pm 24).$ Hence, in the original variables we obtain that the integer solutions to
equation \eqref{2y2ex3p2x} are 
$$
(x,y)=(0,0),(4,\pm 6).
$$

\vspace{10pt}

The next equation we will consider is
\begin{equation}\label{2y2ex3mxp2}
2y^2=x^3-x+2.
\end{equation}
After multiplication by $8$, we obtain
$$
Y^2=X^3-4X+16,
$$ 
where $X=2x$ and $Y=4y$. The SageMath command
${\tt EllipticCurve([0,0,0,-4,16]).integral\_points();}$ returns that the integer solutions to the reduced equation are
$(X,Y)=$ $(-3,\pm 1),$ $(\pm 2,\pm 4),$ $(0,\pm 4),$ $(4, \pm 8),$ $(5, \pm 11),$ $(14, \pm 52),$ $(18, \pm 76),$ $(30, \pm 164),$ $(352,\pm 6604),$ $(992,\pm 31244)$. Hence, in the original variables we obtain that the integer solutions to 
equation \eqref{2y2ex3mxp2} are 
$$
(x,y)=(\pm 1, \pm 1),(0,\pm 1),(2,\pm 2),(7 \pm 13), (9,\pm 19), (15,\pm 41),(176, \pm 1651),(496,\pm 7811).
$$

\vspace{10pt}

The next equation we will consider is
\begin{equation}\label{2y2ex3pxm2}
2y^2=x^3+x-2.
\end{equation}
After multiplication by $8$, we obtain
$$
Y^2=X^3+4X-16
$$ 
where $X=2x$ and $Y=4y$. The SageMath command
${\tt EllipticCurve([0,0,0,4,-16]).integral\_points();}$ returns that the integer solutions to the reduced equation are
$(X,Y)=$ $(2,0),$ $(4,\pm 8),$ $(10,\pm 32)$. Hence, in the original variables we obtain that the integer solutions to 
equation \eqref{2y2ex3pxm2} are 
$$
(x,y)=(1,0),(2,\pm 2),(5, \pm8).
$$

\vspace{10pt}

The next equation we will consider is
\begin{equation}\label{2y2ex3pxp2}
2y^2=x^3+x+2.
\end{equation}
After multiplication by $8$, we obtain
$$
Y^2=X^3+4X+16,
$$ 
where $X=2x$ and $Y=4y$. The SageMath command
${\tt EllipticCurve([0,0,0,4,16]).integral\_points();}$ returns that the integer solutions to the reduced equation are
$(X,Y)=$ $(-2,0),$ $(0,\pm 4),$ $(6,\pm 16),$ $(1056,\pm 34316)$. Hence, in the original variables we obtain that the integer solutions to 
equation \eqref{2y2ex3pxp2} are
$$
(x,y)=(-1,0),(0,\pm 1),(3,\pm 4),(528, \pm8579).
$$

\vspace{10pt}

The next equation we will consider is
\begin{equation}\label{2y2pyex3m2}
2y^2+y=x^3-2.
\end{equation}
After multiplication by $8$, we obtain
$$
Y^2+2Y=X^3-16,
$$ 
where $X=2x$ and $Y=4y$. The SageMath command
${\tt EllipticCurve([0,0,2,0,-16]).integral\_points();}$ returns that the integer solutions to the reduced equation are
$(X,Y)=(4,-8),$ $(4,6)$. Hence, in the original variables we obtain that the integer solutions to 
equation \eqref{2y2pyex3m2} are
$$
(x,y)=(2,-2).
$$

\vspace{10pt}

The next equation we will consider is
\begin{equation}\label{2y2pyex3p2}
2y^2+y=x^3+2.
\end{equation}
After multiplication by $8$, we obtain 
$$
Y^2+2Y=X^3+16,
$$ 
where $X=2x$ and $Y=4y$. The SageMath command
${\tt EllipticCurve([0,0,2,0,16]).integral\_points();}$ returns that the integer solutions to the reduced equation are
$(X,Y)=$ $(-2,2),$ $(-1,-5),$ $(-1,3),$ $(2,-6),$ $\pm(2,4),$ $(4,-10),$ $(4,8),$ $(8,-24),$ $(8,22),$ $(43,-283),$ $(43,281),$ $(52,-376),$ $(52,374),$ $(5234,-378662),$ $(5234,378660)$. Hence, in the original variables we obtain that the integer solutions to 
equation \eqref{2y2pyex3p2} are
$$
(x,y)=\pm(1,1),(2,2),(4,-6),(26,-94),(2617,94665).
$$

\vspace{10pt}

The next equation we will consider is
\begin{equation}\label{2y2pyex3mx}
2y^2+y=x^3-x.
\end{equation}
After multiplication by $8$, we obtain
$$
Y^2+2Y=X^3-4X,
$$ 
where $X=2x$ and $Y=4y$. The SageMath command
${\tt EllipticCurve([0,0,2,-4,0]).integral\_points();}$ returns that the integer solutions to the reduced equation are
$(X,Y)=$ $(-1,-3),$ $(-1,1),$ $(0,-2),$ $(0,0),$ $(\pm 2,-2),$ $(\pm 2,0),$ $(3,-5),$ $(3,3),$ $(4,-8),$ $(4,6),$ $(10,-32),$ $(10,30),$ $(12,-42),$ $(12,40),$ $(20,-90),$ $(20,88),$ $(114,-1218),$ $(114,1216),$ $(1274,-45474),$ $(1274,45472).$ Hence, in the original variables we obtain that the integer solutions to  
equation \eqref{2y2pyex3mx} are 
$$
(x,y)=(0,0),(\pm 1,0),(2,-2),(5,-8),(6,10),(10,22),(57,304),(637,11368).
$$

\vspace{10pt}

The next equation we will consider is
\begin{equation}\label{2y2pyex3px}
2y^2+y=x^3+x.
\end{equation}
After multiplication by $8$, we obtain
$$
Y^2+2Y=X^3+4X,
$$ 
where $X=2x$ and $Y=4y$. The SageMath command
${\tt EllipticCurve([0,0,2,4,0]).integral\_points();}$ returns that the integer solutions to the reduced equation are
$(X,Y)=$ $(0,-2),$ $(0,0),$ $(4,-10),$ $(4,8)$. Hence, in the original variables we obtain that the integer solutions to 
equation \eqref{2y2pyex3px} are
$$
(x,y)=(0,0),(2,2).
$$

\vspace{10pt}

The next equation we will consider is
\begin{equation}\label{2y2p2yex3}
2y^2+2y=x^3.
\end{equation}
After multiplication by $8$, we obtain
$$
Y^2+4Y=X^3,
$$ 
where $X=2x$ and $Y=4y$. The SageMath command
${\tt EllipticCurve([0,0,4,0,0]).integral\_points();}$
$(X,Y)=(0,-4),$ $(0,0)$. Hence, in the original variables we obtain that the integer solutions to 
equation \eqref{2y2p2yex3} are
$$
(x,y)=(0,-1),(0,0).
$$

\vspace{10pt}

The next equation we will consider is
\begin{equation}\label{2y2pxyex3}
2y^2+xy=x^3.
\end{equation}
After  multiplication by $8$ and rearranging, we obtain
$(4y+x)^2=x^2(8x+1)$. If $4y+x=0$ then we must have $(x,y)=(0,0)$. Now assume that $4y+x\neq 0$. This implies that $8x+1$ is the square of an odd integer. Let the integer be $2k+1$, where $x(2k+1)$ has the same sign as $4y+x$. Then $8x+1=(2k+1)^2=4k^2+4k+1$, so, $x=\frac{k^2+k}{2}$, which is an integer for any integer $k$. Also, $(4y+x)^2=x^2(2k+1)^2$, hence,
$$
4y+x = x(2k+1)=\frac{(2k+1)(k^2+k)}{2}
$$ 
which implies that $y=\frac{k^3+k^2}{4}$. To ensure that $y$ is integer we must have $k \neq 1$ modulo $4$. Hence, $k=2u$ or $k=4u+3$ for some integer $u$. After substituting these into $(x,y)=\left(\frac{k^2+k}{2},\frac{k^3+k^2}{4}\right)$, we obtain that all integer solutions to equation \eqref{2y2pxyex3} are
$$
(x,y)=(2u^2+u,2u^3+u^2) \quad \text{or} \quad (8u^2+14u+6 ,16u^3+40u^2+33u+9) \quad \text{for some integer} \quad u.
$$

\vspace{10pt}

The next equation we will consider is
\begin{equation}\label{2y2ex3mx2}
2y^2=x^3-x^2.
\end{equation}
After  multiplication by $8$ and rearranging, we obtain
$(4y)^2=(2x)^2(2x-2)$. If $y=0$ then we have $(x,y)=(0,0)$ or $(1,0)$. Now assume that $y\neq 0$. This implies that $2x-2$ is the square of an even integer. Let the integer be $2u$, where $u$ has the same sign as $y$. Then $2x-2=(2u)^2$, and so $x=2u^2+1$. Also, $(4y)^2=(2x)^2(2u)^2$, hence,
$$
4y=(2x)(2u)=(4u^2+2)(2u)=8u^3+4u,
$$
which implies that $y=2u^3+u$. To conclude, all integer solutions to \eqref{2y2ex3mx2} are 
$$
 (x,y)=(0,0)\quad \text{and} \quad (2u^2+1,2u^3+u) \quad \text{for some integer} \quad u.
$$

\vspace{10pt}

The final equation we will consider is
\begin{equation}\label{2y2ex3px2}
2y^2=x^3+x^2.
\end{equation}
After  multiplication by $8$ and rearranging, we obtain
$(4y)^2=(2x)^2(2x+2)$. If $y=0$ then we have $(x,y)=(-1,0),(0,0)$. Now assume that $y\neq 0$. This implies that $2x+2$ is the square of an even integer. Let the integer be $2u$, where $x(2u)$ has the same sign as $y$. Then $2x+2=(2u)^2$, and so $x=2u^2-1$. Also, $(4y)^2=(2x)^2(2u)^2$, hence,
$$
4y=(2x)(2u)=(4u^2-2)(2u)=8u^3-4u,
$$
which implies that $y=2u^3-u$. To conclude, all integer solutions to \eqref{2y2ex3px2} are 
 $$
 (x,y)=(0,0) \quad \text{and} \quad (2u^2-1,2u^3-u) \quad \text{for some integer} \quad u.
$$

Table \ref{Table:ex3.34} summarises the integer solutions to the equations listed in Table \ref{tab:H20almWei}. 

\begin{center}
\begin{tabular}{|c|c|c|}
\hline
 Equation & Solution $(x,y)$ \\\hline \hline
$2y^2=x^3-1$&$(1,0)$ \\\hline
$2y^2=x^3+1$&$ (-1, 0), (1, \pm1), (23, \pm78)$ \\\hline
$2y^2=x^3+2$&$(0,\pm1)$ \\\hline
$2y^2=x^3-x$&$(0,0),(\pm 1,0)$ \\\hline
$2y^2=x^3+x$&$(0,0),(1, \pm 1)$ \\\hline
$2y^2+y=x^3$&$(0,0),(1,-1)$ \\\hline
$2y^2=x^3-3$&- \\\hline
$2y^2=x^3+3$&$(-1, \pm 1),(5,\pm 8),(4079,\pm 184211)$ \\\hline
$2y^2+y=x^3-1$&$(1,0),(16,45)$ \\\hline
$2y^2+y=x^3+1$&$(-1,0),(0,-1),(3,-4),(20,63)$ \\\hline
$2y^2=x^3-2x$&$(0,0)$ \\\hline
$2y^2=x^3+2x$&$(0,0),(4,\pm 6)$ \\\hline
$2y^2=x^3-x+2$&$(\pm 1, \pm 1),(0,\pm 1),(2,\pm 2),(7 \pm 13),$\\
&$(9,\pm 19), (15,\pm 41),(176, \pm 1651),(496,\pm 7811)$\\\hline
$2y^2=x^3+x-2$&$(1,0),(2,\pm 2),(5,\pm 8)$\\\hline
$2y^2=x^3+x+2$&$(-1,0),(0,\pm 1),(3,\pm 4),(528, \pm8579)$\\\hline
$2y^2+y=x^3-2$&$(2,-2)$ \\\hline
$2y^2+y=x^3+2$&$\pm(1,1),(2,2),(4,-6),(26,-94),(2617,94665)$ \\\hline
$2y^2+y=x^3-x$&$(0,0),(\pm1,0),(2,-2),(5,-8),(6,10),$ \\ & $(10,22),(57,304),(637,11368)$ \\\hline
$2y^2+y=x^3+x$&$(0,0),(2,2)$ \\\hline
$2y^2+2y=x^3$&$(0,-1),(0,0)$ \\\hline
$2y^2+xy=x^3$&$(2u^2+u,2u^3+u^2),$ \\ & $ (8u^2+14u+6 ,16u^3+40u^2+33u+9)$\\\hline
$2y^2=x^3-x^2$&$(0,0),(2u^2+1,2u^3+u)$ \\\hline
$2y^2=x^3+x^2$&$(0,0),(2u^2-1,2u^3-u)$ \\\hline
\end{tabular}
\captionof{table}{Integer solutions to the equations listed in Table \ref{tab:H20almWei}. Assume that $u$ is an arbitrary integer.}\label{Table:ex3.34}
\end{center}

\subsection{Exercise 3.36}\label{ex:H20almWei2}
\textbf{\emph{Solve all equations from Table \ref{tab:H20almWei2}. }}

\begin{center}
\begin{tabular}{ |c|c|c|c|c|c| } 
 \hline
 $H$ & Equation & $H$ & Equation & $H$ & Equation \\ 
 \hline\hline
 $18$ & $x y+3 y = x^3$ & $19$ & $2xy+y=x^3+1$ & $20$ & $xy+4y=x^3$ \\ 
 \hline
 $18$ & $2 x y+y = x^3$ & $20$ & $xy+3y=x^3-2$ & $20$ & $2xy+y=x^3-2$ \\ 
 \hline
 $19$ & $xy+3y=x^3-1$ & $20$ & $xy+3y=x^3+2$ & $20$ & $2xy+y=x^3+2$ \\ 
 \hline
 $19$ & $xy+3y=x^3+1$ & $20$ & $xy+3y=x^3-x$ & $20$ & $2xy+y=x^3-x$ \\ 
 \hline
 $19$ & $2xy+y=x^3-1$ & $20$ & $xy+3y=x^3+x$ & $20$ & $2xy+y=x^3+x$ \\ 
 \hline
\end{tabular}
\captionof{table}{\label{tab:H20almWei2} Almost Weierstrass equations \eqref{eq:Weiformgen} with $f=0$ of size $H\leq 20$.}
\end{center} 

This section solves equations of the form 
$$
axy+cy=gx^3+bx^2+dx+e
$$ 
where $a,b,c,d,e,g$ are integers.

Equation
$$
xy+3y=x^3
$$
is solved in Section 3.2.6 of the book, and its integer solutions are
$$
(x, y) = (-30, 1000), (-12, 192), (-6, 72), (-4, 64), (-2,-8), (0, 0), (6, 24) \quad \text{and} \quad (24, 512).
$$

Equation
$$
2xy+y=x^3
$$
is solved in Section 3.2.6 of the book, and its integer solutions are
$$
(x,y)=(-1,1) \quad \text{and} \quad (0,0).
$$

The first equation we will consider is 
\begin{equation}\label{xyp3yex3m1}
xy+3y=x^3-1,
\end{equation}
which can be factorised as $y(x+3)=x^3-1$. Let $z=x+3$, then $yz=(z-3)^3-1=z^3-9z^2+27z-28$. From this, we can see that $z$ must be a divisor of $28$, and so $z=\pm 28,$ $\pm 14,$ $\pm7,$ $\pm4,$ $\pm2$ or $\pm1$. We can find the integer solutions to \eqref{xyp3yex3m1} using  $x=z-3$ and $y=\frac{(z-3)^3-1}{z}$. To conclude, all integer solutions to equation \eqref{xyp3yex3m1} are 
$$
\begin{aligned}
(x,y)=& (-31,1064),(-17,351),(-10,143),(-7,86),(-5,63),(-4,65),\\ & (-2,-9), (-1,-1),(1,0),(4,9),(11,95),(25,558).
\end{aligned}
$$

\vspace{10pt}

The next equation we will consider is
\begin{equation}\label{xyp3yex3p1}
xy+3y=x^3+1,
\end{equation}
which can be factorised as $y(x+3)=x^3+1$. Let $z=x+3$, then $yz=(z-3)^3+1=z^3-9z^2+27z-26$. From this, we can see that $z$ must be a divisor of $26$, and so $z=\pm 26,$ $\pm 13,$ $\pm2$ or $\pm1$. We can now find the integer solutions to \eqref{xyp3yex3p1} using  $x=z-3$ and $y=\frac{(z-3)^3+1}{z}$. To conclude, all integer solutions to equation \eqref{xyp3yex3p1} are 
$$
(x,y)=(-29,938),(-16,315),(-5,62),(-4,63),(-2,-7),(-1,0),(10,77),(23,468).
$$

\vspace{10pt}

The next equation we will consider is 
\begin{equation}\label{2xypyex3m1}
2xy+y=x^3-1,
\end{equation}
which can be factorised as $y(2x+1)=x^3-1$. Let $z=2x+1$, then $yz=(\frac{z-1}{2})^3-1=\frac{z^3-3z^2+3z-9}{8}$. By multiplying this equation by 8, to get integer coefficients, we have $8yz=z^3-3z^2+3z-9$, from which we can see that $z$ must be a divisor of $9$, and so $z=\pm 9, \pm 3,\pm1$. For each value of $z$, we can find the integer solutions to \eqref{2xypyex3m1} using  $x=\frac{z-1}{2}$ and $y=\frac{(\frac{z-1}{2})^3-1}{z}$. To conclude, all integer solutions to equation \eqref{2xypyex3m1} are 
$$
(x,y)=(-5,14),(-2,3),(-1,2),(0,-1),(1,0),(4,7).
$$

\vspace{10pt}

The next equation we will consider is 
\begin{equation}\label{2xypyex3p1}
2xy+y=x^3+1,
\end{equation}
which can be factorised as $y(2x+1)=x^3+1$. Let $z=2x+1$, then $yz=(\frac{z-1}{2})^3+1=\frac{z^3-3z^2+3z+7}{8}$. By multiplying this equation by 8, to get integer coefficients, we have $8yz=z^3-3z^2+3z+7$, from which, we can see that $z$ must be a divisor of $7$, and so $z=\pm 7, \pm1$. We can now find the integer solutions to \eqref{2xypyex3p1} using  $x=\frac{z-1}{2}$ and $y=\frac{(\frac{z-1}{2})^3+1}{z}$. To conclude, all integer solutions to equation \eqref{2xypyex3p1} are 
$$
(x,y)=(-4,9),(-1,0),(0,1),(3,4).
$$

\vspace{10pt}

The next equation we will consider is 
\begin{equation}\label{xyp3yex3m2}
xy+3y=x^3-2,
\end{equation}
which can be factorised as $y(x+3)=x^3-2$. Let $z=x+3$, then $yz=(z-3)^3-2=z^3-9z^2+27z-29$. From this, we can see that $z$ must be a divisor of $29$. Then for each divisor, we can find the integer solutions to \eqref{xyp3yex3m2} using $x=z-3$ and $y=\frac{(z-3)^3-2}{z}$. To conclude, all integer solutions to equation \eqref{xyp3yex3m2} are 
$$
(x,y)=(-32,1130),(-4,66),(-2,10),(26,606).
$$

\vspace{10pt}

The next equation we will consider is 
\begin{equation}\label{xyp3yex3p2}
xy+3y=x^3+2,
\end{equation}
which can be factorised as $y(x+3)=x^3+2$. Let $z=x+3$, then $yz=(z-3)^3+2=z^3-9z^2+27z-25$. From this, we can see that $z$ must be a divisor of $25$. Then for each divisor, we can find the integer solutions to \eqref{xyp3yex3p2} using $x=z-3$ and $y=\frac{(z-3)^3+2}{z}$. To conclude, all integer solutions to equation \eqref{xyp3yex3p2} are 
$$
(x,y)=(-28,878),(-8,102),(-4,62),(-2,-6),(2,2),(22,426).
$$

\vspace{10pt}

The next equation we will consider is 
\begin{equation}\label{xyp3yex3mx}
xy+3y=x^3-x,
\end{equation}
which can be factorised as $y(x+3)=x^3-x$. Let $z=x+3$, then $yz=(z-3)^3-(z-3)=z^3-9z^2+26z-24$. From this, we can see that $z$ must be a divisor of $24$. Then for each divisor, we can find the integer solutions to \eqref{xyp3yex3mx} using $x=z-3$ and $y=\frac{z^3-9z^2+26z-24}{z}$. To conclude, all integer solutions to equation \eqref{xyp3yex3mx} are 
$$
\begin{aligned}
(x,y)=& (-27,819),(-15,280),(-11,165),(-9,120),(-7,84),(-6,70),(-5,60),\\ & (-4,60),(-2,-6), (0,0),(\pm 1,0),(3,4),(5,15),(9,60),(21,385).
\end{aligned}
$$

\vspace{10pt}

The next equation we will consider is 
\begin{equation}\label{xyp3yex3px}
xy+3y=x^3+x,
\end{equation}
which can be factorised as $y(x+3)=x^3+x$. Let $z=x+3$, then $yz=(z-3)^3+(z-3)=z^3-9z^2+28z-30$. From this, we can see that $z$ must be a divisor of $30$. Then for each divisor, we can find the integer solutions to \eqref{xyp3yex3px} using $x=z-3$ and $y=\frac{z^3-9z^2+28z-30}{z}$. To conclude, all integer solutions to equation \eqref{xyp3yex3px} are 
$$
\begin{aligned}
(x,y)=& (-33,1199),(-18,390),(-13,221),(-9,123),(-8,104),(-6,74),(-5,65),\\ & (-4,68),(-2,-10), (-1,-1),(0,0),(2,2),(3,5),(7,35),(12,116),(27,657).
\end{aligned}
$$

\vspace{10pt}

The next equation we will consider is 
\begin{equation}\label{xyp4yex3}
xy+4y=x^3,
\end{equation}
which can be factorised as $y(x+4)=x^3$. Let $z=x+4$, then $yz=(z-4)^3=z^3-12z^2+48z-64$. From this, we can see that $z$ must be a divisor of $64$. Then for each divisor, we can find the integer solutions to \eqref{xyp4yex3} using $x=z-4$ and $y=\frac{(z-4)^3}{z}$. To conclude, all integer solutions to equation \eqref{xyp4yex3} are 
$$
\begin{aligned}
(x,y)=& (-68,4913),(-36,1458), (-20,500),  (-12,216), (-8,128),(-6,108),(-5,125), \\ & (-3,-27), (-2,-4),  (0,0),(4,8), (12,108), (28,686), (60,3375).
\end{aligned}
$$

\vspace{10pt}

The next equation we will consider is 
\begin{equation}\label{2xypyex3m2}
2xy+y=x^3-2,
\end{equation}
which can be factorised as $y(2x+1)=x^3-2$. Let $z=2x+1$, then $yz=(\frac{z-1}{2})^3-2=\frac{z^3-3z^2+3z-17}{8}$. By multiplying this equation by 8, to get integer coefficients, we obtain $8yz=z^3-3z^2+3z-17$. From this, we can see that $z$ must be a divisor of $17$. Then for each divisor, we can find the integer solutions to \eqref{2xypyex3m2} using  $x=\frac{z-1}{2}$ and $y=\frac{(\frac{z-1}{2})^3-2}{z}$. To conclude, all integer solutions to equation \eqref{2xypyex3m2} are 
$$
(x,y)=(-9,43),(-1,3),(0,-2),(8,30).
$$

\vspace{10pt}

The next equation we will consider is 
\begin{equation}\label{2xypyex3p2}
2xy+y=x^3+2,
\end{equation}
which can be factorised as $y(2x+1)=x^3+2$. Letting $z=2x+1$, then $yz=(\frac{z-1}{2})^3+2=\frac{z^3-3z^2+3z+15}{8}$. By multiplying this equation by 8, to get integer coefficients, we obtain $8yz=z^3-3z^2+3z+15$. From this, we can see that $z$ must be a divisor of $15$. Then for each divisor, we can find the integer solutions to \eqref{2xypyex3p2} using  $x=\frac{z-1}{2}$ and $y=\frac{(\frac{z-1}{2})^3+2}{z}$. To conclude, all integer solutions to equation \eqref{2xypyex3p2} are 
$$
(x,y)=(-8,34),(-3,5),(0,2),\pm (1,1),(\pm 2,2),(7,23).
$$

\vspace{10pt}

The next equation we will consider is 
\begin{equation}\label{2xypyex3mx}
2xy+y=x^3-x,
\end{equation}
which can be factorised as $y(2x+1)=x^3-x$. Let $z=2x+1$, then $yz=(\frac{z-1}{2})^3-\frac{z-1}{2}
=\frac{z^3-3z^2-z+3}{8}$. By multiplying this equation by 8, to get integer coefficients, we obtain $8yz=z^3-3z^2-z+3$. From this, we can see that $z$ must be a divisor of $3$. For each divisor, we can find the integer solutions to \eqref{2xypyex3mx} using  $x=\frac{z-1}{2}$ and $y=\frac{(\frac{z-1}{2})^3-\frac{z-1}{2}}{z}$. To conclude, all integer solutions to equation \eqref{2xypyex3mx} are 
$$
(x,y)=(-2,2),(0,0),(\pm 1,0).
$$

\vspace{10pt}

The next equation we will consider is 
\begin{equation}\label{2xypyex3px}
2xy+y=x^3+x,
\end{equation}
which can be factorised as $y(2x+1)=x^3+x$. Let $z=2x+1$, then $yz=(\frac{z-1}{2})^3+\frac{z-1}{2}
=\frac{z^3-3z^2+7z-5}{8}$. By multiplying this equation by 8, to get integer coefficients, we obtain $8yz=z^3-3z^2+7z-5$. From this, we can see that $z$ must be a divisor of $5$. For each divisor, we can find the integer solutions to \eqref{2xypyex3px} using $x=\frac{z-1}{2}$ and $y=\frac{(\frac{z-1}{2})^3+\frac{z-1}{2}}{z}$. To conclude, all integer solutions to equation \eqref{2xypyex3px} are 
$$
(x,y)=(-3,6),(-1,2),(0,0),(2,2).
$$

Table \ref{Table:ex3.35} summarises all integer solutions to the equations listed in Table \ref{tab:H20almWei2}. 

\begin{center}
\begin{tabular}{|c|c|c|}
\hline
Equation & Solution $(x,y)$ \\\hline \hline
$xy+3y=x^3$&$(-30, 1000), (-12, 192), (-6, 72),  $\\
&$(-4, 64),(-2, -8), (0, 0), (6, 24), (24, 512)$  \\\hline
$2xy+y=x^3$& $(-1, 1), (0, 0)$ \\\hline
$xy+3y=x^3-1$&$(-31,1064),(-17,351),(-10,143),$\\
&$(-7,86),(-5,63),(-4,65),(-2,-9),$\\
&$(-1,-1),(1,0),(4,9),(11,95),(25,558)$ \\\hline
$xy+3y=x^3+1$&$(-29,938),(-16,315),(-5,62),(-4,63),$\\
&$(-2,-7),(-1,0),(10,77),(23,468)$ \\\hline
$2xy+y=x^3-1$&$(-5,14),(-2,3),(-1,2),(0,-1),$\\
&$(1,0),(4,7)$ \\\hline
$2xy+y=x^3+1$&$(-4,9),(-1,0),(0,1),(3,4)$ \\\hline
$xy+3y=x^3-2$&$(-32,1130),(-4,66),(-2,10),(26,606)$ \\\hline
$xy+3y=x^3+2$&$(-28,878),(-8,102),(-4,62),$\\ & $(-2,-6),(2,2),(22,426)$\\\hline
$xy+3y=x^3-x$&$(-27,819),(-15,280),(-11,165),$\\
&$(-9,120),(-7,84),(-6,70),(-5,60),$\\
&$(-4,60),(-2,-6),(0,0),(\pm 1,0),$\\
&$(3,4),(5,15),(9,60),(21,385)$\\\hline
$xy+3y=x^3+x$&$(-33,1199),(-18,390),(-13,221),$\\
&$(-9,123),(-8,104),(-6,74),(-5,65),$\\
&$(-4,68),(-2,-10),(-1,-1),(0,0),$\\
&$(2,2),(3,5),(7,35),(12,116),(27,657)$\\\hline
$xy+4y=x^3$&$(-68,4913),(-36,1458),(-20,500), $\\
&$(-12,216),(-8,128),(-6,108), $\\
&$(-5,125), (-3,-27),(-2,-4),(0,0),   $\\
&$(4,8), (12,108),(28,686),(60,3375) $ \\\hline
$2xy+y=x^3-2$&$(-9,43),(-1,3),(0,-2),(8,30)$\\\hline
$2xy+y=x^3+2$&$(-8,34),(-3,5),(0,2),$\\
&$\pm(1,1),(\pm 2,2),(7,23)$\\\hline
$2xy+y=x^3-x$&$(-2,2),(0,0),(\pm 1,0)$ \\\hline
$2xy+y=x^3+x$&$(-3,6),(-1,2),(0,0),(2,2)$ \\\hline
\end{tabular}
\captionof{table}{Integer solutions to the equations listed in Table \ref{tab:H20almWei2}. }\label{Table:ex3.35}
\end{center}

\subsection{Exercise 3.38}\label{ex:H17nearsquare}
\textbf{\emph{Solve all equations from Table \ref{tab:H17nearsquare}. }}

\begin{center}
\begin{tabular}{ |c|c|c|c|c|c| } 
 \hline
 $H$ & Equation & $H$ & Equation & $H$ & Equation \\ 
 \hline\hline
 $14$ & $y^2+x^2y+x=0$ & $16$ & $y^2+x^2y+2x=0$ & $17$ & $y^2+x^2y+2x-1=0$ \\ 
 \hline
 $15$ & $y^2+x^2y+x-1=0$ & $16$ & $y^2+x^2y-y+x=0$ & $17$ & $y^2+x^2y-y+x-1=0$ \\ 
 \hline
 $15$ & $y^2+x^2y+x+1=0$ & $16$ & $y^2+x^2y+y+x=0$ & $17$ & $y^2+x^2y-y+x+1=0$ \\ 
 \hline
 $16$ & $y^2+x^2y+x-2=0$ & $17$ & $y^2+x^2y+x-3=0$ & $17$ & $y^2+x^2y+y+x-1=0$ \\ 
 \hline
 $16$ & $y^2+x^2y+x+2=0$ & $17$ & $y^2+x^2y+x+3=0$ & $17$ & $y^2+x^2y+y+x+1=0$ \\ 
 \hline
\end{tabular}
\captionof{table}{\label{tab:H17nearsquare} Equations of the form \eqref{eq:genquady} with $a,b,c$ satisfying \eqref{eq:degreecond} of size $H\leq 17$.}
\end{center} 

This section will consider equations quadratic in $x$ or $y$ with near-square discriminants and we denote the discriminant as $D(x)$ or $D(y)$, respectively. For an equation of the form
\begin{equation}\label{eq:genquady}
a(x) y^2 + b(x) y + c(x) = 0,
\end{equation}
where $a(x),b(x),c(x)$ are polynomials with integer coefficients, the discriminant is $D(x)=b(x)^2-4a(x)c(x)$. We will consider equations such that $a(x)$ is not identically $0$ and 
\begin{equation}\label{eq:degreecond}
\text{deg}(b(x)) > \text{deg}(a(x)) + \text{deg}(c(x)),
\end{equation}
where $\text{deg}$ denotes the degree of a polynomial.

Equation
$$
y^2+x^2y+x=0
$$
is solved in Section 3.3.1 of the book and its integer solutions are
$$
(x,y)=(0,0).
$$

Equation
$$
y^2+x^2y+x-1=0
$$
is solved in Section 3.3.1 of the book and its integer solutions are
$$
(x,y)=(-1,-2),(0,\pm 1),\pm (1,-1), \quad \text{and} \quad (1,0).
$$

Equation
$$
y^2+x^2y+x+1=0
$$
is solved in Section 3.3.1 of the book and its integer solutions are
$$
(x,y)=(-1,-1),(-1,0),(2,-3), \quad \text{and} \quad (2,-1).
$$

The first equation we will consider is
\begin{equation}\label{y2px2ypxm2}
	y^2+x^2y+x-2=0.
\end{equation}
Solving this equation as a quadratic in $y$, we have $y_{1,2}=\frac{-x^2\pm\sqrt{x^4-4x+8}}{2}$, so $D(x)=x^4-4x+8$. In order for the equation to have integer solutions, $x^4-4x+8$ must be a perfect square. If $x$ is large, then
$(x^2 + 1)^2 > x^4 - 4x + 8 > (x^2 - 1)^2$, or equivalently, $2x^2 +1 > -4x+8 > -2x^2 +1$, which is true for $|x| \geq 4$. Hence, we are left to check the cases $|x| \leq 3$ and $D(x)=(x^2)^2$. In the first case, we find that $D(x)$ is a perfect square for $x=2$, after substituting this into \eqref{y2px2ypxm2} we obtain the integer solutions
\begin{equation}\label{y2px2ypxm2_sol}
	(x,y)=(2,-4),(2,0).
\end{equation}  
In the second case, $x^4-4x+8=x^4$, hence $-4x+8=0$, or $x=2$, which we have already considered. Therefore, the complete list of integer solutions to equation \eqref{y2px2ypxm2} is \eqref{y2px2ypxm2_sol}.

\vspace{10pt}

The next equation we will consider is
\begin{equation}\label{y2px2ypxp2}
y^2+x^2y+x+2=0.
\end{equation}
Solving this equation as a quadratic in $y$, we have $y_{1,2}=\frac{-x^2\pm\sqrt{x^4-4x-8}}{2}$, so $D(x)=x^4-4x-8$. In order for the equation to have integer solutions, $x^4-4x-8$ must be a perfect square. For large $x$, $(x^2+1)^2 > x^4 - 4x - 8 > (x^2- 1)^2$, or equivalently, $2x^2 +1 > -4x-8 > -2x^2 +1$, which is true for $|x| \geq 4$. Hence, it is left to check the cases $|x| \leq 3$ and $D(x)=(x^2)^2$. In the first case, we find that $D(x)$ is a perfect square for $x=\pm2$, after substituting this into \eqref{y2px2ypxp2} we obtain the integer solutions
\begin{equation}\label{y2px2ypxp2_sol}
(x,y)=(-2,-4),(-2,0),(2,-2).
\end{equation} 
In the second case, $x^4-4x-8=x^4$, hence we must have $x=-2$, which we have already considered. Therefore, the complete list of integer solutions to equation \eqref{y2px2ypxp2} is \eqref{y2px2ypxp2_sol}.

\vspace{10pt}

The next equation we will consider is
\begin{equation}\label{y2px2yp2x}
y^2+x^2y+2x=0.
\end{equation}
Solving this equation as a quadratic in $y$, we have $y_{1,2}=\frac{-x^2\pm\sqrt{x^4-8x}}{2}$, so $D(x)=x^4-8x$. In order for the equation to have integer solutions, $x^4-8x$ must be a perfect square. If $x$ is large, then
$(x^2 + 1)^2 > x^4 - 8x > (x^2 - 1)^2$, or equivalently, $2x^2 + 1 > -8x > -2x^2 + 1$, which is true for $|x| \geq 5$. Hence,  it is left to check the cases $|x| \leq 4$ and $D(x)=(x^2)^2$. These cases return that $D(x)$ is a perfect square for $x=-1,0,2$. Therefore, we find that all integer solutions to equation \eqref{y2px2yp2x} are
$$
(x,y)=(-1, -2), (-1, 1), (0, 0), (2, -2).
$$

\vspace{10pt}

The next equation we will consider is
\begin{equation}\label{y2px2ymypx}
	y^2+x^2y-y+x=0.
\end{equation}
Solving this equation as a quadratic in $y$, we have $y_{1,2}=\frac{-(x^2-1)\pm\sqrt{(x^2-1)^2-4x}}{2}$, so $D(x)=x^4 - 2x^2 - 4x + 1$. In order for the equation to have integer solutions, $x^4-2x^2-4x+1$ must be a perfect square. If $x$ is large, then $(x^2 + 1)^2 > x^4 - 2x^2 - 4x + 1 > (x^2 - 2)^2$, or equivalently, $2x^2 +1 > -2x^2 -4x+1 > -4x^2 +4$, which is true for $|x| \geq 3$. Hence, it is left to check the cases $|x| \leq 2$ and, because we have not used consecutive squares, we must also check the cases $D(x)=(x^2)^2$ and $D(x)=(x^2-1)^2$. These cases return that $D(x)$ is a perfect square for $x=-1,0,2$. Therefore, we find that all integer solutions to equation \eqref{y2px2ymypx} are
$$
	(x, y) = (-1, \pm 1), (0, 0), (0, 1), (2, -2), (2, -1).
$$

\vspace{10pt}

The next equation we will consider is
\begin{equation}\label{y2px2ypypx}
y^2+x^2y+y+x=0.
\end{equation}
Solving this equation as a quadratic in $y$, we have $y_{1,2}=\frac{-(x^2+1)\pm\sqrt{(x^2+1)^2-4x}}{2}$, so $D(x)=x^4+2x^2-4x+1$. In order for the equation to have integer solutions, $x^4+2x^2-4x+1$ must be a perfect square. If $x$ is large, then $(x^2 + 2)^2 > x^4 + 2x^2 - 4x + 1 > (x^2 - 2)^2$, or equivalently, $4x^2 + 4 > 2x^2 -4x + 1 > -4x^2 + 4$, this is true for $|x| \geq 2$. Hence, it is left to check the cases $|x| \leq 1$, $D(x)=(x^2-1)^2$, $D(x)=(x^2)^2$ and $D(x)=(x^2+1)^2$. 
These cases return that $D(x)$ is a perfect square for $x=0,1$. Therefore, we find that all integer solutions to equation \eqref{y2px2ypypx} are
$$
(x,y)=(0,-1),(0,0),(1,-1).
$$

\vspace{10pt}

The next equation we will consider is
\begin{equation}\label{y2px2ypxm3}
	y^2+x^2y+x-3=0.
\end{equation}
Solving this equation as a quadratic in $y$, we have $y_{1,2}=\frac{-x^2\pm\sqrt{x^4-4x+12}}{2}$, so $D(x)=x^4 - 4x +12$. In order for the equation to have integer solutions, $x^4-4x+12$ must be a perfect square. If $x$ is large, then $(x^2 + 1)^2 > x^4 - 4x + 12 > (x^2 - 1)^2$, or equivalently, $2x^2 + 1 > -4x + 12 > -2x^2 + 1$, which is true for $|x|\geq 4$. Hence, it is left to check cases the $|x| \leq 3$ and $D(x)=(x^2)^2$. These cases return that $D(x)$ is a perfect square for $x=-2,1,3$.
 Therefore, we find that all integer solutions to equation \eqref{y2px2ypxm3} are
$$
	(x, y) = (-2, -5), (-2, 1), (1, -2), (1, 1), (3, -9), (3, 0).
$$

\vspace{10pt}

The next equation we will consider is
\begin{equation}\label{y2px2ypxp3}
y^2+x^2y+x+3=0.
\end{equation}
Solving this equation as a quadratic in $y$, we have $y_{1,2}=\frac{-x^2\pm\sqrt{x^4-4x-12}}{2}$, so $D(x)=x^4 - 4x -12$. In order for the equation to have integer solutions, $x^4-4x-12$ must be a perfect square. If $x$ is large, then $(x^2 + 1)^2 > x^4 - 4x - 12 > (x^2 - 1)^2$, or equivalently, $2x^2 + 1 > -4x - 12 > -2x^2 + 1$, which is true for $|x|\geq 4$. Hence, it is left to check the cases $|x| \leq 3$ and $D(x)=(x^2)^2$. These cases return that $D(x)$ is a perfect square for $x=-3$.
Therefore, we find that all integer solutions to equation \eqref{y2px2ypxp3} are 
$$
(x, y) = (-3, -9), (-3, 0).
$$

\vspace{10pt}

The next equation we will consider is
\begin{equation}\label{y2px2yp2xm1}
y^2+x^2y+2x-1=0.
\end{equation}
Solving this equation as a quadratic in $y$, we have $y_{1,2}=\frac{-x^2\pm\sqrt{x^4-8x+4}}{2}$, so $D(x)=x^4-8x+4$. In order for the equation to have integer solutions, $x^4-8x+4$ must be a perfect square. If $x$ is large, then
$(x^2 + 1)^2 > x^4 - 8x+4 > (x^2 - 1)^2$, or equivalently, $2x^2  > -8x+3 > -2x^2 $, which is true for $|x| \geq 5$. Hence, it is left to check the cases $|x| \leq 4$ and $D(x)=(x^2)^2$. These cases return that $D(x)$ is a perfect square for $x=0,\pm2$. Therefore, we find that all integer solutions to equation \eqref{y2px2yp2xm1} are
$$
(x,y)=(-2,-5),(0,\pm 1),(2,-3),\pm(2,-1).
$$

\vspace{10pt}

The next equation we will consider is
\begin{equation}\label{y2px2ymypxm1}
	y^2+x^2y-y+x-1=0.
\end{equation}
Solving this equation as a quadratic in $y$, we have $y_{1,2}=\frac{-(x^2-1)\pm\sqrt{(x^2-1)^2-4x+4}}{2}$, so $D(x)=x^4 - 2x^2 - 4x +5$. In order for the equation to have integer solutions, $x^4-2x^2-4x+5$ must be a perfect square. If $x$ is large, then $(x^2 + 1)^2 > x^4 -2x^2- 4x +5 > (x^2 - 2)^2$, which is true for $|x|\geq 2$. Hence, it is left to check the cases $|x| \leq 1$, $D(x)=(x^2)^2$ and $D(x)=(x^2-1)^2$, which return that $D(x)$ is a perfect square for $x=1$. 
Therefore, we find that all integer solutions to equation \eqref{y2px2ymypxm1} are 
$$
	(x, y) = (1,0).
$$

\vspace{10pt}

The next equation we will consider is
\begin{equation}\label{y2px2ymypxp1}
y^2+x^2y-y+x+1=0.
\end{equation}
Solving this equation as a quadratic in $y$, we have $y_{1,2}=\frac{-(x^2-1)\pm\sqrt{(x^2-1)^2-4x-4}}{2}$, so $D(x)=x^4 - 2x^2 - 4x - 3$. In order for the equation to have integer solutions, $x^4-2x^2-4x-3$ must be a perfect square. If $x$ is large, then $(x^2 + 1)^2 > x^4 -2x^2- 4x -3 > (x^2 - 2)^2$, which is true for $|x|\geq 4$. Hence, it is left to check the cases $|x| \leq 3$, $D(x)=(x^2)^2$ and $D(x)=(x^2-1)^2$. These cases return that $D(x)$ is a perfect square for $x=-1$. 
Therefore, we find that all integer solutions to equation \eqref{y2px2ymypxp1} are 
$$
(x, y) = (-1,0).
$$

\vspace{10pt}

The next equation we will consider is
\begin{equation}\label{y2px2ypypxm1}
y^2+x^2y+y+x-1=0.
\end{equation}
Solving this equation as a quadratic in $y$, we have $y_{1,2}=\frac{-(x^2+1)\pm\sqrt{(x^2+1)^2-4x+4}}{2}$, so $D(x)=x^4 + 2x^2 - 4x +5$. In order for the equation to have integer solutions, $x^4+2x^2-4x+5$ must be a perfect square. If $x$ is large, then $(x^2 + 2)^2 > x^4 +2x^2- 4x +5 > (x^2 - 2)^2$, which is true for $|x|\geq 3$. Hence, it is left to check the cases $|x| \leq 2$, $D(x)=(x^2-1)^2$, $D(x)=(x^2)^2$ and $D(x)=(x^2+1)^2$, which return that $D(x)$ is a perfect square for $x=1$. 
Therefore, we find that all integer solutions to equation \eqref{y2px2ypypxm1} are 
$$
(x, y) = (1,-2),(1,0).
$$

\vspace{10pt}

The final equation we will consider is
\begin{equation}\label{y2px2ypypxp1}
	y^2+x^2y+y+x+1=0.
\end{equation}
Solving this equation as a quadratic in $y$, we have $y_{1,2}=\frac{-(x^2+1)\pm\sqrt{(x^2+1)^2-4x-4}}{2}$, so $D(x)=x^4 + 2x^2 - 4x - 3$. In order for the equation to have integer solutions, $x^4+2x^2-4x-3$ must be a perfect square. If $x$ is large, then $(x^2 + 2)^2 > x^4 +2x^2- 4x -3 > (x^2 - 2)^2$, or equivalently, $4x^2 + 4 > 2x^2 - 4x - 3 > -4x^2 + 4$, this is true for $|x|\geq 2$. Hence, it is left to check the cases $|x| \leq 1$, $D(x)=(x^2-1)^2$, $D(x)=(x^2)^2$ and $D(x)=(x^2+1)^2$, which return that $D(x)$ is a perfect square for $x=-1$. 
Therefore, we find that all integer solutions to equation \eqref{y2px2ypypxp1} are 
$$ 
	(x, y) = (-1,-2),(-1,0).
$$

Table \ref{Table:ex3.37} presents all integer solutions to the equations listed in Table \ref{tab:H17nearsquare}. 

\begin{center}
\begin{tabular}{|c|c|c|}
\hline
 Equation & Solution $(x,y)$ \\\hline \hline
$y^2+x^2y+x=0$&$(0,0)$ \\\hline
$y^2+x^2y+x-1=0$&$(-1, -2), (0, \pm 1), \pm(1, -1), (1, 0) $ \\\hline
$y^2+x^2y+x+1=0$&$(-1, -1), (-1, 0), (2, -3), (2, -1)$ \\\hline
$y^2+x^2y+x+2=0$&$(-2,-4),(-2,0),(2,-2)$ \\\hline
$y^2+x^2y+x-2=0$&$(2,-4),(2,0)$ \\\hline
$y^2+x^2y+2x=0$&$(-1, -2), (-1, 1), (0, 0), (2, -2)$ \\\hline
$y^2+x^2y+y+x=0$&$(0,-1),(0,0),(1,-1)$ \\\hline
$y^2+x^2y-y+x=0$&$(-1, \pm 1), (0, 0), (0, 1), (2, -2), (2, -1)$ \\\hline
$y^2+x^2y+x+3=0$&$(-3, -9), (-3, 0)$\\\hline
$y^2+x^2y+x-3=0$&$(-2, -5), (-2, 1), (1, -2), (1, 1), (3, -9), (3, 0)$\\\hline
$y^2+x^2y+2x-1=0$&$(-2,-5),(0,\pm 1),(2,-3),\pm (2,-1)$ \\\hline
$y^2+x^2y+y+x+1=0$&$(-1,-2),(-1,0)$\\\hline
$y^2+x^2y-y+x+1=0$&$(-1,0)$\\\hline
$y^2+x^2y+y+x-1=0$& $(1,-2),(1,0)$ \\\hline
$y^2+x^2y-y+x-1=0$&$(1,0)$ \\\hline
\end{tabular}
\captionof{table}{Integer solutions to the equations listed in Table \ref{tab:H17nearsquare}.}\label{Table:ex3.37}
\end{center}

\subsection{Exercise 3.39}\label{ex:H20nearsquare}
\textbf{\emph{Solve all equations from Table \ref{tab:H20nearsquare}. }}

\begin{center}
\begin{tabular}{ |c|c|c|c|c|c| } 
 \hline
 $H$ & Equation & $H$ & Equation & $H$ & Equation \\ 
 \hline\hline
 $18$ & $y^2+x^2y-x^2+x=0$ & $20$ & $y^2+x^2y-x^2+x+2=0$ & $20$ & $y^2+x^2y+xy-x^2=0$ \\ 
 \hline
 $18$ & $y^2+x^2y+x^2+x=0$ & $20$ & $y^2+x^2y+x^2+x+2=0$ & $20$ & $y^2+x^2y+xy+x^2=0$ \\ 
 \hline
 $18$ & $xy^2+y-x^3=0$ & $20$ & $y^2+x^2y+y-x^2+x=0$ & $20$ & $y^2+x^2 y+x^3 = 0$ \\ 
 \hline
 $19$ & $y^2+x^2y-x^2+x+1=0$ & $20$ & $y^2+x^2y-y+x^2+x=0$ & $20$ & $xy^2+y-x^3+2=0$ \\ 
 \hline
 $19$ & $y^2+x^2y+x^2+x+1=0$ & $20$ & $y^2+x^2y-x^2+2x=0$ & $20$ & $xy^2+2y-x^3=0$ \\ 
 \hline
 $19$ & $xy^2+y-x^3+1=0$ & $20$ & $y^2+x^2y+x^2+2x=0$ & $20$ & $xy^2+x^2y-y+x=0$ \\ 
 \hline
\end{tabular}
\captionof{table}{\label{tab:H20nearsquare} Equations of the form \eqref{eq:genquady} of size $H\leq 20$ with $D(x)$ satisfying \eqref{eq:detdegreecond} or \eqref{eq:dxs2xrx}.}
\end{center} 

In this exercise we will apply methods from Section 3.3.1 of the book which works for any equation of the form \eqref{eq:genquady} with $a(x) \not \equiv 0$ such that $D(x)=b(x)^2-4a(x)c(x)$ can be represented as
\begin{equation}\label{eq:detdegreecond}
D(x) = P(x)^2 + Q(x), \quad \text{deg}(P(x)) \geq \text{deg}(Q(x)),
\end{equation}
for some polynomials $P(x)$ and $Q(x)$ with integer coefficients, or,
\begin{equation}\label{eq:dxs2xrx}
	D(x)=S(x)^2 R(x),
\end{equation}
where $S(x)$ and $R(x)$ are non-constant polynomials with integer coefficients.

 If the discriminant $D(x)$ is of the form \eqref{eq:detdegreecond} and let $(x,y)$ be any integer solution to \eqref{eq:genquady} with $a(x)\neq 0$, then $D(x)$ must be a perfect square. The degree condition in \eqref{eq:detdegreecond} implies that there exists integers $k_0$ and $x_0=x_0(k_0)$ such that for every $x$ with $|x|> x_0$ we have
\begin{equation}\label{eq:degreecondaux}
|(P(x)+k_0)^2-P(x)^2| \geq 2|k_0||P(x)|-k_0^2 > |Q(x)| = |D(x)-P(x)^2|.
\end{equation}
We must then find a suitable $k_0$ such that the leading coefficient of $2|k_0||P(x)|-k_0^2$ exceeds the leading coefficient of $|Q(x)|$. Then inequality \eqref{eq:degreecondaux} will hold for all $|x| \geq x_0$, where $x_0$ is some constant. We then have a finite number of cases to check, because $D(x)$ may only be a perfect square either if $|x| \leq x_0$ or if $D(x)=(P(x)+k)^2$ for $|k| \leq k_0$.

 If the discriminant $D(x)$ is of the form \eqref{eq:dxs2xrx}, then it is a perfect square if and only if $R(x)$ is a perfect square. The problem is then reduced to solving the equation $R(x)=t^2$ for integer variable $t$.

Equation
$$
y^2+x^2y+x^2+x=0
$$
is solved in Section 3.3.1 of the book, and its integer solutions are
$$
(x,y)=(-1,-1),(-1,0),(0,0).
$$

Equation
$$
xy^2+y-x^3=0
$$
is solved in Section 3.3.1 of the book, and its integer solutions are
$$
(x,y)=(0,0).
$$

Equation
$$
y^2+x^2y+x^3=0
$$
is solved in Section 3.3.1 of the book, and its integer solutions are
$$
(x,y)=(0,0),(4,-8).
$$

The first equation we will consider is
\begin{equation}\label{y2px2ymx2px}
y^2 + x^2y - x^2 + x = 0.
\end{equation}

Solving this equation as a quadratic in $y$, we have $D(x)=x^4+4x^2-4x=(x^2+2)^2-4x-4$. So $D(x)=P(x)^2+Q(x)$ where $P(x)=x^2+2$ and $Q(x)=-4x-4$. So 
$$
|(x^2+2 \pm 1)^2 - (x^2+2)^2| \geq 2|\pm 1||x^2+2| - (\pm 1)^2 > |-4x-4| = |D(x) - (x^2+2)^2|.
$$ 
Here, the first inequality holds for all integer $x$, while the second one  
holds for integers $x<0$ and $x>2$. Hence we need to check the cases $x=0,1,2$ and $D(x)=(x^2+2)^2$. We then obtain that all integer solutions to equation \eqref{y2px2ymx2px} are
$$
(x,y)=(-1,-2),(0,0),\pm (1,-1),(1,0).
$$

\vspace{10pt}

The next equation we will consider is
\begin{equation}\label{y2px2ymx2pxp1}
y^2 + x^2y - x^2 + x +1 = 0.
\end{equation}
Solving this equation as a quadratic in $y$, we have $D(x)=x^4+4x^2-4x-4=(x^2+2)^2-4x-8$. So $D(x)=P(x)^2+Q(x)$ where $P(x)=x^2+2$ and $Q(x)=-4x-8$. So 
$$
|(x^2+2 \pm 1)^2 - (x^2+2)^2| \geq 2|\pm 1||x^2+2| - (\pm 1)^2 > |-4x-8| = |D(x) - (x^2+2)^2|.
$$ 
Here, the first inequality holds for all integer $x$, while the second one holds for integers $x<0$ and $x>2$. 
Hence we need to check the cases $x=0,1,2$ and $D(x)=(x^2+2)^2$. We then obtain that all integer solutions to equation \eqref{y2px2ymx2pxp1} are 
$$
(x,y)=(-2,-5),(-2,1).
$$

\vspace{10pt}

The next equation we will consider is
\begin{equation}\label{y2px2ypx2pxp1}
y^2 + x^2y + x^2 + x +1 = 0.
\end{equation}
Solving this equation as a quadratic in $y$, we have $D(x)=x^4-4x^2-4x-4=(x^2-2)^2-4x-8$. So $D(x)=P(x)^2+Q(x)$ where $P(x)=x^2-2$ and $Q(x)=-4x-8$. So 
$$
|(x^2-2 \pm 1)^2 - (x^2-2)^2| \geq 2|\pm 1||x^2-2| - (\pm 1)^2 > |-4x-8| = |D(x) - (x^2-2)^2|.
$$
Here, the first inequality holds for all integer $x$, while the second one holds for integers $x<-1$ and $x>3$. 
 Hence we need to check the cases $x=-1,0,1,2,3$ and $D(x)=(x^2-2)^2$. We then obtain that all integer solutions to equation \eqref{y2px2ypx2pxp1} are
$$
(x,y)=(-2,-3),(-2,-1).
$$

\vspace{10pt}

The next equation we will consider is
\begin{equation}\label{xy2pymx3p1}
xy^2 + y - x^3 + 1 = 0.
\end{equation}
Solving this equation as a quadratic in $y$, we have $D(x)=4x^4-4x+1=(2x^2)^2-4x+1$. So $D(x)=P(x)^2+Q(x)$ where $P(x)=2x^2$ and $Q(x)=-4x+1$. So
 $$
 |(2x^2 \pm 1)^2 - (2x^2)^2| \geq 2|\pm 1||2x^2| - (\pm 1)^2 > |-4x+1| = |D(x) - (2x^2)^2|.
 $$
 Here, the first inequality holds for all integer $x$, while the second one holds for integers $x<-1$ and $x>1$. 
 Hence we need to check the cases $x=-1,0,1$ and $D(x)=(2x^2)^2$. We then obtain that all integer solutions to equation \eqref{xy2pymx3p1} are
$$
(x,y)=(-1,2),(0,-1),(\pm 1,-1),(1,0).
$$

\vspace{10pt}

The next equation we will consider is
\begin{equation}\label{y2px2ymx2pxp2}
y^2 + x^2y - x^2 + x +2 = 0.
\end{equation}
Solving this equation as a quadratic in $y$, we have $D(x)=x^4+4x^2-4x-8=(x^2+2)^2-4x-12$. So $D(x)=P(x)^2+Q(x)$ where $P(x)=x^2+2$ and $Q(x)=-4x-12$. So 
$$
|(x^2+2 \pm 1)^2 - (x^2+2)^2| \geq 2|\pm 1||x^2+2| - (\pm 1)^2 > |-4x-12| = |D(x) - (x^2+2)^2|.
$$
 Here, the first inequality holds for all integer $x$, while the second one holds for integers $x<-1$ and $x>3$. 
 Hence we need to check the cases $x=-1,0,1,2,3$ and $D(x)=(x^2+2)^2$. We then obtain that all integer solutions to equation \eqref{y2px2ymx2pxp2} are
$$
(x,y)=(-3,-10),(-3,1),(-1,-1),(-1,0),(2,-4),(2,0).
$$

\vspace{10pt}

The next equation we will consider is
\begin{equation}\label{y2px2ypx2pxp2}
y^2 + x^2y + x^2 + x +2 = 0.
\end{equation}
Solving this equation as a quadratic in $y$, we have $D(x)=x^4-4x^2-4x-8=(x^2-2)^2-4x-12$. So $D(x)=P(x)^2+Q(x)$ where $P(x)=x^2-2$ and $Q(x)=-4x-12$. So 
$$
|(x^2-2 \pm 1)^2 - (x^2-2)^2| \geq 2|\pm 1||x^2-2| - (\pm 1)^2 > |-4x-12| = |D(x)- (x^2-2)^2|.
$$
 Here, the first inequality holds for all integer $x$, while the second one holds for integers $x<-2$ and $x>4$. 
 Hence we need to check the cases $-2 \leq x\leq 4$ and $D(x)=(x^2-2)^2$. We then obtain that all integer solutions to equation \eqref{y2px2ypx2pxp2} are
$$
(x,y)=(-3,-8),(-3,-1),(-2,-2),(3,-7),(3,-2).
$$

\vspace{10pt}

The next equation we will consider is
\begin{equation}\label{y2px2ypymx2px}
y^2+x^2y+y-x^2+x=0.
\end{equation}
Solving this equation as a quadratic in $y$, we have $D(x)=x^4+6x^2-4x+1=(x^2+3)^2-4x-8$. So $D(x)=P(x)^2+Q(x)$ where $P(x)=x^2+3$ and $Q(x)=-4x-8$. So 
$$
|(x^2+3 \pm 1)^2 - (x^2+3)^2| \geq 2|\pm 1||x^2+3| - (\pm 1)^2 > |-4x-8| = |D(x) - (x^2+3)^2|.
$$
 Here, the first inequality holds for all integer $x$, while the second one holds for integers $x<0$ and $x>2$. 
 Hence we need to check the cases $x=0,1,2$ and  $D(x)=(x^2+3)^2$. We then obtain that all integer solutions to equation \eqref{y2px2ypymx2px} are
$$
(x,y)=(-2,-6),(-2,1),(0,-1),(0,0),(1,-2),(1,0).
$$

\vspace{10pt}

The next equation we will consider is
\begin{equation}\label{y2px2ymypx2px}
y^2+x^2y-y+x^2+x=0.
\end{equation}
Solving this equation as a quadratic in $y$, we have $D(x)=x^4-6x^2-4x+1=(x^2-3)^2-4x-8$. So $D(x)=P(x)^2+Q(x)$ where $P(x)=x^2-3$ and $Q(x)=-4x-8$. So 
$$
|(x^2-3 + \pm 1)^2 - (x^2-3)^2| \geq 2|\pm 1||x^2-3| - (\pm 1)^2 > |-4x-8| = |D(x)- (x^2-3)^2|.
$$
Here, the first inequality holds for all integer $x$, while the second one holds for integers $x<-1$ and $x>3$. 
 Hence we need to check the cases $-1 \leq x \leq 3$ and $D(x)=(x^2-3)^2$. We then obtain that all integer solutions to equation \eqref{y2px2ymypx2px} are
$$
(x,y)=(-2,-2),(-2,-1),(-1,0),(0,0),(0,1),(3,-6),(3,-2).
$$

\vspace{10pt}

The next equation we will consider is
\begin{equation}\label{y2px2ymx2p2x}
y^2 + x^2y - x^2 + 2x = 0.
\end{equation}
Solving this equation as a quadratic in $y$, we have $D(x)=x^4+4x^2-8x=(x^2+2)^2-8x-4$. So $D(x)=P(x)^2+Q(x)$ where $P(x)=x^2+2$ and $Q(x)=-8x-4$. So 
$$
|(x^2+2 \pm 1)^2 - (x^2+2)^2| \geq 2|\pm 1||x^2+2| - (\pm 1)^2 > |-8x-4| = |D(x)- (x^2+2)^2|.
$$
Here, the first inequality holds for all integer $x$, while the second one holds for integers $x<-2$, $x=-1$ and $x>4$. 
  Hence we need to check the cases $x=-2,0,1,2,3,4$ and $D(x)=(x^2+2)^2$. We then obtain that all integer solutions to equation \eqref{y2px2ymx2p2x} are
$$
(x,y)=(0,0),(2,-4),(2,0).
$$

\vspace{10pt}

The next equation we will solve is
\begin{equation}\label{y2px2ypx2p2x}
y^2 + x^2y + x^2 + 2x  = 0.
\end{equation}
Solving this equation as a quadratic in $y$, we have $D(x)=x^4-4x^2-8x=(x^2-2)^2-8x-4$. So $D(x)=P(x)^2+Q(x)$ where $P(x)=x^2-2$ and $Q(x)=-8x-4$. So 
$$
|(x^2-2 \pm 1)^2 - (x^2-2)^2| \geq 2|\pm 1||x^2-2| - (\pm 1)^2 > |-8x-4| = |D(x) - (x^2-2)^2|.
$$
Here, the first inequality holds for all integer $x$, while the second one holds for integers
 $|x| > 4$. Hence we need to check the cases $|x|\leq 4$ and $D(x)=(x^2-2)^2$. We then obtain that all integer solutions to equation \eqref{y2px2ypx2p2x} are
$$
(x,y)=(-2,-4),(-2,0),(0,0).
$$

\vspace{10pt}

The next equation we will consider is
\begin{equation}\label{y2px2ypxymx2}
y^2 + x^2y + xy - x^2 = 0.
\end{equation}
Solving this equation as a quadratic in $y$, we have $D(x)=x^4+2x^3+5x^2=(x^2+x+2)^2-4x-4$. So $D(x)=P(x)^2+Q(x)$ where $P(x)=x^2+x+2$ and $Q(x)=-4x-4$. So 
$$
\begin{aligned}
|(x^2+x+2 \pm 1)^2 - (x^2+x+2)^2| \geq 2|\pm 1||x^2+x+2| - (\pm 1)^2 > |-4x-4| \\ =  |D(x) - (x^2+x+2)^2|.
\end{aligned}
$$
Here, the first inequality holds for all integer $x$, while the second one holds for integers
 $x < 0$ and $x >1$. Hence we need to check the cases $x=0,1$ and $D(x)=(x^2+x+2)^2$. We then obtain that all integer solutions to equation \eqref{y2px2ypxymx2} are
$$
(x,y)=(-1,\pm 1),(0,0).
$$

\vspace{10pt}

The next equation we will consider is
\begin{equation}\label{y2px2ypxypx2}
y^2 + x^2y + xy + x^2 = 0.
\end{equation}
Solving this equation as a quadratic in $y$, we have $D(x)=x^4+2x^3-3x^2=(x^2+x-2)^2+4x-4$. So $D(x)=P(x)^2+Q(x)$ where $P(x)=x^2+x-2$ and $Q(x)=4x-4$. So 
$$
|(x^2+x-2 \pm 1)^2 - (x^2+x-2)^2| \geq 2|\pm 1||x^2+x-2| - (\pm 1)^2 > |4x-4|= |D(x) - (x^2+x-2)^2|.
$$
Here, the first inequality holds for all integer $x$, while the second one holds for integers $x<-4$ and $x>1$.
  Hence we need to check the cases $x=-4,-3,-2,-1,0,1$ and $D(x)=(x^2+x-2)^2$. We then obtain that all integer solutions to equation \eqref{y2px2ypxypx2} are
$$
(x,y)=(-3,-3),(0,0),(1,-1).
$$

\vspace{10pt}

The next equation we will consider is
\begin{equation}\label{xy2pymx3p2}
xy^2 + y - x^3 + 2 = 0.
\end{equation}
Solving this equation as a quadratic in $y$, we have $D(x)=4x^4-8x+1=(2x^2)^2-8x+1$. So $D(x)=P(x)^2+Q(x)$ where $P(x)=2x^2$ and $Q(x)=-8x+1$. So 
$$
|(2x^2 \pm 1)^2 - (2x^2)^2| \geq 2|\pm 1||2x^2| - (\pm 1)^2 > |-8x+1| = |D(x) - (2x^2)^2|.
$$
Here, the first inequality holds for all integer $x$, while the second one holds for integers $|x|>2$.
  Hence we need to check the cases $|x|\leq 2$ and $D(x)=(2x^2)^2$. We then obtain that all integer solutions to equation \eqref{xy2pymx3p2} are
$$
(x,y)=(0,-2),(\pm 2,-2).
$$

\vspace{10pt}

The next equation we will consider is
\begin{equation}\label{xy2p2ymx3}
xy^2 + 2y - x^3 = 0.
\end{equation}
Solving this equation as a quadratic in $y$, we have $D(x)=4x^4+4=2^2(x^4+1)$. In order for \eqref{xy2p2ymx3} to have integer solutions, we must have that $(x^2)^2+1$ is a perfect square. Because the only consecutive perfect squares are $0$ and $1$, we must have $x=0$. Therefore, the unique integer solution to equation \eqref{xy2p2ymx3} is
$$
(x,y)=(0,0).
$$

\vspace{10pt}

The final equation we will consider is
\begin{equation}\label{xy2px2ymypx}
xy^2+x^2y-y+x=0.
\end{equation}
Solving this equation as a quadratic in $y$, we have $D(x)=x^4-6x^2+1=(x^2-3)^2-8$. In order for \eqref{xy2px2ymypx} to have integer solutions, we must have that $(x^2-3)^2-8$ is a perfect square. This is equivalent to solving the equation $X^2-8=Y^2$, where $X=x^2-3$, which can be solved using the method from Section \ref{ex:PQc}, and we obtain that is only has solutions when $X=\pm 3$. This is possible only if $x=0$. Therefore, the unique integer solution to equation \eqref{xy2px2ymypx} is
$$
(x,y)=(0,0).
$$

Table \ref{Table:ex3.38} summarises all integer solutions to the equations listed in Table \ref{tab:H20nearsquare}. 

\begin{center}
\begin{tabular}{|c|c|c|}
\hline
Equation & Solution $(x,y)$ \\\hline \hline
$y^2 + x^2y - x^2 + x = 0$&$(-1,-2),(0,0),\pm (1,-1),(1,0)$ \\\hline
$y^2+x^2y+x^2+x=0$&$ (-1, -1), (-1, 0), (0, 0)$ \\\hline
$xy^2+y-x^3=0$&$(0,0)$ \\\hline
$y^2 + x^2y - x^2 + x +1 = 0$&$(-2,-5),(-2,1)$\\\hline
$y^2 + x^2y + x^2 + x +1 = 0$&$(-2,-3),(-2,-1)$\\\hline
$xy^2 + y - x^3 + 1 = 0$&$(-1,2),(0,-1),(\pm 1,-1),(1,0)$\\\hline
$y^2 + x^2y - x^2 + x +2 = 0$&$(-3,-10),(-3,1),(-1,-1),(-1,0),(2,-4),(2,0)$\\\hline
$y^2 + x^2y + x^2 + x +2 = 0$&$(-3,-8),(-3,-1),(-2,-2),(3,-7),(3,-2)$\\\hline
$y^2+x^2y+y-x^2+x=0$&$(-2,-6),(-2,1),(0,-1),(0,0),(1,-2),(1,0)$\\\hline
$y^2+x^2y-y+x^2+x=0$&$(-2,-2),(-2,-1),(-1,0),(0,0),(0,1),(3,-6),(3,-2)$\\\hline
$y^2 + x^2y - x^2 + 2x = 0$&$(0,0),(2,-4),(2,0)$\\\hline
$y^2 + x^2y + x^2 + 2x  = 0$&$(-2,-4),(-2,0),(0,0)$\\\hline
$y^2 + x^2y + xy - x^2 = 0$&$(-1,\pm 1),(0,0)$\\\hline
$y^2 + x^2y + xy + x^2 = 0$&$(-3,-3),(0,0),(1,-1)$\\\hline
$y^2+x^2y+x^3=0$&$(0,0),(4,-8)$ \\\hline
$xy^2 + y - x^3 + 2 = 0$&$(0,-2),(\pm 2,-2)$\\\hline
$xy^2 + 2y - x^3 = 0$&$(0,0)$\\\hline
$xy^2+x^2y-y+x=0$&$(0,0)$\\\hline

\end{tabular}
\captionof{table}{Integer solutions to the equations listed in Table \ref{tab:H20nearsquare}.}\label{Table:ex3.38}
\end{center}

\subsection{Exercise 3.40}\label{ex:H21boundedcomb}
\textbf{\emph{Without using computer assistance, solve all equations from Table \ref{tab:H21boundedcomb}.}}

\begin{center}
\begin{tabular}{ |c|c|c|c|c|c| } 
 \hline
 $H$ & Equation & $H$ & Equation & $H$ & Equation \\ 
 \hline\hline
 $18$ & $y^3 = x^3-x$ & $20$ & $y^3=x^3-2x$ & $21$ & $y^3 = x^3+2x+1$ \\ 
 \hline
 $18$ & $y^3 = x^3+x$ & $20$ & $y^3=x^3+2x$ & $21$ & $y^3=x^3+x^2-1$ \\ 
 \hline
 $19$ & $y^3 = x^3-x+1$ & $20$ & $y^3=x^3+x^2$ & $21$ & $y^3=x^3+x^2+1$ \\ 
 \hline
 $19$ & $y^3 = x^3+x+1$ & $20$ & $y^3+xy+x^3=0$ & $21$ & $y^3+xy+x^3-1=0$ \\ 
 \hline
 $20$ & $y^3 = x^3-x+2$ & $21$ & $y^3 = x^3-x+3$ & $21$ & $y^3+xy+x^3+1=0$ \\ 
 \hline
 $20$ & $y^3 = x^3+x+2$ & $21$ & $y^3 = x^3+x+3$ &  & \\ 
 \hline
 $20$ & $y^3-y=x^3+x$ & $21$ & $y^3 = x^3-2x+1$ &  & \\ 
 \hline
\end{tabular}
\captionof{table}{\label{tab:H21boundedcomb} Equations with finite optimal value in \eqref{eq:brealopt2var2}.}
\end{center} 

This exercise solves equations $P(x,y)=0$ which have a finite optimal value in the optimisation problem
\begin{equation}\label{eq:brealopt2var2}
\max_{(x,y) \in \mathbb{R}^2}  \quad \min\{|x|,|y|,|x+y|,|x-y|\} \quad \text{subject to} \quad P(x,y)=0.
\end{equation}

Equation 
$$
y^3=x^3-x
$$
is solved in Section 3.3.2 of the book, and its integer solutions are
$$
(x,y)=(0,0),(\pm 1,0).
$$

Equation 
$$
y^3=x^3+x
$$
is solved in Section 3.3.2 of the book, and its integer solutions are
$$
(x,y)=(0,0).
$$

Let us look at the first equation
\begin{equation}\label{y3ex3mxp1}
y^3=x^3-x+1.
\end{equation}
We observe that perfect cubes $y^3$ and $x^3$ are close to each other. The closest perfect cubes to $x^3$ are $(x+1)^3$ and $(x-1)^3$. Hence, either $y=x$ or
$$
|-x+1| = |y^3 - x^3| \geq |(x \pm 1)^3 - x^3| \geq 3|x^2| - 3|x| - 1,
$$
which is only possible for integers $|x|\leq 1$. After checking the cases $y=x$ and $|x| \leq 1$, 
we can conclude that the integer solutions to equation \eqref{y3ex3mxp1} are
$$
(x,y)=(0,1),(\pm 1,1).
$$

\vspace{10pt}

The next equation we will consider is 
\begin{equation}\label{y3ex3pxp1}
y^3=x^3+x+1.
\end{equation}
We must have either $y=x$ or
$$
|x+1| = |y^3 - x^3| \geq |(x \pm 1)^3 - x^3| \geq 3|x^2| - 3|x| - 1,
$$
which is only possible for integers $|x|\leq 1$. After checking the cases $y=x$ and $|x| \leq 1$, 
we can conclude that the integer solutions to equation \eqref{y3ex3pxp1} are
$$
(x,y)=(-1,-1),(0,1).
$$

\vspace{10pt}

The next equation we will consider is 
\begin{equation}\label{y3ex3mxp2}
y^3=x^3-x+2.
\end{equation}
We must have either $y=x$ or
$$
|-x+2| = |y^3 - x^3| \geq |(x \pm 1)^3 - x^3| \geq 3|x^2| - 3|x| - 1,
$$
which is only possible for integers $|x|\leq 1$. After checking the cases $y=x$ and $|x| \leq 1$, 
we can conclude that the unique integer solution to equation \eqref{y3ex3mxp2} is
$$
(x,y)=(2,2).
$$

\vspace{10pt}

The next equation we will consider is 
\begin{equation}\label{y3ex3pxp2}
y^3=x^3+x+2.
\end{equation}
We must have either $y=x$ or
$$
|x+2| = |y^3 - x^3| \geq |(x \pm 1)^3 - x^3| \geq 3|x^2| - 3|x| - 1,
$$
which is only possible for integers $|x|\leq 1$. After checking the cases $y=x$ and $|x| \leq 1$, 
we can conclude that the integer solutions to equation \eqref{y3ex3pxp2} are
$$
(x,y)=(-2,-2),(-1,0).
$$

\vspace{10pt}

The next equation we will consider is 
\begin{equation}\label{y3myex3px}
y^3-y=x^3+x.
\end{equation}
We have that $|y^3-x^3|=|x+y|$. First, assume that $|x| > |y|$, then $|y^3-x^3|=|x+y|<|2x|$. Hence,  
$$
|2x| > |x+y| = |y^3 - x^3| \geq |(x \pm 1)^3 - x^3| \geq 3|x^2| - 3|x| - 1.
$$
which is only true for integers $|x| \leq 1$. If we assume that $|x| < |y|$, then $|y^3-x^3|=|x+y|< |2y|$. Hence,
$$
|2y| > |x+y| = |y^3 - x^3| \geq |(y \pm 1)^3 - y^3| \geq 3|y^2| - 3|y| - 1.
$$
which is only true for integers $|y| \leq 1$. Hence, either $|x| \leq 1$, or $|y| \leq 1$, or $|x|=|y|$. 
By checking these cases, we can then conclude that the integer solutions to equation \eqref{y3myex3px} are
$$
(x,y)=(0,\pm 1),(0,0).
$$

\vspace{10pt}

The next equation we will consider is
\begin{equation}\label{y3ex3m2x}
y^3=x^3-2x.
\end{equation}
We must have either $y=x$ or
$$
|-2x| = |y^3 - x^3| \geq |(x \pm 1)^3 - x^3| \geq 3|x^2| - 3|x| - 1,
$$
which is only possible for integers $|x|\leq 1$. After checking the cases $y=x$ and $|x|\leq 1$,
we can conclude that the integer solutions to equation \eqref{y3ex3m2x} are
$$
(x,y)=(0,0),\pm (1,-1).
$$

\vspace{10pt}

The next equation we will consider is
\begin{equation}\label{y3ex3p2x}
y^3=x^3+2x.
\end{equation}
We must have either, $y=x$ or
$$
|2x| = |y^3 - x^3| \geq |(x \pm 1)^3 - x^3| \geq 3|x^2| - 3|x| - 1,
$$
which is only possible for integers $|x|\leq 1$. After checking the cases $y=x$ and $|x|\leq 1$,
we can conclude that the unique integer solution to equation \eqref{y3ex3p2x} is
$$
(x,y)=(0,0).
$$

\vspace{10pt}

The next equation we will consider is
\begin{equation}\label{y3ex3px2}
y^3=x^3+x^2.
\end{equation}
We must have either $y=x$ or
$$
|x^2| = |y^3 - x^3| \geq |(x \pm 1)^3 - x^3| \geq 3|x^2| - 3|x| - 1,
$$
which is only possible for integers $|x|\leq 1$. After checking the cases $y=x$ and $|x|\leq 1$,
we can conclude that the integer solutions to equation \eqref{y3ex3px2} are
$$
(x,y)=(-1,0),(0,0).
$$

\vspace{10pt}

The next equation we will consider is
\begin{equation}\label{y3pxypx3}
y^3+xy+x^3=0.
\end{equation}
We have that $|y^3-x^3|=|xy|$. First, assume that $|x| > |y|$, then $|y^3-x^3|=|xy| < |x^2|$. Hence,
$$
|x^2| > |xy| = |y^3 - x^3| \geq |(x \pm 1)^3 - x^3| \geq 3|x^2| - 3|x| - 1.
$$
which is only true for integer $|x| \leq 1$. If we then assume that $|x| < |y|$, then $|y^3-x^3|=|xy|< |y^2|$. Hence, 
$$
|y^2| > |xy| = |y^3 - x^3| \geq |(y \pm 1)^3 - y^3| \geq 3|y^2| - 3|y| - 1.
$$
which is only true for integers $|y| \leq 1$. Hence, either $|x| \leq 1$, or $|y| \leq 1$, or $|x|=|y|$.  By checking these cases, we can then conclude that the unique integer solution to equation \eqref{y3pxypx3} is
$$
(x,y)=(0,0).
$$

\vspace{10pt}

The next equation we will consider is
\begin{equation}\label{y3ex3mxp3}
y^3=x^3-x+3.
\end{equation}
We must have either $y=x$ or
$$
|-x+3| = |y^3 - x^3| \geq |(x \pm 1)^3 - x^3| \geq 3|x^2| - 3|x| - 1,
$$
which is only possible for integers $-2 \leq x \leq 1$. After checking the cases $y=x$ and $-2 \leq x \leq 1$
we can conclude that the unique integer solution to equation \eqref{y3ex3mxp3} is
$$
(x,y)=(3,3).
$$

\vspace{10pt}

The next equation we will consider is
\begin{equation}\label{y3ex3pxp3}
y^3=x^3+x+3.
\end{equation}
We must have either $y=x$ or
$$
|x+3| = |y^3 - x^3| \geq |(x \pm 1)^3 - x^3| \geq 3|x^2| - 3|x| - 1,
$$
which is only possible for integers $-1 \leq x \leq 2$. After checking the cases $y=x$ and $-1 \leq x \leq 2$,
 we can conclude that the integer solutions to equation \eqref{y3ex3pxp3} are
$$
(x,y)=(-3,-3),(-1,1).
$$

\vspace{10pt}

The next equation we will consider is
\begin{equation}\label{y3ex3m2xp1}
y^3=x^3-2x+1.
\end{equation}
We must have either $y=x$ or
$$
|-2x+1| = |y^3 - x^3| \geq |(x \pm 1)^3 - x^3| \geq 3|x^2| - 3|x| - 1,
$$
which is only possible for integers $-2 \leq x \leq 1$. After checking the cases $y=x$ and $-2 \leq x \leq 1$,
we can conclude that the integer solutions to equation \eqref{y3ex3m2xp1} are
$$
(x,y)=(0,1),(1,0).
$$

\vspace{10pt}

The next equation we will consider is
\begin{equation}\label{y3ex3p2xp1}
y^3=x^3+2x+1.
\end{equation}
We must have either $y=x$ or
$$
|2x+1| = |y^3 - x^3| \geq |(x \pm 1)^3 - x^3| \geq 3|x^2| - 3|x| - 1,
$$
which is only possible for integers $-1 \leq x \leq 2$.  After checking the cases $y=x$ and $-1 \leq x \leq 2$,
 we can conclude that the unique integer solution to equation \eqref{y3ex3p2xp1} is
$$
(x,y)=(0,1).
$$

\vspace{10pt}

The next equation we will consider is
\begin{equation}\label{y3ex3px2m1}
y^3=x^3+x^2-1.
\end{equation}
We must have either $y=x$ or
$$
|x^2-1| = |y^3 - x^3| \geq |(x \pm 1)^3 - x^3| \geq 3|x^2| - 3|x| - 1,
$$
which is only possible for integers $|x| \leq 1$. After checking the cases $y=x$ and $|x| \leq 1$,
we can conclude that the integer solutions to equation \eqref{y3ex3px2m1} are
$$
(x,y)=(0,-1),\pm(1,1).
$$

\vspace{10pt}

The next equation we will consider is
\begin{equation}\label{y3ex3px2p1}
y^3=x^3+x^2+1.
\end{equation}
We must have either $y=x$ or
$$
|x^2+1| = |y^3 - x^3| \geq |(x \pm 1)^3 - x^3| \geq 3|x^2| - 3|x| - 1,
$$
which is only possible for integers $|x| \leq 2$. After checking the cases $y=x$ and $|x| \leq 2$,
we can conclude that the integer solutions to equation \eqref{y3ex3px2p1} are
$$
(x,y)=(-1,1),(0,1).
$$

\vspace{10pt}

The next equation we will consider is
\begin{equation}\label{y3pxypx3m1}
y^3+xy+x^3-1=0.
\end{equation}
We have that $|y^3-x^3|=|xy+1|$. First, assume that $|x| > |y|$, then $|y^3-x^3|=|xy+1|< x^2+1$. Hence,
$$
x^2+1 > |xy+1| = |y^3 - x^3| \geq |(x \pm 1)^3 - x^3| \geq 3|x^2| - 3|x| - 1.
$$
which is only true for integers $|x| \leq 1$. If we then assume that $|x| < |y|$, then $|y^3-x^3|=|xy+1|< y^2+1$. Hence, 
$$
y^2+1 > |xy+1| = |y^3 - x^3| \geq |(y \pm 1)^3 - y^3| \geq 3|y^2| - 3|y| - 1.
$$
which is only true for integers $|y| \leq 1$. Hence, either $|x| \leq 1$, or $|y| \leq 1$, or $|x|=|y|$.  By checking these cases, we can then conclude that the integer solutions to equation 
\eqref{y3pxypx3m1} are
$$
(x,y)=(0,1),(1,0).
$$

\vspace{10pt}

The final equation we will consider is
\begin{equation}\label{y3pxypx3p1}
y^3+xy+x^3+1=0.
\end{equation}
We have that $|y^3-x^3|=|xy-1|$. First, assume that $|x| > |y|$, then $|y^3-x^3|=|xy-1|< x^2+1$. Hence, 
$$
x^2+1 > |xy-1| = |y^3 - x^3| \geq |(x \pm 1)^3 - x^3| \geq 3|x^2| - 3|x| - 1.
$$
which is only true for integers $|x| \leq 1$. Next, assume that $|x| < |y|$, then $|y^3-x^3|=|xy-1|< y^2+1$. We also have 
$$
y^2+1 > |xy-1| = |y^3 - x^3| \geq |(y \pm 1)^3 - y^3| \geq 3|y^2| - 3|y| - 1.
$$
which is only true for integers $|y| \leq 1$. Hence, either $|x| \leq 1$, or $|y| \leq 1$, or $|x|=|y|$.  By checking these cases, we can then conclude that the integer solutions to equation 
\eqref{y3pxypx3p1} are
$$
(x,y)=(-1,\pm 1),(-1,0),(0,-1),(1,-1).
$$

Table \ref{Table:ex3.39} summarises the integer solutions to the equations listed in Table \ref{tab:H21boundedcomb}. 

\begin{center}
\begin{tabular}{|c|c||c|c|}
\hline
Equation & Solution $(x,y)$ &  Equation & Solution $(x,y)$ \\\hline \hline
$y^3=x^3-x$ & $ (\pm1,0), (0,0)$ & $y^3+xy+x^3=0$ & $(0,0)$  \\\hline
$y^3=x^3+x$ & $(0,0)$ &  $y^3=x^3-x+3$ & $(3,3)$\\\hline
$y^3=x^3-x+1$ & $(0,1),(\pm 1,1)$ & $y^3=x^3+x+3$ & $(-3,-3),(-1,1)$  \\\hline
$y^3=x^3+x+1$ & $(-1,-1),(0,1)$ &  $y^3=x^3-2x+1$ & $(0,1),(1,0)$ \\\hline
$y^3=x^3-x+2$ & $(2,2)$ &  $y^3=x^3+2x+1$ & $(0,1)$ \\\hline
$y^3=x^3+x+2$ & $(-2,-2),(-1,0)$ &  $y^3=x^3+x^2-1$ & $(0,-1),\pm (1,1)$ \\\hline
$y^3-y=x^3+x$ & $(0,\pm1),(0,0)$ &  $y^3=x^3+x^2+1$ & $(-1,1),(0,1)$ \\\hline
$y^3=x^3-2x$ & $(0,0),\pm(1,-1)$ & $y^3+xy+x^3-1=0$ & $(0,1),(1,0)$  \\\hline
$y^3=x^3+2x$ & $(0,0)$ & $y^3+xy+x^3+1=0$ & $(-1,\pm 1),(-1,0),$ \\ \cline{1-2} 
$y^3=x^3+x^2$ & $(-1,0),(0,0)$ && $(0,-1),(1,-1)$ \\\hline

\end{tabular}
\captionof{table}{Integer solutions to the equations listed in Table \ref{tab:H21boundedcomb}.}\label{Table:ex3.39}
\end{center}

\subsection{Exercise 3.42}\label{ex:H28nonlinsub}
\textbf{\emph{For each equation listed in Table \ref{tab:H28nonlinsub}, use a non-linear substitution to reduce it to an equation describing elliptic curve in Weierstrass form, and then use the {\tt EllipticCurve} command in SageMath to solve the resulting equation and hence the original equation.}}

\begin{center}
\begin{tabular}{ |c|c|c|c|c|c| } 
 \hline
 $H$ & Equation & $H$ & Equation & $H$ & Equation \\ 
 \hline\hline
 $25$ & $x^4+y^3-1=0$ & $27$ & $x^4+y^3+3=0$ & $28$ & $x^4+y^3+y+2=0$ \\ 
 \hline
 $25$ & $x^4+y^3+1=0$ & $27$ & $x^4+y^3-y-1=0$ & $28$ & $x^4+y^3-2y=0$ \\ 
 \hline
 $26$ & $x^4+y^3-2=0$ & $27$ & $x^4+y^3+y-1=0$ & $28$ & $x^4+y^3+2y=0$ \\ 
 \hline
 $26$ & $x^4+y^3+2=0$ & $27$ & $x^4+y^3+y+1=0$ & $28$ & $x^4-x^2+y^3=0$ \\ 
 \hline
 $26$ & $x^4+y^3-y=0$ & $28$ & $x^4+y^3-4=0$ & $28$ & $x^4+x^2+y^3=0$ \\ 
 \hline
 $26$ & $x^4+y^3+y=0$ & $28$ & $x^4+y^3+4=0$ &  & \\ 
 \hline
 $27$ & $x^4+y^3-3=0$ & $28$ & $x^4+y^3+y-2=0$ &  & \\ 
 \hline
\end{tabular}
\captionof{table}{\label{tab:H28nonlinsub} Equations of size $H\leq 28$ solvable by non-linear substitutions.}
\end{center} 

We will next solve equations by reducing equations to the Weierstrass form \eqref{eq:Weiform} using non-linear substitutions. 
The first equation we will consider is 
\begin{equation}\label{x4py3my}
x^4+y^3-y=0.
\end{equation}
This equation can be reduced to 
$$
Y^2=X^3-X
$$ 
by making the substitutions $Y=x^2$ and $X=-y$. We have solved this equation previously, see Table \ref{tab:3.31}, and its integer solutions are $(X,Y)=(0,0), (\pm 1,0)$. Hence, in the original variables $x = \pm \sqrt{Y}$ and $y=-X$, we obtain that the integer solutions to equation \eqref{x4py3my} are 
$$
(x,y)=(0,0),(0,\pm 1).
$$

\vspace{10pt}

The next equation we will consider is
\begin{equation}\label{x4py3py}
x^4+y^3+y=0.
\end{equation}
This equation can be reduced to 
$$
Y^2=X^3+X
$$ 
by making the substitutions $Y=x^2$ and $X=-y$. We have solved this equation previously, see Table \ref{tab:3.31}, and its only integer solution is $(X,Y)=(0,0)$. Hence, in the original variables, we have that the unique integer solution to equation \eqref{x4py3py} is
$$
(x,y)=(0,0).
$$

\vspace{10pt}

The next equation we will consider is
\begin{equation}\label{x4py3p2}
x^4+y^3+2=0.
\end{equation}
This equation can be reduced to 
$$
Y^2=X^3-2
$$ 
by making the substitutions $Y=x^2$ and $X=-y$. We have solved this equation previously, see Table \ref{tab:3.31}, and its integer solutions are $(X,Y)=(3,\pm 5)$.  Then, in the original variables, we obtain solutions $(x,y)=(\pm \sqrt{5},-3)$, which are not integer, hence the equation \eqref{x4py3p2} has no integer solutions.

\vspace{10pt}

The next equation we will consider is
\begin{equation}\label{x4py3m2}
x^4+y^3-2=0.
\end{equation}
This equation can be reduced to 
$$
Y^2=X^3+2
$$ 
by making the substitutions $Y=x^2$ and $X=-y$. We have solved this equation previously, see Table \ref{tab:3.31}, and its integer solutions are $(X,Y)=(-1,\pm 1)$. Hence, in the original variables, we have that the integer solutions to equation \eqref{x4py3m2} are
$$
(x,y)=(\pm 1,1).
$$

\vspace{10pt}

The next equation we will consider is
\begin{equation}\label{x4py3mym1}
x^4+y^3-y-1=0.
\end{equation}
This equation can be reduced to 
$$
Y^2=X^3-X+1
$$ 
by making the substitutions $Y=x^2$ and $X=-y$. We have solved this equation previously, see Table \ref{tab:3.31}, and its integer solutions are $(X,Y)=$ $(0, \pm 1),$ $(\pm 1,\pm 1),$ $(3,\pm 5),$ $(5,\pm 11),$ $(56,\pm 419)$.
 Hence, in the original variables, we have that the integer solutions to equation \eqref{x4py3mym1} are
$$
(x,y)=(\pm 1,\pm1),(\pm 1,0).
$$

\vspace{10pt}

The next equation we will consider is
\begin{equation}\label{x4py3pyp1}
x^4+y^3+y+1=0.
\end{equation}
This equation can be reduced to 
$$
Y^2=X^3+X-1
$$ 
by making the substitutions $Y=x^2$ and $X=-y$. We have solved this equation previously, see Table \ref{tab:3.31}, and its integer solutions are $(X,Y)=(1,\pm 1),$ $(2, \pm 3),$ $(13,\pm 47)$.
Hence, in the original variables, we have that the integer solutions to equation \eqref{x4py3pyp1} are
$$
(x,y)=(\pm 1,-1).
$$

\vspace{10pt}

The next equation we will consider is
\begin{equation}\label{x4py3pym1}
x^4+y^3+y-1=0.
\end{equation}
This equation can be reduced to 
$$
Y^2=X^3+X+1
$$ 
by making the substitutions $Y=x^2$ and $X=-y$. We have solved this equation previously, see Table \ref{tab:3.31} and its integer solutions are $(X,Y)=(0,\pm 1),$ $(72, \pm 611)$.
Hence, in the original variables, we have that the integer solutions to equation \eqref{x4py3pym1} are
$$
(x,y)=(\pm 1,0).
$$

\vspace{10pt}

The next equation we will consider is
\begin{equation}\label{x4py3p3}
x^4+y^3+3=0.
\end{equation}
This equation can be reduced to 
$$
Y^2=X^3-3
$$ 
by making the substitutions $Y=x^2$ and $X=-y$. We have solved this equation previously, see Table \ref{tab:3.31}, and it has no integer solutions, hence, equation \eqref{x4py3p3} has no integer solutions.

\vspace{10pt}

The next equation we will consider is
\begin{equation}\label{x4py3m3}
x^4+y^3-3=0.
\end{equation}
This equation can be reduced to 
$$
Y^2=X^3+3
$$ 
by making the substitutions $Y=x^2$ and $X=-y$. We have solved this equation previously, see Table \ref{tab:3.31} and its integer solutions are $(X,Y)=(1,\pm 2)$.
Hence, in the original variables, we obtain that equation \eqref{x4py3m3} has no integer solutions.

\vspace{10pt}

The next equation we will consider is
\begin{equation}\label{x4py3m4}
x^4+y^3-4=0.
\end{equation}
This equation can be reduced to 
$$
Y^2=X^3+4
$$
by making the substitutions $Y=x^2$ and $X=-y$. We have solved this equation previously, see Table \ref{tab:3.31}, and its integer solutions are $(X,Y)=(0,\pm 2)$.
Hence, in the original variables, we obtain that equation \eqref{x4py3m4} has no integer solutions.

\vspace{10pt}

The next equation we will consider is
\begin{equation}\label{x4py3p4}
x^4+y^3+4=0.
\end{equation}
This equation can be reduced to 
$$
Y^2=X^3-4
$$
by making the substitutions $Y=x^2$ and $X=-y$. We have solved this equation previously, see Table \ref{tab:3.31}, and its integer solutions are $(X,Y)=(2,\pm 2),(5,\pm 11)$.
Hence, in the original variables, we obtain that equation \eqref{x4py3p4} has no integer solutions.

\vspace{10pt}

The next equation we will consider is
\begin{equation}\label{x4py3pym2}
x^4+y^3+y-2=0.
\end{equation}
This equation can be reduced to 
$$
Y^2=X^3+X+2
$$
by making the substitutions $Y=x^2$ and $X=-y$. We have solved this equation previously, see Table \ref{tab:3.31}, and its integer solutions are $(X,Y)=(-1,0),(1,\pm 2)$.
Hence, in the original variables, we obtain that equation \eqref{x4py3pym2} has the unique integer solution
$$
(x,y)=(0,1).
$$

\vspace{10pt}

The next equation we will consider is
\begin{equation}\label{x4py3pyp2}
x^4+y^3+y+2=0.
\end{equation}
This equation can be reduced to 
$$
Y^2=X^3+X-2
$$
by making the substitutions $Y=x^2$ and $X=-y$. We have solved this equation previously, see Table \ref{tab:3.31}, and its only integer solution is $(X,Y)=(1,0)$.
 Hence, in the original variables, we obtain that equation \eqref{x4py3pyp2} has the unique integer solution
$$
(x,y)=(0,-1).
$$

\vspace{10pt}

The next equation we will consider is
\begin{equation}\label{x4py3m2y}
x^4+y^3-2y=0.
\end{equation}
This equation can be reduced to 
$$
Y^2=X^3-2X
$$
by making the substitutions $Y=x^2$ and $X=-y$. We have solved this equation previously, see Table \ref{tab:3.31}, and its integer solutions are $(X,Y)=(-1,\pm 1),$ $(0,0),$ $(2,\pm 2),$ $(338, \pm 6214)$.
Hence, in the original variables, we obtain that equation \eqref{x4py3m2y} has the integer solutions
$$
(x,y)=(0,0),(\pm 1,1).
$$

\vspace{10pt}

The next equation we will consider is
\begin{equation}\label{x4py3p2y}
x^4+y^3+2y=0.
\end{equation}
This equation can be reduced to 
$$
Y^2=X^3+2X
$$
by making the substitutions $Y=x^2$ and $X=-y$. We have solved this equation previously, see Table \ref{tab:3.31}, and its only integer solution is $(X,Y)=(0,0)$.
Hence, in the original variables, equation \eqref{x4py3p2y} has the unique integer solution
$$
(x,y)=(0,0).
$$

\vspace{10pt}

The next equation we will consider is
\begin{equation}\label{x4mx2py3}
x^4-x^2+y^3=0.
\end{equation}
This equation can be reduced to 
$$
Y^2+Y=X^3
$$
by making the substitutions $Y=-x^2$ and $X=-y$. We have solved this equation previously, see Table \ref{tab:3.31}, and its integer solutions are $(X,Y)=(0,-1),(0,0)$.
Hence, in the original variables, equation \eqref{x4mx2py3} has the integer solutions
$$
(x,y)=(0,0),(\pm 1,0).
$$

\vspace{10pt}

The final equation we will consider is
\begin{equation}\label{x4px2py3}
x^4+x^2+y^3=0.
\end{equation}
This equation can be reduced to 
$$
Y^2+Y=X^3
$$
by making the substitutions $Y=x^2$ and $X=-y$.  We have solved this equation previously, see Table \ref{tab:3.31}, and its integer solutions are $(X,Y)=(0,-1),(0,0)$. 
Hence, in the original variables, equation \eqref{x4px2py3} has the unique integer solution
$$
(x,y)=(0,0).
$$

Table \ref{Table:ex3.41} summarises integer solutions to the equations listed in Table \ref{tab:H28nonlinsub}.

\begin{center}
\begin{tabular}{|c|c||c|c|}
\hline
Equation & Solution $(x,y)$ &  Equation & Solution $(x,y)$ \\\hline \hline
$x^4+y^3-1=0$ & $(0, 1), (\pm1, 0)$ & $x^4+y^3-3=0$ & No integer solutions \\\hline
 $x^4+y^3+1=0$ & $(0,-1)$ &  $x^4 + y^3- 4 = 0$ & No integer solutions  \\\hline
  $x^4+y^3-y=0$ & $(0,0),(0,\pm 1)$ &  $x^4 + y^3 + 4 = 0$ & No integer solutions \\\hline
 $x^4+y^3+y=0$ & $(0,0)$ & $x^4+y^3+y-2=0$ & $(0,1)$ \\\hline
 $x^4+y^3+2=0$ & No integer solutions &  $x^4 + y^3 + y + 2 = 0$ & $(0,-1)$ \\\hline
 $x^4+y^3-2=0$ & $(\pm 1,1)$ &  $x^4+y^3-2y=0$ & $(0,0),(\pm 1,1)$ \\\hline
$x^4+y^3-y-1=0$ &  $(\pm 1,\pm1),(\pm 1,0)$ &  $x^4+y^3+2y=0$ & $(0,0)$ \\\hline
$x^4+y^3+y+1=0$ & $(\pm1,-1)$  & $x^4-x^2+y^3=0$ & $(0,0),(\pm 1,0)$ \\\hline
$x^4+y^3+y-1=0$&$(\pm1,0)$ & $x^4+x^2+y^3=0$ & $(0,0)$ \\\hline
$x^4+y^3+3=0$& No integer solutions && \\\hline
\end{tabular}
\captionof{table}{Integer solutions to the equations listed in Table \ref{tab:H28nonlinsub}.}\label{Table:ex3.41}
\end{center}

\subsection{Exercise 3.43}\label{ex:H28powerred}
\textbf{\emph{Solve all the equations from Table \ref{tab:H28powerred}.}}

	\begin{center}
		\begin{tabular}{ |c|c|c|c|c|c| } 
			\hline
			$H$ & Equation & $H$ & Equation & $H$ & Equation \\ 
			\hline\hline
			$25$ & $2y^2=x^4-1$ & $26$ & $y^3=x^4+x$ & $28$ & $y^3+y^2=2x^3$ \\ 
			\hline
			$26$ & $y^3=2x^3-x$ & $27$ & $y^3=2x^3-x+1$ & $28$ & $y^3=2x^3+x^2$ \\ 
			\hline
			$26$ & $y^3=2x^3+x$ & $27$ & $2y^2+y=x^4+1$ & $28$ & $2y^2-2y=x^4$ \\ 
			\hline
			$26$ & $y^3-y=2x^3$ & $28$ & $y^3=2x^3-2x$ & $28$ & $2y^2=x^4+2x$ \\ 
			\hline
			$26$ & $y^3+y=2x^3$ & $28$ & $y^3=2x^3+2x$ & $28$ & $2y^2=x^4+x-2$ \\ 
			\hline
			$26$ & $2y^2=x^4+2$ & $28$ & $y^3-2y=2x^3$ & $28$ & $y^3+y^2=x^4$ \\ 
			\hline
			$26$ & $2y^2=x^4+x$ & $28$ & $y^3+2y=2x^3$ & $28$ & $y^3-y^2=x^4$  \\ 
			\hline
			$26$ & $2y^2+y=x^4$ & $28$ & $y^3+y=2x^3+2$ & $28$ & $y^3=x^4+2x$ \\ 
			\hline
		\end{tabular}
		\captionof{table}{\label{tab:H28powerred} Equations of size $H\leq 28$ of the form $ay^n=P(x)$ or $ax^n=P(y)$ with reducible $P$.}
	\end{center} 
	
In this exercise we will solve equations of the form $ay^n = P(x)$ or $ax^n = P(y)$ with reducible $P$. Assuming we have $ay^n=P(x)$, then as $P$ is reducible, we can factorise $P(x)$, so we will assume that $P(x)=S(x)Q(x)$ where $S$ and $Q$ are polynomials. If $|a|=1$ and $S(x)$ and $Q(x)$ are coprime for every $x$, then we must have $\pm u^n=S(x)$ and $\pm v^n=Q(x)$ for some integers $u,v$, and these equations are often easier than the original one. If $|a|>1$ or $S(x)$ and $Q(x)$ can share some factors, we will combine this method with some case analysis. Many of the equations in Table \ref{tab:H28powerred} reduce to solving a system of equations, which then reduces to solving a Thue equation, see Section \ref{ex:Thue}. For a Thue equation $P(x,y)=0$, the Mathematica command
$$
{\tt Reduce[P(x,y)==0,\{x,y\},Integers]}
$$  
outputs all its integer solutions.

Equation 
$$
2y^2=x^4-1
$$
is solved in Section 3.3.5 of the book and its integer solutions are
$$
(x,y)=(\pm 1,0).
$$

Equation 
$$
y^3=2x^3-x
$$
is solved in Section 3.3.5 of the book and its integer solutions are
$$
(x,y)=(0,0),\pm(1,1).
$$

Equation 
$$
2y^2=x^4+2
$$
is solved in Section 3.3.5 of the book and its integer solutions are
$$
(x,y)=(0, \pm 1), (\pm 2,\pm 3).
$$

Equation 
$$
2y^2=x^4+x
$$
is solved in Section 3.3.5 of the book and its integer solutions are
$$
(x,y)=(-1,0),(0,0),(1,\pm 1),(2, \pm 3).
$$

Let us look at the first equation
\begin{equation}\label{y3e2x3px}
y^3=2x^3+x.
\end{equation}
This equation can be factorised as $y^3=x(2x^2+1)$. If $d$ is the greatest common divisor of $x$ and $2x^2+1$, then $d$ is also a divisor of $1$, hence $x$ and $2x^2+1$ are coprime. This implies that $x$ and $2x^2+1$ are both perfect cubes, say $x=u^3$ and $2x^2+1=v^3$. Then we can solve this system of equations. We obtain $2(u^2)^3+1=v^3$, which can be solved as a Thue equation in variables $u^2$ and $v$. It can be easily verified that the only integer solution is $(u,v)=(0,1)$, which implies that $x=0$. After substituting this into \eqref{y3e2x3px}, we can conclude that the unique integer solution to equation \eqref{y3e2x3px} is
$$
(x,y)=(0,0).
$$

\vspace{10pt}

The next equation we will consider is
\begin{equation}\label{y3mye2x3}
y^3-y=2x^3.
\end{equation}
This equation can be factorised as $y(y^2-1)=2x^3$. If $d$ is the greatest common divisor of $y$ and $y^2-1$, then $d$ is also a divisor of $1$, hence $y$ and $y^2-1$ are coprime. Then we must have either (i) $y=2u^3$ and $y^2-1=v^3$, or (ii) $y=u^3$ and $y^2-1=2v^3$. In case (i), the system reduces to $4(u^2)^3-1=v^3$ which can be solved as a Thue equation in variables $u^2$ and $v$. It can be easily verified that its only integer solution is $(u,v)=(0,-1)$, hence $y=0$ and \eqref{y3mye2x3} implies that $x=0$. In case (ii), the system reduces to $(u^2)^3-1=2v^3$ and it can be easily verified that the only integer solutions are $(u,v)=(\pm 1,0)$, hence $y=\pm 1$ and \eqref{y3mye2x3} implies that $x=0$. Finally, we can conclude that the integer solutions to equation \eqref{y3mye2x3} are
$$
(x,y)=(0,0),(0,\pm 1).
$$

\vspace{10pt}

The next equation we will consider is
\begin{equation}\label{y3pye2x3}
y^3+y=2x^3.
\end{equation}
This equation can be factorised as $y(y^2+1)=2x^3$. If $d$ is the greatest common divisor of $y$ and $y^2+1$, then $d$ is also a divisor of $1$, hence $y$ and $y^2+1$ are coprime. Then we must have either (i) $y=2u^3$ and $y^2+1=v^3$, or (ii) $y=u^3$ and $y^2+1=2v^3$. In case (i), the system reduces to $4(u^2)^3+1=v^3$ which can be solved as a Thue equation, and it can be easily verified that its only integer solution is $(u,v)=(0,1)$ hence $y=0$, and \eqref{y3pye2x3} implies that $x=0$. Similarly, in case (ii), the system reduces to $(u^2)^3+1=2v^3$ whose only integer solutions are $(u,v)=(\pm 1,1)$ hence $y=\pm 1$ and \eqref{y3pye2x3} implies the following solutions $(x,y)=\pm (1,1)$. Finally, we can conclude that the integer solutions to equation \eqref{y3pye2x3} are
$$
(x,y)=(0,0),\pm (1,1).
$$

\vspace{10pt}

The next equation we will consider is
\begin{equation}\label{2y2pyex4}
2y^2+y=x^4.
\end{equation}
This equation can be factorised as $y(2y+1)=x^4$. If $d$ is the greatest common divisor of $y$ and $2y+1$, then $d$ is also a divisor of $1$, hence $y$ and $2y+1$ are coprime. Then we must have either (i) $y=u^4$ and $2y+1=v^4$, or (ii) $y=-u^4$ and $2y+1=-v^4$. In case (i), the system reduces to $2u^4+1=v^4$ and it can be easily verified that its only integer solutions are $(u,v)=(0,\pm1)$, hence $y=0$ and \eqref{2y2pyex4} implies that $x=0$. In case (ii), the system reduces to $-2u^4+1=-v^4$ whose only integer solutions are $(u,v)=(\pm 1,\pm 1),$ hence $y=-1$ and \eqref{2y2pyex4} implies that $x=\pm 1$. Finally, we can conclude that the integer solutions to equation \eqref{2y2pyex4} are
$$
(x,y)=(0,0),(\pm 1,-1).
$$

\vspace{10pt}

The next equation we will consider is
\begin{equation}\label{y3ex4px}
y^3=x^4+x.
\end{equation}
This equation can be factorised as $y^3=x(x^3+1)$. Because $x$ and $x^3+1$ are coprime, we must have that they are both perfect cubes, say, $x=u^3$ and $x^3+1=v^3$. The second equation is a Thue equation, whose only integer solutions are $(x,v)=(-1,0),(0,1)$.  Then \eqref{y3ex4px} implies that all its
integer solutions
are
$$
(x,y)=(-1,0),(0,0).
$$

\vspace{10pt}

The next equation we will consider is
\begin{equation}\label{y3e2x3mxp1}
y^3=2x^3-x+1.
\end{equation}
This equation can be factorised as $y^3=(x+1)(2x^2-2x+1)$. If $d>0$ is the greatest common divisor of $x+1$ and $2x^2-2x+1$, then $d$ is also a divisor of $(2x^2-2x+1)-(x+1)(2x-4)=5$.
So either $d=1$ or $d=5$. If $d=1$ then this implies that $x+1$ and $2x^2-2x+1$ are both perfect cubes, say $x+1=u^3$ and $2x^2-2x+1=v^3$. The last equation is of the form \eqref{eq:Weiformgen} and its only integer solutions are $(v,x)=(1,0),(1,1)$. If $x=1$, then $u^3=x+1=2$, a contradiction, hence we must have $x=0$ which implies that $y=1$.
If $d=5$, then either $2x^2-2x+1=5v^3$ or $2x^2-2x+1=25v^3$. The only integer solutions to the first equation are $(v, x) = (1, -1), (1, 2)$ and the only integer solutions to the second equation are $(v,x)=(1, -3), (1, 4)$. After substituting these values of $x$ into \eqref{y3e2x3mxp1}, we can conclude that all integer solutions to equation \eqref{y3e2x3mxp1} are
$$
(x,y)=(-1, 0), (0, 1),(4,5).
$$

\vspace{10pt}

The next equation we will consider is
\begin{equation}\label{2y2pyex4p1}
2y^2+y=x^4+1.
\end{equation}
This equation can be factorised as $(2y-1)(y+1)=x^4$. If $d>0$ is the greatest common divisor of $y+1$ and $2y-1$, then $d$ is also a divisor of $2(y+1)-(2y-1)=3$. So, either $d=1$ or $d=3$. If $d=1$, then we must have either (a) $2y-1=u^4$ and $y+1=v^4$, or (b) $2y-1=-u^4$ and $y+1=-v^4$. These systems reduce to Thue equations and it can be easily verified that they have no integer solutions. If $d=3$, then we must have either (i) $2y-1=3u^4$ and $y+1=27v^4$, or (ii) $2y-1=-3u^4$ and $y+1=-27v^4$; or (iii) $2y-1=9u^4$ and $y+1=9v^4$, or (iv) $2y-1=-9u^4$ and $y+1=-9v^4$; or (v) $2y-1=27u^4$ and $y+1=3v^4$, or (vi) $2y-1=-27u^4$ and $y+1=-3v^4$.
  These systems reduce to Thue equations. Cases (i) and (iii)-(vi) have no integer solutions. Case (ii) has the only integer solutions $(u,v)=(\pm 1,0)$, hence $y=-1$. In which case \eqref{2y2pyex4p1} implies that $x=0$. Finally, we can conclude that the unique integer solution to equation \eqref{2y2pyex4p1} is
$$
(x,y)=(0,-1).
$$

\vspace{10pt}

The next equation we will consider is
\begin{equation}\label{y3e2x3m2x}
y^3=2x^3-2x.
\end{equation}
This equation can be factorised as $y^3=x(2x^2-2)$. If $d>0$ is the greatest common divisor of $x$ and $2x^2-2$, then $d$ is also a divisor of $2x^2-(2x^2-2)=2$. If $d=1$, then this implies that $x$ and $2x^2-2$ must both be perfect cubes, say $x=u^3$ and $2x^2-2=v^3$. Solving this system reduces to a Thue equation which has the only integer solutions $(u,v)=(\pm 1,0)$ hence $x=\pm 1$, and \eqref{y3e2x3m2x} implies that $y=0$. If $d=2$, then we must have either (i) $x=2u^3$ and $2x^2-2=4v^3$, or (ii) $x=4u^3$ and $2x^2-2=2v^3$.
   Case (i) reduces to $8(u^2)^3-2=4v^3$ which is a Thue equation and it has no integer solutions. Case (ii) reduces to $32(u^2)^3-2=2v^3$ which is a Thue equation and whose only integer solution is $(u,v)=(0,-1)$, hence $x=0$, and \eqref{y3e2x3m2x} implies that $y=0$. Finally, we can conclude that the integer solutions to equation \eqref{y3e2x3m2x} are
$$
(x,y)=(\pm 1,0),(0,0).
$$

\vspace{10pt}

The next equation we will consider is
\begin{equation}\label{y3e2x3p2x}
y^3=2x^3+2x.
\end{equation}
This equation can be factorised as $y^3=x(2x^2+2)$. If $d>0$ is the greatest common divisor of $x$ and $2x^2+2$, then $d$ is also a divisor of $(2x^2+2)-2x^2=2$. If $d=1$, then we must have that $x$ and $2x^2+2$ are both perfect cubes, say $x=u^3$ and $2x^2+2=v^3$. This system reduces to $2(u^2)^3+2=v^3$ which can be solved as a Thue equation and it has no integer solutions. If $d=2$, then we must have either (i) $x=2u^3$ and $2x^2+2=4v^3$, or (ii) $x=4u^3$ and $2x^2+2=2v^3$.
   Case (i) reduces to $8(u^2)^3+2=4v^3$ which has no integer solutions. Case (ii) reduces to $32(u^2)^3+2=2v^3$ which has the only integer solution $(u,v)=(0,1)$ hence $x=0$, and \eqref{y3e2x3p2x} implies that $y=0$. We can then conclude that the unique integer solution to equation \eqref{y3e2x3p2x} is
$$
(x,y)=(0,0).
$$

\vspace{10pt}

The next equation we will consider is
\begin{equation}\label{y3m2ye2x3}
y^3-2y=2x^3.
\end{equation}
This equation can be factorised as $y(y^2-2)=2x^3$. If $d>0$ is the greatest common divisor of $y$ and $y^2-2$, then $d$ is also a divisor of $y^2-(y^2-2)=2$. If $d=1$, then we must have either (i) $y=u^3$ and $y^2-2=2v^3$, or (ii) $y=2u^3$ and $y^2-2=v^3$. Case (i) reduces to $(u^2)^3-2=2v^3$ which has the only integer solution $(u,v)=(0,-1)$ which implies that $y=0$, and \eqref{y3m2ye2x3} implies that $x=0$. Case (ii) reduces to $4(u^2)^3-2=v^3$ which has no integer solutions. If $d=2$, then we have the additional case $y=4u^3$ and $y^2-2=4v^3$ which reduces to $16(u^2)^3-2=4v^3$ which has no integer solutions. Finally, we can conclude that the unique integer solution to equation \eqref{y3m2ye2x3} is
$$
(x,y)=(0,0).
$$

\vspace{10pt}

The next equation we will consider is
\begin{equation}\label{y3p2ye2x3}
y^3+2y=2x^3.
\end{equation}
This equation can be factorised as $y(y^2+2)=2x^3$. If $d>0$ is the greatest common divisor of $y$ and $y^2+2$, then $d$ is also a divisor of $(y^2+2)-y^2=2$. If $d=1$, then we must have either (i) $y=u^3$ and $y^2+2=2v^3$, or (ii) $y=2u^3$ and $y^2+2=v^3$. Case (i) reduces to $(u^2)^3+2=2v^3$ and can be solved as a Thue equation whose only integer solution is $(u,v)=(0,1)$ hence $y=0$ and \eqref{y3p2ye2x3} implies that $x=0$. Case (ii) reduces to $4(u^2)^3+2=v^3$ which has no integer solutions. If $d=2$, then we have the additional case $y=4u^3$ and $y^2+2=4v^3$, which reduces to $16(u^2)^3+2=4v^3$ and it has no integer solutions. Finally, we can conclude that the unique integer solution to equation \eqref{y3p2ye2x3} is
$$
(x,y)=(0,0).
$$

\vspace{10pt}

The next equation we will consider is
\begin{equation}\label{y3pye2x3p2}
y^3+y=2x^3+2.
\end{equation}
This equation can be factorised as $(y-1)(y^2+y+2)=2x^3$. If $d>0$ is the greatest common divisor of $y-1$ and $y^2+y+2$, then $d$ is also a divisor of $y^2+y+2-(y-1)(y+2)=4$. If $d=1$, then we must have either (i) $y-1=u^3$ and $y^2+y+2=2v^3$, or (ii) $y-1=2u^3$ and $y^2+y+2=v^3$. Case (i) can be reduced to $t^2+6t+16=s^3$ by the substitutions $t=2u^3$ and $s=2v$, which is an elliptic curve in Weierstrass form which has integer points $(s,t)=(2,-4),(2,-2),(32,-184),(32,178)$ hence $y=-91,-1,0,90$. Similarly, case (ii) can be reduced to $t^2+3t+4=v^3$, which is an elliptic curve in Weierstrass form, by the substitution $t=2u^3$ hence $y$ as $y=-3,2$. If $d=2$, then we have the additional case $y-1=4u^3$ and $y^2+y+2=4v^3$, which can be reduced to $t^2+12t+64=s^3$, which is an elliptic curve in Weierstrass form, by the substitutions $t=16u^3$ and $s=4v$, hence $y$ as $y=-6,-2,1,5$. If $d=4$, we have no further cases to check. After substituting each value of $y$ into \eqref{y3pye2x3p2}, we can then conclude that all integer solutions to equation \eqref{y3pye2x3p2} are
$$
(x,y)=(-1,0),(0,1),(4,5).
$$

\vspace{10pt}

The next equation we will consider is
\begin{equation}\label{y3py2e2x3}
y^3+y^2=2x^3.
\end{equation}
This equation can be factorised as $y^2(y+1)=2x^3$. Because $y^2$ and $y+1$ are coprime, we must have either (i) $y^2=2u^3$ and $y+1=v^3$ or (ii) $y^2=u^3$ and $y+1=2v^3$. Case (i) reduces to $(v^3-1)^2=2u^3$, which after the substitutions $V=2v^3$ and $U=2u$, can be further reduced to $(V-2)^2=U^3$, whose integer solutions are $(U,V)=(k^2,k^3+2)$ for integer $k$, see \eqref{eq:x3my2sol}. To ensure $u,v$ are integers, $k$ must be even, say $k=2t$, then $(u,v)=(2t^2,\sqrt[3]{4t^3+1})$, which is only integer if $t=0$, hence we obtain the solution $(x,y)=(0,0)$ to \eqref{y3py2e2x3}. Similarly, case (ii) reduces to $(2v^3-1)^2=u^3$, which after the substitution $V=2v^3$, reduces the equation further to $(V-1)^2=u^3$, whose integer solutions are $(u,V)=(k^2,k^3+1)$ for integer $k$, then $2v^3=k^3+1$ which has the integer solutions $(v,k)=(0,-1),(1,1)$, hence $y=\pm 1$. After checking these cases of $y$, we can conclude that the integer solutions to equation \eqref{y3py2e2x3} are 
$$
(x,y)=(0,-1),(0,0),(1,1).
$$

\vspace{10pt}

The next equation we will consider is
\begin{equation}\label{y3e2x3px2}
y^3=2x^3+x^2.
\end{equation}
This equation can be factorised as $y^3=x^2(2x+1)$. Because $x^2$ and $2x+1$ are coprime, we must have that $x^2$ and $2x+1$ are both perfect cubes, say $x^2=u^3$ and $2x+1=v^3$. The first equation has the integer solutions $(x,u)=(k^3,k^2)$ for integer $k$. Then $2k^3+1=v^3$, which is a Thue equation, whose only integer solutions are $(k,v)=(-1,-1),(0,1)$, hence $x=0$ or $x=-1$. After substituting these values of $x$ into \eqref{y3e2x3px2} and solving for $y$, we can conclude that the integer solutions to equation \eqref{y3e2x3px2} are
$$
(x,y)=(-1,-1),(0,0).
$$

\vspace{10pt}

The next equation we will consider is
\begin{equation}\label{2y2m2yex4}
2y^2-2y=x^4.
\end{equation}
This equation can be factorised as $y(2y-2)=x^4$. If $d>0$ is the greatest common divisor of $y$ and $2y-2$, then $d$ is also a divisor of $2y-(2y-2)=2$. If $d=1$, we must have either (a) $y=u^4$ and $2y-2=v^4$ or (b) $y=-u^4$ and $2y-2=-v^4$. By solving these systems, we obtain that $y=1$, and \eqref{2y2m2yex4} implies that $x=0$. If $d=2$ we must have either 
(i) $y=2u^4$ and $2y-2=8v^4$ or (ii) $y=-2u^4$ and $2y-2=-8v^4$;
(iii) $y=4u^4$ and $2y-2=4v^4$ or (iv) $y=-4u^4$ and $2y-2=-4v^4$;
(v) $y=8u^4$ and $2y-2=2v^4$ or (vi) $y=-8u^4$ and $2y-2=-2v^4$.
By solving these systems, we obtain that $y=0$ and \eqref{2y2m2yex4} implies that $x=0$. Hence, all integer solutions to equation \eqref{2y2m2yex4} are
$$
(x,y)=(0,0),(0,1).
$$

\vspace{10pt}

The next equation we will consider is
\begin{equation}\label{2y2ex4p2x}
2y^2=x^4+2x.
\end{equation}
This equation can be factorised as $2y^2=x(x^3+2)$. If $d>0$ is the greatest common divisor of $x$ and $x^3+2$, then $d$ is also a divisor of $2$. If $d=1$ we must have either (i) $x=u^2$ and $x^3+2=2v^2$ or (ii) $x=2u^2$ and $x^3+2=v^2$ or (iii) $x=-u^2$ and $x^3+2=-2v^2$ or (iv) $x=-2u^2$ and $x^3+2=-v^2$. By solving these systems, 
we obtain that $x=0$, and \eqref{2y2ex4p2x} implies $y=0$. If $d=2$ then we have no additional cases to check. Finally, we can conclude that the unique integer solution to equation \eqref{2y2ex4p2x} is
$$
(x,y)=(0,0).
$$

\vspace{10pt}
The next equation we will consider is
\begin{equation}\label{2y2ex4pxm2}
2y^2 = x^4 + x -2.
\end{equation}
This equation can be factorised as $2y^2=(x-1)(x^3+x^2+x+2)$. If $d>0$ is the greatest common divisor of $x-1$ and $x^3+x^2+x+2$, then $d$ is also a divisor of 
$(x^3+x^2+x+2)-(x-1)(x^2+2x+3)=5$, hence $d=1$ or $d=5$. If $d=1$, then $x^3+x^2+x+2$ can be equal to $\pm v^2$ or $\pm 2v^2$ for some integer $v$.
These cases lead to equations in Weierstrass form \eqref{eq:Weiform}, but produce no integer solutions $(x,v)$ for which $y=\pm \sqrt{(x^4+x-2)/2}$ is an integer. If $d=5$, then $x^3+x^2+x+2$ can also be equal to $\pm 5v^2$ or $\pm 10v^2$. The case $x^3+x^2+x+2=5v^2$ has integer solution $(x,v)=(1,\pm 1)$ for which $y=\pm \sqrt{(1^4+1-2)/2}=0$ is an integer, and other cases have no such integer solutions.
Hence, the unique integer solution to equation \eqref{2y2ex4pxm2} is
$$
(x,y)=(1,0).
$$

\vspace{10pt}

The next equation we will consider is
\begin{equation}\label{y3py2ex4}
 y^3 + y^2 = x^4.
\end{equation}
This equation can be factorised as $x^4=y^2(y+1)$. Because $y^2$ and $y+1$ are coprime, we must have either (i) $y^2=u^4$ and $y+1=v^4$ or (ii) $y^2=-u^4$ and $y+1=-v^4$. These systems reduce to Thue equations and we obtain that $y=-1,0$, and \eqref{y3py2ex4} implies that $x=0$. Finally, we can conclude that the integer solutions to equation \eqref{y3py2ex4} are
$$
(x,y)=(0,-1),(0,0).
$$

\vspace{10pt}

The next equation we will consider is
\begin{equation}\label{y3my2ex4}
  y^3 - y^2 = x^4.
\end{equation}
This equation can be factorised as $x^4=y^2(y-1)$. Because $y^2$ and $y-1$ are coprime, we must have either (i) $y^2=u^4$ and $y-1=v^4$ or (ii) $y^2=-u^4$ and $y-1=-v^4$. These systems reduce to Thue equations and we obtain that $y=0,1$ and \eqref{y3my2ex4} implies that $x=0$. Finally, we can conclude that the integer solutions to equation \eqref{y3my2ex4} are
$$
(x,y)=(0,0),(0,1).
$$

\vspace{10pt}

The final equation we will consider is
\begin{equation}\label{y3ex4p2x}
y^3 = x^4 + 2x.
\end{equation}
This equation can be factorised as $y^3=x(x^3+2)$. If $d>0$ is the greatest common divisor of $x$ and $x^3+2$, then $d$ is also a divisor of $2$. If $d=1$ then we must have $x=u^3$ and $x^3+2=v^3$. The second equation is a Thue equation whose only integer solution is $(x,v)=(-1,1)$, and \eqref{y3ex4p2x} implies that $y=-1$. If $d=2$ then we have the additional cases (i) $x=2u^3$ and $x^3+2=4v^3$ or (ii) $x=4u^3$ and $x^3+2=2v^3$ to check. The second equation in case (i) is a Thue equation and it has no integer solutions. The second equation in case (ii) is a Thue equation whose only integer solution is $(x,v)=(0,1)$ and \eqref{y3ex4p2x} implies that $y=0$. Finally, we can conclude that the integer solutions to equation \eqref{y3ex4p2x} are 
$$
(x,y)=(-1,-1),(0,0).
$$

 \begin{center}
\begin{tabular}{|c|c||c|c|c|c|}
		\hline
		Equation & Solution $(x,y)$ & Equation & Solution $(x,y)$  \\\hline \hline
		$2y^2=x^4-1$&$(\pm 1,0)$ & $y^3=2x^3+2x$&$(0,0)$ \\\hline
		$y^3=2x^3-x$&$(0,0),\pm(1,1)$ &  $y^3-2y=2x^3$&$(0,0)$ \\\hline
		$y^3=2x^3+x$&$(0,0)$ & $y^3+2y=2x^3$&$(0,0)$ \\\hline
		$y^3-y=2x^3$ &$(0,0),(0,\pm 1)$ & $y^3+y=2x^3+2$&$(-1,0),(0,1),(4,5)$ \\\hline
		$y^3+y=2x^3$ &$(0,0),\pm(1,1)$ & $y^3+y^2=2x^3$&$(0,-1),(0,0),(1,1)$   \\\hline
		$2y^2=x^4+2$&$(0, \pm1),(\pm 2,\pm 3)$ & $y^3=2x^3+x^2$&$(-1,-1),(0,0)$ \\\hline
		$2y^2=x^4+x$&$(-1,0),(0,0),(1,\pm1),(2,\pm 3)$ &$2y^2-2y=x^4$&$(0,0),(0,1)$ \\\hline
		$2y^2+y=x^4$&$(0,0),(\pm 1,-1)$ & $2y^2=x^4+2x$&$(0,0)$ \\\hline
		$y^3=x^4+x$&$(-1,0),(0,0)$ & $2y^2 = x^4 + x - 2$&$(1,0)$  \\\hline
		$y^3=2x^3-x+1$&    $(-1, 0), (0, 1),(4,5)$ & $ y^3 + y^2 = x^4 $&$(0,-1),(0,0)$  \\\hline
		$2y^2+y=x^4+1$&$(0,-1)$ &$  y^3 - y^2 = x^4 $&$(0,0),(0,1)$ \\\hline
		$y^3=2x^3-2x$&$(\pm1,0),(0,0)$ & $ y^3 = x^4 + 2x $ &$(-1,-1),(0,0)$ \\\hline
	\end{tabular}
	\captionof{table}{Integer solutions to the equations in Table \ref{tab:H28powerred}.}\label{Table:ex3.42}
\end{center}

\subsection{Exercise 3.46}\label{ex:H28liny}
\textbf{\emph{Solve the equations listed in Table \ref{tab:H28liny}.  }}

\begin{center}
\begin{tabular}{ |c|c|c|c|c|c| } 
 \hline
 $H$ & Equation & $H$ & Equation & $H$ & Equation \\ 
 \hline\hline
 $26$ & $x^3+2x^2y-y=0$ & $27$ & $x^4+2xy+y+1=0$ & $28$ & $x^4+xy-3y+x=0$ \\ 
 \hline
 $26$ & $x^3+2x^2y+y=0$ & $28$ & $x^3+2x^2y-y-2=0$ & $28$ & $x^4+xy+3y+x=0$ \\ 
 \hline
 $26$ & $x^4+xy+3y=0$ & $28$ & $x^3+2x^2y-y+2=0$ & $28$ & $x^4+2xy+y-2=0$ \\ 
 \hline
 $26$ & $x^4+2xy+y=0$ & $28$ & $x^3+2x^2y-x+y=0$ & $28$ & $x^4+2xy+y+2=0$ \\ 
 \hline
 $27$ & $x^3+2x^2y-y+1=0$ & $28$ & $x^3+2x^2y+x-y=0$ & $28$ & $x^4+2xy+x-y=0$ \\ 
 \hline
 $27$ & $x^3+2x^2y+y+1=0$ & $28$ & $x^3+2x^2y+2y=0$ & $28$ & $x^4+2xy+x+y=0$ \\ 
 \hline
 $27$ & $x^4+xy+3y-1=0$ & $28$ & $x^4+xy+3y-2=0$ & $28$ & $x^4+x^2y+x+y=0$ \\ 
 \hline
 $27$ & $x^4+xy+3y+1=0$ & $28$ & $x^4+xy+3y+2=0$ &  & \\ 
 \hline
 $27$ & $x^4+2xy+y-1=0$ & $28$ & $x^4+xy+4y=0$ &  & \\ 
 \hline
\end{tabular}
\captionof{table}{\label{tab:H28liny} Equations linear in $y$ of size $H\leq 28$.}
\end{center} 

We will now consider equations which are linear in $y$. To solve these equations, we first need to put the equation in terms of $y$, so 
$$
y=\frac{P(x)}{Q(x)}. 
$$
The cases when $Q(x)=0$ should be considered separately. If $Q(x)\neq 0$, then we can do polynomial division, which can either be done by hand, or by using the Mathematica command 
$$
{\tt PolynomialQuotientRemainder[P(x),Q(x),x]} \quad \text{where } P(x) \text{ is the polynomial, and } Q(x) \text{ is the divisor.}
$$
We will then have an expression of the form,
$$
y=P_1(x)+\frac{P_2(x)}{Q(x)}
$$
If the coefficients of $P_1(x)$ are not integer, we must multiply by a coefficient, say $d$, in order to have integer coefficients. We will then have 
$$
dy=dP_1(x)+\frac{dP_2(x)}{Q(x)}
$$
In order to have integer solutions to the original equation, we must find when $\frac{dP_2(x)}{Q(x)}$ is integer. This will give a finite number of cases for $x$, and then we can substitute these values of $x$ into the original equation, to determine all integer solutions.

Equation
$$
x^3+2x^2y-y=0
$$
is solved in Section 3.3.7 of the book (in the equivalent form $x^3-2x^2y+y=0$), and its integer solutions are
$$
(x,y)=(0,0),\pm(1,-1).
$$ 

Equation
$$
x^4+2xy+y=0
$$
is solved in Section 3.3.7 of the book (in the equivalent form $x^4-2xy-y=0$), and its integer solutions are
$$
(x,y)=(-1,1),(0,0).
$$

The first equation we will consider is
\begin{equation}\label{x3p2x2ypy}
x^3 + 2x^2y + y = 0.
\end{equation}
Because $2x^2+1 \neq 0$, we can rearrange this equation to $y=-\frac{x^3}{2x^2+1}= \frac{x/2}{2x^2+1}-\frac{x}{2}$. 
After multiplying by $2$ to get integer coefficients, we obtain $2y=\frac{x}{2x^2+1}-x$. Hence, we need to determine when $\frac{x}{2x^2+1}$ is integer. As $|x|<2x^2+1$, this is only possible if
$x=0$, and then \eqref{x3p2x2ypy} implies that $y=0$. 
Finally, we can conclude that the unique integer solution to equation \eqref{x3p2x2ypy} is
$$
(x,y)=(0,0).
$$

\vspace{10pt}

The next equation we will consider is
\begin{equation}\label{x4pxyp3y}
x^4 + xy + 3y = 0.
\end{equation}
This equation has no integer solutions with $x=-3$, so we may assume that $x \neq -3$. Then $x+3 \neq 0$, and we can rearrange this equation to 
$y=-\frac{x^4}{x+3}=-x^3+3x^2-9x+27-\frac{81}{x+3}$. 
Hence, we need to determine when $\frac{81}{x+3}$ is integer, so, $x+3$ must be a divisor of 81. Hence, 
$$
x+3=\pm 1, \pm 3, \pm 9, \pm 27, \pm 81,
$$
or equivalently,
$$
x=-84,-30,-12,-6,-4,-2,0,6,24,78.
$$
 After checking these values of $x$, we can conclude that the integer solutions to equation \eqref{x4pxyp3y} are
$$
\begin{aligned}
(x,y)=(-84,614656), (-30,30000), (-12,2304), (-6,432), (-4,256),\\ (-2,-16), (0,0), (6,-144), (24,-12288), (78,-456976).
\end{aligned}
$$

\vspace{10pt}

The next equation we will consider is 
\begin{equation}\label{x3p2x2ymyp1}
x^3+2x^2y-y+1=0.
\end{equation}
Because $x$ is integer, $2x^2-1 \neq 0$. We can rearrange \eqref{x3p2x2ymyp1} to 
 $y=-\frac{x^3+1}{2x^2-1}=-\frac{x}{2}-\frac{1+x/2}{2x^2-1}$. After multiplying by $2$ to get integer coefficients we obtain $2y=-x-\frac{2+x}{2x^2-1}$. To ensure that $y$ is integer,
 $\frac{2+x}{2x^2-1}$ must be integer. Hence, we must have either $|2+x|\geq 2x^2-1$ or $2+x=0$. The inequality is true for integers $|x| \leq 1$, and the equality implies that $x=-2$. After checking these cases for $x$, we can conclude that the integer solutions to equation \eqref{x3p2x2ymyp1} are
$$
(x,y)=(-2,1),(-1,0),(0,1),(1,-2).
$$

\vspace{10pt}

The next equation we will consider is
\begin{equation}\label{x3p2x2ypyp1}
x^3+2x^2y+y+1=0.
\end{equation}
Because $2x^2+1 \neq 0$, we can rearrange this equation to
$y=-\frac{x^3+1}{2x^2+1}=-\frac{x}{2}+\frac{-1+x/2}{2x^2+1}$. After multiplying by $2$ to get integer coefficients, we obtain $2y=-x+\frac{x-2}{2x^2+1}$. To ensure that $y$ is integer, $\frac{x-2}{2x^2+1}$ must be integer. Hence, we must have either $|x-2|\geq 2x^2+1$ or $x-2=0$. The inequality is true for integers $x=-1$ and $0$, and the equality implies that $x=2$. After checking these cases for $x$, we can conclude that the integer solutions to equation \eqref{x3p2x2ypyp1} are
$$
(x,y)=(-1,0),(0,-1),(2,-1).
$$

\vspace{10pt}

The next equation we will consider is 
\begin{equation}\label{x4pxyp3ym1}
x^4 + xy + 3y -1 = 0.
\end{equation}
This equation has no integer solutions with $x=-3$, so we may assume that $x \neq -3$. We can then rearrange the equation to
$y=\frac{1-x^4}{x+3}=-x^3+3x^2-9x+27-\frac{80}{x+3}$. To ensure that $y$ is integer, $\frac{80}{x+3}$ must be an integer. Hence,
 after checking the finite number of values of $x$ such that $x+3$ is a divisor of $80$, we can conclude that the integer solutions to equation \eqref{x4pxyp3ym1} are
$$
\begin{aligned}
(x,y)= (-83,593229), (-43,85470), (-23,13992), (-19,8145), (-13,2856),   (-11,1830), \\  (-8,819), (-7,600), (-5,312), (-4,255), (-2,-15), (\pm 1,0), (2,-3), \\  (5,-78), (7,-240), (13,-1785), (17,-4176), (37,-46854), (77,-439413).
\end{aligned}
$$

\vspace{10pt}

The next equation we will consider is
\begin{equation}\label{x4pxyp3yp1}
	x^4 + xy + 3y +1 = 0.
\end{equation}
Because $x^4+1>0$, we may assume that $x \neq -3$. Then we can rearrange the equation to
$y=-\frac{1+x^4}{x+3}=-x^3+3x^2-9x+27-\frac{82}{x+3}$.  To ensure that $y$ is integer, $\frac{82}{x+3}$ must be integer. Hence, 
after checking the finite number of values of $x$ such that $x+3$ is a divisor of $82$,
we can conclude that the integer solutions to equation \eqref{x4pxyp3yp1} are
$$
\begin{aligned}
	(x,y)=(-85,636593), (-44,91417), (-5,313), (-4,257), \\(-2,-17), (-1,-1), (38,-50857), (79,-475001).
\end{aligned}
$$

\vspace{10pt}

The next equation we will consider is 
\begin{equation}\label{x4p2xypym1}
x^4+2xy+y-1=0.
\end{equation}
Because $x$ is integer, $2x+1 \neq 0$, so we can rearrange the equation to 
 $y=\frac{1-x^4}{2x+1}=-\frac{x^3}{2}+\frac{x^2}{4}-\frac{x}{8}+\frac{1}{16}+\frac{15/16}{2x+1}$. After multiplying by $16$ to get integer coefficients, we obtain
$$
16y=-8x^3+4x^2-2x+1+\frac{15}{2x+1}.
$$
To ensure integer solutions, we must have that $\frac{15}{2x+1}$ is integer. After checking the finite number of values of $x$ such that $2x+1$ is a divisor of $15$,
 we can conclude that the integer solutions to equation \eqref{x4p2xypym1} are
$$
\begin{aligned}
(x,y)=(-8,273),(-3,16),(-2,5),(0,1),(\pm 1,0),(2,-3),(7,-160).
\end{aligned}
$$

\vspace{10pt}

The next equation we will consider is 
\begin{equation}\label{x4p2xypyp1}
x^4+2xy+y+1=0.
\end{equation}
Because $x$ is integer, $2x+1 \neq 0$, so we can rearrange the equation to 
$y=-\frac{1+x^4}{2x+1}=-\frac{x^3}{2}+\frac{x^2}{4}-\frac{x}{8}+\frac{1}{16}-\frac{17/16}{2x+1}$. After multiplying by $16$ to get integer coefficients, we obtain
$$
16y=-8x^3+4x^2-2x+1-\frac{17}{2x+1}.
$$
To ensure integer solutions, we must have that $\frac{17}{2x+1}$ is integer. After checking the finite number of values of $x$ such that $2x+1$ is a divisor of $17$,
 we can conclude that the integer solutions to equation \eqref{x4p2xypyp1} are
$$
\begin{aligned}
(x,y)=(-9,386),(-1,2),(0,-1),(8,-241).
\end{aligned}
$$
 
 \vspace{10pt}

 The next equation we will consider is 
\begin{equation}\label{x3p2x2ymym2}
x^3+2x^2y-y-2=0.
\end{equation}
Because $x$ is integer, $2x^2-1 \neq 0$, so we can rearrange the equation to 
$y=\frac{2-x^3}{2x^2-1}=-\frac{x}{2}+\frac{2-x/2}{2x^2-1}$. After multiplying by $2$ we obtain $2y=-x+\frac{4-x}{2x^2-1}$. To ensure integer solutions, we must have that $\frac{4-x}{2x^2-1}$ is integer. Hence, we must have either $|4-x|\geq 2x^2-1$ or $4-x=0$. The inequality is true for integers $|x|\leq 1$, and the equality implies that $x=4$. After checking these cases of $x$, we can conclude that the integer solutions to equation \eqref{x3p2x2ymym2} are
$$
(x,y)=(-1,3),(0,-2),(1,1),(4,-2).
$$
 
 \vspace{10pt}
 
  The next equation we will consider is 
\begin{equation}\label{x3p2x2ymyp2}
x^3+2x^2y-y+2=0.
\end{equation}
Because $x$ is integer, $2x^2-1 \neq 0$, so we can rearrange the equation to
 $y=-\frac{2+x^3}{2x^2-1}=-\frac{x}{2}-\frac{2+x/2}{2x^2-1}$. After multiplying by $2$ we obtain $2y=-x-\frac{4+x}{2x^2-1}$. To ensure integer solutions, we must have that 
 $\frac{4+x}{2x^2-1}$ is integer. Hence, we must have either $|4+x|\geq 2x^2-1$ or $4+x=0$. The inequality is true for integers $|x|\leq 1$, and the equality implies that $x=-4$. After checking these values of $x$, we can conclude that the integer solutions to equation \eqref{x3p2x2ymyp2} are
$$
(x,y)=(-4,2),(-1,-1),(0,2),(1,-3).
$$
 
 \vspace{10pt}
 
 The next equation we will consider is
 \begin{equation}\label{x3p2x2ymxpy}
 x^3+2x^2y-x+y=0.
 \end{equation}
Because $2x^2+1 \neq 0$, we can rearrange the equation to
 $y=\frac{x-x^3}{2x^2+1}=-\frac{x}{2}+\frac{3x/2}{2x^2+1}$. After multiplying by $2$ we obtain $2y=-x+\frac{3x}{2x^2+1}$. 
 To ensure integer solutions, we must have that $\frac{3x}{2x^2+1}$ is integer. Hence, we must have either $|3x| \geq 2x^2+1$ or $3x=0$. The inequality is true for integers $x=\pm 1$, and the equality implies that $x=0$. After checking these values of $x$, we can conclude that the integer solutions to equation \eqref{x3p2x2ymxpy} are
$$
(x,y)= (0,0),(\pm 1,0).
$$
 
 \vspace{10pt}

 The next equation we will consider is
 \begin{equation}\label{x3p2x2ypxmy}
 x^3+2x^2y+x-y=0.
 \end{equation}
Because $x$ is integer, $2x^2-1 \neq 0$, so we can rearrange the equation to 
  $y=-\frac{x+x^3}{2x^2-1}=-\frac{x}{2}-\frac{3x/2}{2x^2-1}$. After multiplying by $2$ we obtain
  $2y=-x-\frac{3x}{2x^2-1}$. 
   To ensure integer solutions, we must have that $\frac{3x}{2x^2-1}$ is integer. Hence, we must have either $|3x| \geq 2x^2-1$ or $3x=0$. The inequality is true for integers $|x|\leq 1$, and the equality implies that $x=0$. After checking these values of $x$, we can conclude that the integer solutions to equation \eqref{x3p2x2ypxmy} are
$$
(x,y)= (0,0),\pm (1,-2).
$$
 
 \vspace{10pt}
   
   The next equation we will consider is 
\begin{equation}\label{x3p2x2yp2y}
x^3+2x^2y+2y=0.
\end{equation}
Because $2x^2+2 \neq 0$, we can rearrange the equation to
$y=-\frac{x^3}{2x^2+2}=-\frac{x}{2}+\frac{x}{2x^2+2}$. After multiplying by $2$ we obtain 
$2y=-x+\frac{x}{2x^2+2}$. 
To ensure integer solutions, we must have that
 $\frac{x}{2x^2+2}$ is integer. Hence, we must have either $|x|\geq 2x^2+2$ or $x=0$. The inequality is impossible, and the equality and \eqref{x3p2x2yp2y} implies that $y=0$. 
 Therefore, we can conclude that the unique integer solution to equation \eqref{x3p2x2yp2y} is
$$
(x,y)=(0,0).
$$

\vspace{10pt}

The next equation we will consider is 
\begin{equation}\label{x4pxyp3ym2}
x^4+xy+3y-2=0.
\end{equation}
This equation has no integer solutions with $x=-3$, so we may assume that $x \neq -3$. We can then rearrange the equation to 
$y=\frac{2-x^4}{x+3}=-x^3+3x^2-9x+27-\frac{79}{x+3}$. To ensure that $y$ is integer, $\frac{79}{x+3}$ must be integer.
After checking the finite number of values of $x$ such that $x+3$ is a divisor of $79$, 
 we can conclude that the integer solutions to equation \eqref{x4pxyp3ym2} are
$$
\begin{aligned}
(x,y)=(-82,572306),(-4,254),(-2,-14),(76,-422306).
\end{aligned}
$$

\vspace{10pt}

The next equation we will consider is 
\begin{equation}\label{x4pxyp3yp2}
x^4+xy+3y+2=0.
\end{equation}
This equation has no integer solutions with $x=-3$, so we may assume that $x \neq -3$. We can then rearrange the equation to 
$y=-\frac{2+x^4}{x+3}=-x^3+3x^2-9x+27-\frac{83}{x+3}$. To ensure that $y$ is integer, $\frac{83}{x+3}$ must be integer. 
After checking the finite number of values of $x$ such that $x+3$ is a divisor of $83$,
 we can conclude that the integer solutions to equation \eqref{x4pxyp3yp2} are
$$
\begin{aligned}
(x,y)=(-86,659046),(-4,258),(-2,-18),(80,-493494).
\end{aligned}
$$

\vspace{10pt}

The next equation we will consider is 
\begin{equation}\label{x4pxyp4y}
x^4+xy+4y=0.
\end{equation}
This equation has no integer solutions with $x=-4$, so we may assume that $x \neq -4$. We can then rearrange the equation to 
$y=-\frac{x^4}{x+4}=-x^3+4x^2-16x+64-\frac{256}{x+4}$. To ensure that $y$ is integer, $\frac{256}{x+4}$ must be integer.
After checking the finite number of values of $x$ such that $x+4$ is a divisor of $256$, 
 we can conclude that the integer solutions to equation \eqref{x4pxyp4y} are
$$
\begin{aligned}
(x,y)=&(-260,17850625), (-132,2371842), (-68,334084), (-36,52488),\\& (-20,10000),(-12,2592), (-8,1024), (-6,648), (-5,625), \\&(-3,-81), (-2,-8), (0,0), (4,-32), (12,-1296), (28,-19208),\\ &(60,-202500), (124,-1847042), (252,-15752961).
\end{aligned}
$$

\vspace{10pt}

The next equation we will consider is
\begin{equation}\label{x4pxym3ypx}
x^4+xy-3y+x=0.
\end{equation}
This equation has no integer solutions with $x=3$, so we may assume that $x \neq 3$. We can then rearrange the equation to  
$y=-\frac{x^4+x}{x-3}=-x^3-3x^2-9x-28-\frac{84}{x-3}$. To ensure integer solutions, we must have that 
$\frac{84}{x-3}$ is integer.
After checking the finite number of values of $x$ such that $x-3$ is a divisor of $84$, 
 we can conclude that the integer solutions to equation \eqref{x4pxym3ypx} are
$$
\begin{aligned}
(x,y)=&(-81,512460),(-39,55081),(-25,13950), ( -18,4998), ( -11,1045), \\
&( -9,546), (-4,36),(-3,13), (-1,0),(0,0),( 1,1),(2, 18), (4, -260),\\
& ( 5, -315), (6, -434), (7, -602), (9,-1095), (10, -1430), (15,-4220),\\
& (17, -5967), (24,-15800), (31,-32984), (45, -97635), (87, -682022).
\end{aligned}
$$

\vspace{10pt}

The next equation we will consider is 
\begin{equation}\label{x4pxyp3ypx}
x^4+xy+3y+x=0.
\end{equation}
This equation has no integer solutions with $x=-3$, so we may assume that $x \neq -3$. We can then rearrange the equation to 
$y=-\frac{x^4+x}{x+3}=-x^3+3x^2-9x+26-\frac{78}{x+3}$. To ensure integer solutions, we must have that 
$\frac{78}{x+3}$ is integer. 
 After checking the finite number of values of $x$ such that $x+3$ is a divisor of $78$, we can conclude that the integer solutions to equation \eqref{x4pxyp3ypx} are
$$
\begin{aligned}
(x,y)=&(-81,551880), (-42,79786),(-29, 27202),(-16, 5040), (-9,1092),\\
& (-6,430),(-5, 310), (-4,252), (-2 ,-14), (-1, 0), ( 0, 0), ( 3,-14),\\
&(10, -770),(23, -10764), (36,-43068), (75, -405650).
\end{aligned}
$$

\vspace{10pt}

The next equation we will consider is 
\begin{equation}\label{x4p2xypym2}
x^4+2xy+y-2=0.
\end{equation}
Because $x$ is integer, $2x+1 \neq 0$, so, we can rearrange the equation to 
$y=\frac{2-x^4}{2x+1}=-\frac{x^3}{2}+\frac{x^2}{4}-\frac{x}{8}+\frac{1}{16}+\frac{31/16}{2x+1}$. After 
multiplying by $16$, we obtain
$$
16y=-8x^3+4x^2-2x+1+\frac{31}{2x+1}.
$$
To ensure integer solutions, we must have that 
$\frac{31}{2x+1}$ is integer.
 After checking the finite number of values of $x$ such that $2x+1$ is a divisor of 31, we can conclude that the integer solutions to equation \eqref{x4p2xypym2} are
$$
(x,y)=(-16,2114),(-1,-1),(0,2),(15,-1633).
$$

\vspace{10pt}

The next equation we will consider is 
\begin{equation}\label{x4p2xypyp2}
x^4+2xy+y+2=0.
\end{equation}
Because $x$ is integer, $2x+1 \neq 0$, so, we can rearrange the equation to $y=-\frac{2+x^4}{2x+1}=-\frac{x^3}{2}+\frac{x^2}{4}-\frac{x}{8}+\frac{1}{16}-\frac{33/16}{2x+1}$. After multiplying by $16$, we obtain 
$$
16y=-8x^3+4x^2-2x+1-\frac{33}{2x+1}.
$$
To ensure integer solutions, we must have that
$\frac{33}{2x+1}$ is integer. 
 After checking the finite number of values of $x$ such that $2x+1$ is a divisor of 33, we can conclude that the integer solutions to equation \eqref{x4p2xypyp2} are
$$
(x,y)=( -17,2531), (-6,118), (-2 ,6), (-1,3), (0 ,-2), (1, -1),
 (5,-57) , (16, -1986).
$$

\vspace{10pt}

The next equation we will consider is 
\begin{equation}\label{x4p2xypxmy}
x^4+2xy+x-y=0.
\end{equation}
 Because $x$ is integer, $2x-1 \neq 0$, so, we can rearrange the equation to $y=-\frac{x^4+x}{2x-1}=-\frac{x^3}{2}-\frac{x^2}{4}-\frac{x}{8}-\frac{9}{16}-\frac{9/16}{2x-1}$. After multiplying by $16$, we obtain 
$$
16y=-8x^3-4x^2-2x-9-\frac{9}{2x-1}.
$$
To ensure integer solutions, we must have that 
$\frac{9}{2x-1}$ is integer. 
 After checking the finite number of values of $x$ such that $2x-1$ is a divisor of 9, 
 we can conclude that the integer solutions to equation \eqref{x4p2xypxmy} are
$$
(x,y)=(-4,28),(-1,0),(0,0),(1,-2),(2,-6),(5,-70).
$$

\vspace{10pt}

The next equation we will consider is 
\begin{equation}\label{x4p2xypxpy}
x^4+2xy+x+y=0.
\end{equation}
 Because $x$ is integer, $2x+1 \neq 0$, so, we can rearrange the equation to 
 $y=-\frac{x^4+x}{2x+1}=-\frac{x^3}{2}+\frac{x^2}{4}-\frac{x}{8}-\frac{7}{16}+\frac{7/16}{2x+1}$. After multiplying by $16$, we obtain 
$$
16y=-8x^3+4x^2-2x-7+\frac{7}{2x+1}.
$$
To ensure integer solutions, we must have that 
$\frac{7}{2x+1}$ is integer. 
  After checking the finite number of values of $x$ such that $2x+1$ is a divisor of 7, we can conclude that the integer solutions to equation \eqref{x4p2xypxpy} are
$$
(x,y)=(-4,36),(-1,0),(0,0),(3,-12).
$$

\vspace{10pt}

 The final equation we will consider is 
\begin{equation}\label{x4px2ypxpy}
x^4+x^2y+x+y=0.
\end{equation}
Because $x^2+1 \neq 0$, we can rearrange the equation to 
$y=-\frac{x^4+x}{x^2+1}=1-x^2-\frac{1+x}{x^2+1}$. To ensure that $y$ is integer,  
$\frac{1+x}{x^2+1}$ must be integer. Hence, either $|1+x| \geq x^2+1$ or $1+x=0$. The inequality is true for integers $x=0$ and $1$, and the equality implies that $x=-1$. After substituting these values of $x$ into \eqref{x4px2ypxpy} and solving for $y$, we can conclude that the integer solutions to equation \eqref{x4px2ypxpy} are
$$
(x,y)= (-1,0),(0,0),(1,-1).
$$

\begin{center}
\begin{tabular}{|c|c|c|}
\hline
Equation & Solution $(x,y)$ \\\hline \hline
$x^3 + 2x^2y - y = 0$ & $(0, 0), \pm (1, -1) $\\\hline
$x^3 + 2x^2y + y = 0$ & $(0,0)$\\\hline
$x^4 + xy + 3y = 0$&$(-84,614656),(-30,30000),(-12,2304),(-6,432),(-4,256),$\\
&$(-2,-16),(0,0),(6,-144),(24,-12288),(78,-456976)$ \\\hline
$x^4 + 2xy + y = 0$ & $ (-1, 1), (0, 0)$\\\hline
$x^3+2x^2y-y+1=0$&$(-2,1),(-1,0),(0,1),(1,-2)$ \\\hline
$x^3+2x^2y+y+1=0$&$(-1,0),(0,-1),(2,-1)$ \\
\hline
\end{tabular}
\captionof{table}{Integer solutions to the equations listed in Table \ref{tab:H28liny}.}\label{Table:ex3.45}
\end{center}

\subsection{Exercise 3.48}\label{ex:H27noykterms}
\textbf{\emph{Solve the equations listed in Table \ref{tab:H27noykterms}.  }}

\begin{center}
\begin{tabular}{ |c|c|c|c|c|c| } 
 \hline
 $H$ & Equation & $H$ & Equation & $H$ & Equation \\ 
 \hline\hline
 $26$ & $x^3+y-2 x y^2 = 0$ & $26$ & $x^4+y+x y^2 = 0$ & $27$ & $1+x+x y+x^3 y+y^2 = 0$ \\ 
 \hline
 $26$ & $x^3-y+x^2 y-x y^2 = 0$ & $27$ & $1+x^3+y-2 x y^2 = 0$ & $27$ & $-1+x+x y+x^3 y+y^2 = 0$ \\ 
 \hline
 $26$ & $x^3+y+x^2 y-x y^2 = 0$ & $27$ & $1+x^3-y+x^2 y-x y^2 = 0$ & $27$ & $1+x+x^3 y+y^3 = 0$ \\ 
 \hline
 $26$ & $2 x^3+y-x y^2 = 0$ & $27$ & $1+x^3+y+x^2 y-x y^2 = 0$ & $27$ & $-1+x+x^3 y+y^3 = 0$ \\ 
 \hline
 $26$ & $x^3+y+x^2 y^2 = 0$ & $27$ & $1+2 x^3+y-x y^2 = 0$ & $27$ & $1+x^3+y-2 x y+x^2 y^2 = 0$ \\ 
 \hline
 $26$ & $x-x y+x^3 y+y^2 = 0$ & $27$ & $1+x^3+y+x^2 y^2 = 0$ & $27$ & $1+x^4+y+x y^2 = 0$ \\ 
 \hline
 $26$ & $x+x y+x^3 y+y^2 = 0$ & $27$ & $-1+x^3+y+x^2 y^2 = 0$ & $27$ & $1+x^3+y+2 x y+x^2 y^2 = 0$ \\ 
 \hline
 $26$ & $x+x^3 y+y^3 = 0$ & $27$ & $1+x-x y+x^3 y+y^2 = 0$ & $27$ & $-1+x^4+y+x y^2 = 0$ \\ 
 \hline
 $26$ & $x^3+y+x^2 y^2 = 0$ & $27$ & $-1+x-x y+x^3 y+y^2 = 0$ &  & \\ 
 \hline
\end{tabular}
\captionof{table}{\label{tab:H27noykterms} Equations of size $H\leq 27$ with no $cx^k$ or no $cy^k$ terms for $k\geq 2$.}
\end{center} 

We will next look at equations either with no $cx^k$ or no $cy^k$ terms for $k\geq 2$. We will use the method described in Section 3.3.8 of the book, which we summarise below for convenience. Let us assume that we have an equation with no power terms in $y$, then we have an equation of the form 
$$
xQ(x, y) + ay + b = 0,
$$
where $Q$ is a polynomial in $x$ and $y$ with integer coefficients, and assume that $a$ and $b$ are also integers. We can check $x=0$ separately. Otherwise, $\frac{ay+b}{x}=t$ is an integer. After substituting $y=\frac{xt-b}{a}$ into the equation, cancelling $x$ and multiplying by a constant if necessary, we obtain a polynomial $P(x,t)=0$ with integer coefficients. We can then use the Mathematica command 
\begin{equation}\label{command:noykterms}
{\tt IntegerPart[ MaxValue[Min[Abs[x], Abs[t]], \{P(x,t) == 0\}, \{t, x\}]]}
\end{equation}
which outputs either a value or $\infty$. If the command returns $\infty$, we must repeat the above process until a finite value is outputted. If the command outputs a value, say $k$, we need to check the cases $|x| \leq k$ and $|t| \leq k$ in the equation $P(x,t)=0$ to find all values for $x$, and then using the original equation, we can find all integer solutions $(x,y)$.

Equation
$$
x^3+y-2xy^2=0
$$
is solved in Section 3.3.8 of the book, and its only integer solutions are
$$
(x,y)=(0,0),\pm(1,1).
$$

Equations
$$
		x^3+y+x^2y^2=0,
$$
and
		$$
		x^4+y+xy^2=0
		$$
are solved in Section 3.3.8 of the book, and they both only have the unique integer solution
$$
(x,y)=(0,0).
$$

Equation
$$
1+x^4+y+xy^2=0
$$
is solved in Section 3.3.8 of the book, and its integer solutions are
$$
(x,y)=(-1,-1),(-1,2),(0,-1).
$$

Let us look at the first equation
\begin{equation}\label{x3mypx2ymxy2}
x^3-y+x^2y-x y^2=0.
\end{equation}
If $x=0$ then $y=0$. Now, assume that $x \neq 0$. Dividing \eqref{x3mypx2ymxy2} by $x$, we obtain 
\begin{equation}\label{x3mypx2ymxy2red}
	\frac{y}{x}=x^2+xy-y^2.
\end{equation}
 As the right-hand side of the equation is integer, the left-hand side must also be integer. Hence, we can express $y$ as $y=xt$ for some integer $t$. After substituting this into \eqref{x3mypx2ymxy2red}, we obtain
\begin{equation}\label{x3mypx2ymxy2red2}
t=x^2+x^2t-x^2t^2.
\end{equation}
The Mathematica command \eqref{command:noykterms}
$$
{\tt IntegerPart[ MaxValue[Min[Abs[x], Abs[t]], \{x^2 + x^2 t - x^2 t^2 == t\}, \{t, x\}]]}
$$
 outputs $1$. Hence, we must have either $|x| \leq 1$ or $|t| \leq 1$. Substituting these values into \eqref{x3mypx2ymxy2red2} returns integer solutions $(x,t)=(\pm 1, \pm1 ),(0,0)$. Finally, we may conclude that all integer solutions to equation \eqref{x3mypx2ymxy2} are 
$$
(x,y)=(\pm 1, \pm 1),(0,0).
$$

\vspace{10pt}

The next equation we will consider is
\begin{equation}\label{x3pypx2ymxy2}
x^3+y+x^2y-xy^2=0.
\end{equation}
If $x=0$ then $y=0$. Now, assuming $x \neq 0$, $x$ must be a divisor of $y$, 
and we may express $y$ as $y=xt$ for some integer $t$. Substituting this into \eqref{x3pypx2ymxy2} and cancelling $x$, we obtain
$$
x^2+t+x^2t-x^2t^2=0.
$$ 
The Mathematica command \eqref{command:noykterms} for this equation 
outputs $1$. Hence, we must have either $|x|\leq 1$ or $|t| \leq 1$. After checking these values, 
we can conclude that the only integer solution to equation \eqref{x3pypx2ymxy2} is
$$
(x,y)=(0,0).
$$

\vspace{10pt}

The next equation we will consider is
\begin{equation}\label{2x3pymxy2}
2x^3+y-xy^2=0.
\end{equation}
If $x=0$ then $y=0$. Otherwise, $x$ must be a divisor of $y$,
and we can express $y$ as $y=xt$ for some integer $t$. Substituting this into \eqref{2x3pymxy2} and cancelling $x$, we obtain
$$
2x^2+t-x^2t^2=0.
$$
The Mathematica command \eqref{command:noykterms} for this equation
outputs $1$. Hence, we must have either $|x|\leq 1$ or $|t| \leq 1$. After checking these values, 
we can conclude that the integer solutions to equation \eqref{2x3pymxy2} are 
$$
(x,y)=(0,0),\pm (1,-1),\pm (1,2).
$$

\vspace{10pt}

The next equation we will consider is
\begin{equation}\label{x3pypx2y2}
x^3+y+x^2y^2=0.
\end{equation}
If $x=0$ then $y=0$. Otherwise, $x$ must be a divisor of $y$, 
and we can express $y$ as $y=xt$ for some integer $t$. Substituting this into \eqref{x3pypx2y2} and cancelling $x$, we obtain 
$$
x^2+t+x^3t^2=0.
$$
The Mathematica command \eqref{command:noykterms} for this equation outputs $1$.  Hence, we must have either $|x|\leq 1$ or $|t| \leq 1$. After checking these values, 
we can conclude that the unique integer solution to equation \eqref{x3pypx2y2} is 
$$
(x,y)=(0,0).
$$

\vspace{10pt}

The next equation we will consider is
\begin{equation}\label{xmxypx3ypy2}
x-xy+x^3y+y^2=0.
\end{equation}
If $y=0$ then $x=0$. Otherwise, $y$ is a divisor of $x$, and 
we can express $x$ as $x=yt$ for some integer $t$. Substituting this into \eqref{xmxypx3ypy2} and cancelling $y$, we obtain 
$$
t+y-yt+y^3t^3=0.
$$
The Mathematica command \eqref{command:noykterms} for this equation outputs $1$.  Hence, we must have either $|y|\leq 1$ or $|t| \leq 1$. After checking these values, 
 we can conclude that the integer solutions to equation \eqref{xmxypx3ypy2} are 
$$
(x,y)=(-1,\pm 1),(0,0).
$$

\vspace{10pt}

The next equation we will consider is
\begin{equation}\label{xpxypx3ypy2}
x+xy+x^3y+y^2=0.
\end{equation}
If $y=0$ then $x=0$. Otherwise, $y$ is a divisor of $x$, and we can express $x$ as $x=yt$ for some integer $t$. Substituting this into \eqref{xpxypx3ypy2} and cancelling $y$, we obtain
$$
t+y+yt+y^3t^3=0.
$$
The Mathematica command \eqref{command:noykterms} for this equation outputs $1$.  Hence, we must have either $|y|\leq 1$ or $|t| \leq 1$. After checking these values, 
 we can conclude that the integer solutions to equation \eqref{xpxypx3ypy2} are 
$$
(x,y)=(0,0),(1,-1).
$$

\vspace{10pt}

The next equation we will consider is
\begin{equation}\label{xpx3ypy3}
x+x^3y+y^3=0.
\end{equation}
If $y=0$ then $x=0$. Otherwise, $y$ is a divisor of $x$, and we can express $x$ as $x=yt$ for some integer $t$. Substituting this into \eqref{xpx3ypy3} and cancelling $y$, we obtain 
$$
t+y^3t^3+y^2=0.
$$
The Mathematica command \eqref{command:noykterms} for this equation outputs $1$.  Hence, we must have either $|y|\leq 1$ or $|t| \leq 1$. After checking these values, 
we can conclude that the unique integer solution to equation \eqref{xpx3ypy3} is
$$
(x,y)=(0,0).
$$

\vspace{10pt}

The next equation we will consider is
\begin{equation}\label{1px3pym2xy2}
1+x^3+y-2xy^2=0.
\end{equation}
If $x=0$ then $y=-1$. Now, assume that $x\neq 0$. Dividing \eqref{1px3pym2xy2} by $x$, we obtain
\begin{equation}\label{1px3pym2xy2red}
 \frac{y+1}{x}=2y^2-x^2,
\end{equation}
  and because $\frac{y+1}{x}$ must be integer, we can express $y$ as $y=xt-1$ for some integer $t$. After substituting this into \eqref{1px3pym2xy2red}, we obtain
$$
t=2(xt-1)^2-x^2.
$$
The Mathematica command \eqref{command:noykterms} for this equation outputs $1$.  Hence, we must have either $|x|\leq 1$ or $|t| \leq 1$. After checking these values, 
we can conclude that the integer solutions to equation \eqref{1px3pym2xy2} are 
$$
(x,y)=(-3,2),(-1,0),(0,-1).
$$

\vspace{10pt}

The next equation we will consider is
\begin{equation}\label{1px3mypx2ymxy2}
	1+x^3-y+x^2y-xy^2=0.
\end{equation}
If $x=0$ then $y=1$. Otherwise, $x$ must be a divisor of $1-y$, and we can express $y$ as $y=1-xt$ for some integer $t$. After substituting this into \eqref{1px3mypx2ymxy2} and cancelling $x$, we obtain
$$
t+x^2+x(1-xt)-(1-xt)^2=0.
$$
The Mathematica command \eqref{command:noykterms} for this equation outputs $1$.  Hence, we must have either $|x|\leq 1$ or $|t| \leq 1$. After checking these values, we can conclude that the integer solutions to equation \eqref{1px3mypx2ymxy2} are 
$$
(x,y)=(-1,0),(0,1),(2,3),(3,-2).
$$

\vspace{10pt}

The next equation we will consider is
\begin{equation}\label{1px3pypx2ymxy2}
	1+x^3+y+x^2y-xy^2=0.
\end{equation}
If $x=0$ then $y=-1$. Otherwise, $x$ must be a divisor of $1+y$, and we can express $y$ as $y=xt-1$ for some integer $t$. After substituting this into \eqref{1px3pypx2ymxy2} and cancelling $x$, we obtain
$$
t+x^2+x(xt-1)-(xt-1)^2=0.
$$
The Mathematica command \eqref{command:noykterms} for this equation outputs $2$.  Hence, we must have either $|x|\leq 2$ or $|t| \leq 2$. After checking these values, we can conclude that the integer solutions to equation  \eqref{1px3pypx2ymxy2} are 
$$
(x,y)=(-2,1),(-1,-2),(-1,0),(0,-1).
$$

\vspace{10pt}

The next equation we will consider is
\begin{equation}\label{1p2x3pymxy2}
	1+2x^3+y-xy^2=0.
\end{equation}
If $x=0$ then $y=-1$, otherwise, $x$ is a divisor of $y+1$ and 
we can express $y$ as $y=xt-1$ for some integer $t$. After substituting this into \eqref{1p2x3pymxy2} and cancelling $x$, we obtain
$$
2x^2+t-(xt-1)^2=0.
$$
The Mathematica command \eqref{command:noykterms} for this equation outputs $2$.  Hence, we must have either $|x|\leq 2$ or $|t| \leq 2$. After checking these values, we can conclude that the integer solutions to equation  \eqref{1p2x3pymxy2} are 
$$
(x,y)=(-2,-3),(0,-1).
$$

\vspace{10pt}

The next equation we will consider is
\begin{equation}\label{1px3pypx2y2}
	1+x^3+y+x^2y^2=0.
\end{equation}
If $x=0$ then $y=-1$. Otherwise, $x$ is a divisor of $y+1$ and we can
express $y$ as $y=xt-1$ for some integer $t$. After substituting this into \eqref{1px3pypx2y2} and cancelling $x$, we obtain
$$
t+x(xt-1)^2+x^2=0.
$$
The Mathematica command \eqref{command:noykterms} for this equation outputs $1$.  Hence, we must have either $|x|\leq 1$ or $|t| \leq 1$. After checking these values, we can conclude that the integer solutions to equation \eqref{1px3pypx2y2} are 
$$
(x,y)=(-1,-1),(-1,0),(0,-1).
$$

\vspace{10pt}

The next equation we will consider is
\begin{equation}\label{m1px3pypx2y2}
	-1+x^3+y+x^2y^2=0.
\end{equation}
If $x=0$ then $y=1$. Otherwise, $x$ is a divisor of $y-1$ and we can express $y$ as $y=xt+1$ for some integer $t$. After substituting this into \eqref{m1px3pypx2y2} and cancelling $x$, we obtain
$$
t+x(xt+1)^2+x^2=0.
$$
The Mathematica command \eqref{command:noykterms} for this equation outputs $1$.  Hence, we must have either $|x|\leq 1$ or $|t| \leq 1$. After checking these values, we can conclude that the integer solutions to equation \eqref{m1px3pypx2y2} are 
$$
(x,y)=(-1,-2),(0,1),\pm (1,-1),(1,0).
$$

\vspace{10pt}

The next equation we will consider is
\begin{equation}\label{1pxmxypx3ypy2}
	1+x-xy+x^3y+y^2=0.
\end{equation}
If $y=0$ then $x=-1$. Otherwise, $y$ is a divisor of $x+1$, and
we can express $x$ as $x=yt-1$ for some integer $t$. After substituting this into \eqref{1pxmxypx3ypy2} and cancelling $y$, we obtain
$$
t+y+2t y-3t^2 y^2+t^3 y^3=0.
$$
The Mathematica command \eqref{command:noykterms} for this equation outputs $1$.  Hence, we must have either $|y|\leq 1$ or $|t| \leq 1$. After checking these values, we can conclude that the unique integer solution to equation \eqref{1pxmxypx3ypy2} is
$$
(x,y)=(-1,0).
$$

\vspace{10pt}

The next equation we will consider is
\begin{equation}\label{m1pxmxypx3ypy2}
	-1+x-xy+x^3y+y^2=0.
\end{equation}
If $y=0$ then $x=1$. Otherwise, $y$ is a divisor of $x-1$, and 
we can express $x$ as $x=yt+1$ for some integer $t$. After substituting this into \eqref{m1pxmxypx3ypy2} and cancelling $y$, we obtain
$$
t+y+2t y+3t^2 y^2+t^3 y^3=0.
$$
The Mathematica command \eqref{command:noykterms} for this equation outputs $1$.  Hence, we must have either $|y|\leq 1$ or $|t| \leq 1$. After checking these values,  we can conclude that the integer solutions to equation \eqref{m1pxmxypx3ypy2} are 
$$
(x,y)=(0,\pm 1),(1,0).
$$

\vspace{10pt}

The next equation we will consider is
\begin{equation}\label{1pxpxypx3ypy2}
	1+x+xy+x^3y+y^2=0.
\end{equation}
If $y=0$ then $x=-1$. Otherwise, $y$ is a divisor of $x+1$, and we can express $x$ as $x=yt-1$ for some integer $t$. After substituting this into \eqref{1pxpxypx3ypy2} and cancelling $y$, we obtain
$$
-2 + t + y + 4 t y - 3 t^2 y^2 + t^3 y^3=0
$$
The Mathematica command \eqref{command:noykterms} for this equation outputs $1$.  Hence, we must have either $|y|\leq 1$ or $|t| \leq 1$. After checking these values,  we can conclude that the integer solutions to equation \eqref{m1pxpxypx3ypy2} are 
$$
(x,y)=(-1,0),(-1,2).
$$

\vspace{10pt}

The next equation we will consider is
\begin{equation}\label{m1pxpxypx3ypy2}
-1+x+xy+x^3y+y^2=0.
\end{equation}
If $y=0$ then $x=1$. Otherwise, $y$ is a divisor of $x-1$, and we can express $x$ as $x=yt+1$ for some integer $t$. After substituting this into \eqref{m1pxpxypx3ypy2} and cancelling $y$, we obtain
$$
2 + t + y + 4 t y + 3 t^2 y^2 + t^3 y^3=0
$$
The Mathematica command \eqref{command:noykterms} for this equation outputs $1$.  Hence, we must have either $|y|\leq 1$ or $|t| \leq 1$. After checking these values,  we can conclude that the integer solutions to equation \eqref{m1pxpxypx3ypy2} are 
$$
(x,y)=(0,\pm 1),(1,-2),(1,0).
$$

\vspace{10pt}

The next equation we will consider is
\begin{equation}\label{1pxpx3ypy3}
1+x+x^3y+y^3=0.
\end{equation}
If $y=0$ then $x=-1$. Otherwise, $y$ is a divisor of $x+1$, and we can express $x$ as $x=yt-1$ for some integer $t$. After substituting this into \eqref{1pxpx3ypy3} and cancelling $y$, we obtain
$$
-1 + t + 3 t y + y^2 - 3 t^2 y^2 + t^3 y^3=0
$$
The Mathematica command \eqref{command:noykterms} for this equation outputs $1$.  Hence, we must have either $|y|\leq 1$ or $|t| \leq 1$. After checking these values,  we can conclude that the integer solutions to equation \eqref{1pxpx3ypy3} are 
$$
(x,y)=(-1,\pm 1),(-1,0),(0,-1),(1,-1).
$$

\vspace{10pt}

The next equation we will consider is
\begin{equation}\label{m1pxpx3ypy3}
-1+x+x^3y+y^3=0.
\end{equation}
If $y=0$ then $x=1$. Otherwise, $y$ is a divisor of $x-1$, and we can express $x$ as $x=yt+1$ for some integer $t$. After substituting this into \eqref{m1pxpx3ypy3} and cancelling $y$, we obtain
$$
1 + t + 3 t y + y^2 + 3 t^2 y^2 + t^3 y^3=0
$$
The Mathematica command \eqref{command:noykterms} for this equation outputs $1$.  Hence, we must have either $|y|\leq 1$ or $|t| \leq 1$. After checking these values,  we can conclude that the integer solutions to equation \eqref{m1pxpx3ypy3} are 
$$
(x,y)=(-2,3),(0,1),(1,0).
$$

\vspace{10pt}

The next equation we will consider is
\begin{equation}\label{1px3pym2xypx2y2}
1+x^3+y-2xy+x^2y^2=0.
\end{equation}
If $x=0$ then $y=-1$. Otherwise, $x$ is a divisor of $y+1$, and we can express $y$ as $y=xt-1$ for some integer $t$. After substituting this into \eqref{1px3pym2xypx2y2} and cancelling $x$, we obtain
$$
2 + t + x - 2 t x + x^2 - 2 t x^2 + t^2 x^3=0
$$
The Mathematica command \eqref{command:noykterms} for this equation outputs $1$.  Hence, we must have either $|x|\leq 1$ or $|t| \leq 1$. After checking these values,  we can conclude that the integer solutions to equation \eqref{1px3pym2xypx2y2} are 
$$
(x,y)=(-1,-3),(-1,0),(0,-1).
$$

\vspace{10pt}

The next equation we will consider is
\begin{equation}\label{1px3pyp2xypx2y2}
1+x^3+y+2xy+x^2y^2=0.
\end{equation}
If $x=0$ then $y=-1$. Otherwise, $x$ is a divisor of $y+1$, and we can express $y$ as $y=xt-1$ for some integer $t$. After substituting this into \eqref{1px3pyp2xypx2y2} and cancelling $x$, we obtain
$$
-2 + t + x + 2 t x + x^2 - 2 t x^2 + t^2 x^3=0
$$
The Mathematica command \eqref{command:noykterms} for this equation outputs $1$.  Hence, we must have either $|x|\leq 1$ or $|t| \leq 1$. After checking these values, we can conclude that the integer solutions to equation \eqref{1px3pyp2xypx2y2} are 
$$
(x,y)=(-3,2),(-2,-1),(-1,0),(0,-1),(1,-2),\pm(1,-1).
$$

\vspace{10pt}

The final equation we will consider is
\begin{equation}\label{x4m1pypxy2}
	-1+x^4+y+xy^2=0.
\end{equation}
If $x=0$ then $y=1$. Otherwise, $x$ is a divisor of $y-1$, and we can express $y$ as $y=xt+1$ for some integer $t$. After substituting this into \eqref{x4m1pypxy2} and cancelling $x$, we obtain
\begin{equation}\label{x4m1pypxy2red}
1 + t + 2 t x + t^2 x^2 + x^3=0
\end{equation}
The Mathematica command \eqref{command:noykterms} for this equation outputs $\infty$.  So, we must repeat the above process. As $x \neq 0$, \eqref{x4m1pypxy2red} implies that $x$ is a divisor of $t+1$, hence we may express $t$ as $t=xk-1$ for some integer $k$. After substituting this into \eqref{x4m1pypxy2red} and cancelling $x$, we obtain
$$
-2 + k + x + 2 k x + x^2 - 2 k x^2 + k^2 x^3=0.
$$
The Mathematica command \eqref{command:noykterms} for this equation outputs $1$. Hence, we must have either $|x|\leq 1$ or $|k| \leq 1$. After checking these values, we can conclude that the integer solutions to equation \eqref{x4m1pypxy2} are 
$$
(x,y)=(-3,-5),(-2,3),(\pm 1,0),(0,1),\pm(1,-1).
$$

\begin{center}
\begin{tabular}{ |c|c|c| } 
 \hline
 Equation & Solution $(x,y)$ \\ 
 \hline\hline
$x^3+y-2 x y^2 = 0$ &$  (0, 0), \pm (1, 1)$ \\ \hline
$x^3-y+x^2 y-x y^2 = 0$ &$(\pm 1,\pm 1),(0,0)$\\ \hline
$x^3+y+x^2 y-x y^2 = 0$ &$(0,0)$ \\ \hline
$2 x^3+y-x y^2 = 0$ &$(-1,-2),(0,0),\pm (1,-1),(1,2)$ \\ \hline
 $x^3+y+x^2 y^2 = 0$ &$(0,0)$ \\ \hline
$x-x y+x^3 y+y^2 = 0$ &$(-1,\pm 1),(0,0)$\\ \hline
$x+x y+x^3 y+y^2 = 0$ &$(0,0),(1,-1)$ \\ \hline
$x+x^3 y+y^3 = 0$ &$(0,0)$ \\ \hline
 $x^3+y+x^2 y^2 = 0$ &$(0,0)$ \\ \hline
 $x^4+y+x y^2 = 0$ &$(0,0)$ \\ \hline 
$1+x^3+y-2 x y^2 = 0$ &$(-3,2),(-1,0),(0,-1)$ \\ \hline
$1+x^3-y+x^2 y-x y^2 = 0$ &$(-1,0),(0,1),(2,3),(3,-2)$ \\ \hline
$1+x^3+y+x^2 y-x y^2 = 0$ &$(-2,1),(-1,-2),(-1,0),(0,-1)$ \\ \hline
$1+2 x^3+y-x y^2 = 0$ &$(-2,-3),(0,-1)$ \\ \hline
 $1+x^3+y+x^2 y^2 = 0$ &$(-1,-1),(-1,0),(0,-1)$\\ \hline
$-1+x^3+y+x^2 y^2 = 0$ &$(-1,-2),(0,1),\pm (1,-1),(1,0)$ \\ \hline
$1+x-x y+x^3 y+y^2 = 0$ &$(-1,0)$ \\ \hline
$-1+x-x y+x^3 y+y^2 = 0$ &$(0,\pm1),(1,0)$ \\  \hline
$1+x+x y+x^3 y+y^2 = 0$&$(-1,0),(-1,2)$ \\  \hline
$-1+x+x y+x^3 y+y^2 = 0$&$(0,\pm1),(1,-2),(1,0)$ \\ \hline
$1+x+x^3 y+y^3 = 0$&$(-1,\pm 1),(-1,0),(0,-1),(1,-1)$ \\ \hline
$-1+x+x^3 y+y^3 = 0$&$(-2,3),(0,1),(1,0)$ \\  \hline
$1+x^3+y-2 x y+x^2 y^2 = 0$&$(-1,-3),(-1,0),(0,-1)$ \\  \hline
$1+x^4+y+x y^2 = 0$&$(-1, -1), (-1, 2), (0, -1)$ \\  \hline
$1+x^3+y+2 x y+x^2 y^2 = 0$&$(-3,2),(-2,-1),(-1,0),(0,-1),(1,-2),\pm (1,-1)$ \\  \hline
$-1+x^4+y+x y^2 = 0$&$(-3,-5),(-2,3),(\pm 1,0),(0,1),\pm (1,-1)$ \\ \hline

\end{tabular}
\captionof{table}{Integer solutions to the equations listed in Table \ref{tab:H27noykterms}. \label{tab:H27noyktermssol}}
\end{center}

\subsection{Exercise 3.49}\label{ex:H28axbyPxy}
\textbf{\emph{Solve the equations listed in Table \ref{tab:H28axbyPxy}.  }}

\begin{center}
\begin{tabular}{ |c|c|c|c|c|c| } 
 \hline
 $H$ & Equation & $H$ & Equation & $H$ & Equation \\ 
 \hline\hline
 $26$ & $-x+x^3+x^2 y-y^3 = 0$ & $28$ & $-x+x^3+y+x^2 y-y^3 = 0$ & $28$ & $x^3+2 y+x^2 y+y^3 = 0$ \\ 
 \hline
 $26$ & $x^3-y+x^2 y-y^3 = 0$ & $28$ & $x^3-2 y+x^2 y-y^3 = 0$ & $28$ & $x+x^3-y+x^2 y+y^3 = 0$ \\ 
 \hline
 $26$ & $x^3+y+x^2 y-y^3 = 0$ & $28$ & $x^3+2 y+x^2 y-y^3 = 0$ & $28$ & $x+x^3+y+x^2 y+y^3 = 0$ \\ 
 \hline
 $26$ & $x+x^3+x^2 y-y^3 = 0$ & $28$ & $x+x^3-y+x^2 y-y^3 = 0$ & $28$ & $2 x+x^3+x^2 y+y^3 = 0$ \\ 
 \hline
 $26$ & $-x+x^3+x^2 y+y^3 = 0$ & $28$ & $x+x^3+y+x^2 y-y^3 = 0$ &  $28$ & $-x+2 x^3-y+y^3 = 0$ \\ 
 \hline
 $26$ & $x^3-y+x^2 y+y^3 = 0$ & $28$ & $2 x+x^3+x^2 y-y^3 = 0$ &  $28$ & $-x+2 x^3+y+y^3 = 0$ \\ 
 \hline
 $26$ & $x^3+y+x^2 y+y^3 = 0$ & $28$ & $-2 x+x^3+x^2 y+y^3 = 0$ & $28$ & $x+2 x^3-y+y^3 = 0$ \\ 
 \hline
 $26$ & $x+x^3+x^2 y+y^3 = 0$ & $28$ & $-x+x^3-y+x^2 y+y^3 = 0$ & $28$ & $x+2 x^3+y+y^3 = 0$ \\ 
 \hline
 $28$ & $-2 x+x^3+x^2 y-y^3 = 0$ & $28$ & $-x+x^3+y+x^2 y+y^3 = 0$ &  &  \\ 
 \hline
 $28$ & $-x+x^3-y+x^2 y-y^3 = 0$ & $28$ & $x^3-2 y+x^2 y+y^3 = 0$ &  &  \\ 
 \hline
\end{tabular}
\captionof{table}{\label{tab:H28axbyPxy} Equations of size $H\leq 28$ of the form $ax + by = P(x,y)$ with $P$ homogeneous.}
\end{center} 

In this exercise, we will solve equations of the form 
$$
ax + by = P(x,y) \text{ with } P \text{ homogeneous}. 
$$
To solve these equations, we will use the method from Section 3.3.9 of the book, which we summarise below for convenience. First, we reduce the equation to the form
$$
z = P(x,(z-ax)/b) \text{ with } P \text{ homogeneous} 
$$
by making the substitution $z=ax+by$. The case $z=0$ can be checked separately, hence we may assume $z \neq 0$ and $P(x,(z-ax)/b)$ contains a term independent of $z$. By multiplying by a constant, say $A$, if necessary to obtain integer coefficients, we have an equation of the form
$$
Az = Q(x,z) \text{ with } Q \text{ homogeneous}. 
$$
Then we can make substitutions $x=dX$ and $z=dZ$, where $d=\gcd(x,z)$ and $\gcd(X,Z)=1$, and cancelling $d$, to obtain
\begin{equation}\label{AZedkm1QXZ}
AZ =d^{k-1}(Q(X,Z)),
\end{equation}
where $k$ is the degree of $Q$. Because $Q(X,Z)$ contains a term independent of $Z$, which we denote as $cd^{k-1}X^k$, $Z$ must divide this term. But $X$ and $Z$ are coprime, so $Z$ divides the coefficient $cd^{k-1}$. Hence, we may represent $cd^{k-1}$ as $cd^{k-1}=Zt$ for some integer $t$. After multiplying \eqref{AZedkm1QXZ} by $c$, substituting $cd^{k-1}=Zt$ and cancelling $Z$, we obtain 
$$
cA =t(Q(X,Z)).
$$
Then, because $cA$ is an integer and $t$ must be a divisor of this integer, we have reduced the problem to solving a finite number of systems of equations of the form
$$
t=a_1 \quad \text{and} \quad Q(X,Z)=a_2 \quad \text{where} \quad a_1 a_2=cA.
$$
Because $Q(X,Z)$ is homogeneous, $Q(X,Z)=a_2$ is a Thue equation and has a finite number of integer solutions which can be found in Mathematica, see Section \ref{ex:Thue}.

Let us look at the first equation 
\begin{equation}\label{mxpx3px2ymy3}
-x+x^3+x^2 y-y^3 = 0.
\end{equation}
If $(x,y)$ is a solution to this equation, then $(-x,-y)$ is also a solution, so, let us find solutions with $x \geq 0$ and then take both signs. If $x=0$, then $y=0$, so we have the solution $(x,y)=(0,0)$. Next, let $(x,y)$ be any solution with $x>0$, and let $d=\gcd(x,y)$. Then $x=dX$ and $y=dY$ with $X,Y$ coprime and $X >0$. Substituting this into \eqref{mxpx3px2ymy3} and cancelling $d$, we obtain
\begin{equation}\label{mxpx3px2ymy3_x1y1}
d^2(X^3+X^2Y-Y^3)-X=0.
\end{equation}
From this equation, it is clear that $X$ is a divisor of $d^2Y^3$. Because $X$ and $Y$ are coprime, this implies that $X$ is a divisor of $d^2$, so we can write $d^2=Xz$ for some integer $z>0$. Then after substitution and cancelling $X$, \eqref{mxpx3px2ymy3_x1y1} reduces to
$$
z(X^3+X^2Y-Y^3)=1.
$$
Hence, $z>0$ is a divisor of $1$, which implies that $z=1$ and $X^3+X^2Y-Y^3=1$. The last equation is a Thue equation and can be solved in Mathematica using the command
$$
{\tt Reduce[X^3 + X^2 Y - Y^3 == 1, \{X, Y\}, Integers]}
$$
which outputs integer solutions $(X,Y) =$ $(-3,-4),$ $(0,-1),$ $(1,\pm 1)$ and $(1, 0)$. Because $X>0$, the only possible pairs of $(X,Y)$ are $(1,\pm 1)$ and $(1, 0)$. These solutions give $d=1$, hence $x=dX=1$ and $y=dY=\pm 1$ or $y=dY=0$.
In conclusion, all integer solutions to equation \eqref{mxpx3px2ymy3} are 
$$
(x,y)=(\pm1,\pm 1),(\pm 1,0),(0,0).
$$

\vspace{10pt}

Let us consider the next equation
\begin{equation}\label{xpx3px2ymy3}
x+x^3+x^2 y-y^3 = 0.
\end{equation}
If $(x,y)$ is a solution to this equation, then $(-x,-y)$ is also a solution, so, let us find solutions with $x \geq 0$ and then take both signs. If $x=0$ then $y=0$. Next, let $(x,y)$ be any solution with $x>0$, and let $d=\gcd(x,y)$. Then $x=dX$ and $y=dY$ with $X,Y$ coprime and $X >0$. Substituting this into \eqref{xpx3px2ymy3} and cancelling $d$, we obtain
\begin{equation}\label{xpx3px2ymy3_x1y1}
d^2(-X^3-X^2Y+Y^3)-X=0.
\end{equation}
From this equation, it is clear that $X$ is a divisor of $d^2Y^3$. Because $X$ and $Y$ are coprime, this implies that $X$ is a divisor of $d^2$, so we can write $d^2=Xz$ for some integer $z>0$. Then after substitution and cancelling $X$, \eqref{xpx3px2ymy3_x1y1} reduces to
$$
z(-X^3-X^2Y+Y^3)=1.
$$
Hence, $z>0$ is a divisor of $1$, which implies that $z=1$ and $-X^3-X^2Y+Y^3=1$. The last equation is a Thue equation and it can be easily verified that its only integer solutions are $(X,Y)=(-1,\pm 1),$ $(-1,0),$ $(0,-1)$ and $(3,4)$. 
Because $X>0$, the only possible pair of $(X,Y)$ is $(3, 4)$. This solution gives $d=\sqrt{X z}=\sqrt{3}$ which is not integer. In conclusion, the unique integer solution to equation \eqref{xpx3px2ymy3} is 
$$
(x,y)=(0,0).
$$

\vspace{10pt}

Let us consider the next equation
\begin{equation}\label{xpx3px2ypy3}
x+x^3+x^2 y+y^3 = 0.
\end{equation}
If $(x,y)$ is a solution to this equation, then $(-x,-y)$ is also a solution, so, let us find solutions with $x \geq 0$ and then take both signs. If $x=0$ then $y=0$. Next, let $(x,y)$ be any solution with $x>0$, and let $d=\gcd(x,y)$. Then $x=dX$ and $y=dY$ with $X,Y$ coprime and $X >0$. Substituting this into \eqref{xpx3px2ypy3} and cancelling $d$, we obtain
\begin{equation}\label{xpx3px2ypy3_x1y1}
d^2(-X^3-X^2Y-Y^3)-X=0.
\end{equation}
From this equation, it is clear that $X$ is a divisor of $d^2Y^3$. Because $X$ and $Y$ are coprime, this implies that $X$ is a divisor of $d^2$, so we can write $d^2=Xz$ for some integer $z>0$. Then after substitution and cancelling $X$, \eqref{xpx3px2ypy3_x1y1} reduces to
$$
z(-X^3-X^2Y-Y^3)=1.
$$
Hence, $z>0$ is a divisor of $1$, which implies that $z=1$ and $-X^3-X^2Y-Y^3=1$. The last equation is a Thue equation and it can be easily verified that its only integer solutions are $(X,Y)=(-3,2),$ $(-1,0),$ $(0,-1)$ and $(1,-1)$. Because $X>0$, the only possible pair of $(X,Y)$ is $(1,-1)$.
This solution gives $d=\sqrt{X z}=1$, hence $x=dX=1$, and $y=dY=-1$. In conclusion, the integer solutions to \eqref{xpx3px2ypy3} are 
$$
(x,y)=(0,0),\pm (1,-1).
$$

\vspace{10pt}

Let us consider the next equation
\begin{equation}\label{x3mypx2ypy3}
x^3-y+x^2 y+y^3 = 0.
\end{equation}
If $(x,y)$ is a solution to this equation, then $(-x,-y)$ is also a solution, so, let us find solutions with $y \geq 0$ and then take both signs. If $y=0$ then $x=0$. Next, let $(x,y)$ be any solution with $y>0$, and let $d=\gcd (x,y)$. Then $x=dX$ and $y=dY$ with $X,Y$ coprime and $Y >0$. Substituting this into \eqref{x3mypx2ypy3} and cancelling $d$, we obtain
\begin{equation}\label{x3mypx2ypy3_x1y1}
d^2(X^3+X^2Y+Y^3)-Y=0.
\end{equation}
From this equation, it is clear that $Y$ is a divisor of $d^2X^3$. Because $X$ and $Y$ are coprime, this implies that $Y$ is a divisor of $d^2$, so we can write $d^2=Yz$ for some integer $z>0$. Then after substitution and cancelling $Y$, \eqref{x3mypx2ypy3_x1y1} reduces to
$$
z(X^3+X^2Y+Y^3)=1.
$$
Hence, $z>0$ is a divisor of $1$, which implies that $z=1$ and $X^3+X^2Y+Y^3=1$. The last equation is a Thue equation and it can be easily verified that its only integer solutions are $(X,Y)=(-1,1),$ $(0,1),$ $(1,0)$ and $(3,-2)$. 
Because $Y>0$, the possible pairs of $(X,Y)$ are $(-1, 1)$ and $(0,1)$. These solutions give $d=\sqrt{Y z}=1$, hence $y=dY=1$ and $x=dX=-1$ or $x=dX=0$. In conclusion, the integer solutions to \eqref{x3mypx2ypy3} are 
$$
(x,y)=(0,\pm 1),(0,0),\pm(1,-1).
$$

\vspace{10pt}

Let us consider the next equation
\begin{equation}\label{x3pypx2ypy3}
x^3+y+x^2 y+y^3 = 0.
\end{equation}
If $(x,y)$ is a solution to this equation, then $(-x,-y)$ is also a solution, so, let us find solutions with $y \geq 0$ and then take both signs. If $y=0$ then $x=0$. Next, let $(x,y)$ be any solution with $y>0$, and let $d=\gcd (x,y)$. Then $x=dX$ and $y=dY$ with $X,Y$ coprime and $Y >0$. Substituting this into \eqref{x3pypx2ypy3} and cancelling $d$, we obtain
\begin{equation}\label{x3pypx2ypy3_x1y1}
d^2(-X^3-X^2Y-Y^3)-Y=0.
\end{equation}
From this equation, it is clear that $Y$ is a divisor of $d^2X^3$. Because $X$ and $Y$ are coprime, this implies that $Y$ is a divisor of $d^2$, so we can write $d^2=Yz$ for some integer $z>0$. Then after substitution and cancelling $Y$,
\eqref{x3pypx2ypy3_x1y1} reduces to
$$
z(-X^3-X^2Y-Y^3)=1.
$$
Hence, $z>0$ is a divisor of $1$, which implies that $z=1$ and $-X^3-X^2Y-Y^3=1$. The last equation is a Thue equation and it can be easily verified that its only integer solutions are 
$(X,Y)=(-3,2),$ $(-1,0),$ $(0,-1)$ and $(1,-1)$. 
Because $Y>0$, the only possible pair of $(X,Y)$ is $(-3,2)$. This solution gives $d=\sqrt{Y z}=\sqrt{2}$, which is not integer. In conclusion, the unique integer solution to \eqref{x3pypx2ypy3} is 
$$
(x,y)=(0,0).
$$

\vspace{10pt}

Let us look consider next equation
\begin{equation}\label{x3mypx2ymy3}
x^3-y+x^2 y-y^3 = 0.
\end{equation}
If $(x,y)$ is a solution to this equation, then $(-x,-y)$ is also a solution, so, let us find solutions with $y \geq 0$ and then take both signs. If $y=0$ then $x=0$. Next, let $(x,y)$ be any solution with $y>0$, and let $d=\gcd (x,y)$. Then $x=dX$ and $y=dY$ with $X,Y$ coprime and $Y >0$. Substituting this into \eqref{x3mypx2ymy3} and cancelling $d$, we obtain
\begin{equation}\label{x3mypx2ymy3_x1y1}
d^2(X^3+X^2Y-Y^3)-Y=0.
\end{equation}
From this equation, it is clear that $Y$ is a divisor of $d^2X^3$. Because $X$ and $Y$ are coprime, this implies that $Y$ is a divisor of $d^2$, so we can write $d^2=Yz$ for some integer $z>0$. Then after substitution and cancelling $Y$, 
\eqref{x3mypx2ymy3_x1y1} reduces to
$$
z(X^3+X^2Y-Y^3)=1.
$$
Hence, $z>0$ is a divisor of $1$, which implies that $z=1$ and $X^3+X^2Y-Y^3=1$. The last equation is a Thue equation and it can be easily verified that its only integer solutions are $(X,Y)=(-3,-4),$ $(0,-1),$ $(1, \pm 1)$ and $(1,0)$.
Because $Y>0$, the only possible pair of $(X,Y)$ is $(1,1)$. This solution gives $d=\sqrt{Y z}=1$, hence $x=dX=1$, and $y=dY=1$. In conclusion, the integer solutions to \eqref{x3mypx2ymy3} are
$$
(x,y)=(0,0),\pm (1,1).
$$

\vspace{10pt}

Let us consider the next equation
\begin{equation}\label{x3pypx2ymy3}
x^3+y+x^2 y-y^3 = 0.
\end{equation}
If $(x,y)$ is a solution to this equation, then $(-x,-y)$ is also a solution, so, let us find solutions with $y \geq 0$ and then take both signs. If $y=0$ then $x=0$. Next, let $(x,y)$ be any solution with $y>0$, and let $d=\gcd (x,y)$. Then $x=dX$ and $y=dY$ with $X,Y$ coprime and $Y >0$. Substituting this into \eqref{x3pypx2ymy3} and cancelling $d$, we obtain
\begin{equation}\label{x3pypx2ymy3_x1y1}
d^2(-X^3-X^2Y+Y^3)-Y=0.
\end{equation}
From this equation, it is clear that $Y$ is a divisor of $d^2X^3$. Because $X$ and $Y$ are coprime, this implies that $Y$ is a divisor of $d^2$, so we can write $d^2=Yz$ for some integer $z>0$. Then after substitution and cancelling $Y$,
 \eqref{x3pypx2ymy3_x1y1} reduces to
$$
z(-X^3-X^2Y+Y^3)=1.
$$
Hence, $z>0$ is a divisor of $1$, which implies that $z=1$ and $-X^3-X^2Y+Y^3=1$. The last equation is a Thue equation and it can be easily verified that its only integer solutions are 
$(X,Y)=(-1,\pm 1),$ $(-1,0),$ $(0,1)$ and $(3,4)$.
Because $Y>0$, the possible pairs of $(X,Y)$ are $(-1,1), (0,1)$ and $(3,4)$. The first and second solution give $d=\sqrt{Y z}=1$, hence $y=dY=1$ and $x=dX=-1$ or $x=dX=0$. The final solution gives $d=\sqrt{Y z}=2$, hence $x=dX=6$ and $y=dY=8$. In conclusion, the integer solutions to \eqref{x3pypx2ymy3} are 
$$
(x,y)=(0,\pm1),(0,0),\pm (1,-1),\pm(6,8).
$$

\vspace{10pt}

The next equation we will consider is
\begin{equation}\label{xp2x3mypy3}
x+2x^3-y+y^3 = 0.
\end{equation}
Let $z=x-y$ be a new variable. Then $x=y+z$ and \eqref{xp2x3mypy3} reduces to 
\begin{equation}\label{xp2x3mypy3_subx}
z+2(y+z)^3+y^3= z+ 3 y^3 + 6 y^2 z + 6 y z^2 + 2 z^3 = 0.
\end{equation}
If $(z,y)$ is a solution to this equation, then $(-z,-y)$ is also a solution, so, let us find solutions with $z \geq 0$ and then take both signs. If $z=0$ then $(x,y)=(0,0)$. Now, assume that $z>0$. Then with $d=\gcd(y,z)$, $y=dY$ and $z=dZ$ we obtain 
\begin{equation}\label{xp2x3mypy3_YZ}
Z+d^2(2 (Y + Z)^3 + Y^3)=Z+d^2(3 Y^3 + 6 Y^2 Z + 6 Y Z^2 + 2 Z^3) = 0.
\end{equation}
This implies that $Z>0$ is a divisor of $3d^2Y^3$, hence it is also a divisor of $3d^2$, so we can write $3d^2=Zt$ for some integer $t>0$. Then \eqref{xp2x3mypy3_YZ} reduces to 
$$
-3=t(2(Y+Z)^3+Y^3).
$$
Hence $t>0$ is a divisor of $-3$, which implies that either (i) $t=1$ and $2(Y+Z)^3+Y^3=-3$ or (ii) $t=3$ and $2(Y+Z)^3+Y^3=-1$. In both cases, the second equation is a Thue equation, which can be easily solved in Mathematica. In case (i), its only integer solutions are $(Y,Z)=(-1,0)$ and $(5,-9)$.
Because $Z>0$, these pairs $(Y,Z)$ are not possible. In case (ii), its only integer solutions are $(Y,Z)=(-1,1)$ and $(1,-2)$. 
  Because $Z>0$, the only possible pair $(Y,Z)$ is $(-1,1)$. Hence $3d^2=Zt=3$ implies that $d=1$, then $y=dY=-1$ and $z=dZ=1$ and then $x=y+z=-1+1=0$. In conclusion, all integer solutions to \eqref{xp2x3mypy3} are 
$$
(x,y)=(0,0),(0,\pm1).
$$

\vspace{10pt}

The next equation we will consider is
\begin{equation}\label{mxp2x3pypy3}
-x+2x^3+y+y^3 = 0.
\end{equation}
Let $z=y-x$ be a new variable. Then $y=x+z$ and \eqref{mxp2x3pypy3} reduces to 
$$
z+2x^3+(x+z)^3 = 0.
$$
If $(x,z)$ is a solution to this equation, then $(-x,-z)$ is also a solution, so, let us find solutions with $z \geq 0$ and then take both signs. If $z=0$ then $(x,y)=(0,0)$. Now, assume that $z>0$. Then with $d=\gcd(x,z)$, $x=dX$ and $z=dZ$ we obtain 
\begin{equation}\label{mxp2x3pypy3redXZ}
Z+d^2(2X^3+(X+Z)^3) = 0.
\end{equation}
This implies that $Z>0$ is a divisor of $3d^2X^3$, hence it is also a divisor of $3d^2$, so we can write $3d^2=Zt$ for some integer $t>0$. Then \eqref{mxp2x3pypy3redXZ} reduces to 
$$
3=t(-2X^3-(X+Z)^3).
$$
Hence $t>0$ is a divisor of $3$, which implies that either (i) $t=1$ and $-2X^3-(X+Z)^3=3$ or (ii) $t=3$ and $-2X^3-(X+Z)^3=1$. In both cases, the second equation is a Thue equation, which can be easily solved in Mathematica. In case (i), its only integer solutions are $(X,Z)=(-4,9)$ and $(-1,0)$. 
Because $Z>0$, the only possible pair $(X,Z)$ is $(-4,9)$, this solution gives $3d^2=Zt=9$ so $d=\sqrt{3}$ which is not integer. In case (ii), its only integer solutions are $(X,Z)=(-1,2)$ and $(0,-1)$. 
Because $Z>0$, the only possible pair $(Y,Z)$ is $(-1,2)$. Hence $3d^2=Zt=6$ so $d=\sqrt{2}$ which is not integer. In conclusion, the unique integer solution to \eqref{mxp2x3pypy3} is 
$$
(x,y)=(0,0).
$$

\vspace{10pt}

The next equation we will consider is
\begin{equation}\label{mxp2x3mypy3}
-x+2x^3-y+y^3 = 0.
\end{equation}
Let $z=x+y$ be a new variable. Then $y=z-x$ and \eqref{mxp2x3mypy3} reduces to 
$$
-z+2x^3+(z-x)^3 =x^3 - z + 3 x^2 z - 3 x z^2 + z^3 = 0.
$$
If $(x,z)$ is a solution to this equation, then $(-x,-z)$ is also a solution, so, let us find solutions with $z \geq 0$ and then take both signs. If $z=0$ then $(x,y)=(0,0)$. Let us assume that $z>0$. Then with $d=\gcd(x,z)$, $x=dX$ and $z=dZ$ we obtain 
\begin{equation}\label{mxp2x3mypy3_subx}
Z-d^2(X^3 + 3 X^2 Z - 3 X Z^2 + Z^3) = 0.
\end{equation}
This implies that $Z>0$ is a divisor of $d^2X^3$, hence it is also a divisor of $d^2$, so we can write $d^2=Zt$ for some integer $t>0$. Then \eqref{mxp2x3mypy3_subx} reduces to 
$$
1=t(X^3+3X^2Z-3XZ^2+Z^3).
$$
Hence $t>0$ is a divisor of $1$, which implies that $t=1$ and $X^3+3X^2Z-3XZ^2+Z^3=1$. The last equation is a Thue equation and it can be easily verified that its only integer solutions are $(X,Z)=(0,1)$ and $(1,0)$. 
Because $Z>0$, the only possible pair $(X,Z)$ is $(0,1)$. Then $d^2=Zt=1$ so $d=1$, and we obtain $z=dZ=1$, $x=dX=0$ and $y=z-x=1$. In conclusion, all integer solutions to \eqref{mxp2x3mypy3} are 
$$
(x,y)=(0,0),(0,\pm1).
$$

\vspace{10pt}

The next equation we will consider is
\begin{equation}\label{xpx3pypx2ymy3}
x+x^3+y+x^2y-y^3 = 0.
\end{equation}
Let $z=x+y$ be a new variable. Then $x=z-y$ and \eqref{xpx3pypx2ymy3} reduces to 
$$
z+(z-y)^3+y(z-y)^2-y^3 = -y^3 + z + y^2 z - 2 y z^2 + z^3 = 0.
$$
If $(z,y)$ is a solution to this equation, then $(-z,-y)$ is also a solution, so, let us find solutions with $z \geq 0$ and then take both signs. If $z=0$ then $(x,y)=(0,0)$. Let us assume that $z>0$. Then with $d=\gcd(y,z)$, $y=dY$ and $z=dZ$ we obtain  
\begin{equation}\label{xpx3pypx2ymy3_subx}
Z-d^2(Y^3-Y^2Z+2YZ^2-Z^3) = 0.
\end{equation}
This implies that $Z>0$ is a divisor of $d^2Y^3$, hence it is also a divisor of $d^2$, so we can write $d^2=Zt$ for some integer $t>0$. Then \eqref{xpx3pypx2ymy3_subx} reduces to 
$$
1=t(Y^3-Y^2Z+2YZ^2-Z^3).
$$
Hence $t>0$ is a divisor of $1$, which implies that $t=1$ and $Y^3-Y^2Z+2YZ^2-Z^3=1$. The last equation is a Thue equation and it can be easily verified that its only integer solutions are $(Y,Z)=(-1,-2),$ $(0,-1)$, $(1,0)$, $(1,1)$ and $(4,7)$. 
Because $Z>0$, the possible pairs $(Y,Z)$ are $(1,1)$ and $(4,7)$. From the first pair, we obtain 
$d^2=Zt=1$ so $d=1$, $y=dY=1$, $z=dZ=1$ and $x=z-y=1-1=0$. The second pair gives $d^2=Zt=7$ so $d=\sqrt{7}$, which is not integer. In conclusion, all integer solutions to \eqref{xpx3pypx2ymy3} are 
$$
(x,y)=(0,0),(0,\pm1).
$$

\vspace{10pt}

The next equation we will consider is
\begin{equation}\label{xpx3pypx2ypy3}
x+x^3+y+x^2y+y^3 = 0.
\end{equation}
Let $z=x+y$ be a new variable. Then $x=z-y$ and \eqref{xpx3pypx2ypy3} reduces to 
$$
z+(z-y)^3+y(z-y)^2+y^3 = y^3 + z + y^2 z - 2 y z^2 + z^3 =  0.
$$
If $(z,y)$ is a solution to this equation, then $(-z,-y)$ is also a solution, so, let us find solutions with $z \geq 0$ and then take both signs. If $z=0$ then $(x,y)=(0,0)$. Let us assume that $z>0$. Then with $d=\gcd(y,z)$, $y=dY$ and $z=dZ$ we obtain 
\begin{equation}\label{xpx3pypx2ypy3_subx}
Z-d^2(-Y^3-Y^2Z+2YZ^2-Z^3) = 0.
\end{equation}
This implies that $Z>0$ is a divisor of $d^2Y^3$, hence it is also a divisor of $d^2$, so we can write $d^2=Zt$ for some integer $t>0$. Then \eqref{xpx3pypx2ypy3_subx} reduces to 
$$
1=t(-Y^3-Y^2Z+2YZ^2-Z^3).
$$
Hence $t>0$ is a divisor of $1$, which implies that $t=1$ and $-Y^3-Y^2Z+2YZ^2-Z^3=1$. The last equation is a Thue equation and it can be easily verified that its only integer solutions are $(Y,Z)=(-1,-1),$ $(-1,0),$ $(0,-1),$ and $(2,-1)$. 
Because $Z>0$, there are no possible pairs $(Y,Z)$. In conclusion, the unique integer solution to \eqref{xpx3pypx2ypy3} is
$$
(x,y)=(0,0).
$$

\vspace{10pt}

The next equation we will consider is
\begin{equation}\label{xpx3mypx2ypy3}
x+x^3-y+x^2y+y^3 = 0.
\end{equation}
Let $z=x-y$ be a new variable. Then $x=y+z$ and \eqref{xpx3mypx2ypy3} reduces to 
$$
z+(y+z)^3+y(y+z)^2+y^3 = 3 y^3 + z + 5 y^2 z + 4 y z^2 + z^3 = 0.
$$
If $(z,y)$ is a solution to this equation, then $(-z,-y)$ is also a solution, so, let us find solutions with $z \geq 0$ and then take both signs. If $z=0$ then $(x,y)=(0,0)$. Let us assume that $z>0$. Then with $d=\gcd(y,z)$, $y=dY$ and $z=dZ$ we obtain 
\begin{equation}\label{xpx3mypx2ypy3_subx}
Z-d^2(-3Y^3-5Y^2Z-4YZ^2-Z^3) = 0.
\end{equation}
This implies that $Z>0$ is a divisor of $3d^2Y^3$, hence it is also a divisor of $3d^2$, so we can write $3d^2=Zt$ for some integer $t>0$. Then \eqref{xpx3mypx2ypy3_subx} reduces to 
$$
3=t(-3Y^3-5Y^2Z-4YZ^2-Z^3).
$$
Hence $t>0$ is a divisor of $3$, which implies that either (i) $t=1$ and $-3Y^3-5Y^2Z-4YZ^2-Z^3=3$ or (ii) $t=3$ and $-3Y^3-5Y^2Z-4YZ^2-Z^3=1$. In both cases, the second equation is a Thue equation, which can be easily solved in Mathematica. In case (i), its only integer solutions are $(Y,Z)=(-1,0)$ and $(1,-3)$. 
Because $Z>0$, there are no possible pairs $(Y,Z)$. In case (ii), its only integer solutions are $(Y,Z)=(-1,1),$ $(-1,2),$ $(0,-1)$ and $(2,-5)$. 
Because $Z>0$, the possible pairs $(Y,Z)$ are $(-1,1)$ and $(-1,2)$. From the first pair we obtain $3d^2=Zt=3$ so $d=1$, $z=dZ=1$, $y=dY=-1$ and $x=y+z=-1+1=0$. The second pair gives $3d^2=Zt=6$ so $d=\sqrt{6}$ which is not integer. In conclusion, all integer solutions to \eqref{xpx3mypx2ypy3} are 
$$
(x,y)=(0,0),(0,\pm1).
$$

\vspace{10pt}

The next equation we will consider is
\begin{equation}\label{xpx3mypx2ymy3}
x+x^3-y+x^2y-y^3 = 0.
\end{equation}
Let $z=x-y$ be a new variable. Then $x=y+z$ and \eqref{xpx3mypx2ymy3} reduces to 
$$
z+(y+z)^3+y(y+z)^2-y^3 = 0.
$$
If $(z,y)$ is a solution to this equation, then $(-z,-y)$ is also a solution, so, let us find solutions with $z \geq 0$ and then take both signs. If $z=0$ then $(x,y)=(0,0)$. Let us assume that $z>0$. Then with $d=\gcd(y,z)$, $y=dY$ and $z=dZ$ we obtain 
\begin{equation}\label{xpx3mypx2ymy3_subx}
Z-d^2(-Y^3-5Y^2Z-4YZ^2-Z^3) = 0.
\end{equation}
This implies that $Z>0$ is a divisor of $d^2Y^3$, hence it is also a divisor of $d^2$, so we can write $d^2=Zt$ for some integer $t>0$. Then \eqref{xpx3mypx2ymy3_subx} reduces to 
$$
1=t(-Y^3-5Y^2Z-4YZ^2-Z^3).
$$
Hence $t>0$ is a divisor of $1$, which implies that $t=1$ and $-Y^3-5Y^2Z-4YZ^2-Z^3=1$. The last equation is a Thue equation and it can be easily verified that its only integer solutions are $(Y,Z)=(-1,0),$ $(0,-1),$ $(1,-2),$ $(1,-1)$ and $(4,-1)$. 
 Because $Z>0$, there are no possible pairs $(Y,Z)$.  In conclusion, the unique integer solution to \eqref{xpx3mypx2ymy3} is 
$$
(x,y)=(0,0).
$$

\vspace{10pt}

The next equation we will consider is
\begin{equation}\label{mxpx3pypx2ypy3}
-x+x^3+y+x^2y+y^3 = 0.
\end{equation}
Let $z=y-x$ be a new variable. Then $y=x+z$ and \eqref{mxpx3pypx2ypy3} reduces to 
$$
z+x^3+x^2(x+z)+(x+z)^3 = 0.
$$
If $(x,z)$ is a solution to this equation, then $(-x,-z)$ is also a solution, so, let us find solutions with $z \geq 0$ and then take both signs. If $z=0$ then $(x,y)=(0,0)$. Let us assume that $z>0$. Then with $d=\gcd(x,z)$, $x=dX$ and $z=dZ$ we obtain 
\begin{equation}\label{mxpx3pypx2ypy3_subx}
Z-d^2(-3X^3-4X^2Z-3XZ^2-Z^3) = 0.
\end{equation}
This implies that $Z>0$ is a divisor of $3d^2X^3$, hence it is also a divisor of $3d^2$, so we can write $3d^2=Zt$ for some integer $t>0$. Then \eqref{mxpx3pypx2ypy3_subx} reduces to 
$$
3=t(-3X^3-4X^2Z-3XZ^2-Z^3).
$$
Hence $t>0$ is a divisor of $3$, which implies that either (i) $t=1$ and $-3X^3-4X^2Z-3XZ^2-Z^3=3$ or (ii) $t=3$ and $-3X^3-4X^2Z-3XZ^2-Z^3=1$. In both cases, the second equation is a Thue equation, which can be easily solved in Mathematica. In case (i), its only integer solutions are $(X,Z)=(-2,3)$ and $(-1,0)$. 
Because $Z>0$, the only possible pair $(X,Z)$ is $(-2,3)$ then $3d^2=Zt=3$ so $d=1$, $z=dZ=3$, $x=dX=-2$ and $y=x+z=-2+3=1$. In case (ii), its only integer solutions are $(X,Z)=(-3,5)$, $(-1,1)$, $(0,-1)$ and $(1,-2)$. 
 Because $Z>0$, the possible pairs $(X,Z)$ are $(-1,1)$ and $(-3,5)$. The first pair gives $3d^2=Zt=3$ so $d=1$, $z=dZ=1$, $x=dX=-1$ and $y=x+z=-1+1=0$. The second pair gives $3d^2=Zt=15$ so $d=\sqrt{5}$ which is not integer. In conclusion, all integer solutions to \eqref{mxpx3pypx2ypy3} are 
$$
(x,y)=(0,0),(\pm1,0),\pm(2,-1).
$$

\vspace{10pt}

The next equation we will consider is
\begin{equation}\label{mxpx3mypx2ymy3}
-x+x^3-y+x^2y-y^3 = 0.
\end{equation}
Let $z=x+y$ be a new variable. Then $x=z-y$ and \eqref{mxpx3mypx2ymy3} reduces to 
$$
-z+(z-y)^3+y(z-y)^2-y^3 = 0.
$$
If $(z,y)$ is a solution to this equation, then $(-z,-y)$ is also a solution, so, let us find solutions with $z \geq 0$ and then take both signs. If $z=0$ then $(x,y)=(0,0)$. Let us assume that $z>0$. Then with $d=\gcd(y,z)$, $y=dY$ and $z=dZ$ we obtain 
\begin{equation}\label{mxpx3mypx2ymy3_subx}
-Z+d^2(-Y^3+Y^2Z-2YZ^2+Z^3) = 0.
\end{equation}
This implies that $Z>0$ is a divisor of $d^2Y^3$, hence it is also a divisor of $d^2$, so we can write $d^2=Zt$ for some integer $t>0$. Then \eqref{mxpx3mypx2ymy3_subx} reduces to 
$$
1=t(-Y^3+Y^2Z-2YZ^2+Z^3).
$$
Hence $t>0$ is a divisor of $1$, which implies that $t=1$ and $-Y^3+Y^2Z-2YZ^2+Z^3=1$. The last equation is a Thue equation and it can be easily verified that its only integer solutions are $(Y,Z)=(-4,-7),$  $(-1,-1),$ $(-1,0),$ $(0,1)$ and $(1,2)$.
 Because $Z>0$, the possible pairs $(Y,Z)$ are $(0,1)$ and $(1,2)$. The first pair gives $d^2=Zt=1$ so $d=1$ so $y=dY=0$ and $z=dZ=1$ then $x=z-y=1-0=1$. The second pair gives $d^2=Zt=2$ so $d=\sqrt{2}$ which is not integer. In conclusion, all integer solutions to \eqref{mxpx3mypx2ymy3} are 
$$
(x,y)=(0,0),(\pm1,0).
$$

\vspace{10pt}

The next equation we will consider is
\begin{equation}\label{mxpx3mypx2ypy3}
-x+x^3-y+x^2y+y^3 = 0.
\end{equation}
Let $z=x+y$ be a new variable. Then $y=z-x$ and \eqref{mxpx3mypx2ypy3} reduces to 
$$
-z+x^3+x^2(z-x)+(z-x)^3 = 0.
$$
If $(x,z)$ is a solution to this equation, then $(-x,-z)$ is also a solution, so, let us find solutions with $z \geq 0$ and then take both signs. If $z=0$ then $(x,y)=(0,0)$. Let us assume that $z>0$. Then with $d=\gcd(x,z)$, $x=dX$ and $z=dZ$ we obtain 
\begin{equation}\label{mxpx3mypx2ypy3_subx}
Z-d^2(-X^3+4X^2Z-3XZ^2+Z^3) = 0.
\end{equation}
This implies that $Z>0$ is a divisor of $d^2X^3$, hence it is also a divisor of $d^2$, so we can write $d^2=Zt$ for some integer $t>0$. Then \eqref{mxpx3mypx2ypy3_subx} reduces to 
$$
1=t(-X^3+4X^2Z-3XZ^2+Z^3).
$$
Hence $t>0$ is a divisor of $1$, which implies that $t=1$ and $-X^3+4X^2Z-3XZ^2+Z^3=1$. The last equation is a Thue equation and it can be easily verified that its only integer solutions are $(X,Z)=(-1,0),$ $(0,1),$ $(1,1)$ and $(3,1)$. 
Because $Z>0$, the possible pairs $(X,Z)$ are $(0,1),(1,1)$ and $(3,1)$. From these pairs we obtain $d^2=Zt=1$, $d=1$ and $z=dZ=1$. 
Then in the original variables, we obtain that all integer solutions to \eqref{mxpx3mypx2ypy3} are 
$$
(x,y)=(\pm 1,0),(0,0),(0,\pm1),\pm (3,-2).
$$

\vspace{10pt}

The next equation we will consider is
\begin{equation}\label{mxpx3pypx2ymy3}
-x+x^3+y+x^2y-y^3 = 0.
\end{equation}
Let $z=y-x$ be a new variable. Then $y=z+x$ and \eqref{mxpx3pypx2ymy3} reduces to 
$$
z+x^3+x^2(z+x)-(z+x)^3 = 0.
$$
If $(x,z)$ is a solution to this equation, then $(-x,-z)$ is also a solution, so, let us find solutions with $z \geq 0$ and then take both signs. If $z=0$ then $(x,y)=(0,0)$. Let us assume that $z>0$. Then with $d=gcd(x,z)$, $x=dX$ and $z=dZ$ we obtain 
\begin{equation}\label{mxpx3pypx2ymy3_subx}
Z-d^2(-X^3+2X^2Z+3XZ^2+Z^3) = 0.
\end{equation}
This implies that $Z>0$ is a divisor of $d^2X^3$, hence it is also a divisor of $d^2$, so we can write $d^2=Zt$ for some integer $t>0$. Then \eqref{mxpx3pypx2ymy3_subx} reduces to 
$$
1=t(-X^3+2X^2Z+3XZ^2+Z^3).
$$
Hence $t>0$ is a divisor of $1$, which implies that $t=1$ and $-X^3+2X^2Z+3XZ^2+Z^3=1$. The last equation is a Thue equation and it can be easily verified that its only integer solutions are $(X,Z)=(-1,0),$ $(-1,1),$ $(-1,2),$ $(0,1)$ and $(3,1)$. 
 Because $Z>0$, the possible pairs $(X,Z)$ are $(-1,2),(0,1)$ and $(3,1)$. The first pair gives $d^2=Zt=2$ so $d=\sqrt{2}$ which is not integer. The second and third pair give $d^2=Zt=1$ so $d=1$ and $z=dZ=1$. Then in the original variables, we obtain that all
  integer solutions to \eqref{mxpx3pypx2ymy3} are 
$$
(x,y)=(\pm 1,0),(0,0),(0,\pm1),\pm(3,4).
$$

\vspace{10pt}

The next equation we will consider is
\begin{equation}\label{x3p2ypx2ypy3}
x^3+2y+x^2y+y^3=0.
\end{equation}
If $(x,y)$ is a solution to this equation, then $(-x,-y)$ is also a solution, so let us find solutions with $y \geq 0$ and then take both signs. If $y=0$ then $x=0$. Next, let $(x,y)$ be any solution with $y>0$, and let $d=\gcd(x,y)$. Then $x=dX$ and $y=dY$ with $X,Y$ coprime and $Y>0$. Substituting this into \eqref{x3p2ypx2ypy3} and cancelling $d$, we obtain
\begin{equation}\label{x3p2ypx2ypy3_x1y1}
2Y=d^2(-X^3-X^2Y-Y^3).
\end{equation}
From this equation, it is clear that $Y$ is a divisor of $d^2X^3$. Because $X$ and $Y$ are coprime, this implies that $Y$ is a divisor of $d^2$, so we can write $d^2=Yz$ for some integer $z>0$. Then after substitution and cancelling $Y$, \eqref{x3p2ypx2ypy3_x1y1} reduces to 
$$
2=z(-X^3-X^2Y-Y^3).
$$
Hence, $z>0$ is a divisor of $2$, which implies that either (i) $z=1$ and $-X^3-X^2Y-Y^3=2$ or (ii) $z=2$ and $-X^3-X^2Y-Y^3=1$. In both cases, the second equation is a Thue equation, which can be easily solved in Mathematica. In case (i), it has no integer solutions. 
In case (ii), its only integer solutions are $(X,Y)=(-3,2),$ $(-1,0),$ $(0,-1)$ and $(1,-1)$. 
Because $Y>0$, the only possible pair $(X,Y)$ is $(-3,2)$. This solution gives $d^2=Yz=4$ so $d=2$, then $x=dX=-6$ and $y=dY=4$. In conclusion, all integer solutions to \eqref{x3p2ypx2ypy3} are
$$
(x,y)=(0,0),\pm(6,-4).
$$

\vspace{10pt}

The next equation we will consider is
\begin{equation}\label{x3m2ypx2ymy3}
x^3-2y+x^2y-y^3=0.
\end{equation}
If $(x,y)$ is a solution to this equation, then $(-x,-y)$ is also a solution, so let us find solutions with $y \geq 0$ and then take both signs. If $y=0$ then $x=0$. Next, let $(x,y)$ be any solution with $y>0$, and let $d=\gcd(x,y)$. Then $x=dX$ and $y=dY$ with $X,Y$ coprime and $Y>0$. Substituting this into \eqref{x3m2ypx2ymy3} and cancelling $d$, we obtain
\begin{equation}\label{x3m2ypx2ymy3_x1y1}
2Y=d^2(X^3+X^2Y-Y^3).
\end{equation}
From this equation, it is clear that $Y$ is a divisor of $d^2X^3$. Because $X$ and $Y$ are coprime, this implies that $Y$ is a divisor of $d^2$, so we can write $d^2=Yz$ for some integer $z>0$. Then after substitution and cancelling $Y$, \eqref{x3m2ypx2ymy3_x1y1} reduces to 
$$
2=z(X^3+X^2Y-Y^3).
$$
Hence, $z>0$ is a divisor of $2$, which implies that either (i) $z=1$ and $X^3+X^2Y-Y^3=2$ or (ii) $z=2$ and $X^3+X^2Y-Y^3=1$. In both cases, the second equation is a Thue equation, which can be easily solved in Mathematica. In case (i), it has no integer solutions. 
In case (ii), its only integer solutions are $(X,Y)=(-3,-4),$ $(0,-1),$ $(1,\pm 1)$ and $(1,0)$. 
Because $Y>0$, the only possible pair $(X,Y)$ is $(1,1)$. This solution gives $d^2=Yz=2$ so $d=\sqrt{2}$, which is not integer. In conclusion, the unique integer solution to \eqref{x3m2ypx2ymy3} is
$$
(x,y)=(0,0).
$$

\vspace{10pt}

The next equation we will consider is
\begin{equation}\label{x3m2ypx2ypy3}
x^3-2y+x^2y+y^3=0.
\end{equation}
If $(x,y)$ is a solution to this equation, then $(-x,-y)$ is also a solution, so let us find solutions with $y \geq 0$ and then take both signs. If $y=0$ then $x=0$. Next, let $(x,y)$ be any solution with $y>0$, and let $d=\gcd(x,y)$. Then $x=dX$ and $y=dY$ with $X,Y$ coprime and $Y>0$. Substituting this into \eqref{x3m2ypx2ypy3} and cancelling $d$, we obtain
\begin{equation}\label{x3m2ypx2ypy3_x1y1}
2Y=d^2(X^3+X^2Y+Y^3).
\end{equation}
From this equation, it is clear that $Y$ is a divisor of $d^2X^3$. Because $X$ and $Y$ are coprime, this implies that $Y$ is a divisor of $d^2$, so we can write $d^2=Yz$ for some integer $z>0$. Then after substitution and cancelling $Y$, \eqref{x3m2ypx2ypy3_x1y1} reduces to 
$$
2=z(X^3+X^2Y+Y^3).
$$
Hence, $z>0$ is a divisor of $2$, which implies that either (i) $z=1$ and $X^3+X^2Y+Y^3=2$ or (ii) $z=2$ and $X^3+X^2Y+Y^3=1$. In case (i), it has no integer solutions. 
In case (ii), its only integer solutions are $(X,Y)=(-1,1)$, $(0,1)$, $(1,0)$ and $(3,-2)$. 
Because $Y>0$, the only possible pairs $(X,Y)$ are $(-1,1)$ and $(0,1)$. These solutions give $d^2=Yz=2$ so $d=\sqrt{2}$, which is not integer. In conclusion, the unique integer solution to \eqref{x3m2ypx2ypy3} is
$$
(x,y)=(0,0).
$$

\vspace{10pt}

The next equation we will consider is
\begin{equation}\label{x3p2ypx2ymy3}
x^3+2y+x^2y-y^3=0.
\end{equation}
If $(x,y)$ is a solution to this equation, then $(-x,-y)$ is also a solution, so let us find solutions with $y \geq 0$ and then take both signs. If $y=0$ then $x=0$. Next, let $(x,y)$ be any solution with $y>0$, and let $d=\gcd(x,y)$. Then $x=dX$ and $y=dY$ with $X,Y$ coprime and $Y>0$. Substituting this into \eqref{x3p2ypx2ymy3} and cancelling $d$, we obtain
\begin{equation}\label{x3p2ypx2ymy3_x1y1}
2Y=d^2(-X^3-X^2Y+Y^3).
\end{equation}
From this equation, it is clear that $Y$ is a divisor of $d^2X^3$. Because $X$ and $Y$ are coprime, this implies that $Y$ is a divisor of $d^2$, so we can write $d^2=Yz$ for some integer $z>0$. Then after substitution and cancelling $Y$, \eqref{x3p2ypx2ymy3_x1y1} reduces to 
$$
2=z(-X^3-X^2Y+Y^3).
$$
Hence, $z>0$ is a divisor of $2$, which implies that either (i) $z=1$ and $-X^3-X^2Y+Y^3=2$ or (ii) $z=2$ and $-X^3-X^2Y+Y^3=1$. In case (i), it has no integer solutions. 
In case (ii), its only integer solutions are $(X,Y)=(-1,\pm 1),$ $(-1,0),$ $(0,1)$ and $(3,4)$. 
Because $Y>0$, the only possible pairs $(X,Y)$ are $(-1,1),(0,1)$ and $(3,4)$. These solutions do not correspond to integer $d$. 
In conclusion, the unique integer solution to \eqref{x3p2ypx2ymy3} is
$$
(x,y)=(0,0).
$$

\vspace{10pt}

The next equation we will consider is
\begin{equation}\label{x3p2xpx2ypy3}
x^3+2x+x^2y+y^3=0.
\end{equation}
If $(x,y)$ is a solution to this equation, then $(-x,-y)$ is also a solution, so let us find solutions with $x \geq 0$ and then take both signs. If $x=0$ then $y=0$. Next, let $(x,y)$ be any solution with $x>0$, and let $d=\gcd(x,y)$. Then $x=dX$ and $y=dY$ with $X,Y$ coprime and $X>0$. Substituting this into \eqref{x3p2xpx2ypy3} and cancelling $d$, we obtain
\begin{equation}\label{x3p2xpx2ypy3_x1y1}
2X=d^2(-X^3-X^2Y-Y^3).
\end{equation}
From this equation, it is clear that $X$ is a divisor of $d^2Y^3$. Because $X$ and $Y$ are coprime, this implies that $X$ is a divisor of $d^2$, so we can write $d^2=Xz$ for some integer $z>0$. Then after substitution and cancelling $X$, \eqref{x3p2xpx2ypy3_x1y1} reduces to 
$$
2=z(-X^3-X^2Y-Y^3).
$$
Hence, $z>0$ is a divisor of $2$, which implies that either (i) $z=1$ and $-X^3-X^2Y-Y^3=2$ or (ii) $z=2$ and $-X^3-X^2Y-Y^3=1$. In both cases, the second equation is a Thue equation, which can be easily solved in Mathematica. In case (i), it has no integer solutions. 
In case (ii), its only integer solutions are $(X,Y)=(-3,2),$ $(-1,0),$ $(0,-1)$ and $(1,-1)$. 
Because $X>0$, the only possible pair $(X,Y)$ is $(1,-1)$. This solution gives $d^2=Xz=2$ so $d=\sqrt{2}$, which is not integer. In conclusion, the unique integer solution to \eqref{x3p2xpx2ypy3} is
$$
(x,y)=(0,0).
$$

\vspace{10pt}

The next equation we will consider is
\begin{equation}\label{x3p2xpx2ymy3}
x^3+2x+x^2y-y^3=0.
\end{equation}
If $(x,y)$ is a solution to this equation, then $(-x,-y)$ is also a solution, so let us find solutions with $x \geq 0$ and then take both signs. If $x=0$ then $y=0$. Next, let $(x,y)$ be any solution with $x>0$, and let $d=\gcd(x,y)$. Then $x=dX$ and $y=dY$ with $X,Y$ coprime and $X>0$. Substituting this into \eqref{x3p2xpx2ymy3} and cancelling $d$, we obtain
\begin{equation}\label{x3p2xpx2ymy3_x1y1}
2X=d^2(-X^3-X^2Y+Y^3).
\end{equation}
From this equation, it is clear that $X$ is a divisor of $d^2Y^3$. Because $X$ and $Y$ are coprime, this implies that $X$ is a divisor of $d^2$, so we can write $d^2=Xz$ for some integer $z>0$. Then after substitution and cancelling $X$, \eqref{x3p2xpx2ymy3_x1y1} reduces to 
$$
2=z(-X^3-X^2Y+Y^3).
$$
Hence, $z>0$ is a divisor of $2$, which implies that either (i) $z=1$ and $-X^3-X^2Y+Y^3=2$ or (ii) $z=2$ and $-X^3-X^2Y+Y^3=1$. In case (i), it has no integer solutions. 
In case (ii), its only integer solutions are $(X,Y)=(-1,\pm 1),$ $(-1,0),$ $(0,1)$ and $(3,4)$. 
Because $X>0$, the only possible pair $(X,Y)$ is $(3,4)$. This solution gives $d^2=Xz=6$ so $d=\sqrt{6}$, which is not integer. In conclusion, the unique integer solution to \eqref{x3p2xpx2ymy3} is
$$
(x,y)=(0,0).
$$

\vspace{10pt}

The next equation we will consider is
\begin{equation}\label{x3m2xpx2ypy3}
x^3-2x+x^2y+y^3=0.
\end{equation}
If $(x,y)$ is a solution to this equation, then $(-x,-y)$ is also a solution, so let us find solutions with $x \geq 0$ and then take both signs. If $x=0$ then $y=0$. Next, let $(x,y)$ be any solution with $x>0$, and let $d=\gcd(x,y)$. Then $x=dX$ and $y=dY$ with $X,Y$ coprime and $X>0$. Substituting this into \eqref{x3m2xpx2ypy3} and cancelling $d$, we obtain
\begin{equation}\label{x3m2xpx2ypy3_x1y1}
2X=d^2(X^3+X^2Y+Y^3).
\end{equation}
From this equation, it is clear that $X$ is a divisor of $d^2Y^3$. Because $X$ and $Y$ are coprime, this implies that $X$ is a divisor of $d^2$, so we can write $d^2=Xz$ for some integer $z>0$. Then after substitution and cancelling $X$, \eqref{x3m2xpx2ypy3_x1y1} reduces to 
$$
2=z(X^3+X^2Y+Y^3).
$$
Hence, $z>0$ is a divisor of $2$, which implies that either (i) $z=1$ and $X^3+X^2Y+Y^3=2$ or (ii) $z=2$ and $X^3+X^2Y+Y^3=1$. In both cases, the second equation is a Thue equation, which can be easily solved in Mathematica. In case (i), it has no integer solutions. 
In case (ii), its only integer solutions are $(X,Y)=(-1,1),$ $(0,1),$ $(1,0)$ and $(3,-2)$. 
Because $X>0$, the only possible pairs $(X,Y)$ are $(1,0)$ and $(3,-2)$. The first solution gives $d^2=Xz=2$ so $d=\sqrt{2}$, which is not integer. The second solution gives $d^2=Xz=6$ so $d=\sqrt{6}$, which is not integer. In conclusion, the unique integer solution to \eqref{x3m2xpx2ypy3} is
$$
(x,y)=(0,0).
$$

\vspace{10pt}

The final equation we will consider is
\begin{equation}\label{x3m2xpx2ymy3}
x^3-2x+x^2y-y^3=0.
\end{equation}
If $(x,y)$ is a solution to this equation, then $(-x,-y)$ is also a solution, so let us find solutions with $x \geq 0$ and then take both signs. If $x=0$ then $y=0$. Next, let $(x,y)$ be any solution with $x>0$, and let $d=\gcd (x,y)$. Then $x=dX$ and $y=dY$ with $X,Y$ coprime and $X>0$. Substituting this into \eqref{x3m2xpx2ymy3} and cancelling $d$, we obtain
\begin{equation}\label{x3m2xpx2ymy3_x1y1}
2X=d^2(X^3+X^2Y-Y^3).
\end{equation}
From this equation, it is clear that $X$ is a divisor of $d^2Y^3$. Because $X$ and $Y$ are coprime, this implies that $X$ is a divisor of $d^2$, so we can write $d^2=Xz$ for some integer $z>0$. Then after substitution and cancelling $X$, \eqref{x3m2xpx2ymy3_x1y1} reduces to 
$$
2=z(X^3+X^2Y-Y^3).
$$
Hence, $z>0$ is a divisor of $2$, which implies that either (i) $z=1$ and $X^3+X^2Y-Y^3=2$ or (ii) $z=2$ and $X^3+X^2Y-Y^3=1$. In case (i), it has no integer solutions. 
In case (ii), its only integer solutions are $(X,Y)=(-3,-4),$ $(0,-1),$ $(1,\pm 1)$ and $(1,0)$. 
 Because $X>0$, the only possible pairs $(X,Y)$ are $(1,\pm 1)$ and $(-1,0)$. These solutions give $d^2=Xz=2$ so $d=\sqrt{2}$, which is not integer. In conclusion, the unique integer solution to \eqref{x3m2xpx2ymy3} is
$$
(x,y)=(0,0).
$$

\vspace{10pt}

Table \ref{tab:H28axbyPxysol} summarises all integer solutions to the equations listed in Table \ref{tab:H28axbyPxy}.

\begin{center}
\begin{tabular}{ |c|c|c| } 
		\hline
		Equation & Solution $(x,y)$ \\ 
		\hline\hline
		$-x+x^3+x^2 y-y^3 = 0$ &$(\pm 1,\pm 1),(\pm 1,0),(0,0)$ \\\hline
		$x^3-y+x^2 y-y^3 = 0$ &$(0,0),\pm(1,1)$ \\\hline
		$x^3+y+x^2 y-y^3 = 0$ & $(0,\pm 1),(0,0),\pm(1,-1),\pm(6,8)$ \\\hline
		$x+x^3+x^2 y-y^3 = 0$ &$(0,0)$\\\hline
		$-x+x^3+x^2 y+y^3 = 0$ &$(\pm1 ,0),(0,0)$ \\\hline
		$x^3-y+x^2 y+y^3 = 0$ &$(0,\pm 1),(0,0),\pm(1,-1)$ \\\hline
		$x^3+y+x^2 y+y^3 = 0$ &$(0,0)$ \\\hline
		$x+x^3+x^2 y+y^3 = 0$ &$(0,0),\pm(1,-1)$ \\\hline
		$-2 x+x^3+x^2 y-y^3 = 0$ &$(0,0)$\\\hline
		$-x+x^3-y+x^2 y-y^3 = 0$ &$(0,0),(\pm 1,0)$ \\\hline
		$-x+x^3+y+x^2 y-y^3 = 0$ &$(\pm 1,0),(0,0),(0,\pm 1),\pm(3,4)$ \\\hline
		$x^3-2 y+x^2 y-y^3 = 0$ &$(0,0)$ \\\hline
		$x^3+2 y+x^2 y-y^3 = 0$ &$(0,0)$ \\\hline
		$x+x^3-y+x^2 y-y^3 = 0$ &$(0,0)$\\\hline
		$x+x^3+y+x^2 y-y^3 = 0$ & $(0,0),(0,\pm1)$ \\\hline
		$2 x+x^3+x^2 y-y^3 = 0$ &$(0,0)$  \\\hline
		$-2 x+x^3+x^2 y+y^3 = 0$ & $(0,0)$\\\hline    
		$-x+x^3-y+x^2 y+y^3 = 0$ &$(\pm 1,0),(0,0),(0,\pm 1),\pm(3,-2)$\\\hline  
		$-x+x^3+y+x^2 y+y^3 = 0$ & $(\pm 1,0),(0,0),\pm(2,-1)$   \\ \hline
		$x^3-2 y+x^2 y+y^3 = 0$ & $(0,0)$   \\  \hline
		$x^3+2 y+x^2 y+y^3 = 0$&$(0,0),\pm(6,-4)$ \\  \hline
		$x+x^3-y+x^2 y+y^3 = 0$&$(0,0),(0, \pm 1)$ \\  \hline
		$x+x^3+y+x^2 y+y^3 = 0$&$(0,0)$ \\  \hline
		$2 x+x^3+x^2 y+y^3 = 0$&$(0,0)$ \\  \hline
		$-x+2 x^3-y+y^3 = 0$&$(0,0),(0,\pm 1)$ \\  \hline
		$-x+2 x^3+y+y^3 = 0$&$(0,0)$ \\  \hline
		$x+2 x^3-y+y^3 = 0$ &$(0,0),(0,\pm 1)$\\  \hline
		$x+2 x^3+y+y^3 = 0$&$(0,0)$ \\  \hline
	\end{tabular}
	\captionof{table}{\label{tab:H28axbyPxysol} Integer solutions to the equations in Table \ref{tab:H28axbyPxy}.}
\end{center}

\subsection{Exercise 3.54}\label{ex:genus0}
\textbf{\emph{For each equation in Table \ref{tab:H28genus0}, (i) describe all rational solutions in parametric form, and (ii) use this to find all its integer solutions.  }}

\begin{center}
\begin{tabular}{ |c|c|c|c|c|c| } 
 \hline
 $H$ & Equation \\ 
 \hline\hline
 $28$ & $2x^3+y^3+xy=0$ \\ 
 \hline
 $28$ & $x^3+x^2 y-y^3+y^2=0$ \\ 
 \hline
 $28$ & $x^3+x^2 y+y^3+y^2=0$ \\ 
 \hline
 $28$ & $x^3+x^2 y-y^3+xy=0$ \\ 
 \hline
 $28$ & $x^3+x^2 y+y^3+xy=0$ \\ 
 \hline
 $28$ & $x^3+x^2 y-y^3+x^2=0$ \\ 
 \hline
 $28$ & $x^3+x^2 y+y^3+x^2=0$ \\ 
 \hline
\end{tabular}
\captionof{table}{\label{tab:H28genus0} Equations of genus $0$ of size $H\leq 28$.}
\end{center} 

We will now look at equations with genus 0 of size $H \leq 28$. To solve these equations, we will first describe all rational solutions to the equation, and then use this description to find all integer solutions to the equation. 

To solve these equations, we will use the method from Section 3.4.2 of the book, which we summarise below for convenience.
To describe all rational solutions to these equations, we must first find a point on the curve, which we denote as $(x_0,y_0)$. All curves in Table \ref{tab:H28genus0} have the integer point $(x_0,y_0)=(0,0)$. We can then find other rational points on the curve by drawing lines with rational slopes via $(x_0,y_0)$ and compute intersections of these lines with the curve. The general equation of a line has the form $y=kx+b$, as we want $(x_0,y_0)=(0,0)$ to be on the line, we must have $b=y_0-kx_0=0$, so we will consider all lines in the form $y=kx$ where $k$ is a rational parameter. We then substitute $y=kx$ into the equation to find points of intersection of this line with the curve. For each equation in Table \ref{tab:H28genus0} we obtain an equation of the form
$$
x^2 \cdot P(x,k)=0
$$
where $P$ is linear in $x$. Hence either $x=0$ or $P(x,k)=0$. Because $P(x,k)$ is linear in $x$ we can easily express $x$ in terms of $k$. We can then describe all rational solutions $(x,y)$. To use this description to determine all integer solutions to the equation, we represent $k$ as an irreducible fraction $\frac{a}{b}$ and check for which values of $a,b$ the rational solutions $(x,y)$ are in fact integer solutions. 
Let us illustrate this method on a number of examples.  

Equation 
$$
2x^3+y^3+xy=0
$$ 
is solved in Section 3.4.2 of the book, and its rational solutions are of the form
$$
x= -\frac{k}{k^3+2} \quad \text{and} \quad y = - \frac{k^2}{k^3+2} \quad k \in \mathbb{Q},
$$
and all its integer solutions are
$$
(x,y)=(0,0) \quad \text{and} \quad (1,-1).
$$

The first equation we will consider is
\begin{equation}\label{x3px2ymy3py2}
x^3+x^2y-y^3+y^2=0.
\end{equation}
Let us first describe all rational solutions to \eqref{x3px2ymy3py2}. The rational point $(x_0,y_0)=(0,0)$ is on the curve. 
 After substituting $y=kx$, where $k$ is a rational parameter, into \eqref{x3px2ymy3py2}, we obtain
$$
x^3+x^2(kx)-(kx)^3+(kx)^2=0=x^2(k^2+x+kx-k^3x).
$$  
If $x=0$ we find the point $(x_0,y_0)=(0,0)$ that we started with, or $(x,y)=(0,1)$. Otherwise, $k^2+x+kx-k^3x=0$ from which we find that
\begin{equation}\label{x3px2ymy3py2_rational}
(x,y)=(0,1) \quad \text{or} \quad x=\frac{-k^2}{1+k-k^3}, \quad \text{and} \quad y=kx=\frac{-k^3}{1+k-k^3}, \quad k \in \mathbb{Q}.
\end{equation}
We remark that $1+k-k^3\neq 0$ for rational $k$.

We will now use this description to find all integer solutions to \eqref{x3px2ymy3py2}. 
Let us represent $k=\frac{a}{b}$ as an irreducible fraction with $b>0$ and $\gcd(a,b)=1$. Then 
$$
x=\frac{-k^2}{1+k-k^3}=\frac{-(\frac{a}{b})^2}{1+\frac{a}{b}-(\frac{a}{b})^3}=\frac{a^2b}{a^3-ab^2-b^3}, \text{ and } y=\frac{-k^3}{1+k-k^3}=\frac{-(\frac{a}{b})^3}{1+\frac{a}{b}-(\frac{a}{b})^3}=\frac{a^3}{a^3-ab^2-b^3}.
$$
We need to find $a$ and $b$ such that both $x$ and $y$ are integers. Let $p$ be any prime factor of $a^3-ab^2-b^3$. Because $y$ is an integer, $p$ must also be a prime factor of $a^3$, hence $p$ is a prime factor of $a$, but 
in this case $p$ must also be a prime factor of $(a^3-ab^2-b^3)-a(a^2-b^2)=-b^3$, which is a contradiction with the coprimality of $a$ and $b$. Hence, $a^3-ab^2-b^3$ has no prime factors, this is only possible if $a^3-ab^2-b^3=\pm 1$. These are Thue equations and we can easily verify that their only integer solutions are 
 $(a,b)=(\pm1,\pm1),(\pm 1,0),(0,\pm 1),\pm(4,3)$. Because $b>0$, the only possible pairs $(a,b)$ are $(\pm1,1),(0,1)$, and $(4,3)$. Hence, $k=\pm 1,0$ or $\frac{4}{3}$. Then by \eqref{x3px2ymy3py2_rational} we can conclude that the integer solutions to \eqref{x3px2ymy3py2} are
$$
(x,y)=(-1,\pm 1),(0,0),(0,1),(48,64).
$$

\vspace{10pt}

The next equation we will consider is
\begin{equation}\label{x3px2ypy3py2}
x^3+x^2y+y^3+y^2=0.
\end{equation}
Let us first describe all rational solutions to \eqref{x3px2ypy3py2}. 
The rational point $(x_0,y_0)=(0,0)$ is on the curve. 
After substituting $y=kx$, where $k$ is a rational parameter, into \eqref{x3px2ypy3py2}, 
 we obtain
$$
x^3+x^2(kx)+(kx)^3+(kx)^2=0=x^2(k^2+x+kx+k^3x).
$$  
If $x=0$ we find the point $(x_0,y_0)=(0,0)$ that we started with, or $(x,y)=(0,-1)$. Otherwise, $k^2+x+kx+k^3x=0$ from which we find that
\begin{equation}\label{x3px2ypy3py2_rational}
(x,y)=(0,-1) \quad \text{or} \quad x=\frac{-k^2}{1+k+k^3},\quad  \text{and} \quad y=kx=\frac{-k^3}{1+k+k^3}, \quad k \in \mathbb{Q}.
\end{equation}
We remark that $1+k+k^3\neq 0$ for rational $k$.

We will now use this description to find all integer solutions to \eqref{x3px2ypy3py2}. 
Let us represent $k=\frac{a}{b}$ as an irreducible fraction with $b>0$ and $\gcd(a,b)=1$. Then 
$$
x=\frac{-k^2}{1+k+k^3}=\frac{-(\frac{a}{b})^2}{1+\frac{a}{b}+(\frac{a}{b})^3}=\frac{a^2b}{a^3+ab^2+b^3}, \text{ and } y=\frac{-k^3}{1+k+k^3}=\frac{-(\frac{a}{b})^3}{1+\frac{a}{b}+(\frac{a}{b})^3}=\frac{a^3}{a^3+ab^2+b^3}.
$$
Let $p$ be any prime factor of $a^3+ab^2+b^3$. Because $y$ is an integer, $p$ must also be a prime factor of $a^3$, hence $p$ is a prime factor of $a$, but then $p$ is also a prime factor of $(a^3+ab^2+b^3)-a(a^2+b^2)=b^3$, which is a contradiction with the coprimality of $a$ and $b$. Hence, $a^3+ab^2+b^3$ has no prime factors, which is only possible if $a^3+ab^2+b^3=\pm 1$. These are Thue equations and we can easily verify that their only integer solutions satisfying $b>0$ are $(-2,3),(-1,1),(0,1)$. Hence, 
$k=-\frac{2}{3},-1$ or $0$. Then by \eqref{x3px2ypy3py2_rational} we can conclude that the integer solutions to \eqref{x3px2ypy3py2} are
$$
(x,y)=(-12,8),(0,-1),(0,0),(1,-1).
$$

\vspace{10pt}

The next equation we will consider is
\begin{equation}\label{x3px2ymy3pxy}
x^3+x^2y-y^3+xy=0.
\end{equation}
Let us first describe all rational solutions to \eqref{x3px2ymy3pxy}. The rational point $(x_0,y_0)=(0,0)$ is on the curve. After substituting $y=kx$, where $k$ is a rational parameter, into \eqref{x3px2ymy3pxy}, we obtain 
$$
x^3+x^2(kx)-(kx)^3+x(kx)=0=x^2(k+x+kx-k^3x).
$$  
If $x=0$ we find the point $(x_0,y_0)=(0,0)$ that we started with. Otherwise, $k+x+kx-k^3x=0$ from which we find that
\begin{equation}\label{x3px2ymy3pxy_rational}
x=\frac{-k}{1+k-k^3}, \quad \text{and} \quad y=kx=\frac{-k^2}{1+k-k^3}, \quad k \in \mathbb{Q}.
\end{equation}
We remark that $1+k-k^3\neq 0$ for rational $k$.

We will now use this description to find all integer solutions to \eqref{x3px2ymy3pxy}. 
Let us represent $k=\frac{a}{b}$ as an irreducible fraction with $b>0$ and $\gcd(a,b)=1$. Then 
$$
x=\frac{-k}{1+k-k^3}=\frac{-\frac{a}{b}}{1+\frac{a}{b}-(\frac{a}{b})^3}=\frac{ab^2}{a^3-ab^2-b^3}, \text{ and } y=\frac{-k^2}{1+k-k^3}=\frac{-(\frac{a}{b})^2}{1+\frac{a}{b}-(\frac{a}{b})^3}=\frac{a^2b}{a^3-ab^2-b^3}.
$$
Let $p$ be any prime factor of $a^3-ab^2-b^3$. Because $x$ is integer, $p$ must also be a prime factor of $ab^2$, hence either (i) $p$ is a prime factor of $a$ or (ii) $p$ is a prime factor of $b$. In case (i), $p$ is also a prime factor of $a^3-ab^2-b^3-a(a^2-b^2)=-b^3$, hence $p$ must be a prime factor of $b$, but this is a contradiction with $\gcd(a,b)=1$. Hence $a^3-ab^2-b^3$ has no prime factors, which is only possible if $a^3-ab^2-b^3=\pm 1$. These are Thue equations and we can easily verify that their only integer solutions satisfying $b>0$ are 
 $(\pm1,1),(0,1)$, and $(4,3)$. Hence, $k=\pm 1,0$ or $\frac{4}{3}$. In case (ii), $p$ is also a prime factor of $(a^3-ab^2-b^3)+b(ab+b^2)=a^3$, hence $p$ must be a prime factor of $a$, but this is a contradiction with $\gcd(a,b)=1$. Hence 
 $a^3-ab^2-b^3=\pm 1$ which we have already considered. 
 Finally, by \eqref{x3px2ymy3pxy_rational} and $k=\pm 1,0$ or $\frac{4}{3}$, we can conclude that the integer solutions to \eqref{x3px2ymy3pxy} are
$$
(x,y)=(0,0),(\pm 1,-1),(36,48).
$$

\vspace{10pt}

The next equation we will consider is
\begin{equation}\label{x3px2ypy3pxy}
x^3+x^2y+y^3+xy=0.
\end{equation}
Let us first describe all rational solutions to \eqref{x3px2ypy3pxy}. The rational point $(x_0,y_0)=(0,0)$ is on the curve. After substituting $y=kx$, where $k$ is a rational parameter, into \eqref{x3px2ypy3pxy}, we obtain
$$
x^3+x^2(kx)+(kx)^3+x(kx)=0=x^2(k+x+kx+k^3x).
$$  
If $x=0$ we find the point $(x_0,y_0)=(0,0)$ that we started with. Otherwise, $k+x+kx+k^3x=0$ from which we find that
\begin{equation}\label{x3px2ypy3pxy_rational}
x=\frac{-k}{1+k+k^3}, \quad \text{and} \quad y=kx=\frac{-k^2}{1+k+k^3}, \quad k \in \mathbb{Q}.
\end{equation}
We remark that $1+k+k^3\neq 0$ for rational $k$.

We will now use this description to find all integer solutions to \eqref{x3px2ypy3pxy}. 
Let us represent $k=\frac{a}{b}$ as an irreducible fraction with $b>0$ and $\text{gcd}(a,b)=1$. Then 
$$
x=\frac{-k}{1+k+k^3}=\frac{-\frac{a}{b}}{1+\frac{a}{b}+(\frac{a}{b})^3}=\frac{-ab^2}{a^3+ab^2+b^3}, \text{ and } y=\frac{-k^2}{1+k+k^3}=\frac{-(\frac{a}{b})^2}{1+\frac{a}{b}+(\frac{a}{b})^3}=\frac{-a^2b}{a^3+ab^2+b^3}.
$$
Let $p$ be any prime factor of $a^3+ab^2+b^3$. Because $x$ is integer, $p$ must also be a prime factor of $-ab^2$, hence either (i) $p$ is a prime factor of $a$ or (ii) $p$ is a prime factor of $b$. In case (i), $p$ is also a prime factor of $(a^3+ab^2+b^3)-a(a^2+b^2)=b^3$, hence $p$ must be a prime factor of $b$, but this is a contradiction with $\gcd(a,b)=1$. Hence $a^3+ab^2+b^3$ has no prime factors, which is only possible if $a^3+ab^2+b^3=\pm 1$. These are Thue equations and we can easily verify that their only integer solutions satisfying $b>0$ are 
$(-2,3),(-1,1)$, and $(0,1)$. These solutions give $k=-1,-\frac{2}{3}$ or $0$. In case (ii), $p$ is also a prime factor of $a^3+ab^2+b^3-b(ab+b^2)=a^3$, hence $p$ must be a prime factor of $a$, but this is a contradiction with $\gcd(a,b)=1$. Hence 
$a^3+ab^2+b^3=\pm 1$ which we have already considered.
 Finally, by \eqref{x3px2ypy3pxy_rational} and $k=-1,-\frac{2}{3},0$, we can conclude that the integer solutions to \eqref{x3px2ypy3pxy} are
$$
(x,y)=(-1,1),(0,0),(18,-12).
$$

\vspace{10pt}

The next equation we will consider is
\begin{equation}\label{x3px2ymy3px2}
x^3+x^2y-y^3+x^2=0.
\end{equation}
Let us first describe all rational solutions to \eqref{x3px2ymy3px2}. The rational point $(x_0,y_0)=(0,0)$ is on the curve. After substituting $y=kx$, where $k$ is a rational parameter, into \eqref{x3px2ymy3px2}, we obtain
$$
x^3+x^2(kx)-(kx)^3+x^2=0=x^2(1+x+kx-k^3x).
$$  
If $x=0$ we find the point $(x_0,y_0)=(0,0)$ that we started with. Otherwise, $1+x+kx-k^3x=0$ from which we find that
\begin{equation}\label{x3px2ymy3px2_rational}
x=\frac{-1}{1+k-k^3}, \quad \text{and} \quad y=kx=\frac{-k}{1+k-k^3}, \quad k \in \mathbb{Q}.
\end{equation}
We remark that $1+k-k^3\neq 0$ for rational $k$.

We will now use this description to find all integer solutions to \eqref{x3px2ymy3px2}. 
Let us represent $k=\frac{a}{b}$ as an irreducible fraction with $b>0$ and $\gcd(a,b)=1$. Then 
$$
x=\frac{-1}{1+k-k^3}=\frac{-1}{1+\frac{a}{b}-(\frac{a}{b})^3}=\frac{b^3}{a^3-ab^2-b^3}, \text{ and } y=\frac{-k}{1+k-k^3}=\frac{-\frac{a}{b}}{1+\frac{a}{b}-(\frac{a}{b})^3}=\frac{ab^2}{a^3-ab^2-b^3}.
$$
Let $p$ be any prime factor of $a^3-ab^2-b^3$. Because $x$ is integer, $p$ must also be a prime factor of $b^3$, hence $p$ is a prime factor of $b$, but then $p$ is also a prime factor of $(a^3-ab^2-b^3)+b(ab+b^2)=a^3$, which is a contradiction with the coprimality of $a$ and $b$. Hence, $a^3-ab^2-b^3$ has no prime factors, which implies $a^3-ab^2-b^3=\pm 1$. These are Thue equations and we can easily verify that their only integer solutions satisfying $b>0$ are
 $(\pm1,1),(0,1)$, and $(4,3)$. These solutions give $k=\pm 1,0$ or $\frac{4}{3}$. Then by \eqref{x3px2ymy3px2_rational} we can conclude that the integer solutions to \eqref{x3px2ymy3px2} are
$$
(x,y)=(-1,\pm 1),(-1,0),(0,0),(27,36).
$$

\vspace{10pt}

The final equation we will consider is
\begin{equation}\label{x3px2ypy3px2}
x^3+x^2y+y^3+x^2=0.
\end{equation}
Let us first describe all rational solutions to \eqref{x3px2ypy3px2}. The rational point $(x_0,y_0)=(0,0)$ is on the curve. After substituting $y=kx$, where $k$ is a rational parameter, into \eqref{x3px2ypy3px2}, we obtain 
$$
x^3+x^2(kx)+(kx)^3+x^2=0=x^2(1+x+kx+k^3x).
$$  
If $x=0$ we find the point $(x_0,y_0)=(0,0)$ that we started with. Otherwise, $1+x+kx+k^3x=0$ from which we find that
\begin{equation}\label{x3px2ypy3px2_rational}
x=\frac{-1}{1+k+k^3}, \quad \text{and} \quad y=kx=\frac{-k}{1+k+k^3}, \quad k \in \mathbb{Q}.
\end{equation}
We remark that $1+k+k^3\neq 0$ for rational $k$.

We will now use this description to find all integer solutions to \eqref{x3px2ypy3px2}. 
Let us represent $k=\frac{a}{b}$ as an irreducible fraction with $b>0$ and $\gcd(a,b)=1$. Then 
$$
x=\frac{-1}{1+k+k^3}=\frac{-1}{1+\frac{a}{b}+(\frac{a}{b})^3}=\frac{-b^3}{a^3+ab^2+b^3}, \text{ and } y=\frac{-k}{1+k+k^3}=\frac{-\frac{a}{b}}{1+\frac{a}{b}+(\frac{a}{b})^3}=\frac{-ab^2}{a^3+ab^2+b^3}.
$$
Let $p$ be any prime factor of $a^3+ab^2+b^3$. Because $x$ is integer, $p$ must also be a prime factor of $-b^3$, hence $p$ is a prime factor of $b$, but then $p$ is also a prime factor of $(a^3+ab^2+b^3)-b(ab+b^2)=a^3$, which is a contradiction with the coprimality of $a$ and $b$. Hence, $a^3+ab^2+b^3$ has no prime factors, which is only possible if $a^3+ab^2+b^3=\pm 1$. These are Thue equations and we can easily verify that their only integer solutions satisfying $b>0$ are 
 $(-2,3),(-1,1)$, and $(0,1)$. These solutions give $k=- 1,-\frac{2}{3}$ or $0$. Then by \eqref{x3px2ypy3px2_rational} we can conclude that the integer solutions to \eqref{x3px2ypy3px2} are
$$
(x,y)=(-27,18),(-1,0),(0,0),(1,-1).
$$
\vspace{10pt}

Table \ref{tab:H28genus0sol} summarises all rational and integer solutions to the equations listed in Table \ref{tab:H28genus0}.

\begin{center}
\begin{tabular}{ |c|c|c|c|c|c| } 
 \hline
Equation & Rational Solutions $(x,y)$ & Integer Solutions $(x,y)$ \\ 
 \hline\hline
 $2x^3+y^3+xy=0$&$\left(\frac{-k}{k^3+2},\frac{-k^2}{k^3+2}\right)$&$(0,0),(1,-1)$ \\ 
 \hline
$x^3+x^2 y-y^3+y^2=0$&$\left(\frac{-k^2}{1+k-k^3},\frac{-k^3}{1+k-k^3}\right)$, $(0,1)$&$(-1,\pm1),(0,0),(0,1),(48,64)$ \\ 
 \hline
$x^3+x^2 y+y^3+y^2=0$&$\left(\frac{-k^2}{1 + k + k^3}, \frac{-k^3}{1 + k + k^3}\right)$, $(0,-1)$ &$(-12,8),(0,-1),(0,0),(1,-1)$ \\ 
 \hline
 $x^3+x^2 y-y^3+xy=0$ &$\left(\frac{-k}{1 + k - k^3}, \frac{-k^2}{1 + k - k^3}\right)$&$(\pm 1,-1),(0,0),(36,48)$ \\ 
 \hline
 $x^3+x^2 y+y^3+xy=0$ &$\left(\frac{-k}{1 + k + k^3}, \frac{-k^2}{1 + k + k^3}\right)$&$(-1,1),(0,0),(18,-12)$ \\ 
 \hline
$x^3+x^2 y-y^3+x^2=0$ &$\left(\frac{-1}{1 + k - k^3}, \frac{-k}{1 + k - k^3}\right)$ &$(-1,\pm1),(-1,0),(0,0),(27,36)$\\ 
 \hline
$x^3+x^2 y+y^3+x^2=0$ &$\left(\frac{-1}{1 + k + k^3}, \frac{-k}{1 + k + k^3}\right)$&$(-27,18),(-1,0),(0,0),(1,-1)$ \\ 
 \hline
\end{tabular}
\captionof{table}{\label{tab:H28genus0sol}  Rational and integer solutions to the equations in Table \ref{tab:H28genus0}, assume that $k \in \mathbb{Q}$.}
\end{center} 

\subsection{Exercise 3.56}\label{ex:H29quarticdet}
\textbf{\emph{For each equation listed in Table \ref{tab:H29quarticdet}, use the ${\tt IntegralQuarticPoints}$ command in Magma to find all its integer solutions. }}

\begin{center}
\begin{tabular}{ |c|c|c|c|c|c| } 
 \hline
 $H$ & Equation & $H$ & Equation & $H$ & Equation \\ 
 \hline\hline
 $25$ & $1+x^4-2 y^2 = 0$ & $29$ & $1+x^3+x^2 y+y^2-x y^2 = 0$ & $29$ & $-1+x+x^4+y-2 y^2 = 0$ \\ 
 \hline
 $27$ & $-1+x^4+y-2 y^2 = 0$ & $29$ & $-1+x^3+x^2 y+y^2-x y^2 = 0$ & $29$ & $-1+2 x+x^4-2 y^2 = 0$ \\ 
 \hline
 $28$ & $2+x^4+y-2 y^2 = 0$ & $29$ & $-3+x^4+y-2 y^2 = 0$ & $29$ & $-1+x^4+x y-2 y^2 = 0$ \\ 
 \hline
 $28$ & $2+x+x^4-2 y^2 = 0$ & $29$ & $-1+x^4+2 y-2 y^2 = 0$ & $29$ & $1+x^4+x^2 y-y^2 = 0$ \\ 
 \hline
 $28$ & $x+x^4+y-2 y^2 = 0$ & $29$ & $1+x+x^4+y-2 y^2 = 0$ &  $29$ & $-1+x^4+x^2 y-y^2 = 0$ \\ 
 \hline
\end{tabular}
\captionof{table}{\label{tab:H29quarticdet} Equations of size $H\leq 29$ with quartic discriminant.}
\end{center} 

We will now look at solving equations that are quadratic in $y$ with quartic discriminant. In order to solve these equations, we can use the quadratic formula and find all integers $x$ such that the discriminant is a perfect square. Hence, 
we will be solving equations of the form 
$$
ax^4+bx^3+cx^2+dx+e=k^2
$$
where $a,b,c,d,e$ are integer coefficients and $x,k$ are integer variables. If we can find a rational solution $(x_0,k_0)$ to the equation, we can use the Magma command
\begin{equation}\label{intquarMagma}
{\tt IntegralQuarticPoints([a,b,c,d,e],[x_0,k_0]).}
\end{equation}
This will find all integer solutions, up to the sign of $k$.

The equation 
$$
2y^2=x^4+1
$$
is solved in Section 3.4.3 of the book, and its integer solutions are
$$
(x,y)=(\pm 1, \pm 1).
$$

The equation
$$
2y^2-y=x^4-1
$$
is equivalent to $2y^2+y=x^4-1$ which is solved in Section 3.4.3 of the book, and its integer solutions are
$$
(x,y)=(\pm 1,0), (\pm 2,3) \quad \text{and} \quad (\pm 8,-45).
$$

Let us now consider the first equation
\begin{equation}\label{2px4pym2y2}
2+x^4+y-2y^2=0.
\end{equation}
Solving this equation as a quadratic in $y$ results in
\begin{equation}\label{2px4pym2y2_y}
y_{1,2}=\frac{-1\pm \sqrt{8x^4+17}}{-4}.
\end{equation}
To ensure that $y$ is integer (or at least rational), 
the discriminant $8x^4+17$ must be a perfect square, hence we need to find all integer solutions to the equation
$$
8x^4+17=k^2.
$$
This equation has a rational solution $(x,k)=(1,-5)$.  We can use this along with the Magma command \eqref{intquarMagma} 
$$
{\tt IntegralQuarticPoints([8,0,0,0,17],[1,-5])}
$$
to find all integer solutions to the equation. This outputs that 
the discriminant is a perfect square for $x=\pm 1$. Substituting these values of $x$ into \eqref{2px4pym2y2_y} we can find all rational solutions for $y$, and then conclude that the integer solutions to \eqref{2px4pym2y2} are
$$
(x,y)=(\pm 1,-1).
$$

\vspace{10pt}

The next equation we will consider is
\begin{equation}\label{2pxpx4m2y2}
2+x+x^4-2y^2=0.
\end{equation}
After multiplying this equation by $2$ and making the substitution $Y=2y$, the equation is reduced to
$$
Y^2=2x^4+2x+4.
$$
We can see that $(x,Y)=(-1,-2)$ is a rational solution to this equation and the Magma command \eqref{intquarMagma} for this equation returns that 
$2x^4+2x+4$ is a perfect square for $x=-16,-1$ and $0$. After substituting these values of $x$ into \eqref{2pxpx4m2y2} and solving for $y$,
we can conclude that the integer solutions to \eqref{2pxpx4m2y2} are
$$
(x,y)=(-16, \pm 181), (-1, \pm 1),(0, \pm 1).
$$

\vspace{10pt}

Let us next consider the equation
\begin{equation}\label{xpx4pym2y2}
x+x^4+y-2y^2=0.
\end{equation}
Solving this equation as a quadratic in $y$ results in
\begin{equation}\label{xpx4pym2y2_y}
y_{1,2}=\frac{-1\pm \sqrt{8x^4+8x+1}}{-4}.
\end{equation}
To ensure that $y$ is integer (or at least rational), the discriminant $8x^4+8x+1$ must be a perfect square, hence we need to find all integer solutions to the equation
$$
8x^4+8x+1=k^2.
$$
We can see that this equation has a rational solution $(x,k)=(0,1)$ and the Magma command \eqref{intquarMagma} for this equation returns that 
$8x^4+8x+1$ is a perfect square for $x=-3,-1,0$ and $5$. After substituting these values of $x$ into \eqref{xpx4pym2y2_y} and selecting the solutions for which $y$ is an integer, 
we can conclude that the integer solutions to \eqref{xpx4pym2y2} are
$$
(x,y)=(-3,-6),(-1,0),(0,0),(5,18).
$$

\vspace{10pt}

Let us next consider the equation
\begin{equation}\label{1px3px2ypy2mxy2}
1+x^3+x^2y+y^2-xy^2=0.
\end{equation}
If $x=1$ then $y=-2$. 
Next, solving \eqref{1px3px2ypy2mxy2} as a quadratic in $y$ results in
\begin{equation}\label{1px3px2ypy2mxy2_y}
y_{1,2}=\frac{-x^2\pm \sqrt{5x^4-4x^3+4x-4}}{2-2x}.
\end{equation}
To ensure that $y$ is integer (or at least rational), the discriminant must be a perfect square, hence we need to find all integer solutions to the equation
$$
5x^4-4x^3+4x-4=k^2.
$$
 We can see that this equation has a rational solution $(x,k)=(1,1)$ and the Magma command \eqref{intquarMagma} for this equation returns that 
the discriminant is a perfect square for $x=-2$ and $\pm 1$. After substituting these values of $x$ into \eqref{1px3px2ypy2mxy2_y}, 
 we can conclude that the integer solutions to \eqref{1px3px2ypy2mxy2} are
$$
(x,y)=(-2, 1), (-1, 0), (1, -2).
$$

\vspace{10pt}

Let us next consider the equation
\begin{equation}\label{m1px3px2ypy2mxy2}
-1+x^3+x^2y+y^2-xy^2=0.
\end{equation}
If $x=1$ then $y=0$.
Next, solving this equation as a quadratic in $y$ results in
\begin{equation}\label{m1px3px2ypy2mxy2_y}
y_{1,2}=\frac{-x^2\pm \sqrt{5x^4-4x^3-4x+4}}{2-2x}.
\end{equation}
To ensure that $y$ is integer (or at least rational), we need to find all integer solutions to the equation
$$
5x^4-4x^3-4x+4=k^2.
$$
We can see that this equation has a rational solution $(x,k)=(1,1)$ and the Magma command \eqref{intquarMagma} for this equation returns that 
the discriminant is a perfect square for $x=0,1, \pm 3$ and $7$. After substituting these values of $x$ into \eqref{m1px3px2ypy2mxy2_y}, 
we can conclude that the integer solutions to \eqref{m1px3px2ypy2mxy2} are
$$
(x,y)=(-3, -4), (0, \pm 1), (1, 0), (3, -2).
$$

\vspace{10pt}

Let us next consider the equation
\begin{equation}\label{m3px4pym2y2}
-3+x^4+y-2y^2=0.
\end{equation}
Solving this equation as a quadratic in $y$ results in
\begin{equation}\label{m3px4pym2y2_y}
y_{1,2}=\frac{-1\pm \sqrt{8x^4-23}}{-4}.
\end{equation}
To ensure that $y$ is integer (or at least rational), we need to find all integer solutions to the equation
$$
8x^4-23=k^2.
$$
We can see that this equation has a rational solution $(x,k)=(4,-45)$ and the Magma command \eqref{intquarMagma} for this equation returns that
 the discriminant is a perfect square for $x=\pm 3$ and $\pm 4$. After substituting these values of $x$ into \eqref{m3px4pym2y2_y}, 
 we can conclude that the integer solutions to \eqref{m3px4pym2y2} are
$$
(x,y)= (\pm 3, -6), (\pm 4, -11).
$$

\vspace{10pt}

Let us next consider the equation
\begin{equation}\label{m1px4p2ym2y2}
-1+x^4+2y-2y^2=0.
\end{equation}
Solving this equation as a quadratic in $y$ results in
\begin{equation}\label{m1px4p2ym2y2_y}
y_{1,2}=\frac{-2\pm \sqrt{8x^4-4}}{-4}.
\end{equation}
To ensure that $y$ is integer (or at least rational), we need to find all integer solutions to the equation
$$
8x^4-4=k^2.
$$
 We can see that this equation has a rational solution $(x,k)=(-1,-2)$ and the Magma command \eqref{intquarMagma} for this equation returns that 
the discriminant is a perfect square for $x=\pm 1,\pm 13$. After substituting these values of $x$ into \eqref{m1px4p2ym2y2_y}, 
we can conclude that the integer solutions to \eqref{m1px4p2ym2y2} are
$$
(x,y)= (\pm 1, 0), (\pm 1, 1), (\pm 13, -119), (\pm13, 120).
$$

\vspace{10pt}

Let us next consider the equation
\begin{equation}\label{1pxpx4pym2y2}
1+x+x^4+y-2y^2=0.
\end{equation}
Solving this equation as a quadratic in $y$ results in
\begin{equation}\label{1pxpx4pym2y2_y}
y_{1,2}=\frac{-1\pm \sqrt{8x^4+8x+9}}{-4}.
\end{equation}
To ensure that $y$ is integer (or at least rational), we need to find all integer solutions to the equation
$$
8x^4+8x+9=k^2.
$$
We can see that this equation has a rational solution $(x,k)=(-1,3)$ and the Magma command \eqref{intquarMagma} for this equation returns that 
the discriminant is a perfect square for $x=-4,-2,0,\pm1,10$ and $274$. After substituting these values of $x$ into \eqref{1pxpx4pym2y2_y}, 
we can conclude that the integer solutions to \eqref{1pxpx4pym2y2} are
$$
(x,y)=(-4, -11), (-2, 3), (0, 1),  \pm(1, -1), (10, 71),(274, 53087).
$$

\vspace{10pt}

Let us next consider the equation
\begin{equation}\label{m1pxpx4pym2y2}
-1+x+x^4+y-2y^2=0.
\end{equation}
Solving this equation as a quadratic in $y$ results in
\begin{equation}\label{m1pxpx4pym2y2_y}
y_{1,2}=\frac{-1\pm \sqrt{8x^4+8x-7}}{-4}.
\end{equation}
To ensure that $y$ is integer (or at least rational), we need to find all integer solutions to the equation
$$
8x^4+8x-7=k^2.
$$
We can see that this equation has a rational solution $(x,k)=(1,3)$ and the Magma command \eqref{intquarMagma} for this equation returns that 
the discriminant is a only perfect square when $x=1$. After substituting this value of $x$ into \eqref{m1pxpx4pym2y2_y},  
we can conclude that the unique integer solution to \eqref{m1pxpx4pym2y2} is
$$
(x,y)=(1,1).
$$

\vspace{10pt}

The next equation we will consider is
\begin{equation}\label{m1p2xpx4m2y2}
-1+2x+x^4-2y^2=0.
\end{equation}
After multiplying this equation by $2$ and making the substitution $Y=2y$, the equation is reduced to
$$
Y^2=2x^4+4x-2.
$$
We can see that $(x,Y)=(1,-2)$ is a rational solution to this equation and the Magma command \eqref{intquarMagma} for this equation returns that 
the discriminant is only a perfect square when $x=1$. After substituting this value of $x$ into \eqref{m1p2xpx4m2y2} and solving for $y$, 
we can conclude that the integer solutions to \eqref{m1p2xpx4m2y2} are
$$
(x,y)= (1, \pm 1).
$$

\vspace{10pt}

Let us next consider the equation
\begin{equation}\label{m1px4pxym2y2}
-1+x^4+xy-2y^2=0.
\end{equation}
Solving this equation as a quadratic in $y$ results in
\begin{equation}\label{m1px4pxym2y2_y}
y_{1,2}=\frac{-x\pm \sqrt{8x^4+x^2-8}}{-4}.
\end{equation}
To ensure that $y$ is integer (or at least rational), we need to find all integer solutions to the equation 
$$
8x^4+x^2-8=k^2.
$$
We can see that this equation has a rational solution $(x,k)=(1,1)$ and the Magma command \eqref{intquarMagma} for this equation returns that 
the discriminant is a perfect square for $x=\pm 1$. After substituting these values of $x$ into \eqref{m1px4pxym2y2_y}, 
we can conclude that the integer solutions to \eqref{m1px4pxym2y2} are
$$
(x,y)=(\pm 1,0).
$$

\vspace{10pt}

Let us next consider the equation
\begin{equation}\label{1px4px2ymy2}
1+x^4+x^2y-y^2=0.
\end{equation}
Solving this equation as a quadratic in $y$ results in
\begin{equation}\label{1px4px2ymy2_y}
y_{1,2}=\frac{-x^2\pm \sqrt{5x^4+4}}{-2}.
\end{equation}
To ensure that $y$ is integer (or at least rational), we need to find all integer solutions to the equation 
$$
5x^4+4=k^2.
$$
We can see that this equation has a rational solution $(x,k)=(0,-2)$ and the Magma command \eqref{intquarMagma} for this equation returns that 
the discriminant is a perfect square for $x=0,\pm 1$ or $\pm 12$. After substituting these values of $x$ into \eqref{1px4px2ymy2_y},  
we can conclude that the integer solutions to \eqref{1px4px2ymy2} are
$$
(x,y)=(0,\pm 1), (\pm 1, -1), (\pm 1, 2), (\pm 12, -89), (\pm 12, 233).
$$

\vspace{10pt}

The final equation we will consider is
\begin{equation}\label{m1px4px2ymy2}
-1+x^4+x^2y-y^2=0.
\end{equation}
Solving this equation as a quadratic in $y$ results in
\begin{equation}\label{m1px4px2ymy2_y}
y_{1,2}=\frac{-x^2\pm \sqrt{5x^4-4}}{-2}.
\end{equation}
To ensure that $y$ is integer (or at least rational), we need to find all integer solutions to the equation
$$
5x^4-4=k^2.
$$
We can see that this equation has a rational solution $(x,k)=(1,1)$ and the Magma command \eqref{intquarMagma} for this equation returns that 
the discriminant is a perfect square for $x=\pm 1$. After substituting these values of $x$ into \eqref{m1px4px2ymy2_y}, 
we can conclude that the integer solutions to \eqref{m1px4px2ymy2} are
$$
(x,y)=(\pm 1,0),(\pm 1,1).
$$

\vspace{10pt}

Table \ref{tab:H29quarticdetsol} summarises the integer solutions to the equations listed in Table \ref{tab:H29quarticdet}.

\begin{center}
\begin{tabular}{ |c|c|c|c|c|c| } 
 \hline
Equation & Solution $(x,y)$ \\ 
 \hline\hline
$1+x^4-2 y^2 = 0$ & $(\pm1,\pm1)$   \\  \hline
$-1+x^4+y-2 y^2 = 0$ &  $ (\pm1, 0), (\pm 2, 3), (\pm 8, -45)$ \\  \hline
$2+x^4+y-2 y^2 = 0$ & $(\pm 1,-1)$ \\  \hline
$2+x+x^4-2 y^2 = 0$ & $(-16, \pm 181), (-1, \pm 1),(0, \pm 1)$  \\  \hline
 $x+x^4+y-2 y^2 = 0$ & $(-3,-6),(-1,0),(0,0),(5,18)$  \\  \hline
 $1+x^3+x^2 y+y^2-x y^2 = 0$ &$(-2, 1), (-1, 0), (1, -2)$\\  \hline
$-1+x^3+x^2 y+y^2-x y^2 = 0$&$(-3, -4), (0, \pm 1), (1, 0), (3, -2)$ \\  \hline
 $-3+x^4+y-2 y^2 = 0$ &$(\pm 3, -6), (\pm 4, -11)$ \\  \hline
 $-1+x^4+2 y-2 y^2 = 0$ &$(\pm 1, 0), (\pm 1, 1), (\pm 13, -119), (\pm13, 120)$ \\  \hline
$1+x+x^4+y-2 y^2 = 0$ &$(-4, -11), (-2, 3),  (0, 1), \pm (1, -1), (10, 71),(274, 53087)$ \\  \hline
 $-1+x+x^4+y-2 y^2 = 0$&$(1,1)$ \\  \hline
$-1+2 x+x^4-2 y^2 = 0$ &$(1,\pm1)$\\  \hline
 $-1+x^4+x y-2 y^2 = 0$&$(\pm 1,0)$ \\  \hline
 $1+x^4+x^2 y-y^2 = 0$&$(0,\pm 1), (\pm 1, -1), (\pm 1, 2), (\pm 12, -89), (\pm 12, 233)$  \\  \hline
  $-1+x^4+x^2 y-y^2 = 0$&$(\pm 1,0),(\pm 1,1)$  \\  \hline
\end{tabular}
\captionof{table}{\label{tab:H29quarticdetsol} Integer solutions to the equations listed in Table \ref{tab:H29quarticdet}.}
\end{center} 

\subsection{Exercise 3.58}\label{ex:H27rank0ell}
\textbf{\emph{Use the described method to find all integer solutions to the equations
\begin{equation}\label{x3px2ypy3myp1}
x^3+x^2 y+y^3-y+1=0
\end{equation}
and
\begin{equation}\label{x3px2ymy3pyp1}
x^3+x^2 y-y^3+y+1=0
\end{equation}
of size $H=27$. On the other hand, check that the method is not directly applicable to the equations listed in Table \ref{tab:H27open}.}}
\begin{center}
\begin{tabular}{ |c|c|c|c|c|c| } 
 \hline
 $H$ & Equation \\ 
 \hline\hline
 $27$ & $y^3 = 2x^3 + x + 1$ \\ 
 \hline
 $27$ & $x^3+x^2 y+y^3+y+1=0$ \\ 
 \hline
 $27$ & $x^3+x^2 y-y^3-x+1=0$ \\ 
 \hline
 $27$ & $x^3+x^2 y-y^3+x+1=0$ \\ 
 \hline
 $27$ & $x^3+x^2 y+y^3-x+1=0$ \\ 
 \hline
 $27$ & $x^3+x^2 y+y^3+x+1=0$ \\ 
 \hline
 $27$ & $y^3 = x^4+x+1$ \\ 
 \hline
 $27$ & $y^3 = x^4+x-1$ \\ 
 \hline
\end{tabular}
\captionof{table}{\label{tab:H27open} The remaining equations of size $H\leq 27$.}
\end{center} 

First, to solve equations \eqref{x3px2ypy3myp1} and \eqref{x3px2ymy3pyp1}, we will put them into Weierstrass form (equations of the form \eqref{eq:Weiform}), using the Maple Command 
\begin{equation}\label{maple:wei}
{\tt Weierstrassform(P(x,y), x, y, X, Y)}
\end{equation}
where $P(x,y)=0$ is the equation to be transformed. This will output 
$$
{\tt [P'(X,Y),P_1(x,y),P_2(x,y),P_3(X,Y),P_4(X,Y)]}
$$
Then, $P'(X,Y)=0$ is the equation in Weierstrass form, this equation can be found using the substitutions $x=P_3(X,Y)$ and $y=P_4(X,Y)$. Solutions can then be transformed into the original variables using $X=P_1(x,y)$ and $Y=P_2(x,y)$. However, $P'(X,Y)=0$ may not have integer coefficients, so we may need to multiply by a constant and make further substitutions as we did in Section \ref{ex:H20almWei}.

If an equation can be reduced to Weierstrass form we can check the rank of the equation using the Magma command 
$$
{\tt Rank(EllipticCurve([a,b,c,d,e])).}
$$
If the Rank is 0, then the equation has a finite number of rational solutions. If the reduced equation has a finite number of rational solutions, then the original equation has a finite number of rational solutions. All rational points on a curve with Rank 0 are the curves torsion points. We can find torsion points using the SageMath Command
$$
{\tt EllipticCurve([a,b,c,d,e]).torsion\_points();}
$$
where $a,b,c,d,e$ are the coefficients to the equation of the form \eqref{eq:Weiform}. This will output solutions 
$$
{\tt [(X: Y : 1), (0 : 1 : 0)]}.
$$
 The solution $(0 : 1 : 0)$ represents the point at infinity. Solutions of the form $(X: Y : 1)$ are the rational solutions to the equation, we can then use the change of variables found initially to obtain solutions $(x,y)$. 

Let us first solve the equation \eqref{x3px2ypy3myp1}.
The Maple command \eqref{maple:wei}
$$
{\tt Weierstrassform(x^3 + x^2*y + y^3 - y + 1, x, y, X, Y)}
$$
returns that using the rational change of variables 
\begin{equation}\label{x3px2ypy3myp1_XY}
	X=-\frac{-9y + x + 1}{3(x + 1)} \quad \text{and} \quad Y=-\frac{3(3x + 2y - 3)}{2(x + 1)}.
\end{equation}
we can reduce \eqref{x3px2ypy3myp1} to
$$
X^3 - \frac{X}{3} + \frac{721}{108} + Y^2=0.
$$
To make the coefficients integers, we can multiply this equation by $46656$ and make the substitutions $U=-36X$ and $V=216Y$, to obtain
$$
V^2=U^3-432U-311472. 
$$ 
We can then use the Magma command 
$$
{\tt Rank(EllipticCurve([0,0,0,-432,-311472]))}
$$
to check that the Rank is 0. We can now use the SageMath command
$$
{\tt EllipticCurve([0,0,0,-432,-311472]).torsion\_points();}
$$
to find all rational points on the curve, however, this outputs ${\tt [(0 : 1 : 0)]}$, hence the equation has no finite rational solutions. It is clear from \eqref{x3px2ypy3myp1_XY} that the ``point at infinity'' corresponds to $x=-1$, and then \eqref{x3px2ypy3myp1} implies that $y=0$. Therefore we can conclude that the unique rational solution to \eqref{x3px2ypy3myp1} is 
$$
	(x,y)=(-1,0)
$$
and so it is the only integer solution to \eqref{x3px2ypy3myp1}.

\vspace{10pt}

We will now consider the equation \eqref{x3px2ymy3pyp1}.
The Maple command \eqref{maple:wei}
$$
{\tt Weierstrassform(x^3 + x^2*y - y^3 + y + 1, x, y, X, Y)}
$$
returns that using the rational change of variables 
\begin{equation}\label{x3px2ymy3pyp1_XY}
	X=-\frac{5x^2 + 9xy - 6y^2 - 2x + 9y + 5}{3(x+1)^2} \quad \text{and} \quad Y=\frac{5x^3 - 6x^2y + 4xy^2 + 13x^2 - 4y^2 - 13x + 6y - 5}{2(x+1)^3}
\end{equation}
we can reduce \eqref{x3px2ymy3pyp1} to
$$
X^3 - \frac{X}{3} + \frac{521}{108} + Y^2=0.
$$
To make the coefficients integers, we can multiply the equation by $46656$ and make the substitutions $U=-36X$ and $V=216Y$, to obtain 
$$
V^2=U^3-432U-225072. 
$$ 
We can then use the Magma command 
$$
{\tt Rank(EllipticCurve([0,0,0,-432,-225072]))}
$$
to check that the Rank is 0. Now we can use the SageMath command
$$
{\tt EllipticCurve([0,0,0,-432,-225072]).torsion\_points();}
$$
to find all rational points on the curve, however, this outputs $[(0 : 1 : 0)]$, hence the equation has no finite rational solutions. By \eqref{x3px2ymy3pyp1_XY}, it is clear that the ``point at infinity'' corresponds to $x=-1$, and then \eqref{x3px2ymy3pyp1} implies that $y=0$. Therefore we can conclude that the unique rational solution to \eqref{x3px2ymy3pyp1} is 
$$
	(x,y)=(-1,0)
$$
and so it is the only integer solution to \eqref{x3px2ymy3pyp1}.

The equations in Table \ref{tab:H27open} cannot be solved using the same method as equations \eqref{x3px2ypy3myp1} and \eqref{x3px2ymy3pyp1} as they either cannot be reduced to Weierstrass form or the equations they reduce to do not have Rank 0. Table \ref{tab:H27openrank} provides, for each equation, the equation it reduces to (where possible) and its rank.

\begin{center}
\begin{tabular}{ |c|c|c|c|c|c| } 
 \hline
 Equation & Reduced Equation & Rank \\ 
 \hline\hline
 $y^3 = 2x^3 + x + 1$&$Y^2=X^3-29$&$1$ \\ 
 \hline
$x^3+x^2 y+y^3+y+1=0$&$Y^2=X^3-432X-411696$&$1$ \\ 
 \hline
$x^3+x^2 y-y^3-x+1=0$&$Y^2=X^3-16X-304$&$1$ \\ 
 \hline
$x^3+x^2 y-y^3+x+1=0$&$Y^2=X^3-16X-432$&$1$ \\ 
 \hline
 $x^3+x^2 y+y^3-x+1=0$&$Y^2=X^3+16X-432$&$1$ \\ 
 \hline
 $x^3+x^2 y+y^3+x+1=0$&$Y^2=X^3+16X-560$&$1$ \\ 
 \hline
  $y^3 = x^4+x+1$&Cannot be reduced, genus $=3$&- \\ 
 \hline
 $y^3 = x^4+x-1$&Cannot be reduced, genus $=3$& -\\ 
 \hline
\end{tabular}
\captionof{table}{\label{tab:H27openrank} Necessary information to show that the equations in Table \ref{tab:H27open} cannot be solved using the same method as the equations \eqref{x3px2ypy3myp1} and \eqref{x3px2ymy3pyp1}.}
\end{center}

\subsection{Exercise 3.69}\label{ex:H31Runge}
\textbf{\emph{Solve all the equations from Table \ref{tab:H31Runge}. }}
\begin{center}
\begin{tabular}{ |c|c|c|c|c|c| } 
 \hline
 $H$ & Equation & $H$ & Equation & $H$ & Equation \\ 
 \hline\hline
 $28$ & $x^3 y+y^3+x^2=0$ & $30$ & $x^3 y+y^3+x^2-x=0$ & $31$ & $x^3 y+y^3+x^2+y-1=0$ \\ 
 \hline
 $29$ & $x^3 y+y^3+x^2-1=0$ & $30$ & $x^3 y+y^3+x^2+x=0$ & $31$ & $x^3 y+y^3+x^2-x-1=0$ \\ 
 \hline
 $29$ & $x^3 y+y^3+x^2+1=0$ & $31$ & $x^3 y+y^3+x^2-3=0$ & $31$ & $x^3 y+y^3+x^2-x+1=0$ \\ 
 \hline
 $30$ & $x^3 y+y^3+x^2-y=0$ & $31$ & $x^3 y+y^3+x^2+3=0$ & $31$ & $x^3 y+y^3+x^2+x-1=0$ \\ 
 \hline
 $30$ & $x^3 y+y^3+x^2+y=0$ & $31$ & $x^3 y+y^3+x^2-y-1=0$ & $31$ & $x^3 y+y^3+x^2+x+1=0$ \\ 
 \hline
\end{tabular}
\captionof{table}{\label{tab:H31Runge} Equations satisfying strong Runge's condition of size $H\leq 31$.}
\end{center} 

To solve the equations in Table \ref{tab:H31Runge} we will use the method from Section 3.4.7 of the book, which we summarise below for convenience. An equation satisfies Runge's strong condition if the monomials of degree $d$ (referred to as the leading part) of the equation of degree $d$ can be written as a product of two non-constant coprime polynomials. All equations in Table \ref{tab:H31Runge} have degree $4$, and they can be rewritten as
$$
y(x^3+y^2+a)=Q(x),
$$
where $Q(x)$ is some quadratic polynomial, and $a$ is a constant. We can consider the case $y=0$ separately, so we may assume $y \neq 0$. Then $Q(x)$ is divisible by $y$, and we may write $Q(x)=yt$ for some integer $t$. Substituting this into the original equation and cancelling $y$, we obtain an equation of the form
$$
x^3+y^2+a=t.
$$
Because $y$ must divide $Q(x)$ and $x^3+a-t$, $y$ must also divide any expression of the form $P_1 \cdot Q(x)-P_2 \cdot (x^3+a-t)$ and we can find $P_1$ and $P_2$ such that this expression is linear in terms of $x$ and $t$. Let $k$ be this linear expression divided by $y$. Then, when $y$ is large, $|k|$ is bounded, and we have a finite number of cases to consider. Let us illustrate this method on a number of examples.

Equation 
$$
x^3 y+y^3+x^2=0
$$
is solved in Section 3.4.7 of the book, and its integer solutions are
$$
(x,y)=(0,0).
$$

Equation 
$$
x^3 y+y^3+x^2+1=0
$$
is solved in Section 3.4.7 of the book, and its integer solutions are
$$
(x,y)=(-3,5),(0,-1),(1,-1).
$$

Let us first consider the equation
\begin{equation}\label{x3ypy3px2m1}
x^3y+y^3+x^2-1=0.
\end{equation}
If $y=0$ then $x = \pm 1$. Now assume that $y\neq 0$. Then $x^2-1$ is divisible by $y$. So, let $x^2-1=yt$ for some integer $t$. Substituting this into \eqref{x3ypy3px2m1} and cancelling $y$ we obtain $x^3+y^2+t=0$. We can now see that $x^3+t$ is divisible by $y$, but then $x^3+t-x(x^2-1)=t+x$ is divisible by $y$. Let $k$ be an integer such that $k=\frac{x+t}{y}=\frac{x+((x^2-1)/y)}{y}=\frac{xy+x^2-1}{y^2}$. Let $r=\frac{|x|}{|y|}$, and assume that $|y|>5$. Then 
$$
r^3=\frac{|x|^3}{|y|^3}=\frac{|x^3y|}{|y|^4}=\frac{|y^3+x^2-1|}{|y|^4}<\frac{1}{5}+\frac{r^2}{25}+\frac{1}{5^4}=0.04r^2+0.2016.
$$
Then $f(r)<0.2016$, where $f(r)=r^2(r-0.04)$. The function $f(r)$ is increasing when $r\geq 0.6$, and $f(0.6)=(0.6)^2(0.6-0.04)=0.2016$, so $f(r)<0.2016$ implies that $r<0.6$. We then have
$$
|k|=\left|\frac{x}{y} +\frac{x^2}{y^2}-\frac{1}{y^2}\right| \leq r^2+r+\frac{1}{y^2} < (0.6)^2+0.6+0.04=1.
$$
Hence $k=0$. But then $x^2+xy-1=0$ or $x(x+y)=1$ which implies that $x= \pm 1$. Substituting this into \eqref{x3ypy3px2m1} gives $\pm y+y^3+2=0$ which is impossible for $|y| > 5$. Finally, we must consider the case $|y| \leq 5$. We can then conclude that all integer solutions to equation \eqref{x3ypy3px2m1} are
$$
(x,y)=(-2,-3),(-1,\pm 1),(\pm 1,0),(0,1).
$$

\vspace{10pt}

The next equation we will consider is
\begin{equation}\label{x3ypy3px2my}
x^3y+y^3+x^2-y=0.
\end{equation}
If $y=0$ then $x = 0$. Now assume that $y\neq 0$. Then $x^2$ is divisible by $y$. So, let $x^2=yt$ for some integer $t$. Substituting this into \eqref{x3ypy3px2my} and cancelling $y$ we obtain $x^3+y^2+t-1=0$. We can now see that $x^3+t-1$ is divisible by $y$, but then $x^3+t-1-x(x^2)=t-1$ is divisible by $y$. Let $k$ be an integer such that $k=\frac{t-1}{y}=\frac{(x^2/y)-1}{y}=\frac{x^2-y}{y^2}$. Let $r=\frac{|x|}{|y|}$, and assume that $|y|>5$. Then 
$$
r^3=\frac{|x|^3}{|y|^3}=\frac{|x^3y|}{|y|^4}=\frac{|y^3+x^2-y|}{|y|^4}<\frac{1}{5}+\frac{r^2}{25}+\frac{1}{5^3}=0.04r^2+0.208.
$$
Then $f(r)<0.208$, where $f(r)=r^2(r-0.04)$. The function $f(r)$ is increasing when $r\geq 0.606$, and  $f(0.606)=(0.606)^2(0.606-0.04)=0.208$, so $f(r)<0.208$ implies that $r<0.606$. We then have
$$
|k|=\left|\frac{x^2}{y^2}-\frac{1}{y}\right| \leq r^2+\frac{1}{y} < (0.606)^2+0.2=0.567.
$$
Hence $k=0$. But then $x^2-y=0$, hence $(x,y)=(u,u^2)$, where $u$ is an arbitrary integer. After substituting this into \eqref{x3ypy3px2my} we have $u^5+u^6=0$, so $u=-1$ or $0$, which is impossible for $|y| > 5$.

Finally, we must consider the case $|y| \leq 5$. We can then conclude that all integer solutions to equation \eqref{x3ypy3px2my} are
$$
(x,y)=(0,\pm 1),(0,0),\pm(1,-1).
$$

\vspace{10pt}

The next equation we will consider is
\begin{equation}\label{x3ypy3px2py}
x^3y+y^3+x^2+y=0.
\end{equation}
If $y=0$ then $x = 0$. Now assume that $y\neq 0$. Then $x^2$ is divisible by $y$. So, let $x^2=yt$ for some integer $t$. Substituting this into \eqref{x3ypy3px2py} and cancelling $y$, we obtain $x^3+y^2+t+1=0$. We can now see that $x^3+t+1$ is divisible by $y$, but then $x^3+t+1-x(x^2)=t+1$ is divisible by $y$. Let $k$ be an integer such that $k=\frac{t+1}{y}=\frac{(x^2/y)+1}{y}=\frac{x^2+y}{y^2}$. Let $r=\frac{|x|}{|y|}$, and assume that $|y|>5$. Then 
$$
r^3=\frac{|x|^3}{|y|^3}=\frac{|x^3y|}{|y|^4}=\frac{|y^3+x^2+y|}{|y|^4}<\frac{1}{5}+\frac{r^2}{25}+\frac{1}{5^3}=0.04r^2+0.208.
$$
Then $f(r)<0.208$, where $f(r)=r^2(r-0.04)$. The function $f(r)$ is increasing when $r\geq 0.606$, and  $f(0.606)=(0.606)^2(0.606-0.04)=0.208$, so $f(r)<0.208$ implies that $r<0.606$. We then have
$$
|k|=\left|\frac{x^2}{y^2}+\frac{1}{y}\right| \leq r^2+\frac{1}{y} < (0.606)^2+0.2=0.567.
$$
Hence $k=0$. But then $x^2+y=0$, hence $(x,y)=(u,-u^2)$, where $u$ is an arbitrary integer. After substituting this into \eqref{x3ypy3px2py} we have $u^5+u^6=0$, so $u=-1$ or $0$, which is impossible for $|y| > 5$. Finally, we must consider the case $|y| \leq 5$. We can then conclude that all integer solutions to equation \eqref{x3ypy3px2py} are
$$
(x,y)=(-1, -1),(0,0).
$$

\vspace{10pt}

The next equation we will consider is
\begin{equation}\label{x3ypy3px2mx}
x^3 y+y^3+x^2-x=0.
\end{equation}
If $y=0$ then $x=0$ or $x=1$. Now assume that $y \neq 0$. Then $x^2-x$ is divisible by $y$. So, let $x^2-x=yt$ for some integer $t$. Substituting this into \eqref{x3ypy3px2mx} and cancelling $y$, we obtain $x^3+y^2+t=0$. We can now see that $x^3+t$ is divisible by $y$, but then $x^3+t-(x+1)(x^2-x)=t+x$ is divisible by $y$ as well. Let $k$ be an integer such that $k=\frac{x+t}{y}=\frac{x+(x^2-x)/y}{y}=\frac{xy+x^2-x}{y^2}$. Let $r=\frac{|x|}{|y|}$, and assume that $|y|>6$. Then
$$
r^3 = \frac{|x|^3}{|y|^3}=\frac{|x^3y|}{|y|^4}=\left|\frac{1}{y}+\frac{x^2}{y^4}-\frac{x}{y^4}\right|< \frac{1}{6}+\frac{r^2}{36}+\frac{r}{216}.
$$
Then $f(r)<\frac{1}{6}$, where $f(r)=r^3-\frac{r^2}{36}-\frac{r}{216}$. Because $f(r)$ is an increasing function for $r\geq 0.57$, and $f(0.57)=0.173...>\frac{1}{6}$, inequality $f(r)< \frac{1}{6}$ implies that $r<0.57$. We then have 
$$
|k|=\left|\frac{x^2}{y^2}+\frac{x}{y}-\frac{x}{y^2}\right| \leq r^2+r+\frac{r}{y} < (0.57)^2+0.57 + \frac{0.57}{6} < 0.99.
$$ 
Hence $k=0$. But then $xy+x^2-x=0$ or $x(y+x-1)=0$, hence $(x,y)=(u,1-u)$ or $(x,y)=(0,u)$ where $u$ is integer. Substituting the first solution into \eqref{x3ypy3px2mx} we obtain $u=1$, which is impossible for $|y|>6$. Substituting the second solution into \eqref{x3ypy3px2mx} we have $u=0=y$, which is impossible for $|y|>6$. Finally, we must consider the case $|y|\leq 6$. We can then conclude that all integer solutions to equation \eqref{x3ypy3px2mx} are
$$
(x,y)=(0,0),(1,0).
$$

\vspace{10pt}

The next equation we will consider is
\begin{equation}\label{x3ypy3px2px}
x^3 y+y^3+x^2+x=0.
\end{equation}
If $y=0$ then $x=-1$ or $x=0$. Now assume that $y \neq 0$. Then $x^2+x$ is divisible by $y$. So, let $x^2+x=yt$ for some integer $t$. Substituting this into \eqref{x3ypy3px2px} and cancelling $y$, we obtain $x^3+y^2+t=0$. We can now see that $x^3+t$ is divisible by $y$, but then $x^3+t-(x-1)(x^2+x)=t+x$ is divisible by $y$ as well. Let $k$ be an integer such that $k=\frac{x+t}{y}=\frac{x+(x^2+x)/y}{y}=\frac{xy+x^2+x}{y^2}$.
Let $r=\frac{|x|}{|y|}$, and assume that $|y|>6$. Then
$$
r^3 = \frac{|x|^3}{|y|^3}=\frac{|x^3y|}{|y|^4}=\left|\frac{1}{y}+\frac{x^2}{y^4}+\frac{x}{y^4}\right|< \frac{1}{6}+\frac{r^2}{36}+\frac{r}{216}.
$$
Then $f(r)<\frac{1}{6}$, where $f(r)=r^3-\frac{r^2}{36}-\frac{r}{216}$. Because $f(r)$ is an increasing function for $r\geq 0.57$, and $f(0.57)=0.173...>\frac{1}{6}$, inequality $f(r)< \frac{1}{6}$ implies that $r<0.57$. We then have 
$$
|k|=\left|\frac{x^2}{y^2}+\frac{x}{y}+\frac{x}{y^2}\right| \leq r^2+r+\frac{r}{y} < (0.57)^2+0.57 + \frac{0.57}{6} < 0.99.
$$ 
Hence $k=0$. But then $xy+x^2+x=0$ or $x(y+x+1)=0$ hence $(x,y)=(u,-1-u)$ or $(x,y)=(0,u)$ where $u$ is integer. Substituting the first solution into \eqref{x3ypy3px2px} we obtain $u=-1$, which is impossible for $|y|>6$. Substituting the second solution into \eqref{x3ypy3px2px} we have $u=0=y$, which is impossible for $|y|>6$. Finally, we must consider the case $|y|\leq 6$. We can then conclude that all integer solutions to equation \eqref{x3ypy3px2px} are
$$
(x,y)=(-1,0),(-1,1),(0,0),(\pm 1,-1).
$$

\vspace{10pt}

The next equation we will consider is
\begin{equation}\label{x3ypy3px2m3}
	x^3 y+y^3+x^2-3=0.
\end{equation}
If $(x,y)$ is any integer solution, then $y\neq 0$, and $x^2-3$ is divisible by $y$. So, we can write $x^2-3=yt$ for some integer $t$. Substituting this into \eqref{x3ypy3px2m3} and cancelling $y$, we obtain $x^3+y^2+t=0$. Hence $x^3+t$ is divisible by $y$. But then $x^3+t-x(x^2-3)=t+3x$ is divisible by $y$ as well. Let $k$ be an integer such that $k=\frac{t+3x}{y}=\frac{3x+(x^2-3)/y}{y}=\frac{3xy+x^2-3}{y^2}$. Let $r=\frac{|x|}{|y|}$, and assume that $|y|>40$. Then 
$$
r^3=\frac{|x|^3}{|y|^3}=\frac{|x^3y|}{|y|^4}=\frac{|y^3+x^2-3|}{|y|^4}<\frac{1}{40}+\frac{r^2}{40^2}+\frac{3}{40^4}
$$ 
Then $f(r)<\frac{1}{40}+\frac{3}{40^4}<0.0251$ where $f(r)=r^2(r-0.000625)$. The function $f(r)$ is increasing when $r>0.293$, and $f(0.293)=0.0251$. So $f(r)<0.0251$ implies that $r<0.293$. We then have 
$$
|k|=\left| \frac{3x}{y}+\frac{x^2}{y^2}-\frac{3}{y^2}\right| \leq r^2+3r+\frac{3}{y^2}<0.966724.
$$ 
Hence $k=0$. Then $3xy+x^2-3=0$ which has no integer solutions. Finally, we must consider the case $|y| \leq 40$, which can be easily done with a computer search. We can then conclude that the unique integer solution to equation \eqref{x3ypy3px2m3} is
$$
(x,y)=(1,1).
$$

\vspace{10pt}

The next equation we will consider is
\begin{equation}\label{x3ypy3px2p3}
	x^3 y+y^3+x^2+3=0.
\end{equation}
If $(x,y)$ is any integer solution, then $y\neq 0$, and $x^2+3$ is divisible by $y$. So, we can write $x^2+3=yt$ for some integer $t$. Substituting this into \eqref{x3ypy3px2p3} and cancelling $y$ we obtain $x^3+y^2+t=0$. We can now see that $x^3+t$ is divisible by $y$. But then $x^3+t-x(x^2+3)=t-3x$ is divisible by $y$ as well. Let $k$ be an integer such that $k=\frac{t-3x}{y}=\frac{-3x+(x^2+3)/y}{y}=\frac{x^2+3-3xy}{y^2}$. Let $r=\frac{|x|}{|y|}$, and assume that $|y|>40$. Then 
$$
r^3=\frac{|x|^3}{|y|^3}=\frac{|x^3y|}{|y|^4}=\frac{|y^3+x^2+3|}{|y|^4}<\frac{1}{40}+\frac{r^2}{40^2}+\frac{3}{40^4}.
$$ 
Then $f(r)<\frac{1}{40}+\frac{3}{40^4}<0.0251$ where $f(r)=r^2(r-0.000625)$. The function $f(r)$ is increasing when $r>0.293$, and $f(0.293)=0.0251$. So $f(r)<0.0251$ implies that $r<0.293$. We then have 
$$
|k|=\left|\frac{x^2}{y^2}+\frac{3}{y^2}-\frac{3x}{y}\right| \leq r^2+3r+\frac{3}{y^2}<0.966724.
$$ 
Hence $k=0$. Then $x^2+3-3xy=0$ which has no integer solutions. Finally, we must consider the case $|y| \leq 40$, which can be easily done with a computer search. We can then conclude that the unique integer solution to equation \eqref{x3ypy3px2p3} is
$$
(x,y)=(-2,1).
$$

\vspace{10pt}

The next equation we will consider is
\begin{equation}\label{x3ypy3px2mym1}
	x^3 y+y^3+x^2-y-1=0.
\end{equation}
If $y=0$ then $x = \pm 1$. Now assume that $y\neq 0$. Then $x^2-1$ is divisible by $y$. So, let $x^2-1=yt$ for some integer $t$. Substituting this into \eqref{x3ypy3px2mym1} and cancelling $y$, we obtain $x^3+y^2+t-1=0$. We can now see that $x^3+t-1$ is divisible by $y$, but then $x^3+t-1-x(x^2-1)=t-1+x$ is divisible by $y$. Let $k$ be an integer such that $k=\frac{t-1+x}{y}=\frac{((x^2-1)/y)-1+x}{y}=\frac{x^2-1-y+xy}{y^2}$. Let $r=\frac{|x|}{|y|}$, and assume that $|y|>8$. Then 
$$
r^3=\frac{|x|^3}{|y|^3}=\frac{|x^3y|}{|y|^4}=\frac{|y^3+x^2-y-1|}{|y|^4}<\frac{1}{8}+\frac{r^2}{64}+\frac{1}{8^3}+\frac{1}{8^4}=\frac{r^2}{64}+\frac{521}{4096}.
$$
Then $f(r)<\frac{521}{4096}$, where $f(r)=r^2(r-0.015625)$. The function $f(r)$ is increasing when $r\geq 0.51$, and  $f(0.51)=(0.51)^2(0.51-0.015625)=0.1285...>\frac{521}{4096}$, so $f(r)<\frac{521}{4096}$ implies that $r<0.51$. We then have
$$
|k|=\left|\frac{x^2}{y^2}-\frac{1}{y^2}-\frac{1}{y}+\frac{x}{y}\right| \leq r^2+\frac{1}{y^2}+\frac{1}{y}+r <0.910725.
$$ 
Hence $k=0$. But then $x^2-1-y+xy=0$, or $(x-1)(x+1+y)=0$ hence $(x,y)=(u,-1-u)$, or $(x,y)=(1,u)$ where $u$ is an arbitrary integer. Substituting the first solution into \eqref{x3ypy3px2mym1} we obtain $u=-1$, which is impossible for $|y| > 8$. Substituting the second solution into \eqref{x3ypy3px2mym1} we have $u=0=y$, which is impossible for $|y| > 8$. Finally, we must consider the case $|y| \leq 8$. We can then conclude that all integer solutions to equation \eqref{x3ypy3px2mym1} are
$$
(x,y)=(\pm 1,0).
$$

\vspace{10pt}

The next equation we will consider is
\begin{equation}\label{x3ypy3px2pym1}
x^3 y+y^3+x^2+y-1=0.
\end{equation}
If $y=0$ then $x = \pm 1$. Now assume that $y\neq 0$. Then $x^2-1$ is divisible by $y$. So, let $x^2-1=yt$ for some integer $t$. Substituting this into \eqref{x3ypy3px2pym1} and cancelling $y$, we obtain $x^3+y^2+t+1=0$. We can now see that $x^3+t+1$ is divisible by $y$, but then $x^3+t+1-x(x^2-1)=t+1+x$ is divisible by $y$. Let $k$ be an integer such that $k=\frac{t+1+x}{y}=\frac{((x^2-1)/y)+1+x}{y}=\frac{x^2-1+y+xy}{y^2}$. Let $r=\frac{|x|}{|y|}$, and assume that $|y|>8$. Then 
$$
r^3=\frac{|x|^3}{|y|^3}=\frac{|x^3y|}{|y|^4}=\frac{|y^3+x^2+y-1|}{|y|^4}<\frac{1}{8}+\frac{r^2}{64}+\frac{1}{8^3}+\frac{1}{8^4}=\frac{r^2}{64}+\frac{521}{4096}.
$$
Then $f(r)<\frac{521}{4096}$, where $f(r)=r^2(r-0.015625)$. The function $f(r)$ is increasing when $r\geq 0.51$, and  $f(0.51)=(0.51)^2(0.51-0.015625)=0.1285...>\frac{521}{4096}$, so $f(r)<\frac{521}{4096}$ implies that $r<0.51$. We then have 
$$
|k|=\left|\frac{x^2}{y^2}-\frac{1}{y^2}+\frac{1}{y}+\frac{x}{y}\right| \leq r^2+\frac{1}{y^2}+\frac{1}{|y|}+r < (0.51)^2+\frac{1}{8^2}+\frac{1}{8}+0.51 < 1. 
$$ 
Hence $k=0$. But then $x^2-1+y+xy=0$, or $(x+1)(x-1+y)=0$, hence $(x,y)=(u,1-u)$ or $(-1,u)$, where $u$ is an arbitrary integer. Substituting the first solution into \eqref{x3ypy3px2pym1} we obtain $u=1$, which is impossible for $|y| > 8$. Substituting the second solution into \eqref{x3ypy3px2pym1} we have $u=0=y$, which is impossible for $|y| > 8$. Finally, we must consider the case $|y| \leq 8$. We can then conclude that all integer solutions to equation \eqref{x3ypy3px2pym1} are
$$
(x,y)=(\pm 1,0).
$$

\vspace{10pt}

The next equation we will consider is
\begin{equation}\label{x3ypy3px2mxm1}
	x^3 y+y^3+x^2-x-1=0.
\end{equation}
If $(x,y)$ is any integer solution, then $y \neq 0$, and $x^2-x-1$ is divisible by $y$. So, we can write $x^2-x-1=yt$ for some integer $t$. Substituting this into \eqref{x3ypy3px2mxm1} and cancelling $y$, we obtain $x^3+y^2+t=0$. Hence $x^3+t$ is divisible by $y$. But then $(x^3+t)-(x+1)(x^2-x-1)=t+2x+1$ is divisible by $y$ as well. Let $k$ be an integer such that $k=\frac{t+2x+1}{y}=\frac{((x^2-x-1)/y)+2x+1}{y}=\frac{x^2-x-1+2xy+y}{y^2}$. Let $r=\frac{|x|}{|y|}$, and assume that $|y|>18$. Then
$$
r^3=\frac{|x|^3}{|y|^3}=\frac{|x^3y|}{|y|^4}=\frac{|y^3+x^2-x-1|}{|y|^4}<\frac{1}{18}+\frac{r^2}{324}+\frac{r}{18^3}+\frac{1}{18^4}.
$$ 
Then $f(r)<\frac{1}{18^4}+\frac{1}{18}$ where $f(r)=r(r^2-\frac{r}{324}-\frac{1}{18^3})$. The function $f(r)$ is increasing when $r \geq 0.39$, and $f(0.39)=0.0587...>\frac{1}{18^4}+\frac{1}{18}$. So $f(r)<\frac{1}{18^4}+\frac{1}{18}$ implies that $r<0.39$. We then have 
$$
|k|=\left|\frac{x^2}{y^2}-\frac{x}{y^2}-\frac{1}{y^2}+\frac{2x}{y}+\frac{1}{y}\right|\leq r^2+\frac{r}{y}+\frac{1}{y^2}+2r+\frac{1}{y} < 0.9898.
$$ 
Hence $k=0$. Then, $x^2-x-1+2xy+y=0$, whose only integer solutions have $y=1$, which is impossible with $|y|>18$. Finally, we must consider the case $|y| \leq 18$. We can then conclude that the unique integer solution to equation \eqref{x3ypy3px2mxm1} is
$$
(x,y)=(0,1).
$$

\vspace{10pt}

The next equation we will consider is
\begin{equation}\label{x3ypy3px2mxp1}
	x^3 y+y^3+x^2-x+1=0.
\end{equation}
If $(x,y)$ is any integer solution, then $y \neq 0$, and $x^2-x+1$ is divisible by $y$. So, we can write $x^2-x+1=yt$ for some integer $t$. Substituting this into \eqref{x3ypy3px2mxp1} and cancelling $y$, we obtain $x^3+y^2+t=0$. Hence $x^3+t$ is divisible by $y$. But then $(x^3+t)-(x+1)(x^2-x+1)=t-1$ is divisible by $y$ as well. Let $k$ be an integer such that $k=\frac{t-1}{y}=\frac{((x^2-x+1)/y)-1}{y}=\frac{x^2-x+1-y}{y^2}$. Let $r=\frac{|x|}{|y|}$, and assume that $|y|>4$. Then
$$
r^3=\frac{|x|^3}{|y|^3}=\frac{|x^3y|}{|y|^4}=\frac{|y^3+x^2-x+1|}{|y|^4}<\frac{1}{4}+\frac{r^2}{16}+\frac{r}{64}+\frac{1}{256}.
$$ 
Then $f(r)<\frac{65}{256}$ where $f(r)=r(r^2-\frac{r}{16}-\frac{1}{64})$. The function $f(r)$ is increasing when $r \geq 0.67$, and $f(0.67)=0.26622...>\frac{65}{256}$. So $f(r)<\frac{65}{256}$ implies that $r<0.67$. We then have 
$$
|k|=\left|\frac{x^2}{y^2}-\frac{x}{y^2}+\frac{1}{y^2}-\frac{1}{y}\right|\leq r^2+\frac{r}{y}+\frac{1}{y^2}+\frac{1}{y} < 0.918.
$$ 
Hence $k=0$. Then, $x^2+x-1-y=0$ whose integer solutions are $(x,y)=(u,u^2-u+1)$ for some integer $u$. Substituting this into \eqref{x3ypy3px2mxp1} gives $(u^2+2)(1-u+u^2)^2=0$ which has no integer solutions. Finally, we must consider the case $|y| \leq 4$. We can then conclude that all integer solutions to equation \eqref{x3ypy3px2mxp1} are
$$
(x,y)=(-2,1),(0,-1).
$$

\vspace{10pt}

The next equation we will consider is
\begin{equation}\label{x3ypy3px2pxm1}
x^3 y+y^3+x^2+x-1=0.
\end{equation}
If $(x,y)$ is any integer solution, then $y \neq 0$, and $x^2+x-1$ is divisible by $y$. So, we can write $x^2+x-1=yt$ for some integer $t$. Substituting this into \eqref{x3ypy3px2pxm1} and cancelling $y$, we obtain $x^3+y^2+t=0$. Hence $x^3+t$ is divisible by $y$. But then $(x^3+t)-(x-1)(x^2+x-1)=t+2x-1$ is divisible by $y$ as well. Let $k$ be an integer such that $k=\frac{t+2x-1}{y}=\frac{((x^2+x-1)/y)+2x-1}{y}=\frac{x^2+x-1+2xy-y}{y^2}$. Let $r=\frac{|x|}{|y|}$, and assume that $|y|>18$. Then
$$
r^3=\frac{|x|^3}{|y|^3}=\frac{|x^3y|}{|y|^4}=\frac{|y^3+x^2+x-1|}{|y|^4}<\frac{1}{18}+\frac{r^2}{324}+\frac{r}{18^3}+\frac{1}{18^4}.
$$ 
Then $f(r)<\frac{1}{18^4}+\frac{1}{18}$ where $f(r)=r(r^2-\frac{r}{324}-\frac{1}{18^3})$. The function $f(r)$ is increasing when $r \geq 0.39$, and $f(0.39)=0.0587...>\frac{1}{18^4}+\frac{1}{18}$. So $f(r)<\frac{1}{18^4}+\frac{1}{18}$ implies that $r<0.39$. We then have 
$$
|k|=\left|\frac{x^2}{y^2}+\frac{x}{y^2}-\frac{1}{y^2}+\frac{2x}{y}-\frac{1}{y}\right|\leq r^2+\frac{r}{y}+\frac{1}{y^2}+2r+\frac{1}{y} < 0.9898.
$$
Hence $k=0$. Then, $x^2+x-1+2xy-y=0$ whose only integer solutions have $y=-1$, which is impossible with $|y|>18$. Finally, we must consider the case $|y| \leq 18$. We can then conclude that the unique integer solution to equation \eqref{x3ypy3px2pxm1} is
$$
(x,y)=(0,1).
$$

\vspace{10pt}

The final equation we will consider is
\begin{equation}\label{x3ypy3px2pxp1}
	x^3 y+y^3+x^2+x+1=0.
\end{equation}
If $(x,y)$ is any integer solution, then $y \neq 0$, and $x^2+x+1$ is divisible by $y$. So, we can write $x^2+x+1=yt$ for some integer $t$.  Substituting this into \eqref{x3ypy3px2pxp1} and cancelling $y$, we obtain $x^3+y^2+t=0$. Hence $x^3+t$ is divisible by $y$. But then $(x^3+t)-(x-1)(x^2+x+1)=t+1$ is divisible by $y$ as well. Let $k$ be an integer such that $k=\frac{t+1}{y}=\frac{((x^2+x+1)/y)+1}{y}=\frac{x^2+x+1+y}{y^2}$. Let $r=\frac{|x|}{|y|}$, and assume that $|y|>4$. Then
$$
r^3=\frac{|x|^3}{|y|^3}=\frac{|x^3y|}{|y|^4}=\frac{|y^3+x^2+x+1|}{|y|^4}<\frac{1}{4}+\frac{r^2}{16}+\frac{r}{64}+\frac{1}{256}.
$$ 
Then $f(r)<\frac{65}{256}$ where $f(r)=r(r^2-\frac{r}{16}-\frac{1}{64})$. The function $f(r)$ is increasing when $r \geq 0.67$, and $f(0.67)=0.26622...>\frac{65}{256}$. So $f(r)<\frac{65}{256}$ implies that $r<0.67$. We then have 
$$
|k|=\left|\frac{x^2}{y^2}+\frac{x}{y^2}+\frac{1}{y^2}+\frac{1}{y}\right|\leq r^2+\frac{r}{y}+\frac{1}{y^2}+\frac{1}{y} < 0.918.
$$ 
Hence $k=0$. Then, $x^2+x+1+y=0$ whose integer solutions are $(x,y)=(u,-u^2-u-1)$ for some integer $u$. Substituting this into \eqref{x3ypy3px2pxp1} gives $-u(u+2)(1+u+u^2)^2=0$ which has no integer solutions. Finally, we must consider the case $|y| \leq 4$. We can then conclude that all integer solutions to equation \eqref{x3ypy3px2pxp1} are
$$
(x,y)=(-2,-3),(0,-1).
$$

\vspace{10pt}

Table \ref{tab:H31Rungesol} presents the integer solutions to the equations listed in Table \ref{tab:H31Runge}. 

\begin{center}
\begin{tabular}{ |c|c|c|c|c|c| } 
		\hline
		Equation & Solution $(x,y)$ \\ 
		\hline\hline
		$x^3 y+y^3+x^2=0$ &$(0,0)$ \\  \hline
		$x^3 y+y^3+x^2-1=0$ &$(-2,-3),(-1,\pm 1),(\pm1,0),(0,1)$\\  \hline
		$x^3 y+y^3+x^2+1=0$ &$(-3,5),(0,-1),(1,-1)$ \\  \hline
		$x^3 y+y^3+x^2-y=0$ &$(0,\pm 1),(0,0),\pm (1,-1)$\\  \hline
		$x^3 y+y^3+x^2+y=0$ &$(-1,-1),(0,0)$\\  \hline
		$x^3 y+y^3+x^2-x=0$ &$(0,0),(1,0)$\\  \hline
		$x^3 y+y^3+x^2+x=0$ &$(-1,0),(-1,1),(0,0),(\pm 1,-1)$\\  \hline
		$x^3 y+y^3+x^2-3=0$ &$(1,1)$ \\  \hline
		$x^3 y+y^3+x^2+3=0$ &$(-2,1)$ \\  \hline
		$x^3 y+y^3+x^2-y-1=0$ &$(\pm 1,0)$ \\  \hline
		$x^3 y+y^3+x^2+y-1=0$&$(\pm 1,0)$ \\  \hline
		$x^3 y+y^3+x^2-x-1=0$&$(0,1)$ \\ \hline
		$x^3 y+y^3+x^2-x+1=0$&$(-2,1),(0,-1)$ \\  \hline
		$x^3 y+y^3+x^2+x-1=0$&$(0,1)$ \\  \hline
		$x^3 y+y^3+x^2+x+1=0$ &$(-2,-3),(0,-1)$\\  \hline
	\end{tabular}
	\captionof{table}{\label{tab:H31Rungesol} Integer solutions to the equations in Table \ref{tab:H31Runge}.}
\end{center}

\subsection{Exercise 3.70}\label{ex:H36Runge}
\textbf{\emph{Solve all the equations from Table \ref{tab:H36Runge}. }}
\begin{center}
\begin{tabular}{ |c|c|c|c|c|c| } 
 \hline
 $H$ & Equation & $H$ & Equation & $H$ & Equation \\ 
 \hline\hline
 $32$ & $x^3+x^2 y^2+y^3=0$ & $35$ & $x^3+x^2 y^2+y^3-x+1=0$ & $36$ & $x^3+x^2 y^2+y^3+x-y=0$ \\ 
 \hline
 $33$ & $x^3+x^2 y^2+y^3-1=0$ & $35$ & $x^3+x^2 y^2+y^3+x-1=0$ & $36$ & $x^3+x^2 y^2+y^3+x+y=0$ \\ 
 \hline
 $33$ & $x^3+x^2 y^2+y^3+1=0$ & $35$ & $x^3+x^2 y^2+y^3+x+1=0$ & $36$ & $x^3+x^2 y^2+y^3-2x=0$ \\ 
 \hline
 $34$ & $x^3+x^2 y^2+y^3-x=0$ & $36$ & $x^3+x^2 y^2+y^3-x-2=0$ & $36$ & $x^3+x^2 y^2+y^3+2x=0$ \\ 
 \hline
 $34$ & $x^3+x^2 y^2+y^3+x=0$ & $36$ & $x^3+x^2 y^2+y^3-x+2=0$ & $36$ & $x^3+x^2 y^2+y^3-xy=0$ \\ 
 \hline
 $35$ & $x^3+x^2 y^2+y^3-3=0$ & $36$ & $x^3+x^2 y^2+y^3+x-2=0$ & $36$ & $x^3+x^2 y^2+y^3-x^2=0$  \\ 
 \hline
 $35$ & $x^3+x^2 y^2+y^3+3=0$ & $36$ & $x^3+x^2 y^2+y^3+x+2=0$ & $36$ & $x^3+x^2 y^2+y^3+x^2=0$ \\ 
 \hline
 $35$ & $x^3+x^2 y^2+y^3-x-1=0$ & $36$ & $x^3+x^2 y^2+y^3-x-y=0$ & & \\ 
 \hline
\end{tabular}
\captionof{table}{\label{tab:H36Runge} Equations of the form \eqref{eq:Runquartic2} of size $H\leq 36$.}
\end{center} 

We will now look at solving equations of the form
\begin{equation}\label{eq:Runquartic2}
P_1(x,y)\cdot P_2(x,y) + Q(x,y) = 0, 
\end{equation} 
where $P_i(x,y)=a_i x^2 + b_i x y + c_i y^2$, $i=1,2$ are coprime, and $Q$ is a cubic polynomial. To solve these equations we will use the method from Section 3.4.7 in the book, which we summarise below for convenience.
First, we rearrange the equation to 
$$
(P_1(x,y)+ax+by)(P_2(x,y)+cx+dy) = Q_1(x,y), 
$$
where $Q_1$ is a quadratic polynomial, selecting $a,b,c,d$ in order to cancel the cubic terms of $Q(x,y)$. Next, we prove that the ratios $\frac{Q_1(x,y)}{(P_1(x,y)+ax+by)}$ and $\frac{Q_1(x,y)}{(P_2(x,y)+cx+dy)}$ are bounded for large $x$. As the ratios must be integer, there are only a finite number of cases to check. Let us illustrate this with examples.

Equation
$$
x^3+x^2 y^2+y^3=0
$$
is solved in Section 3.4.7 of the book, and its integer solutions are
$$
(x,y)=(-2,-2),(0,0).
$$

Equation
$$
x^3+x^2 y^2+y^3-1=0
$$
is solved in Section 3.4.7 of the book, and its integer solutions are
$$
(x,y)=(-3,-2),(-2,-3),(0,1),\pm(1,-1),(1,0).
$$

The first equation we will consider is    
\begin{equation}\label{x3px2y2py3p1}
x^3+x^2 y^2+y^3+1=0.
\end{equation}
This equation can be written as 
\begin{equation}\label{x3px2y2py3p1red}
(x^2+y)(y^2+x)=xy-1.
\end{equation}
 By swapping $x$ and $y$ if necessary, let us assume that $|x| \leq |y|$. The cases $|y|=|x|$ and $|y|=|x|+1$ can be checked separately, so let us assume that $|y|>|x|+1$. 
Then 
$$
|y^2+x|\geq y^2-|x|>|y|(|x|+1)-|x|=|xy|+|y|-|x|>|xy|+1 \geq |xy-1|,
$$
hence
$$
k=\frac{xy-1}{y^2+x}
$$
satisfies $|k| <1$. But \eqref{x3px2y2py3p1red} implies that $k$ is an integer, so $k=0$. Then we must have $xy-1=0$, which implies that $x= \pm 1$. By checking these cases as well as cases $|y|=|x|$ and $|y|=|x|+1$, we find that integer solutions to \eqref{x3px2y2py3p1} satisfying $|x|\leq |y|$ are $(x,y)=(-1,-1)$ and $(0,-1)$. From symmetry, the case $|x|>|y|$ returns another solution $(-1,0)$. Therefore, we can conclude that all integer solutions to equation \eqref{x3px2y2py3p1} are 
$$
(x,y)=(-1,-1),(-1,0),(0,-1).
$$
   
\vspace{10pt}

The next equation we will consider is    
\begin{equation}\label{x3px2y2py3mx}
x^3+x^2 y^2+y^3-x=0.
\end{equation}
This equation can be written as 
\begin{equation}\label{x3px2y2py3mx_red}
(x^2+y)(y^2+x)=xy+x.
\end{equation}
Let us first assume that $|y| \geq |x|+2$. Then, 
$$
|y^2+x|\geq y^2-|x| \geq |y|(|x|+2)-|x|= |xy|+2|y|-|x| >  |xy|+|x| \geq |xy+x|,
$$
hence 
$$
k=\frac{xy+x}{y^2+x}
$$
satisfies $|k|<1$. But \eqref{x3px2y2py3mx_red} implies that $k$ is an integer, so $k=0$. Then we must have $xy+x=0$ which implies that either $x=0$ or $y=-1$. Substituting this into \eqref{x3px2y2py3mx} we obtain integer solutions $(x,y)=(0,0)$ and $(\pm 1,-1)$. However, these solutions do not satisfy $|y| \geq |x|+2$. 
Let us next consider the case $|x|\geq |y|+2$. Then, 
$$
|x^2+y|\geq x^2-|y| \geq |x|(|y|+2)-|y|= |xy|+2|x|-|y| >  |xy|+|x| \geq |xy+x|,
$$
hence 
$$
k=\frac{xy+x}{y^2+x}
$$
satisfies $|k|<1$. But \eqref{x3px2y2py3mx_red} implies that $k$ is an integer, so $k=0$, or $xy+x=0$, then as before we obtain integer solutions $(x,y)=(0,0)$ and $(\pm 1,-1)$. 
 However, these solutions do not satisfy $|x|\geq |y|+2$. Hence, the possible values of $|x|-|y|$ are $0$ and $\pm 1$. It is easy to check these three cases and 
 conclude that all integer solutions to equation \eqref{x3px2y2py3mx} are
$$
(x,y)=(0,0),(\pm1,-1),(\pm 1,0).
$$

\vspace{10pt}

The next equation we will consider is
\begin{equation}\label{x3px2y2py3px}
x^3+x^2 y^2+y^3+x=0.
\end{equation}
This equation can be written as 
\begin{equation}\label{x3px2y2py3px_red}
(x^2+y)(y^2+x)=xy-x.
\end{equation}
Let us first assume that $|y| \geq |x|+2 $. Then,
$$
|y^2+x|\geq y^2-|x| \geq |y|(|x|+2)-|x|= |xy|+2|y|-|x| >  |xy|+|x| \geq |xy-x|,
$$
hence 
$$
k=\frac{xy-x}{y^2+x}
$$
satisfies $|k|<1$. But \eqref{x3px2y2py3px_red} implies that $k$ is an integer, so $k=0$. Then we must have $xy-x=0$, which implies that $x=0$ or $y=1$. Substituting these into \eqref{x3px2y2py3px} we obtain the integer solutions $(x,y)=(0,0)$ and $(-1,1)$ but these do not satisfy $|y| \geq |x|+2 $. 
Let us next consider the case $|x|\geq |y|+2$. Then, 
$$
|x^2+y|\geq x^2-|y| \geq |x|(|y|+2)-|y|= |xy|+2|x|-|y| >  |xy|+|x| \geq |xy-x|,
$$
hence 
$$
k=\frac{xy-x}{y^2+x}
$$
satisfies $|k|<1$. Hence $k=0$, or $xy-x=0$, then as before, we obtain solutions $(x,y)=(0,0),(- 1,1)$. However, these solutions do not satisfy $|x|\geq |y|+2$. Therefore, the possible values of $|x|-|y|$ are $0$ or $\pm 1$. It is easy to check these cases and conclude that all integer solutions to equation \eqref{x3px2y2py3px} are
$$
(x,y)=(-1,1),(0,0).
$$

\vspace{10pt}

The next equation we will consider is      
\begin{equation}\label{x3px2y2py3m3}
	x^3+x^2 y^2+y^3-3=0.
\end{equation}
This equation can be written as
\begin{equation}\label{x3px2y2py3m3_red}
	(x^2+y)(y^2+x)=xy+3.
\end{equation}
Let us first assume that $|y| \geq |x|+2$. Then,
$$
|y^2+x| \geq y^2-|x| \geq |y|(|x|+2)-|x| = |xy|+2|y|-|x| \geq |xy|+|x|+4 > |xy|+3\geq |xy+3|,
$$
hence
$$
k=\frac{xy+3}{y^2+x}
$$
satisfies $|k|<1$. But \eqref{x3px2y2py3m3_red} implies that $k$ is an integer, so $k=0$. Then we must have $xy+3=0$, so $(x,y)=\pm (-3,1)$ or $\pm (1,-3)$. However, these pairs are not integer solutions to \eqref{x3px2y2py3m3}, hence it has no integer solutions satisfying $|y| \geq |x|+2$. From symmetry, the equations also has no integer solutions such that $|x| \geq |y|+2$. Hence, $|x|-|y|$ can only be $-1,0$ or $1$. By checking these cases, we conclude 
that the unique integer solution to equation \eqref{x3px2y2py3m3} is
$$
(x,y)=(1,1). 
$$

\vspace{10pt}

The next equation we will consider is     
\begin{equation}\label{x3px2y2py3p3}
	x^3+x^2 y^2+y^3+3=0.
\end{equation}
This equation can be written as
\begin{equation}\label{x3px2y2py3p3_red}
	(x^2+y)(y^2+x)=xy-3.
\end{equation}
By swapping $x$ and $y$ if necessary, let us assume that $|x| \leq |y|$. The cases $|y| = |x|$ and $|y| = |x|+1$ can be checked separately, so let us assume that $|y| \geq |x|+2$. Then,
$$
|y^2+x| \geq y^2-|x| \geq |y|(|x|+2)-|x| = |xy|+2|y|-|x| \geq |xy|+|x|+4 > |xy|+3\geq |xy-3|,
$$
hence
$$
k=\frac{xy-3}{y^2+x}
$$
satisfies $|k|<1$. But \eqref{x3px2y2py3p3_red} implies that $k$ is an integer, so $k=0$, or $xy-3=0$, so $(x,y)=\pm(3,1)$ and $\pm(1,3)$. The solutions satisfying $|y| \geq |x|+2$ are $(x,y)=\pm(1,3)$, however, after substituting into \eqref{x3px2y2py3p3} we can see that these are not integer solutions to the equation.
By checking the cases $|y|=|x|+1$ and $|y|=|x|$, we find the integer solution $(x,y)=(1,-2)$. From symmetry, the case $|y| < |x|$ returns the integer solution $(x,y)=(-2,1)$. 
Therefore, we can conclude that all integer solutions to equation \eqref{x3px2y2py3p3} are
$$
(x,y)=(-2,1),(1,-2). 
$$

 \vspace{10pt}

The next equation we will consider is        
\begin{equation}\label{x3px2y2py3mxm1}
	x^3+x^2 y^2+y^3-x-1=0.
\end{equation}
This equation can be written as
\begin{equation}\label{x3px2y2py3mxm1_red}
	(x^2+y)(y^2+x)=xy+x+1.
\end{equation}
Let us first assume that $|y| \geq |x|+2$. Then,
$$
|y^2+x| \geq y^2-|x| \geq |y|(|x|+2)-|x| = |xy|+2|y|-|x| > |xy|+|x|+1 \geq |xy+x+1|,
$$
hence
$$
k=\frac{xy+x+1}{y^2+x}
$$
satisfies $|k|<1$. But \eqref{x3px2y2py3mxm1_red} implies that $k$ is an integer, so $k=0$. Then we must have $xy+x+1=0$, which implies that $(x,y)=(-1,0)$ and $(1,-2)$. However, these solutions do not satisfy $|y| \geq |x|+2$. Let us next consider the case $|x| \geq |y|+2$. Then,
$$
|x^2+y| \geq x^2-|y| \geq |x|(|y|+2)-|y| = |xy|+2|x|-|y| > |xy|+|x|+1 \geq |xy+x+1|,
$$
hence
$$
k=\frac{xy+x+1}{x^2+y}
$$
satisfies $|k|<1$. Hence $k=0$, and as before, we obtain integer solutions $(x,y)=(-1,0)$ and $(1,-2)$. These solutions do not satisfy $|x| \geq |y|+2$. Therefore, the possible values of $|x|-|y|$ are $0$ and $\pm 1$. It is easy to check these cases and conclude that the unique integer solution to equation \eqref{x3px2y2py3mxm1} is
$$
(x,y)=(0,1). 
$$

\vspace{10pt}

The next equation we will consider is    
\begin{equation}\label{x3px2y2py3mxp1}
x^3+x^2 y^2+y^3-x+1=0.
\end{equation}
This equation can be written as
\begin{equation}\label{x3px2y2py3mxp1_red}
(x^2+y)(y^2+x)=xy+x-1.
\end{equation}
Let us first assume that $|y| \geq |x|+2$. Then,
$$
|y^2+x| \geq y^2-|x| \geq |y|(|x|+2)-|x| = |xy|+2|y|-|x| > |xy|+|x|+1 \geq |xy+x-1|,
$$
hence
$$
k=\frac{xy+x-1}{y^2+x}
$$
satisfies $|k|<1$. But \eqref{x3px2y2py3mxp1_red} implies that $k$ is an integer, so $k=0$. Then we must have $xy+x-1=0$ which implies that $(x,y)=(-1,-2)$ and $(1,0)$. However, these solutions do not satisfy $|y| \geq |x|+2$. Let us next consider the case $|x| \geq |y|+2$. Then,
$$
|x^2+y| \geq x^2-|y| \geq |x|(|y|+2)-|y| = |xy|+2|x|-|y| > |xy|+|x|+1 \geq |xy+x-1|,
$$
hence
$$
k=\frac{xy+x-1}{x^2+y}
$$
satisfies $|k|<1$. Hence, $k=0$ and as before, we obtain integer solutions $(x,y)=(-1,-2)$ and $(1,0)$. However, these solutions do not satisfy $|x| \geq |y|+2$. Therefore, the possible values of $|x|-|y|$ are $0$ or $\pm 1$. It is easy to check these cases and conclude that all integer solutions to equation \eqref{x3px2y2py3mxp1} are
$$
(x,y)=(-2,1),(0,-1). 
$$
 
\vspace{10pt}

The next equation we will consider is    
\begin{equation}\label{x3px2y2py3pxm1}
x^3+x^2 y^2+y^3+x-1=0.
\end{equation}
This equation can be written as
\begin{equation}\label{x3px2y2py3pxm1_red}
(x^2+y)(y^2+x)=xy-x+1.
\end{equation}
Let us first assume that $|y| \geq |x|+2$. Then,
$$
|y^2+x| \geq y^2-|x| \geq |y|(|x|+2)-|x| = |xy|+2|y|-|x| > |xy|+|x|+1 \geq |xy-x+1|,
$$
hence
$$
k=\frac{xy-x+1}{y^2+x}
$$
satisfies $|k|<1$. But \eqref{x3px2y2py3pxm1_red} implies that $k$ is an integer, so $k=0$. Then we must have $xy-x+1=0$, which implies that $(x,y)=(-1,2)$ and $(1,0)$. However, these solutions do not satisfy $|y| \geq |x|+2$. Let us next consider the case $|x| \geq |y|+2$. Then,
$$
|x^2+y| \geq x^2-|y| \geq |x|(|y|+2)-|y| = |xy|+2|x|-|y| > |xy|+|x|+1 \geq |xy-x+1|,
$$
hence
$$
k=\frac{xy-x+1}{x^2+y}
$$
satisfies $|k|<1$. Hence, $k=0$ and as before, we obtain integer solutions $(x,y)=(-1,2)$ and $(1,0)$. However, these solutions do not satisfy $|x| \geq |y|+2$. Therefore, the possible values of $|x|-|y|$ are $0$ or $\pm 1$. It is easy to check these cases and 
conclude that the unique integer solution to equation \eqref{x3px2y2py3pxm1} is
$$
(x,y)=(0,1). 
$$

\vspace{10pt}

The next equation we will consider is        
\begin{equation}\label{x3px2y2py3pxp1}
x^3+x^2 y^2+y^3+x+1=0.
\end{equation}
This equation can be written as
\begin{equation}\label{x3px2y2py3pxp1_red}
(x^2+y)(y^2+x)=xy-x-1.
\end{equation}
Let us first assume that $|y| \geq |x|+2$. Then,
$$
|y^2+x| \geq y^2-|x| \geq |y|(|x|+2)-|x| = |xy|+2|y|-|x| > |xy|+|x|+1 \geq |xy-x-1|,
$$
hence
$$
k=\frac{xy-x-1}{y^2+x}
$$
satisfies $|k|<1$. But \eqref{x3px2y2py3pxp1_red} implies that $k$ is an integer, so $k=0$. Then we must have $xy-x-1=0$, which implies that $(x,y)=(-1,0)$ and $(1,2)$. However, these solutions do not satisfy $|y| \geq |x|+2$. Let us next consider the case $|x| \geq |y|+2$. Then,
$$
|x^2+y| \geq x^2-|y| \geq |x|(|y|+2)-|y| = |xy|+2|x|-|y| > |xy|+|x|+1 \geq |xy-x-1|,
$$
hence
$$
k=\frac{xy-x-1}{x^2+y}
$$
satisfies $|k|<1$. Hence $k=0$, and as before, we obtain integer solutions $(x,y)=(-1,0)$ and $(1,2)$. However, these solutions do not satisfy $|x| \geq |y|+2$. Therefore, the possible values of $|x|-|y|$ are $0$ or $\pm 1$. It is easy to check these cases and 
conclude that all integer solutions to equation \eqref{x3px2y2py3pxp1} are
$$
(x,y)=(-2,-3),(0,-1). 
$$
        
\vspace{10pt}
     
The next equation we will consider is         
\begin{equation}\label{x3px2y2py3mxm2}
x^3+x^2 y^2+y^3-x-2=0.
\end{equation}
This equation can be written as
\begin{equation}\label{x3px2y2py3mxm2_red}
(x^2+y)(y^2+x)=xy+x+2.
\end{equation}
Let us first assume that $|y| \geq |x|+2$. Then,
$$
|y^2+x| \geq y^2-|x| \geq |y|(|x|+2)-|x| = |xy|+2|y|-|x| > |xy|+|x|+2 \geq |xy+x+2|,
$$
hence
$$
k=\frac{xy+x+2}{y^2+x}
$$
satisfies $|k|<1$. But \eqref{x3px2y2py3mxm2_red} implies that $k$ is an integer, so $k=0$. Then we must have $xy+x+2=0$ which implies $(x,y)=(-2,0),(-1,1),(1,-3)$ and $(2,-2)$. However, the only solution satisfying $|y| \geq |x|+2$ is $(x,y)=(1,-3)$, after substituting this into \eqref{x3px2y2py3mxm2} we can see that this isn't a solution to the equation. Let us next consider the case $|x| \geq |y|+2$. Then,
$$
|x^2+y| \geq x^2-|y| \geq |x|(|y|+2)-|y| = |xy|+2|x|-|y| > |xy|+|x|+2 \geq |xy+x+2|,
$$
hence
$$
k=\frac{xy+x+2}{x^2+y}
$$
satisfies $|k|<1$. Hence $k=0$, and as before, we obtain integer solutions $(x,y)=(-2,0),(-1,1),(1,-3)$ and $(2,-2)$. The only solution satisfying $|x| \geq |y|+2$ is $(x,y)=(-2,0)$, after substituting this into \eqref{x3px2y2py3mxm2} we can see that this isn't a solution to the equation. Therefore, the possible values of $|x|-|y|$ are $0$ or $\pm 1$. It is easy to check these cases and 
conclude that all integer solutions to equation \eqref{x3px2y2py3mxm2} are
$$
(x,y)=(-2,-2),(\pm 1,1). 
$$

\vspace{10pt}

The next equation we will consider is         
\begin{equation}\label{x3px2y2py3mxp2}
x^3+x^2 y^2+y^3-x+2=0.
\end{equation}
This equation can be written as
\begin{equation}\label{x3px2y2py3mxp2_red}
(x^2+y)(y^2+x)=xy+x-2.
\end{equation}
Let us first assume that $|y| \geq |x|+2$. Then,
$$
|y^2+x| \geq y^2-|x| \geq |y|(|x|+2)-|x| = |xy|+2|y|-|x| > |xy|+|x|+2 \geq |xy+x-2|,
$$
hence
$$
k=\frac{xy+x-2}{y^2+x}
$$
satisfies $|k|<1$. But \eqref{x3px2y2py3mxp2_red} implies that $k$ is an integer, so $k=0$. Then we must have $xy+x-2=0$ which implies that $(x,y)=(-2,-2),(-1,-3),(1,1)$ and $(2,0)$. However, the only solution satisfying $|y| \geq |x|+2$ is $(x,y)=(-1,-3)$, after substituting this into \eqref{x3px2y2py3mxp2} we can see that this isn't a solution to the equation. Let us next consider the case $|x| \geq |y|+2$. Then,
$$
|x^2+y| \geq x^2-|y| \geq |x|(|y|+2)-|y| = |xy|+2|x|-|y| > |xy|+|x|+2 \geq |xy+x-2|,
$$
hence
$$
k=\frac{xy+x-2}{x^2+y}
$$
satisfies $|k|<1$. Hence $k=0$, and as before, we obtain integer solutions $(x,y)=(-2,-2),(-1,-3),(1,1)$ and $(2,0)$. However, the only solution satisfying $|x| \geq |y|+2$ is $(x,y)=(2,0)$, and after substituting this into \eqref{x3px2y2py3mxp2} we can see that this isn't a solution to the equation. Therefore, the possible values of $|x|-|y|$ are $0$ or $\pm 1$. It is easy to check these cases and 
conclude that there are no integer solutions to equation \eqref{x3px2y2py3mxp2}. 
           
           \vspace{10pt}
           
The next equation we will consider is      
\begin{equation}\label{x3px2y2py3pxm2}
x^3+x^2 y^2+y^3+x-2=0.
\end{equation}
This equation can be written as
\begin{equation}\label{x3px2y2py3pxm2_red}
(x^2+y)(y^2+x)=xy-x+2.
\end{equation}
Let us first assume that $|y| \geq |x|+2$. Then,
$$
|y^2+x| \geq y^2-|x| \geq |y|(|x|+2)-|x| = |xy|+2|y|-|x| > |xy|+|x|+2 \geq |xy-x+2|,
$$
hence
$$
k=\frac{xy-x+2}{y^2+x}
$$
satisfies $|k|<1$. But \eqref{x3px2y2py3pxm2_red} implies that $k$ is an integer, so $k=0$. Then we must have $xy-x+2=0$, which implies that $(x,y)=(-2,2),(-1,3),(1,-1)$ and $(2,0)$. However, the only solution satisfying $|y| \geq |x|+2$ is $(x,y)=(-1,3)$, which is not an integer solution to \eqref{x3px2y2py3pxm2}. Let us next consider the case $|x| \geq |y|+2$. Then,
$$
|x^2+y| \geq x^2-|y| \geq |x|(|y|+2)-|y| = |xy|+2|x|-|y| > |xy|+|x|+2 \geq |xy-x+2|,
$$
hence
$$
k=\frac{xy-x+2}{x^2+y}
$$
satisfies $|k|<1$. Hence $k=0$, and as before, we obtain integer solutions $(x,y)=(-2,2),(-1,3),(1,-1)$ and $(2,0)$. The only solution satisfying $|x| \geq |y|+2$ is $(x,y)=(2,0)$, which is not an integer solution to \eqref{x3px2y2py3pxm2}. Therefore, the possible values of $|x|-|y|$ are $0$ or $\pm 1$. It is easy to check these cases and 
 conclude that all integer solutions to equation \eqref{x3px2y2py3pxm2} are
$$
(x,y)=(1,-1),(1,0). 
$$
            
            \vspace{10pt} 
             
The next equation we will consider is        
\begin{equation}\label{x3px2y2py3pxp2}
x^3+x^2 y^2+y^3+x+2=0.
\end{equation}
This equation can be written as
\begin{equation}\label{x3px2y2py3pxp2_red}
(x^2+y)(y^2+x)=xy-x-2.
\end{equation}
Let us first assume that $|y| \geq |x|+2$. Then,
$$
|y^2+x| \geq y^2-|x| \geq |y|(|x|+2)-|x| = |xy|+2|y|-|x| > |xy|+|x|+2 \geq |xy-x-2|,
$$
hence
$$
k=\frac{xy-x-2}{y^2+x}
$$
satisfies $|k|<1$. But \eqref{x3px2y2py3pxp2_red} implies that $k$ is an integer, so $k=0$. Then we must have $xy-x-2=0$ which implies that $(x,y)=(-2,0),(-1,-1),(1,3)$ and $(2,2)$. The only solution satisfying $|y| \geq |x|+2$ is $(x,y)=(1,3)$, which is not an integer solution to \eqref{x3px2y2py3pxp2}. Let us next consider the case $|x| \geq |y|+2$. Then,
$$
|x^2+y| \geq x^2-|y| \geq |x|(|y|+2)-|y| = |xy|+2|x|-|y| > |xy|+|x|+2 \geq |xy-x-2|,
$$
hence
$$
k=\frac{xy-x-2}{x^2+y}
$$
satisfies $|k|<1$. Hence $k=0$, and as before we have $(x,y)=(-2,0),(-1,-1),(1,3)$ and $(2,2)$. However, the only solution satisfying $|x| \geq |y|+2$ is $(x,y)=(-2,0)$, which is not an integer solution to \eqref{x3px2y2py3pxp2}. Therefore, the possible values of $|x|-|y|$ are $0$ or $\pm 1$. It is easy to check these cases and 
conclude that all integer solutions to equation \eqref{x3px2y2py3pxp2} are
$$
(x,y)=(-3,-2),(-2,-2),(-1,-1),(-1,0),(1,-2). 
$$

\vspace{10pt}    

The next equation we will consider is          
\begin{equation}\label{x3px2y2py3mxmy}
	x^3+x^2 y^2+y^3-x-y=0.
\end{equation}
This equation can be written as
\begin{equation}\label{x3px2y2py3mxmy_red}
	(x^2+y)(y^2+x)=xy+x+y.
\end{equation}
By swapping $x$ and $y$ if necessary, let us assume that $|x| \leq |y|$. The cases $|x|=|y|$, $|x|=|y|+1$ and $|x|=|y|+2$ can be checked separately, so let us assume that $|y| \geq |x|+3$. Then,
$$
|y^2+x| \geq y^2-|x| \geq |y|(|x|+3)-|x| = |xy|+3|y|-|x| > |xy|+|x|+|y| \geq |xy+x+y|,
$$
hence
$$
k=\frac{xy+x+y}{y^2+x}
$$
satisfies $|k|<1$. But \eqref{x3px2y2py3mxmy_red} implies that $k$ is an integer, so $k=0$, or $xy+x+y=0$, so $(x,y)=(-2,-2)$ and $(0,0)$. These solutions do not satisfy $|y| \geq |x|+3$. 
Then, checking the cases $|x|=|y|$, $|x|=|y|+1$ and $|x|=|y|+2$, we find that integer solutions $|x| \leq |y|$ are $(x,y)=(0,0)$ and $(0, \pm 1)$. From symmetry, the case $|x| > |y|$ returns solutions $(\pm 1,0)$. 
Therefore, we can conclude that all integer solutions to equation \eqref{x3px2y2py3mxmy} are
$$
(x,y)=(0,\pm 1),(0,0),(\pm 1,0). 
$$

\vspace{10pt}          

The next equation we will consider is        
\begin{equation}\label{x3px2y2py3pxmy}
x^3+x^2 y^2+y^3+x-y=0.
\end{equation}
This equation can be written as
\begin{equation}\label{x3px2y2py3pxmy_red}
(x^2+y)(y^2+x)=xy-x+y.
\end{equation}
Let us first assume that $|y| \geq |x|+3$. Then,
$$
|y^2+x| \geq y^2-|x| \geq |y|(|x|+3)-|x| = |xy|+3|y|-|x| > |xy|+|x|+|y| \geq |xy-x+y|,
$$
hence
$$
k=\frac{xy-x+y}{y^2+x}
$$
satisfies $|k|<1$. But \eqref{x3px2y2py3pxmy_red} implies that $k$ is an integer, so $k=0$. Then we must have $xy-x+y=0$, which implies that $(x,y)=(-2,2)$ and $(0,0)$. However, these solutions do not satisfy $|y| \geq |x|+3$. Let us next consider the case $|x| \geq |y|+3$. Then,
$$
|x^2+y| \geq x^2-|y| \geq |x|(|y|+3)-|y| = |xy|+3|x|-|y| > |xy|+|x|+|y| \geq |xy-x+y|,
$$
hence
$$
k=\frac{xy-x+y}{x^2+y}
$$
satisfies $|k|<1$. Hence $k=0$, and as before $(x,y)=(-2,2)$ and $(0,0)$. However, these solutions do not satisfy $|x| \geq |y|+3$. Therefore, the possible values of $|x|-|y|$ are $0, \pm 1$ or $\pm 2$. It is easy to check these cases and 
conclude that all integer solutions to equation \eqref{x3px2y2py3pxmy} are
$$
(x,y)=(-3,-2),(-2,-2),(0,\pm 1),(0,0),(1,-2). 
$$
  
  \vspace{10pt}

The next equation we will consider is         
\begin{equation}\label{x3px2y2py3pxpy}
x^3+x^2 y^2+y^3+x+y=0.
\end{equation}
This equation can be written as
\begin{equation}\label{x3px2y2py3pxpy_red}
(x^2+y)(y^2+x)=xy-x-y.
\end{equation}
By swapping $x$ and $y$ if necessary, let us assume that $|x| \leq |y|$. The cases $|y|=|x|$, $|y|=|x|+1$, $|y|=|x|+2$ can be checked separately, so let us assume that $|y| \geq |x|+3$. Then,
$$
|y^2+x| \geq y^2-|x| \geq |y|(|x|+3)-|x| = |xy|+3|y|-|x| > |xy|+|x|+|y| \geq |xy-x-y|,
$$
hence
$$
k=\frac{xy-x-y}{y^2+x}
$$
satisfies $|k|<1$. But \eqref{x3px2y2py3pxpy_red} implies that $k$ is an integer, so $k=0$, or $xy-x-y=0$, so $(x,y)=(0,0)$ and $(2,2)$. These solutions do not satisfy $|y| \geq |x|+3$. 
By checking the cases, $|y|=|x|$, $|y|=|x|+1$ and $|y|=|x|+2$, we obtain the integer solution $(x,y)=(0,0)$. From symmetry, the case $|x|>|y|$ returns no further solutions. Therefore, we can conclude that the unique integer solution to equation \eqref{x3px2y2py3pxpy} is
$$
(x,y)=(0,0). 
$$
 
 \vspace{10pt}

The next equation we will consider is
\begin{equation}\label{x3px2y2py3m2x}
	x^3+x^2 y^2+y^3-2x=0.
\end{equation}
This equation can be written as
\begin{equation}\label{x3px2y2py3m2x_red}
	(x^2+y)(y^2+x)=xy+2x.
\end{equation}
Let us first assume that $|y| \geq |x|+3$. Then,
$$
|y^2+x| \geq y^2-|x| \geq |y|(|x|+3)-|x| = |xy|+3|y|-|x| > |xy|+2|x| \geq |xy+2x|,
$$
hence
$$
k=\frac{xy+2x}{y^2+x}
$$
satisfies $|k|<1$. But \eqref{x3px2y2py3m2x_red} implies that $k$ is an integer, so $k=0$. Then we must have $xy+2x=0$, which implies that either $x=0$ or $y=-2$. Substituting this into \eqref{x3px2y2py3m2x} we obtain solutions $(x,y)=(-4,-2)$ and $(0,0)$. However, these solutions do not satisfy $|y| \geq |x|+3$. Let us next consider the case $|x| \geq |y|+3$. Then,
$$
|x^2+y| \geq x^2-|y| \geq |x|(|y|+3)-|y| = |xy|+3|x|-|y| > |xy|+2|x| \geq |xy+2x|,
$$
hence
$$
k=\frac{xy+2x}{x^2+y}
$$
satisfies $|k|<1$. But \eqref{x3px2y2py3m2x_red} implies that $k$ is an integer, so $k=0$, or $xy+2x=0$, then, as before, we obtain integer solutions $(x,y)=(-4,-2)$ and $(0,0)$. However, these solutions do not satisfy $|x| \geq |y|+3$. Hence, the possible values of $|x|-|y|$ are $-2,-1,0,1$ or $2$. It is easy to check these five cases and conclude that the only integer solutions to equation \eqref{x3px2y2py3m2x} are 
$$
	(x,y)=(-4,-2),(0,0). 
$$
 
\vspace{10pt}

The next equation we will consider is          
\begin{equation}\label{x3px2y2py3p2x}
x^3+x^2 y^2+y^3+2x=0.
\end{equation}
This equation can be written as
\begin{equation}\label{x3px2y2py3p2x_red}
(x^2+y)(y^2+x)=xy-2x.
\end{equation}
Let us first assume that $|y| \geq |x|+3$. Then,
$$
|y^2+x| \geq y^2-|x| \geq |y|(|x|+3)-|x| = |xy|+3|y|-|x| > |xy|+2|x| \geq |xy-2x|,
$$
hence
$$
k=\frac{xy-2x}{y^2+x}
$$
satisfies $|k|<1$. But \eqref{x3px2y2py3p2x_red} implies that $k$ is an integer, so $k=0$. Then we must have $xy-2x=0$, which implies that $x=0$ or $y=2$. Substituting these into \eqref{x3px2y2py3p2x} we obtain the integer solutions $(x,y)=(-4,2)$ and $(0,0)$. However, these solutions do not satisfy $|y| \geq |x|+3$. Let us next consider the case $|x| \geq |y|+3$. Then,
$$
|x^2+y| \geq x^2-|y| \geq |x|(|y|+3)-|y| = |xy|+3|x|-|y| > |xy|+2|x| \geq |xy-2x|,
$$
hence
$$
k=\frac{xy-2x}{x^2+y}
$$
satisfies $|k|<1$. Hence $k=0$ and as before we obtain solutions $(x,y)=(-4,2)$ and $(0,0)$. However, these solutions do not satisfy $|x| \geq |y|+3$. Therefore, the possible values of $|x|-|y|$ are $0,\pm 1$ or $\pm 2$. It is easy to check these cases and 
conclude that all integer solutions to equation \eqref{x3px2y2py3p2x} are
$$
(x,y)=(-4,2),(0,0). 
$$
  
  \vspace{10pt}
  
The next equation we will consider is    
\begin{equation}\label{x3px2y2py3mxy}
x^3+x^2 y^2+y^3-xy=0. 
\end{equation}
This equation can be written as
\begin{equation}\label{x3px2y2py3mxy_red}
(x^2+y)(y^2+x)=2xy.
\end{equation}
By swapping $x$ and $y$ if necessary, let us assume that $|x| \geq |y|$. The case $|x| =|y|$ can be checked separately, and returns the integer solution $(x,y)=(0,0)$. So, let us assume that $|x| \geq |y|+1$. Then
$$
|x^2+y| \geq x^2-|y| \geq |x|(|y|+1)-|y| = |xy|+|x|-|y| \geq |xy|+1 > |xy|,
$$
hence
$$
k = \frac{2xy}{x^2+y}
$$
satisfies $|k| < 2$. But \eqref{x3px2y2py3mxy_red} implies that $k$ is an integer, so $k=0$ or $|k|= 1$. If $k=0$ then $xy=0$ so $x=0$ or $y=0$, substituting these solutions into \eqref{x3px2y2py3mxy}, we get the solution $(x,y)=(0,0)$, which is impossible with $|x| \geq |y|+1$. If $|k| =1$ then $|2xy|=|x^2+y|$, which is an easy equation whose solution is $(x,y)=(\pm 1,1),(0,0)$ do not satisfy $|x| \geq |y|+1$. 
From symmetry, the case $|x| < |y|$ returns no further solutions. 
Therefore, we can conclude that the unique integer solution to equation \eqref{x3px2y2py3mxy} is
$$
(x,y)=(0,0).
$$

 \vspace{10pt}

The next equation we will consider is         
\begin{equation}\label{x3px2y2py3mx2}
x^3+x^2 y^2+y^3-x^2=0.
\end{equation}
This equation can be written as
\begin{equation}\label{x3px2y2py3mx2_red}
(x^2+y)(y^2+x-1)=xy-y.
\end{equation}
Let us first assume that $|y| \geq |x|+2$. Then,
$$
|y^2+x-1| \geq y^2-|x|-1 \geq |y|(|x|+2)-|x|-1 = |xy|+2|y|-|x|-1 > |xy|+|y| \geq |xy-y|,
$$
hence
$$
k=\frac{xy-y}{y^2+x-1}
$$
satisfies $|k|<1$. But \eqref{x3px2y2py3mx2_red} implies that $k$ is an integer, so $k=0$. Then we must have $xy-y=0$, which implies that $x=1$ or $y=0$. Substituting these solutions into \eqref{x3px2y2py3mx2} we obtain the solutions $(x,y)=(0,0),(1,-1)$ and $(1,0)$, but these solutions do not satisfy $|y| \geq |x|+2$. Let us next consider the case $|x| \geq |y|+2$. Then,
$$
|x^2+y| \geq x^2-|y| \geq |x|(|y|+2)-|y| = |xy|+2|x|-|y| > |xy|+|y| \geq |xy-y|,
$$
hence
$$
k=\frac{xy-y}{x^2+y}
$$
satisfies $|k|<1$. Hence $k=0$, and as before, we obtain solutions $(x,y)=(0,0),(1,-1)$ and $(1,0)$. However, these solutions do not satisfy $|x| \geq |y|+2$. Therefore, the possible values of $|x|-|y|$ are $0$ or $\pm 1$. It is easy to check these cases and 
  conclude that all integer solutions to equation \eqref{x3px2y2py3mx2} are
$$
(x,y)=(0,0),\pm(1,-1),(1,0). 
$$

\vspace{10pt}

The final equation we will consider is       
\begin{equation}\label{x3px2y2py3px2}
x^3+x^2 y^2+y^3+x^2=0.
\end{equation}
This equation can be written as
\begin{equation}\label{x3px2y2py3px2_red}
(x^2+y)(y^2+x+1)=xy+y.
\end{equation}
Let us first assume that $|y| \geq |x|+2$. Then,
$$
|y^2+x+1| \geq y^2-|x|-1 \geq |y|(|x|+2)-|x|-1 = |xy|+2|y|-|x|-1 > |xy|+|y| \geq |xy+y|,
$$
hence
$$
k=\frac{xy+y}{y^2+x+1}
$$
satisfies $|k|<1$. But \eqref{x3px2y2py3px2_red} implies that $k$ is an integer, so $k=0$. Then we must have $xy+y=0$ which implies that $x=-1$ or $y=0$. Substituting these solutions into \eqref{x3px2y2py3px2} we obtain the solutions $(x,y)=(-1,-1),(-1,0)$ and $(0,0)$. However, these solutions do not satisfy $|y| \geq |x|+2$. Let us next consider the case $|x| \geq |y|+2$. Then,
$$
|x^2+y| \geq x^2-|y| \geq |x|(|y|+2)-|y| = |xy|+2|x|-|y| > |xy|+|y| \geq |xy+y|,
$$
hence
$$
k=\frac{xy+y}{x^2+y}
$$
satisfies $|k|<1$. Hence $k=0$, and as before, this implies $(x,y)=(-1,-1),(-1,0)$ and $(0,0)$. However, these solutions do not satisfy $|x| \geq |y|+2$. Therefore, the possible values of $|x|-|y|$ are $0$ or $\pm 1$. It is easy to check these cases and 
conclude that all integer solutions to equation \eqref{x3px2y2py3px2} are
$$
(x,y)=(-1,-1),(-1,0),(0,0). 
$$

	\begin{center}
	\begin{tabular}{ |c|c|c| } 
		\hline
		Equation & Solution $(x,y)$ \\ 
		\hline\hline
		$x^3+x^2 y^2+y^3=0$ &$(-2,-2),(0,0)$ \\  \hline
		$x^3+x^2 y^2+y^3-1=0$ &$(-3,-2),(-2,-3),(0,1),\pm (1,-1),(1,0)$ \\  \hline
		$x^3+x^2 y^2+y^3+1=0$ &$(-1,-1),(-1,0),(0,-1)$ \\  \hline
		$x^3+x^2 y^2+y^3-x=0$ &$(0,0),(\pm1,-1),(\pm 1,0)$ \\  \hline
		$x^3+x^2 y^2+y^3+x=0$ &$(-1,1),(0,0)$ \\  \hline
		$x^3+x^2 y^2+y^3-3=0$ &$(1,1)$ \\  \hline
		$x^3+x^2 y^2+y^3+3=0$ &$(-2,1),(1,-2)$ \\  \hline
		$x^3+x^2 y^2+y^3-x-1=0$ &$(0,1)$ \\  \hline
		$x^3+x^2 y^2+y^3-x+1=0$ &$(-2,1),(0,-1)$ \\  \hline
		$x^3+x^2 y^2+y^3+x-1=0$ &$(0,1)$ \\  \hline
		$x^3+x^2 y^2+y^3+x+1=0$ & $(-2,-3),(0,-1)$ \\  \hline
		$x^3+x^2 y^2+y^3-x-2=0$ & $(-2,-2),(\pm 1,1)$\\  \hline
		$x^3+x^2 y^2+y^3-x+2=0$ & - \\  \hline
		$x^3+x^2 y^2+y^3+x-2=0$ & $(1,-1),(1,0)$ \\  \hline
		$x^3+x^2 y^2+y^3+x+2=0$ & $(-3,-2),(-2,-2),(-1,-1),(-1,0),(1,-2)$ \\  \hline
		$x^3+x^2 y^2+y^3-x-y=0$ &$(0,\pm1),(0,0),(\pm 1,0)$  \\  \hline
		$x^3+x^2 y^2+y^3+x-y=0$&$(-3,-2),(-2,-2),(0,\pm 1),(0,0),(1,-2)$ \\  \hline
		$x^3+x^2 y^2+y^3+x+y=0$ &$(0,0)$\\  \hline
		$x^3+x^2 y^2+y^3-2x=0$&$(-4,-2),(0,0)$ \\  \hline
		$x^3+x^2 y^2+y^3+2x=0$ &$(-4,2),(0,0)$ \\  \hline
		$x^3+x^2 y^2+y^3-xy=0$ &$(0,0)$\\  \hline 
		$x^3+x^2 y^2+y^3-x^2=0$&$(0,0),\pm(1,-1),(1,0)$  \\  \hline
		$x^3+x^2 y^2+y^3+x^2=0$ &$(-1,-1),(-1,0),(0,0)$\\ \hline
	\end{tabular}
	\captionof{table}{\label{tab:H36Rungesol} Integer solutions to the equations in Table \ref{tab:H36Runge}.}
\end{center}

\subsection{Exercise 3.73}\label{ex:x4paxypy3}
\textbf{\emph{Use the method described in the proof of Proposition 3.72 in the book to solve all equations of the form 
\begin{equation}\label{eq:x4paxypy3}
x^4+axy+y^3=0, \quad a\neq 0
\end{equation} 
with $2 \leq a \leq 10$.  }}

This exercise will solve equations of the form \eqref{eq:x4paxypy3} for $2 \leq a \leq 10$, using the method presented in Section 3.4.8 of the book, which we summarise below for convenience.

If $xy=0$, then $(x,y)=(0,0)$. So, let $(x,y)$ be any solution to \eqref{eq:x4paxypy3} with $xy \neq 0$, and let $d$ be an integer of the same sign as $y$ such that $|d|=\text{gcd}(x,y)$. Then $x=d x_1$ and $y=d y_1$ with $x_1\neq 0$ and $y_1>0$ coprime. We can make the further change of variables $d=k x_1$ for some non-zero integer $k$ and $k^2=y_1 z$ for some integer $z>0$. Then let $u = \text{gcd}(y_1, z) > 0$ and write $y_1/u = v^2$, $z/u = w^2$ for some integers $v,w$, where we may assume that $w > 0$ and $v$ has the same sign as $k$.
Then we need to solve the equations 
\begin{equation}\label{x4paxypy3_pairs}
	u_iw_i^2 x_1^5+a+w_iu_i^3v^5=0, \quad i=1,\dots, N
\end{equation}
where $(u_i,w_i)$ are pairs of positive integers such that $u_iw_i$ is a divisor of $a$, and $N$ is the number of such pairs. 
These equations are Thue-like equations and can be solved with the Magma code:
\begin{equation}\label{thue:magma}
\begin{aligned} 
& {\tt R<x> := PolynomialRing(Integers());}
\\ & {\tt f := u_i*w_i^2*x^5+w_i*u_i^3;}
\\ & {\tt T := Thue(f);}
\\ & {\tt T;}
\\ & {\tt Solutions(T, -a);}
	\end{aligned}
\end{equation}
From the Magma output, we choose the solutions such that $x_1\neq 0$ and $v\neq 0$ and compute solutions $(x,y)\neq (0,0)$ using the formulas 
\begin{equation}\label{eq:kdy1xyformula}
k=uvw, \quad d=kx_1, \quad y_1 =uv^2, \quad x=dx_1 \quad \text{and} \quad y=dy_1.
\end{equation} 
We then add solution $(0,0)$ to the final list.

The first equation we will consider is \eqref{eq:x4paxypy3} with $a=2$, or specifically,
\begin{equation}\label{x4p2xypy3}
x^4+2xy+y^3=0.
\end{equation}
As $a=2$, we obtain the pairs of positive integers 
$(u_i,w_i)=(1,1), (1,2),(2,1)$. By \eqref{x4paxypy3_pairs}, we need to solve the following equations:
$$
x_1^5+v^5= -2, \quad 2x_1^5+v^5=-1, \quad x_1^5+4v^5=-1.
$$ 
These are Thue-like equations and can be solved using the Magma code \eqref{thue:magma}:
\begin{center}
\begin{tabular}{l}
 ${\tt R<x> := PolynomialRing(Integers());}$
\\ ${\tt f := x^5+1;}$
\\ ${\tt T := Thue(f);}$
\\ ${\tt T;}$
\\ ${\tt Solutions(T, -2);}$
\\ ${\tt f := 2*x^5+1;}$
\\ ${\tt T := Thue(f);}$
\\ ${\tt T;}$
\\ ${\tt Solutions(T, -1);}$
\\ ${\tt f := x^5+4;}$
\\ ${\tt T := Thue(f);}$
\\ ${\tt T;}$
\\ ${\tt Solutions(T, -1);}$\\
\end{tabular}
\end{center}
We then obtain the following integer solutions with $x_1 v\neq 0$: 
$(u,w,x_1,v) = (1,1,-1,-1),$ $(1,2,-1,1)$. 
We can then find integer solutions $(x,y) \neq (0,0)$ using the formulas \eqref{eq:kdy1xyformula}. Finally, we can conclude that the integer solutions to equation \eqref{x4p2xypy3} are
$$
(x,y)=(-1, 1), (0,0), (2, -2).
$$

\vspace{10pt}

This next equation we will consider is \eqref{eq:x4paxypy3} with $a=3$, or specifically,
\begin{equation}\label{x4p3xypy3}
x^4+3xy+y^3=0.
\end{equation}
As $a=3$, we obtain the pairs of positive integers 
$(u_i,w_i)=(1,1), (1,3),(3,1)$. By \eqref{x4paxypy3_pairs}, we need to solve the following equations:
$$
x_1^5+v^5=-3, \quad 3 x_1^5+v^5=-1, \quad x_1^5+9v^5=-1.
$$ 
The Magma code \eqref{thue:magma} for these equations outputs that there are no integer solutions with $x_1 v \neq 0$. Therefore, we can conclude that the unique integer solution to equation \eqref{x4p3xypy3} is
$$
(x,y)=(0,0).
$$

\vspace{10pt}

This next equation we will consider is \eqref{eq:x4paxypy3} with $a=4$, or specifically, 
\begin{equation}\label{x4p4xypy3}
x^4+4xy+y^3=0.
\end{equation}
As $a=4$, we obtain the pairs of positive integers 
 $(u_i,w_i)=(1,1),$ $(1,2),$ $(1,4),$ $(2,1),$ $(2,2),$ $(4,1)$. By \eqref{x4paxypy3_pairs}, we need to solve the following equations:
$$
x_1^5+v^5=-4, \quad 2 x_1^5+v^5=-2, \quad 4 x_1^5+v^5=-1, \quad x_1^5+4v^5=-2, 
$$ $$ 2 x_1^5+ 4v^5=-1, \quad  x_1^5+16 v^5=-1.
$$
The Magma code \eqref{thue:magma} for these equations outputs that there are no integer solutions with $x_1 v \neq 0$. Therefore, we can conclude that the unique integer solution to equation \eqref{x4p4xypy3} is
$$
(x,y)=(0,0).
$$

\vspace{10pt}

This next equation we will consider is \eqref{eq:x4paxypy3} with $a=5$, or specifically,
\begin{equation}\label{x4p5xypy3}
x^4+5xy+y^3=0.
\end{equation}
As $a=5$, we obtain the pairs of positive integers 
$(u_i,w_i)=(1,1), (1,5),(5,1)$. By \eqref{x4paxypy3_pairs}, we need to solve the following equations:
$$
x_1^5+v^5=-5, \quad 5 x_1^5+v^5=-1, \quad x_1^5+25v^5=-1.
$$
The Magma code \eqref{thue:magma} for these equations outputs that there are no integer solutions with $x_1 v \neq 0$. Therefore, we can conclude that the unique integer solution to equation \eqref{x4p5xypy3} is
$$
(x,y)=(0,0).
$$

\vspace{10pt}

This next equation we will consider is \eqref{eq:x4paxypy3} with $a=6$, or specifically,
\begin{equation}\label{x4p6xypy3}
x^4+6xy+y^3=0.
\end{equation}
As $a=6$, we obtain the pairs of positive integers 
$(u_i,w_i)=(1,1)$, $(1,2),$ $(1,3),$ $(1,6),$ $(2,1),$ $(2,3),$ $(3,1),$ $( 3,2),$ $( 6,1)$. By \eqref{x4paxypy3_pairs}, we need to solve the following equations:
$$
x_1^5+v^5=-6, \quad 2 x_1^5+v^5=-3, \quad 3 x_1^5+v^5=- 2, \quad 6 x_1^5+v^5=- 1, \quad x_1^5+4v^5=- 3,
$$ 
$$
3 x_1^5+ 4v^5=- 1, \quad  x_1^5+9 v^5=- 2, \quad  2x_1^5+9 v^5=- 1 , \quad  x_1^5+36 v^5=-1.
$$
The Magma code \eqref{thue:magma} for these equations outputs the following integer solutions with $x_1 v\neq 0$: $(u,w,x_1,v)=$ $(1,2,-1,-1),$ $(1,3,-1,1),$ $(2,1,1,-1),$ $(2,3,1,-1)$. We can then find integer solutions $(x,y)\neq (0,0)$ using the formulas \eqref{eq:kdy1xyformula}. Finally, we can conclude that the integer solutions to equation \eqref{x4p6xypy3} are
$$
(x,y)=(-6,-12),(-2,-4),(-2,2),(0,0),(3,-3).
$$

\vspace{10pt}

This next equation we will consider is \eqref{eq:x4paxypy3} with $a=7$, or specifically,
\begin{equation}\label{x4p7xypy3}
x^4+7xy+y^3=0.
\end{equation}
As $a=7$, we obtain the pairs of positive integers 
$(u_i,w_i)=(1,1),$ $(1,7),$ $(7,1)$. By \eqref{x4paxypy3_pairs}, we need to solve the following equations:
$$
x_1^5+v^5=- 7, \quad 7 x_1^5+v^5=- 1, \quad x_1^5+49v^5=- 1.
$$
The Magma code \eqref{thue:magma} for these equations outputs that there are no integer solutions with $x_1 v \neq 0$. 
Therefore, we can conclude that the unique integer solution to equation \eqref{x4p7xypy3} is
$$
(x,y)=(0,0).
$$

\vspace{10pt}

This next equation we will consider is \eqref{eq:x4paxypy3} with $a=8$, or specifically,
\begin{equation}\label{x4p8xypy3}
x^4+8xy+y^3=0.
\end{equation}
As $a=8$, we obtain the pairs of positive integers 
$(u_i,w_i)=(1,1),$ $(1,2),$ $(1,4),$ $(1,8),$ $(2,1),$ $(2,2),$ $(2,4),$ $(4,1),$ $(4,2),$ $(8,1)$. By \eqref{x4paxypy3_pairs}, we need to solve the following equations:
$$
x_1^5+v^5=-8, \quad 2 x_1^5+v^5=- 4, \quad 4 x_1^5+v^5=- 2, \quad 8 x_1^5+v^5=- 1, \quad x_1^5+4v^5=- 4,
$$ $$
x_1^5+ 2v^5=- 1, \quad  4x_1^5+4 v^5=-1, \quad  x_1^5+16 v^5=-2 , \quad  2x_1^5+16 v^5=-1, \quad  x_1^5+64 v^5=- 1.
$$
The Magma code \eqref{thue:magma} for these equations outputs the following integer solution with $x_1 v\neq 0$: 
 $(u,w,x_1,v)= (2,2,1,-1)$. We can then find the integer solution $(x,y)\neq(0,0)$ using the formulas \eqref{eq:kdy1xyformula}. Finally, we can conclude that the integer solutions to equation \eqref{x4p8xypy3} are
$$
(x,y)=(-4,-8),(0,0).
$$

\vspace{10pt}

This next equation we will consider is \eqref{eq:x4paxypy3} with $a=9$, or specifically,
\begin{equation}\label{x4p9xypy3}
x^4+9xy+y^3=0.
\end{equation}
As $a=9$, we obtain the pairs of positive integers 
$(u_i,w_i)=(1,1),$ $(1,3),$ $(1,9),$ $(3,1),$ $(3,3),$ $(9,1)$. By \eqref{x4paxypy3_pairs}, we need to solve the following equations:
$$
x_1^5+v^5=- 9, \quad 3 x_1^5+v^5=- 3, \quad 9 x_1^5+v^5=- 1, \quad x_1^5+9v^5=- 3,
$$ $$
3 x_1^5+ 9v^5=-1, \quad  x_1^5+81 v^5=-1.
$$
The Magma code \eqref{thue:magma} for these equations outputs that there are no integer solutions with $x_1 v \neq 0$.  
Therefore, we can conclude that the unique integer solution to equation \eqref{x4p9xypy3} is
$$
(x,y)=(0,0).
$$

\vspace{10pt}

This next equation we will consider is \eqref{eq:x4paxypy3} with $a=10$, or specifically,
\begin{equation}\label{x4p10xypy3}
x^4+10xy+y^3=0.
\end{equation}
As $a=10$, we obtain the pairs of positive integers 
$(u_i,w_i)=(1,1),$ $(1,2),$ $(1,5),$ $(1,10),$$(2,1),$ $(2,5),$ $(5,1),$ $(5,2),$ $(10,1)$. By \eqref{x4paxypy3_pairs}, we need to solve the following equations:
$$
x_1^5+v^5=-10, \quad 2 x_1^5+v^5=- 5, \quad 5 x_1^5+v^5=- 2, \quad 10 x_1^5+v^5=- 1, \quad x_1^5+4v^5=- 5,  
$$ $$
5 x_1^5+ 4v^5=- 1, \quad x_1^5+25 v^5=-2, \quad  2x_1^5+25 v^5=-1 , \quad  x_1^5+100 v^5=-1.
$$
The Magma code \eqref{thue:magma} for these equations outputs the following integer solution with $x_1 v\neq 0$:  
$(u,w,x_1,v)=(2,1,-1,-1),$ $(2,5,-1,1)$. We can then find integer solutions $(x,y)$ using the formulas \eqref{eq:kdy1xyformula}. Finally, we can conclude that the integer solutions to equation \eqref{x4p10xypy3} are
$$
(x,y)=(-2,4),(0,0),(10,-20).
$$

\begin{center}
\begin{tabular}{ |c|c|c| } 
 \hline
 $a$ & Equation & Solution $(x,y)$ \\ 
 \hline\hline
1&$x^4+xy+y^3=0$&$(0,0)$ \\\hline
2&$x^4+2xy+y^3=0$&$(-1, 1),(0,0), (2, -2)$ \\\hline
3&$x^4+3xy+y^3=0$&$(0,0)$ \\\hline
4&$x^4+4xy+y^3=0$&$(0,0)$ \\\hline
5&$x^4+5xy+y^3=0$&$(0,0)$ \\\hline
6&$x^4+6xy+y^3=0$&$(-6,-12),(-2,-4),(-2,2),(0,0),(3,-3)$ \\\hline
7&$x^4+7xy+y^3=0$&$(0,0)$ \\\hline
8&$x^4+8xy+y^3=0$&$(-4,-8),(0,0)$ \\\hline
9&$x^4+9xy+y^3=0$&$(0,0)$ \\\hline
10&$x^4+10xy+y^3=0$&$(-2,4),(0,0),(10,-20)$ \\\hline
\end{tabular}
\captionof{table}{\label{tab:Ex3.72} Integer solutions to equations of the form \eqref{eq:x4paxypy3} with $1 \leq a \leq 10$.}
\end{center} 

\subsection{Exercise 3.75}\label{ex:H32genus2}
\textbf{\emph{Find all integer solutions to the equations listed in Table \ref{tab:H32genus2}. }}
\begin{center}
\begin{tabular}{ |c|c|c|c|c|c| } 
 \hline
 $H$ & Equation \\ 
 \hline\hline
 $28$ & $x^4+xy+y^3=0$ \\ 
 \hline
 $32$ & $x^4+2xy+y^3=0$  \\ 
 \hline
 $32$ & $x^4+x^2+xy+y^3=0$  \\ 
 \hline
 $32$ & $x^4+xy-y^2+y^3=0$ \\ 
 \hline
 $32$ & $x^4+xy+y^2+y^3=0$ \\ 
 \hline
\end{tabular}
\captionof{table}{\label{tab:H32genus2} Equations of genus $g=2$ and size $H\leq 32$.}
\end{center}

We will now look at equations of genus $g=2$ and size $H \leq 32$. To solve these equations, we can make a rational change of variables to reduce the equation to a Hyperelliptic curve, that is, an equation of the form
$$
ay^m = P(x)
$$
where $m \geq 0$ and $a$ are integers, and $P(x)$ is a polynomial with integer coefficients.
We can find this change of variables using the {\tt Weierstrassform} command in Maple, see \eqref{maple:wei}. In all examples, we will have $a=1$ and $m=2$, so the curve is actually $y^2=P(x)$. 
We can then compute the rank of Jacobian of this curve using the Magma code:
\begin{equation}\label{hypercode}
	\begin{aligned}
& {\tt P<x> := PolynomialRing(Rationals());} \\
& {\tt C := HyperellipticCurve(P(x));} \\
& {\tt J := Jacobian(C);} \\
& {\tt RankBounds(J);}
	\end{aligned}
\end{equation}
This will output two integers, the upper and lower bound of the rank of the Jacobian. If the numbers are equal, then the rank is exactly that integer. If the rank is $0$, then we can compute all rational points on the curve using the code 
$$
{\tt Chabauty0(J)}
$$
 in Magma. This outputs points of the form $(a:b:c)$, if $c=0$ then this point represents a point at infinity. All other points correspond to the rational solutions $x= \frac{a}{c}$ and $y=\frac{b}{c^3}$, to the equation $y^2=P(x)$. If the rank of the Jacobian is $1$, then we can use the Chabauty method implemented in the Magma code:
\begin{equation}\label{hypercode2}
	\begin{aligned}
& {\tt ptsJ := Points(J : Bound := 10);} \\
& {\tt [Order(P) : P\, in \, ptsJ]; }
	\end{aligned}
\end{equation}
This will output the order of points, if the order is 0, then the point has infinite order. Let $n$ be the first point of infinite order, then we can run the following code:
\begin{equation}\label{hypercode3}
	\begin{aligned}
& {\tt P := ptsJ[n];} \\
& {\tt allptsC := Chabauty(P); allptsC;}
	\end{aligned}
\end{equation}
 This outputs all rational solutions to the equation, these solutions are of the same form as those outputted by ${\tt Chabauty0(J)}$. This method is not applicable to the equation $x^4-x^2+xy+y^3=0$ as this equation has rank$=2$ and genus$=2$. 

The equation
$$
x^4+xy+y^3=0
$$
is solved in Section 3.4.8 of the book and its only integer solution is
$$
(x,y)=(0, 0).
$$

The equation
$$
x^4+xy-y^2+y^3=0
$$
is solved in Section 3.4.8 of the book and its integer solutions are
$$
(x,y)=(-1, \pm 1), (0, 0), (0, 1), (2, -2), \quad \text{and} \quad (12, -27).
$$

The equation
$$
x^4+2xy+y^3=0
$$
is the equation \eqref{x4p2xypy3} which is solved in Section \ref{ex:x4paxypy3} and its integer solutions are
$$
(x,y)=(-1,1),(0, 0) \quad \text{and} \quad (2,-2).
$$

Let us first look at the equation
\begin{equation}\label{x4pxypy2py3}
x^4+xy+y^2+y^3=0.
\end{equation}
Using the command 
$$
{\tt Weierstrassform(x^4 + y^3 + x*y + y^2, x, y, X, Y)}
 $$
 in Maple returns
 $${\tt [4*X^5 + 4*X^4 + Y^2 - 1, x/y, (y^2 + 2*x + 2*y)/y^2, 2*X*(1 + X)/(-1 + Y), 2*(1 + X)/(-1 + Y)]}$$
  which means that equation \eqref{x4pxypy2py3} can be reduced to
\begin{equation}\label{x4pxypy2py3_red}
Y^2=-4X^5-4X^4+1,
\end{equation}
after the rational change of variables, 
$$
x=\frac{2X(1+X)}{Y-1}, \quad y=\frac{2(1+X)}{Y-1}.
$$
We can see in the Maple output that there is division by $y$, hence the case $y=0$ needs to be checked separately. By substituting $y=0$ into \eqref{x4pxypy2py3} we get the solution $(x,y)=(0,0)$. We can then compute the rank of the Jacobian of the reduced equation using the following code in Magma \eqref{hypercode}:
$$
\begin{aligned}
& {\tt P<x> := PolynomialRing(Rationals());}
\\ & {\tt C := HyperellipticCurve(-4*x^5 - 4*x^4 + 1); }
\\ & {\tt J := Jacobian(C);}
\\ & {\tt RankBounds(J); }\\
\end{aligned}
$$
This outputs ${\tt 1 \, 1}$, so the rank is exactly $1$. We can then use the Magma code \eqref{hypercode2}, which returns that the first point has finite order, but the second is of infinite order.  We can then use the Magma code \eqref{hypercode3} focussing on point $2$:
$$
\begin{aligned}
& {\tt  P := ptsJ[2]; } \\
& {\tt allptsC := Chabauty(P); allptsC;}
\end{aligned}
$$
which outputs 
$$
{\tt \{ (0 : -1 : 1), (0 : 1 : 1), (1 : 0 : 0), (-1 : 1 : 1), (-6 : -161 : 1), (-6 :161 : 1), (-1 : -1 : 1) \}}
$$
The point $(1:0:0)$ represents the point at infinity, so this can be ignored. The other points correspond to rational points $(X,Y)=(0,\pm 1),$ $(-1, \pm 1)$ and $(-6,\pm 161)$. 
We can then use the change of variables to find the set containing all integer solutions $(x,y)$ and determine which are integer. Therefore, we can conclude that all integer solutions to equation \eqref{x4pxypy2py3} are 
$$
(x,y)=(0,-1),(0,0),(1,-1).
$$

\vspace{10pt}

Let us look at the next equation
\begin{equation}\label{x4pxypx2py3}
x^4+xy+x^2+y^3=0.
\end{equation}
By using the command \eqref{maple:wei} for this equation, 
 \eqref{x4pxypx2py3} can be reduced to
\begin{equation}\label{x4pxypx2py3_red}
Y^2=X^6-4X^4-4,
\end{equation}
after the rational change of variables, 
$$
x=\frac{-2(1+X)}{X^3-Y}, \quad y=\frac{-2X(1+X)}{X^3-Y}.
$$
The case $x=0$ can be checked separately, and we obtain the integer solution 
$(x,y)=(0,0)$. We can then compute the rank of the Jacobian of the reduced equation using the Magma code \eqref{hypercode} 
which outputs that the rank is exactly $1$. We can then use the code \eqref{hypercode2}. 
The first point of infinite order is point $ $. We can then use the Magma code \eqref{hypercode3} focussing on point $2$, 
which outputs 
$$
{ \tt \{ (1 : -1 : 0), (-1 : 1 : 1), (1 : 1 : 0), (-1 : -1 : 1) \}}
$$
The points $(1:-1:0)$ and $(1:1:0)$ represent the points at infinity, so these can be ignored. The other points correspond to rational points $(X,Y)=(-1, \pm 1)$. 
We can then use the change of variables to find the set containing all integer solutions $(x,y)$ and determine which are integer. Therefore, we can conclude that all integer solutions to equation \eqref{x4pxypx2py3} are 
$$
(x,y)=(0,0),(1,-1).
$$

\begin{center}
	\begin{tabular}{ |c|c|c| } 
		\hline
		Equation & Solution $(x,y)$ \\ 
		\hline\hline
		$x^4+xy+y^3=0$&$(0,0)$ \\\hline
		$x^4+2xy+y^3=0$&$(-1, 1),(0,0), (2, -2)$ \\\hline
		$x^4+x^2+xy+y^3=0$&$(0,0),(1,-1)$ \\\hline
		$x^4+xy-y^2+y^3=0$&$(-1,\pm 1),(0,0),(0,1),(2,-2),(12,-27)$ \\\hline
		$x^4+xy+y^2+y^3=0$&$(0,-1),(0,0),(1,-1)$ \\\hline
		
	\end{tabular}
	\captionof{table}{\label{tab:Ex3.74} Integer solutions to the equations in Table \ref{tab:H32genus2}.}
\end{center} 

\subsection{Exercise 3.78}\label{ex:H31open}
\textbf{\emph{For each equation listed in Table \ref{tab:H31Open}, conjecture what its complete set of integer solutions is.}}
\begin{center}
\begin{tabular}{ |c|c|c|c|c|c| } 
 \hline
 $H$ & Equation & $H$ & Equation & $H$ & Equation \\ 
 \hline\hline
 $28$ & $x^4-y^3-x+y=0$ & $30$ & $x^4+y^3+2x-y=0$ & $31$ & $x^4+y^3+2x-y+1=0$ \\ 
 \hline
 $28$ & $x^4-y^3+x-y=0$ & $30$ & $x^4+y^3+2x+y=0$ & $31$ & $x^4+y^3+2x+y+1=0$ \\ 
 \hline
 $29$ & $x^4+y^3+xy+1=0$ & $30$ & $x^4+y^3+y^2+x=0$ & $31$ & $x^4+y^3+2x+y-1=0$ \\ 
 \hline
 $29$ & $x^4+y^3+xy-1=0$ & $30$ & $x^4+y^3+xy-y=0$ & $31$ & $x^4+y^3+xy-y-1=0$ \\ 
 \hline
 $30$ & $x^4-y^3+xy+x=0$ & $30$ & $x^4+y^3+xy+y=0$ & $31$ & $x^4+y^3+xy+3=0$ \\ 
 \hline
 $30$ & $x^4+y^3-y^2+x=0$ & $30$ & $x^4+y^3+xy+x=0$ & $31$ & $x^4+y^3+xy-3=0$  \\ 
 \hline
 $30$ & $x^4+y^3+x-2y=0$ & $31$ & $x^4-y^3+xy+x+1=0$ & $31$ & $x^4+y^3+xy+y-1=0$ \\ 
 \hline
 $30$ & $x^4+y^3+x-y-2=0$ & $31$ & $x^4-y^3+xy+x-1=0$ & $31$ & $x^4+y^3+xy+x+1=0$ \\ 
 \hline
 $30$ & $x^4+y^3+x+y+2=0$ & $31$ & $x^4+y^3+x-2y+1=0$ & $31$ & $x^4+y^3+xy+x-1=0$ \\ 
 \hline
 $30$ & $x^4+y^3+x+y-2=0$ & $31$ & $x^4+y^3+x-2y-1=0$ & & \\ 
 \hline
 $30$ & $x^4+y^3+x+2y=0$ & $31$ & $x^4+y^3+x+2y+1=0$ & & \\ 
 \hline
\end{tabular}
\captionof{table}{\label{tab:H31Open} Equations in $2$ variables of size $H\leq 31$ not solved/excluded in this chapter.}
\end{center} 

This exercise looks at equations of size $H\leq 31$ not solvable by any method considered in this chapter. Table \ref{tab:H31Open} conjectures the complete set of integer solutions to the unsolved equations. To conjecture these solutions, we can use the {\tt Reduce} command in Mathematica with a large bound on the variables. Table \ref{tab:H31Opensol} lists integer solutions to equations from Table \ref{tab:H31Open} in the range $|x|<10000$. Because a heuristic argument suggests that for these equations large integer solutions are unlikely to exist, we conjecture that the found solutions are the only ones. For the equation $x^4+y^3+xy+3=0$, the non-existence of integer solutions is proved in Proposition 7.7 of the book.

\begin{center}
	\begin{tabular}{ |c|c|c|c|c|c| } 
		\hline
		Equation & Solution $(x,y)$ \\ 
		\hline\hline
		$x^4-y^3-x+y=0$ &$(0,\pm 1),(0,0),(1,0),(1,\pm 1)$ \\  \hline
		$x^4-y^3+x-y=0$&$(-1,0),(0,0),(1,1)$ \\  \hline
		$x^4+y^3+xy+1=0$&$(0,-1),(1,-1)$\\  \hline
		$x^4+y^3+xy-1=0$&$(-1,\pm 1),(\pm 1,0),(0,1)$ \\  \hline
		$x^4-y^3+xy+x=0$ & $(-1,0),(0,0)$\\  \hline
		$x^4+y^3-y^2+x=0$ &$(-4,-6),(-1,0),(0,0),(0,1),\pm(1,-1)$\\  \hline
		$x^4+y^3+x-2y=0$ &$(-1,0),(0,0)$\\  \hline
		$x^4+y^3+x+2y=0$ &$(-1,0),(0,0)$\\  \hline
		$x^4+y^3+x-y-2=0$ &$(1,0),(1,\pm 1)$  \\  \hline
		$x^4+y^3+x+y+2=0$ &$(-1,-1),(0,-1)$\\  \hline
		$x^4+y^3+x+y-2=0$ &$(-1,1),(0,1),(1,0)$\\  \hline
		$x^4+y^3+2x-y=0$ &$(-8,-16),(0,0),(0,\pm 1)$\\  \hline
		$x^4+y^3+2x+y=0$ &$(0,0),(8,-16)$\\  \hline
		$x^4+y^3+y^2+x=0$ &$(-1,-1),(-1,0),(0,-1),(0,0),(2,-3)$ \\  \hline
		$x^4+y^3+xy-y=0$ &$(0,\pm 1),(0,0),\pm(1,-1),$\\  \hline
		$x^4+y^3+xy+y=0$ & $(-1,-1),(0,0)$\\  \hline
		$x^4+y^3+xy+x=0$ &$(-1,\pm 1),(-1,0),(0,0),(1,-1)$\\  \hline
		$x^4-y^3+xy+x+1=0$ &$(0,1)$ \\  \hline
		$x^4-y^3+xy+x-1=0$ &$(0,-1)$\\  \hline
		$x^4+y^3+xy+x+1=0$ &$(0,-1)$\\ \hline
		$x^4+y^3+xy+x-1=0$&$(0,1)$ \\  \hline
		$x^4+y^3+x-2y+1=0$ &$(-1,1),(0,1)$\\  \hline
		$x^4+y^3+x-2y-1=0$ &$(0,-1),\pm(1,1)$  \\  \hline
		$x^4+y^3+x+2y+1=0$ & $(1,-1)$ \\  \hline
		$x^4+y^3+2x-y+1=0$ &$(-1,\pm 1),(-1,0)$\\  \hline
		$x^4+y^3+2x+y+1=0$&$(-1,0)$ \\  \hline
		$x^4+y^3+2x+y-1=0$&$\pm (1,-1)$ \\  \hline
		$x^4+y^3+xy-y-1=0$ &$(\pm 1,0)$\\  \hline
		$x^4+y^3+xy+y-1=0$ &$(\pm 1,0),(3,-4)$\\  \hline
		$x^4+y^3+xy+3=0$&- \\ \hline
		$x^4+y^3+xy-3=0$&$(1,1)$  \\  \hline
	\end{tabular}
	\captionof{table}{\label{tab:H31Opensol} Conjectured list of all integer solutions to the equations in Table \ref{tab:H31Open}. }
\end{center}

\subsection{Exercise 3.80}\label{ex:H28abc}
\textbf{\emph{Show that effective abc conjecture reduces equation 
\begin{equation}\label{eq:y3pymx4mx}
y^3+y=x^4+x
\end{equation}
 to a finite computation.}}

We will reduce the equation \eqref{eq:y3pymx4mx} from the previous exercise to a finite computation by assuming the ``effective abc'' conjecture. For any integer $n> 0$, let $\mathrm{rad}(n)$ be the product of distinct prime factors of $n$. 
\begin{proposition}\label{prop:effectiveabc}[Proposition 3.79 in the book]
	The following statements are equivalent.
	\begin{itemize}
		\item[(i)] The effective abc conjecture is true, that is, for every $\epsilon>0$ there is only a finite number of triples of positive integers $(a,b,c)$ satisfying 
$$
a+b=c, \quad \text{gcd}(a,b,c)=1, \quad \text{and} \quad c > \mathrm{rad}(abc)^{1+\epsilon},
$$
		and there is an algorithm that lists them all. 
		\item[(ii)] For every $\epsilon>0$ and $\delta>0$, there is only a finite number of triples of positive integers $(a,b,c)$ satisfying 
		\begin{equation}\label{eq:abcexc2}
			a+b=c, \quad \text{gcd}(a,b,c)=1, \quad \text{and} \quad c > \delta \cdot  \mathrm{rad}(abc)^{1+\epsilon},
		\end{equation}	
		and there is an algorithm that lists them all. 
	\end{itemize}
\end{proposition}

If $y \leq |x|$, then $x^4+x=y^3+y \leq |x|^3+|x|$, which is not true for $|x|>1$. So, we must have $y>|x|$. Let $d=\gcd(x,y)$, then $x=dX$ and $y=dY$ for some integers $X,Y$ such that $\gcd(X,Y)=1$. Substituting this into \eqref{eq:y3pymx4mx} we obtain
$$
d^4X^4+dX=d^3Y^3+dY
$$
then dividing by $d^3$, we have
$$
dX^4=Y^3+\frac{Y-X}{d^2}
$$
as $dX^4$ and $Y^3$ are integers, $\frac{Y-X}{d^2}$ must also be integer. Letting $a=Y^3$, $b=\frac{Y-X}{d^2}$ and $c=dX^4$ we have $a,b,c$ are positive coprime integers. Let us now estimate $\mathrm{rad}(abc)$. We have $\mathrm{rad}(a)\leq Y$, $\mathrm{rad}(b) \leq b$ and $\mathrm{rad}(c)\leq dX$. Then
$$
\mathrm{rad}(abc)=\mathrm{rad}(a)\mathrm{rad}(b)\mathrm{rad}(c) \leq dXY\frac{(Y-X)}{d^2} =\frac{XY(Y-X)}{d}\leq XY(Y-X)
$$
Now, $a + b = c$ implies that $a < c$, or $Y^3 < dX^4$ hence $Y <d^{1/3} X^{4/3}=c^{1/3}$. Also, $|x| < y$ implies that
$|X| < Y$, hence $Y - X \leq Y + |X| < 2Y$. This implies that
$$
\mathrm{rad}(abc) \leq XY(Y-X) \leq 2XY^2 \leq 2X(c^{2/3}) =2\left(\frac{c}{d}\right)^{1/4} (c^{2/3}) \leq 2c^{11/12}.
$$
Hence \eqref{eq:abcexc2} holds with $\epsilon = 1/11$ and $\delta = (1/2)^{12/11}$, and the effective abc conjecture states that there are only finitely many such triples $(a, b, c)$, and there is an algorithm that can list them all. For each given triple $(a, b, c)$, it is easy to check whether the system $a = Y^3$, $b =\frac{ Y-X}{d^2}$, $c =dX^4$ has a solution in integers $X,Y, d$. Hence, effective abc conjecture reduces equation \eqref{eq:y3pymx4mx} to a finite computation.

\section{Chapter 4}
\subsection{Exercise 4.1}\label{ex:H14x1kx2}
\textbf{\emph{Solve all equations listed in Table \ref{tab:H14x1kx2}.}}
	\begin{center}
\begin{tabular}{ |c|c|c|c|c|c| } 
			\hline
			$H$ & Equation & $H$ & Equation & $H$ & Equation \\ 
			\hline\hline
			$13$ & $x^2y+z^2-1=0$ & $14$ & $xyz+t^2+2=0$ & $14$ & $x^2y+z^2-2=0$ \\ 
			\hline
			$13$ & $x^2y+z^2+1=0$ & $14$ & $xyz+t^2+t=0$ & $14$ & $x^2y+z^2+2=0$ \\ 
			\hline
			$14$ & $xyz+t^2-2=0$ & $14$ & $xyz+t^2+x=0$ & $14$ & $x^2y+z^2+z=0$ \\ 
			\hline
		\end{tabular}
		\captionof{table}{\label{tab:H14x1kx2} Equations of the form \eqref{eq:x1kx2} of size $H\leq 14$.}
	\end{center} 

All equation in Table \ref{tab:H14x1kx2} are of the form 
\begin{equation}\label{eq:x1kx2}
P(x_1^kx_2, x_3, \dots, x_n) = 0, \quad k \geq 1
\end{equation}
where $P$ is a polynomial with integer coefficients such that we know how to describe the integer solution set $S$ of the equation $P(y,x_2,\dots,x_n)=0$. We will describe the solutions using the Divisor function \eqref{Divisor Function}. The case $x_1=0$ must be checked separately. 

\begin{center}
\begin{tabular}{ |c|c|c|c|c|c| } 
 \hline
  Equation & Solution $(x,y,z)$ or $(x,y,z,t)$ \\ 
 \hline \hline
$x^2y+z^2-1=0$ & $(0,u,\pm 1), \left(v,\frac{1-u^2}{v^2},u\right), \quad u \in \mathbb{Z}, \quad v \in D_2(1-u^2)$ \\\hline
  $x^2y+z^2+1=0$ &$ \left(v,-\frac{1+u^2}{v^2},u\right), \quad u \in \mathbb{Z}, \quad v \in D_2(1+u^2)$ \\\hline
   $xyz+t^2-2=0$ & $ \left(v,w,\frac{2-u^2}{vw},u\right), \quad u \in \mathbb{Z}, \quad v \in D(2-u^2), \quad w \in D\left(\frac{2-u^2}{v}\right)$ \\\hline
$xyz+t^2+2=0$ & $ \left(v,w,- \frac{2+u^2}{vw},u\right), \quad u \in \mathbb{Z}, \quad v \in D(2+u^2), \quad w \in D\left(\frac{2+u^2}{v}\right)$ \\\hline
 $xyz+t^2+t=0$ & $ \left(v,w,- \frac{u^2+u}{vw},u\right), \quad u \in \mathbb{Z}, \quad v \in D(u^2+u), \quad w \in D\left(\frac{u^2+u}{v}\right),$ \\
 & $(0,u,v,0),(0,u,v,-1),(u,0,v,0),(u,0,v,-1), \quad u,v \in \mathbb{Z}  $ \\\hline
 $xyz+t^2+x=0$ &$ \left(v,w,- \frac{u^2+v}{vw},u\right), \quad u \in \mathbb{Z}, \quad v \in D(u^2), \quad w \in D\left(\frac{u^2+v}{v}\right)$ \\
 & $(0,u,v,0),(-u^2,0,v,u), \quad u,v \in \mathbb{Z}$ \\\hline
  $x^2y+z^2-2=0$ & $\left(v,\frac{2-u^2}{v^2},u\right), \quad u \in \mathbb{Z}, \quad v \in D_2(2-u^2)$ \\\hline
  $x^2y+z^2+2=0$ & $\left(v,- \frac{2+u^2}{v^2},u\right), \quad u \in \mathbb{Z}, \quad v \in D_2(2+u^2)$ \\\hline
    $x^2y+z^2+z=0$ & $(0,u,0),(0,u,-1), \left(v,- \frac{u^2+u}{v^2},u\right), \quad u \in \mathbb{Z}, \quad v \in D_2(u^2+u)$ \\\hline
          \end{tabular}
\captionof{table}{\label{tab:ex:4.1} Integer solutions to the equations listed in Table \ref{tab:H14x1kx2}.}
\end{center}

\subsection{Exercise 4.2}\label{ex:H16x2my2}
\textbf{\emph{Solve all equations listed in Table \ref{tab:H16x2my2}.}}

\begin{center}
\begin{tabular}{ |c|c|c|c|c|c| } 
 \hline
 $H$ & Equation & $H$ & Equation & $H$ & Equation \\ 
 \hline\hline
 $12$ & $x^2+y^2=z^2$ & $15$ & $x^2+y^2=z^2-3$ & $16$ & $x^2+x+y^2=z^2-2$ \\ 
 \hline
 $13$ & $x^2+y^2=z^2-1$ & $15$ & $x^2+y^2=z^2+3$ & $16$ & $x^2+x+y^2=z^2+2$ \\ 
 \hline
 $13$ & $x^2+y^2=z^2+1$ & $15$ & $x^2+x+y^2=z^2-1$ & $16$ & $2x^2+y^2=z^2$ \\ 
 \hline
 $14$ & $x^2+y^2=z^2-2$ & $15$ & $x^2+x+y^2=z^2+1$ & $16$ & $x^2+y^2=z^2+t^2$ \\ 
 \hline
 $14$ & $x^2+y^2=z^2+2$ & $16$ & $x^2+y^2=z^2-4$ & $16$ & $x^2+y^2+z^2=t^2$ \\ 
 \hline
 $14$ & $x^2+x+y^2=z^2$ & $16$ & $x^2+y^2=z^2+4$ & $16$ & $x^3+y^2=z^2$  \\ 
 \hline
\end{tabular}
\captionof{table}{\label{tab:H16x2my2} Equations of the form \eqref{eq:x12mx22mP} of size $H\leq 16$.}
\end{center} 

We will now consider equations of the form 
\begin{equation}\label{eq:x12mx22mP}
x_1^2-x_2^2 = P(x_3,\dots,x_n).
\end{equation}
To solve these equations, we can make the change of variables, $x'_1=x_1+x_2$ and $x'_2=x_1-x_2$, so the equation is reduced to 
$$
x'_1 x'_2=P(x_3,\dots,x_n),
$$
with the condition that $x'_1$ and $x'_2$ have the same parity. We can then do modulo analysis and make further changes of variables in order to remove these conditions and solve the equation with no restrictions. For example, if $P(x_3,\dots,x_n)$ is always even, then $x_1$ and $x_2$ have the same parity, and we may use the change of variables $x'_1=\frac{x_1+x_2}{2}$ and $x'_2=\frac{x_1-x_2}{2}$, after with we have no restrictions on $x'_1$ and $x'_2$. 

We will describe solutions using the Divisor function \eqref{Divisor Function} introduced in Section \ref{ex:divfunct}.

Equation 
$$
x^2+y^2=z^2-1
$$
is solved in Section 4.1.1 of the book and its integer solutions are
$$
(x,y)=\left(2u, \frac{4u^2+1-v^2}{2v}, \frac{4u^2+1+v^2}{2v}\right), \quad u \in {\mathbb Z}, \quad v \in D(4u^2+1).
$$

Equation
$$
x^2+y^2=z^2+t^2
$$
is solved in Section 4.1.1 of the book and its integer solutions are
$$
(x,y,z,t)=\left(wr+uv,vr-uw,wr-uv,vr+uw\right), \quad u,v,w,r \in {\mathbb Z},
$$
as well as the solutions obtained from the above by swapping $z$ and $t$.

Equation
$$
x^2+y^2+z^2=t^2
$$
is solved in Section 4.1.1 of the book and its integer solutions are
$$
(x,y,z,t)=\left(2s(uw+vr),2s(vw-ur),s(w^2+r^2-u^2-v^2),s(w^2+r^2+u^2+v^2)\right), \quad u,v,w,r,s \in {\mathbb Z},
$$
as well as the solutions obtained from the above by permuting $x,y,z$.

Let us now consider the equation
\begin{equation}\label{x2py2mz2m1}
x^2+y^2=z^2+1.
\end{equation}
Modulo $4$ analysis shows that $x$ and $y$ cannot be both even, hence, after swapping $x$ and $y$ if necessary, we may assume that $x$ is odd, that is, 
$x=2x'+1$ for some integer $x'$. Then we have 
$$
(2x'+1)^2-1=4(x')^2+4x'=(z-y)(z+y).
$$  
Because $4(x')^2+4x'$ is always even, $z$ and $y$ have the same parity, hence $y'=\frac{z-y}{2}$ and $z'=\frac{z+y}{2}$ are integers. Then the equation reduces to $(x')^2+x'=y'z'$, which has the integer solutions $(x',y',z')=(-1,0,u),(0,0,u),\left(u,v,\frac{u^2+u}{v}\right)$ with $u \in \mathbb{Z}$, and $v \in D(u^2+u)$. So, in the original variables, we have that the integer solutions to equation \eqref{x2py2mz2m1} are
\begin{equation}\label{x2py2mz2m1_2xp1}
	\begin{aligned}
(x,y,z)=(\pm 1,u,u), \quad \text{or} \quad \left(2u+1,\frac{u^2+u-v^2}{v},\frac{u^2+u+v^2}{v}\right), \quad \text{with} \quad u \in \mathbb{Z}, \\ \text{and} \quad v \in D(u^2+u),
\end{aligned}
\end{equation}
and the ones obtained from \eqref{x2py2mz2m1_2xp1} by swapping $x$ and $y$.

\vspace{10pt}

Let us consider the next equation
\begin{equation}\label{x2py2mz2p2}
x^2+y^2=z^2-2.
\end{equation}
We can rewrite this equation as $x^2+2=y'z'$, where $y'=z-y$ and $z'=z+y$ have the same parity. Modulo 4 analysis shows that $x$ is odd, that is, $x=2x'+1$ for some integer $x'$. After substitution, we then obtain 
$$
4(x')^2+4x'+3=y'z'.
$$  
Because $4(x')^2+4x'+3$ is always odd, its divisors $y'$ and $z'$ are automatically odd, so we may solve this equation without any additional parity restrictions. This equation has the solution $(x',y',z')=\left(u,v,\frac{4u^2+4u+3}{v}\right)$ with $u \in \mathbb{Z}$, and $v \in D(4u^2+4u+3)$. So, in the original variables, we have that the integer solutions to equation \eqref{x2py2mz2p2} are
$$
(x,y,z)=\left(2u+1,\frac{4u^2+4u+3-v^2}{2v},\frac{4u^2+4u+3+v^2}{2v}\right), \quad \text{with} \quad u \in \mathbb{Z}, \quad \text{and} \quad v \in D(4u^2+4u+3).
$$

\vspace{10pt}

Let us consider the next equation 
\begin{equation}\label{x2py2mz2m2}
x^2+y^2=z^2+2.
\end{equation}
We can rewrite this equation as $x^2-2=y'z'$, where $y'=z-y$ and $z'=z+y$ have the same parity. Modulo 4 analysis shows that  $x$ is odd. So let $x=2x'+1$ for some integer $x'$. After substitution, we then obtain 
$$
4(x')^2+4x'-1=y'z'.
$$  
Because $4(x')^2+4x'-1$ is always odd, its divisors $y'$ and $z'$ are automatically odd, so we may solve this equation without any additional parity restrictions. This equation has the solution $(x',y',z')=\left(u,v,\frac{4u^2+4u-1}{v}\right)$ with $u \in \mathbb{Z}$, and $v \in D(4u^2+4u-1)$. So, in the original variables, we have that the integer solutions to equation \eqref{x2py2mz2m2} are
$$
\begin{aligned}
(x,y,z)=\left(2u+1,\frac{4u^2+4u-1-v^2}{2v},\frac{4u^2+4u-1+v^2}{2v}\right), \quad \text{with} \quad u \in \mathbb{Z}, \\ \text{and} \quad v \in D(4u^2+4u-1).
\end{aligned}
$$

\vspace{10pt}

The next equation we will consider is
\begin{equation}\label{x2pxpy2mz2}
x^2+x+y^2=z^2.
\end{equation}
This equation can be rearranged to $x^2+x=(z-y)(z+y)$. Modulo $4$ analysis shows that $x$ is equal to $0$ or $-1$ modulo $4$, while $y$ and $z$ have the same parity, hence $y'=\frac{z-y}{2}$ and $z'=\frac{z+y}{2}$ are integers. The equation then becomes 
\begin{equation}\label{x2pxpy2mz2red}
\frac{x^2+x}{4}=y'z'. 
\end{equation}
We have either $x=4u$ or $x=4u-1$ for integers $u$. In these cases, equation \eqref{x2pxpy2mz2red} reduces to (i) $4u^2+u=y'z'$ and (ii) $4u^2-u=y'z'$, respectively. If $y'\neq 0$, the integer solution to (i) is $\left(u,v,\frac{4u^2+u}{v}\right)$ with $u \in \mathbb{Z}$, and $v \in D(4u^2+u)$, while the integer solution to (ii) is $\left(u,v,\frac{4u^2-u}{v}\right)$ with $u \in \mathbb{Z}$, and $v \in D(4u^2-u)$. The case $y'=0$ returns extra family $(0,0,u)$, $u \in {\mathbb Z}$.  Then, in the original variables, we have that the integer solutions to equation \eqref{x2pxpy2mz2} are
$$
	(x,y,z)=(0,u,u) \quad \text{or} \quad \left(4u,\frac{4u^2+u-v^2}{v},\frac{4u^2+u+v^2}{v}\right), \quad u \in \mathbb{Z}, \quad v \in D(4u^2+u).
$$
and
$$
	(x,y,z)=(-1,u,u) \quad \text{or} \quad \left(4u-1,\frac{4u^2-u-v^2}{v},\frac{4u^2-u+v^2}{v}\right), \quad u \in \mathbb{Z}, \quad v \in D(4u^2-u).
$$

\vspace{10pt}

The next equation we will consider is
\begin{equation}\label{x2py2mz2p3}
x^2+y^2=z^2-3.
\end{equation}
Modulo $4$ anlysis shows that $x$ and $y$ cannot both be even, hence, after swapping $x$ and $y$ if necessary, we may assume that $x$ is odd, that is, 
$x=2x'+1$ for some integer $x'$. Then we have
$$
4((x')^2+x'+1)=(y-z)(y+z).
$$
 Because the left-hand side of this equation is always even, $y$ and $z$ have the same parity, hence 
 $y'=\frac{z-y}{2}$ and $z'=\frac{z+y}{2}$ are integers. Then the equation reduces to 
$(x')^2+x'+1=y'z'$ 
which has the integer solutions $(x',y',z')=\left(u,v,\frac{u^2+u+1}{v} \right)$ with $u \in \mathbb{Z}$, and $v \in D(u^2+u+1)$. So, in the original variables, we have that the integer solutions to equation \eqref{x2py2mz2p3} are
\begin{equation}\label{x2py2mz2p3_2up1}
(x,y,z)=\left(2u+1,\frac{u^2+u+1-v^2}{v},\frac{u^2+u+1+v^2}{v}\right), \quad u \in \mathbb{Z}, \quad v \in D(u^2+u+1),
\end{equation}
and the ones obtained from \eqref{x2py2mz2p3_2up1} by swapping $x$ and $y$.

\vspace{10pt}

The next equation we will consider is
\begin{equation}\label{x2py2mz2m3}
x^2+y^2=z^2+3.
\end{equation}
Modulo $4$ analysis on this equation shows that $x$ must be even, so let 
$x=2x'$ for some integer $x'$. Then we can make the change of variables $y'=z-y$ and $z'=z+y$ which reduces the equation to 
$$
4(x')^2-3=y'z',
$$
with the restriction that $y'$ and $z'$ have the same parity. Because $4(x')^2-3$ is always odd, its divisors $y'$ and $z'$ are automatically odd, so we may solve this equation without any additional parity restrictions. The integer solutions to reduced equation are $(x',y',z')=\left(u,v,\frac{4u^2-3}{v}\right)$ with $u \in \mathbb{Z}$, and $v \in D(4u^2-3)$. Then, in the original variables, we can conclude that the integer solutions to \eqref{x2py2mz2m3} are
\begin{equation}\label{x2py2mz2m3red}
(x,y,z)=\left(2u,\frac{4u^2-3-v^2}{2v},\frac{4u^2-3+v^2}{2v}\right), \quad u \in \mathbb{Z}, \quad v \in D(4u^2-3).
\end{equation}

\vspace{10pt}

The next equation we will consider is
\begin{equation}\label{x2pxpy2mz2p1}
x^2+x+y^2=z^2-1.
\end{equation}
This equation can be rewritten as 
$x^2+x+1=y'z'$  where $y'=z-y$ and $z'=z+y$ have the same parity. Because $x^2+x+1$ is always odd, its divisors $y'$ and $z'$ must both be odd, and we can solve the equation without any additional parity restrictions on the variables. 
The reduced equation has the integer solutions $(x,y',z')=\left(u,v,\frac{u^2+u+1}{v}\right)$ with $u \in \mathbb{Z}$, and $v \in D(u^2+u+1)$. Then, in the original variables, we can conclude that the integer solutions to \eqref{x2pxpy2mz2p1} are
$$
(x,y,z)=\left(u,\frac{u^2+u+1-v^2}{2v},\frac{u^2+u+1+v^2}{2v}\right), \quad u \in \mathbb{Z}, \quad v \in D(u^2+u+1).
$$

\vspace{10pt}

The next equation we will consider is
\begin{equation}\label{x2pxpy2mz2m1}
x^2+x+y^2=z^2+1.
\end{equation}
This equation can be rewritten as 
$x^2+x-1=y'z'$ where $y'=z-y$ and $z'=z+y$ have the same parity. Becuase $x^2+x-1$ is always odd, 
both $y'$ and $z'$ must be odd, and we can solve the equation without any additional parity restrictions on the variables. 
The reduced equation has the integer solutions $(x,y',z')=\left(u,v,\frac{u^2+u-1}{v}\right)$ with $u \in \mathbb{Z}$, and $v \in D(u^2+u-1)$. Then, in the original variables, we have that the integer solutions to \eqref{x2pxpy2mz2m1} are
$$
(x,y,z)=\left(u,\frac{u^2+u-1-v^2}{2v},\frac{u^2+u-1+v^2}{2v}\right), \quad u \in \mathbb{Z}, \quad v \in D(u^2+u-1).
$$

\vspace{10pt}

The next equation we will consider is
\begin{equation}\label{x2py2mz2p4}
x^2+y^2=z^2-4.
\end{equation}
Modulo $4$ analysis shows that $x$ and $y$ cannot both be odd, hence, after swapping $x$ and $y$ if necessary, we may assume that $x$ is even, that is, 
 $x=2x'$ for some integer $x'$. 
We then have $4((x')^2+1)=(z-y)(z+y)$. Because $4((x')^2+1)$ is always even, $y$ and $z$ have the same parity, hence 
$y'=\frac{z-y}{2}$ and $z'=\frac{z+y}{2}$ are integers. Then the equation reduces to 
$$
(x')^2+1=y'z',
$$ 
which has the integer solutions $(x',y',z')=\left(u,v,\frac{u^2+1}{v}\right)$ with $u \in \mathbb{Z}$, and $v \in D(u^2+1)$. So, in the original variables, we have that the integer solutions to equation \eqref{x2py2mz2p4} are
\begin{equation}\label{x2py2mz2p4_2u}
(x,y,z)=\left(2u,\frac{u^2+1-v^2}{v},\frac{u^2+1+v^2}{v}\right), \quad u \in \mathbb{Z}, \quad v \in D(u^2+1),
\end{equation}
and the ones obtained from \eqref{x2py2mz2p4_2u} by swapping $x$ and $y$.

\vspace{10pt}

The next equation we will consider is
\begin{equation}\label{x2py2mz2m4}
x^2+y^2=z^2+4.
\end{equation}
Modulo $4$ analysis shows that $x$ and $y$ cannot both be odd, hence, after swapping $x$ and $y$ if necessary, we may assume that $x$ is even, that is, $x=2x'$ for some integer $x'$. 
Then we have $4((x')^2-1)=(z-y)(z+y)$. Because $4((x')^2-1)$ is always even, $y$ and $z$ have the same parity, hence $y'=\frac{z-y}{2}$ and $z'=\frac{z+y}{2}$ are integers. Then the equation reduces to 
$$
(x')^2-1=y'z'.
$$
which we can solve with no additional parity restrictions on the variables. This equation has the integer solutions $(x',y',z')=(\pm 1,0,u)$ or $\left(u,v,\frac{u^2-1}{v}\right)$ with $u \in \mathbb{Z}$, and $v \in D(u^2-1)$. Then, in the original variables, we obtain the integer solutions to equation \eqref{x2py2mz2m4} are 
\begin{equation}\label{x2py2mz2m4_2u}
(x,y,z)=(\pm 2,v,v) \quad \text{or} \quad \left(2u,\frac{u^2-1-v^2}{v},\frac{u^2-1+v^2}{v}\right), \quad u \in \mathbb{Z}, \quad v \in D(u^2-1),
\end{equation}
and the ones obtained from \eqref{x2py2mz2m4_2u} by swapping $x$ and $y$.

\vspace{10pt}

The next equation we will consider is
\begin{equation}\label{x2pxpy2mz2p2}
x^2+x+y^2=z^2-2.
\end{equation}
Modulo $4$ analysis shows that $x$ is equal to $1$ or $2$ modulo $4$, while $y$ and $z$ have the same parity, hence, $y'=\frac{z-y}{2}$ and $z'=\frac{z+y}{2}$ are integers. The equation then becomes
\begin{equation}\label{x2pxpy2mz2p2red} 
\frac{x^2+x+2}{4}=y'z'.
\end{equation} 
We have either $x=4u+1$ or $x=4u+2$ for integers $u$. In these cases, equation \eqref{x2pxpy2mz2p2red} reduces to (i) $4u^2+3u+1=y'z'$ and (ii) $4u^2+5u+2=y'z'$, respectively. The integer solution to (i) is $\left(u,v,\frac{4u^2+3u+1}{v}\right)$ with $u \in \mathbb{Z}$, and $v \in D(4u^2+3u+1)$, while the integer solution to (ii) is $\left(u,v,\frac{4u^2+5u+2}{v}\right)$ with $u \in \mathbb{Z}$, and $v \in D(4u^2+5u+2)$. Then, in the original variables, we have that the integer solutions to equation \eqref{x2pxpy2mz2p2} are 
$$
(x,y,z)=\left(4u+1,\frac{4u^2+3u+1-v^2}{v},\frac{4u^2+3u+1+v^2}{v}\right), \quad u \in \mathbb{Z}, \quad v \in D(4u^2+3u+1).
$$
and
$$
(x,y,z)=\left(4u+2,\frac{4u^2+5u+2-v^2}{v},\frac{4u^2+5u+2+v^2}{v}\right), \quad u \in \mathbb{Z}, \quad v \in D(4u^2+5u+2).
$$

\vspace{10pt}

The next equation we will consider is
\begin{equation}\label{x2pxpy2mz2m2}
x^2+x+y^2=z^2+2.
\end{equation}
Modulo $4$ analysis shows that $x$ is equal to $1$ or $2$ modulo $4$, while $y$ and $z$ have the same parity, hence, $y'=\frac{z-y}{2}$ and $z'=\frac{z+y}{2}$ are integers. The equation then becomes
\begin{equation}\label{x2pxpy2mz2m2red}
	\frac{x^2+x-2}{4}=y'z'.
\end{equation}
We have either $x=4u+1$ or $x=4u+2$ for integers $u$. In these cases, equation \eqref{x2pxpy2mz2m2red} reduces to (i) $4u^2+3u=y'z'$ and (ii) $4u^2+5u+1=y'z'$, respectively. The integer solution to (i) are $(0,0,u)$ or $\left(u,v,\frac{4u^2+3u}{v}\right)$ with $u \in \mathbb{Z}$, and $v \in D(4u^2+3u)$, while the integer solution to (ii) is $(-1,0,u)$ or $\left(u,v,\frac{4u^2+5u+1}{v}\right)$ with $u \in \mathbb{Z}$, and $v \in D(4u^2+5u+1)$. Then, in the original variables, we have that the integer solutions to equation \eqref{x2pxpy2mz2m2} are 
$$
(x,y,z)=(1,u,u) \quad \text{or} \quad \left(4u+1,\frac{4u^2+3u-v^2}{v},\frac{4u^2+3u+v^2}{v}\right), \quad u \in \mathbb{Z}, \quad v \in D(4u^2+3u).
$$
and
$$
	\begin{aligned}
(x,y,z)=(-2,u,u) \quad \text{or} \quad \left(4u+2,\frac{4u^2+5u+1-v^2}{v},\frac{4u^2+5u+1+v^2}{v}\right), \\ u \in \mathbb{Z}, \quad v \in D(4u^2+5u+1).
\end{aligned}
$$

\vspace{10pt}

The next equation we will consider is
\begin{equation}\label{2x2py2mz2}
2x^2+y^2=z^2.
\end{equation}
Modulo $4$ analysis shows that $x$ must be even and $y$ and $z$ have the same parity, so let $x=2x'$ for integer $x'$. Because integers $z-y$ and $z+y$ are both even, we can make the change of variables $y'=\frac{z-y}{2}$ and $z'=\frac{z+y}{2}$ to reduce the equation to $2(x')^2=y'z'$, and we have no restrictions on the variables. 
It has the integer solution $(x',y',z')=(0,0,u)$ or $\left(u,v,\frac{2u^2}{v}\right)$ with $u \in \mathbb{Z}$, and $v \in D(2u^2)$. Then, in the original variables, we have that the integer solutions to equation \eqref{2x2py2mz2} are
$$
(x,y,z)=(0,u,u) \quad \text{or} \quad \left(2u,\frac{2u^2-v^2}{v},\frac{2u^2+v^2}{v}\right), \quad u \in \mathbb{Z}, \quad v \in D(2u^2).
$$

\vspace{10pt}

The final equation we will consider is
\begin{equation}\label{x3py2mz2}
x^3+y^2=z^2.
\end{equation}
This equation can be rearranged to $x^3=(z-y)(z+y)$. 

Let us first consider the case when $x$ is even, so let $x=2x'$ for some integer $x'$. Then both $z-y$ and $z+y$ must be even as they have the same parity. So we can make the change of variables $y'=\frac{z-y}{2}$ and $z'=\frac{z+y}{2}$ and so the equation reduces to
$$
2(x')^3=y'z',
$$
with no restrictions on the variables, and it has the integer solutions $(x',y',z')=(0,0,u)$ or $\left(u,v,\frac{2u^3}{v}\right)$ with $u \in \mathbb{Z}$, and $v \in D(2u^3)$. Then, in the original variables, we obtain that the integer solutions to equation \eqref{x3py2mz2} with $x$ even are
\begin{equation}\label{x3py2mz2_2u}
(x,y,z)=(0,u,u), \quad \text{or} \quad \left(2u,\frac{2u^3-v^2}{v},\frac{2u^3+v^2}{v}\right), \quad u \in \mathbb{Z}, \quad v \in D(2u^3).
\end{equation}

Let us now consider the case when $x$ is odd. We can reduce the equation to $x^3=y'z'$ where 
integers $y'=z-y$ and $z'=z+y$ have the same parity. 
Then, as $x$ is odd, let $x=2x'+1$ for some integer $x'$. The equation then becomes 
$$
(2x'+1)^3=y'z'
$$ 
and there are no parity restrictions on the variables because the left-hand side is always odd, both $y'$ and $z'$ must be odd. The integer solutions to this equation are $(x',y',z')=\left(u,v,\frac{(2u+1)^3}{v} \right)$ with $u \in \mathbb{Z}$, and $v \in D((2u+1)^3)$. Hence, in the original variables, we have the integer solutions to equation \eqref{x3py2mz2} with odd $x$ are
\begin{equation}\label{x3py2mz2_2up1}
(x,y,z)=\left(2u+1,\frac{(2u+1)^3-v^2}{2v},\frac{(2u+1)^3+v^2}{2v}\right), \quad u \in \mathbb{Z}, \quad v \in D((2u+1)^3).
\end{equation}

Finally, we can conclude that all integer solutions to equation \eqref{x3py2mz2} are of the form \eqref{x3py2mz2_2u} or \eqref{x3py2mz2_2up1}.

\begin{center}
\begin{tabular}{ |c|c|c|c|c|c| } 
 \hline
 Equation & Solution $(x,y,z)$ or $(x,y,z,t)$ \\ 
 \hline\hline
  $x^2+y^2=z^2$&$(2uvw,u(w^2-v^2),u(w^2+v^2))$ or \\
 & $(u(w^2-v^2),2uvw,u(w^2+v^2))$, $u,v,w \in \mathbb{Z}$  \\  \hline
   $x^2+y^2=z^2-1$ & $\left(2u,\frac{4u^2+1-v^2}{2v},\frac{4u^2+1+v^2}{2v}\right)$, $u \in \mathbb{Z}$, $v \in D(4u^2+1)$  \\  \hline
 $x^2+y^2=z^2+1$ 
 &$(\pm 1,u,u),\left(2u+1,\frac{u^2+u-v^2}{v},\frac{u^2+u+v^2}{v}\right)$, $u \in \mathbb{Z}$, $v \in D(u^2+u)$ \\
 & and the ones obtained by swapping $x$ and $y$. \\
   \hline
  $x^2+y^2=z^2-2$ &$\left(2u+1,\frac{4u^2+4u+3-v^2}{2v},\frac{4u^2+4u+3+v^2}{2v}\right)$, $u \in \mathbb{Z}$, $v \in D(4u^2+4u+3)$\\  \hline
  $x^2+y^2=z^2+2$ & $\left(2u+1,\frac{4u^2+4u-1-v^2}{2v},\frac{4u^2+4u-1+v^2}{2v}\right)$, $u \in \mathbb{Z}$, $v \in D(4u^2+4u-1)$\ \\
\hline
\end{tabular}
\captionof{table}{\label{tab:H16x2my2sol} Integer solutions to the equations listed in Table \ref{tab:H16x2my2}.}
\end{center}

\subsection{Exercise 4.3}\label{ex:H16x1Ppax2pb}
\textbf{\emph{Solve all equations listed in Table \ref{tab:H16x1Ppax2pb}.}}

\begin{center}
\begin{tabular}{ |c|c|c|c|c|c| } 
 \hline
 $H$ & Equation & $H$ & Equation & $H$ & Equation \\ 
 \hline\hline
 $14$ & $x^2+xyz+y=0$ & $16$ & $x^2+xyz+x+y=0$ & $17$ & $x^2+xyz+x+y+1=0$ \\ 
 \hline
 $15$ & $x^2+xyz+y-1=0$ & $16$ & $x^2+xyz+2y=0$ & $17$ & $x^2+xyz+2y-1=0$ \\ 
 \hline
 $15$ & $x^2+xyz+y+1=0$ & $17$ & $x^2+xyz+y-3=0$ & $17$ & $x^2+xyz+2y+1=0$ \\ 
 \hline
 $16$ & $x^2+xyz+y-2=0$ & $17$ & $x^2+xyz+y+3=0$ &  & \\ 
 \hline
 $16$ & $x^2+xyz+y+2=0$ & $17$ & $x^2+xyz+x+y-1=0$ &  & \\ 
 \hline
\end{tabular}
\captionof{table}{\label{tab:H16x1Ppax2pb} Equations of the form \eqref{eq:x1Ppax2pb} of size $H\leq 17$.}
\end{center} 
s
In this exercise, we will consider equations of the form 
\begin{equation}\label{eq:x1Ppax2pb}
x_1 P(x_1,x_2,\dots,x_n)+ax_2+b=0
\end{equation}
where $a \neq 0$ and $b$ are integers, and $P$ is a polynomial with integer coefficients. We can see from this equation that $ax_2+b$ is divisible by $x_1$, hence we can write $ax_2+b=x_1 t$ for some integer $t$. The equation is then reduced to 
$x_1( P(x_1,\frac{x_1t-b}{a},\dots,x_n)+ t)=0$. So we have $x_1=0$ or $P(x_1,\frac{x_1t-b}{a},\dots,x_n)+ t=0$. We can multiply this equation by a constant to get integer coefficients. For each equation in Table \ref{tab:H16x1Ppax2pb}, the reduced equation is easier to solve than the original, and can therefore be solved using methods from previous exercises. 

The first equation we will consider is
\begin{equation}\label{x2pxyzpym1}
x^2+xyz+y-1=0.
\end{equation}
We can see from this equation that $y-1$ is divisible by $x$, so we can write $y-1=tx$ for some integer $t$. Then \eqref{x2pxyzpym1} can be written as
$$
x(x+xtz+z+t)=0
$$
so either $x=0$ or $x+xtz+z+t=0$. Up to the names of the variables, this is equation \eqref{xyzpxpypz}, whose integer solutions are given by \eqref{eq:xyzmxmymzsol}. 
 This implies that the integer solutions to equation \eqref{x2pxyzpym1} are
$$
\begin{aligned}
(x,y,z)=&(-1,u,1),(0,1,u),(\pm 1,0,u),(1,u,-1),(u,1-u^2,0),\\&(u,1-u,1),(u,1,-u),(u,u+1,-1), \quad u \in \mathbb{Z}.
\end{aligned}
$$

\vspace{10pt}

The next equation we will consider is
\begin{equation}\label{x2pxyzpxpy}
x^2+xyz+x+y=0.
\end{equation}
We can see from this equation that $y$ is divisible by $x$, so we can write $y=tx$ for some integer $t$. Then \eqref{x2pxyzpxpy} can be written as
$$
x(x+xtz+1+t)=0
$$
so either $x=0$ or $x+xtz+1+t=0$. Up to the names of the variables, this is equation \eqref{xyzpxpyp1}, whose integer solutions are given by \eqref{xyzpxpyp1sol}. 
This implies that the integer solutions to equation \eqref{x2pxyzpxpy} are
$$
\begin{aligned}
(x,y,z)=&(u,-u,1),(u,-u^2-u,0),(-1,0,u),(-1,u,1),(0,0,u), (1,1,-3), \\ &(1,2,-2),(2,2,-2),(2,6,-1),(3,6,-1), \quad u \in \mathbb{Z}.
\end{aligned}
$$

\vspace{10pt}

The next equation we will consider is
\begin{equation}\label{x2pxyzp2y}
x^2+xyz+2y=0.
\end{equation}
We can see from this equation that $2y$ is divisible by $x$, so we can write $2y=tx$ for some integer $t$. Then \eqref{x2pxyzp2y} can be written as
$$
x(2x+xtz+2t)=0
$$
so either $x=0$ or $2x+xtz+2t=0$. This equation has a dominating monomial and can be solved using the method from Section \ref{ex:boundedreal2}, and its integer solutions are 
$(x, t, z) =(0, 0, u),$ $(u, -u, 0),$  $\pm(-1, 2, 1), $  $\pm (1, 1, -4),$ $\pm (1, 2, -3),$ $\pm (2, -1, 1),$ $\pm (2, 1, -3),$ $\pm (2, 2, -2),$ $\pm (3, 6, -1),$ $\pm (4, 4, -1),$ $\pm (6, 3, -1),$  where $u$ is an integer. This implies that the integer solutions to equation \eqref{x2pxyzp2y} are
$$
\begin{aligned}
(x,y,z)=&(2u,-2u^2,0),(0,0,u),(-6,9,1),(-4,8,1),(-3,9,1),(-2,2,2),(-2,1,3), (-2,-1,-1),\\ & (-1,1,3),(-1,-1,1),(1,-1,-1), 
  (1,1,-3),(2,-1,1),(2,1,-3),\\ &(2,2,-2),(3,9,-1), (4,8,-1), (6,9,-1), \quad u \in \mathbb{Z}.
\end{aligned}
$$

\vspace{10pt}

The next equation we will consider is
\begin{equation}\label{x2pxyzpym2}
x^2+xyz+y-2=0.
\end{equation}
We can see from this equation that $y-2$ is divisible by $x$, so we can write $y-2=tx$ for some integer $t$. Then \eqref{x2pxyzpym2} can be written as
$$
x(x+xtz+2z+t)=0
$$
so either $x=0$ or $x+xtz+2z+t=0$. This equation has a dominating monomial and can be solved using the method from Section \ref{ex:boundedreal2}, and its integer solutions are 
$(x,t,z)=$ $(0,2u,-u),$ $(u,-u,0),$ $(2u,0,-u),$ $\pm (1,-3,-2),$ $\pm (3,-1,2),$ where $u$ is an arbitrary integer. This implies that the integer solutions to equation \eqref{x2pxyzpym2} are
$$
\begin{aligned}
(x,y,z)=&(0, 2, u),(u, -u^2 + 2, 0), (2 u, 2, -u),(-3, -1, -2), (-2, 2, 1),  \\
   & (-1, -1, 2), (\pm 1, 1, 0), (1, -1, -2), (3, -1, 2),  \quad u \in \mathbb{Z}.
\end{aligned}
$$

\vspace{10pt}

The next equation we will consider is
\begin{equation}\label{x2pxyzpyp2}
x^2+xyz+y+2=0.
\end{equation}
We can see from this equation that $y+2$ is divisible by $x$, so we can write $y+2=tx$ for some integer $t$. Then \eqref{x2pxyzpyp2} can be written as
$$
x(x+xtz-2z+t)=0
$$
so either $x=0$ or $x+xtz-2z+t=0$. This equation has a dominating monomial and can be solved using the method from Section \ref{ex:boundedreal2}, and its integer solutions are 
$(x,t,z)=$ $(2u,0,u),$ $(0,2u,u),$ $(u,-u,0),$  $\pm (1,1,2),$ $\pm (1,3,-4),$ $\pm (1,5,-2),$ $\pm (2,2,-2),$ $\pm (2,4,-1),$ $\pm (3,1,-4),$ $\pm (4,2,-1),$ $\pm (5,1,-2),$ where $u$ is an integer. This implies that the integer solutions to equation \eqref{x2pxyzpyp2} are
$$
\begin{aligned}
(x,y,z)=&(u, -u^2 - 2, 0), (2 u, -2, u),(0, -2, u), (-5, 3, 2), (-4, 6, 1), \\ &
 (-3, 1, 4), (-2, 2, 2), (-2, 6, 1),(-1, -1, -2),  \\& 
 (-1, 1, 4), (-1, 3, 2),  (1, -1, 2), (1, 1, -4), (1, 3, -2),  \\ & (2, 2, -2), (2, 6, -1), (3, 1, -4), (4, 6, -1), (5, 3, -2), \quad u \in \mathbb{Z}.
\end{aligned}
$$

\vspace{10pt}

The next equation we will consider is
\begin{equation}\label{x2pxyzp2ym1}
x^2+xyz+2y-1=0.
\end{equation}
We can see from this equation that $2y-1$ is divisible by $x$, so we can write $2y-1=tx$ for some integer $t$. Then \eqref{x2pxyzp2ym1} can be written as
$$
x(2x+xtz+z+2t)=0
$$
so either $x=0$ or $2x+xtz+z+2t=0$. This equation has a dominating monomial and can be solved using the method from Section \ref{ex:boundedreal2}, and its integer solutions are 
$(x,t,z)=$  $(0,u,-2u),$ $\pm (1,-1,u),$ $\pm (1,u,-2),$  $(u,-u,0),$ $(u,0,-2u),$ $\pm (u,-1,2),$ $\pm(3,5,-1),$ $\pm (5,3,-1),$ where $u$ is an integer. This implies that the integer solutions to equation \eqref{x2pxyzp2ym1} are
$$
\begin{aligned}
(x,y,z)=&(2u-1,u,-2), (2u+1,-2u^2-2u,0),(2u+1,-u,2),(\pm 1,0,u),\pm (1,u,-2), \\ &
 (-5,8,1),(-3,8,1),(3,8,-1),(5,8,-1) \quad u \in \mathbb{Z}.
\end{aligned}
$$

\vspace{10pt}

The next equation we will consider is
\begin{equation}\label{x2pxyzp2yp1}
x^2+xyz+2y+1=0.
\end{equation}
We can see from this equation that $2y+1$ is divisible by $x$, so we can write $2y+1=tx$ for some integer $t$. Then \eqref{x2pxyzp2yp1} can be written as
$$
x(2x+xtz-z+2t)=0
$$
so either $x=0$ or $2x+xtz-z+2t=0$. This equation has a dominating monomial and can be solved using the method from Section \ref{ex:boundedreal2}, and its integer solutions are 
$(x,t,z)=$ $(u,-u,0),$ $(u,0,2u),$ $(0,u,2u),$ $\pm (1,-3,-1),$ $\pm (1,2,-6),$ $\pm (1,3,-4),$ $\pm (1,5,-3),$ $\pm (2,1,-6),$ $\pm (2,3,-2),$ $\pm (3,-1,1),$ $\pm (3,1,-4),$ $\pm (3,2,-2),$ $\pm (3,7,-1),$ $\pm (5,1,-3),$ $\pm (7,3,-1),$ where $u$ is an integer. This implies that the integer solutions to equation \eqref{x2pxyzp2yp1} are
$$
\begin{aligned}
(x,y,z)=&(2u+1,-2u^2-2u-1,0),(-7,10,1),(-5,2,3),(-3,-2,-1),(-3,1,4),(-3,10,1),\\
  &(-1,-2,1),(-1,1,4),(-1,2,3),(1,-2,-1),(1,1,-4),(1,2,-3),(3,-2,1),\\
  &(3,1,-4),(3,10,-1),(5,2,-3),(7,10,-1),  \quad u \in \mathbb{Z}.
\end{aligned}
$$

\vspace{10pt}

The next equation we will consider is
\begin{equation}\label{x2pxyzpym3}
x^2+xyz+y-3=0.
\end{equation}
We can see from this equation that $y-3$ is divisible by $x$, so we can write $y-3=tx$ for some integer $t$. Then \eqref{x2pxyzpym3} can be written as
$$
x(x+xtz+3z+t)=0
$$
so either $x=0$ or $x+xtz+3z+t=0$. This equation has a dominating monomial and can be solved using the method from Section \ref{ex:boundedreal2}, and its integer solutions are 
$(x,t,z)=(-3u,0,u),$ $(u,-u,0),$ $(0,-3u,u),$ $\pm (-2,1,1),$ $\pm (-1,4,3),$ $\pm (-1,5,2),$  $\pm (1,-2,1),$ $\pm (4,-1,3),$ $\pm (5,-1,2),$ where $u$ is an integer. This implies that the integer solutions to equation \eqref{x2pxyzpym3} are
$$
\begin{aligned}
(x,y,z)=&(-3u,3,u),(0,3,u),(u,-u^2+3,0),(-5,-2,-2),(-4,-1,-3),(-2,1,1),\\
 &(-1,-2,2),(-1,-1,3),(-1,1,-1),(1,-2,-2),(1,-1,-3),\\
 &(1,1,1),(2,1,-1),(4,-1,3),(5,-2,2), \quad u \in \mathbb{Z}.
\end{aligned}
$$

\vspace{10pt}

The next equation we will consider is
\begin{equation}\label{x2pxyzpyp3}
x^2+xyz+y+3=0.
\end{equation}
We can see from this equation that $y+3$ is divisible by $x$, so we can write $y+3=tx$ for some integer $t$. Then \eqref{x2pxyzpyp3} can be written as
$$
x(x+xtz-3z+t)=0
$$
so either $x=0$ or $x+xtz-3z+t=0$. This equation has a dominating monomial and can be solved using the method from Section \ref{ex:boundedreal2}, and its integer solutions are 
$(x,t,z)=$ $(u,-u,0),$ $(3u,0,u),$ $(0,3u,u),$ $\pm (1,1,1),$ $\pm (1,2,3),$ $\pm (1,4,-5),$ $\pm (1,5,-3),$ $\pm (1,7,-2),$ $\pm (2,1,3),$ $\pm (2,2,-4),$ $\pm (2,5,-1),$ $\pm (3,3,-1),$ $\pm (4,1,-5),$ $\pm (5,1,-3),$ $\pm (5,2,-1),$ $\pm (7,1,-2),$  where $u$ is an arbitrary integer. This implies that the integer solutions to equation \eqref{x2pxyzpyp3} are
$$
\begin{aligned}
(x,y,z)=&(0,-3,u),(u,-u^2-3,0),(3u,-3,u),(-7,4,2),(-5,7,1),(-5,2,3),\\
  &(-4,1,5),(-3,6,1),(-2,-1,-3),(-2,1,4),(-2,7,1),(-1,-2,-1),\\
  &(-1,-1,-3),(-1,1,5),(-1,2,3),(-1,4,2),(1,-2,1),(1,-1,3),\\
  &(1,1,-5),(1,2,-3),(1,4,-2),(2,-1,3),(2,1,-4),(2,7,-1),\\
  &(3,6,-1),(4,1,-5),(5,2,-3),(5,7,-1),(7,4,-2), \quad u \in \mathbb{Z}.
\end{aligned}
$$

\vspace{10pt}

The next equation we will consider is
\begin{equation}\label{x2pxyzpxpyp1}
x^2+xyz+x+y+1=0.
\end{equation}
We can see from this equation that $y+1$ is divisible by $x$, so we can write $y+1=tx$ for some integer $t$. Then \eqref{x2pxyzpxpyp1} can be written as
$$
x(x+xtz-z+1+t)=0
$$
so, $x=0$ or $x+xtz-z+1+t=0$. This equation has a dominating monomial and can be solved using the method from Section \ref{ex:boundedreal2}, and its integer solutions are 
$(x,t,z)=$ $(0,u,u+1),$ $(u,-u-1,0),$ $(u,0,u+1),$  $(-2,-2,1),$ $(-2,-1,2),$ $(-1,-2,2),$ $(1,2,-4),$ $(1,4,-2),$ $(2,1,-4),$ $(2,4,-1),$ $(4,1,-2),$ $(4,2,-1),$ where $u$ is an integer. This implies that the integer solutions to equation \eqref{x2pxyzpxpyp1} are
$$
\begin{aligned}
(x,y,z)=&(u,-u^2-u-1,0),(u,-1,u+1),(0,-1,u),(-2,1,2),(-2,3,1),(-1,1,2),\\
 &(1,1,-4),(1,3,-2),(2,1,-4),(2,7,-1),(4,3,-2),(4,7,-1), \quad u \in \mathbb{Z}.
\end{aligned}
$$

\vspace{10pt}

The final equation we will consider is
\begin{equation}\label{x2pxyzpxpym1}
x^2+xyz+x+y-1=0.
\end{equation}
We can see from this equation that $y-1$ is divisible by $x$, so we can write $y-1=tx$ for some integer $t$. Then \eqref{x2pxyzpxpym1} can be written as
$$
x(x+xtz+z+1+t)=0
$$
so, $x=0$ or $x+xtz+z+1+t=0$. This equation has a dominating monomial and can be solved using the method from Section \ref{ex:boundedreal2}, and its integer solutions are 
$(x,t,z)=$ $(0,u,-1-u),$ $(u,0,-u-1),$ $(u,-1-u,0),$ $(-1,2,2),$ $(2,-1,2),$ $(2,2,-1),$  where $u$ is an integer. This implies that the integer solutions to equation \eqref{x2pxyzpxpym1} are
$$
\begin{aligned}
(x,y,z)=& (u,-u^2-u+1,0), (u,1,-u-1),(0,1,u), \\ &  (-1,-1,2),(2,-1,2),(2,5,-1), \quad u \in \mathbb{Z}.
\end{aligned}
$$

\begin{center}
\begin{tabular}{ |c|c|c|c|c|c| } 
	\hline
	Equation & Solution $(x,y,z)$ \\ 
	\hline\hline
	 $x^2+xyz+y=0$ &$(0,0,u), (u,-u,0), \pm(1,1,-2), \pm (2,2,-1)$ \\\hline
	 $x^2+xyz+y+1=0$ &$(0,-1,u),(u,-u^2-1,0),(u,-1,u),(-3,2,2),(-3,5,1),$\\
	&$(-2,1,3),(-2,5,1),(-1,1,3),(-1,2,2),(1,1,-3),$\\
	&$(1,2,-2),(2,1,-3),(2,5,-1),(3,2,-2),(3,5,-1)$ \\\hline	
	$x^2+xyz+y-1=0$ &$(-1,u,1),(0,1,u),(\pm 1,0,u),(1,u,-1),(u,1-u^2,0),$\\
	&$(u,1-u,1),(u,1,-u),(u,u+1,-1)$ \\\hline
	 $x^2+xyz+x+y=0$ &$(u,-u,1),(u,-u^2-u,0),(-1,0,u),(-1,u,1),(0,0,u),$ \\
\hline
\end{tabular}
\captionof{table}{\label{tab:H16x1Ppax2pbsol} Integer solutions to the equations listed in Table \ref{tab:H16x1Ppax2pb}. Assume that $u$ is an integer.}
\end{center}

\subsection{Exercise 4.4}\label{ex:H16subst}
\textbf{\emph{Solve all equations listed in Table \ref{tab:H16subst}.}}

\begin{center}
\begin{tabular}{ |c|c|c|c|c|c| } 
 \hline
 $H$ & Equation & $H$ & Equation & $H$ & Equation \\ 
 \hline\hline
 $14$ & $y^2+z^2=x^2+x$ & $16$ & $y^2+z^2=x^2+x-2$ & $16$ & $y^2+z^2=2x^2$  \\ 
 \hline
 $14$ & $x^2y+y=z^2$ & $16$ & $y^2+z^2=x^2+x+2$ & $16$ & $x^2y-2y=z^2$ \\ 
 \hline
 $14$ & $x^2y-y=z^2$ & $16$ & $y^2+y+z^2=x^2+x$ & $16$ & $x^2y+2y=z^2$ \\ 
 \hline
 $14$ & $x^2y+x=z^2$ & $16$ & $x^2+x+y^2+y=z^2$ & $16$ & $x^2y+2x=z^2$  \\ 
 \hline
 $15$ & $y^2+z^2=x^2+x-1$ & $16$ & $x^2+xy-y^2=z^2$ & $16$ & $x^2y+xy+z^2=0$ \\ 
 \hline
 $15$ & $y^2+z^2=x^2+x+1$ & $16$ & $x^2+xy+y^2=z^2$ & $16$ & $x^2y+xz+z^2=0$ \\ 
 \hline
\end{tabular}
\captionof{table}{\label{tab:H16subst} Equations of size $H\leq 16$ that can be simplified by substitutions.}
\end{center} 

By making substitutions, the equations in Table \ref{tab:H16subst}, can be reduced to equations of the form \eqref{eq:x1kx2} and can therefore be solved using the method in Section \ref{ex:H14x1kx2}.

The first equation we will consider is
\begin{equation}\label{x2ymymz2}
x^2y-y=z^2,
\end{equation}
which can be rewritten as 
$$
y(x^2-1)=z^2.
$$
This suggests to make the non-linear substitution $x'=x^2-1$, and the equation reduces to $x'y=z^2$. Up to the names of variables and parameters, we have solved this equation previously and its solution $(x',y,z)=(rk^2,rw^2,rkw),$  with $r,k,w \in \mathbb{Z}$ can be found in Table \ref{table1.38}. Then, $x^2-1=x'=rk^2$, $x^2-1=rk^2$ is an equation of the form \eqref{eq:x1kx2} and its integer solutions are $(x,k,r)=(\pm 1,0,u),\left(u,v,\frac{u^2-1}{v^2}\right)$ with $u \in \mathbb{Z}$ and $v \in D_2(u^2-1)$.  
Hence, the integer solutions to the original equation \eqref{x2ymymz2} are 
$$
(x,y,z)=(\pm 1,u,0) \quad \text{or} \quad \left(u,\frac{w^2(u^2-1)}{v^2}, \frac{w(u^2-1)}{v}\right), \quad u,w \in {\mathbb Z}, \quad v \in D_2(u^2-1).
$$

\vspace{10pt}

The next equation we will consider is
\begin{equation}\label{y2pz2mx2mxp1}
y^2+z^2=x^2+x-1.
\end{equation}
After multiplying the equation by 4 and making the substitutions $x'=2x+1$, $y'=2y$ and $z'=2z$, the equation is reduced to
$$
4y^2+5=(x')^2-(z')^2=(x'-z')(x'+z').
$$
This suggests to make the substitutions $x''=x'-z'$ and $z''=x'+z'$, which reduces the equation to $4y^2+5=x''z''$ with the restriction that $x=\frac{x''+z''-2}{4}$ and $z=\frac{z''-x''}{4}$ are integers. As $4y^2+5$ is $1$ modulo $4$, we must have either $x'' \equiv z'' \equiv 1 (\text{mod}\, 4)$ or $x'' \equiv z'' \equiv 3 (\text{mod}\, 4)$, and in both cases the listed restrictions are satisfied. The equation $4y^2+5=x''z''$ has the integer solutions 
$$
(x'',y,z'')=\left(v,u,\frac{4u^2+5}{v}\right)\quad u \in {\mathbb Z}, \quad v \in D(4u^2+5).
$$
 Then, in the original variables, the set of integer solutions to equation \eqref{y2pz2mx2mxp1} are
$$
(x,y,z)=\left(\frac{4u^2+5+v^2-2v}{4v},u,\frac{4u^2+5-v^2}{4v}\right), \quad u \in {\mathbb Z}, \quad v \in D(4u^2+5).
$$

\vspace{10pt}

The next equation we will consider is
\begin{equation}\label{y2pz2mx2mxm1}
y^2+z^2=x^2+x+1.
\end{equation}
After multiplying the equation by 4 and making the substitutions $x'=2x+1$, $y'=2y$ and $z'=2z$, the equation is reduced to
$$
(y')^2+(z')^2=(x')^2+3,
$$
with the restriction that $x'$ is odd and $y'$ and $z'$ are even. However, modulo $4$ analysis shows that every integer solution $(x',y',z')$ satisfies this condition, so we can solve it without restrictions.
Up to the names of variables, this is equation \eqref{x2py2mz2m3}, whose integer solutions are given by \eqref{x2py2mz2m3red}. 
 Then, in the original variables, the integer solutions to equation \eqref{y2pz2mx2mxm1} are
$$
(x,y,z)=\left(\frac{4u^2-3+v^2-2v}{4v},u,\frac{4u^2-3-v^2}{4v}\right), \quad u \in {\mathbb Z}, \quad v \in D(4u^2-3).
$$

\vspace{10pt}

The next equation we will consider is
\begin{equation}\label{y2pz2mx2mxp2}
y^2+z^2=x^2+x-2.
\end{equation}
After multiplying the equation by 4 and making the substitutions $x'=2x+1$ and $z'=2z$, the equation is reduced to
$$
4y^2+9=(x')^2-(z')^2=(x'-z')(x'+z').
$$
This suggests to make the substitutions $x''=x'-z'$ and $z''=x'+z'$, which reduces the equation to $4y^2+9=x''z''$ with the restriction that $x=\frac{x''+z''-2}{4}$ and $z=\frac{z''-x''}{4}$ are integers. Therefore, both $x''$ and $z''$ are odd and $x'' \equiv z''(\text{mod}\, 4)$. However, because $4y^2 + 9 \equiv 1(\text{mod}\,4)$, all integer solutions to $4y^2+9=x''z''$ automatically satisfy these conditions. This equation has the integer solutions 
$$
(x'',y,z'')=\left(v,u,\frac{4u^2+9}{v}\right), \quad u \in {\mathbb Z}, \quad v \in D(4u^2+9).
$$
 Then, in the original variables, the set of integer solutions to equation \eqref{y2pz2mx2mxp2} are
$$
(x,y,z)=\left(\frac{4u^2+9+v^2-2v}{4v},u,\frac{v^2-4u^2-9}{4v}\right), \quad u \in {\mathbb Z}, \quad v \in D(4u^2+9).
$$

\vspace{10pt}

The next equation we will consider is
\begin{equation}\label{y2pz2mx2mxm2}
y^2+z^2=x^2+x+2.
\end{equation}
After multiplying the equation by 4 and making the substitutions $x'=2x+1$ and $z'=2z$, the equation is reduced to
$$
4y^2-7=(x')^2-(z')^2=(x'-z')(x'+z').
$$
This suggests to make the substitutions $x''=x'-z'$ and $z''=x'+z'$, which reduces the equation to $4y^2-7=x''z''$ with the restriction that $x=\frac{x''+z''-2}{4}$ and $z=\frac{z''-x''}{4}$ are integers. Because $4y^2-7$ is $1$ modulo $4$, we must have either $x'' \equiv z'' \equiv 1 \, (\text{mod}\, 4)$ or $x'' \equiv z'' \equiv 3 \, (\text{mod}\, 4)$, and in both cases, the listed conditions are satisfied. 
 The equation $4y^2-7=x''z''$ has the integer solutions
$$
(x'',y,z'')=\left(v,u,\frac{4u^2-7}{v}\right),\quad u \in {\mathbb Z}, \quad v \in D(4u^2-7).
$$
 Then, in the original variables, the integer solutions to equation \eqref{y2pz2mx2mxm2} are
$$
(x,y,z)=\left(\frac{4u^2-7+v^2-2v}{4v},u,\frac{v^2-4u^2+7}{4v}\right), \quad u \in {\mathbb Z}, \quad v \in D(4u^2-7).
$$

\vspace{10pt}

 The next equation we will consider is
 \begin{equation}\label{y2pypz2mx2mx}
 y^2+y+z^2=x^2+x,
 \end{equation}
 which can be rewritten as 
 $$
 z^2=x^2+x-y^2-y=(x-y)(x+y+1).
 $$
 This suggests to make the substitutions $x'=x-y$ and $y'=x+y+1$, so $x=\frac{x'+y'-1}{2}$ and $y=\frac{y'-x'-1}{2}$, with the restriction that $x'$ and $y'$ have different parities. Up to the names of variables and parameters, we have solved this equation previously and its solution $(x',y',z)=(lm^2,ln^2,lmn)$ with $l,m,n \in \mathbb{Z}$ can be found in Table \ref{table1.38}. As we have the restriction that $x'$ and $y'$ have different parities, we must have either (i) $x'$ is odd and $y'$ is even or (ii) $y'$ is odd and $x'$ is even. 
 
 In case (i), because $x'$ is odd, both $l$ and $m$ must be odd. Let $l=2u+1$ and $m=2v+1$ for some integers $u$ and $v$. Then because $y'$ is even, $n$ must be even, so let $n=2w$ for some integer $w$. We now have $(x',y',z)=((2u+1)(2v+1)^2, (2u+1)(2w)^2,(2u+1)(2v+1)(2w))$. Hence, in the original variables, we obtain that a set of integer solutions to equation \eqref{y2pypz2mx2mx} is
 \begin{equation}\label{y2pypz2mx2mx_soli}
 \begin{aligned}
 (x,y,z)= \left(\frac{(2u+1)((2v+1)^2+4w^2)-1}{2},\frac{(2u+1)(4w^2-(2v+1)^2)-1}{2},2w(2u+1)(2v+1)\right), \\
  u,v,w \in \mathbb{Z}.
 \end{aligned}
 \end{equation}
 
  In case (ii), because $y'$ is odd, both $l$ and $n$ must be odd. Let $l=2u+1$ and $n=2w+1$ for some integers $u$ and $w$. Because $x'$ is even, $m$ must be even, so let $m=2v$ for some integer $v$. We now have $(x',y',z)=((2u+1)(2v)^2, (2u+1)(2w+1)^2,(2u+1)(2v)(2w+1))$. Hence, in the original variables, we obtain that the other set of integer solutions to equation \eqref{y2pypz2mx2mx} is 
 \begin{equation}\label{y2pypz2mx2mx_solii}
 \begin{aligned}
 (x,y,z)= \left(\frac{(2u+1)(4v^2+(2w+1)^2)-1}{2},\frac{(2u+1)((2w+1)^2-4v^2)-1}{2},2v(2u+1)(2w+1)\right), \\
  u,v,w \in \mathbb{Z}.
 \end{aligned}
 \end{equation}

Finally, we can conclude that integer solutions to \eqref{y2pypz2mx2mx} are of the form \eqref{y2pypz2mx2mx_soli} or \eqref{y2pypz2mx2mx_solii}.

\vspace{10pt}

 The next equation we will consider is
 \begin{equation}\label{x2pxpy2pymz2}
 x^2+x+y^2+y=z^2.
 \end{equation}
 After multiplication by 4 and rearranging, this equation can be rewritten as
 $$
 (2x+1)^2-2=4z^2-4y^2-4y-1=(2z)^2-(2y+1)^2=(2z+2y+1)(2z-2y-1).
 $$
 This suggests to make the substitutions $y'=2z-2y-1$ and $z'=2z+2y+1$, then we have $(2x+1)^2-2=y'z'$. As $y=\frac{z'-y'-2}{4}$ and $z=\frac{y'+z'}{4}$, we have the restrictions that $z'-y'-2 \equiv 0$ (mod $4$) and $y' \equiv -z'$ (mod $4$). Hence, we must have either $y' \equiv 1$ (mod $4$) and $z' \equiv 3$ (mod $4$) or $y' \equiv 3$ (mod $4$) and $z' \equiv 1$ (mod $4$). As the left-hand side of this equation is odd, both $y'$ and $z'$ must be odd. We also have that the left-hand side of this equation is $\equiv 3$ (mod $4$), then $y'z' \equiv 3$ (mod $4$), as $y'$ and $z'$ are odd, we must have $y' \equiv -z'$ (mod $4$), therefore we do not have any restrictions on this equation. The integer solutions to the reduced equation are $(x,y',z')=\left(u,v,\frac{(2u+1)^2-2}{v}\right)$ with $u \in \mathbb{Z}$ and $v \in D((2u+1)^2-2)$. Hence, in the original variables, we have that the integer solutions to \eqref{x2pxpy2pymz2} are
 $$
 (x,y,z)=\left(u,\frac{(2u+1)^2-2}{4v}-\frac{v+2}{4},\frac{v}{4}+\frac{(2u+1)^2-2}{4v}\right), \quad u \in \mathbb{Z}, \quad v \in D((2u+1)^2-2).
$$

\vspace{10pt}

The next equation we will consider is
 \begin{equation}\label{x2pxymy2mz2}
x^2+xy-y^2=z^2.
\end{equation}
After multiplying the equation by 4 and rearranging, this equation can be rewritten as
$$
5y^2=4x^2 + 4xy + y^2 -4z^2 = (2x+y)^2 - (2z)^2 = (2x+y-2z)(2x+y+2z).
$$
This suggests to make the substitutions $x'=2x+y-2z$ and $z'=2x+y+2z$, so $x=\frac{x'+z'-2y}{4}$ and $z=\frac{z'-x'}{4}$. The equation then reduces to 
 \begin{equation}\label{x2pxymy2mz2red}
5y^2=x'z',
\end{equation}
with the restriction that 
\begin{itemize}
	\item[(*)] $x' \equiv z'$ (mod $4$), and $x'+z'-2y \equiv 0$ (mod $4$). 
\end{itemize}

Let us first consider the case where $y$ is odd, so let $y=2y'+1$ for some integer $y'$. Equation \eqref{x2pxymy2mz2red} then becomes  $5(2y'+1)^2=x'z'$. Because $5(2y'+1)^2 \equiv 1$ (mod $4$), we have that $x'z' \equiv 1$ (mod $4$), which is possible only if $x'$ and $z'$ are odd and equivalent modulo $4$. Hence, the restrictions (*) are automatically satisfied and we may solve equation $5(2y'+1)^2=x'z'$ without any restrictions, 
and its integer solutions are $(x',y',z')=\left(v,u,\frac{5(2u+1)^2}{v} \right)$, with $u \in \mathbb{Z}$ and $v \in D(5(2u+1)^2)$. Then, in the original variables, we obtain that the integer solutions to equation \eqref{x2pxymy2mz2} with $y$ odd are
 \begin{equation}\label{x2pxymy2mz2_sol_y_odd}
(x,y,z)=\left(\frac{v}{4}+\frac{5(2u+1)^2}{4v}-u-\frac{1}{2},2u+1,\frac{5(2u+1)^2}{4v}-\frac{v}{4}  \right), \quad u \in \mathbb{Z}, \quad v \in D(5(2u+1)^2).
\end{equation}
Condition (*) guarantees that $(x,y,z)$ are always integers.

Let us now consider the case where $y$ is even. We will do this as two cases, (i) $y=4y'$ and (ii) $y=4y'+2$ for some integer $y'$. Let us first consider case (i), so $y=4y'$ and then \eqref{x2pxymy2mz2red} becomes $80(y')^2=x'z'$, restriction (*) implies that $x' \equiv z' \equiv -z'$ (mod $4$), so both $x'$ and $z'$ are even and equivalent modulo $4$. 
Because $80(y')^2$ is divisible by 16, both $x'$ and $z'$ must be divisible by $4$, so let $x'=4x''$ and $z'=4z''$ for some integers $x''$ and $z''$. Then we have $5(y')^2=x''z''$, which can be solved without restrictions, and its integer solutions are 
$(x'',y',z'')=(0,0,u)$ and $\left(v,u,\frac{5u^2}{v}\right)$ with $u \in \mathbb{Z}$ and $v \in D(5u^2)$. Then, in the original variables, we have that the integer solutions to equation \eqref{x2pxymy2mz2} with $y \equiv 0$ (mod $4$) are
 \begin{equation}\label{x2pxymy2mz2_sol_y_4u}
 (x,y,z)=(u,0,u), \quad  \left( \frac{5u^2}{v}+v-2u,4u,\frac{5u^2}{v}-v\right), \quad u \in \mathbb{Z}, \quad v \in D(5u^2).
 \end{equation}

 We will now consider case (ii), so $y=4y'+2$ for some integer $y'$. Then we have $5(4y'+2)^2=x'z'$. Restriction (*) implies that $x' \equiv z' \equiv -z'$ (mod $4$) so both $x'$ and $z'$ are even and equivalent modulo $4$. 
 As $5(4y'+2)^2$ is even, both $x'$ and $z'$ must be even, so let $x'=2x''$ and $z'=2z''$ for some integers $x''$ and $z''$. The equation is then reduced to $5(2y'+1)^2=x''z''$ with the restriction that $x''$ and $z''$ have the same parity. However, this is not a restriction as $5(2y'+1)^2 \equiv 1$ (mod $4$), so $x''z'' \equiv 1$ (mod $4$), and so we must have $x'' \equiv z''$ modulo $4$. Therefore we can solve the equation without restrictions, and its integer solutions are $(x'',y',z'')=\left(v,u,\frac{5(2u+1)^2}{v}\right)$ with $u \in \mathbb{Z}$ and $v \in D(5(2u+1)^2)$. Then in the original variables, we have that the integer solutions to equation \eqref{x2pxymy2mz2} with $y \equiv 2$ (mod $4$) are
  \begin{equation}\label{x2pxymy2mz2_sol_y_4up2}
 (x,y,z)=\left( \frac{5(2u+1)^2}{2v}+\frac{v}{2}-2u-1,4u+2,\frac{5(2u+1)^2}{2v}-\frac{v}{2}\right), \quad u \in \mathbb{Z}, \quad v \in D(5(2u+1)^2).
 \end{equation}
 Condition (*) guarantees that $(x,y,z)$ are always integers.
 
 Finally, we can conclude that all integer solutions to equation \eqref{x2pxymy2mz2} are of the form \eqref{x2pxymy2mz2_sol_y_odd}, \eqref{x2pxymy2mz2_sol_y_4u}, or \eqref{x2pxymy2mz2_sol_y_4up2}.
 
 \vspace{10pt}

The next equation we will consider is
 \begin{equation}\label{x2pxypy2mz2}
x^2+xy+y^2=z^2.
\end{equation}
After multiplication by 4 and rearranging, we can rewrite the equation as
$$
-3y^2=4x^2+4xy+y^2-4z^2=(2x+y)^2-(2z)^2=(2x+y-2z)(2x+y+2z).
$$
This suggests to make the substitutions $x'=2x+y-2z$ and $z'=2x+y+2z$, so $x=\frac{x'+z'-2y}{4}$ and $z=\frac{z'-x'}{4}$. This reduces the equation to 
$$
-3y^2=x'z',
$$
with the restriction that
\begin{itemize}
\item[(*)] $x' \equiv z'$ (mod $4$), and $x'+z'-2y \equiv 0$ (mod $4$).
\end{itemize}
  
Let us first consider the case where $y$ is odd. So let $y=2y'+1$ for some integer $y'$. Then we obtain $-3(2y'+1)^2=x'z'$. Because $-3(2y'+1)^2 \equiv 1 \, \text{(mod $4$)}$, we have $x'z' \equiv 1 \, \text{(mod $4$)}$, which is possible only if $x'$ and $z'$ are odd and equivalent modulo $4$. 
Hence, the restrictions (*) are automatically satisfied, and we may solve equation $-3(2y'+1)^2=x'z'$ without any restrictions. This equation has the integer solutions $(x',y',z')=\left(v,u,-\frac{3(2u+1)^2}{v} \right)$, with $u \in \mathbb{Z}$ and $v \in D(3(2u+1)^2)$. Then, in the original variables, the set of integer solutions to equation \eqref{x2pxypy2mz2} with $y$ odd are
 \begin{equation}\label{x2pxypy2mz2_sol_y_odd}
(x,y,z)=\left(\frac{v}{4}-\frac{3(2u+1)^2}{4v}-u-\frac{1}{2},2u+1,-\frac{3(2u+1)^2}{4v}-\frac{v}{4}  \right), \quad u \in \mathbb{Z}, \quad v \in D(3(2u+1)^2).
\end{equation}
Condition (*) guarantees that $(x,y,z)$ are always integers.

Let us next consider the case where $y$ is even. We will do this as two cases, (i) $y=4y'$ and (ii) $y=4y'+2$ for some integer $y'$. Let us first consider case (i), so $y=4y'$ and we have $-48(y')^2=x'z'$. As $y$ is even, restriction (*) means that $x' \equiv z' \equiv -z'$ (mod $4$) so both $x'$ and $z'$ are even and equivalent modulo $4$. 
Because $-48(y')^2$ is divisible by 16, both $x'$ and $z'$ must be divisible by $4$, so let $x'=4x''$ and $z'=4z''$ for some integers $x''$ and $z''$. Then we have $-3(y')^2=x''z''$, which can be solved with no restrictions, and its integer solutions are 
$(x'',y',z'')=(0,0,u)$ or $\left(v,u,-\frac{3u^2}{v}\right)$ with $u \in \mathbb{Z}$ and $v \in D(3u^2)$. Then in the original variables, we have that the integer solutions to equation \eqref{x2pxypy2mz2} with $y \equiv 0$ (mod $4$) are
 \begin{equation}\label{x2pxypy2mz2_sol_y_4u}
 (x,y,z)=(u,0,u), \quad \text{or} \quad \left( v-\frac{3u^2}{v}-2u,4u,-\frac{3u^2}{v}-v\right), \quad u \in \mathbb{Z}, \quad v \in D(3u^2).
 \end{equation}
 
 We will now consider case (ii), so $y=4y'+2$ for some integer $y'$. Then we have $-3(4y'+2)^2=x'z'$, restriction (*) means that we must have $x' \equiv z' \equiv -z'$ (mod $4$) so we must have $x'$ and $z'$ even and equivalent modulo $4$.
 However, as $-3(4y'+2)^2=x'z'$ is even, both $x'$ and $z'$ must be even, so let $x'=2x''$ and $z'=2z''$ for some integers $x''$ and $z''$. The equation is then reduced to $-3(2y'+1)^2=x''z''$ with the restriction that $x''$ and $z''$ have the same parity. This is not a restriction as $-3(2y'+1)^2 \equiv 1$ (mod $4$), so $x''z'' \equiv 1$ (mod $4$), and so we must have that $x''$ and $z''$ are odd and equivalent modulo $4$.
 Therefore we can solve the equation without restrictions and its integer solutions are $(x'',y',z'')=\left(v,u,-\frac{3(2u+1)^2}{v}\right)$ with $u \in \mathbb{Z}$ and $v \in D(3(2u+1)^2)$. Then, in the original variables, we have that the integer solutions to \eqref{x2pxypy2mz2} with $y \equiv 2$ (mod $4$) are
  \begin{equation}\label{x2pxypy2mz2_sol_y_4up2}
 (x,y,z)=\left(\frac{v}{2} -\frac{3(2u+1)^2}{2v}-2u-1,4u+2,-\frac{3(2u+1)^2}{2v}-\frac{v}{2}\right), \quad u \in \mathbb{Z}, \quad v \in D(3(2u+1)^2).
 \end{equation}
 Condition (*) guarantees that $(x,y,z)$ are always integers.
 
 Finally, we can conclude that all integer solutions to equation \eqref{x2pxypy2mz2} are of the form \eqref{x2pxypy2mz2_sol_y_odd}, \eqref{x2pxypy2mz2_sol_y_4u} or \eqref{x2pxypy2mz2_sol_y_4up2}.

\vspace{10pt}

The next equation we will consider is
\begin{equation}\label{x2ym2ymz2}
x^2y-2y=z^2,
\end{equation}
which can be rewritten as $y(x^2-2)=z^2$. This suggests to make the non-linear substitution $x'=x^2-2$, and the equation reduces to $x'y=z^2$. Up to the names of variables and parameters, we have solved this equation previously and its integer solutions are $(x',y,z)=(rk^2,rw^2,rkw),$ with $r,k,w \in \mathbb{Z}$, which can be found in Table \ref{table1.38}. We then have $x^2-2=x'=rk^2$, $x^2-2=rk^2$ is an equation of the form \eqref{eq:x1kx2} and its solutions are $(x,k,r)=\left(u,v,\frac{u^2-2}{v^2}\right)$ with $u \in \mathbb{Z}$ and $v \in D_2(u^2-2)$. Hence, the integer solutions to the original equation \eqref{x2ym2ymz2} are 
$$
(x,y,z)=\left(u,\frac{w^2(u^2-2)}{v^2}, \frac{w(u^2-2)}{v}\right), \quad u,w \in {\mathbb Z}, \quad v \in D_2(u^2-2).
$$

\vspace{10pt}

The next equation we will consider is
\begin{equation}\label{x2yp2ymz2}
x^2y+2y=z^2,
\end{equation}
which can be rewritten as $y(x^2+2)=z^2$. This suggests to make the non-linear substitution $x'=x^2+2$, and the equation reduces to $x'y=z^2$. Up to the names of variables, we have solved this equation previously and its integer solutions are $(x',y,z)=(rk^2,rw^2,rkw),$ with $r,k,w \in \mathbb{Z}$, which can be found in Table \ref{table1.38}. We then have $x^2+2=x'=rk^2$, the equation $x^2+2=rk^2$ is of the form \eqref{eq:x1kx2} and its integer solutions are $(x,k,r)=\left(u,v,\frac{u^2+2}{v^2}\right)$ with $u \in \mathbb{Z}$ and $v \in D_2(u^2+2)$. Hence, the integer solutions to the original equation \eqref{x2yp2ymz2} are 
$$
(x,y,z)=\left(u,\frac{w^2(u^2+2)}{v^2}, \frac{w(u^2+2)}{v}\right), \quad u,w \in {\mathbb Z}, \quad v \in D_2(u^2+2).
$$

\vspace{10pt}

The next equation we will consider is
\begin{equation}\label{x2yp2xmz2}
x^2y+2x=z^2,
\end{equation}
which can be rewritten as $x(xy+2)=z^2$. This suggests to make the non-linear substitution $y'=xy+2$, and the equation reduces to $xy'=z^2$. Up to the names of variables, we have solved this equation previously and its integer solutions are $(x,y',z)=(rk^2,rw^2,rkw),$ with $r,k,w \in \mathbb{Z}$, which can be found in Table \ref{table1.38}. We can then deduce that $rw^2=rk^2y+2$. As $r$ must divide $2$, $r=\pm 1$ or $\pm 2$, we then have equations of the form \eqref{eq:x1kx2} and their integer solutions are
$$
(w,k,y)= \left(u,v,\frac{u^2-2/r}{v^2}\right) \quad \text{with} \quad r\in \{\pm2, \pm 1\}, \quad u \in \mathbb{Z} \quad \text{and} \quad v \in D_2(u^2-2/r).
$$ 
and
$$
(w,k,y,r)=(\pm 1, 0, u,2), \quad u \in \mathbb{Z}.
$$
Hence, the integer solutions to the original equation \eqref{x2yp2xmz2} are 
$$
(x,y,z)=(0,u,0), \quad \left(rv^2, \frac{u^2 -2/r}{v^2}, ru v \right), \quad r \in \{\pm1, \pm 2 \}\quad u \in \mathbb{Z}, \quad v \in D_2(u^2 -2/r ).
$$

\vspace{10pt}

 The next equation we will consider is
 \begin{equation}\label{x2ypxypz2}
 x^2y+xy+z^2=0,
 \end{equation}
 which can be rewritten as $y(-x^2-x)=z^2$. This suggests to make the substitution $x'=-x^2-x$, which reduces the equation to $x'y=z^2$. Up to the names of variables and parameters, we have solved this equation previously and its integer solutions are $(x',y,z)=(rk^2,rw^2,rkw)$ with $r,k,w \in \mathbb{Z}$, which can be found in Table \ref{table1.38}. Then we have $x'=-x^2-x=rk^2$, which is an equation of the form \eqref{eq:x1kx2} and its integer solutions are $(x,k,r)=(-1,0,u),(0,0,u),\left(u,v,-\frac{u^2+u}{v^2}\right)$ with $u \in \mathbb{Z}$ and $v \in D_2(u^2+u)$. Hence the integer solutions to the original equation \eqref{x2ypxypz2} are
$$
 (x,y,z)=(-1,u,0),\quad (0,u,0),\quad \left(u,-\frac{w^2(u^2+u)}{v^2},-\frac{w(u^2+u)}{v}\right), \quad u,w \in \mathbb{Z}, \quad v \in D_2(u^2+u).
$$

\vspace{10pt}

The final equation we will solve is
\begin{equation}\label{x2ypxzpz2}
x^2y+xz+z^2=0.
\end{equation}
We can solve this as a quadratic in $z$, so we have
$$
z= \frac{-x \pm \sqrt{x^2-4x^2y}}{2}.
$$
In order to have integer solutions, $x^2-4x^2y=x^2(1-4y)$ must be a prefect square. Hence $1-4y$ must be an odd perfect square. So, $1-4y=(2t+1)^2$ for some integer $t$, then $y=-t^2-t$, and we obtain the integer solutions $(x,y,z)=(u,-t^2-t,ut)$, or $(u,-t^2-t,-ut-u)$ with $ u,t \in \mathbb{Z}$. These families cover the same set of integer solutions, so we can conclude that all integer solutions to \eqref{x2ypxzpz2} are 
$$
(x,y,z)=(u,-t^2-t,ut), \quad \text{with} \quad u,t \in \mathbb{Z}.
$$

\begin{center}
\begin{tabular}{ |c|c|c|c|c|c| } 
 \hline
  Equation & Solution $(x,y,z)$ \\ 
 \hline\hline
   $x^2y+y=z^2$ & $\left(u,\frac{w^2(u^2+1)}{v^2}, \frac{w(u^2+1)}{v}\right), \quad u \in {\mathbb Z}, \quad v \in D_2(u^2+1)$  \\\hline
   $x^2y-y=z^2$ &$(\pm 1,u,0), \quad \text{or} \quad \left(u,\frac{w^2(u^2-1)}{v^2}, \frac{w(u^2-1)}{v}\right), \quad  u \in {\mathbb Z}, \quad v \in D_2(u^2-1)$\\\hline
    $x^2y+x=z^2$ & $\quad (0,u,0) \quad \text{or} \quad \left(\pm v^2,\frac{u^2\mp 1}{v^2},vu\right), $ \\ & $ u \in {\mathbb Z}, \quad v \in D_2(u^2\mp 1)$    \\\hline
     $y^2+z^2=x^2+x$ &$\left(\frac{4u^2+1}{4v}+\frac{v-2}{4}, \frac{4u^2+1}{4v}-\frac{v}{4}, u\right), \quad u \in {\mathbb Z}, \quad v \in D(4u^2+1)$ \\\hline
    $y^2+z^2=x^2+x-1$& $\left(\frac{4u^2+5+v^2-2v}{4v},u,\frac{4u^2+5-v^2}{4v}\right), \quad u \in {\mathbb Z}, \quad v \in D(4u^2+5)$\\\hline
     $y^2+z^2=x^2+x+1$&$\left(\frac{4u^2-3+v^2-2v}{4v},u,\frac{4u^2-3-v^2}{4v}\right), \quad u \in {\mathbb Z}, \quad v \in D(4u^2-3)$ \\
\hline
\end{tabular}
\captionof{table}{\label{tab:H16substsol} Integer solutions to the equations listed in Table \ref{tab:H16subst}.}
\end{center} 

\subsection{Exercise 4.8}\label{ex:H16modularroot}
\textbf{\emph{For each equation listed in Table \ref{tab:H16modularroot}, (a) find an integer solution by trial and error, (b) use Proposition \ref{prop:finlincheck} to conclude that it has infinitely many integer solutions, and (c) present all its integer solutions. }}

\begin{center}
\begin{tabular}{ |c|c|c|c|c|c| } 
 \hline
 $H$ & Equation & $H$ & Equation & $H$ & Equation \\ 
 \hline\hline
 $15$ & $x^2 y-y+z^2-1=0$ & $16$ & $x^2y-y+z^2-2=0$ & $16$ & $x^2y-y+z^2+z=0$ \\ 
 \hline
 $15$ & $x^2 y-y+z^2+1=0$ & $16$ & $x^2 y-y+z^2+2=0$ & $16$ & $x^2y+y+z^2+z=0$ \\ 
 \hline
 $15$ & $x^2 y+y+z^2-1=0$ & $16$ & $x^2 y+y+z^2-2=0$ & $16$ & $x^2y+x+z^2+z=0$ \\ 
 \hline
  $15$ & $x^2 y+y+z^2+1=0$ & $16$ & $x^2 y+y+z^2+2=0$ & $16$ & $x^2y+x-y+z^2=0$ \\ 
 \hline
 $15$ & $x^2 y+x+z^2-1=0$ & $16$ & $x^2 y+x+z^2-2=0$ & $16$ & $x^2y+x+y+z^2=0$  \\ 
 \hline
 $15$ & $x^2 y+x+z^2+1=0$ & $16$ & $x^2 y+x+z^2+2=0$ & $16$ & $xyz+x+y+t^2=0$ \\ 
 \hline
\end{tabular}
\captionof{table}{\label{tab:H16modularroot} Equations of the form \eqref{eq:modularroot} of size $H\leq 16$.}
\end{center} 

In this exercise we will describe the set of integer solutions to equations of the form 
\begin{equation}\label{eq:modularroot}
Px_1=Qx_2^2+Rx_2+T
\end{equation}
where $P,Q,R,T$ are polynomials with integer coefficients in other variables $x_3, \dots, x_n$. In order to solve these equations, we will first prove that the equations have infinite solutions, using the following proposition. 
\begin{proposition}\label{prop:finlincheck}[Proposition 4.7 in the book]
Assume that a Diophantine equation can, possibly after permutation of variables, be written in the form
\begin{equation}\label{eq:finlincheck}
	P(x_3,\dots,x_n) x_1 = Q(x_2,\dots,x_n) 
\end{equation}
where $P,Q$ are polynomials with integer coefficients. If \eqref{eq:finlincheck} has an integer solution, then it has infinitely many of them.
\end{proposition}

To describe the integer solutions to equations in this exercise, we will use the following notation
\begin{equation}\label{S Notation}
S(m,a)=\{r \in \mathbb{Z}: |r| \leq |m|/2 \quad \text{and} \quad r^2 \equiv a \, ( \text{mod} \, |m|)\}.
\end{equation}
The solutions to the congruence $r^2 \equiv a \, ( \text{mod} \, |m|)$ can be computed in Mathematica using the command
$$
{\tt PowerModList[a,1/2,Abs[m]]}.
$$

To solve all equations of the form \eqref{eq:modularroot}, we can use the method from Section 4.1.3 of the book, which we summarise below for convenience. The cases $P=0$ and $Q=0$ can be solved separately, so we may assume that $P \neq 0$ and $Q \neq 0$. Multiplying equation \eqref{eq:modularroot} by $4$, we have
$$
4Px_1=4Qx_2^2+4Rx_2+4T=(2Qx_2+R)^2-R^2+4QT.
$$
Then taking variables $x_3,\dots, x_n$ arbitrary, and let $x_2 = Pu+\frac{r-R}{2Q}$, where $u$ is an arbitrary integer and $r$ is an integer such that $r \equiv R \, (\text{mod} \, |2Q|)$ and $r \in S(4QP,R^2-4QT)$, we guarantee that both $x_2$ and
$$
x_1= \frac{(2Qx_2+R)^2-R^2+4QT}{4QP} = \frac{(2QPu+r)^2-R^2+4QT}{4QP} = QPu^2+ur+\frac{r^2-(R^2-4QT)}{4QP}
$$ 
are integers. 

Table \ref{tab:H16modularroot_sol} shows that each equation in Table \ref{tab:H16modularroot} can be written in the form \eqref{eq:finlincheck} and that each equation has an integer solution. Therefore, by Proposition \ref{prop:finlincheck} we can conclude that every equation in Table \ref{tab:H16modularroot_sol} has infinitely many integer solutions. 

\begin{center}
\begin{tabular}{ |c|c|c|c|c|c| } 
 \hline
  Equation & Form \eqref{eq:finlincheck} & Solution $(x,y,z)$ or $(x,y,z,t)$ \\
 \hline\hline
  $x^2 y-y+z^2-1=0$& $y(1-x^2)=z^2-1$ & $(0,0,1)$  \\ \hline
 $x^2 y-y+z^2+1=0$&$y(1-x^2)=z^2+1$ & $(0,1,0)$  \\ \hline
 $x^2 y+y+z^2-1=0$&$y(-x^2-1)=z^2-1$ & $(0,1,0)$ \\ \hline
  $x^2 y+y+z^2+1=0$ & $y(-x^2-1)=z^2+1 $ & $(-1,-1,-1)$  \\ \hline
 $x^2 y+x+z^2-1=0$ & $y(-x^2)=x+z^2-1$ & $(1,0,0)$\\ \hline
 $x^2 y+x+z^2+1=0$ &$y(-x^2)=x+z^2+1$ &$(-1,-1,-1)$ \\ \hline
 $x^2y-y+z^2-2=0$ & $y(1-x^2)=z^2-2 $ &$(0,-1,-1)$ \\ \hline
  $x^2 y-y+z^2+2=0$ & $y(1-x^2)=z^2+2$ & $(-2,-1,-1)$ \\ \hline
 $x^2 y+y+z^2-2=0$ & $y(-1-x^2)=z^2-2 $ &$(1,1,0)$ \\ \hline
  $x^2 y+y+z^2+2=0$ & $y(-x^2-1)=z^2+2$ & $(-1,-1,0)$ \\ \hline
  $x^2 y+x+z^2-2=0$ & $y(-x^2)=x+z^2-2$ & $(1,1,0)$ \\ \hline
  $x^2 y+x+z^2+2=0$ & $y(-x^2)=x+z^2+2$ &$(-1,-2,-1)$ \\ \hline
  $x^2y-y+z^2+z=0$& $y(1-x^2)=z^2+z$ & $(1,-1,-1)$  \\  \hline
  $x^2y+y+z^2+z=0$ & $y(-x^2-1)=z^2+z$ & $(1,-1,-2)$ \\  \hline
   $x^2y+x+z^2+z=0$ & $y(-x^2)=x+z^2+z$ & $(1,-1,-1)$\\    \hline
    $x^2y+x-y+z^2=0$& $y(1-x^2)=z^2+x$ & $(2,-1,-1)$ \\ \hline
  $x^2y+x+y+z^2=0$& $y(-1-x^2)=z^2+x$ & $(1,-1,-1)$  \\  \hline
   $xyz+x+y+t^2=0$& $z(-xy)=x+y+t^2$ & $(-1,-1,1,-1)$ \\
   \hline
\end{tabular}
\captionof{table}{\label{tab:H16modularroot_sol} An integer solution to each equation listed in Table \ref{tab:H16modularroot}.}
\end{center} 

Equation 
$$
x^2y+y+z^2+1=0
$$
is equivalent to $x^2y+y-z^2-1=0$ which is solved in Section 4.1.3 of the book, and its integer solutions are
$$
(x,y,z) = \left(v, -\left((v^2+1)u^2+2ur+\frac{r^2+1}{v^2+1}\right), (v^2+1)u+r \right), \quad u,v \in {\mathbb Z}, \quad r \in S(v^2+1,-1). 
$$

Equation
$$
x^2y+x+z^2+1=0
$$
is solved in Section 4.1.3 of the book, and its integer solutions are
$$
(x,y,z) = \left(u, -\left(u^2v^2+2vr+\frac{r^2+u+1}{u^2}\right),u^2 v + r \right), \quad u,v \in {\mathbb Z}, \quad r \in S(u^2,-u-1). 
$$

The first equation we will consider is
\begin{equation}\label{x2ymypz2m1}
x^2 y-y+z^2-1=0.
\end{equation}
If $x=\pm 1$, we obtain the integer solutions
\begin{equation}\label{x2ymypz2m1_soli}
	(x,y,z)=\left(\pm 1,u, \pm 1 \right), \quad u \in \mathbb{Z}.
\end{equation}
Now assume that $x \neq \pm 1$. We can then rearrange \eqref{x2ymypz2m1} to
$$
y=\frac{1-z^2}{x^2-1}.
$$
Let $x=v\neq \pm1$ be an integer, and let $m=v^2-1$. Then $z=mu+r$ where $u$ is an integer and $|r| \leq \frac{|m|}{2}$ is the integer remainder. Then we have $y=\frac{1-(mu+r)^2}{m} = -mu^2-2ur-\frac{r^2-1}{m}$. In order for $\frac{r^2-1}{m}$ to be integer, we must have $r^2-1 \equiv 0 \text{ (mod $|m|$)}$, or $r \in S(m,1)$. Substituting $m=v^2-1$ we obtain that the set of integer solutions to equation \eqref{x2ymypz2m1} with $x \neq \pm 1$ is
\begin{equation}\label{x2ymypz2m1_sol}
(x,y,z)=\left(v,-u^2(v^2-1)-2ur-\frac{r^2-1}{v^2-1},u(v^2-1)+r\right), \quad u,v \in \mathbb{Z}, \quad v\neq \pm 1, \quad r \in S(v^2-1,1).
\end{equation}
Therefore, we can conclude that the integer solutions to equation \eqref{x2ymypz2p1} are of the form  \eqref{x2ymypz2m1_soli} or \eqref{x2ymypz2m1_sol}.

\vspace{10pt}

The next equation we will consider is
\begin{equation}\label{x2ymypz2p1}
 x^2 y-y+z^2+1=0.
\end{equation} 
This equation has no integer solutions with $x=\pm 1$, hence we may assume that $x \neq \pm 1$. Let 
$x=v\neq \pm1$ be an integer, and $m=v^2-1$. Then $z=mu+r$ where $u$ is an integer and $|r| \leq \frac{|m|}{2}$ is the integer remainder. Then we have $y=-\frac{1+(mu+r)^2}{m} = -mu^2-2ur-\frac{r^2+1}{m}$. In order for $\frac{r^2+1}{m}$ to be integer, we must have $r^2+1 \equiv 0 \text{ (mod $|m|$)}$, or $r \in S(m,-1)$. Substituting $m=v^2-1$, we can conclude that the set of integer solutions to equation \eqref{x2ymypz2p1} is
$$
 (x,y,z)=\left(v,-u^2(v^2-1)-2ur-\frac{r^2+1}{v^2-1},u(v^2-1)+r\right), \quad u,v \in \mathbb{Z}, \quad v \neq \pm 1, \quad r \in S(v^2-1,-1).
$$

\vspace{10pt}

The next equation we will consider is
 \begin{equation}\label{x2ypypz2m1} 
 x^2 y+y+z^2-1=0.
 \end{equation}
 Let $x=v$ be an integer, and $m=v^2+1>0$. Then $z=mu+r$ where $u$ is an integer and $|r| \leq \frac{m}{2}$ is the integer remainder. Then we have, $y=\frac{1-(mu+r)^2}{m} = -mu^2-2ur+\frac{1-r^2}{m}$. In order for $\frac{1-r^2}{m}$ to be integer, we must have $1-r^2 \equiv 0 \text{ (mod $m$)}$, or $r \in S(m,1)$. Substituting $m=v^2+1$, we can conclude that the set of integer solutions to equation \eqref{x2ypypz2m1} is
$$
   (x,y,z)=\left(v,-u^2(v^2+1)-2ur+\frac{1-r^2}{v^2+1},u(v^2+1)+r\right), \quad u,v \in \mathbb{Z}, \quad r \in S(v^2+1,1).
$$
  
  \vspace{10pt}
  
  The next equation we will consider is
 \begin{equation}\label{x2ypxpz2m1} 
 x^2 y+x+z^2-1=0.
 \end{equation} 
  If $x=0$, we obtain the integer solutions 
 \begin{equation}\label{x2ypxpz2m1_soli} 
 	(x,y,z)=\left(0,u, \pm 1 \right),\quad u \in \mathbb{Z}.
 \end{equation}
 Now assume that $x\neq 0$, and let $x=u\neq 0$ be an integer. Then dividing $z$ by $u^2$ with remainder, we can represent $z$ as $z=u^2v+r$ for some integers $v$ and $r$ such that $|r|\leq \frac{u^2}{2}$. Then 
 $$
 y=\frac{1-x-z^2}{x^2}=\frac{1-u-(u^2v+r)^2}{u^2}=-u^2v^2-2vr+\frac{1-u-r^2}{u^2}.
 $$
To ensure $y$ is integer, $\frac{1-u-r^2}{u^2}$ must be integer, so $1-u-r^2 \equiv 0 \text{ (mod $u^2$)}$, or $r \in S(u^2,1-u)$. We then have that the set of integer solutions to equation \eqref{x2ypxpz2m1} with $x\neq 0$ is
   \begin{equation}\label{x2ypxpz2m1_sol} 
 (x,y,z)=\left(u,-u^2v^2-2vr-\frac{r^2+u-1}{u^2},u^2v+r \right),\quad u,v \in \mathbb{Z}, \quad u \neq 0, \quad r \in S(u^2,1-u).
 \end{equation}

  Therefore, we can conclude that the integer solutions to equation \eqref{x2ypxpz2m1} are of the form \eqref{x2ypxpz2m1_soli} or \eqref{x2ypxpz2m1_sol}.
 
 \vspace{10pt}
 
 The next equation we will consider is
\begin{equation}\label{x2ymypz2m2}
 x^2y-y+z^2-2=0.
\end{equation} 
This equation has no solutions with $x = \pm 1$, so we may assume that $x \neq \pm 1$. We can rearrange the equation to
$$
y=\frac{2-z^2}{x^2-1}.
$$
Let $x=v\neq \pm1$ be an integer, and $m=v^2-1$. Then $z=mu+r$ where $u$ is an integer and $|r| \leq \frac{|m|}{2}$ is the integer remainder. Then we have $y=\frac{2-(mu+r)^2}{m} = -mu^2-2ur+\frac{2-r^2}{m}$. In order for $\frac{2-r^2}{m}$ to be integer, we must have $2-r^2 \equiv 0 \text{ (mod $|m|$)}$, or $r \in S(m,2)$. Substituting $m=v^2-1$, we obtain that the set of integer solutions to equation \eqref{x2ymypz2m2} is
$$
(x,y,z)=\left(v,-u^2(v^2-1)-2ur+\frac{2-r^2}{v^2-1},u(v^2-1)+r\right), \quad u,v \in \mathbb{Z}, \quad v\neq \pm1, \quad r \in S(v^2-1,2).
$$
 
 \vspace{10pt}
 
 The next equation we will consider is
\begin{equation}\label{x2ymypz2p2} 
  x^2 y-y+z^2+2=0.
 \end{equation} 
This equation has no solutions with $x = \pm 1$, so we can assume that $x \neq \pm 1$. We can rearrange the equation to
$$
y=-\frac{2+z^2}{x^2-1}.
$$
Let $x=v\neq \pm1$ be an integer, and $m=v^2-1$. Then $z=mu+r$ where $u$ is an integer and $|r| \leq \frac{|m|}{2}$ is the integer remainder. Then we have $y=-\frac{2+z^2}{x^2-1}=-\frac{2+(mu+r)^2}{m} = -mu^2-2ur-\frac{2+r^2}{m}$. In order for $\frac{2+r^2}{m}$ to be integer, we must have $2+r^2 \equiv 0 \text{ (mod $|m|$)}$, or $r \in S(m,-2)$. Substituting $m=v^2-1$, we obtain that the set of integer solutions to equation \eqref{x2ymypz2p2} is
$$
(x,y,z)=\left(v,-u^2(v^2-1)-2ur-\frac{r^2+2}{v^2-1},u(v^2-1)+r\right), \quad u,v \in \mathbb{Z}, \quad v\neq \pm 1, \quad r \in S(v^2-1,-2).
$$
  
  \vspace{10pt}
  
   The next equation we will consider is
\begin{equation}\label{x2ypypz2m2} 
 x^2 y+y+z^2-2=0.
 \end{equation} 
Let $x=v$ be an integer, and $m=v^2+1>0$. Then $z=mu+r$ where $u$ is an integer and $|r| \leq \frac{m}{2}$ is the integer remainder. Then we have, $y=\frac{2-z^2}{x^2+1}=\frac{2-(mu+r)^2}{m} = -mu^2-2ur+\frac{2-r^2}{m}$. In order for $\frac{2-r^2}{m}$ to be integer, we must have $2-r^2 \equiv 0 \text{ (mod $m$)}$, or $r \in S(m,2)$. Substituting $m=v^2+1$, we obtain that the set of integer solutions to equation \eqref{x2ypypz2m2} is
$$
 (x,y,z)=\left(v,-u^2(v^2+1)-2ur+\frac{2-r^2}{v^2+1},u(v^2+1)+r\right), \quad u,v \in \mathbb{Z}, \quad r \in S(v^2+1,2).
$$
 
 \vspace{10pt}
 
   The next equation we will consider is
\begin{equation}\label{x2ypypz2p2}  
 x^2 y+y+z^2+2=0.
  \end{equation} 
  Let $x=v$ be an integer, and $m=v^2+1>0$. Then $z=mu+r$ where $u$ is an integer and $|r| \leq \frac{m}{2}$ is the integer remainder. Then we have $y=-\frac{z^2+2}{x^2+1}=-\frac{2+(mu+r)^2}{m} = -mu^2-2ur-\frac{2+r^2}{m}$. To ensure that $\frac{2+r^2}{m}$ is an integer, we must have $2+r^2 \equiv 0 \text{ (mod $m$)}$, or $r \in S(m,-2)$. Substituting $m=v^2+1$, we obtain that the set of integer solutions to equation \eqref{x2ypypz2p2} is
$$
 (x,y,z)=\left(v,-\left(u^2(v^2+1)+2ur+\frac{r^2+2}{v^2+1}\right),u(v^2+1)+r\right), \quad u,v \in \mathbb{Z}, \quad r \in S(v^2+1,-2).
$$
 
 \vspace{10pt}
 
  The next equation we will consider is
 \begin{equation}\label{x2ypxpz2m2}  
 x^2 y+x+z^2-2=0.
 \end{equation}
 This equation has no integer solutions with $x=0$, so we can assume that $x\neq 0$. Let $x=u\neq 0$ be an integer. Then dividing $z$ by $u^2$ with remainder, we can represent $z$ as $z=u^2v+r$ for some integer $r$ such that $|r|\leq \frac{u^2}{2}$. Then 
 $$
 y=\frac{2-x-z^2}{x^2}=\frac{2-u-(u^2v+r)^2}{u^2}=-u^2v^2-2vr+\frac{2-u-r^2}{u^2}.
 $$
 To ensure $y$ is an integer, $\frac{2-u-r^2}{u^2}$ must be an integer, so $2-u-r^2 \equiv 0 \text{ (mod $u^2$)}$, or $r \in S(u^2,2-u)$. We can then conclude that the set of integer solutions to equation \eqref{x2ypxpz2m2} is
   $$
 (x,y,z)=\left(u,-u^2v^2-2vr-\frac{r^2+u-2}{u^2},u^2v+r \right),\quad  u,v \in \mathbb{Z}, \quad u \neq 0, \quad r \in S(u^2,2-u).
$$

\vspace{10pt}
  
   The next equation we will consider is
 \begin{equation}\label{x2ypxpz2p2}   
  x^2 y+x+z^2+2=0.
   \end{equation}
   This equation has no integer solutions when $x=0$, so we may assume that $x\neq 0$. Let  $x=u\neq 0$ be an integer. Then dividing $z$ by $u^2$ with remainder, we can represent $z$ as $z=u^2v+r$ for some integer $r$ such that $|r|\leq \frac{u^2}{2}$. Then 
 $$
 y=-\frac{2+x+z^2}{x^2}=-\frac{2+u+(u^2v+r)^2}{u^2}=-u^2v^2-2vr-\frac{2+u+r^2}{u^2}.
 $$
 To ensure $y$ is integer, $\frac{2+u+r^2}{u^2}$ must be an integer, so $2+u+r^2 \equiv 0 \text{ (mod $u^2$)}$, or $r \in S(u^2,-2-u)$. We can then conclude that the set of integer solutions to equation \eqref{x2ypxpz2p2} is
$$
 (x,y,z)=\left(u,-\left(u^2v^2+2vr+\frac{r^2+u+2}{u^2}\right),u^2v+r \right),\quad  u,v \in \mathbb{Z}, \quad u \neq 0, \quad r \in S(u^2,-u-2).
$$
  
  \vspace{10pt}

    The next equation we will consider is
 \begin{equation}\label{x2ymypz2pz}   
x^2y-y+z^2+z=0.
   \end{equation}
   If $x=\pm 1$ then we obtain the integer solutions
       \begin{equation}\label{x2ymypz2pz_solii}  
   	(x,y,z)= (\pm 1,u,-1), (\pm 1,u,0), \quad u \in \mathbb{Z}.
   \end{equation}
   Now assume that $x \neq \pm 1$. 
   We can rearrange \eqref{x2ymypz2pz} to $y(1-x^2)=z^2+z$, which is an equation of the form \eqref{eq:modularroot}, with $P=1-x^2$, $x_1=y$, $x_2=z$, $Q=1=R$ and $T=0$.  
   Then let $x=v \neq \pm1$ be integer, so $P=1-v^2$, and we obtain
   $$
   z=x_2=Pu+\frac{r-R}{2Q}=u(1-v^2)+\frac{r-1}{2}
   $$
  where $u$ is an integer and $r$ is an integer such that $r \equiv 1 \, (\text{mod} \, 2)$ and $r \in S(4(1-v^2),1)$. Then 
  $$
  y=x_1=QPu^2+ur+\frac{r^2-(R^2-4QT)}{4QP}=u^2(1-v^2)+ur+\frac{r^2-1}{4(1-v^2)}.
  $$
  Therefore, we have that the set of integer solutions to equation \eqref{x2ymypz2pz} with $x\neq \pm 1$ is
  \begin{equation}\label{x2ymypz2pz_soli}  
  \begin{aligned}
  (x,y,z)= \left(v,u^2(1-v^2)+ur+\frac{r^2-1}{4(1-v^2)},u(1-v^2)+\frac{r-1}{2}\right),   u,v \in \mathbb{Z},\\  v \neq \pm 1, \quad r \equiv 1\, \text{(mod $2$)}, \quad r \in S(4(1-v^2),1).
  \end{aligned}
  \end{equation}

   Finally, we can conclude that the integer solutions to equation \eqref{x2ymypz2pz} are of the form \eqref{x2ymypz2pz_solii} or \eqref{x2ymypz2pz_soli}.
 
\vspace{10pt}
  
    The next equation we will consider is
 \begin{equation}\label{x2ypypz2pz} 
 x^2y+y+z^2+z=0.
    \end{equation}
    We can rearrange this equation to $y(-1-x^2)=z^2+z$, which is an equation of the form \eqref{eq:modularroot}, with $P=-1-x^2$, $x_1=y$, $x_2=z$, $Q=1=R$ and $T=0$. Then, let $x=v$ be an integer, so $P=-1-v^2$, and we obtain
    $$
    z=x_2=Pu+\frac{r-R}{2Q}=u(-1-v^2)+\frac{r-1}{2}
    $$
     where $u$ is an integer and $r$ is an integer such that $r \equiv 1 \, (\text{mod} \, 2)$ and $r \in S(4(-1-v^2),1)$, and 
  $$
  y=x_1=QPu^2+ur+\frac{r^2-(R^2-4QT)}{4QP}=u^2(-1-v^2)+ur+\frac{r^2-1}{4(-1-v^2)}.
  $$
 Therefore, we can conclude that the set of integer solutions to equation \eqref{x2ypypz2pz} is
  $$
  \begin{aligned}
  (x,y,z)= \left(v,-u^2(1+v^2)+ur-\frac{r^2-1}{4(1+v^2)},-u(1+v^2)+\frac{r-1}{2}\right),   u,v \in \mathbb{Z},\\  r \equiv 1\, \text{(mod $2$)}, \quad r \in S(-4(1+v^2),1).
  \end{aligned}
 $$

 \vspace{10pt} 

    The next equation we will consider is
 \begin{equation}\label{x2ypxpz2pz} 
 x^2y+x+z^2+z=0.
   \end{equation}
   If $x=0$ then we obtain the integer solutions
     \begin{equation}\label{x2ypxpz2pz_solii}  
   	(x,y,z)= (0,u,-1), (0,u,0), \quad u \in \mathbb{Z}.
   \end{equation}
   Now assume that $x \neq 0$. We can rearrange \eqref{x2ypxpz2pz} to $y(-x^2)=z^2+z+x$, which is an equation of the form \eqref{eq:modularroot}, with $P=-x^2$, $x_1=y$, $x_2=z$, $Q=1=R$ and $T=x$. Then, let $x=v \neq 0$ be an integer, so $P=-v^2$, and we obtain
    $$
    z=x_2=Pu+\frac{r-R}{2Q}=-uv^2+\frac{r-1}{2}
    $$
where $u$ is an integer and $r$ is an integer such that $r \equiv 1 \, (\text{mod} \, 2)$ and $r \in S(-4v^2,1-4v)$, and 
  $$
  y=x_1=QPu^2+ur+\frac{r^2-(R^2-4QT)}{4QP}=u^2(-v^2)+ur+\frac{r^2-1+4v}{-4v^2}.
  $$
Therefore, we have the set of integer solutions to equation \eqref{x2ypxpz2pz} with $x \neq 0$ is
  \begin{equation}\label{x2ypxpz2pz_soli}  
  \begin{aligned}
  (x,y,z)= \left(v,-u^2 v^2+ur-\frac{r^2-1+4v}{4v^2},-uv^2+\frac{r-1}{2}\right),   u,v \in \mathbb{Z},\\ v\neq0, \quad r \equiv 1\, \text{(mod $2$)}, \quad r \in S(-4v^2,1-4v).
  \end{aligned}
  \end{equation}

   Therefore, we can conclude that the integer solutions to equation \eqref{x2ypxpz2pz} are of the form \eqref{x2ypxpz2pz_solii} or \eqref{x2ypxpz2pz_soli}.

 \vspace{10pt}
 
  The next equation we will consider is
 \begin{equation}\label{x2ypxmypz2}
 x^2y+x-y+z^2=0.
 \end{equation}
 If $x = \pm 1$ then we obtain the integer solutions
   \begin{equation}\label{x2ypxmypz2_soli}
 	(x,y,z)=(-1,u,\pm 1), \quad u \in \mathbb{Z}.
 \end{equation}
Now assume that $x \neq \pm 1$. We can then express \eqref{x2ypxmypz2} as 
 $$
 y=-\frac{x+z^2}{x^2-1}.
 $$
 Then let $x=u\neq \pm 1$ be an integer, and $m=u^2-1$. Then $z=mv+r$ such that $r$ is an integer where $|r| \leq \frac{|m|}{2}$. Then 
  $$
 y=-\frac{u+(mv+r)^2}{m}=-v^2m-2vr-\frac{u+r^2}{m}.
 $$
 To ensure that $y$ is integer, we must have that $\frac{u+r^2}{m}$ is integer. So $u+r^2\equiv 0 \text{ (mod $|m|$)}$, or $r \in S(m,-u)$. Then substituting $m=u^2-1$, we obtain that the set of integer solutions to equation \eqref{x2ypxmypz2} with $x \neq \pm 1$ is
  \begin{equation}\label{x2ypxmypz2_sol}
(x,y,z)=\left(u,-v^2(u^2-1)-2vr-\frac{u+r^2}{u^2-1},v(u^2-1)+r\right),\quad  u,v \in \mathbb{Z}, \quad u \neq \pm 1, \quad r \in S(u^2-1,-u).
 \end{equation}

 Therefore we can conclude that integer solutions to equation \eqref{x2ypxmypz2} are of the form \eqref{x2ypxmypz2_soli} or \eqref{x2ypxmypz2_sol}.
 
 \vspace{10pt}
 
   The next equation we will consider is
 \begin{equation}\label{x2ypxpypz2}
 x^2y+x+y+z^2=0.
  \end{equation}
Let $x=u$ be integer, and $m=u^2+1>0$. Then $z=mv+r$ such that $r$ is an integer where $|r| \leq \frac{m}{2}$. Then 
  $$
 y=-\frac{z^2+x}{x^2+1}=-\frac{u+(mv+r)^2}{m}=-v^2m-2vr-\frac{u+r^2}{m}.
 $$
 To ensure that $y$ is integer, we must have that $\frac{u+r^2}{m}$ is integer. So $u+r^2\equiv 0 \text{ (mod $m$)}$, or $r \in S(m,-u)$. Then substituting $m=u^2+1$, we obtain that the set of integer solutions to equation \eqref{x2ypxpypz2} is
  $$
(x,y,z)=\left(u,-v^2(u^2+1)-2vr-\frac{u+r^2}{u^2+1},v(u^2+1)+r\right),\quad  u,v \in \mathbb{Z}, \quad r \in S(u^2+1,-u).
$$
 
 \vspace{10pt}
 
 The final equation we will consider is
 \begin{equation}\label{xyzpxpypt2}  
 xyz+x+y+t^2=0.
 \end{equation}
 If $xy=0$ then we obtain the integer solutions 
  \begin{equation}\label{xyzpxpypt2_solii}  
 	(x,y,z,t)=(0,-u^2,v,u) \quad \text{or} \quad (-u^2,0,v,u), \quad u,v \in \mathbb{Z}.
 \end{equation}
 Now we may assume that $xy \neq 0$. 
 We can rearrange \eqref{xyzpxpypt2} to $z(-xy)=x+y+t^2$, which is an equation of the form \eqref{eq:modularroot}, with $P=-xy$, $x_1=z$, $x_2=t$, $Q=1$, $R=0$ and $T=x+y$. Then, let $x=v$ and $y=w$ be integers. So $P=-vw$, and we obtain
   $$
    t=x_2=Pu+\frac{r-R}{2Q}=-uvw+\frac{r}{2}
    $$
where $u$ is an integer and $r$ is an integer such that $r \equiv 0 \, (\text{mod} \, 2)$ and $r \in S(-4vw,-4v-4w)$, and 
  $$
  z=x_1=QPu^2+ur+\frac{r^2-(R^2-4QT)}{4QP}=u^2(-vw)+ur+\frac{r^2+4(v+w)}{-4vw}.
  $$
Therefore, we have that the set of integer solutions to equation \eqref{xyzpxpypt2} with $xy\neq 0$ is
 \begin{equation}\label{xyzpxpypt2_soli}
 \begin{aligned}
 (x,y,z,t)=\left(v,w,-vwu^2+ur-\frac{r^2+4v+4w}{4vw},-uvw+\frac{r}{2}\right),\quad u,v,w \in \mathbb{Z}, 
 \\  r \equiv 0\, \text{(mod $2$)}, \quad vw\neq 0, \quad r \in S(-4vw,-4v-4w)
  \end{aligned}
  \end{equation}
   Therefore we can conclude that integer solutions to equation \eqref{xyzpxpypt2} are of the form \eqref{xyzpxpypt2_solii} or \eqref{xyzpxpypt2_soli}.

\begin{center}
\begin{tabular}{ |c|c|c|c|c|c| } 
 \hline
 Equation & Solution $(x,y,z)$ or $(x,y,z,t)$ \\
 \hline\hline
  $x^2 y-y+z^2-1=0$ & $\left(v,-(v^2-1)u^2-2ur-\frac{r^2-1}{v^2-1},(v^2-1)u+r\right)$, \\
 & $ u,v \in \mathbb{Z}, \quad v\neq \pm 1, \quad r \in S(v^2-1,1)$, \\
 & $\left(\pm 1,u, \pm 1 \right), \quad u \in \mathbb{Z}$ \\ \hline
   $x^2 y-y+z^2+1=0$ & $\left(v,-u^2(v^2-1)-2ur-\frac{r^2+1}{v^2-1},u(v^2-1)+r\right)$, \\
  & $ u,v \in \mathbb{Z}, \quad v \neq \pm1, \quad r \in S(v^2-1,-1)$  \\  \hline
 $x^2 y+y+z^2-1=0$& $\left(v,-u^2(v^2+1)-2ur+\frac{1-r^2}{v^2+1},u(v^2+1)+r\right)$, \\
 & $u,v \in \mathbb{Z}, \quad r \in S(v^2+1,1)$\\
\hline
\end{tabular}
\captionof{table}{\label{tab:H16modularrootsol} Integer solutions to the equations listed in Table \ref{tab:H16modularroot}.}
\end{center} 

\subsection{Exercise 4.10}\label{ex:x3py2mz2}
\textbf{\emph{Write down all integer solutions to equation 
\begin{equation}\label{eq:x3py2mz2}
y^2-z^2=x^3.
\end{equation}
as a finite union of polynomial families.  }}

We can rewrite equation \eqref{eq:x3py2mz2} as 
$$
x^3=(y-z)(y+z)
$$
which suggests to make the substitutions $y'=y-z$ and $z'=y+z$. Then the equation reduces to $x^3=y'z'$ with the restriction that $y'$ and $z'$ have the same parity. Up to the names of variables and parameters, we have solved this equation previously and its integer solutions are given by
$$
(x,y',z')=(v_1v_2v_3v_4, v_2^3 v_3^2 v_4, v_1^3 v_3 v_4^2)
$$
 for some integers $v_1,v_2,v_3$ and $v_4$, which can be found in Table \ref{table1.38}. 
 
 Because $y'$ and $z'$ must have the same parity, we have integer solutions to \eqref{eq:x3py2mz2} when either 
(i) $v_3$ is even, or (ii) $v_4$ is even, or (iii) $v_1$ and $v_2$ have the same parity. 
Hence, in case (i), let $v_3=2u_3$ for some integer $u_3$, and relabel $v_1,v_2$ and $v_4$ as $u_1,u_2$ and $u_4$ respectively, and similarly, in case (ii), let $v_4=2u_4$ for some integer $u_4$, and relabel $v_1,v_2$ and $v_3$ as $u_1,u_2$ and $u_3$ respectively, while in case (iii) we have $v_1=u_1+u_2$ and $v_2=u_1-u_2$ for some integers $u_1$ and $u_2$, and relabel $v_3$ and $v_4$ as $u_3$ and $u_4$ respectively. Then because $y=\frac{y'+z'}{2}$ and $z=\frac{z'-y'}{2}$, we obtain that the integer solutions to equation \eqref{eq:x3py2mz2} are
$$
	(x,y,z)=\left(2u_1u_2u_3u_4, \,\, 2u_2^3 u_3^2 u_4+u_1^3 u_3 u_4^2, \,\, u_1^3 u_3 u_4^2-2u_2^3 u_3^2 u_4\right), \quad u_1,u_2,u_3,u_4 \in \mathbb{Z},
$$
and
$$
	(x,y,z)=\left(2u_1u_2u_3u_4, \,\, u_2^3 u_3^2 u_4+2u_1^3 u_3 u_4^2, \,\, 2u_1^3 u_3 u_4^2-u_2^3 u_3^2 u_4\right), \quad u_1,u_2,u_3,u_4 \in \mathbb{Z},
$$
and
$$
\begin{aligned}
x= & (u_1+u_2)(u_1-u_2)u_3u_4, \\
y= & \frac{(u_1-u_2)^3 u_3^2 u_4 + (u_1+u_2)^3 u_3 u_4^2}{2}, \\
z= & \frac{(u_1-u_2)^3 u_3^2 u_4 - (u_1+u_2)^3 u_3 u_4^2}{2},  \quad u_1,u_2,u_3,u_4 \in \mathbb{Z}.
\end{aligned}
$$
It is easy to check that $y$ and $z$ are always integers.

\subsection{Exercise 4.12}\label{ex:x4py2pz2}
\textbf{\emph{Use the following proposition to describe all integer solutions to equation 
\begin{equation}\label{eq:y2px2ypz2}
y^2+x^2y+z^2=0.
\end{equation}}}

\begin{proposition}\label{prop:x4py2pz2}[Proposition 4.11 in the book]
The integer solutions to 
\begin{equation}\label{eq:x4py2pz2}
x^4=y^2+z^2
\end{equation}
 are, up to exchange of $y$ and $z$,  given by
\begin{equation}\label{eq:x4py2pz2sol}
\begin{split}
& x =  s(u^2+v^2)(w^2+r^2), \\
& y= s^2(u^2+v^2)((u^2-v^2)(r^4-6r^2w^2+w^4)+2uv(4r^3 w-4r w^3)), \\
& z= s^2(u^2+v^2)(2uv(r^4-6r^2w^2+w^4)+(u^2-v^2)(-4r^3 w+4r w^3)),
\end{split}
\end{equation}
where $u,v,w,r,s$ are integer parameters.
\end{proposition}

Equation \eqref{eq:y2px2ypz2} can be solved as a quadratic in $y$, hence
$$
y=\frac{-x^2 \pm \sqrt{x^4-4z^2}}{2}.
$$
In order to have integer solutions to \eqref{eq:y2px2ypz2}, we must have 
$$
y=\frac{-x^2+k}{2},
$$
where $k$ is an integer such that $x^4-4z^2=k^2$. After making the substitution $z'=2z$ and rearranging, we have $x^4=k^2+(z')^2$, with the restriction that $z'$ is even, which, up to the names of variables, is equation \eqref{eq:x4py2pz2}, whose solutions $(x,k,z')$ are given by \eqref{eq:x4py2pz2sol}, up to exchange of $k$ and $z'$. 
Then, in the original variables, we obtain that all integer solutions to \eqref{eq:y2px2ypz2} can be described by
\begin{equation}\label{y2px2ypz2_sol_xyz+}
\begin{aligned}
x= & \,\, s(u^2+v^2)(w^2+r^2), \\ y= & -s^2(u^2+v^2) (r^2 u + 2 r v w - u w^2)^2, \\ z= & \,\, s^2(u^2+v^2)(uv(r^4-6r^2w^2+w^4)+(u^2-v^2)(2r w^3-2r^3 w)),
\end{aligned}
\end{equation}
and
\begin{equation}\label{y2px2ypz2_sol_generic_ii}
\begin{aligned}
x = & \,\, s(u^2+v^2)(w^2+r^2), \\
y= & - \frac{s^2(u^2+v^2)(r^2(u-v)+2rw(u+v)+w^2(v-u))^2}{2} \\
z= & \,\, \frac{s^2(u^2+v^2)((u^2-v^2)(r^4-6r^2w^2+w^4)+2uv(4r^3 w-4r w^3))}{2},
\end{aligned}
\end{equation}
where $s,r,u,v,$ and $w$ are integer parameters. We now need to ensure these solutions are integers. Solution \eqref{y2px2ypz2_sol_xyz+} is always integer, while solution \eqref{y2px2ypz2_sol_generic_ii} is integer if and only if the numerators of both $y$ and $z$ are even.  
Because $2uv(4r^3 w-4r w^3)$ is even, to enusre $z$ is even, we must have that
$s^2(u^2+v^2)(u^2-v^2)(r^4-6r^2w^2+w^4)$ is even. Hence, we have three cases to consider, (i) $s$ is even, (ii) $u^4-v^4$ is even, or (iii) $r^4+w^4$ is even. 

In case (i), let $s=2t$ for some integer $t$. This reduces \eqref{y2px2ypz2_sol_generic_ii} to 
\begin{equation}\label{y2px2ypz2_sol_si}
\begin{aligned}
x= & \,\, 2t(u^2+v^2)(w^2+r^2), \\ y= & -2t^2(u^2+v^2)(((v-u) (w^2-r^2)+2 rw (u + v) )^2), \\ z= &  \,\, 2t^2(u^2+v^2)((u^2-v^2)(r^4-6r^2w^2+w^4)+8uvrw(r^2-w^2)),
\end{aligned}
\end{equation}
for some integer parameters $t,r,u,v,$ and $w$. 

In case (ii), we have that $u^4-v^4$ is even, hence, we must have that $u$ and $v$ have the same parity, then let $u=m-n$ and $v=m+n$ for some integers $m$ and $n$. This reduces \eqref{y2px2ypz2_sol_generic_ii} to
\begin{equation}\label{y2px2ypz2_sol_uv-}
\begin{aligned}
x =& \,\, 2s (m^2 + n^2) (r^2 + w^2), \\
y= & -4s^2 (m^2 + n^2)(2 m r w + n (w^2-r^2))^2, \\
z= & \,\, 4 s^2 ((2 m^4 r w - 2 n^4 r w) (r^2 - w^2) - (m^3 n + m n^3) (r^4 -  6 r^2 w^2 + w^4)),
\end{aligned}
\end{equation}
for some integer parameters $m,n,s,r,$ and $w$. 

Let us now consider case (iii), In order to have $r^4+w^4$ even, we must have that $w$ and $r$ have the same parity. Then let $w=m-n$ and $r=m+n$ for some integers $m$ and $n$. This reduces \eqref{y2px2ypz2_sol_generic_ii} to
\begin{equation}\label{y2px2ypz2_sol_rw-}
\begin{aligned}
x = & \,\, 2s (m^2 + n^2)  (u^2 + v^2), \\
y= & -2 s^2 (u^2 + v^2) ( (u - v)(m^2-n^2) - 2 m n (u + v))^2, \\
z= & \,\, 2s^2 (u^2 + v^2) ((v^2 - u^2) (m^4 - 6 m^2 n^2 + n^4) +8 m n u v(m^2 - n^2)),
\end{aligned}
\end{equation}
for some integer parameters $m,n,s,u,$ and $v$.  After making the change of parameters $s \to t$, $m \to r$, $n \to w$, $v \to -u$ and $u \to -v$ to \eqref{y2px2ypz2_sol_rw-}, we obtain solution \eqref{y2px2ypz2_sol_si}. Hence, solution \eqref{y2px2ypz2_sol_si} is redundant. In all cases, we can easily see that $y$ is integer.

In conclusion, we obtain that all integer solutions to equation \eqref{eq:y2px2ypz2} are of the form \eqref{y2px2ypz2_sol_xyz+},  \eqref{y2px2ypz2_sol_uv-} or \eqref{y2px2ypz2_sol_rw-}.

\subsection{Exercise 4.13}\label{ex:y2px2ymz2}
\textbf{\emph{Solve equation 
		\begin{equation}\label{eq:y2mz2mx4}
			y^2-z^2=x^4.
		\end{equation}
		 Use it to describe all integer solutions to the equation 
\begin{equation}\label{eq:y2px2ymz2}
y^2+x^2y-z^2=0.
\end{equation}}}

Equation \eqref{eq:y2mz2mx4} can be reduced to 
$$
x^4=YZ,
$$
where $Y=y-z$ and $Z=y+z$.
We can solve this equation using the method in Section \ref{ex:2mona}, which gives the solution 
\begin{equation}\label{x4myz_sol}
(x,Y,Z)=\left(u_1 u_2 u_3 u_4 u_5, \, \, u_2 u_3^2 u_4^3 u_5^4, \, \, u_1^4 u_2^3 u_3^2 u_4 \right), \quad \text{with} \quad u_1, u_2, u_3, u_4, u_5 \in \mathbb{Z}.
\end{equation}
Then, by using $y=\frac{Y+Z}{2}$ and $z=\frac{Z-Y}{2}$, we obtain, in the original variables, that all integer solutions to \eqref{eq:y2mz2mx4} can be described by
$$
(x,y,z)=\left(u_1 u_2 u_3 u_4 u_5, \frac{u_2 u_3^2 u_4 (u_1^4 u_2^2 + u_4^2 u_5^4)}{2}, \frac{u_2 u_3^2 u_4 (u_1^4 u_2^2 - u_4^2 u_5^4)}{2}\right), \quad \text{with} \quad u_1, u_2, u_3, u_4, u_5 \in \mathbb{Z},
$$
with parity restrictions on the parameters to ensure that $x,y,z$ are integers.

Let us now solve equation \eqref{eq:y2px2ymz2}.

Let us first assume that $x$ is even, so we can write $x=2x'$ for some integer $x'$. Then 
$$
y=\frac{-4(x')^2 \pm \sqrt{16(x')^4+4z^2}}{2}=-2(x')^2\pm \sqrt{4(x')^4+z^2}.
$$
In order for $y$ to be integer, we must have that $4(x')^4+z^2$ is a perfect square, so $4(x')^4+z^2=k^2$ for some integer $k$. We can rearrange this equation to
$$
4(x')^4=k^2-z^2=(k-z)(k+z).
$$
The left-hand side of this equation is always even, and as $k-z$ and $k+z$ have the same parity, they must both be even, so we can make the change of variables $K=\frac{k-z}{2}$ and $T=\frac{k+z}{2}$, which reduces the equation to $(x')^4=KT$, which has the integer solutions $(x',K,T)=(X,Y,Z)$ where $X,Y,Z$ are given by \eqref{x4myz_sol}. Then, in the original variables, we have the integer solutions to \eqref{eq:y2px2ymz2} where $x$ is even,
\begin{equation}\label{y2px2ymz2_sol_even}
\begin{aligned}
(x,y,z)=& (2u_1 u_2 u_3 u_4 u_5, -2 u_1^2 u_2^2 u_3^2 u_4^2 u_5^2 \pm (u_2 u_3^2 u_4^3 u_5^4 + u_1^4 u_2^3 u_3^2 u_4), u_1^4 u_2^3 u_3^2 u_4 - u_2 u_3^2 u_4^3 u_5^4), \\ & \quad u_1, u_2, u_3, u_4, u_5 \in \mathbb{Z}.
\end{aligned}
\end{equation}

We will now consider the case where $x$ is odd. 
Solving \eqref{eq:y2px2ymz2} as a quadratic in $y$, we have
$$
y=\frac{-x^2 \pm \sqrt{x^4+4z^2}}{2}.
$$
In order to have integer solutions, we must have that $x^4+4z^2$ is a perfect square. So, $x^4+4z^2=k^2$ for some integer $k$. After making the substitution $z'=2z$ and rearranging, we have 
$$
x^4=k^2-(z')^2=(k-z')(k+z'),
$$
with the restriction that $z'$ is even. This suggests to make the substitutions $K=k-z'$ and $T=k+z'$, which then reduces the equation to $x^4=KT$, with solution $(x,K,T)=(X,Y,Z)$ where $X,Y,Z$ are given by \eqref{x4myz_sol}, we also have the restriction that $z=\frac{T-K}{4}$ and $k=\frac{T+K}{2}$ are integers. Then in the original variables we have
$$
(x,y,z)=\left(u_1 u_2 u_3 u_4 u_5, \, \, \frac{-2u_1^2 u_2^2 u_3^2 u_4^2 u_5^2 \pm (u_2 u_3^2 u_4^3 u_5^4 + u_1^4 u_2^3 u_3^2 u_4)}{4}, \, \, \frac{u_1^4 u_2^3 u_3^2 u_4-u_2 u_3^2 u_4^3 u_5^4}{4}\right)
$$

As $x$ is odd, we must have $u_1, u_2, u_3, u_4$ and $u_5$ odd. Then let $u_1=2a+1$, $u_2=2b+1$, $u_1=2c+1$, $u_1=2d+1$, and $u_1=2e+1$, for some integers $a,b,c,d$ and $e$. Then we have the integer solutions to \eqref{eq:y2px2ymz2} where $x$ is odd,
\begin{equation}\label{y2px2ymz2_sol_odd}
\begin{aligned}
x= &(2a+1) (2b+1) (2c+1) (2d+1) (2e+1), \\
y= & (2b+1)(2c+1)^2(2d+1)\left( \frac{-2 (2a+1)^2 (2b+1) (2d+1) (2e+1)^2}{4} \right. \\ & \pm \left. \frac{(2d+1)^2 (2e+1)^4 + (2a+1)^4 (2b+1)^2}{4} \right), \\
z= & (2b+1)(2c+1)^2(2d+1)\left( \frac{(2a+1)^4 (2b+1)^2 - (2d+1)^2 (2e+1)^4)}{4} \right), \\
&  \quad a,b,c,d,e \in \mathbb{Z}.
\end{aligned}
\end{equation}
After modulo $4$ analysis, we can easily see that $y$ and $z$ are integer.

Therefore, we can conclude that all integer solutions to \eqref{eq:y2px2ymz2} are of the form  \eqref{y2px2ymz2_sol_even} or \eqref{y2px2ymz2_sol_odd}.

\subsection{Exercise 4.14}\label{ex:H18vieta}
\textbf{\emph{Solve all equations listed in Table \ref{tab:H18vieta}. }}

\begin{center}
\begin{tabular}{ |c|c|c|c|c|c| } 
 \hline
 $H$ & Equation & $H$ & Equation & $H$ & Equation \\ 
 \hline\hline
 $16$ & $x^2-xyz+y+z=0$ & $17$ & $x^2-xyz+y^2+1=0$ & $18$ & $x^2-xyz-y^2+2=0$ \\ 
 \hline
 $16$ & $x^2-xyz-y^2=0$ & $18$ & $x^2-xyz+y+z-2=0$ & $18$ & $x^2-xyz-y^2+x=0$ \\ 
 \hline
 $16$ & $x^2-xyz+y^2=0$ & $18$ & $x^2-xyz+y+z+2=0$ & $18$ & $x^2-xyz-y^2-z=0$ \\ 
 \hline
 $17$ & $x^2-xyz+y+z-1=0$ & $18$ & $x^2-xyz+x+y-z=0$ & $18$ & $x^2-xyz+y^2-2=0$ \\ 
 \hline
 $17$ & $x^2-xyz+y+z+1=0$ & $18$ & $x^2-xyz+x+y+z=0$ & $18$ & $x^2-xyz+y^2+2=0$ \\ 
 \hline
 $17$ & $x^2-xyz-y^2+1=0$ & $18$ & $x^2-xyz+2y-z=0$ & $18$ &  $x^2-xyz+y^2+x=0$ \\ 
 \hline
 $17$ & $x^2-xyz+y^2-1=0$ & $18$ & $x^2-xyz+yz+y=0$ & $18$ & $x^2-xyz+y^2-z=0$ \\ 
 \hline
\end{tabular}
\captionof{table}{\label{tab:H18vieta} Equations of size $H\leq 18$ solvable by Vieta jumping.}
\end{center} 

To solve equations in Table \ref{tab:H18vieta} we will use a similar method to the one used in Section \ref{ex:Vieta2var}. We define the norm as $N(x,y,z)=|x|+|y|+|z|$ and the lower norm as $N'(x,y,z)=\min\{|x|,|y|,|z|\}$. If we have an equation which is quadratic in one variable, we will be jumping in one variable. Similarly, if the equation is quadratic in two variables, we will be jumping in two variables. Next, we must find a transformation $x'$ such that if we have a solution $(x,y,z)$, we have that $(x',y,z)$ is also a solution, and if we are jumping in two variables, we must also find a transformation $y'$ such that if we have a solution $(x,y,z)$, we have that $(x,y',z)$ is a solution. 
Let us call a solution $(x,y,z)$ minimal if the transformation $(x,y,z) \to (x',y,z)$ does not decrease the norm, if we are jumping in one variable, and for jumping in two variables, we must have that neither transformation $(x,y,z) \to (x',y,z)$ or $(x,y,z) \to (x,y',z)$ decreases the norm. We can then find all integer solutions to the equation by applying these transformations to the minimal solutions.

To find the minimal solutions to the equations, we must first prove that the lower norm of any minimal solution is bounded. To do this, we can solve the optimisation problem
$$
\max_{(x,y,z) \in \mathbb{R}^3} \min\{|x|,|y|,|z|\} \quad \text{subject to} \quad P(x,y,z)=0, |x'|\geq|x|.
$$
if we are jumping in one variable, or, for jumping in two variables,
$$
\max_{(x,y,z) \in \mathbb{R}^3} \min\{|x|,|y|,|z|\} \quad \text{subject to} \quad P(x,y,z)=0, |x'|\geq|x|, |y'|\geq|y|.
$$
This can be done in Mathematica using the command
\begin{equation}\label{eq:vietajump1}
\begin{aligned}
{\tt	N[MaxValue[Min[Abs[x],  Abs[y],Abs[z]], \{P(x,y,z)==0,Abs[x']\geq Abs[x] \},  \{x,y,z\}]]}
\end{aligned}
\end{equation}
for jumping in one variable, and the command
\begin{equation}\label{eq:vietajump2}
\begin{aligned}
{\tt	N[MaxValue[Min[Abs[x],  Abs[y],Abs[z]], \{P(x,y,z)==0,Abs[x']\geq Abs[x],} \phantom{\}} \\ \phantom{\{}  {\tt Abs[y']\geq Abs[y]\},  \{x,y,z\}]]}
\end{aligned}
\end{equation}
for jumping in two variables.

Equation
$$
x^2-x y z+y+z=0
$$
is solved in Section 4.1.5 of the book and its integer solutions are
$$
\begin{aligned}
	(x,y,z)=(0,u,-u), (-u^2,-u,u), (-1,-1,u), (1-u,-1,u),(u,-u^2,0),(1,2,3), \\ (2,1,5), (2,2,2),(3,1,5),(5,2,3), \quad u \in \mathbb{Z},
\end{aligned}
$$
and the solutions obtained by swapping $y$ and $z$.

Equation 
$$
x^2-xyz+y^2=0
$$
is solved in Section 4.1.5 of the book and its integer solutions are
$$
(x,y,z)=(0,0,u),(u,-u.-2),(u,u,2), \quad u \in \mathbb{Z}.
$$

Equation
$$
x^2-xyz+y^2-1=0
$$
is solved in Section 4.1.5 of the book and its integer solutions are
$$
(x,y,z)=(\pm (x_n,y_n),u), \quad \text{or} \quad (\pm (x'_n,y'_n),u), \quad \text{or} \quad (\pm (y_n,x_n),u), \quad \text{or} \quad (\pm (y'_n,x'_n),u), 
$$
where $u$ is an integer and 
\begin{equation}\label{rec:x2mxyzpy2m1}
\begin{aligned}
	(x_0,y_0)=(1,0) \quad \text{and} & \quad (x_{n+1},y_{n+1})=(-x_n+u y_n,-x_n u+y_n(u^2-1)), \quad n=0,1,2,\dots
	\\
	(x'_0,y'_0)=(1,0) \quad \text{and} & \quad (x'_{n+1},y'_{n+1})=(x'_n(u^2-1)-u y'_n,x'_n u-y'_n), \quad n=0,1,2,\dots
\end{aligned}
\end{equation}

Equation 
$$
x^2-xyz+y^2+1=0
$$
is solved in Section 4.1.5 of the book and its integer solutions are
$$
(x,y,z)=(-F_{2n-1},-F_{2n+1},3),(-F_{2n-1},F_{2n+1},-3),(F_{2n-1},-F_{2n+1},-3),(F_{2n-1},F_{2n+1},3), \quad n \geq 0,
$$
and the solutions obtained by swapping $x$ and $y$, where $F_n$ is the sequence of Fibonacci numbers defined by
\begin{equation}\label{eq:fibonacci}
F_0=0, \quad F_1=1, \quad F_{n+2} = F_{n+1} + F_n, \quad n=0,1,2,\dots
\end{equation}

Equation 
$$
x^2-xyz+y^2-z=0
$$
is solved in Section 4.1.5 of the book and its integer solutions are
$$
(x,y,z)=(0,u,u^2), (u,0,u^2),(-2,1,-5),(-1,2,-5),(1,-2,-5),(2,-1,-5), \quad u \in \mathbb{Z},
$$
and the solutions obtained by applying the operations
\begin{equation}\label{eq:x2mxyzpy2mz:tranform}
(a) \,\, (x,y,z) \to (x',y,z), \,\, x'=yz-x, \quad \text{and} \quad (b) \,\, (x,y,z) \to (x,y',z), \,\, y'=xz-y.
\end{equation}

The first equation we will consider is
\begin{equation}\label{x2mxyzmy2}
	x^2-xyz-y^2=0.
\end{equation}
This equation is a quadratic in both $x$ and $y$ and can be solved by Vieta jumping, but we will present an easier argument. Assume that $(x,y,z)$ is a solution to \eqref{x2mxyzmy2} with $xyz\neq 0$. Then $z=\frac{x^2-y^2}{xy}$. Let $d=\text{gcd}(x,y)$, $X=x/d$, $Y=y/d$, so that $\text{gcd}(X,Y)=1$. Then $z=\frac{X^2-Y^2}{XY}$. Let $p$ be any prime factor of $X$. Then $p$ must also divide $X^2-Y^2$, hence $p$ is a divisor of $Y$, a contradiction with $\text{gcd}(X,Y)=1$. Hence, $X$ has no prime divisors, which implies that $|X|=1$. By the same argument, $|Y|=1$, but then $z=\frac{X^2-Y^2}{XY}=0$, a contradiction with $xyz\neq 0$. Hence, $xyz=0$, in which case we easily deduce that
$$
(x,y,z)=(0,0,u),\quad (u,\pm u,0), \quad u \in \mathbb{Z}.
$$  

\vspace{10pt}

The next equation we will consider is
\begin{equation}\label{x2mxyzpypzm1}
x^2-xyz+y+z-1=0.
\end{equation}
As this equation is a quadratic in $x$, we will have one Vieta jumping operation
\begin{equation}\label{transformation_yzmx}
 (x,y,z) \to (x',y,z), \, x'=yz-x.
\end{equation}
We must first prove that the lower norm of any minimal solution is bounded. We can do this in Mathematica, using the command \eqref{eq:vietajump1}
$$
{\tt N[MaxValue[Min[Abs[x], Abs[y], Abs[z]], \{x^2 - x y z + y + z - 1 == 0, Abs[y z - x] \geq Abs[x]\}, \{x, y, z\}]]}
$$
This outputs $1.61803$. Therefore we must check the cases $|x| \leq 1, |y| \leq 1$ and $|z| \leq 1$. This gives the following list of solutions containing all minimal solutions,
$$
(x,y,z)=(0,u,1-u),(u,0,1-u^2), (-2,-1,2),(-1,-2,-2),(1,2,2),(2,1,4),\quad u \in \mathbb{Z}
$$
and the solutions obtained by swapping $y$ and $z$. We can then find all integer solutions by applying transformation \eqref{transformation_yzmx} to the minimal solutions. We can then conclude that all integer solutions to equation \eqref{x2mxyzpypzm1} are
$$
\begin{aligned}
(x,y,z)=&  (0,u,1-u),(u,0,1-u^2),(u,1-u^2,0),(u-u^2,u,1-u), \\ & (-1,-2,-2),(1,2,2),(2,1,4),(2,4,1),(3,2,2),(5,-2,-2),  \quad u \in \mathbb{Z}.
\end{aligned}
$$

\vspace{10pt}

The next equation we will consider is
\begin{equation}\label{x2mxyzpypzp1}
x^2-xyz+y+z+1=0.
\end{equation}
As this equation is a quadratic in $x$, we will have one Vieta jumping operation \eqref{transformation_yzmx}. 
The Mathematica command \eqref{eq:vietajump1} for this equation outputs $2.1479$, hence 
the lower norm of any minimal solution is bounded 
and we must check the cases $|x| \leq 2, |y| \leq 2$ and $|z| \leq 2$. This gives a list of solutions containing all minimal solutions, from which we can determine which are minimal, and then find all integer solutions by applying transformation \eqref{transformation_yzmx} to them. Finally, we can conclude that all integer solutions to equation \eqref{x2mxyzpypzp1} are
$$
\begin{aligned}
(x,y,z)=& (-u-u^2,u,-1-u),(0,u,-1-u),(u,-1-u^2,0),(u,0,-1-u^2),(-2,-1,4), \\ & (-2,4,-1),(1,2,4),(1,4,2),(2,1,6),(2,6,1),(4,1,6),(4,6,1),(7,2,4), (7,4,2),   \quad u \in \mathbb{Z}.
\end{aligned}
$$

\vspace{10pt}

The next equation we will consider is
$$
x^2-xyz-y^2+1=0.
$$
As this is a quadratic equation in both $x$ and $y$, we will have two Vieta jumping operations given by
\begin{equation}\label{transformation_yzmx_mxzmy}
	(a) \,\, (x,y,z) \to (x',y,z), \, x'=yz-x, \quad \text{and} \quad (b) \,\, (x,y,z) \to (x,y',z), \, y'=-xz-y.
\end{equation}
The Mathematica command \eqref{eq:vietajump2} for this equation outputs $0.707107$, hence the lower norm of any minimal solution is equal to $0$ 
and we must check the cases $x=0$, $y=0$ and $z=0$. This gives the minimal solutions $(x,y,z)=(0,\pm 1 ,u)$, where $u$ is an arbitrary integer. We can then find all other solutions by applying operations $(a)$ and $(b)$ given in \eqref{transformation_yzmx_mxzmy}. As a result, we obtain, for every $u \in {\mathbb Z}$, an infinite chain of solutions 
$$
\dots \xleftrightarrow{\text{(b)}} (-u,-1,u) \xleftrightarrow{\text{(a)}}  (0,-1,u) \xleftrightarrow{\text{(b)}} (0,1,u) \xleftrightarrow{\text{(a)}} (u,1,u) \xleftrightarrow{\text{(b)}} (u,-u^2-1,u) \xleftrightarrow{\text{(a)}}   \dots  
$$
which can be summarised as $(x,y,z)= (\pm(x_n,y_n),u)$ or $(\pm(x'_n,y'_n),u)$, where $u$ is an integer and $x_n,y_n$ and $x'_n,y'_n$ are given by
\begin{equation}\label{eq:x2mxyzmy2p1}
\begin{aligned}
	 (x_0,y_0)=(0,1) \quad \text{and} \quad (x_{n+1},y_{n+1})=(uy_n-x_n,-y_n(u^2+1)+ux_n) \\
	  (x'_0,y'_0)=(0,1) \quad \text{and} \quad (x'_{n+1},y'_{n+1})=(-x'_n(u^2+1)-uy'_n,-ux'_n-y'_n).
\end{aligned}
\end{equation}

\vspace{10pt}

The next equation we will consider is
\begin{equation}\label{x2mxyzpypzm2}
x^2-xyz+y+z-2=0.
\end{equation}
As this equation is a quadratic in $x$, we will have one Vieta jumping operation \eqref{transformation_yzmx}. The Mathematica command \eqref{eq:vietajump1} for this equation outputs $1.41421$, hence the lower norm of any minimal solution is bounded 
and we must check the cases $|x| \leq 1, |y| \leq 1$ and $|z| \leq 1$. This gives a list containing all minimal solutions. 
We can then find all integer solutions by applying transformation \eqref{transformation_yzmx} to the minimal solutions and conclude that all integer solutions to equation \eqref{x2mxyzpypzm2} are
$$
\begin{aligned}
(x,y,z)=&(0,u,2-u),(1,1,u),(1,u,1),(2u-u^2,u,2-u),(u,0,2-u^2),(u,2-u^2,0), \\ &(u,1,1+u),(u,1+u,1) (-1,-3,-2),(-1,-2,-3), (7,-3,2),(7,2,-3), \quad u \in \mathbb{Z}.
\end{aligned}
$$

\vspace{10pt}

The next equation we will consider is
\begin{equation}\label{x2mxyzpypzp2}
x^2-xyz+y+z+2=0.
\end{equation}
As this equation is a quadratic in $x$, we will have one Vieta jumping operation \eqref{transformation_yzmx}. The Mathematica command \eqref{eq:vietajump1} for this equation outputs $2.26953$, hence the lower norm of any minimal solution is bounded 
and we must check the cases $|x| \leq 2, |y| \leq 2$ and $|z| \leq 2$. This gives a list of solutions containing all minimal solutions. We can then find all integer solutions by applying transformation \eqref{transformation_yzmx} to the minimal solutions. We can then conclude that all integer solutions to equation \eqref{x2mxyzpypzp2} are
$$
\begin{aligned}
(x,y,z)=& (-2u-u^2,u,-2-u),(u,0,-2-u^2),(0,u,-2-u), (-3,-1,5),(-2,-1,5), \\ & (-1,-2,1),(1,2,5),(1,3,3),(2,1,7), (3,1,6), (5,1,7),(8,3,3),(9,5,2),   \quad u \in \mathbb{Z},
\end{aligned}
$$
and the solutions obtained by swapping $y$ and $z$.

\vspace{10pt}

The next equation we will consider is
\begin{equation}\label{x2mxyzpxpymz}
	x^2-xyz+x+y-z=0.
\end{equation}
As this equation is a quadratic in $x$, we will have one Vieta jumping operation
\begin{equation}\label{transformation_yzmxm1}
	(x,y,z) \to (x',y,z), \, x'=yz-x-1.
\end{equation}
The Mathematica command \eqref{eq:vietajump1} for this equation outputs $1.61803$, hence the lower norm of any minimal solution is bounded 
and we must check the cases $|x| \leq 1, |y| \leq 1$ and $|z| \leq 1$. This gives a list of solutions containing all minimal solutions. 
We can then find all integer solutions by applying transformation \eqref{transformation_yzmxm1} to the list of minimal solutions. Finally, we can conclude that all integer solutions to equation \eqref{x2mxyzpxpymz} are
$$
\begin{aligned}
	(x,y,z)=&(0,u,u), (u,-u^2-u,0),(u,0,u^2+u),(u^2-1,u,u),(-4,2,-2),\\ & (-2,1,-3),(-2,3,-1),(-1,2,-2),(2,-1,-5),(2,5,1),  \quad u \in \mathbb{Z}.
\end{aligned}
$$

\vspace{10pt}

The next equation we will consider is
\begin{equation}\label{x2mxyzpxpypz}
	x^2-xyz+x+y+z=0.
\end{equation}
As this equation is a quadratic in $x$, we will have one Vieta jumping operation \eqref{transformation_yzmxm1}. The Mathematica command \eqref{eq:vietajump1} for this equation outputs $2.1479$, hence the lower norm of any minimal solution is bounded 
and we must check the cases $|x| \leq 2, |y| \leq 2$ and $|z| \leq 2$. This gives a list of solutions containing all minimal solutions. We can then find all integer solutions by applying transformation \eqref{transformation_yzmxm1} to the list of minimal solutions. Finally, we can conclude that all integer solutions to equation \eqref{x2mxyzpxpypz} are
$$
\begin{aligned}
	(x,y,z)=& (0,u,-u),(u,0,-u^2-u),(-u^2-1,u,-u),(-1,-2,-2), \\ & (1,2,4),(2,1,7),(4,-2,-2),(4,1,7),(6,2,4),  \quad u \in \mathbb{Z},
\end{aligned}
$$
and the solutions obtained by swapping $y$ and $z$.

 \vspace{10pt}

The next equation we will consider is
\begin{equation}\label{x2mxyzp2ymz}
x^2-xyz+2y-z=0.
\end{equation}
As this equation is a quadratic in $x$, we will have one Vieta jumping operation \eqref{transformation_yzmx}. The Mathematica command \eqref{eq:vietajump1} for this equation outputs $2.30278$, hence the lower norm of any minimal solution is bounded 
and we must check the cases $|x| \leq 2, |y| \leq 2$ and $|z| \leq 2$. This gives a list of solutions containing all minimal solutions. We can then find all integer solutions by applying transformation \eqref{transformation_yzmx} to the minimal solutions. Finally, we can conclude that all integer solutions to equation \eqref{x2mxyzp2ymz} are
$$
\begin{aligned}
(x,y,z)=& (0,u,2u),(u,0,u^2),(2u^2,u,2u),(2u,-2u^2,0), (-11,4,-3),  \\ & (-9,2,-5),(-7,-2,3),(-7,10,-1),(-4,1,-6), (-4,3,-2), \\ & (-3,10,-1),(-2,1,-6),(-2,3,-2),(-1,-2,-1), (-1,2,-5),\\ &(-1,4,-3),(1,-2,3),(3,-2,-1),(3,8,1), (5,8,1), \quad u \in \mathbb{Z}.
\end{aligned}
$$

 \vspace{10pt}
   
The next equation we will consider is
\begin{equation}\label{x2mxyzpyzpy}
x^2-xyz+yz+y=0.
\end{equation}
As this equation is a quadratic in $x$, we will have one Vieta jumping operation \eqref{transformation_yzmx}. The Mathematica command \eqref{eq:vietajump1} for this equation outputs $2.41421$, hence the lower norm of any minimal solution is bounded 
and we must check the cases $|x| \leq 2, |y| \leq 2$ and $|z| \leq 2$. This gives a list of solutions containing all minimal solutions. We can then find all integer solutions by applying transformation \eqref{transformation_yzmx} to the minimal solutions. Finally, we can conclude that all integer solutions to equation \eqref{x2mxyzpyzpy} are
$$
\begin{aligned}
(x,y,z)= & (1,-1,u),(u,-u,-1),(u,-u^2,0),(u,-1,-1-u),(0,0,u),(0,u,-1), (2,1,5),(2,2,3),   \\ &(2,4,2),(3,1,5), (3,3,2),(3,9,1), (4,2,3),(4,8,1),  (6,4,2),(6,9,1), \quad u \in \mathbb{Z}.
\end{aligned}
$$

\vspace{10pt}

The next equation we will consider is
\begin{equation}\label{x2mxyzmy2p2}
x^2-xyz-y^2+2=0.
\end{equation}
As this equation is a quadratic in both $x$ and $y$, we will have two Vieta jumping operations given by \eqref{transformation_yzmx_mxzmy}. The Mathematica command \eqref{eq:vietajump2} for this equation outputs $1$, hence the lower norm of any minimal solution is bounded 
and we must check the cases $|x|\leq 1$, $|y|\leq1$ and $|z|\leq1$. This gives the minimal solutions $(x,y,z)=(-1,-1,2),(-1,1,-2),(1,-1,-2),(1,1,2)$. We can then find all other solutions by applying operations $(a)$ and $(b)$ given in \eqref{transformation_yzmx_mxzmy}. Note that transformations (a) and (b) do not change $z$. Therefore, all solutions to \eqref{x2mxyzmy2p2} have $|z|=2$, which reduces the equation to $x^2 \pm 2xy-y^2+2=0$, or $(x\pm y)^2 -2y^2+2=0$. As $-2y^2+2$ is even, $x\pm y$ must be even, so let $x \pm y=2t$ for some integer $t$, then the equation is reduced to $y^2-2t^2=1$. Up to the names of variables, we have solved this equation previously, and its integer solutions are $(y,t)=(\pm y_n,\pm t_n)$, with
\begin{equation}\label{x2mxyzmy2p2_rec}
(y_0,t_0)=(1,0) \quad \text{and} \quad (y_{n+1},t_{n+1})=(3y_n+4t_n,2y_n+3t_n) \quad n=0,1,2,\dots,
\end{equation}
 which can be found in Table \ref{table3.10}. We can therefore summarise the integer solutions to equation \eqref{x2mxyzmy2p2} as
$$
 (x,y,z)=(\pm 2t_n +s_1 y_n,s_2 y_n,2 s_3 ), \quad s_i=\pm 1, \,\, \text{such that} \quad s_1 s_2 s_3=1, \,\, \text{and} \quad y_n,t_n \quad \text{are given by \eqref{x2mxyzmy2p2_rec}}.
$$

 \vspace{10pt}

The next equation we will consider is
\begin{equation}\label{x2mxyzmy2px}
	x^2-xyz-y^2+x=0.
\end{equation}
As this equation is a quadratic in both $x$ and $y$, we will have two Vieta jumping operations
\begin{equation}\label{transformation_yzmxm1_mxzmy}
	(a) \,\, (x,y,z) \to (x',y,z), \, x'=yz-x-1, \quad \text{and} \quad (b) \,\, (x,y,z) \to (x,y',z), \, y'=-xz-y.
\end{equation}
 The Mathematica command \eqref{eq:vietajump2} for this equation outputs $1$, hence the lower norm of any minimal solution is bounded
and we must check the cases $|x| \leq 1$, $|y| \leq 1$ and $|z| \leq 1$. This gives the minimal solutions 
\begin{equation}\label{x2mxyzmy2px_min}
	(x,y,z)=(0,0,u), \quad u \in \mathbb{Z}.
\end{equation}
We can then find all other integer solutions by applying operations $(a)$ and $(b)$ given in \eqref{transformation_yzmxm1_mxzmy}. In conclusion, all integer solutions to equation \eqref{x2mxyzmy2px} can be obtained from \eqref{x2mxyzmy2px_min} by applying operations (a) and (b) given by \eqref{transformation_yzmxm1_mxzmy}.

\vspace{10pt}

The next equation we will consider is
\begin{equation}\label{x2mxyzmy2mz}
	x^2-xyz-y^2-z=0.
\end{equation}
As this equation is a quadratic in both $x$ and $y$, we will have two Vieta jumping operations \eqref{transformation_yzmx_mxzmy}. The Mathematica command \eqref{eq:vietajump2} for this equation outputs $1$, hence the lower norm of any minimal solution is bounded 
and we must check the cases $|x| \leq 1$, $|y| \leq 1$ and $|z| \leq 1$. This gives the following list containing all minimal solutions 
$$
	\begin{aligned}
		(x,y,z)=&(1,u,1-u),(0,u,-u^2),(-1,u,u+1),(u,-1,-1-u),\\ & (u,0,u^2),(u,1,u-1),(u,\pm u,0), \quad u \in \mathbb{Z}.
	\end{aligned}
	$$
We can then find all other integer solutions to equation \eqref{x2mxyzmy2mz} by applying operations $(a)$ and $(b)$ given in \eqref{transformation_yzmx_mxzmy}. 

\vspace{10pt}

The next equation we will consider is
\begin{equation}\label{eq:x2mxyzpy2m2}
x^2-xyz+y^2-2=0.
\end{equation}
As this equation is a quadratic in both $x$ and $y$, we will have two Vieta jumping operations
\begin{equation}\label{transformation_yzmx_xzmy}
(a) \,\, (x,y,z) \to (x',y,z), \, x'=yz-x, \quad \text{and} \quad (b) \,\, (x,y,z) \to (x,y',z), \, y'=xz-y.
\end{equation}
The Mathematica command \eqref{eq:vietajump2} for this equation outputs $0.839287$, hence the lower norm of any minimal solution is bounded 
and we must check the cases $x=0$, $y=0$ and $z=0$. This gives the minimal solutions $(x,y,z)=(\pm 1,\pm 1,0)$. We can then find all other solutions by applying operations $(a)$ and $(b)$ given in \eqref{transformation_yzmx_xzmy}, however, this produces no further solutions. Therefore, we can conclude that the integer solutions to equation \eqref{eq:x2mxyzpy2m2} are 
$$
(x,y,z)=(\pm 1,\pm 1,0).
$$

\vspace{10pt}

The next equation we will consider is
\begin{equation}\label{x2mxyzpy2p2}
x^2-xyz+y^2+2=0.
\end{equation}
As this equation is a quadratic in both $x$ and $y$, we will have two Vieta jumping operations given by \eqref{transformation_yzmx_xzmy}. The Mathematica command \eqref{eq:vietajump2} for this equation outputs $2.3593$, hence the lower norm of any minimal solution is bounded 
and we must check the cases $|x|\leq 2$, $|y|\leq2$ and $|z|\leq2$. This gives the minimal solutions $(x,y,z)=(-1,-1,4),(-1,1,-4),(1,-1,-4),(1,1,4)$. We can then find all other solutions by applying operations $(a)$ and $(b)$ given in \eqref{transformation_yzmx_xzmy}. Note that these transformations do not change $z$. Therefore, all integer solutions to \eqref{x2mxyzpy2p2} have $|z|=4$, which reduces the equation to 
$$
x^2\pm 4xyz+y^2+2=0=(x\pm 2y)^2-3y^2+2.
$$
 We can make the substitution $t=x \pm 2y$ for some integer $t$. After rearranging, we can reduce the equation to $t^2-3y^2=-2$, which, up to the names of variables, this is equation \eqref{eq:3.10_ii}, and its integer solutions are $(t,y)=(\pm x_n,\pm y_n)$, where $(x_n,y_n)$ is given by \eqref{eq:x2m3y2p2sol}. 
We can therefore summarise that the integer solutions to equation \eqref{x2mxyzpy2p2} are
$$
 (x,y,z)=(\pm x_n +2 s_1 y_n,s_2 y_n,4 s_3 ), \quad s_i=\pm 1, \,\, \text{such that} \quad s_1 s_2 s_3=1, \,\, \text{and} \quad x_n,y_n \quad \text{are given by \eqref{eq:x2m3y2p2sol}}.
$$

\vspace{10pt}

The final equation we will consider is
\begin{equation}\label{x2mxyzpy2px}
x^2-xyz+y^2+x=0.
\end{equation}
As this equation is a quadratic in both $x$ and $y$, we will have two Vieta jumping operations
\begin{equation}\label{transformation_yzmxm1_xzmy}
(a) \,\, (x,y,z) \to (x',y,z), \, x'=yz-x-1, \quad \text{and} \quad (b) \,\, (x,y,z) \to (x,y',z), \, y'=xz-y.
\end{equation}
The Mathematica command \eqref{eq:vietajump2} for this equation outputs $2.41421$, hence the lower norm of any minimal solution is bounded 
and we must check the cases $|x| \leq 2$, $|y| \leq 2$ and $|z| \leq 2$. This gives the minimal solutions 
$$
(x,y,z)=(0,0,u),(1,-1,-3),(1,1,3).
$$
 We can then find all other integer solutions to equation \eqref{x2mxyzpy2px} by applying operations $(a)$ and $(b)$ given in \eqref{transformation_yzmxm1_xzmy}.

\begin{longtable}[c]{ |c|c|c|c|c|c| } 
 \hline
  Equation & Solution $(x,y,z)$ \\
 \hline\hline
 \endfirsthead
 
 $x^2-xyz+y+z=0$ & $(0,u,-u),(-u^2,-u,u),(-1,-1,u),(1-u,-1,u),$ \\ & $(u,-u^2,0),(1,2,3),(2,1,5),(2,2,2),(3,1,5),(5,2,3),$  \\ & and the solutions obtained by exchanging $y$ and $z$\\\hline
 $x^2-xyz-y^2=0$ &$(0,0,u),(u,\pm u,0)$ \\\hline
 $x^2-xyz+y^2=0$ & $(0,0,u),(u,-u,-2),(u,u,2)$ \\\hline
 $x^2-xyz+y+z-1=0$ &$(0,u,1-u),(u,0,1-u^2),(u-u^2,u,1-u),\pm(1,2,2),(2,1,4)$ \\ & $(3,2,2),(5,-2,-2),$ and the solutions obtained by exchanging $y$ and $z$ \\ \hline
 $x^2-xyz+y+z+1=0$ &  $(-u-u^2,u,-1-u), (u,0,-1-u^2),(0,u,-1-u),$ \\ & $(-2,-1,4),(1,2,4),(2,1,6),(4,1,6),(7,2,4),$ \\ & and the solutions obtained by exchanging $y$ and $z$  \\\hline
 $x^2-xyz-y^2+1=0$ & $(\pm(x_n,y_n),u),$ or $(\pm(x'_n,y'_n),u)$ where $x_n,y_n$ and $x'_n,y'_n$ are given by \eqref{eq:x2mxyzmy2p1}. \\\hline 
 $x^2-xyz+y^2-1=0$ & $(\pm(x_n,y_n),u),$ $(\pm(x'_n,y'_n),u),$ $(\pm(y_n,x_n),u)$ or $(\pm(y'_n,x'_n),u)$ \\
 & where $x_n,y_n$ and $x'_n,y'_n$ are given by \eqref{rec:x2mxyzpy2m1}. \\\hline 
 $x^2-xyz+y^2+1=0$ & $(-F_{2n-1},-F_{2n+1},3),(-F_{2n-1},F_{2n+1},-3),(F_{2n-1},-F_{2n+1},-3),$ \\ & $(F_{2n-1},F_{2n+1},3),$  where $F_n$ is defined by \eqref{eq:fibonacci}, \\ & and the solutions obtained by exchanging $y$ and $z$ \\\hline
 $x^2-xyz+y+z-2=0$ & $(1,1,u),(u,2-u^2,0),(u,1,1+u),(0,u,2-u),$ \\ &$(2u-u^2,u,2-u),(-1,-3,-2),(7,-3,2),$ and the solutions  \\ & obtained by exchanging $y$ and $z$ \\\hline
 $x^2-xyz+y+z+2=0$ &$(0,u,-2-u),(-2u-u^2,u,-2-u),(u,0,-2-u^2),$ \\ & $(-3,-1,5),(-2,-1,5),(-1,-2,1),(1,2,5),(1,3,3),$ \\ & $(2,1,7),(3,1,6),(5,1,7),(8,3,3),(9,5,2)$ \\ & and the solutions obtained by exchanging $y$ and $z$ \\\hline
 $x^2-xyz+x+y-z=0$ & $(0,u,u), (u,-u^2-u,0),(u,0,u^2+u),(u^2-1,u,u),(-4,2,-2),$ \\ & $(-2,1,-3),(-2,3,-1),(-1,2,-2),(2,-1,-5),(2,5,1)$ \\\hline
 $x^2-xyz+x+y+z=0$ & $(0,u,-u),(u,0,-u^2-u),(-u^2-1,u,-u),$ \\ & $(-1,-2,-2),(1,2,4),(2,1,7),(4,-2,-2),(4,1,7),(6,2,4),$ \\ &  and the solutions obtained by exchanging $y$ and $z$ \\\hline
 $x^2-xyz+2y-z=0$ &  $(0,u,2u),(u,0,u^2),(2u^2,u,2u),(2u,-2u^2,0)$ \\ & $(-11,4,-3),(-9,2,-5),(-7,-2,3),(-7,10,-1),(-4,1,-6),$\\ & $(-4,3,-2),(-3,10,-1),(-2,1,-6),(-2,3,-2),(-1,-2,-1),$ \\ & $(-1,2,-5),(-1,4,-3),(1,-2,3),(3,-2,-1),(3,8,1),(5,8,1)$ \\ \hline
 $x^2-xyz+yz+y=0$ & $(0,u,-1),(1,-1,u),(u,-u,-1),(u,-u^2,0),(u,-1,-1-u)$ \\ &$(0,0,u),(2,1,5),(2,2,3),(2,4,2),(3,1,5),(3,3,2),$ \\ & $(3,9,1),(4,2,3),(4,8,1),(6,4,2),(6,9,1)$ \\\hline
 $x^2-xyz-y^2+2=0$ & $ (\pm 2t_n +s_1 y_n,s_2 y_n,2 s_3 ), \quad s_i=\pm 1, \,\, \text{such that}$ \\ & $s_1 s_2 s_3=1, \,\, \text{and} \quad y_n,t_n \quad \text{are given by \eqref{x2mxyzmy2p2_rec}}$.  \\\hline
 $x^2-xyz+y^2-2=0$ &$(\pm 1,\pm 1,0)$ \\\hline  
 $x^2-xyz+y^2+2=0$ &   $(\pm x_n +2 s_1 y_n,s_2 y_n,4 s_3 ), \quad s_i=\pm 1, \,\, \text{such that}$ \\ &$ s_1 s_2 s_3=1, \,\, \text{and} \quad x_n,y_n$ are given by \eqref{eq:x2m3y2p2sol}.   \\\hline  
 $x^2-xyz+y^2+x=0$ & All solutions can be obtained from $(0,0,u),(1,-1,-3),(1,1,3)$ \\ &  by applying operations (a) and (b) given by \eqref{transformation_yzmxm1_xzmy}. \\\hline 
 $x^2-xyz-y^2+x=0$ & All solutions can be obtained from $(0,0,u)$ \\ &  by applying operations (a) and (b) given by \eqref{transformation_yzmxm1_mxzmy}. \\\hline 
 $x^2-xyz+y^2-z=0$ & All solutions can be obtained from \\ & $(0,u,u^2),(u,0,u^2),(-2,1,-5),(-1,2,-5),(1,-2,-5),(2,-1,-5)$ \\ & by applying operations (a) and (b) given by \eqref{eq:x2mxyzpy2mz:tranform}.  \\\hline 
 $x^2-xyz-y^2-z=0$ & All solutions can be obtained from $(1,u,1-u),(0,u,-u^2),$ \\ & $(-1,u,u+1),(u,-1,-1-u), (u,0,u^2),(u,1,u-1),(u,\pm u,0),$ \\ &  by applying operations (a) and (b) given by \eqref{transformation_yzmx_mxzmy}. \\\hline

\caption{\label{tab:H18vietasol} Integer solutions to the equations listed in Table \ref{tab:H18vieta}. Assume that $u$ is an integer.} 
\end{longtable}

\subsection{Exercise 4.15}\label{ex:H21vietasub}
\textbf{\emph{For each equation listed in Table \ref{tab:H21vietasub}, reduce it to an equation of the form
		\begin{equation}\label{eq:vieta}
		a_i x_i^2+Q_ix_i+R_i=0
		\end{equation}
		where $a_i \neq 0$ is an integer and $Q_i$ and $R_i$ are polynomials in other variables, which is
		solvable by Vieta jumping, that is, one for which optimization problem 
		\begin{equation}\label{eq:vietaopt}
		\begin{split}
		& \max_{(x_1,\dots,x_n)\in {\mathbb R}^n} \min\{|x_1|, \dots, |x_n|\}, \\
		& \text{subject to} \quad P=0, \quad |x_i| \leq |-(Q_i/a_i)-x_i|, \,\, i \in S,
		\end{split}
		\end{equation} 
		returns a finite value. You do not need to solve the resulting equations.}}

	\begin{center}
		\begin{tabular}{ |c|c|c|c|c|c| } 
			\hline
			$H$ & Equation & $H$ & Equation & $H$ & Equation \\ 
			\hline\hline
			$19$ & $x^3-xyz+y+1=0$ & $20$ & $x^3-xyz+y-z=0$ & $21$ & $x^3-xyz-x+y+1=0$ \\ 
			\hline
			$20$ & $x^3-xyz+y+2=0$ & $20$ & $x^3-xyz+y+z=0$ & $21$ & $x^3-xyz+x+y+1=0$ \\ 
			\hline
			$20$ & $x^3-xyz-x+y=0$ & $21$ & $x^3-xyz+y+3=0$ & $21$ & $x^3-xyz+y-z+1=0$ \\ 
			\hline
			$20$ & $x^3-xyz+x+y=0$ & $21$ & $x^3-xyz+2y+1=0$ & $21$ & $x^3-xyz+y+z+1=0$ \\ 
			\hline
		\end{tabular}
		\captionof{table}{\label{tab:H21vietasub} Equations of the form \eqref{eq:x1Pplinpb} of size $H\leq 21$ solvable by substitutions in combination with Vieta jumping.}
	\end{center} 

Equations in Table \ref{tab:H21vietasub} are of the form
\begin{equation}\label{eq:x1Pplinpb}
x_i P(x_1,x_2,\dots,x_n)+\sum_{j \neq i} a_j x_j+b=0,
\end{equation}
where $a$ and $b$ are integers. We can then use a change of variable, to reduce the equation to one which is solvable by Vieta jumping. For each equation in Table \ref{tab:H21vietasub}, Table \ref{tab:H21vietasubsol} presents a change of variable and the equation it reduces to which can be solved by Vieta jumping.

\begin{center}
\begin{tabular}{|c|c|c|c|}
		\hline
		Equation & Change of Variable & Reduced Equation \\ \hline \hline
		$x^3-xyz+y+1=0$ &$y_1=-x^2+yz$&$x^2-xy_1z+y_1+z=0$ \\\hline
		$x^3-xyz+y+2=0$ &$y_1=-x^2+yz$&$x^2-xy_1z+y_1+2z=0$ \\\hline
		$x^3-xyz-x+y=0$ & $y_1=-x^2+yz+1$&$x^2-xy_1z+y_1-1=0$ \\\hline
		$x^3-xyz+x+y=0$ &$y_1=-x^2+yz-1$&$x^2-xy_1z+y_1+1=0$ \\\hline
		$x^3-xyz+y-z=0$ &$y_1=-x^2+yz$&$x^2-xy_1z-z^2+y_1=0$ \\\hline
		$x^3-xyz+y+z=0$ & $y_1=-x^2+yz$&$x^2-xy_1z+z^2+y_1=0$  \\\hline
		$x^3-xyz+y+3=0$ &$y_1=-x^2+yz$&$x^2-xy_1z+y_1+3z=0$ \\\hline 
		$x^3-xyz+2y+1=0$ &$y_1=-x^2+yz$&$2x^2-xy_1z+2y_1+z=0$  \\\hline
		$x^3-xyz-x+y+1=0$ &$y_1=-x^2+yz+1$&$x^2-xy_1z+y_1+z-1=0$ \\\hline
		$x^3-xyz+x+y+1=0$ & $y_1=-x^2+yz-1$&$x^2-xy_1z+y_1+z+1=0$ \\\hline
		$x^3-xyz+y-z+1=0$ &  $y_1=-x^2+yz$&$x^2-xy_1z-z^2+y_1+z=0$ \\\hline
		$x^3-xyz+y+z+1=0$ &$y_1=-x^2+yz$&$x^2-xy_1z+z^2+y_1+z=0$ \\\hline
	\end{tabular}
	\captionof{table}{\label{tab:H21vietasubsol} Equations in Table \ref{tab:H21vietasub} reduced to equations solvable by Vieta jumping.}
\end{center}

 \subsection{Exercise 4.16}\label{ex:H22vietasub}
 \textbf{\emph{For each equation listed in Table \ref{tab:H22vietasub}, reduce it to an equation solvable by Vieta jumping, that is, one for which optimization problem \eqref{eq:vietaopt} returns a finite value. You do not need to solve the resulting equations.}}
 
 	\begin{center}
		\begin{tabular}{ |c|c|c|c|c|c| } 
			\hline
			$H$ & Equation & $H$ & Equation & $H$ & Equation \\ 
			\hline\hline
			$21$ & $x^2y+xyz+z^2-1=0$ & $22$ & $x^2y+xyz+z^2+x=0$ & $22$ & $x^2y-yz^2+xz+y=0$ \\ 
			\hline
			$21$ & $x^2y+xyz+z^2+1=0$ & $22$ & $x^2y+xyz+z^2-y=0$ & $22$ & $x^2y-yz^2+x^2-2=0$ \\ 
			\hline
			$21$ & $x^2y-yz^2+x^2-1=0$ & $22$ & $x^2y+xyz+z^2+y=0$ & $22$ & $x^2y-yz^2+x^2+2=0$ \\ 
			\hline
			$21$ & $x^2y-yz^2+x^2+1=0$ & $22$ & $x^2y+xyz+z^2+z=0$ & $22$ & $x^2y-yz^2+x^2+x=0$ \\ 
			\hline
			$21$ & $x^2y-yz^2+xz+1=0$ & $22$ & $x^2y-yz^2+xz+2=0$ & $22$ & $x^2y-yz^2+x^2-y=0$ \\ 
			\hline
			$22$ & $x^2y+xyz+z^2-2=0$ & $22$ & $x^2y-yz^2+xz+x=0$ & $22$ & $x^2y-yz^2+x^2+y=0$ \\ 
			\hline
			$22$ & $x^2y+xyz+z^2+2=0$ & $22$ & $x^2y-yz^2+xz-y=0$ & $22$ & $x^2y-yz^2+x^2+z=0$ \\ 
			\hline
		\end{tabular}
		\captionof{table}{\label{tab:H22vietasub} Equations of size $H\leq 22$ solvable by linear substitutions in combination with Vieta jumping.}
	\end{center} 
 
 The equations listed in Table \ref{tab:H22vietasub} are either quadratic in $x$ or quadratic in both $x$ and $z$. Each equation can be reduced to the form \eqref{eq:vieta}, for which the optimisation problem \eqref{eq:vietaopt} returns a finite value. This can be done in Mathematica using the command
 \begin{equation}\label{eq:vietajumpxz}
 	\begin{aligned}
 	{\tt	N[MaxValue[Min[Abs[x],  Abs[y],Abs[z]], \{P(x,y,z)==0,Abs[x']\geq Abs[x], } \phantom{\}} \\ \phantom{\{}  {\tt Abs[z']\geq Abs[z]\},  \{x,y,z\}]]}
 	\end{aligned}
 \end{equation} 

  Equation 
 $$
 x^2 y+x y z+z^2+1=0
 $$
 is reduced to 
 $$
 x^2-2xz'+xyz'+(z')^2+1=0
 $$
 in Section 4.1.6 of the book, which is solvable by Vieta jumping.
 
 Equation 
 $$
 x^2 y-y z^2+x^2+1=0
 $$
is reduced to 
 $$
 u^2+4uvy+2uv+v^2+4=0
 $$
 in Section 4.1.6 of the book, which is solvable by Vieta jumping.
 
Let us consider a few examples.
 The first equation we will consider is
$$
 x^2y+xyz+z^2-1=0.
$$
 The leading part of this equation is $x^2y+xyz=xy(x+z)$, we can replace $z$ by a new variable $z_1=x+z$, which reduces the equation to
 $$
 x^2+xyz_1-2xz_1+z_1^2-1=0= x^2+xz_1(y-2)+z_1^2-1.
 $$
 Making the substitution $y_1=y-2$, we reduce the equation to $x^2+xy_1z_1+z_1^2-1=0$. We can solve this equation using Vieta jumping as the command \eqref{eq:vietajumpxz}
 $$
 \begin{aligned}
{\tt N[  MaxValue[Min[Abs[x],  Abs[y],Abs[z]], \{x^2 + x y z + z^2 - 1==0,Abs[-y z-x]\geq Abs[x],} \phantom{\}} \\ \phantom{\{}  {\tt Abs[-x y-z]\geq Abs[z]\},  \{x,y,z\}]]}
\end{aligned}
 $$
 outputs $0.618034$. 
 
 \vspace{10pt}
 
 The next equation we will consider is
$$
  x^2y-yz^2+x^2-1=0.
$$
  The leading part of this equation is $x^2y-yz^2=y(x^2-z^2)=y(x+z)(x-z)$, we can make the substitutions $x_1=x+z$ and $z_1=x-z$, so that $x=\frac{x_1+z_1}{2}$ and $z=\frac{x_1-z_1}{2}$. With new variables, and after multiplication by $4$, the equation becomes 
  \begin{equation}\label{x2ymyz2px2m1_red}
  x_1^2+4x_1 y z_1+2x_1z_1+z_1^2-4=0.
  \end{equation}  
   To make sure that $x$ and $z$ are integers, we must have that $x_1$ and $z_1$ have the same parity. It is easy to check that every integer solution to \eqref{x2ymyz2px2m1_red} automatically satisfies this condition. We can solve this equation using Vieta jumping as the command \eqref{eq:vietajumpxz} for this equation
 outputs $1$. 
 
 \vspace{10pt}
 
The final example we will consider is equation
$$
  x^2y-yz^2+xz+2=0.
$$
After multiplying by $4$, we can make the change of variables $x_1=x+z$ and $z_1=x-z$ and reduce the equation to  
  \begin{equation}\label{x2ymyz2pxzp2_red}
  x_1^2+4x_1 y z_1-z_1^2+8=0
  \end{equation}  
   with the restriction that $x_1$ and $z_1$ have the same parity. It is easy to check that every integer solution to \eqref{x2ymyz2pxzp2_red} automatically satisfies this condition. We can solve \eqref{x2ymyz2pxzp2_red} using Vieta jumping as the command \eqref{eq:vietajumpxz} for this equation
 outputs $1.25992$.

The other equations in Table \ref{tab:H22vietasub} can be solved similarly, and for each equation, Table \ref{tab:H22vietasubsol} presents the necessary transformations and the equation they reduce to which is solvable by Vieta jumping. 

\begin{center}

	\captionof{table}{\label{tab:H22vietasubsol} Equations from Table \ref{tab:H22vietasub} reduced by linear substitutions to equations solvable by Vieta jumping.}
\end{center}

\subsection{Exercise 4.19}\label{ex:H18quadred}
\textbf{\emph{Reduce all equations listed in Table \ref{tab:H18quadred} to the form \eqref{eq:genquadform}.  }}

\begin{center}
%
\captionof{table}{\label{tab:H18quadred} Quadratic equations of size $H\leq 18$ not of the form \eqref{eq:genquadform}.}
\end{center} 

All equation in Table \ref{tab:H18quadred} (perhaps after a change of variable) can be written in the form 
\begin{equation}\label{eq:gen_quad}
\sum_{i=1}^n \sum_{j=i}^n a_{ij}x_i x_j +\sum_{i=1}^n b_i x_i +c=0
\end{equation}
where $a_{ij}, b_i$ and $c$ are some integer coefficients and $a_{nn} \neq 0$. By making linear substitution, the equations can be reduced to equations of the form
\begin{equation}\label{eq:genquadform}
\sum_{i=1}^n a_i (x_i ')^2=b
\end{equation}
for some integers $a_1,\dots,a_n,b$. To do this, we will use the method described in Section 4.2.1 of the book, which we summarise below for convenience. Equation \eqref{eq:gen_quad} can be written as
\begin{equation}\label{eq:gen_quadLQ}
a_{nn} x_n^2+x_n L(x_1,\dots,x_{n-1})+Q(x_1,\dots,x_{n-1})=0
\end{equation}
where $L$ and $Q$ are some linear and quadratic polynomials in $x_1,\dots,x_{n-1}$, respectively. Multiplying the equation by $4a_{nn}$ and rewriting as
$$
(2a_{nn} x_n+L(x_1,\dots,x_{n-1}))^2-L^2(x_1,\dots,x_{n-1})+4a_{nn}Q(x_1,\dots,x_{n-1})=0,
$$
suggests to make a new variable $x'_{nn}=2a_{nn} x_n+L(x_1,\dots,x_{n-1})$. By repeating this process if necessary, we will obtain an equation of the form \eqref{eq:genquadform}. 

Let us illustrate this method by considering a number of equations, starting with an easy example.
Equation
\begin{equation}\label{xpx2my2mz2m3}
x+x^2-y^2-z^2-3=0
\end{equation}
is of the form \eqref{eq:gen_quadLQ} with $x_n=x$, $a_{nn}=1$, $L=1$, $Q=-y^2-z^2-3$. Then multiplying \eqref{xpx2my2mz2m3} by 4, we can rewrite this as 
$$
(2x+1)^2-1+4(-y^2-z^2-3)=0,
$$
which suggests to make the new variable $x'=2x+1$. We then obtain
$$
(x')^2-4y^2-4z^2=13,
$$
which is of the form \eqref{eq:genquadform}.

\vspace{10pt}

The next equation we will consider is
\begin{equation}\label{xpx2pymy2pz2m1}
x+x^2+y-y^2+z^2-1=0.
\end{equation}
This equation is of the form \eqref{eq:gen_quadLQ} with
$x_n=x$, $a_{nn}=1$, $L=1$, $Q=y-y^2+z^2-1$. After multiplying \eqref{xpx2pymy2pz2m1} by 4, we can rewrite this as 
$$
(2x+1)^2-1+4(y-y^2+z^2-1)=0,
$$
which suggests to make the new variable $x'=2x+1$. We then obtain
$$
(x')^2+4y-4y^2+4z^2-5=0,
$$
however, this is not of the form \eqref{eq:genquadform}. So, we could repeat the process with $x_n=y$, $a_{nn}=-4$, $L=4$, $Q=(x')^2+4z^2-5$ and then multiply the equation by $4a_{nn}$. However, we can see that $(x')^2+4y-4y^2+4z^2-5=0=(x')^2-(2y-1)^2+1+4z^2-5$, so making the new variable $y'=2y-1$, we obtain 
$$
(x')^2-(y')^2+4z^2=4
$$
which is now of the form \eqref{eq:genquadform}.

\vspace{10pt}

The next equation we will consider is
\begin{equation}\label{xpx2pxypy2mz2}
x+x^2+xy+y^2-z^2=0.
\end{equation}
This equation is of the form \eqref{eq:gen_quadLQ} with
$x_n=x$, $a_{nn}=1$, $L=y+1$, $Q=y^2-z^2$. After multiplying \eqref{xpx2pxypy2mz2} by 4, we can rewrite this as 
$$
(2x+y+1)^2-(y+1)^2+4(y^2-z^2)=0.
$$
This suggests to make the new variable $x'=2x+y+1$, and we obtain 
$$
(x')^2+3y^2-2y-4z^2-1=0.
$$
However, this is not of the form \eqref{eq:genquadform}. So, we could repeat the process with $x_n=y$, $a_{nn}=3$, $L=-2$, $Q=(x')^2-4z^2-1$ and then multiply the equation by $4a_{nn}$. However, by multiplying the equation by $3$, we have $3(x')^2+9y^2-6y-12z^2-3=0=(x')^2+(3y-1)^2-4-12z^2$, so making the new variable $y'=3y-1$, we obtain 
$$
(x')^2+(y')^2-12z^2=4
$$
which is now of the form \eqref{eq:genquadform}.

\vspace{10pt}

The final example we will consider is equation
\begin{equation}\label{eq:x2py2p2tzpz}
	x^2+y^2+2tz+z=0.
\end{equation}
Writing this equation in the form \eqref{eq:gen_quad} is not possible, as we have $a_{nn} =0$ and $a_{in} \neq 0$ for some $i$ and $n$. Instead, let us label the variables as $x_n=z$, $x_3=t$, $x_2=y$ and $x_1=x$. Then we have $a_{nn}=0$ and $a_{3n}=2$. The linear substitution $t=t'+z$ transforms the equation to
$$
x^2+y^2+2(t'+z)z+z=0=x^2+y^2+2t'z+2z^2+z.
$$
This equation is now of the form \eqref{eq:gen_quad} with 
$x_n=z$, $a_{nn}=2$, $L=1+2t'$ and $Q=x^2+y^2$, and multiplying by $4a_{nn}$, it can be rewritten as 
$$
(4z+2t'+1)^2-(2t'+1)^2+8x^2+8y^2=0.
$$
This suggests to make the substitutions $z'=4z+2t'+1$ and $t''=2t'+1$, and \eqref{eq:x2py2p2tzpz} reduces to 
$$
(z')^2-(t'')^2+8x^2+8y^2=0,
$$
which is of the form \eqref{eq:genquadform}.

All other equations in Table \ref{eq:genquadform} can be solved similarly. 
Table \ref{tb:gen_quad} presents, for each equation, the necessary transformations and the equation of the form \eqref{eq:genquadform} it reduces to.

 \begin{center}
\begin{tabular}{|c|c|c|c|}
\hline
Equation & Transformation & Reduced Equation \\ \hline \hline
$x+x^2-y^2-z^2-3=0$ &$x'=2x+1$&$(x')^2-4y^2-4z^2=13$ \\\hline
 $x+x^2-y^2-z^2+3=0$ & $x'=2x+1$ & $(x')^2-4y^2-4z^2=-11$ \\\hline
 $x+x^2+y-y^2+z^2-1=0$ & $x'=2x+1$ & $(x')^2-(y')^2+4z^2=4$ \\
& $y'=2y-1$ & \\\hline
 $x+x^2+y-y^2+z^2+1=0$ & $x'=2x+1$ & $(x')^2-(y')^2+4z^2=-4$\\
& $y'=2y-1$ & \\\hline
 $x+x^2+y+y^2-z^2-1=0$ & $x'=2x+1$ & $(x')^2+(y')^2-4z^2=6$ \\
\hline

 \end{tabular}
 \captionof{table}{Reducing the equations listed in Table \ref{eq:genquadform} to equations of the form \eqref{eq:genquadform}.\label{tb:gen_quad}}
 \end{center}

\subsection{Exercise 4.23}\label{ex:H20hom}
\textbf{\emph{Use formulas \eqref{eq:genhomquadsoln}-\eqref{eq:genhomquadsol}, \eqref{eq:genhomquadsol2} to describe all integer solutions to the equations listed in Table \ref{tab:H20hom}. }}

\begin{center}
\begin{tabular}{ |c|c|c|c|c|c| } 
 \hline
 $H$ & Equation & $H$ & Equation & $H$ & Equation \\ 
 \hline\hline
 $20$ & $x^2+xy+y^2-z^2-t^2=0$ & $20$ & $x^2+2xy+xz+z^2=0$ & $20$ & $2x^2+xy-y^2-z^2=0$ \\ 
 \hline
 $20$ & $x^2+xy-y^2+z^2+t^2=0$ & $20$ & $2x^2-y^2-yz-z^2=0$ & $20$ & $2x^2+xy-y^2+z^2=0$ \\ 
 \hline
 $20$ & $x^2+2xy+xz-z^2=0$ & $20$ & $2x^2+y^2-yz-z^2=0$ & $20$ & $2x^2+xy+y^2-z^2=0$ \\ 
 \hline
\end{tabular}
\captionof{table}{\label{tab:H20hom} Homogeneous quadratic equations in non-diagonal form of size $H\leq 20$.}
\end{center} 

All equations in Table \ref{tab:H20hom} are of the form 
\begin{equation}\label{quad_hom}
\sum_{i=1}^n \sum_{j=i}^n a_{ij} x_i x_j =0,
\end{equation}
where $a_{ij}$ are integer coefficients. To solve these equations, we will use the method presented in Section 4.2.3 of the book, which we summarise below for convenience.
First, we need to find an integer solution with $x_n \neq 0$, which we denote as $(x_1^0,\dots,x_n^0)$. The following theorem allows us to find such a solution.
\begin{theorem}\label{theorem_cassel}[Theorem 4.21 in the book]
	Let $S$ be the set of integer solutions to \eqref{quad_hom} other than $(0,0,\dots,0)$. If $S$ is non-empty, then there exists a solution $(x_1,\dots,x_n) \in S$ such that
	$$
	\max(|x_1|,\dots,|x_n|) \leq (12F)^{(n-1)/2}, \quad \text{where} \quad F=\sum_{i=1}^n \sum_{j=i}^n |a_{ij}|.
	$$
\end{theorem}
We can describe the solutions to these equations $(x_1,\dots,x_n)$ with $x_n \neq 0$, using the formulas
\begin{equation}\label{eq:genhomquadsoln}
x_n = \frac{u}{q} \left( x_n^0 \sum_{i=1}^{n-1} \sum_{j=i}^{n-1} a_{ij} u_i u_j \right),
\end{equation}
and
\begin{equation}\label{eq:genhomquadsol}
x_k = \frac{u}{q} \left( x_k^0 \sum_{i=1}^{n-1} \sum_{j=i}^{n-1} a_{ij} u_i u_j -u_k  \sum_{i=1}^{n-1} \sum_{j=i}^{n-1} a_{ij}(x_i^0 u_j + x_j^0 u_i) - x_n^0 u_k  \sum_{i=1}^{n-1} a_{in} u_i \right), \quad k=1,\dots,n-1,
\end{equation}
where $u,u_1,\dots, u_{n-1}$ are arbitrary parameters and $q$ is any common divisor of the expressions in the parentheses, 
and we may assume that gcd$(u_1, \dots , u_{n-1}) = 1$ and gcd$(u, q) = 1$. In addition, there is a family of solutions
\begin{equation}\label{eq:genhomquadsol2}
	x_k = \frac{x_k^0 w + x_n^0 u_k v}{q}, \quad k=1,\dots,n-1, \quad x_n = \frac{x_n^0 w}{q},
\end{equation}
where $w$ and $v$ are arbitrary integer parameters, $(u_1,\dots u_{n-1})$ is an integer solution to system 
\begin{equation}\label{eq:genhomquadratzero}
	\sum_{i=1}^{n-1} \sum_{j=i}^{n-1} a_{ij} u_i u_j = \sum_{i=1}^{n-1} \sum_{j=i}^{n-1} a_{ij} (x_i^0 u_j + x_j^0 u_i) + x_n^0 \sum_{i=1}^{n-1} a_{in} u_i = 0,
\end{equation}
where we may assume that $\text{gcd}(u_1,\dots, u_{n-1})=1$, and $q$ is any common divisor of the numerators. The second equation in system \eqref{eq:genhomquadratzero} is linear, hence, expressing any variable from it and substituting into the first equation, we obtain an equation in $n-2$ variables, which can be solved using the same method. 

To describe solutions with $x_n=0$, we can apply the same method as above for the resulting equation in $n-1$ variables. 

Let us consider the first equation
\begin{equation}\label{x2pxypy2mz2mt2}
x^2+xy+y^2-z^2-t^2=0.
\end{equation}
This is an equation of the form \eqref{quad_hom} with $n=4$, $x_1=x$, $x_2=y$, $x_3=z$, $x_4=t$, $a_{11}=a_{12}=a_{22}=1$, $a_{33}=a_{44}=-1$ and $a_{13}=a_{14}=a_{23}=a_{24}=a_{34}=0$. We can see that $(x_1^0,x_2^0,x_3^0,x_4^0)=(1,0,0,1)$ is an integer solution to \eqref{x2pxypy2mz2mt2} with $x_4^0 \neq 0$. 
System \eqref{eq:genhomquadratzero} simplifies to
$$
u_1^2+u_1 u_2+u_2^2-u_3^2=2u_1+u_2=0.
$$
From the second equation of this system, $u_2=-2u_1$. Substituting into the first equation gives $3u_1^2-u_3^2=0$, whose only integer solution is $(u_1,u_3)=(0,0)$. Therefore all solutions with $t \neq 0$ are given by \eqref{eq:genhomquadsoln}-\eqref{eq:genhomquadsol}, which for the given equation reduces to
\begin{equation}\label{x2pxypy2mz2mt2_sol}
\begin{aligned}
x=&\frac{u}{q}(-u_1^2+u_2^2-u_3^2), \\
y=&\frac{u}{q}(-u_2^2-2u_1 u_2), \\
z=&\frac{u}{q}(-2u_1 u_3 -u_2 u_3), \\
t=&\frac{u}{q}(u_1^2+u_1 u_2+u_2^2-u_3^2), \\
\end{aligned}
\end{equation}
where $u,u_1,u_2,u_3$ are integer parameters and $q$ is any common divisor of the expressions in the parentheses. 
To find integer solutions to \eqref{x2pxypy2mz2mt2} with $t=0$, we can repeat this method for the resulting equation $x^2+xy+y^2-z^2=0$ in three variables. In this case, $x_1=x$, $x_2=y$, $x_3=z$, $a_{11}=a_{12}=a_{22}=1$, $a_{13}=a_{23}=0$ and $a_{33}=-1$. If $x_3=0$, then we have the integer solution $(x,y,z,t)=(0,0,0,0)$. Assuming $x_3 \neq 0$, we can see that $(x_1^0,x_2^0,x_3^0)=(0,1,1)$ is an integer solution with $x_3^0 \neq 0$. System \eqref{eq:genhomquadratzero} has no non-zero integer solutions, hence all integer solutions are given by \eqref{eq:genhomquadsoln}-\eqref{eq:genhomquadsol}. Then applying formulas \eqref{eq:genhomquadsoln} and \eqref{eq:genhomquadsol}, we obtain that the integer solutions with $t=0$ are
\begin{equation}\label{x2pxypy2mz2mt2_sol2}
\begin{aligned}
x=&\frac{u}{q}(-u_1^2-2u_1 u_2), \\
y=&\frac{u}{q}(u_1^2-u_2^2), \\
z=&\frac{u}{q}(u_1^2+u_1u_2+ u_2^2), \\
t=& \, \, 0,\\
\end{aligned}
\end{equation}  
where $u,u_1,u_2,u_3$ are integer parameters and $q$ is any common divisor of the expressions in the parentheses. 
In conclusion, all integer solutions to equation \eqref{x2pxypy2mz2mt2} are given by either \eqref{x2pxypy2mz2mt2_sol} or \eqref{x2pxypy2mz2mt2_sol2}. 

\vspace{10pt}

The next equation we will consider is
\begin{equation}\label{x2pxymy2pz2pt2}
x^2+xy-y^2+z^2+t^2=0.
\end{equation}
This is an equation of the form \eqref{quad_hom} with $n=4$, $x_1=x$, $x_2=z$, $x_3=t$, $x_4=y$, $a_{11}=a_{14}=a_{22}=a_{33}=1$, $a_{44}=-1$ and $a_{12}=a_{13}=a_{23}=a_{24}=a_{34}=0$. If $x_4=y=0$ then we have the solution $(x,y,z,t)=(0,0,0,0)$. Now, let us assume that $x_4 \neq 0$. We can see that $(x_1^0,x_2^0,x_3^0,x_4^0)=(0,1,0,1)$ is an integer solution to \eqref{x2pxymy2pz2pt2} with $x_4^0 \neq 0$. System \eqref{eq:genhomquadratzero} has no non-zero integer solutions, hence all integer solutions are given by \eqref{eq:genhomquadsoln}-\eqref{eq:genhomquadsol}, which in our case reduces to 
\begin{equation}\label{x2pxymy2pz2pt2_sol}
\begin{aligned}
x=&\frac{u}{q}(-2u_1 u_2 -u_1^2), \\
y=&\frac{u}{q}(u_1^2+u_2^2+u_3^2), \\
z=&\frac{u}{q}(u_1^2-u_2^2+u_3^2-u_1 u_2), \\
t=&\frac{u}{q}(-2u_2 u_3 -u_1 u_3), \\
\end{aligned}
\end{equation}
where $u,u_1,u_2,u_3$ are integer parameters and $q$ is any common divisor of the expressions in the parentheses. The solution $(x,y,z,t)=(0,0,0,0)$ is included in \eqref{x2pxymy2pz2pt2_sol} with $u=0$. In conclusion, all integer solutions to equation \eqref{x2pxymy2pz2pt2} are given by \eqref{x2pxymy2pz2pt2_sol}.

\vspace{10pt}

The next equation we will consider is
\begin{equation}\label{x2p2xypxzmz2}
x^2+2xy+xz-z^2=0.
\end{equation}
This is an equation of the form \eqref{quad_hom} with $n=3$, $x_1=x$, $x_2=z$, $x_3=y$, $a_{11}=a_{12}=1$, $a_{22}=-1$, $a_{13}=2$ and $a_{23}=a_{33}=0$. If $x_3=y=0$ then we have the solution $(x,y,z)=(0,0,0)$. Now, let us assume that $x_3 \neq 0$. We can see that $(x_1^0,x_2^0,x_3^0)=(0,0,1)$ is an integer solution to \eqref{x2p2xypxzmz2} with $x_3^0 \neq 0$. System \eqref{eq:genhomquadratzero} has no non-zero integer solutions, hence all integer solutions are given by \eqref{eq:genhomquadsoln}-\eqref{eq:genhomquadsol}, which in our case reduces to 
\begin{equation}\label{x2p2xypxzmz2_sol}
\begin{aligned}
x=&\frac{u}{q}(-2u_1^2), \\
y=&\frac{u}{q}(u_1^2+u_1 u_2 -u_2^2), \\
z=&\frac{u}{q}(-2u_1 u_2), \\
\end{aligned}
\end{equation}
where $u,u_1,u_2$ are integer parameters and $q$ is any common divisor of the expressions in the parentheses. The solution $(x,y,z)=(0,0,0)$ is included in \eqref{x2p2xypxzmz2_sol} with $u=0$. Let us now attempt to simplify this solution. As $q$ divides $-2u_1^2$, $u_1^2+u_1 u_2 -u_2^2$ and $-2u_1 u_2$, then
$$
q \mid (-2u_1^2)+2(u_1^2+u_1 u_2 -u_2^2)+(-2u_1 u_2)=-2u_2^2.
$$
As $u_1$ and $u_2$ are coprime, $q=2$ or $q=1$. If $q=2$ then at least one of $u_1$ and $u_2$ are odd, however, $u_1^2+u_1 u_2 -u_2^2$ is then always odd, and so not divisible by $q=2$. Hence, $q=1$ and therefore we can simplify \eqref{x2p2xypxzmz2_sol} to 
\begin{equation}\label{x2p2xypxzmz2_soli}
\begin{aligned}
x=& u(-2u_1^2), \\
y=& u(u_1^2+u_1 u_2 -u_2^2), \\
z=& u(-2u_1 u_2), \quad u,u_1,u_2 \in \mathbb{Z}. \\
\end{aligned}
\end{equation}
In conclusion, all integer solutions to equation \eqref{x2p2xypxzmz2} are given by \eqref{x2p2xypxzmz2_soli}.

\vspace{10pt}

The next equation we will consider is
\begin{equation}\label{x2p2xypxzpz2}
x^2+2xy+xz+z^2=0.
\end{equation}
This is an equation of the form \eqref{quad_hom} with $n=3$, $x_1=x$, $x_2=z$, $x_3=y$, $a_{11}=a_{12}=a_{22}=1$, $a_{13}=2$ and $a_{23}=a_{33}=0$. If $x_3=y=0$ then we have the solution $(x,y,z)=(0,0,0)$. Now, let us assume that $x_3 \neq 0$. We can see that $(x_1^0,x_2^0,x_3^0)=(0,0,1)$ is an integer solution to \eqref{x2p2xypxzpz2} with $x_3^0 \neq 0$. 
System \eqref{eq:genhomquadratzero} has no non-zero integer solutions, hence all integer solutions are given by \eqref{eq:genhomquadsoln}-\eqref{eq:genhomquadsol}, which in our case reduces to
\begin{equation}\label{x2p2xypxzpz2_sol}
\begin{aligned}
x=&\frac{u}{q}(-2u_1^2), \\
y=&\frac{u}{q}(u_1^2+u_1 u_2 +u_2^2), \\
z=&\frac{u}{q}(-2u_1 u_2), \\
\end{aligned}
\end{equation}
where $u,u_1,u_2$ are integer parameters and $q$ is any common divisor of the expressions in the parentheses. The solution $(x,y,z)=(0,0,0)$ is included in \eqref{x2p2xypxzpz2_sol} with $u=0$. Let us simplify this solution. As $q$ divides $-2u_1^2$, $u_1^2+u_1 u_2 +u_2^2$ and $-2u_1 u_2$, then
$$
q \mid (-2u_1^2)+2(u_1^2+u_1 u_2 +u_2^2)+(-2u_1 u_2)=2u_2^2.
$$
As $u_1$ and $u_2$ are coprime, $q=2$ or $q=1$. If $q=2$ then at least one of $u_1$ and $u_2$ are odd, however, $u_1^2+u_1 u_2 +u_2^2$ is then always odd, and so not divisible by $q=2$. Hence, $q=1$ and therefore we can simplify \eqref{x2p2xypxzpz2_sol} to 
\begin{equation}\label{x2p2xypxzpz2_soli}
\begin{aligned}
x=& u(-2u_1^2), \\
y=& u(u_1^2+u_1 u_2 +u_2^2), \\
z=& u(-2u_1 u_2), \quad u,u_1,u_2 \in \mathbb{Z}. \\
\end{aligned}
\end{equation}
In conclusion, all integer solutions to equation \eqref{x2p2xypxzpz2} are given by \eqref{x2p2xypxzpz2_soli}.

\vspace{10pt}

The next equation we will consider is
\begin{equation}\label{2x2my2myzmz2}
2x^2-y^2-yz-z^2=0.
\end{equation}
This is an equation of the form \eqref{quad_hom}. We can see that $(x,y,z)=(0,0,0)$ is a solution to this equation. We can use Theorem \ref{theorem_cassel} to find a non-trivial integer solution to this equation. In this case, $F=5$ and $n=3$, so we must check the cases $\max(|x|,|y|,|z|) \leq 60$. We can do this with a computer search, however, this returns no non-trivial integer solutions. Therefore, we can conclude that the only integer solution to equation \eqref{2x2my2myzmz2} is $(x,y,z)=(0,0,0)$.

\vspace{10pt}

The next equation we will consider is
\begin{equation}\label{2x2py2myzmz2}
2x^2+y^2-yz-z^2=0.
\end{equation}
This is an equation of the form \eqref{quad_hom}. We can see that $(x,y,z)=(0,0,0)$ is a solution to this equation. We can use Theorem \ref{theorem_cassel} to find a non-trivial integer solution to this equation. In this case, $F=5$ and $n=3$, so we must check the cases $\max(|x|,|y|,|z|) \leq 60$. We can do this with a computer search, however, this returns no non-trivial integer solutions. Therefore, we can conclude that the only integer solution to equation \eqref{2x2py2myzmz2} is $(x,y,z)=(0,0,0)$.

\vspace{10pt}

The next equation we will consider is
\begin{equation}\label{2x2pxymy2mz2}
	2x^2+xy-y^2-z^2=0.
\end{equation}
This is an equation of the form \eqref{quad_hom} with $n=3$, $x_1=z$, $x_2=y$, $x_3=x$, $a_{11}=a_{22}=-1$, $a_{33}=2$, $a_{12}=a_{13}=0$ and $a_{23}=1$. If $x_3=x=0$ then we have the solution $(x,y,z)=(0,0,0)$. Now, let us assume that $x_3 \neq 0$. We can see that $(x_1^0,x_2^0,x_3^0)=(0,-2,-1)$ is an integer solution to \eqref{2x2pxymy2mz2} with $x_3^0 \neq 0$. System \eqref{eq:genhomquadratzero} has no non-zero integer solutions, hence all integer solutions are given by \eqref{eq:genhomquadsoln}-\eqref{eq:genhomquadsol}, which in our case reduces to 
\begin{equation}\label{2x2pxymy2mz2_sol}
	\begin{aligned}
		x=&\frac{u}{q}(u_1^2+u_2^2), \\
		y=&\frac{u}{q}(2u_1^2 - u_2^2), \\
		z=&\frac{u}{q}(-3 u_1 u_2), \\
	\end{aligned}
\end{equation}
where $u,u_1,u_2$ are integer parameters and $q$ is any common divisor of the expressions in the parentheses. The solution $(x,y,z)=(0,0,0)$ is included in \eqref{2x2pxymy2mz2_sol} with $u=0$. Let us simplify this solution. As $q$ divides $u_1^2+u_2^2$ and $2u_1^2-u_2^2$, then
$$
q \mid 2(u_1^2+u_2^2)-(2u_1^2-u_2^2)=3u_2^2, \quad \text{and} \quad q \mid (u_1^2+u_2^2)+(2u_1^2-u_2^2)=3u_1^2.
$$
But $u_1$ and $u_2$ are coprime so $q=3$ or $q=1$. If $q=3$ then $3$ must divide $u_1^2+u_2^2$, however, $u_1^2+u_2^2 \not\equiv 0$ $($mod $3)$. Hence, $q=1$ and therefore we can simplify \eqref{2x2pxymy2mz2_sol} to 
\begin{equation}\label{2x2pxymy2mz2_soli}
	\begin{aligned}
		x=& u(u_1^2+u_2^2), \\
		y=& u(2u_1^2 - u_2^2), \\
		z=& u(-3 u_1 u_2), \quad u,u_1,u_2 \in \mathbb{Z}. \\
	\end{aligned}
\end{equation}
In conclusion, all integer solutions to equation \eqref{2x2pxymy2mz2} are given by \eqref{2x2pxymy2mz2_soli}.

\vspace{10pt}

The next equation we will consider is
\begin{equation}\label{2x2pxymy2pz2}
2x^2+xy-y^2+z^2=0.
\end{equation}
This is an equation of the form \eqref{quad_hom} with $n=3$, $x_1=x$, $x_2=z$, $x_3=y$, $a_{11}=2$, $a_{12}=a_{23}=0$, $a_{22}=a_{13}=1$, and $a_{33}=-1$. If $x_3=y=0$ then we have the solution $(x,y,z)=(0,0,0)$. Now, let us assume that $x_3 \neq 0$. We can see that $(x_1^0,x_2^0,x_3^0)=(0,1,1)$ is an integer solution to \eqref{2x2pxymy2pz2} with $x_3^0 \neq 0$. System \eqref{eq:genhomquadratzero} has no non-zero integer solutions, hence all integer solutions are given by \eqref{eq:genhomquadsoln}-\eqref{eq:genhomquadsol}, which in our case reduces to 
\begin{equation}\label{2x2pxymy2pz2_sol}
\begin{aligned}
x=&\frac{u}{q}(-u_1^2-2 u_1 u_2), \\
y=&\frac{u}{q}(2u_1^2+u_2^2), \\
z=&\frac{u}{q}(2u_1^2 -u_2^2-u_1 u_2), \\
\end{aligned}
\end{equation}
where $u,u_1,u_2$ are integer parameters and $q$ is any common divisor of the expressions in the parentheses. The solution $(x,y,z)=(0,0,0)$ is included in \eqref{2x2pxymy2pz2_sol} with $u=0$. In conclusion, all integer solutions to equation \eqref{2x2pxymy2pz2} are given by \eqref{2x2pxymy2pz2_sol}.

\vspace{10pt}

The final equation we will consider is
\begin{equation}\label{2x2pxypy2mz2}
	2x^2+xy+y^2-z^2=0.
\end{equation}
This is an equation of the form \eqref{quad_hom} with $n=3$, $x_1=x$, $x_2=z$, $x_3=y$, $a_{11}=2$, $a_{12}=a_{23}=0$, $a_{22}=-1$, and $a_{13}=a_{33}=1$. If $x_3=y=0$ then we have the solution $(x,y,z)=(0,0,0)$. Now, let us assume that $x_3 \neq 0$. We can see that $(x_1^0,x_2^0,x_3^0)=(0,1,1)$ is an integer solution to \eqref{2x2pxypy2mz2} with $x_3^0 \neq 0$. System \eqref{eq:genhomquadratzero} has no non-zero integer solutions, hence all integer solutions are given by \eqref{eq:genhomquadsoln}-\eqref{eq:genhomquadsol}, which in our case reduces to 
\begin{equation}\label{2x2pxypy2mz2_sol}
	\begin{aligned}
		x=&\frac{u}{q}(-u_1^2+2 u_1 u_2), \\
		y=&\frac{u}{q}(2u_1^2-u_2^2), \\
		z=&\frac{u}{q}(2u_1^2 +u_2^2-u_1 u_2), \\
	\end{aligned}
\end{equation}
where $u,u_1,u_2$ are integer parameters and $q$ is any common divisor of the expressions in the parentheses. The solution $(x,y,z)=(0,0,0)$ is included in \eqref{2x2pxypy2mz2_sol} with $u=0$. In conclusion, all integer solutions to equation \eqref{2x2pxypy2mz2} are given by \eqref{2x2pxypy2mz2_sol}.

\begin{center}

	\captionof{table}{Integer solutions to the equations in Table \ref{tab:H20hom}.  Assume that $u,u_1,u_2,u_3 \in \mathbb{Z}$, and $q$ is any common divisor of the expressions in brackets.  \label{tb:homogeneous_quad}}
\end{center}

\subsection{Exercise 4.24}\label{ex:H24hom}
\textbf{\emph{Describe all integer solutions to the equations listed in Table \ref{tab:H24hom}. }}

\begin{center}
%
\captionof{table}{\label{tab:H24hom} Homogeneous quadratic equations in diagonal form \eqref{eq:quadhomdiag} of size $H\leq 24$.}
\end{center} 

All equations in Table \ref{tab:H24hom} are 
homogeneous quadratic equations in diagonal form, that is, equations of the form
\begin{equation}\label{eq:quadhomdiag}
\sum_{i=1}^n a_i x_i^2=0.
\end{equation}
We will solve these by a similar method to the one in Section \ref{ex:H20hom}.  
For equations of the form \eqref{eq:quadhomdiag}, formulas \eqref{eq:genhomquadsoln} and \eqref{eq:genhomquadsol} simplify to 
\begin{equation}\label{quad_hom_diag_xn}
x_k = \frac{u}{q}\left( x_k^0 \sum_{i=1}^{n-1} a_i u_i ^2 - 2u_k \sum_{i=1}^{n-1} a_i x_i^0 u_i \right), \quad k=1,\dots,n-1, \quad x_n = \frac{u}{q}\left( x_n^0 \sum_{i=1}^{n-1} a_i u_i^2 \right),
\end{equation}
and system \eqref{eq:genhomquadratzero} simplifies to
\begin{equation}\label{eq:genhomquadratzero_diag}
	\sum_{i=1}^{n-1}  a_{i} u_i^2 = 2 \sum_{i=1}^{n-1}  a_{i} x_i^0 u_i  = 0.
\end{equation}

Equation
$$
2x^2-y^2-z^2-t^2=0
$$
is solved in Section 4.2.3 of the book, and its integer solutions are
$$
\begin{aligned}
	x=&\frac{u}{q}(u_1^2+u_2^2+u_3^2), \\
	y=&\frac{u}{q}(u_1^2+u_2^2+u_3^2-2u_1(u_1+u_2)), \\
	z=&\frac{u}{q}(u_1^2 +u_2^2+u_3^2-2u_2(u_1 +u_2)), \\
	t=&\frac{u}{q}(-2u_3(u_1+u_2)),
\end{aligned}
$$
where $u,u_1,u_2$ are integer parameters and $q$ is any common divisor of the expressions in the parentheses. 

Equation
$$
2x^2+2y^2-z^2-t^2=0
$$
is solved in Section 4.2.3 of the book, and its integer solutions are
$$
\begin{aligned}
(x,y,z,t)=& (\pm u(v^2+2vw-w^2),u(w^2+2vw-v^2),2u(w^2+v^2),0), \quad \text{or} \\
& (w+v,\pm v,w+2v,w), \quad u,v,w \in \mathbb{Z},
\end{aligned}
$$
or
$$
\begin{aligned}
	x=&\frac{u}{q}(-2u_1^2+2u_1 u_3+2u_2^2-u_3^2), \\
	y=&\frac{u}{q}(-4u_1 u_2+2u_2 u_3), \\
	z=&\frac{u}{q}(2u_1^2 +2u_2^2+u_3^2-4u_1 u_3), \\
	t=&\frac{u}{q}(2u_1^2+2u_2^2-u_3^2),
\end{aligned}
$$
where $u,u_1,u_2$ are integer parameters and $q$ is any common divisor of the expressions in the parentheses. 

Let us consider the first equation
\begin{equation}\label{3x2my2mz2}
3x^2-y^2-z^2=0.
\end{equation}
We can see that this equation has the trivial solution $(x,y,z)=(0,0,0)$. Let us find a non-trivial integer solution using Theorem \ref{theorem_cassel}. In this case, $F=5$ and $n=3$, so we must check the cases $\max(|x|,|y|,|z|)\leq  60$. We can do this with a computer search, however, this returns no non-trivial integer solutions. Therefore, we can conclude that the only integer solution to equation \eqref{3x2my2mz2} is $(x,y,z)=(0,0,0)$. 

\vspace{10pt}

The next equation we will consider is
\begin{equation}\label{3x2m2y2mz2}
3x^2-2y^2-z^2=0.
\end{equation}
This is an equation of the form \eqref{eq:quadhomdiag} with $n=3$, $x_1=y$, $x_2=z$, $x_3=x$, $a_{1}=2$, $a_{2}=1$ and $a_{3}=-3$. If $x_3=x=0$ then we have the solution $(x,y,z)=(0,0,0)$. Now, let us assume that $x_3 \neq 0$. We can see that $(x_1^0,x_2^0,x_3^0)=(1,1,1)$ is an integer solution to \eqref{3x2m2y2mz2} with $x_3^0 \neq 0$. System \eqref{eq:genhomquadratzero_diag} has no non-zero integer solutions, hence all integer solutions are given by \eqref{quad_hom_diag_xn}, which in our case reduces to 
\begin{equation}\label{3x2m2y2mz2_sol}
\begin{aligned}
x=&\frac{u}{q}(2u_1^2+u_2^2), \\
y=&\frac{u}{q}(-2u_1^2-2 u_1 u_2+u_2^2), \\
z=&\frac{u}{q}(2u_1^2 -u_2^2-4u_1 u_2), \\
\end{aligned}
\end{equation}
where $u,u_1,u_2$ are integer parameters and $q$ is any common divisor of the expressions in the parentheses. The solution $(x,y,z)=(0,0,0)$ is included in \eqref{3x2m2y2mz2_sol} with $u=0$. In conclusion, all integer solutions to equation \eqref{3x2m2y2mz2} are given by \eqref{3x2m2y2mz2_sol}.

\vspace{10pt}

The next equation we will consider is
\begin{equation}\label{3x2m2y2pz2}
3x^2-2y^2+z^2=0.
\end{equation}
We can see that this equation has the trivial solution $(x,y,z)=(0,0,0)$. Let us find a non-trivial solution using Theorem \ref{theorem_cassel}. In this case, $F=6$ and $n=3$, so we must check the cases $\max(|x|,|y|,|z|)\leq  72$. We can do this with a computer search, however, this returns no non-trivial integer solutions. Therefore, we can conclude that the only integer solution to equation \eqref{3x2m2y2pz2} is $(x,y,z)=(0,0,0)$. 

\vspace{10pt}

The next equation we will consider is
\begin{equation}\label{3x2p2y2mz2}
3x^2+2y^2-z^2=0.
\end{equation}
We can see that this equation has the trivial solution $(x,y,z)=(0,0,0)$. Let us find a non-trivial integer solution using Theorem \ref{theorem_cassel}. In this case, $F=6$ and $n=3$, so we must check the cases $\max(|x|,|y|,|z|)\leq  72$. We can do this with a computer search, however, this returns no non-trivial integer solutions. Therefore, we can conclude that the only integer solution to equation \eqref{3x2p2y2mz2} is $(x,y,z)=(0,0,0)$. 

\vspace{10pt}

The next equation we will consider is
\begin{equation}\label{2x2m2y2pz2pt2}
2x^2-2y^2+z^2+t^2=0.
\end{equation}
This is an equation of the form \eqref{eq:quadhomdiag} with $n=4$, $x_1=x$, $x_2=z$, $x_3=t$, $x_4=y$, $a_{1}=2$, $a_{2}=a_3=1$ and $a_{4}=-2$. If $x_4=y=0$ then we have the solution $(x,y,z,t)=(0,0,0,0)$. We can see that $(x_1^0,x_2^0,x_3^0,x_4^0)=(1,0,0,-1)$ is an integer solution to \eqref{2x2m2y2pz2pt2} with $x_4^0 \neq 0$. System \eqref{eq:genhomquadratzero_diag} has no non-zero integer solutions, hence all integer solutions are given by \eqref{quad_hom_diag_xn}, which in our case reduces to 
\begin{equation}\label{2x2m2y2pz2pt2_sol}
\begin{aligned}
x=&\frac{u}{q}(-2u_1^2+u_2^2+u_3^2), \\
y=&\frac{u}{q}(-2u_1^2-u_2^2-u_3^2), \\
z=&\frac{u}{q}(-4u_1 u_2), \\
t=&\frac{u}{q}(-4 u_1 u_3), \\
\end{aligned}
\end{equation}
where $u,u_1,u_2,u_3$ are integer parameters and $q$ is any common divisor of the expressions in the parentheses.  The solution $(x,y,z,t)=(0,0,0,0)$ is included in \eqref{2x2m2y2pz2pt2_sol} with $u=0$. In conclusion, all integer solutions to equation \eqref{2x2m2y2pz2pt2} are given by \eqref{2x2m2y2pz2pt2_sol}.

\vspace{10pt}

The next equation we will consider is
\begin{equation}\label{3x2my2mz2mt2}
3x^2-y^2-z^2-t^2=0.
\end{equation}
This is an equation of the form \eqref{eq:quadhomdiag} with $n=4$, $x_1=y$, $x_2=z$, $x_3=t$, $x_4=x$, $a_{1}=a_2=a_3=1$, and $a_{4}=-3$. If $x_4=x=0$ then we have the solution $(x,y,z,t)=(0,0,0,0)$. Now, let us assume that $x_4 \neq 0$. We can see that $(x_1^0,x_2^0,x_3^0,x_4^0)=(1,1,1,1)$ is an integer solution to \eqref{3x2my2mz2mt2} with $x_4^0 \neq 0$. System \eqref{eq:genhomquadratzero_diag} has no non-zero integer solutions, hence all integer solutions are given by \eqref{quad_hom_diag_xn}, which in our case reduces to 
\begin{equation}\label{3x2my2mz2mt2_sol}
\begin{aligned}
x=&\frac{u}{q}(u_1^2+u_2^2+u_3^2), \\
y=&\frac{u}{q}(-u_1^2-2 u_1 u_2+u_2^2+u_3^2-2u_1 u_3), \\
z=&\frac{u}{q}(u_1^2 -u_2^2+u_3^2-2u_1 u_2-2u_2 u_3), \\
t=&\frac{u}{q}(u_1^2+u_2^2-u_3^2-2u_1 u_3 - 2u_2 u_3), \\
\end{aligned}
\end{equation}
where $u,u_1,u_2,u_3$ are integer parameters and $q$ is any common divisor of the expressions in the parentheses. The solution $(x,y,z,t)=(0,0,0,0)$ is included in \eqref{3x2my2mz2mt2_sol} with $u=0$. In conclusion, all integer solutions to equation \eqref{3x2my2mz2mt2} are given by \eqref{3x2my2mz2mt2_sol}.

\vspace{10pt}

The final equation we will consider is
\begin{equation}\label{2x52mx12mx22mx32mx42}
2x_5^2-x_1^2-x_2^2-x_3^2-x_4^2=0.
\end{equation}
This is an equation of the form \eqref{eq:quadhomdiag} with $n=5$, $a_{1}=a_2=a_3=a_4=1$, and $a_{5}=-2$. If $x_5=0$ then we have the solution $(x_1,x_2,x_3,x_4,x_5)=(0,0,0,0,0)$. Now, let us assume that $x_5 \neq 0$. We can see that $(x_1^0,x_2^0,x_3^0,x_4^0,x_5^0)=(1,1,0,0,1)$ is an integer solution to \eqref{2x52mx12mx22mx32mx42} with $x_5^0 \neq 0$. System \eqref{eq:genhomquadratzero_diag} has no non-zero integer solutions, hence all integer solutions are given by \eqref{quad_hom_diag_xn}, which in our case reduces to 
\begin{equation}\label{2x52mx12mx22mx32mx42_sol}
\begin{aligned}
x_1=&\frac{u}{q}(-u_1^2+u_2^2+u_3^2+u_4^2-2u_1 u_2), \\
x_2=&\frac{u}{q}(u_1^2 -u_2^2+u_3^2+u_4^2-2u_1 u_2), \\
x_3=&\frac{u}{q}(-2u_1 u_3 - 2u_2 u_3), \\
x_4=&\frac{u}{q}(-2u_1 u_4 - 2u_2 u_4), \\
x_5=&\frac{u}{q}(u_1^2+u_2^2+u_3^2+u_4^2), \\
\end{aligned}
\end{equation}
where $u,u_1,u_2,u_3,u_4$ are integer parameters and $q$ is any common divisor of the expressions in the parentheses. The solution $(x_1,x_2,x_3,x_4,x_5)=(0,0,0,0,0)$ is included in \eqref{2x52mx12mx22mx32mx42_sol} with $u=0$. Therefore we can describe all integer solutions to equation \eqref{2x52mx12mx22mx32mx42}  by \eqref{2x52mx12mx22mx32mx42_sol}. 

 \begin{center}
\begin{tabular}{|c|c|c|}
\hline
Equation & Solution $(x,y,z)$ or $(x,y,z,t)$ \\ \hline \hline
 $3x^2-y^2-z^2=0$ &$(0,0,0)$ \\\hline
 $2x^2-y^2-z^2-t^2=0$ & $x=\frac{u}{q}(u_1^2+u_2^2+u_3^2)$ \\
& $y=\frac{u}{q}(-u_1^2+u_2^2+u_3^2-2u_1 u_2)$ \\
& $z=\frac{u}{q}(u_1^2-u_2^2+u_3^2-2u_1 u_2)$ \\
& $t=\frac{u}{q}(-2u_1 u_3 -2 u_2 u_3)$ \\\hline
 $3x^2-2y^2-z^2=0$ & $x=\frac{u}{q}(2u_1^2+u_2^2)$ \\
& $y=\frac{u}{q}(-2u_1^2+u_2^2-2u_1 u_2)$ \\
& $z=\frac{u}{q}(2u_1^2-4u_1 u_2 -u_2^2)$ \\\hline
 $3x^2-2y^2+z^2=0$ & $(0,0,0)$ \\\hline
 $3x^2+2y^2-z^2=0$ & $(0,0,0)$ \\\hline
 $2x^2+2y^2-z^2-t^2=0$ & $x=\frac{u}{q}(-2u_1^2+2u_2^2-u_3^2+2u_1 u_3)$ \\
& $y=\frac{u}{q}(-4u_1 u_2+2u_2 u_3)$ \\
& $z=\frac{u}{q}(2u_1^2-4u_1 u_3 +2u_2^2+u_3^2)$ \\
& $t=\frac{u}{q}(2u_1^2+2u_2^2-u_3^2)$ \\
 & $(\pm u(u_1^2+2u_1 u_2-u_2^2),u(u_2^2+2 u_1 u_2-u_1^2),2u(u_2^2+u_1^2),0)$ \\
 & $(w+v,\pm v,w+2v,w), \quad v,w \in \mathbb{Z}$ \\\hline
 $2x^2-2y^2+z^2+t^2=0$ & $x=\frac{u}{q}(-2u_1^2+u_2^2+u_3^2)$ \\
& $y=\frac{u}{q}(-2u_1^2-u_2^2-u_3^2)$ \\
& $z=\frac{u}{q}(-4u_1 u_2)$ \\
& $t=\frac{u}{q}(-4u_1 u_3)$ \\\hline
 $3x^2-y^2-z^2-t^2=0$ &  $x=\frac{u}{q}(u_1^2+u_2^2+u_3^2)$ \\
& $y=\frac{u}{q}(-u_1^2+u_2^2+u_3^2-2u_1 u_2-2u_1 u_3)$ \\
& $z=\frac{u}{q}(u_1^2-u_2^2+u_3^2-2u_1 u_2-2u_2 u_3)$ \\
& $t=\frac{u}{q}(u_1^2+u_2^2-u_3^2-2u_1 u_3-2u_2 u_3)$\\\hline
$2x_5^2-x_1^2-x_2^2-x_3^2-x_4^2=0$ & $x_1=\frac{u}{q}(-u_1^2+u_2^2+u_3^2+u_4^2-2u_1 u_2)$ \\
& $x_2=\frac{u}{q}(u_1^2-u_2^2+u_3^2+u_4^2-2u_1 u_2)$ \\
& $x_3=\frac{u}{q}(-2u_1 u_3-2u_2 u_3)$ \\
& $x_4=\frac{u}{q}(-2u_1 u_4-2u_2 u_4)$ \\
& $x_5=\frac{u}{q}(u_1^2+u_2^2+u_3^2+u_4^2)$ \\\hline

 \end{tabular}
 \captionof{table}{Integer solutions to the equations listed in Table \ref{tab:H24hom}. Assume that $u,u_1,\dots,u_{n-1} \in \mathbb{Z}$, and $q$ is any common divisor of the expressions in brackets.  \label{tb:homogeneous_quad_diagonal}}
 \end{center}

\subsection{Exercise 4.26}\label{ex:H18quadsol}
\textbf{\emph{Solve all equations listed in Table \ref{tab:H18quadred}. For every equation, present the answer in as nice a form as you can.}}
 
For each equation in Table \ref{tab:H18quadred}, we will describe solutions either using $S$ notation, a finite union of polynomial families, or showing that solutions can be written as a finite union of polynomial families subject to some conditions. By showing solutions in a variety of ways, we can let the reader decide which way they think is the ``nicest'' description of solutions. 

We will first look at describing solutions with $S$ notation as defined in \eqref{S Notation}. Let us consider the first equation \eqref{xpx2my2mz2m3}.
We can rearrange this equation to 
$$
x+x^2=y^2+z^2+3.
$$
To make the equation linear in $x$ we can use the substitution $y=y_1-x$, which reduces the equation to 
$$
x(2y_1+1)=z^2+y_1^2+3.
$$
This equation is of the form $P(x_3,\dots,x_n)=Q(x_2,\dots,x_n)$ with $x_1=x$, $x_2=z$, $P=2y_1+1$, $Q=1$, $R=0$, $T=y_1^2+3$. This equation has an integer solution $(x,y_1,z)=(4,0,1)$ and therefore by Proposition \ref{prop:finlincheck} this equation has infinitely many integer solutions. Then, we can solve this equation using the method in Section \ref{ex:H16modularroot}, to obtain that the integer solutions to equation \eqref{xpx2my2mz2m3} are 
\begin{equation}\label{xpx2my2mz2m3_sol}
\begin{aligned}
(x,y,z)=& \left( (2y_1+1)u^2+ur+\frac{r^2+4y_1^2+12}{4(2y_1+1)}, y_1-(2y_1+1)u^2-ur-\frac{r^2+4y_1^2+12}{4(2y_1+1)}, 
 (2y_1+1)u+\frac{r}{2}\right) \\  & r \in S(4(2y_1+1),-4(y_1^2+3)), \quad u,y_1 \in \mathbb{Z}, \\     
\end{aligned}
 \end{equation}
subject to the additional condition that $ r \equiv \, 0 \, \text{(mod 2)}$. However, as $r^2+4y_1^2+12$ must be even to be divisible by $4(2y_1+1)$, then $r^2$ must be even and so $r$ is even, therefore this restriction can be omitted. In summary, all integer solutions to equation \eqref{xpx2my2mz2m3} are given by \eqref{xpx2my2mz2m3_sol}.
 
 \vspace{10pt}
 
 Let us consider the next equation
\begin{equation}\label{xpx2my2mz2p3}
x+x^2-y^2-z^2+3=0.
\end{equation}
We can rearrange this equation to 
$$
x+x^2=y^2+z^2-3.
$$
To make the equation linear in $x$ we can use the substitution $y=y_1-x$, which reduces the equation to 
$$
x(2y_1+1)=z^2+y_1^2-3.
$$
This equation is of the form $P(x_3,\dots,x_n)=Q(x_2,\dots,x_n)$ with $x_1=x$, $x_2=z$, $P=2y_1+1$, $Q=1$, $R=0$, $T=y_1^2-3$. This equation has an integer solution $(x,y_1,z)=(-2,0,1)$ and therefore by Proposition \ref{prop:finlincheck} this equation has infinitely many integer solutions. Then, we can solve this equation using the method in Section \ref{ex:H16modularroot}, to obtain that the integer solutions to equation \eqref{xpx2my2mz2p3} are 
\begin{equation}\label{xpx2my2mz2p3_sol}
\begin{aligned}
x=& \, (2y_1+1)u^2+ur+\frac{r^2+4y_1^2-12}{4(2y_1+1)}, \\
  y= & \, y_1-(2y_1+1)u^2-ur-\frac{r^2+4y_1^2-12}{4(2y_1+1)}, \\
  z= & \, (2y_1+1)u+\frac{r}{2}, \\
  & r \in S(4(2y_1+1),-4(y_1^2-3)), \quad y_1,u \in \mathbb{Z},
  \end{aligned}
 \end{equation}
 subject to the restriction that $r \equiv \, 0 \, \text{(mod 2)}$. However, as $r^2+4y_1^2-12$ must be even to be divisible by $4(2y_1+1)$, then $r^2$ must be even and so $r$ is even, hence, we do not need the restriction that $r \equiv 0 \, \text{(mod 2)}$. Therefore, we can describe all integer solutions to equation \eqref{xpx2my2mz2p3} by \eqref{xpx2my2mz2p3_sol}.

\vspace{10pt}
 
Let us consider the next equation \eqref{xpx2pymy2pz2m1}.
We can rearrange this equation to 
$$
x+x^2=y^2-y-z^2+1.
$$
To make the equation linear in $x$ we can use the substitution $y=y_1-x$, which reduces the equation to 
$$
x(2y_1)=y_1^2-y_1-z^2+1.
$$
This equation is of the form $P(x_3,\dots,x_n)=Q(x_2,\dots,x_n)$ with $x_1=x$, $x_2=z$, $P=2y_1$, $Q=-1$, $R=0$, $T=y_1^2-y_1+1$. This equation has an integer solution $(x,y_1,z)=(0,0,1)$ and therefore by Proposition \ref{prop:finlincheck} it has infinitely many integer solutions. Then, we can solve this equation using the method in Section \ref{ex:H16modularroot}, to obtain that the integer solutions to equation \eqref{xpx2pymy2pz2m1} are 
\begin{equation}\label{xpx2pymy2pz2m1_sol}
\begin{aligned}
x=& -2y_1 u^2+ur-\frac{r^2-4y_1^2+4y_1-4}{8y_1}, \\
  y= & \, y_1+2y_1u^2-ur+\frac{r^2-4y_1^2+4y_1-4}{8y_1}, \\
  z=& \, 2y_1 u-\frac{r}{2}, \\
  & r \in S(-8y_1,4(y_1^2-y_1+1)), \quad y_1,u \in \mathbb{Z},
  \end{aligned}
   \end{equation}
subject to the condition that $r \equiv  \, 0 \, \text{(mod 2)}$. However, as $r^2-4(y_1^2-y_1+1)$ must be even to be divisible by $8y_1$, then $r^2$ must be even and so $r$ is even, therefore we do not need the restriction that $r \equiv 0 \, \text{(mod 2)}$. 

We must now consider the case where $P=0$, or eqiuvalently, $2y_1=0$ or $y=-x$. This reduces equation \eqref{xpx2pymy2pz2m1} to $z^2-1=0$, and we obtain the integer solutions 
\begin{equation}\label{xpx2pymy2pz2m1_sol2}
(x,y,z)=(u,-u,\pm 1), \quad u \in \mathbb{Z}.
\end{equation}
Therefore all integer solutions to equation \eqref{xpx2pymy2pz2m1} are given by \eqref{xpx2pymy2pz2m1_sol} or \eqref{xpx2pymy2pz2m1_sol2}.
 
 \vspace{10pt}
 
   Let us consider the next equation
\begin{equation}\label{xpx2pymy2pz2p1}
x+x^2+y-y^2+z^2+1=0.
\end{equation}
We can rearrange this equation to 
$$
x+x^2=y^2-y-z^2-1.
$$
To make the equation linear in $x$ we can use the substitution $y=y_1-x$, which reduces the equation to 
$$
x(2y_1)=y_1^2-y_1-z^2-1.
$$
This equation is of the form $P(x_3,\dots,x_n)=Q(x_2,\dots,x_n)$ with $x_1=x$, $x_2=z$, $P=2y_1$, $Q=-1$, $R=0$, $T=y_1^2-y_1-1$. This equation has an integer solution $(x,y_1,z)=(0,1,0)$ and therefore by Proposition \ref{prop:finlincheck} it has infinitely many integer solutions. Then, we can solve this equation using the method in Section \ref{ex:H16modularroot}, and obtain that the integer solutions to equation \eqref{xpx2pymy2pz2p1} are 
\begin{equation}\label{xpx2pymy2pz2p1_sol}
\begin{aligned}
x=& -2y_1 u^2+ur-\frac{r^2-4y_1^2+4y_1+4}{8y_1}, \\
  y=& \, y_1+2y_1u^2-ur+\frac{r^2-4y_1^2+4y_1+4}{8y_1}, \\
  z=& \, 2y_1 u-\frac{r}{2}, \\
  & r \in S(-8y_1,4(y_1^2-y_1-1)), \quad y_1,u \in \mathbb{Z},
  \end{aligned}
 \end{equation}
 subject to the condition that $r \equiv  \, 0 \, \text{(mod 2)}$. However, as $r^2-4(y_1^2-y_1-1)$ must be even to be divisible by $8y_1$, then $r^2$ must be even and so $r$ is even, therefore, we do not need the restriction that $r \equiv 0 \, \text{(mod 2)}$.
 
 We must now consider the case where $P=0$ or equivalently $2y_1=0$ or $y=-x$. This reduces equation \eqref{xpx2pymy2pz2p1} to 
 $$
 z^2+1=0,
 $$
 which has no integer solutions.
  Therefore, we can describe all integer solutions to equation \eqref{xpx2pymy2pz2p1} by \eqref{xpx2pymy2pz2p1_sol}.

 \vspace{10pt}
 
    Let us consider the next equation
\begin{equation}\label{xpx2pypy2mz2m1}
x+x^2+y+y^2-z^2-1=0.
\end{equation}
We can rearrange this equation to 
$$
x+x^2=z^2-y^2-y+1.
$$
To make the equation linear in $x$, we can use the substitution $z=z_1-x$, which reduces the equation to 
$$
x(1+2z_1)=z_1^2-y^2-y+1.
$$
This equation is of the form $P(x_3,\dots,x_n)=Q(x_2,\dots,x_n)$ with $x_1=x$, $x_2=y$, $P=2z_1+1$, $Q=-1$, $R=-1$, $T=z_1^2+1$. This equation has an integer solution $(x,y,z_1)=(-2,-1,-1)$ and therefore by Proposition \ref{prop:finlincheck} it has infinitely many integer solutions. Then, we can solve this equation using the method in Section \ref{ex:H16modularroot}, to obtain that the integer solutions to equation \eqref{xpx2pypy2mz2m1} are 
\begin{equation}\label{xpx2pypy2mz2m1_sol}
\begin{aligned}
x=& \, -(1+2z_1) u^2+ur-\frac{r^2-4z_1^2-5}{4(1+2z_1)}, \\
y=& \, (1+2z_1) u-\frac{r+1}{2}, \\
  z=& \, z_1+(1+2z_1)u^2-ur+\frac{r^2-4z_1^2-5}{4(1+2z_1)}, \\
  & r \in S(-4(1+2z_1),4z_1^2+5), \quad z_1,u \in \mathbb{Z},
  \end{aligned}
 \end{equation}
 subject to the condition that $r \equiv  \, 1 \, \text{(mod 2)}$. However, as $r^2-4z_1^2-5$ must be even to be divisible by $4(1+2z_1)$, $r$ must be odd, hence, we do not need the restriction that $r \equiv 1 \, \text{(mod 2)}$. Therefore, we can describe all integer solutions to equation \eqref{xpx2pypy2mz2m1} by
\eqref{xpx2pypy2mz2m1_sol}.
 
 \vspace{10pt}
 
Let us consider the next equation
\begin{equation}\label{xpx2pypy2mz2p1}
x+x^2+y+y^2-z^2+1=0.
\end{equation}
We can rearrange this equation to 
$$
x+x^2=z^2-y^2-y-1.
$$
To make the equation linear in $x$, we can use the substitution $z=z_1-x$, which reduces the equation to 
$$
x(1+2z_1)=z_1^2-y^2-y-1.
$$
This equation is of the form $P(x_3,\dots,x_n)=Q(x_2,\dots,x_n)$ with $x_1=x$, $x_2=y$, $P=2z_1+1$, $Q=-1$, $R=-1$, $T=z_1^2-1$. This equation has an integer solution $(x,y,z_1)=(-2,-3,1)$ and therefore by Proposition \ref{prop:finlincheck} it has infinitely many integer solutions. Then, we can solve this equation using the method in Section \ref{ex:H16modularroot}, to obtain that the integer solutions to equation \eqref{xpx2pypy2mz2p1} are
\begin{equation}\label{xpx2pypy2mz2p1_sol}
\begin{aligned}
x=& \, -(1+2z_1) u^2+ur-\frac{r^2-4z_1^2+3}{4(1+2z_1)}, \\
y=& \, (1+2z_1) u-\frac{r+1}{2}, \\
  z=& \, z_1+(1+2z_1)u^2-ur+\frac{r^2-4z_1^2+3}{4(1+2z_1)}, \\
 &  r \in S(-4(1+2z_1),4z_1^2-3), \quad u,z_1 \in \mathbb{Z},
  \end{aligned}
 \end{equation}
 subject to the additional condition that $r \equiv  \, 1 \, \text{(mod 2)}$. However, as $r^2-4z_1^2+3$ must be even to be divisible by $4(1+2z_1)$, $r$ must be odd, hence, we do not need the restriction that $r \equiv 1 \, \text{(mod 2)}$. Therefore, we can describe all integer solutions to equation \eqref{xpx2pypy2mz2p1} by \eqref{xpx2pypy2mz2p1_sol}.
 
 \vspace{10pt}
 
      Let us consider the next equation
\begin{equation}\label{x2pxymy2pz2m1}
x^2+xy-y^2+z^2-1=0.
\end{equation}
We can rearrange this equation to 
$$
xy+x^2=y^2-z^2+1.
$$
To make the equation linear in $x$, we can use the substitutions $z=z_1-x$ and $y=y_1-x$, which reduces the equation to 
$$
x(3y_1-2z_1)=y_1^2-z_1^2+1.
$$
We can then introduce a new integer variables $t$ and $k$, such that $t=y_1-z_1$ and $k=y_1+2t$, which reduces the equation to 
\begin{equation}\label{x2pxymy2pz2m1_red}
xk=2kt-5t^2+1.
\end{equation}
This equation is of the form $P(x_3,\dots,x_n)=Q(x_2,\dots,x_n)$ with $x_1=x$, $x_2=t$, $P=k$, $Q=-5$, $R=2k$, $T=1$. This equation has an integer solution $(x,k,t)=(1,1,0)$ and therefore by Proposition \ref{prop:finlincheck} it has infinitely many integer solutions. Then, we can solve this equation using the method in Section \ref{ex:H16modularroot}, to obtain that the integer solutions to equation \eqref{x2pxymy2pz2m1} are
\begin{equation}\label{x2pxymy2pz2m1_sol}
\begin{aligned}
x=& \, -5k u^2+ur-\frac{r^2-4k^2-20}{20k}, \\
y=& \, k-2ku+5ku^2-ur+\frac{r-2k}{5}+\frac{r^2-4k^2-20}{20k}, \\
  z=& \, k-3ku+5ku^2-ur+\frac{3(r-2k)}{10}+\frac{r^2-4k^2-20}{20k}, \\
  r \equiv & \, 2k \, \text{(mod 10)}, \quad r \in S(-20k,4k^2+20).
  \end{aligned}
\end{equation}
We must now consider the case where $P=0$, or equivalently $k=0$. This reduces equation \eqref{x2pxymy2pz2m1_red} to 
$$
-5t^2+1=0,
$$
which has no integer solutions. 
Therefore, we can describe all integer solutions to equation \eqref{x2pxymy2pz2m1} by \eqref{x2pxymy2pz2m1_sol}. 

\vspace{10pt}

      Let us consider the next equation
\begin{equation}\label{x2pxymy2pz2p1}
x^2+xy-y^2+z^2+1=0.
\end{equation}
We can rearrange this equation to 
$$
xy+x^2=y^2-z^2-1.
$$
To make the equation linear in $x$, we can use the substitutions $z=z_1-x$ and $y=y_1-x$, which reduces the equation to 
$$
x(3y_1-2z_1)=y_1^2-z_1^2-1.
$$
We can then introduce a new integer variables $t$ and $k$, such that $t=y_1-z_1$ and $k=y_1+2t$, which reduces the equation to 
\begin{equation}\label{x2pxymy2pz2p1_red}
xk=2kt-5t^2-1.
\end{equation}
This equation is of the form $P(x_3,\dots,x_n)=Q(x_2,\dots,x_n)$ with $x_1=x$, $x_2=t$, $P=k$, $Q=-5$, $R=2k$, $T=-1$. This equation has an integer solution $(x,k,t)=(-1,1,0)$ and therefore by Proposition \ref{prop:finlincheck} it has infinitely many integer solutions. Then, we can solve this equation using the method in Section \ref{ex:H16modularroot}, to obtain that the integer solutions to equation \eqref{x2pxymy2pz2p1} are
\begin{equation}\label{x2pxymy2pz2p1_sol}
\begin{aligned}
x=& \, -5k u^2+ur-\frac{r^2-4k^2+20}{20k}, \\
y=& \, k-2ku+5ku^2-ur+\frac{r-2k}{5}+\frac{r^2-4k^2+20}{20k}, \\
  z=& \, k-3ku+5ku^2-ur+\frac{3(r-2k)}{10}+\frac{r^2-4k^2+20}{20k}, \\
  r \equiv & \, 2k \, \text{(mod 10)}, \quad r \in S(-20k,4k^2-20).
  \end{aligned}
\end{equation}
We must now consider the case where $P=0$, or equivalently $k=0$. This reduces equation \eqref{x2pxymy2pz2p1_red} to 
$$
-5t^2-1=0,
$$
which has no integer solutions. 
Therefore, we can describe all integer solutions to equation \eqref{x2pxymy2pz2p1} by \eqref{x2pxymy2pz2p1_sol}. 

\vspace{10pt}

      Let us consider the next equation
\begin{equation}\label{x2pxypy2mz2m1}
x^2+xy+y^2-z^2-1=0.
\end{equation}
We can rearrange this equation to 
$$
xy+x^2=z^2-y^2+1.
$$
To make the equation linear in $x$, we can use the substitutions $z=z_1-x$ and $y=y_1-x$, which reduces the equation to 
$$
x(2z_1-y_1)=z_1^2-y_1^2+1.
$$
We can then introduce a new integer variable $t$ such that $t=2z_1-y_1$, which reduces the equation to 
\begin{equation}\label{x2pxypy2mz2m1_red}
xt=-3z_1^2+4z_1 t-t^2+1.
\end{equation}
This equation is of the form $P(x_3,\dots,x_n)=Q(x_2,\dots,x_n)$ with $x_1=x$, $x_2=z_1$, $P=t$, $Q=-3$, $R=4t$, $T=-t^2+1$. This equation has an integer solution $(x,z_1,t)=(0,0,1)$ and therefore by Proposition \ref{prop:finlincheck} it has infinitely many integer solutions. Then, we can solve this equation using the method in Section \ref{ex:H16modularroot}, to obtain that the integer solutions to equation \eqref{x2pxypy2mz2m1} are 
\begin{equation}\label{x2pxypy2mz2m1_sol}
\begin{aligned}
x=& \, -3t u^2+ur-\frac{r^2-4t^2-12}{12t}, \\
y=& \, 2tu-t+3tu^2-ur-\frac{r-4t}{3}+\frac{r^2-4t^2-12}{12t}, \\
  z=& \, tu+3tu^2-ur-\frac{r-4t}{6}+\frac{r^2-4t^2-12}{12t}, \\
  r \equiv & \, 4t \, \text{(mod 6)}, \quad r \in S(12t,4t^2+12).
  \end{aligned}
\end{equation}
We must now consider the case where $P=0$, or equivalently $t=0$. This reduces equation \eqref{x2pxypy2mz2m1_red} to 
$$
-3z_1^2+1=0,
$$
which has no integer solutions. 
Therefore, we can describe all integer solutions to equation \eqref{x2pxypy2mz2m1} by \eqref{x2pxypy2mz2m1_sol}. 

\vspace{10pt}

Let us consider the next equation
\begin{equation}\label{x2pxypy2mz2p1}
x^2+xy+y^2-z^2+1=0.
\end{equation}
We can rearrange this equation to 
$$
xy+x^2=z^2-y^2-1.
$$
To make the equation linear in $x$, we can use the substitutions $z=z_1-x$ and $y=y_1-x$, which reduces the equation to 
$$
x(2z_1-y_1)=z_1^2-y_1^2-1.
$$
We can then introduce a new integer variable $t$ such that $t=2z_1-y_1$, which reduces the equation to 
\begin{equation}\label{x2pxypy2mz2p1_red}
xt=-3z_1^2+4z_1 t-t^2-1.
\end{equation}
This equation is of the form $P(x_3,\dots,x_n)=Q(x_2,\dots,x_n)$ with $x_1=x$, $x_2=z_1$, $P=t$, $Q=-3$, $R=4t$, $T=-t^2-1$. This equation has an integer solution $(x,z_1,t)=(1,-1,-1)$ and therefore by Proposition \ref{prop:finlincheck} it has infinitely many integer solutions. Then, we can solve this equation using the method in Section \ref{ex:H16modularroot}, to obtain that the integer solutions to equation \eqref{x2pxypy2mz2p1} are
\begin{equation}\label{x2pxypy2mz2p1_sol}
\begin{aligned}
x=& \, -3t u^2+ur-\frac{r^2-4t^2+12}{12t}, \\
y=& \, 2tu-t+3tu^2-ur-\frac{r-4t}{3}+\frac{r^2-4t^2+12}{12t}, \\
  z=& \, tu+3tu^2-ur-\frac{r-4t}{6}+\frac{r^2-4t^2+12}{12t}, \\
  r \equiv & \, 4t \, \text{(mod 6)}, \quad r \in S(12t,4t^2-12).
  \end{aligned}
\end{equation}
We must now consider the case where $P=0$, or equivalently $t=0$. This reduces \eqref{x2pxypy2mz2p1_red} to 
$$
-3z_1^2-1=0,
$$
which has no integer solutions. 
Therefore, we can describe all integer solutions to equation \eqref{x2pxypy2mz2p1} by \eqref{x2pxypy2mz2p1_sol}. 

\vspace{10pt}

We have now considered the equations not of the form \eqref{eq:genquadform} of size $H \leq 17$. Many equations with $H =18$ are very similar to the equations we have already solved in this section, so we will now consider equations which are not similar to the ones we have previously solved. 

The first equation we will consider is 
\begin{equation}\label{xpx2my2mz2mt2}
x+x^2-y^2-z^2-t^2=0.
\end{equation}
We can rearrange this equation to 
$$
x+x^2=y^2+z^2+t^2.
$$
To make the equation linear in $x$, we can use the substitution $y=y_1+x$, which reduces the equation to 
$$
x(1-2y_1)=z^2+y_1^2+t^2.
$$
This equation is of the form $P(x_3,\dots,x_n)=Q(x_2,\dots,x_n)$ with $x_1=x$, $x_2=z$, $P=1-2y_1$, $Q=1$, $R=0$, $T=y_1^2+t^2$. This equation has an integer solution $(x,y_1,z,t)=(1,-2,-1,0)$ and therefore by Proposition \ref{prop:finlincheck} it has infinitely many integer solutions. Then, we can solve this equation using the method in Section \ref{ex:H16modularroot}, to obtain that the integer solutions to equation \eqref{xpx2my2mz2mt2} are
\begin{equation}\label{xpx2my2mz2mt2_sol}
\begin{aligned}
x=& (1-2y_1)u^2+ur+\frac{r^2+4y_1^2+4v^2}{4(1-2y_1)}, \\
y= & y_1+(1-2y_1)u^2+ur+\frac{r^2+4y_1^2+v^2}{4(1-2y_1)}, \\
z=& (1-2y_1)u-\frac{r}{2},\\
t=& v, \\  
& r \in S(4(1-2y_1),-4(y_1^2+v^2)), \quad v,u,y_1 \in \mathbb{Z}, \\     
\end{aligned}
 \end{equation}
subject to the additional condition that $ r \equiv \, 0 \, \text{(mod 2)}$. However, as $r^2+4y_1^2+4v^2$ must be even to be divisible by $4(1-2y_1)$, $r$ must be even, hence, this restriction can be omitted. In summary, all integer solutions to equation \eqref{xpx2my2mz2mt2} are given by \eqref{xpx2my2mz2mt2_sol}.

\vspace{10pt}

The next equation we will consider is 
\begin{equation}\label{xpx2py2pypzmz2}
x+x^2+y^2+y+z-z^2=0.
\end{equation}
We can rearrange this equation to 
$$
x+x^2y^2+y=z^2-z.
$$
To make the equation linear in $z$, we can use the substitution $x=x'-z$, which reduces the equation to 
$$
2x'z=(x')^2+x'+y^2+y.
$$
This equation is of the form $P(x_3,\dots,x_n)=Q(x_2,\dots,x_n)$ with $x_1=z$, $x_2=y$, $P=2x'$, $Q=1$, $R=1$, $T=(x')^2+x'$. This equation has an integer solution $(x',y,z)=(1,1,2)$ and therefore by Proposition \ref{prop:finlincheck} it has infinitely many integer solutions. Then, we can solve this equation using the method in Section \ref{ex:H16modularroot}, to obtain that the integer solutions to equation \eqref{xpx2py2pypzmz2} are
\begin{equation}\label{xpx2py2pypzmz2_sol}
\begin{aligned}
x=& x'-2x' u^2-ur-\frac{r^2-1+4(x')^2+4x'}{8x'}, \\
y= & \, 2x' u+\frac{r-1}{2}, \\
z=& \, 2x' u^2+u r+\frac{r^2-1+4x'+4(x')^2}{8x'},\\
& r \in S(8x',1-4(x')^2-4x'), \quad u,x' \in \mathbb{Z}, \\     
\end{aligned}
 \end{equation}
subject to the additional condition that $ r \equiv \, 1 \, \text{(mod 2)}$. However, as $r^2-1+4x'+4(x')^2$ must be even to be divisible by $8x'$, $r$ must be odd, hence, this restriction can be omitted. 
We must now consider the case where $P=0$, or equivalently $x'=0$ or $x=-z$. This reduces \eqref{xpx2py2pypzmz2} to 
$$
y^2+y=0,
$$
and we obtain the integer solutions 
\begin{equation}\label{xpx2py2pypzmz2_sol2}
(x,y,z)=(u,-1,-u), (u,0,-u), \quad u \in \mathbb{Z}.
 \end{equation}
In summary, all integer solutions to equation \eqref{xpx2py2pypzmz2} are given by either \eqref{xpx2py2pypzmz2_sol} or \eqref{xpx2py2pypzmz2_sol2}.

\vspace{10pt}

The next equation we will consider is 
\begin{equation}\label{x2pxymy2pzpz2}
x^2+xy-y^2+z+z^2=0.
\end{equation}
We can rearrange this equation to 
$$
x^2+xy=y^2-z-z^2.
$$
To make the equation linear in $x$, we can use the substitutions $y=y_1-x$ and $z=z_1-x$, which reduces the equation to 
$$
x(3y_1-2z_1-1)=y_1^2-z_1^2-z_1.
$$
By making the further substitutions $w=y_1-z_1$ and $k=3w+z_1$, the equation is reduced to 
\begin{equation}\label{x2pxymy2pzpz2_red}
x(k-1)=-5w^2+w(2k+3)-k.
\end{equation}
This equation is of the form $P(x_3,\dots,x_n)=Q(x_2,\dots,x_n)$ with $x_1=x$, $x_2=w$, $P=k-1$, $Q=-5$, $R=2k+3$, $T=-k$. This equation has an integer solution $(x,k,w)=(-2,2,0)$ and therefore by Proposition \ref{prop:finlincheck} it has infinitely many integer solutions. Then, we can solve this equation using the method in Section \ref{ex:H16modularroot}, and obtain that the integer solutions to equation \eqref{x2pxymy2pzpz2} are
\begin{equation}\label{x2pxymy2pzpz2_sol}
\begin{aligned}
x = & -5 (k - 1) u^2 +   u r - \frac{r^2 - (2 k + 3)^2 + 20 k}{20 (k - 1)}, \\
y= & \, k - 2 u (k - 1) + \frac{r - 2 k - 3}{5} + 5 (k - 1) u^2 - u r + \frac{r^2 - (2 k + 3)^2 + 20 k}{20 (k - 1)}, \\
 z=& \, k - 3 u (k - 1) + \frac{3 (r - 2 k - 3)}{10} + 5 (k - 1) u^2 -   u r + \frac{r^2 - (2 k + 3)^2 + 20 k}{20 (k - 1)}, \\
  r \equiv & \, 2k+3 \, \text{(mod 10)}, \quad r \in S(20(k-1),(2k+3)^2-20k), \quad u,k \in \mathbb{Z}.   
\end{aligned}
 \end{equation}
 
 We must now consider the case where $P=0$, or equivalently $k-1=0$ or $k=1$. This reduces \eqref{x2pxymy2pzpz2_red} to 
$$
5w^2-5w+1=0
$$
which has no integer solutions. Therefore all integer solutions to equation \eqref{x2pxymy2pzpz2} are given by \eqref{x2pxymy2pzpz2_sol}.

\vspace{10pt}

The final equation we will consider is 
\begin{equation}\label{2x2pymy2mz2}
2x^2+y-y^2-z^2=0.
\end{equation}
We can rearrange this equation to 
$$
y^2-y=2x^2-z^2.
$$
To make the equation linear in $y$, we can use the substitutions $x=x'-y$ and $z=z'-y$, which reduces the equation to 
$$
y(4x'-2z'-1)=2(x')^2-(z')^2.
$$
By making the further substitution $w=2x'-z'$, the equation is reduced to 
$$
y(2w-1)=-2(x')^2+4x'w-w^2.
$$
This equation is of the form $P(x_3,\dots,x_n)=Q(x_2,\dots,x_n)$ with $x_1=y$, $x_2=x'$, $P=2w-1$, $Q=-2$, $R=4t$, $T=-w^2$. This equation has an integer solution $(x',y,w)=(1,1,1)$ and therefore by Proposition \ref{prop:finlincheck} it has infinitely many integer solutions. Then, we can solve this equation using the method in Section \ref{ex:H16modularroot}, and obtain that the integer solutions to equation \eqref{2x2pymy2mz2} are
\begin{equation}\label{2x2pymy2mz2_sol}
\begin{aligned}
x = &  \,  (2w - 1)(2u^2+u) - \frac{r }{4}+w  - u r + \frac{r^2 -8 t^2}{8 (2w - 1)}, \\ 
 y= & -2u^2 (2w - 1) +  u r - \frac{r^2 - 8 w^2}{8 (2w - 1)}, \\
z = & \, 2(u^2+u) (2w - 1) - \frac{r }{2} +w  -    u r + \frac{r^2 -8 w^2}{8 (2w - 1)}, \\
& r \in S(8(2w-1),8w^2) \quad u,w \in \mathbb{Z},  
\end{aligned}
 \end{equation}
subject to the additional condition that $ r \equiv \, 4w \, \text{(mod 4)}$. This implies that $r \equiv \, 0 \, \text{(mod 4)}$. However, as $r^2-8w^2$ must be divisible by $8$, $r^2$ must be divisible by $8$, but because $r^2$ is a perfect square, to be divisible by $8$, it must also be divisible by $16$, hence $r$ must be divisible by $4$. Therefore, this restriction can be omitted. In summary, all integer solutions to equation \eqref{2x2pymy2mz2} are given by \eqref{2x2pymy2mz2_sol}.



We will now attempt to solve equations by describing the solutions using a finite union of polynomial families. To solve equations in this way, we may need to use the following proposition.
\begin{proposition}\label{prop:x2py2pzt}[Proposition 2.23 in the book]
Integers $x,y,z,t$ satisfy equation 
\begin{equation}\label{eq:x2py2pzt}
x^2+y^2=zt
\end{equation}
if and only if 
\begin{equation}\label{eq:x2py2pztsol}
(x,y,z,t)=(s(uw+vr),s(vw-ur),s(u^2+v^2),s(w^2+r^2)) \quad \text{for some integers} \quad u,v,w,r,s.
\end{equation}
\end{proposition} 

The first equation we will consider is
\begin{equation}\label{x2py2p2tzpz}
x^2+y^2+2tz+z=0.
\end{equation}
We can rearrange this equation to 
$$
x^2+y^2=z(-2t-1),
$$
which suggests to make the substitution $k=-2t-1$, reducing the equation to 
$$
x^2+y^2=kz.
$$
Up to the names of variables, this is equation \eqref{eq:x2py2pzt} whose integer solutions are $(x,y,z,k)=(x,y,z,t)$ where $x,y,z,t$ are given by \eqref{eq:x2py2pztsol}. Then, in the original variables, we have that the integer solutions to equation \eqref{x2py2p2tzpz} are
\begin{equation}\label{x2py2p2tzpz_sol}
(x,y,z,t)=\left(s(uw+vr), s(vw-ur),s(u^2+v^2),\frac{-s(w^2+r^2)-1}{2}\right) \quad \text{for some integers} \quad u,v,w,r,s,
\end{equation}
Because $t$ must be an integer, we must have that $s$ and $w^2+r^2$ are odd, so $s$ must be odd and $w$ and $r$ must have different parities. By swapping $(w,v)$ with $(r,u)$ if necessary, we may assume that $w$ is odd and $r$ is even. Then $s=2s'+1$, $w=2w'$ and $r=2r'+1$, for some integers $s',w',r'$. After substituting these into \eqref{x2py2p2tzpz_sol}, we have that the integer solutions to equation \eqref{x2py2p2tzpz} are
$$
\begin{aligned}
x=& (2 s' + 1) (2 u w' +  2v r' + v),\\ 
y=& (2 s' + 1) (2 v w' - 2 ur' - u), \\
z=& (2 s' +  1) (u^2 + v^2), \\
t= & -(4s' (w')^2+2(w')^2+4r's'+2r'+s'+4s'(r')^2+2(r')^2+1), \\ & \text{for some integers} \quad u,v,w',r',s'.
\end{aligned}
$$

\vspace{10pt}

The next equation we will consider is
\begin{equation}\label{x2pxypy2pzmz2}
x^2+xy+y^2+z-z^2=0.
\end{equation}
We can make this equation linear in $x$ by making the substitutions $y=y_1-x$ and $z=z_1-x$. This reduces the equation to
$$
x(2z_1-y_1-1)=z_1^2-z_1-y_1^2.
$$
We can then make the further substitution $t=2z_1-y_1-1$ which reduces the equation to
$$
xt=-3z_1^2+4z_1t+3z_1-t^2-2t-1.
$$
Multiplying the equation by $12$ and making the substitutions $X=-12x+4t$ and $Z=-6z_1+4t+3$, the equation is reduced to 
$$
Xt=Z^2+3,
$$
which up to the names of variables, is equation \eqref{eq:yzmx2m3} whose integer solutions are given by \eqref{eq:yzmx2m3solb} or \eqref{eq:yzmx2m3sola}. Then, in the original variables, we have that the integer solutions to equation \eqref{x2pxypy2pzmz2} are
\begin{equation}\label{x2pxypy2pzmz2_sola}
\begin{aligned}
x=& \frac{8u^2+8ru+8r^2-2w^2-2vw-2v^2}{12}, \\
y= & \frac{w^2+vw+v^2}{6} - \frac{2uw+uv+wr+2vr}{3}, \\
z= & \frac{2u^2+2r u+2r^2}{3}+\frac{w^2+v w+v^2-2u w-u v-w r-2v r+3}{6}, \quad \text{or,} \\ \\
\end{aligned}
\end{equation}
\begin{equation}\label{x2pxypy2pzmz2_solb}
	\begin{aligned}
x=& \frac{-8u^2+8ru-8r^2+2w^2-2vw+2v^2}{12}, \\
y= &\frac{2 r v -  u v-  r w + 2 u w}{3}- \frac{  v^2  - v w + w^2}{6}, \\
z= &\frac{- 2 r^2 - 2 u^2 + 2 u r + v r + u w }{3} + \frac{3 - u v - v^2 - wr   + v w -    w^2}{6} , \quad \text{or,} \\ \\
\end{aligned}
\end{equation}
\begin{equation}\label{x2pxypy2pzmz2_solc}
	\begin{aligned}
x=& \frac{4u^2+12r^2-3v^2-w^2}{12}, \\
y= &-v r+ \frac{ 3 v^2  - 4 u w + w^2}{12}, \\
z= & r^2+ \frac{6+4u^2-6vr+ 3 v^2 -2uw+ w^2}{12} , \quad \text{or,}  \\ \\
\end{aligned}
\end{equation}
\begin{equation}\label{x2pxypy2pzmz2_sold}
	\begin{aligned}
x=& \frac{-4u^2-12r^2+3v^2+w^2}{12}, \\
y= &\frac{ -u^2-3r^2+ u w+3 vr}{3} + \frac{4u^2+12r^2-3v^2-w^2}{12}, \\
z= & - r^2+ \frac{6-4u^2+6vr- 3 v^2 +2uw- w^2}{12}, \\
\end{aligned}
\end{equation}
with the restriction that $uv-wr=1$ and the congruence conditions modulo $12$ to ensure $x,y,z$ are integers.

\vspace{10pt}

The next equation we will consider is \eqref{xpx2pxypy2mz2}.
We can make this equation linear in $x$ by making the substitutions $y=y_1-x$ and $z=z_1-x$, which reduces the equation to
$$
x(1+2z_1-y_1)=z_1^2-y_1^2.
$$
We can then make the further substitution $t=1+2z_1-y_1$, which reduces the equation to
$$
xt=-3z_1^2+4z_1t-4z_1-t^2+2t-1.
$$
Multiplying the equation by $12$ and making the substitutions $X=-12x+4t-8$ and $Z=-6z_1+4t-4$, the equation is reduced to 
$$
Xt=Z^2-4,
$$
which up to the names of variables, is equation \eqref{eq:yzmx2p4}, whose integer solutions are given by \eqref{eq:yzmx2p4sola}, \eqref{eq:yzmx2p4solb} and \eqref{eq:yzmx2p4solc}. Then, in the original variables, we have that the integer solutions to equation \eqref{xpx2pxypy2mz2} are
\begin{equation}\label{xpx2pxypy2mz2_sol}
\begin{aligned}
(x,y,z)=&\left(\frac{-4u^2+16r^2+w^2-4v^2-8}{12},\frac{4 - 16 r v + 4 v^2 + 4 u w - w^2}{12},  \frac{16 r^2 - 4 u^2 - 8 r v + 4 v^2 + 2 u w - w^2}{12}\right), \\
& \left(\frac{4 u^2 - 16 r^2 - w^2 + 4 v^2 - 8}{12},\frac{4 + 16 r v - 4 v^2 - 4 u w + w^2}{12},\frac{-16 r^2 + 4 u^2 + 8 r v - 4 v^2 - 2 u w + w^2}{12}\right), \\
& \left( \frac{8 u^2 - 8 r^2 - 2 w^2 + 2 v^2 - 8}{12}, \frac{2 + 4 r v - v^2 - 4 u w + w^2}{6}, \frac{-4 r^2 + 4 u^2 + 2 r v - v^2 - 2 u w + w^2}{6} \right), \\
& \left(\frac{-8 u^2 + 8 r^2 + 2 w^2 - 2 v^2 - 8}{12}, \frac{2 - 4 r v + v^2 + 4 u w - w^2}{6},  \frac{4 r^2 - 4 u^2 - 2 r v + v^2 + 2 u w - w^2}{6} \right), \\
& \left(\frac{16 r u - 8 - 4 v w}{12}, \frac{1 - 2 u v - 2 r w + v w}{3},  \frac{4 r u - u v - r w + v w}{3}\right),
\end{aligned}
\end{equation}
with the restriction that $uv-wr=1$ and the congruence conditions modulo $12$ to ensure $x,y,z$ are integers.

For equations of the form \eqref{eq:x2pyzpc}, the solution set is a finite set of polynomial families, the number of polynomial families is dependent on the number of factors of $c$. In order to write solutions more concisely, we will now accept solutions written in parameters $X,Y,Z$ where $X,Y,Z$ are solutions to $XY-Z^2=c$ for some integer $c$. 

The first equation we will consider is
\begin{equation}\label{xpx2pymy2pz2m2}
x+x^2+y-y^2+z^2-2=0.
\end{equation}
We can make this equation linear in $x$ by making the substitution $y=y'-x$. This reduces the equation to
$$
2xy'=(y')^2-y'-z^2+2.
$$
After multiplying this equation by $4$ and making the substitutions $x'=-8x+4y'-4$ and $z'=2z$ we reduce the equation to 
$$
x'y'=(z')^2-8,
$$
which can be solved using the method in Section \ref{ex:x2pyzpc}. Therefore, we can describe all integer solutions to equation \eqref{xpx2pymy2pz2m2} by
$$
(x,y,z)=\left( \frac{4y'-x'-4}{8}, y' -\frac{4y' -x'-4}{8},\frac{z'}{2} \right), \quad \text{where} \,\, x',y' ,z' \,\, \text{are solutions to} \,\, x'y'=(z')^2-8,
$$
subject to the condition that $z'$ is even and $4y'-x'-4 \equiv 0 \, (\text{mod 8})$.

\vspace{10pt}

Let us consider the next equation
\begin{equation}\label{x2pxymy2pz2m2}
x^2+xy-y^2+z^2-2=0.
\end{equation}
We can make this equation linear in $x$ by making the substitutions $y=y_1-x$ and $z=z_1-x$. This reduces the equation to
$$
x(3y_1-2z_1)=y_1 ^2-z_1 ^2+2.
$$
After the further substitutions $t=y_1-z_1$ and $k=3t+z_1$, we reduce the equation to
$$
xk=-5t^2+2tk+2.
$$
Multiplying the equation by $20$ we can make the change of variables $x'=-20x+4k$ and $z'=2k-10t$, which reduces the equation to 
$$
x'k=(z')^2-40,
$$
which can be solved using the method in Section \ref{ex:x2pyzpc}. Therefore, we can describe all integer solutions to equation \eqref{x2pxymy2pz2m2} by
$$
(x,y,z)=\left( \frac{4k-x'}{20}, \frac{8k+4z'+x'}{20},\frac{4k+6z'+x'}{20} \right), \quad \text{where} \,\, x',k ,z' \,\, \text{are solutions to} \,\, x'k=(z')^2-40,
$$
subject to the condition that $4k-x' \equiv 8k+4z'+x' \equiv 4k+6z'+x' \equiv 0 \, (\text{mod 20})$.

\vspace{10pt}

The final equation we will consider is
\begin{equation}\label{2x2pymy2pz2}
2x^2+y-y^2+z^2=0.
\end{equation}
We can make this equation linear in $y$ by making the substitution $z=z_1-y$. This reduces the equation to
$$
y(2z_1-1)=2x_1 ^2+z_1 ^2.
$$
After multiplying by 4 and making the further substitution $k=2z_1-1$, we reduce the equation to
$$
4yk=8x^2+k^2+2k+1.
$$
Then, multiplying the equation by $32$ and the making the further change of variables $X=128y'-32k-64$ and $Z=16x$, we reduce the equation to 
$$
Xk=Z^2+32,
$$
which can be solved using the method in Section \ref{ex:x2pyzpc}. Therefore, we can describe all integer solutions to equation \eqref{2x2pymy2pz2} by
$$
(x,y,z)=\left( \frac{Z}{16}, \frac{X+32k+64}{128},\frac{k+1}{2} - \frac{X+32k+64}{128} \right), \quad \text{where} \,\, X,k,Z \,\, \text{are solutions to} \,\, Xk=Z^2+32,
$$
subject to the conditions that  $Z \equiv 0 \, (\text{mod 16})$, $X+32k+64 \equiv 0 \, (\text{mod 128})$ and $k$ is odd.

 \begin{center}
\begin{tabular}{|c|c|c|c|c|c|}
\hline
Equation &  Solution $(x,y,z)$  \\ \hline \hline   
   $x^2+y^2+2tz+z=0$ & $x = (2 s + 1) (  2 uw + u + 2v r),$ \\
  &$ y = (2 s + 1) (    2 v w + v - 2u r),$ \\
  &$ z = (2 s + 1) (u^2 +     v^2),$ \\
   & $ t = -(4sw^2+2w^2+4ws+2w+s+4sr^2+2r^2+1), $  \\
  & or,  $x=(2 s + 1) (2 u w +  2v r + v),$ \\
  & $ y = (2 s + 1) (2 v w - 2 ur - u),$ \\
  &$ z = (2 s +  1) (u^2 + v^2),$ \\
  & $ t = -(4sw^2+2w^2+4rs+2r+s+4sr^2+2r^2+1)$   \\ \hline
 $x+x^2+y-y^2+z^2-2=0$ & $\left( \frac{4y'-x'-4}{8}, y' -\frac{4y' -x'-4}{8}, \frac{z'}{2} \right) $, \\
 & where $x',y' ,z'$ are solutions to \\
& $x'y'=(z')^2-8=0,$ subject to the condition that  \\
& $z'$ is even and $4y'-x'-4 \equiv 0 \, (\text{mod 8})$.\\  \hline
  $x^2+xy-y^2+z^2-2=0$  &$\left( \frac{4k-x'}{20}, \frac{8k+4z'+x'}{20},\frac{4k+6z'+x'}{20} \right)$ \\
  & where $x',y' ,z'$ are solutions to \\
  & $ x'k=(z')^2-40$, subject to the condition that  \\
  & $4k-x' \equiv 4z'+3x' \equiv 6z'+2x' \equiv 0 \, (\text{mod 20})$. \\  \hline
 $x^2+xy+y^2+z-z^2=0$& \eqref{x2pxypy2pzmz2_sola}, \eqref{x2pxypy2pzmz2_solb}, \eqref{x2pxypy2pzmz2_solc}, or \eqref{x2pxypy2pzmz2_sold}  \\  \hline
$x+x^2+xy+y^2-z^2=0$ & \eqref{xpx2pxypy2mz2_sol}  \\  \hline
  $2x^2+y-y^2+z^2=0$ & $\left( \frac{Z}{16}, \frac{X+32k+64}{128},\frac{k+1}{2} - \frac{X+32k+64}{128} \right)$, \\
 & where $X,k,Z$ are solutions to \\
 & $Xk=Z^2+32$, subject to the conditions that \\
 &  $Z \equiv 0 \, (\text{mod 16})$, $X+32k+64 \equiv 0 \, (\text{mod 128})$ \\ &  and $k$ is odd. \\  \hline
 \end{tabular}
 \captionof{table}{\label{tab:H18quadred_para} Integer solutions presented as a finite union of polynomial families to some of the equations listed in Table \ref{tab:H18quadred}.} 
 \end{center}

Table \ref{tab:H18quadred_excluded} contains quadratic equations of size $H\leq 18$ not of the form \eqref{eq:genquadform} which are similar to equations in Table \ref{tab:H18quadredsol} and Table \ref{tab:H18quadred_para}. As these equations are similar, they have been excluded, and therefore, have not been solved in this section. 
The integer solutions to all other equations in Table \ref{tab:H18quadred} have been described as a finite union of polynomial families (or showing they can be written as a finite union of polynomial families) and using $S$ notation. We will let the reader decide which answer form they think is the ``nicest.''

 \begin{center}
\begin{tabular}{|c|c|c|c|c|c|}
\hline
$H$ & Equation & $H$ &  Equation  \\ \hline \hline
      $18$ & $x+x^2-y^2-z^2-4=0$   &  $18$ & $x+x^2-y^2-z^2+4=0$      \\  \hline 
 $18$ & $x+x^2+y-y^2+z^2+2=0$ & $18$ & $x+x^2+y+y^2-z^2-2=0$ \\  \hline 
$18$ & $x^2+xy-y^2+z^2+2=0$   & $18$ & $x+x^2+y+y^2-z^2+2=0$   \\  \hline 
  $18$ & $x^2+xy+y^2-z^2-2=0$   &    $18$ & $x^2+xy+y^2-z^2+2=0$     \\  \hline  
 $18$ & $x+x^2+xy-y^2+z^2=0$ &  $18$ & $2x^2+y+y^2-z^2=0$ \\  \hline
 \end{tabular}
 \captionof{table}{\label{tab:H18quadred_excluded} Equations listed in Table \ref{tab:H18quadred} which are similar to equations in Table \ref{tab:H18quadredsol} and Table \ref{tab:H18quadred_para}.} 
\end{center}

\subsection{Exercise 4.29}\label{ex:H24nonhom}
\textbf{\emph{Solve equations of the form 
		\begin{equation}\label{eq:anisotropic}
			Q(x_1,\dots,x_n)=b, 
		\end{equation}
	 where $Q$ is given by 
	 \begin{equation}\label{eq:anisotropiclist}
	 	Q(x,y,z)=3x^2-y^2-z^2, \quad 3x^2-2y^2+z^2 \quad \text{and} \quad 3x^2+2y^2-z^2,
	 \end{equation}
  and $b$ is an integer satisfying $1 \leq |b|\leq 3$. }}

To solve these equations, we will use the method from Section 4.2.5 of the book, which we summarise below for convenience.

For an arbitrary quadratic form $Q$ of the form $Q(x_1,\dots,x_n)=\sum_{i=1}^n \sum_{j=i}^n a_{i j} x_i x_j$ we have the transformation
\begin{equation}\label{eq:hwdef2}
h_w(x)=x-\frac{Q(x,w)}{Q(w)} w, \quad \text{where} \quad Q(x,w)=Q(x+w)-Q(x)-Q(w)=\sum_{i=1}^n \sum_{j=i}^n a_{ij}(x_i w_j+x_j w_i)
\end{equation}
To generate transformations we should select $w$ such that $h_w(x)$ has integer coefficients. Transformations should also satisfy $h_w(h_w(x))=x$ and $Q(x)=Q(h_w(x))$.

We will find a finite set $S$ of linear transformations and a finite set $X$ of integer solutions to \eqref{eq:anisotropic} such that every integer solution to \eqref{eq:anisotropic} can be obtained from some solution $x \in X$ by applying transformations from $S$ in some order. 

The solutions to all equations in this exercise are preserved after transformations
\begin{equation}\label{3x2my2mz2mb_trivial}
	(x,y,z) \to (-x,y,z), \quad (x,y,z) \to (x,-y,z), \quad (x,y,z)\to (x,y,-z).
\end{equation}
If we apply any of these transformations twice, we return to the solution we started with. 
We will call transformations \eqref{3x2my2mz2mb_trivial} ``trivial.'' 

The first equations we will consider are
\begin{equation}\label{3x2my2mz2mb}
	3x^2-y^2-z^2=b,
\end{equation}
where $b$ is equal to $\pm 1$, $\pm 2$ or $\pm 3$. 
For these equations, we have an additional trivial transformation 
\begin{equation}\label{3x2my2mz2mb_trivial2}
	(x,y,z)\to (x,z,y).
\end{equation}
Let us try to find non-trivial transformations of the solutions using formula \eqref{eq:hwdef2} with $Q(x,y,z)=3x^2-y^2-z^2$. 
To ensure that coefficients in \eqref{eq:hwdef2} are integers, let us select $w=(w_1,w_2,w_3)$ such that $Q(w)=Q(w_1,w_2,w_3)=2$. A simple vector satisfying this condition is $(w_1,w_2,w_3)=(1,1,0)$. In this case, \eqref{eq:hwdef2} implies
$$
h_w(x) = (x,y,z)-\frac{6x-2y}{2}(1,1,0)=(-2x+y,-3x+2y,z),
$$
resulting in the transformation 
\begin{equation}\label{3x2my2mz2mb_transform}
	(x,y,z) \to (-2x+y,-3x+2y,z).
\end{equation}
Substituting this into \eqref{3x2my2mz2mb}, we verify that the transformation \eqref{3x2my2mz2mb_transform} indeed maps integer solutions of \eqref{3x2my2mz2mb} into integer solutions of the same equation. Also, the identities
$$
-2(-2x+y)+(-3x+2y)=x \quad \text{and} \quad -3(-2x+y)+2(-3x+2y)=y
$$
verify that if we apply \eqref{3x2my2mz2mb_transform} twice, we obtain the solution we started with.

Let us say that two integer solutions to \eqref{3x2my2mz2mb} are in the same orbit if they can be obtained from each other by a sequence of transformations \eqref{3x2my2mz2mb_trivial}, \eqref{3x2my2mz2mb_trivial2} and \eqref{3x2my2mz2mb_transform}. 

In any orbit, choose a solution $(x,y,z)$ with $|x|$ minimal. Swapping $y$ and $z$ and the signs of $x$ and $y$ if necessary, we may assume that $x \geq 0$ and $y \geq |z| \geq 0$. Because we have chosen the solution with $|x|$ minimal in the orbit, we must have $|-2x+y| \geq x$. If $-2x+y>0$, then $-2x+y \geq x$, or $3x \leq y$. We then have
$$
9x^2 \leq y^2 \leq y^2+z^2 = 3x^2-b.
$$
This inequality is possible only if $b<0$ and $|x| \leq \sqrt{\frac{|b|}{6}}$. Because $|b|\leq 3$ and $x$ is an integer, this implies that $x=0$ and $|y| \leq \sqrt{|b|}$. This together with $y \geq |z| \geq 0$ implies that $(x_0,y_0,z_0)=(0,1,0)$ for $b=-1$, $(x_0,y_0,z_0)=(0,1,1)$ for $b=-2$ and that there are no solutions for $b=-3$. 

Now assume that $-2x+y \leq 0$. Then $|-2x+y| \geq x$ reduces to $2x-y \geq x$ or $x \geq y$, and we have
$$
2x^2 \geq 2y^2 \geq y^2  +z^2 = 3x^2-b.
$$
This implies that $b>0$ and $|x| \leq \sqrt{b}$. Because $x$ is a non-negative integer, and $1\leq b \leq 3$, we have either $x=0$ or $x=1$. Substituting these values into the original equation, we find that its only solutions satisfying $0\leq x \leq 1$ and $y \geq |z| \geq 0$ are $(x_0,y_0,z_0)=(1,1,1)$ for $b=1$, $(x_0,y_0,z_0)=(1,1,0)$ for $b=2$, and $(x_0,y_0,z_0)=(1,0,0)$ for $b=3$. 

In conclusion, equation \eqref{3x2my2mz2mb} has no integer solutions for $b=-3$, while for $b=-2,-1,1,2,3$ all solutions are in the same orbit as $(x_0,y_0,z_0)$, and therefore can be obtained from it by a sequence of transformations  \eqref{3x2my2mz2mb_trivial}, \eqref{3x2my2mz2mb_trivial2} and \eqref{3x2my2mz2mb_transform}.

		\vspace{10pt}
		
		The next equations we will consider are
		\begin{equation}\label{3x2p2y2mz2pb}
			3x^2+2y^2-z^2=b,
		\end{equation}
		where $b$ is equal to $\pm 1$, $\pm 2$ or $\pm 3$. Let us try to find non-trivial transformations of the solutions using formula \eqref{eq:hwdef2} with $Q(x,y,z)=3x^2+2y^2-z^2$. 
		To ensure that the coefficients in \eqref{eq:hwdef2} are integers, let us select $w=(w_1,w_2,w_3)$ such that $Q(w)=Q(w_1,w_2,w_3)$ is equal to $1$ or $2$. A simple vector satisfying the condition $Q(w)=1$ is $(w_1,w_2,w_3)=(0,1,1)$. In this case, \eqref{eq:hwdef2} implies that 
		$$
		h_w(x) = (x,y,z)-\frac{4y-2z}{1}(0,1,1)=(x,-3y+2z,-4y+3z),
		$$ 
		resulting in the transformation
		\begin{equation}\label{3x2p2y2mz2pb_transformi}
			(x,y,z) \to (x,-3y+2z,-4y+3z).
		\end{equation}
		Next, a simple solution to $Q(w)=2$ is $(w_1,w_2,w_3)=(1,0,1)$. In this case, \eqref{eq:hwdef2} implies that 
		$$
		h_w(x) = (x,y,z)-\frac{6x-2z}{2}(1,0,1)=(-2x+z,y,-3x+2z),
		$$ 
		resulting in the transformation
		\begin{equation}\label{3x2p2y2mz2pb_transformii}
			(x,y,z) \to (-2x+z,y,-3x+2z).
		\end{equation}
		Substituting the transformations \eqref{3x2p2y2mz2pb_transformi} and \eqref{3x2p2y2mz2pb_transformii} into \eqref{3x2p2y2mz2pb}, 
		we verify that these transformations indeed map integer solutions of \eqref{3x2p2y2mz2pb} into integer solutions of the same equation.
		
		Also, the identities
		$$
		-2(-2x+z)+(-3x+2z)=x \quad \text{and} \quad -3(-2x+z)+2(-3x+2z)=z, 
		$$
		and 
		$$
		-3(-3y+2z)+2(-4y+3z)=y \quad \text{and} \quad -4(-3y+2z)+3(-4y+3z)=z
		$$
		verify that if we apply either of the transformations \eqref{3x2p2y2mz2pb_transformi} or \eqref{3x2p2y2mz2pb_transformii} twice, we obtain the solution we have started with.
		
		Let us say that two integer solutions to \eqref{3x2p2y2mz2pb} are in the same orbit if they can be obtained from each other by a sequence of \eqref{3x2my2mz2mb_trivial}, \eqref{3x2p2y2mz2pb_transformi} and \eqref{3x2p2y2mz2pb_transformii}. 
		In any orbit, choose a solution $(x,y,z)$ with $|z|$ minimal. Swapping signs of $x$, $y$ and $z$ if necessary, we may assume that $x \geq 0$, $y \geq 0$ and $z \geq 0$. Because we have chosen the solution with $|z|$ minimal in the orbit, we must have $|-3x+2z| \geq z$ and $|-4y+3z| \geq z$. This gives us four cases to consider, (i) $x \geq z$ and $y \geq z$,
		(ii) $x \geq z$ and $z \geq 2y$,
		(iii) $z \geq 3x$ and $y \geq z$,
		and (iv) $z \geq 3x$ and $z \geq 2y$.
		We can consider cases (i) and (ii) simultaneously by considering the case where $x \geq z$. Then 
		$$
		3z^2 \leq 3x^2 \leq 3x^2+2y^2 =z^2+b,
		$$
		which implies that $b>0$ and $|z| \leq \sqrt{\frac{b}{2}}$. If $b=1$ then $z=0$, and if $b=2$ or $b=3$, $z$ is either $0$ or $1$. By using this together with $x \geq z \geq 0$, we obtain that there are no solutions for $b=1$, $(x_0,y_0,z_0)=(0,1,0)$ and $(x'_0,y'_0,z'_0)=(1,0,1)$ for $b=2$, and $(x_0,y_0,z_0)=(1,0,0)$ for $b=3$.
		
		In case (iii) we have $z^2 \geq 9x^2$ and $y^2 \geq z^2$, which implies that
		$$
		2z^2 \leq 2y^2 \leq 2y^2+3x^2 = z^2+b.
		$$
		This is possible only if $b>0$ and $|z| \leq \sqrt{b}$. Now, $1 \leq b \leq 3$ implies that $z=0$ or $z=1$. This together with $z \geq 3x \geq 0$ and $y \geq z \geq 0$ this implies that $(x_0,y_0,z_0)=(0,1,1)$ for $b=1$ and no further solutions for $b=2$ and $b=3$.
		
		In case (iv) we have $z^2 \geq 9x^2$ and $z^2 \geq 4y^2$, which implies that
		$$
		5z^2 = 2z^2 + 3z^2 \geq 18x^2+12y^2=6(z^2+b).
		$$
		This implies that $b<0$ and $|z| \leq \sqrt{-6b}$. However, $-3 \leq b \leq -1$, together with $z \geq 3x \geq 0$ and $z \geq 2y \geq 0$ implies that $(x_0,y_0,z_0)=(0,0,1)$ for $b=-1$, $(x_0,y_0,z_0)=(0,1,2)$ for $b=-2$, and no solutions for $b=-3$.
		
		Hence, all integer solutions to \eqref{3x2p2y2mz2pb} with $1\leq |b| \leq 3$, $b\neq 2$ are in the same orbit as $(x_0,y_0,z_0)$ and therefore can be obtained from it by a sequence of transformations \eqref{3x2my2mz2mb_trivial}, \eqref{3x2p2y2mz2pb_transformi} and \eqref{3x2p2y2mz2pb_transformii}. 
		
		All integer solutions to \eqref{3x2p2y2mz2pb} with $b = 2$ are in the same orbits as $(x_0,y_0,z_0)=(1,0,1)$ or $(x'_0,y'_0,z'_0)=(0,1,0)$ and therefore can be obtained from these solutions by a sequence of transformations \eqref{3x2my2mz2mb_trivial}, \eqref{3x2p2y2mz2pb_transformi} and \eqref{3x2p2y2mz2pb_transformii}. The first orbit corresponds to integer solutions with $y$ even, and the second orbit corresponds to integer solutions with $y$ odd.

			\vspace{10pt}

			The next equations we will consider are
			\begin{equation}\label{3x2m2y2pz2mb}
				3x^2-2y^2+z^2=b
			\end{equation}
			where $b$ is equal to $\pm 1,\pm 2$ or $\pm 3$. Let us find non-trivial transformations of the solutions using formula \eqref{eq:hwdef2} with $Q(x,y,z)=3x^2-2y^2+z^2$. To ensure the coefficients in \eqref{eq:hwdef2} are integers, let us select $w=(w_1,w_2,w_3)$ such that  $Q(w)=Q(w_1,w_2,w_3)=1$. Simple vectors satisfying this condition are $(w_1,w_2,w_3)=(0,1,1)$, $(1,1,0)$ and $(1,2,2)$. In these cases, \eqref{eq:hwdef2} implies that
			$$
			h_w(x) = (x,y,z)-\frac{-4y+2z}{-1}(0,1,1)= (x,-3y+2z,-4y+3z),
			$$ 
			$$
			h_w(x) = (x,y,z)-\frac{6x-4y}{1}(1,1,0)=(-5x+4y,-6x+5y,z),
			$$ 
			and 
			$$
			h_w(x) = (x,y,z)-\frac{6 x - 8 y + 4 z}{-1}(1,2,2)=(7 x - 8 y + 4 z, 12 x - 15 y + 8 z, 12 x - 16 y + 9 z),
			$$ 
			respectively, resulting in the transformations
			\begin{equation}\label{3x2m2y2pz2mb_tri}
				(x,y,z) \to (x,-3y+2z,-4y+3z),
			\end{equation}
			\begin{equation}\label{3x2m2y2pz2mb_trii}
				(x,y,z) \to (-5x+4y,-6x+5y,z),
			\end{equation}
			and 
			\begin{equation}\label{3x2m2y2pz2mb_triii}
				(x,y,z) \to (7x-8y+4z,12x-15y+8z,12x-16y+9z).
			\end{equation}
			Substituting the transformations \eqref{3x2m2y2pz2mb_tri}, \eqref{3x2m2y2pz2mb_trii} and \eqref{3x2m2y2pz2mb_triii} into \eqref{3x2m2y2pz2mb}, 
			we verify that these transformations indeed map integer solutions of \eqref{3x2m2y2pz2mb} into integer solutions of the same equation. Also, the identities
			$$
			-3(-3y+2z)+2(-4y+3z)=y \quad \text{and} \quad -4(-3y+2z)+3(-4y+3z)=z,
			$$
			$$
			-5(-5x+4y)+4(-6x+5y)=x \quad \text{and} \quad -6(-5x+4y)+5(-6x+5y)=y,
			$$
			and
			$$
			\begin{aligned}
				7(7x-8y+4z)-8(12x-15y+8z)+4(12x-16y+9z)=x, \\ 12(7x-8y+4z)-15(12x-15y+8z)+8(12x-16y+9z)=y, \\
				12(7x-8y+4z)-16(12x-15y+8z)+9(12x-16y+9z)=z,
			\end{aligned}
			$$
			verify that if we apply any of the listed transformations twice, we obtain the solution we have started with.
			
			Let us say that two integer solutions to \eqref{3x2m2y2pz2mb} are in the same orbit if they can be obtained from each other by a sequence of transformations \eqref{3x2my2mz2mb_trivial}, \eqref{3x2m2y2pz2mb_tri}, \eqref{3x2m2y2pz2mb_trii} and \eqref{3x2m2y2pz2mb_triii}.
			
			In any orbit, choose a solution $(x,y,z)$ with $|y|$ minimal. By changing the signs of $x$, $y$ and $z$ if necessary, we may assume that $x \geq 0$, $y \geq 0$ and $z \geq 0$. Because we have chosen the solution with $|y|$ minimal in the orbit, we must have $|-3y+2z| \geq y$, $|-6x+5y| \geq y$ and $|12x-15y+8z| \geq y$. 
			Using Mathematica, we derive that the only integer solutions to \eqref{3x2m2y2pz2mb} satisfying these conditions are  $(x_0,y_0,z_0)=(0,0,1)$ and $(x'_0,y'_0,z'_0)=(1,1,0)$ for $b=1$, $(x_0,y_0,z_0)=(0,1,1)$ and $(x'_0,y'_0,z'_0)=(1,2,2)$ for $b=-1$, $(x_0,y_0,z_0)=(0,1,2)$ and $(x'_0,y'_0,z'_0)=(1,1,1)$ for $b=2$, $(x_0,y_0,z_0)=(0,1,0)$ and $(x'_0,y'_0,z'_0)=(2,3,2)$ for $b=-2$, $(x_0,y_0,z_0)=(2,3,3)$ and $(x'_0,y'_0,z'_0)=(1,0,0)$ for $b=3$ and no solutions for $b=-3$. 
			
			Hence, equation \eqref{3x2m2y2pz2mb} has no integer solutions for $b=-3$, while for $b=-2,-1,1,2,3$ all integer solutions are in the same orbits as $(x_0,y_0,z_0)$ or $(x'_0,y'_0,z'_0)$ and therefore can be obtained from these solutions by a sequence of transformations \eqref{3x2my2mz2mb_trivial}, \eqref{3x2m2y2pz2mb_tri}, \eqref{3x2m2y2pz2mb_trii}, and \eqref{3x2m2y2pz2mb_triii}. For $b=-1,1,2,3$, the orbit of $(x_0,y_0,z_0)$ corresponds to integer solutions to \eqref{3x2m2y2pz2mb} with $x$ even, and the orbit of $(x'_0,y'_0,z'_0)$ corresponds to integer solutions to \eqref{3x2m2y2pz2mb} with $x$ odd. For $b = -2$, the orbit of $(x_0,y_0,z_0)$ corresponds to integer solutions with $x \equiv 0 \, (\text{mod} \, 4)$ and the orbit of $(x'_0,y'_0,z'_0)$ corresponds to integer solutions with $x \equiv 2 \, ( \text{mod } \, 4)$.
			
			Table \ref{tab:ex:4.28} summarises the solutions to the equations solved in this exercise.
			
			\begin{table}[htbp]
				\begin{center}

					\caption{\label{tab:ex:4.28} Integer solutions to equations $Q(x,y,z)=b$ with $b$ satisfying $1\leq |b| \leq 3$ and $Q(x,y,z)$ given by \eqref{eq:anisotropiclist}. In all cases, all integer solutions can be obtained from the listed initial solutions by applying the listed transformations in any order.} 
				\end{center}
			\end{table}

\subsection{Exercise 4.32}\label{ex:H18mon3}
\textbf{\emph{For each of the equations listed in Table \ref{tab:H18mon3}, determine whether for any of the three ways of writing it in the form \eqref{eq:gen3mon2}, both systems \eqref{eq:gen3monsyst1} and \eqref{eq:gen3monsyst2} are solvable in non-negative integers. For equations satisfying this condition, describe all their integer solutions using formula \eqref{eq:gen3mon2sol}.}}

	\begin{center}
		\begin{tabular}{ |c|c|c|c|c|c| } 
			\hline
			$H$ & Equation & $H$ & Equation & $H$ & Equation \\ 
			\hline\hline
			$18$ & $x+x^2 y-y z^2 = 0$ & $18$ & $3 x+x^2 y+z^2 = 0$ & $18$ & $x+x^2 y+y^2 z = 0$ \\ 
			\hline
			$18$ & $x^2 y+z+2 y z = 0$ & $18$ & $-y+x^2 y+2 z^2 = 0$ & $18$ & $x^3+y+2 y z = 0$ \\ 
			\hline
			$18$ & $x+x^2 y+2 y z = 0$ & $18$ & $y+x^2 y+2 z^2 = 0$ & $18$ & $x^3+y+x y z = 0$ \\ 
			\hline
			$18$ & $-y+x^2 y+2 x z = 0$ & $18$ & $x+x^2 y+2 z^2 = 0$ & $18$ & $x^3-z+y^2 z = 0$ \\ 
			\hline
			$18$ & $x^2 y+z+2 x z = 0$ & $18$ & $x^2 y+z+t y z = 0$ & $18$ & $x^3+z+y^2 z = 0$ \\ 
			\hline
			$18$ & $y+x^2 y+2 x z = 0$ & $18$ & $x^2 y+z+t x z = 0$ & $18$ &  $x^3+y+y^2 z = 0$ \\ 
			\hline
			$18$ & $-3 y+x^2 y+z^2 = 0$ & $18$ & $-y+x^2 y+x z^2 = 0$ & $18$ & $2 t x+y+x y z = 0$ \\ 
			\hline
			$18$ & $3 y+x^2 y+z^2 = 0$ & $18$ & $y+x^2 y+x z^2 = 0$ & & \\ 
			\hline
		\end{tabular}
		\captionof{table}{\label{tab:H18mon3} Three-monomial equations of size $H\leq 18$.}
	\end{center} 

In this exercise we will attempt to describe all integer solutions to three-monomial equations listed in Table \ref{tab:H18mon3}, using the method described in Section 4.3.1 of the book, which we summarise below for convenience. Three-monomial equations are of the form
\begin{equation}\label{eq:gen3mon2}
a \prod_{i=1}^n x_i^{\alpha_i} + b \prod_{i=1}^n x_i^{\beta_i} = c \prod_{i=1}^n x_i^{\gamma_i},
\end{equation} 
where $x_1,\dots,x_n$ are variables, $\alpha_i$, $\beta_i$ and $\gamma_i$ are non-negative integers, and $a,b,c$ are non-zero integer coefficients. If the systems
\begin{equation}\label{eq:gen3monsyst1}
	\sum_{i=1}^n \alpha_i z_i = \sum_{i=1}^n \beta_i z_i = \sum_{i=1}^n \gamma_i z_i - 1, \quad \quad z_i \geq 0, \quad i=1,\dots, n, 
\end{equation}
and
\begin{equation}\label{eq:gen3monsyst2}
	\sum_{i=1}^n \alpha_i t_i = \sum_{i=1}^n \beta_i t_i = \sum_{i=1}^n \gamma_i t_i + 1, \quad \quad t_i \geq 0, \quad i=1,\dots, n,
\end{equation}
are solvable in non-negative integers, we can describe all solutions using the following formula.
\begin{equation}\label{eq:gen3mon2sol}
	x_i = \left(a \prod_{j=1}^n u_j^{\alpha_j} + b \prod_{j=1}^n u_j^{\beta_j}\right)^{z_i} \left(c \prod_{j=1}^n u_j^{\gamma_j}\right)^{t_i} w^{-z_i-t_i} \cdot u_i, \quad i=1,\dots, n, 
\end{equation} 
where
$$
u_i \in {\mathbb Z}, \, i=1,\dots, n, \quad 
w \in D\left(a \prod_{j=1}^n u_j^{\alpha_j} + b \prod_{j=1}^n u_j^{\beta_j}\right) \cap D\left(c \prod_{j=1}^n u_j^{\gamma_j}\right).
$$

Let us consider, for example, the equation
$$
x^2 y+z+2 y z = 0.
$$
We may write this equation in the form \eqref{eq:gen3mon2} as
$$
z+2yz=-x^2y.
$$
Then the corresponding equation \eqref{eq:gen3monsyst1} is 
$$
z_3=z_2+z_3=2z_1+z_2-1,
$$ 
which has a non-negative integer solution $(z_1,z_2,z_3)=(2,0,3)$, and the corresponding equation \eqref{eq:gen3monsyst2} is 
$$
t_3=t_2+t_3=2t_1+t_2+1,
$$ 
which has a non-negative integer solution $(t_1,t_2,t_3)=(1,0,3)$. Applying formula \eqref{eq:gen3mon2sol} we obtain
\begin{equation}\label{eq:x2ypzp2yzsol}
\begin{aligned}
(x,y,z)=\left( \frac{-u_1^3 u_2 (u_3+2u_2 u_3)^2}{w^3},u_2, \frac{-u_1^6 u_2^3 u_3 (u_3+2u_2 u_3)^3}{w^6}\right), \\ u_1,u_2,u_3 \in \mathbb{Z}, \quad w \in D(u_3+2u_2 u_3) \cap D(-u_1^2 u_2).
\end{aligned}
\end{equation}
As formula \eqref{eq:gen3mon2sol} only describes non-trivial solutions, we must also add the trivial solutions (those with at least one variable equal to $0$). These solutions are $(x,y,z)=(0,u,0)$ and $(u,0,0)$, for integer $u$. The first solution is included in description \eqref{eq:x2ypzp2yzsol} with $u_1=0$ and $u_2=u$. However, the second solution is not.
Finally, we conclude that all integer solutions are of the form \eqref{eq:x2ypzp2yzsol} or $(x,y,z)=(u,0,0)$ for integer $u$. 

Table \ref{tab:H18mon3systsol} lists equations from Table \ref{tab:H18mon3} that can be written in the form \eqref{eq:gen3mon2} in such a way that systems \eqref{eq:gen3monsyst1} and \eqref{eq:gen3monsyst2} are solvable in non-negative integers. The simplest integer solutions to these systems are included in the last two columns. Then all non-trivial solutions to the corresponding equations are given by  \eqref{eq:gen3mon2sol}, while trivial solutions are easy to find. For every equation in Table \ref{tab:H18mon3systsol}, Table \ref{tab:H18mon3sol} presents the integer solutions obtained by formulas \eqref{eq:gen3mon2sol}.

	\begin{center}

		\captionof{table}{\label{tab:H18mon3systsol} Equations in Table \ref{tab:H18mon3} that can be written in the form \eqref{eq:gen3mon2} such that systems \eqref{eq:gen3monsyst1} and \eqref{eq:gen3monsyst2} are solvable in non-negative integers.}
	\end{center}

	\begin{center}
		%
		\captionof{table}{\label{tab:H18mon3sol} Integer solutions to the equations listed in Table \ref{tab:H18mon3systsol}. Assume that  $u_i \in {\mathbb Z}, \, i=1,\dots, n$.}
	\end{center} 

\subsection{Exercise 4.36}\label{ex:H18mon32}
\textbf{\emph{For every equation listed in Table \ref{tab:H18mon3}, except
$$
x+x^2y-yz^2=0,
$$
(i) check that it satisfies the condition of Proposition \ref{prop:3monsuffcond}, (ii) reduce the equation to the equations with independent monomials, and (iii) solve the resulting equations and hence the original equation.}}

In this exercise, we will solve most of the equations in Table \ref{tab:H18mon3} except $x+x^2y-yz^2=0$. Although many of these equations have been solved in Section \ref{ex:H18mon3}, we will attempt to provide a ``nicer'' solution to the ones given in Table \ref{tab:H18mon3sol}.  Let us first introduce the following Proposition. 
\begin{proposition}\label{prop:3monsuffcond}[Proposition 4.35 in the book]
There is an algorithm for describing all integer solutions to any three-monomial Diophantine equation that can be written in the form \eqref{eq:gen3mon2} such that system \eqref{eq:gen3monsyst1} is solvable in non-negative integers. 
\end{proposition}
To solve these equations in this exercise, we will first check that the equation satisfies the condition in Proposition \ref{prop:3monsuffcond}, we will then reduce the equations to ones with independent monomials and solve the resulting equation, and hence the original equation.

Table \ref{tab:H18mon3syst1} shows that the condition in Proposition \ref{prop:3monsuffcond} is satisfied for the equations in Table \ref{tab:H18mon3} except $x+x^2y-yz^2=0$. 

	\begin{center}
		\begin{tabular}{ |c|c|c|c|c|c| } 
			\hline
			Equation & Form \eqref{eq:gen3mon2} & Solution to \eqref{eq:gen3monsyst1}  \\ 
			\hline\hline
			$x+x^2 y+y^2 z = 0$ & $x+x^2 y=-y^2 z$ & $(0,0,1)$  \\ 	\hline
			  $3 x+x^2 y+z^2 = 0$ & $3x+z^2=-x^2y$ & $(0,1,0)$ \\\hline
			 $x^2 y+z+2 y z = 0$ & $z+2yz=-x^2y$ & $( 2, 0, 3)$   \\\hline
			 $-y+x^2 y+2 z^2 = 0$ & $-y+x^2 y=-2 z^2$& $(0,3,2)$ \\\hline
			 $y+x^2 y+2 z^2 = 0$ & $y+x^2 y=-2 z^2 $ & $(0,3,2)$  \\\hline
			  $x^3+y+2 y z = 0$& $-x^3=y+2 y z $ & $(2,5,0)$  \\ 		\hline
			 $x+x^2 y+2 y z = 0$ & $x+x^2 y=-2 y z $ & $(0,0,1)$ \\\hline
			 $x^3+y+x y z = 0$ & $x^3+y=-x y z $ & $(0,0,1)$  \\ 	\hline
			 $-y+x^2 y+2 x z = 0$ & $-y+x^2 y=-2 x z $ & $(0,1,2)$  \\\hline
			$x+x^2 y+2 z^2 = 0$ & $x+2 z^2 = -x^2 y$  & $(0,1,0)$ \\\hline
			 $x^2 y+z+2 x z = 0$ & $-x^2 y=z+2 x z$ &$(0,2,1)$  \\\hline
			 $x^2 y+z+t y z = 0$ & $-x^2 y=z+t y z $ & $(1,0,1,0)$  \\\hline
			 $x^3-z+y^2 z = 0$& $-x^3=-z+y^2 z $ & $(2,0,5)$  \\ 	\hline
			 $x^3+z+y^2 z = 0$ & $-x^3=z+y^2 z $ & $(2,0,5)$  \\ \hline
			 $y+x^2 y+2 x z = 0$ & $y+x^2 y=-2 x z $ & $(0,1,2)$  \\\hline
			$x^2 y+z+t x z = 0$ & $-x^2 y=z+t x z $ &$(0,2,1,0)$  \\\hline
			  $x^3+y+y^2 z = 0$ & $x^3+y=-y^2 z $& $(0,0,1)$  \\ 	\hline
			$-3 y+x^2 y+z^2 = 0$ & $-3 y+x^2 y=-z^2$ & $(0,3,2)$  \\\hline
			 $3 y+x^2 y+z^2 = 0$ & $3 y+x^2 y=-z^2$ & $(0,3,2)$ \\\hline
			 $-y+x^2 y+x z^2 = 0$ & $-y+x^2 y=-x z^2$  & $(0,3,2)$ \\\hline
			  $y+x^2 y+x z^2 = 0$ & $ y+x^2 y=-x z^2$   & $(0,3,2)$  \\ 		\hline
			 $2 t x+y+x y z = 0$& $-2 t x=y+x y z$ & $(0,1,0,2)$  \\ 
			\hline
		\end{tabular}
		\captionof{table}{\label{tab:H18mon3syst1} Equations listed in Table \ref{tab:H18mon3} (excluding $x+x^2 y -y z^2 = 0$) written in form \eqref{eq:gen3mon2} and their solutions to \eqref{eq:gen3monsyst1} in non-negative integers. }
	\end{center}

We will now attempt to solve the equations in Table \ref{tab:H18mon3syst1}, for the equations which have previously been solved in Table \ref{tab:H18mon3sol} we will attempt to provide a ``nicer'' solution. 

The first equation we will consider is
\begin{equation}\label{xpx2ypy2z}
x+x^2 y+y^2 z = 0.
\end{equation}
Assume that $y\neq 0$. From the equation it is clear that $x=-x^2 y-y^2 z$ is divisible by $y$. But then $x^2$ is divisible by $y^2$, hence $x=-x^2 y-y^2 z$ is divisible by $y^2$. Let $x=y^2u$ for some integer $u$. Substituting this into the equation and cancelling $y^2$, we obtain $u+y^3u^2+z = 0$, hence we can take $y=v$ arbitrary, and then $z=-u-v^3u^2$. Hence,
\begin{equation}\label{xpx2ypy2zsol}
(x,y,z)=(v^2u, v, -u-v^3u^2),  \quad u,v \in \mathbb{Z}. 
\end{equation}
If $y=0$, then \eqref{xpx2ypy2z} implies that $x=0$ and $z$ any. This is covered by \eqref{xpx2ypy2zsol} with $v=0$, $u=-z$.

\vspace{10pt}

The next equation we will consider is
\begin{equation}\label{3xpx2ypz2}
 3 x+x^2 y+z^2 = 0.
 \end{equation}
 This equation can be rearranged as $-x(-3-xy)+z^2=0$ which suggests to make the substitution $t=-3-xy$. This reduces the equation to $-xt+z^2=0$, which, up to the names of variables, we have solved previously, and its integer solutions are 
 $$
 (x,t,z)=(uv^2,uw^2,uvw), \quad u,v,w \in \mathbb{Z},
 $$
 which can be found in Table \ref{table1.38}. This solution with $t=-3-xy$ implies that $uw^2=-3-uv^2y$. We can see that $u$ must divide $3$, so $u\in \{\pm1,\pm 3\}$. If $v=0$, we have the solution $(x,y,z)=(0,w,0)$ for integer $w$. Otherwise, $v \neq 0$ and $y=-\frac{w^2+3/u}{v^2}$ with $u \in  \{\pm1,\pm 3\}$ and $v \in D_2(w^2+3/u)$. Therefore, we have that the integer solutions to equation \eqref{3xpx2ypz2} are
	$$
(x,y,z)=(0,w,0),\left(uv^2, - \frac{w^2+3/u}{v^2},uvw \right), \quad w \in \mathbb{Z}, \quad u \in  \{\pm1,\pm 3\}, \quad v \in D_2(w^2+3/u).
$$

 \vspace{10pt}

 The next equation we will consider is
\begin{equation}\label{x2ypzp2yz} 
x^2 y+z+2 y z = 0.
 \end{equation}
 We can rearrange this equation as $y(x^2+2z)+z=0$. This suggests to make the substitution $X=x^2+2z$ so $z=\frac{X-x^2}{2}$. Then after substituting and multiplying by $2$, we obtain $X(2y+1)-x^2=0$. After making the further substitution $Y=2y+1$, we have $XY-x^2=0$, which, up to the names of variables, we have solved previously, and its integer solutions are
 $$
 (X,Y,x)=(Uv^2,UW^2,UvW), \quad U,v,W \in \mathbb{Z},
 $$
 which can be found in Table \ref{table1.38}. As $Y=2y+1$ we must have that $UW^2$ is odd, so $U$ and $W$ must be odd. Let $U=2u+1$ and $W=2w+1$. Then, in the original variables, we have that the integer solutions to equation \eqref{x2ypzp2yz} are
$$
 \begin{aligned}
(x,y,z)=(v(2 u + 1)  (2 w + 1), 2 w (1 + w) + u (1 + 2 w)^2,  -(1 + 2 u) v^2 (2 w (1 + w) + u (1 + 2 w)^2)), \\ u,w,v \in \mathbb{Z}.
\end{aligned}
 $$
 
 \vspace{10pt}

 The next equation we will consider is	
 \begin{equation}\label{mypx2yp2z2}
-y+x^2 y+2 z^2 = 0
\end{equation} 
We can rearrange this equation as $-y(1-x^2)+2z^2=0$ which suggests to make the substitution ${X=1-x^2}$, and the equation is reduced to $-yX+2z^2=0$ which up to the names of variables, we have solved previously, and its integer solutions are
$$
(z,y,X)=(uvw,2uv^2,uw^2), \quad \text{or} \quad (uvw, uv^2, 2uw^2) \quad u,v,w \in \mathbb{Z},
$$
which can be found in Table \ref{tab:H16mon2bsol}. The first solution with $X=1-x^2$ gives $uw^2=1-x^2$, which, up to the names of variables, we have solved this equation previously, and its integer solutions are 
$$
(w,u,x)=(0,u,\pm 1), \quad \left(n,\frac{1-k^2}{n^2},k\right), \quad n \in D_2(1-k^2)
$$
which can be found in Table \ref{tab:ex:4.1}. Then, in the original variables, we have obtain that a set of integer solutions to equation \eqref{mypx2yp2z2} are
 \begin{equation}\label{mypx2yp2z2_sol1}
(x,y,z)=(\pm 1,u,0), \quad \left(k,\frac{2v^2(1-k^2)}{n^2},\frac{v(1-k^2)}{n}\right), \quad n \in D_2(1-k^2), \quad k,v,u \in \mathbb{Z}.
\end{equation} 

The solution $(z,y,X)=(uvw, uv^2, 2uw^2)$ with $X=1-x^2$ implies that $2uw^2=1-x^2$. Because the left-hand side is even, the right-hand side must be even, hence, $x$ must be odd. Let $x=2k+1$ then $u=\frac{-2k-2k^2}{w^2}$ with $w \in D_2(-2k-2k^2)$. Then, in the original variables, we obtain that a set of integer solutions to equation \eqref{mypx2yp2z2} are
\begin{equation}\label{mypx2yp2z2_sol2}	
(x,y,z)=\left(2k+1,-\frac{v^2(2k+2k^2)}{w^2},\frac{v(2k+2k^2)}{w}\right), \quad w \in D_2(2k+2k^2), \quad k,v \in \mathbb{Z}.
\end{equation}
Therefore all integer solutions to equation \eqref{mypx2yp2z2} are either of the form \eqref{mypx2yp2z2_sol1} or \eqref{mypx2yp2z2_sol2}.
	
	\vspace{10pt}

The next equation we will consider is		 
\begin{equation}\label{ypx2yp2z2}		 
y+x^2 y+2 z^2 = 0. 
\end{equation}		 
We can rearrange this equation as $-y(-1-x^2)+2z^2=0$, which suggests to make the substitution ${X=-1-x^2}$. The equation then reduces to $-yX+2z^2=0$, which, up to the names of variables, we have solved previously, and its integer solutions are
$$
(z,y,X)=(uvw,2uv^2,uw^2), \quad \text{or} \quad (uvw, uv^2, 2uw^2) \quad u,v,w \in \mathbb{Z},
$$
which can be found in Table \ref{tab:H16mon2bsol}. The first solution with $X=-1-x^2$ implies that $uw^2=-1-x^2$. So, $u=\frac{-1-x^2}{w^2}$ where $w \in D_2(-1-x^2)$. Then, in the original variables, we have that a set of integer solutions to equation \eqref{ypx2yp2z2} are
\begin{equation}\label{ypx2yp2z2_sol1}	
(x,y,z)=\left(k,-\frac{2v^2(1+k^2)}{w^2},\frac{v(1+k^2)}{w}\right), \quad w \in D_2(1+k^2), \quad k,v \in \mathbb{Z}.
\end{equation}
The solution $(z,y,X)=(uvw, uv^2, 2uw^2)$ with $X=-1-x^2$ implies that $2uw^2=-1-x^2$. Because the left-hand side is even, the right-hand side must be even, and so $x$ must be odd. Let $x=2k+1$, so $u=\frac{-1-2k-2k^2}{w^2}$ with $w \in D_2(-1-2k-2k^2)$. Then, in the original variables, we have that the other set of integer solutions to equation \eqref{ypx2yp2z2} are
\begin{equation}\label{ypx2yp2z2_sol2}	
(x,y,z)=\left(2k+1,-\frac{v^2(1+2k+2k^2)}{w^2}, \frac{v(1+2k+2k^2)}{w}\right), \quad w \in D_2(1+2k+2k^2), \quad k,v \in \mathbb{Z}.
\end{equation}
Therefore, all integer solutions to equation \eqref{ypx2yp2z2} are either of the form \eqref{ypx2yp2z2_sol1} or \eqref{ypx2yp2z2_sol2}.

\vspace{10pt}
		   
The next equation we will consider is		 
\begin{equation}\label{x3pyp2yz}		 
x^3+y+2 y z = 0.
\end{equation}
We can rearrange this equation as $x^3-y(-1-2z)=0$, which suggests to make the substitution $Z=-1-2z$, and reduces the equation to  $x^3-yZ=0$, which, up to the names of variables, we have solved previously, and its integer solutions are
$$
(x,y,Z)=(u_1 u_2 u_3 u_4,  u_2^3 u_3^2 u_4,u_1^3 u_3 u_4^2), \quad u_1,u_2,u_3,u_4 \in \mathbb{Z},
$$
which can be found in Table \ref{table1.38}. We then have $u_1^3 u_3 u_4^2=-1-2z$. Because the right-hand side is odd, we must have that $u_1,u_3$ and $u_4$ are odd. Let $u_1=2a+1$, $u_2=b$, $u_3=2c+1$ and $u_4=2d+1$ for integers $a,b,c,$ and $d$. Then, in the original variables, we have that the integer solutions to equation \eqref{x3pyp2yz} are
$$		 
(x,y,z)=\left(b(2a+1)(2c+1)(2d+1),b^3(2c+1)^2(2d+1),- \frac{(2a+1)^3(2c+1)(2d+1)^2+1}{2}\right).
$$
	
	\vspace{10pt}
		 
The next equation we will consider is		 
\begin{equation}\label{xpx2yp2yz}		 
x+x^2 y+2 y z = 0.
\end{equation}
The equation implies that $x=yt$ for $t=-x^2-2z=-(yt)^2-2z$, hence $2z=-y^2t^2-t$. This implies that either (i) $t$ is even, or (ii) $y$ and $t$ are both odd. In case (i), $t=2w$ for some integer $w$, $y=u$ is arbitrary, and then $z=(-y^2t^2-t)/2=-2u^2w^2-w$, resulting in
\begin{equation}\label{xpx2yp2yz_sol1}
	(x,y,z)=(2uw,u,-2u^2w^2-w), \quad u,w \in \mathbb{Z}.	
\end{equation}
In case (ii), $t=2w+1$ and $y=2u+1$ for some integers $u,w$, resulting in
\begin{equation}\label{xpx2yp2yz_sol2}
	(x,y,z)=\left((2u+1)(2w+1),2u+1,-\frac{(2w+1)(1+(2u+1)^2(2w+1))}{2}\right), \quad u,w \in \mathbb{Z}.	
\end{equation}
Therefore, all integer solutions to equation \eqref{xpx2yp2yz} are given by either \eqref{xpx2yp2yz_sol1} or \eqref{xpx2yp2yz_sol2}.

\vspace{10pt}

The next equation we will consider is				 
\begin{equation}\label{x3pypxyz}			  
	x^3+y+x y z = 0.
\end{equation}
If $x=0$, then \eqref{x3pypxyz} implies that $x=y=0$, and then $z=u$ can be arbitrary. 
Now assume that $x\neq 0$. Then the equation can be rearranged as $x^3=y(-1-xz)$. Because $-1-xz$ is coprime with $x^3$, this implies that $|-1-xz|=1$, hence $1+xz=1$ or $1+xz=-1$. In the first case, $xz=0$, and, because $x\neq 0$, we must have $z=0$, resulting in $(x,y,z)=(u,-u^3,0)$ for any integer $u$. In the second case, $xz=-2$, which results in the solutions $(x,y,z)=\pm(1,1,-2)$ and $\pm (2,8,-1)$. In conclusion, all integer solutions to equation \eqref{x3pypxyz} are
$$
(x,y,z)=(0,0,u),(u,-u^3,0),\pm(1,1,-2),\pm (2,8,-1), \quad u \in \mathbb{Z}.
$$

	\vspace{10pt}
	
The next equation we will consider is			
\begin{equation}\label{mypx2yp2xz}			 
-y+x^2 y+2 x z = 0.
\end{equation}
If $x=0$ then \eqref{mypx2yp2xz} implies that $x=y=0$, and then $z=v$ can be arbitrary. Now assume that $x\neq 0$. From the equation it is clear that $y=x^2y+2xz$ is divisible by $x$, so let $y=xt$ for some integer $t$. Substituting this into \eqref{mypx2yp2xz} and cancelling $x$, we obtain $-t+x^2t+2z=0$, hence $z=\frac{t-x^2t}{2}$, which is only integer if either (i) $t$ is even, or (ii) $x$ is odd. In case (i), $t=2v$ and $x=u$ for integers $u,v$, which results in the integer solutions $(x,y,z)=(u,2uv,v-u^2v)$, note that this description with $u=0$ also covers the solution with $(x,y,z)=(0,0,v)$. In case (ii), $t=v$ and $x=2u+1$ for integers $u,v$, which results in $(x,y,z)=(2u+1,v(2u+1),-2 u v - 2 u^2 v)$. 
In conclusion, all integer solutions to equation \eqref{mypx2yp2xz} are
$$
(x,y,z)=(u,2uv,v-u^2v) \quad \text{or} \quad (2u+1,v(2u+1),-2 u v - 2 u^2 v), \quad u,v \in \mathbb{Z}.
$$  

	\vspace{10pt}
		
The next equation we will consider is		
\begin{equation}\label{xpx2p2z2}		 
x+x^2 y+2 z^2 = 0.
\end{equation}
We can rearrange this equation to $-x(1+xy)=2z^2$. Because $x$ and $1+xy$ are coprime, we must have that $-x=u_1 v^2$ and $1+xy=u_2 w^2$ where $u_1 u_2 =2$ and $v,w$ are integers. These systems then reduce to solving the four equations (i) $1-v^2y=2w^2$, (ii) $1+v^2y=-2w^2$, (iii) $1-2v^2y=w^2$ and (iv) $1+2v^2y=-w^2$. In all cases, the equations are linear in $y$ and can be easily solved using the divisor function. To ensure that $y$ is integer in cases (iii) and (iv), note that $w$ must be odd, so let $w=2u+1$. In conclusion, all integer solutions to \eqref{xpx2p2z2} are described by
\begin{equation}\label{xpx2p2z2_soli}	
	(x,y,z)=(0,u,0),\left(2v^2,- \frac{2u^2+2u+1}{v^2}, -v(2u+1)\right), \quad v \in D_2(2u^2+2u+1), \quad u \in \mathbb{Z}, 
\end{equation}	 
or,
\begin{equation}\label{xpx2p2z2_sol}	
	(x,y,z)=\left(-2v^2,- \frac{2u^2+2u}{v^2}, v(2u+1)\right), \quad v \in D_2(2u^2+2u), \quad u \in \mathbb{Z},  
\end{equation}
or,
\begin{equation}\label{xpx2p2z2_solii}	
	(x,y,z)=\left(uv^2,- \frac{u+2w^2}{v^2}, vw\right), \quad u \in \{\pm 1\}, \quad v \in D_2(u+2w^2), \quad w \in \mathbb{Z}.  
\end{equation} 

		\vspace{10pt}

The next equation we will consider is		 
\begin{equation}\label{x2ypzp2xz}		 
x^2 y+z+2 x z = 0.
\end{equation}
Assume that $x \neq 0$. From the equation it is clear that $z=-(x^2y+2xz)$ is divisible by $x$. But then $xz$ is divisible by $x^2$, therefore $z=-(x^2y+2xz)$ is divisible by $x^2$. Let $z=x^2v$ for some integer $v$. Substituting this into the equation and cancelling $x^2$, we obtain $ y+v+2 x v = 0$, then taking $x=u$ arbitrary, $y=-v-2uv$. Hence, 
\begin{equation}\label{x2ypzp2xzsol}		 
	(x,y,z)=(u,-v-2uv,u^2v), \quad u,v \in \mathbb{Z}.
\end{equation}
If $x=0$, then \eqref{x2ypzp2xz} implies that $z=0$ and $y$ any. This is covered by \eqref{x2ypzp2xzsol} with $u=0$ and $v=-y$.

\vspace{10pt}
			
The next equation we will consider is
\begin{equation}\label{x2ypzptyz}		 
x^2 y+z+t y z = 0.
\end{equation}
Assume that $y \neq 0$. Then, we can rearrange this equation to $x^2y=z(-1-ty)$. Because $y$ and $-1-ty$ are coprime, $y$ must divide $z$. Let $z=yk$ for some integer $k$. Substituting this into the equation and cancelling $y$, we obtain $x^2=k(-1-ty)$, which suggests to make the substitution $T=-1-ty$, reducing the equation to $x^2=kT$, which, up to the names of variables, we have solved previously, and its integer solutions $(x,k,T)=(uvw,uv^2,uw^2)$ for integers $u,v,w$ can be found in Table \ref{table1.38}. This solution together with $T=-1-ty$, implies that $-1-ty=uw^2$, then taking $y=q$, and $u,w$ arbitrary, then $t=-\frac{uw^2+1}{q}$ with $ q \in D(uw^2+1)$. Then in the original variables, we obtain that the integer solutions to equation \eqref{x2ypzptyz} with $y \neq 0$ are
	$$
	(x,y,z,t)=\left(uvw,q,quv^2,-\frac{uw^2+1}{q}\right), \quad u,v,w \in \mathbb{Z}, \quad q \in D(uw^2+1).
	$$
If $y=0$, then \eqref{x2ypzptyz} implies that $z=0$, hence we also have the integer solutions
$$
(x,y,z,t)=(u,0,0,w), \quad u,w \in \mathbb{Z}.
$$	

	\vspace{10pt}
			 
The next equation we will consider is		 
\begin{equation}\label{x3mzpy2z}			  
x^3-z+y^2 z = 0.
\end{equation}
If $y=\pm 1$, then we obtain the integer solutions 
\begin{equation}\label{x3mzpy2z_sol2}			  
	(x,y,z)=(0,\pm 1,u), \quad u \in \mathbb{Z}.
\end{equation}
 Otherwise, assume that $y \neq \pm 1$. The equation can be rearranged as $x^3-z(1-y^2)=0$, which suggests to make the substitution $Y=1-y^2$ and reduces the equation to $x^3-Yz=0$. Up to the names of variables, we have solved previously, and its integer solutions are
$$
(x,Y,z)=\left(u_1 u_2 u_3 u_4, u_2^3 u_3^2 u_4, u_1^3 u_3 u_4^2\right), \quad u_1,u_2,u_3,u_4 \in \mathbb{Z},
$$ 
which can be found in Table \ref{table1.38}. This solution with $Y=1-y^2$ gives $u_2^3 u_3^2 u_4=1-y^2$, then taking $y=u$ arbitrary then $u_4=\frac{1-u^2}{u_2^3 u_3^2}$ with $u_3 \in D_2(1-u^2)$ and $u_2 \in D_3\left(\frac{1-u^2}{u_3^2}\right)$. Then in the original variables we have the solution to \eqref{x3mzpy2z} is
\begin{equation}\label{x3mzpy2z_sol1}			  
(x,y,z)=\left(\frac{u_1(1-u^2)}{u_2^2 u_3},u,\frac{u_1^3(1-u^2)^2}{u_2^6 u_3^3}\right), \quad  u_3 \in D_2(1-u^2), \quad u_2 \in D_3\left(\frac{1-u^2}{u_3^2}\right), \quad u_1,u \in \mathbb{Z}.
\end{equation}
Therefore, all integer solutions to equation \eqref{x3mzpy2z} are of the form \eqref{x3mzpy2z_sol2} or \eqref{x3mzpy2z_sol1}. 

\vspace{10pt}
			
The next equation we will consider is				 
\begin{equation}\label{x3pzpy2z}			 
x^3+z+y^2 z = 0.
\end{equation}
This equation can be rearranged as $x^3-z(-1-y^2)=0$, which suggests to make the substitution $Y=-1-y^2$ which reduces the equation to $x^3-Yz=0$ which up to the names of variables we have solved previously, and its solution
$$
(x,Y,z)=\left(u_1 u_2 u_3 u_4, u_2^3 u_3^2 u_4, u_1^3 u_3 u_4^2\right), \quad u_1,u_2,u_3,u_4 \in \mathbb{Z}
$$ 
can be found in Table \ref{table1.38}. This solution with $Y=-1-y^2$ gives $u_2^3 u_3^2 u_4=-1-y^2$, then taking $y=u$ arbitrary then $u_4=\frac{-1-u^2}{u_2^3 u_3^2}$ with $u_3 \in D_2(-1-u^2)$ and $u_2 \in D_3\left(\frac{-1-u^2}{u_3^2}\right)$. Then, in the original variables, we have that the integer solutions to equation \eqref{x3pzpy2z} are
$$		  
(x,y,z)=\left(-\frac{u_1(1+u^2)}{u_2^2 u_3},u,\frac{u_1^3(1+u^2)^2}{u_2^6 u_3^3}\right), \quad  u_3 \in D_2(1+u^2), \quad u_2 \in D_3\left(\frac{1+u^2}{u_3^2}\right), \quad u_1,u \in \mathbb{Z}.
$$
	
\vspace{10pt}

The next equation we will consider is			
\begin{equation}\label{ypx2yp2xz}			 
	y+x^2 y+2 x z = 0.
\end{equation}
If $x=0$ then \eqref{mypx2yp2xz} implies that $x=y=0$, and then $z=v$ can be arbitrary. Now assume that $x\neq 0$. From the equation it is clear that $y=-(x^2y+2xz)$ is divisible by $x$, so let $y=xt$ for some integer $t$. Substituting this into \eqref{ypx2yp2xz} and cancelling $x$, we obtain $t+x^2t+2z=0$, hence $z=-\frac{t(1+x^2)}{2}$, which is only integer if either (i) $t$ is even, or (ii) $x$ is odd. In case (i), $t=2v$ and $x=u$ for integers $u,v$, which results in the integer solutions $(x,y,z)=(u,2uv,-v-u^2v)$, note that this description with $u=0$ also covers the solution with $(x,y,z)=(0,0,v)$. In case (ii), $t=v$ and $x=2u+1$ for integers $u,v$, which results in $(x,y,z)=\left(2u+1,v(2u+1),-\frac{v(1+(2u+1)^2}{2}\right)$. 
In conclusion, all integer solutions to equation \eqref{mypx2yp2xz} are given by 
$$
(x,y,z)=(u,2uv,-v-u^2v) \quad \text{or} \quad \left(2u+1,v(2u+1),-v(2u^2+2u+1) \right), \quad u,v \in \mathbb{Z}.
$$  
 
	 \vspace{10pt}

The next equation we will consider is	 
\begin{equation}\label{x2ypzptxz}	 
x^2 y+z+t x z = 0.
\end{equation}
Assume that $x \neq 0$. From the equation it is clear that $z=-(x^2y+txz)$ is divisible by $x$. But then $txz$ is divisible by $x^2$, which implies that $z=-(x^2y+txz)$ is in fact divisible by $x^2$. Let $z=x^2v$ for some integer $v$. Substituting this into the equation and cancelling $x^2$, we obtain $ y+ v+t x v = 0$, then taking $x=u$ and $t=r$ arbitrary, $y=-v-uvr$. Hence, 
\begin{equation}\label{x2ypzptxzsol}
(x,y,z,t)=(u,-v-uvr,u^2v,r), \quad u,v,r \in \mathbb{Z}.
\end{equation}
If $x=0$, then \eqref{x2ypzptxz} implies that $z=0$, and $y,t$ are arbitrary. This description is included in \eqref{x2ypzptxzsol} with $u=0$, $v=-y$ and $r=t$. 

\vspace{10pt}
	 
The next equation we will consider is			
\begin{equation}\label{x3pypy2z}			 
x^3+y+y^2 z = 0.
\end{equation}
If $y=0$ then \eqref{x3pypy2z} implies that $x=0$, and we obtain the integer solutions $(x,y,z)=(0,0,u)$ for integer $u$. Now assume that $y \neq 0$.  
We can rearrange the equation to $x^3=y(-1-yz)$. Because $y$ and $-1-yz$ are coprime, as their product is a perfect cube, they must both be perfect cubes. Hence, $y=u^3$ and $-1-yz=v^3$ for some integers $u,v$. Then, after substituting and rearranging, $z=-\frac{v^3+1}{u^3}$, and we obtain that the integer solutions to equation \eqref{x3pypy2z} with $y \neq 0$ are
\begin{equation}\label{x3pypy2z_sol}
(x,y,z)=\left(uv,u^3,-\frac{v^3+1}{u^3}\right), \quad v \in \mathbb{Z}, \quad u \in D_3(v^3+1).
\end{equation}
Finally, we can conclude that all integer solutions to equation \eqref{x3pypy2z} are of the form \eqref{x3pypy2z_sol} or $(x,y,z)=(0,0,u)$ where $u$ is an arbitrary integer.

\vspace{10pt}

The next equation we will consider is			  
\begin{equation}\label{m3ypx2ypz2}			  
-3 y+x^2 y+z^2 = 0.
\end{equation}
We can rearrange this equation to $y(3-x^2)=z^2$, which suggests to make the substitution $X=3-x^2$, and reduces the equation to $Xy=z^2$, which, up to the names of variables, we have solved previously, and its integer solutions are $(X,y,z)=(uw^2,uv^2,uvw)$, which can be found in Table \ref{table1.38}. Then $X=uw^2=3-x^2$. This equation can be solved by taking $x=k$ arbitrary, then taking any $w \in D_2(3-k^2)$, and then expressing $u$ as $u=\frac{3-k^2}{w^2}$. Then, in the original variables, we have that the integer solutions to equation \eqref{m3ypx2ypz2} are 
$$			  
(x,y,z)=\left(k,\frac{v^2(3-k^2)}{w^2},\frac{v(3-k^2)}{w} \right), \quad k,v \in \mathbb{Z}, \quad w \in D_2(3-k^2).
$$

\vspace{10pt}

The next equation we will consider is			 
\begin{equation}\label{3ypx2ypz2}			  
3 y+x^2 y+z^2 = 0.
\end{equation}
We can rearrange this equation to $-y(3+x^2)=z^2$, which suggests to make the substitution $X=3+x^2$ and $Y=-y$, and reduces the equation to $XY=z^2$, which, up to the names of variables, we have solved previously, and its integer solutions are $(X,Y,z)=(uw^2,uv^2,uvw)$, which can be found in Table \ref{table1.38}. Then $X=uw^2=3+x^2$. This equation can be solved by taking $x=k$ arbitrary, then taking any $w \in D_2(3+k^2)$, and then expressing $u$ as $u=\frac{3+k^2}{w^2}$. Then, in the original variables, we have that the integer solutions to equation \eqref{3ypx2ypz2} are 
$$		  
(x,y,z)=\left(k,-\frac{v^2(3+k^2)}{w^2},\frac{v(3+k^2)}{w} \right), \quad k,v \in \mathbb{Z}, \quad w \in D_2(3+k^2).
$$

\vspace{10pt}

The next equation we will consider is			 
\begin{equation}\label{mypx2ypxz2}			  
-y+x^2 y+x z^2 = 0.
\end{equation}
If $x=0$ then \eqref{mypx2ypxz2} implies that $y=0$, and we obtain the integer solutions $(x,y,z)=(0,0,w)$ for integer $w$. Assume that $x \neq 0$. 
From the equation it is clear that $y=x^2y+xz^2$ is divisible by $x$, so let $y=xt$ for some integer $t$. After substituting this into \eqref{mypx2ypxz2}, cancelling $x$ and rearranging, we obtain
$t(1-x^2)=z^2$. Up to the names of variables and equivalence, we have solved this equation previously, and its integer solutions are $(x,t,z)=(\pm 1,u,0)$ or $\left(u,-\frac{w^2(u^2-1)}{v^2},\frac{w(u^2-1)}{v}\right)$, where $u \in \mathbb{Z}$ and $v \in D_2(u^2-1)$, which can be found in Table \ref{tab:H16substsol}. This implies that the integer solutions to the original equation \eqref{mypx2ypxz2} with $x \neq 0$ are
$$
(x,y,z)=(\pm 1,u,0) \quad \text{or} \quad \left(u,-\frac{uw^2(u^2-1)}{v^2},\frac{w(u^2-1)}{v}\right), \quad u,w \in \mathbb{Z}, \quad v \in D_2(u^2-1).
$$
Note that the integer solution $(x,y,z)=(0,0,w)$ is included in the above description with $u=0$ and $v=1$. 

\vspace{10pt}

The next equation we will consider is			
\begin{equation}\label{ypx2ypxz2}			 
y+x^2 y+x z^2 = 0.
\end{equation}
If $x=0$ then \eqref{ypx2ypxz2} implies that $y=0$, and we obtain the integer solutions $(x,y,z)=(0,0,u)$ for integer $u$. Assume that $x \neq 0$. 
From the equation it is clear that $y=-(x^2y+xz^2)$ is divisible by $x$, so let $y=xt$ for some integer $t$. After substituting this into \eqref{ypx2ypxz2}, cancelling $x$ and rearranging, we obtain $-t(1+x^2)=z^2$. Up to the names of variables, we have solved this equation previously, and its integer solutions are $(x,t,z)=\left(u,-\frac{w^2(u^2+1)}{v^2},\frac{w(u^2-1)}{v}\right)$ for integers $u,v$ where $v \in D_2(u^2+1)$, which can be found in Table \ref{tab:H16substsol}. This solution implies that the integer solutions to equation \eqref{ypx2ypxz2} with $x \neq 0$ are
$$
(x,y,z)= \left(u,-\frac{u w^2(1+u^2)}{v^2},\frac{w(1+u^2)}{v}\right), \quad u,w \in \mathbb{Z}, \quad v \in D_2(1+u^2).
$$
Note that the integer solution $(x,y,z)=(0,0,u)$ is included in the above description with $r=0$ and $w=u$ and $v=1$.
	
\vspace{10pt}

The final equation we will consider is			
\begin{equation}\label{2txpypxyz}			
2 t x+y+x y z = 0.
\end{equation}
If $x=0$, then \eqref{2txpypxyz} implies that $y=0$. Now assume that $x \neq 0$. 
From this equation it is clear that $y$ is divisible by $x$, so let $y=xk$ for some integer $k$. After substituting this into \eqref{2txpypxyz} and cancelling $x$, we obtain $2t+k(1+xz)=0$, hence $t=-\frac{k(1+xz)}{2}$, which is only integer if either (i) $x$ and $z$ are both odd, or (ii) $k$ is even. In case (i), $x=2u+1$, $z=2v+1$ and $k=w$, which results in the integer solutions
\begin{equation}\label{2txpypxyzsol1}		
(x,y,z,t)=(2u+1,w(2u+1),2v+1,-w(2uv+u+v+1)), \quad u,v,k \in \mathbb{Z}.
\end{equation}
In case (ii), $k=2w$, $x=u$ and $z=v$, which results in integer solutions 
\begin{equation}\label{2txpypxyzsol2}		
(x,y,z,t)=(u,2uw,v,-w(1+uv)), \quad u,v,w \in \mathbb{Z}.
\end{equation}
Hence, all integer solutions to equation \eqref{2txpypxyz} are given by either \eqref{2txpypxyzsol1} or \eqref{2txpypxyzsol2}. 

	\begin{center}

		\captionof{table}{\label{tab:H18mon3syst1sol1} Integer solutions to the equations listed in Table \ref{tab:H18mon3} which have a solution to \eqref{eq:gen3monsyst1} in non-negative integers. Unless otherwise stated, assume that all parameters are arbitrary integers. }
	\end{center}

 \subsection{Exercise 4.37}\label{ex:H26mon3mixed}
 \textbf{\emph{Solve all equations from Table \ref{tab:H26mon3mixed}.}}
 
 	\begin{center}
		%
		\captionof{table}{\label{tab:H26mon3mixed} Equations of the forms \eqref{eq:axpbx2ypcyz2} and \eqref{eq:atpbty2pcz2} of size $H\leq 26$.}
	\end{center} 

We will solve the equations in Table \ref{tab:H26mon3mixed}	using the method in Section 4.3.3 of the book, which we summarise below for convenience.
As explained in the book, equations of the form
 \begin{equation}\label{eq:axpbx2ypcyz2}
 ax+bx^2y+cyz^2=0,
 \end{equation}
can be reduced to the form
  \begin{equation}\label{eq:atpbty2pcz2}
  Ax+Bx^2y^2+Cz^2=0.
 \end{equation}
 Equations of the form \eqref{eq:atpbty2pcz2}, can be solved by analysing when the product of $x$ and $A+Bxy^2$ can be a multiple of a perfect square.
 
Equation 
$$
x+x^2y-yz^2=0
$$
is solved in Section 4.3.3 of the book, and its integer solutions are
$$
(x,y,z)=(0,0,w),(0,w,0),\pm (-1,1,0), \quad w \in \mathbb{Z}.
$$

Equation 
$$
2x+x^2y-yz^2=0
$$
is solved in Section 4.3.3 of the book, and its integer solutions are
$$
(x,y,z)=(0,0,w),(0,w,0),\pm (-2,1,0),\pm(-1,2,0), \quad w \in \mathbb{Z}.
$$

Equation 
$$
x+2x^2y-yz^2=0
$$
is solved in Section 4.3.3 of the book, and its integer solutions are
$$
\begin{aligned}
(x,y,z)=& (0,0,w),\quad w \in \mathbb{Z}, \\
& \left(\frac{y_n^2}{t},t,\pm \frac{x_n y_n}{t}\right),  t \in D(y_n), \quad \text{or} \quad \left(-\frac{(y'_n)^2}{t},t,\pm \frac{x'_n y'_n}{t}\right), t \in D(y'_n) ,
\end{aligned}
$$
where $(x_n,y_n)$ is given by \eqref{eq:x2m2y2m1sol} and $(x'_n,y'_n)$ is given by \eqref{eq:x2m2y2p1sol}.

The first equation we will consider is
\begin{equation}\label{3xpx2ymyz2}
3 x+x^2 y-y z^2 = 0.
\end{equation}
Assume that $xyz \neq 0$. This equation implies that $3x$ is divisible by $y$, hence $3x=ty$ for some integer $t \neq 0$. Substituting $x=ty/3$ into \eqref{3xpx2ymyz2}, cancelling $y$ and multiplying by $9$, we obtain
$$
t(9+ty^2)=(3z)^2.
$$
Because any common factor of $t$ and $9+ty^2$ is a divisor of $9=3^2$, their product can be a non-zero perfect square only if $t=eu^2$ and $9+ty^2=ev^2$, where $u,v$ are non-zero integers and $e$ is a divisor of $3$. Then $9+eu^2y^2=ev^2$ or $v^2-(uy)^2=9/e$. This equation is easy to solve, and we obtain non-zero integer solutions 
$(e,uy,v)=(-3,\pm 2, \pm 1),$ $(-1,\pm 5,\pm 4),$ $(1,\pm 4,\pm 5),$ $(3, \pm 1,\pm 2)$.
To obtain integer solutions to \eqref{3xpx2ymyz2} we must have that $x=ty/3$ and $z =\pm euv/3$ are integers, hence the only suitable solutions are
$$
(e,uy,v)=(-3,\pm 2, \pm 1),(3, \pm 1,\pm 2),
$$
which correspond to
$$
\begin{aligned}
(e,u,y,v,t)=& (-3,\pm 2,\pm 1, \pm 1,-12), (-3,\pm 1, \pm 2, \pm 1,-3),(3, \pm 1,\pm 1,\pm 2,3).
\end{aligned}
$$
Then, in the original variables, we obtain that the non-zero integer solutions to equation \eqref{3xpx2ymyz2} are
 \begin{equation}\label{3xpx2ymyz2_sol1}
(x,y,z)= \pm(1,1,\pm 2), \pm(4,-1,\pm 2), \pm(2,-2,\pm 1).
 \end{equation}
We must now consider the case $xyz=0$. We then obtain the integer solutions to \eqref{3xpx2ymyz2}
\begin{equation}\label{3xpx2ymyz2_sol2}
(x,y,z)=(0,0,w),(0,w,0),\pm(-3,1,0),\pm(1,-3,0), \quad w \in \mathbb{Z}.
\end{equation}
Finally, we can conclude that all integer solutions to equation \eqref{3xpx2ymyz2} are described by \eqref{3xpx2ymyz2_sol1} or \eqref{3xpx2ymyz2_sol2}.
 
 \vspace{10pt}

 The next equation we will consider is
\begin{equation}\label{xpx2y2mz2}
x+x^2 y^2-z^2 = 0.
\end{equation}
This equation can be rearranged to $x(1+xy^2)=z^2$. Because $x$ and $1+xy^2$ are coprime, their product can be a perfect square only if $x=eu^2$ and $1+xy^2=ev^2$ where $e$ is a divisor of $1$. Then $1+eu^2y^2=ev^2$ or $v^2-(uy)^2=e$. This equation is easy to solve, and we obtain the integer solutions
$$
(e,uy,v)=(-1,\pm 1,0),(1,0,\pm 1),
$$
which correspond to
 $$
 (e,u,y,v)=(-1,\pm 1,\pm 1,0),(1,0,w,\pm 1),(1,w,0,\pm 1), \quad w \in \mathbb{Z}.
 $$
 Then with $x=eu^2$ and $z=\pm euv$, we can conclude that all the integer solutions to equation \eqref{xpx2y2mz2} are
$$
(x,y,z)=(w^2,0,w),(0,w,0),(-1,\pm 1,0), \quad w \in \mathbb{Z}.
$$

\vspace{10pt}

The next equation we will consider is
\begin{equation}\label{4xpx2ymyz2}
4 x+x^2 y-y z^2 = 0.
\end{equation}
Assume that $xyz \neq 0$. This equation implies that $4x$ is divisible by $y$, hence $4x=ty$ for some integer $t \neq 0$. Substituting $x=ty/4$ into \eqref{4xpx2ymyz2}, cancelling $y$ and multiplying by $16$, we obtain
$$
t(16+ty^2)=(4z)^2.
$$
Because any common factor of $t$ and $16+ty^2$ is a divisor of $16=2^4$, their product can be a non-zero perfect square only if $t=eu^2$ and $16+ty^2=ev^2$, where $u,v$ are non-zero integers and $e$ is a divisor of $2$. Then $16+eu^2y^2=ev^2$ or $v^2-(uy)^2=16/e$. This equation is easy to solve, and we obtain non-zero integer solutions
$$
(e,uy,v)=(-2,\pm 3,\pm 1),(-1,\pm 5,\pm 3),(1,\pm 3, \pm 5),(2,\pm 1, \pm 3).
$$
For all the listed solutions, either $x=ty/4$ or $z =\pm euv/4$ is not an integer, hence \eqref{4xpx2ymyz2} has no integer solutions satisfying $xyz\neq 0$.

We must now consider the case $xyz=0$. We then obtain that all the integer solutions to equation \eqref{4xpx2ymyz2} are
$$
	(x,y,z)=(0,0,w),(0,w,0),\pm(1,-4,0),\pm(2,-2,0),\pm(4,-1,0), \quad w \in \mathbb{Z}.
$$

\vspace{10pt}

The next equation we will consider is
\begin{equation}\label{2xpx2y2mz2}
2 x+x^2 y^2-z^2 = 0.
\end{equation}
This equation can be rearranged to $x(2+xy^2)=z^2$. Because any common factor of $x$ and $2+xy^2$ is a divisor of $2$, their product can be a perfect square only if $x=eu^2$ and $2+xy^2=ev^2$ where $e$ is a divisor of $2$. Then $2+eu^2y^2=ev^2$ or $v^2-(uy)^2=2/e$. This equation is easy to solve, and we obtain the integer solutions 
$$
(e,uy,v)=(-2,\pm 1,0),(2,0,\pm 1),
$$
 which correspond to 
 $$
 (e,u,y,v)=(-2,\pm 1,\pm 1,0),(2,0,w,\pm 1),(2,w,0,\pm 1), \quad w \in \mathbb{Z}.
 $$
 Then with $x=eu^2$ and $z=\pm euv$, we can conclude that all the integer solutions to equation \eqref{2xpx2y2mz2} are
$$
(x,y,z)=(2w^2,0,2w),(0,w,0),(-2,\pm 1,0), \quad w \in \mathbb{Z}.
$$

\vspace{10pt}

The next equation we will consider is
\begin{equation}\label{5xpx2ymyz2}
5 x+x^2 y-y z^2 = 0.
\end{equation}
Assume that $xyz \neq 0$. This equation implies that $5x$ is divisible by $y$, hence $5x=ty$ for some integer $t \neq 0$. Substituting $x=ty/5$ into \eqref{5xpx2ymyz2}, cancelling $y$ and multiplying by $25$, we obtain
$$
t(25+ty^2)=(5z)^2.
$$
Because any common factor of $t$ and $25+ty^2$ is a divisor of $25=5^2$, their product can be a non-zero perfect square only if $t=eu^2$ and $25+ty^2=ev^2$, where $u,v$ are non-zero integers and $e$ is a divisor of $5$. Then $25+eu^2y^2=ev^2$ or $v^2-(uy)^2=25/e$. This equation is easy to solve, and we obtain the non-zero integer solutions $(e,uy,v)=(-5,\pm 3, \pm 2),$ $(-1,\pm 13,\pm 12),$ $(1,\pm 12, \pm 13),$ $(5,\pm 2,\pm 3)$. To obtain integer solutions to \eqref{5xpx2ymyz2} we must have that $x=ty/5$ and $z =\pm euv/5$ are integers, hence the only suitable solutions are
$$
(e,uy,v)=(-5,\pm 3, \pm 2),(5,\pm 2,\pm 3).
$$
Then, in the original variables, we can conclude that all the non-zero integer solutions to equation \eqref{5xpx2ymyz2} are
\begin{equation}\label{5xpx2ymyz2_sol1}
(x,y,z)=\pm(2,2,\pm 3),\pm(3,-3,\pm2), \pm(4,1,\pm 6), \pm(9,-1,\pm 6).
\end{equation}
We must now consider the case $xyz=0$. We then obtain the additional integer solutions to \eqref{5xpx2ymyz2}
\begin{equation}\label{5xpx2ymyz2_sol2}
(x,y,z)=(0,0,w),(0,w,0),\pm(-5,1,0),\pm(-1,5,0),\quad w \in \mathbb{Z}.
\end{equation}
Finally, we can conclude that all integer solutions to equation \eqref{5xpx2ymyz2} are described by \eqref{5xpx2ymyz2_sol1} and \eqref{5xpx2ymyz2_sol2}.

\vspace{10pt}

The next equation we will consider is
\begin{equation}\label{3xpx2y2mz2}
3 x+x^2 y^2-z^2 = 0.
\end{equation}
This equation can be rearranged to $x(3+xy^2)=z^2$. Because any common factor of $x$ and $3+xy^2$ is a divisor of $3$, their product can be a perfect square only if $x=eu^2$ and $3+xy^2=ev^2$ where $e$ is a divisor of $3$. Then $3+eu^2y^2=ev^2$ or $v^2-(uy)^2=3/e$. This equation is easy to solve, and we obtain the integer solutions 
$$
(e,uy,v)=(-3,\pm 1,0),(-1,\pm 2, \pm 1), (1,\pm 1, \pm 2),(3,0,\pm 1).
$$
 Then, in the original variables, we can conclude that all integer solutions to equation \eqref{3xpx2y2mz2} are
$$
\begin{aligned}
(x,y,z)=&(3w^2,0,3w),(0,w,0),(-4,\pm 1,\pm 2),(-3,\pm 1,0), \\ &(-1,\pm 2,\pm 1),(1,\pm 1,\pm 2), \quad w \in \mathbb{Z}.
\end{aligned}
$$

\vspace{10pt}

The next equation we will consider is
\begin{equation}\label{2x2ypzmyz2}
2 x^2 y+z-y z^2 = 0.
 \end{equation}
 Assume that $xyz \neq 0$. This equation implies that $z$ is divisible by $y$, hence $z=ty$ for some integer $t \neq 0$. Substituting this into \eqref{2x2ypzmyz2} and cancelling $y$, we obtain
 $$
 t(ty^2-1)=2x^2.
 $$
 Because $t$ and $ty^2-1$ are coprime, their product can be two times a perfect square only if $t=e_1u^2$  and $ty^2-1=e_2v^2$, where $e_1 e_2=2$. When $e_1 = \pm 2$ we have $v^2-2(uy)^2=\pm 1$ and when $e_1= \pm 1$ we have $(uy)^2-2v^2=\pm 1$. These are Pell's equations which, up to the names of variables, are equations \eqref{eq:x2m2y2m1} and \eqref{eq:x2m2y2p1}, with integer solutions \eqref{eq:x2m2y2m1sol} and \eqref{eq:x2m2y2p1sol}, respectively.
 Then, using $t=e_1 u^2=e_1 (uy)^2/y^2$, $x=\pm u v = \pm (uy) v/y$ and $z=t y = e_1 (uy)^2/y$, we obtain the integer solutions \begin{equation}\label{2x2ypzmyz2_sol1}
 	\begin{aligned}
 		(x,y,z)=\,& \left(\pm\frac{ x'_n y'_n}{w}, w,  \frac{2 (y'_n)^2}{w} \right), \quad w \in D(y'_n), \\
 		& \left(\pm \frac{ x_n y_n}{w}, w,  - \frac{2 (y_n)^2}{w} \right), \quad w \in D(y_n), \\
 		& \left(\pm \frac{x_n y_n}{w}, w,   \frac{x_n^2}{w} \right), \quad w \in D(x_n), \\
 		& \left(\pm \frac{x'_n y'_n}{w}, w,  - \frac{(x'_n)^2}{w} \right), \quad w \in D(x'_n),
 	\end{aligned}
 \end{equation}
 where $(x_n,y_n)$ is given by \eqref{eq:x2m2y2m1sol}, and $(x'_n,y'_n)$ is given by \eqref{eq:x2m2y2p1sol}.
 We must also consider the case $xyz=0$, and we obtain the additional integer solutions
 \begin{equation}\label{2x2ypzmyz2_sol2}
 	(x,y,z)=(0,w,0),(w,0,0),\pm(0,1,1), \quad w \in \mathbb{Z}.
 \end{equation}
 Finally, we can conclude that all integer solutions to \eqref{2x2ypzmyz2} are described by \eqref{2x2ypzmyz2_sol1} or \eqref{2x2ypzmyz2_sol2}.

\vspace{10pt}

 The final equation we will consider is
\begin{equation}\label{xpx2y2m2z2}
x+x^2 y^2-2 z^2 = 0.
\end{equation}
Assume that $y \neq 0$.
This equation can be rearranged to $x(1+xy^2)=2z^2$. Because $x$ and $1+xy^2$ are coprime, their product can be two times a perfect square only if $x=e_1u^2$  and $1+xy^2=e_2v^2$, where $e_1 e_2=2$. These systems reduce to solving Pell's equations \eqref{eq:x2m2y2m1} and \eqref{eq:x2m2y2p1} whose solutions are given by \eqref{eq:x2m2y2m1sol} and \eqref{eq:x2m2y2p1sol}, respectively. 
Then, using $x=e_1 u^2=e_1 (uy)^2/y^2$ and $z=\pm uv = \pm (uy)v/y$, we obtain that the integer solutions to equation \eqref{xpx2y2m2z2} are 
 \begin{equation}\label{xpx2y2m2z2_sol1}
 \begin{aligned}
 (x,y,z)=& \left( \frac{(x'_n)^2}{w^2}, w, \pm\frac{ x'_n y'_n}{w} \right), \quad w \in D(x'_n), \\
& \left(- \frac{x_n^2}{w^2}, w, \pm\frac{ x_n y_n}{w} \right), \quad w \in D(x_n), \\
& \left(\frac{2 y_n^2}{w^2}, w, \pm \frac{y_n x_n}{w} \right), \quad w \in D(y_n), \\
& \left(-\frac{2 (y'_n)^2}{w^2}, w, \pm \frac{y'_n x'_n}{w} \right), \quad w \in D(y'_n), \\
 \end{aligned}
 \end{equation} 
 where $(x_n,y_n)$ are given by \eqref{eq:x2m2y2m1sol} and $(x'_n,y'_n)$ is given by \eqref{eq:x2m2y2p1sol}. We must also consider the case $y=0$, and we obtain the additional integer solutions
 \begin{equation}\label{xpx2y2m2z2_sol2}
 (x,y,z)=(2u^2,0,u) \quad u \in \mathbb{Z}.
 \end{equation}
 To summarise, all integer solutions to equation \eqref{xpx2y2m2z2} are either of the form \eqref{xpx2y2m2z2_sol1} or \eqref{xpx2y2m2z2_sol2}.
 
 \begin{center}
\begin{tabular}{ |c|c|c|c|c|c| } 
\hline
 Equation & Solution $(x,y,z)$ \\ \hline\hline
$x+x^2 y-y z^2 = 0$ & $(0,0,w),(0,w,0),\pm(-1,1,0)$  \\ \hline
 $2 x+x^2 y-y z^2 = 0$ & $(0,0,w),(0,w,0),\pm(-2,1,0),\pm(-1,2,0)$ \\ \hline
 $3 x+x^2 y-y z^2 = 0$ & $(0,0,w),(0,w,0),\pm(-3,1,0),\pm(1,-3,0),$ \\ & $ \pm(1,1,\pm 2), \pm(2,-2,\pm 1),\pm(4,-1,\pm 2)$  \\ \hline
 $x+x^2 y^2-z^2 = 0$ & $(w^2,0,w),(0,w,0),(-1,\pm 1,0)$ \\ \hline
 $4 x+x^2 y-y z^2 = 0$ & $(0,0,w),(0,w,0),\pm(1,-4,0),\pm(2,-2,0),\pm(4,-1,0)$  \\ \hline
 $2 x+x^2 y^2-z^2 = 0$ & $(2w^2,0,2w),(0,w,0),(-2,\pm 1,0)$  \\ \hline
 $5 x+x^2 y-y z^2 = 0$ & $(0,0,w),(0,w,0),\pm(-5,1,0),\pm(-1,5,0),$\\ &$\pm(2,2,\pm 3),\pm(3,-3,\pm2),\pm(4,1,\pm 6),\pm(9,-1,\pm 6)$ \\ \hline
 $3 x+x^2 y^2-z^2 = 0$ & $(3w^2,0,3w),(0,w,0),(-4,\pm 1,\pm 2),$\\ & $(-3,\pm 1,0),(-1,\pm 2,\pm 1),(1,\pm 1,\pm 2)$ \\ \hline
 $2 x^2 y+z-y z^2 = 0$ & $(0,w,0),(w,0,0),\pm(0,1,1)$    \\ 
& $\left(\pm\frac{ x'_n y'_n}{t}, t,  \frac{2 (y'_n)^2}{t} \right), \quad t \in D(y'_n),$ \\
& $\left(\pm \frac{ x_n y_n}{t}, t,  - \frac{2 (y_n)^2}{t} \right), \quad t \in D(y_n),$ \\
& $\left(\pm \frac{x_n y_n}{t}, t,   \frac{x_n^2}{t} \right), \quad t \in D(x_n),$ \\
& $\left(\pm \frac{x'_n y'_n}{t}, t,  - \frac{(x'_n)^2}{t} \right), \quad t \in D(x'_n),$ \\	 \hline
 $x+2 x^2 y-y z^2 = 0$ & $(0,0,w),\left(\frac{y_n^2}{t},t,\pm \frac{x_n y_n}{t}\right),  t \in D(y_n)$ \\ & $\left(-\frac{(y'_n)^2}{t},t,\pm \frac{x'_n y'_n}{t}\right), t \in D(y'_n)$ \\ 
\hline
 $x+x^2 y^2-2 z^2 = 0$ & $(2w^2,0,w),\left( \frac{(x'_n)^2}{t^2}, t, \pm\frac{ x'_n y'_n}{t} \right), \quad t \in D(x'_n),$ \\
& $\left(- \frac{x_n^2}{t^2}, t, \pm\frac{ x_n y_n}{t} \right), \quad t \in D(x_n),$ \\
& $\left(\frac{2 y_n^2}{t^2}, t, \pm \frac{x_n y_n}{t} \right), \quad t \in D(y_n),$ \\
& $\left(-\frac{2 (y'_n)^2}{t^2}, t, \pm \frac{x'_n y'_n}{t} \right), \quad t \in D(y'_n),$ \\ \hline		
\end{tabular}
\captionof{table}{\label{tab:H26mon3mixedsol} Integer solutions to the equations in Table \ref{tab:H26mon3mixed}. In all solutions $w$ is an arbitrary integer, and $(x_n,y_n)$ is given by \eqref{eq:x2m2y2m1sol} and $(x'_n,y'_n)$ is given by \eqref{eq:x2m2y2p1sol}.}
\end{center}

 \subsection{Exercise 4.38}\label{ex:mon3linyhom}
 \textbf{\emph{Find all integer solutions to the equations 
 		$$x^2+x^2 y-y z^2 = 0,$$ $$x^2 y+2 x z-y z^2 = 0$$ and $$2 x^2+x^2 y-y z^2 = 0.$$}}
 
 The first equation we will consider is 
\begin{equation}\label{eq:x2px2ymyz2}
 x^2+x^2y-yz^2=0.
\end{equation}
  Assume that $xyz \neq 0$. Then 
  $|x| \neq |z|$ and $y=\frac{x^2}{z^2-x^2}$ is an integer. Let $d=\text{gcd}(x,z)$, and $x=dx_1$, $z=dz_1$ for some non-zero coprime integers $x_1,z_1$. Then $y=\frac{x_1^2}{z_1^2-x_1^2}$. If $p$ is any prime factor of $z_1^2-x_1^2$, then $p$ must also be a prime factor of $x_1^2$ (otherwise $y$ would not be an integer), but then it is also a factor of $(z_1^2-x_1^2)+x_1^2=z_1^2$, which is a contradiction with $\text{gcd}(x_1,z_1)=1$. So we must have $z_1^2-x_1^2 = \pm 1$ but this impossible with non-zero $x_1$ and $z_1$. Therefore, \eqref{eq:x2px2ymyz2} has no integer solutions satisfying $xyz \neq 0$. We may easily find all solutions with $xyz=0$, and conclude that all integer solutions to equation \eqref{eq:x2px2ymyz2} are
   $$
   (x,y,z)=(u,-1,0),(0,u,0), \quad \text{or} \quad (0,0,u).
   $$
 
 \vspace{10pt}
 
  The next equation we will consider is 
  \begin{equation}\label{eq:x2yp2xzmyz2}
  	x^2y+2xz-yz^2=0.
  \end{equation}
  Assume that $xyz \neq 0$. Then 
  $|x| \neq |z|$ and $y=\frac{2xz}{x^2-z^2}$ is an integer. 
  Let $d=\gcd (x,z)$, and $x=dx_1$, $z=dz_1$ for some non-zero coprime integers $x_1,z_1$. Then $y=\frac{2x_1z_1}{x_1^2-z_1^2}$. If $p \neq 2$ is any prime factor of $x_1^2-z_1^2$, then $p$ must also be a prime factor of $x_1z_1$ (otherwise $y$ would not be an integer), but then $p$ must be a common factor of $x_1$ and $z_1$, which is a contradiction with $\gcd (x_1,z_1)=1$. Hence, $x_1^2-z_1^2$ does not have any prime factors $p\neq 2$, which is possible only if $x_1^2-z_1^2 = \pm 2^k$ for some integer $k \geq 0$. Equations $x_1^2-z_1^2 = \pm 1$ and $x_1^2-z_1^2 = \pm 2$ do not have non-zero integer solutions, hence $k\geq 2$. In this case $x_1^2-z_1^2$ is divisible by $4$, hence $x_1$ and $z_1$ have the same parity. Because $x_1$ and $z_1$ are coprime, they cannot be both even, hence they must be both odd. But then $2x_1z_1$ is not divisible by $4$, hence $y=\frac{2x_1z_1}{x_1^2-z_1^2}$ is not an integer. This contradiction shows that 
  equation \eqref{eq:x2yp2xzmyz2} has no integer solutions satisfying $xyz \neq 0$. We may easily find all solutions with $xyz=0$, and conclude that all integer solutions to \eqref{eq:x2yp2xzmyz2} are
  $$
  (x,y,z)=(u,0,0), (0,u,0), \quad \text{or} \quad (0,0,u) \quad u \in \mathbb{Z}.
  $$

 \vspace{10pt}
 
  The final equation we will consider is 
\begin{equation}\label{eq:2x2px2ymyz2}
 2x^2+x^2y-yz^2=0.
\end{equation}
Assume that $xyz \neq 0$. Then 
$|x| \neq |z|$ and $y=\frac{2x^2}{z^2-x^2}$ is an integer. Let $d=\text{gcd}(x,z)$, and $x=dx_1$, $z=dz_1$ for some non-zero coprime integers $x_1,z_1$. Then $y=\frac{2x_1^2}{z_1^2-x_1^2}$. If $p\neq 2$ is any prime factor of $z_1^2-x_1^2$, then $p$ must also be a prime factor of $x_1^2$ (otherwise $y$ would not be an integer), but then $p$ must also be a common factor of $(z_1^2-x_1^2)+x_1^2=z_1^2$, which is a contradiction with $\text{gcd}(x_1,z_1)=1$. Hence, $x_1^2-z_1^2$ does not have any prime factors $p \neq 2$, which is possible only if 
 $z_1^2-x_1^2 = \pm 2^k$ for some integer $k \geq 0$. Equations $z_1^2-x_1^2 = \pm 1$ and $z_1^2-x_1^2 = \pm 2$ do not have non-zero integer solutions, hence $k \geq 2$. In this case $z_1^2-x_1^2$ is divisible by $4$, hence $x_1$ and $z_1$ must have the same parity. Because they are coprime, they cannot both be even, so they must both be odd. But then $2x_1^2$ is not divisible by $4$, hence, $y=\frac{2x_1^2}{z_1^2-x_1^2}$ is not an integer. This contradiction shows that equation \eqref{eq:2x2px2ymyz2} has no integer solutions satisfying $xyz \neq 0$. We may easily find all solutions with $xyz=0$, and conclude that all integer solutions to \eqref{eq:2x2px2ymyz2} are
 $$
 (x,y,z)=(u,-2,0),(0,u,0), \quad \text{or} \quad (0,0,u), \quad u \in \mathbb{Z}.
 $$

 Table \ref{tab:ex4.37} summarises the integer solutions to equations solved in this exercise.
  \begin{center}
 \begin{tabular}{ |c|c|c|c|c|c| } 
			\hline
			 Equation & Solution $(x,y,z)$ \\ \hline\hline
	 $ x^2+x^2y-yz^2=0$& $(u,-1,0)$, $(0,u,0)$, or $(0,0,u)$ \\\hline
	 $x^2y+2xz-yz^2=0$ & $(u,0,0)$, $(0,u,0)$, or $(0,0,u)$ \\\hline
	 $2x^2+x^2y-yz^2=0$ & $(u,-2,0)$, $(0,u,0)$, or $(0,0,u)$ \\\hline		
\end{tabular}
\captionof{table}{\label{tab:ex4.37} Integer solutions to the equations solved in this exercise. In all solutions $u$ is an arbitrary integer.}
	\end{center} 
 
 \subsection{Exercise 4.43}\label{ex:H32hom3mon}
 \textbf{\emph{Solve all equations listed in Table \ref{tab:H32hom3mon}. For equations of rank $r>0$, present the answer in the form of a list of generators.}}
 
 	\begin{center}
		\begin{tabular}{ |c|c|c|c|c|c| } 
			\hline
			$H$ & Equation & $H$ & Equation & $H$ & Equation \\ 
			\hline\hline
			$24$ & $x^2 y+y^2 z+x z^2 = 0$ & $32$ & $x^3-2 x y^2+y z^2 = 0$ & $32$ & $2 x^3+x y^2+y z^2 = 0$  \\ 
			\hline
			$24$ & $x^3-x y^2+y z^2 = 0$ & $32$ & $x^3-x y^2+2 y z^2 = 0$ & $32$ & $2 x^3-x y^2+z^3 = 0$ \\ 
			\hline
			$24$ & $x^3+y^2 z+y z^2 = 0$ & $32$ & $x^3+x y^2+2 y z^2 = 0$ & $32$ & $2 x^3+x y^2+z^3 = 0$ \\ 
			\hline
			$24$ & $x^3+x y^2+y z^2 = 0$ & $32$ & $x^3+2 x y^2+y z^2 = 0$ & $32$ & $2 x^3+y^3-y z^2 = 0$ \\ 
			\hline
			$24$ & $x^3-x y^2+z^3 = 0$ & $32$ & $x^3-2 x y^2+z^3 = 0$ & $32$ & $2 x^3+y^3+y z^2 = 0$ \\ 
			\hline
			$24$ & $x^3+x y^2+z^3 = 0$ & $32$ & $x^3+2 x y^2+z^3 = 0$ & $32$ & $2 x^3+y^3+z^3 = 0$ \\ 
			\hline	
			$24$ & $x^3+y^3+z^3 = 0$ & $32$ & $2 x^3-x y^2+y z^2 = 0$ &  & \\ 
			\hline
			$32$ & $2 x^2 y+y^2 z+x z^2 = 0$ & $32$ & $2 x^3+y^2 z+y z^2 = 0$ &  & \\ 
			\hline		
		\end{tabular}
		\captionof{table}{\label{tab:H32hom3mon} Homogeneous three-monomial equations of size $H\leq 32$.}
	\end{center} 

To solve the equations in Table \ref{tab:H32hom3mon} we will use the method in Section 4.3.4 of the book, which we summarise below for convenience.	
 We will reduce the equations to ones in Weierstrass form \eqref{eq:Weiform} using a rational change of variables. We can find the rank of this equation using the Magma command
 \begin{equation}\label{cmd:rank}
 {\tt Rank(EllipticCurve([a,b,c,d,e]));}
 \end{equation}
 If the rank of the resulting curve is $r=0$, it has a finite number of rational points, and we can use these points to describe all integer solutions to the original equation. All the rational points can be found using the SageMath command 
\begin{equation}\label{cmd:torsion}
 {\tt EllipticCurve([a,b,c,d,e]).torsion\_points();}
 \end{equation}
 
 By using a rational change of variables, we may always reduce the equations to short Weierstrass form
$$
 y^2=x^3+Ax+B.
$$
Moreover, we can reduce the equations to an elliptic curve in homogeneous Weierstrass form
 \begin{equation}\label{rank1_3var}
	y^2z=x^3+Axz^2+Bz^3
\end{equation}
by making an integer change of variables.

If the rank of the resulting elliptic curve is $r>0$, the curve contains infinitely many rational points, which can be described by their generators, which can be found by using the Magma command
\begin{equation}\label{cmd:generators}
{\tt Generators(EllipticCurve([0,0,0,A,B]));}
\end{equation}
 This outputs solutions of the form $[ (x : y : 1)]$ for some rational numbers $x,y$. By multiplying them by a common denominator and cancelling any common factor, we can obtain a solution to \eqref{rank1_3var} in the form $[(x : y : z)]$, where $x,y,z$ are integers such that $\text{gcd}(x,y,z)=1$. In this case, notation $[(x : y : z)]$ means a family of solutions of the form $(xu,yu,zu)$, $u \in {\mathbb Z}$. 
 Then we can then find all integer points to \eqref{rank1_3var} using the following operation
$$
P_1+P_2=	(x_1: y_1: z_1) + (x_2:y_2:z_2) = (x_3: y_3: z_3)=P_3,
$$
where  
 \begin{equation}\label{eq:ellprod1}
 \begin{aligned}
x_3 = & (x_2 z_1 -  x_1 z_2) ((y_2 z_1 - y_1 z_2)^2 z_1 z_2 - (x_2 z_1 - x_1 z_2)^2 (x_2 z_1 +  x_1 z_2)), \\ 
y_3 = & (y_2 z_1 -    y_1 z_2) ((x_2 z_1 - x_1 z_2)^2 (x_2 z_1 +  2 x_1 z_2) - (y_2 z_1 - y_1 z_2)^2 z_1 z_2) - (x_2 z_1 -
      x_1 z_2)^3 y_1 z_2, \\
      z_3 = & (x_2 z_1 - x_1 z_2)^3 z_1 z_2,
\end{aligned}
\end{equation}
provided that $x_1\neq x_2$, and
\begin{equation}\label{eq:ellprod2}
\begin{aligned}
x_3 = & 2 y_1 z_1 (-8 x_1 y_1^2 z_1 + (3 x_1^2 + A z_1^2)^2) \\
y_3 =&  -8 y_1^4 z_1^2 + (3 x_1^2 + A z_1^2) (12 x_1 y_1^2 z_1 - (3 x_1^2 + A z_1^2)^2) \\
z_3 = & 8 y_1^3 z_1^3
\end{aligned}
\end{equation}
provided that $(x_1: y_1: z_1) = (x_2:y_2:z_2)$ and $y_1 \neq 0$. We also have $O+P=P+O=P$ for any $P=(x : y : z)$, and $P+(-P)=O$ where $-P=(x : -y : z)$, where $O=(0 : 1 : 0)$. 

If $P$ is a rational point on  \eqref{rank1_3var} then we can compute the sequence
\begin{equation}\label{eq:ntimesP}
	P, \quad 2P := P+P, \quad 3P := 2P+P, \quad \dots, \quad nP := (n-1)P+P, \quad \dots
\end{equation}
Also note that by definition $0P=(0:1:0)$ and $(-n)P$ as $-(nP)$.

Equation 
$$
x^2y+y^2z+z^2x=0
$$
is solved in Section 4.3.4 of the book, and its integer solutions are
$$
(x,y,z)=(0,0,u), (0,u,0), \quad \text{and} \quad (u,0,0), \quad u \in \mathbb{Z}.
$$

Equation
$$
x^3+y^3+z^3=0
$$
is solved in Section 4.3.4 of the book, and its integer solutions are
$$
(x,y,z)=(0,u,-u),(-u,0,u), \quad \text{and} \quad (u,-u,0), \quad u \in \mathbb{Z}.
$$

Equation
$$
2x^3+y^3+yz^2=0
$$
is equivalent to $x^3+xy^2+2z^3=0$ which is solved in Section 4.3.4 of the book, and its integer solutions are
$$
\begin{aligned}
(x,y,z)= & \left(- \frac{u x_{2n+1}}{2},u z_{2n+1}, \frac{u y_{2n+1}}{2} \right), \quad \text{or} \\ &  \left(- u x_{2n},2u z_{2n}, u y_{2n} \right) \quad \text{where} \quad (x_n,y_n,z_n)=n(2,-2,1), \quad \text{and} \quad u \in \mathbb{Z}.
\end{aligned}
$$

 The first equation we will consider is
 \begin{equation}\label{x3mxy2pyz2}
 x^3-x y^2+y z^2 = 0.
  \end{equation}
  Assuming that $y \neq 0$, we can divide this equation by $y^3$ and by making the change of variables $X=-x/y$ and $Z=z/y$, the equation is reduced to 
 $$
 X^3-X=Z^2,
 $$ 
 which is an equation in Weierstrass form. The Magma command \eqref{cmd:rank}
 $$
{\tt  Rank(EllipticCurve([0, 0, 0, -1, 0]))}
 $$
outputs that the rank of this equation is $0$, hence the its only rational points are torsion points. We can find these torsion points with the SageMath command \eqref{cmd:torsion}
$$
{\tt EllipticCurve([0,0,0,-1,0]).torsion\_points()}
$$
which outputs $(X,Z)=(\pm 1,0),(0,0)$. Then using $X=-x/y$ and $Z=z/y$, we have that the integer solutions to equation \eqref{x3mxy2pyz2} with $y \neq 0$ are
  \begin{equation}\label{x3mxy2pyz2_sol1}
(x,y,z)=(\pm u,u,0),(0,u,0), \quad u \in \mathbb{Z}.
  \end{equation}
  We must also consider the case $y=0$, which gives the additional integer solutions
    \begin{equation}\label{x3mxy2pyz2_sol2}
(x,y,z)=(0,0,u), \quad u \in \mathbb{Z}.
  \end{equation}
  Finally, all integer solutions to equation \eqref{x3mxy2pyz2} are either of the form \eqref{x3mxy2pyz2_sol1} or \eqref{x3mxy2pyz2_sol2}.
  
 \vspace{10pt}
  
   The next equation we will consider is
 \begin{equation}\label{x3pyz2pyz2}
 x^3+y^2 z+y z^2 = 0.
  \end{equation}
    Assuming that $y \neq 0$, we can divide this equation by $y^3$ and by making the change of variables $X=-x/y$ and $Z=z/y$, the equation is reduced to 
 $$
 X^3=Z^2+Z,
 $$ 
 which is an equation in Weierstrass form. The Magma command \eqref{cmd:rank} for this equation
outputs that the rank is $0$, hence the only rational points are torsion points. We can find these torsion points with the SageMath command \eqref{cmd:torsion} 
which outputs $(X,Z)=(0,0),(0,-1)$. Then using $X=-x/y$ and $Z=z/y$, we have that the integer solutions to equation \eqref{x3pyz2pyz2} with $y \neq 0$ are
  \begin{equation}\label{x3pyz2pyz2_sol1}
(x,y,z)=(0,u,0),(0,u,-u), \quad u \in \mathbb{Z}.
  \end{equation}
  We must also consider the case $y=0$, which gives the additional integer solutions
    \begin{equation}\label{x3pyz2pyz2_sol2}
(x,y,z)=(0,0,u), \quad u \in \mathbb{Z}.
  \end{equation}
  Finally, all integer solutions to equation \eqref{x3pyz2pyz2} are either of the form \eqref{x3pyz2pyz2_sol1} or \eqref{x3pyz2pyz2_sol2}.

\vspace{10pt}

 The next equation we will consider is 
 \begin{equation}\label{x3pxy2pyz2}
 x^3+x y^2+y z^2 = 0.
  \end{equation}
    Assuming that $y \neq 0$, we can divide this equation by $y^3$ and by making the change of variables $X=-x/y$ and $Z=z/y$, the equation is reduced to 
 $$
 X^3+X=Z^2,
 $$ 
 which is an equation in Weierstrass form. The Magma command \eqref{cmd:rank} for this equation
outputs that the rank is $0$, hence the only rational points are torsion points. The SageMath command \eqref{cmd:torsion} for this equation
outputs $(X,Z)=(0,0)$. Then using $X=-x/y$ and $Z=z/y$, we have that the integer solutions to equation \eqref{x3pxy2pyz2} with $y \neq 0$ are
  \begin{equation}\label{x3pxy2pyz2_sol1}
(x,y,z)=(0,u,0), \quad u \in \mathbb{Z}.
  \end{equation}
  We must also consider the case $y=0$, which gives the additional integer solutions
    \begin{equation}\label{x3pxy2pyz2_sol2}
(x,y,z)=(0,0,u), \quad u \in \mathbb{Z}.
  \end{equation}
  Finally, all integer solutions to equation \eqref{x3pxy2pyz2} are either of the form \eqref{x3pxy2pyz2_sol1} or \eqref{x3pxy2pyz2_sol2}.
  
  \vspace{10pt}

 The next equation we will consider is 
 \begin{equation}\label{x3mxy2pz3}
 x^3-x y^2+z^3 = 0.
  \end{equation}
      Assuming that $x \neq 0$, we can divide this equation by $x^3$ and by making the change of variables $Y=y/x$ and $Z=z/x$, the equation is reduced to 
 $$
 Z^3+1=Y^2
 $$ 
 which is an equation in Weierstrass form. The Magma command \eqref{cmd:rank} for this equation
outputs that the rank is $0$, hence the only rational points are torsion points. The SageMath command \eqref{cmd:torsion} for this equation
outputs $(Z,Y)=(-1,0),(0,\pm 1),(2,\pm 3)$. Then using $Y=y/x$ and $Z=z/x$, we have that the integer solutions to equation  \eqref{x3mxy2pz3} with $x \neq 0$ are
  \begin{equation}\label{x3mxy2pz3_sol1}
(x,y,z)=(u,\pm u,0),(u,0,-u),(u,\pm 3 u,2u), \quad u \in \mathbb{Z}.
  \end{equation}
  We must also consider the case $x=0$, which gives the additional integer solutions
    \begin{equation}\label{x3mxy2pz3_sol2}
(x,y,z)=(0,u,0), \quad u \in \mathbb{Z}.
  \end{equation}
  Finally, all integer solutions to equation \eqref{x3mxy2pz3} are either of the form \eqref{x3mxy2pz3_sol1} or \eqref{x3mxy2pz3_sol2}.
 
 \vspace{10pt}
 
  The next equation we will consider is
 \begin{equation}\label{x3pxy2pz3}
 x^3+x y^2+z^3 = 0.
  \end{equation}
        Assuming that $x \neq 0$, we can divide this equation by $x^3$ and by making the change of variables $Y=y/x$ and $Z=-z/x$, the equation is reduced to 
 $$
 Z^3-1=Y^2,
 $$ 
 which is an equation in Weierstrass form. The Magma command \eqref{cmd:rank} for this equation
outputs that the rank is $0$, hence the only rational points are torsion points. The SageMath command \eqref{cmd:torsion} for this equation
outputs $(Z,Y)=(1,0)$. Then using $Y=y/x$ and $Z=-z/x$, we have that the integer solutions to equation \eqref{x3pxy2pz3} with $x \neq 0$ are
  \begin{equation}\label{x3pxy2pz3_sol1}
(x,y,z)=(u,0,-u), \quad u \in \mathbb{Z}.
  \end{equation}
  We must also consider the case $x=0$, which gives the additional integer solutions
    \begin{equation}\label{x3pxy2pz3_sol2}
(x,y,z)=(0,u,0), \quad u \in \mathbb{Z}.
  \end{equation}
  Finally, all integer solutions to equation \eqref{x3pxy2pz3} are either of the form \eqref{x3pxy2pz3_sol1} or \eqref{x3pxy2pz3_sol2}.
 
 \vspace{10pt}
 
  The next equation we will consider is
 \begin{equation}\label{2x2ypy2zpxz2}
 2 x^2 y+y^2 z+x z^2 = 0.
  \end{equation}
    Assume that $xyz \neq 0$. Dividing by $x^3$ and by making the change of variables $Y=y/x$ and $Z=z/x$ the equation is reduced to
  $$
  2Y+Y^2Z+Z^2=0.
  $$
  We can use the Maple command
  $$
  {\tt Weierstrassform(2Y+Y^2Z+Z^2, Y,Z, w, v)}
  $$
  which shows that we can reduce this equation to
  $$
  w^3 +1=v^2,
  $$
  by using the change of variables $Y=\frac{1-v}{w}$ and $Z=-w$.
  The Magma command \eqref{cmd:rank} for this equation
outputs that the rank is $0$, hence the only rational points are torsion points. The SageMath command \eqref{cmd:torsion} for this equation
outputs $(w,v)=(-1,0),(0,\pm 1),(2,\pm 3)$, the second solution is a contradiction with $xyz \neq 0$. Then using $Y=\frac{1-v}{w}$, $Z=-w$, $Y=y/x$ and $Z=z/x$, we have that the integer solutions to equation \eqref{2x2ypy2zpxz2} with $xyz \neq 0$ are
  \begin{equation}\label{2x2ypy2zpxz2_sol1}
(x,y,z)=(u,-u,u),(u,2u,-2u),(u,-u,-2u), \quad u \in \mathbb{Z}.
  \end{equation}
   We must also consider the case $xyz=0$, which gives the additional integer solutions
    \begin{equation}\label{2x2ypy2zpxz2_sol2}
(x,y,z)=(0,0,u),(0,u,0),(u,0,0), \quad u \in \mathbb{Z}.
  \end{equation}
  Finally, all integer solutions to equation \eqref{2x2ypy2zpxz2} are either of the form \eqref{2x2ypy2zpxz2_sol1} or \eqref{2x2ypy2zpxz2_sol2}.
 
 \vspace{10pt}
 
  The next equation we will consider is
\begin{equation}\label{eq:x3m2xy2pyz2}
 x^3-2 x y^2+y z^2 = 0.
\end{equation}
  By making the change of variables $y \leftrightarrow z$ and $x \to -x$, we obtain
$$
  y^2z=x^3-2xz^2.
  $$
The Magma command \eqref{cmd:rank} for this equation 
 outputs that the rank is $1$. Therefore, this equation has infinitely many rational points. The Magma command \eqref{cmd:generators}
 $$
{\tt Generators(EllipticCurve([0, 0, 0, -2, 0]))}
 $$
 outputs $[ (0 : 0 : 1), (-1 : -1 : 1) ]$ which means that all rational points are generated by the points $(x,y,z)=(0,0,1)$ and $(-1,-1,1)$. We have 
 $$
 2(0,0,1) = (0,0,1) + (0,0,1) = (0,1,0),
 $$ 
 so point $(0,0,1)$ has order $2$. Because the rank is $1$, the order of $(-1,-1,1)$ must be infinite. Hence, all rational points on this curve are of the form
 $$
 	P = e (0,0,1) + n (-1,-1,1), \quad e \in \{0,1\}, \quad n \in {\mathbb Z}.
 $$
 Let us denote $(x_n,y_n,z_n)=n (-1,-1,1)$ and $(x'_n,y'_n,z'_n)=(0,0,1) + n (-1,-1,1)$ for $n \in {\mathbb Z}$. Then all the integer solutions to equation \eqref{eq:x3m2xy2pyz2} are
 $$
  	(x,y,z)=(-u x_n,u z_n,u y_n) \quad \text{and} \quad (-u x'_n,u z'_n,u y'_n), \quad u, n \in {\mathbb Z}.
 $$

 \vspace{10pt}

  The next equation we will consider is
 \begin{equation}\label{x3mxy2p2yz2}
 x^3-x y^2+2 y z^2 = 0.
  \end{equation}
      Assuming that $y \neq 0$, we can divide this equation by $y^3$ and by making the change of variables $X=-x/y$ and $Z=z/y$, the equation is reduced to 
 $$
 X^3-X=2Z^2.
 $$ 
 After multiplying the equation by $8$ and making the change of variables $X'=2X$ and $Z'=4Z$, the equation reduces to 
 $$
 (X')^3-4X'=(Z')^2,
 $$
 which is an equation in Weierstrass form.
 The Magma command \eqref{cmd:rank} for this equation outputs 
that the rank is $0$, hence the only rational points are torsion points, and the SageMath command \eqref{cmd:torsion} 
outputs that these points are $(X',Z')=(0,0),(\pm1,0)$. Then using $X'=2X$, $Z'=4Z$, $X=-x/y$ and $Z=z/y$, we have that the integer solutions to equation \eqref{x3mxy2p2yz2} with $y \neq 0$ are
  \begin{equation}\label{x3mxy2p2yz2_sol1}
(x,y,z)=(0,u,0),(\pm u,u,0) \quad u \in \mathbb{Z}.
  \end{equation}
  We must also consider the case $y=0$, which gives the additional integer solutions
    \begin{equation}\label{x3mxy2p2yz2_sol2}
(x,y,z)=(0,0,u), \quad u \in \mathbb{Z}.
  \end{equation}
  Finally, all integer solutions to equation \eqref{x3mxy2p2yz2} are either of the form \eqref{x3mxy2p2yz2_sol1} or \eqref{x3mxy2p2yz2_sol2}.
 
 \vspace{10pt}
  
  The next equation we will consider is
 \begin{equation}\label{x3pxy2p2yz2}
 x^3+x y^2+2 y z^2 = 0.
  \end{equation}
        Assuming that $y \neq 0$, we can divide this equation by $y^3$ and by making the change of variables $X=-x/y$ and $Z=z/y$, the equation is reduced to 
 $$
 X^3+X=2Z^2.
 $$ 
 After multiplying by $8$ and making the change of variables $X'=2X$ and $Z'=4Z$ the equation reduces to 
 $$
 (X')^3+4X'=(Z')^2,
 $$
 which is an equation in Weierstrass form.
 The Magma command \eqref{cmd:rank} for this equation outputs that the rank is $0$, 
hence the only rational points are torsion points, and the SageMath command \eqref{cmd:torsion} outputs that these points are 
$(X',Z')=(0,0),(2,\pm 4)$. Then using $X'=2X$, $Z'=4Z$, $X=-x/y$ and $Z=z/y$, we have that the integer solutions to equation \eqref{x3pxy2p2yz2} with $y \neq 0$ are
  \begin{equation}\label{x3pxy2p2yz2_sol1}
(x,y,z)=(0,u,0),(-u,u,\pm u) \quad u \in \mathbb{Z}.
  \end{equation}
  We must also consider the case $y=0$, which gives the additional integer solutions
    \begin{equation}\label{x3pxy2p2yz2_sol2}
(x,y,z)=(0,0,u), \quad u \in \mathbb{Z}.
  \end{equation}
  Finally, all integer solutions to equation \eqref{x3pxy2p2yz2} are either of the form \eqref{x3pxy2p2yz2_sol1} or \eqref{x3pxy2p2yz2_sol2}.
 
\vspace{10pt}

  The next equation we will consider is
 \begin{equation}\label{x3p2xy2pyz2}
 x^3+2 x y^2+y z^2 = 0.
  \end{equation}
   Assuming that $y \neq 0$, we can divide by $y^3$ and make the change of variables $X=-x/y$ and $Z=z/y$, which reduces the equation to 
  $$
  X^3+2X=Z^2,
  $$
  which is an equation in Weierstrass form. The Magma command \eqref{cmd:rank} for this equation outputs that the rank is $0$, 
 hence the only rational points are torsion points, and the SageMath command \eqref{cmd:torsion} outputs that the only torsion point is 
$(X,Z)=(0,0)$. Then using $X=-x/y$ and $Z=z/y$, we have that the integer solutions to equation \eqref{x3p2xy2pyz2} with $y \neq 0$ are
 \begin{equation}\label{x3p2xy2pyz2_sol1}
(x,y,z)=(0,u,0), \quad u \in \mathbb{Z}.
\end{equation}
We must now consider the case with $y =0$, we then obtain the additional integer solutions
 \begin{equation}\label{x3p2xy2pyz2_sol2}
(x,y,z)=(0,0,u), \quad u \in \mathbb{Z}.
\end{equation}
Finally, we can conclude that all integer solutions to equation \eqref{x3p2xy2pyz2} are of the form \eqref{x3p2xy2pyz2_sol1} or \eqref{x3p2xy2pyz2_sol2}.
 
 \vspace{10pt}

  The next equation we will consider is
\begin{equation}\label{eq:x3m2xy2pz3}
 x^3-2 x y^2+z^3 = 0.
\end{equation}
  Multiplying this equation by $8$ and making the change of variables $X=2z$, $Y=4y$ and $Z=x$, the equation is reduced to
  $$
   Y^2Z=X^3+8Z^3,
  $$
  which is an equation in Weierstrass form.
The Magma command \eqref{cmd:rank} for this equation outputs that the rank is $1$. 
 Therefore, this equation has infinitely many rational points, and the command \eqref{cmd:generators} 
returns $[ (-2 : 0 : 1), (2 : -4 : 1) ]$, which means that all points are generated by $(X,Y,Z)=(-2,0,1)$ and $(2,-4,1)$. We have 
 $$
 2(-2,0,1) = (-2,0,1) + (-2,0,1) = (0,1,0),
 $$ 
  so the point $(-2,0,1)$ has order $2$. Because the rank is $1$, the order of $(2,-4,1)$ must be infinite. Hence, all rational points on this curve are of the form
 $$
 	P = e (-2,0,1) + n (2,-4,1), \quad e \in \{0,1\}, \quad n \in {\mathbb Z}.
 $$
 Let us denote $(x_n,y_n,z_n)=n (2,-4,1)$ and $(x'_n,y'_n,z'_n)=(-2,0,1) + n (2,-4,1)$ for $n \in {\mathbb Z}$. Then all the integer solutions to equation \eqref{eq:x3m2xy2pz3} are 
 $$
(x,y,z)=\left(u z_{2n+1},\frac{u y_{2n+1}}{4},\frac{u x_{2n+1}}{2}\right) \quad \text{and} \quad \left(u z'_{2n},\frac{u y'_{2n}}{4},\frac{u x'_{2n}}{2}\right), 
$$
$$
\text{and} \quad \left(4 u z_{2n}, u y_{2n},2 u x_{2n}\right)  \quad \text{and} \quad \left(4u z'_{2n+1},u y'_{2n+1},2u x'_{2n+1} \right), \quad u, n \in {\mathbb Z}.
$$
Modulo $4$ analysis shows that $y_{2n+1}/4$, $x_{2n+1}/2$, $y'_{2n}/4$ and $x'_{2n}/2$ are integers.

\vspace{10pt}
 
  The next equation we will consider is
 \begin{equation}\label{x3p2xy2pz3}
 x^3+2 x y^2+z^3 = 0.
  \end{equation}
  Assume that $xyz \neq 0$. Dividing by $z^3$ and by making the change of variables $X=x/z$ and $Y=y/z$, the equation is reduced to
  $$
  X^3+2XY^2+1=0
  $$
  The Maple command
  $$
{\tt  Weierstrassform(X^3+2XY^2+1,X, Y, u, v)}
  $$
  shows that we can reduce this equation to $u^3 -8=v^2$, which is in Weierstrass form, using the change of variables $X=-\frac{2}{u}$ and $Y=\frac{v}{2u}$.
  The Magma command \eqref{cmd:rank} for this equation outputs that the rank is $0$, 
hence the only rational points are torsion points, and the SageMath command \eqref{cmd:torsion} outputs that the only torsion point is 
$(u,v)=(2,0)$. This point corresponds to $(X,Y)=(-1,0)$ which is a contradiction with the assumption $xyz\neq 0$. Therefore, all the integer solutions to equation \eqref{x3p2xy2pz3} must have $xyz=0$, and we easily obtain that all integer solutions to equation \eqref{x3p2xy2pz3} are
$$
(x,y,z)=(0,u,0),(u,0,-u), \quad u \in \mathbb{Z}.
$$
 
 \vspace{10pt}

  The next equation we will consider is
 \begin{equation}\label{2x3pxy2pyz2}
 2 x^3+x y^2+y z^2 = 0.
  \end{equation}
  Assuming that $y \neq 0$, we can divide this equation by $y^3$ and by making the substitutions $X=-x/y$ and $Z=z/y$, the equation is reduced to 
 $$
 2X^3+X=Z^2.
 $$  
 After multiplying the equation by $4$ and making the change of variables $X'=2X$ and $Z'=2Z$ the equation reduces to 
 $$
 (X')^3+2X'=(Z')^2,
 $$
 which is an equation in Weierstrass form.
 The Magma command \eqref{cmd:rank} for this equation outputs that the rank is $0$, 
hence the only rational points are torsion points, and the SageMath command \eqref{cmd:torsion} outputs the only torsion point is 
$(X',Z')=(0,0)$. Then using $X'=2X$, $Z'=2Z$, $X=-x/y$ and $Z=z/y$, we have that the integer solutions to equation \eqref{2x3pxy2pyz2} with $y \neq 0$ are
  \begin{equation}\label{2x3pxy2pyz2_sol1}
(x,y,z)=(0,u,0) \quad u \in \mathbb{Z}.
  \end{equation}
  We must also consider the case $y=0$, which gives the additional integer solutions
    \begin{equation}\label{2x3pxy2pyz2_sol2}
(x,y,z)=(0,0,u), \quad u \in \mathbb{Z}.
  \end{equation}
  Finally, all integer solutions to equation \eqref{2x3pxy2pyz2} are either of the form \eqref{2x3pxy2pyz2_sol1} or \eqref{2x3pxy2pyz2_sol2}.
 
 \vspace{10pt}

The next equation we will consider is
$$
 2 x^3-x y^2+z^3 = 0.
$$
        By making the change of variables $x \leftrightarrow z$, we obtain
$$
  y^2z=x^3+2z^3.
 $$
The Magma command \eqref{cmd:rank} outputs that the rank is $1$. 
 Therefore, this equation has infinitely many rational points. The Magma command \eqref{cmd:generators} 
 outputs $[ (-1 : -1 : 1) ]$, which means that all points are generated by $(x,y,z)=(-1,-1,1)$. Because the rank is $1$, the order of $(-1,-1,1)$ must be infinite. Hence, all rational points on this curve are of the form
 $$
 	P =  n (-1,-1,1),  \quad n \in {\mathbb Z}.
 $$
  Let us denote $(x_n,y_n,z_n)=n (-1,-1,1)$ for $n \in {\mathbb Z}$. Then all integer solutions to the original equation are
 $$
  	(x,y,z)=\left(u z_n,u y_n, u x_n \right) \quad u,n \in {\mathbb Z}.
 $$

\vspace{10pt}
 
  The next equation we will consider is
\begin{equation}\label{eq:2x3pxy2pz3}
 2 x^3+x y^2+z^3 = 0.
\end{equation}
      By making the change of variables $z \rightarrow -x$ and $x \rightarrow z$, we obtain
$$
  y^2z=x^3-2z^3.
  $$
The Magma command \eqref{cmd:rank}
 outputs that the rank of this equation is $1$. Therefore, this equation has infinitely many rational points. The command \eqref{cmd:torsion}
 outputs $[ (3 : -5 : 1) ]$, which means that all points are generated by the points $(x,y,z)=(3,-5,1)$. Because the rank is $1$, the order of $(3,-5,1)$ must be infinite. Hence, all rational points on this curve are of the form
 $$
 	P =  n (3,-5,1),  \quad n \in {\mathbb Z}.
 $$
  Let us denote $(x_n,y_n,z_n)=n (3,-5,1)$ for $n \in {\mathbb Z}$. Then all the integer solutions to equation \eqref{eq:2x3pxy2pz3} are
 $$
  	(x,y,z)=\left(u z_n,u y_n, -u x_n \right), \quad u,n \in {\mathbb Z}.
 $$

 \vspace{10pt}

  The next equation we will consider is
 \begin{equation}\label{2x3py3myz2}
 2 x^3+y^3-y z^2 = 0.
  \end{equation}
   Assuming that $y \neq 0$, we can divide this equation by $y^3$ and by making the change of variables $X=x/y$ and $Z=z/y$, the equation is reduced to 
 $$
 2X^3+1=Z^2.
 $$ 
 To reduce this equation to Weierstrass form, we can 
 After multiplying by $4$ and making the change of variables $X'=2X$ and $Z'=2Z$ the equation reduces to 
 $$
 (X')^3+4=(Z')^2,
 $$
 which is an equation in Weierstrass form.
 The Magma command \eqref{cmd:rank}
outputs that the rank of this equation is $0$, hence the only rational points are torsion points, and
the SageMath command \eqref{cmd:torsion} outputs that the only torsion points are
$(X',Z')=(0,\pm 2)$. Then using $X'=2X$, $Z'=2Z$, $X=x/y$ and $Z=z/y$, we have that the integer solutions to equation \eqref{2x3py3myz2} with $y \neq 0$ are
  \begin{equation}\label{2x3py3myz2_sol1}
(x,y,z)=(0,u,\pm u) \quad u \in \mathbb{Z}.
  \end{equation}
  We must also consider the case $y=0$, which gives the additional integer solutions
    \begin{equation}\label{2x3py3myz2_sol2}
(x,y,z)=(0,0,u), \quad u \in \mathbb{Z}.
  \end{equation}
  Finally, all integer solutions to equation \eqref{2x3py3myz2} are either of the form \eqref{2x3py3myz2_sol1} or \eqref{2x3py3myz2_sol2}. 

\vspace{10pt}

The next equation we will consider is
\begin{equation}\label{eq:2x3py3pyz2}
 2 x^3+y^3+y z^2 = 0.
\end{equation}
      Multiplying the equation by $4$ and by making the substitutions $X=-2x$, $Y=2z$ and $Z=y$, the equation is reduced to
$$
  Y^2Z=X^3-4Z^3.
$$
The Magma command \eqref{cmd:rank}
 outputs that the rank of this equation is $1$. Therefore, this equation has infinitely many rational points. The command \eqref{cmd:generators}
 outputs $[ (2 : -2 : 1) ]$, which means that all points are generated by the points $(X,Y,Z)=(2,-2,1)$. Because the rank is $1$, the order of $(2,-2,1)$ must be infinite. Hence, all rational points on this curve are of the form
 $$
 	P =  n (2,-2,1),  \quad n \in {\mathbb Z}.
 $$
  Let us denote $(x_n,y_n,z_n)=n (2,-2,1)$ for $n \in {\mathbb Z}$. Then all the integer solutions to equation \eqref{eq:2x3py3pyz2} are
 $$
  	(x,y,z)=\left(- \frac{u x_{2n+1}}{2},u z_{2n+1}, \frac{u y_{2n+1}}{2} \right), \quad \text{and} \quad \left(- u x_{2n},2u z_{2n}, u y_{2n} \right), \quad u,n \in {\mathbb Z}.
 $$
 Modulo $2$ analysis shows that $x_{2n+1}/2$, $y_{2n+1}/2$ are integers.

\vspace{10pt}

  The next equation we will consider is
 \begin{equation}\label{2x3py3pz3}
 2 x^3+y^3+z^3 = 0.
  \end{equation}
  Assume that $xyz \neq 0$. Dividing by $x^3$ and by making the change of variables $Y=y/x$ and $Z=z/x$, the equation is reduced to
  $$
  Y^3+Z^3+2=0
  $$
  We can use the Maple command
  $$
  {\tt Weierstrassform(Y^3 + Z^3 + 2, Y, Z, u, v)}
  $$
  which returns that we can reduce this equation to $u^3 -27=v^2$, which is an equation in Weierstrass form, using the change of variables $Y=-\frac{9-v}{3u}$ and $Z=-\frac{9+v}{3u}$. The Magma command \eqref{cmd:rank} outputs that the rank is $0$, 
hence the only rational points are torsion points, and the SageMath command \eqref{cmd:torsion} outputs that the only torsion point is
$(u,v)=(3,0)$. Then using $Y=-\frac{9-v}{3u}$, $Z=-\frac{9+v}{3u}$, $Y=y/x$ and $Z=z/x$, we have that the integer solutions to equation \eqref{2x3py3pz3} with $y \neq 0$ are
  \begin{equation}\label{2x3py3pz3_sol1}
(x,y,z)=(-u,u,u) \quad u \in \mathbb{Z}.
  \end{equation}
  We must also consider the case $xyz=0$, which gives the additional integer solutions
    \begin{equation}\label{2x3py3pz3_sol2}
(x,y,z)=(0,u,-u), \quad u \in \mathbb{Z}.
  \end{equation}
  Finally, all integer solutions to equation \eqref{2x3py3pz3} are either of the form \eqref{2x3py3pz3_sol1} or \eqref{2x3py3pz3_sol2}. 
  
  \vspace{10pt}

 The next equation we will consider is	
\begin{equation}\label{eq:2x3mxy2pyz2}
 2 x^3-x y^2+y z^2 = 0.
 \end{equation}
  Multiplying this equation by $4$ and making the change of variables $X=-2x$, $Y=2z$ and $Z=y$, the equation reduces to
$$
   Y^2Z=X^3-2XZ^2.
  $$
The Magma command \eqref{cmd:rank}
 outputs that the rank of this equation is $1$. Therefore, this equation has infinitely many rational points. The command \eqref{cmd:generators}
 outputs $[ (0 : 0 : 1), (-1 : -1 : 1) ]$ which means that all rational points are generated by $(X,Y,Z)=(0,0,1),(-1,-1,1)$. We have 
 $$
 2(0,0,1) = (0,0,1) + (0,0,1) = (0,1,0),
 $$ 
  so the point $(0,0,1)$ has order $2$. Because the rank is $1$, the order of $(-1,-1,1)$ must be infinite. Hence, all rational points on this curve are of the form
 $$
 	P = e (0,0,1) + n (-1,-1,1), \quad e \in \{0,1\}, \quad n \in {\mathbb Z}.
 $$
 Let us denote $(x_n,y_n,z_n)=n (-1,-1,1)$ and $(x'_n,y'_n,z'_n)=(0,0,1) + n (-1,-1,1)$ for $n \in {\mathbb Z}$. Then all the integer solutions to equation \eqref{eq:2x3mxy2pyz2} are
 $$
  	(x,y,z)=\left(-u x_n,2u z_n,u y_n \right) \quad \text{and} \quad \left(-\frac{u x'_n}{2},u z'_n,\frac{u y'_n}{2} \right), \quad u, n \in {\mathbb Z}.
 $$
 Modulo $2$ analysis shows that $x'_n/2$, $y'_n/2$ are integers.
 
 \vspace{10pt}
	
  The final equation we will consider is
 \begin{equation}\label{2x3py2zpyz2}
 2 x^3+y^2 z+y z^2 = 0.
 \end{equation}
    Assuming that $y \neq 0$, we can divide this equation by $y^3$ and by making the change of variables $X=-x/y$ and $Z=z/y$, the equation is reduced to 
 $$
 2X^3=Z^2+Z.
 $$ 
 After multiplying by $4$ and making the change of variables $X'=2X$ and $Z'=2Z$, the equation reduces to 
 $$
 (X')^3=(Z')^2+2Z',
 $$
 which is an equation in Weierstrass form.
 The Magma command \eqref{cmd:rank}
outputs that the rank is $0$, hence the only rational points are torsion points, and the SageMath command \eqref{cmd:torsion} outputs that these points are
$(X',Z')=(-1,-1),(0,-2),(0,0),(2,-4),(2,2)$. Then using $X'=2X$, $Z'=2Z$, $X=-x/y$ and $Z=z/y$, we have that the integer solutions to equation \eqref{2x3py2zpyz2} with $y \neq 0$ are
  \begin{equation}\label{2x3py2zpyz2_sol1}
(x,y,z)=(0,u,-u),(0,u,0),(u,-u,2u),(-u,u,u),(u,2u,-u), \quad u \in \mathbb{Z}.
  \end{equation}
  We must also consider the case $y=0$, which gives the additional integer solutions
    \begin{equation}\label{2x3py2zpyz2_sol2}
(x,y,z)=(0,0,u), \quad u \in \mathbb{Z}.
  \end{equation}
  Finally, all the integer solutions to equation \eqref{2x3py2zpyz2} are either of the form \eqref{2x3py2zpyz2_sol1} or \eqref{2x3py2zpyz2_sol2}. 
  
  Table \ref{tab:H32hom3monsol} summarises the integer solutions to the equations solved in this exercise.
 
 \begin{center}

\captionof{table}{\label{tab:H32hom3monsol} Integer solutions to the equations in Table \ref{tab:H32hom3mon}. In all solutions, $u$ is an arbitrary integer, and notation $n(x,y,z)$ is defined in \eqref{eq:ntimesP}, where the addition is defined in \eqref{eq:ellprod1} and \eqref{eq:ellprod2}.}
\end{center}

 \subsection{Exercise 4.55}\label{ex:2mon3var}
 \textbf{\emph{Find all integer solutions to the equations listed in Table \ref{tab:H562mon3var}. }}
 
 	\begin{center}
		%
		\captionof{table}{\label{tab:H562mon3var} Three-monomial two-variable equations of size $H\leq 56$.}
	\end{center} 
 
The equations in Table \ref{tab:H562mon3var} are of the form
 $$
	a x^n y^q + b x^k y^l + c x^r y^m = 0
$$
 where $a,b,c$ are integers and $n,q,k,l,r,m$ are non-negative integers, and we will solve them using the method in Section 4.3.7 of the book, which we summarise below for convenience. The cases $x=0$ and $y=0$ are easy to check, and so we can assume that $xy\neq 0$. By cancelling the common term $x^{\min\{n,k,r\}} y^{\min\{q,l,m\}}$ we will have equations of the form
 \begin{equation}\label{eq:axnpbxkylpcym}
	a x^n + b x^k y^l + c y^m = 0.
 \end{equation}
 To solve these equations, we will need to use the following proposition.
 \begin{proposition}\label{le:main}[Proposition 4.33 in the book]
	Let $A,B,C$ be any integers such that $A+B+C=0$. Let $p$ be any prime, and let $A_p, B_p$ and $C_p$ be the exponents of $p$ in the prime factorizations of $A,B$ and $C$, respectively. Then the smallest two of the integers $A_p, B_p, C_p$ must be equal.	
\end{proposition}

There is a finite number of pairs $(X_i,Y_i)$, $i=1,\dots,N$ dependent on the divisors of $abc$, such that every solution to \eqref{eq:axnpbxkylpcym} with $xy\neq 0$ satisfies 
$$
x = X_i u^{l'} v^{m'}, \quad y = Y_i u^{n'} v^{k'},
$$ 
for some $i=1,\dots, N$ and some non-zero integers $u$ and $v$, where 
$$
l'=\frac{l}{\text{gcd}(l,n-k)}\,; \quad n'=\frac{n-k}{\text{gcd}(l,n-k)}\,; \quad m'=\frac{m-l}{\text{gcd}(k,m-l)}\,; \quad k'= \frac{k}{\text{gcd}(k,m-l)}.
$$
Note that the prime factors of $abc$ are the only prime factors of $X_i$ and $Y_i$. Let $K$ denote the set of prime divisors of $abc$, then
$$
X_i = \prod_{k \in K} k^{x_k} \quad \text{and} \quad Y_i = \prod_{k \in K} k^{y_k},
$$
where $x_k$ and $y_k$ are the exponents of $k$ in the prime factorisations of $x$ and $y$, respectively. Then equation \eqref{eq:axnpbxkylpcym} can be reduced to 
$$
	a X_i^n v^{nm'-km'-lk'} + b X_i^k Y_i^l + c Y_i^m u^{mn'-kl'-ln'} = 0,
$$
which reduces the equation to solving a finite number of Thue equations. These have a finite number of integer solutions and can be solved using the Magma commands \eqref{thue:magma}. 

The equation
$$
x^4+xy+2y^3=0
$$
is solved in Section 4.3.7 of the book and its integer solutions are
$$
(x,y)=(-1,-1) \quad \text{and} \quad (0,0).
$$
 
 The first family of equations we will consider are the equations of the form
 \begin{equation}\label{ax4pbxypcy3}
 ax^4+bxy+cy^3=0,
 \end{equation}
 where $a,b$ and $c$ are integers. Assume that $xy\neq 0$. 
In this equation, we have $n=4$, $k=l=1$ and $m=3$, so $n-k=3$ and $m-l=2$. Then 
$$
	l'=\frac{1}{\text{gcd}(1,3)}=1 \,; \quad n'=\frac{3}{\text{gcd}(1,3)}=3 \,; \quad m'=\frac{2}{\text{gcd}(1,2)}=2 \,; \quad k'= \frac{1}{\text{gcd}(1,2)}=1 .
$$
So 
\begin{equation}\label{ax4pbxypcy3lmnk}
x = X_i u^{l'} v^{m'}=X_i uv^2, \quad y = Y_i u^{n'} v^{k'}= Y_i u^{3} v,
\end{equation}
for some non-zero integers $u,v$.

\vspace{10pt}

 The next family of equations we will consider are the equations of the form
 \begin{equation}\label{ax5pbxypcy3}
 ax^5+bxy+cy^3=0,
 \end{equation}
 where $a,b$ and $c$ are integers. Assume that $xy\neq 0$. 
In this equation, we have $n=5$, $k=l=1$ and $m=3$, so $n-k=4$ and $m-l=2$. Then 
$$
	l'=\frac{1}{\text{gcd}(1,4)}=1 \,; \quad n'=\frac{4}{\text{gcd}(1,4)}=4 \,; \quad m'=\frac{2}{\text{gcd}(1,2)}=2 \,; \quad k'= \frac{1}{\text{gcd}(1,2)}=1 .
$$
So 
\begin{equation}\label{ax5pbxypcy3lmnk}
x = X_i u^{l'} v^{m'}=X_i uv^2, \quad y = Y_i u^{n'} v^{k'}= Y_i u^{4} v,
\end{equation}
for some non-zero integers $u,v$. 

\vspace{10pt}

 The next family of equations we will consider are the equations of the form
 \begin{equation}\label{ax5pbxypcy4}
 ax^5+bxy+cy^4=0,
 \end{equation}
 where $a,b$ and $c$ are integers. Assume that $xy\neq 0$. 
In this equation, we have $n=5$, $k=l=1$ and $m=4$, so $n-k=4$ and $m-l=3$. Then 
$$
	l'=\frac{1}{\text{gcd}(1,4)}=1 \,; \quad n'=\frac{4}{\text{gcd}(1,4)}=4 \,; \quad m'=\frac{3}{\text{gcd}(1,3)}=3 \,; \quad k'= \frac{1}{\text{gcd}(1,3)}=1 .
$$
So 
\begin{equation}\label{ax5pbxypcy4lmnk}
x = X_i u^{l'} v^{m'}=X_i uv^3, \quad y = Y_i u^{n'} v^{k'}= Y_i u^{4} v,
\end{equation}
 for some non-zero integers $u,v$.
 
 \vspace{10pt}
 
  The next family of equations we will consider are the equations of the form
 \begin{equation}\label{ax4pbxypcy4}
 ax^4+bxy+cy^4=0,
 \end{equation}
 where $a,b$ and $c$ are integers. Assume that $xy\neq 0$. 
 In this equation, we have $n=4=m$ and $k=l=1$, so $n-k=m-l=3$. Then 
$$
	l'=\frac{1}{\text{gcd}(1,3)}=1 \,; \quad n'=\frac{3}{\text{gcd}(1,3)}=3 \,; \quad m'=\frac{3}{\text{gcd}(1,3)}=3 \,; \quad k'= \frac{1}{\text{gcd}(1,3)}=1 .
$$
So 
\begin{equation}\label{ax4pbxypcy4lmnk}
x = X_i u^{l'} v^{m'}=(-1)^{e_x} |X_i| uv^3, \quad y = Y_i u^{n'} v^{k'}=(-1)^{e_y} |Y_i| u^{3} v,
\end{equation} 
 for some non-zero integers $u,v$, and where $e_x=0$ if $x>0$ and $e_x=1$ if $x<0$, and $e_y$ is defined similarly.
 
 \vspace{10pt}
 
The next family of equations we will consider are the equations of the form
 \begin{equation}\label{ax5pbxy2pcy4}
 ax^5+bxy^2+cy^4=0,
 \end{equation}
 where $a,b$ and $c$ are integers. Assume that $xy\neq 0$. 
In this equation, we have $n=5$, $k=1$, $l=2$ and $m=4$, so $n-k=4$ and $m-l=2$. Then 
$$
	l'=\frac{2}{\text{gcd}(2,4)}=1 \,; \quad n'=\frac{4}{\text{gcd}(2,4)}=2 \,; \quad m'=\frac{2}{\text{gcd}(1,2)}=2 \,; \quad k'= \frac{1}{\text{gcd}(1,2)}=1 .
$$
So
\begin{equation}\label{ax5pbxy2pcy4lmnk}
x = X_i u^{l'} v^{m'}=X_i uv^2, \quad y = Y_i u^{n'} v^{k'}= Y_i u^{2} v,
\end{equation}
for some non-zero integers $u,v$.
 
 \vspace{10pt}
 
The first specific equation we will consider is
\begin{equation}\label{x4p2xyp2y3}
x^4+2xy+2y^3=0,
\end{equation}
which is equation \eqref{ax4pbxypcy3} with $a=1$ and $b=c=2$. We have an obvious solution $(x,y)=(0,0)$ and otherwise $x\neq 0$ and $y\neq 0$. The only prime divisor of $abc=2$ is $p=2$. Let $x_2$ and $y_2$ be the exponents of $2$ in the prime factorization of $x$ and $y$, respectively. Then $2$ appears in the prime factorisation of the monomials $x^4$, $2xy$, and $2y^3$ with the exponents $4x_2$, $x_2+y_2+1$ and $3y_2+1$, respectively. Let $M_1=\min\{4x_2,x_2+y_2+1,3y_2+1\}$, and $M_1'$ be the second-smallest among these integers. Then, by Proposition \ref{le:main}, we must have $M_1=M_1'$. The case $x_2+y_2+1 \geq 4x_2 = 3y_2+1$ is clearly impossible for non-negative $x_2,y_2$, hence we must have either (i) $x_2+y_2+1 = 4x_2$ or (ii) $x_2+y_2+1 = 3y_2 +1 $. 

Let us consider these cases separately. In case (i), $y_2+1=3x_2$, so $(x_2,y_2)=(w+1,3w+2)$ for some integer $w\geq 0$. Then substitution \eqref{ax4pbxypcy3lmnk} takes the form
$$
x = 2^{x_2} u v^2 = 2Uv^2, \quad 
y=2^{y_2} u^3 v =2^2 (2^{w})^3 u^3 v = 4U^3 v, \quad \text{where} \quad U=2^{w} u,
$$ 
and $u,v$ are non-zero integers. 
Substituting this into equation \eqref{x4p2xyp2y3} and cancelling $16U^4v^3$, we obtain 
$$
v^5 + 8U^5 = -1.
$$
This is a Thue equation, whose only integer solution is $(v,U)=(-1,0)$. This solution is impossible because $U\neq 0$ by definition.

In case (ii), $x_2= 2y_2$, hence substitution \eqref{ax4pbxypcy3lmnk} takes the form
$$
x = 2^{x_2} u v^2 =  (2^{y_2})^2 u v^2 = uV^2, \quad 
y = 2^{y_2} u^3 v = u^3 V, \quad \text{where} \quad V=2^{y_2} v,
$$ 
and $u,v$ are non-zero integers. 
Substituting this into equation \eqref{x4p2xyp2y3} and cancelling $u^4V^3$, we obtain  
$$
V^5 + 2u^5 = -2.
$$
This is a Thue equation, whose only integer solution is $(V,u)=(0,-1)$, which is impossible because $V\neq 0$ by definition. Therefore, the unique integer solution to equation \eqref{x4p2xyp2y3} is
$$
(x,y)=(0,0).
$$

\vspace{10pt}

 The next equation we will consider is 
 \begin{equation}\label{x4p3xyp2y3}
x^4+3xy+2y^3=0,
 \end{equation}
 which is equation \eqref{ax4pbxypcy3} with $a=1$, $b=3$ and $c=2$.
 We have an obvious solution $(x,y)=(0,0)$ and otherwise $x\neq 0$ and $y\neq 0$. The prime divisors of $abc$ are ${p=2}$ or $p=3$. Let $x_p$ and $y_p$ be the exponents of $p$ in the prime factorization of $x$ and $y$, respectively. 
 Then $2$ appears in the prime factorisation of the monomials $x^4$, $3xy$, and $2y^3$ with the exponents $4x_2$, $x_2+y_2$ and $3y_2+1$, respectively. Let $M_1=\min\{4x_2,x_2+y_2,3y_2+1 \}$, and $M_1'$ be the second-smallest among these integers. Then, by Proposition \ref{le:main}, we must have $M_1=M_1'$. The case $x_2+y_2 \geq 4x_2 = 3y_2+1$ is impossible for non-negative $x_2,y_2$, hence we must have either (2i) $x_2+y_2 = 4x_2$, or (2ii) $x_2+y_2= 3y_2+1$.
Also $3$ appears in the prime factorisation of the monomials $x^4$, $3xy$, and $2y^3$ with the exponents $4x_3$, $x_3+y_3+1$ and $3y_3$, respectively. Let $M_2=\min\{4x_3,x_3+y_3+1,3y_3 \}$, and $M_2'$ be the second-smallest among these integers. Then, by Proposition \ref{le:main}, we must have $M_2=M_2'$. The case $x_3+y_3+1 \geq 4x_3 = 3y_3$ is only possible for $(x_3,y_3)=(0,0)$, hence we must have either (3i) $x_3=y_3=0$, or (3ii) $4x_3=x_3+y_3+1$, or (3iii) $x_3+y_3+1= 3y_3$.
We now have $6$ cases to consider. Let us start with the case where both (2i) and (3i) hold, so, $y_2 = 3x_2$ and $x_3=y_3=0$. Then substitution \eqref{ax4pbxypcy3lmnk} takes the form
$$
x = 2^{x_2} 3^{x_3} u v^2 =2^{x_2} u v^2=Uv^2 , \quad 
y=2^{y_2} 3^{y_3} u^3 v = (2^{x_2})^3 u^3 v=U^3v, \quad \text{where} \quad U=2^{x_2} u,
$$ 
and $u,v$ are non-zero integers.
Substituting this into equation \eqref{x4p3xyp2y3} and cancelling $U^4 v^3$, we obtain 
 $$
2U^5+v^5=-3.
$$
This is a Thue equation whose only integer solution is $(U,v)=(-1,-1)$, from which we can obtain the integer solution 
$(x,y)=(Uv^2,U^3v)=(-1,1)$ to \eqref{x4p3xyp2y3}.

Let us now consider the case where both (2i) and (3ii) hold, so, $y_2=3x_2$ and $3x_3=y_3+1$ or $(x_3,y_3)=(t+1,3t+2)$ for some integer $t\geq 0$. Then substitution \eqref{ax4pbxypcy3lmnk} takes the form
$$
\begin{aligned}
x = & \,\, 2^{x_2} 3^{x_3} u v^2 =2^{x_2} (3)3^{t} u v^2=3Uv^2 , \\
y=& \,\, 2^{y_2} 3^{y_3} u^3 v = (2^{x_2})^3 (3^2)(3^{t})^3 u^3 v=9U^3v, \quad \text{where} \quad U=2^{x_2} 3^t u,
\end{aligned}
$$ 
and $u,v$ are non-zero integers.
Substituting this into equation \eqref{x4p3xyp2y3} and cancelling $81U^4 v^3$, we obtain 
 $$
18U^5+v^5=-1.
$$
This is a Thue equation whose only integer solution is $(U,v)=(0,-1)$, which is impossible because $U \neq 0$ by definition.

Let us now consider the case where both (2i) and (3iii) holds. So, $y_2=3x_2$ and $x_3+1= 2y_3$ or $(x_3,y_3)=(2t+1,t+1)$ for some integer $t \geq 0$. Then substitution \eqref{ax4pbxypcy3lmnk} takes the form
$$
\begin{aligned}
x = & \,\, 2^{x_2} 3^{x_3} u v^2 =2^{x_2} (3)(3^{t})^2 u v^2=3UV^2 , \\
y=& \,\, 2^{y_2} 3^{y_3} u^3 v = (2^{x_2})^3 (3)(3^{t}) u^3 v=3U^3V, \quad \text{where} \quad U=2^{x_2} u, \quad \text{and} \quad V=3^t v,
\end{aligned}
$$ 
and $u,v$ are non-zero integers.
Substituting this into equation \eqref{x4p3xyp2y3} and cancelling $27U^4 V^3$, we obtain 
 $$
2U^5+3V^5=-1.
$$
This is a Thue equation whose only integer solution is $(U,V)=(1,-1)$, from which we obtain the integer solution $(x,y)=(3UV^2,3U^3V)=(3,-3)$ to \eqref{x4p3xyp2y3}.

Let us now consider the case where both (2ii) and (3i) hold, so, $x_2= 2y_2+1$ and $x_3=y_3=0$. Then substitution \eqref{ax4pbxypcy3lmnk} takes the form
$$
x = 2^{x_2} 3^{x_3} u v^2 =2(2^{y_2})^2 u v^2=2uV^2 , \quad 
y=2^{y_2} 3^{y_3} u^3 v = 2^{y_2} u^3 v=u^3V, \quad \text{where} \quad V=2^{y_2} v,
$$ 
and $u,v$ are non-zero integers.
Substituting this into equation \eqref{x4p3xyp2y3} and cancelling $2u^4 V^3$, we obtain 
 $$
u^5+8V^5=-3.
$$
This is a Thue equation and it has no integer solutions.

Let us now consider the case where both (2ii) and (3ii) hold, so, $x_2= 2y_2+1$ and $(x_3,y_3)=(t+1,3t+2)$ for some integer $t \geq 0$. Then substitution \eqref{ax4pbxypcy3lmnk} takes the form
$$ 
\begin{aligned}
x = & \,\, 2^{x_2} 3^{x_3} u v^2 =2(2^{y_2})^2 (3)3^{t} u v^2=6UV^2 , \\
y=& \,\, 2^{y_2} 3^{y_3} u^3 v = 2^{y_2} (3^2)(3^{t})^3 u^3 v=9U^3V, \quad \text{where} \quad U=3^t  u \quad \text{and} \quad V=2^{y_2} v,
\end{aligned}
$$ 
and $u,v$ are non-zero integers.
Substituting this into equation \eqref{x4p3xyp2y3} and cancelling $162U^4 V^3$, we obtain 
 $$
9U^5+8V^5=-1.
$$
This is a Thue equation whose only integer solution is $(U,V)=(-1,1)$, from which we obtain the integer solution $(x,y)=(6UV^2,9U^3V)=(-6,-9)$ to \eqref{x4p3xyp2y3}.

Let us now consider the case where both (2ii) and (3iii) hold, so, $x_2= 2y_2+1$ and $(x_3,y_3)=(2t+1,t+1)$ for some integer $t \geq 0$. Then substitution \eqref{ax4pbxypcy3lmnk} takes the form
$$
\begin{aligned}
x = & \,\, 2^{x_2} 3^{x_3} u v^2 =2(2^{y_2})^2 (3)(3^{t})^2 u v^2=6uV^2 , \\
y=& \,\, 2^{y_2} 3^{y_3} u^3 v = 2^{y_2} (3)3^{t} u^3 v=3u^3V,  \quad \text{where} \quad V=2^{y_2} 3^t v,
\end{aligned}
$$ 
and $u,v$ are non-zero integers.
Substituting this into equation \eqref{x4p3xyp2y3} and cancelling $54u^4 V^3$, we obtain 
 $$
u^5+24V^5=-1.
$$
This is a Thue equation whose only integer solution is $(u,V)=(-1,0)$, which is impossible because $V \neq 0$ by definition.

To summarise, the integer solutions to equation \eqref{x4p3xyp2y3} are
$$
(x,y)=(-6,-9), \quad (-1,1), \quad (3,-3) \quad \text{and} \quad (0,0).
$$

\vspace{10pt} 

 The next equation we will consider is
 \begin{equation}\label{x4pxyp3y3}
x^4+xy+3y^3=0,
\end{equation}
which is equation \eqref{ax4pbxypcy3} with $a=b=1$ and $c=3$.
We have an obvious solution $(x,y)=(0,0)$ and otherwise $x\neq 0$ and $y\neq 0$. The only prime divisor of $abc$ is $p=3$. Let $x_3$ and $y_3$ be the exponents of $3$ in the prime factorization of $x$ and $y$, respectively. Then $3$ appears in the prime factorisation of the monomials $x^4$, $xy$, and $3y^3$ with the exponents $4x_3$, $x_3+y_3$ and $3y_3+1$, respectively. Let $M_1=\min\{4x_3,x_3+y_3,3y_3+1\}$, and $M_1'$ be the second-smallest among these integers. Then we must have $M_1=M_1'$. The case $x_3+y_3 \geq 4x_3 = 3y_3+1$ is clearly impossible for non-negative $x_3,y_3$, hence we must have either (i) $x_3+y_3 = 4x_3$ or (ii) $x_3+y_3= 3y_3 +1 $. 

Let us consider these cases separately. In case (i), $y_3=3x_3$. Then substitution \eqref{ax4pbxypcy3lmnk} takes the form
$$
x = 3^{x_3} u v^2 = Uv^2, \quad 
y=3^{y_3} u^3 v =(3^{x_3})^3 u^3 v = U^3 v, \quad \text{where} \quad U=3^{x_3} u,
$$ 
and $u,v$ are non-zero integers.
Substituting this into equation \eqref{x4pxyp3y3} and cancelling $U^4v^3$, we obtain
$$
v^5 + 3U^5 = -1.
$$
This is a Thue equation whose only integer solution is $(v,U)=(-1,0)$, which is impossible because $U\neq 0$ by definition.

In case (ii), $x_3= 2y_3+1$, hence substitution \eqref{ax4pbxypcy3lmnk} takes the form
$$
x = 3^{x_3} u v^2 =  3(3^{y_2})^2 u v^2 = 3uV^2, \quad 
y = 3^{y_3} u^3 v = u^3 V, \quad \text{where} \quad V=3^{y_3} v,
$$ 
and $u,v$ are non-zero integers.
Substituting this into equation \eqref{x4pxyp3y3} and cancelling $3u^4V^3$, we obtain
$$
27V^5 + u^5 = -1.
$$
This is a Thue equation, whose only integer solution is $(V,u)=(0,-1)$. Therefore the unique integer solution to equation \eqref{x4pxyp3y3} is
$$
(x,y)=(0,0).
 $$

\vspace{10pt}

The next equation we will consider is
\begin{equation}\label{2x4pxypy3}
2x^4+xy+y^3=0,
\end{equation}
which is equation \eqref{ax4pbxypcy3} with $a=2$ and $b=c=1$.
We have an obvious solution $(x,y)=(0,0)$ and otherwise $x\neq 0$ and $y\neq 0$. The only prime divisor of $abc=2$ is $p=2$. Let $x_2$ and $y_2$ be the exponents of $2$ in the prime factorisation of $x$ and $y$, respectively. Then $2$ appears in the prime factorisation of the monomials $2x^4$, $xy$, and $y^3$ with the exponents $4x_2+1$, $x_2+y_2$ and $3y_2$, respectively. Let $M_1=\min\{4x_2+1,x_2+y_2,3y_2 \}$, and $M_1'$ be the second-smallest among these integers. Then we must have $M_1=M_1'$. The case $x_2+y_2 \geq 4x_2+1 = 3y_2$ is clearly impossible for non-negative $x_2,y_2$, hence we must have either (i) $x_2+y_2 = 4x_2+1$ or (ii) $x_2+y_2 = 3y_2 $. Let us consider these cases separately. In case (i), $y_2=3x_2+1$. Then substitution \eqref{ax4pbxypcy3lmnk} takes the form
$$
x = 2^{x_2} u v^2 = Uv^2, \quad 
y=2^{y_2} u^3 v =2 (2^{x_2})^3 u^3 v =2 U^3 v, \quad \text{where} \quad U=2^{x_2} u,
$$ 
and $u,v$ are non-zero integers.
Substituting this into equation \eqref{2x4pxypy3} and cancelling $2U^4v^3$, we obtain
$$
v^5 + 4U^5 = -1.
$$
This is a Thue equation, whose only integer solution is $(v,U)=(-1,0)$, which is impossible because $U\neq 0$ by definition.
In case (ii), $x_2= 2y_2$, hence substitution \eqref{ax4pbxypcy3lmnk} takes the form
$$
x = 2^{x_2} u v^2 =  (2^{y_2})^2 u v^2 = uV^2, \quad 
y = 2^{y_2} u^3 v = u^3 V, \quad \text{where} \quad V=2^{y_2} v,
$$ 
and $u,v$ are non-zero integers.
Substituting this into equation \eqref{2x4pxypy3} and cancelling $u^4V^3$, we obtain
$$
2V^5 + u^5 = -1.
$$
This is a Thue equation whose only integer solutions are $(V,u)=(0,-1),(-1,1)$. The first solution impossible because $V\neq 0$ by definition, while the second one results in the integer solution $(x,y)=(uV^2,u^3V)=(1,-1)$ to \eqref{2x4pxypy3}. 
Therefore the integer solutions to equation \eqref{2x4pxypy3} are
$$
(x,y)=(0,0) \quad \text{and} \quad (1,-1).
$$

\vspace{10pt}

The next equation we will consider is
\begin{equation}\label{x5pxypy3}
x^5+xy+y^3=0,
\end{equation}
which is equation \eqref{ax5pbxypcy3} with $a=b=c=1$.
We have an obvious solution $(x,y)=(0,0)$ and otherwise $x \neq 0$ and $y \neq 0$. There are no prime divisors of $abc=1$, so $X_i=1$ and $Y_i=1$. Therefore $x=uv^2$ and $y=u^4v$, where $u,v$ are non-zero integers. 
Substituting these into equation \eqref{x5pxypy3} and cancelling $u^5 v^3$, we obtain
$$
v^7+u^7=-1.
$$
This is a Thue equation whose only integer solutions are $(v,u)=(-1,0),(0,-1)$, which are impossible as $v\neq 0$ and $u \neq 0$ by definition. Therefore, the unique integer solution to equation \eqref{x5pxypy3} is
$$
(x,y)=(0,0).
$$

\vspace{10pt}

The next equation we will consider is
\begin{equation}\label{x4p4xyp2y3}
x^4+4xy+2y^3=0,
\end{equation}
which is equation \eqref{ax4pbxypcy3} with $a=1$, $b=4$ and $c=2$.
We have an obvious solution $(x,y)=(0,0)$ and otherwise $x\neq 0$ and $y\neq 0$. The only prime divisor of $abc$ is $p=2$. Let $x_2$ and $y_2$ be the exponents of $2$ in the prime factorisation of $x$ and $y$, respectively.  Then $2$ appears in the prime factorisation of the monomials $x^4$, $4xy$, and $2y^3$ with the exponents $4x_2$, $x_2+y_2+2$ and $3y_2+1$, respectively. Let $M_1=\min\{4x_2,x_2+y_2+2,3y_2+1 \}$, and $M_1'$ be the second-smallest among these integers. Then we must have $M_1=M_1'$. We must have either (i) $x_2+y_2+2 = 4x_2$ or (ii) $x_2+y_2+2 = 3y_2+1 $ or (iii) $4x_2 = 3y_2+1$. Let us consider these cases separately. In case (i), $y_2+2=3x_2$, so $(x_2,y_2)=(w+1,3w+1)$ for some integer $w \geq 0$. Then substitution \eqref{ax4pbxypcy3lmnk} takes the form
$$
x = 2^{x_2} u v^2 =2(2^{w}) u v^2 = 2Uv^2, \quad 
y=2^{y_2} u^3 v =2 (2^{w})^3 u^3 v =2 U^3 v, \quad \text{where} \quad U=2^{w} u,
$$ 
and $u,v$ are non-zero integers. 
Substituting this into equation \eqref{x4p4xyp2y3} and cancelling $16 U^4v^3$, we obtain
$$
v^5 + U^5 = -1.
$$
This is a Thue equation whose only integer solutions are $(v,U)=(-1,0),(0,-1)$, which are impossible because $U\neq 0$ and $v\neq 0$ by definition.
In case (ii), $x_2+1= 2y_2$, so $(x_2,y_2)=(2w+1,w+1)$ for some integer $w \geq 0$. Then substitution \eqref{ax4pbxypcy3lmnk} takes the form
$$
x = 2^{x_2} u v^2 = 2 (2^{w})^2 u v^2 = 2 uV^2, \quad 
y = 2^{y_2} u^3 v = 2(2^{w}) u^3 v =2 u^3 V, \quad \text{where} \quad V=2^{w} v,
$$ 
and $u,v$ are non-zero integers. 
Substituting this into equation \eqref{x4p4xyp2y3} and cancelling $16u^4V^3$, we obtain
$$
V^5 + u^5 = -1.
$$
This is a Thue equation whose only integer solutions are $(V,u)=(0,-1),(-1,0)$, which are impossible because $V\neq 0$ and $u\neq 0$ by definition.
In case (iii), $4x_2 = 3y_2+1$, so $(x_2,y_2)=(3w+1,4w+1)$ for some integer $w \geq 0$. Then substitution \eqref{ax4pbxypcy3lmnk} takes the form
$$
x = 2^{x_2} u v^2 = 2^{2w+1} 2^w u v^2= 2^{2w+1}U v^2, \quad 
y = 2^{y_2} u^3 v =2^{w+1}(2^{w})^3 u^3 v =2^{w+1}U^3 v, \quad \text{where} \quad U=2^{w} u,
$$ 
and $u,v$ are non-zero integers. 
Substituting this into equation \eqref{x4p4xyp2y3} and cancelling $2^{3w+3}U^4v^3$, we obtain
$$
2^{5w+1}v^5 +4 U^5 = -1.
$$
It is clear that this equation has no integer solutions with $w \geq 0$, therefore this case produces no integer solutions.
Therefore, the unique integer solution to equation \eqref{x4p4xyp2y3} is
$$
(x,y)=(0,0).
$$

\vspace{10pt}

The next equation we will consider is
\begin{equation}\label{x4p2xyp3y3}
x^4+2xy+3y^3=0,
\end{equation}
which is equation \eqref{ax4pbxypcy3} with $a=1$, $b=2$ and $c=3$.
We have an obvious solution $(x,y)=(0,0)$ and otherwise $x\neq 0$ and $y\neq 0$. For this equation, the prime divisors of $abc$ are $p=2$ or $p=3$. Let $x_p$ and $y_p$ be the exponents of $p$ in the prime factorisation of $x$ and $y$, respectively. Then $2$ appears in the prime factorisation of the monomials $x^4$, $2xy$, and $3y^3$ with the exponents $4x_2$, $x_2+y_2+1$ and $3y_2$, respectively. Let $M_1=\min\{4x_2,x_2+y_2+1,3y_2 \}$, and $M_1'$ be the second-smallest among these integers. Then we must have $M_1=M_1'$. The case $x_2+y_2+1 \geq 4x_2 = 3y_2$ is only possible for $(x_2,y_2)=(0,0)$, hence we must have either (2i) $x_2=y_2=0$, or (2ii) $x_2+y_2 +1= 4x_2$, or (2iii) $x_2+y_2+1= 3y_2$.
Also $3$ appears in the prime factorisation of the monomials $x^4$, $2xy$, and $3y^3$ with the exponents $4x_3$, $x_3+y_3$ and $3y_3+1$, respectively. Let $M_2=\min\{4x_3,x_3+y_3,3y_3+1 \}$, and $M_2'$ be the second-smallest among these integers. Then we must have $M_2=M_2'$. The case $x_3+y_3 \geq 4x_3 = 3y_3+1$ is impossible for non-negative $x_3,y_3$, hence we must have either (3i) $4x_3=x_3+y_3$, or (3ii) $x_3+y_3= 3y_3+1$.
We now have $6$ cases to consider. We will first consider the case where both (2i) and (3i) hold. So, $y_2=x_2=0$ and $3x_3=y_3$. Then substitution \eqref{ax4pbxypcy3lmnk} takes the form
$$
x = 2^{x_2} 3^{x_3} u v^2 =3^{x_3} u v^2=Uv^2 , \quad 
y=2^{y_2} 3^{y_3} u^3 v = (3^{x_3})^3 u^3 v=U^3v, \quad \text{where} \quad U=3^{x_3} u,
$$ 
and $u,v$ are non-zero integers.
Substituting this into equation \eqref{x4p2xyp3y3} and cancelling $U^4 v^3$, we obtain
 $$
3U^5+v^5=-2.
$$
This is a Thue equation whose only integer solution is $(U,v)=(-1,1)$, from which we obtain the integer solution $(x,y)=(Uv^2,U^3v)=(-1,-1)$ to \eqref{x4p2xyp3y3}.

Let us now consider the case where both (2i) and (3ii) hold, so, $y_2=x_2=0$ and $x_3=2y_3+1$. Then substitution \eqref{ax4pbxypcy3lmnk} takes the form
$$
x = 2^{x_2} 3^{x_3} u v^2 =3(3^{y_3})^2 u v^2=3uV^2 , \quad 
y=2^{y_2} 3^{y_3} u^3 v = 3^{y_3} u^3 v=u^3V, \quad \text{where} \quad V=3^{y_3},
$$ 
and $u,v$ are non-zero integers.
Substituting this into equation \eqref{x4p2xyp3y3} and cancelling $3u^4 V^3$, we obtain
 $$
u^5+27V^5=-2.
$$
This is a Thue equation and it has no integer solutions.

Let us now consider the case where both (2ii) and (3i) hold, so, $y_2 +1= 3x_2$ or $(x_2,y_2)=(w+1,3w+2)$ for some integer $w \geq 0$ and $3x_3=y_3$. Then substitution \eqref{ax4pbxypcy3lmnk} takes the form
$$
\begin{aligned}
x =& \,\, 2^{x_2} 3^{x_3} u v^2 =2(2^w)3^{x_3} u v^2=2Uv^2 , \\
y=& \,\, 2^{y_2} 3^{y_3} u^3 v = 2^2(2^w)^3 (3^{x_3})^3 u^3 v=4U^3v, \quad \text{where} \quad U=2^w 3^{x_3} u,
\end{aligned}
$$ 
and $u,v$ are non-zero integers.
Substituting this into equation \eqref{x4p2xyp3y3} and cancelling $16U^4 v^3$, we obtain
 $$
12U^5+v^5=-1.
$$
This is a Thue equation whose only integer solution is $(U,v)=(0,-1)$, which is impossible as $U \neq 0$ by definition.

Let us now consider the case where both (2ii) and (3ii) hold, so, $(x_2,y_2)=(w+1,3w+2)$ for some integer $w \geq 0$ and $x_3=2y_3+1$. Then substitution \eqref{ax4pbxypcy3lmnk} takes the form
$$
\begin{aligned}
x = & \,\, 2^{x_2} 3^{x_3} u v^2 =2(2^w)(3)(3^{x_3})^2 u v^2=6UV^2 , \\
y=& \,\, 2^{y_2} 3^{y_3} u^3 v = 2^2(2^w)^3 3^{y_3} u^3 v=4U^3V, \quad \text{where} \quad U=2^w u, \quad \text{and} \quad V=3^{x_3} v,
\end{aligned}
$$ 
and $u,v$ are non-zero integers.
Substituting this into equation \eqref{x4p2xyp3y3} and cancelling $48U^4 V^3$, we obtain
 $$
4U^5+27V^5=-1.
$$
This is a Thue equation and it has no integer solutions. 

Let us now consider the case where both (2iii) and (3i) hold, so, $x_2+1= 2y_2$ or $(x_2,y_2)=(2w+1,w+1)$ for some integer $w \geq 0$ and $3x_3=y_3$. Then substitution \eqref{ax4pbxypcy3lmnk} takes the form
$$
\begin{aligned}
x =& \,\, 2^{x_2} 3^{x_3} u v^2 =2(2^w)^2 3^{x_3} u v^2=2UV^2 , \\
y= & \,\, 2^{y_2} 3^{y_3} u^3 v = 2(2^w) (3^{x_3})^3 u^3 v=2U^3V, \quad \text{where} \quad U=3^{x_3} u \quad \text{and} \quad V=2^w v,
\end{aligned}
$$ 
and $u,v$ are non-zero integers.
Substituting this into equation \eqref{x4p2xyp3y3} and cancelling $8U^4 V^3$, we obtain
 $$
3U^5+2V^5=-1.
$$
This is a Thue equation whose only integer solution is $(U,V)=(-1,1)$, from which we obtain the integer solution $(x,y)=(2UV^2,2U^3V)=(-2,-2)$ to \eqref{x4p2xyp3y3}.

Let us now consider the case where both (2iii) and (3ii) hold, so, $(x_2,y_2)=(2w+1,w+1)$ for some integer $w \geq 0$ and $x_3=2y_3+1$. Then substitution \eqref{ax4pbxypcy3lmnk} takes the form
$$
\begin{aligned}
x =&\,\, 2^{x_2} 3^{x_3} u v^2 =2(2^w)^2 (3)(3^{y_3})^2 u v^2=6uV^2 , \\
y=& \,\,2^{y_2} 3^{y_3} u^3 v = 2(2^w) 3^{y_3}) u^3 v=2u^3V, \quad \text{where} \quad  V=2^w 3^{y_3} v,
\end{aligned}
$$ 
and $u,v$ are non-zero integers.
Substituting this into equation \eqref{x4p2xyp3y3} and cancelling $24u^4 V^3$, we obtain
 $$
u^5+54V^5=-1.
$$
This is a Thue equation whose only integer solution is $(u,V)=(0,-1)$, which is impossible as $u \neq 0$.

To summarise, the integer solutions to equation \eqref{x4p2xyp3y3} are
$$
(x,y)=(-2,-2), \quad (-1,-1), \quad \text{and} \quad (0,0).
$$

\vspace{10pt}

The next equation we will consider is
\begin{equation}\label{2x4p2xypy3}
2x^4+2xy+y^3=0,
\end{equation}
which is equation \eqref{ax4pbxypcy3} with $a=b=2$ and $c=1$.
We have an obvious solution $(x,y)=(0,0)$ and otherwise $x\neq 0$ and $y\neq 0$. The only prime divisor of $abc$ is $p=2$. Let $x_2$ and $y_2$ be the exponents of $2$ in the prime factorisation of $x$ and $y$, respectively.  Then $2$ appears in the prime factorisation of the monomials $2x^4$, $2xy$, and $y^3$ with the exponents $4x_2+1$, $x_2+y_2+1$ and $3y_2$, respectively. Let $M_1=\min\{4x_2+1,x_2+y_2+1,3y_2 \}$, and $M_1'$ be the second-smallest among these integers. Then we must have $M_1=M_1'$. The case $x_2+y_2+1 \geq 4x_2+1 = 3y_2$ is clearly impossible for non-negative $x_2,y_2$, hence we must have either (i) $x_2+y_2+1 = 4x_2+1$ or (ii) $x_2+y_2+1 = 3y_2$. Let us consider these cases separately. In case (i), $y_2=3x_2$. Then substitution \eqref{ax4pbxypcy3lmnk} takes the form
$$
x = 2^{x_2} u v^2 =2^{x_2} u v^2 = Uv^2, \quad 
y=2^{y_2} u^3 v =(2^{x_2})^3 u^3 v =U^3 v, \quad \text{where} \quad U=2^{x_2} u,
$$ 
and $u,v$ are non-zero integers.
Substituting this into equation \eqref{2x4p2xypy3} and cancelling $U^4 v^3$, we obtain
 $$
2v^5+U^5=-2.
$$
This is a Thue equation whose only integer solution is $(v,U)=(-1,0)$, which is impossible because $U\neq 0$ by definition. In case (ii), $x_2+1=2y_2$ or $(x_2,y_2)=(2w+1,w+1)$ for some integer $w \geq 0$. Then substitution \eqref{ax4pbxypcy3lmnk} takes the form
$$
x = 2^{x_2} u v^2 =2(2^{w})^2 u v^2 = 2uV^2, \quad 
y=2^{y_2} u^3 v =2(2^{w}) u^3 v =2u^3 V, \quad \text{where} \quad V=2^{w} v,
$$
and $u,v$ are non-zero integers.
Substituting this into equation \eqref{2x4p2xypy3} and cancelling $8u^4 V^3$, we obtain
 $$
4V^5+u^5=-1.
$$ 
This is a Thue equation whose only integer solution is $(V,u)=(0,-1)$, which is impossible because $V \neq 0$ by definition. Therefore, the unique integer solution to equation \eqref{2x4p2xypy3} is
$$
(x,y)=(0,0).
$$

 \vspace{10pt}
 
The next equation we will consider is
\begin{equation}\label{x5p2xypy3}
x^5+2xy+y^3=0,
\end{equation}
which is equation \eqref{ax5pbxypcy3} with $a=c=1$ and $b=2$.
We have an obvious solution $(x,y)=(0,0)$ and otherwise $x\neq 0$ and $y\neq 0$. The only prime divisor of $abc=2$ is $p=2$. Let $x_2$ and $y_2$ be the exponents of $2$ in the prime factorisation of $x$ and $y$, respectively.  Then $2$ appears in the prime factorisation of the monomials $x^5$, $2xy$, and $y^3$ with the exponents $5x_2$, $x_2+y_2+1$ and $3y_2$, respectively. Let $M_1=\min\{5x_2,x_2+y_2+1,3y_2 \}$, and $M_1'$ be the second-smallest among these integers. Then we must have $M_1=M_1'$. The case $x_2+y_2+1 \geq 5x_2 = 3y_2$ only holds for $(x_2,y_2)=(0,0)$ hence we must have either (i) $x_2=y_2=0$ or (ii) $x_2+y_2+1 = 5x_2$ or (iii) $x_2+y_2+1 = 3y_2$. Let us consider these cases separately. In case (i),  $x_2=y_2=0$, then substitution \eqref{ax5pbxypcy3lmnk} takes the form
$$
x = 2^{x_2} u v^2 = u v^2  \quad 
y=2^{y_2} u^4 v = u^4 v,
$$
where $u,v$ are non-zero integers.
Substituting this into equation \eqref{x5p2xypy3} and cancelling $u^5 v^3$, we obtain
 $$
v^7+u^7=-2.
$$
This is a Thue equation whose only integer solution is $(u,v)=(-1,-1)$, from which we obtain the integer solution $(x,y)=(uv^2,u^4v)=(-1,-1)$ to \eqref{x5p2xypy3}. In case (ii), $y_2+1=4x_2$, or $(x_2,y_2)=(w+1,4w+3)$ for some integer $w \geq 0$. Then substitution \eqref{ax5pbxypcy3lmnk} takes the form
$$
x = 2^{x_2} u v^2 =2(2^{w}) u v^2 = 2Uv^2, \quad 
y=2^{y_2} u^4 v =2^3(2^{w})^4 u^4 v =8U^4 v, \quad \text{where} \quad U=2^{w},
$$
and $u,v$ are non-zero integers.
Substituting this into equation \eqref{x5p2xypy3} and cancelling $32U^5 v^3$, we obtain
 $$
v^7+16U^7=-1.
$$
This is a Thue equation whose only integer solution is $(U,v)=(0,-1)$, which is impossible because $U \neq 0$ by definition. In case (iii), $x_2+1 = 2y_2$, or $(x_2,y_2)=(2w+1,w+1)$ for some integer $w \geq 0$. Then substitution \eqref{ax5pbxypcy3lmnk} takes the form
$$
x = 2^{x_2} u v^2 =2(2^{w})^2 u v^2 = 2uV^2, \quad 
y=2^{y_2} u^4 v =2(2^{w}) u^4 v =2u^4 V, \quad \text{where} \quad V=2^{w} v,
$$
and $u,v$ are non-zero integers.
Substituting this into equation \eqref{x5p2xypy3} and cancelling $8u^5 V^3$, we obtain
 $$
4V^7+u^7=-1.
$$ 
This is a Thue equation whose only integer solution is $(V,u)=(0,-1)$, which is impossible because $V \neq 0$ by definition. Therefore, the integer solutions to equation \eqref{x5p2xypy3} are
$$
(x,y)=(-1,-1) \quad \text{and} \quad (0,0).
$$

\vspace{10pt}

The next equation we will consider is
\begin{equation}\label{x4p5xyp2y3}
x^4+5xy+2y^3=0,
\end{equation}
which is equation \eqref{ax4pbxypcy3} with $a=1$, $b=5$ and $c=2$. 
We have an obvious solution $(x,y)=(0,0)$ and otherwise $x\neq 0$ and $y\neq 0$. The prime divisors of $abc$ are $p=2$ or $p=5$. Let $x_p$ and $y_p$ be the exponents of $p$ in the prime factorisation of $x$ and $y$, respectively. 
 Then $2$ appears in the prime factorisation of the monomials $x^4$, $5xy$, and $2y^3$ with the exponents $4x_2$, $x_2+y_2$ and $3y_2+1$, respectively. Let $M_1=\min\{4x_2,x_2+y_2,3y_2+1 \}$, and $M_1'$ be the second-smallest among these integers. Then we must have $M_1=M_1'$. The case $x_2+y_2 \geq 4x_2 = 3y_2+1$ is impossible for non-negative $x_2,y_2$, hence we must have either (2i) $x_2+y_2 = 4x_2$, or (2ii) $x_2+y_2= 3y_2+1$.
Also $5$ appears in the prime factorisation of the monomials $x^4$, $5xy$, and $2y^3$ with the exponents $4x_5$, $x_5+y_5+1$ and $3y_5$, respectively. Let $M_2=\min\{4x_5,x_5+y_5+1,3y_5 \}$, and $M_2'$ be the second-smallest among these integers. Then we must have $M_2=M_2'$. The case $x_5+y_5+1 \geq 4x_5 = 3y_5$ is only possible for $(x_5,y_5)=(0,0)$, hence we must have either (5i) $x_5=y_5=0$, or (5ii) $4x_5=x_5+y_5+1$, or (5iii) $x_5+y_5+1= 3y_5$.
We now have $6$ cases to consider. We will first consider the case where both (2i) and (5i) hold, so, $y_2 = 3x_2$ and $x_5=y_5=0$. Then substitution \eqref{ax4pbxypcy3lmnk} takes the form
$$
x = 2^{x_2} 5^{x_5} u v^2 =2^{x_2} u v^2=Uv^2 , \quad 
y=2^{y_2} 5^{y_5} u^3 v = (2^{x_2})^3 u^3 v=U^3v, \quad \text{where} \quad U=2^{x_2} u,
$$ 
and $u,v$ are some non-zero integers. 
Substituting this into equation \eqref{x4p5xyp2y3} and cancelling $U^4 v^3$, we obtain
 $$
2U^5+v^5=-5.
$$
This is a Thue equation and it has no integer solutions.

Let us now consider the case where both (2i) and (5ii) hold, so, $y_2=3x_2$ and $3x_5=y_5+1$ or $(x_5,y_5)=(t+1,3t+2)$ for some integer $t \geq 0$. Then substitution \eqref{ax4pbxypcy3lmnk} takes the form
$$
\begin{aligned}
x = & \,\, 2^{x_2} 5^{x_5} u v^2 =2^{x_2} (5)5^{t} u v^2=5Uv^2 , \\
y=& \,\, 2^{y_2} 5^{y_5} u^3 v = (2^{x_2})^3 (5^2)(5^{t})^3 u^3 v=25U^3v, \quad \text{where} \quad U=2^{x_2} 5^t u,
\end{aligned}
$$ 
and $u,v$ are some non-zero integers. 
Substituting this into equation \eqref{x4p5xyp2y3} and cancelling $625 U^4 v^3$, we obtain
 $$
50U^5+v^5=-1.
$$
This is a Thue equation whose only integer solution is $(U,v)=(0,-1)$, which is impossible as $U \neq 0$ by definition.

Let us now consider the case where both (2i) and (5iii) hold, so, $y_2=3x_2$ and $x_5+1= 2y_5$ or $(x_5,y_5)=(2t+1,t+1)$ for some integer $t \geq 0$. Then substitution \eqref{ax4pbxypcy3lmnk} takes the form
$$
\begin{aligned}
x = & \,\, 2^{x_2} 5^{x_5} u v^2 =2^{x_2} (5)(5^{t})^2 u v^2=5UV^2 , \\
y=& \,\, 2^{y_2} 5^{y_5} u^3 v = (2^{x_2})^3 (5)(5^{t}) u^3 v=5U^3V, \quad \text{where} \quad U=2^{x_2} u, \quad \text{and} \quad V=5^t v,
\end{aligned}
$$ 
and $u,v$ are non-zero integers.
Substituting this into equation \eqref{x4p5xyp2y3} and cancelling $125U^4 V^3$, we obtain
 $$
2U^5+5V^5=-1.
$$
This is a Thue equation and it has no integer solutions.

Let us now consider the case where both (2ii) and (5i) hold, so, $x_2= 2y_2+1$ and $x_5=y_5=0$. Then substitution \eqref{ax4pbxypcy3lmnk} takes the form
$$
x = 2^{x_2} 5^{x_5} u v^2 =2(2^{y_2})^2 u v^2=2uV^2 , \quad 
y=2^{y_2} 5^{y_5} u^3 v = 2^{y_2} u^3 v=u^3V, \quad \text{where} \quad V=2^{y_2} v,
$$ 
and $u,v$ are non-zero integers.
Substituting this into equation \eqref{x4p5xyp2y3} and cancelling $2u^4 V^3$, we obtain
 $$
u^5+8V^5=-5.
$$
This is a Thue equation and it has no integer solutions.

Let us now consider the case where both (2ii) and (5ii) hold, so, $x_2= 2y_2+1$ and $(x_5,y_5)=(t+1,3t+2)$ for some integer $t \geq 0$. Then substitution \eqref{ax4pbxypcy3lmnk} takes the form
$$
\begin{aligned}
x = & \,\, 2^{x_2} 5^{x_5} u v^2 =2(2^{y_2})^2 (5)5^{t} u v^2=10UV^2 , \\
y=& \,\, 2^{y_2} 5^{y_5} u^3 v = 2^{y_2} (5^2)(5^{t})^3 u^3 v=25U^3V, \quad \text{where} \quad U=5^t  u \quad \text{and} \quad V=2^{y_2} v,
\end{aligned}
$$ 
and $u,v$ are non-zero integers.
Substituting this into equation \eqref{x4p5xyp2y3} and cancelling $1250U^4 V^3$, we obtain
 $$
25U^5+8V^5=-1.
$$
This is a Thue equation and it has no integer solutions.

Let us now consider the case where both (2ii) and (5iii) hold, so, $x_2= 2y_2+1$ and $(x_5,y_5)=(2t+1,t+1)$ for some integer $t \geq 0$. Then substitution \eqref{ax4pbxypcy3lmnk} takes the form
$$
\begin{aligned}
x = & \,\, 2^{x_2} 5^{x_5} u v^2 =2(2^{y_2})^2 (5)(5^{t})^2 u v^2=10uV^2 , \\
y=& \,\, 2^{y_2} 5^{y_5} u^3 v = 2^{y_2} (5)5^{t} u^3 v=5u^3V,  \quad \text{where} \quad V=2^{y_2} 5^t v,
\end{aligned}
$$ 
and $u,v$ are non-zero integers.
Substituting this into equation \eqref{x4p5xyp2y3} and cancelling $250u^4 V^3$, we obtain
 $$
u^5+40V^5=-1.
$$
This is a Thue equation whose only integer solution is $(u,V)=(-1,0)$, which is impossible because $V \neq 0$ by definition.

To summarise, the unique integer solution to equation \eqref{x4p5xyp2y3} is
$$
(x,y)=(0,0).
$$

\vspace{10pt}
 
The next equation we will consider is
\begin{equation}\label{x4p3xy3y3}
x^4+3xy+3y^3=0,
\end{equation}
which is equation \eqref{ax4pbxypcy3} with $a=1$ and $b=c=3$.
We have an obvious solution $(x,y)=(0,0)$ and otherwise $x\neq 0$ and $y\neq 0$. The only prime divisor of $abc$ is $p=3$. Let $x_3$ and $y_3$ be the exponents of $3$ in the prime factorisation of $x$ and $y$, respectively.  Then $3$ appears in the prime factorisation of the monomials $x^4$, $3xy$, and $3y^3$ with the exponents $4x_3$, $x_3+y_3+1$ and $3y_3+1$, respectively. Let $M_1=\min\{4x_3,x_3+y_3+1,3y_3+1 \}$, and $M_1'$ be the second-smallest among these integers. Then we must have $M_1=M_1'$. The case $x_3+y_3+1 \geq 4x_3 = 3y_3+1$ is clearly impossible for non-negative $x_3,y_3$, hence we must have either (i) $x_3+y_3+1=4x_3$ or (ii) $x_3+y_3+1 = 3y_3+1$. Let us consider these cases separately. In case (i),  $y_3+1=3x_3$, or $(x_3,y_3)=(w+1,3w+2)$ for some integer $w \geq 0$. Then substitution \eqref{ax4pbxypcy3lmnk} takes the form
$$
x = 3^{x_3} u v^2 =3(3^{w}) u v^2 = 3Uv^2, \quad 
y=3^{y_3} u^3 v =3^2(3^{w})^3 u^3 v =9U^3 v, \quad \text{where} \quad U=3^{w} u,
$$ 
and $u,v$ are non-zero integers.
Substituting this into equation \eqref{x4p3xy3y3} and cancelling $3^4 U^4 v^3$, we obtain
 $$
27U^5+v^5=-1.
$$
This is a Thue equation whose only integer solution is $(U,v)=(0,-1)$, which is impossible because $U \neq 0$ by definition. In case (ii), $x_3= 2y_3$. Then substitution \eqref{ax4pbxypcy3lmnk} takes the form
$$
x = 3^{x_3} u v^2 =(3^{y_3})^2 u v^2=uV^2 , \quad 
y=3^{y_3} u^3 v =3^{y_3} u^3 v =u^3 V, \quad \text{where} \quad V=3^{y_3} v,
$$
and $u,v$ are non-zero integers.
Substituting this into equation \eqref{x4p3xy3y3} and cancelling $u^4 V^3$, we obtain
 $$
 3u^5+V^5=-3.
$$
This is a Thue equation whose only integer solution is $(u,V)=(-1,0)$, which is impossible because $V \neq 0$ by definition. Therefore, the unique integer solution to equation \eqref{x4p3xy3y3} is
$$
(x,y)=(0,0).
$$

\vspace{10pt}
 
The next equation we will consider is
\begin{equation}\label{x4pxyp4y3}
x^4+xy+4y^3=0,
\end{equation}
which is equation \eqref{ax4pbxypcy3} with $a=b=1$ and $c=4$.
We have an obvious solution $(x,y)=(0,0)$ and otherwise $x\neq 0$ and $y\neq 0$. The only prime divisor of $abc$ is $p=2$. Let $x_2$ and $y_2$ be the exponents of $2$ in the prime factorisation of $x$ and $y$, respectively. Then $2$ appears in the prime factorisation of the monomials $x^4$, $xy$, and $4y^3$ with the exponents $4x_2$, $x_2+y_2$ and $3y_2+2$, respectively. Let $M_1=\min\{4x_2,x_2+y_2,3y_2+2 \}$, and $M_1'$ be the second-smallest among these integers. Then we must have $M_1=M_1'$. The case $x_2+y_2 \geq 4x_2 = 3y_2+2$ is clearly impossible for non-negative $x_2,y_2$, hence we must have either (i) $x_2+y_2 = 4x_2$ or (ii) $x_2+y_2= 3y_2+2$. Let us consider these cases separately. In case (i), $y_2=3x_2$. Then substitution \eqref{ax4pbxypcy3lmnk} takes the form
$$
x = 2^{x_2} u v^2 =2^{x_2} u v^2 = Uv^2, \quad 
y=2^{y_2} u^3 v =(2^{x_2})^3 u^3 v =U^3 v, \quad \text{where} \quad U=2^{x_2} u,
$$ 
and $u,v$ are non-zero integers.
Substituting this into equation \eqref{x4pxyp4y3} and cancelling $U^4 v^3$, we obtain
 $$
v^5+4U^5=-1.
$$
This is a Thue equation whose only integer solution is $(v,U)=(-1,0)$, which is impossible because $U \neq 0$ by definition. In case (ii), $x_2=2y_2+2$. Then substitution \eqref{ax4pbxypcy3lmnk} takes the form
$$
x = 2^{x_2} u v^2 =2^2(2^{y_2})^2 u v^2 = 4uV^2, \quad 
y=2^{y_2} u^3 v =2^{y_2} u^3 v =u^3 V, \quad \text{where} \quad V=2^{y_2} v,
$$ 
and $u,v$ are non-zero integers.
Substituting this into equation \eqref{x4pxyp4y3} and cancelling $4u^4 V^3$, we obtain
 $$
64V^5+u^5=-1.
$$
This is a Thue equation whose only integer solution is $(V,u)=(0,-1)$, which is impossible because $V \neq 0$ by definition. Therefore, the unique integer solution to equation \eqref{x4pxyp4y3} is
$$
(x,y)=(0,0).
$$

\vspace{10pt}
 
The next equation we will consider is
\begin{equation}\label{2x4pxymy4}
2x^4+xy-y^4=0,
\end{equation}
which is equation \eqref{ax4pbxypcy4} with $a=2$, $b=1$ and $c=-1$.
We have an obvious solution $(x,y)=(0,0)$ and otherwise $x\neq 0$ and $y\neq 0$. The only prime divisor of $abc$ is $p=2$. Let $x_2$ and $y_2$ be the exponents of $2$ in the prime factorisation of $x$ and $y$, respectively.  Then $2$ appears in the prime factorisation of the monomials $2x^4$, $xy$, and $y^4$ with the exponents $4x_2+1$, $x_2+y_2$ and $4y_2$, respectively. Let $M_1=\min\{4x_2+1,x_2+y_2,4y_2 \}$, and $M_1'$ be the second-smallest among these integers. Then we must have $M_1=M_1'$. The case $x_2+y_2 \geq 4x_2 +1= 4y_2$ is clearly impossible for non-negative $x_2,y_2$, hence we must have either (i) $x_2+y_2 = 4x_2+1$ or (ii) $x_2+y_2= 4y_2$. Let us consider these cases separately. In case (i), $y_2=3x_2+1$. Then substitution \eqref{ax4pbxypcy4lmnk} takes the form
$$
\begin{aligned}
x=(-1)^{e_x} 2^{x_2} u v^3 =(-1)^{e_x} U v^3, \quad
y=(-1)^{e_y} 2^{y_2} u^3 v =(-1)^{e_y} 2(2^{x_2})^3 u^3 v =(-1)^{e_y} 2U^3 v  \\ \text{where} \quad U=2^{x_2} u,
\end{aligned}
$$
and $u,v$ are non-zero integers.
Substituting this into equation \eqref{2x4pxymy4} and cancelling $2U^4 v^4$, we obtain
$$
v^8  -8U^8 =(-1)^{e_x+e_y+1}.
$$
This is a Thue equation whose integer solutions are $(v,U,e_x,e_y)=(\pm 1,0,0,1),(\pm 1,0,1,0)$, which are impossible because $U \neq 0$ by definition.  In case (ii), $x_2=3y_2$. Then substitution \eqref{ax4pbxypcy4lmnk} takes the form
$$
\begin{aligned}
x=(-1)^{e_x} 2^{x_2} u v^3 =(-1)^{e_y} (2^{y_2})^3 u v^3 =(-1)^{e_x} u V^3, \quad
y=(-1)^{e_y} 2^{y_2} u^3 v  =(-1)^{e_y} u^3 V \\ \text{where} \quad V=2^{y_2} v,
\end{aligned}
$$
and $u,v$ are non-zero integers.
Substituting this into equation \eqref{2x4pxymy4} and cancelling $u^4 V^4$, we obtain
$$
2 V^8-u^8=(-1)^{e_x+e_y+1}.
$$
This is a Thue equation whose integer solutions are $(u,V,e_x,e_y)=(\pm 1,\pm 1,0,1),(\pm 1,\pm 1,1,0),(\pm1,0,0,0)$. The final solution is impossible as $V \neq 0$ by definition. The first two solutions correspond the solutions $(x,y)=((-1)^{e_x} u V^3,(-1)^{e_y} u^3 V)=\pm(1,-1)$. Therefore, the integer solutions to equation \eqref{2x4pxymy4} are
$$
(x,y)=(0,0) \quad \text{and} \quad \pm(1,-1).
$$

\vspace{10pt}

The next equation we will consider is
\begin{equation}\label{2x4p3xypy3}
2x^4+3xy+y^3=0,
\end{equation}
which is equation \eqref{ax4pbxypcy3} with $a=2$, $b=3$ and $c=1$.
We have an obvious solution $(x,y)=(0,0)$ and otherwise $x\neq 0$ and $y\neq 0$. The prime divisors of $abc$ are $p=2$ or $p=3$. Let $x_p$ and $y_p$ be the exponents of $p$ in the prime factorisation of $x$ and $y$, respectively. Then $2$ appears in the prime factorisation of the monomials $2x^4$, $3xy$, and $y^3$ with the exponents $4x_2+1$, $x_2+y_2$ and $3y_2$, respectively. Let $M_1=\min\{4x_2+1,x_2+y_2,3y_2 \}$, and $M_1'$ be the second-smallest among these integers. Then we must have $M_1=M_1'$. The case $x_2+y_2 \geq 4x_2+1 = 3y_2$ impossible for non-negative $x_2,y_2$, hence we must have either (2i) $x_2+y_2 = 4x_2+1$ or (2ii) $x_2+y_2= 3y_2$.
Also $3$ appears in the prime factorisation of the monomials $2x^4$, $3xy$, and $y^3$ with the exponents $4x_3$, $x_3+y_3+1$ and $3y_3$, respectively. Let $M_2=\min\{4x_3,x_3+y_3+1,3y_3 \}$, and $M_2'$ be the second-smallest among these integers. Then we must have $M_2=M_2'$. The case $x_3+y_3+1 \geq 4x_3 = 3y_3$ is only possible for $(x_3,y_3)=(0,0)$, hence we must have either (3i) $x_3=y_3=0$, or (3ii) $4x_3=x_3+y_3+1$, or (3iii) $x_3+y_3+1= 3y_3$.
We now have $6$ cases to consider. We will first consider the case where both (2i) and (3i) hold, so, $y_2= 3x_2+1$ and $x_3=y_3=0$. Then substitution \eqref{ax4pbxypcy3lmnk} takes the form
$$
x = 2^{x_2} 3^{x_3} u v^2 =2^{x_2} u v^2=Uv^2 , \quad 
y=2^{y_2} 3^{y_3} u^3 v = 2(2^{x_2})^3 u^3 v=2U^3v, \quad \text{where} \quad U=2^{x_2} u,
$$ 
and $u,v$ are non-zero integers.
Substituting this into equation \eqref{2x4p3xypy3} and cancelling $2U^4 v^3$, we obtain
 $$
4U^5+v^5=-3.
$$
This is a Thue equation whose only integer solution is $(U,v)=(-1,1)$, from which we obtain the integer solution $(x,y)=(Uv^2,2U^3v)=(-1,-2)$ to \eqref{2x4p3xypy3}.

Let us now consider the case where both (2i) and (3ii) hold, so $y_2= 3x_2+1$, and $3x_3=y_3+1$ or $(x_3,y_3)=(t+1,3t+2)$ for some integer $t \geq 0$. Then substitution \eqref{ax4pbxypcy3lmnk} takes the form
$$
\begin{aligned}
x = & \,\, 2^{x_2} 3^{x_3} u v^2 =2^{x_2}(3)3^{t} u v^2=3Uv^2 , \\
y=&  \,\, 2^{y_2} 3^{y_3} u^3 v = 2(2^{x_2})^3(3^2)(3^t)^3 u^3 v=18U^3v, \quad \text{where} \quad U=2^{x_2}3^t u,
\end{aligned}
$$ 
and $u,v$ are non-zero integers.
Substituting this into equation \eqref{2x4p3xypy3} and cancelling $162 U^4 v^3$, we obtain
 $$
36U^5+v^5=-1.
$$
This is a Thue equation whose only integer solution is $(U,v)=(0,-1)$, which is impossible as $U \neq 0$ by definition.

Let us now consider the case where both (2i) and (3iii) hold, so $y_2= 3x_2+1$, and $x_3+1= 2y_3$ or $(x_3,y_3)=(2t+1,t+1)$ for some integer $t \geq 0$. Then substitution \eqref{ax4pbxypcy3lmnk} takes the form
$$
\begin{aligned}
x = & \,\, 2^{x_2} 3^{x_3} u v^2 =2^{x_2}(3)(3^{t})^2 u v^2=3UV^2 , \\
y=&  \,\, 2^{y_2} 3^{y_3} u^3 v = 2(2^{x_2})^3(3)3^t u^3 v=6U^3V, \quad \text{where} \quad U=2^{x_2} u \quad \text{and} \quad V= 3^t v,
\end{aligned}
$$ 
and $u,v$ are non-zero integers.
Substituting this into equation \eqref{2x4p3xypy3} and cancelling $54 U^4 V^3$, we obtain
 $$
4U^5+3v^5=-1.
$$
This is a Thue equation whose only integer solution is $(U,V)=(-1,1)$, from which we obtain the integer solution  $(x,y)=(3UV^2,6U^3V)=(-3,-6)$ to \eqref{2x4p3xypy3}.

 We will now consider the case where both (2ii) and (3i) hold, so, $x_2= 2y_2$ and $x_3=y_3=0$. Then substitution \eqref{ax4pbxypcy3lmnk} takes the form
$$
x = 2^{x_2} 3^{x_3} u v^2 =(2^{y_2})^2 u v^2=uV^2 , \quad 
y=2^{y_2} 3^{y_3} u^3 v = 2^{y_2} u^3 v=u^3V, \quad \text{where} \quad V=2^{y_2} v,
$$ 
and $u,v$ are non-zero integers.
Substituting this into equation \eqref{2x4p3xypy3} and cancelling $u^4 V^3$, we obtain
 $$
u^5+2V^5=-3.
$$
This is a Thue equation whose only integer solution is $(u,V)=(-1,-1)$, from which we obtain the integer solution $(x,y)=(uV^2,u^3V)=(-1,1)$ to \eqref{2x4p3xypy3}.

We will now consider the case where both (2ii) and (3ii) hold, so, $x_2= 2y_2$ and $(x_3,y_3)=(t+1,3t+2)$ for some integer $t \geq 0$. Then substitution \eqref{ax4pbxypcy3lmnk} takes the form
$$
\begin{aligned}
x = & \,\, 2^{x_2} 3^{x_3} u v^2 =(2^{y_2})^2 (3)3^{t} u v^2=3UV^2 , \\
y=& \,\, 2^{y_2} 3^{y_3} u^3 v = 2^{y_2} 3^2 (3^t)^3 u^3 v=9U^3V, \quad \text{where} \quad U=3^t u \quad \text{and} \quad V= 2^{y_2}v,
\end{aligned}
$$ 
and $u,v$ are non-zero integers.
Substituting this into equation \eqref{2x4p3xypy3} and cancelling $81U^4 V^3$, we obtain
 $$
9U^5+2V^5=-1.
$$
This is a Thue equation and it has no integer solutions. 

 We will now consider the case where both (2ii) and (3iii) hold, so, $x_2= 2y_2$ and $(x_3,y_3)=(2t+1,t+1)$ for some integer $t \geq 0$. Then substitution \eqref{ax4pbxypcy3lmnk} takes the form
$$
\begin{aligned}
x = & \,\, 2^{x_2} 3^{x_3} u v^2 =(2^{y_2})^2(3)(3^t)^2 u v^2=3uV^2 , \\
y= & \,\, 2^{y_2} 3^{y_3} u^3 v = 2^{y_2} (3) 3^t u^3 v=3u^3V, \quad \text{where} \quad V=2^{y_2} 3^t v,
\end{aligned}
$$ 
and $u,v$ are non-zero integers.
Substituting this into equation \eqref{2x4p3xypy3} and cancelling $27u^4 V^3$, we obtain
 $$
u^5+6V^5=-1.
$$
This is a Thue equation and its only integer solution is $(u,V)=(-1,0)$, which is impossible as $V \neq 0$ by definition.

To summarise, the integer solutions to equation \eqref{2x4p3xypy3} are
$$
(x,y)=(-3,-6), \quad (-1,-2), \quad (-1,1), \quad \text{and} \quad (0,0).
$$

\vspace{10pt}
 
The next equation we will consider is
\begin{equation}\label{2x4pxyp2y3}
2x^4+xy+2y^3=0,
\end{equation}
which is equation \eqref{ax4pbxypcy3} with $a=c=2$ and $b=1$.
We have an obvious solution $(x,y)=(0,0)$ and otherwise $x\neq 0$ and $y\neq 0$. The only prime divisor of $abc$ is $p=2$. Let $x_2$ and $y_2$ be the exponents of $2$ in the prime factorisation of $x$ and $y$, respectively. Then $2$ appears in the prime factorisation of the monomials $2x^4$, $xy$, and $2y^3$ with the exponents $4x_2+1$, $x_2+y_2$ and $3y_2+1$, respectively. Let $M_1=\min\{4x_2+1,x_2+y_2,3y_2+1 \}$, and $M_1'$ be the second-smallest among these integers. Then we must have $M_1=M_1'$. The case $x_2+y_2 \geq 4x_2+1 = 3y_2+1$ is clearly impossible for non-negative $x_2,y_2$, hence we must have either (i) $x_2+y_2 = 4x_2+1$ or (ii) $x_2+y_2= 3y_2+1$. Let us consider these cases separately. In case (i), $y_2=3x_2+1$, then substitution \eqref{ax4pbxypcy3lmnk} takes the form
$$
x = 2^{x_2} u v^2 =2^{x_2} u v^2 = Uv^2, \quad 
y=2^{y_2} u^3 v =2(2^{x_2})^3 u^3 v =2U^3 v, \quad \text{where} \quad U=2^{x_2} u,
$$ 
and $u,v$ are non-zero integers.
Substituting this into equation \eqref{2x4pxyp2y3} and cancelling $2U^4 v^3$, we obtain
 $$
v^5+8U^5 =-1.
$$
This is a Thue equation whose only integer solution is $(v,U)=(-1,0)$, which is impossible because $U \neq 0$ by definition. In case (ii), $x_2=2y_2+1$, then substitution \eqref{ax4pbxypcy3lmnk} takes the form
$$
x = 2^{x_2} u v^2 =2(2^{y_2})^2 u v^2 = 2uV^2, \quad 
y=2^{y_2} u^3 v =2^{y_2} u^3 v =u^3 V, \quad \text{where} \quad V=2^{y_2} v,
$$ 
and $u,v$ are non-zero integers.
Substituting this into equation \eqref{2x4pxyp2y3} and cancelling $2u^4 V^3$, we obtain
 $$
16V^5+u^5=-1.
$$
This is a Thue equation whose only integer solution is $(V,u)=(0,-1)$, which is impossible because $V \neq 0$ by definition. Therefore, the unique integer solution to equation \eqref{2x4pxyp2y3} is
$$
(x,y)=(0,0).
$$

\vspace{10pt}

The next equation we will consider is
\begin{equation}\label{x5p3xypy3}
x^5+3xy+y^3=0,
\end{equation}
which is equation \eqref{ax5pbxypcy3} with $a=c=1$ and $b=3$.
We have an obvious solution $(x,y)=(0,0)$ and otherwise $x\neq 0$ and $y\neq 0$. The only prime divisor of $abc$ is $p=3$. Let $x_3$ and $y_3$ be the exponents of $3$ in the prime factorisation of $x$ and $y$, respectively. Then $3$ appears in the prime factorisation of the monomials $x^5$, $3xy$, and $y^3$ with the exponents $5x_3$, $x_3+y_3+1$ and $3y_3$, respectively. Let $M_1=\min\{5x_3,x_3+y_3+1,3y_3 \}$, and $M_1'$ be the second-smallest among these integers. Then we must have $M_1=M_1'$. The case $x_3+y_3+1 \geq 5x_3 = 3y_3$ only holds for $(x_3,y_3)=(0,0)$ hence we must have either (i) $x_3=y_3=0$ or (ii) $x_3+y_3+1 = 5x_3$ or (iii) $x_3+y_3+1 = 3y_3$. Let us consider these cases separately. In case (i),  $x_3=y_3=0$. Then substitution \eqref{ax5pbxypcy3lmnk} takes the form
$$
x = 3^{x_3} u v^2 = u v^2  \quad 
y=3^{y_3} u^4 v = u^4 v,
$$
where $u,v$ are non-zero integers.
Substituting this into equation \eqref{x5p3xypy3} and cancelling $u^5 v^3$, we obtain
 $$
v^7+u^7=-3.
$$
This is a Thue equation and it has no integer solutions. In case (ii), $y_3+1=4x_3$, or $(x_3,y_3)=(w+1,4w+3)$ for some integer $w \geq 0$. Then substitution \eqref{ax5pbxypcy3lmnk} takes the form
$$
x = 3^{x_3} u v^2 =3(3^{w}) u v^2 = 3Uv^2, \quad 
y=3^{y_3} u^4 v =3^3(3^{w})^4 u^4 v =27U^4 v, \quad \text{where} \quad U=3^{w} u,
$$
and $u,v$ are non-zero integers.
Substituting this into equation \eqref{x5p3xypy3} and cancelling $3^5U^5 v^3$, we obtain
 $$
v^7+81U^7=-1.
$$
This is a Thue equation, whose only integer solution is $(U,v)=(0,-1)$, which is impossible because $U \neq 0$ by definition. In case (iii), $x_3+1 = 2y_3$, or $(x_3,y_3)=(2w+1,w+1)$ for some integer $w \geq 0$. Then substitution \eqref{ax5pbxypcy3lmnk} takes the form
$$
x = 3^{x_3} u v^2 =3(3^{w})^2 u v^2 = 3uV^2, \quad 
y=3^{y_3} u^4 v =3(3^{w}) u^4 v =3u^4 V, \quad \text{where} \quad V=3^{w} v,
$$
and $u,v$ are non-zero integers.
Substituting this into equation \eqref{x5p3xypy3} and cancelling $27u^5 V^3$, we obtain
 $$
9V^7+u^7=-1.
$$ 
This is a Thue equation whose only integer solution is $(V,u)=(0,-1)$, which is impossible because $V \neq 0$ by definition. Therefore, the unique integer solution to equation \eqref{x5p3xypy3} is
$$
(x,y)= (0,0).
$$ 

\vspace{10pt}

The next equation we will consider is
\begin{equation}\label{x5pxyp2y3}
x^5+xy+2y^3=0,
\end{equation}
which is equation \eqref{ax5pbxypcy3} with $a=b=1$ and $c=2$.
We have an obvious solution $(x,y)=(0,0)$ and otherwise $x\neq 0$ and $y\neq 0$. The only prime divisor of $abc$ is $p=2$. Let $x_2$ and $y_2$ be the exponents of $2$ in the prime factorisation of $x$ and $y$, respectively.  Then $2$ appears in the prime factorisation of the monomials $x^5$, $xy$, and $2y^3$ with the exponents $5x_2$, $x_2+y_2$ and $3y_2+1$, respectively. Let $M_1=\min\{5x_2,x_2+y_2,3y_2+1 \}$, and $M_1'$ be the second-smallest among these integers. Then we must have $M_1=M_1'$. The case $x_2+y_2 \geq 5x_2 = 3y_2+1$ is clearly impossible for non-negative $x_2,y_2$, hence we must have either (i) $x_2+y_2 = 5x_2$ or (iii) $x_2+y_2 = 3y_2+1$. Let us consider these cases separately. In case (i), $y_2=4x_2$. Then substitution \eqref{ax5pbxypcy3lmnk} takes the form
$$
x = 2^{x_2} u v^2 =2^{x_2} u v^2 = Uv^2, \quad 
y=2^{y_2} u^4 v =(2^{x_2})^4 u^4 v =U^4 v, \quad \text{where} \quad U=2^{x_2} u,
$$
and $u,v$ are non-zero integers.
Substituting this into equation \eqref{x5pxyp2y3} and cancelling $U^5 v^3$, we obtain
 $$
v^7+2U^7=-1.
$$
This is a Thue equation whose integer solutions are $(U,v)=(0,-1),(-1,1)$. The first solution is impossible because $U \neq 0$ by definition. The second solution gives the integer solution $(x,y)=(Uv^2,U^4v)=(-1,1)$ to \eqref{x5pxyp2y3}. In case (ii), $x_2 = 2y_2+1$. Then substitution \eqref{ax5pbxypcy3lmnk} takes the form
$$ 
x = 2^{x_2} u v^2 =2(2^{y_2})^2 u v^2 = 2uV^2, \quad 
y=2^{y_2} u^4 v =(2^{y_2}) u^4 v =u^4 V, \quad \text{where} \quad V=2^{y_2} v,
$$
and $u,v$ are non-zero integers.
Substituting this into equation \eqref{x5pxyp2y3} and cancelling $8u^5 V^3$, we obtain
 $$
16V^7+u^7=-1.
$$ 
This is a Thue equation whose only integer solution is $(V,u)=(0,-1)$, which is impossible because $V \neq 0$ by definition. Therefore, the integer solutions to equation \eqref{x5pxyp2y3} are
$$
(x,y)=(-1,1) \quad \text{and} \quad (0,0).
$$

\vspace{10pt}

The next equation we will consider is
\begin{equation}\label{x5pxypy4}
x^5+xy+y^4=0,
\end{equation}
which is equation \eqref{ax5pbxypcy4} with $a=b=c=1$.
We have an obvious solution $(x,y)=(0,0)$ and otherwise $x \neq 0$ and $y \neq 0$. There are no prime divisors of $abc$, so $X_i=1$ and $Y_i=1$. Therefore substitution \eqref{ax5pbxypcy4lmnk} takes the form $x=uv^3$ and $y=u^4v$, where $u,v$ are non-zero integers. Substituting these into equation \eqref{x5pxypy4} and cancelling $u^5 v^4$, we obtain
$$
u^{11}+v^{11}=-1.
$$
This is a Thue equation whose only integer solutions are $(v,u)=(-1,0),(0, -1)$, which are impossible as $u$ and $v$ are non-zero by definition. Therefore, the unique integer solution to equation \eqref{x5pxypy4} is 
$$
(x,y)=(0,0).
$$

\vspace{10pt}

The next equation we will consider is
\begin{equation}\label{x4p6xyp2y3}
x^4+6xy+2y^3=0,
\end{equation}
which is equation \eqref{ax4pbxypcy3} with $a=1$, $b=6$ and $c=2$.
We have an obvious solution $(x,y)=(0,0)$ and otherwise $x\neq 0$ and $y\neq 0$. The prime divisors of $abc$ are $p=2$ or $p=3$. Let $x_p$ and $y_p$ be the exponents of $p$ in the prime factorisation of $x$ and $y$, respectively.  Then $2$ appears in the prime factorisation of the monomials $x^4$, $6xy$, and $2y^3$ with the exponents $4x_2$, $x_2+y_2+1$ and $3y_2+1$, respectively. Let $M_1=\min\{4x_2,x_2+y_2+1,3y_2+1 \}$, and $M_1'$ be the second-smallest among these integers. Then we must have $M_1=M_1'$. The case $x_2+y_2+1 \geq 4x_2 = 3y_2+1$ impossible for non-negative $x_2,y_2$, hence we must have either (2i) $x_2+y_2 +1= 4x_2$ or (2iii) $x_2+y_2+1= 3y_2+1$.
Also $3$ appears in the prime factorisation of the monomials $x^4$, $6xy$, and $2y^3$ with the exponents $4x_3$, $x_3+y_3+1$ and $3y_3$, respectively. Let $M_2=\min\{4x_3,x_3+y_3+1,3y_3 \}$, and $M_2'$ be the second-smallest among these integers. Then we must have $M_2=M_2'$. The case $x_3+y_3+1 \geq 4x_3 = 3y_3$ is only possible for $(x_3,y_3)=(0,0)$, hence we must have either (3i) $x_3=y_3=0$, or (3ii) $4x_3=x_3+y_3+1$, or (3iii) $x_3+y_3+1= 3y_3$.
We now have $6$ cases to consider. We will first consider the case where both (2i) and (3i) hold, so, $x_2+y_2 +1= 4x_2$ or $(x_2,y_2)=(w+1,3w+2)$ for some integer $w \geq 0$ and $x_3=y_3=0$. Then substitution \eqref{ax4pbxypcy3lmnk} takes the form
$$
x = 2^{x_2} 3^{x_3} u v^2 =2(2^{w}) u v^2=2Uv^2 , \quad 
y=2^{y_2} 3^{y_3} u^3 v = 2^2(2^{w})^3 u^3 v=4U^3v, \quad \text{where} \quad U=2^{w} u,
$$ 
and $u,v$ are non-zero integers.
Substituting this into equation \eqref{x4p6xyp2y3} and cancelling $16U^4 v^3$, we obtain
 $$
8U^5+v^5=-3.
$$
This is a Thue equation and it has no integer solutions.

Let us now consider the case where both (2i) and (3ii) hold, so $(x_2,y_2)=(w+1,3w+2)$ for some integer $w \geq 0$, and $3x_3=y_3+1$ or $(x_3,y_3)=(t+1,3t+2)$ for some integer $t \geq 0$. Then substitution \eqref{ax4pbxypcy3lmnk} takes the form
$$
\begin{aligned}
x = & \,\, 2^{x_2} 3^{x_3} u v^2 =2(2^{w})(3)3^{t} u v^2=6Uv^2 , \\
y=&  \,\, 2^{y_2} 3^{y_3} u^3 v = 2^2(2^{w})^3(3^2)(3^t)^3 u^3 v=36U^3v, \quad \text{where} \quad U=2^{w}3^t u,
\end{aligned}
$$ 
and $u,v$ are non-zero integers.
Substituting this into equation \eqref{x4p6xyp2y3} and cancelling $1296 U^4 v^3$, we obtain
 $$
72U^5+v^5=-1.
$$
This is a Thue equation whose only integer solution is $(U,v)=(0,-1)$, which is impossible as $U \neq 0$ by definition.

We will now consider the case where both (2i) and (3iii) hold, so, $(x_2,y_2)=(w+1,3w+2)$ for some integer $w \geq 0$ and  $x_3+1= 2y_3$ or $(x_3,y_3)=(2t+1,t+1)$ for some integer $t \geq 0$. Then substitution \eqref{ax4pbxypcy3lmnk} takes the form
$$
\begin{aligned}
x = & \,\, 2^{x_2} 3^{x_3} u v^2 =2(2^w)(3) (3^{t})^2 u v^2=6U V^2, \\
y=& \,\, 2^{y_2} 3^{y_3} u^3 v =2^2(2^w)^3 (3)3^{t} u^3 v=12U^3 V, \quad \text{where} \quad U=2^w u \quad \text{and} \quad V=3^t  v,
\end{aligned}
$$
and $u,v$ are non-zero integers. 
Substituting this into equation \eqref{x4p6xyp2y3} and cancelling $432U^4 V^3$, we obtain
 $$
8U^5+3v^5=-1.
$$
This is a Thue equation and it has no integer solutions. 

We will now consider the case where both (2ii) and (3i) hold, so, $x_2= 2y_2$ and $x_3=y_3=0$. Then substitution \eqref{ax4pbxypcy3lmnk} takes the form
$$
x = 2^{x_2} 3^{x_3} u v^2 =(2^{y_2})^2 u v^2=uV^2 , \quad 
y=2^{y_2} 3^{y_3} u^3 v = 2^{y_2} u^3 v=u^3V, \quad \text{where} \quad V=2^{y_2} v,
$$ 
and $u,v$ are non-zero integers.
Substituting this into equation \eqref{x4p6xyp2y3} and cancelling $u^4 V^3$, we obtain
 $$
2u^5+V^5=-6.
$$
This is a Thue equation and it has no integer solutions.

Let us now consider the case where both (2ii) and (3ii) hold, so $x_2= 2y_2$, and $(x_3,y_3)=(t+1,3t+2)$ for some integer $t \geq 0$. Then substitution \eqref{ax4pbxypcy3lmnk} takes the form
$$
\begin{aligned}
x =& \,\,  2^{x_2} 3^{x_3} u v^2 =(2^{y_2})^2(3)3^{t} u v^2=3UV^2 , \\ 
y=& \,\, 2^{y_2} 3^{y_3} u^3 v = 2^{y_2}(3^2)(3^t)^3 u^3 v=9U^3V, \quad \text{where} \quad U=3^t u \quad \text{and} \quad V=2^{y_2}v,
\end{aligned}
$$ 
and $u,v$ are non-zero integers.
Substituting this into equation \eqref{x4p6xyp2y3} and cancelling $81 U^4 V^3$, we obtain
 $$
18U^5+V^5=-2.
$$
This is a Thue equation and it has no integer solutions.

We will now consider the case where both (2ii) and (3iii) hold, so, $x_2= 2y_2$ and  $(x_3,y_3)=(2t+1,t+1)$ for some integer $t \geq 0$. Then substitution \eqref{ax4pbxypcy3lmnk} takes the form
$$
\begin{aligned}
x = & \,\, 2^{x_2} 3^{x_3} u v^2 =(2^{y_2})^2(3) (3^{t})^2 u v^2=3u V^2, \\
y=&\,\, 2^{y_2} 3^{y_3} u^3 v =2^{y_2} (3)3^{t} u^3 v=3u^3 V, \quad \text{where}  \quad V=2^{y_2} 3^t  v,
\end{aligned}
$$ 
and $u,v$ are non-zero integers.
Substituting this into equation \eqref{x4p6xyp2y3} and cancelling $27u^4 V^3$, we obtain
 $$
2u^5+3v^5=-2.
$$
This is a Thue equation whose only integer solution is $(u,V)=(-1,0)$, which is impossible as $V \neq 0$ by definition. 

To summarise, the unique integer solution to equation \eqref{x4p6xyp2y3} is
$$
(x,y)=(0,0).
$$

\vspace{10pt}

The next equation we will consider is
\begin{equation}\label{x4p4xyp3y3}
x^4+4xy+3y^3=0,
\end{equation}
which is equation \eqref{ax4pbxypcy3} with $a=1$, $b=4$ and $c=3$.
We have an obvious solution $(x,y)=(0,0)$ and otherwise $x\neq 0$ and $y\neq 0$. The prime divisors of $abc$ are $p=2$ or $p=3$. Let $x_p$ and $y_p$ be the exponents of $p$ in the prime factorisation of $x$ and $y$, respectively. Then $2$ appears in the prime factorisation of the monomials $x^4$, $4xy$, and $3y^3$ with the exponents $4x_2$, $x_2+y_2+2$ and $3y_2$, respectively. Let $M_1=\min\{4x_2,x_2+y_2+2,3y_2 \}$, and $M_1'$ be the second-smallest among these integers. Then we must have $M_1=M_1'$. The case $x_2+y_2+2 \geq 4x_2 = 3y_2$ is only possible for $(x_2,y_2)=(0,0)$, hence we must have either (2i) $x_2=y_2=0$, or (2ii) $x_2+y_2 +2= 4x_2$ or (2iii) $x_2+y_2+2= 3y_2$.
Also $3$ appears in the prime factorisation of the monomials $x^4$, $4xy$, and $3y^3$ with the exponents $4x_3$, $x_3+y_3$ and $3y_3+1$, respectively. Let $M_2=\min\{4x_3,x_3+y_3,3y_3+1 \}$, and $M_2'$ be the second-smallest among these integers. Then we must have $M_2=M_2'$. The case $x_3+y_3 \geq 4x_3 = 3y_3+1$ is impossible for non-negative $x_3,y_3$, hence we must have either (3i) $4x_3=x_3+y_3$, or (3ii) $x_3+y_3= 3y_3+1$.
We now have $6$ cases to consider. We will first consider the case where both (2i) and (3i) hold, so, $x_2=y_2=0$ and $3x_3=y_3$. Then substitution \eqref{ax4pbxypcy3lmnk} takes the form
$$
x = 2^{x_2} 3^{x_3} u v^2 =3^{x_3} u v^2=U v^2, \quad 
y=2^{y_2} 3^{y_3} u^3 v = (3^{x_3})^3 u^3 v=U^3 v, \quad \text{where} \quad U=3^{x_3} u,
$$ 
and $u,v$ are non-zero integers.
Substituting this into equation \eqref{x4p4xyp3y3} and cancelling $U^4 v^3$, we obtain
 $$
3U^5+v^5=-4.
$$
This is a Thue equation whose only integer solution is $(U,v)=(-1,-1)$, from which we obtain the integer solution $(x,y)=(Uv^2,U^3v)=(-1,1)$ to \eqref{x4p4xyp3y3}. 

Let us now consider the case where both (2i) and (3ii) hold, so $x_2=y_2=0$, and  $x_3= 2y_3+1$. Then substitution \eqref{ax4pbxypcy3lmnk} takes the form
$$
x = 2^{x_2} 3^{x_3} u v^2 = 3(3^{y_3})^2 u v^2=3uV^2, \quad 
y=2^{y_2} 3^{y_3} u^3 v  =3^{y_3} u^3 v =u^3V, \quad \text{where} \quad V=3^{y_3} v,
$$ 
and $u,v$ are non-zero integers.
Substituting this into equation \eqref{x4p4xyp3y3} and cancelling $3 u^4 V^3$, we obtain
 $$
27V^5+u^5=-4.
$$
This is a Thue equation and it has no integer solutions.

We will now consider the case where both (2ii) and (3i) hold, so, $y_2 +2= 3x_2$ or $(x_2,y_2)=(w+1,3w+1)$ for some integer $w \geq 0$ and $3x_3=y_3$. Then substitution \eqref{ax4pbxypcy3lmnk} takes the form
$$
x = 2^{x_2} 3^{x_3} u v^2 =2(2^w) 3^{x_3} u v^2=2U v^2, \quad 
y=2^{y_2} 3^{y_3} u^3 v =2(2^w)^3 (3^{x_3})^3 u^3 v=2U^3 v, \quad \text{where} \quad U=2^w 3^{x_3} u,
$$ 
and $u,v$ are non-zero integers.
Substituting this into equation \eqref{x4p4xyp3y3} and cancelling $8U^4 v^3$, we obtain
 $$
3U^5+2v^5=-2.
$$
This is a Thue equation whose only integer solution is $(U,v)=(0,-1)$, which is impossible because $U \neq 0$ by definition. 

Let us now consider the case where both (2ii) and (3ii) hold, so $(x_2,y_2)=(w+1,3w+1)$ for some integer $w \geq 0$, and  $x_3= 2y_3+1$. Then substitution \eqref{ax4pbxypcy3lmnk} takes the form
$$
\begin{aligned}
x = &\,\,2^{x_2} 3^{x_3} u v^2 = 2(2^w) (3)(3^{y_3})^2 u v^2=6UV^2, \\
y=&\,\,2^{y_2} 3^{y_3} u^3 v  =2(2^w)^3 3^{y_3} u^3 v =2U^3V, \\ &\text{where} \quad V=3^{y_3} v, \quad \text{and} \quad U=2^w u,
\end{aligned}
$$ 
and $u,v$ are non-zero integers.
Substituting this into equation \eqref{x4p4xyp3y3} and cancelling $24 U^4 V^3$, we obtain
 $$
54V^5+U^5=-2.
$$
This is a Thue equation and it has no integer solutions.

We will now consider the case where both (2iii) and (3i) hold, so, $x_2+2= 2y_2$ or $(x_2,y_2)=(2w,w+1)$ for some integer $w \geq 0$ and $3x_3=y_3$. Then substitution \eqref{ax4pbxypcy3lmnk} takes the form
$$
\begin{aligned}
x = & \,\,2^{x_2} 3^{x_3} u v^2 =(2^w)^2 3^{x_3} u v^2=U V^2, \\ 
y=& \,\, 2^{y_2} 3^{y_3} u^3 v =2(2^w) (3^{x_3})^3 u^3 v=2U^3 V, \\ & \text{where} \quad U=3^{x_3} u \quad \text{and} \quad V=2^{w} v,
\end{aligned}
$$ 
and $u,v$ are non-zero integers.
Substituting this into equation \eqref{x4p4xyp3y3} and cancelling $U^4 V^3$, we obtain
 $$
24U^5+V^5=-8.
$$
This is a Thue equation whose only integer solution is $(U,V)=(1,-2)$, from which we obtain the integer solution $(x,y)=(UV^2,2U^3V)=(4,-4)$ to \eqref{x4p4xyp3y3}. 

Finally, let us now consider the case where both (2iii) and (3ii) hold, so $(x_2,y_2)=(2w,w+1)$ for some integer $w \geq 0$, and  $x_3= 2y_3+1$. Then substitution \eqref{ax4pbxypcy3lmnk} takes the form
$$
x = 2^{x_2} 3^{x_3} u v^2 = (2^w)^2 (3)(3^{y_3})^2 u v^2=3uV^2, \quad 
y=2^{y_2} 3^{y_3} u^3 v  =2(2^w) 3^{y_3} u^3 v =2u^3V, \quad \text{where} \quad V=2^w 3^{y_3} v,
$$ 
and $u,v$ are non-zero integers.
Substituting this into equation \eqref{x4p4xyp3y3} and cancelling $3 u^4 V^3$, we obtain
 $$
27V^5+8u^5=-8.
$$
This is a Thue equation whose only integer solution is $(V,u)=(0,-1)$, which is impossible because $V \neq 0$ by definition.

To summarise, the integer solutions to equation \eqref{x4p4xyp3y3} are
$$
(x,y)=(-1,1), \quad (0,0) \quad \text{and} \quad (4,-4).
$$

\vspace{10pt}

The next equation we will consider is
\begin{equation}\label{x4p2xyp4y3}
x^4+2xy+4y^3=0,
\end{equation}
which is equation \eqref{ax4pbxypcy3} with $a=1$, $b=2$ and $c=4$.
We have an obvious solution $(x,y)=(0,0)$ and otherwise $x\neq 0$ and $y\neq 0$. The only prime divisor of $abc$ is $p=2$. Let $x_2$ and $y_2$ be the exponents of $2$ in the prime factorisation of $x$ and $y$, respectively.  Then $2$ appears in the prime factorisation of the monomials $x^4$, $2xy$, and $4y^3$ with the exponents $4x_2$, $x_2+y_2+1$ and $3y_2+2$, respectively. Let $M_1=\min\{4x_2,x_2+y_2+1,3y_2+2 \}$, and $M_1'$ be the second-smallest among these integers. Then we must have $M_1=M_1'$. The case $x_2+y_2+1 \geq 4x_2 = 3y_2+2$ is clearly impossible for non-negative $x_2,y_2$, hence we must have either (i) $x_2+y_2 +1= 4x_2$ or (ii) $x_2+y_2+1= 3y_2+2$. Let us consider these cases separately. In case (i), $y_2+1=3x_2$, or $(x_2,y_2)=(w+1,3w+2)$ for some integer $w \geq 0$. Then substitution \eqref{ax4pbxypcy3lmnk} takes the form
$$
x = 2^{x_2} u v^2 =2(2^{w}) u v^2 = 2Uv^2, \quad 
y=2^{y_2} u^3 v =2^2(2^{w})^3 u^3 v =4U^3 v, \quad \text{where} \quad U=2^{w} u,
$$ 
and $u,v$ are non-zero integers.
Substituting this into equation \eqref{x4p2xyp4y3} and cancelling $2^4 U^4 v^3$, we obtain
 $$
v^5+16U^5=-1.
$$
This is a Thue equation whose only integer solution is $(v,U)=(-1,0)$, which is impossible because $U \neq 0$ by definition. In case (ii), $x_2=2y_2+1$. Then substitution \eqref{ax4pbxypcy3lmnk} takes the form
$$
x = 2^{x_2} u v^2 =2(2^{y_2})^2 u v^2 = 2uV^2, \quad 
y=2^{y_2} u^3 v =2^{y_2} u^3 v =u^3 V, \quad \text{where} \quad V=2^{y_2} v,
$$ 
and $u,v$ are non-zero integers.
Substituting this into equation \eqref{x4p2xyp4y3} and cancelling $4u^4 V^3$, we obtain
 $$
4V^5+u^5=-1.
$$
This is a Thue equation whose only integer solution is $(V,u)=(0,-1)$, which is impossible because $V \neq 0$ by definition. Therefore, the unique integer solution to equation \eqref{x4p2xyp4y3} is
$$
(x,y)=(0,0).
$$

\vspace{10pt}

The next equation we will consider is
\begin{equation}\label{2x4p2xymy4}
2x^4+2xy-y^4=0,
\end{equation}
which is equation \eqref{ax4pbxypcy4} with $a=b=2$, and $c=-1$.
We have an obvious solution $(x,y)=(0,0)$ and otherwise $x\neq 0$ and $y\neq 0$. The only prime divisor of $abc$ is $p=2$. Let $x_2$ and $y_2$ be the exponents of $2$ in the prime factorisation of $x$ and $y$, respectively.  Then $2$ appears in the prime factorisation of the monomials $2x^4$, $2xy$, and $y^4$ with the exponents $4x_2+1$, $x_2+y_2+1$ and $4y_2$, respectively. Let $M_1=\min\{4x_2+1,x_2+y_2+1,4y_2 \}$, and $M_1'$ be the second-smallest among these integers. Then we must have $M_1=M_1'$. The case $x_2+y_2+1 \geq 4x_2 +1= 4y_2$ is clearly impossible for non-negative $x_2,y_2$, hence we must have either (i) $x_2+y_2+1 = 4x_2+1$ or (ii) $x_2+y_2+1= 4y_2$. Let us consider these cases separately. In case (i), $y_2=3x_2$. Then substitution \eqref{ax4pbxypcy4lmnk} takes the form
$$
\begin{aligned}
x=(-1)^{e_x} 2^{x_2} u v^3 =(-1)^{e_x} U v^3, \quad
y=(-1)^{e_y} 2^{y_2} u^3 v =(-1)^{e_y} (2^{x_2})^3 u^3 v =(-1)^{e_y} U^3 v  \\ \text{where} \quad U=2^{x_2} u,
\end{aligned}
$$
and $u,v$ are non-zero integers.
Substituting this into equation \eqref{2x4p2xymy4} and cancelling $U^4 v^4$, we obtain
$$
2v^8  -U^8 =2(-1)^{e_x+e_y+1}.
$$
This is a Thue equation whose only integer solutions are $(v,U,e_x,e_y)=(\pm 1,0,0,1),(\pm 1,0,1,0)$, which are impossible because $U \neq 0$ by definition.  In case (ii), $x_2+1=3y_2$ or $(x_2,y_2)=(3w+2,w+1)$ for some integer $w \geq 0$. Then substitution \eqref{ax4pbxypcy4lmnk} takes the form
$$
\begin{aligned}
x=& (-1)^{e_x} 2^{x_2} u v^3 =(-1)^{e_y} 2^2 (2^{w})^3 u v^3 =(-1)^{e_x} 4u V^3, \\
y= &(-1)^{e_y} 2^{y_2} u^3 v  =(-1)^{e_y} 2(2^w) u^3 v=(-1)^{e_y} 2 u^3 V \quad \text{where} \quad V=2^{w} v,
\end{aligned}
$$
and $u,v$ are non-zero integers.
Substituting this into equation \eqref{2x4p2xymy4} and cancelling $16u^4 V^4$, we obtain
$$
32 V^8-u^8=(-1)^{e_x+e_y+1}.
$$
which is a Thue equation whose integer solutions are $(u,V,e_x,e_y)=(\pm1,0,0,0)$, which are impossible as $V \neq 0$ by definition. Therefore, the unique integer solution to equation \eqref{2x4p2xymy4} is
$$
(x,y)=(0,0).
$$  

\vspace{10pt}

The next equation we will consider is
\begin{equation}\label{2x4p4xypy3}
2x^4+4xy+y^3=0,
\end{equation}
which is equation \eqref{ax4pbxypcy3} with $a=2$, $b=4$ and $c=1$.
We have an obvious solution $(x,y)=(0,0)$ and otherwise $x\neq 0$ and $y\neq 0$. The only prime divisor of $abc$ is $p=2$. Let $x_2$ and $y_2$ be the exponents of $2$ in the prime factorisation of $x$ and $y$, respectively.  Then $2$ appears in the prime factorisation of the monomials $2x^4$, $4xy$, and $y^3$ with the exponents $4x_2+1$, $x_2+y_2+2$ and $3y_2$, respectively. Let $M_1=\min\{4x_2+1,x_2+y_2+2,3y_2 \}$, and $M_1'$ be the second-smallest among these integers. Then we must have $M_1=M_1'$. The case $x_2+y_2+2 \geq 4x_2+1 = 3y_2$ is only possible for $(x_2,y_2)=(0,0)$, hence we must have either (i) $x_2=y_2=0$, or (ii) $x_2+y_2 +2= 4x_2+1$ or (iii) $x_2+y_2+2= 3y_2$. Let us consider these cases separately. In case (i), $x_2=y_2=0$. Then substitution \eqref{ax4pbxypcy3lmnk} takes the form
$$
x = 2^{x_2} u v^2 =u v^2, \quad 
y=2^{y_2} u^3 v = u^3 v,
$$ 
where $u,v$ are non-zero integers.
Substituting this into equation \eqref{2x4p4xypy3} and cancelling $u^4 v^3$, we obtain
 $$
2v^5+u^5=-4.
$$
This is a Thue equation and it has no integer solutions. In case (ii), $y_2 +1= 3x_2$, or $(x_2,y_2)=(w+1,3w+2)$ for some integer $w \geq 0$. Then substitution \eqref{ax4pbxypcy3lmnk} takes the form
$$
x = 2^{x_2} u v^2 =2(2^{w}) u v^2 = 2Uv^2, \quad 
y=2^{y_2} u^3 v =2^2(2^{w})^3 u^3 v =4U^3 v, \quad \text{where} \quad U=2^{w} u,
$$ 
and $u,v$ are non-zero integers.
Substituting this into equation \eqref{2x4p4xypy3} and cancelling $2^5 U^4 v^3$, we obtain
 $$
v^5+2U^5=-1.
$$
This is a Thue equation whose integer solutions are $(v,U)=(-1,0),(1,-1)$. The first solution is impossible because $U \neq 0$ by definition. From the second solution we obtain the integer solution $(x,y)=(2Uv^2,4U^3 v)=(-2,-4)$ to \eqref{2x4p4xypy3}. In case (iii), $x_2+2=2y_2$, or $(x_2,y_2)=(2w,w+1)$ for some integer $w \geq 0$. Then substitution \eqref{ax4pbxypcy3lmnk} takes the form
$$
x = 2^{x_2} u v^2 =(2^{w})^2 u v^2 = uV^2, \quad 
y=2^{y_2} u^3 v =2(2^{w}) u^3 v =2u^3 V, \quad \text{where} \quad V=2^{w} v,
$$ 
and $u,v$ are non-zero integers.
Substituting this into equation \eqref{2x4p4xypy3} and cancelling $2u^4 V^3$, we obtain
 $$
V^5+4u^5=-4.
$$
This is a Thue equation whose only integer solution is $(V,u)=(-1,0)$, which is impossible because $u \neq 0$ by definition. Therefore, the integer solutions to equation \eqref{2x4p4xypy3} are
$$
(x,y)=(0,0) \quad \text{and} \quad (-2,-4).
$$

\vspace{10pt}

The next equation we will consider is
\begin{equation}\label{x5p4xypy3}
x^5+4xy+y^3=0,
\end{equation}
which is equation \eqref{ax5pbxypcy3} with $a=c=1$ and $b=4$.
We have an obvious solution $(x,y)=(0,0)$ and otherwise $x\neq 0$ and $y\neq 0$. The only prime divisor of $abc$ is $p=2$. Let $x_2$ and $y_2$ be the exponents of $2$ in the prime factorisation of $x$ and $y$, respectively.  Then $2$ appears in the prime factorisation of the monomials $x^5$, $4xy$, and $y^3$ with the exponents $5x_2$, $x_2+y_2+2$ and $3y_2$, respectively. Let $M_1=\min\{5x_2,x_2+y_2+2,3y_2 \}$, and $M_1'$ be the second-smallest among these integers. Then we must have $M_1=M_1'$. The case $x_2+y_2+2 \geq 5x_2 = 3y_2$ is clearly impossible for non-negative $x_2,y_2$, hence we must have either (i) $x_2+y_2+2 = 5x_2$ or (ii) $x_2+y_2+2 = 3y_2$. Let us consider these cases separately. In case (i), $y_2+2=4x_2$, or $(x_2,y_2)=(w+1,4w+2)$ for some integer $w \geq 0$. Then substitution \eqref{ax5pbxypcy3lmnk} takes the form
$$
x = 2^{x_2} u v^2 =2(2^{w}) u v^2 = 2Uv^2, \quad 
y=2^{y_2} u^4 v =2^2(2^{w})^4 u^4 v =4U^4 v, \quad \text{where} \quad U=2^{w} u,
$$
and $u,v$ are non-zero integers.
Substituting this into equation \eqref{x5p4xypy3} and cancelling $U^5 v^3$, we obtain
 $$
v^7 +2U^7 =-1.
$$
This is a Thue equation whose integer solutions are $(U,v)=(0,-1),(-1,1)$. The first solution is impossible as $U \neq 0$ by definition. From the second solution we obtain the integer solution $(x,y)=(2Uv^2,4U^4 v)=(-2,4)$ to \eqref{x5p4xypy3}. In case (ii), $x_2+2 = 2y_2$, or $(x_2,y_2)=(2w,w+1)$ for some integer $w \geq 0$. Then substitution \eqref{ax5pbxypcy3lmnk} takes the form
$$ 
x = 2^{x_2} u v^2 =(2^{w})^2 u v^2 = uV^2, \quad 
y=2^{y_2} u^4 v =2(2^{w}) u^4 v =2u^4 V, \quad \text{where} \quad V=2^{w} v,
$$
and $u,v$ are non-zero integers.
Substituting this into equation \eqref{x5p4xypy3} and cancelling $u^5 V^3$, we obtain
 $$
V^7+8u^7=-8.
$$ 
This is a Thue equation whose only integer solution is $(V,u)=(0,-1)$, which is impossible because $V \neq 0$ by definition. Therefore, the integer solutions to equation \eqref{x5p4xypy3} are
$$
(x,y)= (0,0) \quad \text{and} \quad (-2,4).
$$

\vspace{10pt}

The next equation we will consider is
\begin{equation}\label{x5p2xyp2y3}
x^5+2xy+2y^3=0,
\end{equation}
which is equation \eqref{ax5pbxypcy3} with $a=1$ and $b=c=2$.
We have an obvious solution $(x,y)=(0,0)$ and otherwise $x\neq 0$ and $y\neq 0$. The only prime divisor of $abc$ is $p=2$. Let $x_2$ and $y_2$ be the exponents of $2$ in the prime factorisation of $x$ and $y$, respectively.  Then $2$ appears in the prime factorisation of the monomials $x^5$, $2xy$, and $2y^3$ with the exponents $5x_2$, $x_2+y_2+1$ and $3y_2+1$, respectively. Let $M_1=\min\{5x_2,x_2+y_2+1,3y_2+1 \}$, and $M_1'$ be the second-smallest among these integers. Then we must have $M_1=M_1'$. The case $x_2+y_2+1 \geq 5x_2 = 3y_2+1$ is clearly impossible for non-negative $x_2,y_2$, hence we must have either (i) $x_2+y_2+1 = 5x_2$ or (ii) $x_2+y_2+1 = 3y_2+1$. Let us consider these cases separately. In case (i), $y_2+1=4x_2$, or $(x_2,y_2)=(w+1,4w+3)$ for some integer $w \geq 0$. Then substitution \eqref{ax5pbxypcy3lmnk} takes the form
$$
x = 2^{x_2} u v^2 =2(2^{w}) u v^2 = 2Uv^2, \quad 
y=2^{y_2} u^4 v =2^3(2^{w})^4 u^4 v =8U^4 v, \quad \text{where} \quad U=2^{w} u,
$$
and $u,v$ are non-zero integers.
Substituting this into equation \eqref{x5p2xyp2y3} and cancelling $2^5U^5 v^3$, we obtain
 $$
v^7+32 U^7 =-1.
$$
This is a Thue equation whose only integer solution is $(U,v)=(0,-1)$, which is impossible as $U \neq 0$ by definition.  In case (ii), $x_2= 2y_2$. Then substitution \eqref{ax5pbxypcy3lmnk} takes the form
$$ 
x = 2^{x_2} u v^2 =(2^{y_2})^2 u v^2 = uV^2, \quad 
y=2^{y_2} u^4 v = u^4 V, \quad \text{where} \quad V=2^{y_2} v,
$$
and $u,v$ are non-zero integers.
Substituting this into equation \eqref{x5p2xyp2y3} and cancelling $u^5 V^3$, we obtain
 $$
V^7+2u^7=-2.
$$ 
This is a Thue equation whose only integer solution is $(V,u)=(0,-1)$, which is impossible because $V \neq 0$ by definition. Therefore, the unique integer solution to equation \eqref{x5p2xyp2y3} is
$$
(x,y)= (0,0).
$$

\vspace{10pt}

The next equation we will consider is
\begin{equation}\label{x5p2xypy4}
x^5+2xy+y^4=0,
\end{equation}
which is equation \eqref{ax5pbxypcy4} with $a=c=1$ and $b=2$.
We have an obvious solution $(x,y)=(0,0)$ and otherwise $x\neq 0$ and $y\neq 0$. The only prime divisor of $abc$ is $p=2$. Let $x_2$ and $y_2$ be the exponents of $2$ in the prime factorisation of $x$ and $y$, respectively.  Then $2$ appears in the prime factorisation of the monomials $x^5$, $2xy$, and $y^4$ with the exponents $5x_2$, $x_2+y_2+1$ and $4y_2$, respectively. Let $M_1=\min\{5x_2,x_2+y_2+1,4y_2 \}$, and $M_1'$ be the second-smallest among these integers. Then we must have $M_1=M_1'$. The case $x_2+y_2+1 \geq 5x_2 = 4y_2$ is only possible for $(x_2,y_2)=(0,0)$, hence we must have either (i) $x_2=y_2=0$, or  (ii) $x_2+y_2+1 = 5x_2$, or (iii) $x_2+y_2+1=4y_2$. Let us consider these cases separately. In case (i), $x_2=y_2=0$. Then substitution \eqref{ax5pbxypcy4lmnk} takes the form
$$
x = 2^{x_2} u v^3 =u v^3, \quad 
y=2^{y_2} u^4 v = u^4 v,
$$ 
where $u,v$ are non-zero integers.
Substituting this into equation \eqref{x5p2xypy4} and cancelling $u^5 v^4$, we obtain
 $$
v^{11}+u^{11}=-2.
$$
This is a Thue equation whose only integer solution is $(u,v)=(-1,-1)$, from which we obtain the integer solution $(x,y)=(uv^3,u^4v)=(1,-1)$ to \eqref{x5p2xypy4}. 
In case (ii), $y_2+1 = 4x_2$, or $(x_2,y_2)=(w+1,4w+3)$ for some integer $w \geq 0$. Then substitution \eqref{ax5pbxypcy4lmnk} takes the form
$$
x = 2^{x_2} u v^3 =2(2^{w}) u v^3=2U v^3, \quad 
y=2^{y_2} u^4 v = 2^3(2^{w})^4 u^4 v=8U^4 v, \quad \text{where} \quad U=2^{w} u,
$$ 
and $u,v$ are non-zero integers.
Substituting this into equation \eqref{x5p2xypy4} and cancelling $2^5 U^5 v^4$, we obtain
 $$
v^{11}+128U^{11}=-1.
$$
This is a Thue equation whose only integer solution is $(v,U)=(-1,0)$, which is impossible as $U \neq 0$ by definition. In case (iii), $x_2+1=3y_2$, or $(x_2,y_2)=(3w+2,w+1)$ for some integer $w \geq 0$. Then substitution \eqref{ax5pbxypcy4lmnk} takes the form 
$$
x = 2^{x_2} u v^3 =2^2(2^{w})^3 u v^3=4u V^3, \quad 
y=2^{y_2} u^4 v = 2(2^{w}) u^4 v=2u^4 V, \quad \text{where} \quad V=2^{w} v,
$$ 
and $u,v$ are non-zero integers.
Substituting this into equation \eqref{x5p2xypy4} and cancelling $16 u^5 V^4$, we obtain
 $$
64 V^{11}+u^{11}=-1.
$$
This is a Thue equation whose only integer solution is $(V,u)=(0,-1)$, which is impossible as $V \neq 0$ by definition. Therefore, the integer solutions to equation \eqref{x5p2xypy4} are
$$
(x,y)=(0,0) \quad \text{and} \quad (1,-1).
$$

\vspace{10pt}

The next equation we will consider is
\begin{equation}\label{x5mxy2py4}
x^5-xy^2+y^4=0,
\end{equation}
which is equation \eqref{ax5pbxy2pcy4} with $a=c=1$ and $b=-1$.
 We have an obvious solution $(x,y)=(0,0)$ and otherwise $x\neq 0$ and $y\neq 0$. There are no prime divisors of $abc$, so $X_i=1$ and $Y_i=1$. Therefore, substitutions \eqref{ax5pbxy2pcy4lmnk} take the form $x=uv^2$ and $y=u^2v$, where $u,v$ are non-zero integers. Substituting these into \eqref{x5mxy2py4} and cancelling $u^5 v^4$ we obtain
 $$
 v^6-1+u^3=0=V^3-1+u^3
 $$
 where $V=v^2$. This is a Thue equation whose only integer solutions are $(V,u)=(1,0),(0, 1)$, which are impossible as $u$ and $v$ (and $V$) are non-zero by definition. Therefore, the unique integer solution to equation \eqref{x5mxy2py4} is 
$$
(x,y)=(0,0).
$$

\vspace{10pt}

The next equation we will consider is
\begin{equation}\label{x5pxy2py4}
x^5+xy^2+y^4=0,
 \end{equation}
 which is equation \eqref{ax5pbxy2pcy4} with $a=b=c=1$.
 We have an obvious solution $(x,y)=(0,0)$ and otherwise $x\neq 0$ and $y\neq 0$. There are no prime divisors of $abc$, so $X_i=1$ and $Y_i=1$. Therefore, substitutions \eqref{ax5pbxy2pcy4lmnk} take the form $x=uv^2$ and $y=u^2v$, where $u,v$ are non-zero integers. Substituting these into \eqref{x5pxy2py4} and cancelling $u^5 v^4$, we obtain
 $$
 v^6+1+u^3=0=V^3+1+u^3
 $$
 where $V=v^2$. This is a Thue equation whose only integer solutions are $(V,u)=(0,-1),(-1,0)$, which are impossible as $u,v$ (and $V$) are non-zero by definition. Therefore, the unique integer solution to equation \eqref{x5pxy2py4} is 
 $$
(x,y)=(0,0).
$$
 
\vspace{10pt}
 
 The final equation we will consider is
 \begin{equation}\label{x5px2ypy4}
x^5+x^2y+y^4=0.
\end{equation}
We have an obvious solution $(x,y)=(0,0)$, and otherwise $x \neq 0$ and $y \neq 0$.  So now, assume that $x$ and $y$ are non-zero. This is equation \eqref{eq:axnpbxkylpcym} with $a=b=c=1$, $n=5$, $k=2$, $l=1$ and $m=4$, so $n-k=3$ and $m-l=3$. Then
$$
l'=\frac{1}{\text{gcd}(1,3)}=1 \,; \quad n'=\frac{3}{\text{gcd}(1,3)}=3 \,; \quad m'=\frac{3}{\text{gcd}(2,3)}=3 \,; \quad k'= \frac{2}{\text{gcd}(2,3)}=2 .
$$
Because $abc$ has no prime divisors, $X_i=Y_i=1$.
So 
$$
x = X_i u^{l'} v^{m'}=uv^3, \quad y = Y_i u^{n'} v^{k'}=  u^{3} v^2
$$
for some non-zero integers $u,v$. 
Substituting these into equation \eqref{x5px2ypy4} and cancelling $u^5 v^8$, we obtain
$$
v^7+u^7=-1.
$$
This is a Thue equation whose only integer solutions are $(v,u)=(-1,0),(0,-1)$, which are impossible as $v\neq 0$ and $u \neq 0$ by definition. Therefore, the unique integer solution to equation \eqref{x5px2ypy4} is
 $$
(x,y)=(0,0).
$$
 
 Table \ref{tab:H562mon3varsol} summarises the integer solutions to the equations solved in this section.
 
 \begin{center}
\begin{tabular}{ |c|c||c|c|c| } 
\hline
 Equation & Solution $(x,y)$&  Equation & Solution $(x,y)$ \\ 
\hline\hline
$x^4+xy+2y^3=0$ & $(-1,-1),(0,0)$ &  $2x^4+xy+2y^3=0$ & $(0,0)$ \\ \hline
 $x^4+2xy+2y^3=0$ & $(0,0)$ &  $x^5+3xy+y^3=0$ & $(0,0)$  \\ \hline
 $x^4+3xy+2y^3=0$ & $(-6,-9),(-1,1),$ & $x^5+xy+2y^3=0$ & $(-1,1),(0,0)$  \\ 
& $(0,0),(3,-3)$ && \\\hline
 $x^4+xy+3y^3=0$ & $(0,0)$ &  $x^5+xy+y^4=0$ & $(0,0)$  \\ \hline
 $2x^4+xy+y^3=0$ & $(0,0),(1,-1)$ & $x^4+6xy+2y^3=0$ & $(0,0)$   \\ \hline
 $x^5+xy+y^3=0$ & $(0,0)$ &  $x^4+4xy+3y^3=0$ & $(-1,1),(0,0),(4,-4)$ \\ \hline
 $x^4+4xy+2y^3=0$ & $(0,0)$ & $x^4+2xy+4y^3=0$& $(0,0)$  \\ \hline
 $x^4+2xy+3y^3=0$ & $(-2,-2),(-1,-1),(0,0)$ & $2x^4+2xy-y^4=0$ & $(0,0)$    \\ \hline
 $2x^4+2xy+y^3=0$ & $(0,0)$ & $2x^4+4xy+y^3=0$& $(-2,-4),(0,0)$  \\ \hline
 $x^5+2xy+y^3=0$ & $(-1,-1),(0,0)$ &  $x^5+4xy+y^3=0$& $(-2,4),(0,0)$  \\ \hline
 $x^4+5xy+2y^3=0$ & $(0,0)$ & $x^5+2xy+2y^3=0$ & $(0,0)$   \\ \hline
 $x^4+3xy+3y^3=0$ & $(0,0)$ & $x^5+2xy+y^4=0$ & $(0,0),(1,-1)$  \\ \hline
 $x^4+xy+4y^3=0$ & $(0,0)$ & $x^5-xy^2+y^4=0$ & $(0,0)$  \\ \hline
 $2x^4+xy-y^4=0$ & $(0,0),\pm(1,-1)$ & $x^5+xy^2+y^4=0$ & $(0,0)$  \\ \hline
 $2x^4+3xy+y^3=0$ & $(-3,-6),(-1,-2),$ &$x^5+x^2y+y^4=0$ & $(0,0)$ \\
  & $(-1,1),(0,0)$ &&  \\ \hline
\end{tabular}
\captionof{table}{\label{tab:H562mon3varsol} Integer solutions to the equations in Table \ref{tab:H562mon3var}.}
\end{center}

\subsection{Exercise 4.56}\label{ex:x2py2pz2mxyzme}
\textbf{\emph{Solve all equations of the form
	$$
	x^2 + y^2 + z^2 = xyz + e
	$$
	for integer $|e|\leq 5$.}}
	
To describe all integer solutions to
\begin{equation}\label{eq:x2py2pz2mxyzme}
x^2+y^2+z^2=xyz+e,
\end{equation}
for $|e| \leq 5$, we will use the Vieta jumping method as we did in Sections \ref{ex:Vieta2var} and \ref{ex:H18vieta}, however, we will now be jumping in three variables. 
As we are jumping in three variables, we will have three Vieta jumping operations
$$
\begin{aligned}
(i) \,\, (x,y,z) \to (x',y,z), \quad x'=yz-x, \\
(ii) \,\, (x,y,z) \to (x,y',z), \quad y'=xz-y, \\
(iii) \,\, (x,y,z) \to (x,y,z'), \quad z'=xy-z. \\
\end{aligned}
$$
Recall that a minimal solution is one for which the value $N=|x|+|y|+|z|$ is not decreased by any of the transformations (i), (ii) or (iii).
We must first prove that the lower norm of any minimal solution is bounded, by solving the optimisation problem
$$
\max_{(x,y,z) \in \mathbb{R}^3} 
\min\{|x|,|y|,|z|\} \quad \text{subject to} \quad P(x,y,z)=0,  |x'|\geq|x|, |y'|\geq |y|, |z'|\geq |z|.
$$
We can do this in Mathematica using the command
\begin{equation}\label{maxval_4.54}
\begin{aligned}
{\tt N[MaxValue[Min[Abs[x],Abs[y],Abs[z]], \{x^2 + y^2 + z^2 == x y z + e , Abs[z y - x] \geq Abs[x],} \phantom{\}} \\  \phantom{\{} {\tt Abs[z x - y] \geq Abs[y],    Abs[x y - z] \geq Abs[z]\}, \{x, y, z\}, Reals]]}.
   \end{aligned}
   \end{equation}
   
   As \eqref{eq:x2py2pz2mxyzme} is symmetric, we only need to check the case $|x| \leq t$, where $t$ is the optimal value. Checking these values outputs a set of solutions containing all minimal solutions, from which we can determine the minimal solutions. We can find all integer solutions to the equation from the minimal ones by applying a sequence of the transformations (i)-(iii). To reduce the number of initial minimal solutions to start with, we also consider elementary transformations of permuting the variables and changing signs of two variables, that is, we consider the following transformations: 
   \begin{equation}\label{3var_vieta_sym}
   	\begin{aligned}
(a) & \,\, (x,y,z) \to (yz-x,y,z), & \quad (b) & \,\, (x,y,z) \to (x,z,y), \\
(c) & \,\, (x,y,z) \to (y,x,z), & \quad  (d) & \,\, (x,y,z) \to (-x,-y,z).
\end{aligned}
\end{equation}
We will call minimal solutions equivalent if they can be transformed to each other by transformations \eqref{3var_vieta_sym}. Then it is sufficient to select only one minimal solution from each equivalence class.

Let us first consider the equation \eqref{eq:x2py2pz2mxyzme} with $e=-5$, or,
\begin{equation}\label{x2py2pz2mxyzp5}
x^2+y^2+z^2=xyz-5.
\end{equation}
The MaxValue command \eqref{maxval_4.54} outputs $3.42599$. As \eqref{x2py2pz2mxyzp5} is symmetric, we only need to consider the case $|x| \leq 3$. By checking these values of $x$, we do not obtain any integer solutions to \eqref{x2py2pz2mxyzp5}. Therefore, equation \eqref{x2py2pz2mxyzp5} has no integer solutions.

\vspace{10pt}

The next equation we will consider is \eqref{eq:x2py2pz2mxyzme} with $e=-4$, or,
\begin{equation}\label{x2py2pz2mxyzp4}
x^2+y^2+z^2=xyz-4.
\end{equation}
The MaxValue command \eqref{maxval_4.54} outputs $3.3553$. 
Hence, we only need to consider the case $|x| \leq 3$. By checking these values of $x$, we do not obtain integer solutions to \eqref{x2py2pz2mxyzp4}. Therefore, equation \eqref{x2py2pz2mxyzp4} has no integer solutions.

\vspace{10pt}

The next equation we will consider is \eqref{eq:x2py2pz2mxyzme} with $e=-3$, or,
\begin{equation}\label{x2py2pz2mxyzp3}
x^2+y^2+z^2=xyz-3.
\end{equation}
The MaxValue command \eqref{maxval_4.54} outputs $3.27902$. 
By checking the case $|x| \leq 3$, we do not obtain integer solutions to \eqref{x2py2pz2mxyzp3}. Therefore, equation \eqref{x2py2pz2mxyzp3} has no integer solutions.

\vspace{10pt}

The next equation we will consider is \eqref{eq:x2py2pz2mxyzme} with $e=-2$, or,
\begin{equation}\label{x2py2pz2mxyzp2}
x^2+y^2+z^2=xyz-2.
\end{equation}
The MaxValue command \eqref{maxval_4.54} outputs $3.19582$. 
By checking the case $|x| \leq 3$, we obtain that, up to equivalence, the only minimal solution is $(x,y,z)=(3,3,4)$. Therefore we can obtain all integer solutions to \eqref{x2py2pz2mxyzp2} by applying a sequence of transformations \eqref{3var_vieta_sym} to $(x,y,z)=(3,3,4)$.

\vspace{10pt}

The next equation we will consider is \eqref{eq:x2py2pz2mxyzme} with $e=-1$, or,
\begin{equation}\label{x2py2pz2mxyzp1}
x^2+y^2+z^2=xyz-1.
\end{equation}
The MaxValue command \eqref{maxval_4.54} outputs $3.1038$. 
 By checking the case $|x| \leq 3$, we do not obtain integer solutions to \eqref{x2py2pz2mxyzp1}. Therefore, equation \eqref{x2py2pz2mxyzp1} has no integer solutions.

\vspace{10pt}

The next equation we will consider is \eqref{eq:x2py2pz2mxyzme} with $e=0$, or,
\begin{equation}\label{x2py2pz2mxyz}
x^2+y^2+z^2=xyz.
\end{equation}
The MaxValue command \eqref{maxval_4.54} outputs $3$. 
By checking the case $|x| \leq 3$, we obtain that, up to equivalence, the minimal solutions are
$$
(x,y,z)=(0,0,0),(3,3,3).
$$
Applying transformations \eqref{3var_vieta_sym} to $(0,0,0)$, we only obtain $(0,0,0)$.
To summarise, excluding $(x,y,z)=(0,0,0)$, we can obtain all integer solutions to \eqref{x2py2pz2mxyz} by applying a sequence of transformations \eqref{3var_vieta_sym} to $(x,y,z)=(3,3,3)$.

\vspace{10pt}

The next equation we will consider is \eqref{eq:x2py2pz2mxyzme} with $e=1$, or,
\begin{equation}\label{x2py2pz2mxyzm1}
x^2+y^2+z^2=xyz+1.
\end{equation}
The MaxValue command \eqref{maxval_4.54} outputs $2.87939$. 
By checking the case $|x| \leq 2$, we obtain that, up to equivalence, the only minimal solution is $(x,y,z)=(1,0,0)$. 
 Therefore, we can obtain all integer solutions to \eqref{x2py2pz2mxyzm1} by applying a sequence of transformations \eqref{3var_vieta_sym} to $(x,y,z)=(1,0,0)$. By applying these transformations, we see that all integer solutions to equation \eqref{x2py2pz2mxyzm1} are 
$$
(x,y,z)=(\pm1,0,0),(0,\pm 1,0), \quad \text{and} \quad (0,0,\pm 1).
$$

\vspace{10pt}

The next equation we will consider is \eqref{eq:x2py2pz2mxyzme} with $e=2$, or,
\begin{equation}\label{x2py2pz2mxyzm2}
x^2+y^2+z^2=xyz+2.
\end{equation}
The MaxValue command \eqref{maxval_4.54} outputs $2.73205$. 
By checking the case $|x| \leq 2$, we obtain that, up to equivalence, the only minimal solution is $(x,y,z)=(0,1,1)$. 
 Therefore, all integer solutions to equation \eqref{x2py2pz2mxyzm2} can be obtained by applying a sequence of transformations \eqref{3var_vieta_sym} to $(x,y,z)=(0,1,1)$.

\vspace{10pt}

The next equation we will consider is \eqref{eq:x2py2pz2mxyzme} with $e=3$, or,
\begin{equation}\label{x2py2pz2mxyzm3}
x^2+y^2+z^2=xyz+3.
\end{equation}
The MaxValue command \eqref{maxval_4.54} outputs $2.53209$. 
By checking the case $|x| \leq 2$, we do not obtain integer solutions to \eqref{x2py2pz2mxyzm3}. Therefore, equation \eqref{x2py2pz2mxyzm3} has no integer solutions.

\vspace{10pt}

The next equation we will consider is \eqref{eq:x2py2pz2mxyzme} with $e=4$, or,
\begin{equation}\label{x2py2pz2mxyzm4}
x^2+y^2+z^2=xyz+4.
\end{equation}
The MaxValue command \eqref{maxval_4.54} outputs $2$. 
By checking the case $|x| \leq 2$, we obtain that, up to equivalence, the only minimal solution is $(x,y,z)=(2,u,u)$ for integer $u$. 
Therefore, all integer solutions to equation \eqref{x2py2pz2mxyzm4} can be obtained by applying a sequence of transformations \eqref{3var_vieta_sym} to $(x,y,z)=(2,u,u)$.

\vspace{10pt}

The final equation we will consider is \eqref{eq:x2py2pz2mxyzme} with $e=5$, or,
\begin{equation}\label{x2py2pz2mxyzm5}
x^2+y^2+z^2=xyz+5.
\end{equation}
The MaxValue command \eqref{maxval_4.54} outputs $1.1038$. 
By checking the case $|x| \leq 1$, we obtain that, up to equivalence, the only minimal solution is $(x,y,z)=(0,1,2)$.
 To summarise, we can obtain all integer solutions to equation \eqref{x2py2pz2mxyzm5} by applying transformations \eqref{3var_vieta_sym} to $(x,y,z)=(0,1,2)$.

Table \eqref{tab:x2py2pz2mxyzme} summarises the list of minimal solutions, up to equivalence, to the equations solved in this exercise.

\begin{center}
\begin{tabular}{ |c|c|c|c|c|c| } 
\hline
 Equation &  Minimal Solutions $(x,y,z)$ \\ 
\hline\hline
$x^2+y^2+z^2=xyz-5$ & - \\\hline
 $x^2+y^2+z^2=xyz-4$ & - \\\hline
 $x^2+y^2+z^2=xyz-3$ & -  \\\hline
 $x^2+y^2+z^2=xyz-2$ & $(3,3,4)$ \\\hline
 $x^2+y^2+z^2=xyz-1$ & - \\\hline
 $x^2+y^2+z^2=xyz$ & $(0,0,0),(3,3,3)$ \\\hline
 $x^2+y^2+z^2=xyz+1$ & $(1,0,0)$ \\\hline
 $x^2+y^2+z^2=xyz+2$ & $(0, 1,1)$ \\\hline
 $x^2+y^2+z^2=xyz+3$ & - \\\hline
 $x^2+y^2+z^2=xyz+4$ & $(2,u,u), \quad u \in \mathbb{Z}$\\\hline
 $x^2+y^2+z^2=xyz+5$ & $(0,1,2)$ \\\hline

\end{tabular}
\captionof{table}{\label{tab:x2py2pz2mxyzme} Minimal solutions, up to equivalence, to equation \eqref{eq:x2py2pz2mxyzme} with $|e|\leq 5$. All integer solutions to these equations can be found by applying transformations \eqref{3var_vieta_sym} to the listed minimal solutions.}
\end{center}

\subsection{Exercise 4.60}\label{ex:H325sym2var}
\textbf{\emph{Solve all equations listed in Table \ref{H325sym2var} of size $H=324$. For the listed equations of size $H=325$, prove that they have a finite number of integer solutions and that there is an algorithm for listing these solutions.}}

	\begin{center}
		\begin{tabular}{ |c|c|c|c| } 
			\hline
			$H$ & Equation &  $H$ & Equation \\ 
			\hline\hline
			$324$ & $x^6+x^5 y+x^3 y^3+x y^5+y^6-xy=0$ &  $325$ & $x^6+x^5 y+x^3 y^3+x y^5+y^6+xy-1=0$ \\ 
			\hline
			$324$ & $x^6+x^5 y+x^3 y^3+x y^5+y^6+xy=0$ &  $325$ & $x^6-x^4 y^2+x^3 y^3-x^2 y^4+y^6-xy-1=0$ \\ 
			\hline
			$324$ & $x^6-x^4 y^2+x^3 y^3-x^2 y^4+y^6-xy=0$ &  $325$ & $x^6-x^4 y^2+x^3 y^3-x^2 y^4+y^6+xy+1=0$ \\ 
			\hline
			$324$ & $x^6-x^4 y^2+x^3 y^3-x^2 y^4+y^6+xy=0$ &  $325$ & $x^6-x^4 y^2+x^3 y^3-x^2 y^4+y^6+xy-1=0$ \\ 
			\hline
			$324$ & $x^6+3 x^3 y^3+y^6-xy=0$ &  $325$ & $x^6+3 x^3 y^3+y^6-xy-1=0$ \\ 
			\hline
			$324$ & $x^6+3 x^3 y^3+y^6+xy=0$ &  $325$ & $x^6+3 x^3 y^3+y^6+xy-1=0$ \\ 
			\hline
			$325$ & $x^6+x^5 y+x^3 y^3+x y^5+y^6-xy-1=0$ & & \\ 
			\hline
		\end{tabular}
		\captionof{table}{\label{H325sym2var} Some symmetric two-variable equations of size $H\leq 325$.}
	\end{center}

To find the integer solutions to the equations of size $H=324$ in Table \ref{H325sym2var}, we will reduce them using the substitutions
\begin{equation}\label{symsub}
w=(x+y)^2 \quad \text{and} \quad v=xy.
\end{equation}
In all examples, the resulting equation is of the form
\begin{equation}\label{eq:axpbypPxy}
	aw + bv = P(w,v),
\end{equation}
for homogeneous ${P(w,v)}$, and can be reduced to solving a finite number of Thue equations, as explained in Section \ref{ex:H28axbyPxy}. We can then solve the resulting Thue equations as we did in Section \ref{ex:Thue}.

The equations of size $H=325$ can be reduced by the substitutions \eqref{symsub} to equations of genus $g=1$, hence they have finitely many integer solutions. The Maple command
$$
\begin{aligned}
	& {\tt with(algcurves) } \\
	& {\tt genus(P(x,y), x, y) }
\end{aligned}
$$
can be used to find the genus of a curve.  The following theorem of Baker and Coates \cite{baker1970integer} can then be used to show that there exists an algorithm that can list these solutions. 

\begin{theorem}\label{th:genus1}[Theorem 3.61 in the book]
	There is an algorithm that, given an absolutely irreducible polynomial $P(x,y)$ of genus $g=1$ with integer coefficients, determines all integer solutions to $P(x,y)=0$.
\end{theorem}

The following equalities will be useful for reducing the equations of size $H=324$ and $H=325$ that we are considering,
$$
x^2+y^2 = (x+y)^2 - 2xy = w - 2v,
$$
$$
x^4+y^4 = (x+y)^4 - 4xy(x^2+y^2)-6x^2y^2 = w^2-4v(w - 2v) - 6v^2 = w^2 - 4wv+2v^2,
$$
$$
x^6+y^6 = (x^2+y^2)(x^4+y^4-x^2y^2)=(w - 2v)(w^2 - 4wv+2v^2-v^2)=w^3-6 w^2 v+9 w v^2-2 v^3,
$$
where $w,v$ are defined in \eqref{symsub}.

The first equation we will consider is
\begin{equation}\label{x6px5ypx3y3pxy5py6mxy}
x^6+x^5 y+x^3 y^3+x y^5+y^6-xy=0.
\end{equation}
Using the substitutions \eqref{symsub}, we can reduce this equation to 
\begin{equation}\label{x6px5ypx3y3pxy5py6mxy_red}
w^3-5w^2v+5wv^2+v^3-v=0.
\end{equation}
This is an equation of the form \eqref{eq:axpbypPxy}, and it can be solved by the method of Section \ref{ex:H28axbyPxy}. If $v=0$, then we have the solution $(v,w)=(0,0)$, leading to $(x,y)=(0,0)$. If $(v,w)$ is a solution to \eqref{x6px5ypx3y3pxy5py6mxy_red}, then $(-v,-w)$ is also a solution. So, we may first find solutions with $v >0$, and then take both signs. Let $d=\text{gcd}(v,w)$ and $w=d w_1$ and $v=d v_1$, where $w_1$ and $v_1>0$ are coprime. Substituting this into \eqref{x6px5ypx3y3pxy5py6mxy_red} and cancelling $d$, we obtain
\begin{equation}\label{x6px5ypx3y3pxy5py6mxy_red1}
d^2(w_1^3-5w_1^2v_1+5w_1v_1^2+v_1^3)=v_1.
\end{equation}
From this equation it is clear that $v_1$ is a divisor of $d^2 w_1^3$. Because $v_1$ and $w_1$ are coprime, this implies that $v_1$ is a divisor of $d^2$, so we can write $d^2=v_1z$ for some integer $z >0$. Then, substituting this into \eqref{x6px5ypx3y3pxy5py6mxy_red1} and cancelling $v_1$, we obtain
$$
z(w_1^3-5w_1^2v_1+5w_1v_1^2+v_1^3)=1.
$$
Hence $z>0$ is a divisor of $1$, so we must have $z=1$ and $w_1^3-5w_1^2v_1+5w_1v_1^2+v_1^3=1$. This is a Thue equation and its integer solutions with $v_1 >0$ are $(v_1,w_1)=(1,0),(2,7),(45,76)$. Because $d=\sqrt{v_1 z}$, the second and third solutions do not provide integer solutions. The first solution implies that $d=1$, so $v=d v_1=1$ and $w=dw_1=0$. Then using \eqref{symsub}, we have $(x+y)^2=0$ and $xy=1$. Hence, all integer solutions to equation \eqref{x6px5ypx3y3pxy5py6mxy} are
$$
(x,y)=\pm(1,-1),(0,0).
$$

\vspace{10pt}

The next equation we will consider is
\begin{equation}\label{x6px5ypx3y3pxy5py6pxy}
x^6+x^5 y+x^3 y^3+x y^5+y^6+xy=0.
\end{equation}
Using the substitutions \eqref{symsub}, we can reduce this equation to 
\begin{equation}\label{x6px5ypx3y3pxy5py6pxy_red}
w^3-5w^2v+5wv^2+v^3+v=0.
\end{equation}
If $v=0$, then we have the solution $(v,w)=(0,0)$, leading to $(x,y)=(0,0)$. If $(v,w)$ is a solution to \eqref{x6px5ypx3y3pxy5py6pxy_red}, then $(-v,-w)$ is also a solution. So, we may first find solutions with $v >0$, and then take both signs. Let $d=\text{gcd}(v,w)$ and $w=d w_1$ and $v=d v_1$, where $w_1$ and $v_1>0$ are coprime. Substituting this into \eqref{x6px5ypx3y3pxy5py6pxy_red} and cancelling $d$, we obtain
\begin{equation}\label{x6px5ypx3y3pxy5py6pxy_red1}
d^2(w_1^3-5w_1^2v_1+5w_1v_1^2+v_1^3)=-v_1.
\end{equation}
From this equation it is clear that $v_1$ is a divisor of $d^2 w_1^3$. Because $v_1$ and $w_1$ are coprime, this implies that $v_1$ is a divisor of $d^2$, so we can write $d^2=v_1z$ for some integer $z >0$. Then, substituting this into \eqref{x6px5ypx3y3pxy5py6pxy_red1} and cancelling $v_1$, we obtain
$$
z(w_1^3-5w_1^2v_1+5w_1v_1^2+v_1^3)=-1.
$$
Hence $z>0$ is a divisor of $1$, so we must have $z=1$ and $w_1^3-5w_1^2v_1+5w_1v_1^2+v_1^3=-1$. This is a Thue equation and its only integer solution with $v_1 >0$ is $(v_1,w_1)=(1,2)$. So, $d=\sqrt{v_1 z}=1$, so $v=d v_1=1$ and $w=dw_1=2$, but then $x+y=\sqrt{2}$, which is not integer. 
Therefore, the unique integer solution to equation \eqref{x6px5ypx3y3pxy5py6pxy} is
$$
(x,y)=(0,0).
$$

\vspace{10pt}

The next equation we will consider is
\begin{equation}\label{x6mx4y2px3y3mx2y4py6mxy}
x^6-x^4 y^2+x^3 y^3-x^2 y^4+y^6-xy=0.
\end{equation}
Using the substitutions \eqref{symsub}, we can reduce this equation to 
\begin{equation}\label{x6mx4y2px3y3mx2y4py6mxy_red}
w^3-6w^2v+8wv^2+v^3=v.
\end{equation}
If $v=0$, then we have the solution $(v,w)=(0,0)$, leading to $(x,y)=(0,0)$. If $(v,w)$ is a solution to \eqref{x6mx4y2px3y3mx2y4py6mxy_red}, then $(-v,-w)$ is also a solution. So, we may first find solutions with $v >0$, and then take both signs. Let $d=\text{gcd}(v,w)$ and $w=d w_1$ and $v=d v_1$, where $w_1$ and $v_1>0$ are coprime. Substituting this into \eqref{x6mx4y2px3y3mx2y4py6mxy_red} and cancelling $d$, we obtain
\begin{equation}\label{x6mx4y2px3y3mx2y4py6mxy_red1}
d^2(w_1^3-6w_1^2v_1+8w_1v_1^2+v_1^3)=v_1.
\end{equation}
From this equation it is clear that $v_1$ is a divisor of $d^2 w_1^3$. Because $v_1$ and $w_1$ are coprime, this implies that $v_1$ is a divisor of $d^2$, so we can write $d^2=v_1z$ for some integer $z >0$. Then, substituting this into \eqref{x6mx4y2px3y3mx2y4py6mxy_red1} and cancelling $v_1$, we obtain
$$
z(w_1^3-6w_1^2v_1+8w_1v_1^2+v_1^3)=1.
$$
Hence $z>0$ is a divisor of $1$, so we must have $z=1$ and $w_1^3-6w_1^2v_1+8w_1v_1^2+v_1^3=1$. This is a Thue equation and its integer solutions with $v_1 >0$ are $(v_1,w_1)=(1,2),(1,4),(4,9),(273,1054)$. Then using $d=\sqrt{v_1 z}$, $v=d v_1=1$ and $w=dw_1$, we obtain the integer solutions $(\sqrt{w},v)=(\pm 2,1),(0,1)$. Finally, we can find all integer solutions to \eqref{x6mx4y2px3y3mx2y4py6mxy} using \eqref{symsub}. Hence, all integer solutions to equation \eqref{x6mx4y2px3y3mx2y4py6mxy} are
$$
(x,y)=(\pm 1,\pm 1),(0,0).
$$

\vspace{10pt}

The next equation we will consider is
\begin{equation}\label{x6mx4y2px3y3mx2y4py6pxy}
x^6-x^4 y^2+x^3 y^3-x^2 y^4+y^6+xy=0.
\end{equation}
Using the substitutions \eqref{symsub}, we can reduce this equation to 
\begin{equation}\label{x6mx4y2px3y3mx2y4py6pxy_red}
w^3-6w^2v+8wv^2+v^3=-v.
\end{equation}
If $v=0$, then we have the solution $(v,w)=(0,0)$, leading to $(x,y)=(0,0)$. If $(v,w)$ is a solution to \eqref{x6mx4y2px3y3mx2y4py6pxy_red}, then $(-v,-w)$ is also a solution. So, we may first find solutions with $v >0$, and then take both signs. Let $d=\text{gcd}(v,w)$ and $w=d w_1$ and $v=d v_1$, where $w_1$ and $v_1>0$ are coprime. Substituting this into \eqref{x6mx4y2px3y3mx2y4py6pxy_red} and cancelling $d$, we obtain
\begin{equation}\label{x6mx4y2px3y3mx2y4py6pxy_red1}
d^2(w_1^3-6w_1^2v_1+8w_1v_1^2+v_1^3)=-v_1.
\end{equation}
From this equation it is clear that $v_1$ is a divisor of $d^2 w_1^3$. Because $v_1$ and $w_1$ are coprime, this implies that $v_1$ is a divisor of $d^2$, so we can write $d^2=v_1z$ for some integer $z >0$. Then, substituting this into \eqref{x6mx4y2px3y3mx2y4py6pxy_red1} and cancelling $v_1$, we obtain
$$
z(w_1^3-6w_1^2v_1+8w_1v_1^2+v_1^3)=-1.
$$
Hence $z>0$ is a divisor of $1$, so we must have $z=1$ and $w_1^3-6w_1^2v_1+8w_1v_1^2+v_1^3=-1$. This is a Thue equation and it has no solutions satisfying $v_1 >0$. Therefore, the unique integer solution to equation \eqref{x6mx4y2px3y3mx2y4py6pxy} is
$$
(x,y)=(0,0).
$$

\vspace{10pt}

The next equation we will consider is
\begin{equation}\label{x6p3x3y3py6mxy}
x^6+3x^3 y^3+y^6-xy=0.
\end{equation}
Using the substitutions \eqref{symsub}, we can reduce this equation to 
\begin{equation}\label{x6p3x3y3py6mxy_red}
w^3-6w^2v+9wv^2+v^3=v.
\end{equation}
If $v=0$, then we have the solution $(v,w)=(0,0)$, leading to $(x,y)=(0,0)$. If $(v,w)$ is a solution to \eqref{x6p3x3y3py6mxy_red}, then $(-v,-w)$ is also a solution. So, we may first find solutions with $v >0$, and then take both signs. Let $d=\text{gcd}(v,w)$ and $w=d w_1$ and $v=d v_1$, where $w_1$ and $v_1>0$ are coprime. Substituting this into \eqref{x6p3x3y3py6mxy_red} and cancelling $d$, we obtain
\begin{equation}\label{x6p3x3y3py6mxy_red1}
d^2(w_1^3-6w_1^2v_1+9w_1v_1^2+v_1^3)=v_1.
\end{equation}
From this equation it is clear that $v_1$ is a divisor of $d^2 w_1^3$. Because $v_1$ and $w_1$ are coprime, this implies that $v_1$ is a divisor of $d^2$, so we can write $d^2=v_1z$ for some integer $z >0$. Then, substituting this into \eqref{x6p3x3y3py6mxy_red1} and cancelling $v_1$, we obtain
$$
z(w_1^3-6w_1^2v_1+9w_1v_1^2+v_1^3)=1.
$$
Hence $z>0$ is a divisor of $1$, so we must have $z=1$ and $w_1^3-6w_1^2v_1+9w_1v_1^2+v_1^3=1$. This is a Thue equation and its integer solutions with $v_1 >0$ are $(v_1,w_1)=(1,0),(1,3)$. As $v_1=1$, $d=\sqrt{v_1 z}=1$, the second solution gives $w=dw_1=3$, so $x+y=\sqrt{3}$ which is not integer. The first solution with $v=d v_1=1$ and $w=dw_1$ gives the integer solution $(v,w)=(1,0)$. Finally, we can find all integer solutions to \eqref{x6p3x3y3py6mxy} using \eqref{symsub}. Hence, all integer solutions to equation \eqref{x6p3x3y3py6mxy} are
$$
(x,y)=\pm ( 1,-1),(0,0).
$$

\vspace{10pt}

The next equation we will consider is
\begin{equation}\label{x6p3x3y3py6pxy}
x^6+3x^3 y^3+y^6+xy=0.
\end{equation}
Using the substitutions \eqref{symsub}, we can reduce this equation to 
\begin{equation}\label{x6p3x3y3py6pxy_red}
w^3-6w^2v+9wv^2+v^3=-v.
\end{equation}
If $v=0$, then we have the solution $(v,w)=(0,0)$, leading to $(x,y)=(0,0)$. If $(v,w)$ is a solution to \eqref{x6p3x3y3py6pxy_red}, then $(-v,-w)$ is also a solution. So, we may first find solutions with $v >0$, and then take both signs. Let $d=\text{gcd}(v,w)$ and $w=d w_1$ and $v=d v_1$, where $w_1$ and $v_1>0$ are coprime. Substituting this into \eqref{x6p3x3y3py6pxy_red} and cancelling $d$, we obtain
\begin{equation}\label{x6p3x3y3py6pxy_red1}
d^2(w_1^3-6w_1^2v_1+9w_1v_1^2+v_1^3)=-v_1.
\end{equation}
From this equation it is clear that $v_1$ is a divisor of $d^2 w_1^3$. Because $v_1$ and $w_1$ are coprime, this implies that $v_1$ is a divisor of $d^2$, so we can write $d^2=v_1z$ for some integer $z >0$. Then, substituting this into \eqref{x6p3x3y3py6pxy_red1} and cancelling $v_1$, we obtain
$$
z(w_1^3-6w_1^2v_1+9w_1v_1^2+v_1^3)=-1.
$$
Hence $z>0$ is a divisor of $1$, so we must have $z=1$ and $w_1^3-6w_1^2v_1+9w_1v_1^2+v_1^3=-1$. This is a Thue equation and it has no solutions satisfying $v_1 >0$.
Therefore, the unique integer solution to equation \eqref{x6p3x3y3py6pxy} is
$$
(x,y)=(0,0).
$$

\begin{center}
\begin{tabular}{ |c|c|c|c| } 
\hline
 Equation & Solution $(x,y)$ \\ 
\hline\hline
 $x^6+x^5 y+x^3 y^3+x y^5+y^6-xy=0$ & $\pm(1,-1),(0,0)$  \\\hline
$x^6+x^5 y+x^3 y^3+x y^5+y^6+xy=0$ & $(0,0)$ \\\hline
 $x^6-x^4 y^2+x^3 y^3-x^2 y^4+y^6-xy=0$ & $(\pm 1, \pm 1),(0,0)$ \\\hline
 $x^6-x^4 y^2+x^3 y^3-x^2 y^4+y^6+xy=0$ & $(0,0)$ \\\hline
 $x^6+3 x^3 y^3+y^6-xy=0$ &$(0,0),\pm(1,-1)$ \\\hline
 $x^6+3 x^3 y^3+y^6+xy=0$ &$(0,0)$ \\\hline

\end{tabular}
\captionof{table}{\label{H324sym2varsol} Integer solutions to the equations in Table \ref{H325sym2var} of size $H=324$.}
\end{center}

We will now consider the equations in Table \ref{H325sym2var} of size $H=325$.
The first equation we will consider is
\begin{equation}\label{x6px5ypx3y3pxy5py6mxym1}
x^6+x^5 y+x^3 y^3+x y^5+y^6-xy-1=0.
\end{equation}
Using the substitutions \eqref{symsub}, equation \eqref{x6px5ypx3y3pxy5py6mxym1} is reduced to 
\begin{equation}\label{x6px5ypx3y3pxy5py6mxym1_red}
w^3-5w^2v+5wv^2+v^3-v-1=0.
\end{equation}
The Maple command 
$$
{\tt genus(w^3-5w^2v+5wv^2+v^3-v-1, w,v)}
$$
outputs $1$, while the Mathematica command
$$
{\tt IrreduciblePolynomialQ[w^3-5w^2v+5wv^2+v^3-v-1, Extension \to All]}
$$
outputs True. Hence, \eqref{x6px5ypx3y3pxy5py6mxym1_red} is an absolutely irreducible polynomial with genus $1$. Therefore, using Theorem \ref{th:genus1}, \eqref{x6px5ypx3y3pxy5py6mxym1_red} has finitely many integer solutions. Then for each solution $(w,v)$, we can check whether $w$ is a perfect square and then check whether equations \eqref{symsub} are solvable in integers. Hence equation \eqref{x6px5ypx3y3pxy5py6mxym1} has finitely many integer solutions.

The other equations of size $H=325$ in Table \ref{H325sym2var} can be analysed similarly and for each equation, Table \ref{H325sym2vargenus} presents the absolutely irreducible equation it reduces to after substitutions \eqref{symsub}, and the genus of that equation. Because the genus is equal to $1$, these equations have finitely many integer points, and they can all be listed by Theorem \ref{th:genus1}. For each solution $(w,v)$, we can check whether $w$ is a perfect square and then find $u = \pm \sqrt{w}$ and check whether the system of equations $u=x+y$ and $v=xy$ is solvable in integers $x,y$.

\begin{center}
\begin{tabular}{ |c|c|c|c| } 
\hline
 Equation & Reduced Equation & Genus   \\ \hline \hline

 $x^6+x^5 y+x^3 y^3+x y^5+y^6-xy-1=0$ & $w^3-5w^2v+5wv^2+v^3-v-1=0$& $1$  \\ 
\hline
 $x^6+x^5 y+x^3 y^3+x y^5+y^6+xy-1=0$ & $w^3-5w^2v+5wv^2+v^3+v-1=0$& $1$ \\ 
\hline
 $x^6-x^4 y^2+x^3 y^3-x^2 y^4+y^6-xy-1=0$ &  $w^3-6w^2v+8wv^2+v^3-v-1=0$& $1$ \\ 
\hline
 $x^6-x^4 y^2+x^3 y^3-x^2 y^4+y^6+xy+1=0$ & $w^3-6w^2v+8wv^2+v^3+v+1=0$&  $1$ \\ 
\hline
 $x^6-x^4 y^2+x^3 y^3-x^2 y^4+y^6+xy-1=0$ & $w^3-6w^2v+8wv^2+v^3+v-1=0$&  $1$ \\ 
\hline
 $x^6+3 x^3 y^3+y^6-xy-1=0$ & $w^3-6w^2v+9wv^2+v^3-v-1=0$&  $1$ \\ 
\hline
 $x^6+3 x^3 y^3+y^6+xy-1=0$ & $w^3-6w^2v+9wv^2+v^3+v-1=0$&  $1$ \\ 
\hline
\end{tabular}
\captionof{table}{\label{H325sym2vargenus} Equations from Table \ref{H325sym2var} of size $H=325$. In all equations, $v=xy$ and $w=(x+y)^2$.}
\end{center}

\subsection{Exercise 4.66}\label{ex:x3py3pz3pm2xyz}
\textbf{\emph{Use the following proposition to find all integer solutions to the equations
	$$
		x^3+y^3+z^3\pm 2xyz=0.
	$$
	\begin{proposition}\label{prop:x3py3pcz3mdxyz}[Proposition 4.65 in the book]
	Let $c,d$ be integers. If integers $(x,y,z)$ with $z\neq 0$ solve equation
	\begin{equation}\label{eq:x3py3pcz3mdxyz}
		x^3+y^3+c z^3-d xyz=0
	\end{equation} 
	then the rational numbers 
	$$
	X=-4\frac{xy}{z^2}, \quad Y=4\frac{x^3-y^3}{z^3} 
	$$
	solve equation
	\begin{equation}\label{eq:x3py3pcz3mdxyza}
		Y^2 = X^3 + (d X + 4 c)^2.
	\end{equation} 
\end{proposition}}}
	
The equations 
\begin{equation}\label{eq:x3py3pz3pm2xyz}
	(i) \quad x^3+y^3+z^3 + 2xyz=0, \quad (ii) \quad x^3+y^3+z^3 - 2xyz=0 
\end{equation} 
are of the form \eqref{eq:x3py3pcz3mdxyz} with $c=1$ and $d=\pm 2$. In this case, equation \eqref{eq:x3py3pcz3mdxyza} takes the form
$$
Y^2 = X^3 + (\pm 2X + 4)^2 = X^3 + 4X^2 \pm 16 X + 16.
$$
These are equations in Weierstrass form \eqref{eq:Weiform} with $(a,b,c,d,e)=(0,4,0,\pm 16,16)$, and their rational points can be computed as explained in Section \ref{ex:H27rank0ell}. In this example, the Magma commands
$$
{\tt Rank(EllipticCurve([0,4,0,-16,16]));}
$$
$$
{\tt Rank(EllipticCurve([0,4,0,16,16]));}
$$
return ``0 true'', hence both equations define curves of rank $0$, and their only rational points are the torsion points. Using the SageMath commands
$$
{\tt EllipticCurve([0,4,0,-16,16]).torsion\_points()}
$$
$$
{\tt EllipticCurve([0,4,0,16,16]).torsion\_points()}
$$
we can easily compute that the only rational solutions in both cases are $(X,Y)=(0,\pm 4)$. In particular, there are no rational solutions with $X\neq 0$, hence, by Proposition \ref{prop:x3py3pcz3mdxyz}, all integer solutions to equations \eqref{eq:x3py3pz3pm2xyz} must satisfy $xyz=0$. From this, it is easy to deduce that all integer solutions to equations \eqref{eq:x3py3pz3pm2xyz} are
\begin{equation}\label{sol:x3py3pz3p2xyz}
(x,y,z)=(0,u,-u), \quad (u,0,-u),\quad \text{and} \quad (u,-u,0), \quad u \in \mathbb{Z}.
\end{equation}

\subsection{Exercise 4.71}\label{ex:x2ypy2zpz2xpm2xyz}	
\textbf{\emph{Describe all integer solutions to the equations
	$$
	x^2y+y^2z+z^2x \pm 2xyz=0.
	$$}}
	
With $a=x^2y$, $b=y^2z$ and $c=z^2x$, we may observe that 
\begin{equation} \label{eq:x2ypy2zpz2xpm2xyz}
x^2y+y^2z+z^2x \mp 2xyz=0
\end{equation}
 implies that $a+b+c=\pm 2xyz$ and 
\begin{equation} \label{eq:8abcapbpc3}
8abc=8(x^2y)(y^2z)(z^2x)=(2xyz)^3=\pm(a+b+c)^3.
\end{equation}
Let us prove that any integer solution $(a,b,c)$ to equation \eqref{eq:8abcapbpc3} must satisfy $abc=0$. By contradiction, assume that there is a solution with $abc\neq 0$. By cancelling any common factors if necessary, we may assume that $\gcd(a,b,c)=1$. If, for example, $a$ and $b$ have any common prime factor $p$, then $p$ is a factor of $abc$, hence $p$ is a factor of $(a+b+c)^3$, which implies that $p$ is a factor of $c$, a contradiction with $\text{gcd}(a,b,c)=1$. Hence, $a,b,c$ are pairwise coprime. We have that $8abc$ must be a perfect cube which implies that $a,b,c$ are all perfect cubes, that is $a=u^3$, $b=v^3$ and $c=w^3$ for some integers $u,v,w$. But then $\pm(u^3+v^3+w^3)^3=8u^3v^3w^3=(2uvw)^3$ simplifies to $u^3+v^3+w^3=\pm 2uvw$. Up to the names of variables these are equations \eqref{eq:x3py3pz3pm2xyz} whose solutions \eqref{sol:x3py3pz3p2xyz}  satisfy $uvw=0$. This implies that $abc=0$, hence $xyz=0$, and therefore all integer solutions to equation \eqref{eq:x2ypy2zpz2xpm2xyz} are
$$
(x,y,z)=(0,u,0),(0,0,u),(u,0,0), \quad u \in \mathbb{Z}.
$$
	
\section{Chapter 5}
In the next exercises, we will solve the following problem. 
\begin{problem}\label{prob:fin}[Problem 5.1 in the book]
	Given a Diophantine equation, either list all its integer solutions, or prove that there are infinitely many of them. 
\end{problem}

\subsection{Exercise 5.2}\label{ex:H19trivinf}
\textbf{\emph{For each of the equations listed in Table \ref{tab:H19trivinf}, find the values of two variables such that the equation is true for any value of the remaining variable. }}

	\begin{center}

		\captionof{table}{\label{tab:H19trivinf} Some equations with trivially infinite solution sets of size $H\leq 19$.}
	\end{center} 

Table \ref{tab:H19trivinfsol} presents, for each equation in Table \ref{tab:H19trivinf}, an infinite family of integer solutions.

	\begin{center}
		%
		\captionof{table}{\label{tab:H19trivinfsol} An example of an infinite set of solutions for the equations listed in Table \ref{tab:H19trivinf}.  Assume that $u$ is an integer.}
	\end{center} 	
	
\subsection{Exercise 5.3}\label{ex:H20trivinf}
\textbf{\emph{For each of the equations listed in Table \ref{tab:H20trivinf}, find the values of some variables that reduce the equation to a two-monomial two-variable equation with trivially infinite set of integer solutions. }}
	\begin{center}
		%
		\captionof{table}{\label{tab:H20trivinf} Some equations with trivially infinite solution sets of size $H\leq 20$.}
	\end{center} 

For each equation in Table \ref{tab:H20trivinf}, we will reduce the equation to a two-monomial two-variable equation with infinitely many solutions. Table \ref{tab:H20trivinfsol} presents, for each equation in Table \ref{tab:H20trivinf}, a suitable value to substitute into the equation, the two-variable two-monomial equation it reduces to, and the resulting infinite integer solution to the original equation.

	\begin{center}
		%
		\captionof{table}{\label{tab:H20trivinfsol} Trivially infinite solution sets to the equations listed in Table \ref{tab:H20trivinf}. Assume that $u$ is an integer.}
	\end{center}

\subsection{Exercise 5.5}\label{ex:H19Gauss}
\textbf{\emph{For each of the equations listed in Table \ref{tab:H19Gauss}, find the value of one of the variables, such that the resulting two-variable equation is quadratic and has infinitely many integer solutions by either of the following propositions. 
		\begin{proposition}[Gauss's Theorem]\label{prop:Gaussquad}[Proposition 3.14 in the book]
			Assume that the general two-variable quadratic equation with integer coefficients \eqref{eq:quad2vargen2} is such that $D=b^2-4ac > 0$, $D$ is not a perfect square, and $\Delta = 4acf + bde - ae^2 - cd^2 - fb^2 \neq 0$. Then if equation \eqref{eq:quad2vargen2} has at least one integer solution then it has infinitely many of them. Moreover, if $b\leq 0$ and \eqref{eq:quad2vargen2} has any integer solution $(x_0,y_0)$ satisfying $Dy_0 + b d - 2 a e > 0$ and $2 a x_0 + b y_0 + d > 0$, then it has infinitely many solutions in \textbf{positive} integers $(x,y)$.
\end{proposition}
\begin{proposition}\label{prop:Gaussquad2}[Proposition 5.4 in the book]
	Let $a,b,c$ be integers such that $a>0$, $a$ is not a perfect square, $b^2-4ac\neq 0$, and there exists an integer $x$ such that $P(x)=ax^2+bx+c$ is a perfect square. Then $P(x)$ is a perfect square for infinitely many integers $x$.
\end{proposition} }}

	\begin{center}

		\captionof{table}{\label{tab:H19Gauss} Equations of size $H\leq 19$ having infinitely many integer solutions by Gauss's Theorem.}
	\end{center}

For each equation in Table \ref{tab:H19Gauss}, there exists a value of variable $y$ such that the resulting quadratic equation in $(x,z)$ has an integer solution $(x_0,z_0)$ satisfying the conditions in Proposition \ref{prop:Gaussquad}, see Table \ref{tab:H19Gausssol}. This implies that the original equations have infinitely many integer solutions.

\begin{center}
%
\captionof{table}{\label{tab:H19Gausssol} Reducing the equations in Table \ref{tab:H19Gauss} to two-variable quadratic equations with infinitely many integer solutions by Proposition \ref{prop:Gaussquad}.}
\end{center}

\subsection{Exercise 5.6}\label{ex:H24Gausssub}
\textbf{\emph{Use a suitable linear substitution in combination with Proposition \ref{prop:Gaussquad} or \ref{prop:Gaussquad2} to prove that each of the equations listed in Table \ref{tab:H24Gausssub} has infinitely many integer solutions. }}

	\begin{center}
		%
		\captionof{table}{\label{tab:H24Gausssub} Equations of size $H\leq 24$ solvable by Gauss's Theorem after a linear substitution.}
	\end{center}

Let us consider, for example, the equation
\begin{equation}\label{eq:x3py3pz2p2}
	x^3+y^3+z^2+2=0.
\end{equation}
Using a well-known factorization for the sum of cubes, we can rewrite this equation as
\begin{equation}\label{eq:x3py3pz2p2red}
	(x+y)(x^2-xy+y^2)+z^2+2=0.
\end{equation}
This representation suggests to look for solutions such that $x+y$ is equal to some fixed constant $c$. In this case, $x=c-y$, and \eqref{eq:x3py3pz2p2red} reduces to
\[
c((c-y)^2-(c-y)y+y^2)+z^2+2=0,
\]
or equivalently
\begin{equation}\label{eq:x3py3pz2p2red2}
	z^2 = -3cy^2 +3c^2y -c^3 -2.
\end{equation}
For every fixed $c$, this is a quadratic equation in $y$ and $z$, so we may hope to prove that it has infinitely many integer solutions using Propositions \ref{prop:Gaussquad} or \ref{prop:Gaussquad2}. In this example, Proposition \ref{prop:Gaussquad2} is applicable provided that $c<0$, $3|c|$ is not a perfect square, $(3c^2)^2-4(-3c)(-c^3-2)\neq 0$, and \eqref{eq:x3py3pz2p2red2} has \emph{some} integer solution $(y_0,z_0)$. A quick computer search with small values of $c$ returns that these conditions are satisfied for $c=-11$, in which case \eqref{eq:x3py3pz2p2red2} reduces to
$$
z^2=33y^2+363y+1329.
$$
Indeed, $33>0$, not a perfect square, $363^2-4 \cdot 33\cdot 1329 \neq 0$, and the equation has a solution $(y_0,z_0)=(5,63)$. Hence, it has infinitely many integer solutions by Proposition \ref{prop:Gaussquad2}. This implies that the original equation \eqref{eq:x3py3pz2p2} has infinitely many integer solutions as well.

The other equations in Table \ref{tab:H24Gausssub} can be solved similarly. Table \ref{tab:H24Gausssubsol} presents, for each such equation, a linear substitution in the form $x=c-y$ for the appropriately chosen constant $c$, the two-variable quadratic equation it reduces to, and an integer solution $(y_0,z_0)$ to this equation. Then Proposition \ref{prop:Gaussquad} or  \ref{prop:Gaussquad2} guarantees that the original equation has infinitely many integer solutions.

\begin{center}
\begin{tabular}{ |c|c|c|c|c|c| } 
\hline
 Equation & Assumption & Reduced Quadratic & Solution $(y_0,z_0)$ \\ 
\hline\hline
 $x^3+y^3+z^2-2=0$ & $x=-1-y$ & $z^2=3y^2+3y+3$ & $(-2,3)$  \\\hline
 $x^3+y^3+z^2+2=0$ & $x=-11-y$ & $z^2=33y^2+363y+1329$ & $(5,63)$ \\\hline
 $x^3+y^3+z^2-3=0$ & $x=-1-y$ & $z^2=3y^2+3y+4$ & $(-1,2)$ \\\hline
 $x^3+y^3+z^2+3=0$ & $x=-1-y$ & $z^2=3y^2+3y-2$ & $(-2,2)$ \\\hline
 $x^3-x+y^3+z^2-1=0$ & $x=-1-y$ & $z^2=3y^2+2y+1$ & $(-2,3)$ \\\hline
$x^3-x+y^3+z^2+1=0$ & $x=-1-y$ & $z^2=3y^2+2y-1$ & $(-1,0)$ \\\hline
 $x^3+x+y^3+z^2-1=0$ & $x=-2-y$ & $z^2=6y^2+13y+11$ & $(-2,3)$ \\\hline
 $x^3+x+y^3+z^2+1=0$ &  $x=-1-y$ & $z^2=3y^2+4y+1$ & $(-1,0)$ \\\hline
 $x^3+y^3+z^2+z+1=0$ & $x=-1-y$ & $z^2+z=3y^2+3y$ & $(-1,0)$ \\\hline
 $x^3+y^3+z^2+z-1=0$ & $x=-1-y$ & $z^2+z=3y^2+3y+2$ & $(-1,1)$ \\\hline
 $x^3-x+y^3-y+z^2=0$ & $x=-1-y$ & $z^2=3y^2+3y$ & $(0,0)$ \\\hline
 $x^3-x+y^3+y+z^2=0$ & $x=-1-y$ & $z^2=3y^2+y$ & $(0,0)$ \\\hline
 $x^3+x+y^3-y+z^2=0$ & $x=-1-y$ & $z^2=3y^2+5y+2$ & $(-1,0)$ \\\hline
 $x^3+x+y^3+y+z^2=0$  & $x=-2-y$ & $z^2=6y^2+12y+10$ & $(-1,2)$ \\\hline
 $x^3+y^3+z^2+z+2=0$ & $x=-2-y$ & $z^2+z=6y^2+12y+6$ & $(-1,0)$ \\ \hline
 $x^3+y^3+z^2+4=0$ & $x=-2-y$ & $z^2=6y^2+12y+4$ & $(-2,2)$  \\ \hline
 $x^3-x+y^3+z^2-2=0$ & $x=-2-y$ & $z^2=6y^2+11y+8$ & $(-7,15)$ \\ \hline
  $x^3-x+y^3+z^2+2=0$ & $x=-2-y$ & $z^2=6y^2+11y+4$ & $(-3,5)$ \\ \hline
 $x^3-x+y^3+z^2+z=0$ & $x=-1-y$ & $z^2+z=3y^2+2y$ & $(0,0)$  \\ \hline
 $x^3+x+y^3+z^2+z=0$ & $x=-1-y$ & $z^2+z=3y^2+4y+2$ & $(0,1)$ \\ 
\hline

\end{tabular}
\captionof{table}{\label{tab:H24Gausssubsol} Equations in Table \ref{tab:H24Gausssub} and the necessary information to prove these equations have infinitely many integer solutions by Proposition \ref{prop:Gaussquad} or \ref{prop:Gaussquad2}.}
\end{center}

\subsection{Exercise 5.8}\label{ex:H21linfinchecksub}
\textbf{\emph{For each of the equations listed in Table \ref{tab:H21linfinchecksub}, use a suitable linear substitution in combination with Proposition \ref{prop:finlincheck} or \ref{prop:finlincheck2} to prove that it has infinitely many integer solutions. }}

	\begin{center}
		\begin{tabular}{ |c|c|c|c|c|c| } 
			\hline
			$H$ & Equation & $H$ & Equation & $H$ & Equation \\ 
			\hline\hline
			$19$ & $x^3+y^2+y-z^2-1=0$ & $21$ & $x^3-x+y^2+y-z^2-1=0$ &  $21$ & $x^3+x+y^2+y-z^2-1=0$  \\ 
			\hline
			$19$ & $x^3+y^2+y-z^2+1=0$ & $21$ & $x^3-x+y^2+y-z^2+1=0$ &  $21$ & $x^3+x+y^2+y-z^2+1=0$  \\ 
			\hline
			$20$ & $x^3+y^2+y-z^2+2=0$ & $21$ & $x^3+y^2+y-z^2-3=0$ &  $21$ & $x^3+xy+y^2-z^2-1=0$ \\ 
			\hline
			$20$ & $x^3-x+y^2+y-z^2=0$ & $21$ & $x^3+y^2+y-z^2+3=0$ &  $21$ & $x^3+x y+y^2-z^2+1=0$ \\ 
			\hline
			$20$ & $x^3+x+y^2+y-z^2=0$ & $21$ & $x^3+y^2+y-z^2+z-1=0$ &   &   \\ 
			\hline
			$20$ & $x^3+y^2+y-z^2+z=0$ & $21$ & $x^3+y^2+y-z^2+z+1=0$ &   &   \\ 
			\hline
		\end{tabular}
		\captionof{table}{\label{tab:H21linfinchecksub} Equations of size $H\leq 21$ having infinitely many integer solutions by Proposition \ref{prop:finlincheck} or \ref{prop:finlincheck2}.}
	\end{center} 

To prove that equations in Table \ref{tab:H21linfinchecksub} have infinitely many solutions, it will be convenient to use an extended version of Proposition \ref{prop:finlincheck}.

\begin{proposition}\label{prop:finlincheck2}[Proposition 5.7 in the book]
	Assume that a Diophantine equation in $n\geq 3$ variables can, possibly after permutation of variables, be written in the form
	\begin{equation}\label{eq:finlincheck2}
		(a x_2 + P(x_3,\dots,x_n)) x_1 = Q(x_2,\dots,x_n) 
	\end{equation}
	where $a$ is an integer and $P,Q$ are polynomials with integer coefficients. If \eqref{eq:finlincheck2} has an integer solution, then it has infinitely many of them.
\end{proposition}

Each equation in Table \ref{tab:H21linfinchecksub} can be reduced, after substitution $y=t+z$, to an equation of the form \eqref{eq:finlincheck2} that has an integer solution, see Table \ref{tab:H21linfinchecksubsol}. By Proposition \ref{prop:finlincheck2} this implies that the original equations have infinitely many integer solutions.

\begin{center}
\begin{tabular}{ |c|c|c|c|c|c| } 
\hline
Equation & Form \eqref{eq:finlincheck2}  & Solution $(z_0,t_0,x_0)$ \\ 
\hline\hline
$x^3+y^2+y-z^2-1=0$ & $z(2t+1)=-x^3 - t - t^2 + 1$ & $(0,0,1)$  \\  \hline
$x^3+y^2+y-z^2+1=0$  & $z(2t+1)=-x^3-t^2-t-1$ & $(-1,0,0)$ \\  \hline
$x^3+y^2+y-z^2+2=0$ &  $z(2t+1)=-x^3-t^2-t-2$ & $(-1,1,-1)$ \\  \hline
$x^3-x+y^2+y-z^2=0$ & $z(2t+1)=-x^3+x-t^2-t$  & $(0,0,1)$ \\  \hline
$x^3+x+y^2+y-z^2=0$ &  $z(2t+1)=-x^3-x-t^2-t$  & $(0,-1,0)$\\  \hline
$x^3+y^2+y-z^2+z=0$ &  $z(2+2t)=-t-t^2-x^3$& $(1,-2,0)$ \\  \hline
$x^3-x+y^2+y-z^2-1=0$ &  $z(2t+1)=-x^3+x-t^2-t+1$ & $(1,0,0)$ \\  \hline
$x^3-x+y^2+y-z^2+1=0$  & $z(2t+1)=-x^3+x-t^2-t-1$ & $(-1,-2,-2)$ \\  \hline
$x^3+y^2+y-z^2-3=0$ &  $z(2t+1)=-x^3 - t^2 - t + 3$ & $(0,1,1)$ \\  \hline
$x^3+y^2+y-z^2+3=0$ & $z(2t+1)=-x^3 - t^2 - t - 3$ & $(-1,-2,-2)$ \\  \hline
$x^3+y^2+y-z^2+z-1=0$ &  $z(2+2t)=-t-t^2-x^3+1$& $(1,-2,1)$   \\ \hline
$x^3+y^2+y-z^2+z+1=0$ &  $z(2+2t)=-t-t^2-x^3-1$& $(-1,-1,-1)$      \\ \hline
$x^3+x+y^2+y-z^2-1=0$ &  $z(2t+1)=-x^3 - x - t^2 - t + 1$ & $(1,0,0)$ \\ 	\hline
$x^3+x+y^2+y-z^2+1=0$ &  $z(2t+1)=-x^3 - x - t^2 - t - 1$ & $(-1,-1,-1)$  \\ \hline
$x^3+xy+y^2-z^2-1=0$ & $z(2t+x)=1-t^2-tx-x^3$ & $(0,0,1)$ \\ \hline
$x^3+x y+y^2-z^2+1=0$ &  $z(2t+x)=-1-t^2-tx-x^3$ & $(0,0,-1)$ \\  \hline

\end{tabular}
\captionof{table}{\label{tab:H21linfinchecksubsol} Equations in Table \ref{tab:H21linfinchecksub} and the necessary information to prove these equations have infinitely many integer solutions by Proposition \ref{prop:finlincheck2}. All equations are reduced to form \eqref{eq:finlincheck2} using $y=t+z$.}
\end{center}

\subsection{Exercise 5.12}\label{ex:H20sumsquares1}
\textbf{\emph{By using either human argument or computer search, investigate for which equations from Table \ref{tab:H20sumsquares} you can find non-constant polynomials $Q_1(u), Q_2(u), Q_3(u)$ satisfying 
$$
P(Q_1(u), \dots, Q_n(u)) = 0 \quad \text{for all} \quad u \in {\mathbb Z},
$$
and conclude that the corresponding equations have infinitely many integer solutions.}}

	\begin{center}
		\begin{tabular}{ |c|c|c|c|c|c| } 
			\hline
			$H$ & Equation & $H$ & Equation & $H$ & Equation \\ 
			\hline\hline
			$17$ & $y^2+z^2=x^3-1$ & $19$ & $y^2+z^2=x^3+x-1$ &  $20$ & $y^2+z^2=x^3-2x$  \\ 
			\hline
			$18$ & $y^2+z^2=x^3-2$ & $19$ & $y^2+z^2=x^3+x+1$ &  $20$ & $y^2+z^2=x^3+2x$  \\ 
			\hline
			$18$ & $y^2+z^2=x^3+2$ & $19$ & $y^2+y+z^2=x^3-1$ &  $20$ & $y^2+z^2=x^3-x^2$ \\ 
			\hline
			$18$ & $y^2+z^2=x^3-x$ & $19$ & $y^2+y+z^2=x^3+1$ &  $20$ & $y^2+z^2=x^3+x^2$ \\ 
			\hline
			$18$ & $y^2+z^2=x^3+x$ & $20$ & $y^2+z^2=x^3-4$ & $20$ & $y^2+y+z^2=x^3-2$ \\ 
			\hline
			$19$ & $y^2+z^2=x^3-3$ & $20$ & $y^2+z^2=x^3-x-2$ & $20$ & $y^2+y+z^2=x^3-x$  \\ 
			\hline
			$19$ & $y^2+z^2=x^3+3$ & $20$ & $y^2+z^2=x^3-x+2$ &  $20$ & $y^2+y+z^2=x^3+x$  \\ 
			\hline
			$19$ & $y^2+z^2=x^3-x-1$ & $20$ & $y^2+z^2=x^3+x-2$ &  $20$ & $y^2+y+z^2+z=x^3$ \\ 
			\hline
			$19$ & $y^2+z^2=x^3-x+1$ & $20$ & $y^2+z^2=x^3+x+2$ &   &  \\ 
			\hline
		\end{tabular}
		\captionof{table}{\label{tab:H20sumsquares} Equations involving sums of two squares of size $H\leq 20$.}
	\end{center} 

For a number of equations in Table \ref{tab:H20sumsquares} we can find either a human argument or use a computer search to find polynomials $Q_1(u), Q_2(u), Q_3(u)$, satisfying 
$$
P(Q_1(u), Q_2(u), Q_3(u))=0, \quad \text{for all} \quad u \in \mathbb{Z}.
$$
This shows that the corresponding equation has infinitely many solutions. 

Let 
$$
S = \{n \in {\mathbb Z} : \,\, n=y^2+z^2 \,\, \text{for some} \,\, y,z \in {\mathbb Z}\}.
$$
Identity
\begin{equation}\label{eq:sumsquaresid}
	(a^2+b^2)(c^2+d^2) = (ac-bd)^2 + (ad+bc)^2
\end{equation}
implies that set $S$ is closed under multiplication, that is
$$
	n \in S \quad \text{and} \quad m \in S \quad \Rightarrow \quad nm \in S.
$$
In particular, if all prime factors of an integer $n$ belong to $S$, then so does $n$. 

For some equations, we can find use a short argument, by hand, to find an infinite family of integer solutions, and the details of these are presented below.

Let us begin with the easy equation
\begin{equation}\label{eq:y2pz2mx3mx}
y^2+z^2=x^3+x.
\end{equation}
We can factorise the right-hand side as $x(x^2+1)$, $x^2+1$ is a sum of two squares, so we need to find when $x$ is a sum of two squares, and then we will have that their product is a sum of two squares. Let $x=u^2+0^2=u^2$ for $u$ arbitrary. 
Then we have 
$$
y^2+z^2=x^3+x=(u^3)^2+u^2,
$$
and hence, we obtain that an infinite family of integer solutions to equation \eqref{eq:y2pz2mx3mx} is
$$
(x,y,z)=(u^2,u^3,u), \quad u \in \mathbb{Z}.
$$

\vspace{10pt}

Another easy equation is
\begin{equation}\label{eq:y2pz2mx3px2}
y^2+z^2=x^3-x^2.
\end{equation}
We can factorise the right-hand side to $x^2(x-1)$, $x^2=x^2+0^2$ is a sum of two squares, so we need to find when $x-1$ is a sum of two squares, and then we will have that their product is a sum of two squares. Let $x=u^2+1$ for $u$ arbitrary, so we have that $x-1=u^2=u^2+0^2$ is a sum of two squares. Then we have 
$$
y^2+z^2=x^3-x^2=(u^3+u)^2+0^2,
$$
and hence, we obtain that an infinite family of integer solutions to equation \eqref{eq:y2pz2mx3px2} is 
$$
(x,y,z)=(u^2+1,u^3+u,0), \quad u \in \mathbb{Z}.
$$

\vspace{10pt}

Another easy equation is
\begin{equation}\label{eq:y2pz2mx3mx2}
y^2+z^2=x^3+x^2.
\end{equation}
We can factorise the right-hand side as $x^2(x+1)$, $x^2=x^2+0^2$ is a sum of two squares, so we need to find when $x+1$ is a sum of two squares, and then we will have that their product is a sum of two squares. Let $x=u^2$ for $u$ arbitrary, so we have that $x+1=u^2+1$ which is a sum of two squares. Then we have 
$$
y^2+z^2=x^3+x^2=(u^3)^2+u^2,
$$
and hence, we obtain that an infinite family of integer solutions to equation \eqref{eq:y2pz2mx3mx2} is 
$$
(x,y,z)=(u^2,u^3,u^2), \quad u \in \mathbb{Z}.
$$

\vspace{10pt}

The next equation we will consider is
\begin{equation}\label{eq:y2pypz2pzmx3}
y^2+y+z^2+z=x^3.
\end{equation}
Multiplying the equation by $4$ and rearranging, we obtain
$$
(2y+1)^2+(2z+1)^2=4x^3+2.
$$
As the left-hand side is a sum of two squares, we must now find when $4x^3+2 \in S$.  This can be factorised as $2(2x^3+1)$, as $2=1^2+1^2$ is a sum of two squares, it is sufficient to find when $2x^3+1 \in S$. For this, it is sufficient to choose $x$ such that $2x^3$ is a perfect square. To ensure this, let $x=2u^2$ for some integer $u$. Then we have
$$
(2y+1)^2+(2z+1)^2=(1^2+1^2)((4u^3)^2+1^2),
$$
and using identity \eqref{eq:sumsquaresid} with $(a,b,c,d)=(4u^3,1,1,1)$ we obtain that an infinite family of integer solutions to equation \eqref{eq:y2pypz2pzmx3} is 
$$
(x,y,z)=(2u^2,2u^3-1,2u^3), \quad u \in \mathbb{Z}.
$$

\vspace{10pt}

The next equation we will consider is
\begin{equation}\label{eq:y2pz2mx3mxp2}
y^2+z^2=x^3+x-2.
\end{equation}
The right-hand side can be factorised as $(x-1)(x^2+x+2)$. Let $x$ be of the form $x=t^2-2$, for $t$ arbitrary, then $x^2+x+2=x^2+t^2$. Then, we need to find when $x-1=t^2-3 \in S$, so we are looking for $t^2-3=a^2+b^2$ for some integers $a,b$. We can rearrange this to $(t-b)(t+b)=a^2+3$, then letting $b=t-1$ we obtain $2t-1=a^2+3$, and then $t=\frac{a^2+4}{2}$. Let $a=2u$ for an arbitrary integer $u$, and we have that $t=2u^2+2$. Then
$$
y^2+z^2=((2u)^2+(2u^2+1)^2)((4u^4+8u^2+2)^2+(2u^2+2)^2),
$$
and using identity \eqref{eq:sumsquaresid} with $(a,b,c,d)=(2u,2u^2+1,4u^4+8u^2+2,2u^2+2)$, we obtain that an infinite family of integer solutions to equation \eqref{eq:y2pz2mx3mxp2} is 
$$
(x,y,z)=(4 u^4+ 8 u^2+2, 8 u^5- 4 u^4+ 16 u^3- 6 u^2+ 4 u-2 , 8 u^6+ 20 u^4+ 4 u^3+ 12 u^2+ 4 u+2), \quad u \in \mathbb{Z}.
$$

\vspace{10pt}

The next equation we will consider is
\begin{equation}\label{eq:y2pz2mx3p2x}
y^2+z^2=x^3-2x.
\end{equation}
The right-hand side of the equation can be factorised as $x(x^2-2)$. Let us first determine when $x^2-2 \in S$. We are looking for $x^2-2=a^2+b^2$ for some integers $a,b$. We can rearrange this to $(x-b)(x+b)=a^2+2$. Let $b=x-1$ and $a=2t+1$, so $x=2(t^2+t+1)$. Now let us select $t$ such that $x \in S$. 
As $2$ is a sum of two squares, we need $t^2+t+1$ to be a sum of two squares. Letting $t=u^2-1$, we obtain
$$
x=2(t^2+t+1)=2((u^2-1)^2+u^2)=(u^2-u-1)^2 + (u^2+u-1)^2,
$$
which is a sum of two squares,
and 
$$
x^3-2x=x(a^2+b^2)=((u^2-u-1)^2 + (u^2+u-1)^2)((2u^2-1)^2+(2u^4-2u^2+1)^2).
$$
Then, using identity \eqref{eq:sumsquaresid} with $(a,b,c,d)=(u^2-u-1,u^2+u-1,2u^2-1,2u^4-2u^2+1)$, we obtain that an infinite family of integer solutions to equation \eqref{eq:y2pz2mx3p2x} is
$$
(x,y,z)=(2(u^4-u^2+1),2(u^6-u^5-u^4+2u^3-u),2(u^6+u^5-3u^4+3u^2-1)), \quad u \in \mathbb{Z}.
$$

Equations $y^2+z^2=x^3-1$,  $y^2+z^2=x^3+2$ and $y^2+z^2=x^3-x$ can be solved using a similar method to the equations already solved in this section. For some other equations from Table \ref{tab:H20sumsquares}, like, for example, $y^2+z^2=x^3-x+1$, we have found an infinite family of solutions by a simple computer search (searching for solutions where $x$ is quadratic and $y,z$ are cubic polynomials in $u$ with small coefficients). 

Table \ref{tab:H20sumsquaressol} presents an infinite family of solutions in parametric form for some of the equations in Table \ref{tab:H20sumsquares}. These families were either found by a computer search or by a human argument like the ones used above. 
Table \ref{tab:H20sumsquaresunsolved} lists the equations from Table \ref{tab:H20sumsquares} for which we did not find an infinite family of solutions in parametric form.

\begin{center}
\begin{tabular}{ |c|c|c|c|c|c| } 
\hline
 Equation & Solution $(x,y,z)$  \\ 
\hline\hline
 $y^2+z^2=x^3-1$ & $(4u^4 +8u^3 +12u^2 +8u+3,$ \\ & $8u^6 +24u^5 +44u^4 +44u^3 +28u^2 +8u+1, $ \\ & $8u^5 +24u^4 +40u^3 +38u^2 +20u+5)$ \\ \hline
 $y^2+z^2=x^3+2$ & $(2u^2,2u^3-1,2u^3+1)$ \\ \hline
 $y^2+z^2=x^3-x$ & $(2u^4 +4u^3 +2u^2 +1,2u^6 +8u^5 +8u^4 +2u^3 +2u^2 +2u,$ \\ & $2u^6 +4u^5 +2u^4 +2u^3 +2u^2)$ \\ \hline
 $y^2+z^2=x^3+x$ & $(u^2,u^3,u)$ \\ \hline
 $y^2+z^2=x^3-x+1$ & $(u^2+3,-u^2-5,u^3+4u)$ \\ \hline
 $y^2+z^2=x^3+x-1$ & $(u^2+2,u^3+u,2u^2+3)$   \\ \hline
 $y^2+y+z^2=x^3+1$ & $(u^2-1,-u^2,u^3-2u)$\\ \hline
 $y^2+z^2=x^3-x+2$ & $(2 u^2+2 u-1,2 u^3+4 u^2-1,2u^3 + 2 u^2 - 2 u - 1)$ \\ \hline
 $y^2+z^2=x^3+x-2$ & $(4u^4+ 8 u^2 +2,8u^5-4u^4+16u^3-6u^2+4u-2,$ \\ & $  
8u^6+20u^4+4u^3+12u^2+4u+2)$ \\ \hline
 $y^2+z^2=x^3+x+2$ & $(u^2-1,u^2,-u^3+2u)$  \\ \hline
 $y^2+z^2=x^3-2x$ & $(2u^4-2u^2+2,2u^6-2u^5-6u^4+6u^2-2,$ \\ & $2u^6+2u^5-2u^4-4u^3+2u)$  \\ \hline
 $y^2+z^2=x^3-x^2$ & $(u^2+1,u^3+u,0)$ \\ \hline
 $y^2+z^2=x^3+x^2$ & $(u^2,u^3,u^2)$  \\ \hline
$y^2+y+z^2=x^3-x$ & $(u^2+1,u^2,u^3+u)$ \\ \hline
 $y^2+y+z^2=x^3+x$ & $(u^2+1,u^2+1,u^3+u)$ \\ \hline
 $y^2+y+z^2+z=x^3$ & $(2u^2,2u^3-1,2u^3)$ \\ \hline
\end{tabular}
\captionof{table}{\label{tab:H20sumsquaressol} An infinite family of solutions for some equations in Table \ref{tab:H20sumsquares}. In all solutions $u$ is an arbitrary integer. }
\end{center} 

\begin{center}
\begin{tabular}{ |c||c||c|c||c|c| } 
\hline
 Equation & Equation  &  Equation \\ 
\hline\hline

 $y^2+z^2=x^3-2$ & $y^2+z^2=x^3-x-1$  &  $y^2+y+z^2=x^3-2$  \\ \hline
 $y^2+y+z^2=x^3-1$ & $y^2+z^2=x^3+x+1$ & $y^2+z^2=x^3+2x$\\ \hline
 $y^2+z^2=x^3+3$ & $y^2+z^2=x^3-4$    &     \\ \hline
 $y^2+z^2=x^3-3$ & $y^2+z^2=x^3-x-2$  &\\ \hline

\end{tabular}
\captionof{table}{\label{tab:H20sumsquaresunsolved} The equations from Table \ref{tab:H20sumsquares} not solved in Section \ref{ex:H20sumsquares1}. }
\end{center} 
	
\subsection{Exercise 5.15}\label{ex:H20sumsquares2}
\textbf{\emph{Solve Problem \ref{prob:fin} for all equations from Table \ref{tab:H20sumsquares} that you have not solved in Section \ref{ex:H20sumsquares1}.}}

The equations from Table \ref{tab:H20sumsquares} for which an infinite family was not found in Section \ref{ex:H20sumsquares1} are those in Table \ref{tab:H20sumsquaresunsolved}. 

In order to prove these equations have infinitely many solutions, we will use the following propositions.
\begin{proposition}\label{prop:sosproduct}[Proposition 5.13 in the book]
	If the product $ab$ of two positive integers $a$ and $b$ is a sum of two squares, and $a,b$ do not share any prime factor of the form $p=4k+3$, then each of $a,b$ is the sum of two squares.
\end{proposition}

\begin{proposition}\label{prop:Gaussquadpos}[Proposition 5.14 in the book]
	Let $a,b,c,d$ be integers such that $ad>0$, $ad$ is not a perfect square, and $b^2-4ac\neq 0$. If equation $ax^2+bx+c=dy^2$ has an integer solution $(x_0,y_0)$ satisfying $x_0>-\frac{b}{2a}$, then it has infinitely many solutions in positive integers $(x,y)$.
\end{proposition}

We remark that the equations $y^2+z^2=x^3-2$, $y^2+z^2=x^3+3$ and $y^2+z^2=x^3-x-1$ from Table \ref{tab:H20sumsquaresunsolved} are not included in this section because they have been solved in Section 5.2.3 of the book.

Let us begin with the equation 
$$
y^2+z^2=x^3-3.
$$
Multiplying $P(x)=x^3-3$ by $x$, we obtain that for any integer $k$
$$
x(x^3-3) = x^4-3x=(x^2-k)^2+(2kx^2-3x-k^2).
$$ 
Let us choose $k$ such that $Q(x)=2kx^2-3x-k^2$ is a perfect square (say $t^2$) infinitely often by Proposition \ref{prop:Gaussquadpos}. A computer search returns that the conditions of this proposition are satisfied for $k=1$. Indeed, in this case $Q(x)=2x^2-3x-1=t^2$, so in Proposition \ref{prop:Gaussquadpos} $a=2$, $b=-3$, $c=-1$, $d=1$, hence $ad=2>0$, $b^2-4ac=9-4(-2) \neq 0$ and a solution satisfying $x_0 > - \frac{-3}{4}=0.75$ is $(x_0,t_0)=(2,1)$. Hence, by Proposition \ref{prop:Gaussquadpos}, there are infinitely many positive integers $x$ such that $2x^2-3x-1$ is a perfect square, which implies that $x(x^3-3)$ is a sum of two squares infinitely often.

Let $a=x$ and ${b=x^3-3}$. If $p$ is a common prime factor of $a$ and $b$ then $p$ is also a prime factor of $a^3-b=3$, so $p=3$. However, if $x$ were divisible by $3$, then we would have $Q(x)=2x^2-3x-1 \equiv 2 \text{ (mod 3)}$, which is impossible for a perfect square.
Therefore, $a$ and $b$ are coprime positive integers. Hence by Proposition \ref{prop:sosproduct}, $ab$ can be the sum of two squares only if $a$ and $b$ are. Hence, $b=x^3-3$ is a sum of two squares infinitely often.

\vspace{10pt}

The next equation we will consider is
\begin{equation}\label{eq:y2pz2mx3mxm1}
y^2+z^2=x^3+x+1.
\end{equation}
Multiplying $P(x)=x^3+x+1$ by $x$, we obtain that for any integer $k$
$$
x(x^3+x+1)=x^4+x^2+x=(x^2-k)^2+((2k+1)x^2+x-k^2).
$$ 
Let us choose $k$ such that $Q(x)=(2k+1)x^2+x-k^2$ is a perfect square (say $t^2$) infinitely often by Proposition \ref{prop:Gaussquadpos}. A computer search returns that the conditions of this proposition are satisfied for $k=2$. Indeed, in this case $Q(x)=5x^2+x-4=t^2$, so in Proposition \ref{prop:Gaussquadpos} $a=5$, $b=1$, $c=-4$, $d=1$, hence ${ad=5>0}$, ${b^2-4ac=1-4(5)(-4) \neq 0}$ and a solution satisfying $x_0 > - \frac{1}{10}=-0.1$ is $(x_0,t_0)=(8,18)$. Hence, by Proposition \ref{prop:Gaussquadpos}, there are infinitely many positive integers $x$ such that $5x^2+x-4$ is a perfect square, which implies that $x(x^3+x+1)$ is a sum of two squares infinitely often.

Let $a=x$ and ${b=x^3+x+1}$. If $p$ is a common prime factor of $a$ and $b$ then $p$ is also a prime factor of $b-a^3-a=1$, a contradiction.  Hence by Proposition \ref{prop:sosproduct}, $ab$ can be the sum of two squares only if $a$ and $b$ are. Hence by Proposition \ref{prop:sosproduct}, $b=x^3+x+1$ is a sum of two squares infinitely often.

\vspace{10pt}

The next equation we will consider is
$$
y^2+y+z^2=x^3-1.
$$
Multiplying this equation by $4$ and rearranging, we obtain
\begin{equation}\label{eq:y2pypz2mx3p1red}
(2y+1)^2+(2z)^2=4x^3-3.
\end{equation}
Multiplying $P(x)=4x^3-3$ by $x$, we obtain that for any integer $k$
$$
x(4x^3-3)=4x^4-3x=(2x^2+k)^2+(-4kx^2-3x-k^2).
$$ 
Let us choose $k$ such that $Q(x)=-4kx^2-3x-k^2$ is a perfect square (say $t^2$) infinitely often by Proposition \ref{prop:Gaussquadpos}. A computer search returns that the conditions of this proposition are satisfied for $k=-2$. Indeed, in this case $Q(x)=8x^2-3x-4=t^2$, so in Proposition \ref{prop:Gaussquadpos} $a=8$, $b=-3$, $c=-4$, $d=1$, hence ${ad=8>0}$, ${b^2-4ac=9-4(8)(-4) \neq 0}$ and a solution satisfying $x_0 > \frac{3}{16}$ is $(x_0,t_0)=(1,1)$. Hence, by Proposition \ref{prop:Gaussquadpos}, there are infinitely many positive integers $x$ such that $8x^2-3x-4$ is a perfect square, which implies that $x(4x^3-3)$ is a sum of two squares infinitely often.

Let $a=x$ and $b=4x^3-3$.  If $p$ is a common prime factor of $a$ and $b$ then $p$ is also a prime factor of $4a^3-b=3$, so $p=3$. However, if $x$ were divisible by $3$, then we would have $8x^2-3x-4 \equiv 2 \text{ (mod 3)}$, which is impossible for a perfect square.
Therefore, $a$ and $b$ are coprime positive integers. Hence by Proposition \ref{prop:sosproduct}, $ab$ can be the sum of two squares only if $a$ and $b$ are. Therefore, $b=4x^3-3$ is a sum of two square infinitely often. Because $b$ is odd, one of these squares must be even, and another one odd, hence equation \eqref{eq:y2pypz2mx3p1red} has infinitely many integer solutions.

\vspace{10pt}

The next equation we will consider is
$$
y^2+z^2=x^3-4.
$$
Multiplying $P(x)=x^3-4$ by $x+2$, we obtain that for any integer $k$
$$
(x+2)(x^3-4)=x^4+2x^3-4x-8=(x^2+x+k)^2+((-2k-1)x^2-4x-2kx-8-k^2).
$$ 
Let us choose $k$ such that $Q(x)=(-2k-1)x^2-4x-2kx-8-k^2$ is a perfect square (say $t^2$) infinitely often by Proposition \ref{prop:Gaussquadpos}. A computer search returns that the conditions of this proposition are satisfied for $k=-6$. Indeed, in this case $Q(x)=11x^2+8x-44=t^2$, so in Proposition \ref{prop:Gaussquadpos} $a=11$, $b=8$, $c=-44$, $d=1$, hence $ad=11>0$, $b^2-4ac=64-4(11)(-44) \neq 0$ and a solution satisfying $x_0 >- \frac{8}{22}$ is $(x_0,t_0)=(2,4)$. Hence, by Proposition \ref{prop:Gaussquadpos}, there are infinitely many positive integers $x$ such that $11x^2+8x-44$ is a perfect square, which implies that $(x+2)(x^3-4)$ is a sum of two squares infinitely often.

Let $a=x+2$ and $b=x^3-4$. If $p$ is a common prime factor of $a$ and $b$ then $p$ is also a prime factor of $a(x^2-2x+4)-b=12$, so $p=2$ or $p=3$. However, if $x+2$ is divisible by $3$, or $x \equiv 1 \text{ (mod 3)}$, then $11x^2+8x-44 \equiv 2 \text{ (mod 3)}$, which is impossible for a perfect square. Hence, $p\neq 3$, and the only option is $p=2$.
In particular, positive integers $a$ and $b$ cannot share prime factors of the form $p=4k+3$. Hence by Proposition \ref{prop:sosproduct}, $ab$ can be the sum of two squares only if $a$ and $b$ are. Therefore, $b=x^3-4$ is a sum of two squares infinitely often.

\vspace{10pt}

The next equation we will consider is
\begin{equation}\label{eq:y2pz2mx3pxp2}
y^2+z^2=x^3-x-2.
\end{equation}
Multiplying $P(x)=x^3-x-2$ by $x$, we obtain that for any integer $k$
$$
x(x^3-x-2)=x^4-x^2-2x=(x^2-k)^2+((2k-1)x^2-2x-k^2).
$$ 
Let us choose $k$ such that $Q(x)=(2k-1)x^2-2x-k^2$ is a perfect square (say $t^2$) infinitely often by Proposition \ref{prop:Gaussquadpos}. A computer search returns that the conditions of this proposition are satisfied for $k=2$. Indeed, in this case $Q(x)=3x^2-2x-4=t^2$, so in Proposition \ref{prop:Gaussquadpos} $a=3$, $b=-2$, $c=-4$, $d=1$, hence ${ad=3>0}$, ${b^2-4ac=4-4(3)(-4) \neq 0}$ and a solution satisfying $x_0 > -\frac{-2}{6}=\frac{1}{3}$ is $(x_0,t_0)=(2,2)$. Hence, by Proposition \ref{prop:Gaussquadpos}, there are infinitely many positive integers $x$ such that $3x^2-2x-4$ is a perfect square, which implies that $(x+2)(x^3+2x)$ is a sum of two squares infinitely often.

Let $a=x$ and $b=x^3-x-2$. If $p$ is a common prime factor of $a$ and $b$ then $p$ is also a prime factor of $a^3-a-b=2$. 
Hence positive integers $a$ and $b$ cannot share factors of the form $p=4k+3$. Hence by Proposition \ref{prop:sosproduct}, $ab$ can be the sum of two squares only if $a$ and $b$ are. Therefore, $b=x^3-x-2$ is a sum of two squares infinitely often. 

\vspace{10pt}

The next equation we will consider is
$$
y^2+z^2=x^3+2x.
$$
Multiplying $P(x)=x^3+2x$ by $x+2$, we obtain that for any integer $k$
$$
(x+2)(x^3+2x) = x^4+2x^3+2x^2+4x=(x^2+x-k)^2+((2k+1)x^2+(2k+4)x-k^2).
$$ 
Let us choose $k$ such that $Q(x)=(2k+1)x^2+(2k+4)x-k^2$ is a perfect square (say $t^2$) infinitely often by Proposition \ref{prop:Gaussquadpos}. A computer search returns that the conditions of this proposition are satisfied for $k=2$. However, in this case $a=x+2$ and $b=x^3+2x$ can be both divisible by $3$, and we cannot apply Proposition \ref{prop:sosproduct}. The same problem occurs with any $k$ equal to $2$ modulo $3$.

The smallest $k$ not equal to $2$ modulo $3$ for which we can apply Proposition \ref{prop:Gaussquadpos} is $k=48$. In this case $Q(x)=97x^2+100x-2304=t^2$, so in Proposition \ref{prop:Gaussquadpos} $a=97$, $b=100$, $c=-2304$, $d=1$, hence ${ad=97>0}$, ${b^2-4ac \neq 0}$ and a solution satisfying $x_0 >- \frac{-100}{194}$ is $(x_0,t_0)=(7232, 71232)$. Hence, by Proposition \ref{prop:Gaussquadpos}, there are infinitely many positive integers $x$ such that $Q(x)$ is a perfect square, which implies that $(x+2)(x^3+2x)$ is a sum of two squares infinitely often.

Let $a=x+2$ and $b=x^3+2x$. If an odd prime $p$ is a common factor of $a$ and $b$, then $p$ would also be a factor of $b-a^3+6a^2-14a=-12$, hence $p=3$. But if $a$ is divisible by $3$, then $Q(x)=97x^2+100x-2304$ would be $2$ modulo $3$ which is impossible for a perfect square. Therefore $a$ and $b$ do not share any odd prime factors. Hence, by Proposition \ref{prop:sosproduct}, $ab$ can be the sum of two squares only if $a$ and $b$ are. Therefore, $b=x^3+2x$ is a sum of two squares infinitely often.

\vspace{10pt}

The next equation we will consider is
\begin{equation}\label{eq:y2pypz2mx3p2}
	y^2+y+z^2=x^3-2.
\end{equation}
Multiplying this equation by $4$ and rearranging, we obtain
\begin{equation}\label{eq:y2pypz2mx3p2red}
(2y+1)^2+(2z)^2=4x^3-7.
\end{equation}
Multiplying $P(x)=4x^3-7$ by $x$, we obtain that for any integer $k$
$$
x(4x^3-7) = 4x^4-7x=(2x^2+k)^2+(-4kx^2-7x-k^2).
$$ 
Let us choose $k$ such that $-4kx^2-7x-k^2$ is a perfect square (say $t^2$) infinitely often by Proposition \ref{prop:Gaussquadpos}. A computer search returns that the conditions of this proposition are satisfied for $k=-3$. Indeed, in this case $-4kx^2-7x-k^2=12 x^2-7 x -9=t^2$, so in Proposition \ref{prop:Gaussquadpos} $a=12$, $b=-7$, $c=-9$, $d=1$, hence ${ad=12>0}$, ${b^2-4ac=49-4\cdot 12(-9) \neq 0}$ and a solution satisfying $x_0 >- \frac{-7}{24}$ is $(x_0,t_0)=(2,5)$. Hence, by Proposition \ref{prop:Gaussquadpos}, there are infinitely many positive integers $x$ such that $12 x^2-7 x -9$ is a perfect square, which implies that $x(4x^3-7)$ is a sum of two squares infinitely often.

Let $a=x$ and $b=4x^3-7$. If a prime $p$ is a common factor of $a$ and $b$, then $p$ would also be a factor of $4a^3-b=7$, hence $p=7$. But if $x$ is divisible by $7$, then $12 x^2-7 x -9$ would be $5$ modulo $7$, which is impossible for a perfect square. 
This contradiction shows that $a$ and $b$ are coprime positive integers. Hence by Proposition \ref{prop:sosproduct}, $ab$ can be the sum of two squares only if $a$ and $b$ are. Hence, $b=4x^3-7$ is a sum of two squares infinitely often. Because $b$ is odd, one of these squares must be even and another one odd, hence equation \eqref{eq:y2pypz2mx3p2red} has infinitely many integer solutions.

\subsection{Exercise 5.16}\label{ex:x2ymypy2pz2p4}
\textbf{\emph{Prove that the equation
\begin{equation}\label{eq:y2mx2ypypz2p4a}
	y^2-x^2y+y+z^2+4=0
\end{equation}
has no integer solutions.}}

Let us first consider the case where $x$ is odd. Modulo $4$ analysis determines that if $x$ is odd, both $y$ and $z$ must be even. Let us make the substitutions $x=2k+1$, $y=2m$ and $z=2n$ for some integers $k,m,n$. Then the equation is reduced to 
$$
m^2-2k^2m-2km+n^2+1=0.
$$
This equation has no solutions modulo $4$, and therefore there are no integer solutions to equation \eqref{eq:y2mx2ypypz2p4a} where $x$ is odd.

We can now assume that $x$ is even. Multiplying equation \eqref{eq:y2mx2ypypz2p4a} by $4$ and rearranging, we obtain
$$
(2y-x^2+1)^2+(2z)^2=(x^2-5)(x^2+3).
$$
We now need to determine for what values of $x$ the right-hand side is the sum of two squares. It is easy to check that this is not the case for $|x|\leq 2$, so we may assume that $|x|\geq 3$. Then $a=x^2-5$ and $b=x^2+3$ are positive integers. If $p$ is any common prime factor of $a$ and $b$, then $p$ is also a factor of $b-a=8$, so $p=2$, but this is impossible for even $x$. Therefore, $a$ and $b$ are coprime positive integers. Thus, by Proposition \ref{prop:sosproduct}, $ab$ can be the sum of two squares only if both $a$ and $b$ are. However, because $x$ is even, we have $b=x^2+3 \equiv 3$ modulo $4$, and therefore it cannot be a sum of two squares.

\subsection{Exercise 5.24}\label{ex:quadres}
\textbf{\emph{Compute $\left(\frac{a}{p}\right)$ for $a=\pm 4, \pm 5, \pm 6$ and $\pm 7$ and all odd primes $p$ coprime to $a$.}}

Let us recall the following property of the Legendre symbol from Section 5.3.2 of the book. We have 
for any integers $a,b$
\begin{equation}\label{eq:legsymprod}
	\left(\frac{a b}{p}\right) = \left(\frac{a}{p}\right) \left(\frac{b}{p}\right).
\end{equation}

The following theorem will also be useful in computing the Legendre symbol in this exercise.
\begin{theorem}[Law of quadratic reciprocity]\label{th:lawquadre}[Theorem 5.23 in the book]
	 Let $p$ and $q$ be different odd primes. Then
	\begin{equation}\label{eq:repos}
		\left(\frac{p}{q}\right) \left(\frac{q}{p}\right) = (-1)^{\frac{p-1}{2}\frac{q-1}{2}}.
	\end{equation}  
\end{theorem}

The values of $\left(\frac{a}{p}\right)$ for $|a| \leq 3$ are computed in Section 5.3.2 of the book, and the answers are
$$
\left(\frac{1}{p}\right)=1,
$$
\begin{equation}\label{leg:m1}
\left(\frac{-1}{p}\right)= \begin{cases}
	1, & \text{if $p \equiv 1(\text{mod}\, 4)$,}\\
	-1, & \text{if $p \equiv 3(\text{mod}\, 4)$,}
\end{cases}
\end{equation}
\begin{equation}\label{leg:p2}
	\left(\frac{2}{p}\right)= \begin{cases}
		1, & \text{if $p \equiv \pm 1 \, (\text{mod}\, 8)$,}\\
		-1, & \text{if $p \equiv \pm 3 \, (\text{mod}\, 8)$,}
	\end{cases}
\end{equation}
$$
	\left(\frac{-2}{p}\right)= \begin{cases}
	1, & \text{if $p \equiv 1 \,\, \text{or} \,\, 3 \, (\text{mod}\, 8)$,}\\
-1, & \text{if $p \equiv -1 \,\, \text{or} \,\, -3 \, (\text{mod}\, 8)$,}
	\end{cases}
$$
\begin{equation}\label{leg:p3}
	\left(\frac{3}{p}\right)= \begin{cases}
	1, & \text{if $p \equiv \pm 1 \, (\text{mod}\, 12)$,}\\
	-1, & \text{if $p \equiv \pm 5 \, (\text{mod}\, 12)$,}
	\end{cases}
\end{equation}
\begin{equation}\label{leg:m3}
	\left(\frac{-3}{p}\right)= \begin{cases}
	1, & \text{if $p \equiv 1 \, (\text{mod}\, 3)$,}\\
-1, & \text{if $p \equiv -1 \, (\text{mod}\, 3)$.}
\end{cases}
\end{equation}

The symbols $\left(\frac{\pm 4}{p}\right)$ are easy to compute. For odd $p$, \eqref{eq:legsymprod} implies that
\begin{equation}\label{leg:p4}
\left(\frac{4}{p}\right)=\left(\frac{2}{p}\right)^2=1,
\end{equation}
and
$$
\left(\frac{-4}{p}\right)=\left(\frac{-1}{p}\right) \left(\frac{4}{p}\right)=\left(\frac{-1}{p}\right)= \begin{cases}
	1, & \text{if $p \equiv 1 \, (\text{mod}\, 4)$,}\\
	-1, & \text{if $p \equiv 3 \, (\text{mod}\, 4)$,}
\end{cases}
$$
where the first equality follows from \eqref{leg:p4} and the final equality follows from \eqref{leg:m1}.

For odd $p\neq 5$, \eqref{eq:repos} implies that $\left(\frac{p}{5}\right) \left(\frac{5}{p}\right) = (-1)^{\frac{p-1}{2}\frac{5-1}{2}}=1$. Obviously, $\left(\frac{p}{5}\right)=1$ if and only if $p \equiv \pm 1 (\text{mod } 5)$. From this, it is easy to deduce that
\begin{equation}\label{leg:p5}
	\left(\frac{5}{p}\right) = 	\begin{cases}
		1, & \text{if $p \equiv \pm 1 \, (\text{mod}\, 5)$,}\\
		-1, & \text{if $p \equiv \pm 2 \, (\text{mod}\, 5)$.}
	\end{cases}
\end{equation}

For odd $p\neq 5$, \eqref{eq:legsymprod} implies that $\left(\frac{-5}{p}\right) =\left(\frac{5}{p}\right) \left(\frac{-1}{p}\right) $. From this, \eqref{leg:m1} and \eqref{leg:p5}, we have $\left(\frac{-5}{p}\right)=1$ if and only if (i) $p \equiv \pm 1 (\text{mod } 5)$ and $p \equiv  1 (\text{mod } 4)$ or (ii) $p \equiv \pm 2 (\text{mod } 5)$ and $p \equiv  -1 (\text{mod } 4)$ . From this, it is easy to deduce that
\begin{equation}\label{leg:m5}
	\left(\frac{-5}{p}\right) = 	\begin{cases}
		1, & \text{if $p \equiv 1,3,7,9 \, (\text{mod}\, 20)$,}\\
		-1, & \text{if $p \equiv 11,13,17,19 \, (\text{mod}\, 20)$.}
	\end{cases}
\end{equation}

For odd $p\neq 3$, \eqref{eq:legsymprod} implies that $\left(\frac{6}{p}\right) =\left(\frac{2}{p}\right) \left(\frac{3}{p}\right) $. From this, \eqref{leg:p2} and \eqref{leg:p3}, we have $\left(\frac{6}{p}\right)=1$ if and only if (i) $p \equiv \pm 1 (\text{mod } 8)$ and $p \equiv  \pm 1 (\text{mod } 12)$ or (ii) $p \equiv \pm 3 (\text{mod } 8)$ and $p \equiv  \pm 5 (\text{mod } 12)$ . From this, it is easy to deduce that
\begin{equation}\label{leg:p6}
	\left(\frac{6}{p}\right) = 	\begin{cases}
		1, & \text{if $p \equiv \pm 1, \pm 5 \, (\text{mod}\, 24)$,}\\
		-1, & \text{if $p \equiv \pm 7, \pm 11 \, (\text{mod}\, 24)$.}
	\end{cases}
\end{equation}

For odd $p\neq 3$, \eqref{eq:legsymprod} implies that $\left(\frac{-6}{p}\right) =\left(\frac{-1}{p}\right) \left(\frac{6}{p}\right) $. From this, \eqref{leg:m1} and \eqref{leg:p6}, we have $\left(\frac{-6}{p}\right)=1$ if and only if (i) $p \equiv  1 (\text{mod } 4)$ and $p \equiv  \pm 1,\pm 5 (\text{mod } 24)$ or (ii) $p \equiv -1 (\text{mod } 4)$ and $p \equiv  \pm 7, \pm 11 (\text{mod } 24)$. From this, it is easy to deduce that
$$
	\left(\frac{-6}{p}\right) = 	\begin{cases}
		1, & \text{if $p \equiv  1, 5,7,11 \, (\text{mod}\, 24)$,}\\
		-1, & \text{if $p \equiv 13,17,19,23 \, (\text{mod}\, 24)$.}
	\end{cases}
$$

For odd $p\neq 7$, \eqref{eq:repos} implies that $\left(\frac{p}{7}\right) \left(\frac{7}{p}\right) = (-1)^{\frac{p-1}{2}\frac{7-1}{2}}=(-1)^{\frac{p-1}{2}}=\left(\frac{-1}{p}\right)$. Obviously, $\left(\frac{p}{7}\right)=1$ if and only if $p \equiv 1,2,4 \, (\text{mod } 7)$. From this and \eqref{leg:m1}, it is easy to deduce that
\begin{equation}\label{leg:p7}
	\left(\frac{7}{p}\right) = 	\begin{cases}
		1, & \text{if $p \equiv  \pm 1,\pm 3, \pm 9 \, (\text{mod}\, 28)$,}\\
		-1, & \text{if $p \equiv  \pm 5,\pm 11,\pm 13 \, (\text{mod}\, 28)$.}
	\end{cases}
\end{equation}

For odd $p\neq 7$, \eqref{eq:legsymprod} implies that $\left(\frac{-7}{p}\right) =\left(\frac{-1}{p}\right) \left(\frac{7}{p}\right) $. 
We previously had $\left(\frac{p}{7}\right) \left(\frac{7}{p}\right)  =(-1)^{\frac{p-1}{2}}=\left(\frac{-1}{p}\right)$, hence  
$\left(\frac{-7}{p}\right) =\left(\frac{p}{7}\right) \left(\frac{7}{p}\right)\left(\frac{7}{p}\right)=\left(\frac{p}{7}\right)$. From this, it is easy to deduce that
$$
	\left(\frac{-7}{p}\right) = 	\begin{cases}
		1, & \text{if $p \equiv  1,2,4 \, (\text{mod}\, 7)$,}\\
		-1, & \text{if $p \equiv  3,5,6 \, (\text{mod}\, 7)$.}
	\end{cases}
$$

\begin{center}

	\captionof{table}{\label{tab:legendresymbol} Legendre symbol for $1 \leq |a| \leq 7$, and prime $p$ coprime to $a$. }
\end{center}

\subsection{Exercise 5.30}\label{ex:quadtriv}
\textbf{\emph{For the equations listed in Table \ref{tab:H25}, find all its integer solutions (if any exist) and, produce a detailed human proof that they have no other integer solutions.}}

\begin{center}
		%
		\captionof{table}{\label{tab:H25} Equations of size $H\leq 25$ solvable by Algorithm 5.29 in the book.}
	\end{center} 

	We will use the following results.
	
	\begin{proposition}\label{prop:z2p1div}[Proposition 5.17 in the book]
		For any integer $z$, all positive divisors of $z^2+1$ must be $1$ or $2$ modulo $4$. In particular, if $p$ is a prime divisor of $z^2+1$, then either $p=2$ or $p \equiv 1(\text{mod}\, 4)$.
	\end{proposition}
	
	\begin{proposition}\label{prop:z2m2div}[Proposition 5.18 in the book]
		For any integer $z$, all divisors of $z^2-2$ must be $\pm 1$ or $\pm 2$ modulo $8$. In particular, if $p$ is a prime divisor of $z^2-2$, then either $p=2$ or $p \equiv \pm 1(\text{mod}\, 8)$.
	\end{proposition}
	
	\begin{proposition}\label{prop:z2p2div}[Proposition 5.19 in the book]
		For any integer $z$, all positive divisors of $z^2+2$ must be $1$, $3$ or $\pm 2$ modulo $8$. In particular, if $p$ is a prime divisor of $z^2+2$, then either $p=2$ or $p \equiv 1(\text{mod}\, 8)$ or $p \equiv 3(\text{mod}\, 8)$.
	\end{proposition}

\begin{proposition}\label{cor:quadform}[Proposition 5.26 in the book]
	Let $x,y$ be integers and let $p$ be an odd prime that is not a common divisor of $x$ and $y$. Then:
	\begin{center}
		\begin{tabular}{ |c|c|c|c|c|c|c| } 
			\hline
			& if $p$ is a prime factor of & then \\ 
			\hline
			(a) & $x^2+y^2$ & $p \equiv 1 (\text{mod } 4)$ \\ 
			\hline
			(b) & $x^2+2y^2$ & $p \equiv 1 \text{ or } 3 (\text{mod } 8)$ \\ 
			\hline
			(c) & $x^2-2y^2$ & $p \equiv 1 \text{ or } 7 (\text{mod } 8)$ \\
			\hline
			(d) & $x^2+3y^2$ & $p=3$ or $p \equiv 1 (\text{mod } 3)$ \\  
			\hline
			(e) & $x^2-3y^2$ & $p=3$ or $p \equiv 1 \text{ or } 11 (\text{mod } 12)$ \\ 
			\hline
		\end{tabular}
	\end{center}
\end{proposition}

The equation 
$$
x^2y+4y-z^2-2=0
$$
is solved in Section 5.3.1 of the book, and it has no integer solutions.

The equation 
$$
y^2-x^2y+y+z^2+4=0
$$
is solved in Section 5.3.3 of the book, and it has no integer solutions.

The equation 
$$
y^2-x^2y+4y+z^2+1=0
$$
is solved in Section 5.3.3 of the book, and its integer solutions are
$$
(x,y,z)=(\pm 1,-2,\pm 1) \quad \text{and} \quad (\pm 1,-1,\pm 1).
$$

The equation 
$$
y^2-x^2 y+z^2+1=0
$$
can be rewritten as $(x^2-y) y=z^2+1$. Then Proposition \ref{prop:z2p1div} implies that the right-hand side can only have positive divisors equal to $1,2$ modulo $4$. Because $z^2+1>0$, factors $x^2-y$ and $y$ must have the same sign. If both are negative, then their sum is also negative, however $x^2-y+y=x^2 \geq 0$, which is a contradiction. So both $x^2-y$ and $y$ must be positive, and therefore must be equal to $1$ or $2$ modulo $4$. If they are both even, then $z^2+1$ would be divisible by $4$ which is a contradiction. Hence either $x^2-y$ or $y$ is odd, but then $x^2=(x^2-y)+y=1+1=2$ or $1+2=3$ modulo $4$, in both cases a contradiction. 

\vspace{10pt}

The equation 
$$
x^2 y+3y-z^2-1=0
$$
can be rewritten as $(x^2+3) y=z^2+1$. Then Proposition \ref{prop:z2p1div} implies that the right-hand side can only have positive divisors equal to $1,2$ modulo $4$. Because $x^2+3>0$, it must be $1$ or $2$ modulo $4$, but then $x^2=2$ or $3$ modulo $4$, which is a contradiction. 

\vspace{10pt}

The equation 
$$
x^2 y+2y-2z^2+1=0
$$
can be rewritten as $(x^2+2) y=2z^2-1$. 
Hence, $x^2+2$ is a divisor of $2z^2-1$, which implies that $x^2+2$ is odd, and is also a divisor of $2(2z^2-1)=(2z)^2-2$, hence, by Proposition \ref{prop:z2m2div}, it must be $\pm 1$ modulo $8$. But then $x^2$ is equal to $5$ or $7$ modulo $8$, which is a contradiction.

\vspace{10pt}

Let us now consider the equation 
$$
x^2 y+4y-z^2+2=0.
$$
If $x$ is even, then $z^2$ is $2$ modulo $4$, which is a contradiction. Hence, $x$ is odd. The equation can be rewritten as $(x^2+4) y=z^2-2$. Then Proposition \ref{prop:z2m2div} implies that the right-hand side can only have odd divisors equal to $\pm 1$ modulo $8$. Because $x^2+4$ is odd, $x^2+4=\pm 1$ modulo $8$ but then $x^2=\pm 3$ modulo $8$, which is a contradiction. 

\vspace{10pt}

Let us now consider the equation
$$
x^2 y-4y-z^2+2=0.
$$
If $x$ is even, then $z^2$ is $2$ modulo $4$, which is a contradiction. Hence, $x$ is odd. The equation can be rewritten as $(x^2-4) y=z^2-2$. Then Proposition \ref{prop:z2m2div} implies that the right-hand side can only have odd divisors equal to $\pm 1$ modulo $8$. Because $x^2-4$ is odd, $x^2-4=\pm 1$ modulo $8$ but then $x^2=\pm 3$ modulo $8$, which is, again, a contradiction. 

\vspace{10pt}

The equation 
$$
x^2 y+4y-2z^2-1=0.
$$
Modulo $2$ analysis shows that $x$ is odd. After multiplication by $2$, the equation can be rewritten as $2(x^2+4)y=(2z)^2+2$. Because $x$ is odd, positive integer $x^2+4$ is equal to $5$ modulo $8$, but Proposition \ref{prop:z2p2div} implies that $(2z)^2+2$ cannot have such positive divisors.

\vspace{10pt}

The equation
$$
x^2 y+4y-2z^2+1=0
$$
can be rewritten as $(x^2+4) y=2z^2-1$. Hence, $x^2+4$ is a divisor of $2z^2-1$, which implies that $x^2+4$ is odd, and is also a divisor of $2(2z^2-1)=(2z)^2-2$, hence, by Proposition \ref{prop:z2m2div}, it must be $\pm 1$ modulo $8$.
But then $x^2+4=\pm 1$ modulo $8$ so $x^2=\pm 3$ modulo $8$, which is a contradiction. 

\vspace{10pt}

The equation 
$$
x^2 y-4y-2z^2+1=0
$$
can be rewritten as $(x^2-4) y=2z^2-1$. Hence, $x^2-4$ is a divisor of $2z^2-1$, which implies that $x^2-4$ is odd, and is also a divisor of $2(2z^2-1)=(2z)^2-2$, hence, by Proposition \ref{prop:z2m2div}, it must be $\pm 1$ modulo $8$.
But then $x^2-4=\pm 1$ modulo $8$ so $x^2=\pm  3$ modulo $8$, which is a contradiction.

\vspace{10pt}

The equation 
$$
x^2 y+2y-2z^2-2z-1=0
$$
can be rewritten as $(x^2+2) y=2z^2+2z+1=z^2+(z+1)^2$. Because the right-hand side is odd, and $z$ and $z+1$ are coprime, Proposition \ref{cor:quadform} (a) implies that all its prime factors are equal to $1$ modulo $4$. Hence, all prime factors of $x^2+2>0$ are $1$ modulo $4$, but then $x^2+2$ is $1$ modulo $4$, hence $x^2$ is $3$ modulo $4$, which is a contradiction. 

\vspace{10pt}

The equation 
$$
y^2-x^2 y-4y+z^2+1=0
$$
can be rewritten as $(x^2+4-y) y=z^2+1$. Then Proposition \ref{prop:z2p1div} implies that the right-hand side can only have positive divisors equal to $1$ or $2$ modulo $4$. Because $z^2+1>0$, the factors $x^2+4-y$ and $y$ must have the same sign. If they are both negative, then their sum must also be negative, however, $x^2+4-y+y=x^2+4>0$, which is a contradiction. Hence, the factors must both be positive. If the factors $x^2+4-y$ and $y$ are both even, then $z^2+1$ would be divisible by $4$ which is a contradiction. Hence either $t=x^2+4-y$ or $y$ is odd, but then $y+t=x^2+4$ is either $1+1=2$ or $1+2=3$ modulo $4$, which is, again, a contradiction. 

\vspace{10pt}

The equation 
$$
x^3+x y^2+2x-z^2-1=0
$$
can be rewritten as $x(x^2+y^2+2)=z^2+1$. Then Proposition \ref{prop:z2p1div} implies that the right-hand side can only have positive divisors equal to $1$ or $2$ modulo $4$. Because $z^2+1>0$, the factors $x^2+y^2+2$ and $x$ must have the same sign. Because $x^2+y^2+2>0$, $x$ must be positive. If they are both even, then $z^2+1$ would be divisible by $4$, which is a contradiction. Hence either $x^2+y^2+2$ or $x$ is odd, but then $y^2=(x^2+y^2+2)-x^2-2$ is equal to either $1-1^2-2=2$, $2-1^2-2=3$ or $1-2^2-2=3$ modulo $4$, so, in all cases we obtain a contradiction. 

\vspace{10pt}

The equation 
\begin{equation}\label{eq:x3pxy2m2xmz2m1}
x^3+x y^2-2x-z^2-1=0
\end{equation}
can be rewritten as $x(x^2+y^2-2)=z^2+1$. Then Proposition \ref{prop:z2p1div} implies that the right-hand side can only have positive divisors equal to $1$ or $2$ modulo $4$. Because $z^2+1>0$, both $x^2+y^2-2$ and $x$ must have the same sign. If $x^2+y^2-2<0$ and $x<0$ then $x=-1$, and we obtain the integer solution $(x,y,z)=(-1,0,0)$ to \eqref{eq:x3pxy2m2xmz2m1}. Let us now consider the case where both $x^2+y^2-2$ and $x$ are positive, and therefore equal to $1$ or $2$ modulo $4$. If these numbers are both even, then $z^2+1$ would be divisible by $4$ which is a contradiction. Hence either $x^2+y^2-2$ or $x$ is odd. But then $y^2=(x^2+y^2-2)-x^2+2$ is equal to either $1-1^2+2=2$, $2-1^2+2=3$ or $1-2^2+2=3$ modulo $4$, so, in all cases we obtain a contradiction. Therefore, the only integer solution to equation \eqref{eq:x3pxy2m2xmz2m1} is
$$
(x,y,z)=(-1,0,0).
$$

\subsection{Exercise 5.37}\label{ex:ay2pbyzpcz2mPx}
\textbf{\emph{Prove that all equations listed in Table \ref{tab:H24quadforms} have infinitely many integer solutions.}}
	\begin{center}
		\begin{tabular}{ |c|c|c|c|c|c| } 
			\hline
			$H$ & Equation &  $H$ & Equation & $H$ & Equation \\ 
			\hline\hline
			$22$ & $y^2+yz+z^2=x^3-2$ &  $23$ & $2y^2+y+z^2=x^3-1$ & $24$ & $2y^2+y+z^2=x^3-2$ \\ 
			\hline
			$22$ & $y^2+yz+z^2=x^3-x$ &  $23$ & $y^2+yz+z^2=x^3-3$ & $24$ & $y^2+yz+z^2=x^3-x-2$ \\ 
			\hline
			$22$ & $y^2+yz+z^2=x^3+x$ &  $24$ & $2y^2+z^2=x^3-x-2$ & $24$ & $y^2+yz+z^2=x^3+x-2$ \\ 
			\hline
			$23$ & $2y^2+z^2=x^3-x-1$ &  $24$ & $2y^2+z^2=x^3-x+2$ & $24$ & $y^2+yz+z^2=x^3+2x$ \\ 
			\hline
			$23$ & $2y^2+z^2=x^3+x+1$ &  $24$ & $2y^2+z^2=x^3+x-2$ & $24$ & $y^2+yz+z^2+y=x^3-2$ \\ 
			\hline
			$23$ & $2y^2+z^2=x^3-3$ &  $24$ & $2y^2+z^2+z=x^3-2$ & $24$ & $y^2+yz+z^2+y=x^3+2$ \\ 
			\hline			
		\end{tabular}
		\captionof{table}{\label{tab:H24quadforms} Equations of the form \eqref{eq:ay2pbyzpcz2mPx} with $(a,b,c)\neq (1,0,1)$ of size $H\leq 24$.}
	\end{center} 

The equations in Table \ref{tab:H24quadforms} are of the form 
\begin{equation}\label{eq:ay2pbyzpcz2mPx}
ay^2+byz+cz^2=P(x)
\end{equation}
with $(a,b,c)\neq (1,0,1)$. 

We will first recall some notation and necessary propositions from Section 5.4.1 of the book.
Let us define the set $S_2$ as
\begin{equation}\label{S2_notation}
S_2 = \{n \in {\mathbb Z} : \,\, n=2y^2+z^2 \,\, \text{for some} \,\, y,z \in {\mathbb Z}\},
\end{equation}
and the set $S_3$ as
\begin{equation}\label{S3_notation}
S_3 = \{n \in {\mathbb Z} : \,\, n=3y^2+z^2 \,\, \text{for some} \,\, y,z \in {\mathbb Z}\}= \{n \in {\mathbb Z} : \,\, n=y^2+yz+z^2 \,\, \text{for some} \,\, y,z \in {\mathbb Z}\}.
\end{equation}

\begin{proposition}\label{prop:2y2pz2product}[Proposition 5.32 in the book]
	If the product $ab$ of two positive integers $a$ and $b$ belongs to $S_2$, and $a,b$ do not share any prime factor equal to $5$ or $7$ modulo $8$, then $a \in S_2$ and $b \in S_2$.
\end{proposition}

\begin{proposition}\label{prop:3y2pz2product}[Proposition 5.35 in the book]
	If the product $ab$ of two positive integers $a$ and $b$ belongs to $S_3$, and $a,b$ do not share any prime factor of the form $p=3k+2$, then $a \in S_3$ and $b \in S_3$.
\end{proposition}

Note that the sets $S_2$ and $S_3$ are closed under multiplication as implied by the identity 
$$
	(ka^2 + b^2)(kc^2+d^2) = k(ad + bc)^2 + (bd-kac)^2
$$
with $k=2$ and $k=3$, respectively.
	
Let us begin with the equation 
\begin{equation}\label{eq:y2pyzpz2mx3px}
y^2+yz+z^2=x^3-x.
\end{equation}
Multiplying $P(x)=x^3-x$ by $x+2$, we obtain that for any integer $k$,
$$
(x+2)(x^3-x)=(x^2+x-k)^2+((2k-2)x^2+(2k-2)x-k^2).
$$ 
Let us try to select an integer $k$ such that the equation $(2k-2)x^2+(2k-2)x-k^2=3t^2$ has infinitely many integer solutions with $x$ being positive and odd. With $x=2u-1$, the equation reduces to
$$
8ku^2-8u^2-4ku+4u-k^2=3t^2.
$$ 
A computer search returns that for $k=2$, the resulting equation $8u^2-4u-4=3t^2$ has solution $(u,t)=(4,6)$. Because we also have $8 \cdot 3=24 >0$, $24$ is not a perfect square, $(-4)^2-4 \cdot 8 (-4) \neq 0$, and $4 > \frac{1}{4}$, Proposition \ref{prop:Gaussquadpos} implies that it has infinitely many solutions in positive integers. Hence, for infinitely many odd positive integers $x$ we have $(x+2)(x^3-x)=(x^2+x-2)^2+3t^2 \in S_3$, where $S_3$ is defined in \eqref{S3_notation}. If $p$ is a common prime factor of $a=x+2$ and $b=x^3-x$ then $p$ is odd, and $p$ is also a prime factor of $a(x^2-2x+3)-b=6$. Hence $p=3$, 
which is not of the form $p=3k+2$, and Proposition \ref{prop:3y2pz2product} implies that $b=x^3-x \in S_3$. Therefore, equation  \eqref{eq:y2pyzpz2mx3px} has infinitely many integer solutions.

\vspace{10pt}

The next equation we consider is
\begin{equation}\label{eq:y2pyzpz2mx3mx}
y^2+yz+z^2=x^3+x.
\end{equation}
Modulo $4$ analysis shows that $y$ and $z$ must be even, and $x$ is divisible by $4$. Substituting $y=2Y$, $z=2Z$, $x=4X$ and cancelling $4$, we obtain the equation
$$
Y^2+YZ+Z^2=16 X^3+X.
$$
Multiplying $P(X)=16 X^3+X$ by $X-2$, we obtain that for any integer $k$,
$$
(X-2)(16 X^3+X)=(4 X^2-4 X-k)^2+((8k-15)X^2+(-8k-2)X-k^2).
$$ 
Let us try to select an integer $k$ such that the equation $(8k-15)X^2+(-8k-2)X-k^2=3t^2$ has infinitely many integer solutions with $X$ being positive and odd. With $X=2u-1$, the equation reduces to
$$
-13 + 16 k - k^2 + 56 u - 48 k u - 60 u^2 + 32 k u^2 =3t^2.
$$ 
A computer search returns that for $k=6$, the resulting equation $132 u^2 - 232 u + 47=3t^2$ has solution $(u,t)=(5,27)$. Because we also have $3 \cdot 132=396 >0$, $396$ is not a perfect square, $(-232)^2-4 \cdot 132 \cdot 47 = 29008 \neq 0$, and $5 > \frac{232}{2\cdot 132}$, Proposition \ref{prop:Gaussquadpos} implies that it has infinitely many solutions in positive integers. Hence, for infinitely many odd positive integers $X$ we have $(X-2)(16 X^3+X)=(4 X^2-4 X-6)^2+3t^2 \in S_3$, where $S_3$ is defined in \eqref{S3_notation}. If $p$ is a common prime factor of $a=X-2$ and $b=16 X^3+X$ then $p$ is odd, and $p$ is also a prime factor of $b-a(16X^2+32X+65)=130$, so $p=5$ or $p=13$. However, modulo $5$ analysis of the equation $132 u^2 - 232 u + 47=3t^2$ shows that $u$ can only be $0$ or $1$ modulo $5$, hence $a=X-2=2u-3$ is not divisible by $5$. Hence, 
$p=13$, which is not of the form $p=3k+2$, and Proposition \ref{prop:3y2pz2product} implies that $b=16 X^3+X \in S_3$. Therefore, equation \eqref{eq:y2pyzpz2mx3mx} has infinitely many integer solutions.

\vspace{10pt}

The next equation we will consider is
\begin{equation}\label{eq:2y2pz2mx3pxp1}
2y^2+z^2=x^3-x-1.
\end{equation}
Multiplying $P(x)=x^3-x-1$ by $2x$, we obtain that for any integer $k$, 
$$
2x(x^3-x-1)=2x^4-2x^2-2x=2(x^2-k)^2+((4 k - 2) x^2 -2 x -2 k^2).
$$
Let us try to select an integer $k$ such that the equation $(4 k - 2) x^2 -2 x -2 k^2=t^2$ has infinitely many integer solutions with $x$ being positive. A computer search returns that for $k=2$ the resulting equation $6x^2-2x-8=t^2$ has solution $(x,t)=(6,14)$. Because we also have $6 \cdot 1=6 >0$, $6$ is not a perfect square, $(-2)^2-4 \cdot 6 (-8) \neq 0$, and $6 > \frac{2}{12}$, Proposition \ref{prop:Gaussquadpos} implies that it has infinitely many solutions in positive integers. 
Hence, for infinitely many positive integers $x$ we have $2x(x^3-x-1)=2(x^2-2)^2+t^2 \in S_2$, where $S_2$ is defined in \eqref{S2_notation}.
If $p$ is a common prime factor of $a=2x$ and $b=x^3-x-1$ then $p$ is also a divisor of $-ax^2+a+2b=2$, so $p=2$. Because $p=2$ is not equal to $5$ or $7$ modulo $8$, Proposition \ref{prop:2y2pz2product} implies that $b=x^3-x-1 \in S_2$. Therefore, equation \eqref{eq:2y2pz2mx3pxp1} has infinitely many integer solutions.

\vspace{10pt}

The next equation we will consider is
\begin{equation}\label{eq:2y2pz2mx3mxm1}
2y^2+z^2=x^3+x+1.
\end{equation}
Multiplying $P(x)=x^3+x+1$ by $2x$, we obtain that for any integer $k$, 
$$
2x(x^3+x+1)=2x^4+2x^2+2x=2(x^2-k)^2+((4 k + 2) x^2 +2 x -2 k^2).
$$
Let us try to select an integer $k$ such that the equation $(4 k + 2) x^2 +2 x -2 k^2=t^2$ has infinitely many integer solutions with $x$ being positive. A computer search returns that for $k=2$ the resulting equation $10x^2+2x-8=t^2$ has solution $(x,t)=(1,2)$. Because we also have $10 \cdot 1=10 >0$, $10$ is not a perfect square, $(2)^2-4 \cdot 10 (-8) \neq 0$, and $1 > -\frac{2}{20}$, Proposition \ref{prop:Gaussquadpos} implies that it has infinitely many solutions in positive integers. Hence, for infinitely many positive integers $x$ we have $2x(x^3+x+1)=2(x^2-2)^2+t^2 \in S_2$, where $S_2$ is defined in \eqref{S2_notation}. If $p$ is a common prime factor of $a=2x$ and $b=x^3+x+1$ then $p$ is also a divisor of $2b-ax^2-a=2$, so $p=2$. Because $p=2$ is not equal to $5$ or $7$ modulo $8$, Proposition \ref{prop:2y2pz2product} implies that $b=x^3+x+1 \in S_2$. Therefore, equation \eqref{eq:2y2pz2mx3mxm1} has infinitely many integer solutions.

\vspace{10pt}

The next equation we will consider is
\begin{equation}\label{eq:2y2pz2mx3p3}
2y^2+z^2=x^3-3.
\end{equation}
Multiplying $P(x)=x^3-3$ by $x$, we obtain that for any integer $k$, 
$$
x(x^3-3)=(x^2-k)^2+(2kx^2-3x-k^2).
$$
Let us try to select an integer $k$ such that the equation $2kx^2-3x-k^2=2t^2$ has infinitely many integer solutions with $x$ being positive. A computer search returns that for $k=3$ the resulting equation $6x^2-3x-9=2t^2$ has solution $(x,t)=(9,15)$. Because we also have $6 \cdot 2=12 >0$, $12$ is not a perfect square, $(-3)^2-4 \cdot 6 (-9) \neq 0$, and $9 > \frac{3}{12}$, Proposition \ref{prop:Gaussquadpos} implies that it has infinitely many solutions in positive integers. Hence, for infinitely many positive integers $x$ we have $x(x^3-3)=(x^2-3)^2+2t^2 \in S_2$, where $S_2$ is defined in \eqref{S2_notation}. If $p$ is a common prime factor of $a=x$ and $b=x^3-3$ then $p$ is also a divisor of $a^3-b=3$, so $p=3$. Because $p=3$ is not equal to $5$ or $7$ modulo $8$, Proposition \ref{prop:2y2pz2product} implies that $b=x^3-3 \in S_2$. Therefore, equation \eqref{eq:2y2pz2mx3p3} has infinitely many integer solutions.

\vspace{10pt}

The next equation we will consider is
\begin{equation}\label{eq:2y2pypz2mx3p1}
2y^2+y+z^2=x^3-1.
\end{equation}
Multiplying this equation by $8$ and rearranging, we obtain
\begin{equation}\label{eq:2y2pypz2mx3p1red}
(4y+1)^2+2(2z)^2=8x^3-7.
\end{equation}
Multiplying $P(x)=8x^3-7$ by $x$, we obtain that for any integer $k$,
$$
x(8x^3-7)=2(2x^2-k)^2+(8kx^2-7x-2k^2).
$$
Let us try to select an integer $k$ such that the equation $8kx^2-7x-2k^2=t^2$ has infinitely many integer solutions with $x$ being positive. A computer search returns that for $k=1$ the resulting equation $8x^2-7x-2=t^2$ has solution $(x,t)=(2,4)$. Let $A$ be the set of positive integers such that $8x^2-7x-2$ is a perfect square. Because we also have $8 \cdot 1=8 >0$, $8$ is not a perfect square, $(-7)^2-4 \cdot 8(-2) \neq 0$, and $2 > \frac{7}{16}$, Proposition \ref{prop:Gaussquadpos} implies that the set $A$ is infinite.

For every $x\in A$, we have $x(8x^3-7) = 2(2x^2-1)^2+t^2 \in S_2$, where $S_2$ is defined in \eqref{S2_notation}. If $p$ is a common prime factor of $a=x$ and $b=8x^3-7$ then $p$ is also a divisor of $8a^3-b=7$, so $p=7$. Because $p=7$ is equal to $7$ modulo $8$, it appears that Proposition \ref{prop:2y2pz2product} is not applicable. However, if $x$ is divisible by $7$, then $8x^2-7x-2 \equiv 5 \, \, \text{(mod 7)}$, and therefore it cannot be a perfect square, thus all integers in the set $A$ are coprime to $7$. Hence, for every $x \in A$, the integers $a$ and $b$ are coprime, and Proposition \ref{prop:2y2pz2product} implies that $b \in S_2$. By the definition of $S_2$, $8x^3-7=2u^2+v^2$ for some integers $u,v$. Because $8x^3-7$ is $1$ modulo $8$, this is possible only if $u$ is even and $v$ is odd. By replacing $v$ with $-v$ if necessary, we may assume that $v$ is equal to $1$ modulo $4$.  Therefore, equation \eqref{eq:2y2pypz2mx3p1red} has infinitely many integer solutions.

\vspace{10pt}
	 
The next equation we will consider is
\begin{equation}\label{eq:y2pyzpz2mx3p3}
y^2+yz+z^2=x^3-3.
\end{equation}
Multiplying $P(x)=x^3-3$ by $x$, we obtain that for any integer $k$, 
$$
x(x^3-3)=(x^2-k)^2+(2kx^2-3x-k^2).
$$
Let us try to select an integer $k$ such that the equation $2kx^2-3x-k^2=3t^2$ has infinitely many integer solutions with $x$ being positive. A computer search returns that for $k=2$ the resulting equation $4x^2-3x-4=3t^2$ has solution $(x,t)=(4,4)$. Because we also have $4 \cdot 3=12 >0$, $12$ is not a perfect square, $(-3)^2-4 \cdot 4 (-4) \neq 0$, and $4 > \frac{3}{8}$, Proposition \ref{prop:Gaussquadpos} implies that it has infinitely many solutions in positive integers. Hence, for infinitely many positive integers $x$ we have $x(x^3-3)=(x^2-2)^2+3t^2 \in S_3$, where $S_3$ is defined in \eqref{S3_notation}. If $p$ is a common prime factor of $a=x$ and $b=x^3-3$ then $p$ is also a divisor of $a^3-b=3$. Hence $p=3$, and Proposition \ref{prop:3y2pz2product} implies that $b=x^3-3 \in S_3$. Therefore, equation \eqref{eq:y2pyzpz2mx3p3} has infinitely many integer solutions.	 

\vspace{10pt}

The next equation we will consider is
\begin{equation}\label{eq:2y2pz2mx3pxp2}
2y^2+z^2=x^3-x-2.
\end{equation}
Multiplying $P(x)=x^3-x-2$ by $x$, we obtain that for any integer $k$, 
$$
x(x^3-x-2)=(x^2-k)^2+((2k-1)x^2-2x-k^2).
$$
Let us try to select an integer $k$ such that the equation $(2k-1)x^2-2x-k^2=2t^2$ has infinitely many integer solutions with $x$ being positive. A computer search returns that for $k=4$ the resulting equation $7x^2-2x-16=2t^2$ has solution $(x,t)=(2,2)$. Because we also have $7 \cdot 2=14 >0$, $14$ is not a perfect square, $(-2)^2-4 \cdot 7 (-16) \neq 0$, and $2 > \frac{2}{14}$, Proposition \ref{prop:Gaussquadpos} implies that it has infinitely many solutions in positive integers. Hence, for infinitely many positive integers $x$ we have $x(x^3-x-2)=(x^2-4)^2+2t^2 \in S_2$, where $S_2$ is defined in \eqref{S2_notation}. If $p$ is a common prime factor of $a=x$ and $b=x^3-x-2$ then $p$ is also a divisor of $a^3-a-b=2$. Hence $p=2$, and Proposition \ref{prop:2y2pz2product} implies that $b=x^3-x-2 \in S_2$. Therefore, equation \eqref{eq:2y2pz2mx3pxp2} has infinitely many integer solutions.

\vspace{10pt}

The next equation we will consider is
\begin{equation}\label{eq:2y2pz2mx3pxm2}
2y^2+z^2=x^3-x+2.
\end{equation}
Multiplying $P(x)=x^3-x+2$ by $x$, we obtain that for any integer $k$, 
$$
x(x^3-x+2)=(x^2-k)^2+((2k-1)x^2+2x-k^2).
$$
Let us try to select an integer $k$ such that the equation $(2k-1)x^2+2x-k^2=2t^2$ has infinitely many integer solutions with $x$ being positive. A computer search returns that for $k=1$ the resulting equation $x^2+2x-1=2t^2$ has solution $(x,t)=(1,1)$. Because we also have $1 \cdot 2=2 >0$, $2$ is not a perfect square, $(2)^2-4 \cdot 1 (-1) \neq 0$, and $1 > \frac{-2}{2}=-1$, Proposition \ref{prop:Gaussquadpos} implies that it has infinitely many solutions in positive integers. Hence, for infinitely many positive integers $x$ we have $x(x^3-x+2)=(x^2-1)^2+2t^2 \in S_2$, where $S_2$ is defined in \eqref{S2_notation}. If $p$ is a common prime factor of $a=x$ and $b=x^3-x+2$ then $p$ is also a divisor of $b-a^3+a=2$. Hence $p=2$, and Proposition \ref{prop:2y2pz2product} implies that $b=x^3-x+2 \in S_2$. Therefore, equation \eqref{eq:2y2pz2mx3pxm2} has infinitely many integer solutions.

\vspace{10pt}

The next equation we will consider is
\begin{equation}\label{eq:2y2pz2mx3mxp2}
2y^2+z^2=x^3+x-2.
\end{equation}
Multiplying $P(x)=x^3-x+2$ by $x$, we obtain that for any integer $k$, 
$$
x(x^3+x-2)=(x^2-k)^2+((2k+1)x^2-2x-k^2).
$$
Let us try to select an integer $k$ such that the equation $(2k-1)x^2-2x-k^2=2t^2$ has infinitely many integer solutions with $x$ being positive. A computer search returns that for $k=1$ the resulting equation $3x^2-2x-1=2t^2$ has solution $(x,t)=(1,0)$. Because we also have $3 \cdot 2=6 >0$, $6$ is not a perfect square, $(-2)^2-4 \cdot 3 (-1) \neq 0$, and $1 > \frac{2}{6}$, Proposition \ref{prop:Gaussquadpos} implies that it has infinitely many solutions in positive integers. Hence, for infinitely many positive integers $x$ we have $x(x^3+x-2)=(x^2-1)^2+2t^2 \in S_2$, where $S_2$ is defined in \eqref{S2_notation}. If $p$ is a common prime factor of $a=x$ and $b=x^3+x-2$ then $p$ is also a divisor of $a^3+a-b=2$. Hence $p=2$, and Proposition \ref{prop:2y2pz2product} implies that $b=x^3+x-2 \in S_2$. Therefore, equation \eqref{eq:2y2pz2mx3mxp2} has infinitely many integer solutions.

\vspace{10pt}

The next equation we will consider is
\begin{equation}\label{eq:2y2pz2pzmx3p2}
2y^2+z^2+z=x^3-2.
\end{equation}
Multiplying this equation by $4$ and rearranging, we obtain
\begin{equation}\label{eq:2y2pz2pzmx3p2red}
2(2y)^2+(2z+1)^2=4x^3-7.
\end{equation}
Multiplying $P(x)=4x^3-7$ by $2x$, we obtain that for any integer $k$,
$$
2x(4x^3-7)=2(2x^2-k)^2+(8kx^2-14x-2k^2).
$$
Let us try to select an integer $k$ such that the equation $8kx^2-14x-2k^2=t^2$ has infinitely many integer solutions with $x$ being positive. A computer search returns that for $k=10$ the resulting equation $80x^2-14x-200=t^2$ has solution $(x,t)=(4,32)$. Let $A$ be the set of positive integers such that $80x^2-14x-200$ is a perfect square. Because we also have $80 \cdot 1=80 >0$, $80$ is not a perfect square, $(-14)^2-4 \cdot 80 (200) \neq 0$, and $4 > \frac{14}{160}$, Proposition \ref{prop:Gaussquadpos} implies that the set $A$ is infinite. 

For every $x \in A$, we have $2x(4x^3-7)=2(2x^2-10)^2+t^2 \in S_2$, where $S_2$ is defined in \eqref{S2_notation}. If $p$ is a common prime factor of $a=2x$ and $b=4x^3-7$, then $p$ is also a divisor of $xa^2-b=7$. Because $p=7$ is equal to $7$ modulo $8$, it appears that Proposition \ref{prop:2y2pz2product} is not applicable. However, if $x$ is divisible by $7$, then $80x^2-14x-200 \equiv 3 \, \, \text{(mod 7)}$, and therefore it cannot be a perfect square, thus all integers in the set $A$ are coprime to $7$. Hence, for every $x \in A$, integers $a$ and $b$ are coprime, and Proposition \ref{prop:2y2pz2product} implies that $b \in S_2$. By the definition of $S_2$, $4x^3-7=2u^2+v^2$ for some integers $u,v$. Because $4x^3-7$ is $1$ modulo $4$, this is possible only if $u$ is even and $v$ is odd. 
Therefore, equation \eqref{eq:2y2pz2pzmx3p2red} has infinitely many integer solutions. 

\vspace{10pt}

The next equation we will consider is
\begin{equation}\label{eq:2y2pypz2mx3p2}
2y^2+y+z^2=x^3-2.
\end{equation}
Multiplying this equation by $8$ and rearranging, we obtain
\begin{equation}\label{eq:2y2pypz2mx3p2red}
(4y+1)^2+2(2z)^2=8x^3-15.
\end{equation}
Multiplying $P(x)=8x^3-15$ by $x$, we obtain that for any integer $k$,
$$
x(8x^3-15)=2(2x^2-k)^2+(8kx^2-15x-2k^2).
$$
Let us try to select an integer $k$ such that the equation $8kx^2-15x-2k^2=t^2$ has infinitely many integer solutions with $x$ being positive. A computer search returns that for $k=1$ the resulting equation $8x^2-15x-2=t^2$ has solution $(x,t)=(2,0)$. Let $A$ be the set of positive integers such that $8x^2-15x-2$ is a perfect square. Because we also have $8 \cdot 1=8 >0$, $8$ is not a perfect square, $(-15)^2-4 \cdot 8(-2) \neq 0$, and $2 > \frac{15}{16}$, Proposition \ref{prop:Gaussquadpos} implies that the set $A$ is infinite.

For every $x \in A$, we have $x(8x^3-15)=2(2x^2-1)^2+t^2 \in S_2$, where $S_2$ is defined in \eqref{S2_notation}. If $p$ is a common prime factor of $a=x$ and $b=8x^3-15$ then $p$ is also a divisor of $8a^3-b=15$, so $p=3$ or $p=5$. Because $p=5$ is equal to $5$ modulo $8$, it appears that Proposition \ref{prop:2y2pz2product} is not satisfied. However, if $x$ is divisible by $5$, then $8x^2-15x-2 \equiv 3 \, \, \text{(mod 5)}$, and therefore it cannot be a perfect square, thus all integers in the set $A$ are coprime to $5$. Because $x$ is coprime to $5$, the only possible common prime factor of integers $a$ and $b$ is $p=3$, which is not equal to $5$ or $7$ modulo $8$. Hence, Proposition \ref{prop:2y2pz2product} implies that $b \in S_2$. By the definition of $S_2$, $8x^3-15=2u^2+v^2$ for some integers $u,v$. Because $8x^3-15$ is $1$ modulo $8$, this is possible only if $u$ is even and $v$ is odd. By replacing $v$ with $-v$ if necessary, we may assume that $v$ is equal to $1$ modulo $4$.  Therefore, equation \eqref{eq:2y2pypz2mx3p2red} has infinitely many integer solutions.

\vspace{10pt}

The next equation we consider is
\begin{equation}\label{eq:y2pyzpz2mx3pxp2}
	y^2+yz+z^2=x^3-x-2.
\end{equation}
Modulo $4$ analysis shows that $y$ and $z$ must be even, and $x$ is $2$ modulo $4$. Substituting $y=2Y$, $z=2Z$, $x=4X+2$ and cancelling $4$, we obtain the equation
$$
	Y^2+YZ+Z^2=16 X^3+24X^2+11X+1.
$$
Multiplying $P(X)=16 X^3+24X^2+11X+1$ by $X$, we obtain that for any integer $k$,
$$
X(16 X^3+24X^2+11X+1)=(4 X^2+3 X-k)^2+((8k+2)X^2+(6k+1)X-k^2).
$$ 
Let us try to select an integer $k$ such that the equation $(8k+2)X^2+(6k+1)X-k^2=3t^2$ has infinitely many integer solutions with $X$ being positive. A computer search returns that for $k=2$, the resulting equation $18X^2+13X-4=3t^2$ has solution $(X,t)=(1,3)$. Because we also have $3 \cdot 18=54 >0$, $54$ is not a perfect square, $13^2-4 \cdot 18(-4) = 457 \neq 0$, and $1 > -\frac{13}{2\cdot 18}$, Proposition \ref{prop:Gaussquadpos} implies that it has infinitely many solutions in positive integers. Hence, for infinitely many positive integers $X$ we have $X(16 X^3+24X^2+11X+1)=(4 X^2+3 X-k)^2+3t^2 \in S_3$, where $S_3$ is defined in \eqref{S3_notation}. If $p$ is a common prime factor of $a=X$ and $b=16 X^3+24X^2+11X+1$ then $p$ is also a prime factor of $1$, a contradiction, so $a$ and $b$ are coprime. Hence, Proposition \ref{prop:3y2pz2product} implies that $b=16 X^3+24X^2+11X+1 \in S_3$. Therefore, equation \eqref{eq:y2pyzpz2mx3pxp2} has infinitely many integer solutions.

\vspace{10pt}

The next equation we consider is
\begin{equation}\label{eq:y2pyzpz2mx3mxp2}
y^2+yz+z^2=x^3+x-2.
\end{equation}
Modulo $4$ analysis shows that $y$ and $z$ must be even. Let us prove that the equation has infinitely many integer solutions with $x$ being positive and odd. Substituting $y=2Y$, $z=2Z$, $x=2X+1$ and cancelling $4$, we obtain the equation
$$
Y^2+YZ+Z^2=2 X^3+3X^2+2X.
$$
Let us prove that $P(X)=2 X^3+3X^2+2X$ belongs to $S_3$, where $S_3$ is defined in \eqref{S3_notation}, for infinitely many positive integers $X$. Multiplying $P(X)$ by $2X+1$, we obtain that for any integer $k$,
$$
(2X+1)(2 X^3+3X^2+2X)=(2X^2+2X-k)^2+((4k+3)X^2+(4k+2)X-k^2).
$$
A computer search returns that for $k=1$, the resulting equation $7X^2+6X-1=3t^2$ has solution $(X,t)=(1,2)$. Because we also have $3 \cdot 7=21 >0$, $21$ is not a perfect square, $6^2-4 \cdot 7(-1) \neq 0$, and $1 > -\frac{6}{2\cdot 7}$, Proposition \ref{prop:Gaussquadpos} implies that it has infinitely many solutions in positive integers. Hence, for infinitely many positive integers $X$ we have $(2X+1)(2 X^3+3X^2+2X)=(2 X^2+2 X-k)^2+3t^2 \in S_3$. If $p$ is a common prime factor of $a=2X+1$ and $b=2 X^3+3X^2+2X$ then $p$ is also a prime factor of $a(2X^2+2X+1)-2b = 1$, a contradiction, so $a$ and $b$ are coprime. Hence, Proposition \ref{prop:3y2pz2product} implies that $b=2 X^3+3X^2+2X \in S_3$. Therefore, equation \eqref{eq:y2pyzpz2mx3mxp2} has infinitely many integer solutions.

\vspace{10pt}

The next equation we will consider is
\begin{equation}\label{eq:y2yzpz2mx3m2x}
y^2+yz+z^2=x^3+2x.
\end{equation} 
Let us prove that the equation has infinitely many integer solutions with all variables being even. Substituting $x=2X$, $y=2Y$ and $z=2Z$ and cancelling $4$, we obtain the equation
$$
Y^2+YZ+Z^2=2X^3+X.
$$
Let us prove that $P(X)=2X^3+X$ belongs to $S_3$, where $S_3$ is defined in \eqref{S3_notation}, for infinitely many positive odd integers $X$. First, let $X=2u+1$ for integer $u$. 
Then multiplying $P(X)$ by $u$, we obtain for any integer $k$,
$$
u(2X^3+X)=16u^4+24u^3+14u^2+3u=(4u^2+3u-k)^2+((8k+5)u^2+(3+6k)u-k^2).
$$
Let us try to select an integer $k$ such that the equation $(8k+5)u^2+(3+6k)u-k^2=3t^2$ has infinitely many integer solutions with $u$ being positive. A computer search returns that a suitable $k$ is $k=1$ and the equation reduces to 
$$
13u^2+9u-1=3t^2.
$$
This equation has solution $(u,t)=(4,9)$. Because we also have $3 \cdot 13=39 >0$, $39$ is not a perfect square, $9^2-4 \cdot 13 \cdot (-1) = 133 \neq 0$, and $4 > -\frac{9}{2\cdot 13}$, Proposition \ref{prop:Gaussquadpos} implies that it has infinitely many solutions in positive integers. Hence, for infinitely many positive integers $u$ we have $16u^4+24u^3+14u^2+3u=(4u^2+3u-k)^2+3t^2 \in S_3$. If $p$ is a common prime factor of $a=u$ and $b=2X^3+X=16u^3+24u^2+14u+3$ then $p$ is also a prime factor of $3$. As $3$ is not of the form $p=3k+2$, Proposition \ref{prop:3y2pz2product} implies that $b=2X^3+X \in S_3$. Therefore, equation \eqref{eq:y2yzpz2mx3m2x} has infinitely many integer solutions.

\vspace{10pt}

The final equation we will consider is
\begin{equation}\label{eq:y2pyzpz2pymx3m2}
	y^2+yz+z^2+y=x^3+2.
\end{equation}
Multiplying the equation by $12$, and rearranging, we obtain
\begin{equation}\label{eq:y2pyzpz2pymx3m2aux}
3(2y+1+z)^2+(3z-1)^2=4(3x^3+7).
\end{equation}
Let us prove that $P(x)=3x^3+7$ belongs to $S_3$, where $S_3$ is defined in \eqref{S3_notation}, for infinitely many positive integers $x$. Multiplying $P(x)$ by $x$, we obtain that for any $k$,
$$
x(3x^3+7)=3(x^2-k)^2+(6kx^2+7x-3k^2).
$$
A computer search returns that for $k=1$, the resulting equation $6x^2+7x-3=t^2$ has solution $(x,t)=(4,11)$. Because we also have $1 \cdot 6=6 >0$, $6$ is not a perfect square, $7^2-4 \cdot 6(-3) = 121 \neq 0$, and $4 > -\frac{7}{2\cdot 6}$, Proposition \ref{prop:Gaussquadpos} implies that it has infinitely many solutions in positive integers. Hence, for infinitely many positive integers $x$ we have $x(3x^3+7)=3(x^2-1)^2+t^2 \in S_3$. If $p$ is a common prime factor of $a=x$ and $b=3x^3+7$ then $p$ is also a prime factor of $7$. As $p=7$ is not of the form $3k+2$, Proposition \ref{prop:3y2pz2product} implies that $b=3x^3+7 \in S_3$. Because $4 = 3\cdot 1^2 + 1^2 \in S_3$, and set $S_3$ is multiplicative, we have that $4(3x^3+7) \in S_3$. By definition of $S_3$, this implies that $3u^2+v^2=4(3x^3+7)$ for some integers $u$ and $v$. Modulo $3$ analysis shows that $v$ is not divisible by $3$, so, by replacing $v$ with $-v$ if necessary, we may assume that $v=3z-1$ for some integer $z$. Also, modulo $2$ analysis shows that $u$ and $v$ have the same parity, hence $u$ and $z$ have opposite parity, which implies that $y=\frac{u-z-1}{2}$ in an integer. 
Then $(x,y,z)$ is a solution to \eqref{eq:y2pyzpz2pymx3m2aux}. 
Therefore, equation \eqref{eq:y2pyzpz2pymx3m2} has infinitely many integer solutions.

Finally, the equations $y^2+yz+z^2=x^3-2$ and $y^2+yz+z^2+y=x^3-2$ are solved in Section 5.4.1 of the book.

\subsection{Exercise 5.42}\label{ex:claasnumber1}
\textbf{\emph{Let $D=-3,-4,-7,-8,-11,-12,-16,-19,-27,-28,-43,-67$ or $-163$. Prove that $h(D)=1$.}}
 
 To prove this statement, we must recall some information from Section 5.4.2 of the book. 
 Binary quadratic forms $f(y,z)=a y^2 + b y z + c z^2$ with $a>0$ and $D=b^2-4ac<0$ are called \emph{positive definite}, and a form $f(y,z)=a y^2 + b y z + c z^2$ is called \emph{primitive} if $\text{gcd}(a,b,c)=1$.
\begin{definition}\label{def:formreduced}[Definition 5.40 in the book]
	A primitive positive definite form $ay^2 + byz + cz^2$ will be called \emph{reduced} if
	(i) $|b|\leq a \leq c$, and (ii) if either $|b|=a$ or $a=c$, then $b\geq 0$.
\end{definition}

For any reduced form of discriminant $D<0$, inequalities $b^2 \leq a^2$ and $a\leq c$ imply that
$$
	-D = 4ac - b^2 \geq 4a^2 - a^2 = 3a^2,
$$
thus
\begin{equation}\label{eq:Dboundforab}
	\sqrt{-D/3} \geq a \geq |b|.
\end{equation}

The class number, denoted as $h(D)$, is the number of equivalence classes of primitive positive definite forms of discriminant $D$. This number is known to be finite, and it is equal to the number of reduced forms of discriminant $D$.

\begin{proof}
Let $D=-3,-7,-11$. Then \eqref{eq:Dboundforab} implies that $2> a \geq |b|$. Because $D$ and $b$ have the same parity, $b$ must be odd, and therefore $|b|=1$. Then we must also have $a=1$. As $1=|b|\leq a =1$, condition (ii) in Definition \ref{def:formreduced} implies that $b \geq 0$, and we obtain the triple $(a,b,c)=\left(1,1,\frac{1^2-D}{4}\right),$ which is determined uniquely, thus $h(-3)=h(-7)=h(-11)=1$. 

Let $D=-4,-8$. Then \eqref{eq:Dboundforab} implies that $2> a \geq |b|$. Because $D$ and $b$ have the same parity, $b$ must be even, and therefore $b=0$. Then we must have $ac=-D/4$, hence $a=1$, and we obtain the triple $(a,b,c)=\left(1,0,-\frac{D}{4}\right)$, which is determined uniquely, thus $h(-4)=h(-8)=1$. 

Let $D=-12,-16$. Then \eqref{eq:Dboundforab} implies that $3> a \geq |b|$. Because $D$ and $b$ have the same parity, $b$ must be even, and therefore $b=0$ or $|b|=2$. Then we must have $a=1$ or $a=2$, respectively. Then $c=-\frac{b^2-D}{4a}$, and we obtain the triples $(a,b,c)=(2,2,2)$ and $(1,0,3)$ for $D=-12$ and the triples $(a,b,c)=(2,0,2)$ and $(1,0,4)$ for $D=-16$. For a primitive form, we must have $\gcd(a,b,c)=1$, however, the first triples for both $D=-12$ and $D=-16$ have $\gcd(a,b,c)=2$, so these forms are not primitive. Therefore, we obtain the triple $(a,b,c)=\left(1,1,-\frac{D}{4} \right)$, which is determined uniquely, thus $h(-12)=h(-16)=1$. 

Let $D=-19$. Then \eqref{eq:Dboundforab} implies that $3> a \geq |b|$. Because $D$ and $b$ have the same parity, $b$ must be odd, and therefore $|b|=1$. Then, $a=1$ or $a=2$. If $a=2$, we must have $c=\frac{b^2-D}{4a}=\frac{20}{8}$, hence $c$ is not integer. So $a=1$, then $|b|=a=1$ and condition (ii) in Definition \ref{def:formreduced} implies that $b \geq 0$, and we obtain the triple $(a,b,c)=(1,1,5)$, which is determined uniquely, 
thus $h(-19)=1$. 

Let $D=-27,-43$. Then \eqref{eq:Dboundforab} implies that $4> a \geq |b|$. Because $D$ and $b$ have the same parity, $b$ must be odd, and therefore $|b|=1$ or $|b|=3$. In the first case, we have $a=1,2$ or $3$, and in the second case, we must have $a=3$. When $|b|=3$, then $ac=\frac{3^2-(-27)}{4}=9$, hence $(|a|,|b|,|c|)=(3,3,3)$ we do not obtain integer triples with $\gcd(a,b,c)=1$, or $ac=\frac{3^2-(-43)}{4}=14$, which is impossible for $4>a \geq 3$. So $|b|=1$. Then $ac=\frac{1^2-D}{4}$, hence 
 $a=1$, and condition (ii) in Definition \ref{def:formreduced} implies that $b \geq 0$, and we obtain the triple $(a,b,c)=\left(1,1,\frac{1^2-D}{4}\right)$, which is determined uniquely, 
thus $h(-27)=h(-43)=1$. 

Let $D=-28$. Then \eqref{eq:Dboundforab} implies that $4> a \geq |b|$. Because $D$ and $b$ have the same parity, $b$ must be even, and therefore $b=0$ or $|b|=2$. Then we must have $ac=\frac{b^2-D}{4}$. When $|b|=2$, $ac=32$ so the only possibility is $a=2$, and we obtain the integer triple $(a,b,c)=(2,2,16)$, which is a contradiction with $\gcd(a,b,c)=1$, hence this triple is not a primitive form.  So we must have $b=0$, then $ac=-D/4$, hence $a=1$, and we obtain the triple $(a,b,c)=(1,0,7)$, which is determined uniquely, thus $h(-28)=1$. 

Let $D=-67$. Then \eqref{eq:Dboundforab} implies that $5> a \geq |b|$. Because $D$ and $b$ have the same parity, $b$ must be odd, and therefore $|b|=1$ or $3$. When $|b|=3$, then $ac=\frac{3^2-D}{4}=19$, which is impossible for $5>a \geq 3$. 
So $|b|=1$, then we must have $ac=\frac{1^2-D}{4}=17$, hence $a=1$. As $|b| =a =1$, condition (ii) in Definition \ref{def:formreduced} implies that $b \geq 0$, and we obtain the triple $(a,b,c)=(1,1,17)$, which is determined uniquely, thus $h(-67)=1$. 

Let $D=-163$. Then \eqref{eq:Dboundforab} implies that $8> a \geq |b|$. Because $D$ and $b$ have the same parity, $b$ must be odd, and therefore $|b|=1,3,5$ or $7$. Then we must have $ac=\frac{b^2-D}{4}$. If $|b|=3,5$ or $7$, then $a$ is a divisor of $43,47$ or $53$, which have no divisors satisfying $3 \leq a \leq 7$, hence $|b|=1$. Then we must have $ac=\frac{1-D}{4}=41$, hence $a=1$. As $|b| =a =1$, condition (ii) in Definition \ref{def:formreduced} implies that $b \geq 0$, and we obtain the triple $(a,b,c)=(1,1,41)$, which is determined uniquely, thus $h(-163)=1$. 
\end{proof}

\subsection{Exercise 5.48}\label{ex:H25Tomita}
\textbf{\emph{Prove, using computer assistance, that all equations listed in Table \ref{tab:H25Tomita} have infinitely many integer solutions.}}

	\begin{center}
		\begin{tabular}{ |c|c|c|c|c|c| } 
			\hline
			$H$ & Equation &  $H$ & Equation & $H$ & Equation \\ 
			\hline\hline
			$22$ & $z^2+y^2z+x^3+2 = 0$ &  $24$ & $z^2+y^2z+z+x^3+2=0$ & $25$ & $z^2+y^2z+x^3+2x+1 = 0$ \\ 
			\hline
			$22$ & $y^2 + x y z + x^3 -2=0$ &  $25$ & $z^2+y^2z-z+x^3-x+1=0$ & $25$ & $y^2+xyz+x^3+5=0$ \\ 
			\hline
			$23$ & $z^2+y^2z+x^3-x+1=0$ &  $25$ & $z^2+y^2z+z+x^3-x+1=0$ & $25$ & $y^2+xyz+x^3+x^2+1=0$ \\ 
			\hline
			$23$ & $y^2+xyz+x^3+3=0$ &  $25$ & $z^2+y^2z+z+x^3+3=0$ &  &  \\ 
			\hline			
		\end{tabular}
		\captionof{table}{\label{tab:H25Tomita} Equations solvable by Tomita's method of size $H\leq 25$.}
	\end{center} 

We first present a useful special case of Proposition \ref{prop:Gaussquad}.
\begin{proposition}\label{prop:Ax2pBxpCmDy2}[Proposition 9.1 in the book]
If the equation 
$$
Ax^2+Bx+C=Dy^2,
$$
where $A,B,C,D$ are integer coefficients, such that $AD >0$ and not a perfect square and $B^2-4AC \neq 0$, has an integer solution, then it has infinitely many. 
\end{proposition}

To prove that the equations listed in Table \ref{tab:H25Tomita} have infinitely many integer solutions, we will use the method described in Section 5.5.1.
All the equations are quadratic in one of the variables, and we need to prove that the discriminant of this equation is a perfect square for infinitely many values of the two remaining variables. The idea is to find substitutions of these variables as polynomials in a new variable (say $u$), such that the discriminant is a perfect square for infinitely many values of $u$ by Proposition \ref{prop:Ax2pBxpCmDy2}.

The first equation we will consider is
\begin{equation}\label{eq:z2py2zpx3p2}
z^2+y^2z+x^3+2 = 0.
\end{equation}
Solving this equation as a quadratic in $z$, we obtain that
$$
z_{1,2} = \frac{-y^2\pm \sqrt{y^4-4 (x^3+2)}}{2}
$$
are integers if and only if the discriminant $D(x,y)=y^4-4 (x^3+2)$ is a perfect square. Let us prove that $D(x,y)$ is a perfect square for infinitely many integer pairs $(x,y)$. A computer search suggests to look for solutions in the form $x=-3u^2+2u-2$ and $y=3u-1$ for some integer $u$. Substituting these values into the discriminant, we obtain
$$
D(x,y)=(3 u^2-2u+1)^2 (12 u^2-8u+25).
$$
The equation $12 u^2-8u+25=t^2$ has a solution $(u,t)=(0,5)$. Because we also have $12 \cdot 1=12 >0$, $12$ is not a perfect square and $(-8)^2-4 \cdot 12 \cdot 25 = -1136 \neq 0$, Proposition \ref{prop:Ax2pBxpCmDy2} implies that $12 u^2-8u+25$ is a perfect square infinitely often. This implies that equation \eqref{eq:z2py2zpx3p2} has infinitely many integer solutions.

\vspace{10pt}

The next equation we will consider is
\begin{equation}\label{eq:y2pxyzpx3m2}
y^2 + x y z + x^3 -2=0.
\end{equation}
Solving this equation as a quadratic in $y$, the discriminant is $D(x,z)=x^2 z^2-4 (x^3-2)$. Let us prove that it is a perfect square infinitely often. A computer search suggests to look for solutions in the form $x=6u^2+1$ and $z=6u$ for some integer $u$. Substituting these values into the discriminant, we obtain
$$
D(x,z)=(3 u^2+1) (12 u^2-2)^2.
$$
The equation $3 u^2+1=t^2$ has a solution $(u,t)=(0,1)$. Because we also have $3 \cdot 1=3 >0$, $3$ is not a perfect square and $0^2-4 \cdot 1 \cdot 3 = -12 \neq 0$, Proposition \ref{prop:Ax2pBxpCmDy2} implies that $3 u^2+1$ is a perfect square infinitely often. This implies that equation \eqref{eq:y2pxyzpx3m2} has infinitely many integer solutions.

\vspace{10pt}

The next equation we will consider is
\begin{equation}\label{eq:z2py2zpx3mxp1}
z^2+y^2z+x^3-x+1=0.
\end{equation}
Solving this equation as a quadratic in $z$, the discriminant is $D(x,y)=y^4 -4 (x^3-x+1)$. Let us prove that it is a perfect square infinitely often. A computer search suggests to look for solutions in the form $x=-2u^2$ and $y=2u$ for some integer $u$. Substituting these values into the discriminant, we obtain
$$
D(x,y)=(2 u^2-1) (4 u^2+2)^2.
$$
The equation $2 u^2-1=t^2$ has a solution $(u,t)=(1,1)$. Because we also have $2 \cdot 1=2 >0$, $2$ is not a perfect square and $0^2-4 \cdot 2 (-1) = 8 \neq 0$, Proposition \ref{prop:Ax2pBxpCmDy2} implies that $2u^2-1$ is a perfect square infinitely often. This implies that equation \eqref{eq:z2py2zpx3mxp1} has infinitely many integer solutions.

\vspace{10pt}

The next equation we will consider is
\begin{equation}\label{eq:y2pxyzpx3p3}
 y^2+xyz+x^3+3=0.
\end{equation}
Solving this equation as a quadratic in $y$, the discriminant is $D(x,z)= x^2 z^2-4 (x^3+3)$. Let us prove that it is a perfect square infinitely often. A computer search suggests to look for solutions in the form $x=7u^2 - 8u + 1$, and $z= 4 - 7u$ for some integer $u$. Substituting these values into the discriminant, we obtain
$$
D(x,z)=( 21 u^2-24u)(7u^2-8u+3)^2.
$$
The equation $21 u^2-24u=t^2$ has a solution $(u,t)=(0,0)$. Because we also have $5 \cdot 1=5 >0$, $5$ is not a perfect square and $(-30)^2-4 \cdot 5 (-35) = 1600 \neq 0$, Proposition \ref{prop:Ax2pBxpCmDy2} implies that $21 u^2-24u$ is a perfect square infinitely often. This implies that equation \eqref{eq:y2pxyzpx3p3} has infinitely many integer solutions.

\vspace{10pt}

The next equation we will consider is
\begin{equation}\label{eq:z2py2zpzpx3p2}
z^2+y^2z+z+x^3+2=0.
\end{equation}
Solving this equation as a quadratic in $z$, the discriminant is $D(x,y)=(1 + y^2)^2-4 ( x^3+2)$. Let us prove that it is a perfect square infinitely often. A computer search suggests to look for solutions in the form $x=-3u^2-2$ and $y=3u$ for some integer $u$. Substituting these values into the discriminant, we obtain
$$
D(x,y)=(3 u^2+1)^2 (12 u^2+25).
$$
The equation $12 u^2+25=t^2$ has a solution $(u,t)=(0,5)$. Because we also have $12 \cdot 1=12 >0$, $12$ is not a perfect square and $0^2-4 \cdot 12 \cdot 25 = -1200 \neq 0$, Proposition \ref{prop:Ax2pBxpCmDy2} implies that $12u^2+25$ is a perfect square infinitely often. This implies that equation \eqref{eq:z2py2zpzpx3p2} has infinitely many integer solutions.

\vspace{10pt}

The next equation we will consider is
\begin{equation}\label{eq:z2py2zmzpx3mxp1}
	z^2+y^2z-z+x^3-x+1=0.
\end{equation}
Solving this equation as a quadratic in $z$, the discriminant is $D(x,y)=(-1 + y^2)^2-4 (x^3-x+1)$. Let us prove that it is a perfect square infinitely often. A computer search suggests to look for solutions in the form $x=-2u^2+2u-2$ and $y=2u-1$ for some integer $u$. Substituting these values into the discriminant, we obtain
$$
D(x,y)=(4 u^2-4u+2)^2 (2 u^2-2u+5).
$$
The equation $2 u^2-2u+5=t^2$ has a solution $(u,t)=(-1,3)$. Because we also have $2 \cdot 1=2 >0$, $2$ is not a perfect square and $(-2)^2-4 \cdot 2 \cdot 5 = -36 \neq 0$, Proposition \ref{prop:Ax2pBxpCmDy2} implies that $2u^2-2u+5$ is a perfect square infinitely often. This implies that equation \eqref{eq:z2py2zmzpx3mxp1} has infinitely many integer solutions.

\vspace{10pt}

The next equation we will consider is
\begin{equation}\label{eq:z2py2zpzpx3mxp1}
	z^2+y^2z+z+x^3-x+1=0.
\end{equation}
Solving this equation as a quadratic in $z$, the discriminant is $D(x,y)=(1 + y^2)^2-4 (x^3-x+1)$. Let us prove that it is a perfect square infinitely often. A computer search suggests to look for solutions in the form $x=-2u^2+2u-1$ and $y=2u-1$ for some integer $u$. Substituting these values into the discriminant, we obtain
$$
D(x,y)=(-2u +2 u^2) (4u^2-4u+4)^2.
$$
The equation $2u^2-2u=t^2$ has a solution $(u,t)=(0,0)$. Because we also have $2 \cdot 1=2 >0$, $2$ is not a perfect square and $(-2)^2-4 \cdot 2 \cdot 0 = 4 \neq 0$, Proposition \ref{prop:Ax2pBxpCmDy2} implies that $2u^2-2u=t^2$ is a perfect square infinitely often. This implies that equation \eqref{eq:z2py2zpzpx3mxp1} has infinitely many integer solutions.

\vspace{10pt}

The next equation we will consider is
\begin{equation}\label{eq:z2py2zpzpx3p3}
	z^2+y^2z+z+x^3+3=0.
\end{equation}
Solving this equation as a quadratic in $z$, the discriminant is $D(x,y)=(y^2+1)^2-4 (x^3+3)$. Let us prove that it is a perfect square infinitely often. A computer search suggests to look for solutions in the form $x=-6u^2+4u-3$ and $y=3u-1$ for some integer $u$. Substituting these values into the discriminant, we obtain 
$$
D(x,y)=(3 u^2-2u+2)^2 (96 u^2-64u+25).
$$
The equation $96 u^2-64u+25=t^2$ has a solution $(u,t)=(0,5)$. Because we also have $96 \cdot 1=96 >0$, $96$ is not a perfect square and $(-64)^2-4 \cdot 96 \cdot 25 =-5504 \neq 0$, Proposition \ref{prop:Ax2pBxpCmDy2} implies that $96 u^2-64u+25$ is a perfect square infinitely often. This implies that equation \eqref{eq:z2py2zpzpx3p3} has infinitely many integer solutions.

\vspace{10pt}

The next equation we will consider is
\begin{equation}\label{eq:z2py2zpx3p2xp1}
	z^2+y^2z+x^3+2x+1 = 0.
\end{equation}
Solving this equation as a quadratic in $z$, the discriminant is $D(x,y)=y^4-4 (x^3+2x+1)$. Let us prove that it is a perfect square infinitely often. A computer search suggests to look for solutions in the form $x=-2u^2-1$ and $y=2u$ for some integer $u$. Substituting these values into the discriminant, we obtain
$$
D(x,y)=(2u^2+2) (4 u^2+2)^2.
$$
The equation $2 u^2+2=t^2$ has a solution $(u,t)=(1,2)$. Because we also have $2 \cdot 1=2 >0$, $2$ is not a perfect square and $0^2-4 \cdot 2 \cdot 2 = -16 \neq 0$, Proposition \ref{prop:Ax2pBxpCmDy2} implies that $2u^2+2$ is a perfect square infinitely often. This implies that equation \eqref{eq:z2py2zpx3p2xp1} has infinitely many integer solutions.

\vspace{10pt}

The next equation we will consider is
\begin{equation}\label{eq:y2pxyzpx3p5}
	y^2+xyz+x^3+5=0.
\end{equation}
Solving this equation as a quadratic in $y$, the discriminant is $D(x,z)= x^2 z^2-4 (x^3+5)$. Let us prove that it is a perfect square infinitely often. A computer search suggests to look for solutions in the form $x=-u^2+6u+6$ and $z=u-3$ for some integer $u$. Substituting these values into the discriminant, we obtain
$$
D(x,z)=(5 u^2-30u-35)(u^2-6u-4)^2.
$$
The equation $5 u^2-30u-35=t^2$ has a solution $(u,t)=(-1,0)$. Because we also have $5 \cdot 1=5 >0$, $5$ is not a perfect square and $(-30)^2-4 \cdot 5 (-35) = 1600 \neq 0$, Proposition \ref{prop:Ax2pBxpCmDy2} implies that $5 u^2-30u-35=t^2$ is a perfect square infinitely often. This implies that equation \eqref{eq:y2pxyzpx3p5} has infinitely many integer solutions.

\vspace{10pt}

The final equation we will consider is
\begin{equation}\label{eq:y2pxyzpx3px2p1}
y^2+xyz+x^3+x^2+1=0.
\end{equation}
Solving this equation as a quadratic in $y$, the discriminant is $D(x,z)=x^2 z^2-4 (x^3+x^2+1)$. Let us prove that it is a perfect square infinitely often. A computer search suggests to look for solutions in the form $x=4-u^2$ and $z=2u$ for some integer $u$. Substituting these values into the discriminant, we obtain
$$
D(x,z)= (2u^2-6)^2 (2u^2 -9).
$$
The equation $2 u^2-9=t^2$ has a solution $(u,t)=(3,3)$. Because we also have $2 \cdot 1=2 >0$, $2$ is not a perfect square and $0^2-4 \cdot 2 (-9) = 72 \neq 0$, Proposition \ref{prop:Ax2pBxpCmDy2} implies that $2 u^2-9$ is a perfect square infinitely often. This implies that equation \eqref{eq:y2pxyzpx3px2p1} has infinitely many integer solutions.

Table \ref{tab:H25Tomitasol} lists, for each equation in Table \ref{tab:H25Tomita}, a suitable substitution such that the discriminant of the quadratic equation in the remaining variable is a perfect square infinitely often by Proposition \ref{prop:Ax2pBxpCmDy2}. This then implies that the original equation has infinitely many integer solutions. 

\begin{center}
\begin{tabular}{ |c|c|c|c|c|c| } 
\hline
 Equation &  Substitution \\ 
\hline\hline
 $z^2+y^2z+x^3+2 = 0$ & $(x,y)=(-3u^2+2u-2,3u-1)$  \\\hline
 $y^2 + x y z + x^3 -2=0$ & $(x,z)=(6u^2+1,6u)$ \\\hline
 $z^2+y^2z+x^3-x+1=0$ & $(x,y)=(-2u^2,2u)$ \\\hline
 $y^2+xyz+x^3+3=0$ & $(x,z)=(7u^2-8u+1,4-7u)$ \\\hline
 $z^2+y^2z+z+x^3+2=0$ & $(x,y)=(-3u^2-2,3u)$ \\\hline
  $z^2+y^2z-z+x^3-x+1=0$ & $(x,y)=(-2u^2+2u-2,2u-1)$  \\ \hline
 $z^2+y^2z+z+x^3-x+1=0$ & $(x,y)=(-2u^2+2u-1,2u-1)$  \\ \hline
 $z^2+y^2z+z+x^3+3=0$ & $(x,y)=(-6u^2+4u-3,3u-1)$   \\ \hline	
 $z^2+y^2z+x^3+2x+1 = 0$ & $(x,y)=(-2u^2-1,2u)$ \\ \hline
$y^2+xyz+x^3+5=0$ & $(x,z)=(-u^2+6u+6,u-3)$ \\\hline 
 $y^2+xyz+x^3+x^2+1=0$ & $(x,z)=(4-u^2,2u)$ \\\hline	
\end{tabular}
\captionof{table}{\label{tab:H25Tomitasol} Equations from Table \ref{tab:H25Tomita} and a substitution such that the discriminant of the equation in the remaining variable is a perfect square infinitely often.}
\end{center} 

\subsection{Exercise 5.51}\label{ex:x3y2az2b}
\textbf{\emph{Solve Problem \ref{prob:fin} for all equations of the form $x^3y^2 = az^2 + b$ of size $H \leq 42$.}}

The equation
$
x^3y^2=z^2-1
$
has infinitely many solutions in $(x,y)$ for $z=\pm 1$. Similarly, the equations $x^3y^2=z^2-4$ and $x^3y^2=2z^2-2$ have infinitely many integer solutions with $z=\pm 2$ and $z=\pm 1$, respectively.

The other equations have only finitely many integer solutions for any fixed $z$, but we will prove that they have infinitely many integer solutions for some fixed $x$ using Proposition \ref{prop:Ax2pBxpCmDy2}. For each equation, Table \ref{tab:H42powerful} lists a suitable value of $x$. In all cases, we obtain a quadratic equation of the form $dy^2 = az^2 + c$ where $d=x^3$. In all cases, we have $ad >0$, which is also not a perfect square, the discriminant $-4ac \neq 0$, and for each equation an integer solution $(y_0,z_0)$ is presented. Hence, by Proposition \ref{prop:Ax2pBxpCmDy2}, the equations have infinitely many integer solutions.

\begin{center}
\begin{tabular}{ |c|c|c|c|c|c| } 
\hline
Equation &  $x$ & Quadratic & Solution $(y_0,z_0)$  \\ 
\hline\hline
 $x^3y^2=z^2+1$ & $5$ & $125y^2=z^2+1$ & $(61,682)$ \\\hline
 $x^3y^2=z^2+2$ & $3$ & $27y^2=z^2+2$ & $(1,5)$ \\\hline
 $x^3y^2=z^2-2$ & $7$ & $343 y^2=z^2-2$ & $(617,11427)$ \\\hline
 $x^3y^2=z^2+3$ & $7$ & $343y^2=z^2+3$ & $(2,37)$ \\\hline
 $x^3y^2=z^2-3$ & $13$ & $2197y^2=z^2-3$ & $(15503069909027, 726662475293296)$ \\\hline
 $x^3y^2=z^2+4$ & $2$ & $8y^2=z^2+4$ & $(1,2)$ \\\hline
 $x^3y^2=z^2+5$ & $29$ & $24389y^2=z^2+5$ & $(2244043487947648410385982128562847,$ \\ &&& $ 350451776492276095430077495449157936)$ \\\hline
 $x^3y^2=z^2-5$ & $11$ & $1331y^2=z^2-5 $ & $(2,73)$ \\\hline
 $x^3y^2=2z^2+1$ & $1$ & $y^2=2z^2+1$ & $(1,0)$ \\\hline
 $x^3y^2=2z^2-1$ & $1$ & $y^2=2z^2-1$ & $(1,1)$ \\\hline
 $x^3y^2=z^2+6$ & $7$ & $343y^2=z^2+6$ & $(25,463)$ \\\hline
 $x^3y^2=z^2-6$ & $19$ & $6859y^2=z^2-6$ & $(755031379, 62531004125)$ \\\hline
 $x^3y^2=2z^2+2$ & $1$ & $y^2=2z^2+2$ & $(2,1)$ \\\hline

\end{tabular}
\captionof{table}{\label{tab:H42powerful} Equations of the form $x^3 y^2=a z^2+b$ of size $H\leq 42$ which can be reduced to two-variable quadratic equations with infinitely many integer solutions for some fixed integer $x$. }
\end{center}

In the next exercises, we will solve the following problem.
\begin{problem}\label{prob:large}[Problem 5.54 in the book]
	Given a Diophantine equation $P(x_1,\dots,x_n)=0$, determine whether for any $k\geq 0$ it has a solution such that
	$$
	\min(|x_1|, \dots, |x_n|) \geq k.
	$$
	If yes, the problem is solved. If not, then describe all its integer solutions.
\end{problem}

\subsection{Exercise 5.55}\label{ex:H21linz}

\textbf{\emph{Solve Problem \ref{prob:large} for the equations listed in Table \ref{tab:H21linz}. }}

	\begin{center}
		\begin{tabular}{ |c|c|c|c|c|c| } 
			\hline
			$H$ & Equation &  $H$ & Equation &  $H$ & Equation \\ 
			\hline\hline
			$20$ & $x^3 + y^2 z - z + 2 = 0$ & $21$ & $x^3 + y^2 z + 2y + 1 = 0$ & $21$ & $1 - x + x^3 + y + 2yz = 0$ \\ 
			\hline
			$20$ & $x^3 + y^2 z + z + 2 = 0$ & $21$ & $x^3 + y^2 z + y + 3 = 0$ & $21$ & $1 + x + x^3 + y + 2yz = 0$ \\ 
			\hline
			$20$ & $x^3 + y^2 z + y + 2 = 0$ & $21$ & $x^3 + y^2 z - z + 3 = 0$ & $21$ & $3 + x^3 + y + 2yz=0$ \\ 
			\hline
			$21$ & $x^3 + y^2 z + y - x + 1 = 0$ & $21$ & $x^3 + y^2 z + z + 3 = 0$ & $21$ & $-1 + x^3 + y + z + 2yz = 0$ \\ 
			\hline
			$21$ & $x^3 + y^2 z + y + x + 1 = 0$ & $21$ & $x^3 + y^2 z + y - z + 1 = 0$ &  &  \\ 
			\hline
		\end{tabular}
		\captionof{table}{\label{tab:H21linz} Some equations linear in $z$ of size $H\leq 21$.}
	\end{center} 

For each equation, we will find non-constant polynomials $X(u)$, $Y(u)$, $Z(u)$ with integer coefficients such that $(x,y,z)=(X(u),Y(u),Z(u))$ is a solution to the equation for every $u$. We will do this by applying the method in Section 5.6.1 of the book, which we summarise below for convenience.
Each equation in Table \ref{tab:H21linz} is of the form
$$
Q(y) z = P(x,y)
$$
for some polynomials $P,Q$ with integer coefficients. Hence, we need to investigate when the ratio
$$
\frac{P(x,y)}{Q(y)},
$$
is an integer. We will start with finding a solution of the form $(x,y,z)=(F(u),G(u),c)$, where $F$, $G$ are non-constant polynomials and $c$ is a constant. Then we observe that if $k=Q(y)$ is a divisor of $P(x,y)$, then it is also a divisor of $P(x+k,y)$. Hence, we may replace $x$ by $x+k$ and $z=c$ by $z=\frac{P(x+k,y)}{k}$, which results in a solution of the form  $(x,y,z)=(X(u),Y(u),Z(u))$ with non-constant $X,Y,Z$.

Let us now consider, for example, the equation
$$
x^3+y^2z+y+2=0.
$$
For $z=0$ this equation reduces to a trivial equation $x^3+y+2=0$ which has infinitely many integer solutions of the form $(x,y,z)=(u,-u^3-2,0)$. 
Now, observe that if $x^3+y+2$ is divisible by $k=y^2$, then so is $(x+k)^3+y+2$. We can then use this to transform the solution $(x,y)=(u,-u^3-2)$ into the solution
 $$
 (x,y)=(x+k,-u^3-2)=(x+y^2,-u^3-2)=(u^6+4u^3+u+4,-u^3-2).
 $$
 Then $z=-\frac{x^3+y+2}{y^2}=-(u^{12}+8u^9+3u^7+24u^6+12u^4+32u^3+3u^2+12u+16)$, and we obtain the solution 
 $$
 (x,y,z)=\left(u^6+4u^3+u+4,-u^3-2,-\left(u^{12}+8u^9+3u^7+24u^6+12u^4+32u^3+3u^2+12u+16\right)\right), \quad u \in \mathbb{Z},
 $$
  which solves Problem \ref{prob:large} for this equation.

The other equations in Table \ref{tab:H21linz} can be solved similarly, and for each equation Table \ref{tab:H21linzsol} presents a solution of the form $(x,y,z)=(X(u),Y(u),Z(u))$ with non-constant $X,Y,Z$ and integer $u$.

\begin{center}
\begin{tabular}{ |c|c|c|c|c|c| } 
\hline
Equation &  Solution $(x,y,z)$ \\ 
\hline\hline
$x^3 + y^2 z - z + 2 = 0$ & $( 4 u^6- 2 u^2  - 1,  2 u^3,  -16u^{12}+24u^8+8u^6-12u^4-6u^2+1)$ \\ \hline
 $x^3 + y^2 z + z + 2 = 0$ & $( 4 u^6+ 2 u^2 + 1,  2 u^3, -(16u^{12}+24u^8+8u^6+12u^4+6u^2+3))$  \\ \hline
  $x^3 + y^2 z + y + 2 = 0$ & $(u^6+4u^3+u+4,-u^3-2,$ \\ & $-(u^{12}+8u^9+3u^7+24u^6+12u^4+32u^3+3u^2+12u+16))$ \\ \hline
 $x^3 + y^2 z + y - x + 1 = 0$ & $(u^6-2u^4+2u^3+u^2-u+1,-u^3+u-1, $ \\ & $-u^{12}+4u^{10}-4u^9-6u^8+9u^7-2u^6-6u^5+5u^4-3u^3-3u^2+u)$ \\ \hline
 $x^3 + y^2 z + y + x + 1 = 0$ & $(u^6+2u^4+2u^3+u^2+3u+1,-(u^3+u+1), $ \\ & $-(u^{12}+4u^{10}+4u^9+6u^8+15u^7+10u^6+$ \\ & $18u^5+19u^4+11u^3+15u^2+7u+2))$ \\ \hline
 $x^3 + y^2 z + 2y + 1 = 0$ & $(16u^6+48u^5+60u^4+44u^3+21u^2+8u+2,-(4u^3+6u^2+3u+1),$ \\ & $-(256u^{12}+1536u^{11}+4224u^{10}+7168u^9+8496u^8+7584u^7+$ \\ & $5400u^6+3168u^5+1533u^4+598u^3+189u^2+48u+7))$ \\ \hline
 $x^3 + y^2 z + y + 3 = 0$ & $(u^6+6u^3+u+9,-u^3-3,$ \\ & $-(u^{12}+12u^9+3u^7+54u^6+18u^4+108u^3+3u^2+27u+81))$  \\ \hline
 $x^3 + y^2 z - z + 3 = 0$ & $(-9u^6-3u^2+1,3u^3,81u^{12}+81u^8-18u^6+27u^4-9u^2+4)$ \\ \hline
 $x^3 + y^2 z + z + 3 = 0$ & $(9u^6+3u^2+1,3u^3,-(81u^{12}+81u^8+18u^6+27u^4+9u^2+4))$ \\ \hline
$x^3 + y^2 z + y - z + 1 = 0$ & $(u^6+2u^3+u,-u^3-1,-(u^{12}+4u^9+3u^7+4u^6+6u^4+3u^2))$ \\ \hline
 $1 - x + x^3 + y + 2yz = 0$ & $(2u^3-u+2,-u^3+u-1,4u^6-2u^4+8u^3+u^2-2u+3)$ \\ \hline
 $1 + x + x^3 + y + 2yz = 0$ & $(2u^3+3u+2,-u^3-u-1,4u^6+14u^4+8u^3+13u^2+14u+5)$ \\ \hline
 $3 + x^3 + y + 2yz=0$ & $(2u^3+u+6,-u^3-3,4u^6+6u^4+24u^3+3u^2+18u+36)$ \\ \hline
 $-1 + x^3 + y + z + 2yz = 0$ & $(-2u^3+u+3,-u^3+1,-4u^6+6u^4+12u^3-3u^2-9u-9)$ \\ \hline
\end{tabular}
\captionof{table}{\label{tab:H21linzsol} Non-constant polynomial solutions to the equations in Table \ref{tab:H21linz}, in all solutions $u$ is an arbitrary integer.}
\end{center}

\subsection{Exercise 5.56}\label{ex:H22param}
\textbf{\emph{Find non-constant parametric solutions (in the form
$$
	x_i = Q_i(u), \quad i=1,\dots, n,
$$
where $Q_i$ are \emph{non-constant} polynomials with integer coefficients) for the equations
	$$
	x^3+y^2z-x\pm z+2=0, \quad x^3+y^2z+yz+2=0, \quad x^3+xy^2+z^2+y=0
	$$	
	of size $H=22$.}}

The first equation we will consider is
\begin{equation}\label{eq:x3py2zmxpzp2}	
x^3+y^2z-x+z+2=0.
\end{equation}
A computer search discovered that this equation is solvable in integers $(x,z)$ for $y$ of the form $y=4u^3+u$ for every small $u$. To guess parametric solutions for $x,z$, we will substitute the powers of $10$ for $u$. For $u=10$, the simplest solution is $(x,z)= (-80200, 32079998)$, and for $u=100$, the simplest solution is $(x,z)=(-800020000 , 32000799999998)$. An obvious candidate is $x=-(8u^4+2u^2)$, then we can deduce that $z=32u^6+8u^4-2$.
Therefore, a family
$$
(x,y,z)=(-8u^4 -2u^2,4u^3 +u,32u^6+8u^4-2)
$$
is a candidate solution to equation \eqref{eq:x3py2zmxpzp2}. A direct substitution verifies that it is indeed a solution. 

\vspace{10pt}

The next equation we will consider is
\begin{equation}\label{eq:x3py2zmxmzp2}	
x^3+y^2z-x-z+2=0.
\end{equation}
A computer search discovered that this equation is solvable in integers $(x,z)$ for $y$ of the form $y=4u^3-u$ for every small $u$.  To guess parametric solutions for $x,z$, we will substitute the powers of $10$ for $u$. For $u=10$, the simplest solution is $(x,z)=( -79800,31920002)$, and for $u=100$, the simplest solution is $(x,z)=(-799980000 , 31999200000002)$. An obvious candidate is $x=-8u^4 +2u^2$, then we can deduce that $z=32u^6-8u^4+2$.
Therefore, a family
$$
(x,y,z)=(-8u^4 +2u^2,4u^3 -u,32u^6-8u^4+2)
$$
is a candidate solution to equation \eqref{eq:x3py2zmxmzp2}. A direct substitution verifies that it is indeed a solution. 

\vspace{10pt}

The next equation we will consider is
\begin{equation}\label{eq:x3py2zpyzp2}	
x^3+y^2z+yz+2=0.
\end{equation}
A computer search discovered that this equation is solvable in integers $(x,z)$ for $y$ of the form $y=4u^3+1$ for every small $u$.  To guess parametric solutions for $x,z$, we will substitute the powers of $10$ for $u$. For $u=10$, the simplest solution is $(x,z)=(880240,-42594906001)$, and for $u=100$, the simplest solution is $(x,z)=(80800020400,-32969632244826000001)$. An obvious candidate is $x=8u^5 +8u^4 +2u^2 +4u$, then we can deduce that $z=-(32u^9+96u^8+96u^7+32u^6+24u^5 +48u^4 +26u^3 +1 )$. 
Therefore, a family 
$$
(x,y,z)=(8u^5 +8u^4 +2u^2 +4u,4u^3 +1, -(32u^9+96u^8+96u^7+32u^6+24u^5 +48u^4 +26u^3 +1 ))
$$
is a candidate solution to equation \eqref{eq:x3py2zpyzp2}. A direct substitution verifies that it is indeed a solution. 

\vspace{10pt}

The final equation we will consider is
\begin{equation}\label{eq:x3pxy2pz2py}	
x^3+xy^2+z^2+y=0.
\end{equation}
A computer search discovered that this equation is solvable in integers $(y,z)$ for $x$ of the form $x=-u^2$ for every small $u$. To guess parametric solutions for $y,z$, we will substitute the powers of $10$ for $u$. For $u=10$, the simplest solution is $(y,z)=(1000000,10000000)$, and for $u=100$, the simplest solution is $(y,z)=(1000000000000, 100000000000000)$. Then we have the obvious candidates $y=u^6$ and $z=u^7$.
Therefore, a family
$$
(x,y,z)=(-u^2,u^6,u^7)
$$
is a candidate solution to equation \eqref{eq:x3pxy2pz2py}. A direct substitution verifies that it is indeed a solution. 

\subsection{Exercise 5.57}\label{ex:H22Gauss}
\textbf{\emph{Solve Problem \ref{prob:large} for the equations
	$$
	x^3+y^3+z^2\pm 1 = 0, \quad x^2 y + x y^2 - z^2 + 1 = 0
	$$
	of size $H=21$, as well as the equations
	$$
	x^2y + 2y^2 \pm z^2 + 2 = 0, \quad x^2y + xy^2 - z^2 + 2 =0, \quad x^3 - xy^2 + z^2 - 2 = 0 
	$$
	of size $H=22$.}}
	
In order to solve these equations we will need to use the following theorem of Siegel \cite{siegel1929uber}.	
\begin{theorem}\label{th:genus1finite}[Theorem 3.59 in the book]
Let $P(x,y)$ be an absolutely irreducible polynomial with integer (or rational) coefficients of genus $g \geq 1$. Then equation
$$
P(x,y)=0
$$
has (at most) a finite number of integer solutions.
\end{theorem}

To solve the equations in this exercise, we will use the methods presented in Section 5.6.2 of the book, which we summarise below for convenience. First, if necessary, make a substitution to reduce the equation to one that is quadratic in two of the variables. Then select a value of the other variable, so that the quadratic two-variable equation has infinitely many solutions by Gauss's Theorem. We can then select a large solution to the reduced equation and use this to show that the original equation also has a large solution.

The first equation we will consider is
\begin{equation}\label{eq:x3py3pz2p1}	
	x^3+y^3+z^2+1=0.
\end{equation}
Making the substitution $x=t-y$, the equation is reduced to
\begin{equation}\label{eq:x3py3pz2p1red}	
	1+t^3-3t^2y+3ty^2+z^2=0.
\end{equation}
This equation is quadratic in $y,z$ for every fixed $t$. For $t=-1$, it reduces to $z^2=3y^2+3y$. This equation has solution $(y,z)=(3,6)$ and all conditions of Proposition \ref{prop:Ax2pBxpCmDy2} are satisfied, hence the equation has infinitely many solutions in integers $(y,z)$. For any large solution $(y_0,z_0)$, $x=-1-y_0$ is also large, hence the problem is solved.

\vspace{10pt}

The next equation we will consider is
\begin{equation}\label{eq:x3py3pz2m1}	
	x^3+y^3+z^2-1=0.
\end{equation}
Making the substitution $x=t-y$, the equation is reduced to
\begin{equation}\label{eq:x3py3pz2m1red}	
	-1+t^3-3t^2y+3ty^2+z^2=0.
\end{equation}
This equation is quadratic in $y,z$ for every fixed $t$. For $t=-2$, it reduces to $z^2=6y^2+12y+9$. This equation has solution $(y,z)=(0,3)$ and all conditions of Proposition \ref{prop:Ax2pBxpCmDy2} are satisfied, hence the equation has infinitely many solutions in integers $(y,z)$. For any large solution $(y_0,z_0)$, $x=-2-y_0$ is also large, hence the problem is solved.

\vspace{10pt}

The next equation we will consider is
\begin{equation}\label{eq:x2ypxy2mz2p1}	
	x^2y+xy^2-z^2+1=0.
\end{equation}
Consider this equation as a quadratic in $(y,z)$ for fixed $x$. A computer search returns that for $x=2$ the resulting equation $z^2=2y^2+4y+1$ has a solution $(y,z)=(4,7)$, and the other conditions of Proposition \ref{prop:Gaussquadpos} are satisfied, hence \eqref{eq:x2ypxy2mz2p1} has infinitely many solutions with $x=2$ and $(y,z)$ positive integers. Further, there are only finitely many solutions with $x=2$ and $y$ a perfect square, because if $y=t^2$ then the resulting equation $z^2=2t^4+4t^2+1$ defines a curve with genus $1$ that can have at most finitely many integer points by Theorem \ref{th:genus1finite}. Now choose a solution $(x,y,z)=(2,y_0,z_0)$ with $y_0>0$ large and not a perfect square. Then consider the original equation as a quadratic in $(x,z)$ with $y=y_0$ and apply Proposition \ref{prop:Gaussquad2} with $(a,b,c)=(y_0,y_0^2,1)$. Then we have $a=y_0>0$, $a=y_0$ is not a perfect square, $b^2-4ac=y_0^4-4y_0 \neq 0$ and the equation has a solution $(x,z)=(2,z_0)$. Hence it has infinitely many integer solutions $(x,z)$ and we can choose one with both $|x|$ and $|z|$ being as large as we like.

\vspace{10pt}

The next equation we will consider is
\begin{equation}\label{eq:x2yp2y2pz2p2}	
	x^2y+2y^2+z^2+2=0.
\end{equation}
We will prove that it has arbitrarily large integer solutions with $x$ and $z$ both even. Then  substituting  $x=2X$ and $z=2Z$ into \eqref{eq:x2yp2y2pz2p2} and cancelling $2$, we obtain
\begin{equation}\label{eq:x2yp2y2pz2p2red}	
	1+2X^2y+y^2+2Z^2=0,
\end{equation}
reducing the coefficient of $y^2$ to $1$.
By choosing $y=-89$, we have
\[
-2(-3961+89X^2-Z^2)=0,
\]
which has the solution $(X,Z)=(7,-20)$, hence by Proposition \ref{prop:Gaussquad} it has infinitely many. Let us choose a large solution $(X_0,Z_0)$ to $-3961+89X^2-Z^2=0$, then $(x,y,z)=(2X_0,-2X_0^2+89,2Z_0)$ is a large integer solution to equation \eqref{eq:x2yp2y2pz2p2}.

\vspace{10pt}

The next equation we will consider is
\begin{equation}\label{eq:x2yp2y2mz2p2}	
	x^2y+2y^2-z^2+2=0.
\end{equation}
We will prove that is has arbitrarily large integer solutions with $x$ and $z$ both even. Then  substituting  $x=2X$ and $z=2Z$ into \eqref{eq:x2yp2y2mz2p2} and cancelling $2$, we obtain 
\begin{equation}\label{eq:x2yp2y2mz2p2red}	
	1+2X^2y+y^2-2z^2=0,
\end{equation}
reducing the coefficient of $y^2$ to $1$.
By choosing $y=7$, we have $2(25+7X^2-Z^2)=0$, which has the solution $(X,Z)=(0,5)$, and by Proposition \ref{prop:Gaussquad} it has infinitely many integer solutions. Let us choose a large solution $(X_0,Z_0)$ to $25+7X^2-Z^2=0$, then $(x,y,z)=(2X_0,-2X_0^2-7,2Z_0)$ is a large solution to equation \eqref{eq:x2yp2y2mz2p2}.

\vspace{10pt}

The next equation we will consider is
\begin{equation}\label{eq:x2ypxy2mz2p2}	
	x^2y+xy^2-z^2+2=0.
\end{equation}
Consider this equation as a quadratic in $(y,z)$ for fixed $x$. A computer search returns that for $x=194$ the resulting equation $z^2=194y^2+37636y+2$ has a solution $(y,z)=(23,984)$, and the other conditions of Proposition \ref{prop:Gaussquadpos} are satisfied, hence \eqref{eq:x2ypxy2mz2p2} has infinitely many solutions with $x=194$ and $(y,z)$ positive integers. Further, there are only finitely many solutions with $x=194$ and $y$ a perfect square, because if $y=t^2$ then the resulting equation $z^2=194t^4+37636t^2+2$ defines a curve with genus $1$ that can have at most finitely many integer points by Theorem \ref{th:genus1finite}. Now choose a solution $(x,y,z)=(194,y_0,z_0)$ with $y_0>0$ large and not a perfect square. Then consider the original equation as a quadratic in $(x,z)$ with $y=y_0$ and apply Proposition \ref{prop:Gaussquad2} with $(a,b,c)=(y_0,y_0^2,2)$ then we have $a=y_0>0$, which is not a perfect square, $b^2-4ac=y_0^4-8y_0 \neq 0$ for $y_0>2$, and the equation has a solution $(x,z)=(194,z_0)$. Hence it has infinitely many integer solutions $(x,z)$ and we can choose one with both $|x|$ and $|z|$ being large.

\vspace{10pt}

The final equation we will consider is
\begin{equation}\label{eq:x3mxy2pz2m2}	
	x^3-xy^2+z^2-2=0.
\end{equation} This equation can be rewritten as $x(x-y)(x+y)+z^2-2=0$. Let us look for solutions such that $y=x+k$, then we have $x(-k)(2x+k)+z^2-2=0$ or $z^2=2kx^2+k^2x+2$. For $k=167$, this equation has solution $(x,z)=(8401,154295)$ hence it has infinitely many integer solutions by Proposition \ref{prop:Ax2pBxpCmDy2}. For any large solution $(x_0,z_0)$, $y=x_0+k$ is also large, hence the problem is solved.

\subsection{Exercise 5.59}\label{ex:H20finxiu}
\textbf{\emph{Use Proposition \ref{prop:finxiu} to solve Problem \ref{prob:large} for the equations
	$$
	y^2+z^2 = x^3-x-1, \quad y^2+z^2=x^3+x+1, \quad y^2+z^2=x^3+3
	$$
	of size $H=19$, as well as the equations
	$$
	y^2+z^2 = x^3-x-2, \quad y^2+y+z^2 = x^3-2 
	$$
	of size $H=20$.}}

In order to solve these equations, we first cite the following theorems and propositions. 
\begin{theorem}\label{th:BakerElliptic} \cite{baker1968diophantine} [Theorem 3.30 in the book]
For any integers $a,b$ satisfying 
$$
4a^3+27b^2 \neq 0,
$$
equation 
\begin{equation}\label{eq:x3my2paxpb}
y^2=x^3+ax+b,
\end{equation} 
 has finitely many integer solutions, and there is an algorithm that can list them all.
\end{theorem}

\begin{theorem}\label{th:BakerEllipticGen}[Theorem 3.31 in the book] 
Let $a,b,c,d,e$ be integers such that integers $f=a^2+4b$, $g=ac+2d$ and $h=c^2+4e$ satisfy \eqref{eq:Weiformcond}. Then equation \eqref{eq:Weiform} has a finitely many integer solutions, and there is an algorithm that can list them all.
\end{theorem}

\begin{proposition}\label{prop:finxiu}[Proposition 5.58 in the book]
	Assume that the equation $P(x_1,\dots,x_n)=0$ has infinitely many integer solutions. Assume further that for any $i=1,\dots, n$ and any integer $u$, the equation has at most finitely many integer solutions such that $x_i=u$. Then, for any $k$, the equation has a solution such that $\text{min}(|x_1|, \dots, |x_n|) \geq k$. 
\end{proposition} 

The first equation we will consider is
\begin{equation}\label{eq:y2pz2mx3pxp1}	
y^2+z^2 =x^3-x-1.
\end{equation}
This equation is in Table \ref{tab:H20sumsquaresunsolved} and it is proven in Section 5.2.3 of the book that it has infinitely many integer solutions. 
Let us check the conditions of Proposition \ref{prop:finxiu}.
 This equation obviously has at most finitely many integer solutions $(y,z)$ for any fixed $x=u$. For a fixed $z=u$, the equation can be rewritten as $y^2=x^3-x-(u^2+1)$. This is an equation of the form \eqref{eq:x3my2paxpb} with $(a,b)=(-1,-u^2-1)$, and, by Theorem \ref{th:BakerElliptic}, it has at most finitely many integer solutions provided that $4a^3+27b^2 \neq 0$, which in our case reduces to $-4+27(-u^2-1)^2\neq 0$. The case of fixed $y=u$ is similar. Hence, by Proposition \ref{prop:finxiu}, the fact that equation \eqref{eq:y2pz2mx3pxp1} has infinitely many integer solutions immediately implies that it has solutions with absolute values of all variables arbitrarily large.

\vspace{10pt}

The next equation we will consider is
$$
y^2+z^2=x^3+x+1,
$$
which is equation \eqref{eq:y2pz2mx3mxm1}, and it 
has infinitely many integer solutions. 
This equation obviously has at most finitely many integer solutions $(y,z)$ for any fixed $x=u$. For a fixed $z=u$, the equation can be rewritten as $y^2=x^3+x-u^2+1$. This is an equation of the form \eqref{eq:x3my2paxpb} with $(a,b)=(1,1-u^2)$, and, by Theorem \ref{th:BakerElliptic}, it has at most finitely many integer solutions provided that $4a^3+27b^2 \neq 0$, which in our case reduces to $4+27(1-u^2)^2\neq 0$. The case of fixed $y=u$ is similar. Hence, by Proposition \ref{prop:finxiu}, the fact that equation \eqref{eq:y2pz2mx3mxm1} has infinitely many integer solutions immediately implies that it has solutions with absolute values of all variables arbitrarily large.

\vspace{10pt}

The next equation we will consider is
\begin{equation}\label{eq:y2pz2mx3m3}	
y^2+z^2 =x^3+3. 
\end{equation}
This equation is in Table \ref{tab:H20sumsquaresunsolved} and it is proven in Section 5.2.3 of the book that it has infinitely many integer solutions. 
This equation obviously has at most finitely many integer solutions $(y,z)$ for any fixed $x=u$. For a fixed $z=u$, the equation can be rewritten as $y^2=x^3+3-u^2$. This is an equation of the form \eqref{eq:x3my2paxpb} with $(a,b)=(0,3-u^2)$, and, by Theorem \ref{th:BakerElliptic}, it has at most finitely many integer solutions provided that $4a^3+27b^2 \neq 0$, which in our case reduces to $27(3-u^2)^2\neq 0$. The case of fixed $y=u$ is similar. Hence, by Proposition \ref{prop:finxiu}, the fact that equation \eqref{eq:y2pz2mx3m3} has infinitely many integer solutions immediately implies that it has solutions with absolute values of all variables arbitrarily large.

\vspace{10pt}

The next equation we will consider is 
$$
y^2+z^2=x^3-x-2,
$$
which is equation \eqref{eq:y2pz2mx3pxp2}, 
and it has infinitely many integer solutions. 
This equation obviously has at most finitely many integer solutions $(y,z)$ for any fixed $x=u$. For a fixed $z=u$, the equation can be rewritten as $y^2=x^3-x-(u^2+2)$. This is an equation of the form \eqref{eq:x3my2paxpb} with $(a,b)=(-1,-u^2-2)$, and, by Theorem \ref{th:BakerElliptic}, it has at most finitely many integer solutions provided that $4a^3+27b^2 \neq 0$, which in our case reduces to $-4+27(-u^2-2)^2\neq 0$. The case of fixed $y=u$ is similar. Hence, by Proposition \ref{prop:finxiu}, the fact that equation \eqref{eq:y2pz2mx3pxp2} has infinitely many integer solutions immediately implies that it has solutions with absolute values of all variables arbitrarily large.

\vspace{10pt}

The final equation we will consider is 
$$
y^2+y+z^2=x^3-2,
$$
which is equation \eqref{eq:y2pypz2mx3p2}, and it has infinitely many integer solutions.  
For a fixed $y=u$, the equation can be rewritten as $z^2 =x^3-(2+u+u^2)$. This is an equation of the form \eqref{eq:x3my2paxpb} with $(a,b)=(0,-2-u-u^2)$, and, by Theorem \ref{th:BakerElliptic}, it has at most finitely many integer solutions provided that $4a^3+27b^2 \neq 0$, which in our case reduces to $27(-2-u-u^2)^2\neq 0$. 
Let us now consider the case of fixed $z=u$.
The equation can be rewritten as $y^2+y =x^3-2-u^2$. This is an equation of the form \eqref{eq:Weiform} with $(a,b,c,d,e)=(0,0,1,0,-2-u^2)$, then by Theorem \ref{th:BakerEllipticGen}, $f=a^2+4b=0$, $g=ac+2d=0$ and $h=c^2+4e=-7-4u^2$, it has at most finitely many integer solutions provided that $-f^2 g^2 + 32 g^3 + f^3 h - 36 f g h + 108 h^2 \neq 0$, which in our case reduces to $5292 + 6048 u^2 + 1728 u^4 \neq 0$. 
The equation also obviously has at most finitely many integer solutions $(y,z)$ for any fixed $x=u$. Hence, by Proposition \ref{prop:finxiu}, the fact that equation \eqref{eq:y2pypz2mx3p2} has infinitely many integer solutions immediately implies that it has solutions with absolute values of all variables arbitrarily large.

\subsection{Exercise 5.60}\label{ex:H22quadforms}
\textbf{\emph{Solve Problem \ref{prob:large} for the equations
	$$
	y^2+z^2 = x^3+x^2+1, \quad y^2 + xy + z^2=x^3+1
	$$
	of size $H=21$, as well as the equations
	$$
	y^2 + xy + z^2 + z = x^3, \quad 2y^2 + z^2 + z = x^3
	$$
	of size $H=22$.}}

Let us first consider the equation
\begin{equation}\label{eq:y2pz2mx3mx2m1}
y^2+z^2=x^3+x^2+1.
\end{equation}
Let us start by finding families of solutions with one of the variables fixed. There are obviously finitely many integer solutions $(y,z)$ for any fixed $x=u$.
For a fixed $y=u$, the equation can be rewritten as $z^2 =x^3+x^2+(1-u^2)$. This is an equation of the form \eqref{eq:Weiform} with $(a,b,c,d,e)=(0,1,0,0,1-u^2)$, then $f=a^2+4b=4$, $g=ac+2d=0$ and $h=c^2+4e=4-4u^2$, and, by Theorem \ref{th:BakerEllipticGen}, it has at most finitely many integer solutions provided that $-f^2 g^2 + 32 g^3 + f^3 h - 36 f g h + 108 h^2 \neq 0$, which in our case reduces to $4^3 (4-4u^2) + 108 (4-4u^2)^2 \neq 0$, which is true provided that $u \neq \pm 1$. For $y=\pm 1$, equation \eqref{eq:y2pz2mx3mx2m1} has solutions $(x,y,z)=(u^2-1,\pm 1,u(u^2-1))$. By the same argument, \eqref{eq:y2pz2mx3mx2m1} has at most finitely many integer solutions for any fixed $z \neq \pm 1$, while for $z = \pm 1$ it has solutions $(u^2-1,u(u^2-1),\pm 1)$. In conclusion, equation \eqref{eq:y2pz2mx3mx2m1} has only the following infinite families of solutions with one of the variables fixed:  
$$
(x,y,z)=(u^2-1,\pm 1,u(u^2-1)), \quad (u^2-1,u(u^2-1),\pm 1).
$$ 
In all these families, $x$ is of the form $u^2-1$ for integer $u$. 
 
Let us use the method introduced in Section \ref{ex:H20sumsquares2} to prove that the equation is solvable for infinitely many values of $x$ not of the form $u^2-1$. 

We have
$$
(x^3+x^2+1)(x+1)=x^4+2x^3+x^2+x+1=(x^2+x-k)^2+(2kx^2+(2k+1)x-k^2+1),
$$
for an arbitrary integer $k$. Letting $k=1$, we have that $2kx^2+(2k+1)x-k^2+1=2x^2+3x$. Because $2x^2+3x=u^2$ has a solution $(x,u)=(12,18)$ and all other conditions in Proposition \ref{prop:Gaussquadpos} are satisfied, we can conclude that there is an infinite set $S$ of positive integers $x$ such that $2x^2+3x$ is a perfect square.

For all $x \in S$, $x^3+x^2+1$ and $x+1$ are coprime positive integers, and their product is a sum of two squares. By Proposition \ref{prop:sosproduct}, this implies that $x^3+x^2+1$ is a sum of two squares for every $x \in S$. It is now left to prove that at most finitely many $x \in S$ can be of the form $x=u^2-1$. If $x=u^2-1$, then $2x^2+3x=2u^4-u^2-1=t^2$ for some integer $t$. Equation $2u^4-u^2-1=t^2$ is absolutely irreducible and defines a curve of genus $1$, therefore, by Proposition \ref{th:genus1finite} it has only finitely many integer solutions. This proves that the original equation \eqref{eq:y2pz2mx3mx2m1} has infinitely many integer solutions such that $x \neq u^2-1$. However, for each fixed $y$ or $z$ there can be at most finitely many such solutions. This implies that the equation must have solutions $(x, y, z)$ with absolute values of all variables arbitrarily large. 

\vspace{10pt}

The next equation we will consider is
\begin{equation}\label{eq:y2pxypz2mx3m1}
y^2+xy+z^2=x^3+1.
\end{equation}
Let us start by finding families of solutions with one of the variables fixed. There are obviously finitely many integer solutions $(y,z)$ for any fixed $x=u$.
For a fixed $y=u$, the equation can be rewritten as $z^2 =x^3-xu+(1-u^2)$. This is an equation of the form \eqref{eq:Weiform} with $(a,b,c,d,e)=(0,0,0,-u,1-u^2)$, then $f=a^2+4b=0$, $g=ac+2d=-2u$ and $h=c^2+4e=4-4u^2$, and, by Theorem \ref{th:BakerEllipticGen}, it has at most finitely many integer solutions provided that $-f^2 g^2 + 32 g^3 + f^3 h - 36 f g h + 108 h^2 \neq 0$, which in our case reduces to $32 (-2u)^3  + 108 (4-4u^2)^2 \neq 0$, which is always true. For a fixed $z=u$, the equation can be rewritten as $y^2+xy =x^3+(1-u^2)$. This is an equation of the form \eqref{eq:Weiform} with $(a,b,c,d,e)=(1,0,0,0,1-u^2)$, then $f=a^2+4b=1$, $g=ac+2d=0$ and $h=c^2+4e=4-4u^2$, and, by Theorem \ref{th:BakerEllipticGen}, it has at most finitely many integer solutions provided that $-f^2 g^2 + 32 g^3 + f^3 h - 36 f g h + 108 h^2 \neq 0$, which in our case reduces to $ 109 (4-4u^2)^2 \neq 0$, which is true provided that $u \neq \pm 1$. For $z=\pm 1$, equation \eqref{eq:y2pxypz2mx3m1} has solutions $(x,y,z)=(u^2+u,u^3+u^2,\pm 1)$. In conclusion, equation \eqref{eq:y2pxypz2mx3m1} has only the following infinite families of solutions with one of the variables fixed:  
$$
(x,y,z)=(u^2+u,u^3+u^2,\pm 1).
$$ 
In all these families, $x$ is of the form $u^2+u$ for integer $u$. 

Let us now prove that the equation is solvable for infinitely many values of $x$ not of the form $u^2+u$.
Multiplying \eqref{eq:y2pxypz2mx3m1} by $4$ and rearranging, we obtain
\begin{equation}\label{eq:y2pxypz2mx3m1a}
(2y+x)^2+(2z)^2=4x^3+x^2+4.
\end{equation}
We have
$$
16x(4x^3+x^2+4)=64x^4+16x^3+64x=(8x^2+x-k)^2+((16k-1)x^2+(2k+64)x-k^2)
$$
for an arbitrary integer $k$. Letting $k=9$, we have that $(16k-1)x^2+(2k+64)x-k^2=143x^2+82x-81$. Because $143x^2+82x-81=u^2$ has a solution $(x,u)=(1,12)$ and all other conditions in Proposition \ref{prop:Gaussquadpos} are satisfied, we can conclude that there exists an infinite set $S$ of positive integers $x$ such that $143x^2+82x-81$ is a perfect square. 

Every $x \in S$ is odd, as otherwise $143x^2+82x-81 \equiv 3 \,\, (\text{mod} \,\, 4)$. 
Hence for all $x \in S$, $16x$ and $4x^3+x^2+4$ are coprime positive integers, and their product is a sum of two squares. By Proposition \ref{prop:sosproduct}, this implies that 
for every $x \in S$, $4x^3+x^2+4$ is the sum of two squares, say $Y^2+Z^2$. Because $x$ is odd, $Y$ and $Z$ must have opposite parity. By swapping $Y$ and $Z$ if necessary, we may assume that $Y$ is odd and $Z$ is even. Then $y=\frac{Y-x}{2}$ and $z=\frac{Z}{2}$ are integers, and $(x,y,z)$ is a solution to \eqref{eq:y2pxypz2mx3m1a}, and therefore is a solution to \eqref{eq:y2pxypz2mx3m1}.
It is now left to prove that only finitely many of $x \in S$ are of the form $x=u^2+u$. If $x=u^2+u$, then $143x^2+82x-81=143u^4+286u^3+225u^2+82u-81=t^2$ for some integer $t$. The last equation is absolutely irreducible and defines a curve of genus $1$,  therefore, by Proposition \ref{th:genus1finite}, it has only finitely many integer solutions. This proves that the original equation \eqref{eq:y2pxypz2mx3m1} has infinitely many integer solutions such that $x \neq u^2+u$. However, for each fixed $y$ or $z$ there can be at most finitely many such solutions. This implies that equation \eqref{eq:y2pxypz2mx3m1} must have solutions $(x, y, z)$ with absolute values of all variables arbitrarily large.

\vspace{10pt}

The next equation we will consider is
\begin{equation}\label{eq:y2pxypz2pzmx3}
y^2+xy+z^2+z=x^3.
\end{equation}
Let us start by finding families of solutions with one of the variables fixed. There are obviously finitely many integer solutions $(y,z)$ for any fixed $x=u$.
For a fixed $y=u$, the equation can be rewritten as $z^2+z =x^3-xu-u^2$. This is an equation of the form \eqref{eq:Weiform} with $(a,b,c,d,e)=(0,0,1,-u,-u^2)$, then $f=a^2+4b=0$, $g=ac+2d=-2u$ and $h=c^2+4e=1-4u^2$, and, by Theorem \ref{th:BakerEllipticGen}, it has at most finitely many integer solutions provided that $-f^2 g^2 + 32 g^3 + f^3 h - 36 f g h + 108 h^2 \neq 0$, which in our case reduces to $ 32 (-2u)^3  + 108 (1-4u^2)^2 \neq 0$, which is always true. For a fixed $z=u$, the equation can be rewritten as $y^2+xy =x^3+(-u-u^2)$. This is an equation of the form \eqref{eq:Weiform} with $(a,b,c,d,e)=(1,0,0,0,-u-u^2)$, then $f=a^2+4b=1$, $g=ac+2d=0$ and $h=c^2+4e=-4u-4u^2$, and, by Theorem \ref{th:BakerEllipticGen}, it has at most finitely many integer solutions provided that $-f^2 g^2 + 32 g^3 + f^3 h - 36 f g h + 108 h^2 \neq 0$, which in our case reduces to $109(-4u-4u^2) ^2 \neq 0$, which is true provided that $u \neq -1$ or $u \neq 0$. For $z=0$ or $z=-1$, equation \eqref{eq:y2pxypz2pzmx3} has solutions $(x,y,z)=(u^2+u,u^3+u^2,0)$ or $(u^2+u,u^3+u^2,-1)$. In conclusion, equation \eqref{eq:y2pxypz2pzmx3} has only the following infinite families of solutions with one of the variables fixed:  
$$
(x,y,z)=(u^2+u,u^3+u^2,0),(u^2+u,u^3+u^2,-1).
$$ 
In all these families, $x$ is of the form $u^2+u$ for integer $u$. 

Let us prove that the equation is solvable for infinitely many values of $x$ not of the form $u^2+u$. 
Multiplying \eqref{eq:y2pxypz2pzmx3} by $4$ and rearranging, we obtain
\begin{equation}\label{eq:y2pxypz2pzmx3a}
(2y+x)^2+(2z+1)^2=4x^3+x^2+1.
\end{equation}
We have
$$
16x(4x^3+x^2+1)=64x^4+16x^3+16x=(8x^2+x-k)^2+((16k-1)x^2+(2k+16)x-k^2)
$$
for an arbitrary integer $k$. Letting $k=6$, we have that $(16k-1)x^2+(2k+16)x-k^2=95x^2+28x-36$. 
Because $95x^2+28x-36=u^2$ has a solution $(x,u)=(2,20)$ and all other conditions in Proposition \ref{prop:Gaussquadpos} are satisfied, we can conclude that there exists an infinite set $S$ of positive integers $x$ such that $95x^2+28x-36$ is a perfect square. 

Every $x \in S$ is even, as otherwise $95x^2+28x-36 \equiv 3 \,\, (\text{mod} \,\, 4)$. 
Hence, for all $x \in S$, $16x$ and $4x^3+x^2+1$ are coprime positive integers and their product is a sum of two squares. By Proposition \ref{prop:sosproduct} 
this implies that for every $x \in  S$, $4x^3+x^2+1$ is the sum of two squares, say $Y^2+Z^2$. Because $x$ is even, $Y$ and $Z$ must have opposite parity. By swapping $Y$ and $Z$ if necessary, we may assume that $Y$ is even and $Z$ is odd. Then $y=\frac{Y-x}{2}$ and $z=\frac{Z-1}{2}$ are integers, and $(x,y,z)$ is a solution to \eqref{eq:y2pxypz2pzmx3a}, and therefore is a solution to \eqref{eq:y2pxypz2pzmx3}.
It is now left to prove that only finitely many of integers $x \in S$ are of the form $x=u^2+u$. If $x=u^2+u$, then $95x^2+28x-36=95u^4+190u^3+123u^2+28u-36=t^2$ for some integer $t$. Equation $95u^4+190u^3+123u^2+28u-36=t^2$ is absolutely irreducible and defines a curve of genus $1$, therefore, by Proposition \ref{th:genus1finite}, it has only finitely many integer solutions. This proves that the original equation \eqref{eq:y2pxypz2pzmx3} has infinitely many integer solutions such that $x \neq u^2+u$. However, for each fixed $y$ or $z$ there can be at most finitely many such solutions. This implies that the equation must have solutions $(x, y, z)$ with absolute values of all variables arbitrarily large.

\vspace{10pt}

The final equation we will consider is
\begin{equation}\label{eq:2y2pz2pzmx3}
2y^2+z^2+z=x^3.
\end{equation}
Let us start by finding families of solutions with one of the variables fixed. There are obviously finitely many integer solutions $(y,z)$ for any fixed $x=u$.
For a fixed $y=u$, the equation can be rewritten as $z^2+z =x^3-2u^2$. This is an equation of the form \eqref{eq:Weiform} with $(a,b,c,d,e)=(0,0,1,0,-2u^2)$, then $f=a^2+4b=0$, $g=ac+2d=0$ and $h=c^2+4e=1-4u^2$, and, by Theorem \ref{th:BakerEllipticGen}, it has at most finitely many integer solutions provided that $-f^2 g^2 + 32 g^3 + f^3 h - 36 f g h + 108 h^2 \neq 0$, which in our case reduces to $  108 (1-4u^2)^2 \neq 0$, which is always true. For a fixed $z=u$, the equation can be rewritten as $2y^2=x^3+(-u-u^2)$, after multiplying by $8$ and the change of variables $X=2x$ and $Y=4y$, we obtain $Y^2=X^3+(-8u-8u^2)$. This is an equation of the form \eqref{eq:Weiform} with $(a,b,c,d,e)=(0,0,0,0,-8u-8u^2)$, then $f=a^2+4b=0$, $g=ac+2d=0$ and $h=c^2+4e=-32u-32u^2$, and, by Theorem \ref{th:BakerEllipticGen}, it has at most finitely many integer solutions provided that $-f^2 g^2 + 32 g^3 + f^3 h - 36 f g h + 108 h^2 \neq 0$, which in our case reduces to $ 108 (-8u-8u^2)^2\neq 0$, which is true provided that $u \neq -1$ or $u \neq 0$. For $z=0$ or $z=-1$, equation \eqref{eq:2y2pz2pzmx3} has solutions $(x,y,z)=(2u^2,2u^3,0)$ or $(2u^2,2u^3,-1)$. In conclusion, equation \eqref{eq:2y2pz2pzmx3} has only the following infinite families of solutions with one of the variables fixed:  
$$
(x,y,z)=(2u^2,2u^3,0),(2u^2,2u^3,-1).
$$  
In all these families, $x$ is of the form $2u^2$ for integer $u$. 

Let us prove that the equation is solvable for infinitely many values of $x$ not of the form $2u^2$. 
Multiplying \eqref{eq:2y2pz2pzmx3} by $4$ and rearranging, we obtain
\begin{equation}\label{eq:2y2pz2pzmx3a}
2(2y)^2+(2z+1)^2=4x^3+1.
\end{equation}
We have
$$
2x(4x^3+1)=8x^4+2x=2(2x^2-k)^2+(8kx^2+2x-2k^2)
$$
for an arbitrary integer $k$. Letting $k=4$, we have that $8kx^2+2x-2k^2=32x^2+2x-32$. Because $32x^2+2x-32=u^2$ has a solution $(x,u)=(2,10)$ and all other conditions in Proposition \ref{prop:Gaussquadpos} are satisfied, we conclude that there exists an infinite set $S$ of positive integers $x$ such that $32x^2+2x-32$ is a perfect square. 

Every $x \in S$ is even, as otherwise $32x^2+2x-32 \equiv 2 \,\, (\text{mod} \,\, 4)$.  Hence, for all $x \in S$, $2x$ and $4x^3+1$ are coprime positive integers and their product is $\in S_2$, where $S_2$ is defined in \eqref{S2_notation}. By Proposition \ref{prop:2y2pz2product} 
this implies that $4x^3+1 \in S_2$ for every $x \in S$. By definition of $S_2$, we have $4x^3+1=2Y^2+Z^2$ for some integers $Y,Z$. Because $4x^3+1$ is $1$ modulo $4$, $Z$ must be odd, and $Y$ must be even. Then $y=\frac{Y}{2}$ and $z=\frac{Z-1}{2}$ are integers, and $(x,y,z)$ is a solution to \eqref{eq:2y2pz2pzmx3a}, and therefore is a solution to \eqref{eq:2y2pz2pzmx3}. It is now left to prove that only finitely many of integers $x \in S$ are of the form $x=2u^2$. If $x=2u^2$, then $32x^2+2x-32=128u^4+4u^2-32=t^2$ for some integer $t$. Equation $128u^4+4u^2-32=t^2$ defines a curve of genus $1$, therefore, by Proposition \ref{th:genus1finite} it has only finitely many solutions in integers $(u,t)$. This proves that the original equation \eqref{eq:2y2pz2pzmx3} has infinitely many integer solutions such that $x \neq 2u^2$. However, for each fixed $y$ or $z$ there can be at most finitely many such solutions. This implies that the equation must have solutions $(x, y, z)$ with absolute values of all variables arbitrarily large.

\subsection{Exercise 5.64}\label{ex:x2m2y2}
\textbf{\emph{Prove the analogues of Theorem \ref{th:5y2mz2} and Proposition \ref{prop:5y2mz2product} for integers representable as $y^2-2z^2$.}}

Let us first introduce some notation and a necessary proposition from the book. Let us define the set $S_{-2}$ as
$$
S_{-2}=\{ n \in \mathbb{Z}: n = 2y^2-z^2 \, \text{for some} \, y,z \in \mathbb{Z}\}.
$$
This set is closed under multiplication as implied by the identity 
\begin{equation}\label{eq:multidentity}
(ka^2 + b^2)(kc^2+d^2) = k(ad + bc)^2 + (bd-kac)^2
\end{equation}
 with $k=-2$, and it can be equivalently written as 
$$
S_{-2}=\{ n \in \mathbb{Z}: n = y^2-2z^2 \, \text{for some} \, y,z \in \mathbb{Z}\},
$$
because $n=2Y^2-Z^2=(2Y-Z)^2-2(Z-Y)^2=y^2-2z^2$ where $y=2Y-Z$ and $z=Z-Y$. Similarly, if $n=y^2-2z^2$, then $n=2Y^2-Z^2$ for $Y=y+z$ and $Z=y+2z$.

\begin{proposition}\label{prop:quadform}[Proposition 5.25 in the book]
	Let $a,b,c$ be integers, and let $D=b^2-4ac$. Let $p$ be a prime factor of $ax^2+bxy+cy^2$ for some integers $x,y$. Then either $p$ is a divisor of $2D$, or $\left(\frac{D}{p}\right)=1$, or $p$ is a common divisor of $x$ and $y$.
\end{proposition}

We may now state the following theorem. 

\begin{theorem}\label{th:2y2mz2}
An integer $n$ belongs to $S_{-2}$ if and only if every prime equal to $3$ or $5$ modulo $8$ enters the prime factorisation of $n$ with even exponent.
\end{theorem}

\begin{proof}
First assume that $n=2y^2-z^2$ for some integers $y,z$, and let $p$ be any prime factor of $n$ that appears in its prime factorisation with an odd exponent. Let $d=\gcd(y,z)$, $y_1=y/d$ and $z_1=z/d$. Then $n$ is divisible by $d^2$, and $n/d^2=2y_1^2-z_1^2$. Then $p$ is a prime factor of $n/d^2$, but it is not a common divisor of $y_1$ and $z_1$. Hence, by Proposition \ref{prop:quadform} with $(a,b,c)=(2,0,-1)$, either $p$ is a divisor of $2D$ where $D=0^2-4 \cdot 2 \cdot (-1)=8$, or $\left(\frac{8}{p}\right)=1$. In the first case, $p=2$. In the last case, $\left(\frac{2}{p}\right)=1$ which implies that $p$ is equal to $\pm 1$ modulo $8$ by \eqref{leg:p2}. In any case, $p$ cannot be equal to $3$ or $5$ modulo $8$. 

Conversely, assume that every prime $p$ equal to $3$ or $5$ modulo $8$ enters the prime factorisation of $n$ with an even exponent. Then $n$ can be written as $n=ek^2 \cdot \prod_{i=1}^m p_i$, where $e$ is equal to $1$ or $-1$, $k$ is an integer, and each $p_i$ is a prime that is either equal to $2$ or is $\pm 1$ modulo $8$. We need to prove that $p \in S_{-2}$.  
Because $1=2(1)^2-1^2 \in S_{-2}$ and $-1=2(0)^2-1^2 \in S_{-2}$ and $k^2=2(k)^2-k^2 \in S_{-2}$, $2=2(1)^2-0^2 \in S_{-2}$ and $S_{-2}$ is closed under multiplication, it suffices to prove that $p \in S_{-2}$ for every prime $p$ equal to $\pm 1$ modulo $8$. For every such prime, we have $\left( \frac{2}{p} \right)=1$ by \eqref{leg:p2}, hence, by the definition of the Legendre symbol, there exists an integer $x$ such that $x^2-2$ is divisible by $p$. Because $p$ is odd, we conclude that $p \in S_{-2}$ by Proposition \ref{prop:divsofx2pc}.
\end{proof}

Theorem \ref{th:2y2mz2} implies the following result.

\begin{proposition}\label{prop:2y2mz2}
If the product $ab$ of two non-zero integers $a$ and $b$ belongs to $S_{-2}$ and $a,b$ do not share any prime factor equal to $3$ or $5$ modulo $8$, then $a \in S_{-2}$ and $b \in S_{-2}$.
\end{proposition}

\subsection{Exercise 5.66}\label{ex:y2pyzmz2mx3px}
\textbf{\emph{Prove that the equation
	\begin{equation}\label{eq:y2pyzmz2mx3px}
		y^2+yz-z^2 = x^3-x
	\end{equation}
	is solvable in $(y,z)$ for infinitely many values of $x$.}}

We first need to recall the following proposition.
\begin{proposition}\label{prop:y2pyzmz2set}[Proposition 5.61 in the book]
	An integer $n$ can be represented as $n=y^2+yz-z^2$ for some integers $y,z$ if and only if it can be represented as $n=5Y^2-Z^2$ for some integers $Y,Z$.
\end{proposition}

By Proposition \ref{prop:y2pyzmz2set}, we can solve equation \eqref{eq:y2pyzmz2mx3px} by showing that $x^3-x \in S_{-5}$ for infinitely many $x$, where
$$
S_{-5}=\{ n \in \mathbb{Z}: n = 5y^2-z^2 \, \text{for some} \, y,z \in \mathbb{Z}\} = \{ n \in \mathbb{Z}: n = y^2-5z^2 \, \text{for some} \, y,z \in \mathbb{Z}\},
$$
see Section 5.6.3 of the book. 
Identity \eqref{eq:multidentity} with $k=-5$ implies that set $S_{-5}$ is closed under multiplication. Hence, to prove that $x^3-x=(x-1)x(x+1) \in S_{-5}$ for infinitely many $x$, it is sufficient to prove that the three integers $x-1,x,x+1$ are in $S_{-5}$ infinitely often. Let us choose $x=5t^2$ for some integer $t$. In this case, we have $x-1=5t^2-1 \in S_{-5}$ and $x=5t^2-0 \in S_{-5}$, so it is left to determine whether $x+1=5t^2+1 \in S_{-5}$, or, equivalently, whether $5t^2+1=u^2-5v^2$ for some integers $u,v$. By selecting $v=0$, we obtain equation $5t^2+1=u^2$, which has the integer solution $(t,u)=(0,1)$. By using Proposition \ref{prop:Ax2pBxpCmDy2} we have that this equation has infinitely many integer solutions, and we can conclude that $x+1 \in S_{-5}$ infinitely often. Therefore, equation \eqref{eq:y2pyzmz2mx3px} is solvable in $(y,z)$ for infinitely many values of $x$. 

\section{Chapter 6}
In the next exercises, we will solve the following problem. 
\begin{problem}\label{prob:finhom}[Problem 6.2 in the book]
	Given a homogeneous polynomial $P$ with integer (or rational) coefficients, determine whether the equation $P=0$ has a non-trivial integer (equivalently, rational) solution.
\end{problem} 

\subsection{Exercise 6.3}\label{ex:H48Hasse}
\textbf{\emph{Solve Problem \ref{prob:finhom} for the equations listed in Table \ref{tab:H48nosol}, except for the last 4 equations in the second column (marked in bold).}}

	\begin{center}
		\begin{tabular}{ |c|c|c|c| } 
			\hline
			$H$ & Equation &  $H$ & Equation \\ 
			\hline\hline
			$40$ & $2 x^3+y^3+y^2 z-z^3 = 0$ &  $48$ & $2 x^3+x y^2+y^3+y^2 z-z^3 = 0$ \\ 
			\hline
			$40$ & $2 x^3+y^3+y^2 z+z^3 = 0$ &  $48$ & $3 x^3+y^3+y^2 z-z^3=0$ \\ 
			\hline
			$48$ & $x^4+y^4+z^4=0$ &  $48$ & $\mathbf{x^3+x^2 y-y^3+x^2 z+y^2 z+z^3=0}$ \\ 
			\hline
			$48$ & $x^3+x^2 y-y^3+y^2 z+x z^2-z^3=0$ &  $48$ & $\mathbf{2 x^3+x^2 y-y^3+y^2 z-z^3=0}$ \\ 
			\hline
			$48$ & $x^3+x^2 y-y^3+y^2 z-x z^2+z^3=0$ &  $48$ & $\mathbf{2 x^3+x y^2+y^3-y z^2+z^3=0}$ \\ 
			\hline
			$48$ & $x^3+x^2 y+y^3+y^2 z-x z^2-z^3=0$ &  $48$ & $\mathbf{2 x^3+y^3+x y z-y^2 z-z^3=0}$ \\ 
			\hline
			$48$ & $x^3+x^2 y+y^3+y^2 z+x z^2+z^3=0$ &  &  \\ 
			\hline	
		\end{tabular}
		\captionof{table}{\label{tab:H48nosol} Homogeneous equations of size $H\leq 48$ with no non-trivial solutions}
	\end{center} 
	
In this exercise we will use the following proposition to prove that the equations listed in Table \ref{tab:H48nosol}, excluding those in bold, have no non-trivial integer solutions.
	
\begin{proposition}\label{prop:homequiv}[Proposition 6.1 in the book]
Let $P(x_1,\dots,x_n)$ be a homogeneous polynomial. Then the following statements are equivalent.
	\begin{itemize}
		\item[(a)] The equation $P=0$ has a non-trivial rational solution.
		\item[(b)] The equation $P=0$ has a non-trivial integer solution.
		\item[(c)] The equation $P=0$ has a non-trivial integer solution $(x_1,\dots,x_n)$ such that $\text{gcd}(x_1,\dots,x_n)=1$.
		\item[(d)] The equation $P=0$ has infinitely many integer solutions.
	\end{itemize}
\end{proposition}

Modulo $4$ analysis of equations 
$$
\begin{aligned}
x^3+x^2 y-y^3+y^2 z+x z^2-z^3=0, & \quad &  x^3+x^2 y-y^3+y^2 z-x z^2+z^3=0, \\ 
x^3+x^2 y+y^3+y^2 z-x z^2-z^3=0, & \quad \text{and} \quad & x^3+x^2 y+y^3+y^2 z+x z^2+z^3=0,
\end{aligned}
$$
suggests that all variables must be even. Hence, these equations have no integer solutions satisfying $\gcd(x,y,z)=1$, and therefore, by Proposition \ref{prop:homequiv}, they have no non-trivial integer solutions. 
Similarly, modulo $9$ analysis of the equation 
$$
2 x^3+x y^2+y^3+y^2 z-z^3 = 0
$$
shows that all variables must be divisible by $3$. Hence, this equation has no integer solutions 
satisfying $\gcd(x,y,z)=1$, and therefore, by Proposition \ref{prop:homequiv}, it has no non-trivial integer solutions. 

Finally, the equations $2x^3+y^3+y^2z \pm z^3=0$, $x^4+y^4+z^4=0$ and $3x^3+y^3+y^2-z^3=0$ are solved in Section 6.1.1 of the book.

\subsection{Exercise 6.5}\label{ex:H48Hassecount}
\textbf{\emph{Check that the equations $x^3+x^2 y-y^3+x^2 z+y^2 z+z^3=0$ and $2 x^3+x y^2+y^3-y z^2+z^3=0$ from Table \ref{tab:H48nosol} are also solvable everywhere locally but do not have any non-trivial integer/rational solutions.}}

The equations in Table \ref{tab:H48nosol} are of the form 
\begin{equation}\label{eq:homcubgen}
	a_1x^3+a_2y^3+a_3z^3+a_4x^2y+a_5x^2z+a_6y^2x+a_7y^2z+a_8z^2x+a_9z^2y+a_{10}xyz=0,
\end{equation}  
with integer coefficients. We can check that an equation \eqref{eq:homcubgen} is solvable everywhere locally by using the Magma code
\begin{align}\label{cmd:local}
	&{\tt > C := GenusOneModel(3,\, "a_1 \, a_2 \, a_3 \, a_4 \, a_5 \, a_6 \, a_7 \, a_8 \, a_9 \, a_{10}");}\nonumber\\
	&{\tt > IsLocallySoluble(C);}\nonumber\\
	&{\tt > E := Jacobian(C);}\nonumber\\
	&{\tt > MW,iso := MordellWeilGroup(E);}\nonumber\\
	&{\tt > G,quomap := quo<MW|3*MW>;}\\
	&{\tt > PP := [iso(g @@ quomap): g\, in \, G];}\nonumber\\
	&{\tt > cubics := [P \, eq \, E!0 \, select \, GenusOneModel(3,E) \, else \, GenusOneModel(3,P): \, P \, in \, PP];}\nonumber\\
	&{\tt > exists\{ deltaP : deltaP \, in \, cubics | IsEquivalent(C,deltaP) \};}\nonumber
\end{align}
described in Section 6.1.2 of the book. This code outputs ``true/false true/false''. If the first ``true/false'' is ``true'' then the equation is solvable everywhere locally. If the final ``true/false'' is ``false'' then the equation does not have any non-trivial integer/rational solutions. 

The first equation we will consider is
\begin{equation}\label{eq:x3px2ymy3px2zpy2zpz3}
x^3+x^2y-y^3+x^2z+y^2z+z^3=0.
\end{equation}
This is equation \eqref{eq:homcubgen} with $(a_1,...,a_{10})=(1, -1 ,1 ,1 ,1, 0 ,1 ,0 ,0 ,0)$. We can therefore use the Magma code \eqref{cmd:local} with the first line 
$$
{\tt > C := GenusOneModel(3,\, "1 \, -1 \, 1 \, 1 \, 1 \, 0 \, 1 \, 0 \, 0 \,0");}
$$ 
which outputs 
{\tt true false}.
This means that equation \eqref{eq:x3px2ymy3px2zpy2zpz3} is solvable everywhere locally, but it has no non-trivial integer solutions.  

\vspace{10pt}

Let us now consider the equation
\begin{equation}\label{eq:2x3pxy2py3myz2pz3}
2x^3+xy^2+y^3-yz^2+z^3=0.
\end{equation}
This is equation \eqref{eq:homcubgen} with $(a_1,...,a_{10})=(2, 1 ,1 ,0,0,1,0,0,-1,0)$. So, in this case we can use the code \eqref{cmd:local} with the first line 
$$
{\tt > C := GenusOneModel(3,\, "2 \, 1 \, 1 \, 0 \, 0 \, 1 \, 0 \, 0 \, -1 \,0");}
$$ 
which outputs 
 {\tt true false}.
This means that equation \eqref{eq:2x3pxy2py3myz2pz3} is solvable everywhere locally, but it has no non-trivial integer solutions.   

\subsection{Exercise 6.7} \label{ex:H96hom}
\textbf{\emph{Observe that each of the equations in Table \ref{tab:H96homa} has a variable that appears in the equation with even exponents only. Use this to reduce each such equation to a genus $1$ curve. If it has no rational points, or is an elliptic curve of rank $0$, solve the equation. }}\\
	\begin{center}
	\begin{tabular}{ |c|c| }
		\hline
		$H$ & Equation \\
		\hline\hline
		$96$ & $x^4+x^3 y-y^4+x^2 z^2+y^2 z^2+z^4=0$  \\
		\hline
		$96$ & $x^4+x^3 y+y^4+x^2 z^2+y^2 z^2-z^4=0$  \\
		\hline
		$96$ & $2 x^4-y^4+x y^2 z+y^2 z^2+z^4 = 0$ \\
		\hline
		$96$ & $2 x^4+x^2 y^2-y^4+y z^3+z^4=0$ \\
		\hline
		$96$ & $2 x^4+y^4+x y^2 z+y^2 z^2-z^4=0$ \\
		\hline
		$96$ & $2 x^4+x^3 y+2 y^4-z^4=0$ \\
		\hline		
	\end{tabular}
	\captionof{table}{\label{tab:H96homa} Three-variable homogeneous equations of size $H=96$}
\end{center}

Each of the equations listed in Table \ref{tab:H96homa} have a variable that enters the equation in even exponents only, and we can reduce these to genus $1$ curve, as shown in Section 6.1.3 of the book. If this curve has an easy-to-find rational point, then it is an elliptic curve, for which we can compute the rank in Magma, and, if the rank is $0$, solve the equation. If we cannot find a rational point on our genus $1$ curve, we can try to prove their non-existence by using the following proposition.

\begin{proposition}\label{prop:gen1quarjab}[Proposition 6.6 in the book]
	Let $a,b,c,d,e$ be rational numbers, such that the equation 
	$$
		Y^2 = aX^4 + bX^3+cX^2+dX+e
	$$
	has a rational solution with $Y\neq 0$. Then the equation
	$$
		Y^2 = X^3 + i X + j,
$$
	where
	$$
		i = -27(12ae-3bd+c^2), \quad j=-27(72ace-27ad^2-27b^2e+9bcd-2c^3),
	$$
	has a rational solution as well.
\end{proposition}

Let us first consider the equation 
\begin{equation}\label{eq:x4px3ymy4px2z2py2z2pz4}
x^4+x^3y-y^4+x^2z^2+y^2z^2+z^4=0.
\end{equation}
If $y=0$ then the equation is reduced to $x^4+x^2z^2+z^4=0$, whose only integer solution is $(x,z)=(0,0)$. Now assume that $y \neq 0$. Then by dividing by $y^4$ and making the rational change of variables $X=\frac{x}{y}$ and $Z=\frac{z}{y}$, the equation is reduced to
$$
X^4+X^3-1+X^2Z^2+Z^2+Z^4=0.
$$
We can then make the further change of variable $T=Z^2$, and the equation reduces to 
$$
X^4+X^3-1+X^2T+T+T^2=0.
$$
Multiplying by $4$ and rearranging, we obtain,
\begin{equation}\label{ell:x4px3ymy4px2z2py2z2pz4}
(2T+X^2+1)^2=-3X^4-4X^3+2X^2+5.
\end{equation}
Then letting $Y=2T+X^2+1$ and using Proposition \ref{prop:gen1quarjab}, we have that for any rational solution to \eqref{ell:x4px3ymy4px2z2py2z2pz4}, the equation 
$$
Y^2=X^3+4752X+117072
$$
has a rational solution as well. This equation is an elliptic curve with rank $0$, and all its rational points are torsion points and the SageMath command \eqref{cmd:torsion} returns $[(0 : 1 : 0)]$. Hence, the only integer solution to equation \eqref{eq:x4px3ymy4px2z2py2z2pz4} is the trivial solution $(x,y,z)=(0,0,0)$.

\vspace{10pt}

Let us next consider the equation 
\begin{equation}\label{eq:2x4my4pxy2zpy2z2pz4}
2x^4-y^4+xy^2z+y^2z^2+z^4=0.
\end{equation}
If $z=0$ then the equation is reduced to $2x^4-y^4=0$, whose only integer solution is $(x,y)=(0,0)$. Now assume that $z \neq 0$. Then by dividing by $z^4$ and making the rational change of variables $X=\frac{x}{z}$ and $Y=\frac{y}{z}$, the equation is reduced to
$$
2X^4-Y^4+XY^2+Y^2+1=0.
$$
We can then make the further change of variable $T=Y^2$, and the equation reduces to 
$$
2X^4-T^2+XT+T+1=0.
$$
This curve has the rational point $(X,T)=(1,-1)$, which implies that it is an elliptic curve. As explained in Section \ref{ex:H27rank0ell}, we can compute its rank, which is $1$, hence it has infinitely many rational points.

\vspace{10pt}

Let us next consider the equation 
\begin{equation}\label{eq:x4px3ypy4px2z2py2z2mz4}
x^4+x^3y+y^4+x^2z^2+y^2z^2-z^4=0.
\end{equation}
If $y=0$ then the equation is reduced to $x^4+x^2z^2-z^4=0$, whose only integer solution is $(x,z)=(0,0)$. Now assume that $y \neq 0$. Then by dividing by $y^4$ and making the rational change of variables $X=\frac{x}{y}$ and $Z=\frac{z}{y}$, the equation is reduced to
$$
X^4+X^3+1+X^2Z^2+Z^2-Z^4=0.
$$
We can then make the further change of variable $T=Z^2$, and the equation reduces to 
$$
X^4+X^3+1+X^2T+T-T^2=0.
$$
This curve has the rational point $(X,T)=(1,-1)$ which implies that it is an elliptic curve. As explained in Section \ref{ex:H27rank0ell}, we can compute its rank, which is $1$, hence it has infinitely many rational points.

\vspace{10pt}

Let us next consider the equation 
$$
2x^4+x^2y^2-y^4+yz^3+z^4=0.
$$
If $z=0$ then the equation is reduced to $2x^4+x^2y^2-y^4=0$, whose only integer solution is $(x,y)=(0,0)$. Now assume that $z \neq 0$. Then by dividing by $z^4$ and making the rational change of variables $X=\frac{x}{z}$ and $Y=\frac{y}{z}$, the equation is reduced to
$$
2X^4+X^2Y^2-Y^4+Y+1=0.
$$
We can then make the further change of variable $T=X^2$, and the equation reduces to 
$$
2T^2+TY^2-Y^4+Y+1=0.
$$
This curve has the rational point $(Y,T)=(-1,-1)$ which implies that it is an elliptic curve. As explained in Section \ref{ex:H27rank0ell}, we can compute its rank, which is $1$, hence it has infinitely many rational points.

\vspace{10pt}

The next equation we will consider is
\begin{equation}\label{eq:2x4px3yp2y4mz4}
2x^4+x^3y+2y^4-z^4=0.
\end{equation}
If $y=0$ then the equation is reduced to $2x^4-z^4=0$, whose only integer solution is $(x,z)=(0,0)$. Now assume that $y \neq 0$. Then by dividing by $y^4$ and making the rational change of variables $X=\frac{x}{y}$ and $Z=\frac{z}{y}$, the equation is reduced to
$$
2X^4+X^3+2-Z^4=0.
$$
We can then make the further change of variable $T=Z^2$, and the equation reduces to 
\begin{equation}\label{ell:2x4px3yp2y4mz4}
2X^4+X^3+2-T^2=0.
\end{equation} 
This curve has the rational point $(X,T)=(1/2,3/2)$ which implies that it is an elliptic curve. As explained in Section \ref{ex:H27rank0ell}, we can compute its rank, which is $1$, hence it has infinitely many rational points.

Finally, the equation $2x^4+y^4+xy^2z+y^2z^2-z^4=0$ is solved in Section 6.1.3 of the book.

\subsection{Exercise 6.12}\label{ex:homdia34}
\textbf{\emph{Prove that all equations listed in Table \ref{tab:H800hom4} have no non-trivial integer solutions. }}

	\begin{center}
		\begin{tabular}{ |c|c|c|c|c|c| } 
			\hline
			$H$ & Equation &  $H$ & Equation &  $H$ & Equation \\ 
			\hline\hline
			$320$ & $9x^4+7y^4-4z^4=0$ & $608$ & $21x^4+14y^4-3z^4=0$ & $768$ & $29x^4+11y^4-8z^4=0$ \\ 
			\hline
			$320$ & $11x^4+7y^4-2z^4=0$ & $640$ & $23x^4+9y^4-8z^4=0$ & $800$ & $28x^4+21y^4-z^4=0$ \\ 
			\hline
			$352$ & $17x^4-4y^4-z^4=0$ & $704$ & $23x^4+14y^4-7z^4=0$ & $800$ & $32x^4-11y^4-7z^4=0$ \\ 
			\hline
			$368$ & $14x^4-7y^4+2z^4=0$ & $720$ & $38x^4-6y^4+z^4=0$ & $800$ & $43x^4+6y^4-z^4=0$ \\ 
			\hline
			$480$ & $19x^4+8y^4-3z^4=0$ & $736$ & $29x^4-14y^4+3z^4=0$ & $800$ & $46x^4+3y^4-z^4=0$  \\ 
			\hline	
			$592$ & $17x^4-16y^4-4z^4=0$ & $736$ & $39x^4-4y^4-3z^4=0$ &  & \\ 
			\hline	
		\end{tabular}
		\captionof{table}{\label{tab:H800hom4} Equations of the form \eqref{eq:homdia34} of size $H\leq 800$ having elliptic curve  \eqref{eq:threecurves} of rank $0$}
	\end{center}

Assuming that an integer solution exists to an equation of the form
\begin{equation}\label{eq:homdia34}
ax^4+by^4+cz^4=0
\end{equation}
 such that $(x,y,z) \neq (0,0,0)$, the equation can be reduced to the three elliptic curves 
\begin{equation}\label{eq:threecurves}
	Y^2 = X^3 + abc^2 X, \quad Y^2 = X^3 + ab^2c X \quad \text{and} \quad Y^2 = X^3 + a^2bc X,
\end{equation}
as shown in Section 6.2.2 of the book.
If any of these curves has rank $0$, we can determine all rational points on the elliptic curve and then check if any of these correspond to integer solutions to the original equation. 

For every equation in Table \ref{tab:H800hom4}, Table \ref{tab:H800hom4sol} presents an elliptic curve it reduces to, which has rank $0$, and the rational points on the elliptic curve. For each equation, there are no rational points with $XY \neq 0$, hence the equations have no non-trivial integer solutions.

\begin{center}

\captionof{table}{\label{tab:H800hom4sol} Equations in Table \ref{tab:H800hom4} reduced to an elliptic curve \eqref{eq:threecurves} of rank $0$, and the rational points on this curve.}
\end{center} 
	
\subsection{Exercise 6.16}\label{ex:homdia34all}
\textbf{\emph{Prove that all equations listed in Table \ref{tab:H1000hom4} have no non-trivial integer solutions. }}

	\begin{center}
		%
		\captionof{table}{\label{tab:H1000hom4} Equations of the form \eqref{eq:homdia34} of size $H\leq 1000$ with all elliptic curves \eqref{eq:threecurves} of positive rank}
	\end{center} 

We first recall from Section 6.2.2 of the book that if $(x,y,z)$ is a non-trivial solution to equation \eqref{eq:homdia34} then the triples $(X,Y,Z)$ given by
$$
	(X,Y,Z)=(-b c x y^2, b c^2 y z^2, x^3), \,\, (-a b z x^2, a b^2 x y^2, z^3) \, \text{and} \, (-c a y z^2, c a^2 z x^2, y^3) \,\, 
$$
are the non-trivial solutions to the homogeneous equations
$$
	Y^2 Z = X^3 + abc^2 X Z^2, \quad Y^2 Z = X^3 + ab^2c X Z^2 \quad \text{and} \quad Y^2 Z = X^3 + a^2bc X Z^2,
	$$
respectively.

To prove that equations in Table \ref{tab:H1000hom4} have no non-trivial integer solutions, we will use the method in Section 6.2.2 of the book, which we summarise below for convenience. We first reduce the equation to an elliptic curve with positive rank. Then all the rational points on the curve are obtained from its generators. We then need to find a suitable prime to obtain a contradiction and conclude that the original equation has no non-trivial integer solutions. We will illustrate this method below with a number of examples.

Let us consider, for example, the equation
\begin{equation}\label{eq:38x4p7y4m6z4}
38 x^4+7 y^4-6 z^4 = 0.
\end{equation}
This equation can be reduced to the curve
$$
Y^2Z=X^3+9576XZ^2,
$$
where $X=42xy^2$, $Y=252yz^2$ and $Z=x^3$.
This curve has generators $(0 : 0 : 1)$, $(57 : 855 : 1)$, and $(25 : -505 : 1)$. By adding these points modulo $13$ using Theorem \ref{ref:ellmodpred}, we conclude that all integer points on this elliptic curve must be in one of the families
$$
(0: 1: 0), (0: 0: 1), (12: 2: 1), (5: 10: 1), (12: 11: 1), (5: 3:1),(4:5:1),(1:12:1) 
$$
modulo $13$. On the other hand, solving the original equation \eqref{eq:38x4p7y4m6z4} modulo $13$, applying transformation $(X,Y,Z)=(42xy^2,252yz^2,x^3)$ and removing duplicates, we conclude that $(X,Y,Z)$ must be in one of the families
$$
(3: 5: 1),(3: 8: 1),(10: 1: 1), (10: 12: 1)
$$
modulo $13$. We can see that there is no intersection between these two lists and therefore no non-trivial integer solutions can exist for \eqref{eq:38x4p7y4m6z4}.

\vspace{10pt}

Let us consider another equation, for example,	
\begin{equation}\label{eq:37x4m4y4mz4}
37 x^4-4 y^4-z^4 = 0.
\end{equation}
Using the change of variables $(X,Y,Z)=(-4xy^2,-4yz^2,x^3)$ and $(X,Y,Z)=(148zx^2,592xy^2,z^3)$
we can obtain the elliptic curves in Weierstrass form
\begin{equation}\label{eq:37x4m4y4mz4curves}
Y^2Z=X^3-148XZ^2 \quad \text{and} \quad  Y^2Z=X^3-592XZ^2,
\end{equation}
respectively.
The first curve has generators $(0:0:1)$ and $(2775:-27935:27)$. By adding these generators modulo $23$ using Theorem \ref{ref:ellmodpred}, we conclude that integer points on this elliptic curve must belong to one of the families
\begin{equation}\label{eq:37x4m4y4mz4list}
(0:1:0),(0:0:1),(21:14:1),(5:12:1),(6:15:1),(6:8:1),(5:11:1),(21:9:1)
\end{equation}
modulo $23$. By solving the original equation modulo $23$ and making the transformation of solutions $(X,Y,Z)=(-4xy^2,-4yz^2,x^3)$, we conclude that $(X,Y,Z)$ belong to the families \eqref{eq:37x4m4y4mz4list} only if $(x,y,z)$ belong to one of the following families modulo $23$:
$$
(2: 5: 1), (2: 18: 1),(3: 11: 1),(3: 12: 1),(20: 11: 1),(20: 12: 1),(21: 5: 1),(21: 18: 1).
$$
Then making the transformations $(X,Y,Z)=(148zx^2,592xy^2,z^3)$ to these solutions and removing duplicates, we obtain that $(X,Y,Z)$ must be in one of the families 
$$
(17: 22: 1),(21: 7: 1),(21: 16: 1),(17: 1: 1)
$$ 
modulo $23$. On the other hand, the second curve in \eqref{eq:37x4m4y4mz4curves} has generators: $(0:0:1)$ and $(-4:-48:-1)$. By adding these generators modulo $23$ using Theorem \ref{ref:ellmodpred}, we conclude that integer points on this elliptic curve must belong to one of the families
$$
(0:1:0),(0:0:1),(19:21:1),(10:18:1),(12:11:1),(12:12:1),(10:5:1),(19:2:1)
$$
modulo $23$. We can see that there is no intersection between the last two lists and therefore we have a contradiction and no non-trivial integer solutions can exist for \eqref{eq:37x4m4y4mz4}.

The other equations from Table \ref{tab:H1000hom4} can be solved similarly. 
Table \ref{tab:H1000hom4sol} shows reduction to what elliptic curve we used, and modulo which prime we have obtained a contradiction.

\begin{center}
\begin{tabular}{ |c|c|c|c|c|c| } 
\hline
 Equation &  Elliptic Curves & Prime \\ 
\hline\hline
 $19x^4-4y^4-3z^4=0$ & $Y^2Z = X^3 - 684XZ^2$ and $Y^2Z = X^3 -912XZ^2$ & $11$ \\ \hline	
 $13x^4+11y^4-8z^4 = 0$ & $Y^2Z=X^3 +3146XZ^2$ & $17$ \\ \hline	
 $37 x^4-4 y^4-z^4 = 0$ & $Y^2Z=X^3-148XZ^2$ and $Y^2Z=X^3-592XZ^2$ & $23$ \\ \hline	
 $31 x^4-13 y^4-2 z^4 = 0$ &$Y^2Z=X^3-10478XZ^2$ and $Y^2Z=X^3+24986XZ^2$ & $41$ \\ \hline	
 $31 x^4-11 y^4-5 z^4 = 0$ & $Y^2Z=X^3-18755XZ^2$ and $Y^2Z=X^3+52855XZ^2$ & $3$ \\ \hline	
 $26 x^4+22 y^4-z^4 = 0$ & Reduces to: $13x^4+11y^4-8z^4=0$ &  \\ \hline	
 $38 x^4+7 y^4-6 z^4 = 0$ & $Y^2Z=X^3+9576XZ^2$ & $13$ \\ \hline	
 $35 x^4-10 y^4+7 z^4 = 0$ & $Y^2Z=X^3+24500XZ^2$ & $13$ \\ \hline	
 $34 x^4-19 y^4-2 z^4 = 0$ & $Y^2Z=X^3-24548XZ^2$ & $11$ \\ \hline	
 $43 x^4-11 y^4+z^4 = 0$ & $Y^2Z=X^3-473XZ^2$& $41$ \\ \hline	
 $37 x^4-16 y^4-4 z^4 = 0$ & $Y^2Z=X^3-9472XZ^2$ and $Y^2Z=X^3-37888XZ^2$ & $23$ \\ \hline	
 $39 x^4-13 y^4+6 z^4 = 0$ & $Y^2Z=X^3+39546XZ^2$& $19$ \\ \hline	
 $53 x^4-4 y^4-z^4 = 0$ & $Y^2Z=X^3-848XZ^2$ and $Y^2Z=X^3+11236XZ^2$ & $17$ \\ \hline	
 $29 x^4+17 y^4-14 z^4 = 0$ & $Y^2Z=X^3-117334XZ^2$ and $Y^2Z=X^3-200158XZ^2$ & $11$ \\ \hline	
 $41 x^4-12 y^4+7 z^4 = 0$ & $Y^2Z=X^3+41328XZ^2$& $61$ \\ \hline
 $29 x^4-23 y^4-9 z^4 = 0$ &$Y^2Z=X^3-138069XZ^2$ and $Y^2Z=X^3-174087XZ^2$ & $17$  \\ \hline
 $38 x^4-23 y^4+z^4 = 0$ & $Y^2Z=X^3-874XZ^2$ & $31$ \\ \hline
$43 x^4+12 y^4-7 z^4 = 0$ & $Y^2Z=X^3-43344XZ^2$ and $Y^2Z=X^3-155316XZ^2$ & $13$ \\ \hline
$47 x^4-11 y^4-4 z^4 = 0$ & $Y^2Z=X^3-8272XZ^2$ & $29$  \\ \hline	
\end{tabular}
\captionof{table}{\label{tab:H1000hom4sol} Necessary information to conclude that the equations from Table \ref{tab:H1000hom4} have no non-trivial integer solutions.}
\end{center}

\subsection{Exercise 6.18}\label{ex:H96homb}
\textbf{\emph{Solve all three-variable equations of size $H=96$ from Table \ref{tab:H96homa} that are not covered by Section \ref{ex:H96hom}.  }}

The remaining equations to solve are those in Table \ref{tab:H96hom}. By using the Mordell-Weil Sieve we will show that these equations have no non-trivial integer solutions. To do this, we will use the following theorem.

\begin{theorem}\label{ref:ellmodpred}[Theorem 6.13 in the book]
	Let $a,b$ be integers such that $D=4a^3+27b^2 \neq 0$, let $P_1$ and $P_2$ be rational points on the elliptic curve 
	$$
	Y^2 Z = X^3 + a X Z^2 + b Z^3,
$$
 and let $p$ be a prime that is not a divisor of $D$. Then
	$$
	\psi(P_1 + P_2) = \psi(P_1) + \psi(P_2),
	$$
	where the $+$ on the left-hand side is the addition defined in Section \ref{ex:H32hom3mon}, see formulas \eqref{eq:ellprod1}-\eqref{eq:ellprod2}, while the $+$ on the right-hand side is defined by the same formulas but with arithmetic operations modulo $p$.
\end{theorem}

	\begin{center}
		\begin{tabular}{ |c|c| } 
			\hline
			$H$ & Equation \\ 
			\hline\hline
			$96$ & $x^4+x^3 y+y^4+x^2 z^2+y^2 z^2-z^4=0$  \\ 
			\hline
			$96$ & $2 x^4-y^4+x y^2 z+y^2 z^2+z^4 = 0$ \\
			\hline
			$96$ & $2 x^4+x^3 y+2 y^4-z^4=0$ \\
			\hline		
		\end{tabular}
		\captionof{table}{\label{tab:H96hom} Homogeneous equations of size $H\leq 100$ left open by Section \ref{ex:H96hom}.}
	\end{center}

	Let us first consider the equation $2 x^4+x^3 y+2 y^4-z^4=0$ which is \eqref{eq:2x4px3yp2y4mz4}.
This equation is reducible to
$$
X^3 - 16X Z^2 + 2Z^3 - Y^2 Z = 0
$$
after the change of variables
\begin{equation}\label{eq:2x4px3yp2y4mz4change}
\begin{aligned}
X = & 12 x^3 - 4 x^2 y + 31 x y^2 - 16 y^3 + 24 x z^2 - 12 y z^2,\\ 
Y = & 60 x^3 + 18 x^2 y + 96 y^3 + 14 x z^2 + 65 y z^2, \\
Z = & 8 x^3 - 12 x^2 y + 6 x y^2 - y^3.
\end{aligned}
\end{equation}
The resulting elliptic curve has rank $1$ and generator $P = (-4596 : 8767 : 1728)$. Therefore, all integer solutions to the equation are of the form
$$
P_n=nP, \quad n \in \mathbb{Z}.
$$
Because the prime divisors of $D=4(-16)^3+27(2^2)=-16276$ are $2,13$ and $313$, Theorem \ref{ref:ellmodpred} is applicable to all other primes. If $p=17$, then we conclude that $(X,Y,Z)$ belong to one of the families
$$
(0: 1: 0), (1: 15: 1), (16: 0: 1), (1: 2: 1).
$$
However, solving \eqref{eq:2x4px3yp2y4mz4} modulo $17$ and applying the substitutions \eqref{eq:2x4px3yp2y4mz4change}, we conclude that $(X,Y,Z)$ belong to one of the families
$$
(0: 6: 1),(12: 12: 1),(10: 3: 1),(4: 11: 1),(4: 6: 1),(3: 10: 1),(15: 14: 1),(10: 14: 1).
$$ 
Because these lists have no intersection, the original equation \eqref{eq:2x4px3yp2y4mz4} has no integer solutions. 

\vspace{10pt}

The next equation we will consider is $x^4+x^3y+y^4+x^2z^2+y^2z^2-z^4$ which is \eqref{eq:x4px3ypy4px2z2py2z2mz4}.
This equation is reducible to
$$
X^3 - 131328X Z^2 - 2460672 Z^3 - Y^2 Z=0
$$
after substitutions
$$
X = 816x^3 - 576x^2y + 432xy^2 - 672y^3 - 576xz^2 + 576yz^2,
$$
$$
Y = 28512x^3 + 19872x^2y + 11232xy^2 + 23328y^3 - 15552xz^2 - 12096yz^2,
$$
$$
Z = x^3 - 3x^2y + 3xy^2 - y^3.
$$
In this example, the rank is $1$ and the generator is $P=(-327 : -2349 : 1)$. Therefore, all integer solutions to the equation are of the form
$$
P_n=nP, \quad n \in \mathbb{Z}.
$$
Because the prime divisors of $D=4(-131328)^3+27(-2460672)^2=-8896596478525440$ are $2,3,5,31$ and $103$, Theorem \ref{ref:ellmodpred} is applicable to all other primes. In this example, we were not able to find a single prime to obtain a contradiction, so we need to use several primes, as explained in the proof of Proposition 6.17 in the book. Analysis modulo $p=61$ shows that $(X:Y:Z)=nP$ for $n$ equal to $0,3,4,8,11$ or $14$ modulo $16$. Analysis modulo $p=67$ shows that $n$ should be $2$ or $5$ modulo $8$, that is, $2,5,10$ or $13$ modulo $16$. This gives the desired contradiction.

\vspace{10pt}

Finally, we will consider the equation $2x^4-y^4+xy^2z+y^2z^2+z^4=0$, which is \eqref{eq:2x4my4pxy2zpy2z2pz4}.
This equation is reducible to
$$
X^3 - 851116032 X Z^2 - 912311451648 Z^3 - Y^2 Z=0
$$
after substitutions
$$
X = 37632x^3 - 36864xy^2 - 11520x^2z + 36864y^2z + 20736xz^2 - 46848z^3,
$$
$$
Y = 14155776x^3 - 7962624xy^2 + 4866048x^2z - 6193152y^2z + 10616832xz^2 + 12828672z^3,
$$
$$
Z = x^3 - 3x^2z + 3xz^2 - z^3.
$$
In this example, the rank is $1$ and the generator is $P=(-28608 : -152064 : 1)$. Therefore, all integer solutions to the equation are of the form
$$
P_n=nP, \quad n \in \mathbb{Z}.
$$
Because the prime divisors of 
$$
D=4(-851116032)^3+27(-912311451648)^2=-2443716278390110909485809664
$$
are $2,3$ and $31907$, Theorem \ref{ref:ellmodpred} is applicable to all other primes. Analysis modulo $p=67$ shows that $(X:Y:Z)=nP$ for $n$ equal to $2$ or $4$ modulo $8$. Analysis modulo $p=109$ shows that $n$ should be $0,1$ or $3$ modulo $8$. This gives us the desired contradiction.

\subsection{Exercise 6.19}\label{ex:H704hom5}
\textbf{\emph{Solve all equations listed in Table \ref{tab:H704hom5}. }}

	\begin{center}
		\begin{tabular}{ |c|c|c|c|c|c| } 
			\hline
			$H$ & Equation &  $H$ & Equation &  $H$ & Equation \\ 
			\hline\hline
			$352$ & $6x^5+3y^5+2z^5=0$ & $512$ & $7x^5+5y^5+4z^5=0$ & $608$ & $14x^5+3y^5+2z^5=0$ \\ 
			\hline
			$352$ & $7x^5+3y^5+z^5=0$ & $512$ & $9x^5+5y^5+2z^5=0$ & $608$ & $15x^5+3y^5+z^5=0$ \\ 
			\hline
			$384$ & $7x^5+4y^5+z^5=0$ & $512$ & $13x^5+2y^5+z^5=0$ & $640$ & $12x^5+7y^5+z^5=0$ \\ 
			\hline
			$384$ & $9x^5+2y^5+z^5=0$ & $544$ & $10x^5+6y^5+z^5=0$ & $640$ & $13x^5+5y^5+2z^5=0$ \\ 
			\hline
			$416$ & $6x^5+5y^5+2z^5=0$ & $544$ & $14x^5+2y^5+z^5=0$ & $640$ & $15x^5+4y^5+z^5=0$ \\ 
			\hline
			$416$ & $7x^5+5y^5+z^5=0$ & $576$ & $7x^5+6y^5+5z^5=0$ & $672$ & $9x^5+7y^5+5z^5=0$ \\
			\hline 
			$416$ & $10x^5+2y^5+z^5=0$ & $576$ & $10x^5+7y^5+z^5=0$ & $672$ & $10x^5+9y^5+2z^5=0$ \\
			\hline 
			$448$ & $6x^5+5y^5+3z^5=0$ & $576$ & $12x^5+5y^5+z^5=0$ & $672$ & $12x^5+5y^5+4z^5=0$ \\ 
			\hline
			$448$ & $10x^5+3y^5+z^5=0$ & $576$ & $14x^5+3y^5+z^5=0$ & $672$ & $13x^5+5y^5+3z^5=0$ \\ 
			\hline
			$480$ & $10x^5+3y^5+2z^5=0$ & $608$ & $12x^5+4y^5+3z^5=0$ & $672$ & $14x^5+5y^5+2z^5=0$ \\ 
			\hline	
		\end{tabular}
		\captionof{table}{\label{tab:H704hom5} Homogeneous three-monomial equations of degree $d\geq 5$ and size $H<704$.}
	\end{center}

Table \ref{tab:H700hom5sol} reduces equations from Table \ref{tab:H704hom5} to hyperelliptic curves using a rational change of variables. If the original equation has solution $(x,y,z)\neq (0,0,0)$, then the reduced one has a rational point with $X\neq 0$. We then calculate the bounds for the rank of the Jacobian of the reduced equation. If the upper bound is exactly $0$, we can use the \emph{Chabauty0} command, if the upper bound of the rank is $1$, we can use the \emph{Chabauty} command in Magma (see Section \ref{ex:H32genus2} for example code).
Table \ref{tab:H700hom5sol} displays the finite rational points found by the \emph{Chabauty0} and \emph{Chabauty} commands. There are four equations which have a solution with $X \neq 0$, however, when transforming these solutions to $(x,y,z)$, we have not obtained rational/integer solutions. Hence, all equations in Table \ref{tab:H704hom5} have no non-trivial integer solutions.

\begin{center}

\captionof{table}{\label{tab:H700hom5sol} 
	Equations from Table \ref{tab:H704hom5}, the corresponding hyperelliptic curves, bounds for their ranks, and Chabauty output.}
\end{center}

In the next exercises we will solve the following problem.

\begin{problem}\label{prob:homlarge}[Problem 6.22 in the book]
	Given a homogeneous polynomial $P(x_1,\dots,x_n)$ with integer coefficients, determine whether equation $P=0$ has an integer solution such that all variables are different from $0$.
\end{problem}

\subsection{Exercise 6.23}\label{ex:H48genus1}

\textbf{\emph{Prove that the equations listed in Table \ref{tab:H48genus1} have no integer solutions satisfying $xyz\neq 0$.}}

	\begin{center}
		\begin{tabular}{ |c|c|c|c|c|c| } 
			\hline
			$H$ & Equation & $H$ & Equation & $H$ & Equation \\ 
			\hline\hline
			$48$ & $x^3 y+x^2 z^2-y^2 z^2=0$ & $48$ & $x^4+x^2 y^2-z^4=0$ & $48$ & $x^4+x y^3+y^2 z^2=0$ \\ 
			\hline
			$48$ & $x^3 y+x^2 z^2+y^2 z^2=0$ & $48$ & $x^4+y^3 z-y^2 z^2=0$  & $48$ & $x^4-x^2 y^2+z^4=0$ \\ 
			\hline
			$48$ & $x^4-x^2 y^2-z^4=0$ & $48$ & $x^4+y^3 z+y^2 z^2=0$ & $48$ & $x^3y + x^3z + y^2z^2=0$ \\ 
			\hline		
		\end{tabular}
		\captionof{table}{\label{tab:H48genus1} Quartic equations of size $H\leq 48$ defining genus $1$ curves}
	\end{center} 
	
	To prove that the equations listed in Table \ref{tab:H48genus1} have no integer solutions satisfying $xyz \neq 0$, we will use the method from Section 6.3.2 of the book, which we summarise below for convenience. First, we reduce the equations to elliptic curves with rank $0$. As these curves have a finite number of rational points, we can use these to show that none of these points correspond to integer solutions to the original equation satisfying $xyz \neq 0$.
	
	We remark that the equation $x^3y+x^2z^2+y^2z^2=0$ is not included as it is solved in Section 6.3.2 of the book.

Let us try to find an integer solution to the equation
\begin{equation}\label{eq:x3ypx2z2my2z2}
x^3y+x^2z^2-y^2z^2=0,
\end{equation}
with $xyz \neq 0$.
We can divide the equation by $z^4$ and make the change of variables $X=x/z$ and $Y=y/z$, which reduces the equation to
\begin{equation}\label{red:x3ypx2z2my2z2}
X^3Y+X^2-Y^2=0,
\end{equation}
and we are interested in finding solutions with $XY \neq 0$.
The Maple command \eqref{maple:wei}
$$
{\tt Weierstrassform(X^3* Y + X^2 - Y^2, X, Y, u, v)}
$$
returns that this equation reduces to 
\begin{equation}\label{ell:x3ypx2z2my2z2}
v^2=u^3-u
\end{equation}
with the change of variables $u=\frac{X^3-Y}{X}$ and $v=X^3-Y$. This is an elliptic curve and the SageMath commands \eqref{cmd:rank}
$$
{\tt EllipticCurve([0, 0, 0,-1,0]).rank()}
$$
and \eqref{cmd:torsion}
$$
{\tt EllipticCurve([0, 0, 0,-1,0]).torsion\_points()}
$$
return ${\tt 0}$ and ${\tt [(-1 : 0 : 1), (0 : 0 : 1), (0 : 1 : 0), (1 : 0 : 1)]}$, respectively, so the finite rational solutions to equation \eqref{ell:x3ypx2z2my2z2} are $(u,v)=(\pm 1,0),(0,0)$. The first solution is impossible for $(X,Y)$. The second solution corresponds to $X^3-Y=v=0$ or $Y=X^3$. Substituting this into \eqref{red:x3ypx2z2my2z2} results in $X^6+X^2-X^6=0$ or $X^2=0$ which is impossible for $X \neq 0$. Hence equation \eqref{eq:x3ypx2z2my2z2} has no integer solutions with $xyz \neq 0$.

\vspace{10pt}

Let us try to find a solution with $xyz \neq 0$ to
\begin{equation}\label{eq:x4mx2y2mz4}
x^4-x^2y^2-z^4=0.
\end{equation}
We can divide the equation by $z^4$ and make the change of variables $X=x/z$ and $Y=y/z$, which reduces the equation to
\begin{equation}\label{red:x4mx2y2mz4}
X^4-X^2Y^2-1=0,
\end{equation}
and we are interested in finding solutions with $XY \neq 0$.
The Maple command \eqref{maple:wei} for this equation 
returns that the equation can be reduced to 
\begin{equation}\label{ell:x4mx2y2mz4}
v^2=u^3+4u
\end{equation}
with the change of variables $u=2X^2 - 2XY$ and $v=4 X^3-4 X^2 Y$. This is an elliptic curve and the SageMath commands \eqref{cmd:rank} and \eqref{cmd:torsion} for this equation 
return ${\tt 0}$ and ${\tt [(0 : 0 : 1), (0 : 1 : 0), (2 : -4 : 1), (2 : 4 : 1)]}$ respectively, so the finite rational solutions to \eqref{ell:x4mx2y2mz4} are $(u,v)=(0,0),(2,\pm 4)$. The first solution corresponds to $2X^2 - 2XY=u=0$ so either $X=0$ or $Y=X$. Substituting $Y=X$ into \eqref{red:x4mx2y2mz4} results in $X^4-X^4-1=0$ or $1=0$ which is impossible. The second solution $(u,v)$ corresponds to $2X(X-Y)=2$ and $4X^2(X-Y)=\pm 4$ so $X= \pm 1$. But then $(X,Y)=(\pm 1,0)$.  Hence equation \eqref{eq:x4mx2y2mz4} has no integer solutions with $xyz \neq 0$.

\vspace{10pt}

Let us try to find a solution with $xyz \neq 0$ to
\begin{equation}\label{eq:x4px2y2mz4}
x^4+x^2y^2-z^4=0.
\end{equation}
We can divide the equation by $z^4$ and make the change of variables $X=x/z$ and $Y=y/z$, which reduces the equation to
$$
X^4+X^2Y^2-1=0,
$$
and we are interested in finding solutions with $XY \neq 0$.
The Maple command \eqref{maple:wei} for this equation 
returns that the equation reduces to \eqref{ell:x4mx2y2mz4} with the change of variables $u=-\frac{2(X-1)}{X+1}$ and $v=\frac{4XY}{(X+1)^2}$, and its finite rational solutions are $(u,v)=(0,0),(2,\pm 4)$. The first solution corresponds to $X=0$ or $Y=0$. The second solution corresponds to $-\frac{2(X-1)}{X+1}=2$ so $X= 0$. Hence equation \eqref{eq:x4px2y2mz4} has no integer solutions with $xyz \neq 0$.

\vspace{10pt}

Let us try to find a solution with $xyz \neq 0$ to
\begin{equation}\label{eq:x4py3zmy2z2}
x^4+y^3z-y^2z^2=0.
\end{equation}
We can divide the equation by $z^4$ and make the change of variables $X=x/z$ and $Y=y/z$, which reduces the equation to
$$
X^4+Y^3-Y^2=0,
$$
and we are interested in finding solutions with $XY \neq 0$.
The Maple command \eqref{maple:wei} for this equation 
returns that the equation reduces to \eqref{ell:x3ypx2z2my2z2} with the change of variables $u=\frac{Y(Y-1)}{X^2}$ and $v=X$, and its finite rational solutions are $(u,v)=(\pm 1,0),(0,0)$. These solutions have $v=X=0$. Hence equation \eqref{eq:x4py3zmy2z2} has no integer solutions with $xyz \neq 0$.

\vspace{10pt}

Let us try to find a solution with $xyz \neq 0$ to
\begin{equation}\label{eq:x4py3zpy2z2}
x^4+y^3z+y^2z^2=0.
\end{equation}
We can divide the equation by $z^4$ and make the change of variables $X=x/z$ and $Y=y/z$, which reduces the equation to
$$
X^4+Y^3+Y^2=0,
$$
and we are interested in finding solutions with $XY \neq 0$.
The Maple command \eqref{maple:wei} for this equation 
returns that the equation reduces to 
$$
v^2=u^3+u
$$
with the change of variables $u=\frac{Y(Y+1)}{X^2}$ and $v=X$, and its only finite rational solution is $(u,v)=(0,0)$. This solution has $v=X=0$. Hence equation \eqref{eq:x4py3zpy2z2} has no integer solutions with $xyz \neq 0$.

\vspace{10pt}

Let us try to find a solution with $xyz \neq 0$ to
\begin{equation}\label{eq:x4pxy3py2z2}
x^4+xy^3+y^2z^2=0.
\end{equation}
We can divide the equation by $x^4$ and make the change of variables $X=-\frac{y}{x}$ and $Y=\frac{yz}{x^2}$, which reduces the equation to
\begin{equation}\label{red:x4pxy3py2z2}
Y^2=X^3-1
\end{equation}
and we are interested in finding solutions with $XY \neq 0$.
This is an elliptic curve and the SageMath commands \eqref{cmd:rank} and \eqref{cmd:torsion} for this equation 
return ${\tt 0}$ and ${\tt [(0 : 1 : 0), (1 : 0 : 1)]}$, respectively, so the only finite rational solution to \eqref{red:x4pxy3py2z2} is $(X,Y)=(1,0)$, which is impossible with $Y \neq 0$. Hence equation \eqref{eq:x4pxy3py2z2} has no integer solutions with $xyz \neq 0$.

\vspace{10pt}

Let us try to find a solution with $xyz \neq 0$ to
\begin{equation}\label{eq:x4mx2y2pz4}
x^4-x^2y^2+z^4=0
\end{equation}
We can divide the equation by $z^4$ and make the change of variables $X=x/z$ and $Y=y/z$, which reduces the equation to
\begin{equation}\label{red:x4mx2y2pz4}
X^4-X^2Y^2+1=0,
\end{equation}
and we are interested in finding solutions with $XY \neq 0$.
The Maple command \eqref{maple:wei} for this equation 
returns that the equation reduces to 
\begin{equation}\label{ell:x4mx2y2pz4}
v^2=u^3-4u
\end{equation}
with the change of variables $u=2X^2 - 2XY$ and $v=4X^3 - 4X^2Y$. This is an elliptic curve and the SageMath commands \eqref{cmd:rank} and \eqref{cmd:torsion} for this equation 
return ${\tt 0}$ and ${\tt [(-2 : 0 : 1), (0 : 0 : 1), (0 : 1 : 0), (2 : 0 : 1)]}$, respectively, so the finite rational solutions to \eqref{ell:x4mx2y2pz4} are $(u,v)=(0,0),(\pm 2,0)$. Both solutions have $v=0$ so $X=0$ or $X=Y$. Substituting $Y=X$ into \eqref{red:x4mx2y2pz4} we obtain $1=0$ which is impossible. Hence equation \eqref{eq:x4mx2y2pz4} has no integer solutions with $xyz \neq 0$.

\vspace{10pt}

Let us try to find a solution with $xyz \neq 0$ to
\begin{equation}\label{eq:x3ypx3zpy2z2}
x^3y+x^3z+y^2z^2=0.
\end{equation}
We can divide the equation by $z^4$ and make the change of variables $X=x/z$ and $Y=y/z$, which reduces the equation to
\begin{equation}\label{red:x3ypx3zpy2z2}
X^3Y+X^3+Y^2=0,
\end{equation}
and we are interested in finding solutions with $XY \neq 0$.
This equation reduces to 
\begin{equation}\label{ell:x3ypx3zpy2z2}
v^2=u^3+16
\end{equation}
with the change of variables $u=\frac{4(X^3+Y)}{X}$ and $v=8X^3+8Y-4$. This is an elliptic curve and the SageMath commands \eqref{cmd:rank} and \eqref{cmd:torsion} 
return ${\tt 0}$ and ${\tt [(0 : -4 : 1), (0 : 1 : 0), (0 : 4 : 1)]}$, respectively, so the finite rational solutions to \eqref{ell:x3ypx3zpy2z2} are $(u,v)=(0,\pm 4)$. This solution corresponds to $u=0$ so $X^3+Y=0$ or $Y=-X^3$. Substituting $Y=-X^3$ into \eqref{red:x3ypx3zpy2z2} we obtain $-X^6+X^3+X^6=0$ so $X^3=0$ which is impossible with $X \neq 0$. Hence equation \eqref{eq:x3ypx3zpy2z2} has no integer solutions with $xyz\neq 0$.

\subsection{Exercise 6.24}\label{ex:H48genus2rank0}

\textbf{\emph{Reduce Problem \ref{prob:homlarge} for equations
	$$
		x^4-y^4+x y z^2 = 0, \quad x^4-x^2 y^2+y z^3 = 0, \quad \text{and} \quad x^4+x y^2 z+y z^3 = 0
	$$
	to the problem of finding rational points with non-zero coordinates on a hyperelliptic curve with rank of Jacobian $0$. Then use either the {\tt Chabauty0} command or LMFDB to check that there are no such points. }}
	
The first equation we will consider is 
\begin{equation}\label{eq:x4my4pxyz2}
x^4-y^4+xyz^2=0.
\end{equation}
We are looking for solutions with $xyz \neq 0$, so by dividing the equation by $z^4$ and making the rational change of variables $X=x/z$ and $Y=y/z$, we can reduce the equation to
$$
X^4-Y^4+XY=0.
$$
The Maple command \eqref{maple:wei} for this equation 
suggests to use the change of variables $u=Y/X$ and $v=Y/X^2$ to reduce the equation to $v^2=u^5-u.$ This equation has rank of the Jacobian $0$ and using the {\tt Chabauty0} command from Section \ref{ex:H32genus2} we find that all rational points are $(u,v)=(0,0),(\pm 1,0)$. 
As there are no finite rational points with $v \neq 0$, there are no rational points with $Y \neq 0$. This implies that equation \eqref{eq:x4my4pxyz2} has no integer solutions with $xyz \neq 0$.

\vspace{10pt}

The next equation we will consider is 
\begin{equation}\label{eq:x4mx2y2pyz3}
x^4-x^2 y^2+yz^3=0.
\end{equation}
We are looking for solutions with $xyz \neq 0$, so by dividing the equation by $z^4$ and making the rational change of variables $X=x/z$ and $Y=y/z$, we can reduce the equation to
$$
X^4-X^2 Y^2+Y=0.
$$
The Maple command \eqref{maple:wei} for this equation 
suggests to use the change of variables $u=X$ and $v=-2X^2Y+1$ to reduce the equation to $v^2=4u^6+1$. This equation has rank of the Jacobian $0$ and using the {\tt Chabauty0} command from Section \ref{ex:H32genus2} we find that all finite rational points are $(u,v)=(0,\pm 1)$. As there are no rational points with $u \neq 0$, there are no rational points with $X \neq 0$. This implies that equation \eqref{eq:x4mx2y2pyz3} has no integer solutions with $xyz \neq 0$.

\vspace{10pt}

The next equation we will consider is 
\begin{equation}\label{eq:x4pxy2zpyz3}
x^4+x y^2 z+yz^3=0.
\end{equation}
We are looking for solutions with $xyz \neq 0$, so dividing by $z^4$ and making the rational change of variables $X=x/z$ and $Y=y/z$ we can reduce this equation to
$$
X^4+X Y^2+Y=0.
$$
The Maple command \eqref{maple:wei} for this equation 
suggests to use the change of variables $u=-X$ and $v=2XY+1$ to reduce the equation to $v^2=4u^5+1$. This equation has rank of the Jacobian $0$ and using the {\tt Chabauty0} command from Section \ref{ex:H32genus2} we find that all finite rational points are $(u,v)=(0,\pm 1)$. As there are no rational points with $u \neq 0$, there are no rational points with $X \neq 0$. This implies that equation \eqref{eq:x4pxy2zpyz3} has no integer solutions with $xyz \neq 0$.

\subsection{Exercise 6.25}\label{ex:H64genus2}
\textbf{\emph{By using either the {\tt Chabauty} command described in Section \ref{ex:H32genus2} or LMFDB, check that the equations listed in Table \ref{tab:H64genus2} have no integer solutions such that $xyz \neq 0$.}}

	\begin{center}
		\begin{tabular}{ |c|c|c|c|c|c| } 
			\hline
			$H$ & Equation & $H$ & Equation & $H$ & Equation \\ 
			\hline\hline
			$48$ & $x^4+y^4+xyz^2=0$ & $64$ & $x^4+x^3y-y^4+xyz^2=0$ & $64$ & $x^4-x^2y^2+y^2z^2+z^4=0$ \\ 
			\hline
			$48$ & $x^4+x^2y^2+yz^3=0$ & $64$ & $x^4-x^2y^2+x^2yz+yz^3=0$ & $64$ & $x^4+x^3y+y^2z^2+z^4=0$ \\ 
			\hline
			$64$ & $x^3y+x^2y^2+x^3z-yz^3=0$ & $64$ & $x^4+x^3y+xy^2z-yz^3=0$ &  &  \\ 
			\hline		
		\end{tabular}
		\captionof{table}{\label{tab:H64genus2} Equations of size $H\leq 64$ defining genus $2$ curves with rank of the Jacobian $1$.}
	\end{center}

We remark that the equation $x^4+y^4+xyz^2=0$ is not included in this section as it has been solved in Section 6.3.2 of the book.

The first equation we will consider is
\begin{equation}\label{eq:x4px2y2pyz3}
 x^4+x^2y^2+yz^3=0.
 \end{equation}
 We are looking for solutions with $xyz\neq 0$ so by dividing by $x^4$, this equation can be reduced to $XY^3 + X^2 + 1=0$ with the rational change of variables $X=y/x$ and $Y=z/x$, such that $XY \neq 0$. Then the Maple command \eqref{maple:wei} for this equation
 $$
 {\tt Weierstrassform(X*Y^3 + X^2 + 1, X, Y, u, v);}
 $$
 suggests to use the change of variables $u=Y$ and $v=Y^3+2X$ to reduce the equation to
\begin{equation}\label{red:x4px2y2pyz3}
v^2=u^6-4.
 \end{equation}
 The Magma code
 \begin{equation}\label{hypercode:6.25}
 \begin{aligned}
 & {\tt > P < u >:= PolynomialRing(Rationals());}
\\ & {\tt > C := HyperellipticCurve(u^6-4); }
\\ & {\tt > J := Jacobian(C);}
\\ & {\tt > RankBounds(J);}
\\ & {\tt > ptsJ := Points(J : Bound := 10); }
\\ & {\tt > [Order(P) : P\,\, in\,\, ptsJ];}
\end{aligned}
 \end{equation}
 outputs ${\tt 1 \,\, 1 \quad [ 1, 3, 3, 0, 0, 0, 0, 0, 0 ]}$. This suggests to use $n=4$ in the code 
\begin{equation}\label{hypercode3:6.25}
 \begin{aligned}
 & {\tt > P := ptsJ[n];}
\\ & {\tt > allptsC := Chabauty(P); allptsC;}
\end{aligned}
 \end{equation}
 which outputs that \eqref{red:x4px2y2pyz3} has no finite rational points. Therefore equation \eqref{eq:x4px2y2pyz3} has no integer points with $xyz \neq 0$. 
 
 \vspace{10pt}
 
 The next equation we will consider is
\begin{equation}\label{eq:x3ypx2y2px3zmyz3}
 x^3y+x^2y^2+x^3z-yz^3=0. 
 \end{equation}
 We are looking for solutions with $xyz\neq 0$ so dividing by $z^4$, this equation can be reduced to $X^3Y+X^2Y^2+X^3-Y=0$ with the rational change of variables $X=x/z$ and $Y=y/z$, such that $XY \neq 0$. Then the Maple command \eqref{maple:wei} for this equation 
  suggests to use the change of variables $u=X$ and $v=X^3 + 2X^2Y - 1$ to reduce the equation to
\begin{equation}\label{red:x3ypx2y2px3zmyz3}
v^2=u^6 - 4u^5 - 2u^3 +1.
 \end{equation}
 The Magma code \eqref{hypercode:6.25} with the second line 
 $$
 {\tt C := HyperellipticCurve(u^6 - 4*u^5 - 2*u^3 + 1); }
 $$
outputs ${\tt 1\,\, 1 \quad [ 1, 0, 0, 0, 0 ]}$. Then the code \eqref{hypercode3:6.25} with $n=2$ 
outputs 
$$
{\tt \{ (0 : -1 : 1), (0 : 1 : 1), (1 : -1 : 0), (1 : 1 : 0) \}}
$$ 
so the finite rational points to \eqref{red:x3ypx2y2px3zmyz3} are $(u,v)=(0,\pm 1)$, which implies $X=0$. Therefore equation \eqref{eq:x3ypx2y2px3zmyz3} has no integer points with $xyz \neq 0$.

\vspace{10pt} 
 
 The next equation we will consider is
\begin{equation}\label{eq:x4px3ymy4pxyz2}
 x^4+x^3y-y^4+xyz^2=0.
 \end{equation}
  We are looking for solutions with $xyz\neq 0$ so dividing by $z^4$, this equation can be reduced to $X^4+X^3Y-Y^4+XY=0$ with the rational change of variables $X=x/z$ and $Y=y/z$, such that $XY \neq 0$. 
  Then the Maple command \eqref{maple:wei} for this equation 
   suggests to use the change of variables $u=Y/X$ and $v=Y/X^2$ to reduce the equation to
\begin{equation}\label{red:x4px3ymy4pxyz2}
v^2=u^5 - u^2-u.
 \end{equation}
  The Magma code \eqref{hypercode:6.25} with the second line 
  $$
  {\tt > C := HyperellipticCurve(u^5 -u^2-u); }
  $$
outputs that the only finite rational solution to \eqref{red:x4px3ymy4pxyz2} is $(u,v)=(0,0)$, which implies $Y=0$. Therefore equation \eqref{eq:x4px3ymy4pxyz2} has no integer points with $xyz \neq 0$. 

\vspace{10pt}

 The next equation we will consider is
\begin{equation}\label{eq:x4mx2y2px2yzpyz3}
 x^4-x^2y^2+x^2yz+yz^3=0. 
 \end{equation}
   We are looking for solutions with $xyz\neq 0$ so dividing by $z^4$, this equation can be reduced to $X^4 - X^2Y^2 + X^2Y + Y=0$ with the rational change of variables $X=x/z$ and $Y=y/z$, such that $XY \neq 0$. 
  Then the Maple command \eqref{maple:wei} for this equation 
   suggests to use the change of variables $u=X$ and $v=-2X^2Y + X^2 + 1$ to reduce the equation to
\begin{equation}\label{red:x4mx2y2px2yzpyz3}
v^2=4u^6+u^4+2u^2+1.
 \end{equation}
  The Magma code \eqref{hypercode:6.25} with the second line
  $$
  {\tt > C := HyperellipticCurve(4*u^6+u^4+2*u^2+1); }
  $$
outputs ${\tt 1\,\, 1 \quad [ 1, 0, 0 ]}$. Then the code \eqref{hypercode3:6.25} with $n=2$ 
outputs that the only finite rational solution to \eqref{red:x4mx2y2px2yzpyz3} is $(u,v)=(0,\pm 1)$, which implies $X=0$. Therefore equation \eqref{eq:x4mx2y2px2yzpyz3} has no integer points with $xyz \neq 0$.
 
 \vspace{10pt}
 
 The next equation we will consider is
\begin{equation}\label{eq:x4px2y2my2z2pz4}
x^4+x^2y^2-y^2z^2+z^4=0. 
\end{equation}
 We are looking for solutions with $xyz\neq 0$ so dividing by $y^4$, this equation can be reduced to 
 \begin{equation}\label{wei:x4px2y2my2z2pz4}
 X^4 + Z^4 + X^2 - Z^2=0
 \end{equation}
 with the rational change of variables $X=x/y$ and $Z=z/y$, such that $XZ \neq 0$. 
  Then the Maple command \eqref{maple:wei} for this equation 
   suggests to use the change of variables $u=-\frac{X}{X-Z}$, and $v= -\frac{Z + X}{(X-Z)^2}$ to reduce the equation to
\begin{equation}\label{red:x4px2y2my2z2pz4}
v^2=4u^5+10u^4+16u^3+14u^2+6u+1.
 \end{equation}
   The Magma code \eqref{hypercode:6.25} with the second line
   $$
   {\tt > C := HyperellipticCurve(4*u^5+10*u^4+16*u^3+14*u^2+6*u+1);  }
   $$
outputs ${\tt 1\,\, 1 \quad [ 1, 0, 0 ]}$. Then the code \eqref{hypercode3:6.25} with $n=2$ 
outputs that the only finite rational solutions to \eqref{red:x4px2y2my2z2pz4} are $(u,v)=(0,\pm 1),(-1/2,0)$, the first solution implies $X=0$. The second solution has $v=0$ which implies $X=-Z$, substituting this into \eqref{wei:x4px2y2my2z2pz4} we must have $2X^4=0$ or $X=0$. Therefore equation \eqref{eq:x4px2y2my2z2pz4} has no integer points with $xyz \neq 0$.

\vspace{10pt}

The next equation we will consider is
\begin{equation}\label{eq:x4px3ypxy2zmyz3}
x^4+x^3y+xy^2z-yz^3=0.
\end{equation}
 We are looking for solutions with $xyz\neq 0$ so dividing by $z^4$, this equation can be reduced to 
 $$
 X^4 + X^3Y + XY^2 - Y=0
$$
 with the rational change of variables $X=x/z$ and $Y=y/z$, such that $XY \neq 0$. 
  Then the Maple command \eqref{maple:wei} for this equation 
    suggests to use the change of variables $u=X$, and $v=X^3 + 2XY - 1$ to reduce the equation to
\begin{equation}\label{red:x4px3ypxy2zmyz3}
v^2=u^6-4u^5-2u^3+1.
 \end{equation}
   The Magma code \eqref{hypercode:6.25} with the second line
   $$
   {\tt > C := HyperellipticCurve(u^6-4*u^5-2*u^3+1);}
   $$
outputs ${\tt 1\,\, 1 \quad [ 1, 0, 0, 0, 0 ]}$. Then the code \eqref{hypercode3:6.25} with $n=2$ 
outputs that the only finite rational solutions to \eqref{red:x4px3ypxy2zmyz3} are $(u,v)=(0,\pm 1)$, which implies $X=0$. Therefore equation \eqref{eq:x4px3ypxy2zmyz3} has no integer points with $xyz \neq 0$.

\vspace{10pt}

The final equation we will consider is
\begin{equation}\label{eq:x4px3ypy2z2pz4}
x^4+x^3y+y^2z^2+z^4=0.
\end{equation}
 We are looking for solutions with $xyz\neq 0$ so dividing by $z^4$, this equation can be reduced to 
 $$
 X^4 + X^3Y + Y^2 + 1=0
$$
 with the rational change of variables $X=x/z$ and $Y=y/z$, such that $XY \neq 0$. 
  Then the Maple command \eqref{maple:wei} for this equation 
    suggests to use the change of variables $u=X$, and $v=X^3 + 2Y$ to reduce the equation to
\begin{equation}\label{red:x4px3ypy2z2pz4}
v^2=u^6-4u^4-4.
 \end{equation}
   The Magma code \eqref{hypercode:6.25} with the second line
   $$
   {\tt > C := HyperellipticCurve(u^6-4*u^4-4); }
   $$
outputs ${\tt 1\,\, 1 \quad [ 1, 0, 0 ]}$. Then the code\eqref{hypercode3:6.25} with $n=2$ 
outputs that \eqref{red:x4px3ypy2z2pz4} has no finite rational solutions. Therefore equation \eqref{eq:x4px3ypy2z2pz4} has no integer points with $xyz \neq 0$.

\subsection{Exercise 6.26}\label{ex:H48x2t}
\textbf{\emph{Prove that the equations
	$$
		x^4+y^3z-z^4=0, \quad x^4+y^3z-yz^3=0, \quad \text{and} \quad x^4+xy^3+z^4=0
	$$ 
	have no integer solutions with $xyz\neq 0$.}}

The first equation we will consider is
\begin{equation}\label{eq:x4py3zmz4}
x^4+y^3z-z^4=0.
\end{equation}
We are looking for solutions such that $xyz \neq 0$, or equivalently, solutions to
\begin{equation}\label{eq:x4py3zmz4red}
	X^4+Y^3-1=0
\end{equation}
in rational variables $X,Y$ such that $XY \neq 0$. This is a genus $3$ curve. However, the substitution $X^2=T$ reduces the equation to
$$
	T^2=(-Y)^3+1.
$$
This equation defines an elliptic curve, and it is easy to check with Magma (see Section \ref{ex:H32hom3mon}), that is has rank $0$ and its only finite torsion points are $(Y,T)=(1,0),(0,\pm 1),(-2,\pm 3)$. Hence, the only rational points on \eqref{eq:x4py3zmz4red} are $(X,Y)=(0,1),(\pm 1,0)$, hence the original equation \eqref{eq:x4py3zmz4} has no integer solutions with $xyz\neq 0$.

\vspace{10pt}

The next equation we will consider is 
\begin{equation}\label{eq:x4py3zmyz3}
	x^4+y^3 z-y z^3=0.
\end{equation}
We are looking for solutions such that $xyz \neq 0$, or equivalently, solutions to
\begin{equation}\label{eq:x4py3zpyz3red}
	X^4+Y^3-Y=0
\end{equation}
in rational variables $X,Y$ such that $XY \neq 0$. This is a genus $3$ curve. However, the substitution $X^2=T$ reduces the equation to
$$
	T^2=(-Y)^3-(-Y).
$$
This equation defines an elliptic curve, and it is easy to check with Magma, that is has rank $0$ and its only finite torsion points are $(Y,T)=(0,0),(\pm 1,0)$. Hence, the only rational points on \eqref{eq:x4py3zpyz3red} are $(X,Y)=(0,0),(0,\pm 1)$, hence the original equation \eqref{eq:x4py3zmyz3} has no integer solutions with $xyz\neq 0$.

\vspace{10pt}

The final equation we will consider is 
\begin{equation}\label{eq:x4pxy3pz4}
	x^4+xy^3+ z^4=0.
\end{equation}
We are looking for solutions such that $xyz \neq 0$, or equivalently, solutions to
\begin{equation}\label{eq:x4pxy3pz4red}
	Z^4+Y^3+1=0
\end{equation}
in rational variables $Y,Z$ such that $YZ \neq 0$. This is a genus $3$ curve. However, the substitution $Z^2=T$ reduces the equation to
$$
	T^2=(-Y)^3-1.
$$
This equation defines an elliptic curve, and it is easy to check with Magma, that is has rank $0$ and its only finite torsion points are $(Y,T)=(-1,0)$. Hence, the only rational point on \eqref{eq:x4pxy3pz4red} is $(Z,Y)=(0,-1)$, hence the original equation \eqref{eq:x4pxy3pz4} has no integer solutions with $xyz\neq 0$.

\subsection{Exercise 6.29}\label{ex:H48hom4var}
\textbf{\emph{Prove that the equations listed in Table \ref{tab:H64redto3var} have no integer solutions with $xyzt\neq 0$.}}

	\begin{center}
		\begin{tabular}{ |c|c|c|c|c|c| } 
			\hline
			$H$ & Equation & $H$ & Equation & $H$ & Equation \\ 
			\hline\hline
			$48$ & $x^4-y^4- z^2t^2=0$ & $64$ & $x^2 y^2+x^2 z t-y^2 z t+z^2 t^2=0$ & $64$ & $2 x^3 y+x y^3+z^2 t^2=0$ \\ 
			\hline
			$48$ & $x^4+y^4- z^2t^2=0$ & $64$ & $x^2 y^2+x^2 z t+ y^2 z t-t^2 z^2=0$ & $64$ & $x^4+y z t^2+2 y^2 z^2=0$ \\ 
			\hline
			$48$ & $x^4+t^2 y z-y^2 z^2 =0$ & $64$ & $x^3 y+x y z t+ y^2 z t-z^2 t^2=0$ & $64$ &  $x^4+x^2 y^2-y^4+z^2 t^2=0$ \\ 
			\hline
			$48$ & $x^4+t^2 y z+y^2 z^2 =0$ & $64$ & $x^3 y+x^2 y^2-x y^3+z^2 t^2=0$ & $64$ & $x^4+2 y z t^2-y^2 z^2=0$ \\ 
			\hline		
			$48$ & $x^3y-xy^3 + z^2t^2=0$ & $64$ & $x^3 y+x^2 y^2+x y^3-z^2 t^2=0$ & $64$ & $x^4+x^2 y^2+y^2 z t+z^2 t^2=0$ \\ 
			\hline
			$48$ & $x^3y+xy^3 + z^2t^2=0$ & $64$ & $x^4-y z t^2+x^2 y z+y^2 z^2=0$ & $64$ & $x^4+x^2 y^2+x y^3+z^2 t^2=0$ \\ 
			\hline
			$48$ & $x^4+xy^3+z^2t^2=0$ & $64$ & $x^4-x^2 y^2-y^4-z^2 t^2=0$ & $64$ & $x^4+x^3 y-x y^3+z^2 t^2=0$ \\ 
			\hline
			$64$ & $x^4-y^4-2 z^2 t^2=0$ & $64$ & $x^2 y^2+2 x y z t-z^2 t^2=0$ & $64$ & $x^4+x^2 y^2+y^4-z^2 t^2=0$ \\ 
			\hline	
			$64$ & $x^4+x y^3+2 z^2 t^2=0$  & $64$ & $x^2 y^2+x y z^2+x y t^2-z^2 t^2=0$ & $64$ & $2 x^4+y z t^2+y^2 z^2=0$ \\ 
			\hline
			$64$ & $x^4+2 x y^3-z^2 t^2=0$  & $64$ & $x^2 y^2-x y t^2+x y z^2+z^2 t^2=0$ & $64$ & $x^3 y-x y^3+2 z^2 t^2=0$ \\ 
			\hline
			$64$ & $2 x^4+y^4-z^2 t^2=0$ & $64$ & $x^4+y z t^2+x^2 y z-y^2 z^2=0$ & &  \\ 
			\hline						
		\end{tabular}
		\captionof{table}{\label{tab:H64redto3var} Four-variable equations of size $H\leq 64$ reducible to three-variable ones.}
	\end{center} 

We first consider equation
\begin{equation}\label{eq:x4my4mz2t2}
	x^4 - y^4 - z^2t^2 = 0.
\end{equation}
Making the substitution $s=zt$, we can reduce \eqref{eq:x4my4mz2t2} to the three-variable equation
$$
x^4-y^4=s^2.
$$
Dividing this equation by $y^4$ and making the rational change of variables $X=x/y$ and $Y=s/y^2$ where $XY \neq 0$, we obtain
$$
X^4-1=Y^2.
$$ 
The Maple command
$$
{\tt Weierstrassform(X^4 - 1 - Y^2, X, Y, u, v)}
$$
suggests to use the substitutions $u=\frac{2(X - 1)}{X + 1}$ and $v=-\frac{4Y}{X^2 + 2X + 1}$ to reduce the equation to
$$
v^2=u^3 + 4u.
$$
This is an elliptic curve in Weierstrass form with rank $0$ and its only finite rational solutions are $(u,v)=(0,0),(2,\pm 4)$.
The first solution implies that $Y=0$. The second solution $u=\frac{2(X - 1)}{X + 1}=2$ is impossible in $X$ as this reduces to $2=0$. Therefore, \eqref{eq:x4my4mz2t2} has no integer solutions with $xyzt \neq 0$. 

All other equations listed in Table \ref{tab:H64redto3var} can be solved similarly: we can (i) reduce the original equation to a three-variable equation in integer variables after a suitable substitution; (ii) reduce the resulting equation to an equation in two rational variables $(X,Y)$; (iii) reduce the equation further to Weierstrass form in variables $(u,v)$, and (iv) check that the resulting equation has rank $0$ and find all its rational points. Table \ref{tab:H64redto3varsol} presents, for each equations from Table \ref{tab:H64redto3var}, all steps in this algorithm. 
Excluding the equation $v^2=u^3+4u$ we have solved above, all of the equations in the final column of Table \ref{tab:H64redto3varsol} have no non-trivial rational points, hence the original equations have no non-trivial integer solutions.

\begin{sidewaystable}
	\begin{center}

		\captionof{table}{\label{tab:H64redto3varsol} Four-variable equations of size $H\leq 64$ reducible to three-variable ones. The two-variable equations are produced either by solving when the reduced equation has a discriminant which is a perfect square, or by dividing by the fourth power of a variable. Excluding the highlighted equation, the equations in Weierstrass form have Rank $0$ and have no non-zero torsion points.  }
	\end{center} 
\end{sidewaystable}

\subsection{Exercise 6.30}\label{ex:H64hom4var}
\textbf{\emph{Prove that the equations listed in Table \ref{tab:H64red2varrat} have no integer solutions with $xyzt\neq 0$.}}

	\begin{center}
		\begin{tabular}{ |c|c|c|c|c|c| } 
			\hline
			$H$ & Equation & $H$ & Equation & $H$ & Equation \\ 
			\hline\hline
			$48$ & $x^2 y^2+y z t^2-x^2 z^2 =0$ & $64$ & $x^2 y^2+x^2 z^2+y^2t^2-z^2 t^2=0$ & $64$ & $x^2 y^2-y z t^2+x y z t+ x z^2 t=0$ \\ 
			\hline
			$48$ & $x^2 y^2+ y z t^2+x^2 z^2 =0$ & $64$ & $x^2 y z+xy^2t+x z^2 t - yzt^2=0$ & $64$ & $x^2 y^2+x^2 y z-y z t^2+z^2 t^2=0$ \\ 
			\hline
			$64$ & $x^2 y^2+2 y z t^2 -x^2 z^2=0$ & $64$ & $x^2 y^2+x^2 y z + y z t^2 -x^2 z^2=0$ & $64$ &  $x^2 y^2+x^2 y z+y^2t^2+z^2 t^2=0$ \\ 
			\hline
			$64$ & $2 x^2 y^2+y z t^2+x^2 z^2=0$ & $64$ & $x^2 y^2+x^2 y z+y z t^2-z^2 t^2=0$ & $64$ & $x^2 y^2+x^2 y z+x^2 z^2-y z t^2=0$ \\ 
			\hline					
		\end{tabular}
		\captionof{table}{\label{tab:H64red2varrat} Four-variable equations of size $H\leq 64$ reducible to curves.}
	\end{center} 

For each equation  listed in Table \ref{tab:H64red2varrat}, Table \ref{tab:H64red2varratsol} presents the necessary substitutions to reduce the equation to two-variables, and this equation in Weierstrass form. Excluding the highlighted equations, all of the equations in the final column of Table \ref{tab:H64red2varratsol} have no non-trivial rational points, hence the original equations have no non-trivial integer solutions.

Let us consider the two equations that have a highlighted Weierstrass form equation in more detail.
The first equation we will consider is
\begin{equation}\label{eq:x2y2px2z2py2t2mz2t2}
x^2y^2+x^2z^2+y^2t^2-z^2t^2=0.
\end{equation}
As we are looking for solutions with $xyzt \neq 0$ we can divide by $x^2y^2$ and make the rational change of variables $X=t/x$ and $Y=z/y$ to reduce the equation to
\begin{equation}\label{eq:x2y2px2z2py2t2mz2t2red}
1+Y^2+X^2-X^2Y^2=0.
\end{equation}
The Maple command 
$$
{\tt Weierstrassform(1+Y^2+X^2-X^2Y^2, X, Y, u, v)}
$$
suggests to make the change of variables $u=\frac{2(X+1)}{X-1}$ and $v=\frac{4Y(X+1)}{X-1}$ to reduce the equation to
$$
v^2=u^3+4u.
$$
This equation is an elliptic curve with rank $0$ and its only finite rational points are $(u,v)=(0,0),(2,\pm 4)$. The first point suggests $X=-1$, substituting this into \eqref{eq:x2y2px2z2py2t2mz2t2red} we see that this does not lead to rational solutions. The second point $u=\frac{2(X+1)}{X-1}=2$, leads to $2=0$ which is impossible. Therefore, \eqref{eq:x2y2px2z2py2t2mz2t2} has no integer solutions with $xyzt \neq 0$.

\vspace{10pt}

The next equation we will consider is
\begin{equation}\label{eq:x2y2myzt2pxyztpxz2t}
x^2y^2-yzt^2+xyzt+xz^2t=0.
\end{equation}
As we are looking for solutions with $xyzt \neq 0$ we can divide by $x^2y^2$ and make the rational change of variables $X=t/x$ and $Y=z/y$ to reduce the equation to
$$
1-X^2Y+XY+XY^2=0.
$$
The Maple command 
$$
{\tt Weierstrassform(1-X^2Y+XY+XY^2, X, Y, u, v)}
$$
suggests to make the change of variables $u=\frac{X-12}{12X}$ and $v=\frac{2Y-X+1}{2X}$ to reduce the equation to
$$
v^2=u^3+\frac{23u}{48}+\frac{181}{864}.
$$
This equation is an elliptic curve with rank $0$ and its only finite rational points are $(u,v)=\left(\frac{1}{12},\pm \frac{1}{2} \right)$. However, $u=\frac{X-12}{12X}=\frac{1}{12}$, leads to $12=0$ which is impossible. Therefore, \eqref{eq:x2y2myzt2pxyztpxz2t} has no integer solutions with $xyzt \neq 0$.

	\begin{center}
		\begin{tabular}{ |c|c|c|c|c|c|c| } 
			\hline
			 Equation & Two-Variable & X & Y & Weierstrass form  \\ 
			\hline\hline
			 $x^2 y^2+y z t^2-x^2 z^2 =0$ & $Y^2+YX^2-1=0$ & $t/x$ & $y/z$&  $v^2=u^3-u$ \\\hline
			 $x^2 y^2+ y z t^2+x^2 z^2 =0$ & $Y^2+YX^2+1=0$ & $t/x$ & $y/z$&  $v^2=u^3+u$ \\\hline
			 $x^2 y^2+2 y z t^2 -x^2 z^2=0$ & $Y^2+2X^2Y-1=0$ & $t/x$ & $y/z$ & $v^2=u^3-4u$ \\\hline
			 $2 x^2 y^2+y z t^2+x^2 z^2=0$ & $2+X^2Y+Y^2=0$ & $t/x$ & $z/y$ & $v^2=u^3+2u$ \\\hline
			$x^2 y^2+x^2 z^2+y^2t^2-z^2 t^2=0$ & $1+Y^2+X^2-X^2Y^2=0$ & $t/x$ & $z/y$ & \hl{$v^2=u^3+4u$} \\\hline
			 $x^2 y z+xy^2t+x z^2 t - yzt^2=0$ & $Y+X+XY^2-X^2Y=0$ & $t/x$ & $z/y$ & $v^2=u^3-u$ \\\hline
			  $x^2 y^2+x^2 y z + y z t^2 -x^2 z^2=0$ & $Y^2+Y+X^2 Y-1=0$ & $t/x$ & $y/z$ & $v^2=u^3 - \frac{4u}{3}  -\frac{11}{27} $ \\\hline
			  $x^2 y^2+x^2 y z+y z t^2-z^2 t^2=0$ & $1+Y+X^2Y-X^2Y^2=0$ & $t/x$ & $z/y$ & $v^2=u^3-u$ \\\hline
			 $x^2 y^2-y z t^2+x y z t+ x z^2 t=0$ & $1-X^2Y+XY+XY^2=0$ & $t/x$ & $z/y$ & \hl{$v^2=u^3 + \frac{23u}{48} +\frac{181}{864} $} \\ 
			\hline
			  $x^2 y^2+x^2 y z-y z t^2+z^2 t^2=0$ & $X^2 Y^2+X^2Y-Y+1=0$ & $x/t$ & $y/z$ & $v^2=u^3-u$ \\ 
			\hline
			 $x^2 y^2+x^2 y z+y^2t^2+z^2 t^2=0$ & $X^2 Y^2+X^2Y+Y^2+1=0$ & $x/t$ & $y/z$& $v^2=u^3 + \frac{32u}{3} - \frac{1280}{27}$\\ \hline 
			$x^2 y^2+x^2 y z+x^2 z^2-y z t^2=0$ & $1+Y+Y^2-X^2Y=0$ & $t/x$ & $z/y$ & $v^2=u^3 + \frac{2u}{3}+ \frac{7}{27}$  \\ 
			\hline					
		\end{tabular}
		\captionof{table}{\label{tab:H64red2varratsol} Equations listed in Table \ref{tab:H64red2varrat} reduced to elliptic curves in Weierstrass form.}
	\end{center} 
	
\subsection{Exercise 6.31} \label{ex:H64hom}
\textbf{\emph{By either direct search or any kind of reasoning, find an integer solution with $xyzt\neq 0$ to each of the equations 
	$$
		x^4-y^4+z^3t+zt^3 = 0, \quad x^4+x^2 y^2+y z^3+t^4 = 0, \quad x^4-x^2 y^2+z^4+t^4 = 0,
	$$
	as well as an integer solution with $xyzts \neq 0$ to each of the equations
	$$
		x^4+y^3z-yz^3+t^2s^2 =0, \quad x^4+y^4+z^2t^2+zts^2 = 0.
	$$}}
	
	For each equation listed above, Table \ref{tab:H64red2varratintsol} presents an integer solution with $xyzt \neq 0$ for four-variable equations, and $xyzts \neq 0$ for five-variable equations.
	\begin{center}
		\begin{tabular}{ |c|c|c|c|c|c|c| } 
			\hline
			 Equation & Solution $(x,y,z,t)$ or $(x,y,z,t,s)$  \\ 
			\hline\hline
 $x^4-y^4+z^3t+zt^3 = 0$ & $(22,124,17,240)$ \\\hline
 $x^4+x^2 y^2+y z^3+t^4 = 0$ & $(3,146,-11,-7)$ \\\hline
 $x^4-x^2 y^2+z^4+t^4 = 0$ & $(144,481,180,240)$ \\\hline
$x^4+y^3z-yz^3+t^2s^2 =0$ & $(6,-17,1,1,60)$ \\\hline
$x^4+y^4+z^2t^2+zts^2 = 0$ & $(20,28,1,-512,45)$ \\\hline						
		\end{tabular}
		\captionof{table}{\label{tab:H64red2varratintsol} An integer solution with $xyzt \neq 0$ or $xyzts \neq 0$ for each equation listed in Exercise \ref{ex:H64hom}.}
	\end{center} 

\subsection{Exercise 6.38}\label{ex:dif4power}
\textbf{\emph{\begin{itemize}
		\item[(a)] Produce a complete list of positive integers up to $n\leq 200$ that are the differences of two non-zero rational forth powers.
		\item[(b)] For each of the equations listed in Table \ref{tab:kx4mky4mz4}, determine whether it has an integer solution such that $xyz\neq 0$.
		\item[(c)] List all equations of the form 
		\begin{equation}\label{eq:kx4mky4mz4}
			kx^4-ky^4=mz^4
		\end{equation}
	 	of size $H \leq 608$ that have an integer solution such that $xyz \neq 0$.
	\end{itemize}}}
	
		\begin{center}

		\captionof{table}{\label{tab:kx4mky4mz4} Non-trivial equations of the form \eqref{eq:kx4mky4mz4} of size $H\leq 608$.}
	\end{center}

	\subsubsection{Exercise 6.38 (a)}

The following positive integers $n \leq 200$ can be expressed as a difference of two non-zero rational fourth powers: 
$$
	5,15,34,39,65,80,84,111,145,150,175.
$$ 
The representations for the listed integers can be found in Table \ref{Tab:x4my4maz4sols}.

\begin{center}
%
\captionof{table}{\label{Tab:x4my4maz4sols} Equations of the form \eqref{eq:x4my4mnz4a} with non-trivial solutions.}
\end{center}

To prove that other positive integers $n\leq 200$ are not representable in this way, we need to prove that the corresponding equations 
\begin{equation}\label{eq:x4my4mnz4a}
x^4-y^4=nz^4
\end{equation} 
have no integer solutions such that $z\neq 0$. For most values of $n$, at least one of the elliptic curves 
\begin{equation}\label{eq:threecurvesa}
Y^2=X^3-n^2 X, \quad Y^2=X^3-nX, \quad \text{and} \quad Y^2=X^3+nX
\end{equation}
 has rank $0$, and the statement is easy to check, see Section 6.3.4 of the book for details. All values of $n\leq 200$ for which this method does not work are listed in Tables \ref{Tab:x4my4maz4nosolsell} and \ref{Tab:x4my4maz4nosolsmws}. 
In this case, as explained in Section 6.3.4 of the book, equation \eqref{eq:x4my4mnz4a} reduces to a finite number of equations of the form
\begin{equation}\label{eq:a2u8mb2v8mcw4}
a^2u^8+b^2v^8+cw^4=0,
\end{equation}
where $a,b,c$ are pairwise coprime positive integers such that $abc=2n$, and we are looking for pairwise coprime solutions $(u,v,w)$. Most of the resulting equations \eqref{eq:a2u8mb2v8mcw4} have local obstructions, that is, for some prime $p$ and positive integer $r$ all solutions to the equation modulo $p^r$ have at least two of $u,v,w$ divisible by $p$, which is a contradiction to the pairwise co-primality of $u,v,w$.

All equations \eqref{eq:a2u8mb2v8mcw4} with no local obstructions are listed in Tables \ref{Tab:x4my4maz4nosolsell} and \ref{Tab:x4my4maz4nosolsmws}. 
Some of these equations can be solved by considering the corresponding homogeneous equations $a^2U^4+b^2V^4+cw^4=0$, where $(U,V)=(u^2,v^2)$ and proving that at least one of the elliptic curves \eqref{eq:threecurves} has rank $0$. In our case, the curves \eqref{eq:threecurves} are
\begin{equation}\label{eq:x4my4mnz4threecurvesabc}
Y^2=X^3+a^2b^2c^2X, \quad Y^2=X^3+a^2b^4cX, \quad Y^2=X^3+a^4b^2cX.
\end{equation} 
Values of $n$ for which all equations \eqref{eq:a2u8mb2v8mcw4} are solvable by this method are collected in Table \ref{Tab:x4my4maz4nosolsell}.

For some equations \eqref{eq:a2u8mb2v8mcw4}, all the ranks of the curves \eqref{eq:x4my4mnz4threecurvesabc} are positive. These equations can be solved using the Mordell-Weil Sieve, as in the proof of Proposition 6.36 in the book. Values of $n$ for which the Mordell-Weil Sieve was used are collected in Table \ref{Tab:x4my4maz4nosolsmws}.

\begin{center}

	\captionof{table}{\label{Tab:x4my4maz4nosolsell} Equations of the form \eqref{eq:x4my4mnz4a} with no non-trivial solutions, for which all equations \eqref{eq:a2u8mb2v8mcw4} with no local obstructions reduce to a curve \eqref{eq:x4my4mnz4threecurvesabc} of rank $0$.}
\end{center}

\begin{center}
%
	\captionof{table}{\label{Tab:x4my4maz4nosolsmws} Equations of the form \eqref{eq:x4my4mnz4a} with no non-trivial solutions, for which the proof requires the Mordell-Weil sieve.}
\end{center}

Let us consider two examples in detail, one from each of the Tables \ref{Tab:x4my4maz4nosolsell} and \ref{Tab:x4my4maz4nosolsmws}.

The first equation we will consider is \eqref{eq:x4my4mnz4a} with $n=20$, or, 
$$
x^4-y^4=20z^4.
$$
For this equation, the elliptic curves \eqref{eq:threecurvesa} take the form
$$
Y^2=X^3-400X, \quad Y^2=X^3-20X, \quad Y^2=X^3+20X,
$$
and they all have rank $1$. Hence, this method does not solve the equation, and we need to consider equations \eqref{eq:a2u8mb2v8mcw4}. 
All triples $(a,b,c)$ in \eqref{eq:a2u8mb2v8mcw4} are $(1,40,1),(1,1,40),(8,5,1),(8,1,5),(5,1,8)$. For triples $(1,1,40)$, $(8,5,1)$, $(8,1,5)$ and  $(5,1,8)$ analysis modulo $16$ shows that at least two of $u,v$ and $w$ must be even, which contradicts the pairwise coprimality assumption. The final case to consider is $(1,40,1)$, that is, equation
\begin{equation}\label{red:x4my4m20z4}
u^8+40^2v^8=w^4.
\end{equation}
In this case, we have not discovered a contradiction modulo any prime $p$. Instead, we can reduce this equation to 
$$
w^4-U^4=100V^4
$$
where $U=u^2$ and $V=2v^2$. This is equation \eqref{eq:x4my4mnz4a} with $n=100$, and one of the elliptic curves \eqref{eq:x4my4mnz4threecurvesabc}  
$$
Y^2=X^3-100X
$$
is a rank $0$ curve. This curve has no finite rational points such that $XY \neq0$, and so the equation \eqref{red:x4my4m20z4} has no non-trivial integer solutions. 

\vspace{10pt}

The next equation we will consider is \eqref{eq:x4my4mnz4a} with $n=21$, or, 
\begin{equation}\label{eq:x4my4m21z4}
x^4-y^4=21z^4.
\end{equation}
For this equation, the elliptic curves \eqref{eq:threecurvesa} take the form
$$
Y^2=X^3-441X, \quad Y^2=X^3-21X, \quad Y^2=X^3+21X
$$
and they all have rank $1$. Hence, this method does not solve the equation, and we need to consider equations \eqref{eq:a2u8mb2v8mcw4}.

All triples $(a,b,c)$ in \eqref{eq:a2u8mb2v8mcw4} are $(2,3,7)$, $(2,7,3)$, $(3,7,2)$, $(6,1,7)$, $(7,6,1)$, $(7,1,6)$, $(1,1,42)$, $(42,1,1)$, $(21,1,2)$, $(1,2,21)$, and $(21,2,1)$. For triples $(2,3,7)$, $(2,7,3)$, $(3,7,2)$, $(6,1,7)$, $(7,1,6)$, $(1,1,42)$, $(21,1,2)$ and $(21,2,1)$ analysis modulo 16 shows that at least two of $u,v$ and $w$ must be even, which contradicts the pairwise coprimality assumption. For the triple $(1,2,21)$, analysis modulo $3$ shows that $u$ and $v$ must be divisible by $3$, again a contradiction. The final cases to consider are $(42,1,1)$ and $(7,6,1)$, these are equations 
\begin{equation}\label{red:x4my4m21z4}
(a) \quad 42^2u^8+v^8=w^4, \quad \text{and} \quad (b) \quad 7^2 u^8+6^2v^8=w^4.
\end{equation}
In these cases, we have not discovered a contradiction modulo any prime $p$, hence, we must use an alternative method. For equation \eqref{red:x4my4m21z4} (a) the second elliptic curve \eqref{eq:x4my4mnz4threecurvesabc} takes the form 
$$
Y^2=X^3-1764X
$$
which has rank $0$ and no torsion points satisfying $XY \neq 0$. However, for equation \eqref{red:x4my4m21z4} (b), the three elliptic curves \eqref{eq:x4my4mnz4threecurvesabc} take the form
$$
Y^2=X^3+1764X, \quad Y^2=X^3-63504X, \quad Y^2=X^3-86436X,
$$ 
which all have positive rank. Therefore, we need to use the Mordell-Weil sieve, as explained in the proof of Proposition 6.36 in the book. Our computer search returns that it works modulo $p=157$, which proves there are no non-zero integer points to the equation \eqref{red:x4my4m21z4} (b). This finishes the proof that the original equation \eqref{eq:x4my4m21z4} has no integer solutions with $z \neq 0$.

\subsubsection{Exercise 6.38 (b)} 

We will show that all equations in Table \ref{tab:kx4mky4mz4} have no integer solutions with $xyz \neq 0$. This can be proved exactly as in part (a). The details are presented in Tables \ref{Tab:kx4mky4mmz4nosolsell} and \ref{Tab:kx4mky4mmz4nosolsmws}.

\begin{center}

\captionof{table}{\label{Tab:kx4mky4mmz4nosolsell} Equations from Table \ref{tab:kx4mky4mz4} with no non-trivial solutions, for which all equations \eqref{eq:a2u8mb2v8mcw4} with no local obstructions reduce to a curve \eqref{eq:x4my4mnz4threecurvesabc} of rank $0$.}
\end{center}

\begin{center}
%
\captionof{table}{\label{Tab:kx4mky4mmz4nosolsmws} Equations from Table \ref{tab:kx4mky4mz4} with no non-trivial solutions, for which the proof requires the Mordell-Weil sieve.}
\end{center}

\subsubsection{Exercise 6.38 (c)} 

Table \ref{Tab:kx4mky4mmz4sols} lists equations of the form \eqref{eq:kx4mky4mz4} of size $H\leq 608$ that have integer solutions with $xyz \neq 0$. This happens if and only if rational number $\frac{m}{k}$ can be represented as a difference of two non-zero rational fourth powers. 

\begin{center}
	\begin{tabular}{|c|c|c|c|}
		\hline
		Equation & Solution & Representation \\\hline \hline
		$x^4-y^4= 5z^4$ & $(3,1,2)$ & $5=\left(\frac{3}{2}\right)^4 - \left(\frac{1}{2}\right)^4$ \\\hline
		$x^4-y^4= 15z^4$ & $(2,1,1)$ & $15=2^4-1^4$  \\\hline
		$5(x^4-y^4)=11 z^4$ & $(13,9,10)$ & $\frac{11}{5}=1.3^4 - 0.9^4$  \\\hline
		$4(x^4-y^4)= 21z^4$ & $(93,51,60)$ & $\frac{21}{4}=\left(\frac{93}{60}\right)^4 - \left(\frac{51}{60}\right)^4$ \\\hline
		$8(x^4-y^4)= 17z^4$ & $(10,6,8)$ & $\frac{17}{8}=\left(\frac{5}{4}\right)^4 - \left(\frac{3}{4}\right)^4$  \\\hline
		$x^4-y^4=34 z^4$ & $(5,3,2)$ & $34=\left(\frac{5}{2}\right)^4 - \left(\frac{3}{2}\right)^4$ \\\hline
		$16(x^4-y^4)=5 z^4$ & $(6,2,8)$ & $\frac{5}{16}=\left(\frac{3}{4}\right)^4 - \left(\frac{1}{4}\right)^4$  \\\hline
		
	\end{tabular}
	\captionof{table}{\label{Tab:kx4mky4mmz4sols} Equations of the form \eqref{eq:kx4mky4mz4} of size $H \leq 608$ that have integer solutions with $xyz \neq 0$.}
\end{center}

\subsection{Exercise 6.39}\label{ex:dif6pow}
\textbf{\emph{Solve Problem \ref{prob:homlarge} for all equations of the form 
		\begin{equation}\label{eq:kx6mky6mz6}
			kx^6-ky^6=mz^6
		\end{equation}
	of size $H\leq 640$.}}

Equations of the form \eqref{eq:kx6mky6mz6} have size $H\leq 640$ if and only if $2|k|+|m|\leq 10$. By multiplying \eqref{eq:kx6mky6mz6} by $-1$ and swapping $x$ and $y$ if necessary, we may assume that $k$ and $m$ are coprime positive integers, so that $0<m \leq 10-2k$, which implies that $0<k\leq 4$. In each case, equation \eqref{eq:kx6mky6mz6} reduces to 
\begin{equation}\label{eq:x6my6mnz6}
	x^6-y^6=nz^6
\end{equation}
with $n = mk^5$. As a result, we need to solve \eqref{eq:x6my6mnz6} for 
\begin{equation}\label{eq:nvalues}
	n=1,2,3,4,5,6,7,8,32,96,160,243,486,972,\, \text{and} \, 1024.
\end{equation}

We can reduce equations of the form \eqref{eq:x6my6mnz6} to the hyperelliptic curves
\begin{equation}\label{hyper:x6my6mnz6}
	(a) \,\, Y^2=4X^6+n^2  \quad \text{and} \quad (b) \,\, Y^2 = 4 n X^6 + 1
\end{equation}
where $(X, Y) = (xy/z^2, 2(x/z)^6 -n)$ in (a) and $(X,Y)=(xz/y^2, -2n(z/y)^6-1)$ in (b), see Section 6.3.4 of the book.   
For $1\leq n \leq 8$ and $n=32,160, 972, 1024$, the rank of the Jacobian of the curve (a) is at most $1$, and its rational points can be computed by the {\tt Chabauty0} or {\tt Chabauty} Magma commands, as explained in Section \ref{ex:H32genus2}. For the remaining values $n=96,243,486$, the rank of the Jacobian of the curve (b) is at most $1$, and {\tt Chabauty} Magma commands is applicable to it. In all cases, the command outputs no finite rational points with $X \neq 0$. Therefore equations of the form \eqref{eq:x6my6mnz6} with values of $n$ given by \eqref{eq:nvalues} have no integer solutions with $xyz \neq 0$.

 \begin{center}
\begin{tabular}{|c|c|c|c|c|}
	\hline
	Equation & \eqref{hyper:x6my6mnz6} & Rank Bounds & Rational Points  \\\hline \hline
 $x^6-y^6=z^6$ & $Y^2=4X^6+1$ & $0 \,\, 0$ & $(0,\pm 1)$ \\\hline
 $x^6-y^6=2z^6$ & $Y^2=4X^6+4$ & $0\,\, 0$ & $(0,\pm 2)$ \\\hline
 $x^6-y^6=3z^6$ & $Y^2=4X^6+9$ & $0\,\, 0$ & $(0,\pm 3)$ \\\hline
 $x^6-y^6=4z^6$ & $Y^2=4X^6+16$ & $0\,\, 0$ &$(0,\pm 4)$ \\\hline
 $x^6-y^6=5z^6$ & $Y^2=4X^6+25$ & $0\,\, 0$ & $(0,\pm 5)$ \\\hline
$x^6-y^6=6z^6$ & $Y^2=4X^6+36$ & $1 \,\, 1$ & $(0,\pm 6)$ \\\hline
$x^6-y^6=7z^6$ & $Y^2=4X^6+49$ & $0 \,\, 1$ & $(0,\pm 7)$ \\\hline
 $x^6-y^6=8z^6$ & $Y^2=4X^6+64$ & $0\,\, 0$ & $(0,\pm 8)$ \\\hline
\end{tabular}
 \captionof{table}{The necessary information to conclude that equations of the form \eqref{eq:x6my6mnz6} with $1\leq n \leq 8$ have no integer solutions satisfying $xyz \neq 0$.}
\end{center}

\section{Chapter 7}
In this chapter we are interested in determining whether integer solutions to an equation exist, or more formally, we will be solving the following problem.
\begin{problem}\label{prob:main}[Problem 7.1 in the book]
	Given a Diophantine equation, determine whether it has an integer solution.
\end{problem}

\subsection{Exercise 7.3}\label{ex:H17inthasse}
\textbf{\emph{Using methods described in Section 1.3 or in Section 1.4 in the book, show that all equations from Table \ref{tab:H17inthasse} that are linear in $y$ have no integer solutions.}}

	\begin{center}
		\begin{tabular}{ |c|c|c|c|c|c| } 
			\hline
			$H$ & Equation & $H$ & Equation & $H$ & Equation \\ 
			\hline\hline
			$13$ & $x^2y+2y+1=0$ & $17$ & $x^2y+4y+1=0$ & $17$ & $y^2=x^3-5$ \\ 
			\hline
			$15$ & $x^2y+3y+1=0$ & $17$ & $x^2y+2y+5=0$  & $17$ & $y^2+y=x^3-3$ \\ 
			\hline
			$15$ & $y^2=x^3-3$ & $17$ & $x^2y+2y+x^2+1=0$ & $17$ & $y^2=x^3-x^2-1$ \\ 
			\hline
			$16$ & $x^2y+3y+2=0$ & $17$ & $y^2=x^3-x-3$ & $17$ & $y^2 - x^2 y + z^2 + 1 = 0$ \\ 
			\hline
			$17$ & $x^2y-4y+1=0$ & $17$ & $y^2=x^3+x-3$ &  &  \\ 
			\hline			
		\end{tabular}
		\captionof{table}{\label{tab:H17inthasse} Counterexamples to the integral Hasse principle of size $H\leq 17$}
	\end{center}

Equation 
$$
x^2y+2y+x^2+1=0
$$ 
implies that $y=- \frac{x^2+1}{x^2+2}=-1+\frac{1}{x^2+2}$, hence $0 < y+1 \leq 1/2$, which implies that the equation has no integer solutions. For all other equations listed in Table \ref{tab:H17inthasse} that are linear in $y$, Table \ref{tab:H17inthassesol} rewrites them in the form \eqref{eq:PQc}, specifically in the form $y(x^2+a)=c$, which implies that $x^2+a$ is a divisor of $c$, but easy case analysis shows that this is impossible.

\begin{center}

 \captionof{table}{\label{tab:H17inthassesol} Equations of the form \eqref{eq:PQc} from Table \ref{tab:H17inthasse} with no integer solutions.}
\end{center}

\subsection{Exercise 7.5}\label{ex:H17wei}
\textbf{\emph{Use SageMath or Magma to verify that all equations listed in Table \ref{tab:H17wei} have no integer solutions. Investigate for which equations you can prove the non-existence of integer solutions by a direct argument, without computer assistance.}}

	\begin{center}
		%
		\captionof{table}{\label{tab:H17wei} Almost Weierstrass equations \eqref{eq:Weiformgen} of size $H\leq 19$ with no integer solutions}
	\end{center} 

Let us start with the computer-assisted solution. To solve an equation of the form
$$
y^2 + a x y + c y = x^3 + b x^2 + d x + e ,
$$
we will use the SageMath command 
$$
{\tt EllipticCurve([a,b,c,d,e]).integral\_points();}
$$ 
The values of $a,b,c,d,e$ are presented in Table \ref{tab:H17weiellcommand}. In all cases, the command outputs ${\tt [ \, \, ]}$ which shows that these equations do not have integer solutions.

	\begin{center}
		%
		\captionof{table}{\label{tab:H17weiellcommand} Almost Weierstrass equations \eqref{eq:Weiformgen} of size $H\leq 19$ with no integer solutions.}
	\end{center} 

Let us now present, for all equations in Table \ref{tab:H17wei} except the ones in bold, proofs, without computer assistance, that they have no integer solutions. In some of these proofs, we will need the following result.

\begin{proposition}\label{prop:755}[Proposition 7.62 in the book]
Let $a$, $b$ be positive integers such that $ab = x^2 + y^2$ for some integers $x$, $y$. If either $a$ or $b$ is equal to $3$ modulo $4$, then all integers $a$, $b$, $x$, $y$ must have a common prime divisor $p$ equal to $3$ modulo $4$.
\end{proposition}
	
The non-existence of integer solutions to the equation
$$
y^2=x^3-3
$$
is proven in Section 7.1.2 of the book.
	
The first equation we will consider is	
$$
y^2=x^3-x-3.
$$
If $x$ is odd, then $x^3=x$ modulo $8$, hence $x^3-x-3=5$ modulo $8$, and cannot be a perfect square. Hence, $x$ is even. Next, we can rewrite the equation as
$$
(x+2) (x^2-2x+3)=y^2+9 
$$
Because $x$ is even, $x^2-2x+3>0$ is $3$ modulo $4$, hence, by Proposition \ref{prop:755}, integers $x+2$, $x^2-2x+3$, $y$ and $3$ must have a common prime divisor $p$. Because $p$ is a divisor of $3$, we must have $p=3$, but it is easy to check that both $x+2$ and $x^2-2x+3$ cannot be divisible by $3$, a contradiction.   
 
\vspace{10pt}

The next equation we will consider is
$$
y^2=x^3+x-3.
$$
Because $y^2\geq 0$, we must have $x>1$. Also, the only solutions to this equation modulo $4$ are $(x,y)=(0,1)$ and $(0,3)$. Next, we can rewrite the equation as
$$
(x-1) (x^2+x+2)=y^2+1.
$$
Because $x>1$ and divisible by $4$, $x-1$ is a positive integer equal to $3$ modulo $4$. However, by Proposition \ref{prop:z2p1div}, all positive divisors of $y^2 + 1$ must be $1$ or $2$ modulo $4$, a contradiction.  

\vspace{10pt}

The next equation we will consider is
$$
y^2=x^3-5.
$$
The only solutions to this equation modulo $4$ are $(x,y)=(1,0)$ and $(1,2)$. 
We can rewrite this equation as
$$
(x-1) (x^2+x+1)=y^2+4. 
$$
Because $x$ is equal to $1$ modulo $4$, $x^2+x+1$ is a positive integer equal to $3$ modulo $4$. However, by Proposition \ref{prop:divsofx2pc}, all positive odd divisors of $y^2+4$ are of the form $u^2 + 4r^2$ for coprime integers $u,r$ and therefore cannot be equal to $3$ modulo $4$.

\vspace{10pt}

The next equation we will consider is
\begin{equation}\label{eq:y2mx3px2p1}
y^2=x^3-x^2-1.
\end{equation}
If $x$ is even, then $x^3-x^2-1$ is $3$ modulo $4$ and is not a perfect square. Also, if $x$ is $0$ or $1$ modulo $3$, then $x^3-x^2-1$ is $2$ modulo $3$ and is not a perfect square. Hence, $x$ is odd and it is equal to $2$ modulo $3$. Then we can 
rewrite the equation as
$$
(x+1) (x^2-2x+2)=y^2+3. 
$$
Then $x^2-2x+2$ is an odd positive integer equal to $2$ modulo $3$. However, by Proposition \ref{prop:divsofx2pc}, all positive odd divisors of $y^2+3$ are of the form $u^2+3r^2$ for coprime integers $u,r$ and therefore cannot be equal to $2$ modulo $3$. 

\vspace{10pt}

The next equation we will consider is
$$
y^2=x^3-2x-2.
$$
From the equation, it is clear that $x>1$. Also, the only solutions to this equation modulo $4$ are $(x,y)=(1,1)$ or $(1,3)$. We can then rewrite the equation as
$$
(x+1) (x^2-x-1)=y^2+1. 
$$
Because $x$ is equal to $1$ modulo $4$, $x^2-x-1$ is a positive integer equal to $3$ modulo $4$. However, by Proposition \ref{prop:z2p1div}, all positive divisors of $y^2 + 1$ must be $1$ or $2$ modulo $4$, a contradiction.   

\vspace{10pt}

The next equation we will consider is
$$
y^2=x^3-6.
$$
Solving this equation modulo $8$, we see that $x$ must be $7$ modulo $8$. 
We can rewrite the equation as
$$
(x-2) (x^2+2x+4)=y^2-2.
$$
Then $x-2$ is $5$ modulo $8$. However, by Proposition \ref{prop:z2m2div}, all divisors of $y^2-2$ must be $\pm 1$ or $\pm 2$ modulo $8$, a contradiction.

\vspace{10pt}

The next equation we will consider is
$$
y^2=x^3+6.
$$
From the equation, it is clear that $x>-2$. 
Solving this equation modulo $8$, we see that $x$ must be $3$ modulo $8$. We can rewrite the equation as
$$
(x+2) (x^2-2x+4)=y^2+2.
$$
Because $x$ is equal to $3$ modulo $8$, $x+2$ is an positive integer equal to $5$ modulo $8$. However, by Proposition \ref{prop:z2p2div}, all positive odd divisors of $y^2+2$ must be $1$ or $3$ modulo $8$, a contradiction. 

\vspace{10pt}

The next equation we will consider is
$$
y^2=x^3+7.
$$
It is clear from this equation that $x>-2$. 
Solving this equation modulo $4$, we see that solutions are $x$ is $1$ modulo $4$. We can rewrite the equation as
$$
(x+2) (x^2-2x+4)=y^2+1. 
$$
Because $x$ is equal to $1$ modulo $4$, $x+2$ is a positive integer equal to $3$ modulo $4$. However, by Proposition \ref{prop:z2p1div}, all positive divisors of $y^2+1$ must be $1$ or $2$ modulo $4$, a contradiction.

\vspace{10pt}

The next equation we will consider is
$$
y^2=x^3+x+5.
$$
We can rewrite this equation as
$$
(x-1) (x^2+x+2)=y^2-7. 
$$
Modulo $8$ analysis shows that $x$ is equal to $4$ modulo $8$. Substituting $x=8k+4$ for integer $k$, we have
$$
2(8k+3) (32k^2+36k+11)=y^2-7. 
$$
Modulo $7$ analysis shows that $k=1,4,5$ or $6$ modulo $7$. Then $8k+3=11,7,15$ or $23$ modulo $28$. On the other hand, Proposition \ref{prop:quadform} with $(a,b,c,x)=(-7,0,1,1)$ implies that if $p$ is a prime divisor of $y^2-7$, then either $p=2$, or $p=7$, or $\left(\frac{28}{p}\right)=1$. The last condition is equivalent to $\left(\frac{7}{p}\right)=1$, and, by \eqref{leg:p7}, this happens if and only if $p \equiv \pm 1, \pm 3, \pm 9$ modulo $28$. Because the product of primes equal to $p \equiv \pm 1, \pm 3, \pm 9$ modulo $28$ is again $p \equiv \pm 1, \pm 3, \pm 9$ modulo $28$, it follows that $y^2-7$ cannot have positive divisors equal to $11,15$ or $23$ modulo $28$. Because $8k+3$ is a divisor of $y^2-7$, the only possibility is that $8k+3=7$ modulo $28$, or $k=4$ modulo $7$. But in this case $32k^2+36k+11$ is a divisor of $y^2-7$ which is equal to $23$ modulo $28$, which is impossible as explained above. 

\vspace{10pt}

The next equation we will consider is
$$
y^2=x^3+3x-1.
$$
We can rewrite this equation as
$$
(x+1) (x^2-x+4)=y^2+5.
$$ 
Modulo $8$ analysis shows that $x$ is equal to $6$ modulo $8$. Substituting $x=8k+6$ for integer $k$, we have
$$
2(8k+7) (32k^2+44k+17)=y^2+5. 
$$
Modulo $5$ analysis shows that $k=1,3$ or $4$ modulo $5$. Then $8k+7=15,11$ or $19$ modulo $20$. On the other hand, Proposition \ref{prop:quadform} with $(a,b,c,x)=(5,0,1,1)$ implies that if $p$ is a prime divisor of $y^2+5$, then either $p=2$, or $p=5$, or $\left(\frac{-20}{p}\right)=1$. The last condition is equivalent to $\left(\frac{-5}{p}\right)=1$, and, by \eqref{leg:m5}, happens if and only if $p \equiv 1, 3, 7, 9$ modulo $20$. Because the product of primes equal to $p \equiv 1, 3, 7, 9$ modulo $20$ is again $p \equiv 1, 3, 7, 9$ modulo $20$, it follows that $y^2+5$ cannot have positive divisors equal to $11$ or $19$ modulo $20$. Because $8k+7$ is a divisor of $y^2+5$, the only possibility is that $8k+7=15$ modulo $20$, or $k=1$ modulo $5$. But in this case $32k^2+44k+17$ is a divisor of $y^2+5$ which is equal to $13$ modulo $20$, which is impossible as explained above.
 
\vspace{10pt}

The next equation we will consider is
\begin{equation}\label{eq:y2mxymx3p3}
y^2-xy=x^3-3.
\end{equation}
This equation only has real solutions for $x>1$. 
If $x$ is divisible by $3$, then so is $y$, but then $y^2-xy-x^3$ is divisible by $9$ and cannot be equal to $-3$. Hence, $x$ is not divisible by $3$. Then modulo $3$ analysis shows that $(x,y)=(2,1)$ modulo $3$. 
We can rewrite this equation as
$$
 x (x^2+y)=y^2+3.
$$
Because $x>1$ and $y^2+3>0$, $x$ and $x^2+y$ are positive integers, 
equal to $2$ modulo $3$. It is clear from \eqref{eq:y2mxymx3p3} that $x$ and $y$ are not both even, hence either $x$ or $x^2+y$ is odd. However, as explained in our solution to equation \eqref{eq:y2mx3px2p1} above, Proposition \ref{prop:divsofx2pc} implies that $y^2+3$ cannot have positive odd divisors equal to $2$ modulo $3$.  

\vspace{10pt}

The next equation we will consider is
$$
2y^2=x^3-3.
$$
From this equation, it is clear that $x>1$. 
Solving this equation modulo $8$, we see that we only have solutions for $x=3$ or $x=5$ modulo $8$. We can rewrite this equation as
$$
(x+1) (x^2-x+1)=2 y^2+4.
$$
Because $x=3$ or $5$ modulo $8$, $x^2-x+1$ is an odd positive integer equal to $7$ or $5$ modulo $8$. However, Proposition \ref{prop:z2p2div} implies that all odd positive divisors of $2(y^2+2)$ are equal to $1$ or $3$ modulo $8$. 

\vspace{10pt}

The next equation we will consider is
$$
y^2=x^3+x^2+3.
$$
From this equation it is clear that $x>-2$.
Solving this equation modulo $4$, we see that $x$ must be equal to $1$ modulo $4$. We can rewrite the equation as
$$
(x+2) (x^2-x+2)=y^2+1. 
$$
Because $x=1$ modulo $4$, $x+2$ is a positive integer equal to $3$ modulo $4$. However, by Proposition \ref{prop:z2p1div}, all positive divisors of $y^2+1$ must be $1$ or $2$ modulo $4$, a contradiction.

\vspace{10pt}

The next equation we will consider is
$$
y^2=x^3+x^2+x-1.
$$
From this equation, it is clear that $x>0$.
Solving this equation modulo $8$, we see that $x$ is equal to $6$ modulo $8$. We can rewrite this equation as
$$
(x+1) (x^2+1)=y^2+2.
$$
Because $x=6$ modulo $8$, $x+1$ is a positive integer equal to $7$ modulo $8$. However, by Proposition \ref{prop:z2p2div}, all odd positive divisors of $y^2+2$ must be $1$ or $3$ modulo $8$, a contradiction.

\subsection{Exercise 7.6}\label{ex:H30gen1norat}
\textbf{\emph{Check that the equations listed in Table \ref{tab:H30gen1norat} have no rational solutions and therefore no integer solutions.}}

	\begin{center}
	\begin{tabular}{ |c|c|c|c|c|c| } 
		\hline
		$H$ & Equation & $H$ & Equation & $H$ & Equation \\ 
		\hline\hline
		$28$ & $x^3+x^2 y-y^3-x+2=0$ & $29$ & $y^3=2x^3-x+3$ & $30$ & $y^3=2x^3-x+4$ \\ 
		\hline
		$29$ & $x^3+x^2y-y^3+x+3=0$ & $30$ & $x^3+x^2y-y^3-y+4=0$  & $30$ & $x^4-2y^2+y-4=0$ \\ 
		\hline
		$29$ & $x^3+x^2y+y^3-x+3=0$ & $30$ & $x^3+x^2y-y^3+y^2+2=0$ & $30$ & $x^4-2y^2+x-4=0$ \\ 
		\hline
		$29$ & $x^3+x^2y+y^3-x+y+1=0$ & $30$ & $x^3+x^2y+y^3-x+4=0$ & $30$ & $x^4-2y^2+x+y+2=0$ \\ 
		\hline
		$29$ & $x^3+x^2y+y^3+y+3=0$ & $30$ & $x^3+x^2y+y^3+xy-2=0$ &  &  \\ 
		\hline	
		$29$ & $x^3+x^2y+y^3+x+3=0$ & $30$ & $2x^3+y^3-x-y+2=0$ &  &  \\ 
		\hline			
	\end{tabular}
	\captionof{table}{\label{tab:H30gen1norat} Genus $1$ equations of size $H\leq 30$ with no rational solutions.}
\end{center}

In this section, we will show that some genus $1$ equations of size $H \leq 30$ have no rational solutions and therefore no integer solutions. 
We start with general two-variable cubic equations of the form
\begin{equation}\label{cubic,a1-10}
a_1 x^3 + a_2 y^3 + a_3 + a_4 x^2 y + a_5 x^2  + a_6 y^2 x + a_7 y^2  + a_8 x + a_9 y + a_{10} x y  = 0.
\end{equation}
For such equations we will use the Magma code from Section \ref{ex:H48Hassecount}, which we present below for convenience.
\begin{align}\label{cubic,a1-10_magma}
	&{\tt > C := GenusOneModel(3,\, ``a_1 \, a_2 \, a_3 \, a_4 \, a_5 \, a_6 \, a_7 \, a_8 \, a_9 \, a_{10}");}\nonumber\\
	&{\tt > IsLocallySoluble(C);}\nonumber\\
	&{\tt > E := Jacobian(C);}\nonumber\\
	&{\tt > MW,iso := MordellWeilGroup(E);}\nonumber\\
	&{\tt > G,quomap := quo<MW|3*MW>;}\\
	&{\tt > PP := [iso(g @@ quomap): g\, in \, G];}\nonumber\\
	&{\tt > cubics := [P \, eq \, E!0 \, select \, GenusOneModel(3,E) }\nonumber\\
	&\quad\quad\quad\quad\quad\quad{\tt else \, GenusOneModel(3,P): \, P \, in \, PP];}\nonumber\\
	&{\tt > exists\{ deltaP : deltaP \, in \, cubics | IsEquivalent(C,deltaP) \};}\nonumber
\end{align}
If the output is {\tt true, false} then the equation is locally solvable but it has no rational solutions.

For example, equation
$$
x^3+x^2 y-y^3-x+2=0,
$$
is of the form \eqref{cubic,a1-10} with $(a_1,\dots,a_{10}) =(1,-1,2,1,0,0,0,-1,0,0)$. The Magma code \eqref{cubic,a1-10_magma} with the first line
$$
{\tt > C := GenusOneModel(3, ``1 \, -1 \, 2 \, 1 \, 0 \, 0 \, 0 \, -1 \, 0 \, 0 "); }
$$
outputs {\tt true, false} so the equation is locally solvable and it has no non-trivial integer solutions.

All other cubic equations from Table \ref{tab:H30gen1norat} are similar: they are of the form \eqref{cubic,a1-10} with $(a_1,\dots,a_{10})$ given in Table \ref{tab:H30gen1noratcube}. For all equations, the output of \eqref{cubic,a1-10_magma} is {\tt true, false} so the equations are solvable everywhere locally but they have no rational solutions, and therefore no integer solutions.

	\begin{center}
	\begin{tabular}{ |c|c|c|c|c|c|c|c| } 
		\hline
		 Equation & $(a_1, \dots, a_{10})$  \\ 
		\hline\hline
		 $x^3+x^2 y-y^3-x+2=0$ & $(1,-1,2,1,0,0,0,-1,0,0)$  \\ \hline
		 $x^3+x^2y-y^3+x+3=0$ & $(1,-1,3,1,0,0,0,1,0,0)$   \\ \hline
		 $x^3+x^2y+y^3-x+3=0$ & $(1,1,3,1,0,0,0,-1,0,0)$   \\ \hline
		 $x^3+x^2y+y^3-x+y+1=0$ &$(1,1,1,1,0,0,0,-1,1,0)$  \\ \hline
		 $x^3+x^2y+y^3+y+3=0$ & $(1,1,3,1,0,0,0,0,1,0)$   \\ \hline
		$x^3+x^2y+y^3+x+3=0$ & $(1,1,3,1,0,0,0,1,0,0)$  \\ \hline
		$y^3=2x^3-x+3$ & $(2,-1,3,0,0,0,0,-1,0,0)$  \\ \hline
		 $x^3+x^2y-y^3-y+4=0$  & $(1,-1,4,1,0,0,0,0,-1,0)$  \\ \hline
		 $x^3+x^2y-y^3+y^2+2=0$ & $(1,-1,2,1,0,0,1,0,0,0)$  \\ \hline
		$x^3+x^2y+y^3-x+4=0$ & $(1,1,4,1,0,0,0,-1,0,0)$ \\ \hline
		$x^3+x^2y+y^3+xy-2=0$ & $(1,1,-2,1,0,0,0,0,0,1)$   \\ \hline
		$2x^3+y^3-x-y+2=0$ & $(2,1,2,0,0,0,0,-1,-1,0)$   \\ \hline
		$y^3=2x^3-x+4$ & $(2,-1,4,0,0,0,0,-1,0,0)$   \\ \hline
	\end{tabular}
	\captionof{table}{\label{tab:H30gen1noratcube} The coefficients $(a_1,\dots,a_{10})$ in \eqref{cubic,a1-10} for the cubic equations listed in Table \ref{tab:H30gen1norat}. }
\end{center}

We next consider equations of the form
\begin{equation}\label{eq:genus1quartic}
	y^2 + y(fx^2 + gx + h) - (ax^4 + bx^3 + cx^2 + dx + e) = 0,
\end{equation}
where $a,b,c,d,e,f,g,h$ are integer coefficients. 

To solve this type of equation, we will use the Magma code
\begin{align*}\label{eq:genusonemodel}
	&{\tt > C := GenusOneModel(2,``f \, g \, h\, a \, b \, c \, d \, e ");}\\
	&{\tt > IsLocallySoluble(C);}\\
	&{\tt > E := Jacobian(C);}\\
	&{\tt > MW,iso := MordellWeilGroup(E);}\\
	&{\tt > G,quomap := quo<MW|3*MW>;}\\
	&{\tt > PP := [iso(g @@ quomap): g\, in \, G];}\\
	&{\tt > quartics := [P \, eq \, E!0 \, select \, GenusOneModel(2,E) \quad else \, GenusOneModel(2,P): \, P \, in \, PP];}\\
	&{\tt > exists\{ deltaP : deltaP \, in \, quartics | IsEquivalent(C,deltaP) \};}
\end{align*}
In all examples (see Table \ref{tab:H30gen1noratsols}), the code returns {\tt true, false}. This means that the equation is solvable everywhere locally but it has no rational solutions and therefore no integer solutions.

	\begin{center}
		\begin{tabular}{ |c|c|c|c|c|c| } 
			\hline
			Equation & $(f, g , h, a , b , c ,d , e)$ \\ 
			\hline\hline

			  $x^4-2y^2+y-4=0$ & $(0 , 0 , -1 , 2 , 0 , 0 , 0 , -8 )$ \\ \hline
			  $x^4-2y^2+x-4=0$ & $(0 , 0 , 0 , 2 , 0 , 0 , 2 , -8 )$ \\ \hline
			 $x^4-2y^2+x+y+2=0$&  $(0 , 0 , -1 , 2 , 0 , 0 , 2 , 4 )$ \\ \hline		
		\end{tabular}
		\captionof{table}{\label{tab:H30gen1noratsols} Information to prove that the equations listed in Table \ref{tab:H30gen1norat} of the form \eqref{eq:genus1quartic} have no integer solutions.}
	\end{center}

\subsection{Exercise 7.10}\label{ex:vietasolv}
\textbf{\emph{Check that the equations listed in Table $\ref{tab:vietasolv}$ have no integer solutions with $\min(|x|,|y|,|z|)\leq t^*$ and conclude that they have no integer solutions.}}

		\begin{center}
	\begin{tabular}{ |c|c|c|c|c|c| }
		\hline
		$H$ & Equation & $t^*$ & $H$ & Equation & $t^*$\\ 
		\hline\hline
		$22$ & $2+x^2+y^2+x y z-z^2=0$ & $0$ & $29$ & $3+x^2+y+x^2 y+y^2+y z^2=0$ & $1.40$ \\ 
		\hline
		$24$ & $4+x^2+y^2+x y z+z^2=0$ & $3.36$ & $29$ & $1+x^2+2 y+x^2 y-y^2+y z^2=0$ & $1$\\ 
		\hline
		$25$ & $1+y^2+x^2 y z+z^2=0$ & $1.55$ & $30$ & $6+y^2+x^2 y z+z^2=0$ & $1.91$ \\ 
		\hline
		$27$ & $3+x^2+x^2 y-y^2+y z^2=0$ & $1.14$ & $30$ & $6+x^2+x^2 y+y^2+y z^2=0$ & $1.86$\\ 
		\hline
		$27$ & $1+x^2+y+x^2 y-y^2+y z^2=0$ & $1$ & $30$ & $6-x^2+x^2 y+y^2+y z^2=0$ & $1.44$\\ 
		\hline
		$28$ & $4+y^2+x^2 y z+z^2=0$ & $1.80$ & $30$ & $2+x^2+y^2+x^2 y z+z^2=0$ & $1.89$\\ 
		\hline
		$28$ & $8+x^2+y^2+x y z+z^2=0$ & $3.61$ & $30$ & $2+x^2+y^2+x^2 y z-z^2=0$ & $1.41$\\ 
		\hline
		$28$ & $2+3 x+x^2 y+y^2+y z^2=0$ & $1.71$ & $30$ & $4+x^2+y+x^2 y-y^2+y z^2=0$ & $1.13$\\ 
		\hline
		$29$ & $5+y^2+x^2 y z+z^2=0$ & $1.86$ & $30$ & $2+x^2+2 y+x^2 y-y^2+y z^2=0$ & $1$\\ 
		\hline
	\end{tabular}
	\captionof{table}{\label{tab:vietasolv} Equations of size $H\leq 30$ solvable by Vieta jumping.}
\end{center} 

The Mathematica command 
$$
\begin{aligned}
{\tt N[MaxValue[Min[Abs[x],Abs[y],Abs[z]],\{P(x,y,z)==0,Abs[x']\geq Abs[x],Abs[y']\geq Abs[y], } \phantom{\}} \\ \phantom{\{} {\tt Abs[z']\geq Abs[z]\},\{x,y,z\}]]}
\end{aligned}
$$
outputs $t^*$. The theory in Section 7.1.4 of the book implies that if the equation has no integer solutions with $|x| \leq t^*$, $|y| \leq t^*$ or $|z| \leq t^*$, then we can conclude that it has no integer solutions.

For each equation in Table \ref{tab:vietasolv}, the final column of Table \ref{tab:vietasolvsols} lists the whether the output to
\begin{equation}\label{reduce:vietaxyz}
{\tt Reduce[\{ P(x,y,z) == 0,   Abs[x] \leq t^* || Abs[y] \leq t^* || Abs[z] \leq t^* \}, \{x, y, z\}, Integers]}
\end{equation}
outputs {\tt False}.

		\begin{center}
			\begin{tabular}{ |c|c|c|c|c|c|c| }
				\hline
				 Equation & $t^*$ & $x'$ & $y'$ & $z'$ & \eqref{reduce:vietaxyz} output {\tt False}? \\ 
				\hline\hline
				 $2+x^2+y^2+x y z-z^2=0$ & $0$ & $-yz-x$ & $-xz-y$ & $xy-z$ & Yes \\\hline
				$4+x^2+y^2+x y z+z^2=0$ & $3.36$ & $-yz-x$ & $-xz-y$ & $-xy-z$ & Yes \\\hline
				$1+y^2+x^2 y z+z^2=0$ & $1.55$ & & $-x^2 z-y$ & $-x^2y-z$  & Yes \\\hline
				 $3+x^2+x^2 y-y^2+y z^2=0$ & $1.14$ & & $x^2+z^2-y$ & & No \\\hline
				 $1+x^2+y+x^2 y-y^2+y z^2=0$ & $1$ & & $x^2+z^2+1-y$ &&  No \\\hline
				 $4+y^2+x^2 y z+z^2=0$ & $1.80$ & & $-x^2 z-y$ & $-x^2y-z$  & Yes \\\hline
				$8+x^2+y^2+x y z+z^2=0$ & $3.61$ & $-yz-x$ & $-xz-y$ & $-xy-z$ & Yes \\\hline
				$2+3 x+x^2 y+y^2+y z^2=0$ & $1.71$ & & $-x^2-z^2-y$ & & Yes \\\hline
				 $5+y^2+x^2 y z+z^2=0$ & $1.86$ & & $-x^2 z-y$ & $-x^2y-z$ & Yes \\\hline
				 $3+x^2+y+x^2 y+y^2+y z^2=0$ & $1.40$ && $-x^2-1-y$ & & Yes \\ \hline
				  $1+x^2+2 y+x^2 y-y^2+y z^2=0$ & $1$ & & $x^2+z^2+2-y$ & & No \\ \hline
				$6+y^2+x^2 y z+z^2=0$ & $1.91$ & & $-x^2 z-y$ & $-x^2y-z$ & Yes  \\ \hline
				 $6+x^2+x^2 y+y^2+y z^2=0$ & $1.86$ & & $-x^2-z^2-y$ & & No \\ \hline
				$6-x^2+x^2 y+y^2+y z^2=0$ & $1.44$ & &  $-x^2-z^2-y$ & & No \\ \hline
				  $2+x^2+y^2+x^2 y z+z^2=0$ & $1.89$ &&  $-x^2 z-y$ & $-x^2y-z$ & No  \\\hline
				 $2+x^2+y^2+x^2 y z-z^2=0$ & $1.41$ &&  $-x^2 z-y$ & $x^2y-z$ & No \\ \hline
				 $4+x^2+y+x^2 y-y^2+y z^2=0$ & $1.13$ && $x^2+z^2+1-y$ & & No \\ \hline
				 $2+x^2+2 y+x^2 y-y^2+y z^2=0$ & $1$ && $x^2+z^2+2-y$ & & No \\ \hline
			\end{tabular}
			\captionof{table}{\label{tab:vietasolvsols} Equations of size $H\leq 30$ solvable by Vieta jumping.}
		\end{center} 

For the equations for which command \eqref{reduce:vietaxyz} outputs {\tt False}, we are done. 
For the other equations, we now provide more details. For these equations, we run the three commands 
$$
{\tt Reduce[\{ P(x,y,z) == 0,   Abs[x] \leq t^* \}, \{x, y, z\}, Integers]}
$$
$$
{\tt Reduce[\{ P(x,y,z) == 0, Abs[y] \leq t^* \}, \{x, y, z\}, Integers]}
$$
$$
{ \tt Reduce[\{ P(x,y,z) == 0,  Abs[z] \leq t^* \}, \{x, y, z\}, Integers]}
$$
and only consider the cases which do not output {\tt False}. 

Let us first consider the equation
\begin{equation}\label{eq:3px2px2ymy2py}
3+x^2+x^2y-y^2+yz^2=0.
\end{equation}
For this equation, $t^*=1.14$, so it suffices to prove that it has no integer solutions such that $\min(|x|,|y|,|z|)\leq 1$. The {\tt Reduce} command outputs false for $|x|\leq 1$ and $|y| \leq 1$, so we only need to consider the case $|z|\leq 1$. 
 
For $z=\pm 1$, the equation is reduced to $3+x^2+x^2y-y^2+y=0$. Let us consider this equation as quadratic in $y$ with parameter $x$. To have integer solutions, the discriminant $D$ must be a perfect square, so $D=x^4+6x^2+13=d^2$ for some integer $d$. However, we have that
$$
(x^2+3)^2< x^4+6x^2+13 < (x^2+4)^2
$$
is true for all integers $x$, therefore the discriminant cannot be a perfect square. 
For $z=0$, the equation is reduced to $3+x^2+x^2y-y^2=0$. To have integer solutions, the discriminant must be a perfect square, so $D=x^4+4x^2+12=d^2$ for some integer $d$. However, we have that
$$
(x^2+2)^2<x^4+4x^2+12 < (x^2+3)^2
$$
is true for all $|x|>1$. As \eqref{eq:3px2px2ymy2py} has no solutions with $|x|\leq 1$, we can conclude that \eqref{eq:3px2px2ymy2py} has no integer solutions. 

\vspace{10pt}

We will now consider the equation
\begin{equation}\label{eq:1px2pypx2ymy2pyz2}
1+x^2+y+x^2y-y^2+yz^2=0.
\end{equation}
For this equation, $t^*=1$, so it suffices to prove that it has no integer solutions such that $\min(|x|,|y|,|z|)\leq 1$. 
The reduce command outputs {\tt False} for $|x|\leq 1$, $|y| \leq 1$ and $|z|=1$.  For $z=0$, the equation is reduced to $1+x^2+y+x^2y-y^2=0$. Let us consider this equation as quadratic in $y$ with parameter $x$. To have integer solutions, the discriminant $D$ must be a perfect square, so $D=x^4+6x^2+5=d^2$ for some integer $d$. However, we have that
$$
(x^2+2)^2 < x^4+6x^2+5 < (x^2+3)^2
$$
is true for all $x$, therefore the discriminant cannot be a perfect square. Therefore, equation  \eqref{eq:1px2pypx2ymy2pyz2} has no integer solutions. 

\vspace{10pt}

We will now consider the equation
\begin{equation}\label{eq:1px2p2ypx2ymy2pyz2}
1+x^2+2y+x^2y-y^2+yz^2=0.
\end{equation}
For this equation, $t^*=1$, so it suffices to prove that it has no integer solutions such that $\min(|x|,|y|,|z|)\leq 1$.
The {\tt Reduce} command outputs false for $|x|\leq1$, $|y| \leq 1$ and $z=0$. For $|z|=1$, the equation is reduced to $1+x^2+3y+x^2y-y^2=0$. Let us consider this equation as quadratic in $y$ with parameter $x$. To have integer solutions, the discriminant $D$ must be a perfect square, so $D=x^4+10x^2+13=d^2$ for some integer $d$. We have that
$$
(x^2+4)^2 < x^4+10x^2+13<(x^2+5)^2
$$
is true for all $|x| > 1$. As \eqref{eq:1px2p2ypx2ymy2pyz2} has no solutions with $|x|\leq 1$, we can conclude that \eqref{eq:1px2p2ypx2ymy2pyz2} has no integer solutions.

\vspace{10pt}

We will now consider the equation
\begin{equation}\label{eq:6px2px2ypy2pyz2}
6+x^2+x^2y+y^2+yz^2=0.
\end{equation}
For this equation, $t^*=1.86$, so it suffices to prove that it has no integer solutions such that $\min(|x|,|y|,|z|)\leq 1$.
The {\tt Reduce} command outputs false for $|x|\leq1$, $|y| \leq 1$ and $z=0$. For $|z|=1$, the equation is reduced to $6+x^2+x^2y+y^2+y=0$. Let us consider this equation as quadratic in $y$ with parameter $x$. To have integer solutions, the discriminant $D$ must be a perfect square, so $D=x^4-4x^2-24=d^2$ for some integer $d$. We have that
$$
(x^2-3)^2 < x^4-4x^2-24<(x^2-2)^2
$$
is true for all $|x| > 4$. By checking values $|x|\leq 3$ we obtain no integer solutions, and therefore \eqref{eq:6px2px2ypy2pyz2} has no integer solutions.

\vspace{10pt}

We will now consider the equation
\begin{equation}\label{eq:6mx2px2ypy2pyz2}
6-x^2+x^2y+y^2+yz^2=0.
\end{equation}
For this equation, $t^*=1.44$, so it suffices to prove that it has no integer solutions such that $\min(|x|,|y|,|z|)\leq 1$.
The {\tt Reduce} command outputs false for $|x|\leq1$ and $|y| \leq 1$. For $|z|=1$, the equation is reduced to $6-x^2+x^2y+y^2+y=0$. Let us consider this equation as quadratic in $y$ with parameter $x$. To have integer solutions, the discriminant $D$ must be a perfect square, so $D=x^4+6x^2-23=d^2$ for some integer $d$. We have that
$$
(x^2+2)^2 < x^4+6x^2-23<(x^2+3)^2
$$
is true for all $|x| > 4$. By checking values $|x|\leq 3$ we obtain no integer solutions. For $z=0$, the equation is reduced to $6-x^2+x^2y+y^2=0$. Let us consider this equation as quadratic in $y$ with parameter $x$. To have integer solutions, the discriminant $D$ must be a perfect square, so $D=x^4+4x^2-24=d^2$ for some integer $d$. We have that
$$
(x^2+1)^2 < x^4+6x^2-23<(x^2+2)^2
$$
is true for all $|x| > 5$. By checking values $|x|\leq 4$ we obtain no integer solutions, and therefore \eqref{eq:6mx2px2ypy2pyz2} has no integer solutions.

\vspace{10pt}

We will now consider the equation
\begin{equation}\label{eq:2px2py2px2yzpz2}
2+x^2+y^2+x^2yz+z^2=0.
\end{equation}
For this equation, $t^*=1.89$, so it suffices to prove that it has no integer solutions such that $\min(|x|,|y|,|z|)\leq 1$.
This equation is symmetric in $y,z$. The {\tt Reduce} command outputs false for $|x|\leq1$ and $y=0$. For $|y|=1$, the equation is reduced to $3+x^2 \pm x^2z+z^2=0$. Let us consider this equation as quadratic in $z$ with parameter $x$. To have integer solutions, the discriminant $D$ must be a perfect square, so $D=x^4-4x^2-12=d^2$ for some integer $d$. We have that
$$
(x^2-3)^2 < x^4-4x^2-12<(x^2-2)^2
$$
is true for all $|x| > 5$. By checking values $|x|\leq 4$ we obtain no integer solutions, and therefore \eqref{eq:2px2py2px2yzpz2} has no integer solutions.

\vspace{10pt}

We will now consider the equation
\begin{equation}\label{eq:2px2py2px2yzmz2}
2+x^2+y^2+x^2yz-z^2=0.
\end{equation}
For this equation, $t^*=1.41$, so it suffices to prove that it has no integer solutions such that $\min(|x|,|y|,|z|)\leq 1$.
The {\tt Reduce} command outputs false for $|x|\leq1$, $y=0$ and $|z| \leq 1$. For $|y|=1$, the equation is reduced to $3+x^2 \pm x^2z-z^2=0$. Let us consider this equation as quadratic in $z$ with parameter $x$. To have integer solutions, the discriminant $D$ must be a perfect square, so $D=x^4+4x^2+12=d^2$ for some integer $d$. We have that
$$
(x^2+2)^2 < x^4+4x^2+12 < (x^2+3)^2
$$
is true for all $|x| > 1$. By checking values $|x|\leq 1$ we obtain no integer solutions, and therefore \eqref{eq:2px2py2px2yzmz2} has no integer solutions.

\vspace{10pt}

We will now consider the equation
\begin{equation}\label{eq:4px2pypx2ymy2pyz2}
4+x^2+y+x^2y-y^2+yz^2=0.
\end{equation}
For this equation, $t^*=1.13$, so it suffices to prove that it has no integer solutions such that $\min(|x|,|y|,|z|)\leq 1$.
The {\tt Reduce} command outputs false for $|x|\leq1$, $|y|\leq1$ and $|z| = 1$. For $z=0$, the equation is reduced to $4+x^2+y+x^2y-y^2=0$. Let us consider this equation as quadratic in $y$ with parameter $x$. To have integer solutions, the discriminant $D$ must be a perfect square, so $D=x^4+6x^2+17=d^2$ for some integer $d$. We have that
$$
(x^2+3)^2 < x^4+6x^2+17 < (x^2+4)^2
$$
is true for all $x \neq 0$. By checking the value $x=0$ we obtain no integer solutions, and therefore \eqref{eq:4px2pypx2ymy2pyz2} has no integer solutions.

\vspace{10pt}

We will now consider the equation
\begin{equation}\label{eq:2px2p2ypx2ymy2pyz2}
2+x^2+2y+x^2y-y^2+yz^2=0.
\end{equation}
For this equation, $t^*=1$, so it suffices to prove that it has no integer solutions such that $\min(|x|,|y|,|z|)\leq 1$.
The {\tt Reduce} command outputs false for $|x|\leq1$, $|y|\leq1$ and $|z| = 1$. For $z=0$, the equation is reduced to $2+x^2+2y+x^2y-y^2=0$. Let us consider this equation as quadratic in $y$ with parameter $x$. To have integer solutions, the discriminant $D$ must be a perfect square, so $D=x^4+8x^2+12=d^2$ for some integer $d$. We have that
$$
(x^2+3)^2 < x^4+8x^2+12 < (x^2+4)^2
$$
is true for all $x$. Therefore \eqref{eq:2px2p2ypx2ymy2pyz2} has no integer solutions.

\subsection{Exercise 7.13}\label{ex:H31f}
\textbf{\emph{Show that equation
		\begin{equation}\label{eq:h31f}
			3+x^2+x^2 y+2 y^2-y z^2 = 0
		\end{equation}
		has no integer solutions.}}

We can prove that this equation has no integer solutions using Vieta jumping. After multiplying \eqref{eq:h31f} by $4$, and making the substitutions $x = (u + v)/2$, $z= (u - v)/2$ and $y= w/2$, \eqref{eq:h31f} is reduced to 
\begin{equation}\label{eq:h31fred}
	12 + u^2 + 2 u v + v^2 + 2 u v w + 2 w^2 =0 
\end{equation}
where $u,v$ have the same parity and $w$ is even. If $(u,v,w)$ is any integer solution to \eqref{eq:h31fred} then so are
$$
(a) \,\, (-2 v - 2 v w - u,v,w), \quad (b) \,\, (u,-2 u - 2 u w - v,w), \quad \text{and} \quad \,\, (c) \,\, (u,v,-u v - w).
$$

Recall that $(u,v,w)$ is a minimal solution to \eqref{eq:h31fred} if none of the operations (a), (b) and (c) decrease $|u|+|v|+|w|$. The Mathematica command
$$
\begin{aligned}
	& {\tt N[MaxValue[
		Min[Abs[u], Abs[v], 
		Abs[w]], \{12 + u^2 + 2 u v + v^2 + 2 u v w + 2 w^2 == 0, } \\ &
	{\tt Abs[-2 v - 2 v w - u] \geq Abs[u], Abs[-2 u - 2 u w - v] \geq Abs[v], 
		Abs[-u v - w] \geq Abs[w]\}, \{u, v, w\}]]}
\end{aligned}
$$
returns $3.49203$, hence, any minimal solution to \eqref{eq:h31fred} must have $\min\{|u|,|v|,|w|\} \leq 3$. Direct analysis does not return any integer solutions to \eqref{eq:h31fred} satisfying this condition, so the equation has no minimal solutions. Therefore no integer solutions exist to equations \eqref{eq:h31fred} or \eqref{eq:h31f}.

\subsection{Exercise 7.15}\label{ex:H30}
\textbf{\emph{Prove that the equations listed in Table \ref{tab:H30} have no integer solutions. }}

		\begin{center}
	\begin{tabular}{ |c|c|c|c|c| }
		\hline
		$H$ & Equation & Representation & $m$ \\ 
		\hline\hline
		$26$ & $2+x^2 y+y^2+x^2 z+z^2=0$ & $(x^2+y-z)(y+z)=-2-2 z^2$ & $8$ \\ 
		\hline
		$27$ & $3+x^2+x^2 y+y^2-2 z^2=0$ & $(1+y) (-1+x^2+y)=-4+2 z^2$ & $8$ \\ 
		\hline
		$27$ & $3-x^2+x^2 y+y^2+2 z^2=0$ & $(-1+y) (1+x^2+y)=-4-2 z^2$ & $8$ \\ 
		\hline
		$29$ & $5-2 x^2+x^2 y+y^2+z^2=0$ & $(-2+y) (2+x^2+y)=-9-z^2$ & $24$ \\ 
		\hline
		$29$ & $5+x^3+y^2+x y z+z^2=0$ & $(2+x)(4-2 x+x^2+y z)=3-(y-z)^2$ & $36$ \\
		\hline
		$30$ & $6+x^2+x^2 y+2 y^2+z^2=0$ & $(1+y)(-2+x^2+2 y)=-8-z^2$ & $32$ \\  
		\hline
		$30$ & $2-x^2-x y+x^3 y+z^2=0$ & $(-1+x) (1+x) (-1+x y)=-1-z^2$ & $4$ \\ 
		\hline
		$30$ & $2+x^2+2 y+x^2 y-2 y^2+z^2=0$ & $(4+x^2-2 y) (1+y)=2-z^2$ & $8$ \\ 
		\hline
		$30$ & $2+2 x+x^3+y^2+x y^2-z^2=0$ & $(1+x) (3-x+x^2+y^2)=1+z^2$ & $4$ \\ 
		\hline
		$30$ & $2-2 x+x^3-y^2+x y^2+z^2=0$ & $(-1+x)(-1+x+x^2+y^2)=-1-z^2$ & $4$ \\ 
		\hline
	\end{tabular}
	\captionof{table}{\label{tab:H30} Equations of size $H\leq 30$ solvable by Algorithm 5.29 in the book after a linear transformation.}
\end{center} 

Equations $2+x^2 y+y^2+x^2 z+z^2=0$, $5+x^3+y^2+x y z+z^2=0$ and $2+2 x+x^3+y^2+x y^2-z^2=0$ have no integer solutions, which is proved in Section 7.2.1 of the book. To prove that other equations from Table \ref{tab:H30} have no integer solutions, we can use the third column of the table to rewrite the equations as
$$
Q(x,y)=a(z^2+b),
$$
where $Q(x,y)$ is a factorizable polynomial and $a>0$ and $b$ are integers. We can then apply Proposition \ref{prop:quadform} to find the prime divisors of $a(z^2+b)$, and obtain a contradiction with the divisors of the factors of $Q(x,y)$.
 
The first equation we will consider is
\begin{equation}\label{eq:3px2px2ypy2m2z2}
3+x^2+x^2y+y^2-2z^2=0,
\end{equation}
which we can rewrite as
$$
(1+y)(-1+x^2+y)=2(z^2-2).
$$
Modulo $8$ analysis shows that $x$ is always odd, therefore $x^2=1$ modulo $8$. The left-hand side modulo $8$ is then $y(y+1)$. Also, modulo $8$ analysis shows that $y$ can only be $2,3,4$ or $5$ modulo $8$, hence either $y$ or $y+1$ must be $3$ or $5$ modulo $8$. However, by Proposition \ref{prop:z2m2div}, all odd divisors of $2(z^2-2)$ must be $1$ or $7$ modulo $8$, which is a contradiction. In conclusion, equation \eqref{eq:3px2px2ypy2m2z2} has no integer solutions. 

\vspace{10pt}

The next equation we will consider is 
\begin{equation}\label{eq:3mx2px2ypy2p2z2}
3-x^2+x^2y+y^2+2z^2=0,
\end{equation}
which we can rewrite as 
$$
-(y-1)(1+x^2+y)=2(z^2+2).
$$
This representation suggests to make the change of variable $w=-(y-1)$, which reduces the equation to
\begin{equation}\label{red:3mx2px2ypy2p2z2}
	w(x^2-w+2)=2(z^2+2),
\end{equation}
and we can see that the right-hand side is positive, so the integers $w$ and $x^2-w+2$ must have the same sign. Because equation \eqref{red:3mx2px2ypy2p2z2} has no real solutions with $w \leq 0$, this implies that $w$ and $x^2-w+2$ are both positive.
Modulo $8$ analysis on \eqref{red:3mx2px2ypy2p2z2} shows $x$ is always odd, therefore $x^2=1$ modulo $8$. The left-hand side modulo $8$ is then $w(3-w)$. Also, modulo $8$ analysis shows that $w$ can only be $4,5,6$ or $7$ modulo $8$, hence either $w$ or $3-w$ must be $5$ or $7$ modulo $8$. However, by Proposition \ref{prop:z2p2div}, all positive odd divisors of $2(z^2+2)$ must be $1$ or $3$ modulo $8$, which is a contradiction. In conclusion, equation \eqref{eq:3mx2px2ypy2p2z2} has no integer solutions.

\vspace{10pt}

The next equation we will consider is 
\begin{equation}\label{eq:5m2x2px2ypy2pz2}
5-2x^2+x^2y+y^2+z^2=0,
\end{equation}
which we can rewrite as
$$
-(y-2)(2+x^2+y)=z^2+9.
$$
This representation suggests to make the change of variable $w=-(y-2)$, which reduces the equation to
\begin{equation}\label{red:5m2x2px2ypy2pz2}
	w(x^2-w+4)=z^2+9.
\end{equation}
This equation has no real solutions with $w \leq 0$, hence $w>0$ and $x^2-w+4>0$.  Modulo $8$ analysis on \eqref{red:5m2x2px2ypy2pz2} shows $x$ is always odd, so $x^2$ is always equal to $1$ modulo $8$. The left-hand side modulo $8$ is then $w(5-w)$. Also modulo $8$ analysis shows that $w$ can only be $6$ or $7$ modulo $8$, hence either $w$ or $5-w$ is $7$ modulo $8$, which is $3$ modulo $4$. By Proposition \ref{prop:abx2y2} integers $w$, $x^2-w+4$, $z$ and $3$ must have a common prime divisor $p$, and it is obvious that $p=3$. However, if $w$ and $x^2-w+4$ are both divisible by $3$, then $x^2$ must be equal to $2$ modulo $3$, which is a contradiction. In conclusion, equation \eqref{eq:5m2x2px2ypy2pz2} has no integer solutions.

\vspace{10pt}

The next equation we will consider is
\begin{equation}\label{eq:6px2px2yp2y2pz2}
6+x^2+x^2y+2y^2+z^2=0,
\end{equation}
which we can rewrite as
$$
-(1+y)(-2+x^2+2y)=z^2+8.
$$
This representation suggests to make the change of variable $w=-(1+y)$, which reduces the equation to
$$
w(x^2-2w-4)=z^2+8.
$$
This equation has no real solutions with  $w \leq 0$, hence $w>0$ and $x^2-2w-4>0$.
Proposition \ref{prop:quadform} with $(a,b,c,y)=(1,0,8,1)$ implies that if $p$ is a prime divisor of $z^2+8$, then either $p=2$ or $\left(\frac{-32}{p}\right)=1$. The last condition is equivalent to $\left(\frac{-2}{p}\right)=1$ which by \eqref{eq:reposm2} happens if and only if $p \equiv 1,3$ modulo $8$. Note that the product of primes equal to $p \equiv 1,3$ modulo $8$ is again $\equiv 1,3$ modulo $8$. Hence, all odd positive divisors of $z^2+8$ must be $1$ or $3$ modulo $8$. Modulo $8$ analysis shows that if $w$ is odd, then $w$ must be equal to $7$ modulo $8$, which is a contradiction. Hence, $w$ must be even, which implies that $z$ is even as well. Let us write $w=2u$ and $z=2t$ for some integers $u,t$. Then the equation simplifies to
$$
u(x^2-4 u-4)=2(t^2+2).
$$
If $x$ is odd, then $u$ is even, but then $x^2-4 u-4$ is $5$ modulo $8$. However, by Proposition \ref{prop:z2p2div} all positive odd divisors of $2(t^2+2)$ are $1$ or $3$ modulo $8$, which is a contradiction. Hence, $x$ is even, which implies that $t$ is even. Let $x=2v$ and $t=2s$ for integers $v,t$, then the equation reduces to
$$
u(v^2-u-1)=2s^2+1.
$$
However, this equation has no solutions modulo $8$ and therefore has no integer solutions. Hence, equation \eqref{eq:6px2px2yp2y2pz2} has no integer solutions.

\vspace{10pt}

The next equation we will consider is 
\begin{equation}\label{eq:2mx2mxypx3ypz2}
2-x^2-x y+x^3 y+z^2=0,
\end{equation}
which we can rewrite as 
$$
(x-1)(x+1)(x y-1)=-(z^2+1).
$$
Modulo $4$ analysis shows $x$ is $2$ modulo $4$. Hence, $|x-1|$ and $|x+1|$ are odd positive integers that differ by $2$, thus one of them is equal to $3$ modulo $4$. However, by Proposition \ref{prop:z2p1div}, all positive odd divisors of $z^2+1$ are equal to $1$ modulo $4$, a contradiction. In conclusion, equation \eqref{eq:2mx2mxypx3ypz2} has no integer solutions.

\vspace{10pt}

The next equation we will consider is
\begin{equation}\label{eq:2px2p2ypx2ym2y2pz2}
2+x^2+2y+x^2y-2y^2+z^2=0,
\end{equation}
which we can rewrite as 
\begin{equation}\label{eq:2px2p2ypx2ym2y2pz2a}
(1+y)(4+x^2-2y)=-(z^2-2).
\end{equation}
Modulo $8$ analysis on \eqref{eq:2px2p2ypx2ym2y2pz2} shows $x$ is odd, so $x^2=1$ modulo $8$. The left-hand side of \eqref{eq:2px2p2ypx2ym2y2pz2a} is then $(y+1)(5-2y)$ modulo $8$. Also, modulo $8$ analysis shows that $y$ can only be $1,4$ or $5$ modulo $8$, hence $5-2y$ is equal to $3$ or $5$ modulo $8$. However, by Proposition \ref{prop:z2m2div}, all odd divisors of $z^2-2$ must be $\pm 1$ modulo $8$. In conclusion, equation \eqref{eq:2px2p2ypx2ym2y2pz2} has no integer solutions.

\vspace{10pt}

The next equation we will consider is
\begin{equation}\label{eq:2p2xpx3py2pxy2mz2}
	2+2 x+x^3+y^2+x y^2-z^2=0,
\end{equation}
which we can rewrite as 
$$
(x+1)(y^2+x^2-x+3)=z^2+1
$$
This representation suggests to make the change of variable $w=x+1$, which reduces the equation to
$$
w(y^2+w^2-3w+5)=z^2+1
$$
This equation has no real solutions with $w \leq 0$, so $w>0$ and $y^2+w^2-3w+5>0$. Modulo $4$ analysis  shows that $w$ is $2$ or $3$ modulo $4$, and $y$ and $w$ have the same parity. Then, $w$ or $y^2+w^2-3w+5$ is equal to $3$ modulo $4$. However, by Proposition \ref{prop:z2p1div}, all positive odd divisors of $z^2+1$ are equal to $1$ modulo $4$, a contradiction. In conclusion, equation \eqref{eq:2p2xpx3py2pxy2mz2} has no integer solutions.

\vspace{10pt}

The final equation we will consider is
\begin{equation}\label{eq:2m2xpx3my2pxy2pz2}
2-2x+x^3-y^2+x y^2+z^2=0,
\end{equation}
which we can rewrite as 
$$
-(x-1)(y^2+x^2+x-1)=z^2+1.
$$
This representation suggests to make the change of variable $w=-(x-1)$, which reduces the equation to
\begin{equation}\label{red:2m2xpx3my2pxy2pz2}
w(y^2+w^2-3w+1)=z^2+1.
\end{equation}
This equation has no real solutions with $w \leq 0$, so $w>0$ and $y^2+w^2-3w+1>0$. Modulo $4$ analysis on \eqref{red:2m2xpx3my2pxy2pz2} shows that $w$ is $2$ or $3$ modulo $4$, and $y$ and $w$ have the same parity. Then, $w$ or $y^2+w^2-3w+1$ is equal to $3$ modulo $4$. However, by Proposition \ref{prop:z2p1div}, all positive odd divisors of $z^2+1$ are equal to $1$ modulo $4$, a contradiction. In conclusion, equation \eqref{eq:2m2xpx3my2pxy2pz2} has no integer solutions.

\subsection{Exercise 7.18}\label{ex:a2_2b2}
\textbf{\emph{Prove that the equations listed in Table \ref{tab:H33} have no integer solutions. }}

		\begin{center}
	\begin{tabular}{ |c|c|c| }
		\hline
		$H$ & Equation & Representation \eqref{eq:quadform2}  \\ 
		\hline\hline
		$26$ & $2+x^3 y+y^2+z^2=0$ & $(-2+x^2)(4+2 x^2+x^4)=(x^3 + 2y)^2 + 4z^2$ \\ 
		\hline
		$28$ & $4-x^4+y^2+z^2=0$ & $(-2+x^2)(2+x^2)=y^2+z^2$\\ 
		\hline
		$29$ & $1+x^3 y+2 y^2+z^2=0$ & $(-2+x^2)(4+2 x^2+x^4)=(x^3+4 y)^2+8 z^2$ \\ 
		\hline
		$30$ & $2+x^3 y+y^2+2 z^2=0$ & $(-2+x^2)(4+2 x^2+x^4)=(x^3+2 y)^2+8 z^2$ \\
		\hline
		$33$ & $1+x^3 y+2 y^2+2 z^2 =0$ & $(-2+x^2) (4+2 x^2+x^4)=(x^3+4 y)^2+16 z^2$ \\
		\hline
		$33$ & $1+2 y+x^2 y-2 y^2+2 z-2 z^2=0$ & $(4+2x+x^2)(4-2x+x^2)=(4y-x^2-2)^2+(4z-2)^2$ \\
		\hline
		$33$ & $1+4 x-4 y+x^2 y+y^2+z^2=0$ & $(2-4 x+x^2) (6+4 x+x^2)=(-4+x^2+2 y)^2+4 z^2$ \\
		\hline
		$33$ & $1+4 y+x y+x^2 y-y^2+2 z^2=0$ & $(4+x^2) (5+2 x+x^2)=(4+x+x^2-2 y)^2- 8 z^2$ \\ 
		\hline
		$33$ & $1+4 y+x y+x^2 y-y^2-2 z^2=0$ & $(4+x^2) (5+2 x+x^2)=(4+x+x^2-2 y)^2+8 z^2$ \\ 
		\hline
		$33$ & $3+3 y+x y+x^2 y-y^2+y z+z^2=0$ & $-4 (4-x+x^2) (6+3 x+x^2)=$ \\ && $(-3-x-x^2-5 z)^2-5 (3+x+x^2-2 y+z)^2$ \\ 
		\hline
		$33$ & $5+x+x^3-y+x^2 y-y^2-z^2=0$ & $(3-2 x+x^2)(7+6 x+x^2)=(-1+x^2-2 y)^2+4 z^2$ \\ 
		\hline
		$33$ & $1+2 x+x^3+y+x y-y^2+z+x z-z^2=0$ & $2(1+2 x)(3+x^2)=(1+x-2 y)^2+(-1-x+2 z)^2$ \\ 
		\hline
	\end{tabular}
	\captionof{table}{\label{tab:H33} Equations of size $H\leq 33$ reducible to the form \eqref{eq:quadform2}, which have no integer solutions.}
\end{center} 

All equations in Table \ref{tab:H33} are of the form
\begin{equation}\label{eq:quadform2}
	QT = aR^2+bRS+cS^2,
\end{equation}
where $a,b,c$ are integers, $R$ and $S$ are arbitrary non-constant polynomials, and $Q,T$ are non-constant polynomials in one variable. Recall the notation
$$
S_{k}=\{ n \in \mathbb{Z}: n = y^2+kz^2 \, \text{for some} \, y,z \in \mathbb{Z}\},
$$
and recall that $S_{-2}$ and $S_{-5}$ are symmetric, that is
$$
S_{-2}=\{ n \in \mathbb{Z}: n = 2y^2-z^2 \, \text{for some} \, y,z \in \mathbb{Z}\}=\{ n \in \mathbb{Z}: n = y^2-2z^2 \, \text{for some} \, y,z \in \mathbb{Z}\},
$$
and
$$
S_{-5}=\{ n \in \mathbb{Z}: n = 5y^2-z^2 \, \text{for some} \, y,z \in \mathbb{Z}\}=\{ n \in \mathbb{Z}: n = y^2-5z^2 \, \text{for some} \, y,z \in \mathbb{Z}\},
$$
see Section \ref{ex:x2m2y2} and Section 5.6.3 of the book. 

To prove that the equations in Table \ref{tab:H33} have no integer solutions, we will use representation \eqref{eq:quadform2} and a relevant proposition describing the divisors of $S_k$. For some equations, we will need the following theorem and proposition.

\begin{theorem}\label{th:5y2mz2}[Theorem 5.62 in the book]
	An integer $n\neq 0$ belongs to $S_{-5}$ if and only if every prime $p$ equal to $2$ or $3$ modulo $5$ appears in the prime factorization of $n$ with an even exponent.
\end{theorem}

\begin{proposition}\label{prop:5y2mz2product}[Proposition 5.63 in the book]
	If the product $ab$ of two non-zero integers $a$ and $b$ belongs to $S_{-5}$, and $a,b$ do not share any prime factor equal to $2$ or $3$ modulo $5$, then $a \in S_{-5}$ and $b \in S_{-5}$.
\end{proposition}

Equations $2+x^3 y+y^2+z^2=0$ and $1+x^3 y+2 y^2+z^2=0$ have been solved in Section 7.2.1 of the book.

The first equation we will consider is
$$
4-x^4+y^2+z^2=0.
$$
Modulo $16$ analysis shows that $x$ is odd. The equation can be rewritten as
$$
(x^2-2)(x^2+2)=y^2+z^2.
$$
Because $x$ is odd, and we may assume that $|x| > 1$, positive integers $a = x^2 -2$ and $b = x^2 + 2$ are odd and coprime. Hence, by Proposition \ref{prop:sosproduct}, $ab$ can be the sum of two squares only if both $a$ and $b$ are. However, $a = x^2 - 2$ is congruent to $3$ modulo $4$ for every odd $x$, and therefore cannot be the sum of two squares. 

\vspace{10pt}

The next equation we will consider is
\begin{equation}\label{eq:2px3ypy2p2z2}
2+x^3y+y^2+2z^2=0.
\end{equation} 
Modulo $16$ analysis shows that $x$ is odd. Next, after multiplying \eqref{eq:2px3ypy2p2z2} by $4$, it can be rewritten as
$$
(x^2-2)(x^4+2x^2+4)=(2y+x^3)^2+2(2z)^2.
$$
Because $x$ is odd and we may assume that $|x|>1$, the positive integers $a=x^2-2$ and $b=x^4 + 2 x^2 + 4$ are odd and relatively prime, because if they had a common (odd) prime factor $p$, it would also divide $(x^2-2)(x^2+4)-(x^4 + 2 x^2 + 4)=-12$, hence $p$ must be $3$, but $x^2-2$ is never divisible by $3$. Hence, by Proposition \ref{prop:2y2pz2product}, $ab \in S_2$ only if $a \in S_2$ and $b \in S_2$. However, $b=x^2-2$ is congruent to $7$ modulo $8$, while modulo $8$ analysis shows that $S_2$ contains no integers equal to $7$ modulo $8$.

\vspace{10pt}

The next equation we will consider is
\begin{equation}\label{eq:1px3yp2y2p2z2}
1+x^3y+2y^2+2z^2=0.
\end{equation}
Modulo $2$ analysis shows that $x$ is odd. Next, after multiplying \eqref{eq:1px3yp2y2p2z2} by $8$, it can be rewritten as
$$
(x^2-2)(x^4+2x^2+4)=(x^3+4y)^2+(4z)^2.
$$
Because $x$ is odd and we may assume that $|x| > 1$, positive integers $a = x^2 -2$ and $b =x^4+2x^2+4$ are odd and relatively coprime as shown in the solution of \eqref{eq:2px3ypy2p2z2}. Hence, by Proposition \ref{prop:sosproduct}, $ab$ can be the sum of two squares only if both $a$ and $b$ are. However, $a = x^2 - 2$ is congruent to $3$ modulo $4$ for every odd $x$, and therefore cannot be the sum of two squares. 

\vspace{10pt}

The next equation we will consider is
\begin{equation}\label{eq:1p2ypx2ym2y2p2zm2z2}
1+2y+x^2y-2y^2+2z-2z^2=0.
\end{equation}
Modulo $2$ analysis shows that $x$ is odd. Next, after multiplying \eqref{eq:1p2ypx2ym2y2p2zm2z2} by $8$, it can be rewritten as
$$
(x^2+2x+4)(x^2-2x+4)=(4y-x^2-2)^2+(4z-2)^2.
$$
Because $x$ is odd, positive integers $a = x^2+2x+4$ and $b =x^2-2x+4$ are odd and relatively coprime as if $p|a$ and $p|b$ then $p|a-b=4x$ so either $p|4$ or $p|x$ but if $p|x$ then $p|a-x(x+2)=4$, so in either case $p=2$, a contradiction with $a,b$ odd. Hence, by Proposition \ref{prop:sosproduct}, $ab$ can be the sum of two squares only if both $a$ and $b$ are. However, $a = x^2 +2x+4$ is congruent to $3$ modulo $4$ for every odd $x$, and therefore cannot be the sum of two squares. 

\vspace{10pt}

The next equation we will consider is
$$
1+4x-4y+x^2y+y^2+z^2=0.
$$
Modulo $4$ analysis shows that $x$ must be odd. Next, after multiplying the equation by $4$, it can be rewritten as
$$
(x^2-4x+2)(x^2+4x+6)=(2y+x^2-4)^2+(2z)^2.
$$
Because $x$ is odd, and we may assume that $|x|>3$, positive integers $a=x^2-4x+2$ and $b=x^2+4x+6$ are odd. If they have a common (odd) prime factor $p$, it would also divide $(x^2+4x+6)-(x^2-4x+2)=4(2x+1)$, hence $p$ must divide $2x+1$, but then $p$ also divides  $4(x^2-4x+2)-(2x+1)(2x-9)=17$, hence $p=17$. Thus, $a,b$ do not share any prime factor of the form $p=4k+3$. Hence, by Proposition \ref{prop:sosproduct}, $ab$ can be the sum of two squares only if both $a$ and $b$ are. However,  $b=x^2+4x+6$ is congruent to $3$ modulo $4$ for every odd $x$, and therefore cannot be the sum of two squares.

\vspace{10pt}

The next equation we will consider is
$$
1+4y+xy+x^2y-y^2+2z^2=0.
$$
Multiplying this equation by $4$ and rearranging, we obtain
$$
(x^2+4)(x^2+2x+5) = (4 + x + x^2 - 2y)^2 -2(2z)^2.
$$
If the positive integers $a=x^2+4$ and $b=x^2+2x+5$ had a common prime factor $p$, it would also be a factor of $(x^2+2x+5)-(x^2+4)=2x+1$, so $p$ must divide $4(x^2+2x+5)-(2x+1)(2x+3)=17$, hence $p=17$. Thus, $p$ is not equal to $3$ or $5$ modulo $8$. By Proposition \ref{prop:2y2mz2}, $ab \in S_{-2}$  only if $a \in S_{-2}$ and $b \in S_{-2}$. If $x$ is odd, then $a=x^2+4$ is congruent to $5$ modulo $8$, while if $x$ is even, then $b=x^2+2x+5$ is congruent to $5$ modulo $8$. 
However, modulo $8$ analysis shows that $S_{-2}$ contains no integers equal to $5$ modulo $8$.

\vspace{10pt}

The next equation we will consider is
$$
1+4y+xy+x^2y-y^2-2z^2=0.
$$
Multiplying this equation by $4$ and rearranging, we obtain
$$
(x^2+4)(x^2+2x+5)=(x^2+x+4-2y)^2+2(2z)^2.
$$
We have that for all integers $x$, $a=x^2+4$ and $b=x^2+2x+5$ are positive integers. If they have a common prime factor, it would also divide $(x^2+2x+5)-(x^2+4)=2x+1$, so $p$ must divide $2x+1$ but then $4(x^2+2x+5)-(2x+1)(2x+3)=17$, hence $p=17$. Thus $p$ is not equal to $5$ or $7$ modulo $8$. By Proposition \ref{prop:2y2pz2product}, $ab \in S_2$ only if $a \in S_2$ and $b \in S_2$. If $x$ is odd, then $a=x^2+4$ is congruent to $5$ modulo $8$, while if $x$ is even, then $b=x^2+2x+5$ is congruent to $5$ modulo $8$. However modulo $8$ analysis shows that $S_2$ contains no integers equal to $5$ modulo $8$.

\vspace{10pt}

The next equation we will consider is
$$
3+3y+xy+x^2y-y^2+yz+z^2=0.
$$
Multiplying this equation by $20$ and rearranging, we obtain
$$
4(x^2-x+4)(x^2+3x+6) = 5(3 + x + x^2 - 2y + z)^2 - (3+x + x^2 + 5z)^2.
$$
For all $x$, $a=\frac{x^2-x+4}{2}$ and $b=\frac{x^2+3x+6}{2}$ are positive integers. Then the equation can be rewritten as
$$
16ab=5Z^2-Y^2,
$$
where $Y=3+x + x^2 + 5z$ and $Z=3 + x + x^2 - 2y + z$.
Hence, $16 a b \in S_{-5}$. By Theorem \ref{th:5y2mz2}, this is equivalent to $ab \in S_{-5}$. If $a$ and $b$ have any common prime factor $p$, then $p$ would also divide $b-a=2x+1$ and then $p$ would divide $8a-(2x+1)(2x-3)=19$, hence $p=19$. Thus $p$ is not equal to $2$ or $3$ modulo $5$, and by Proposition \ref{prop:5y2mz2product}, $ab \in S_{-5}$ only if $a \in S_{-5}$ and $b \in S_{-5}$. However, if $x=0$ or $1$ modulo $5$ then $a$ is congruent to $2$ modulo $5$. If $x=2$ or $4$ modulo $5$, then $a$ is congruent to $3$ modulo $5$, and if $x=3$ modulo $5$, then $b$ is congruent to $2$ modulo $5$. 
On the other hand, modulo $5$ analysis shows that no integer equal to $2$ or $3$ modulo $5$ can belong to $S_{-5}$.

\vspace{10pt}

The next equation we will consider is
$$
5+x+x^3-y+x^2y-y^2-z^2=0.
$$
Modulo $4$ analysis shows that $x$ must be even. Next, after multiplying the equation by $4$, it can be rewritten as 
$$
(3-2x+x^2)(7+6x+x^2)=(-1+x^2-2y)^2+(2z)^2.
$$
Direct substitution shows that equation has no integer solutions with $|x|\leq 4$. 
Because $x$ is even and $|x|>4$, integers $a=x^2-2x+3$ and $b=x^2+6x+7$ are both positive and odd. If they have a common (odd) prime factor $p$, it would also divide $(x^2+6x+7)-(x^2-2x+3)=8x+4$, hence $p$ must divide $2x+1$. But then, $p$ must also divide $4(x^2+6x+7)-(11+2x)(2x+1)=17$, hence $p=17$. Thus, $a,b$ do not share any prime factor of the form $p=4k+3$. Hence, by Proposition \ref{prop:sosproduct}, $ab$ can be the sum of two squares only if both $a$ and $b$ are. However, for any even $x$, $b=x^2+6x+7$ is congruent to $3$ modulo $4$, and therefore cannot be the sum of two squares.

\vspace{10pt}

The final equation we will consider is
$$
1+2x+x^3+y+xy-y^2+z+xz-z^2=0.
$$
Modulo $2$ analysis shows that $x$ must be odd. Next, after multiplying the equation by $4$, it can be rewritten as
$$
2(x^2+3)(2x+1)=(1+x-2y)^2+(-1-x+2z)^2.
$$ 
Because the right-hand side and $a=x^2+3>0$ are positive, this implies that odd integer $b=2x+1$ is positive. If integers $a$ and $b$ have a common (odd) prime factor, it would also divide $4 (x^2 + 3) - (2 x + 1) (2 x - 1)=13$, hence $p=13$. Thus, $a,b$ do not share any prime factor of the form $p=4k+3$. Hence, by Proposition \ref{prop:sosproduct}, $ab$ can be the sum of two squares only if both $a$ and $b$ are. However, for any odd $x$, $b=2x+1$ is congruent to $3$ modulo $4$, and therefore cannot be the sum of two squares. 

\subsection{Exercise 7.19}\label{ex:H32sol}
\textbf{\emph{Find an integer solution for each of the equations listed in Table \ref{tab:H32sol}.}}

	\begin{center}

		\captionof{table}{\label{tab:H32sol} Some equations with large integer solutions.}
	\end{center} 

Table \ref{tab:H32sols} presents an integer solution for each equation in Table \ref{tab:H32sol} found by a computer search, see the discussion in Section 7.2.2 of the book for details how this computer search has been performed.

		\begin{center}
			%
			\captionof{table}{\label{tab:H32sols} Examples of an integer solution for each of the equations listed in Table \ref{tab:H32sol}.}
		\end{center}

\subsection{Exercise 7.24}\label{ex:2varsolsearch}
\textbf{\emph{Check that the equations listed in Table \ref{tab:2varH34} have no integer solutions with $|y|\leq 10^6$.}}

 	\begin{center}
	%
	\captionof{table}{\label{tab:2varH34} Two-variable equations of size $H\leq 34$ with open Problem \ref{prob:main}.}
\end{center}

Notice that for each equation in Table \ref{tab:2varH34}, when $x$ is large, we must have $y>0$, as otherwise the left-hand side will be negative, whilst the right-hand side is a large positive. In fact, every equation in Table \ref{tab:2varH34} has no real solutions with $y<-2$. Hence, it is sufficient to search the region $-1 \leq y \leq 10^6$. 

To show these equations have no integer solutions with $|y| \leq 10^6$, we can use the Mathematica code
$$
\begin{aligned}
{\tt For[ }& {\tt y = -1, y \leq 10^6, y++, } \\
& {\tt p = Reduce[P(x,y)==0, x, Integers]; } \\
& {\tt If[Length[p] > 0, Print[``the \, solution \, is", p]]]}
\end{aligned}
$$
For each equation $P(x,y)=0$ in Table \ref{tab:2varH34}, the Mathematica code does not output any solutions.

\subsection{Exercise 7.25}\label{ex:h33cubica}
\textbf{\emph{By analysing separately positive and negative values of $x$, prove that the equation
	\begin{equation}\label{eq:1px2px3py2pxy2mxz2}
		1+x^2+x^3+y^2+x y^2-x z^2 = 0 
	\end{equation}
	has no integer solutions.}}
	
It is easy to see that this equation has no solutions with $x=0$.
Now, we may consider the cases $x>0$ and $x<0$ separately. Let us begin with $x<0$. 
We can rearrange \eqref{eq:1px2px3py2pxy2mxz2} to
\begin{equation}\label{red1:1px2px3py2pxy2mxz2}
	x(y^2-z^2+x^2+x)=-y^2-1.
\end{equation}
We then make the change of variable $x \to -x$, giving the equation
\begin{equation}\label{red2:1px2px3py2pxy2mxz2}
	x(y^2-z^2+x^2-x)=y^2+1
\end{equation}
with $x>0$, hence $y^2-z^2+x^2-x>0$. Proposition \ref{prop:z2p1div} states that $y^2+1$ can only have positive divisors equal to $1$ or $2$ modulo $4$. If $x$ is $1$ modulo $4$, then $y^2-z^2+x^2-x$ reduces to $y^2-z^2$ modulo $4$, which can only be $1$ or $2$ modulo $4$ if $y$ is odd and $z$ is even. However the left-hand side of \eqref{red2:1px2px3py2pxy2mxz2} is then odd, whilst the right-hand side is even, a contradiction. If $x$ is $2$ modulo $4$, then $y^2-z^2+x^2-x$ reduces to $y^2-z^2+2$ modulo $4$, which can only be $1$ or $2$ modulo $4$ if $y$ is even, or both $y$ and $z$ are odd. If $y$ is even, then the left-hand side of \eqref{red2:1px2px3py2pxy2mxz2} is even, whilst the right-hand side is odd, a contradiction. If both $y$ and $z$ are odd, then the left-hand side of \eqref{red2:1px2px3py2pxy2mxz2} is divisible by $4$ whilst the right-hand side is only divisible by $2$, a contradiction.

Let us now consider the case $x>0$ in \eqref{eq:1px2px3py2pxy2mxz2}. The substitution $x=w-1$ reduces the equation to
\begin{equation}\label{red2b:1px2px3py2pxy2mxz2}
	w(-1+2w-w^2-y^2+z^2) = 1 + z^2.
\end{equation}
By Proposition \ref{prop:z2p1div}, $w$ must be either $1$ or $2$ modulo $4$. 
If $w$ is $1$ modulo $4$, then $-1+2w-w^2-y^2+z^2$ reduces to $z^2-y^2$ modulo $4$. It is impossible for a difference of squares to be $2$ modulo $4$, hence $z^2-y^2$ must be $1$ modulo $4$, which then implies that $z$ is odd. However the left-hand side of \eqref{red2b:1px2px3py2pxy2mxz2} is then odd, whilst the right-hand side is even, a contradiction. If $w$ is $2$ modulo $4$, then $-1+2w-w^2-y^2+z^2$ reduces to $z^2-y^2-1$ modulo $4$, which can only be $1$ or $2$ modulo $4$ if $z$ is even. However the left-hand side of \eqref{red2b:1px2px3py2pxy2mxz2} is then even, whilst the right-hand side is odd, a contradiction. Therefore, equation \eqref{eq:1px2px3py2pxy2mxz2} has no integer solutions.

\subsection{Exercise 7.32}\label{ex:h33formeropen}
\textbf{\emph{Prove that equation 
		\begin{equation}\label{eq:h33cubicopen2}
			2 x^2 - xyz + 2 z^2 = y^2 - 9 y + 23
		\end{equation}
		 has no integer solutions.}}

To solve this equation, we will need to use the following propositions.
\begin{proposition}\label{prop:Jacobi}[Proposition 7.26 in the book]
	For any integers $a,b$ and positive odd integers $m,n$, the following properties hold.
	\begin{itemize}
		\item[(J1)] If $a \equiv b\,(\text{mod}\, n)$ then $\left(\frac{a}{n}\right) = \left(\frac{b}{n}\right)$.
		\item[(J2)] The multiplicative properties:
		$
		\left(\frac{ab}{n}\right) = \left(\frac{a}{n}\right)\left(\frac{b}{n}\right)
		$
		and
		$
		\left(\frac{a}{nm}\right) = \left(\frac{a}{n}\right)\left(\frac{a}{m}\right).
		$
		\item[(J3)] 
		$
		\left(\frac{-1}{n}\right) = (-1)^{\frac{n-1}{2}}
		$
		and
		$
		\left(\frac{2}{n}\right) = (-1)^{\frac{n^2-1}{8}}.
		$
		\item[(J4)] Law of quadratic reciprocity: if $m$ and $n$ are odd positive coprime integers, then
		$$
		\left(\frac{m}{n}\right) \left(\frac{n}{m}\right) = (-1)^{\frac{m-1}{2}\frac{n-1}{2}}.
		$$
		\item[(J5)] If $a$ is a quadratic residue modulo $n$ and $\gcd(a,n) = 1$, then $\left(\frac{a}{n}\right)=1$.
	\end{itemize} 
\end{proposition}

\begin{proposition}\label{prop:evenmult}[Proposition 5.27 in the book]
	Let $a,b,c$ be integers, and let $p$ be an odd prime not dividing $D=b^2-4ac$, such that $\left(\frac{D}{p}\right)=-1$. Then, for any integers $x$ and $y$, $p$ appears in the prime factorization of $m=ax^2+bxy+cy^2$ with an even exponent.
\end{proposition}

\begin{proposition}\label{prop:h33open}[Proposition 7.29 in the book]
	Assume that \eqref{eq:h33cubicopen2} has an integer solution, and let $f_1 = y^2-16$, $f_2=y^2-9y+23$, and $d=\gcd(f_1,f_2)$. Then $f_1,f_2$ are odd positive integers, and $d$ is a divisor of $75$. If $d=1,3,5$ or $15$, then $\gcd(f_1/d,8f_2)=1$ and $\left(\frac{8f_2}{f_1/d}\right)=1$. If $d=25$ or $75$, then $\gcd(f_1/d,8f_2/25)=1$ and $\left(\frac{8f_2/25}{f_1/d}\right)=1$.
\end{proposition}

The properties $(J2)$, $(J3)$ and $(J4)$ in Proposition \ref{prop:Jacobi} imply that
\begin{equation}\label{eq:reposm2}
	\left(\frac{-2}{n}\right)=\left(\frac{-1}{n}\right)\left(\frac{2}{n}\right) = (-1)^{\frac{n-1}{2}} (-1)^{\frac{n^2-1}{8}} =  
	\begin{cases}
		1, & \text{if $n \equiv 1 \, \text{or} \, 3 \, (\text{mod}\, 8)$,}\\
		-1, & \text{if $n \equiv -1 \, \text{or} \, -3 \, (\text{mod}\, 8)$.}
	\end{cases}
\end{equation}

We now recall the following information from Section 7.3.2 of the book for convenience. Consider any equation of the form
\begin{equation}\label{eq:jacobynosol}
	A(x) y^2 + B(x) yz + C(x) z^2 = P_1(x)P_2(x), 
\end{equation}
where $A(x),B(x),C(x),P_1(x),P_2(x)$ are polynomials in $x$ with integer coefficients. Let $D(x) = B^2(x) - 4 A(x) C(x)$. Then for each $i=1,2$, we have the following result. Let $x$ be any integer such that 
\begin{itemize}
	\item[(i)] $P_i(x)$ is positive and odd;
	\item[(ii)] $\text{gcd}(P_i(x),D(x))=1$;
	\item[(iii)] $P_1(x)$ and $P_2(x)$ have no common prime factor $p$ such that $\left(\frac{D(x)}{p}\right)=-1$.
\end{itemize}
Then Proposition \ref{prop:evenmult} states that any prime factor $p$ of $P_i(x)$ such that $\left(\frac{D(x)}{p}\right)=-1$ appears in the prime factorization of $A(x) y^2 + B(x) yz + C(x) z^2$ (and hence of $P_i(x)$) with an even exponent. By definition of Jacobi symbol, this implies that
\begin{equation}\label{eq:jacobidpi1}
	\left(\frac{D(x)}{P_i(x)}\right) = 1.
\end{equation}

To prove that an equation has no integer solutions, if we can use the law of quadratic reciprocity to show that $\left(\frac{D(x)}{P_i(x)}\right) = -1$, this is a desired contradiction.

Therefore to prove that \eqref{eq:h33cubicopen2} has no integer solutions, by Proposition \ref{prop:h33open} and the above theory, it suffices to prove that for every odd integer $y\geq 5$ we have $\left(\frac{8f_2}{f_1/d}\right)=-1$ if $d=1,3,5$ or $15$, and $\left(\frac{8f_2/25}{f_1/d}\right)=-1$ if $d=25$ or $75$, where  $d=\gcd(f_1,f_2)$, and $f_1=y^2-16$ and $f_2=y^2-9y+23$. The cases $d=1$ and $d=3$ have been proven in Section 7.3.2 of the book. So, the first case we need to consider is $d=5$.

\subsubsection*{The case $d=5$}

Let us first consider the case with $d=5$. In this case, both $y^2-9y+23$ and $y^2-16$ are divisible by $5$, hence $y=1$ modulo $5$. Also, $y^2-9y+23$ and $y^2-16$ do not share factor $3$, which is only possible if $y=0$ modulo $3$. 
We have
$$
\left(\frac{8f_2}{f_1/d}\right) = \left(\frac{8(y^2-9y+23)}{(y^2-16)/5}\right) = \left(\frac{2}{(y^2-16)/5}\right)^3\cdot \left(\frac{y^2-9y+23}{y-4}\right) \cdot \left(\frac{y^2-9y+23}{(y+4)/5}\right),
$$
where we have used (J2) and the fact that $y$ is odd and $1$ modulo $5$. Because $y^2-16$ is $1$ modulo $8$, $(y^2-16)/5$ must be $5$ modulo $8$, and Proposition \ref{prop:Jacobi} (J3) implies that $\left(\frac{2}{(y^2-16)/5}\right)=-1$. Further, 
\begin{equation}\label{eq:h33cubicyaux}
	y^2-9y+23 = (y-4)(y-5)+3 = (y+4)(y-13)+75
\end{equation}
implies that $y^2-9y+23$ is equal to either $3$ or $75=3\cdot 5^2$ modulo $y-4$ and modulo $(y+4)/5$. Hence, by Proposition \ref{prop:Jacobi} (J1)
$$
\left(\frac{8f_2}{f_1/5}\right) = (-1)^3 \left(\frac{3}{(y+4)/5}\right) \cdot \left(\frac{3}{y-4}\right).
$$
Because $y=0$ modulo $3$, we have $\text{gcd}(3,(y+4)/5)=\text{gcd}(3,y-4)=1$. Then we can apply Proposition \ref{prop:Jacobi} (J4) and conclude that
$$
\left(\frac{8f_2}{f_1/5}\right) = - \left(\frac{(y+4)/5}{3}\right) (-1)^{\frac{-1+(y+4)/5}{2}} \cdot \left(\frac{y-4}{3}\right) (-1)^{\frac{-1+(y-4)}{2}}.
$$
The integers $\frac{-1+(y+4)/5}{2}$ and $\frac{-1+(y-4)}{2}$ have the same parity, thus we have
$$
(-1)^{\frac{-1+(y+4)/5}{2}} (-1)^{\frac{-1+(y-4)}{2}}=1.
$$
Hence,
$$
\left(\frac{8f_2}{f_1/5}\right) = - \left(\frac{(y+4)/5}{3}\right)  \cdot \left(\frac{y-4}{3}\right) = - \left(\frac{2}{3}\right)  \cdot \left(\frac{-4}{3}\right) =-(-1)(-1)=-1.
$$

\subsubsection*{The case $d=15$}

Let us now consider the case with $d=15$. In this case, $y=1$ or $11$ modulo $15$. 
We will first consider the case when $y=1$ modulo $15$. Then we have
$$
\left(\frac{8f_2}{f_1/d}\right) = \left(\frac{8(y^2-9y+23)}{(y^2-16)/15}\right) = \left(\frac{2}{(y^2-16)/15}\right)^3\cdot \left(\frac{y^2-9y+23}{(y-4)/3}\right) \cdot \left(\frac{y^2-9y+23}{(y+4)/5}\right),
$$
where we have used Proposition \ref{prop:Jacobi} (J2) and the fact that $y$ is odd and $1$ modulo $15$. Because $y^2-16$ is $1$ modulo $8$, $(y^2-16)/15$ must be $7$ modulo $8$, and Proposition \ref{prop:Jacobi} (J3) implies that $\left(\frac{2}{(y^2-16)/15}\right)=1$.
Further, \eqref{eq:h33cubicyaux} implies that $y^2-9y+23$ is equal to either $3$ or $75=3\cdot 5^2$ modulo $(y+4)/5$ and modulo $(y-4)/3$. Hence, by Proposition \ref{prop:Jacobi} (J1)
$$
\left(\frac{8f_2}{f_1/15}\right) = \left(\frac{3}{(y+4)/5}\right) \cdot \left(\frac{3}{(y-4)/3}\right).
$$
As $d=15$, we have $\text{gcd}(f_1/15,8f_2)=1$ by Proposition \ref{prop:h33open}. Because $f_2$ is divisible by $3$, $f_1/15$ is not. Hence, $\gcd(3,(y+4)/5)=\gcd(3,(y-4)/3)=1$. Then we can apply Proposition \ref{prop:Jacobi} (J4) and obtain that
$$
\left(\frac{8f_2}{f_1/15}\right) = \left(\frac{(y+4)/5}{3}\right) (-1)^{\frac{-1+(y+4)/5}{2}} \cdot \left(\frac{(y-4)/3}{3}\right) (-1)^{\frac{-1+(y-4)/3}{2}}.
$$
Then because $y$ is odd and $y=1$ modulo $15$, we can make the change of variable $y=30k+1$ for integer $k$, so,
$$
(-1)^{\frac{-1+(y+4)/5}{2}}  (-1)^{\frac{-1+(y-4)/3}{2}}=(-1)^{8k-1} = -1.
$$
Hence,
$$
\left(\frac{8f_2}{f_1/15}\right) = - \left(\frac{(y+4)/5}{3}\right) \cdot \left(\frac{(y-4)/3}{3}\right) = - \left(\frac{1}{3}\right) \cdot \left(\frac{k-1}{3}\right) = -  \left(\frac{k-1}{3}\right).
$$
Because $\gcd(3,(y-4)/3)=1$, $k-1$ cannot be divisible by $3$. If $k$ is divisible by $3$, then $y=1$ modulo $9$, but the original equation has no solutions with $y=1$ modulo $9$. Hence, $k=2$ modulo $3$, and
$$
\left(\frac{8f_2}{f_1/15}\right) =  -  \left(\frac{k-1}{3}\right)=-\left(\frac{1}{3}\right)=-1.
$$
We must now consider the case where $y=11$ modulo $15$. Then we have
$$
\left(\frac{8f_2}{f_1/d}\right) = \left(\frac{8(y^2-9y+23)}{(y^2-16)/15}\right) = \left(\frac{2}{(y^2-16)/15}\right)^3\cdot \left(\frac{y^2-9y+23}{y-4}\right) \cdot \left(\frac{y^2-9y+23}{(y+4)/15}\right),
$$
where we have used Proposition \ref{prop:Jacobi} (J2) and the fact that $y$ is odd and $11$ modulo $15$. Because $y^2-16$ is $1$ modulo $8$, $(y^2-16)/15$ must be $7$ modulo $8$, and Proposition \ref{prop:Jacobi} (J3) implies that $\left(\frac{2}{(y^2-16)/15}\right)=1$.
Further, \eqref{eq:h33cubicyaux} implies that $y^2-9y+23$ is equal to either $3$ or $75=3\cdot 5^2$ modulo $(y+4)/15$ and modulo $y-4$. Hence, by Proposition \ref{prop:Jacobi} (J1)
$$
\left(\frac{8f_2}{f_1/15}\right) = \left(\frac{3}{(y+4)/15}\right) \cdot \left(\frac{3}{y-4}\right).
$$
As $d=15$, we have $\text{gcd}(f_1/15,8f_2)=1$ by Proposition \ref{prop:h33open}. Because $f_2$ is divisible by $3$, $f_1/15$ is not. Hence, $\gcd(3,(y+4)/15)=\gcd(3,y-4)=1$. Then we can apply Proposition \ref{prop:Jacobi} (J4) and obtain that
$$
\left(\frac{8f_2}{f_1/15}\right) = \left(\frac{(y+4)/15}{3}\right) (-1)^{\frac{-1+(y+4)/15}{2}} \cdot \left(\frac{y-4}{3}\right) (-1)^{\frac{-1+y-4}{2}}.
$$
Then because $y$ is odd and $y=11$ modulo $15$, we can make the change of variable $y=30k+11$ for integer $k$, so,
$$
(-1)^{\frac{-1+(y+4)/15}{2}} (-1)^{\frac{y-5}{2}}= (-1)^{16k+3} = -1.
$$
Hence,
$$
\left(\frac{8f_2}{f_1/15}\right) = - \left(\frac{(y+4)/15}{3}\right) \cdot \left(\frac{y-4}{3}\right) = - \left(\frac{2k+1}{3}\right) \cdot \left(\frac{1}{3}\right) = -  \left(\frac{2k+1}{3}\right).
$$

We can also deduce that $k=0$ modulo $3$, because otherwise 
$y=5$ or $8$ modulo $9$ for which the original equation has no solutions. 
Finally, we may conclude that
$$
\left(\frac{8f_2}{f_1/15}\right) =  -  \left(\frac{2k+1}{3}\right)=-\left(\frac{1}{3}\right)=-1.
$$

\subsubsection*{The case $d=25$}

Let us now consider the case with $d=25$. In this case, $y=21$ modulo $25$. We have
$$
\begin{aligned}
	\left(\frac{8f_2/25}{f_1/d}\right) & = \left(\frac{8(y^2-9y+23)/25}{(y^2-16)/25}\right) \\ & = \left(\frac{2}{(y^2-16)/25}\right)^3\cdot \left(\frac{(y^2-9y+23)/25}{y-4}\right) \cdot \left(\frac{(y^2-9y+23)/25}{(y+4)/25}\right),
\end{aligned}
$$
where we have used Proposition \ref{prop:Jacobi} (J2) and the fact that $y$ is odd and $21$ modulo $25$. Because $y^2-16$ is $1$ modulo $8$, $(y^2-16)/25$ must be $1$ modulo $8$, and Proposition \ref{prop:Jacobi} (J3) implies that $\left(\frac{2}{(y^2-16)/25}\right)=1$. Proposition \ref{prop:h33open} implies that $\gcd(25,y-4)=\gcd(25,(y+4)/25)=1$. Using this and Proposition \ref{prop:Jacobi} (J2),
$$
\left(\frac{y^2-9y+23}{y-4}\right)=\left(\frac{(y^2-9y+23)/25}{y-4}\right) \cdot \left(\frac{25}{y-4}\right)=\left(\frac{(y^2-9y+23)/25}{y-4}\right) \cdot \left(\frac{5}{y-4}\right)^2
$$
and
$$
\left(\frac{y^2-9y+23}{(y+4)/25}\right) =\left(\frac{(y^2-9y+23)/25}{(y+4)/25}\right) \cdot \left(\frac{5}{(y+4)/25}\right)^2.
$$
Then
$$
\left(\frac{8f_2/25}{f_1/d}\right) = 1^3\cdot \left(\frac{y^2-9y+23}{y-4}\right) \cdot \left(\frac{y^2-9y+23}{(y+4)/25}\right).
$$
Further, \eqref{eq:h33cubicyaux} and Proposition \ref{prop:Jacobi} (J1) implies
$$
\left(\frac{8f_2/25}{f_1/d}\right) = \left(\frac{3}{y-4}\right) \cdot \left(\frac{3}{(y+4)/25}\right).
$$
As $d=25$, we have $\text{gcd}(f_1/25,8f_2/25)=1$ by Proposition \ref{prop:h33open}. Then we can apply Proposition \ref{prop:Jacobi} (J4) and conclude that
$$
\left(\frac{8f_2/25}{f_1/25}\right) = \left(\frac{(y+4)/25}{3}\right) (-1)^{\frac{-1+(y+4)/25}{2}} \cdot \left(\frac{y-4}{3}\right) (-1)^{\frac{-1+y-4}{2}}.
$$

Then because $y$ is odd and $y=21$ modulo $25$, we can make the change of variable $y=50k+21$ for integer $k$, so,
$$
(-1)^{\frac{-1+(y+4)/25}{2}} (-1)^{\frac{y-5}{2}}= (-1)^{2(13k+4)} = 1.
$$
Then,
$$
\begin{aligned}
	\left(\frac{8f_2/25}{f_1/25}\right) = & \left(\frac{(y+4)/25}{3}\right) \cdot \left(\frac{y-4}{3}\right) =\left(\frac{2k+1}{3}\right) \cdot \left(\frac{50k+17}{3}\right) \\ = & \left(\frac{2k+1}{3}\right)   \cdot \left(\frac{2}{3}\right) \cdot \left(\frac{k+1}{3}\right)= - \left(\frac{2k+1}{3}\right) \cdot \left(\frac{k+1}{3}\right).
\end{aligned}$$
Because $\gcd(3,(y+4)/25)=\gcd(3,y-4)=1$, $k+1$ and $2k+1$ cannot be divisible by $3$ which means that $k=0$ modulo $3$. Hence
$$
\left(\frac{8f_2/25}{f_1/25}\right) = - \left(\frac{2k+1}{3}\right) \cdot \left(\frac{k+1}{3}\right) = - \left(\frac{1}{3}\right) \cdot \left(\frac{1}{3}\right) = -1.
$$

\subsubsection*{The case $d=75$}

Let us now consider the case with $d=75$. In this case, $y=46$ or $71$ modulo $75$. We will first consider the case when $y=46$ modulo $75$. Then we have
$$
\begin{aligned}
	\left(\frac{8f_2/25}{f_1/d}\right) & = \left(\frac{8(y^2-9y+23)/25}{(y^2-16)/75}\right) \\ & = \left(\frac{2}{(y^2-16)/75}\right)^3\cdot \left(\frac{(y^2-9y+23)/25}{(y-4)/3}\right) \cdot \left(\frac{(y^2-9y+23)/25}{(y+4)/25}\right),
\end{aligned}
$$
where we have used Proposition \ref{prop:Jacobi} (J2) and the fact that $y$ is odd and $46$ modulo $75$. Because $y^2-16$ is $1$ modulo $8$, $(y^2-16)/75$ must be $3$ modulo $8$, and Proposition \ref{prop:Jacobi} (J3) implies that $\left(\frac{2}{(y^2-16)/75}\right)=-1$. Proposition \ref{prop:h33open} implies $\gcd(25,(y-4)/3)=\gcd(25,(y+4)/25)=1$. Using this and Proposition \ref{prop:Jacobi} (J3)
$$
\left(\frac{8f_2/25}{f_1/d}\right) =-\cdot \left(\frac{y^2-9y+23}{(y-4)/3}\right) \cdot \left(\frac{y^2-9y+23}{(y+4)/25}\right) 
=- \left(\frac{3}{(y-4)/3}\right) \cdot \left(\frac{3}{(y+4)/25}\right).
$$
As $d=75$, we have $\gcd(f_1/75,8f_2/25)=1$ by Proposition \ref{prop:h33open}. Then we can apply Proposition \ref{prop:Jacobi} (J4) and obtain that
$$
\left(\frac{8f_2/25}{f_1/75}\right) =- \left(\frac{y+4}{3}\right) (-1)^{\frac{-1+(y+4)/25}{2}} \cdot \left(\frac{(y-4)/3}{3}\right) (-1)^{\frac{-1+(y-4)/3}{2}}.
$$

Then because $y$ is odd and $y=46$ modulo $75$, we can make the change of variable $y=75(2k-1)+46$ for integer $k$, so,
$$
(-1)^{\frac{-1+(y+4)/25}{2}} (-1)^{\frac{-1+(y-4)/3}{2}}= (-1)^{28k-7} = -1.
$$
Then,
$$
\begin{aligned}
	\left(\frac{8f_2/25}{f_1/75}\right) = & \left(\frac{75(2k-1)+46+4}{3}\right) \cdot \left(\frac{(75(2k-1)+46-4)/3}{3}\right) \\  = &  \left(\frac{150k-25}{3}\right) \cdot \left(\frac{50k-11}{3}\right) = \left(\frac{2}{3}\right) \cdot \left(\frac{2k+1}{3}\right) = - \left(\frac{2k+1}{3}\right).
\end{aligned}
$$
Because $\gcd(3,(y-4)/3)=1$, $2k+1$ cannot be divisible by $3$, so $k \neq 1$ modulo $3$. If $k=2$ modulo $3$, then $y=1$ modulo $9$, but the original equation has no solutions with $y=1$ modulo $9$. Hence, $k=0$ modulo $3$, and
$$
\left(\frac{8f_2/25}{f_1/75}\right) =  -  \left(\frac{2k+1}{3}\right)=-\left(\frac{1}{3}\right)=-1.
$$

We will next consider the case when $y=71$ modulo $75$. Then we have
$$
\left(\frac{8f_2/25}{f_1/d}\right) = \left(\frac{2}{(y^2-16)/75}\right)^3\cdot \left(\frac{(y^2-9y+23)/25}{y-4}\right) \cdot \left(\frac{(y^2-9y+23)/25}{(y+4)/75}\right),
$$
where we have used Proposition \ref{prop:Jacobi} (J2). Because $y^2-16$ is $1$ modulo $8$, $(y^2-16)/75$ must be $3$ modulo $8$, and Proposition \ref{prop:Jacobi} (J3) implies that $\left(\frac{2}{(y^2-16)/75}\right)=-1$. Proposition \ref{prop:h33open} implies $\gcd(25,y-4)=\gcd(25,(y+4)/75)=1$. Using this and Proposition \ref{prop:Jacobi} (J3)
$$
\left(\frac{8f_2/25}{f_1/d}\right) =(-1)^3\cdot \left(\frac{y^2-9y+23}{y-4}\right) \cdot \left(\frac{y^2-9y+23}{(y+4)/75}\right) =- \left(\frac{3}{y-4}\right) \cdot \left(\frac{3}{(y+4)/75}\right).
$$
As $d=75$, we have $\text{gcd}(f_1/75,8f_2/25)=1$ by Proposition \ref{prop:h33open}. Then we can apply Proposition \ref{prop:Jacobi} (J4) and obtain that
$$
\left(\frac{8f_2/25}{f_1/75}\right) =- \left(\frac{(y+4)/75}{3}\right) (-1)^{\frac{-1+(y+4)/75}{2}} \cdot \left(\frac{y-4}{3}\right) (-1)^{\frac{-1+y-4}{2}}.
$$
Then because $y$ is odd and $y=71$ modulo $75$, we can make the change of variable $y=150k+71$ for integer $k$, so,
$$
(-1)^{\frac{-1+(y+4)/75}{2}} (-1)^{\frac{-1+y-4}{2}}= (-1)^{76k+33} = -1.
$$

Then,
$$
\left(\frac{8f_2/25}{f_1/75}\right) = \left(\frac{2k+1}{3}\right) \cdot \left(\frac{150k+67}{3}\right) =  \left(\frac{2k+1}{3}\right) \cdot \left(\frac{1}{3}\right)  = \left(\frac{2k+1}{3}\right).
$$
Because $\gcd(3,(y+4)/75)=1$, $2k+1$ cannot be divisible by $3$, so $k \neq 1$ modulo $3$. If $k=0$ modulo $3$, then $y=8$ modulo $9$, but the original equation has no solutions with $y=8$ modulo $9$. Hence, $k=2$ modulo $3$, and
$$
\left(\frac{8f_2/25}{f_1/75}\right) = \left(\frac{2k+1}{3}\right)=\left(\frac{5}{3}\right)=-1.
$$

\subsection{Exercise 7.33}\label{ex:h34cubic}
\textbf{\emph{Prove that the equations
	$$
	2+x^2+2 x^2 y+y^2+x z^2 = 0 
	$$
	and
	$$
		2+2x^2+x^2y-y^2+x y z+z^2=0
	$$
	of size $H=34$ have no integer solutions.}}

To prove that the equations in this exercise have no integer solutions, we will use the Jacobi symbol to obtain a contradiction, see Section \ref{ex:h33formeropen}.

Let us first prove that the equation
\begin{equation}\label{eq:2px2p2x2ypy2pxz2}
	2+x^2+2x^2y+y^2+x z^2=0
\end{equation}
has no integer solutions.

The equation can be rewritten as $(y+x^2)^2+xz^2=x^4-x^2-2$, or after substitution $y+x^2 \to y$ as $y^2+xz^2=x^4-x^2-2=(x^2+1)(x^2-2)$. Note that $x^2+1$ and $x^2-2$ are coprime because the only potential common prime factor is $3$, but $x^2-2$ is never divisible by $3$. Modulo $8$ analysis on equation \eqref{eq:2px2p2x2ypy2pxz2} shows that $x$ must be $2,5,$ or $6$ modulo $8$, hence either $x=4a+2$ or $x=8a+5$ for some integer $a$.

Let us first consider the case where $x=4a+2$.
Then
\eqref{eq:jacobidpi1} implies that
\begin{equation}\label{eq:prop:2px2p2x2ypy2pxz2}
	\left(\frac{-4(4a+2)}{(4a+2)^2+1}\right)=\left(\frac{-(4a+2)}{(4a+2)^2+1}\right)=1.
\end{equation}
Let us compute this Jacobi symbol. First assume that $a\geq 0$. Then
\begin{equation}\label{eq:2px2p2x2ypy2pxz2_4ap2}
	\left(\frac{-(4a+2)}{(4a+2)^2+1}\right)=\left(\frac{-2}{(4a+2)^2+1}\right) \cdot \left(\frac{2a+1}{(4a+2)^2+1}\right).
\end{equation}
As $(4a+2)^2+1 \equiv 5$ modulo $8$, by \eqref{eq:reposm2}, $\left(\frac{-2}{(4a+2)^2+1}\right)=-1$.
Using Proposition \ref{prop:Jacobi} $(J4)$  we have
$$
\left(\frac{2a+1}{(4a+2)^2+1}\right) \cdot \left(\frac{(4a+2)^2+1}{2a+1}\right) = (-1)^{\frac{2a}{2} \frac{(4a+2)^2}{2}} = 1.
$$
Hence,
$$
\left(\frac{2a+1}{(4a+2)^2+1}\right) = \left(\frac{(4a+2)^2+1}{2a+1}\right).
$$
As $(4a+2)^2+1=1+4(2a+1)^2$,
$$
\left(\frac{(4a+2)^2+1}{2a+1}\right)= \left(\frac{1}{2a+1}\right) = 1.
$$
Therefore, \eqref{eq:2px2p2x2ypy2pxz2_4ap2} reduces to $\left(\frac{-(4a+2)}{(4a+2)^2+1}\right)=-1 \cdot 1 =-1$, a contradiction with \eqref{eq:prop:2px2p2x2ypy2pxz2}.

Now assume that $a<0$. Then by $(J2)$,
\begin{equation}\label{eq:2px2p2x2ypy2pxz2_4ap2i}
	\left(\frac{-(4a+2)}{(4a+2)^2+1}\right)=\left(\frac{2}{16|a|^2-16|a|+5}\right) \cdot \left(\frac{2|a|-1}{16|a|^2-16|a|+5}\right).
\end{equation}
Using Proposition \ref{prop:Jacobi} $(J3)$ we obtain that $\left(\frac{2}{16|a|^2-16|a|+5}\right)=-1$.
Using Proposition \ref{prop:Jacobi} $(J4)$  we have
$$
\left(\frac{2|a|-1}{16|a|^2-16|a|+5}\right) \cdot \left(\frac{16|a|^2-16|a|+5}{2|a|-1}\right) = (-1)^{\frac{2|a|-2}{2} \frac{16|a|^2-16|a|+4}{2}} = 1.
$$
Hence,
$$
\left(\frac{2|a|-1}{16|a|^2-16|a|+5}\right) = \left(\frac{16|a|^2-16|a|+5}{2|a|-1}\right).
$$
As $16|a|^2-16|a|+5=1+4(2|a|-1)^2$,
$$
\left(\frac{16|a|^2-16|a|+5}{2|a|-1}\right)= \left(\frac{1}{2|a|-1}\right) = 1.
$$
Therefore, \eqref{eq:2px2p2x2ypy2pxz2_4ap2i} reduces to $\left(\frac{-(4a+2)}{(4a+2)^2+1}\right)=-1 \cdot 1 =-1$, a contradiction with \eqref{eq:prop:2px2p2x2ypy2pxz2}.

Now consider the case $x=8a+5$ for some integer $a$.
In this case, 
\eqref{eq:jacobidpi1}
implies that
\begin{equation}\label{eq:prop:2px2p2x2ypy2pxz2a}
	\left(\frac{-4(8a+5)}{((8a+5)^2+1)/2}\right)=\left(\frac{-(8a+5)}{32a^2+40a+13}\right)=1.
\end{equation}
Let us compute this Jacobi symbol. First assume that $a\geq 0$. Then
\begin{equation}\label{eq:2px2p2x2ypy2pxz2_8ap5}
	\left(\frac{-(8a+5)}{32a^2+40a+13}\right)=\left(\frac{-1}{32a^2+40a+13}\right) \cdot \left(\frac{8a+5}{32a^2+40a+13}\right).
\end{equation}
As $32a^2+40a+13 \equiv 1$ modulo $4$, Proposition \ref{prop:Jacobi} $(J3)$ shows that $\left(\frac{-1}{32a^2+40a+13}\right)=1$.
Using Proposition \ref{prop:Jacobi} $(J4)$ we have
$$
\left(\frac{8a+5}{32a^2+40a+13}\right) \cdot \left(\frac{32a^2+40a+13}{8a+5}\right) = (-1)^{\frac{8a+4}{2} \frac{32a^2+40a+12}{2}} = 1.
$$
Hence,
$$
\left(\frac{8a+5}{32a^2+40a+13}\right) = \left(\frac{32a^2+40a+13}{8a+5}\right).
$$
As $32a^2+40a+13=(4a+2)(8a+5)+4a+3$, then $\left(\frac{32a^2+40a+13}{8a+5}\right)=\left(\frac{4a+3}{8a+5}\right)$. Then by Proposition \ref{prop:Jacobi} $(J4)$,
$$
\left(\frac{4a+3}{8a+5}\right) \cdot \left(\frac{8a+5}{4a+3}\right) = (-1)^{\frac{4a+2}{2} \frac{8a+4}{2}} = 1.
$$
So,
$$
\left(\frac{4a+3}{8a+5}\right) = \left(\frac{8a+5}{4a+3}\right),
$$
and because $8a+5=2(4a+3)-1$,
$$
\left(\frac{8a+5}{4a+3}\right)= \left(\frac{-1}{4a+3}\right) = -1
$$
by Proposition \ref{prop:Jacobi} $(J3)$. Therefore, \eqref{eq:2px2p2x2ypy2pxz2_8ap5} reduces to $\left(\frac{-(8a+5)}{32a^2+40a+13}\right)=1 \cdot -1 = -1$, a contradiction with \eqref{eq:prop:2px2p2x2ypy2pxz2a}.

Now assume that $a<0$. 
Then
\begin{equation}\label{eq:2px2p2x2ypy2pxz2_8ap5i}
	\left(\frac{-(8a+5)}{32a^2+40a+13}\right)= \left(\frac{8|a|-5}{32|a|^2-40|a|+13}\right).
\end{equation}
Using Proposition \ref{prop:Jacobi} $(J4)$ we have
$$
\left(\frac{8|a|-5}{32|a|^2-40|a|+13}\right) \cdot \left(\frac{32|a|^2-40|a|+13}{8|a|-5}\right) = (-1)^{\frac{8|a|-6}{2} \frac{32|a|^2-40|a|+12}{2}} = 1.
$$
Hence,
\begin{equation}\label{eq:2px2p2x2ypy2pxz2_8ap5iii}
	\left(\frac{8|a|-5}{32|a|^2-40|a|+13}\right) = \left(\frac{32|a|^2-40|a|+13}{8|a|-5}\right).
\end{equation}
Because $32|a|^2-40|a|+13=(8|a|-5)(4|a|-3)+4|a|-2$,
\begin{equation}\label{eq:2px2p2x2ypy2pxz2_8ap5ii}
	\left(\frac{32|a|^2-40|a|+13}{8|a|-5}\right)=\left(\frac{4|a|-2}{8|a|-5}\right)=\left(\frac{2}{8|a|-5}\right) \cdot \left(\frac{2|a|-1}{8|a|-5}\right).
\end{equation}
$8|a|-5 \equiv 3$ modulo $4$, so $\left(\frac{2}{8|a|-5}\right)=-1$ by Proposition \ref{prop:Jacobi} (J3). Hence, it is left to determine $ \left(\frac{2|a|-1}{8|a|-5}\right)$.  Using Proposition \ref{prop:Jacobi} $(J4)$ we have
$$
\left(\frac{2|a|-1}{8|a|-5}\right) \cdot  \left(\frac{8|a|-5}{2|a|-1}\right) = (-1)^{\frac{2|a|-2}{2} \frac{8|a|-6}{2}} =
\begin{cases}
	1, & \text{if $a$ is odd}, \\
	-1, & \text{if $a$ is even}. \\
\end{cases}
$$
Therefore, if $a$ is odd, we have $\left(\frac{2|a|-1}{8|a|-5}\right) =  \left(\frac{8|a|-5}{2|a|-1}\right)$, and when $a$ is even, we have $\left(\frac{2|a|-1}{8|a|-5}\right)=-  \left(\frac{8|a|-5}{2|a|-1}\right)$. Because $8|a|-5=4(2|a|-1)-1$,
$$
\left(\frac{8|a|-5}{2|a|-1}\right)=\left(\frac{-1}{2|a|-1}\right)=
\begin{cases}
	1, & \text{if $a$ is odd}, \\
	-1, & \text{if $a$ is even}. \\
\end{cases}
$$
Hence, in any case, $\left(\frac{2|a|-1}{8|a|-5}\right)=1$. Then \eqref{eq:2px2p2x2ypy2pxz2_8ap5ii} implies that $\left(\frac{32|a|^2-40|a|+13}{8|a|-5}\right)=-1$, and by \eqref{eq:2px2p2x2ypy2pxz2_8ap5i} and \eqref{eq:2px2p2x2ypy2pxz2_8ap5iii}  $\left(\frac{-(8a+5)}{32a^2+40a+13}\right)=-1$, a contradiction with \eqref{eq:prop:2px2p2x2ypy2pxz2a}.

Finally, we can conclude that equation \eqref{eq:2px2p2x2ypy2pxz2} has no integer solutions.

\vspace{10pt}

We will now consider the equation
\begin{equation}\label{eq:2p2x2px2ymy2pxyzpz2}
	2+2x^2+x^2y-y^2+x y z+z^2=0.
\end{equation}
It is easy to check that this equation has no integer solutions with $|y| \leq 1$.
Consider the equation as a quadratic in $x,z$ with parameter $y$.

First assume that $y$ is odd. We may rewrite the equation as
\begin{equation}\label{eq:2p2x2px2ymy2pxyzpz2red}
	x^2(y+2)+y \cdot x z +z^2=y^2-2.
\end{equation}
The discriminant of the quadratic form on the left-hand side of \eqref{eq:2p2x2px2ymy2pxyzpz2red} is $y^2-4y-8$. Letting $P_1(y)=y^2-2$ and $P_2(y)=1$ in \eqref{eq:jacobynosol}, conditions (i), (ii)  and (iii) listed after \eqref{eq:jacobynosol} are satisfied when $y$ is odd. Hence, the Jacobi symbol $\left( \frac{y^2-4y-8}{y^2-2} \right)$ must be equal to $1$ by \eqref{eq:jacobidpi1}. Let us prove that it is equal to $-1$ for every odd $y$, which will give us the desired contradiction.
First assume that $y>1$. Then,
$$
\left( \frac{y^2-4y-8}{y^2-2} \right)=\left( \frac{-4y-6}{y^2-2} \right)=\left( \frac{-2}{y^2-2} \right) \cdot \left( \frac{2y+3}{y^2-2} \right)= - \left( \frac{y^2-2}{2y+3} \right)
$$
because $\left( \frac{2y+3}{y^2-2} \right) \cdot \left( \frac{y^2-2}{2y+3} \right)=(-1)^{\frac{y^2-3}{2}\frac{2y+2}{2}}=1$ by Proposition \ref{prop:Jacobi} (J4), and $\left( \frac{-2}{y^2-2} \right)=-1$ by \eqref{eq:reposm2} as $y^2-2=7$ modulo $8$ when $y$ is odd. Also,
$$
\left( \frac{y^2-2}{2y+3} \right) =  \left( \frac{2}{2y+3} \right)^2 \cdot \left( \frac{y^2-2}{2y+3} \right) = \left( \frac{(2y+3)(2y-3)+1}{2y+3} \right) = \left( \frac{1}{2y+3} \right) = 1.
$$
Hence,
$$
\left( \frac{y^2-4y-8}{y^2-2} \right)=- \left( \frac{y^2-2}{2y+3} \right)=-1.
$$
Now assume that $y<-1$. Then
$$
\left( \frac{y^2-4y-8}{y^2-2} \right) =\left( \frac{4|y|-6}{y^2-2} \right) 
=\left( \frac{2}{y^2-2} \right) \cdot \left( \frac{2|y|-3}{y^2-2} \right) =  \left( \frac{2|y|-3}{y^2-2} \right)
$$
because $\left( \frac{2}{y^2-2} \right)=1$ by Proposition \ref{prop:Jacobi} (J3) as $y^2-2=7$ modulo $8$ when $y$ is odd. Also,
$$
\left( \frac{2|y|-3}{y^2-2} \right) = - \left( \frac{y^2-2}{2|y|-3} \right)  
= - \left( \frac{4y^2-8}{2|y|-3} \right) = - \left( \frac{1}{2|y|-3} \right) = -1
$$
by Proposition \ref{prop:Jacobi} (J2) and (J4).

Now assume that $y$ is even. Then modulo $4$ analysis on equation \eqref{eq:2p2x2px2ymy2pxyzpz2} shows that $y$ is divisible by $4$, and $z$ is even. After making the substitutions $y=4Y$, $z=2Z$ and cancelling $2$, we obtain
$$
x^2(2Y+1)+4Y \cdot x Z +2Z^2=8Y^2-1.
$$
The discriminant of the quadratic form on the left-hand side is $8 (2Y^2- 2 Y-1)$. Letting $P_1(Y)=8Y^2-1$ and $P_2(Y)=1$ in \eqref{eq:jacobynosol}, conditions (i), (ii)  and (iii) listed after \eqref{eq:jacobynosol} are satisfied for all $Y$. Hence, the Jacobi symbol $\left( \frac{8 (2Y^2- 2 Y-1)}{8Y^2-1} \right)$ must be equal to $1$ by \eqref{eq:jacobidpi1}. Let us prove that it is equal to $-1$ for every $Y$, which will give us the desired contradiction.
First assume that $Y>0$. Then,
$$
\left( \frac{8 (2Y^2- 2 Y-1)}{8Y^2-1} \right)=\left( \frac{2}{8Y^2-1} \right) \cdot \left( \frac{8 Y^2- 8 Y-4}{8Y^2-1} \right) = \left( \frac{8 Y^2- 8 Y-4}{8Y^2-1} \right) =  \left( \frac{-8 Y-3}{8Y^2-1} \right)
$$
because $\left( \frac{2}{8Y^2-1} \right)=1$ by Proposition \ref{prop:Jacobi} (J3). Then,
$$
\left( \frac{-8 Y-3}{8Y^2-1} \right) = \left( \frac{-1}{8Y^2-1} \right) \cdot \left( \frac{8 Y+3}{8Y^2-1} \right) =  \left( \frac{8Y^2-1}{8 Y+3} \right) = \left( \frac{-3Y-1}{8 Y+3} \right) = \left( \frac{2}{8 Y+3} \right)=-1
$$
because $\left( \frac{-1}{8Y^2-1} \right)=-1$ and $\left( \frac{2}{8 Y+3} \right)=-1$ by Proposition \ref{prop:Jacobi} (J3) and $\left( \frac{8 Y+3}{8Y^2-1} \right) \cdot  \left( \frac{8Y^2-1}{8 Y+3} \right)=-1$ by Proposition \ref{prop:Jacobi} (J4).
Now assume that $Y<0$. Then,
$$
\left( \frac{8 (2Y^2- 2 Y-1)}{8Y^2-1} \right) 
=\left( \frac{2}{8Y^2-1} \right) \cdot \left( \frac{8 Y^2+ 8 |Y|-4}{8Y^2-1} \right) = \left( \frac{8 Y^2+ 8 |Y|-4}{8Y^2-1} \right)
$$
because $\left( \frac{2}{8Y^2-1} \right)=1$ by Proposition \ref{prop:Jacobi} (J3). Then,
$$
\left( \frac{8 Y^2+ 8 |Y|-4}{8Y^2-1} \right) = \left( \frac{8 |Y|-3}{8Y^2-1} \right) =  \left( \frac{8Y^2-1}{8 |Y|-3} \right) =  \left( \frac{16(8Y^2-1)}{8 |Y|-3} \right)=\left( \frac{2}{8 |Y|-3} \right)=-1
$$
because $\left( \frac{2}{8 |Y|-3} \right)=-1$ by Proposition \ref{prop:Jacobi} (J3), $\left( \frac{8 |Y|-3}{8Y^2-1} \right) \cdot  \left( \frac{8Y^2-1}{8 |Y|-3} \right)=1$ by Proposition \ref{prop:Jacobi} (J4) and $16(8Y^2-1)=(8Y-3)(16Y+6)+2$.

Finally, we can conclude that equation \eqref{eq:2p2x2px2ymy2pxyzpz2} has no integer solutions.

\subsection{Exercise 7.42}\label{ex:cubicfactors}
\textbf{\emph{Prove that the equation
	\begin{equation}\label{eq:z(y2pyp3)mx3mx2pxm1}
	z(y^2+y+3) = x^3+x^2-x+1
	\end{equation}
	has no integer solutions.}}

To solve this equation, we will follow the method in Section 7.3.3 of the book.  
This equation is linear in $z$ and its integer solvability is equivalent to the question whether a number of the form $x^3+x^2-x+1$ can have a divisor of the form $y^2+y+3$. 
Let
	$$
		 f(x)=x^3 + ax^2 + bx + c,
$$
	where $a,b,c$ are integer coefficients, be an irreducible cubic monic polynomial. 
Recall from Section 7.3.3 of the book that $P(f)$ is the set of all primes that are divisors of $f(x)$ for some integer $x$. In other words, $P(f)$ is the set of all primes $p$ for which the congruence 
	$$
		f(x) \equiv 0 \,(\,\text{mod}\, p)
	$$
	has at least one solution.

	The structure of set $P(f)$ depends on the discriminant of $f$, which is given by
	\begin{equation}\label{eq:disccubicmon}
		D(f)=a^2b^2 + 18abc - 4b^3 - 4a^3c - 27c^2.
	\end{equation}

	We next cite the following well-known result, see e.g. \cite[Theorem 1.1]{MR1994473}.
	
	\begin{proposition}\label{prop:cubic013law}[Proposition 7.35 in the book]
		Let $p>3$ be a prime which is not a divisor of $D(f)$. Then the equation $f(x)=0$ has $1$ solution modulo $p$ if $\left(\frac{D}{p}\right)=-1$ and $0$ or $3$ solutions modulo $p$ if $\left(\frac{D}{p}\right)=1$.
	\end{proposition}	

We also state an important theorem of Spearman and Williams \cite{spearman1992cubic} that characterises sets $P(f)$ for cubic monic polynomials $f$ with non-square discriminant.
\begin{theorem}\label{th:spwi}[Theorem 7.37 in the book]
		Let $a,b,c$ be integers such that $f(x)=x^3+ax^2+bx+c$ is an irreducible polynomial with nonsquare discriminant $D$. Then there exists a unique subgroup $J = J(a,b,c)$ of $H(D)$ having $h(D)/3$ elements with the following property: if $p>3$ is a prime which is not a divisor of $D$ such that $\left(\frac{D}{p}\right)=1$, then the equation $f(x)=0$ has $3$ solutions modulo $p$ if and only if $p$ is represented by one of the forms in $J(a,b,c)$.
	\end{theorem} 

We will also use the following proposition.
\begin{proposition}\label{prop:h33cubicepart2}[Proposition 7.39 in the book]
 	There is no integer $y$ such that all prime factors $p$ of $y^2+y+3$ are of the form $p=u^2+11v^2$.
 \end{proposition}

Let us now state and prove the following proposition.
\begin{proposition}\label{prop:x3px2mxp1factors}
	A prime $p$ is a factor of $f(x)=x^3+x^2-x+1$ for some integer $x$ if and only if $p=2$, or   $\left(\frac{-44}{p}\right)=-1$, or $p$ can be represented as $p=u^2+11v^2$ for some integers $u,v$.
\end{proposition}	
\begin{proof}
The discriminant \eqref{eq:disccubicmon} of $f$ is $D=-44$, its prime factors are $2$ and $11$. Let us consider cases $p=2,3,11$ separately. It is easy to see that $f(1)=2$ is divisible by $p=2$, while $f(2)=11$ is divisible by $p=11$. However, an easy case analysis shows that $f(x)$ is never divisible by $p=3$. Now assume that $p>3$ is not a divisor of $D$. If $\left(\frac{-44}{p}\right)=-1$, then $p$ is a factor of $f(x)$ for some integer $x$ by Proposition \ref{prop:cubic013law}. 

Now assume that $\left(\frac{-44}{p}\right)=1$. Then equation $f(x)=0$ has $0$ or $3$ solutions modulo $p$, and we will use Theorem \ref{th:spwi} to distinguish these cases. For $D=-44$, the group $H(-44)$ has three elements: the identity elements consisting of the forms equivalent to $\langle 1,0,11\rangle$, and two other elements consisting of the forms equivalent to $\langle 3,\pm 2,4\rangle$. Hence, $h(-44)=3$. The subgroup $J$ in Theorem \ref{th:spwi} has $h(-44)/3=1$ element, hence $J$ consists of the identity element only. Thus, by Theorem \ref{th:spwi}, the equation $f(x)=0$ has $3$ solutions modulo $p$ if and only if $p=u^2+11v^2$ for some integers $u,v$. Because $p=11$ is of this form with $(u,v)=(0,1)$, we do not need to list it separately in the statement of the proposition.
\end{proof}

We next use Proposition \ref{prop:x3px2mxp1factors} to study equation \eqref{eq:z(y2pyp3)mx3mx2pxm1}. Assume that it has an integer solution $(x,y,z)$, and let $p$ be any prime factor of $y^2+y+3$. Then \eqref{eq:z(y2pyp3)mx3mx2pxm1} implies that $p$ is also a factor of $x^3+x^2-x+1$. Because $p$ is a factor of $y^2+y+3$, we must have $p\neq 2$, and, moreover, $\left(\frac{-11}{p}\right)\neq -1$ by Proposition \ref{prop:quadform}. Because $\left(\frac{-44}{p}\right)=\left(\frac{-11}{p}\right)= 1$, Proposition \ref{prop:x3px2mxp1factors} implies that $p$ can be represented as $p=u^2+11v^2$ for some integers $u,v$. However, by Proposition \ref{prop:h33cubicepart2}, there is no integer $y$ such that all prime factors $p$ of $y^2+y+3$ are of the form $u^2+11v^2$ and so \eqref{eq:z(y2pyp3)mx3mx2pxm1} has no integer solutions.

\subsection{Exercise 7.47}\label{ex:H34cubic}
\textbf{\emph{For each of the equations listed in Table \ref{tab:H34cubic}, determine whether is has any integer solutions.}}
		\begin{center}
			\begin{tabular}{ |c|c|c|c| } 
				\hline
				$H$ & Equation & $H$ & Equation \\ 
				\hline\hline
				$34$ & $4+x^2 y+x y^2+z+x y z+z^2 = 0$ & $34$ & $6+x^2+x^3+y^2+x y z+z^2 = 0$ \\ 
				\hline
				$34$ & $2+x^2+x^2 y+2 y^2+z^2-y z^2 = 0$ & $34$ &  $6+x^2+x^3-y^2+x y z-z^2  = 0$ \\ 
				\hline	
				$34$ &  $8+x^2+y+x^2 y+y^2-y z^2 = 0$   &&   \\
				\hline		
			\end{tabular}
			\captionof{table}{\label{tab:H34cubic} Some cubic equations of size $H\leq 34$.}
		\end{center}

The first equation we will consider is 
	$$
		4+x^2 y+x y^2+z+x y z+z^2=0.
$$
The substitution $z=t-x-y$ for a new variable $t$, reduces the equation to
\begin{equation}\label{eq:4px2ypxy2pzpxyzpz2red}
4 + t + t^2 - x - 2 t x + x^2 - y - 2 t y + 2 x y + t x y + y^2=0.
\end{equation}
We can solve this equation by using the Vieta jumping method, see Section \ref{ex:H18vieta}. For this equation we have the following Vieta jumping operations
$$
\begin{aligned}
	& (x,y,t) \to (x',y,t), \quad x'=1 - 2 y + 2 t - y t - x, \\
	& (x,y,t) \to (x,y',t), \quad y'= 1 - 2 x + 2 t - x t - y, \\
	& (x,y,t) \to (x,y,t'), \quad t'= -1 + 2 x + 2 y - x y - t. 
\end{aligned}
$$
Recall that an integer solution to \eqref{eq:4px2ypxy2pzpxyzpz2red} is called minimal if none of these operations decreases $\min\{|x|,|y|,|t|\}$. 
The Mathematica command 
$$
\begin{aligned}
	{\tt N[MaxValue[Min[Abs[t], Abs[x], Abs[y]], \{4 + t + t^2 - x - 2 t x + x^2 - y - 2 t y + 2 x y + t x y + y^2 == 0, } \phantom{\}} \\ \phantom{\{} {\tt Abs[-1 + 2 x + 2 y - x y - t] \geq Abs[t], Abs[1 - 2 y + 2 t - y t - x] \geq Abs[x], } \phantom{\}} \\ \phantom{\{} {\tt Abs[1 - 2 x + 2 t - x t - y] \geq Abs[y]\}, \{t, x, y \}]]}
\end{aligned}
$$
outputs ${\tt 8.70825}$. Hence, any minimal solution to \eqref{eq:4px2ypxy2pzpxyzpz2red} must satisfy $\min\{|x|,|y|,|t|\}\leq 8$.  Therefore, we must check the cases $|x| \leq 8$,  $|y| \leq 8$ and  $|t| \leq 8$. After checking these cases, we obtain no minimal solutions, therefore the equation has no integer solutions. 

\vspace{10pt}

The next equation we will consider is
$$
	2+x^2+x^2y+2y^2+z^2-y z^2=0.
$$
After multiplication by $2$, the equation can be rearranged to
$$
	4+u^2+v^2+w^2+uvw=0,
$$
where $w=2y$, $u=-(x+z)$ and $v=z-x$. We have solved this equation in Section \ref{ex:vietasolv} and it has no integer solutions, therefore, the original equation does not have any integer solutions either.

\vspace{10pt}

The next equation we will consider is
$$
	8+x^2+y+x^2y+y^2-y z^2=0.
$$
After multiplication by $4$, the equation can be rearranged to
\begin{equation}\label{eq:8px2pypx2ypy2myz2red}
32+u^2+2uv+4uvy+v^2+4y+4y^2 = 0
\end{equation}
where $u=x+z$ and $v=x-z$. 
We can solve this equation by using the Vieta jumping method, see Section \ref{ex:H18vieta}. For this equation we have the following Vieta jumping operations
$$
\begin{aligned}
	& (u,v,y) \to (u',v,y), \quad u'=-2v-4vy-u, \\
	& (u,v,y) \to (u,v',y), \quad v'= -2u-4uy-v, \\
	& (u,v,y) \to (u,v,y'), \quad y'= -1-uv-y. 
\end{aligned}
$$
Recall that an integer solution to \eqref{eq:8px2pypx2ypy2myz2red} is called minimal if none of these operations decreases $\min\{|u|,|v|,|w|\}$. 
The Mathematica command 
$$
\begin{aligned}
& {\tt N[MaxValue[Min[Abs[u], Abs[v], Abs[y]], \{4 y^2 + 4 y + v^2 + 2 u v + u^2 + 32 + 4 u v y == 0,} \phantom{\}} \\ & \phantom{\{\}}
	{\tt Abs[-2 v - 4 v y - u] \geq Abs[u], Abs[-2 u - 4 u y - v] \geq Abs[v], 
	Abs[-1 - u v - y] \geq Abs[y]\}, \{u, v, y\}]]}
\end{aligned}
$$
outputs ${\tt 2.71619}$. Hence, any minimal solution to \eqref{eq:8px2pypx2ypy2myz2red} must satisfy $\min\{|u|,|v|,|w|\}\leq 2$.  Therefore, we must check the cases $|u| \leq 2$,  $|v| \leq 2$ and  $|y| \leq 2$. After checking these cases, we obtain no minimal solutions, therefore the equation \eqref{eq:8px2pypx2ypy2myz2red} has no integer solutions. 
		
\vspace{10pt}
		
The next equation we will consider is
\begin{equation}\label{eq:6px2px3py2pxyzpz2}
	6+x^2+x^3+y^2+xyz+z^2=0.
\end{equation}
Solving this equation modulo $8$, we obtain that $x=0,1,5$ and $yz$ has opposite parity to $x$. We can rearrange equation \eqref{eq:6px2px3py2pxyzpz2} to
\begin{equation}\label{eq:6px2px3py2pxyzpz2a}
	(x+2)(2-x+x^2+yz)=-(2+(y-z)^2),
\end{equation}
and
\begin{equation}\label{eq:6px2px3py2pxyzpz2b}
	(x-2)(6+3x+x^2+yz)=-(2(3)^2+(y+z)^2).
\end{equation}

Proposition \ref{cor:quadform} (b) implies that odd positive divisors of the right-hand sides of \eqref{eq:6px2px3py2pxyzpz2a} and \eqref{eq:6px2px3py2pxyzpz2b} are $1$ or $3$ modulo $8$. Hence, if $x$ is odd then both $|x+2|$ and $|x-2|$ must be $1$ or $3$ modulo $8$. However, if $|x|\geq 2$, then these numbers differ by $4$ and therefore cannot be both $1$ or $3$ modulo $8$, a contradiction. The case $|x|<2$ can be excluded by direct substitution.

The final case to consider is $x=0$ modulo $8$ and $yz$ is odd. Because $x=0$ is not a solution, $x+2$ and $x-2$ have the same sign. Let $e=-1$ if they are both positive, and $e=1$ if they are both negative. Then integers $e(2-x+x^2+yz)$ and $e(6+3x+x^2+yz)$ are odd and positive, hence they must be $1$ or $3$ modulo $8$. Because $x=0$ modulo $8$, we obtain that both $e(2+yz)$ and $e(6+yz)$ are $1$ or $3$ modulo $8$. However, these numbers differ by $4$ and therefore cannot be both $1$ or $3$ modulo $8$, a contradiction. 
Therefore \eqref{eq:6px2px3py2pxyzpz2} has no integer solutions. 		
	
	\vspace{10pt}
		
Let us consider the equation
\begin{equation}\label{eq:6px2px3my2pxyzmz2}
	6+x^2+x^3-y^2+xyz-z^2=0.
\end{equation}
We may rearrange \eqref{eq:6px2px3my2pxyzmz2} to
$$
(x+2)(2-x+x^2+yz)=(y+z)^2-2,
$$
and
$$
(x-2)(6+3x+x^2+yz)=(y-z)^2-2(3)^2.
$$
Modulo $9$ analysis shows that $x=2$ modulo $9$ and $(y-z)/3$ is an integer. Hence, $x=9t+2$ for some integer $t$. Substituting these into the above representations, we obtain
$$
(9t+4)(2-x+x^2+yz)=(y+z)^2-2,
$$
and
$$
t(6+3x+x^2+yz)=((y-z)/3)^2-2.
$$
By Proposition \ref{prop:z2m2div}, all odd divisors of the right-hand sides are $\pm 1$ modulo $8$. If $t$ is odd, then $t$ and $9t+4$ are odd but cannot both be $\pm 1$ modulo $8$, a contradiction. If $t$ is even, then $x=9t+2$ is even and modulo $4$ analysis on \eqref{eq:6px2px3my2pxyzmz2} shows that $y,z$ are both odd. Then integers $6+3x+x^2+yz$ and $2-x+x^2+yz$ are both odd, but their difference $4+4x$ is equal to $4$ modulo $8$, hence these integers cannot both be $\pm 1$ modulo $8$, a contradiction. Therefore equation \eqref{eq:6px2px3my2pxyzmz2} has no integer solutions. 

\subsection{Exercise 7.48}\label{ex:h34cubicopen}
	\textbf{\emph{Prove that if $(x,y,z)$ is any integer solution to 
	\begin{equation}\label{eq:h34cubicopen}
	2+2 x+x^3-y^2-x y^2+x z^2 = 0,
\end{equation} 
then $x<0$ and $x \equiv 2\, (\text{mod}\, 16)$.}}

First, modulo $16$ analysis of \eqref{eq:h34cubicopen} shows that $x$ can only be equal to $2,5,6$ or $13$ modulo $16$. In particular, $x$ is $1$ or $2$ modulo $4$, and $\pm 2$ or $5$ modulo $8$. We may represent \eqref{eq:h34cubicopen} as
$$
(x+1)(3-x+x^2-y^2+z^2)=z^2+1.
$$
Because the right-hand side is positive, the integers $x+1$ and $3-x+x^2-y^2+z^2$ must have the same sign. By Proposition \ref{prop:z2p1div}, positive divisors of $z^2+1$ must be equal to $1$ or $2$ modulo $4$. If $x>0$ then $x=1$ modulo $4$, but then $z$ must be odd, and $3-x+x^2-y^2+z^2$ is equal to $3-1+1^2-y^2+1=-y^2=0$ or $3$ modulo $4$, which is a contradiction. Hence, $x<0$.

We may also represent \eqref{eq:h34cubicopen} as
$$
x(2+x^2-y^2+z^2)=y^2-2.
$$
By Proposition \ref{prop:z2m2div}, all odd divisors of $y^2-2$ are $\pm 1$ modulo $8$.
If $x$ is odd, then $x=5$ modulo $8$, and cannot be a divisor of $y^2-2$. Hence, $x$ is even. Then $x$ is equal to $2$ or $6$ modulo $16$, that is, $x=2t$ for an integer $t$ equal to $1$ or $3$ modulo $8$. Because $t$ is a divisor of $y^2-2$, $t \equiv 1\, (\text{mod}\, 8)$, hence $x \equiv 2\, (\text{mod}\, 16)$.

\subsection{Exercise 7.50}\label{ex:H34quartic}
\textbf{\emph{For each of the equations listed in Table \ref{tab:H34quartic}, determine whether is has any integer solutions.}}

	\begin{center}
		\begin{tabular}{ |c|c|c|c| } 
			\hline
			$H$ & Equation & $H$ & Equation \\ 
			\hline\hline
			$33$ & $1+2 y+2 y^2+x^2 y z+z^2 = 0$ & $34$ & $2+2 y+y^2+x^2 y z+2 z^2 = 0$ \\ 
			\hline
			$33$ & $1+x-x^2+y^2-x^2 y^2+z+z^2 = 0$ & $34$ & $2+x^2+y^2+x^2 y z+2 z^2 = 0$ \\ 
			\hline	
			$34$ & $2+x^4+y^2+y^2 z+z^2=0$ & $34$ & $2+x^2+y^2+x^2 y^2+x z^2 = 0$ \\ 
			\hline
			$34$ & $2+x^4-y^2+y^2 z+z^2=0$ &  &  \\ 
			\hline			
		\end{tabular}
		\captionof{table}{\label{tab:H34quartic} Some quartic equations of size $33 \leq H\leq 34$.}
	\end{center} 
	
	For each equation in Table \ref{tab:H34quartic}, we either present an integer solution to the equation, or a proof that no integer solutions exist.
	For some equations, we may need to use the following theorem.
	\begin{theorem}\label{th:sos}[Theorem 5.11 in the book]
		A positive integer $n$ can be represented as a sum of two squares if and only if every prime $p$ of the form $p=4k+3$ appears in the prime factorization of $n$ with an even exponent.
	\end{theorem}
	
	The first equation we will consider is
	\begin{equation}\label{eq:1p2yp2y2px2yzpz2}
		1+2 y+2 y^2+x^2 y z+z^2 = 0.	
	\end{equation}
	Modulo $4$ analysis on this equation shows that $x$ and $z$ are both odd while $y$ is $2$ modulo $4$. After substituting $y=2Y$ for some integer $Y$ into this equation, we obtain
	$$
	1+4Y+8 Y^2+2x^2 Y z+z^2 = 0
	$$
	where $Y$ is odd. This equation is equivalent to
	$$
	(1+2Y)^2+(x^2Y+z)^2=Y^2(x^4-4),
	$$
	or after the substitution $Z=x^2Y+z$,
	$$
	(1+2Y)^2+Z^2=Y^2(x^2-2)(x^2+2).
	$$
	Modulo $4$ analysis of this equation shows that $x$ and $Y$ are odd and $Z$ is even. Hence, $x^2-2$ is  equal to $3$ modulo $4$, and it must have a prime divisor $p$ equal to $3$ modulo $4$ appearing its prime factorization with odd exponent. Because $p$ must appear in the prime factorization of the sum of squares $(1+2Y)^2+Z^2$ with even exponent by Theorem \ref{th:sos}, $p$ must also be a divisor of $x^2+2$. However, $x^2+2$ and $x^2-2$ cannot share a prime divisor $p \neq 2$, a contradiction.
	
	\vspace{10pt}
	
	The next equation we will consider is
	\begin{equation}\label{eq:1pxmx2py2mx2y2pzpz2}
		1+x-x^2+y^2-x^2 y^2+z+z^2 = 0.
	\end{equation}
	After the substitution $z=xy+t$ for a new variable $t$, we obtain
	\begin{equation}\label{eq:1pxmx2py2mx2y2pzpz2red}
		1 + t + t^2 + x - x^2 + x y + 2 t x y + y^2=0.
	\end{equation}
	We can prove this equation has no integer solutions using Vieta Jumping. Assume that $(x,y,t)$ is an integer solution to \eqref{eq:1pxmx2py2mx2y2pzpz2red}, then
	$$
	{\rm (a)} \,\, (1 + y + 2 t y - x,y,t), \quad {\rm (b)} \,\, (x,-x - 2 t x - y,t), \quad \text{and} \quad {\rm (c)} \,\, (x,y,-1 - 2 x y - t)
	$$
	are also solutions to the same equation. Recall that a minimal solution to \eqref{eq:1pxmx2py2mx2y2pzpz2red} is one where the value $|x|+|y|+|t|$ does not decrease with any of the above operations. The Mathematica command
	\begin{align*}
		& {\tt N[MaxValue[Min[Abs[t], Abs[x], Abs[y]], \{1 + t + t^2 + x - x^2 + x y + 2 t x y + y^2 == 0, \phantom{\}}} \\
		& {\tt Abs[-1 - 2 x y - t] \geq Abs[t], Abs[1 + y + 2 t y - x]\geq Abs[x], Abs[-x - 2 t x - y] \geq Abs[y]\}, \{t, x, y\}]]}
	\end{align*}
	outputs $1.21825$, hence any minimal solution must have $\min\{|x|,|y|,|t|\} \leq 1$. Direct analysis shows that no minimal solutions to \eqref{eq:1pxmx2py2mx2y2pzpz2red} exist, therefore it has no integer solutions and hence equation \eqref{eq:1pxmx2py2mx2y2pzpz2} has no integer solutions as well.
	
	\vspace{10pt}
	A computer search returns that the equation
	\begin{equation}
		2+x^4+y^2+y^2 z+z^2 = 0
	\end{equation}
	has an integer solution
	$$
	(x,y,z)=(526303, 842989, -132778118093),
	$$
	and the equation
	\begin{equation}
		2+x^4-y^2+y^2 z+z^2 = 0
	\end{equation}
	has an integer solution
	$$
	(x,y,z)=(3013, 4533, -5463342).
	$$	
	
	\vspace{10pt}
	
	\textbf{\emph{The following equations in this exercise have been solved by Denis Shatrov.}}
	
	\smallskip
	The next equation we will consider is
	\begin{equation}\label{eq:2p2ypy2px2yzp2z2}
		2+2y+y^2+x^2yz+2z^2=0.
	\end{equation}
	After making the transformation $z \to -z$ we obtain
	$$
	2+2y+y^2-x^2yz+2z^2=0
	$$
	and $y$ and $z$ must have the same sign. We can then rearrange this equation as
	$$
	(y+1)^2+1=z(x^2y-2z),
	$$
	and also as
	$$
	2(z^2+1)=y(x^2z-y-2).
	$$
	The integers $y,z,x^2y-2z$ and $x^2z-y-2$ are either all positive or all negative. If they are all positive, then none of them can be $3$ modulo $4$ by Proposition \ref{prop:z2p1div}. If they are all negative, then their absolute values are positive and cannot be $3$ modulo $4$ by Proposition \ref{prop:z2p1div}, therefore the integers themselves cannot be $1$ modulo $4$. On the other hand, modulo $8$ analysis shows that at least one of these four integers is congruent to $3$ modulo $4$ and at least one of them is congruent to $1$ modulo $4$, which is a contradiction.
	
	\vspace{10pt}
	
	The next equation we will consider is
	\begin{equation}\label{eq:2px2py2px2yzp2z2}
		2+x^2+y^2+x^2 y z+2 z^2 = 0.
	\end{equation}
	We can rearrange equation \eqref{eq:2px2py2px2yzp2z2} to
	$$
	2+y^2+2z^2=x^2(-1-yz),
	$$
	from which it follows that $-1-yz>0$, hence, in particular, $z\neq 0$.
	Then multiplying by $z^2$ and rearranging we obtain
	$$
	(z^2)^2+(z^2+1)^2=(-1-yz)(z^2x^2+yz-1).
	$$
	Modulo $16$ analysis of \eqref{eq:2px2py2px2yzp2z2} shows that $z$ must be divisible by $4$, hence  $-1-yz=3$ modulo $4$. Because $-1-yz>0$, Proposition \ref{prop:abx2y2} implies that $z^2$, $z^2+1$, $-1-yz$ and $z^2x^2+yz-1$ must have a common prime divisor $p$ equal to $3$ modulo $4$. However, $z^2$ and $z^2+1$ are coprime, which is a contradiction.
	
	\vspace{10pt}
	
	The final equation we will consider is
	\begin{equation}\label{eq:2px2py2px2y2pxz2}
		2+x^2+y^2+x^2y^2+xz^2=0.
	\end{equation}
 It is clear that $x<0$.
	After making the transformation $x \to -x$ and rearranging, we obtain
	\begin{equation}\label{eq:2px2py2px2y2pxz2'}
		(y^2+2)(x^2+1)=x(z^2+x),
	\end{equation}
	where $x>0$. Let $p$ be any odd prime divisor of $x^2+1$, then $p$ is a divisor of $x(z^2+x)$. Because $\text{gcd}(x^2+1,x)=1$, $p$ is not a divisor of $x$, hence $p$ is a divisor of $z^2+x$. But then $p$ is also a divisor of $(z^2+x)(z^2-x)+(x^2+1)=z^4+1$. Because $z^4+1=(z^2-1)^2+2z^2$, Proposition \ref{cor:quadform} (b) implies that $p$ is equal to $1$ or $3$ modulo $8$. Because $p$ is a positive divisor of $x^2+1$, by Proposition \ref{prop:z2p1div} it cannot be $3$ modulo $8$, hence it must be $1$ modulo $8$. Because this is true for every odd prime factor of $x^2+1$, this implies that $x^2+1$ can only be $1$, $2$ or $9$ modulo $16$.
	
	On the other hand, modulo $16$ analysis of \eqref{eq:2px2py2px2y2pxz2'} shows that $x=3,5,6,10,11$ or $13$ modulo $16$, but then $x^2+1$ is congruent to $5$ or $10$ modulo $16$, which is a contradiction.

\subsection{Exercise 7.55}\label{ex:H45sym}
\textbf{\emph{Find an integer solution to each of the equations listed in Table \ref{tab:H45sym}.}}
		\begin{center}
			\begin{tabular}{ |c|c|c|c| } 
				\hline
				$H$ & Equation & $H$ & Equation \\ 
				\hline\hline
				$41$ & $2(x^2+y^2+z^2)=xyz+9$ & $44$ & $x^3+y^3+z^3=2xyz+4$ \\ 
				\hline
				$41$ & $2(x^2+y^2+z^2)=xyz-9$ & $44$ & $x^2 y+y^2 z+z^2x + 3(x+y+z) = 2$ \\ 
				\hline
				$42$ & $x^3+y^3+z^3+x y+y z+z x=6$ & $45$ & $x^3+y^3+z^3+xyz=13$ \\ 
				\hline
				$43$ & $x^3+y^3+z^3=xyz+11$ & $45$ & $x^3+y^3+z^3+x y+y z+z x=9$  \\ 
				\hline
				$44$ & $x^3+y^3+z^3-(x y+y z+z x)-(x+y+z)=2$ & $45$ & $x^3+y^3+z^3=2xyz+5$  \\ 
				\hline		
			\end{tabular}
			\captionof{table}{\label{tab:H45sym} Cyclic equations of size $H\leq 45$ with large integer solutions.}
		\end{center} 

Table \ref{tab:H45sym_sol} presents an integer solution for each equation in Table \ref{tab:H45sym} found by a computer search, see the discussion in Section 7.3.5 of the book for details on how this computer search has been performed.

		\begin{center}
			\begin{tabular}{ |c|c|c|c| } 
				\hline
				 Equation &Solution $(x,y,z)$ \\ 
				\hline\hline
				 $2(x^2+y^2+z^2)=xyz+9$ & $(15,-1701,-231)$  \\ 
				\hline
				$2(x^2+y^2+z^2)=xyz-9$ & $(75,-21753,-815157)$ \\ 
				\hline
				 $x^3+y^3+z^3+x y+y z+z x=-6$ & $(-245,-245,309)$  \\ 
				\hline
				 $x^3+y^3+z^3=xyz+11$ & $(-106179597,-530403428,567068662)$  \\ 
				\hline
				 $x^3+y^3+z^3-(x y+y z+z x)-(x+y+z)=2$ & $(64,228,-230)$  \\ 
				\hline
				 $x^3+y^3+z^3=2xyz+4$ & $(17097,-35353,22718)$ \\ 
				\hline
				 $x^2 y+y^2 z+z^2x + 3(x+y+z) = 2$ & $(18,310,-5338)$ \\ 
				\hline
				 $x^3+y^3+z^3+xyz=13$ & $(2389570,5322418,-4710359)$\\ 
				\hline
				 $x^3+y^3+z^3+x y+y z+z x=9$ & $(-8360,-8426,10575)$ \\ 
				\hline
				 $x^3+y^3+z^3=2xyz+5$ & $(-41542008,-45693347,77205280)$\\ 
				\hline	
			\end{tabular}
			\captionof{table}{\label{tab:H45sym_sol} Examples of integer solutions for the equations in Table \ref{tab:H45sym}.}
		\end{center} 

In the next exercises, we will solve the following problem. 
\begin{problem}\label{prob:pos}[Problem 7.59 in the book]
	Given a Diophantine equation, determine whether it has a solution in \textbf{positive} integers.
\end{problem}

\subsection{Exercise 7.60}\label{ex:H26conssq}
\textbf{\emph{Solve Problem \ref{prob:pos} for the equations
	$$
		1 - x + x^2 y^2 - z - z^2 = 0, \quad -1 + x + y + x^2 y^2 - z^2 = 0 \quad \text{and} \quad 1 - x - y + x^2 y^2 - z^2 = 0
	$$
	of size $H=25$, as well as the equations 
	$$
	x - y + x^2 y^2 - z - z^2 = 0 \quad \text{and} \quad  x - y + x^2 y^2 + z - z^2 = 0 
	$$
	of size $H=26$.}}

Let us consider the first equation
$$
1 - x + x^2 y^2 - z - z^2 = 0.
$$
If $z \geq xy$, then
$$
1 - x + x^2 y^2 = z^2+z \geq x^2y^2 + xy, 
$$
hence $1-x \geq xy$, which is impossible for positive integers $x,y$. On the other hand, if $z \leq xy-1$, then
$$
1 - x + x^2 y^2 = z^2+z \leq x^2y^2 - xy, 
$$
hence $1-x+xy \leq 0$, which is again impossible for positive integers $x,y$.

\vspace{10pt}

The next equation we will consider is 
$$
-1+x+y+x^2y^2-z^2=0.
$$
Because $x,y,z$ are positive integers, the equation implies that $x^2y^2<z^2$, we then have $xy<z$ or $xy+1 \leq z$. Then
$$
0=-1+x+y+x^2y^2-z^2 \leq -1+x+y+x^2y^2-(xy+1)^2=-2+x+y-2xy=-(x-1)(y-1)-xy-1,
$$
which is impossible for positive integers $x,y$.

\vspace{10pt}

The next equation we will consider is 
$$
1-x-y+x^2y^2-z^2=0.
$$
Because $x,y,z$ are positive integers, the equation implies that $z^2<x^2y^2$, so we have $z<xy$ or $z \leq xy-1$. Then
$$
0=-(1-x-y+x^2y^2-z^2) \leq (xy-1)^2+x+y-1-x^2y^2=x+y-2xy,
$$
or $2xy \leq x+y$. The only positive integer solution to this inequality is $(x,y)=(1,1)$, however, in this case $z \leq xy-1=0$, a contradiction.

\vspace{10pt}

The next equation we will consider is 
$$
x-y+x^2y^2-z-z^2=0.
$$
If $z \geq xy$, then 
$$
x-y+x^2y^2=z^2+z \geq xy+x^2y^2,
$$
hence $x-y-xy \geq 0$, which is impossible for positive integers $x,y$. On the other hand, if $z \leq xy-1$, then 
$$
x-y+x^2y^2=z^2+z \leq x^2y^2-xy,
$$
hence $x-y+xy \leq 0$ which is again impossible for positive integers $x,y$.

\vspace{10pt}

The final equation we will consider is 
$$
x-y+x^2y^2+z-z^2=0.
$$
If $z \leq xy$, then
$$
x-y+x^2y^2=z^2-z \leq x^2y^2-xy,
$$
hence $x-y+xy\leq 0$, which is impossible for positive integers $x,y$. On the other hand, if $z \geq xy+1$ then 
$$
x-y+x^2y^2=z^2-z \geq x^2y^2+xy,
$$
hence $x-y-xy \geq 0$, which is again impossible for positive integers $x,y$.

\subsection{Exercise 7.61}\label{ex:H26factorcheck}
\textbf{\emph{Use Proposition \ref{prop:z2p1div} to solve Problem \ref{prob:pos} for the equations
	$$
	-1 - 4 y + x^2 y - y^2 - z^2 = 0, \quad -1 + x^2 y + 3 y^2 - z^2 = 0 \quad \text{and} \quad -1 - 2 x + x^3 + x y^2 - z^2 = 0
	$$
	of size $H=25$, as well as the equations 
	$$
	-2 + x^3 - x y^2 - 2 z^2 = 0 \quad \text{and} \quad  2 + x^3 - y^2 - x y^2 + z^2 = 0 
	$$
	of size $H=26$.}}

The first equation we will consider is
$$
-1 - 4 y + x^2 y - y^2 - z^2 = 0,
$$
which we can rewrite as
$$
y(- 4  + x^2  - y)  = z^2+1.
$$
As $y>0$ and $z^2+1>0$, $s=-4+x^2-y>0$. By Proposition \ref{prop:z2p1div}, $y$ must be either $1$ or $2$ modulo $4$. If $y$ is odd, as $x^2=0$ or $1$ modulo $4$, $s=0$ or $3$ modulo $4$ which is a contradiction with Proposition \ref{prop:z2p1div}. On the other hand, if $y$ is even, then $s=2$ or $3$ modulo $4$. However, if $s=2$ and $y$ is even then $z^2+1$ would be divisible by $4$, a contradiction, while if $s=3$ then this is a contradiction with Proposition \ref{prop:z2p1div}. 

\vspace{10pt}

The next equation we will consider is
$$
-1 + x^2 y + 3 y^2 - z^2 = 0,
$$
which we can rewrite as
$$
y(x^2  +3 y)  = z^2+1.
$$
As $y>0$ and $z^2+1>0$, $s=x^2+3y>0$. By Proposition \ref{prop:z2p1div}, $y$ must be either $1$ or $2$ modulo $4$. If $y$ is odd, as $x^2=0$ or $1$ modulo $4$, $s=0$ or $3$ modulo $4$ which is a contradiction with Proposition \ref{prop:z2p1div}. On the other hand, if $y$ is even, then $s=2$ or $3$ modulo $4$, however, if $s=2$ and $y$ is even then $z^2+1$ would be divisible by $4$, a contradiction, while if $s=3$ then this is a contradiction with Proposition \ref{prop:z2p1div}. 

\vspace{10pt}

The next equation we will consider is
$$
-1 - 2 x + x^3 + x y^2 - z^2 = 0,
$$
which we can rewrite as
$$
x(-2+x^2  + y^2)  = z^2+1.
$$
As $x>0$ and $z^2+1>0$, $s=-2+x^2+y^2>0$. By Proposition \ref{prop:z2p1div}, $x$ must be either $1$ or $2$ modulo $4$.
If $x$ is odd, as $y^2=0$ or $1$ modulo $4$, $s=0$ or $3$ modulo $4$ which is a contradiction with Proposition \ref{prop:z2p1div}. On the other hand, if $x$ is even, then $s=2$ or $3$ modulo $4$, however, if $s=2$ and $x$ is even then $z^2+1$ would be divisible by $4$, a contradiction, while if $s=3$ then is a contradiction with Proposition \ref{prop:z2p1div}. 

\vspace{10pt}

The next equation we will consider is
$$
-2 + x^3 - x y^2 - 2 z^2 = 0,
$$
which we can rewrite as
$$
x(x^2-y^2)=2(z^2+1).
$$
As $x>0$ and $2(z^2+1)>0$, $s=x^2-y^2>0$. By Proposition \ref{prop:z2p1div}, $x$ must be either $1$ or $2$ modulo $4$. If $x$ is odd, as $y^2=0$ or $1$ modulo $4$, $s=1$ or $3$ modulo $4$, however, if $s=1$ then the left-hand side is $1$ modulo $4$, whilst the right-hand side is either $0$ or $2$, a contradiction, if $s=3$ then this is a contradiction with Proposition \ref{prop:z2p1div}. On the other hand, if $x$ is even, then $s=0$ or $3$ modulo $4$, which is a contradiction with Proposition \ref{prop:z2p1div}.

\vspace{10pt}

The final equation we will consider is
$$
2 + x^3 - y^2 - x y^2 + z^2 = 0,
$$
which we can rewrite as
$$
(x+1)(y^2-x^2+x-1)=z^2+1.
$$
Modulo $4$ analysis shows that $x$ is $1$ or $2$ modulo $4$. The second case is impossible as $x+1$ would be $3$ modulo $4$, a contradiction with Proposition \ref{prop:z2p1div}. In the first case, $x+1=2$ modulo $4$, and as both the divisors $x+1$ and $y^2-x^2+x-1$ of $z^2+1$ cannot be equal to $2$ modulo $4$, we must have $y^2-x^2+x-1=1$ modulo $4$. This implies that $y^2=2$ modulo $4$, a contradiction.

\subsection{Exercise 7.64}\label{ex:abx2y2gen}
\textbf{\emph{Prove the following proposition.
\begin{proposition}\label{prop:abx2y2gen}[Proposition 7.63 in the book]
	Let $a,b,x,y$ be integers satisfying one of the following conditions
	\begin{itemize}
		\item[(i)] $a>0$, $b>0$, $ab=x^2+2y^2$, and either $a$ or $b$ is $5$ or $7$ modulo $8$;
		\item[(ii)] $ab=x^2-2y^2$, and either $a$ or $b$ is $3$ or $5$ modulo $8$;
		\item[(iii)] $a>0$, $b>0$, $ab=x^2+3y^2$ or $ab=x^2\pm xy + y^2$ and either $a$ or $b$ is $2$ modulo $3$;
		\item[(iv)] $ab=x^2-3y^2$ and either $a$ or $b$ is $5$ or $7$ modulo $12$.
	\end{itemize}
	Then all integers $a,b,x,y$ must have a common prime divisor $p$ satisfying the same congruence condition.
\end{proposition}}}

We will need to use the following proposition.
\begin{proposition}\label{prop:y2pyzpz2set}[Proposition 5.33 in the book]
	An integer $n$ can be represented as $n=y^2+yz+z^2$ for some integers $y,z$ if and only if it can be represented as $n=3Y^2+Z^2$ for some integers $Y,Z$.
\end{proposition}

\begin{proof}[Proof of Proposition \ref{prop:abx2y2gen}]
	
	We will first prove (i). By symmetry, we may assume that $a$ is equal to $5$ or $7$ modulo $8$. Then $a$ has a prime factor $p$ equal to $5$ or $7$ modulo $8$ that appears in its prime factorization with an odd exponent. Then Proposition \ref{cor:quadform} implies that $p$ is a common divisor of $x$ and $y$, while Proposition \ref{prop:evenmult} implies that $p$ appears in the prime factorization of $x^2+2y^2$ with an even exponent. This is possible only if $p$ is also a divisor of $b$.
	
	Let us next prove (ii). By symmetry, we may assume that $a$ is equal to $3$ or $5$ modulo $8$. Then $a$ has a prime factor $p$ equal to $3$ or $5$ modulo $8$ that appears in its prime factorization with an odd exponent. Then Proposition \ref{cor:quadform} implies that $p$ is a common divisor of $x$ and $y$, while Proposition \ref{prop:evenmult} implies that $p$ appears in the prime factorization of $x^2-2y^2$ with an even exponent. This is possible only if $p$ is also a divisor of $b$.
	
	Let us next prove (iii). We first use Proposition \ref{prop:y2pyzpz2set} to state that representations $n=x^2+3y^2$ and $n=x^2\pm xy+y^2$ cover the same set of integers for integers $x,y$.  By symmetry, we may assume that $a$ is equal to $2$ modulo $3$. Then $a$ has a prime factor $p$ equal to $2$ modulo $3$ that appears in its prime factorization with an odd exponent.  Then Proposition \ref{cor:quadform} implies that $p$ is a common divisor of $x$ and $y$, while Proposition \ref{prop:evenmult} implies that $p$ appears in the prime factorization of $x^2+3y^2$ with an even exponent. This is possible only if $p$ is also a divisor of $b$.
	
	Let us finally prove (iv). By symmetry, we may assume that $a$ is equal to $5$ or $7$ modulo $12$. Then $a$ has a prime factor $p$ equal to $5$ or $7$ modulo $12$ that appears in its prime factorization with an odd exponent. Then Proposition \ref{cor:quadform} implies that $p$ is a common divisor of $x$ and $y$, while Proposition \ref{prop:evenmult} implies that $p$ appears in the prime factorization of $x^2-3y^2$ with an even exponent. This is possible only if $p$ is also a divisor of $b$.
\end{proof}

\subsection{Exercise 7.65}\label{ex:H26pos}
\textbf{\emph{Use Propositions \ref{prop:abx2y2gen} and \ref{prop:abx2y2} to solve Problem \ref{prob:pos} for all equations listed in Table \ref{tab:H26pos}.}}

	\begin{center}
		\begin{tabular}{ |c|c|c|c| } 
			\hline
			$H$ & Equation & $H$ & Equation \\ 
			\hline\hline
			$23$ & $-1 - x + x^2 y - x y^2 - z^2 = 0$  & $26$ & $x^2+y+x^2 y-y^2 z+z^2=0$ \\ 
			\hline
			$25$ & $-1 + x^3 + y^2 - x y z + z^2 = 0$   & $26$ & $-x + x^2 y - y^2 - x^2 z - z^2=0$ \\ 
			\hline
			$25$ & $1 + x - y + x^2 y + y^2 - x z^2 = 0$  & $26$ & $x+x^2 y+y^2-x^2 z+z^2 = 0$ \\ 
			\hline
			$25$ & $-1 - x - y + x^2 y - y^2 + x z^2=0$ & $26$ & $-x+x^2 y-y^2+x^2 z-z^2 = 0$  \\ 
			\hline
			$25$ & $1 - x - y + x^2 y - y^2 - x z^2=0$ & $26$ & $-x - 2 x^2 + x^3 - y^2 - z^2 = 0$  \\ 
			\hline	
			$26$ & $-x + x^2 y - y^2 - 3 z^2 = 0$ & $26$ & $-2 - 2 x + x^2 + x^3 - y^2 - z^2=0$  \\ 
			\hline	
			$26$ & $-x^2 + y + x^2 y - 2 y^2 - z^2 = 0$ & $26$ & $-x - x^2 + x^3 - y^2 - y z - z^2 = 0$  \\ 
			\hline	
			$26$ & $-x^2 - y + x^2 y - y^2 - 2 z^2 = 0$  & $26$ & $-x - x^2 + x^3 - y^2 + y z - z^2 = 0$	\\
			\hline
			$26$ & $-x^2 - y + x^2 y - y^2 z - z^2 = 0$  & $26$ & $x + x^3 + y^2 - x y z + z^2 = 0$ \\
			\hline
			$26$ & $-x^2-y+x^2 y+y^2 z-z^2 = 0$ & $26$ & $x + x^3 + y^2 - x y^2 + z^2 = 0$ \\
			\hline	
			$26$ & $2 - x^2 + x^3 - x y^2 + z^2 = 0$ & &  \\
			\hline				
		\end{tabular}
		\captionof{table}{\label{tab:H26pos} Equations of size $H\leq 26$ solvable by Propositions \ref{prop:abx2y2gen} and \ref{prop:abx2y2}.}
	\end{center} 

\begin{proposition}\label{prop:abx2y2}[Proposition 7.62 in the book]
	Let $a,b$ be positive integers such that $ab=x^2+y^2$ for some integers $x,y$. If either $a$ or $b$ is equal to $3$ modulo $4$, then all integers $a,b,x,y$ must have a common prime divisor $p$ equal to $3$ modulo $4$.
\end{proposition}

The proofs that the equations $-1 - x + x^2 y - x y^2 - z^2 = 0$ and $-1 + x^3 + y^2 - x y z + z^2 = 0$ do not have positive integer solutions are in Section 7.5.2 of the book.

The first equation we will consider is
\begin{equation}\label{eq:1pxmypx2ypy2mxz2}
1 + x - y + x^2 y + y^2 - x z^2 = 0,
\end{equation}
which can be rearranged to
$$
(x+1)(z^2-xy+y-1)=y^2+z^2
$$
Because $x>0$ and $y^2+z^2>0$, we must have $z^2-xy+y-1>0$. Modulo $4$ analysis on \eqref{eq:mxm2x2px3my2mz2} shows that $x$ must be odd while $y$ and $z$ are even. Then postive integer $z^2-xy+y-1=3$ modulo $4$ and Proposition \ref{prop:abx2y2} implies that $x+1$, $z^2-xy+y-1$, $y$ and $z$ have a common prime factor $p=3$ modulo $4$. However, $p$ is then a factor of $(z^2-xy+y-1)-z^2+xy-y=1$, which is a contradiction.

\vspace{10pt}

The next equation we will consider is
\begin{equation}\label{eq:m1mxmypx2ymy2pxz2}
-1 - x - y + x^2 y - y^2 + x z^2=0.
\end{equation}
After the substitution $y \to -y$, the equation is reduced to \eqref{eq:1pxmypx2ypy2mxz2}, which has no positive integer solutions. As the solution of \eqref{eq:1pxmypx2ypy2mxz2} does not assume that $y>0$, this implies that equation \eqref{eq:m1mxmypx2ymy2pxz2} also has no solutions in positive integers. 

\vspace{10pt}

The next equation we will consider is
\begin{equation}\label{eq:1mxmypx2ymy2mxz2}
1 - x - y + x^2 y - y^2 - x z^2=0
\end{equation}
which can be rearranged to
$$
(x-1)(xy+y-z^2-1)=y^2+z^2.
$$
Because $x-1>0$ and $y^2+z^2>0$, we must have $xy+y-z^2-1>0$.
Modulo $4$ analysis on \eqref{eq:1mxmypx2ymy2mxz2} shows that $x$ is odd and while $y$ and $z$ are even. Then positive integer $xy+y-z^2-1=3$ modulo $4$ and Proposition \ref{prop:abx2y2} implies that $x-1$, $xy+y-z^2-1$, $y$ and $z$ have a common prime factor $p=3$ modulo $4$. However, $p$ is then a factor of $-z^2+xy+y-(xy+y-z^2-1)=1$, which is a contradiction.

\vspace{10pt}

The next equation we will consider is
\begin{equation}\label{eq:mxpx2ymy2m3z2}
-x + x^2 y - y^2 - 3 z^2 = 0
\end{equation}
which can be rearranged to
$$
x(-1 + x y) = y^2 + 3 z^2.
$$
Because $x>0$ and $y^2+3z^2>0$, we must have $-1+xy>0$. Modulo $9$ analysis on \eqref{eq:mxpx2ymy2m3z2} shows that both $x$ and $y$ must be divisible by $3$. Then $-1+xy=2$ modulo $3$ and Proposition \ref{prop:abx2y2gen} (iii) implies that positive integers $x$, $-1+xy$, $y$ and $z$ must have a common prime divisor $p=2$ modulo $3$. However, $x$ and $-1+xy$ are coprime, which is a contradiction.                

\vspace{10pt}

The next equation we will consider is
\begin{equation}\label{eq:mx2pypx2ym2y2mz2}
-x^2+y+x^2y-2y^2-z^2=0.
\end{equation}
Modulo $4$ analysis shows that $x$ is even and $y=0$ or $3$ modulo $4$. We can rewrite equation \eqref{eq:mx2pypx2ym2y2mz2} as
$$
y(1+x^2 -2 y) = x^2 + z^2.
$$
Because $y>0$ and $x^2+z^2>0$, we must have $1+x^2-2y>0$. If $y$ is $3$ modulo $4$, Proposition \ref{prop:abx2y2} implies that $y$, $1+x^2 -2 y$, $x$ and $z$ have a common divisor $p=3$ modulo $4$, but then $p$ is a divisor of $1$, a contradiction. 
If $y=0$ modulo $4$, write $y=t+1$ where $t$ is $3$ modulo $4$, and rewrite the equation \eqref{eq:mx2pypx2ym2y2mz2} as 
$$
t(x^2-2t-3)=z^2+1.
$$
By Proposition \ref{prop:z2p1div}, $z^2+1$ cannot have positive divisors equal to $3$ modulo $4$, a contradiction. 

\vspace{10pt}

The next equation we will consider is
$$
-x^2-y+x^2y-y^2-2z^2=0.
$$
The substitution $y \to t+1$ reduces the equation to
\begin{equation}\label{eq:mx2mypx2ymy2m2z2red}
t(x^2-t-3)=2(z^2+1).
\end{equation}
Because $t>-1$ and $2(z^2+1)>0$, we must have $x^2-t-3>0$. Modulo $8$ analysis on \eqref{eq:mx2mypx2ymy2m2z2red} shows that $x$ must be even. By Proposition \ref{prop:z2p1div}, the positive divisors of $z^2+1$ are $1,2$ modulo $4$. If $t=1$ modulo $4$ then $x^2-t-3=2$ modulo $4$ to ensure both sides of the equation have the same parity. However, $x^2$ would then equal $2$ modulo $4$, a contradiction. If $t=2$ modulo $4$, then $x^2=2,3$ modulo $4$ which is impossible.

\vspace{10pt}

The next equation we will consider is
\begin{equation}\label{eq:mx2mypx2ymy2zmz2}
-x^2-y+x^2y-y^2z-z^2=0,
\end{equation}
which can be rearranged to
$$
y(-1+x^2-y z)=x^2+z^2.
$$
Because $y>0$ and $x^2+z^2>0$, we must have $-1+x^2-yz>0$. Modulo $4$ analysis on \eqref{eq:mx2mypx2ymy2zmz2} shows that $x,y,z$ must all be even. Then $-1+x^2-y z$ is $3$ modulo $4$ and Proposition \ref{prop:abx2y2} implies that $y$, $-1+x^2-yz$, $x$ and $z$ have a common prime factor $p=3$ modulo $4$. However, $p$ is then a factor of $x^2+ z(y)-(-1+x^2-y z)=1$, which is a contradiction.

\vspace{10pt}

The next equation we will consider is
\begin{equation}\label{eq:mx2mypx2ypy2zmz2}
-x^2-y+x^2y+y^2z-z^2=0,
\end{equation}
which can be rearranged to
$$
y(-1+x^2+y z)=x^2+z^2.
$$
Because $y>0$ and $x^2+z^2>0$, we must have $-1+x^2+yz>0$. Modulo $4$ analysis on \eqref{eq:mx2mypx2ypy2zmz2} shows that $x,y,z$ must all be even. Then $-1+x^2+y z$ is $3$ modulo $4$ and Proposition \ref{prop:abx2y2} implies that $y$, $-1+x^2+yz$, $x$ and $z$ have a common prime factor $p=3$ modulo $4$. However, $p$ is then a factor of $x^2-z(y)-(-1+x^2-y z)=1$, which is a contradiction.

\vspace{10pt}

The next equation we will consider is
$$
	2-x^2+x^3-xy^2+z^2=0.
$$
Making the change of variable $x \to t-1$, the equation can be rearranged to
\begin{equation}\label{eq:2mx2px3mxy2pz2red}
t(-5+4t+y^2-t^2)=y^2+z^2.
\end{equation}
Because $y^2+z^2>0$ and $x>0$, we must have $t>0$ and $-5+4t+y^2-t^2>0$. Modulo $4$ analysis on \eqref{eq:2mx2px3mxy2pz2red} shows that $t$ is either $0$ or $3$ modulo $4$, and $y$ has the same parity as $t$. If $t$ is $3$ modulo $4$, then Proposition \ref{prop:abx2y2} implies that $t$, $-5+4t+y^2-t^2$, $y$ and $z$ have a common prime factor $p=3$ modulo $4$. However, $p$ is then a factor of $y^2-t^2+4t-(-5+4t+y^2-t^2)=5$, which is a contradiction. If $t$ is $0$ modulo $4$ then $y$ is even, hence $-5+4t+y^2-t^2=y^2-5$ is $3$ modulo $4$, and we have the same contradiction as in the previous case.

\vspace{10pt}

The next equation we will consider is
\begin{equation}\label{eq:x2pypx2ymy2zpz2}
x^2+y+x^2y-y^2z+z^2=0,
\end{equation}
which can be rearranged to
$$
y(-1-x^2+y z)=x^2+z^2.
$$
Because $y>0$ and $x^2+z^2>0$, we must have $-1-x^2+yz>0$. Modulo $4$ analysis on \eqref{eq:x2pypx2ymy2zpz2} shows that $x,y,z$ must all be even. Then $-1-x^2+y z$ is $3$ modulo $4$ and Proposition \ref{prop:abx2y2} implies that $y$, $-1-x^2+yz$, $x$ and $z$ have a common prime factor $p=3$ modulo $4$. However, $p$ is then a factor of $-x^2+z(y)-(-1-x^2+y z)=1$, which is a contradiction.

\vspace{10pt}

The next equation we will consider is
\begin{equation}\label{eq:mxpx2ymy2mx2zmz2}
-x+x^2y-y^2-x^2 z-z^2=0,
\end{equation}
which can be rearranged to
$$
x(-1+xy-x z)=y^2+z^2.
$$
Because $x>0$ and $y^2+z^2>0$, we must have $-1+xy-xz>0$. Modulo $4$ analysis on \eqref{eq:mxpx2ymy2mx2zmz2} shows that $x$ is even and $y$ and $z$ have the same parity. Then $-1+xy-x z$ is $3$ modulo $4$ and Proposition \ref{prop:abx2y2} implies that $x$, $-1+xy-xz$, $y$ and $z$ have a common prime factor $p=3$ modulo $4$. However $p$ is then a factor of $xy-xz-(-1+xy-x z)=1$, which is a contradiction.

\vspace{10pt}

The next equation we will consider is
\begin{equation}\label{eq:xpx2ypy2mx2zpz2}
x+x^2y+y^2-x^2 z+z^2=0,
\end{equation}
which can be rearranged to
$$
x(-1-xy+x z)=y^2+z^2.
$$
Because $x>0$ and $y^2+z^2>0$, we must have $-1-xy+xz>0$. Modulo $4$ analysis on \eqref{eq:xpx2ypy2mx2zpz2} shows that $x$ is even and $y$ and $z$ have the same parity. Then $-1-xy+x z$ is $3$ modulo $4$ and Proposition \ref{prop:abx2y2} implies that $x$, $-1-xy+xz$, $y$ and $z$ have a common prime factor $p=3$ modulo $4$. However, $p$ is then a factor of $-xy+xz-(-1-xy+x z)=1$, which is a contradiction.

\vspace{10pt}

The next equation we will consider is
\begin{equation}\label{eq:mxpx2ymy2px2zmz2}
-x+x^2y-y^2+x^2 z-z^2=0,
\end{equation}
which can be rearranged to
$$
x(-1+xy+x z)=y^2+z^2.
$$
Because $x>0$ and $y^2+z^2>0$, we must have $-1+xy+xz>0$. Modulo $4$ analysis on \eqref{eq:mxpx2ymy2px2zmz2} shows that $x$ is even and $y$ and $z$ have the same parity. Then $-1+xy+x z$ is $3$ modulo $4$ and Proposition \ref{prop:abx2y2} implies that $x$, $-1+xy+xz$, $y$ and $z$ have a common prime factor $p=3$ modulo $4$. However, $p$ is then a factor of $xy+xz-(-1+xy+x z)=1$, which is a contradiction.

\vspace{10pt}

The next equation we will consider is
\begin{equation}\label{eq:mxm2x2px3my2mz2}
-x-2x^2+x^3-y^2-z^2=0,
\end{equation}
which can be rearranged to
$$
x(-1-2x+x^2)=y^2+z^2.
$$
Because $x>0$ and $y^2+z^2>0$, we must have $-1-2x+x^2>0$. Modulo $8$ analysis of \eqref{eq:mxm2x2px3my2mz2} shows that $x$ is even. Then $-1-2x+x^2$ is $3$ modulo $4$ and Proposition \ref{prop:abx2y2} implies that $x$, $-1-2x+x^2$, $y$ and $z$ have a common prime factor $p=3$ modulo $4$. However, $p$ is then a factor of $x(-2+x)-(-1-2x+x^2)=1$, which is a contradiction.

\vspace{10pt}

The next equation we will consider is
\begin{equation}\label{eq:m2m2xpx2px3my2mz2}
-2-2x+x^2+x^3-y^2-z^2=0.
\end{equation}
Making the change of variable $x \to t-1$, we reduce the equation to
$$
t(-1-2t+t^2)=y^2+z^2.
$$
which is equivalent to \eqref{eq:mxm2x2px3my2mz2} with $x \to t$. As $x>0$, we must have $t>0$, hence equation \eqref{eq:m2m2xpx2px3my2mz2} has no positive integer solutions. 

\vspace{10pt}

The next equation we will consider is
\begin{equation}\label{eq:mxmx2px3my2myzmz2}
-x - x^2 +x^3 - y^2 - yz - z^2 = 0,
\end{equation}
which can be rearranged to
$$
x(-1 -x+x^2) = y^2 + yz+ z^2.
$$
Because $x>0$ and $y^2+yz+z^2>0$, $-1-x+x^2>0$. Modulo $3$ analysis on \eqref{eq:mxmx2px3my2myzmz2} shows that $x$ must be divisible by $3$. Then $-1-x+x^2=2$ modulo $3$ and Proposition \ref{prop:abx2y2gen} (iii) implies that positive integers $x$, $-1-x+x^2$, $y$ and $z$ must have a common prime divisor $p=2$ modulo $3$. However, $x$ and $-1-x+x^2$ are coprime, which is a contradiction.      

\vspace{10pt}

The next equation we will consider is
\begin{equation}\label{eq:mxmx2px3my2pyzmz2}
-x - x^2 +x^3 - y^2 + yz - z^2 = 0,
\end{equation}
After the substitution $y \to -y$, the equation is reduced to \eqref{eq:mxmx2px3my2myzmz2} which has no positive integer solutions. As the solution of \eqref{eq:mxmx2px3my2myzmz2} does not assume that $y>0$, this implies that equation \eqref{eq:mxmx2px3my2pyzmz2} also has no solutions in positive integers. 

\vspace{10pt}

The next equation we will consider is
\begin{equation}\label{eq:xpx3py2mxyzpz2}
x+x^3+y^2-xyz+z^2=0,
\end{equation}
which can be rearranged to
$$
x(-1-x^2+yz)=y^2+z^2.
$$
Because $x>0$ and $y^2+z^2>0$, we must have $-1-x^2+yz>0$. Modulo $4$ analysis on \eqref{eq:xpx3py2mxyzpz2} shows that $x,y,z$ must all be even. Then $-1-x^2+y z$ is $3$ modulo $4$ and Proposition \ref{prop:abx2y2} implies that $x$, $-1-x^2+yz$, $y$ and $z$ have a common prime factor $p=3$ modulo $4$. However, $p$ is then a factor of $-x^2+z(y)-(-1-x^2+y z)=1$, which is a contradiction.

\vspace{10pt}

The final equation we will consider is
\begin{equation}\label{eq:xpx3py2mxy2pz2}
x+x^3+y^2-xy^2+z^2=0,
\end{equation}
which can be rearranged to
$$
x(-1-x^2+y^2)=y^2+z^2.
$$
Because $x>0$ and $y^2+z^2>0$, we must have $-1-x^2+y^2>0$. Modulo $4$ analysis on \eqref{eq:xpx3py2mxy2pz2} shows that $x$ is either $0$ or $3$ modulo $4$ and $y$ has the same parity as $x$. Hence, either $x$ or $-1-x^2+y^2$ is $3$ modulo $4$ and Proposition \ref{prop:abx2y2} implies that $x$, $-1-x^2+y^2$, $y$ and $z$ have a common prime factor $p=3$ modulo $4$. However, if $p$ divides $x$, $y$ and $-1-x^2+y^2$, then $p$ divides $-1$, which is a contradiction. 

\section{Chapter 8}
\subsection{Exercise 8.4}\label{ex:l8nocases}
	\textbf{\emph{Describe all integer solutions to the equations
	$$
		\pm x + x^3 y^2 + z^2 = 0.
	$$}}
The first equation we will consider is
\begin{equation}\label{eq:xpx3y2pz2}
x+x^3y^2+z^2=0,
\end{equation}
which can be rearranged to 
$$
-x(1+x^2y^2)=z^2,
$$
as $-x$ and $1+x^2y^2>0$ are coprime, they must both be perfect squares, so let us denote $-x=s^2$ and $1+x^2y^2=t^2$, then $1+(s^2)^2y^2=t^2$ whose integer solutions are $(s,y,t)=(0,u,\pm 1)$ or $(u,0,\pm 1)$. Hence the only integer solutions to equation \eqref{eq:xpx3y2pz2} are $(x,y,z)=(0,u,0)$ and $(-u^2,0,u)$ where $u$ is an arbitrary integer.

\vspace{10pt}

The next equation we will consider is
\begin{equation}\label{eq:mxpx3y2pz2}
-x+x^3y^2+z^2=0,
\end{equation}
which can be rearranged to 
$$
-x(-1+x^2y^2)=z^2.
$$
Assuming $xy \neq 0$, then $-x$ and $-1+x^2y^2>0$ are coprime and they must both be perfect squares. Let us denote $-x=s^2$ and $-1+x^2y^2=t^2$, then $-1+(s^2)^2y^2=t^2$ whose integer solutions are $(s,y,t)=(\pm 1, \pm 1, 0)$. In the original variables, equation \eqref{eq:mxpx3y2pz2} has integer solutions $(x,y,z)=(\pm 1, \pm 1,0)$. Otherwise, if $xy =0$ then we have the integer solutions $(x,y,z)=(0,u,0)$ and $(u^2,0,u)$ for arbitrary integer $u$. To summarise, all integer solutions to equation \eqref{eq:mxpx3y2pz2} are $(\pm 1, \pm 1,0), (0,u,0)$ and $(u^2,0,u)$ for arbitrary integer $u$. 

\subsection{Exercise 8.5}\label{ex:2x2ymyz2pz}

\textbf{\emph{Solve Problem \ref{prob:large} for the equations
	$$
	1 + 2 x^2 y + z - y z^2 = 0, \quad	x + 2 x^2 y - y z^2 = 0 \quad \text{and} \quad 1 + x + 2 x^2 y - y z^2 = 0. 
	$$}}

From Section \ref{ex:genPell} we know that there is an infinite sequence of pairs of integers $(z_n,x_n)$, $n=0,1,2,\dots$ satisfying 
\begin{equation}\label{eq:zn2m2xn2m1}
	z_n^2-2x_n^2=1,
\end{equation}
as this is equation \eqref{eq:x2m2y2m1}. Of course, their absolute values can be arbitrarily large, and we will use these solutions to solve the equations in this exercise.

The first equation we will consider is
$$
1 + 2 x^2 y + z - y z^2 = 0,
$$
which can be rearranged to
$$
y(z^2-2x^2)=z+1.
$$
Then this equation has integer solutions $(x,y,z)=(x_n,z_n+1,z_n)$ with $(z_n,x_n)$, $n=0,1,2,\dots$ satisfying \eqref{eq:zn2m2xn2m1}, and their absolute values can be arbitrarily large.

\vspace{10pt}

The next equation we will consider is
$$
x + 2 x^2 y - y z^2 = 0,
$$
which can be rearranged to
$$
y(z^2-2x^2)=x.
$$
Then this equation has integer solutions $(x,y,z)=(x_n,x_n,z_n)$ with $(z_n,x_n)$, $n=0,1,2,\dots$ satisfying \eqref{eq:zn2m2xn2m1}, and their absolute values can be arbitrarily large.

\vspace{10pt}

The final equation we will consider is
$$
1 + x + 2 x^2 y - y z^2 = 0,
$$
which can be rearranged to
$$
y(z^2-2x^2)=x+1.
$$
Then this equation has integer solutions $(x,y,z)=(x_n,x_n+1,z_n)$ with $(z_n,x_n)$, $n=0,1,2,\dots$ satisfying \eqref{eq:zn2m2xn2m1}, and their absolute values can be arbitrarily large.

\subsection{Exercise 8.6}\label{ex:1pxpx2y2m2z2}

\textbf{\emph{Solve Problem \ref{prob:large} for the equations
	$$
	-1 + x^2 + x^2 y^2 - z^2 = 0, \quad -1 + x^2 y^2 + x z - z^2 = 0, \quad -1 + y + x^3 y^2 + z^2 = 0, \quad 1 + x + x^3 y^2 + z^2 = 0. 
	$$}}

The first equation we will consider is
$$
-1 + x^2 + x^2 y^2 - z^2 = 0.
$$ 
A simple computer search returns the solution $(x,y,z)=(4u^2+1,u,4u^3+3u)$, where $u$ is an arbitrary integer.

\vspace{10pt}

The next equation we will consider is
$$
-1 + x^2 y^2 + x z - z^2 = 0.
$$
A simple computer search returns the solution $(x,y,z)=(-32u^2-2,u,-(32u^3+16u^2+6u+1))$, where $u$ is an arbitrary integer.

\vspace{10pt}

The next equation we will consider is
$$
-1 + y + x^3 y^2 + z^2 = 0.
$$
Take any negative $x$ such that $|x|$ is not a perfect square. The resulting equation is quadratic in $(y,z)$ and has a solution $(y,z)=(0,1)$, hence it has infinitely many integer solutions by Proposition \ref{prop:Gaussquad}.

\vspace{10pt}

The final equation we will consider is
$$
1 + x + x^3 y^2 + z^2 = 0.
$$
Let us choose $x$ in the form $x=-u^2-1$ for non-zero integer $u$. Then equation reduces to $z^2-(u^2+1)^3y^2 = u^2$. This equation has a solution $(y,z)=(0,u)$, hence it has infinitely many integer solutions by Proposition \ref{prop:Gaussquad}.

\subsection{Exercise 8.8}\label{ex:Perron}
	\textbf{\emph{Solve Problem \ref{prob:large} for the equations
	$$
	x^3 y^2 + z + z^2 = 0, \quad y + x^3 y^2 + z^2 = 0, \quad 1 + x^3 y^2 + 2 z^2 = 0, \quad -2 + x^3 y^2 + z^2 = 0,
	$$
	$$
	2 + x^3 y^2 + z^2 = 0, \quad -1 + 2 x^3 y^2 + z^2 = 0 \quad \text{and} \quad 1 + 2 x^3 y^2 + z^2 = 0.
	$$}}
	
For each equation in this exercise we can use a combination of the following proposition and theorems and Proposition \ref{prop:Gaussquad} to solve Problem \ref{prob:large}.

\begin{proposition}\label{prop:Pellexists}[Proposition 3.6 in the book]
	Let $d$ be a positive integer that is not a perfect square. Then equation \eqref{eq:x2mdy2m1} has a solution in positive integers $(x,y)$.
\end{proposition}

\begin{theorem}\label{th:Dirichletprimes}[Theorem 8.3 in the book]
	For any coprime integers $a\neq 0$ and $b$ there are infinitely many primes $p$ equal to $b$ modulo $a$.
\end{theorem}

\begin{theorem}\label{th:Perron}[Theorem 8.7 in the book]
	Assume that $d$ is a positive integer, not a perfect square, that can be written in the form $d=p^n$ or $d=2p^n$ for an odd prime $p$ and integer $n\geq 1$. Then out of the three equations
	\begin{equation}\label{eq:Perron}
		(a) \quad x^2-dy^2=-1, \quad (b) \quad x^2-dy^2=2, \quad (c) \quad x^2-dy^2=-2
	\end{equation} 
	exactly one is solvable in integers $(x,y)$. In particular,
	\begin{itemize}
		\item[(i)] if $p\equiv 5\,(\text{mod}\,8)$, then equation \eqref{eq:Perron} (a) is solvable;
		\item[(ii)] if $p\equiv 7\,(\text{mod}\,8)$, then equation \eqref{eq:Perron} (b) is solvable;
		\item[(iii)] if $p\equiv 3\,(\text{mod}\,8)$, then equation \eqref{eq:Perron} (c) is solvable.
	\end{itemize} 
\end{theorem}

The first equation we will consider is
\begin{equation}\label{eq:x3y2pzpz2}
	x^3 y^2 + z + z^2 = 0.
\end{equation}
Let $p$ be a prime equal to $5$ modulo $8$. By Theorem \ref{th:Dirichletprimes}, there are infinitely many such primes, and we may select one as large as we like. Then Theorem \ref{th:Perron} (i) states that the equation $X^2 - p^3 Y^2 = -1$ is solvable in integers $(X,Y)$. Take $z=X^2$, then $z+1=X^2+1=p^3Y^2$, then $z^2+z=z(z+1)=p^3(XY)^2$. Hence, equation \eqref{eq:x3y2pzpz2} has a solution $(x,y,z)=(-p,XY,X^2)$. 
Then it has infinitely many solutions in $(y,z)$ with $x=-p$ by Proposition \ref{prop:Gaussquad}, and we can choose one with $|y|,|z|$ arbitrarily large.

\vspace{10pt}

The next equation we will consider is
\begin{equation}\label{eq:ypx3y2pz2}
	y + x^3 y^2 + z^2 = 0,
\end{equation}
which we can rearrange to
$$
y(1+x^3y)=-z^2.
$$
This equation is satisfied if $y=u^2$ and $1+x^3y=-v^2$, so that
\begin{equation}\label{eq:ypx3y2pz2red}
	1+x^3u^2=-v^2.
\end{equation}
Let $p$ be a prime equal to $5$ modulo $8$. By Theorem \ref{th:Dirichletprimes}, there are infinitely many such primes, and we may select one as large as we like. Then Theorem \ref{th:Perron} (i) states that the equation $X^2 - p^3 Y^2 = -1$ is solvable in integers $(X,Y)$. This means that equation \eqref{eq:ypx3y2pz2red} is solvable in integers $(u,v)$ for $x=-p$. Then \eqref{eq:ypx3y2pz2} has infinitely many such solutions by Proposition \ref{prop:Gaussquad}, and we can choose one with $|y|,|z|$ arbitrarily large.

\vspace{10pt}

The next equation we will consider is
\begin{equation}\label{eq:1px3y2p2z2}
	1 + x^3 y^2 + 2 z^2 = 0.
\end{equation}
Let $p$ be a prime equal to $3$ modulo $8$. By Theorem \ref{th:Dirichletprimes}, there are infinitely many such primes, and we may select one as large as we like. Then Theorem \ref{th:Perron} (iii) states that the equation $X^2 - 2 p^3 y^2 = -2$ is solvable in integers $(X,y)$. Obviously, $X$ is even, so let us write $X=2z$ for integer $z$. Then $4z^2- 2 p^3 y^2 + 2=0$, or $2z^2- p^3 y^2 + 1=0$. This means that equation \eqref{eq:1px3y2p2z2} is solvable in integers $(y,z)$ for $x=-p$. Then it has infinitely many such solutions by Proposition \ref{prop:Gaussquad}, and we can choose one with $|y|,|z|$ arbitrarily large.

\vspace{10pt}

The next equation we will consider is
\begin{equation}\label{eq:m2px3y2pz2}
	-2 + x^3 y^2 + z^2 = 0.
\end{equation}
Let $p$ be a prime equal to $7$ modulo $8$. By Theorem \ref{th:Dirichletprimes}, there are infinitely many such primes, and we may select one as large as we like. Then Theorem \ref{th:Perron} (ii) states that the equation $z^2 -  p^3 y^2 = 2$ is solvable in integers $(z,y)$. This means that equation \eqref{eq:m2px3y2pz2} is solvable in integers $(y,z)$ for $x=-p$. Then it has infinitely many such solutions by Proposition \ref{prop:Gaussquad}, and we can choose one with $|y|,|z|$ arbitrarily large.

\vspace{10pt}

The next equation we will consider is
\begin{equation}\label{eq:2px3y2pz2}
	2 + x^3 y^2 + z^2 = 0.
\end{equation}
Let $p$ be a prime equal to $3$ modulo $8$. By Theorem \ref{th:Dirichletprimes}, there are infinitely many such primes, and we may select one as large as we like. Then Theorem \ref{th:Perron} (iii) states that the equation $z^2 -  p^3 y^2 =-2$ is solvable in integers $(z,y)$. This means that equation \eqref{eq:2px3y2pz2} is solvable in integers $(y,z)$ for $x=-p$. Then it has infinitely many such solutions by Proposition \ref{prop:Gaussquad}, and we can choose one with $|y|,|z|$ arbitrarily large.

\vspace{10pt}

The next equation we will consider is
\begin{equation}\label{eq:m1p2x3y2pz2}
	-1 + 2 x^3 y^2 + z^2 = 0.
\end{equation}
Let $x$ be any integer, with $|x|$ as large as we like, such that $-2x^3$ is positive and not a perfect square. Then the equation
$$
z^2-(-2x^3)y^2=1
$$
has infinitely many positive integer solutions $(y,z)$ by Proposition \ref{prop:Pellexists}, and we can choose one with $|y|,|z|$ arbitrarily large.

\vspace{10pt}

The final equation we will consider is
\begin{equation}\label{eq:1p2x3y2pz2}
	1 + 2 x^3 y^2 + z^2 = 0.
\end{equation}
Let $p$ be a prime equal to $5$ modulo $8$. By Theorem \ref{th:Dirichletprimes}, there are infinitely many such primes, and we may select one as large as we like. Then Theorem \ref{th:Perron} (i) states that the equation $z^2 - 2 p^3 y^2 = -1$ is solvable in integers $(z,y)$. This means that equation \eqref{eq:1p2x3y2pz2} is solvable in integers $(y,z)$ for $x=-p$. Then it has infinitely many such solutions by Proposition \ref{prop:Gaussquad}, and we can choose one with $|y|,|z|$ arbitrarily large.

\subsection{Exercise 8.11}\label{ex:m1px3y2pzpz2}
\textbf{\emph{Solve Problem \ref{prob:large} for the equation
	$$
	-1 + x^3 y^2 + z + z^2 = 0.
	$$}}

To solve Problem \ref{prob:large} for the equation
\begin{equation}\label{eq:m1px3y2pzpz2}
	-1 + x^3 y^2 + z + z^2 = 0,
\end{equation}
we can begin by multiplying by $4$ and rearranging to obtain
$$
(2z+1)^2-(-x)^3(2y)^2=5.
$$
Now, the problem is reduced to studying the solvability of equation
$$
X^2-dY^2=5
$$
in integers $X,Y$ such that $X$ is odd, $Y$ is even and $d$ is a perfect cube.

Let $m\geq 0$ be integer and let $d=11^{2m+1}$. Let us prove, by induction in $m$, that the equation
$$
x^2-dy^2=1
$$
has a solution $(x,y)=(A,B)$ such that $B$ is even and coprime with $11$.

For $m=0$, we have $(A,B)=(199,60)$. Assume that it is true for $m$, that is, $A^2-dB^2=1$ for $d=11^{2m+1}$, and integers $A,B$ such that $B$ is even and coprime with $11$. Then let
$$
A'+B'\sqrt{d}=(A+B\sqrt{d})^{11},
$$
that is
$$
B' = 11 A^{10} B + 165 A^8 d B^3 + 462 A^6 d^2 B^5 + 330 A^4 B^7 d^3+ 55 A^2 d^4 B^9 + d^5 B^{11}.
$$
Because $B$ is even, $B'$ is also even, and because $A^{10}B$ is coprime with $11$, integer $B_1 = B'/11$ is not divisible by $11$. Now,
\begin{align*}
	(A')^2-121d(B_1)^2 &= (A'-11\sqrt{d}B_1)(A'+11\sqrt{d}B_1) = (A'-\sqrt{d}B')(A'+\sqrt{d}B') \\
	&= (A-B\sqrt{d})^{11} (A+B\sqrt{d})^{11} = (A^2-dB^2)^{11} = 1^{11} = 1.
\end{align*}
Hence, $(A',B_1)$ is a solution to $x^2-121dy^2=1$ such that $B_1$ is even and coprime with $121$, which completes the proof by induction.

Now let us prove, by induction in $m$, that for $d=11^{2m+1}$ the equation
$$
x^2-dy^2=5
$$
has a solution $(x,y)=(a,b)$ such that $a$ is odd and $b$ is even. For $m=0$, we have that the statement is true with $(a,b)=(7,2)$. Now, assume that it is true for $m$. Consider
$$
a_k + b_k\sqrt{d} = (a+b\sqrt{d})(A+B\sqrt{d})^k,
$$
where $(A,B)$ is a solution to $x^2-dy^2=1$ such that $B$ is even and coprime with $11$.\linebreak

Then we have $a_k^2-db_k^2=5$ for all $k$. It can be easily seen by induction in $k$ that $b_k$ is even for all $k$.
Now,
$$
b_k = bA^k + a k A^{k-1}B + \cdots = A^{k-1}(bA+kaB) + \cdots
$$
where the terms $\cdots$ are divisible by $d$. Because $aB$ is coprime with $11$, we may select $k$ such that $bA+kaB$ is divisible by $11$. Then $b_k$ is divisible by $11$ and $(a_k,b_k/11)$ solves $x^2-121dy^2=5$, which completes the proof by induction.

Finally, let $x$ be any odd power of $11$. Then $x^3$ is again an odd power of $11$, and we have proved that there exist integers $z,y$ such that
$$
(2z+1)^2 - (-x)^3 (2y)^2 = 5,
$$
or equivalently $-1 + x^3 y^2 + z + z^2=0$.

\subsection{Exercise 8.13}\label{ex:l9vieta}
\textbf{\emph{List all integer solutions to each of the equations
	$$
	2 - x^2 + x^2 y^2 - z^2=0, \quad 1 - x^2 + x^2 y^2 + z - z^2=0, \quad -1 + x^2 + x^2 y^2 + z - z^2=0.
	$$}}

The substitution $z=xy+u$ reduces each equation in this exercise to one that is solvable by Vieta jumping, see Section \ref{ex:H18vieta}. We now provide details to list the full set of integer solutions to the original equations.

Equation
\begin{equation}\label{eq:2mx2px2y2mz2}
2-x^2+ x^2 y^2  - z^2 = 0 
\end{equation}
has only finitely many integer solutions, and they are easy to find. Let us make the substitution $z=xy+u$ for a new variable $u$, which reduces the equation to $-2 + u^2 + x^2 +2 u x y=0$. This equation is solvable by Vieta jumping. The optimization problem \eqref{eq:vietaopt} takes the form
$$
	\begin{aligned}
		& \max_{(x,y,u)\in {\mathbb R}^3} \min\{|x|,|y|,|u|\}, \\
		& \text{subject to} \quad -2 + u^2 + x^2 +2 u x y=0, \quad |x| \leq |-2yu-x|, \quad |u| \leq |-2xy-u|,
	\end{aligned}
$$
and its optimal value is $0.754878...<1$. Hence, any minimal integer solution must satisfy $\min\{|x|,|y|,|u|\}=0$. The case $x=0$ results in $u^2-2=0$, which is impossible in integers. Similarly, $u=0$ results in $x^2-2=0$, again a contradiction. Hence, we must have $y=0$. Because the Vieta jumps do not change $y$, we must have $y=0$ for all integer solutions, not only the minimal ones. Substituting $y=0$ into the original equation \eqref{eq:2mx2px2y2mz2}, we obtain $2-x^2  - z^2 = 0 $, which is trivial to solve, and its integer solution $(x,z)=(\pm 1,\pm 1)$ is listed in Table \ref{tab:H12boundedrealsol}. From which we easily deduce the final answer: $(x,y,z)=(\pm1,0,\pm 1)$.

\vspace{10pt}

Equation
\begin{equation}\label{eq:1mx2px2y2pzmz2}
1 - x^2 + x^2 y^2 + z - z^2= 0 
\end{equation}
has only finitely many integer solutions, and they are easy to find.  Let us make the substitution $z=xy+u$ for a new variable $u$, which reduces the equation to $1 - u^2 - x^2 -2 u x y+xy+u=0$. This equation is solvable by Vieta jumping. The optimization problem \eqref{eq:vietaopt} takes the form
$$
	\begin{aligned}
		& \max_{(x,y,u)\in {\mathbb R}^3} \min\{|x|,|y|,|u|\}, \\
		& \text{subject to} \quad 1 - u^2 - x^2 -2 u x y+xy+u=0, \quad |x| \leq |-2yu+y-x|, \quad |u| \leq |-2xy+1-u|,
	\end{aligned}
$$
and its optimal value is $0.829484<1$.  Hence, any minimal integer solution must satisfy $\min\{|x|,|y|,|u|\}=0$. The case $x=0$ results in $ 1 - u^2 +u=0$, which is impossible in integers. Similarly, $u=0$ results in $1 - x^2 +xy=0$, which has integer solutions $(x,y)=(\pm 1,0)$. Let us now consider the case when $y=0$. Because the Vieta jumps do not change $y$, we must have $y=0$ for all integer solutions, not only the minimal ones. Substituting $y=0$ into the original equation \eqref{eq:1mx2px2y2pzmz2}, we obtain $1 - x^2 + z - z^2= 0$, which is equivalent to an equation in Table \ref{tab:H12boundedrealsol}, and its integer solutions are $(x,z)=(\pm 1, 1),(\pm 1,0)$. From which we easily deduce the final answer: $(x,y,z)=(\pm 1,0, 1), (\pm 1,0,0)$.

\vspace{10pt}

Finally, the equation
\begin{equation}\label{eq:1px2px2y2pzmz2}
-1 + x^2 + x^2 y^2 + z - z^2= 0 
\end{equation}
has only finitely many integer solutions, and they are easy to find.  Let us make the substitution $z=xy+u$ for a new variable $u$, which reduces the equation to $-1 - u^2 + x^2 -2 u x y+xy+u=0$. This equation is solvable by Vieta jumping. The optimization problem \eqref{eq:vietaopt} takes the form
$$
	\begin{aligned}
		& \max_{(x,y,u)\in {\mathbb R}^3} \min\{|x|,|y|,|u|\}, \\
		& \text{subject to} \quad -1 - u^2 + x^2 -2 u x y+xy+u=0, \quad |x| \leq |2yu-y-x|, \quad |u| \leq |-2xy+1-u|,
	\end{aligned}
$$
and its optimal value is $0.738984...<1$. Hence, any minimal integer solution must satisfy $\min\{|x|,|y|,|u|\}=0$. The case $x=0$ results in $-1 - u^2 +u=0$, which is impossible in integers. Similarly, $u=0$ results in $-1  + x^2 +xy=0$, which has integer solutions $(x,y)=(\pm 1,0)$. Let us now consider the case when $y=0$. Because the Vieta jumps do not change $y$, we must have $y=0$ for all integer solutions, not only the minimal ones. Substituting $y=0$ into the original equation \eqref{eq:1px2px2y2pzmz2}, we obtain $-1 + x^2  + z - z^2= 0$, which is equivalent to equation $x^2+x-y^2+1=0$ from Table \ref{exercise1.12}, and its integer solutions are $(x,z)=(\pm 1, 1),(\pm 1,0)$. From which we easily deduce the final answer: $(x,y,z)=(\pm1,0, 1), (\pm 1,0,0)$.

\subsection{Exercise 8.14}\label{ex:l9dommon}
	\textbf{\emph{\begin{itemize}
		\item[(a)] List all integer solutions of equation 
		\begin{equation}\label{eq:1px3yp2zpxyz}
			1 + x^3 y + 2 z + x y z = 0. 
		\end{equation}
		\item[(b)] Solve the similar equation
\begin{equation}\label{eq:m1px3yp2zpxyz}
	-1 + x^3 y + 2 z + x y z = 0. 
\end{equation}
	\end{itemize}}}
\subsubsection{Exercise 8.14 (a)}
The case $xyz=0$ can be analysed easily, and we obtain solutions $(x,y,z)=\pm(1,-1,0)$. Now let us assume that $xyz\neq 0$. Then \eqref{eq:1px3yp2zpxyz} implies that $2z+1$ is divisible by $xy$. After substituting $2z+1=txy$ into \eqref{eq:1px3yp2zpxyz} and cancelling $xy$, it reduces to $-1 + 2 t + 2 x^2 + t x y=0$. This implies that $2t-1$ is divisible by $x$. After further substitution $2t-1=kx$ and cancelling $x$ the equation reduces to $2 k + 4 x + y + k x y = 0$, which is an equation with a dominating monomial. The conditions $2z+1=txy$ and $2t-1=kx$ imply that $k,x,y$ are odd. The optimisation command
$$
{\tt MaxValue[Min[Abs[x], Abs[y], Abs[k]], 
 2 k + 4 x + y + k x y == 0,\{k, x, y\}, Reals]}
$$
outputs
$\sqrt{5}$. Hence, we must search for integer solutions with $\min\{|x|,|y|,|k|\}\leq 2$. Excluding the solutions with $kxy$ even, these solutions are
$$
(x,y,k)=\pm(-1,1,3),\pm(1,-3,1), \pm(3,-1,11), \pm(3,5,-1), \pm(5,1,-3),\pm(9,-1,5).
$$
We can then use the transformation $z=\frac{txy-1}{2}=\frac{-1+xy(kx+1)/2}{2}=\frac{kx^2y+xy-2}{4}$, and conclude that the integer solutions to equation \eqref{eq:1px3yp2zpxyz} are  
$$
\begin{aligned}
(x,y,z)= & \,\, (-9,1,-104),(-5,-1,-18),(-3,-5,-8),(-3,1,-26),(-1,3,-2),\\
& \,\, (1,-3,-2), \pm(1,-1,0), (3,-1,-26),(3,5,-8), (5,1,-18), (9,-1,-104).
\end{aligned}
$$

\subsubsection{Exercise 8.14 (b)}
The case $xyz=0$ can be analysed easily, and we obtain solutions $(x,y,z)=\pm(1,1,0)$. Now let us assume that $xyz\neq 0$.  Then \eqref{eq:m1px3yp2zpxyz} implies that $2z-1$ is divisible by $xy$. Substituting $2z-1=txy$ into \eqref{eq:m1px3yp2zpxyz} and cancelling $xy$, it reduces to $1+2t+2x^2+t x y=0$. This implies that $2t+1$ is divisible by $x$. After the further substitution $2t+1=kx$ and cancelling $x$, the equation is reduced to $2k+4x-y+k x y=0$, which is an equation with a dominating monomial. The conditions $2z-1=txy$ and $2t+1=kx$ imply that $k,x,y$ are odd. The optimisation command
$$
{\tt MaxValue[Min[Abs[x], Abs[y], Abs[k]], 
 2 k + 4 x - y + k x y == 0,\{k, x, y\}, Reals]}
$$
outputs $\sqrt{7}$.  Hence, we must search for integer solutions with $\min\{|x|,|y|,|k|\}\leq 2$. Excluding the solutions with $kxy$ even, these solutions are
$$
\begin{aligned}
(x,y,k)=&\pm(-1,1,5),\pm(1,-5,3),\pm(1,-3,7),\pm(1,1,-1), \pm(3,-7,1),\pm(3,-1,13), \\ & \pm(5,-1,7),\pm(5,3,-1),\pm(7,-5,1),\pm(7,1,-3),\pm(11,-1,5).
\end{aligned}
$$
We can then use the transformation $z=\frac{txy+1}{2}=\frac{1+xy(kx-1)/2}{2}=\frac{kx^2y-xy+2}{4}$, and conclude that the integer solutions to equation \eqref{eq:m1px3yp2zpxyz} are
$$
\begin{aligned}
(x,y,z)=& (-11,1,-148),(-7,-1,-38),(-7,5,-52),(-5,-3,-22),(-5,1,-42),\\ & (-3,1,-28),  (-3,7,-10),(-1,1,2),(-1,3,-4),(-1,5,-2),(1,-5,-2),\\ & (1,-3,-4),(1,-1,2),  \pm(1,1,0),(3,-7,-10),(3,-1,-28),(5,-1,-42), \\ & (5,3,-22),  (7,-5,-52),(7,1,-38),(11,-1,-148).
\end{aligned}
$$

\subsection{Exercise 8.15}\label{ex:l9infinite}
\textbf{\emph{Prove that the equation
	\begin{equation}\label{eq:1px2y2pz2pxz2}
		1 + x^2 y^2 + z^2 + x z^2 = 0
	\end{equation}
	has infinitely many integer solutions.}}

This equation is quadratic in $(y,z)$ for any fixed $x$. If $|x|>1$, then Proposition \ref{prop:Gaussquad} is applicable, hence it is enough to find one integer solution to conclude there are infinitely many of them. A computer search returns a solution with $x=-314$. The solution is 
$$
\begin{aligned}
x =&  -314, \\
y =& 
9379720965297275947816439847522981744539618335794117144937168010154158\\ &
8014754381217479016004010778515334476645458604244904428420004895696073\\ &
5119805509807876173985241036970594595328291615091717930653678982215761\\ &
4283604197441283424712641075801866574960778171164428037268766905909928\\ &
597053037242648372293936810300569052485431088, \\
z = &
1664743769486740643084686001668270915842497122740295450339914593767148\\ &
2435865518985132840500878237446705397349753405038393608727035061347025\\ &
9116206338537063348664390026231305344194514585842501822147048866727058\\ &
2058117187634122975105587594739057145091319956080243430837190155953197\\ &
29189493541695069608566023203361773805777308285.
\end{aligned}
$$
Hence, Proposition \ref{prop:Gaussquad} implies that equation \eqref{eq:1px2y2pz2pxz2} has infinitely many integer solutions. 	

\subsection{Exercise 8.16}\label{ex:l9sos}	
\textbf{\emph{Prove that the equation
	\begin{equation}\label{eq:2y2pz2mx4p1}
	2y^2+z^2 = x^4-1
	\end{equation}
	has infinitely many integer solutions.}}

To prove that the equation \eqref{eq:2y2pz2mx4p1} has infinitely many integer solutions, equivalently, we can prove that there are infinitely many integers $x$ such that $x^4-1 \in S_2$, where $S_2$ is defined in \eqref{S2_notation}.
We have
$$
8(x^4-1)=2 (2 x^2 + k)^2 + ( -8 k x^2- 2 k^2-8),
$$
where $k$ is an arbitrary integer. 
We must now find an integer $k$ such that $-8 k x^2- 2 k^2-8$ is a perfect square infinitely often. A computer search shows that a satisfactory $k$ is $k=-82$, as the equation $ -8 (-82) x^2- 2 (-82)^2-8=656 x^2-13456 =t^2$ has an integer solution $(x_0,t_0)=(145,3712)$ and therefore by Proposition \ref{prop:Gaussquad2} there exists an infinite set ${\cal X}$ such that $656 x^2-13456$ is a perfect square for every $x \in {\cal X}$. Hence, $8(x^4-1) \in S_2$ for every $x \in {\cal X}$. 
Because the only common prime factor of integers $8$ and $x^4-1$ can be $2$, Proposition \ref{prop:2y2pz2product} implies that $x^4-1 \in S_2$ for every $x \in {\cal X}$. Hence, equation \eqref{eq:2y2pz2mx4p1} has infinitely many integer solutions.

\subsection{Exercise 8.22}\label{ex:l12cubic}
	\textbf{\emph{For each of the equations listed in Table \ref{tab:l12cubic}, determine whether it has any integer solutions.}}
		\begin{center}
		\begin{tabular}{ |c|c|c|c| }
	\hline
	$l$ & Equation & $l$ & Equation \\
	\hline\hline
	$11$ & $1+x^2 y-4 y^2-4 z^2=0$ & $\approx 11.7$ & $1-2 x^2-2 y^2+7 x y z=0$ \\
	\hline
	$11$ & $1+2 x^2 y+y^2 z+4 z^2=0$ & $12$ & $1-16 x^2+x y^2+2 z^2=0$ \\
	\hline
	$11$ & $2+x^2 y-8 y^2-z^2=0$ & $12$ & $1+2 x^2 y-16 y^2+z^2=0$ \\
	\hline
	$11$ & $2+x^2 y-8 y^2+z^2=0$ & $12$ & $1+x y+2 x^2 y+y^2+2 x z^2=0$ \\
	\hline
	$\approx 11.6$ & $1-y+x^2 y-6 y^2-2 z^2=0$ & $12$ & $2+8 x^2 y+x y^2+z^2=0$ \\
	\hline
	$\approx 11.6$ & $1+2 x^2 y-4 y^2-3 z^2=0$ & $12$ & $1+8 x^2 y+y^2+2 x z^2=0$ \\
	\hline
	$\approx 11.6$ & $1+3 x^2 y+y^2 z+4 z^2=0$ & $12$ & $1+2x+x^2y+4y^2+2z^2=0$ \\
	\hline
	$\approx 11.7$ & $1+2 x^2 y+13 y z+z^2=0$ & $12$ & $1+4x^3+xy^2+2yz^2=0$ \\
	\hline
\end{tabular}
		\captionof{table}{\label{tab:l12cubic} Some cubic equations of length $11\leq l\leq 12$.}
	\end{center} 

We first remark that formulas \eqref{leg:p3} and \eqref{leg:m3} for the Legendre symbols remain correct for the Jacobi symbols: for any positive odd integer $n$, we have
\begin{equation}\label{eq:repos3Jacobi}
	\begin{aligned}
		&{\rm (a)} \quad 	\left(\frac{3}{n}\right) = 	\begin{cases}
			1, & \text{if $n \equiv \pm 1\, (\text{mod}\, 12)$,}\\
			-1, & \text{if $n \equiv \pm 5\, (\text{mod}\, 12)$,}
		\end{cases}
		\quad \text{and}\\	
		&{\rm (b)} \quad	
		\left(\frac{-3}{n}\right) =
		\begin{cases}
			1, & \text{if $n \equiv 1\, (\text{mod}\, 3)$,}\\
			-1, & \text{if $n \equiv -1\, (\text{mod}\, 3)$.}
		\end{cases}
	\end{aligned}
\end{equation}

The first equation we will consider is
\begin{equation}\label{eq:1px2ym4y2m4z2}
	1+x^2 y-4 y^2-4 z^2=0,
\end{equation}
which can be rewritten as
\begin{equation}\label{eq:1px2ym4y2m4z2red}
	yx^2-4z^2=(2y+1)(2y-1).
\end{equation}
Solving the equation modulo $4$, we obtain that $y=3$ modulo $4$.
It is clear that this equation has no solutions with $y\leq 0$. Let us treat $y>0$ as a parameter. The discriminant of the quadratic form on the left-hand side of \eqref{eq:1px2ym4y2m4z2red} is $16y$. Note that $16y$ and $2y+1$ are coprime, and the odd integers $2y+1$ and $2y-1$ are also coprime. Hence, as explained in Section \ref{ex:h33formeropen}, the equation may only have integer solutions if $\left(\frac{16y}{2y+1}\right)=1$, see the discussion below \eqref{eq:jacobynosol}. On the other hand,
$$
\left(\frac{16y}{2y+1}\right)= \left(\frac{8(2y+1)-8}{2y+1}\right)=\left(\frac{-8}{2y+1}\right)= \left(\frac{-2}{2y+1}\right)=-1, 
$$
where the last equality follows from \eqref{eq:reposm2} and the fact that $y=3$ modulo $4$.

\vspace{10pt}

The next equation we will consider is
\begin{equation}\label{eq:1p2x2ypy2zp4z2}
1+2 x^2 y+y^2 z+4 z^2=0.
\end{equation}
We can easily see that $y$ is odd. After multiplying the equation by $16$, we can rewrite it as
$$
32y x^2 + (y^2 + 8z)^2 = y^4 - 16 =(y^2 - 4)(y^2 + 4),
$$
or after the change of variables $Z=8z+y^2$ and $X=4x$, the equation is reduced to
$$
2yX^2+Z^2=(y^2-4)(y^2+4).
$$
Note that the odd integers $y^2-4$ and $y^2+4$ are coprime. Hence, by \eqref{eq:jacobidpi1}, we must have $\left(\frac{-8y}{y^2+4}\right)=1$. However,
$$
\left(\frac{-8y}{y^2+4}\right) = \left(\frac{-2e}{y^2+4}\right) \left(\frac{|y|}{y^2+4}\right) = -\left(\frac{|y|}{y^2+4}\right),
$$
where $e=y/|y|$, and we have used that $\left(\frac{2}{y^2+4}\right)=\left(\frac{-2}{y^2+4}\right)=-1$ by Proposition \ref{prop:Jacobi} (J3) and because $y^2+4=5$ modulo $8$. Also, by Proposition \ref{prop:Jacobi} (J4),
$$
\left(\frac{|y|}{y^2+4}\right) \cdot \left(\frac{y^2+4}{|y|}\right) = (-1)^{\frac{y^2+4-1}{2} \frac{|y|-1}{2}} = 1
$$
because $y^2+3$ is $4$ modulo $8$. Hence, 
$$
\left(\frac{-8y}{y^2+4}\right) =  -\left(\frac{|y|}{y^2+4}\right) =-\left(\frac{y^2+4}{|y|}\right)=-\left(\frac{4}{|y|}\right)= -\left(\frac{2}{|y|}\right) ^2 =-1,
$$
a contradiction.
Thus, equation \eqref{eq:1p2x2ypy2zp4z2} has no integer solutions.

\vspace{10pt}

The next equation we will consider is
\begin{equation}\label{eq:2px2ym8y2mz2}
2+x^2 y-8 y^2-z^2=0.
\end{equation}
It is clear that $y>0$. Modulo $4$ analysis shows that $y=2$ or $3$ modulo $4$. We can rewrite the equation as
$$
yx^2-z^2=8y^2-2=2(2y-1)(2y+1).
$$
This is an equation of the form \eqref{eq:jacobynosol} and has $D(y)=4y$. Note that $4y$ and $2y+1$ are coprime, and the odd integers $2y-1$ and $2y+1$ are also coprime. Hence, by \eqref{eq:jacobidpi1}, we must have $\left(\frac{4y}{2y+1}\right)=1$. On the other hand,
$$
\left(\frac{4y}{2y+1}\right)=\left(\frac{-2y}{2y+1}\right) \cdot \left(\frac{-2}{2y+1}\right) = \left(\frac{1}{2y+1}\right) \left(\frac{-2}{2y+1}\right) = -1,
$$
where the last equality follows from \eqref{eq:reposm2} and the fact that $y=2$ or $3$ modulo $4$. 
This is a contradiction.

\vspace{10pt}

The next equation we will consider is
\begin{equation}\label{eq:2px2ym8y2pz2}
2+x^2 y-8 y^2+z^2=0.
\end{equation}
Modulo $4$ analysis on \eqref{eq:2px2ym8y2pz2} shows that $y$ must be equal to $1$ or $2$ modulo $4$.
We may rewrite the equation as
$$
yx^2+z^2=2(4y^2-1)=2(2|y|-1)(2|y|+1),
$$
which is an equation of the form \eqref{eq:jacobynosol} with $D(y)=-4y$. Note that $-4y$ and $2|y|+1$ are coprime, and the odd integers $2|y|-1$ and $2|y|+1$ are also coprime. Hence, by \eqref{eq:jacobidpi1}, we must have $\left(\frac{-4y}{2|y|+1}\right)=1$. Let us compute this Jacobi symbol
$$
\left(\frac{-4y}{2|y|+1}\right)=\left(\frac{-2|y|}{2|y|+1}\right) \cdot \left(\frac{2(y/|y|)}{2|y|+1}\right) = \left(\frac{2e}{2|y|+1}\right),
$$
where $e=y/|y|$.
If $y>0$, then
$$
\left(\frac{2e}{2|y|+1}\right)=\left(\frac{2}{2y+1}\right) =
\begin{cases}
1, \quad \text{if} \,\, y=0,3 \,\, \text{modulo} \,\, 4, \\
-1, \quad \text{if} \,\, y=1,2 \,\, \text{modulo} \,\, 4,
\end{cases}
$$
where we have used Proposition \ref{prop:Jacobi} (J3), a contradiction.
If $y<0$,
$$
\left(\frac{2e}{2|y|+1}\right)=\left(\frac{-2}{2|y|+1}\right) =
\begin{cases}
1, \quad \text{if} \,\, |y|=0,1 \,\, \text{modulo} \,\, 4, \\
-1, \quad \text{if} \,\, |y|=2,3 \,\, \text{modulo} \,\, 4,
\end{cases}
$$
where we have used Proposition \ref{prop:Jacobi} (J3).
If $y<0$ and $|y|=2,3$ modulo $4$, then, $y=1,2$ modulo $4$, a contradiction. Finally, equation \eqref{eq:2px2ym8y2pz2} has no integer solutions with $y=0$, hence, the equation has no integer solutions.

\vspace{10pt}

The next equation we will consider is
\begin{equation}\label{eq:1mypx2ym6y2m2z2}
1-y+x^2y-6y^2-2z^2=0.
\end{equation}	
It is clear that $y$ is odd and $x$ is even. After the substitutions $y=2Y+1$ and $x=2X$ for integers $X,Y$ and dividing by $2$, the equation can be rewritten as
$$
(4 Y + 2) X^2 - z^2 = (4 Y + 3) (3 Y + 1),
$$
and it is clear that $Y>0$. The above equation is of the form \eqref{eq:jacobynosol} with $D(Y)=4(4Y+2)$. Note that $4(4Y+2)$ and $4Y+3$ are coprime, hence condition (ii) after equation \eqref{eq:jacobynosol} works. Let us now check condition (iii), which states that $4 Y + 3$ and $3 Y + 1$ cannot share a prime factor $p$ such that $\left(\frac{4(4Y+2)}{p}\right)=-1$. Because $3(4Y+3)-4(3Y+1)=5$, the only possible common prime factor is $5$. If $3Y+1$ is divisible by $5$, then $Y=5t+3$ for some integer $t$. But then
$$
-1 = \left(\frac{4(4Y+2)}{p}\right) = \left(\frac{4(4(5t+3)+2)}{5}\right) = \left(\frac{56}{5}\right) = 1,
$$
a contradiction. Hence, all conditions (i)--(iii) hold, and, by \eqref{eq:jacobidpi1}, we must have $\left(\frac{4(4Y+2)}{4Y+3}\right)=1$. However, the Jacobi symbol
$$
\left(\frac{4(4Y+2)}{4Y+3}\right)=\left(\frac{4Y+2}{4Y+3}\right)=\left(\frac{-1}{4Y+3}\right)=-1
$$
for all $Y>0$ by Proposition \ref{prop:Jacobi} (J3). Hence \eqref{eq:1mypx2ym6y2m2z2} has no integer solutions.

\vspace{10pt}

The next equation we will consider is
\begin{equation}\label{eq:1p2x2ym4y2m3z2}
1 + 2x^2y - 4y^2 - 3z^2 = 0.
\end{equation}	
It is clear that $y>0$.
After adding $x^2$ to both sides and rearranging, the equation can be rewritten as
$$
(2y+1)(x^2-2y+1) = x^2+3z^2,
$$
where both factors $2y+1$ and $x^2-2y+1$ are positive.
Modulo $3$ analysis on \eqref{eq:1p2x2ym4y2m3z2} shows that $x$ is divisible by $3$, while $y$ is not. This implies that either $2y+1$ or $x^2-2y+1$ is equal to $2$ modulo $3$. Then by Proposition \ref{prop:abx2y2gen} (iii), the integers $2y+1$, $x^2-2y+1$, $x$ and $z$ must have a common prime factor $p$. But then $p$ is also a factor of $(x^2-2y+1)-x^2+(2y+1)=2$, hence $p=2$, but $2$ is not a factor of $2y+1$. This contradiction proves that equation \eqref{eq:1p2x2ym4y2m3z2} has no integer solutions.

\vspace{10pt}

The next equation we will consider is
\begin{equation}\label{eq:1p3x2ypy2zp4z2}
1+3 x^2 y+y^2 z+4 z^2=0	.
\end{equation}	
It is clear that $y$ must be odd. Multiplying the equation by $16$ and rearranging, we obtain
$$
(8 z + y^2)^2 + 48 y x^2 = (y^2-4)(y^2+4).
$$
Let $t=8 z + y^2$, then the above equation is of the form \eqref{eq:jacobynosol} with $D(y)=-4(48y)$. Note that the odd integers $y^2-4$ and $y^2+4$ are coprime. Also, any common prime factor $p$ of integers $-4(48y)$ and $y^2+4$ would also divide $-4y(48y)+192(y^2+4)=768 = 2^8\cdot 3$, so $p=2$ or $p=3$. However, because $y$ is odd, $y^2+4$ is not divisible by $2$ or $3$, hence, the integers $-4(48y)$ and $y^2+4$ are coprime. Then by \eqref{eq:jacobidpi1} we must have $\left(\frac{-4(48y)}{y^2+4}\right)=1$. However, the Jacobi symbol
$$
\left(\frac{-4(48y)}{y^2+4}\right)=\left(\frac{-3y}{y^2+4}\right)= \left(\frac{-3(y/|y|)}{y^2+4}\right) \cdot \left(\frac{|y|}{y^2+4}\right).
$$
Because $y$ is odd, $\left(\frac{|y|}{y^2+4}\right) \cdot \left(\frac{y^2+4}{|y|} \right)=(-1)^{\frac{|y|-1}{2} \frac{y^2+3}{2}} = 1$, and $\left(\frac{y^2+4}{|y|} \right)=\left(\frac{4}{|y|} \right)=1$. But then
$$
1 = \left(\frac{-4(48y)}{y^2+4}\right) = \left(\frac{-3(y/|y|)}{y^2+4}\right) = \left(\frac{- 3}{y^2+4}\right) = \begin{cases}
1, \quad \text{if} \,\, y=0 \,\, \text{modulo} \,\, 3, \\
-1, \quad \text{if} \,\, y=1,2 \,\, \text{modulo} \,\, 3,
\end{cases}
$$
where the third and the last equations follow from $\left(\frac{-1}{y^2+4}\right)=1$ by Proposition \ref{prop:Jacobi} (J3), and \eqref{eq:repos3Jacobi} (b), respectively.  Therefore, $y$ must be divisible by $3$. However, from the original equation, we must then have that $4z^2+1$ is divisible by $3$, which is impossible.

\vspace{10pt}

The next equation we will consider is
\begin{equation}\label{eq:1p2x2yp13yzpz2}
1+2 x^2 y+13 y z+z^2=0.
\end{equation}	
If $(x,y,z)$ is an integer solution to \eqref{eq:1p2x2yp13yzpz2}, then $v=2x^2+13z$ is an integer, such that $v=0,\pm 2,\pm 5,\pm 6$ modulo $13$. Then $z=\frac{v-2x^2}{13}$, and substituting this into \eqref{eq:1p2x2yp13yzpz2} and multiplying by $169$, we obtain
\begin{equation}\label{eq:1p2x2yp13yzpz2red}
169 + v^2 - 4 v x^2 + 4 x^4 + 169 v y=0.
\end{equation}
From this, it is clear that $169+4x^4$ is divisible by $v$. We can represent $4x^4+169$ as
\begin{equation}\label{eq:1p2x2yp13yzpz2red2}
4x^4+169=(2x^2+13)^2-13(2x)^2.
\end{equation}
Then applying Proposition \ref{prop:quadform} with $(a,b,c)=(1,0,-13)$, we can conclude that all odd prime divisors $p$ of $(2x^2+13)^2-13(2x)^2$ must be either divisors of $13$, or common divisors of $2x^2+13$ and $2x$, or such that $13$ is a quadratic residue modulo $p$. Hence $p=\pm 1,\pm 3, \pm 4$ modulo $13$ or $p=13$, so $v=0, \pm 1,\pm 3, \pm 4$ modulo $13$. Because we observed above that $v=0,\pm 2,\pm 5,\pm 6$ modulo $13$, the only possible case is when $v$ is divisible by $13$, which implies $x$ is also divisible by $13$. We can then make the change of variables $x=13X$ and $V=13v$ for integers $X,V$. Substituting this into \eqref{eq:1p2x2yp13yzpz2red} and cancelling $169$, we obtain
$$
1 + V^2 - 52 V X^2 + 676 X^4 + 13 V y=0.
$$
Any integer solution $(X,y,V)$ to this equation must have (i) $V^2+1$ is divisible by $13$ and (ii) $676X^4+1$ divisible by $V$. Condition (i) implies that $V=\pm 5$ modulo $13$. To analyse condition (ii), we can represent $676X^4+1$ as $(26X^2+1)^2-13(2X)^2$, then by Proposition \ref{prop:quadform}, we can conclude that all odd prime divisors $p$ of $(26X^2+1)^2-13(2X)^2$ must be such that $13$ is a quadratic residue modulo $p$. Therefore condition (ii) implies $V=\pm 1,\pm 3, \pm 4$ modulo $13$. Because this contradicts condition (i), equation \eqref{eq:1p2x2yp13yzpz2} has no integer solutions.

\vspace{10pt}

The next equation we will consider is
\begin{equation}\label{eq:1m2x2m2y2p7xyz}
1-2 x^2-2 y^2+7 x y z=0	.
\end{equation}	
This equation has an integer solution
$$
\begin{aligned}
&(x,y,z)=\\
& (24018909353284539760071841456040984563736520914233717685321, \\
&
99458886949844324979096721605263193799582258676037843673,
69).
\end{aligned}
$$

\vspace{10pt}

The next equation we will consider is
\begin{equation}\label{eq:1m16x2pxy2p2z2}
1-16 x^2+x y^2+2 z^2=0.
\end{equation}	
It is clear that both $x$ and $y$ are odd. After multiplication by $16$, we may rewrite the equation as
$$
(32 x - y^2)^2 - 2 (8 z)^2 = y^4 + 64=(y^2 - 4 y + 8) (y^2 + 4 y + 8).
$$
Because $y$ is odd, the integers $a = y^2 - 4 y + 8$ and $b = y^2 + 4 y + 8$ are odd. Hence, if $p$ is any common prime factor of the integers $a$ and $b$, then $p$ is a divisor of $b-a=8y$, so $p$ is a divisor of $y$. However, $p$ is then also a divisor of $b-y(y+4)=8$, a contradiction. Hence, the integers $a$ and $b$ are coprime, and their product is of the form $X^2-2Y^2$, hence, by Proposition \ref{prop:2y2mz2}, they must both be of this form. However, $y^2 \pm 4 y + 8=5$ modulo $8$ when $y$ is odd, hence these integers cannot be of the form $X^2-2Y^2$.

\vspace{10pt}

The next equation we will consider is
\begin{equation}\label{eq:1p2x2ym16y2pz2}
1+2 x^2 y-16 y^2+z^2=0.
\end{equation}	
Modulo $4$ analysis shows that all variables must be odd. Multiplying the equation by $16$ and rearranging, we obtain
$$
(x^2 - 16y- 4z)(x^2 - 16y + 4z) = x^4 + 16=(x^2+4)^2-2(2x)^2.
$$
As the right-hand side is odd and $x^2+4$ and $2x$ are coprime, then all prime factors of $x^4 + 16$ are $\pm 1$ modulo $8$ by Proposition \ref{cor:quadform} (c).
However, for all odd $x$ and $z$, $x^2 - 16y + 4z=5$ modulo $8$, a contradiction.

\vspace{10pt}

The next equation we will consider is
\begin{equation}\label{eq:1pxyp2x2ypy2p2xz2}
1+x y+2 x^2 y+y^2+2 x z^2=0.
\end{equation}	
We can rewrite this equation as
$$
x(-y-2xy-2z^2)=y^2+1.
$$
Modulo $4$ analysis shows that $x$ is $2$ modulo $4$ and $y$ is odd, hence we may let $x=4u+2$ for integer $u$. Then,
$$
2(2u+1)(-y-2(4u+2)y-2z^2)=y^2+1.
$$
Because $|2u+1|$ is a positive odd divisor of $y^2+1$, it must be $1$ modulo $4$. Hence, either $u$ is positive and even, or it is negative and odd. Write $y=2v-4(2u+1)^2-(2u+1)$ for a new integer variable $v$. After substituting into the above equation, dividing by $4$ and rearranging, we obtain
$$
v^2 + (2u+1)z^2 = (8u^2+9u+2)(8u^2+9u+3).
$$
Modulo $4$ analysis on this equation shows that $u \neq 3$ modulo $4$.
If $u$ is positive and even, then $8u^2+9u+3$ is odd. Note that $8u^2+9u+2$ and $8u^2+9u+3$ are coprime, and $8u^2+9u+3$ and $2u+1$ are coprime, hence, by \eqref{eq:jacobidpi1}, the Jacobi symbol $\left(\frac{-4(2u+1)}{8u^2+9u+3}\right)$ should be equal to $1$. Let us compute this. We have
$$
\left(\frac{-4(2u+1)}{8u^2+9u+3}\right) =\left(\frac{-1}{8u^2+9u+3}\right) \cdot \left(\frac{2u+1}{8u^2+9u+3}\right).
$$
Next,
$$
\left(\frac{2u+1}{8u^2+9u+3}\right) = \left(\frac{8u^2+9u+3}{2u+1}\right) = \left(\frac{u+1}{2u+1}\right) = \left(\frac{2u+1}{u+1}\right) = \left(\frac{-1}{u+1}\right),
$$
because $\left(\frac{8u^2+9u+3}{2u+1}\right) \cdot \left(\frac{2u+1}{8u^2+9u+3}\right)=-1^{u (\frac{8u^2+9u+2}{2})}=1$, and $\left(\frac{u+1}{2u+1}\right) \cdot \left(\frac{2u+1}{u+1}\right)=-1^{\frac{u^2}{2}}=1$ for even $u$, and $(8 u^2 + 9 u + 3) - (4u+2) (2 u + 1)=u+1$. Hence,
$$
\left(\frac{-4(2u+1)}{8u^2+9u+3}\right) = \left(\frac{-1}{8u^2+9u+3}\right) \cdot \left(\frac{-1}{u+1}\right).
$$
For even $u>0$, one of $8u^2+9u+3$ and $u+1$ is equal $1$ modulo $4$, and the other is equal to $3$ modulo $4$. Hence, by Proposition \ref{prop:Jacobi} (J3),
$$
\left(\frac{-1}{8u^2+9u+3}\right) \cdot \left(\frac{-1}{u+1}\right)=1(-1)=-1,
$$
a contradiction.

If $u$ is negative and odd, then we must have $u=-4w+1$ for positive integer $w$. The equation then reduces to
$$
v^2 - (8w-3)z^2 = (128w^2-100w+19)(128w^2-100w+20).
$$
Note that $128w^2-100w+19$ and $128w^2-100w+20$ are coprime and $128w^2-100w+19$ and $8w-3$ are coprime, hence, by \eqref{eq:jacobidpi1}, the Jacobi symbol
\[
J = \left(\frac{4(8w-3)}{128w^2-100w+19}\right)
\]
should be equal to $1$. Let us compute $J$. We have
$$
J =\left(\frac{8w-3}{128w^2-100w+19}\right)=\left(\frac{128w^2-100w+19}{8w-3}\right)=\left(\frac{4w-2}{8w-3}\right)
$$
because
$$
\left(\frac{8w-3}{128w^2-100w+19}\right) \cdot \left(\frac{128w^2-100w+19}{8w-3}\right)=-1^{\frac{128w^2-100w+18}{2} \frac{8w-4}{2}}=1,
$$
and $(128 w^2 - 100 w + 19) - (16 w - 7) (8 w - 3)=4w-2$.
Further, we have
$$
\left(\frac{4w-2}{8w-3}\right)=\left(\frac{2}{8w-3}\right) \cdot \left(\frac{2w-1}{8w-3}\right) = - \left(\frac{2w-1}{8w-3}\right) = - \left(\frac{8w-3}{2w-1}\right) = - \left(\frac{1}{2w-1}\right) =-1
$$
where we have used $\left(\frac{2w-1}{8w-3}\right) \cdot \left(\frac{8w-3}{2w-1}\right)=-1^{\frac{8w-4}{2} \frac{2w-2}{2}}=1$, $8w-3=-3$ modulo $8$, Proposition \ref{prop:Jacobi} (J3), and $8w-3-4(2w-1)=1$. In conclusion, $J=-1$, which is a contradiction. Hence, equation \eqref{eq:1pxyp2x2ypy2p2xz2} has no integer solutions.

\vspace{10pt}

	The next equation we will consider is
	\begin{equation}\label{eq:2p8x2ypxy2pz2}
		2+8 x^2 y+x y^2+z^2=0.
	\end{equation}	
	Modulo $4$ analysis shows that $x$ is $1$ or $2$ modulo $4$. We can rewrite the equation as
	$$
	x(y+4x)^2+z^2=16x^3-2=2(2x-1)(4x^2+2x+1).
	$$
	Considering the above equation as a quadratic in $y+4x,z$ with parameter $x$, it is of the form \eqref{eq:jacobynosol} and its discriminant $D(x)=-4x$, which is also coprime with $4x^2+2x+1$. Note that the odd integers $2x-1$ and $4x^2+2x+1>0$ are coprime. Hence, by \eqref{eq:jacobidpi1}, we must have $\left(\frac{-4x}{4x^2+2x+1}\right)=1$. Let us now compute the Jacobi symbol
	$$
	J = \left(\frac{-4x}{4x^2+2x+1}\right)=\left(\frac{-x}{4x^2+2x+1}\right).
	$$
	We will first consider the case where $x \geq 0$ and $x$ is $1$ modulo $4$. Then
	$$
	J =\left(\frac{-1}{4x^2+2x+1}\right)\left(\frac{x}{4x^2+2x+1}\right)=-\left(\frac{x}{4x^2+2x+1}\right)
	$$
	using Proposition \ref{prop:Jacobi} (J3) and because $4x^2+2x+1=3$ modulo $4$. We also have
	$$
	\left(\frac{x}{4x^2+2x+1}\right) = \left(\frac{4x^2+2x+1}{x}\right)
	$$
	because $\left(\frac{x}{4x^2+2x+1}\right) \cdot \left(\frac{4x^2+2x+1}{x}\right) = -1^{\frac{4x^2+2x}{2} \frac{x-1}{2}}=1$ when $x$ is $1$ modulo $4$. Then,
	$$
	J =-\left(\frac{x}{4x^2+2x+1}\right)=- \left(\frac{4x^2+2x+1}{x}\right)=-\left(\frac{1}{x}\right)=-1,
	$$
	a contradiction.
	
	Let us now consider the case where $x \geq 0$ and $x$ is $2$ modulo $4$. Let us make the substitution $x=4t+2$ for positive integer $t$. Then we have
	$$
	J =\left(\frac{-2}{64t^2+72t+21}\right) \cdot \left(\frac{2t+1}{64t^2+72t+21}\right)=-\left(\frac{2t+1}{64t^2+72t+21}\right)
	$$
	by Proposition \ref{prop:Jacobi} (J3) and because $64t^2+72t+21=5$ modulo $8$. Then
	$$
	\left(\frac{2t+1}{64t^2+72t+21}\right)=\left(\frac{64t^2+72t+21}{2t+1}\right)
	$$
	because $\left(\frac{2t+1}{64t^2+72t+21}\right) \cdot \left(\frac{64t^2+72t+21}{2t+1}\right)=-1^{\frac{2t}{2} \frac{64t^2+72t+20}{2}}=1$. Also, $64 t^2 + 72 t + 21 - (20 + 32 t) (2 t + 1)=1$, hence
	$$
	J =-\left(\frac{2t+1}{64t^2+72t+21}\right)=-\left(\frac{64t^2+72t+21}{2t+1}\right)=-\left(\frac{1}{2t+1}\right)=-1,
	$$
	a contradiction.
	
	Let us now consider the case where $x < 0$ and $x$ is $1$ modulo $4$. Then
	$$
	J = \left(\frac{-x}{4x^2+2x+1}\right)=\left(\frac{|x|}{4x^2-2|x|+1}\right)=-\left(\frac{4x^2-2|x|+1}{|x|}\right) = -\left(\frac{1}{|x|}\right)=-1,
	$$
	because $\left(\frac{|x|}{4x^2-2|x|+1}\right) \cdot \left(\frac{4x^2-2|x|+1}{|x|}\right) =-1^{\frac{|x|-1}{2} \frac{4x^2-2|x|}{2}}=-1$, as $|x|=3$ modulo $4$. This is a contradiction.
	
	Let us finally consider the case where $x < 0$ and $x$ is $2$ modulo $4$. Let us make the substitution $x=4t+2$ with $t<0$ an integer. We want to compute
	$$
	J = \left(\frac{-x}{4x^2+2x+1}\right)=\left(\frac{-4t-2}{64t^2+72t+21}\right)=\left(\frac{4|t|-2}{64t^2-72|t|+21}\right).
	$$
	Then,
	$$
	J =
	\left(\frac{2}{64t^2-72|t|+21}\right) \cdot \left(\frac{2|t|-1}{64t^2-72|t|+21}\right)
	=-\left(\frac{2|t|-1}{64t^2-72|t|+21}\right)
	$$
	by Proposition \ref{prop:Jacobi} (J3) and because $64t^2-72|t|+21=5$ modulo $8$. Also,
	$$
	\left(\frac{2|t|-1}{64t^2-72|t|+21}\right)=\left(\frac{64t^2-72|t|+21}{2|t|-1}\right)
	$$
	because $\left(\frac{2|t|-1}{64t^2-72|t|+21}\right) \cdot \left(\frac{64t^2-72|t|+21}{2|t|-1}\right)=-1^{\frac{2|t|-2}{2} \frac{64t^2-72|t|+20}{2}}=1$. We also have that $64t^2-72|t|+21-(2|t|-1)(32|t|-20)=1$, hence,
	$$
	J =-\left(\frac{2|t|-1}{64t^2-72|t|+21}\right)
	=-\left(\frac{64t^2-72|t|+21}{2|t|-1}\right)=-\left(\frac{1}{2|t|-1}\right)=-1,
	$$
	a contradiction.
	
		\vspace{10pt}
		
		The next equation we will consider is
		\begin{equation}\label{eq:1p8x2ypy2p2xz2}
			1+8 x^2 y+y^2+2 x z^2=0.
		\end{equation}	
		Modulo $4$ analysis shows that $x$ is odd. We can rewrite the equation as
		$$
		(4x^2+y)^2+2x z^2=16x^4-1=(4x^2+1)(4x^2-1).
		$$
		Considering the above equation as a quadratic in $4x^2+y,z$ with parameter $x$, it is of the form \eqref{eq:jacobynosol} and its discriminant $D(x)=-8x$  is also coprime with $4x^2+1$. Note that the odd integers $4x^2+1$ and $4x^2-1$ are coprime. Hence, by \eqref{eq:jacobidpi1}, we must have $\left(\frac{-8x}{4x^2+1}\right)=1$. Let us now compute this Jacobi symbol. We have
		$$
		\left(\frac{-8x}{4x^2+1}\right)=\left(\frac{-2x}{4x^2+1}\right)=\left(\frac{-2e}{4x^2+1}\right) \cdot \left(\frac{|x|}{4x^2+1}\right)=- \left(\frac{|x|}{4x^2+1}\right),
		$$
		where $e=x/|x|$, and we have used that $\left(\frac{2}{4x^2+1}\right)=\left(\frac{-2}{4x^2+1}\right)=-1$ by Proposition \ref{prop:Jacobi} (J3) because $4x^2+1$ is $5$ modulo $8$. Also,
		$$
		\left(\frac{|x|}{4x^2+1}\right) \cdot \left(\frac{4x^2+1}{|x|}\right)=-1^{\frac{4x^2+1-1}{2} \frac{|x|-1}{2}}=1,
		$$
		so, $\left(\frac{|x|}{4x^2+1}\right) = \left(\frac{4x^2+1}{|x|}\right)$, and
		$$
		\left(\frac{4x^2+1}{|x|}\right)=\left(\frac{1}{|x|}\right)=1.
		$$
		Hence,
		$$
		\left(\frac{-8x}{4x^2+1}\right)=- \left(\frac{|x|}{4x^2+1}\right)=-1,
		$$
		a contradiction.
		
		The next equation we will consider is
		\begin{equation}\label{eq:1p2xpx2yp4y2p2z2}
			1+2 x+x^2 y+4 y^2+2 z^2 = 0.
		\end{equation}
		This equation was originally solved by Denis Shatrov on the MathOverflow website,\footnote{ \url{https://mathoverflow.net/questions/459877}} with quite a long proof. Here, we present a much simpler argument. It is clear that $y$ is negative and odd. First multiply the equation by $y$ and rewrite it as
		$$
		y-1+(1+xy)^2+4y^3+2yz^2=0.
		$$
		Now let $X=1+xy$ and $t=-y$ to obtain the equation
		$$
		X^2-2tz^2 = 4t^3+t+1 = (2t+1)(2t^2-t+1),
		$$
		where $t$ is now positive and odd. Then \eqref{eq:jacobidpi1} implies that $\left(\frac{8 t}{2 t+1}\right) = 1$. On the other hand,
		$$
		\left(\frac{8 t}{2 t+1}\right) = \left(\frac{2 t}{2 t+1}\right) = \left(\frac{-1}{2 t+1}\right) = -1
		$$
		for every odd positive $t$, which is a contradiction.
		
		The final equation we consider is the equation
		\begin{equation}\label{eq:1p4x3pxy2p2yz2}
			1+4x^3+x y^2+2y z^2=0.
		\end{equation}
		This equation was solved by Denis Shatrov on the MathOverflow  website.\footnote{\url{https://mathoverflow.net/questions/460003}} We can easily see that $x$ must be odd. Multiplying the equation by $x$ and introducing $Y=xy+z^2$, the equation can be rewritten as
		\begin{equation}\label{eq:1p4x3pxy2p2yz2red}
			x+Y^2=(z^2-2x^2)(z^2+2x^2).
		\end{equation}
		Let us prove that \eqref{eq:1p4x3pxy2p2yz2red} has no integer solutions with odd $x$. Note that \eqref{eq:1p4x3pxy2p2yz2red} implies that $-x$ is a quadratic residue modulo $|z^2-2x^2|$ and $z^2+2x^2$.
		
		Modulo $4$ analysis of \eqref{eq:1p4x3pxy2p2yz2red} shows that for odd $x$ the equation can only have integer solutions when $z$ is odd and $x=1$ modulo $4$, or when $z$ is even and $x=3$ modulo $4$. Also modulo $8$ analysis shows that in the second case $x$ is in fact equal to $3$ modulo $8$.
		
		We now have $6$ cases to consider.
		
		\smallskip
		\textbf{Case 1: $x>0$ and $z$ is odd.} Because $z$ is odd, $x=1$ modulo $4$. Also, $x>0$ and \eqref{eq:1p4x3pxy2p2yz2red} implies that $z^2-2x^2>0$. As $-x$ is a quadratic residue modulo $z^2-2x^2$, we have $\left(\frac{-x}{z^2-2x^2}\right)=1$. Let us calculate this Jacobi symbol. We have
		\begin{align*}
			\left(\frac{-x}{z^2-2x^2}\right) &=\left(\frac{-1}{z^2-2x^2}\right) \cdot \left(\frac{x}{z^2-2x^2}\right)\\ &=-\left(\frac{x}{z^2-2x^2}\right)=-\left(\frac{z^2-2x^2}{x}\right)=-\left(\frac{z^2}{x}\right)=-1.
		\end{align*}
		Here, we have used that $z^2-2x^2=3$ modulo $4$ as $x$ and $z$ are both odd, and that $\left(\frac{x}{z^2-2x^2}\right) \cdot \left(\frac{z^2-2x^2}{x}\right)=(-1)^{\frac{x-1}{2} \frac{z^2-2x^2-1}{2}}=1$ as $x=1$ modulo $4$. Therefore, $\left(\frac{-x}{z^2-2x^2}\right)=-1$, a contradiction.
		
		\smallskip
		\textbf{Case 2: $-Y^2<x<0$ and $z$ is odd.} Let $x=-x'$, so $x'>0$ and because $z$ is odd, $x'=3$ modulo $4$. Then \eqref{eq:1p4x3pxy2p2yz2red} reduces to
		$$
		Y^2-x'=(z^2-2x'^2)(z^2+2x'^2).
		$$
		Because $Y^2-x'>0$, we have $z^2-2x'^2>0$. Also, $x'$ is a quadratic residue modulo $z^2-2x'^2$, so $\left(\frac{x'}{z^2-2x'^2}\right)=1$. Let us calculate this Jacobi symbol. We have
		$$
		\left(\frac{x'}{z^2-2x'^2}\right) =(-1)^{\frac{x'-1}{2} \frac{z^2-2x'^2-1}{2}} \left(\frac{z^2-2x'^2}{x'}\right) 
		= -\left(\frac{z^2-2x'^2}{x'}\right)=-\left(\frac{z^2}{x'}\right)=-1,
		$$
		which is a contradiction.
		
		\smallskip
		\textbf{Case 3: $x<-Y^2$ and $z$ is odd.} Let $x=-x'$, so $x'>0$ and because $z$ is odd, $x'=3$ modulo $4$. Then \eqref{eq:1p4x3pxy2p2yz2red} reduces to
		$$
		x'-Y^2=(2x'^2-z^2)(z^2+2x'^2),
		$$
		and we have that $x'$ is a quadratic residue modulo $2x'^2-z^2$, so $\left(\frac{x'}{2x'^2-z^2}\right)=1$. Let us calculate this Jacobi symbol. We have
		\begin{align*}
			\left(\frac{x'}{2x'^2-z^2}\right) &=(-1)^{\frac{x'-1}{2} \frac{2x'^2-z^2-1}{2}}\left(\frac{2x'^2-z^2}{x'}\right)\\
			&=\left(\frac{2x'^2-z^2}{x'}\right)=\left(\frac{-z^2}{x'}\right)=\left(\frac{-1}{x'}\right)\left(\frac{z}{x'}\right)^2=-1,
		\end{align*}
		again a contradiction.
		
		\smallskip
		\textbf{Case 4: $x>0$ and $z$ is even.} Because $x>0$ and $z$ is even, $x=3$ modulo $8$. Let $z=2z'$ for integer $z'$, then \eqref{eq:1p4x3pxy2p2yz2red} reduces to
		$$
		x+Y^2=4(2z'^2-x^2)(2z'^2+x^2),
		$$
		and we have that $-x$ is a quadratic residue modulo $2z'^2-x^2$, so $\left(\frac{-x}{2z'^2-x^2}\right)=1$. Let us calculate this Jacobi symbol. We have
		$$
		\left(\frac{-x}{2z'^2-x^2}\right)=\left(\frac{-1}{2z'^2-x^2}\right) \cdot \left(\frac{x}{2z'^2-x^2}\right) =\left(\frac{2z'^2-x^2}{x}\right) = \left(\frac{2z'^2}{x}\right)= \left(\frac{2}{x}\right)=-1,
		$$
		where we have used that $x=3$ modulo $8$,
		\begin{align*}
			&\left(\frac{-1}{2z'^2-x^2}\right)= \begin{cases}
				1,& \quad \text{if} \,\, z' \,\, \text{is odd}, \\
				-1,& \quad \text{if} \,\, z' \,\, \text{is even,}
			\end{cases}
			\quad \text{and} \\
			&\left(\frac{2z'^2-x^2}{x}\right) \cdot \left(\frac{x}{2z'^2-x^2}\right)= \begin{cases}
				1, & \quad \text{if} \,\, z' \,\, \text{is odd}, \\
				-1, &  \quad \text{if} \,\, z' \,\, \text{is even.}
			\end{cases}
		\end{align*}
		Hence, $\left(\frac{-x}{2z'^2-x^2}\right)=-1$, a contradiction.
		
		\smallskip
		\textbf{Case 5: $-Y^2<x<0$ and $z$ is even.} Let $x=-x'$, so $x'>0$ and because $z$ is even, $x'=5$ modulo $8$. Let $z=2z'$ for integer $z'$, then \eqref{eq:1p4x3pxy2p2yz2red} reduces to
		$$
		Y^2-x'=4(2z'^2-x'^2)(2z'^2+x'^2),
		$$
		and we have that $x'$ is a quadratic residue modulo $2z'^2-x'^2$, so $\left(\frac{x'}{2z'^2-x'^2}\right)=1$. Let us calculate this Jacobi symbol. We have
		\begin{align*}
			\left(\frac{x'}{2z'^2-x'^2}\right) &=(-1)^{\frac{2z'^2-x'^2-1}{2} \frac{x'-1}{2}}\left(\frac{2z'^2-x'^2}{x'}\right)\\
			&=\left(\frac{2z'^2-x'^2}{x'}\right)=\left(\frac{2z'^2}{x'}\right)=\left(\frac{2}{x'}\right)=-1,
		\end{align*}
		where we have used that $x'=5$ modulo $8$. Hence, $\left(\frac{x'}{2z'^2-x'^2}\right)=-1$, a contradiction.
		
		\smallskip
		\textbf{Case 6: $x<-Y^2$ and $z$ is even.}	Let $x=-x'$, so $x'>0$ and because $z$ is even, $x'=5$ modulo $8$. Let $z=2z'$ for integer $z'$, then \eqref{eq:1p4x3pxy2p2yz2red} reduces to
		$$
		x'-Y^2=4(x'^2-2z'^2)(2z'^2+x'^2),
		$$
		and we have that $x'$ is a quadratic residue modulo $x'^2-2z'^2$, so $\left(\frac{x'}{x'^2-2z'^2}\right)=1$. Let us calculate this Jacobi symbol.
		\begin{align*}
			\left(\frac{x'}{x'^2-2z'^2}\right) &=(-1)^{\frac{x'^2-2z'^2-1}{2} \frac{x'-1}{2}}\left(\frac{x'^2-2z'^2}{x'}\right)\\
			&=\left(\frac{x'^2-2z'^2}{x'}\right)=\left(\frac{-2z'^2}{x'}\right)=\left(\frac{-2}{x'}\right)=-1,
		\end{align*}
		again a contradiction.
		
		The only case not considered is $x=-Y^2$, but it is easy to see that this case also does not lead to integer solutions. Therefore equation \eqref{eq:1p4x3pxy2p2yz2red} has no integer solutions with odd $x$, which implies that the original equation \eqref{eq:1p4x3pxy2p2yz2} has no integer solutions.

\newpage

\section{List of solved equations}
In the below table, we list all equations solved in the book and in this document, ordered by $H$. 
Similar to Section 8.5 of the book, in the ``Question'' column we use the following shortcuts. For each equation, we select the ``best'' description of the solution set for column ``Question''. For example, if we have solved both ``Infinite'' and ``Parametrize'', we will state ``Parametrize''.

\begin{itemize}
	\item ``No solutions'' means that it has been proven that the equation has no integer solutions.
	\item ``List'' means that the solution set is finite and all integer solutions are listed.
	\item ``Parametrize'' means that the solution set to the equation can be represented as a finite union of polynomial families.
	\item ``Explicit'' means that it is known that there is an algorithm that can solve the equation in a finite time.
	\item ``Describe all'' means that all integer solutions to the equation can be explicitly described in one of the ways (i)-(v) listed in Section 4.1.7 of the book.
	\item ``Primitive'' means that all primitive integer solutions to the equation have been found.
	\item ``Infinite'' means that it has been proven that the equation has infinitely many integer solutions.
	\item ``Large'' means that it has been proven that the equation has integer solutions with absolute values of all variables arbitrarily large.
	\item ``Non-zero'' means that it is determined that a homogeneous equation has any integer solution with \textbf{all} variables non-zero, see Problem \ref{prob:homlarge}.
	\item ``Existence'' means that an integer solution to the equation has been found.
	\item ``No positive'' means that it is determined that the equation has no solutions in positive integers.
\end{itemize}



\appendix

\section{Progress on open questions}
\subsection{Open question 3.78} \label{open:3.78}

In this section we provide the complete list of integer solutions to four two-variable equations left open in \cite{mainbook}.

Recall the notation $v_p(n)$, for a prime $p$ and integer $n$, which denotes the exponent with which $p$ enters the prime factorisation of $n$. 

First we consider equations
\begin{equation}\label{eq:x4my3pxypx}
x^4-y^3+xy+x=0,
\end{equation}
and 
\begin{equation}\label{eq:x4py3pxypx}
x^4+y^3+xy+x=0.
\end{equation}
First consider the case $xy=0$. For \eqref{eq:x4my3pxypx}, if $x=0$, then $y=0$. If $y=0$, then $x^4+x=0$, hence $x=0$ or $x=-1$. This gives two solutions $(x,y)=(-1,0)$ and $(0,0)$. For \eqref{eq:x4py3pxypx}, we obtain the same two solutions $(x,y)=(-1,0)$ and $(0,0)$.

Assume that $xy\neq 0$. Then, all the monomials of the equation are non-zero.
Let $p$ be any prime factor of $x$, which implies that $p$ is also a factor of $y$. Then $p$ enters the prime factorisations of the monomials of the equation with exponents $4v_p(x), 3v_p(y), v_p(x)+v_p(y), v_p(x)$, respectively. Because $4v_p(x)>v_p(x)$ and $v_p(x)+v_p(y)>v_p(x)$, we must have $3v_p(y)=v_p(x)$. Hence, $x=t^3$ and $y=ts$ for some non-zero integers $t$ and $s$.

Let us solve Diophantine equation \eqref{eq:x4my3pxypx}. 
Substituting $x=t^3$ and $y=ts$ for some non-zero integers $t$ and $s$ into the equation and cancelling $t^3$, we obtain
$$
t^9-s^3+ts+1=0, \quad ts\neq 0.
$$  
The optimisation problem
\begin{flushleft}
\texttt{N[Maximize[Min[Abs[t], Abs[s], Abs[t\^{}3 - s]], t\^{}9 - s\^{}3 + t s + 1 == 0, \{t, s\}, Reals]]}
\end{flushleft}
returns optimal objective value less than $1$. Hence in any real solution we must have either $|t|<1$, or $|s|<1$, or $|t^3-s|<1$; because $t$ and $s$ are integers with $ts\neq 0$, this implies that $s=t^3$. But then the equation reduces to $t^4+1=0$, which is impossible.

Next, we will solve \eqref{eq:x4py3pxypx} in the same way. 
Substituting $x=t^3$ and $y=ts$ for some non-zero integers $t$ and $s$ into the equation and cancelling $t^3$, we obtain
$$
t^9+s^3+ts+1=0, \quad ts\neq 0.
$$  
The optimisation problem
\begin{flushleft}
\texttt{N[Maximize[Min[Abs[t], Abs[s], Abs[t\^{}3 + s], Abs[t\^{}3 - s]], t\^{}9 + s\^{}3 + t s + 1 == 0, \{t, s\}, Reals]]}
\end{flushleft}
returns optimal objective value less than $1$; because $t$ and $s$ are integers with $ts\neq 0$, this implies that $s=\pm t^3$. Substituting these into $t^9+s^3+ts+1=0$ returns $(s,t)=\pm(-1,1),(-1,-1)$. Then in the original variables, we obtain the solutions $(x,y)=(-1,-1), \pm(1,-1)$. 

Therefore, the complete set of integer solutions to equation \eqref{eq:x4my3pxypx} is $(x,y)=(-1,0)$ and $(0,0)$, while the complete set of integer solutions to equation \eqref{eq:x4py3pxypx} is $(x,y)=(-1,-1)$, $(-1,0)$, $(0,0)$, $\pm(1,-1)$. 

Now we consider the equations 
\begin{equation}\label{eq:x4py3pxymy}
x^4+y^3+xy-y=0
\end{equation}
and
\begin{equation}\label{eq:x4py3pxypy}
x^4+y^3+xy+y=0.
\end{equation}
Assume that $xy\neq 0$. Then all monomials of the equation are non-zero.
Let $p$ be any prime factor of $y$, which implies that $p$ is also a factor of $x$. Then $p$ enters the prime factorisations of the monomials of the equation with exponents $4v_p(x), 3v_p(y), v_p(x)+v_p(y), v_p(y)$, respectively. Because $3v_p(y)>v_p(y)$ and $v_p(x)+v_p(y)>v_p(y)$, we must have $4v_p(x)=v_p(y)$. Hence, $y=et^4$ and $x=ets$ for some non-zero integers $t$ and $s$, where $e=\pm 1$. 

We now solve the equation \eqref{eq:x4py3pxymy}. If $x=0$ then $y=\pm 1,0$. If $y=0$ then $x=0$. 

Now assume that $xy\neq 0$. Substituting $y=et^4$ and $x=ets$ into the equation and cancelling $t^4$, we obtain
$$
et^8+s^4+ts-e=0, \quad ts\neq 0, \quad e = \pm 1.
$$  
This equation is easy. The optimisation problems
\begin{flushleft}
\texttt{N[Maximize[Min[Abs[t], Abs[s], Abs[t\^{}2 - s],Abs[t\^{}2+s]], et\^{}8 + s\^{}4 + t s - e == 0, \{t, s\}, Reals]]}
\end{flushleft}
returns optimal objective value less than $1$. Because $t$ and $s$ are integers with $ts\neq 0$, this implies that $s=\pm t^2$, the equation then reduces to $2t^8\pm t^3-1=0$ which returns solutions $(s,t)=\pm(1,-1)$, which in the original variables correspond to solutions $(x,y)=\pm(1,-1)$. 

Hence all integer solutions to \eqref{eq:x4py3pxymy} are $(x,y)=(0,0),(0,\pm 1),\pm(1,-1)$.

We finally solve the equation \eqref{eq:x4py3pxypy}.
If $x=0$ then $y=0$. If $y=0$ then $x=0$. 

Substituting $y=et^4$ and $x=ets$ into the equation and cancelling $t^4$, we obtain
$$
et^8+s^4+ts+e=0, \quad ts\neq 0.
$$  
When $e=1$ the equation has no real solutions. When $e=-1$, the optimisation problem
\begin{flushleft}
\texttt{N[Maximize[Min[Abs[t], Abs[s], Abs[t\^{}2 - s],Abs[t\^{}2+s]], -t\^{}8 + s\^{}4 + t s -1 == 0, \{t, s\}, Reals]]}
\end{flushleft}
returns optimal objective value less than $1$. Because $t$ and $s$ are integers with $ts\neq 0$, this implies that $s=\pm t^2$, the equation then reduces to $\pm t^3-1=0$ which returns solutions $(s,t)=\pm(1,1)$, which in the original variables corresponds to the solution $(x,y)=(-1,-1)$. 

Hence all integer solutions to \eqref{eq:x4py3pxypy} are $(x,y)=(0,0),(-1,-1)$.

\subsection{Open question 4.52} \label{open:4.52}

In this section we are considering equations
\begin{equation}\label{eq:genfermat}
	ax^p+by^q=cz^r,
\end{equation} 
where $x,y,z$ are integer variables, while positive integers $p,q,r$ and integers $a,b,c$ are parameters. These equations are referred to as generalised Fermat equations. For such equations, we determine them as solved if we describe all \emph{primitive} integer solutions, that is, solutions $(x,y,z)$ with $\gcd(x,y,z)=1$. Specifically, we wish to progress on Open Question 4.52 of \cite{mainbook}, which we state here for convenience. 
\begin{problem}{\cite[Open Question 4.52]{mainbook}}
	Find all primitive solutions to the equations listed in Table \ref{tab:H60main}. Investigate for which of the listed equations this problem has already been solved in the literature.
\end{problem}

	\begin{center}

	\captionof{table}{\label{tab:H60main} Equations \eqref{eq:genfermat} with $|a|2^p+|b|2^q+|c|2^r \leq 60$, which are uninvestigated in \cite[Section 4.3.6]{mainbook}.}
\end{center}

We begin with the second part of the open question. 
Our literature review \cite{wilcox2024generalisedfermatequationsurvey} shows that four equations from Table \ref{tab:H60main} have been solved in the literature. These equations are listed in Table \ref{tab:prim}.

	\begin{center}
	%
	\captionof{table}{\label{tab:prim} Primitive integer solutions to equations in Table \ref{tab:H60main} solved in the literature.}
\end{center} 

 Removing these equations from Table \ref{tab:H60main}, we obtain Table \ref{tab:H60Fermat2} which contains the equations that remain to be solved.

	\begin{center}
		%
		\captionof{table}{\label{tab:H60Fermat2} Remaining equations \eqref{eq:genfermat} with $|a|2^p+|b|2^q+|c|2^r \leq 60$. Those in bold are solved in this section.}
	\end{center}

\begin{proposition}
The only primitive integer solutions to 
\begin{equation}\label{eq:x4my4m2z3}
x^4-y^4=2z^3
\end{equation}
are $(x,y,z)=(\pm 1,\pm 1,0)$.
\end{proposition}

\begin{proof}
We first note that if $(x,y,z)$ is a primitive solution to \eqref{eq:x4my4m2z3}, then $x,y,z$ are pairwise coprime, and both $x$ and $y$ are odd. Then $S=(x^2-y^2)/2$ and $T=(x^2+y^2)/2$ are integers, with $S$ being even, and $T$ odd. If $p$ is a common prime divisor of $S$ and $T$, then $p$ divides $S+T=x^2$ and $S-T=y^2$ hence $p$ divides both $x$ and $y$, which is a contradiction. Hence, $S$ and $T$ are coprime, and  equation \eqref{eq:x4my4m2z3} reduces to    
\begin{equation}\label{red:x4my4m2z3}
	2ST=z^3 \quad \text{with} \quad \gcd(S,T)=1.
\end{equation} 

The integer solutions to \eqref{red:x4my4m2z3} are $(S,T,z)=(4u^3,v^3,2uv)$ or $(u^3,4v^3,2uv)$, where $u,v$ are arbitrary coprime integers. Note here that because $S$ is even and $T$ is odd, the only possible solution is $(S,T,z)=(4u^3,v^3,2uv)$.

Then in the original variables, we obtain the system of equations $8u^3=x^2-y^2$ and $2v^3=x^2+y^2$. We can factorise the first equation and write it as
\begin{equation}\label{red2:x4my4m2z3}
(x+y)(x-y)=8u^3.
\end{equation}
Because $x$ and $y$ are coprime and both odd, $a=(x+y)/2$ and $b=(x-y)/2$ are coprime integers, and \eqref{red2:x4my4m2z3} reduces to
$$
ab=2u^3 \quad \text{with} \quad \gcd(a,b)=1
$$
whose integer solutions are $(a,b,u)=(d^3,2e^3,de)$ or $(2d^3,e^3,de)$ for integers $d,e$. Then either $(x,y)=(d^3+2e^3, d^3-2e^3)$ or $(x,y)=(2d^3+e^3, 2d^3-e^3)$. Substituting this into $2v^3=x^2+y^2$ and cancelling $2$, we obtain that either 
\begin{equation}\label{eq:v3d64e6}
v^3=d^6+4e^6
\end{equation}
or $v^3=4d^6+e^6$. 

If $d=0$, then \eqref{eq:v3d64e6} implies that $v=e=0$, but then $x=y=0$, a contradiction with $\text{gcd}(x,y)=1$. Hence, $d\neq 0$, and, after dividing \eqref{eq:v3d64e6} by $d^6$ and using the change of variables $X=\frac{4v}{d^2}$ and $Y=16\frac{e^3}{d^3}$, the problem is reduced to describing all rational solutions to
$$
Y^2=X^3-64.
$$
This is an elliptic curve with rank $0$, hence all its rational points are torsion points \cite[Section 3.4.4]{mainbook}. Using SageMath we can verify that the only torsion point is $(X,Y)=(4,0)$. Therefore $Y=16\frac{e^3}{d^3}=0$, so $e=0$, which implies that $x=y$ and leads to primitive solutions $(x,y,z)=\pm (1, 1, 0)$. Similarly, equation $v^3=4d^6+e^6$ implies that $d=0$, and results in primitive solutions $(x,y,z)=\pm (1, -1, 0)$. In conclusion, the only primitive integer solutions to \eqref{eq:x4my4m2z3} are $(x,y,z)=(\pm 1, \pm 1,0)$.
\end{proof}

\begin{proposition}
The only primitive integer solutions to 
\begin{equation}\label{eq:x4my4m3z3}
	x^4-y^4=3z^3
\end{equation}
are $(x,y,z)=(\pm 1,\pm 1,0)$.
\end{proposition}

\begin{proof}
We first note that if $(x,y,z)$ is a primitive solution to \eqref{eq:x4my4m3z3}, then $x,y,z$ are pairwise coprime. Equation \eqref{eq:x4my4m3z3} can be rewritten as
$$
st=3z^3
$$
where $s=x^2-y^2$ and $t=x^2+y^2$ are integers. If $s$ and $t$ share a prime factor $p$ then it must divide $s+t=2x^2$ and $s-t=2y^2$, hence $p=2$ or $p$ divides both $x$ and $y$, which is a contradiction, hence $\gcd(s,t)$ is at most $2$. Because $x$ and $y$ cannot both be even, they must either both be odd or have different parities. If $x$ and $y$ have different parities, then variables $s$ and $t$ are both odd and therefore coprime. If $x$ and $y$ are both odd, then $s$ and $t$ are both even and the change of variables $S=s/2$ and $T=t/2$ ensures that $S$ and $T$ are coprime.  Hence we have two equations to solve:
\begin{equation}\label{reda:x4my4m3z3}
	st=3z^3 \quad \text{with} \quad s,t \,\,\, \text{both odd}, \quad \gcd(s,t)=1, 
\end{equation} 
and
\begin{equation}\label{redb:x4my4m3z3}
	4ST=3z^3 \quad \text{with} \quad \gcd(S,T)=1.
\end{equation} 
Let us first consider equation \eqref{reda:x4my4m3z3}. In this case, we have $(s,t,z)=(u^3,3v^3,uv)$ or $(3u^3,v^3,uv)$ for odd integers $u,v$. 

First, consider solution $(s,t,z)=(u^3,3v^3,uv)$. Then in the original variables, we obtain the system of equations $u^3=x^2-y^2$ and $3v^3=x^2+y^2$. We can factorise the first equation and write it as
\begin{equation}\label{eq:sqmu3}
	(x+y)(x-y)=u^3.	
\end{equation}	
Because $x$ and $y$ have different parity, $a=x+y$ and $b=x-y$ are odd and coprime, and the equation reduces to
$$
ab=u^3 \quad \text{with} \quad \gcd(a,b)=1
$$
whose integer solutions are $(a,b,u)=(d^3,e^3,de)$ for odd integers $d,e$. Then $x=\frac{d^3+e^3}{2}$ and $y=\frac{d^3-e^3}{2}$. Substituting this into $3v^3=x^2+y^2$ and multiplying by $2$, we obtain equation
\begin{equation}\label{eq:6v3d6e6}
	6v^3=d^6+e^6,
\end{equation}
for which it is left to determine all its integer solutions. If $e=0$, then \eqref{eq:6v3d6e6} implies that $d=v=0$, but then $x=y=0$, a contradiction with $\text{gcd}(x,y)=1$. Hence, $e\neq 0$, and after dividing equation \eqref{eq:6v3d6e6} by $e^6$ and making the change of variables $X=\frac{6v}{e^2}$ and $Y=\frac{6d^3}{e^3}$, the problem is reduced to finding all rational points on elliptic curve
$$
Y^2=X^3-36.
$$
This is a rank $0$ curve with no torsion points, and hence we obtain no primitive solutions $(x,y,z)$ in this case. 

Next, consider solution $(s,t,z)=(3u^3,v^3,uv)$. Then in the original variables, we obtain the system of equations $3u^3=x^2-y^2$ and $v^3=x^2+y^2$. We can factorise the first equation and write it as
\begin{equation}\label{eq:sqm3u3}
	(x+y)(x-y)=3u^3.	
\end{equation}	
Because $x$ and $y$ have different parities, $a=x+y$ and $b=x-y$ are odd and coprime, and the equation reduces to
$$
ab=3u^3 \quad \text{with} \quad \gcd(a,b)=1
$$
whose integer solutions are $(a,b,u)=(d^3,3e^3,de)$ or $(3d^3,e^3,de)$ for odd integers $d,e$. 
Then substituting $x=(a+b)/2$ and $y=(a-b)/2$ into $v^3=x^2+y^2$ and multiplying by $2$, we conclude that either
\begin{equation}\label{eq:2v3d69e6}
	2v^3=9d^6+e^6,
\end{equation}
or $2v^3 = d^6+9e^6$. Let us solve \eqref{eq:2v3d69e6}. It is clear that $d\neq 0$. After dividing the equation by $d^6$ and making the change of variables $X=\frac{2v}{d^2}$ and $Y=\frac{2e^3}{d^3}$, the problem is reduced to finding all rational points on elliptic curve
$$
Y^2=X^3-36.
$$
This is a rank $0$ curve with no torsion points, and hence we obtain no primitive solutions $(x,y,z)$ in this case. The case $2v^3 = d^6+9e^6$ is the same up to exchange of $d$ and $e$.

We must now consider the equation \eqref{redb:x4my4m3z3}. Note that because $x$ and $y$ are both odd, $S=\frac{x^2-y^2}{2}$ is even, while $T=\frac{x^2+y^2}{2}$ is odd. We can see that $z$ is even, and make the change of variable $z=2Z$, and after cancelling $4$, we obtain equation
\begin{equation}\label{redb2:x4my4m3z3}
	ST=6Z^3 \quad \text{with} \quad \gcd(S,T)=1,
\end{equation} 
whose solutions are $(S,T,Z)=(6u^3,v^3,uv)$ and $(2u^3,3v^3,uv)$ for integers $u,v$ (note that solution $(S,T,Z)=(3u^3,2v^3,uv)$ is impossible due to $T$ being odd). 

First, consider solution $(S,T,Z)=(6u^3,v^3,uv)$. Then in the original variables, we obtain the system of equations $12u^3=x^2-y^2$ and $2v^3=x^2+y^2$. We can factorise the first equation and write it as
\begin{equation}\label{eq:sqm12u3}
(x+y)(x-y)=12u^3.	
\end{equation}	
Because $x$ and $y$ are both odd, $a=(x+y)/2$ and $b=(x-y)/2$ are coprime and \eqref{eq:sqm12u3} reduces to
$$
ab=3u^3 \quad \text{with} \quad \gcd(a,b)=1
$$
whose integer solutions are $(a,b,u)=(3d^3,e^3,de)$ and $(d^3,3e^3,de)$ for integers $d,e$. Then, either $(x,y)=(3d^3+e^3,3d^3-e^3)$ or $(x,y)=(d^3+3e^3,d^3-3e^3)$. Substituting this into
$2v^3=x^2+y^2$ and cancelling $2$, we obtain equation
\begin{equation}\label{eq1:v3d69e6}
v^3=9d^6+e^6,
\end{equation}
or $v^3=d^6+9e^6$.
First consider \eqref{eq1:v3d69e6}. If $d=0$, then $x=e^3$ and $y=-e^3$ which are only coprime if $e=\pm 1$ and results in primitive solutions $(x,y,z)=\pm(1,-1,0)$. If $d \neq 0$, after dividing equation \eqref{eq1:v3d69e6} by $d^6$ and making the change of variables $X=\frac{v}{d^2}$ and $Y=\frac{e^3}{d^3}$, the problem is reduced to finding all rational points on elliptic curve
$$
Y^2=X^3-9.
$$
This is a rank $0$ curve with no torsion points, and hence we obtain no primitive solutions $(x,y,z)$ in this case. Similarly, equation $v^3=d^6+9e^6$ implies that $e=0$ and results in primitive solutions $(x,y,z)=\pm(1,1,0)$. 

Next, consider solution $(S,T,Z)=(2u^3,3v^3,uv)$. Then in the original variables, we obtain the system of equations $4u^3=x^2-y^2$ and $6v^3=x^2+y^2$. 
We can factorise the first equation and write it as
\begin{equation}\label{eq:sqm4u3}
	(x+y)(x-y)=4u^3.	
\end{equation}	
Because $x$ and $y$ are both odd, $a=(x+y)/2$ and $b=(x-y)/2$ are coprime and \eqref{eq:sqm4u3} reduces to
$$
ab=u^3 \quad \text{with} \quad  \gcd(a,b)=1
$$
whose integer solutions are $(a,b,u)=(d^3,e^3,de)$ for integers $d,e$. Then, $x=d^3+e^3$ and $y=d^3-e^3$. Substituting this into $6v^3=x^2+y^2$ and cancelling $2$, we obtain equation
\begin{equation}\label{eq:3v3d6e6}
	3v^3=d^6+e^6,
\end{equation}
for which it is left to determine all its integer solutions. If $d=0$, then \eqref{eq:3v3d6e6} implies that $v=e=0$, but then $x=y=0$, a contradiction. Hence, $d \neq 0$ and after dividing the equation by $d^6$ and making the change of variables $X=\frac{3v}{d^2}$ and $Y=\frac{3e^3}{d^3}$, the problem is reduced to finding all rational points on elliptic curve
$$
Y^2=X^3-9.
$$
This is a rank $0$ curve with no torsion points, and hence we obtain no primitive solutions $(x,y,z)$ in this case. 

Finally, we conclude that equation \eqref{eq:x4my4m3z3} has the only primitive integer solutions $(x,y,z)=(\pm 1,\pm 1,0)$. 
 \end{proof}
 
 \begin{proposition}
 The only primitive integer solutions to 
  \begin{equation}\label{eq:x4p2y3p2z3}
 	x^4+2y^3+2z^3=0.
 \end{equation}
 are $(x,y,z)=\pm (0,1,-1)$.
 \end{proposition}
 
 \begin{proof}
 First note that if $(x,y,z)$ is a primitive solution to \eqref{eq:x4p2y3p2z3} then $x,y$ and $z$ must be pairwise coprime.
 We can easily see that $x$ must be even, hence $y$ and $z$ are both odd. After the change of variables $x=2t$, $y\to -y$, and cancelling $2$, we obtain
 \begin{equation}\label{eq:8t4my3pz3}
 8t^4=y^3-z^3=(y-z)(y^2+yz+z^2).
 \end{equation}
 If prime $p$ is a common factor of $y-z$ and $y^2+yz+z^2$, then $p$ also divides $(y^2+yz+z^2)-(y-z)^2=3yz$.  
 Because variables $y$ and $z$ are coprime, the only option is $p=3$, and we conclude that integers $y-z$ and $y^2+yz+z^2$ are either coprime or have greatest common divisor equal to $3$. If they are coprime, then, because $y^2+yz+z^2$ is odd and positive, we must have
 \begin{equation}\label{sys1:x4p2y3p2z3}
 	y-z=8u^4 \quad \text{and} \quad y^2+yz+z^2=v^4
 \end{equation}
 where $u,v$ are coprime integers.
 
 Substituting the equation on the left-hand side of \eqref{sys1:x4p2y3p2z3} into the one on the right-hand side, we obtain
 $$
 3z^2+24u^4z+64u^8=v^4.
 $$ 
 If $u=0$ then $y=z$ and \eqref{sys1:x4p2y3p2z3} implies that $y=v=0$, a contradiction with $\gcd(y,z)=1$. Hence $u \neq 0$ and after dividing by $u^8$, and the change of variables $X=\frac{z}{u^4}$ and $Y=\frac{v}{u^2}$, the equation reduces to
 $$
 3X^2+24X+64=Y^4.
 $$
 This is a genus $1$ curve and assuming $X \neq -4$, this can be transformed to the Weierstrass equation $T^2= S^3+ 576S$ after change of variables 
  \begin{equation}\label{eq:ST}
  \begin{multlined}
 S=\frac{8(4Y^3 + 3X^2 + 8Y^2 + 24X + 16Y + 80)}{(X+4)^2},  \\ 
 T= -32 \left( \frac{3X^2Y^2 + 12X^2Y + 24XY^2 + 32Y^3 + 36X^2}{(X+4)^3} \right.  \\  \left.   + \frac{96XY + 112Y^2 + 288X + 320Y + 832}{(X+4)^3}  \right) 
 \end{multlined}
 \end{equation}
 The Weierstrass equation has rank $0$ and all its rational points are its torsion points and its only one is $(S,T)=(0,0)$. In this case, \eqref{eq:ST} implies that $(X,Y)=(-4,-2)$. However, if $X=-4$, then $z=-4u^4$ and because $y-z=8u^4$, we have $y=4u^4$. Because $y$ and $z$ are coprime, there is no suitable value for $u$, hence we obtain no non-trivial primitive solutions in this case.

 Now consider case $\text{gcd}(y-z,y^2+yz+z^2)=3$. Then $x$ is divisible by $3$, and $\text{gcd}(x,y,z)=1$ implies that $y$ and $z$ are not both divisible by $3$. Then modulo $27$ analysis of \eqref{eq:x4p2y3p2z3} implies that $y-z$ is divisible by $9$, while $y^2+yz+z^2$ is divisible by $3$ but not by $9$. Together with the fact that $y^2+yz+z^2$ is odd and positive, this implies \eqref{eq:8t4my3pz3} is possible only if
 \begin{equation}\label{sys2:x4p2y3p2z3}
 	y-z=8\cdot 27 \cdot u^4 \quad \text{and} \quad y^2+yz+z^2=3v^4
 \end{equation}
for some coprime integers $u$ and $v$.
Substituting the equation on the left-hand side of \eqref{sys2:x4p2y3p2z3} into the one on the right-hand side, we obtain
$$
3z^2+648u^4z+46656u^8=3v^4.
$$ 
If $u=0$ then $z=\pm v^2$ and \eqref{sys2:x4p2y3p2z3} implies that $y=z=\pm 1$, which leads to the primitive solutions to \eqref{eq:x4p2y3p2z3} $(x,y,z)=\pm(0,1,-1)$. If $u \neq 0$, after dividing by $3u^8$, and the change of variables $X=\frac{v}{u^2}$ and $Y=\frac{z}{u^4}+108$, we obtain
$$
Y^2 = X^4 - 3888.
$$
This equation in homogeneous form is 
$$
\tilde{Y}^2Z^2=\tilde{X}^4-3888Z^4
$$
where $\tilde{X}=XZ$, $\tilde{Y}=YZ$ and $Z\neq 0$ are integers. 
It is clear that $Y \neq 0$ and then we can reduce this to the rank $0$ elliptic curve
\begin{equation}\label{eq:T2mS3m1259712S}
T^2=S^3+1259712S
\end{equation}
where $S=\frac{139968\tilde{X}^2}{\tilde{Y}^2}$ and $T=\frac{419904(\tilde{X}^5+3888\tilde{X}Z^4)}{\tilde{Y}^3Z^2}$. The only rational points of \eqref{eq:T2mS3m1259712S} are its torsion points, and its only one is $(S,T)=(0,0)$. This implies that $X=0$, hence $v=0$, and \eqref{sys2:x4p2y3p2z3} implies that $(y,z)=(0,0)$, a contradiction with $\gcd(y,z)=1$. 
Hence we obtain no primitive solutions in this case.

Therefore, we can conclude that equation \eqref{eq:x4p2y3p2z3} has the only primitive integer solutions $(x,y,z)=\pm (0,1,-1)$. 

\end{proof}

\subsection{Open question 8.12} \label{open:8.12}

\subsection{Arbitrarily large solutions}
Some equations have trivial infinite families of integer solutions, with, for example, one of the variables equal to $0$. For such equations, a more interesting question to ask is whether it has solutions with the absolute values of all variables being arbitrarily large. 
We now formulate this problem as in \cite[Problem 5.54]{mainbook}.
Given a Diophantine equation $P(x_1,\dots,x_n)=0$, determine whether for any $k \geq 0$ it has a solution such that 
\begin{equation}\label{cond:arblarge}
\min(|x_1|,\dots,|x_n|)\geq k.
\end{equation}
If so, the problem is solved, if not, then describe all its integer solutions.

In this section we solve this problem for two equations left open in \cite{mainbook}.

The equations we will solve are 
\begin{equation}\label{eq:x2y2pxz2pypc}
x^2y^2+xz^2+y + c=0, \quad c \in \{-1,1\}.
\end{equation}
These equations have obvious infinite families of integer solutions, for example, we can take $x=0$, $y=-c$, and $z$ arbitrary.  However, finding arbitrarily large integer solutions is much more difficult. 

\begin{theorem}\label{thm:arblarge}
For any integer $k \geq 0$, the equations \eqref{eq:x2y2pxz2pypc} have an integer solution such that $\min(|x|,|y|,|z|)\geq k$.  
\end{theorem}

Let us define $a$ as $a=r^2+1$, where $r>0$ is a solution of 
\begin{equation}\label{eq:Pell2}
s^2-2r^2=1.
\end{equation}
Infinitely many integers $(s,r)$ satisfying \eqref{eq:Pell2} can be obtained by applying the recursive formula 
$$
(s_0,r_0)= (1,0), \,\, \text{and} \,\, (s_{n+1},r_{n+1})= (3s_n + 4r_n,2s_n + 3r_n),
$$
see \cite[Section 3.1.3]{mainbook}, and these solutions satisfy $s_n \to \infty$ and $r_n \to \infty$ as $n \to \infty$.  Modulo $4$ analysis of equation \eqref{eq:Pell2} returns that $s$ is odd and $r$ is even, hence $a$ is odd. So, by definition and by \eqref{eq:Pell2} we have that $a=\frac{s^2+1}{2}$. Also, by definition, $a$ is not a perfect square, and $a \to \infty$ as $r \to \infty$.

\begin{lemma}\label{lm:cond}
Let $a$ be defined as above. Then the equations
\begin{equation}\label{eq:redpell}
		u^2-av^2=1- 4a^2c, \quad c\in \{-1,1\}
\end{equation} 
have infinitely many integer solutions such that $v \equiv 0 \pmod{2a}$ and $u \equiv 1 \pmod{2a^2}$.  
\end{lemma}

\begin{proof}

	If $c=-1$ then it is easy to verify that for all $a$, $(u,v)=(u_0,v_0):=(2a+1,2)$ is an integer solution to \eqref{eq:redpell}, and $\gcd(u_0,a)=\gcd(2a+1,a)=1$. While if $c=1$ then it is easy to verify that for all $a$, $(u,v)=(u_0,v_0):=(1-2a,2s)$ is an integer solution to \eqref{eq:redpell} and $\gcd(u_0,a)=\gcd(2a-1,a)=1$. 
	
Recall that $s,r$ are integers satisfying \eqref{eq:Pell2}. 
Integer solutions to \eqref{eq:redpell} can be generated by the recurrence relations $(u,v)=(u_n,v_n)$ where $(u_0,v_0)$ is defined above and $(u_{n+1},v_{n+1})$ is given by\footnote{Given a solution of $(a,b)\neq (1,0)$ to $a^2-db^2=1$, $n \neq 0$ and $d>0$ not a perfect square, an infinite family of integer solutions to $x^2-dy^2=m$ with $m \neq 0$ can be represented using the recursive formula, $(x,y)=\pm (x_n,y_n)$ where $(x_n,y_n)$, $n \geq 0$ is given by  
$
(x_0,y_0)=(x^*,y^*) \text{ and } (x_{n+1},y_{n+1})=(ax_n+bdy_n,ay_n+bx_n)$, where $(x^*,y^*)$ is an integer solution to $x^2-dy^2=m$, see Section \ref{ex:genPell}.}
$$
(u_{n+1},v_{n+1})=(s^2 u_n+2ra v_n, 2ru_n+s^2v_n), \quad n=0,1,2,\dots.
$$

Because $u_0=1 \pm 2a$, we have $u_0 \equiv 1 \pmod{2a}$. Then, considering the recurrence relations modulo $2a$, we obtain that $u_{n+1}=-u_n \pmod{2a}$, hence $|u_n| \equiv 1  \pmod{2a}$ for all $n$.

Next, $v_0$ is even. Then the recurrence relation preserves that $v_n$ is even for all $n$. Also, because $u_n \equiv (-1)^n \pmod{2a}$ we have that $2r u_{2n+1} \equiv -2r \pmod{2a}$ and $2r u_{2n} \equiv 2r \pmod{2a}$, $s^2v_n \equiv - v_n \pmod{2a}$ hence, $v_{2n+1}\equiv 2r-v_{2n} \pmod{2a}$ and $v_{2n}\equiv -2r-v_{2n-1} \pmod{2a}$. Now, $v_{2n+2}\equiv v_{2n}-4r \pmod{2a}$, and the even subsequence satisfies $v_{2n} \equiv v_{0}-4nr \pmod{2a}$. Because $a$ is odd and fixed, $r$ is even and fixed and $\gcd(4r,2a)=2$, the sequence visits every even residue modulo $2a$. Therefore $v_{2n} \equiv 0 \pmod{2a}$ for infinitely many $n$. Hence, we can generate infinitely many pairs $(u,v)$ satisfying \eqref{eq:redpell} such that $(|u|,v) \equiv (1, 0) \pmod{2a}$. By changing $u$ to $-u$ if necessary, we can in fact ensure that $(u,v) \equiv (1, 0) \pmod{2a}$ infinitely often.

Because $u \equiv 1 \pmod{2a}$ we can write $u=1+2ka$ where $k$ is an integer, and because $v \equiv 0 \pmod{2a}$, we can write $v=2am$ for some integer $m$. Substituting these into \eqref{eq:redpell} and dividing by $4a$, we obtain $a(-c-k^2+am^2)=k$, so we can write $k=at$ where $t$ is an integer. Substituting these into $u=1+2ka$ we obtain that $u=1+2a^2t$ hence $u \equiv 1 \pmod{2a^2}$.
\end{proof}

We can now prove Theorem \ref{thm:arblarge}. 
\begin{proof}[Proof of Theorem \ref{thm:arblarge}]

For any $k>0$, let $(s,r)$ be an integer solution to \eqref{eq:Pell2} such that $a:=r^2+1>k$. By Lemma \ref{lm:cond} there are infinitely many pairs $(u,v)$ satisfying \eqref{eq:redpell} such that $u \equiv 1 \pmod{2a^2}$ and $v \equiv 0 \pmod{2a}$. For each such $(u,v)$, $(x,y,z)=\left(-a,\frac{u-1}{2a^2}, \frac{v}{2a}\right)$ is an integer solution to \eqref{eq:x2y2pxz2pypc}, which can be verified by direct substitution, and we can clearly select $(u,v)$ such that $|y|>k$ and $|z|>k$. 
\end{proof}

The final equation we will consider in this section is $x^3y^2+z^2+y+1=0$. 

\begin{theorem}
There is a sequence \((x_n,y_n,z_n)\in\mathbb{Z}^3\) such that
\[
 x_n^3y_n^2+z_n^2+y_n+1=0,
 \qquad
 \min\{|x_n|,|y_n|,|z_n|\}\longrightarrow\infty.
\]
\end{theorem}

\begin{proof}
For \(n\geq0\), define integers \(P_n,Q_n\) by
\[
 P_n+Q_n\sqrt6=(4321+1764\sqrt6)(5+2\sqrt6)^{6n}.
\]
Since \(4321^2-6\cdot1764^2=865\) and \(5^2-6\cdot2^2=1\),
we have \(P_n^2-6Q_n^2=865\).
The coefficient-wise congruence
\[
 (5+2\sqrt6)^6\equiv1+36\sqrt6\pmod{216}
\]
shows inductively that \(P_n\equiv1\pmod{216}\) and
\(Q_n\equiv0\pmod{36}\). Consequently
\[
 r_n=\frac{P_n-1}{12},\qquad t_n=\frac{Q_n}{2}
\]
are positive multiples of \(18\), starting with \((r_0,t_0)=(360,882)\).
The norm equation gives
\begin{equation}\label{eq:parameters}
 t_n^2=6r_n^2+r_n-36.
\end{equation}
Both parameters tend to infinity. Thus \(q_n=r_n^2-6\) tends to infinity;
also \((r_n-1)^2<q_n<r_n^2\), so \(q_n\) is nonsquare.

Fix one pair \(r=r_n,\ t=t_n\), and put \(q=r^2-6\).
Then \eqref{eq:parameters} says \(t^2=6q+r\). In \(\mathbb{Z}[\sqrt q]\), set
\[
 A=\frac{r^2}{3}-1,\qquad B=\frac r3,\qquad
 \eta=A+B\sqrt q,\qquad
 \alpha=-(2qr+1)+2t\sqrt q.
\]
Here \(A,B\) are positive integers, \(B\) is even, and \(\eta>1\).
For the norm \(N(a+b\sqrt q)=a^2-qb^2\), direct expansion gives
\[
 N(\eta)=1,\qquad
 N(\alpha)=(2qr+1)^2-4qt^2=1+4q^3.
\]
For every integer \(h\geq0\), define
\begin{equation}\label{eq:inner}
 K_h=rt+2qh,\qquad
 U_h+W_h\sqrt q=\alpha\eta^{K_h}.
\end{equation}
These coefficients are integers and satisfy
\begin{equation}\label{eq:norm}
 U_h^2-qW_h^2=1+4q^3.
\end{equation}

All congruences in this step are coefficient-wise modulo \(2q\).
Since \(A^2-qB^2=1\) and \(B\) is even,
\[
 A^2\equiv1,\qquad (B\sqrt q)^2\equiv0.
\]
The integer \(K_h\) is even. The binomial theorem therefore yields
\[
 \eta^{K_h}\equiv1+K_hAB\sqrt q.
\]
Moreover,
\[
 K_hAB\equiv rtAB\equiv2t\pmod{2q},
 \qquad
 rtAB-2t=\frac{qt(r^2+3)}9\in2q\mathbb{Z},
\]
where the last divisibility uses \(18\mid t\). It follows that
\[
 \alpha\eta^{K_h}
 \equiv(-1+2t\sqrt q)(1+2t\sqrt q)
 \equiv-1\pmod{2q}.
\]
Thus \(U_h\equiv-1\pmod{2q}\) and \(W_h\equiv0\pmod{2q}\).
Set
\[
 u_h=\frac{U_h+1}{2},\qquad z_h=\frac{W_h}{2q};
\]
these are integers and \(q\mid u_h\).
Equation \eqref{eq:norm} becomes
\[
 u_h(u_h-1)=q^3(z_h^2+1).
\]
Since \(\gcd(u_h-1,q)=1\), it follows that \(q^3\mid u_h\).
Hence
\begin{equation}\label{eq:solution}
 x=-q,\qquad y_h=\frac{U_h+1}{2q^3},\qquad
 z_h=\frac{W_h}{2q}
\end{equation}
are integers. Substituting \(u_h=q^3y_h\) into the preceding identity gives
\(q^3y_h^2-y_h=z_h^2+1\), or
\[
 (-q)^3y_h^2+z_h^2+y_h+1=0.
\]

For fixed \(r,t,q\), write
\(\alpha_\pm=-(2qr+1)\pm2t\sqrt q\).
Their product is \(1+4q^3\), so \(\alpha_+\neq0\).
As \(N(\eta)=1\), its conjugate is \(\eta^{-1}\), and
\[
 U_h=\frac{\alpha_+\eta^{K_h}+\alpha_-\eta^{-K_h}}2,
 \qquad
 W_h=\frac{\alpha_+\eta^{K_h}-\alpha_-\eta^{-K_h}}{2\sqrt q}.
\]
Since \(\eta>1\) and \(K_h\to\infty\), both \(|U_h|\) and \(|W_h|\)
tend to infinity. Therefore, with \(q\) fixed,
\(|y_h|,|z_h|\to\infty\).

For each \(n\geq1\), use the outer pair \((r_n,t_n)\) and choose \(h_n\)
large enough that the corresponding solution satisfies
\(|y_n|,|z_n|\geq n\). Then \(x_n=-q_n\), with \(q_n\to\infty\).
\end{proof}

\bibliography{references}

@misc{Z2018,
   title  = {What is the smallest unsolved {D}iophantine equation?},
   author = {Zidane},
   howpublished = "\url{https://mathoverflow.net/questions/316708/}",
   year   = {2018 (accessed June 12, 2020)}
}

@article{baker1970integer,
  title={Integer points on curves of genus 1},
  author={Baker, Alan and Coates, John},
  journal={Mathematical Proceedings of the Cambridge Philosophical Society},
  volume={67},
  number={3},
  pages={595--602},
  year={1970},
  organization={Cambridge University Press}
}

@article{vaserstein2010polynomial,
  title={Polynomial parametrization for the solutions of {D}iophantine equations and arithmetic groups},
  author={Vaserstein, Leonid},
  journal={Annals of mathematics},
  pages={979--1009},
  year={2010},
  publisher={JSTOR}
}

@article{siegel1929uber,
  title={{\"U}ber einige anwendungen diophantischer approximationen},
  author={Siegel, Carl L},
  journal={Sitzungsberichte der Preussischen Akademie der Wissenschaften},
  year={1929},
  volume={1},
  pages={14--72}
}

@article{spearman1992cubic,
  title={The cubic congruence $x^3+ Ax^2+ Bx+ C\equiv 0\, (\text{mod}\, p)$ and binary quadratic forms},
  author={Spearman, Blair K and Williams, Kenneth S},
  journal={Journal of the London Mathematical Society},
  volume={2},
  number={3},
  pages={397--410},
  year={1992},
  publisher={Oxford Academic}
}

@book{dujella2021number,
  title={Number {T}heory},
  author={Dujella, Andrej},
  year={2021},
  publisher={{\v{S}}kolska knjiga Zagreb}
}

@article{baker1968diophantine,
  title={The {D}iophantine equation $y^2= ax^3+ bx^2+ cx+ d$},
  author={Baker, Alan},
  journal={Journal of the London Mathematical Society},
  volume={1},
  number={1},
  pages={1--9},
  year={1968},
  publisher={Oxford Academic}
}

@article {MR1994473,
	AUTHOR = {Sun, Zhi-Hong},
	TITLE = {Cubic and quartic congruences modulo a prime},
	JOURNAL = {J. Number Theory},
	FJOURNAL = {Journal of Number Theory},
	VOLUME = {102},
	YEAR = {2003},
	NUMBER = {1},
	PAGES = {41--89},
	ISSN = {0022-314X},
	MRCLASS = {11A07 (11A15 11B39 11B50)},
	MRNUMBER = {1994473},
	MRREVIEWER = {Charles Helou},
	DOI = {10.1016/S0022-314X(03)00067-2},
	URL = {https://doi.org/10.1016/S0022-314X(03)00067-2},
}

@book{mainbook,
	title={Polynomial {D}iophantine {E}quations. {A} {S}ystematic {A}pproach},
	author={Grechuk, Bogdan},
	year={2024},
	publisher={Springer}
}

@article{wilcox2024generalisedfermatequationsurvey,
title={Generalised Fermat equation: a survey of solved cases}, 
      author={Ashleigh Wilcox and Bogdan Grechuk},
      year={2024},
	journal={arXiv preprint arXiv:2412.11933},
}
\bibliographystyle{plain}

\end{document}